\documentclass{article}
\usepackage{amsmath}
\usepackage{amssymb}
\usepackage{amsthm}
\usepackage{cite}
\usepackage[arrow,curve,matrix,arc]{xy}
\begin{document}
\def\e#1\e{\begin{equation}#1\end{equation}}
\def\ea#1\ea{\begin{align}#1\end{align}}
\def\eq#1{{\rm(\ref{#1})}}
\theoremstyle{plain}
\newtheorem{thm}{Theorem}[section]
\newtheorem{qthm}[thm]{``Theorem''}
\newtheorem{lem}[thm]{Lemma}
\newtheorem{prop}[thm]{Proposition}
\newtheorem{cor}[thm]{Corollary}
\newtheorem{cla}[thm]{Claim}
\newtheorem{conj}[thm]{Conjecture}
\theoremstyle{definition}
\newtheorem{dfn}[thm]{Definition}
\newtheorem{ex}[thm]{Example}
\newtheorem{rem}[thm]{Remark}
\newtheorem{alg}[thm]{Algorithm}
\newtheorem{hypo}[thm]{Hypothesis}
\newtheorem{cond}[thm]{Condition}
\newtheorem{conv}[thm]{Convention}
\newtheorem{princ}[thm]{Principle}
\newtheorem{property}[thm]{Property}
\numberwithin{figure}{section}
\def\dim{\mathop{\rm dim}\nolimits}
\def\vdim{\mathop{\rm vdim}\nolimits}
\def\Ker{\mathop{\rm Ker}}
\def\Coker{\mathop{\rm Coker}}
\def\Diff{\mathop{\rm Diff}}
\def\Im{\mathop{\rm Im}}
\def\Stab{\mathop{\rm Stab}\nolimits}
\def\rank{\mathop{\rm rank}}
\def\Hom{\mathop{\rm Hom}\nolimits}
\def\Aut{\mathop{\rm Aut}}
\def\Fix{\mathop{\rm Fix}}
\def\sign{\mathop{\rm sign}}
\def\supp{\mathop{\rm supp}}
\def\Pd{\mathop{\rm Pd}\nolimits}
\def\id{\mathop{\rm id}\nolimits}
\def\Iso{{\rm Iso}}
\def\na{{\rm na}}
\def\ef{{\rm ef}}
\def\ec{{\rm ec}}
\def\ecb{{\rm ecb}}
\def\cs{{\rm cs}}
\def\rcs{{\rm rcs}}
\def\eo{{\rm eo}}
\def\eb{{\rm eb}}
\def\tr{{\rm tr}}
\def\trc{{\rm trc}}
\def\eac{{\rm eac}}
\def\eca{{\rm eca}}
\def\Kh{{\rm Kh}}
\def\Kch{{\rm Kch}}
\def\rKh{{\rm rKh}}
\def\rKch{{\rm rKch}}
\def\rsi{{\rm si}}
\def\rrsi{{\rm rsi}}
\def\orb{{\rm orb}}
\def\nov{{\rm nov}}
\def\ev{\mathop{\rm ev}\nolimits}
\def\evb{\mathop{\rm evb}\nolimits}
\def\bev{\mathop{\bf ev}\nolimits}
\def\bevb{\mathop{\bf evb}\nolimits}
\def\ac{{\rm ac}}
\def\ca{{\rm ca}}
\def\Kb{{\rm Kb}}
\def\Kcb{{\rm Kcb}}
\def\bo{{\rm bo}}
\def\tr{{\rm tr}}
\def\ma{{\rm ma}}
\def\bs{\boldsymbol}
\def\ge{\geqslant}
\def\le{\leqslant\nobreak}
\def\R{{\mathbin{\mathbb R}}}
\def\Z{{\mathbin{\mathbb Z}}}
\def\Q{{\mathbin{\mathbb Q}}}
\def\N{{\mathbin{\mathbb N}}}
\def\C{{\mathbin{\mathbb C}}}
\def\CP{{\mathbin{\mathbb{CP}}}}
\def\RP{{\mathbin{\mathbb{RP}}}}
\def\M{{\mathbin{\mathcal M}}}
\def\oM{{\mathbin{\smash{\,\,\overline{\!\!\mathcal M\!}\,}}}}
\def\uP{\underline{P\!}\,}
\def\umu{\underline{\smash{\mu}}}
\def\uG{{\underline{G\!}\,}}
\def\uH{{\underline{H\!}\,}}
\def\utG{{\underline{\ti G\!}\,}{}}
\def\uaG{{\underline{\acute G\!}\,}{}}
\def\utH{{\underline{\ti H\!}\,}{}}
\def\ucG{{\underline{\check G\!}\,}{}}
\def\uC{{\underline{C\!}\,}}
\def\utC{{\underline{\ti C\!}\,}{}}
\def\ucC{{\underline{\check C\!}\,}{}}
\def\uD{{\underline{D\!}\,}}
\def\ubG{{\underline{\boldsymbol G\!}\,}}
\def\ubC{{\underline{\boldsymbol C\!}\,}}
\def\ubtG{{\underline{\boldsymbol{\ti G}\!}\,}{}}
\def\ubaG{{\underline{\boldsymbol{\acute G}\!}\,}{}}
\def\ubcG{{\underline{\boldsymbol{\check G}\!}\,}{}}
\def\ubtC{{\underline{\boldsymbol{\ti C}\!}\,}{}}
\def\ubcC{{\underline{\boldsymbol{\check C}\!}\,}{}}
\def\ubH{{\underline{\boldsymbol H\!}\,}}
\def\ubD{{\underline{\boldsymbol D\!}\,}}
\def\ubtH{{\underline{\boldsymbol{\ti H}\!}\,}{}}
\def\ubtD{{\underline{\boldsymbol{\ti D}\!}\,}{}}
\def\m{{\mathfrak m}}
\def\U{\mathbin{\rm U}}
\def\ul{\underline}
\def\al{\alpha}
\def\be{\beta}
\def\ga{\gamma}
\def\de{\delta}
\def\io{\iota}
\def\ep{\epsilon}
\def\la{\lambda}
\def\ka{\kappa}
\def\th{\theta}
\def\ze{\zeta}
\def\up{\upsilon}
\def\vp{\varphi}
\def\si{\sigma}
\def\om{\omega}
\def\De{\Delta}
\def\La{\Lambda}
\def\Om{\Omega}
\def\Up{\Upsilon}
\def\Ga{\Gamma}
\def\Si{\Sigma}
\def\Th{\Theta}
\def\pd{\partial}
\def\ts{\textstyle}
\def\sst{\scriptscriptstyle}
\def\w{\wedge}
\def\sm{\setminus}
\def\bu{\bullet}
\def\op{\oplus}
\def\ot{\otimes}
\def\ov{\overline}
\def\bigop{\bigoplus}
\def\bigot{\bigotimes}
\def\iy{\infty}
\def\es{\emptyset}
\def\ra{\rightarrow}
\def\Ra{\Rightarrow}
\def\ab{\allowbreak}
\def\longra{\longrightarrow}
\def\hookra{\hookrightarrow}
\def\t{\times}
\def\ci{\circ}
\def\ti{\tilde}
\def\gr{\grave}
\def\d{{\rm d}}
\def\ha{{\ts\frac{1}{2}}}
\def\md#1{\vert #1 \vert}
\def\bmd#1{\big\vert #1 \big\vert}
\def\ms#1{\vert #1 \vert^2}
\title{Kuranishi homology and Kuranishi cohomology}
\author{Dominic Joyce}
\date{}
\maketitle
\begin{abstract}
A {\it Kuranishi space\/} is a topological space equipped with a
{\it Kuranishi structure}, defined by Fukaya and Ono. Kuranishi
structures occur naturally on many moduli spaces in differential
geometry. If $(M,\om)$ is a symplectic manifold with compatible
almost complex structure $J$, then moduli spaces of stable
$J$-holomorphic curves in $M$ have Kuranishi structures. Kuranishi
spaces are important in symplectic geometry, for defining
Gromov--Witten invariants, Lagrangian Floer cohomology, and so on.

Let $Y$ be an orbifold, and $R$ a commutative ring or $\Q$-algebra.
We shall define two new homology and cohomology theories of $Y$:
{\it Kuranishi (co)homology} $KH_*,KH^*(Y;R)$, for $R$ a
$\Q$-algebra, and {\it effective Kuranishi (co)homology}
$KH_*^\ef,KH^*_\ec(Y;R)$, for $R$ a commutative ring, such as
$R=\Z$. Our main result is that $KH_*,KH_*^\ef(Y;R)$ are isomorphic
to singular homology, and $KH^*,KH^*_\ec(Y;R)$ are isomorphic to
compactly-supported cohomology. We also define five different kinds
of {\it Kuranishi (co)bordism\/} $KB_*,KB^*(Y;R)$. These are new
topological invariants of $Y$ which turn out to be very large, far
larger than the (co)homology groups.

The (co)bordism theories $KB_*,KB^*(Y;R)$ are spanned over $R$ by
isomorphism classes $[X,\bs f]$, where $X$ is a compact oriented
Kuranishi space without boundary, and $\bs f:X\ra Y$ is a strongly
smooth map. The (co)homology theories are defined using (co)chain
complexes $KC_*,\ab\ldots,\ab KC^*_\ec(Y;R)$ spanned by isomorphism
classes $[X,\bs f,\bs G]$, where $X$ is a compact oriented Kuranishi
space with boundary and corners, $\bs f:X\ra Y$ is strongly smooth,
and $\bs G$ is some extra {\it gauge-fixing data\/} upon $X$. The
main purpose of $\bs G$ is to ensure the automorphism groups
$\Aut(X,\bs f,\bs G)$ are finite, which is necessary to get a
well-behaved (co)homology theory.

These theories are powerful new tools in symplectic geometry. We
define new Gromov--Witten type invariants in Kuranishi bordism,
containing more information than conventional Gromov--Witten
invariants. In a sequel we hope to use these to prove the
integrality conjecture for Gopakumar--Vafa invariants. We also
explain how to reformulate Lagrangian Floer cohomology using
Kuranishi cohomology rather than singular homology; this will lead
to significant technical simplifications, and improved results.

The results of this book are briefly surveyed in \cite{Joyc1}.
\end{abstract}

\newpage

\setcounter{tocdepth}{2}
\tableofcontents

\section{Introduction}
\label{kh1}

The goal of this book is to develop new tools for use in areas of
symplectic geometry concerned with `counting' $J$-holomorphic curves
--- Gromov--Witten theory, Lagrangian Floer cohomology, contact
homology, Symplectic Field Theory and so on --- and elsewhere, such
as in the `string topology' of Chas and Sullivan \cite{ChSu}. Moduli
spaces $\oM$ of (stable) $J$-holomorphic curves $\Si$ in a
symplectic manifold $(M,\om)$ are generally not smooth manifolds,
but {\it Kuranishi spaces}, topological spaces equipped with a {\it
Kuranishi structure}, as in Fukaya and Ono \cite{FuOn1}. To `count'
the number of curves in $\oM$ one generally perturbs $\oM$ to
construct a singular homology class called a {\it virtual moduli
cycle\/} or {\it chain}. Several versions of this virtual moduli
cycle construction
exist~\cite{FuOn1,FOOO,HWZ3,LiTi2,LiuTi,LuTi,McDu,Ruan,Sieb1}.

A great deal of the technical complexity of these areas is due to
the uneasy interplay between Kuranishi spaces, the raw geometric
data, and virtual moduli cycles and chains in singular homology, the
abstract perturbations. This is particularly true in areas like
Lagrangian Floer cohomology in which the moduli spaces $\oM$ have
boundaries and corners, and one must consider compatibility of
virtual moduli chains at the boundaries. Anyone who has seriously
tried to read the mammoth Fukaya et al.\ \cite{FOOO} will appreciate
this.

Let $Y$ be an orbifold, and $R$ a commutative ring or $\Q$-algebra.
The basic idea of this book is to define a new homology theory, the
{\it Kuranishi homology\/} $KH_*(Y;R)$ of $Y$ with coefficients in
$R$. There are two variations, {\it Kuranishi homology} $KH_*(Y;R)$
for $R$ a $\Q$-algebra, and {\it effective Kuranishi homology}
$KH_*^\ef(Y;R)$ for $R$ a commutative ring, such as $R=\Z$.
Kuranishi homology $KH_*(Y;R)$ is better behaved at the chain level.
Both are isomorphic to singular homology $H_*^\rsi(Y;R)$, so they
are not new topological invariants; the point is to use them as a
replacement for singular homology, in situations where one wants to
construct virtual chains for moduli spaces.

Kuranishi homology is the homology of a chain complex
$\bigl(KC_*(Y;R),\pd\bigr)$ of {\it Kuranishi chains}. Elements of
$KC_k(Y;R)$ are finite sums $\sum_{a\in A}\rho_a[X_a,\bs f_a,\bs
G_a]$ for $\rho_a\in R$, where $[X_a,\bs f_a,\bs G_a]$ is the
isomorphism class of a triple $(X_a,\bs f_a,\bs G_a)$, with $X_a$ a
compact oriented Kuranishi space with boundary and corners, of
virtual dimension $\vdim X_a=k$, $\bs f_a:X_a\ra Y$ a strongly
smooth map, and $\bs G_a$ some {\it gauge-fixing data} for $(X_a,\bs
f_a)$. The main point of $\bs G_a$ is to ensure that $\Aut(X_a,\bs
f_a,\bs G_a)$ is finite, since it turns out that if we allow
infinite automorphism groups then the resulting homology groups are
always zero.

We also define Poincar\'e dual theories of {\it (effective)
Kuranishi cohomology\/} $KH^*(Y;R)$, $KH^*_\ec(Y;R)$, which are
isomorphic to compactly-supported cohomology $H^*_\cs(Y;R)$, using a
cochain complex $\bigl(KC^*(Y;R),\d\bigr)$ of {\it Kuranishi
cochains}. Elements of $KC^k(Y;R)$ are finite sums $\sum_{a\in
A}\rho_a[X_a,\bs f_a,\bs C_a]$, where $X_a$ is a compact Kuranishi
space with $\vdim X_a=\dim Y-k$, $\bs f_a:X\ra Y$ is a {\it strong
submersion} with a relative orientation ({\it coorientation}), and
$\bs C_a$ is {\it co-gauge-fixing data} for~$(X_a,\bs f_a)$.

Here is an important difference between manifolds and Kuranishi
spaces. If $f:X\ra Y$ is a submersion of manifolds, then $\dim
X\ge\dim Y$. Thus, if we defined $KC^k(Y;R)$ using manifolds $X_a$
we would have $KC^k(Y;R)=0$ for $k>0$, which would be no use at all.
However, if $X$ is a Kuranishi space and $\bs f:X\ra Y$ is a strong
submersion, we need not have $\vdim X\ge\dim Y$. In fact, if $X$ is
a manifold or Kuranishi space and $\bs f:X\ra Y$ is (strongly)
smooth, we can always change the Kuranishi structure of $X$ to get a
strong submersion $\bs f{}^Y:X{}^Y\ra Y$. So strong submersions are
easy to produce.

If $X,\ti X$ are Kuranishi spaces and $\bs f:X\ra Y$, $\bs{\ti
f}:\ti X\ra Y$ are strong submersions, there is a {\it fibre
product\/} Kuranishi space $X\t_Y\ti X$, with a strong submersion
$\bs\pi_Y:X\t_Y\ti X\ra Y$. Using this we define a {\it cup
product\/} on $KC^*(Y;R)$ by
\e
[X,\bs f,\bs C]\cup[\ti X,\bs{\ti f},\bs{\ti C}]=
\smash{\bigl[X\t_Y\ti X,\bs\pi_Y,\bs C\t_Y \bs{\ti C}\bigr]},
\label{kh1eq1}
\e
where $\bs C\t_Y\bs{\ti C}$ is a fibre product of co-gauge-fixing
data. Then $\cup$ is associative, supercommutative, and compatible
with $\d$, and induces $\cup$ on $KH^*(Y;R)$ which is identified
with the usual $\cup$ on $H^*_\cs(Y;R)$ under $H_\cs^*(Y;R)\cong
KH^*(Y;R)$. There is also a {\it cap product\/} $\cap:KC_*(Y;R)\t
KC^*(Y;R)\ra KC_*(Y;R)$.

This illustrates the fact that {\it Kuranishi homology and
cohomology are very well behaved at the (co)chain level}. For
comparison, if one works with chains in singular homology, to take
the intersection of two singular chains one must perturb them so
they are transverse, and then triangulate the intersection by
simplices, making arbitrary choices. Working with Kuranishi
cochains, we just take cup products in a canonical way. So, {\it
using Kuranishi (co)chains can eliminate problems with
transversality}.

Here is how all this helps in the formation of virtual cycles and
chains. Consider for example the moduli space $\oM_{g,m}(M,J,\be)$
of genus $g$ stable $J$-holomorphic curves with $m$ marked points in
a symplectic manifold $M$. Then $\oM_{g,m}(M,J,\be)$ is a compact
oriented Kuranishi space without boundary, and has {\it evaluation
maps\/} $\bev_i:\oM_{g,m}(M,J,\be)\ra M$ for $i=1,\ldots,m$, strong
submersions sending a curve to its $i^{\rm th}$ marked point. We
define a virtual cycle for $\oM_{g,m}(M,J,\be)$ to be
$\bigl[\oM_{g,m}(M,J,\be),\bev_1\t\cdots\bev_m,\bs C\bigr]$ in
$KC^*(M^m;\Q)$, where $\bs C$ is any choice of co-gauge-fixing data
for~$\bigl(\oM_{g,m}(M,J,\be),\bev_1\t\cdots\bev_m\bigr)$.

Thus, {\it the moduli space is its own virtual cycle}. There is no
need to perturb the moduli space, or triangulate it with simplices;
we only have to choose (co-)gauge-fixing data, which is always
possible, and is a much milder and less disruptive process. Choosing
(co-)gauge-fixing data to satisfy prescribed boundary conditions or
other compatibilities is much easier than trying to do the same with
perturbations and singular chains. As a result, we can generally
translate relationships between moduli spaces into {\it exact
algebraic identities between virtual chains, at the (co)chain level,
using cup and cap products}. This will lead to huge simplifications
in the theory of Lagrangian Floer cohomology~\cite{AkJo}.

We also define {\it Kuranishi bordism\/} $KB_*(Y;R)$, and Poincar\'e
dual {\it Kuranishi cobordism\/} $KB^*(Y;R)$. Elements of
$KB_k(Y;R)$ are finite sums $\sum_{a\in A}\rho_a[X_a,\bs f_a]$ for
$\rho_a\in R$, where $[X_a,\bs f_a]$ is the isomorphism class of
$(X_a,\bs f_a)$, for $X_a$ a compact oriented Kuranishi space
without boundary with $\vdim X_a=k$, and $\bs f_a:X_a\ra Y$ is
strongly smooth. The relations in $KB_k(Y;R)$ are that $[X,\bs
f]+[X',\bs f]=[X\amalg X',\bs f\amalg\bs f']$ and $[\pd W,\bs
e\vert_{\pd W}]=0$, for $W$ compact oriented Kuranishi space with
boundary but without corners, with $\vdim W=k+1$, and $\bs e:W\ra Y$
strongly smooth. The definition for $KB^k(Y;R)$ is the same except
that $\vdim X_a=\dim Y-k$ and $\bs f_a:X_a\ra Y$, $\bs e:W\ra Y$ are
cooriented strong submersions. We define a {\it cup product\/} on
$KB^*(Y;R)$ by
\e
[X,\bs f]\cup[\ti X,\bs{\ti f}]=[X\t_Y\ti X,\bs\pi_Y],
\label{kh1eq2}
\e
as for \eq{kh1eq1}. Generally, to define or prove something for
Kuranishi (co)bordism, we just take the analogue in Kuranishi
(co)homology, restrict to Kuranishi spaces without boundary, and
omit (co-)gauge-fixing data.

In fact we define five variations on this idea, by imposing extra
conditions on or adding extra structure to $X_a,\bs f_a$. One
version, {\it trivial stabilizers Kuranishi (co)bordism\/}
$KB_*^\tr,KB^*_\trc(Y;R)$ has $KB_*^\tr(Y;\Z)\cong MSO_*(Y)$ and
$KB^*_\trc(Y;\Z)\cong MSO^*_\cs(Y)$ for manifolds $Y$, where
$MSO_*(Y)$ and $MSO^*_\cs(Y)$ are the conventional bordism and
compactly-supported cobordism groups defined by Atiyah \cite{Atiy}
using algebraic topology. However, the other four Kuranishi
(co)bordism theories are new topological invariants of $Y$. We show
that they are {\it huge\/}: $KB_k(Y;R)$ has infinite rank over $R$
for all nonempty orbifolds $Y$, commutative rings $R$ without
torsion (such as $R=\Z$ or $\Q$), and even $k\in\Z$, for instance.
So it is probably not feasible to compute them and write them down.

The author intends these Kuranishi (co)bordism theories as new tools
for studying (closed) Gromov--Witten theory. Since moduli spaces of
(closed) $J$-holomorphic curves are compact, oriented Kuranishi
spaces without boundary, they define classes in Kuranishi
(co)bordism groups. So we can define new Gromov--Witten type
invariants $GW_{g,m}^\Kb(M,\om, \be)$ in $KB_*(M^m;\Z)$ or
$KB_*(M^m\ab\t\oM_{g,m};\Z)$, which lift conventional Gromov--Witten
invariants $GW_{g,m}^\rsi(M,\om,\be)$ in $H^\rsi_*(M^m;\Q)$ or
$H_*^\rsi(M^m\!\t\!\oM_{g,m};\Q)$ in singular homology up to
Kuranishi bordism.

There are two important points about these new invariants
$GW_{g,m}^\Kb(M,\om,\be)$. Firstly, they take values in very large
groups, so they {\it contain more information than conventional
Gromov--Witten invariants}. Much of this extra information concerns
the {\it orbifold strata\/} of the moduli spaces, that is, for each
finite group $\Ga$ the $GW_{g,m}^\Kb(M,\om,\be)$ include data which
`counts' $J$-holomorphic curves in $M$ with symmetry group $\Ga$.
Secondly, they are defined in groups over $\Z$ rather than $\Q$.
This makes them a good context for studying {\it integrality
properties} of Gromov--Witten invariants. In \S\ref{kh63} we will
outline a method for using them to prove the Gopakumar--Vafa
Integrality Conjecture for Gromov--Witten invariants of Calabi--Yau
3-folds, which the author hopes to complete in~\cite{Joyc4}.

There will also be applications of Kuranishi (co)homology or
Kuranishi (co)bordism which are not related to $J$-holomorphic
curves. For example, in \cite{Joyc3} we will exploit the
particularly good behaviour of Kuranishi (co)chains to give a nice
formulation of the {\it String Topology} programme of Chas and
Sullivan \cite{ChSu}, in which the String Topology operations
satisfy the desired identities {\it at the chain level}. Also, we
show in \S\ref{kh54} that for compact manifolds $Y,$ trivial
stabilizers Kuranishi cobordism $KH^*_\trc(Y;\Z)$ is isomorphic to
cobordism $MSO^*(Y)$. This is the only differential-geometric
description of $MSO^*(Y)$ that the author knows, and it may be
useful to somebody.

We begin in Chapter \ref{kh2} by defining Kuranishi structures and
Kuranishi spaces, operations on them such as strongly smooth maps
and fibre products, and geometric structures such as
(co)orientations and (co-)almost complex structures. This is partly
taken from Fukaya and Ono \cite{FuOn1,FOOO}, partly modifying
\cite{FuOn1,FOOO}, and partly new. Chapter \ref{kh3} discusses {\it
gauge-fixing data} and {\it co-gauge-fixing data}, vital ingredients
of our Kuranishi (co)homology theories. Their main purpose is to
make automorphism groups finite. If $X$ is a compact Kuranishi
space, $\bs f:X\ra Y$ is strongly smooth, and $\bs G$ is
gauge-fixing data for $(X,\bs f)$, then $\Aut(X,\bs f,\bs G)$ is
finite, although $\Aut(X,\bs f)$ may be infinite.

Chapter \ref{kh4} studies {\it Kuranishi (co)homology\/}
$KH_*,KH^*(Y;R)$ for $R$ a $\Q$-algebra, and {\it effective
Kuranishi (co)homology\/} $KH_*^\ef,KH^*_\ec(Y;R)$ for $R$ a
commutative ring. We define pushforwards $h_*$ on chains and
homology, and pullbacks $h^*$ on cochains and cohomology, and cup
and cap products $\cup,\cap$ --- all the usual package for homology
and cohomology theories, and all work at the (co)chain level. Our
main results are Theorems \ref{kh4thm1} and \ref{kh4thm2}, which
show that $KH_*(Y;R)$ and $KH_*^\ef(Y;R)$ are isomorphic to singular
homology $H^\rsi_*(Y;R)$. Their proofs, deferred to Appendices
\ref{khA}--\ref{khC}, take up more than one third of the book. We
deduce using Poincar\'e duality that $KH^*(Y;R)$ and $KH^*_\ec(Y;R)$
are isomorphic to compactly-supported cohomology $H^*_\cs(Y;R)$. We
show by example that if we omitted (co-)gauge-fixing data $\bs G$
from the chains $[X,\bs f,\bs G]$ then we would have
$KH_*(Y;R)=KH^*(Y;R)=\{0\}$ for all $Y$, because of problems when
$\Aut(X,\bs f)$ is infinite.

Chapter \ref{kh5} defines two kinds of classical bordism groups
$B_*,B_*^\eo(Y;R)$, five kinds of Kuranishi bordism groups
$KB_*,KB^\eb_*,KB^\tr_*,KB^\ac_*,KB^\eac_*(Y;R)$, and their
Poincar\'e dual cobordism groups $KB^*,KB_\ecb^*,KB_\trc^*,
KB_\ca^*,KB_\eca^*(Y;R)$. We define pushforwards, cup and cap
products, and morphisms between the various (co)bordism groups, and
from Kuranishi (co)bordism to (effective) Kuranishi (co)homology. We
also define the {\it orbifold strata\/} $X^{\Ga,\rho}$ a of
Kuranishi space $X$. Restricting to orbifold strata induces
morphisms $\Pi^{\Ga,\rho}:KB_*(Y;R)\ra KB_{*-\dim\rho}(Y;R)$ of
Kuranishi (co)bordism groups, which change degree. Using these
composed with projection to homology, we prove that Kuranishi
(co)bordism groups are very large.

In Chapter \ref{kh6} we discuss applications to symplectic geometry.
Our aim is not to give a complete treatment, but only to demonstrate
the potential of our theories with a few examples, which may perhaps
inspire readers to learn about Kuranishi (co)homology or Kuranishi
(co)bordism, and use it in their own research. Sections
\ref{kh61}--\ref{kh63} discuss closed Gromov--Witten invariants
using Kuranishi bordism, and outline a method of attack on the
Gopakumar--Vafa Integrality Conjecture.

Sections \ref{kh64}--\ref{kh67} discuss moduli spaces of
$J$-holomorphic curves in $M$ with boundary in a Lagrangian $L$, and
open Gromov--Witten theory, using Kuranishi cochains and Kuranishi
cohomology. In \S\ref{kh66} we associate a gapped filtered $A_\iy$
algebra $\bigl(\widehat{KC}{}^*(L; \La_\nov^*)[1],\m\bigr)$ to $L$
in a much simpler way than Fukaya et al.\ \cite{FOOO}, and so define
{\it bounding cochains\/} $b$ and {\it Lagrangian Floer
cohomology\/} $HF^*(L,b;\La_\nov^*)$ for one Lagrangian $L$. Then
\S\ref{kh67} gives a new definition of {\it open Gromov--Witten
invariants} $OGW(M,L,\la)$ for a Lagrangian $L$ with $[L]=0\in
H_n(M;\Q)$, together with a choice of bounding cochain $b$. Both
$HF^*(L,b;\La_\nov^*)$ and the $OGW(M,L,\la)$ are independent of the
choice of almost complex structure $J$, in a certain sense.

Finally, Appendices \ref{khA}--\ref{khD} prove the main results of
Chapter \ref{kh4}, Theorems \ref{kh4thm1}, \ref{kh4thm2} and
\ref{kh4thm6}. The proofs of Theorems \ref{kh4thm1} and
\ref{kh4thm2} in Appendices \ref{khB} and \ref{khC} draw some
inspiration from results of Fukaya and Ono on multisections and
virtual cycles \cite[\S 3, \S 6]{FuOn1}, \cite[\S A1]{FOOO}, and
also from Fukaya et al.\ \cite[\S 30]{FOOO}, who consider a
situation in which transverse multisections must be carefully chosen
on a collection of Kuranishi spaces, with compatibility conditions
at the boundaries.

Significant parts of this book may read like a kind of commentary on
Fukaya, Ono at al.\ \cite{FuOn1,FOOO}, which includes suggesting
improvements and correcting errors. Lest this create the wrong
impression, I would like to record here my great admiration for
\cite{FuOn1,FOOO}, and my gratitude to Fukaya, Oh, Ohta and Ono for
the creativity and long labour involved in writing them.

In \cite{Joyc1} we provide a short `User's Guide' for Kuranishi
(co)homology and Kuranishi (co)bordism, consisting of the most
important definitions and results from this book, but omitting the
proofs. Some readers may find it helpful to begin with \cite{Joyc1}
rather than this book.
\medskip

\noindent{\it Acknowledgements.} I am grateful to Mohammed Abouzaid,
Manabu Akaho, Kenji Fukaya, Ezra Getzler, Shinichiroh Matsuo,
Yong-Geun Oh, Hiroshi Ohta, Kauru Ono, Paul Seidel, Ivan Smith and
Dennis Sullivan for useful conversations. This research was
partially supported by EPSRC grant EP/D07763X/1.

\section{Kuranishi structures and Kuranishi spaces}
\label{kh2}

We begin by discussing {\it Kuranishi spaces}, that is, topological
spaces with Kuranishi structures. Kuranishi structures were
introduced by Fukaya and Ono \cite{FuOn1} as a tool for studying
moduli problems, and further developed by Fukaya, Oh, Ohta and Ono
\cite{FOOO}. We take ideas freely from both, and also give some new
definitions and (in the author's opinion) some improvements of those
in \cite{FuOn1,FOOO}. Since Kuranishi structures are so important in
what follows, and are not yet in widespread use, we give plenty of
detail.

\subsection{Manifolds with corners and generalized corners}
\label{kh21}

We shall define four classes of manifolds. The first three are
entirely standard: {\it manifolds without boundary}, which are
locally modelled on $\R^n$, {\it manifolds with boundary}, locally
modelled on $\R^n$ or $[0,\iy)\t \R^{n-1}$, and {\it manifolds with
corners}, locally modelled on $[0,\iy)^k\t\R^{n-k}$
for~$k=0,\ldots,n$.

The fourth class, which we call {\it manifolds with generalized
corners}, or {\it g-corners} for short, is new. It is similar to,
but more general than, the notion of {\it locally polygonal set\/}
used by Fukaya et al.\ \cite[Def.~35.4]{FOOO}. We begin by defining
the local models for manifolds with g-corners in $\R^n$, which we
call {\it regions with g-corners}.

\begin{dfn} Let $W$ be a nonempty, connected, open set in $\R^n$,
and $f_1,\ldots,\ab f_N:W\ra\R$ be smooth maps, such that:
\begin{itemize}
\setlength{\itemsep}{0pt}
\setlength{\parsep}{0pt}
\item[(a)] For every subset $\{i_1,\ldots,i_l\}$ of
$\{1,\ldots,N\}$, the subset $S_{\{i_1,\ldots,i_l\}}=\bigl\{{\bf
x}\in W:f_{i_1}({\bf x})=\cdots=f_{i_l}({\bf x})=0\bigr\}$ is an
embedded submanifold of $\R^n$, without boundary.
\item[(b)] In (a), for all ${\bf x}\in S_{\{i_1,\ldots,i_l\}}$, we
have $\big\langle\d f_{i_1}\vert_{\bf x},\ldots,\d f_{i_l}\vert_{\bf
x}\big\rangle_\R=\bigl\{\al\in(\R^n)^*:\al\vert_{T_{\bf
x}S_{\{i_1,\ldots,i_l\}}}\equiv 0\bigr\}$ as vector subspaces of
$(\R^n)^*$.
\end{itemize}
Define $U=\bigl\{{\bf x}\in W:f_i({\bf x})\ge 0$,
$i=1,\ldots,N\bigr\}$. Then $U$ is a closed subset of $W$. Let
$U^\ci$ be the {\it interior\/} of $U$ in $W$, and
$\overline{U^\ci}$ the closure of $U^\ci$ in $W$. Suppose $U=
\overline{U^\ci}$. Then we call $U$ a {\it region with generalized
corners}, or {\it g-corners}, in~$\R^n$.
\label{kh2def1}
\end{dfn}

\begin{dfn} Let $X$ be a paracompact Hausdorff topological space.
\begin{itemize}
\setlength{\itemsep}{0pt}
\setlength{\parsep}{0pt}
\item[(i)] An {\it $n$-dimensional chart on\/ $X$ without boundary\/} is
a pair $(U,\phi)$, where $U$ is an open subset in $\R^n$, and
$\phi:U\ra X$ is a homeomorphism with a nonempty open set $\phi(U)$
in~$X$.
\item[(ii)] An {\it $n$-dimensional chart on\/ $X$ with boundary\/} for
$n\ge 1$ is a pair $(U,\phi)$, where $U$ is an open subset in $\R^n$
or in $[0,\iy)\t\R^{n-1}$, and $\phi:U\ra X$ is a homeomorphism with
a nonempty open set $\phi(U)$.
\item[(iii)] An {\it $n$-dimensional chart on\/ $X$ with corners\/}
for $n\ge 1$ is a pair $(U,\phi)$, where $U$ is an open subset in
$[0,\iy)^k\t\R^{n-k}$ for some $0\le k\le n$, and $\phi:U\ra X$ is a
homeomorphism with a nonempty open set~$\phi(U)$.
\item[(iv)] An {\it $n$-dimensional chart on\/ $X$ with generalized
corners}, or {\it g-corners}, is a pair $(U,\phi)$, where $U$ is
a region with g-corners in $\R^n$, and $\phi:U\ra X$ is a
homeomorphism with a nonempty open set~$\phi(U)$.
\end{itemize}
These are increasing order of generality, that is, (i) $\Rightarrow$
(ii) $\Rightarrow$ (iii) $\Rightarrow$ (iv).

Let $A,B$ be subsets of $\R^n$ and $\al:A\ra B$ be continuous. We
call $\al$ {\it smooth\/} if it extends to a smooth map between open
neighbourhoods of $A,B$, that is, if there exists an open subset
$A'$ of $\R^n$ with $A\subseteq A'$ and a smooth map
$\al':A'\ra\R^n$ with $\al'\vert_A\equiv\al$. If $A$ is open we can
take $A'=A$ and $\al'=\al$. We call $\al:A\ra B$ a {\it
diffeomorphism\/} if it is a homeomorphism and $\al:A\ra B$,
$\al^{-1}:B\ra A$ are smooth.

Let $(U,\phi),(V,\psi)$ be $n$-dimensional charts on $X$, which may
be without boundary, or with boundary, or with corners, or with
g-corners. We call $(U,\phi)$ and $(V,\psi)$ {\it compatible\/} if
$\psi^{-1}\ci\phi:\phi^{-1}\bigl(\phi(U)\cap\psi(V)\bigr)\ra
\psi^{-1}\bigl(\phi(U)\cap\psi(V)\bigr)$ is a diffeomorphism between
subsets of $\R^n$, in the sense above.

An $n$-{\it dimensional manifold structure\/} on $X$, which may be
{\it without boundary}, or {\it with boundary}, or {\it with
corners}, or {\it with generalized corners} ({\it with g-corners}),
is a system of pairwise compatible $n$-dimensional charts
$\bigl\{(U^i,\phi^i):i\in I\bigr\}$ on $X$ with $X=\bigcup_{i\in
I}\phi^i(U^i)$, where the $(U^i,\phi^i)$ are without boundary, or
with boundary, or with corners, or with g-corners, respectively.
When we just refer to a manifold, we will usually mean a manifold
with g-corners.
\label{kh2def2}
\end{dfn}

\begin{rem}{\bf (a)} Note that in Definition \ref{kh2def1}(a)
we do not prescribe the dimension of $S_{\{i_1,\ldots,i_l\}}$, and
in particular, we do not require $\dim S_{\{i_1,\ldots,i_l\}}=n-l$.
Definition \ref{kh2def1}(b) implies that $\dim
S_{\{i_1,\ldots,i_l\}}\ge n-l$. If $\dim S_{\{i_1,\ldots,i_l\}}=n-l$
then $\d f_{i_1},\ldots,\d f_{i_l}$ are linearly independent near
$S_{\{i_1,\ldots,i_l\}}$, but if $\dim S_{\{i_1,\ldots,i_l\}}>n-l$
then $\d f_{i_1},\ldots,\d f_{i_l}$ are linearly dependent along
$S_{\{i_1,\ldots,i_l\}}$.

If $\dim S_{\{i_1,\ldots,i_l\}}=n-l$ for all $\{i_1,\ldots,i_l\}$
then $U$ is an $n$-manifold with corners (not g-corners), locally
modelled on $[0,\iy)^l\t\R^{n-l}$. Thus, we only generalize
manifolds with corners by allowing~$\dim
S_{\{i_1,\ldots,i_l\}}>n-l$.

\noindent{\bf(b)} The purpose of requiring $U=\overline{U^\ci}$ in
Definition \ref{kh2def1} is to ensure $U$ has no pieces of dimension
less than $n$. For example, we could take $n=N=2$, $W=\R^2$ and
$f_1(x_1,x_2)=x_1$, $f_2(x_1,x_2)=-x_1$, which gives
$U=\bigl\{(0,x_2): x_2\in\R\bigr\}$, with dimension 1 rather than 2.
But then $U^\ci=\emptyset$, so $\overline{U^\ci}=\emptyset$,
and~$U\ne\overline{U^\ci}$.

\noindent{\bf(c)} An important class of examples of regions with
g-corners is given by taking $W=\R^n$ and each $f_i:\R^n\ra\R$ to be
linear. Then Definition \ref{kh2def1}(a),(b) hold automatically,
with each $S_{\{i_1,\ldots,i_l\}}$ a linear subspace of $\R^n$, and
$U$ is a {\it polyhedral cone} in $\R^n$, which satisfies $U=
\overline{U^\ci}$ if and only if $U^\ci\ne\emptyset$.

More generally, we can take $W=\R^n$ and each $f_i:\R^n\ra\R$ to be
affine, that is, of the form $f_i(x_1,\ldots,x_n)=a_1x_1+
\cdots+a_nx_n+b$. Then $U$ is a {\it convex polyhedron} in $\R^n$,
not necessarily compact, which satisfies $U=\overline{U^\ci}$ if and
only if $U^\ci\ne\emptyset$, that is, if and only if~$\dim U=n$.

For example, the solid octahedron in $\R^3$
\e
O=\bigl\{(x_1,x_2,x_3)\in\R^3:\md{x_1}+\md{x_2}+\md{x_3}\le 1\bigr\}
\label{kh2eq1}
\e
is a region with g-corners, which may be defined using the eight
affine functions
\begin{align*}
f_1(x_1,x_2,x_3)&=x_1+x_2+x_3+1,& f_2(x_1,x_2,x_3)&=x_1+x_2-x_3+1,\\
f_3(x_1,x_2,x_3)&=x_1-x_2+x_3+1,& f_4(x_1,x_2,x_3)&=x_1-x_2-x_3+1,\\
f_5(x_1,x_2,x_3)&=-x_1+x_2+x_3+1,&f_6(x_1,x_2,x_3)&=-x_1+x_2-x_3+1,\\
f_7(x_1,x_2,x_3)&=-x_1-x_2+x_3+1,& f_8(x_1,x_2,x_3)&=-x_1-x_2-x_3+1.
\end{align*}
This is not a manifold with corners, since four 2-dimensional faces
of $O$ meet at the vertex $(1,0,0)$, but three 2-dimensional faces
of $[0,\iy)^3$ meet at the vertex $(0,0,0)$, so $O$ near $(1,0,0)$
is not locally modelled on $[0,\iy)^3$ near~$(0,0,0)$.

\noindent{\bf(d)} To a first approximation, regions with g-corners
and manifolds with corners look locally like convex polyhedrons in
$\R^n$. Let $U$ be a region with g-corners defined using
$W,f_1,\ldots,f_N$ in Definition \ref{kh2def1}, and let ${\bf y}\in
U$. Define $\hat W=\R^n$ and $\hat f_1,\ldots,\hat f_N:\R^n\ra\R$ by
$\hat f_i({\bf x})=f_i({\bf y})+\d f_i\vert_{\bf y}({\bf x}-{\bf
y})$, so that $\hat f_1,\ldots,\hat f_N$ are the unique affine
functions on $\R^n$ with $\hat f_i({\bf x})=f_i({\bf x})+O(\ms{{\bf
x}-{\bf y}})$.

Let $\hat U$ be the region with g-corners defined using $\hat W,\hat
f_1,\ldots,\hat f_N$. Then $\hat U$ is a convex polyhedron in $\R^n$
(not necessarily compact), and $U,\hat U$ are basically isomorphic
near $\bf y$. In particular, from Definition \ref{kh2def1}(b) it
follows that for all subsets $\{i_1,\ldots,i_l\}\subseteq
\{1,\ldots,N\}$, the submanifolds $S_{\{i_1,\ldots,i_l\}}$ from
$f_{i_1},\ldots,f_{i_l}$ and $\hat S_{\{i_1,\ldots,i_l\}}$ from
$\hat f_{i_1},\ldots,\hat f_{i_l}$ have the same tangent space at
$\bf y$, and the same dimension near $\bf y$, and this implies that
$U$and $\hat U$ have the same stratification into faces and corners
of different dimensions near~$\bf y$.

\noindent{\bf(e)} However, despite (d), regions with g-corners in
$\R^n$ are a {\it more general class of local models} than convex
polyhedrons in $\R^n$, and yield a more general notion of manifolds
with g-corners than we would have obtained by taking charts to be
modelled on open sets in convex polyhedra in~$\R^n$.

To see this, note that polyhedral cones $C$ in $\R^n$ are classified
up to the action of GL$(n,\R)$ by both {\it discrete parameters}
--- the number of $k$-dimensional faces for $0<k<n$, and how they
intersect --- and {\it continuous parameters} --- if $C$ has at
least $n+2$ codimension one faces then the family of polyhedral
cones with the same discrete structure as $C$ up to GL$(n,\R)$ is in
general positive-dimensional. Classifying polyhedral cones in $\R^n$
up to GL$(n,\R)$ is important, as this is what is preserved by our
notion of diffeomorphism of subsets of $\R^n$ in
Definition~\ref{kh2def2}.

If $U$ is a convex polyhedron in $\R^n$, and $F$ is a connected,
open dimension $l$ face of $U$, then $U$ is modelled on a fixed
polyhedral cone $C_F$ in $\R^n$ near $\bf y$ for all $\bf y$ in $F$.
However, if $U$ is a region with g-corners in $\R^n$, and $F$ is a
connected, open, dimension $l$ face of $U$ for $l>0$, then $U$ is
modelled on a polyhedral cone $C_{\bf y}$ in $\R^n$ near $\bf y$ for
each $\bf y$ in $F$, where the discrete parameters of $C_{\bf y}$ up
to GL$(n,\R)$ are independent of $\bf y$, but {\it the continuous
parameters can vary with\/}~$\bf y$.

This can only happen when $n\ge 4$, since the local models for
manifolds with g-corners are unique near strata of codimension 1 and
2, so we need both $l\ge 3$ and $\dim F=n-l>0$ to have nontrivial
continuous variation of $C_{\bf y}$ along $F$. The simplest example
would be a region with g-corners in $\R^4$ in which five
3-dimensional faces come together along a 1-dimensional edge.

The point of the rather general, non-explicit local models of
Definition \ref{kh2def1} is to allow this local variation of the
polyhedral cone model.

\noindent{\bf(f)} The only places in this book where the details of
the definition of manifolds with g-corners matter is in the material
on {\it tent functions\/} in Appendix \ref{khA}, and their
application in Appendices \ref{khB} and \ref{khC} to prove Theorems
\ref{kh4thm1} and \ref{kh4thm2}. For the rest of the book, all we
need is that manifolds with g-corners have a well-behaved notion of
boundary.
\label{kh2rem1}
\end{rem}

Next we define the {\it boundary\/} $\pd X$ of a manifold $X$ with
boundary or \hbox{(g-)}\ab corners. To motivate the definition,
consider $[0,\iy)^2$ in $\R^2$, regarded as a manifold with corners.
If we took $\pd\bigl([0,\iy)^2\bigr)$ to be the subset
$[0,\iy)\t\{0\}\ab\cup\ab\{0\}\t[0,\iy)$ of $[0,\iy)^2$, then
$\pd\bigl([0,\iy)^2\bigr)$ would not be a manifold with corners near
$(0,0)$. Instead, we take $\pd\bigl([0,\iy)^2\bigr)$ to be the {\it
disjoint union\/} of the two boundary strata $[0,\iy)\t\{0\}$ and
$\{0\}\t[0,\iy)$. This is a manifold with boundary, but now
$\pd\bigl([0,\iy)^2\bigr)$ is {\it not a subset of\/} $[0,\iy)^2$,
since two points in $\pd\bigl([0,\iy)^2\bigr)$ correspond to the
vertex $(0,0)$ in~$[0,\iy)^2$.

\begin{dfn} Let $U$ be a region with g-corners in $\R^n$, defined
using $W,f_1,\ab\ldots,\ab f_N$ as in Definition \ref{kh2def1}. If
$u\in U$, let $\{i_1,\ldots,i_l\}$ be the subset of $i$ in
$\{1,\ldots,N\}$ for which $f_i(u)=0$, and let
$S_{\{i_1,\ldots,i_l\}}$ be as in Definition \ref{kh2def1}. Then
$S_{\{i_1,\ldots,i_l\}}$ is an embedded submanifold of $W$
containing $u$. Define the {\it codimension of\/ $u$ in\/} $U$ to
be~$n-\dim S_{\{i_1,\ldots,i_l\}}$.

It is not difficult to see that if $u$ has codimension $m$, then
$n-m$ is the maximum dimension of embedded submanifolds $S$ in
$\R^n$ with $u\in S\subseteq U\subseteq\R^n$, and this maximum is
achieved by $S_{\{i_1,\ldots,i_l\}}$. Therefore the codimension of
$u$ depends only on $u$ and the set $U$ up to diffeomorphism (in the
sense of Definition \ref{kh2def2}), and is independent of the choice
of $W,f_1,\ldots,f_N$ used to define~$U$.

Let $X$ be an $n$-manifold with g-corners. Then $X$ is covered by a
system of charts $\bigl\{(U^i,\phi^i):i\in I\bigr\}$, with each
$U^i$ a region with g-corners in $\R^n$. Let $x\in X$, so that
$x=\phi^i(u)$ for some $i\in I$ and $u\in U^i$. Define the {\it
codimension of\/ $x$ in\/} $X$ to be the codimension of $u$ in
$U^i$. This is independent of the choice of chart $(U^i,\phi^i)$,
since if $x=\phi^i(u)=\phi^j(u')$, then $U^i$ near $u$ is
diffeomorphic to $U^j$ near $u'$ using $\phi_j\ci(\phi_i)^{-1}$. For
each $k=0,\ldots,n$, define $S_k(X)$ to be the subset of $x\in X$ of
codimension $k$, the {\it codimension $k$ stratum\/} of $X$. Then
$X=\coprod_{k=0}^nS_k(X)$. One can show that each $S_k(X)$ has the
structure of an $(n-k)$-manifold without boundary, with closure
$\overline{S_k(X)}=\bigcup_{l=k}^nS_l(X)$.

Let $x\in\overline{S_1(X)}=\bigcup_{l=1}^nS_l(X)$. A {\it local
boundary component\/ $B$ of\/ $X$ at\/} $x$ is a local choice of
connected component of $S_1(X)$ near $x$. That is, for each
sufficiently small open neighbourhood $U$ of $x$ in $X$, $B$ gives a
choice of connected component $V$ of $U\cap S_1(X)$ with $x\in\bar
V$, and any two such choices $U,V$ and $U',V'$ must be compatible in
the sense that $x\in\overline{(V\cap V')}$. It is easy to see that
if $x\in S_1(X)$ then $x$ has a unique local boundary component $B$,
and if $x\in S_2(X)$ then $x$ has exactly two local boundary
components $B_1,B_2$, and if $x\in S_k(X)$ for $k\ge 2$ then $x$ has
at least $k$ local boundary components.

As a set, define the {\it boundary\/} $\pd X$ of $X$ by
\begin{equation*}
\pd X=\bigl\{(x,B):\text{$x\in\overline{S_1(X)}$, $B$ is a local
boundary component for $X$ at $x$}\bigr\}.
\end{equation*}
Then $\pd X$ naturally has the structure of an $(n-1)$-{\it manifold
with g-corners}. To see this, let $(x,B)\in\pd X$, and suppose
$(U^i,\phi^i)$ is a chart on $X$ with $x\in\phi^i(U^i)$. Then $U^i$
is a region with g-corners in $\R^n$ defined using some
$W,f_1,\ldots,f_N$, with $x=\phi^i(u)$ for some unique $u\in U^i$.
Now $(\phi^i)^{-1}(B)$ is a local boundary component of $U^i$ at
$u$. From Definition \ref{kh2def1} it follows that
$(\phi^i)^{-1}(B)$ coincides with $f_j^{-1}(0)$ near $u$ for
some~$j=1,\ldots,N$.

Definition \ref{kh2def1}(b) with $l=1$ and $i_1=j$ then implies that
$\d f_j$ is nonzero on $f_j^{-1}(0)=S_{\{j\}}$, so $f_j^{-1}(0)$ is
a closed, embedded $(n-1)$-submanifold of $U^j$ containing $u$. Thus
we can choose a diffeomorphism $\psi:\ti W\ra f_j^{-1}(0)\subset
U^i$ from an open subset $\ti W\subset\R^{n-1}$ to an open
neighbourhood $\psi(\ti W)$ of $u$ in $f_j^{-1}(0)$. Define
functions $\ti f_1,\ldots,\ti f_{N-1}:\ti W\ra\R$ to be $\ti
f_a=f_a\ci \psi$ for $a=1,\ldots,j-1$ and $\ti f_a=f_{a+1}\ci\psi$
for $a=j,\ldots,N-1$. Then $\ti W,\ti f_1,\ldots,\ti f_{N-1}$
satisfy Definition \ref{kh2def1} and define a region with corners
$\ti U$ in $\R^{n-1}$, with $\psi(\ti U)$ an open neighbourhood of
$u$ in $U^i\cap f_j^{-1}(0)$. Make $W,\ti U$ smaller if necessary to
ensure that $\psi(\ti U)\subseteq\overline{(\phi^i)^{-1}(B)}$.
Define $\ti\phi:\ti U\ra\pd X$ by $\ti\phi(\ti
u)=\bigl(\phi^i\ci\psi(\ti u),B\bigr)$. Then $(\ti U,\ti\phi)$ is an
$(n-1)$-dimensional chart on $\pd X$ with g-corners, with
$(x,B)\in\ti\phi(\ti U)$. Any two such charts are compatible, by
compatibility of the $(U^i,\phi^i)$. Thus, covering $\pd X$ with a
system of these charts makes it into an $(n-1)$-manifold with
g-corners.

If $X$ is a manifold without boundary then $\pd X=\es$, and if $X$
is a manifold with boundary then $\pd X$ is a manifold without
boundary. In this case we may regard $\pd X$ as a subset of $X$,
since each $x\in S_1(X)$ has a unique local boundary component $B$,
so $(x,B)\mapsto x$ is injective. If $X$ has corners (not g-corners)
then $\pd X$ has corners (not g-corners). If $X$ is compact then
$\pd X$ is compact.
\label{kh2def3}
\end{dfn}

We define {\it smooth maps\/} and {\it submersions\/} of manifolds
with g-corners.

\begin{dfn} Let $X,Y$ be manifolds with g-corners, with $\dim X=m$,
$\dim Y=n$, and $f:X\ra Y$ be a continuous map. We call $f$ {\it
smooth\/} if whenever $(U,\phi),(V,\psi)$ are charts on $X,Y$
respectively compatible with the manifold structures on $X,Y$, where
$U\subseteq\R^m$ and $V\subseteq\R^n$ are regions with g-corners,
then
\e
\psi^{-1}\ci f\ci\phi:(f\ci\phi)^{-1}(\psi(V))\longra
\psi^{-1}(f\ci\phi(U))
\label{kh2eq2}
\e
is a smooth map from $(f\ci\phi)^{-1}(\psi(V))\subseteq
U\subseteq\R^m$ to $\psi^{-1}(f\ci\phi(U))\subseteq V\subseteq\R^n$.
That is, as in Definition \ref{kh2def2}, $\psi^{-1}\ci f\ci\phi$
should extend to a smooth map between open subsets of $\R^m,\R^n$
containing $(f\ci\phi)^{-1}(\psi(V))$ and~$\psi^{-1}(f\ci\phi(U))$.

We will define when a smooth map $f:X\ra Y$ is a {\it submersion}.
When $\pd X=\es=\pd Y$, the definition is well-known. But when $X,Y$
have boundary and (g-)corners, we will impose complicated extra
conditions on $f$ over $\pd^kX$ and $\pd^lY$ for $k,l\ge 0$. Suppose
first that $f:X\ra Y$ is a smooth map and $\pd X=\es=\pd Y$. We call
$f$ a {\it submersion\/} if for all $(U,\phi),(V,\psi)$ as above and
all $u\in (f\ci\phi)^{-1}(\psi(V))$, the linear map
$\d\bigl(\psi^{-1}\ci f\ci\phi\bigr)\vert_u:\R^m\ra\R^n$ is
surjective. This implies $m\ge n$. More geometrically, $X,Y$ have
tangent bundles $TX,TY$, and the derivative $\d f$ is a morphism of
vector bundles $\d f:TX\ra f^*(TY)$, and $f$ is a submersion if $\d
f$ is surjective.

Next, suppose that $f:X\ra Y$ is a smooth map and $\pd Y=\es$, but
do not assume that $\pd X=\es$. We call $f$ a {\it submersion\/} if
for all $k\ge 0$, the map $f\vert_{\pd^kX}:\pd^kX\ra Y$ satisfies
$\d\bigl(f\vert_{\pd^kX}\bigr):T(\pd^kX)\ra f\vert_{\pd^kX}^*(TY)$
is surjective. This implies that $\pd^kX=\es$ for all $k\ge m-n$.
Also, if $f:X\ra Y$ is a submersion then $f\vert_{\pd X}:X\ra Y$ is
a submersion.

Finally, suppose that $f:X\ra Y$ is a smooth map, but do not assume
that $\pd X=\es$ or $\pd Y=\es$. We give an {\it inductive
definition\/} of what it means for $f$ to be a submersion. We call
$f$ a submersion if:
\begin{itemize}
\setlength{\itemsep}{0pt}
\setlength{\parsep}{0pt}
\item[(a)] $\d f:TX\ra f^*(TY)$ is a surjective morphism of vector
bundles over $X$;
\item[(b)] there is a unique decomposition $\pd
X=\pd_+^fX\amalg\pd_-^fX$, where $\pd_+^fX,\pd_-^fX$ are open and
closed subsets of $\pd X$, such that $f_+=f\vert_{\pd_+^fX}$ maps
$(\pd_+^fX)^\ci\ra Y^\ci$, and $f\vert_{\pd_-^fX}:\pd_-^fX\ra Y$
factors uniquely as $f\vert_{\pd_-^fX}=f_-\ci\io$, where
$f_-:\pd_-^fX\ra\pd Y$ is a smooth map and $\io:\pd Y\ra Y$ is the
natural immersion, and $f_-$ maps $(\pd_+^fX)^\ci\ra(\pd Y)^\ci$;
and
\item[(c)] $f_+:\pd_+^fX\ra Y$ and $f_-:\pd_-^fX\ra\pd Y$ are both
submersions.
\end{itemize}
This definition is recursive, as the definition of $f$ being a
submersion involves $f_\pm$ being submersions in (c), but since
$\dim\pd_\pm X=\dim X-1$, the definition makes sense by induction
on~$\dim X=0,1,2,\ldots$.

The following example may help readers to understand the definition.
Let $W,Y$ be manifolds with g-corners, and set $X=W\t Y$, with
$f:X\ra Y$ the projection $W\t Y\ra Y$. Then $f$ is a submersion in
the sense above. We have $\pd X=\pd(W\t Y)=(\pd W\t Y)\amalg(W\t\pd
Y)$. Then $\pd_+^fX=\pd W\t Y$, with $f_+:\pd_+^fX\ra Y$ the
projection $\pd W\t Y\ra Y$, and $\pd_-^fX=W\t\pd Y$, with
$f_-:\pd_-^fX\ra\pd Y$ the projection $W\t\pd Y\ra\pd Y$. In
general, one can think of submersions $f:X\ra Y$ as being {\it
locally\/} like a projection $W\t Y\ra Y$ for $W,Y$ manifolds with
g-corners. If $Y$ is a manifold with g-corners then $\id_Y:Y\ra Y$
is a submersion, taking $\pd_+^{\id_Y}Y=\es$ and~$\pd_-^{\id_Y}Y=\pd
Y$.
\label{kh2def4}
\end{dfn}

One justification for this complicated definition of submersions is
that if $f:X\ra Y$, $f':X'\ra Y$ are smooth maps of manifolds with
g-corners with at least one of $f,f'$ a submersion, there is a
well-behaved notion of {\it fibre product\/} $X\t_{f,Y,f'}X'$, which
is a manifold with g-corners.

\begin{prop} Suppose\/ $X,X',Y$ are manifolds with g-corners
and\/ $f:X\ra Y,$ $f':X'\ra Y$ are smooth maps, with at least one
of\/ $f,f'$ a submersion. Define the \begin{bfseries}fibre
product\end{bfseries} $X\t_{f,Y,f'}X'$ or $X\t_YX'$ by
\e
X\t_{f,Y,f'}X'=X\t_YX'=\bigl\{(p,p')\in X\t X':f(p)=f'(p')\bigr\}.
\label{kh2eq3}
\e
Then $X\t_YX'$ is a closed, embedded submanifold of\/ $X\t X',$ and
so $X\t_YX'$ is a manifold with g-corners, which is compact if\/
$X,X'$ are compact.

Define maps $\pi_X:X\t_YX'\ra X,$ $\pi_{X'}:X\t_YX'\ra X',$ and\/
$\pi_Y:X\t_YX'\ra Y$ by $\pi_X:(p,p')\mapsto p,$
$\pi_{X'}:(p,p')\mapsto p',$ and\/ $\pi_Y:(p,p')\mapsto
f(p)=f'(p')$. Then $\pi_X,\pi_{X'},\pi_Y$ are smooth, and\/ $\pi_X$
is a submersion if\/ $f'$ is a submersion, and\/ $\pi_{X'}$ is a
submersion if\/ $f$ is a submersion, and\/ $\pi_Y$ is a submersion
if both\/ $f',f'$ are submersions.

If\/ $\pd Y=\es$ then we have a natural diffeomorphism
\e
\pd\bigl(X\t_{f,Y,f'}X'\bigr)\cong \bigl(\pd X\t_{f\vert_{\pd
X},Y,f'}X'\bigr)\amalg \bigl(X\t_{f,Y,f'\vert_{\pd X}}\pd X'\bigr),
\label{kh2eq4}
\e
where if\/ $f$ is a submersion then $f\vert_{\pd X}:\pd X\ra Y$ is a
submersion, and the same for $f'$, so the fibre products in
\eq{kh2eq4} make sense.

In general, allowing $\pd Y\ne\es,$ if\/ $f$ is a submersion then
\e
\pd\bigl(X\t_{f,Y,f'}X'\bigr)\cong\bigl(\pd_+^fX\t_{f_+,Y,f'}X'
\bigr)\amalg\bigl(X\t_{f,Y,f'\vert_{\pd X}}\pd X'\bigr),
\label{kh2eq5}
\e
and if\/ $f'$ is a submersion then
\e
\pd\bigl(X\t_{f,Y,f'}X'\bigr)\cong\bigl(\pd X\t_{f\vert_{\pd X},
Y,f'}X'\bigr)\amalg \bigl(X\t_{f,Y,f'_+}\pd_+^{f'} X'\bigr),
\label{kh2eq6}
\e
and if both\/ $f,f'$ are submersions then both\/
\eq{kh2eq5}--\eq{kh2eq6} hold, and also
\ea
\begin{split}
\pd\bigl(X\t_{f,Y,f'}X'\bigr)&\cong\bigl(\pd_+^fX\t_{f_+,Y,f'}X'\bigr)
\amalg\bigl(X\t_{f,Y,f'_+}\pd_+^{f'}X'\bigr)\\
&\qquad\qquad\qquad\amalg\bigl( \pd_-^fX\t_{f_-,\pd
Y,f_-'}\pd_-^{f'}X'\bigr),
\end{split}
\label{kh2eq7}\\
\text{with\/}\quad
\pd_+^{\pi_Y}\bigl(X\t_{f,Y,f'}X'\bigr)&\cong\bigl(\pd_+^fX\t_{f_+,Y,f'}X'\bigr)
\amalg\bigl(X\t_{f,Y,f'_+}\pd_+^{f'}X'\bigr),
\label{kh2eq8}\\
\text{and\/}\quad \pd_-^{\pi_Y}\bigl(X\t_{f,Y,f'}X'\bigr)&\cong
\pd_-^fX\t_{f_-,\pd Y,f_-'}\pd_-^{f'}X'.
\label{kh2eq9}
\ea
\label{kh2prop1}
\end{prop}

Note that Proposition \ref{kh2prop1} {\it would not be true\/} if in
defining submersions $f:X\ra Y$ we had not imposed conditions over
$\pd^kX$ and $\pd^lY$ for all $k,l\ge 0$, since otherwise $X\t_YX'$
might fail to be a manifold with g-corners on strata coming from
$\pd^kX$ or $\pd^{k'}X'$. Thus we do need the complex, inductive
definition in Definition \ref{kh2def4}. The proof of Proposition
\ref{kh2prop1} is straightforward. The main thing is to check that
$X\t_YX'$ is a {\it submanifold\/} of $X\t X'$, by working with
coordinate charts on $X,X'$. The submersion conditions on
$\pd^kX,\pd^{k'}X'$ ensure that $X\t_YX'$ intersects the parts of
the boundary and corners of $X\t X'$ coming from $\pd^kX,\pd^{k'}X'$
suitably transversely.

A manifold $X$ with (g-)corners has a natural
involution~$\si:\pd^2X\ra\pd^2X$.

\begin{dfn} Let $X$ be an $n$-manifold with g-corners. Then
$\pd^2X=\pd(\pd X)$ is an $(n-2)$-manifold with g-corners. Points of
$\pd^2X$ are of the form $\bigl((x,B_1),B'\bigr)$, where $(x,B_1)\in
\pd X$, so that $x\in X$ and $B_1$ is a local boundary component for
$X$ at $x$, and $B'$ is a local boundary component for $\pd X$ at
$(x,B_1)$. We would like to understand $B'$ in terms of local
boundary components for $X$ at $x$. Let $U$ be a sufficiently small
open neighbourhood of $x\in X$. Then $B_1$ determines a connected
component $V_1$ of $U\cap S_1(X)$ with $x\in\bar V_1$, and $B'$
determines a connected component $V'$ of $U\cap\bar V_1\cap S_2(X)$.

From above, if $y\in S_2(X)$ then there are exactly two local
boundary components for $X$ at $y$. Applying this to $y\in V'$ with
the same $U$, we find that there exists a unique second connected
component $V_2$ of $U\cap S_1(X)$ with $V_1\ne V_2$ such that $V'$
is a subset of $U\cap\bar V_2\cap S_2(X)$. This $V_2$ determines a
unique local boundary component $B_2$ of $X$ at $x$, with $B_1\ne
B_2$. Conversely, $B_1,B_2$ determine $B'$, since $V'$ is the unique
connected component of $U\cap\bar V_1\cap\bar V_2\cap S_2(X)$ with
$x\in\bar V'$. Let us write $B_1\cap B_2$ for $B'$. Then each point
in $\pd^2X$ may be uniquely written as $\bigl((x,B_1),B_1\cap
B_2\bigr)$, where $x\in X$ and $B_1,B_2$ are distinct local boundary
components of $X$ at $x$ which intersect in codimension 2 in~$X$.

Define $\si:\pd^2X\ra\pd^2X$ by $\si:\bigl((x,B_1),B_1\cap
B_2\bigr)\mapsto\bigl((x,B_2),B_2\cap B_1\bigr)$. Then $\si$ is a
{\it natural, smooth, free involution\/} of $\pd^2X$. If $X$ is
oriented, so that $\pd X$ and $\pd^2X$ are oriented, then $\si$ is
orientation-reversing.
\label{kh2def5}
\end{dfn}

These involutions $\si$ are important in problems of extending data
from $\pd X$ to $X$. We call this the {\it Extension Principle}, and
we state it for arbitrary maps, for continuous functions or
sections, and for smooth functions and sections.

\begin{princ}{\bf$\!\!$(Extension Principle)(a)} Let $X$ be a manifold
with g-corners, $P$ a nonempty set, and $H:\pd X\ra P$ an arbitrary
map. Then a necessary and sufficient condition for there to exist
$G:X\ra P$ with $G\vert_{\pd X}\equiv H$ is that $H\vert_{\pd^2X}$
is invariant under $\si:\pd^2X\ra\pd^2X$ in Definition
\ref{kh2def5}, that is,~$H\vert_{\pd^2X}\ci\si\equiv
H\vert_{\pd^2X}$.

\noindent{\bf(b)} Let $X$ be a manifold with g-corners, $E\ra X$ a
vector bundle, $\ze:\pd X\ra\R$ a continuous map, and $t$ a
continuous section of $E\vert_{\pd X}$. Then necessary and
sufficient conditions for there to exist continuous $\eta:X\ra\R$
with $\eta\vert_{\pd X}\equiv\ze$, and a continuous section $s$ of
$E$ with $s\vert_{\pd X}\equiv t$, are that $\ze\vert_{\pd^2X}$ and
$t\vert_{\pd^2X}$ should be invariant under $\si:\pd^2X\ra\pd^2X$.
Here $\si:\pd^2X\ra\pd^2X$ has a natural lift
$\hat\si:E\vert_{\pd^2X}\ra E\vert_{\pd^2X}$, and $t\vert_{\pd^2X}$
is invariant under $\si$ means that~$t\vert_{\pd^2X}\ci\si\equiv
\hat\si\ci t\vert_{\pd^2X}$.

\noindent{\bf(c)} Let $X$ be a manifold with corners ({\it not\/}
g-corners), $E\ra X$ a vector bundle, $\ze:\pd X\ra\R$ be smooth,
and $t$ a smooth section of $E\vert_{\pd X}$. Then necessary and
sufficient conditions for there to exist smooth $\eta:X\ra\R$ with
$\eta\vert_{\pd X}\equiv\ze$, and a smooth section $s$ of $E$ with
$s\vert_{\pd X}\equiv t$, are that $\ze\vert_{\pd^2X}$ and
$t\vert_{\pd^2X}$ should be invariant under~$\si:\pd^2X\ra\pd^2X$.
\label{kh2pri}
\end{princ}

The point is that for manifolds with (g-)corners, the `inclusion'
map $\io:\pd X\ra X$ taking $\io:(x,B)\mapsto x$ is {\it not
injective} over $S_k(X)$ for $k\ge 2$. By the `restriction'
$G\vert_{\pd X}$, we mean $G\ci\io$. Clearly, in (a) there can only
exist $G:X\ra P$ with $G\vert_{\pd X}=H$ if whenever $p\ne q\in\pd
X$ with $\io(p)=\io(q)$ we have $H(p)=H(q)$. Then $p,q$ lie over
$\pd^2X$, and it is easy to show that $H(p)=H(q)$ for all such $p,q$
if and only if $H\vert_{\pd^2X}$ is $\si$-invariant.

This establishes necessity in (a), and for (b),(c) the argument is
the same. Sufficiency in (a) is more-or-less trivial, we can take
$G\vert_{X^\ci}\equiv p\in P$. For (b),(c), we prove local
sufficiency near any point of $X$ that lifts to $\pd X$, then join
local choices of $\eta,s$ with a partition of unity to get a global
choice.

Note that if $X$ has g-corners, not corners, then in (c) for
$\ze\vert_{\pd^2X},t\vert_{\pd^2X}$ to be $\si$-invariant is a
necessary but {\it not\/} sufficient for $\ze,t$ to extend smoothly
to $X$. As an example, take $X$ to be the octahedron $O$ in
\eq{kh2eq1}, and let $\ze:\pd O\ra\R$ be smooth with
$\ze\vert_{\pd^2X}\equiv\ze\vert_{\pd^2X}\ci\si$. Then at the vertex
$(1,0,0)$ of $O$, we can prescribe the derivatives of $\ze$ along
the four edges of $O$ meeting at $(1,0,0)$ independently. However,
if $\ze=\eta\vert_{\pd O}$ with $\eta$ smooth then these four
derivatives would be determined by $\d\eta\vert_{(1,0,0)}$, which
has only three parameters, a contradiction.

That is, {\it the smooth Extension Principle fails for manifolds
with g-corners}. This will cause problems in proofs in Appendices
\ref{khA}--\ref{khC}. To get round them, we develop a notion of {\it
piecewise smooth extensions\/} in \S\ref{khA14}, which do work for
manifolds and orbifolds with g-corners.

Some very similar ideas can be found in Fukaya et al.\ \cite{FOOO},
which is not surprising, as they also have to solve essentially the
same problems that we need the Extension Principle for. They prove a
result close to Principle \ref{kh2pri}(c) in
\cite[Lem.~30.121]{FOOO}, and in \cite[Rem.~30.122]{FOOO} they note
that this result does not hold for polygons, which is their way of
saying that the smooth Extension Principle fails for manifolds with
g-corners, and that instead one can work with piecewise smooth
multisections, just as we use piecewise smooth extensions for
manifolds with g-corners.

\subsection{Orbifolds and orbibundles}
\label{kh22}

Next we introduce {\it orbifolds}. Some useful references are Adem,
Leida and Ruan \cite{ALR}, the appendix in Chen and Ruan
\cite{ChRu1}, and Satake's original paper~\cite{Sata}.

\begin{dfn} Let $X$ be a paracompact Hausdorff topological space.
For orbifolds {\it without boundary}, an {\it $n$-dimensional
orbifold chart on\/} $X$ is a triple $(U,\Ga,\phi)$, where $\Ga$ is
a finite group with a linear representation on $\R^n$, $U$ is a
$\Ga$-invariant open subset in $\R^n$, and $\phi:U/\Ga\ra X$ is a
homeomorphism with a nonempty open set $\phi(U/\Ga)$ in~$X$.

For orbifolds {\it with boundary\/} we also allow orbifold charts of
the form $(U,\Ga,\phi)$, where $\Ga$ is a finite group with a linear
representation on $\R^n$ preserving the subset $[0,\iy)\t\R^{n-1}$,
$U$ is a $\Ga$-invariant open subset in $[0,\iy)\t\R^{n-1}$, and
$\phi:U/\Ga\ra X$ is a homeomorphism with a nonempty open set
$\phi(U/\Ga)$ in~$X$.

For orbifolds {\it with corners\/} we also allow orbifold charts
$(U,\Ga,\phi)$, where $\Ga$ is a finite group with a linear
representation on $\R^n$ preserving the subset $[0,\iy)^k\t\R^{n-k}$
for some $0\le k\le n$, $U$ is a $\Ga$-invariant open subset in
$[0,\iy)^k\t\R^{n-k}$, and $\phi:U/\Ga\ra X$ is a homeomorphism with
a nonempty open set $\phi(U/\Ga)$ in~$X$.

For orbifolds {\it with generalized corners\/} ({\it g-corners}) we
also allow orbifold charts $(U,\Ga,\phi)$, where $\Ga$ is a finite
group with a linear representation on $\R^n$, $U$ is a
$\Ga$-invariant region with g-corners in $\R^n$, and $\phi:U/\Ga\ra
X$ is a homeomorphism with a nonempty open set $\phi(U/\Ga)$ in~$X$.

For such a chart $(U,\Ga,\phi)$, let $p\in\Im\phi\subseteq X$,
choose $u_p\in U$ with $\phi(\Ga u_p)=p$, and define
$\Ga_p=\{\ga\in\Ga:\ga\,u_p=u_p\}$. Choose a $\Ga_p$-invariant open
neighbourhood $U_p$ of $u_p$ in $U$ such that if $\ga\in\Ga$ and
both $u,\ga\,u$ lie in $U_p$ then $\ga\in\Ga_p$; this holds if $U_p$
is sufficiently small. Then the map $U_p/\Ga_p\ra(\Ga
U_p)/\Ga\subseteq U/\Ga$ given by $\Ga_p u\mapsto\Ga u$ is a
homeomorphism. Define $\phi_p:U_p/\Ga_p\ra X$ by $\phi_p(\Ga_pu)=
\phi(\Ga u)$. Then $\phi_p$ is a homeomorphism with its image, so
$(U_p,\Ga_p,\phi_p)$ is an orbifold chart on $X$. We call
$(U_p,\Ga_p,\phi_p)$ an {\it induced chart\/} of $(U,\Ga,\phi)$
for~$p$.

Let $(U,\Ga,\phi)$ and $(V,\De,\psi)$ be $n$-dimensional orbifold
charts on $X$. We say $(U,\Ga,\phi)$, $(V,\De,\psi)$ are {\it
compatible\/} if for each $p\in\phi(U/\Ga)\cap\psi(V/\De)$ there
exist induced charts $(U_p,\Ga_p,\phi_p)$ and $(V_p,\De_p,\psi_p)$
for $p$, an isomorphism of finite groups $\rho_p:\Ga_p\ra\De_p$, and
a $\rho_p$-equivariant diffeomorphism $\si_p:U_p\ra V_p$, such that
$\phi_p=\psi_p\ci(\si_p)_*$, where $(\si_p)_*:U_p/\Ga_p\ra
V_p/\De_p$ is the induced homeomorphism.

An $n$-{\it dimensional orbifold structure\/} on $X$ is a system of
$n$-dimensional orbifold charts $\bigl\{(U^i,\Ga^i,\phi^i):i\in
I\bigr\}$ on $X$ such that $X=\bigcup_{i\in I}\phi^i(U^i/\Ga^i)$ and
$(U^i,\Ga^i,\phi^i),(U^j,\Ga^j,\phi^j)$ are pairwise compatible for
all $i,j\in I$. We call $X$ with its $n$-dimensional orbifold
structure an $n$-{\it dimensional orbifold\/} (either {\it without
boundary}, or {\it with boundary}, or {\it with corners}, or {\it
with g-corners}, according to which classes of orbifold charts above
we allow). When we just refer to an orbifold, we usually mean an
orbifold with g-corners.

If $X$ is an $n$-dimensional orbifold and $p\in X$ then for some
orbifold chart $(U,\Ga,\phi)$ on $X$ we have $p=\phi(\Ga u_p)$ for
$u_p\in U$. Define the {\it stabilizer group\/} $\Stab(p)$ of $p$ to
be the finite group $\{\ga\in\Ga:\ga\,u_p=u_p\}$. Up to isomorphism,
this is independent of choices of $(U,\Ga,\phi)$ and $u_p$. The
representation of $\Stab(p)$ on $T_{u_p}U\cong\R^n$ is also
independent of choices, up to isomorphisms of $\Stab(p)$ and
automorphisms of~$\R^n$.

Generalizing Definition \ref{kh2def3}, we can define the {\it
boundary\/} $\pd X$ of an orbifold with (g-)corners, which is also
an orbifold with (g-)corners. Following Definition \ref{kh2def5} we
can define a natural involution $\si:\pd^2X\ra\pd^2X$. This $\si$
may not be free in the orbifold case, but if $X$ is oriented then
$\si$ is orientation-reversing. Principle \ref{kh2pri} holds, with
orbifolds in place of manifolds.
\label{kh2def6}
\end{dfn}

\begin{rem}{\bf(a)} In defining orbifolds, some authors require the
group $\Ga$ to act {\it effectively} on $\R^n$ in orbifold charts
$(U,\Ga,\phi)$, that is, every non-identity element in $\Ga$ must
act nontrivially on $\R^n$. Our definition is more general, and
allows orbifolds such as $\{0\}/\Ga$ in which generic points have
nontrivial stabilizer groups. We will call orbifolds in which $\Ga$
acts effectively in orbifold charts $(U,\Ga,\phi)$ {\it effective
orbifolds}. They will be discussed in \S\ref{kh39}, and will be
important in defining effective Kuranishi (co)homology.

One reason for allowing non-effective orbifolds is that some
important examples of moduli spaces are not effective. For example,
the moduli space $\oM_{1,1}$ of Deligne--Mumford stable Riemann
surfaces of genus one with one marked point, to be discussed in
\S\ref{kh61}, is a non-effective orbifold, whose generic points have
stabilizer group~$\Z_2$.

\noindent{\bf(b)} Suppose $(U,\Ga,\phi)$ is an orbifold chart on $X$
with corners or g-corners, and $u\in S_k(U)$ for $k\ge 2$ so that
there are $l\ge k$ local boundary components $B_1,\ldots,B_l$ of $U$
at $u$, and $\De\subset\Ga$ is the subgroup of $\Ga$ fixing $u$.
Then we allow $\De$ to permute $B_1,\ldots,B_l$, we do not require
$\De$ to fix each local boundary component of $U$ at $u$. Thus, the
stabilizer group of $\phi(\Ga u)$ in $X$ is $\De$, but the
stabilizer group of $\phi\bigl(\Ga(u,B_i)\bigr)$ in $\pd X$ is the
subgroup of $\De$ fixing $B_i$, for $i=1,\ldots,l$. So, under
$\io:\pd X\ra X$ mapping $(x,B)\mapsto x$, the induced morphisms of
stabilizer groups $\io_*:\Stab\bigl((x,B)\bigr)\ra \Stab(x)$ are
injective, but need not be isomorphisms.

In Definition \ref{kh5def9}, for $\Ga$ a finite group and $\rho$ an
isomorphism class of nontrivial $\Ga$-representations, we will
define the {\it orbifold stratum} $X^{\Ga,\rho}$ of an orbifold $X$.
It is an immersed suborbifold of $X$, with points $(x,\la)$ for
$x\in X$ and $\la:\Ga\ra\Stab_X(x)$ an injective group morphism. The
fact above has the  important consequence that {\it for orbifolds
with corners or g-corners, restricting to orbifold strata does not
commute with taking boundaries}, that is,~$\pd(X^{\Ga,\rho})\ne(\pd
X)^{\Ga,\rho}$.

In Chapter \ref{kh5} we will discuss {\it Kuranishi (co)bordism
groups} $KB_*,KB^*(Y;R)$. We will show that $KB_*,KB^*(Y;R)$ are
{\it very large}, as for `chains' $[X,\bs f]$ there is invariant
information in orbifold strata $X^{\Ga,\rho}$ for all $\Ga,\rho$
with $\md{\Ga}$ odd, which can be recovered from $KB_*,KB^*(Y;R)$.
However, this does not apply to Kuranishi (co)homology. Since
restricting to orbifold strata does not commute with boundaries for
orbifolds and Kuranishi spaces with (g-)corners, for Kuranishi
(co)chains $[X,\bs f,\bs G]$ or $[X,\bs f,\bs C]$, no information
from orbifold strata $X^\Ga$ survives in Kuranishi (co)homology.
This is implicit in the proof of the isomorphism of (effective)
Kuranishi homology with singular homology in Appendices
\ref{khA}--\ref{khC}, see in particular~\S\ref{khA2}.
\label{kh2rem2}
\end{rem}

\begin{dfn} Let $X,Y$ be orbifolds, and $f:X\ra Y$ a continuous
map. We call $f$ {\it smooth\/} if for all $p\in X$ and $q=f(p)$ in
$Y$ there are induced orbifold charts $(U_p,\Ga_p,\phi_p)$,
$(V_q,\De_q,\psi_q)$ on $X,Y$ for $p,q$, a morphism of finite groups
$\rho_{pq}:\Ga_p\ra\De_q$, and a $\rho_{pq}$-equivariant smooth map
$\si_{pq}:U_p\ra V_q$, such that $f\ci\phi_p(\Ga_pu)=\psi_q(\De_q
\si_{pq}(u))$ for all $u\in U_p$. ({\bf Warning:} this definition is
oversimplified and requires modification, as in Remark \ref{kh2rem3}
below.)

We call $f$ a {\it submersion\/} if all the $\si_{pq}$ are
submersions of manifolds, in the sense of Definition \ref{kh2def4}.
We call $f$ an {\it embedding\/} if $f$ is injective, all the
$\rho_{pq}$ are isomorphisms, and all the $\si_{pq}$ are embeddings.
We call $f$ a {\it diffeomorphism\/} if it is a homeomorphism, all
the $\rho_{pq}$ are isomorphisms, and all the $\si_{pq}$ are
diffeomorphisms with their images. If $f$ is a diffeomorphism then
$f^{-1}$ is a diffeomorphism.

If $X,Y,Z$ are orbifolds and $f:X\ra Y$, $g:Y\ra Z$ are smooth then
$g\ci f:X\ra Z$ is smooth. If $f$ is smooth and $p\in X$ with
$q=f(p)$ in $Y$ there is a natural group homomorphism
$f_*:\Stab(p)\ra\Stab(q)$, unique up to conjugation in
$\Stab(p),\Stab(q)$, given by $f_*=\rho_{pq}$ above under the
isomorphisms $\Stab(p)\cong\Ga_p$, $\Stab(q)\cong\De_q$. Embeddings
and diffeomorphisms induce isomorphisms $f_*$ on stabilizer groups.
\label{kh2def7}
\end{dfn}

\begin{rem} This definition of smooth maps of orbifolds is based on
Satake \cite[p.~361]{Sata}, who introduced orbifolds (then called
{\it V-manifolds\/}). However, there is a problem with it, discussed
by Adem et al.\ \cite[p.~23-4, p.~47-50]{ALR}: smooth maps of
orbifolds defined this way do not have all the properties one wants,
in particular, if $f:X\ra Y$ is a smooth map of orbifolds and $E\ra
Y$ is an orbifold vector bundle (see Definition \ref{kh2def8}), then
the pullback orbifold vector bundle $f^*(E)\ra X$ may not be
well-defined.

The problem is this: if $X,Y,f,p,q$ and some choices of orbifold
charts $(U_p,\Ga_p,\phi_p)$, $(V_q,\De_q,\psi_q)$ are as above, then
the $\rho_{pq}$-equivariant smooth map $\si_{pq}:U_p\ra V_q$
satisfying $f\ci\phi_p(\Ga_pu)=\psi_q(\De_q\si_{pq}(u))$ for $u\in
U_p$ {\it may not be uniquely defined}. For instance, if $\De_q$ is
abelian and $\de\in\De_q$ then $\si_{pq}'=\de\ci\si_{pq}$ is also a
$\rho_{pq}$-equivariant smooth map $\si_{pq}':U_p\ra V_q$ with
$f\ci\phi_p(\Ga_pu)=\psi_q(\De_q\si_{pq}'(u))$ for all $u\in U_p$,
but we may have $\si_{pq}'\not\equiv\si_{pq}$. Thus, if we write
down some construction involving the $\si_{pq}$, and the result is
not independent of the choice of $\si_{pq}$, then the construction
is not well-defined. This holds for pullbacks $f^*(E)$ of orbifold
vector bundles, and also for fibre products of orbifolds below.

There is a solution to this problem. Using their `groupoid' point of
view on orbifolds (for them, a {\it groupoid\/} is roughly speaking
equivalent to an orbifold with a particular choice of atlas of
charts $(U,\Ga,\phi)$), Moerdijk and Pronk \cite[\S 5]{MoPr} define
{\it strong maps\/} of orbifolds. Also, Chen and Ruan \cite[\S
4.4]{ChRu1} define {\it good maps} of orbifolds, also with choices
of atlases. Both strong and good maps have well-behaved pullbacks of
orbifold vector bundles, and fibre products. Lupercio and Uribe
\cite[Prop.~5.1.7]{LuUr} prove that these notions of strong maps and
good maps coincide. Adem at al.~\cite[\S 2.4]{ALR} call strong/good
maps {\it morphisms of orbifolds}. One can also use the theory of
{\it stacks\/} from algebraic geometry: as in \cite[\S 8.4]{LuUr},
an orbifold is basically a Deligne--Mumford stack over smooth
manifolds, and then strong and good maps correspond to morphisms of
stacks.

In fact the notion of smooth map of orbifolds we want to work with
is that of strong maps \cite{MoPr}, good maps \cite{ChRu1} or
orbifold morphisms \cite{ALR}, rather than the weaker notion of
Satake \cite{Sata}. But this is not a central issue for us, and to
work in one of the frameworks of \cite{ALR,ChRu1,MoPr} would
increase notation, and sacrifice clarity and geometric intuition. So
we will have our cake and eat it, by making the following
convention: strong/good maps are equivalent to Satake's smooth maps
together with some {\it extra data}, and this extra data enables us
to make unique choices of $\si_{pq}$ in Definition \ref{kh2def7}, in
a functorial way. We will leave this extra data implicit, but use it
to choose $\si_{pq}$ uniquely whenever we need to, such as in
Definition \ref{kh2def9} below. This convention entitles us to
ignore the problem of non-uniqueness of~$\si_{pq}$.
\label{kh2rem3}
\end{rem}

Our definition of orbifold vector bundles is equivalent to the {\it
good orbifold vector bundles} of Chen and Ruan~\cite[\S 4.3]{ChRu1}.

\begin{dfn} Let $X$ be an orbifold. An {\it orbifold vector
bundle\/}, or {\it orbibundle}, $E\ra X$ with {\it fibre\/} $\R^k$
is an orbifold $E$ called the {\it total space\/} and a smooth map
$\pi:E\ra X$ called the {\it projection}, such that:
\begin{itemize}
\setlength{\itemsep}{0pt}
\setlength{\parsep}{0pt}
\item[(a)] for each orbifold chart $(U,\Ga,\phi)$ on $X$ there is
an orbifold chart $(E_U,\Ga,\hat\phi)$ on $E$ with
$\hat\phi(E_U/\Ga)=\pi^{-1}(\phi(U/\Ga))$, and a $\Ga$-equivariant
submersion $\pi_U:E_U\ra U$ with $\pi\ci\hat\phi(\Ga\hat
e)=\phi(\Ga\pi_U(\hat e))$ for all $\hat e\in E_U$, such that $E_U$
has the structure of a vector bundle over $U$ with fibre $\R^k$ and
projection $\pi_U:E_U\ra U$, and the $\Ga$-action on $E_U$ preserves
the vector bundle structure.
\item[(b)] Let $(U,\Ga,\phi)$, $(V,\De,\psi)$, $p\in\phi(U/\Ga)\cap
\psi(V/\De)$, $(U_p,\Ga_p,\phi_p)$, $(V_p,\De_p,\psi_p)$,
$\rho_p:\Ga_p\ra\De_p$ and $\si_p:U_p\ra V_p$ be as in Definition
\ref{kh2def6}, and $(U,\Ga,\phi)$, $(V,\De,\psi)$ lift to
$(E_U,\Ga,\hat\phi)$, $(E_V,\De,\hat\psi)$ on $E$ with projections
$\pi_U:E_U\ra U$, $\pi_V:E_V\ra V$ as in (a). Set
$E_{U_p}=\pi_U^{-1}(U_p)$, $E_{V_p}=\pi_V^{-1}(V_p)$, so that
$E_{U_p},E_{V_p}$ are vector bundles over $U_p,V_p$ with projections
$\pi_U:E_{U_p}\ra U_p$, $\pi_V:E_{V_p}\ra V_p$. Define
$\hat\phi_p:E_{U_p}/\Ga_p\ra E$, $\hat\psi_p:E_{V_p}/\De_p\ra E$ by
$\hat\phi_p(\Ga_p\hat e)=\hat\phi(\Ga\hat e)$, $\hat\psi_p(\De_p\hat
f)=\hat\psi(\De\hat f)$, so that $(E_{U_p},\Ga_p,\hat\phi_p)$,
$(E_{V_p},\De_p,\hat\psi_p)$ are induced charts for $E$. Then there
must exist a $\rho_p$-equivariant isomorphism of vector bundles
$\hat\si_p:E_{U_p}\ra E_{V_p}$ over $U_p,V_p$, such that
$\pi_V\ci\hat\si_p=\si_p\ci\pi_U$ and~$\hat\phi_p=\hat\psi_p
\ci(\hat\si_p)_*$.
\end{itemize}

A {\it section\/} $s$ of an orbibundle $\pi:E\ra X$ is a smooth map
$s:X\ra E$ with $\pi\ci s=\id_X$, the identity on $X$. For any
$n$-dimensional orbifold $X$ we may define the {\it tangent
bundle\/} $TX$, an orbibundle with fibre $\R^n$. Morphisms (or
isomorphisms, embeddings, \ldots) of orbibundles may be defined in
the obvious way, as smooth maps (or diffeomorphisms, embeddings,
\ldots) of the total spaces of the bundles, which lift to morphisms
(or isomorphisms, embeddings, \ldots) of vector bundles $E_U\ra U$
on orbifold charts~$(E_U,\Ga,\hat\phi)$.
\label{kh2def8}
\end{dfn}

In Proposition \ref{kh2prop1} we defined {\it fibre products\/} of
manifolds. We now extend this to orbifolds. There are some
subtleties to do with stabilizer groups which mean that even as a
set, the orbifold fibre product $X\t_YX'$ is not given by
\eq{kh2eq3}. For orbifolds without boundary our definition is
equivalent to that of Moerdijk and Pronk \cite[\S 5]{MoPr}, but
because of our convention in Remark \ref{kh2rem3} we can present it
in a more explicit, differential-geometric way.

\begin{dfn} Let $X,X',Y$ be orbifolds and $f:X\ra Y$, $f':X'\ra Y$
be smooth maps, at least one of which is a submersion. Then for
$p\in X$ and $p'\in X'$ with $f(p)=q=f(p')$ in $Y$ we have
homomorphisms of stabilizer groups $f_*:\Stab(p)\ra\Stab(q)$,
$f'_*:\Stab(p')\ra\Stab(q)$. Thus we can form the double coset space
\begin{align*}
&f_*(\Stab(p))\backslash\Stab(q)/f'_*(\Stab(p'))\\
&=\bigl\{ \{f_*(\ga)\de f_*(\ga'):\ga\in\Stab(p),\;\>
\ga'\in\Stab(p')\}:\de\in\Stab(q)\bigr\}.
\end{align*}
As a set, we define
\e
\begin{split}
X\t_YX'=\bigl\{(p,p',\De):\,&\text{$p\in X$, $p'\in X'$,
$f(p)=f'(p')$,}\\
&\De\in f_*(\Stab(p))\backslash\Stab(f(p))/f'_*(\Stab(p'))\bigr\}.
\end{split}
\label{kh2eq10}
\e

We give this the structure of an orbifold as follows. Suppose $p\in
X$, $p'\in X'$ and $q\in Y$ with $f(p)=y=f'(p')$, and let
$(U_p,\Ga_p,\phi_p)$, $(U'_{p'},\Ga'_{p'},\phi'_{p'})$,
$(V_q,\De_q,\psi_q)$ be induced orbifold charts, so that
$\Ga_p\cong\Stab(p)$, $\Ga'_{p'}\cong\Stab(p')$,
$\De_q\cong\Stab(q)$. Making $U_p,U'_{p'}$ smaller if necessary,
there exist homomorphisms $\rho_{pq}:\Ga_p\ra\De_q$ and
$\rho'_{p'q}:\Ga'_{p'}\ra\De_q$ identified with
$f_*:\Stab(p)\ra\Stab(q)$ and $f'_*:\Stab(p') \ra\Stab(q)$, and
$\rho_{pq},\rho'_{p'q}$-equivariant smooth maps $\si_{pq}:U_p\ra
V_q$ and $\si'_{p'q}:U'_{p'}\ra V_q$ at least one of which is a
submersion, inducing $f,f'$ on $\Im\phi_p$ and $\Im\phi'_{p'}$.
Define an orbifold chart on $X\t_YX'$ by
\e
\bigl((U_p\t U_{p'})\t_{\si_{pq}\t\si'_{p'q},V_p\t
V_p,\pi_q}(V_q\t\De_q),\Ga_p\t\Ga'_{p'},\psi_{pp'q}\bigr).
\label{kh2eq11}
\e
Here the first term in \eq{kh2eq11} is a fibre product of smooth
manifolds, as in Proposition \ref{kh2prop1}. Define
$\pi_q:(V_q\t\De_q)\ra V_q\t V_q$ by $\pi_q(v,\de)\mapsto(v,\de\cdot
v)$. The action of $\Ga_p\t\Ga'_{p'}$ on the first term in
\eq{kh2eq11} is induced by the obvious action on $U_p\t U_{p'}$, and
the action on $V_q\t\De_q$
by~$(\ga,\ga'):(v,\de)\mapsto\bigl(\rho_{pq}(\ga)\cdot
v,\rho'_{p'q}(\ga')\de\rho_{pq}(\ga)^{-1}\bigr)$.

Defining the map $\psi_{pp'q}$ from the first term of \eq{kh2eq11}
to the set \eq{kh2eq10} is tricky, because
$\Stab(p),\Stab(p'),\Stab(q)$ and the morphisms $f_*,f'_*$ between
them are only really defined up to conjugation. (Actually, one needs
the ideas of Remark \ref{kh2rem3} at this point.) Let $(u,u',v,\de)$
lie in the first term of \eq{kh2eq11}, with $\phi_p(\Ga_pu)=x$,
$\phi'_{p'}(\Ga'_{p'}u')=x'$ and $\psi_q(\De_qv)=y$. Then we have
natural identifications $\Stab(x)\cong\Stab_{\Ga_p}(u)$,
$\Stab(x')\cong\Stab_{\Ga_{p'}}(u')$ and
$\Stab(y)\cong\Stab_{\De_q}(v)$. But $f_*:\Stab(x)\ra\Stab(y)$ is
naturally identified with $\rho'_{p'q}:\Stab_{\Ga_{p'}}(u')\ra
\Stab_{\De_q}(\de\cdot v)$. To identify $\Stab_{\De_q}(\de\cdot v)$
with our model $\Stab_{\De_q}(v)$ for $\Stab(y)$ we must choose
$\ep\in\De_q$ such that $\ep\cdot v=\de\cdot v$. (Note that just
fixing $\ep=\de$ is not a good choice; we want the choice of $\ep$
to depend only on the points $v,\de\cdot v$, not on $\de$.) With
these identifications, we can define $\psi_{pp'q}$ by
\begin{align*}
&\psi_{pp'q}:(u,u',v,\de)\longmapsto\\
&\bigl(\psi_p(\Ga_pu),\psi'_{p'}(\Ga'_{p'}u'),\rho_{pq}
(\Stab_{\Ga_p}(u))\de\ep^{-1}\rho'_{p'q}(\Stab_{\Ga'_{p'}}(u'))\bigr).
\end{align*}

The orbifold charts \eq{kh2eq11} are compatible, and make $X\t_YX'$
into an orbifold. The projections $\pi_X:X\t_YX'\ra X$,
$\pi_{X'}:X\t_YX'\ra X'$ and $\pi_Y:X\t_YX'\ra Y$  are smooth maps,
and the analogue of Proposition \ref{kh2prop1} holds. If $Y$ is a
manifold, or more generally if all the maps
$f_*:\Stab(x)\ra\Stab(f(x))$ or $f'_*:\Stab(x')\ra\Stab(f'(x'))$ are
surjective, then as a topological space $X\t_YX'$ is the topological
space fibre product of $X,X'$ over $Y$ in \eq{kh2eq3}. But in
general the orbifold and topological space fibre products do not
coincide.
\label{kh2def9}
\end{dfn}

\subsection{Kuranishi neighbourhoods and coordinate changes}
\label{kh23}

Let $X$ be a paracompact Hausdorff topological space.

\begin{dfn} A {\it Kuranishi neighbourhood\/}
$(V_p,E_p,s_p,\psi_p)$ of $p\in X$ satisfies:
\begin{itemize}
\setlength{\itemsep}{0pt}
\setlength{\parsep}{0pt}
\item[(i)] $V_p$ is an orbifold, which may have boundary, corners,
or g-corners;
\item[(ii)] $E_p\ra V_p$ is an orbifold vector bundle over $V_p$;
\item[(iii)] $s_p:V_p\ra E_p$ is a smooth section; and
\item[(iv)] $\psi_p$ is a homeomorphism from $s_p^{-1}(0)$ to
an open neighbourhood of $p$ in $X$, where $s_p^{-1}(0)$ is the
subset of $V_p$ where the section $s_p$ is zero.
\end{itemize}
We call $E_p$ the {\it obstruction bundle}, and $s_p$ the {\it
Kuranishi map}.
\label{kh2def10}
\end{dfn}

\begin{rem} Definition \ref{kh2def10} follows \cite[Def.~5.1]{FuOn1},
save that in (i) we allow the $V_p$ to have boundary or (g-)corners
which Fukaya and Ono do not, though in \cite[Def.~17.3]{FuOn1} they
allow the $V_p$ to have boundary, and in (iii) following \cite{FOOO}
we require $s_p$ to be {\it smooth}, whereas Fukaya and Ono allow
the $s_p$ to be only {\it continuous}. We discuss smoothness of
$s_p$ in Remark \ref{kh2rem5} below.

In \cite[Def.~A1.1]{FOOO} Fukaya et al.\ define a Kuranishi
neighbourhood to be $(V_p,E_p,\Ga_p,\psi_p,s_p)$, where $V_p$ is a
manifold possibly with boundary or corners, $E_p\ra V_p$ a trivial
vector bundle, $\Ga_p$ a finite group acting effectively on $V_p$
and $E_p$, $s_p$ a smooth $\Ga_p$-equivariant section of $E_p$, and
$\psi_p:s_p^{-1}(0)/\Ga_p\ra X$ a homeomorphism with its image. Then
$(V_p/\Ga_p,E_p/\Ga_p,s_p/\Ga_p,\psi_p)$ is a Kuranishi
neighbourhood in our sense. One reason why we opted for the
$(V_p,E_p,s_p,\psi_p)$ form in this book is that much of Chapter
\ref{kh3} would not work for Kuranishi neighbourhoods of the
form~$(V_p,E_p,\Ga_p,\psi_p,s_p)$.
\label{kh2rem4}
\end{rem}

\begin{dfn} Let $(V_p,E_p,s_p,\psi_p)$, $(\ti V_p,\ti
E_p,\ti s_p,\ti\psi_p)$ be two Kuranishi neighbourhoods of $p\in X$.
We call $(V_p,\ldots,\psi_p)$ and $(\ti V_p,\ldots,\ti\psi_p)$ {\it
isomorphic\/} with {\it isomorphism\/} $(\al,\hat\al)$ if there
exist a diffeomorphism $\al:V_p\ra\ti V_p$ and an isomorphism of
orbibundles $\hat\al:E_p\ra\al^*(\ti E_p)$, such that $\ti
s_p\ci\al\equiv\hat\al\ci s_p$ and~$\ti\psi_p\ci\al\equiv\psi_p$.

We call $(V_p,\ldots,\psi_p),(\ti V_p,\ldots,\ti\psi_p)$ {\it
equivalent\/} if there exist open neighbourhoods $U_p\!\subseteq\!
V_p$, $\ti U_p\!\subseteq\!\ti V_p$ of $\psi_p^{-1}(p),
\ti\psi_p^{-1}(p)$ such that
$(U_p,E_p\vert_{U_p},s_p\vert_{U_p},\psi_p\vert_{U_p})$ and $(\ti
U_p,\ti E_p\vert_{\ti U_p},\ti s_p\vert_{\ti U_p},\ti\psi_p
\vert_{\ti U_p})$ are isomorphic.
\label{kh2def11}
\end{dfn}

\begin{dfn} Let $(V_p,E_p,s_p,\psi_p)$ and $(V_q,E_q,s_q,
\psi_q)$ be Kuranishi neighbourhoods of $p\in X$ and
$q\in\psi_p(s_p^{-1}(0))$ respectively. We call a pair
$(\phi_{pq},\hat\phi_{pq})$ a {\it coordinate change\/} from
$(V_q,\ldots,\psi_q)$ to $(V_p,\ldots,\psi_p)$ if:
\begin{itemize}
\setlength{\itemsep}{0pt}
\setlength{\parsep}{0pt}
\item[(a)] $\phi_{pq}:V_q\ra V_p$ is a smooth embedding of
orbifolds. When $V_p,V_q$ have boundaries or (g-)corners,
$\phi_{pq}$ maps $S_k(V_q)\ra S_k(V_p)$ in the notation of
Definition \ref{kh2def3}, for all $k\ge 0$. Thus $\phi_{pq}$ is
compatible with boundaries and (g-)corners, and induces embeddings
$\phi_{pq}:\pd^kV_q\ra\pd^kV_p$ for all~$k\ge 0$.
\item[(b)] $\hat\phi_{pq}:E_q\ra\phi_{pq}^*(E_p)$ is an embedding of
orbibundles over~$V_q$;
\item[(c)] $\hat\phi_{pq}\ci s_q\equiv s_p\ci\phi_{pq}$;
\item[(d)] $\psi_q\equiv \psi_p\ci\phi_{pq}$; and
\item[(e)] Choose an open neighbourhood $W_{pq}$ of
$\phi_{pq}(V_q)$ in $V_p$, and an orbifold vector subbundle $F_{pq}$
of $E_p\vert_{W_{pq}}$ with $\phi_{pq}^*(F_{pq})=\hat\phi_{pq}
(E_q)$, as orbifold vector subbundles of $\phi_{pq}^*(E_p)$ over
$V_q$. Write $\hat s_p:W_{pq}\ra E_p/F_{pq}$ for the projection of
$s_p\vert_{W_{pq}}$ to the quotient bundle $E_p/F_{pq}$. Now
$s_p\vert_{\phi_{pq}(V_q)}$ lies in $F_{pq}$ by (c), so $\hat
s_p\vert_{\phi_{pq}(V_q)}\equiv 0$. Thus there is a well-defined
derivative
\begin{equation*}
\d\hat s_p:N_{\phi_{pq}(V_q)}V_p\ra
(E_p/F_{pq})\vert_{\phi_{pq}(V_q)},
\end{equation*}
where $N_{\phi_{pq}(V_q)}V_p$ is the normal orbifold vector bundle
of $\phi_{pq}(V_q)$ in $V_p$. Pulling back to $V_q$ using
$\phi_{pq}$, and noting that $\phi_{pq}^*(F_{pq})=\hat\phi_{pq}
(E_q)$, gives a morphism of orbifold vector bundles over~$V_q$:
\e
\d\hat s_p:\frac{\phi_{pq}^*(TV_p)}{(\d\phi_{pq})(TV_q)}\longra
\frac{\phi_{pq}^*(E_p)}{\hat\phi_{pq}(E_q)}\,.
\label{kh2eq12}
\e

We require that \eq{kh2eq12} should be an {\it isomorphism\/} over
$s_q^{-1}(0)$, and hence on an open neighbourhood of $s_q^{-1}(0)$
in $V_q$. Note that this forces $\dim V_p-\rank E_p=\dim V_q-\rank
E_q$. Now since $s_q\equiv 0$ on $s_q^{-1}(0)$, we have $s_p\equiv
0$ on $\phi_{pq}(s_q^{-1}(0))$ by (c), and because of this, the
restriction of $\d\hat s_p$ in \eq{kh2eq12} to $s_q^{-1}(0)$ is {\it
independent\/} of the choices of $W_{pq}$ and $F_{pq}$. So the
condition that \eq{kh2eq12} is an isomorphism over $s_q^{-1}(0)$ is
independent of choices.
\end{itemize}
\label{kh2def12}
\end{dfn}

\begin{rem} Definition \ref{kh2def12} is based on Fukaya, Ono et
al.\ \cite[Def.~5.3]{FuOn1}, \cite[Def.~A1.3]{FOOO}, with two
important modifications. Firstly, the condition in (a) on
$\phi_{pq}$ mapping $S_k(V_q)\ra S_k(V_p)$ is new; it is necessary
for the notion of a boundary or codimension $k$ corner point to be
preserved by coordinate changes, so that boundaries of Kuranishi
spaces are well-behaved. Secondly, (e) is also new, and it replaces
the notion in \cite[Def.~5.6]{FuOn1}, \cite[Def.~A1.14]{FOOO} that a
Kuranishi structure {\it has a tangent bundle}.

There is a possibility of confusion in the literature here. As well
as the 1999 paper of Fukaya and Ono \cite{FuOn1}, there are
currently three different versions of Fukaya et al.\ \cite{FOOO} in
circulation, dated 2000, 2006 and 2008. Each of these four uses
different definitions of Kuranishi space. In this book, references
are to the 2008 version of \cite{FOOO}. In \cite{FuOn1} and the 2000
and 2006 versions of \cite{FOOO}, by saying that a Kuranishi
structure `has a tangent bundle', Fukaya et al.\ mean that for each
coordinate change $(\phi_{pq},\hat\phi_{pq})$ in the Kuranishi
structure, we are {\it given\/} an isomorphism of orbibundles over
$V_q$
\e
\chi_{pq}:\frac{\phi_{pq}^*(TV_p)}{(\d\phi_{pq})(TV_q)}\longra
\frac{\phi_{pq}^*(E_p)}{\hat\phi_{pq}(E_q)}\,,
\label{kh2eq13}
\e
satisfying compatibilities under composition of coordinate changes.
But these $\chi_{pq}$ need have no particular relation to $\d\hat
s_p$ in~\eq{kh2eq12}.

In the 2008 version of \cite{FOOO}, the definition
\cite[Def.~A1.14]{FOOO} that a Kuranishi structure `has a tangent
bundle' is basically equivalent to Definition \ref{kh2def12}(e).
That is, Fukaya et al.\ identify the maps $\chi_{pq}$ in
\eq{kh2eq13} with $\d\hat s_p$, which they write as $d_{\rm
fiber}s_p$. This change is necessary, and important. With the
definitions of `Kuranishi space with a tangent bundle' used in
\cite{FuOn1} and the 2000, 2006 versions of \cite{FOOO}, in which
$\chi_{pq}$ is unrelated to $\d\hat s_p$, the construction of
virtual cycles and virtual classes is invalid. This is remarked upon
in the 2008 version \cite[A1.7(iv)]{FOOO}, which also gives an
example \cite[Ex.~A1.64]{FOOO} of bad behaviour using the previous
definition of tangent bundles.

Our Definition \ref{kh2def12}(e), and the corresponding definition
\cite[Def.~A1.14]{FOOO} in the 2008 version of \cite{FOOO}, require
the Kuranishi maps $s_p$ to be {\it smooth}, as in
\cite[Def.~A1.1]{FOOO}, or at least $C^1$, rather than merely
continuous as in \cite[Def.~5.1]{FuOn1}. This is a significant issue
because in constructing Kuranishi structures on moduli spaces of
$J$-holomorphic curves, it is straightforward to make the $s_p$
continuous at singular curves, but one has to do much more work to
arrange for the $s_p$ to be smooth. We will return to this issue in
Remark \ref{kh6rem1} and in~\S\ref{kh68}.
\label{kh2rem5}
\end{rem}

\subsection{Kuranishi structures}
\label{kh24}

We can now define {\it Kuranishi structures\/} on topological
spaces.

\begin{dfn} Let $X$ be a paracompact Hausdorff topological space. A
{\it germ of Kuranishi neighbourhoods at\/} $p\in X$ is an
equivalence class of Kuranishi neighbourhoods $(V_p,E_p,s_p,\psi_p)$
of $p$, for equivalence as in Definition~\ref{kh2def11}.

Suppose $(V_p,E_p,s_p,\psi_p)$ lies in such a germ. Then for any
open neighbourhood $U_p$ of $\psi_p^{-1}(p)$ in $V_p$,
$(U_p,E_p\vert_{U_p},s_p\vert_{U_p},\psi_p\vert_{U_p})$ also lies in
the germ. As a shorthand, we say that some condition on the germ
{\it holds for sufficiently small\/} $(V_p,\ldots,\psi_p)$ if
whenever $(V_p,\ldots,\psi_p)$ lies in the germ, the condition holds
for $(U_p,\ldots,\psi_p\vert_{U_p})$ for all sufficiently small
$U_p$ as above.

A {\it Kuranishi structure\/} $\ka$ on $X$ assigns a germ of
Kuranishi neighbourhoods for each $p\in X$ and a {\it germ of
coordinate changes\/} between them in the following sense: for each
$p\in X$, for all sufficiently small $(V_p,\ldots,\psi_p)$ in the
germ at $p$, for all $q\in\Im\psi_p$, and for all sufficiently small
$(V_q,\ldots,\psi_q)$ in the germ at $q$, we are given a coordinate
change $(\phi_{pq},\hat\phi_{pq})$ from $(V_q,\ldots,\psi_q)$ to
$(V_p,\ldots,\psi_p)$. These coordinate changes should be compatible
with equivalence in the germs at $p,q$ in the obvious way, and
satisfy:
\begin{itemize}
\setlength{\itemsep}{0pt}
\setlength{\parsep}{0pt}
\item[(i)] $\dim V_p-\rank E_p$ is independent of $p$\/; and
\item[(ii)] if $q\in\Im\psi_p$ and $r\in\Im\psi_q$ then
$\phi_{pq}\ci\phi_{qr}=\phi_{pr}$ and~$\hat\phi_{pq}\ci
\hat\phi_{qr}=\hat\phi_{pr}$.
\end{itemize}
We call $\vdim X=\dim V_p-\rank E_p$ the {\it virtual dimension\/}
of the Kuranishi structure. Note that Definition \ref{kh2def12}(e)
already implies that $\dim V_p-\rank E_p$ is invariant under
coordinate changes, so it is locally constant on $X$; part (i) says
it is globally constant.

A {\it Kuranishi space\/} $(X,\ka)$ is a topological space $X$ with
a Kuranishi structure $\ka$. Usually we will just refer to $X$ as
the Kuranishi space, suppressing $\ka$, except when we need to
distinguish between different Kuranishi structures $\ka,\ti\ka$ on
the same topological space~$X$.

We say that a Kuranishi space $(X,\ka)$ is {\it without boundary},
or {\it with boundary}, or {\it with corners}, or {\it with
g-corners}, if for all $p\in X$ and all sufficiently small
$(V_p,\ldots,\psi_p)$ in the germ at $p$ in $\ka$, $V_p$ is an
orbifold without boundary, or with boundary, or with corners, or
with g-corners, respectively. When we just refer to a Kuranishi
space, we usually mean a Kuranishi space with g-corners.
\label{kh2def13}
\end{dfn}

\begin{rem} Definition \ref{kh2def13} is based on
\cite[Def.~5.3]{FuOn1}, \cite[Def.~A1.5]{FOOO}, but since our
definition of coordinate change includes the extra condition
Definition \ref{kh2def12}(e), as in Remark \ref{kh2rem5}, our notion
of Kuranishi structure is {\it stronger\/} than that of Fukaya et
al., and roughly corresponds to their definition of Kuranishi
structure {\it with a tangent bundle} \cite[Def.~5.6]{FuOn1},
\cite[Def.~A1.14]{FOOO}. Note too that \cite{FuOn1} and the 2000 and
2006 versions of \cite{FOOO} define Kuranishi structures in terms of
{\it germs\/} of Kuranishi neighbourhoods, as above, but the 2008
version of \cite{FOOO} does not use germs.
\label{kh2rem6}
\end{rem}

An {\it orbifold\/} $X$ is equivalent to a Kuranishi space with
$E_p$ the zero vector bundle and $s_p\equiv 0$ in all Kuranishi
neighbourhoods $(V_p,E_p,s_p,\psi_p)$, and can be covered by the
single Kuranishi neighbourhood~$(X,X,0,\id_X)$.

\begin{dfn} Let $(X,\ka)$ be a Kuranishi space. A Kuranishi
neighbourhood $(V,E,s,\psi)$ on $X$ is called {\it compatible with
the Kuranishi structure\/} $\ka$ if for all $p\in\Im\psi$ and for
all sufficiently small Kuranishi neighbourhoods
$(V_p,\ldots,\psi_p)$ in the germ at $p$ in $\ka$, there exists a
coordinate change $\smash{(\phi_p,\hat\phi_p)}$ from
$(V_p,\ldots,\psi_p)$ to~$(V,E,s,\psi)$.

Actually, when we use compatible Kuranishi neighbourhoods, we shall
usually implicitly suppose that for all sufficiently small
$(V_p,\ldots,\psi_p)$ we are given a {\it particular choice\/} of
coordinate change $(\phi_p,\hat\phi_p)$ from $(V_p,\ldots,\psi_p)$
to $(V,E,s,\psi)$, and that these choices are compatible with
coordinate changes $(\phi_{pq},\hat\phi_{pq})$ in the germ of
coordinate changes in $\ka$ in the obvious way. So probably these
choices of $(\phi_p,\hat\phi_p)$ should be part of the definition.

Definition \ref{kh2def13} implies that for $p\in X$, all
sufficiently small $(V_p,\ldots,\psi_p)$ in the germ at $p$ in $\ka$
are compatible with $\ka$. Thus, every Kuranishi space can be
covered by compatible Kuranishi neighbourhoods (with particular
choices of coordinate changes $(\phi_p,\hat\phi_p)$ as above). If
$X$ is compact we can choose a finite cover by compatible Kuranishi
neighbourhoods.
\label{kh2def14}
\end{dfn}

We can also define {\it stabilizer groups\/} of points in a
Kuranishi space.

\begin{dfn} Let $X$ be a Kuranishi space, $p\in X$, and
$(V_p,\ldots,\psi_p)$ be a Kuranishi neighbourhood in the germ at
$p$. Define the {\it stabilizer group} $\Stab(p)$ of $p$ to be the
stabilizer group $\Stab(\psi_p^{-1}(p))$ of the point
$\psi_p^{-1}(p)$ in the orbifold $V_p$. By definition of equivalence
in Definition \ref{kh2def11}, up to isomorphism this is independent
of the choice of $(V_p,\ldots,\psi_p)$, and so is well-defined.
\label{kh2def15}
\end{dfn}

\subsection{Strongly smooth maps and strong submersions}
\label{kh25}

The next definition is based on Fukaya and Ono
\cite[Def.~6.6]{FuOn1}. The equivalent definition in
\cite[Def.~A1.13]{FOOO} instead uses {\it good coordinate systems},
as in \S\ref{kh31}, but we prefer not to choose a good coordinate
system unnecessarily.

\begin{dfn} Let $X$ be a Kuranishi space, and $Y$ a topological
space. A {\it strongly continuous map\/} $\bs f:X\ra Y$ consists of,
for all $p\in X$ and all sufficiently small $(V_p,E_p,s_p,\psi_p)$
in the germ of Kuranishi neighbourhoods at $p$, a choice of
continuous map $f_p:V_p\ra Y$, such that for all $q\in\Im\psi_p$ and
sufficiently small $(V_q,\ldots,\psi_q)$ in the germ at $q$ with
coordinate change $(\phi_{pq},\hat\phi_{pq})$ from
$(V_q,\ldots,\psi_q)$ to $(V_p,\ldots,\psi_p)$ in the germ of
coordinate changes, we have $f_p\ci\phi_{pq}=f_q$. Then $\bs f$
induces a continuous map $f:X\ra Y$ in the obvious way.

If $Y$ is an orbifold and all $f_p$ are smooth, we call $\bs f$ {\it
strongly smooth}. We call $\bs f$ a {\it strong submersion\/} if all
the $f_p$ are submersions, in the sense of Definitions \ref{kh2def4}
and \ref{kh2def7} which includes conditions on $f_p$ over $\pd^kV_p$
and $\pd^lY$ for $k,l\ge 0$. In fact, in this book we will usually
restrict to the case $\pd Y=\es$ for simplicity, except in
\S\ref{kh45}. Note that Fukaya and Ono \cite{FuOn1,FOOO} use the
notation {\it weakly submersive}, rather than strong submersion.

We call two strongly continuous or strongly smooth maps $\bs f:X\ra
Y$, $\bs f':X\ra Y$ {\it equal\/} if for all $p\in X$ and all
sufficiently small $(V_p,\ldots,\psi_p)$ in the germ at $p$, the
maps $f_p:V_p\ra Y$ and $f_p':V_p\ra Y$ are the same. Note that this
allows there to exist $(V_p,\ldots,\psi_p)$ in the germ at $p$ such
that only one of $f_p:V_p\ra Y$ and $f_p':V_p\ra Y$ are given, or
both are given but are not equal; if both are given, the definition
of `sufficiently small' requires only that $f_p\vert_{U_p}\equiv
f'_p\vert_{U_p}$ for some open $\psi_p^{-1}(p)\in U_p\subseteq V_p$.
Equal $\bs f:X\ra Y$, $\bs f':X\ra Y$ induce the same continuous
maps $f:X\ra Y$, $f':X\ra Y$.
\label{kh2def16}
\end{dfn}

We shall need to generalize this to the case when $Y$ is a Kuranishi
space. As far as the author knows, the next two definitions are new.

\begin{dfn} Let $(V_p,E_p,s_p,\psi_p),(W_r,F_r,t_r,\xi_r)$ be
Kuranishi neighbourhoods in Kuranishi spaces $X,Y$. A {\it smooth
map\/} from $(V_p,\ldots,\psi_p)$ to $(W_r,\ldots,\xi_r)$ is a pair
$(f_{rp},\hat f_{rp})$ such that
\begin{itemize}
\setlength{\itemsep}{0pt}
\setlength{\parsep}{0pt}
\item[(a)] $f_{rp}:V_p\ra W_r$ is a smooth map of orbifolds;
\item[(b)] $\hat f_{rp}:E_p\ra f_{rp}^*(F_r)$ is a morphism of
orbibundles over $V_p$; and
\item[(c)] $\hat f_{rp}\ci s_p\equiv t_r\ci f_{rp}$.
\end{itemize}

We call $(f_{rp},\hat f_{rp})$ a {\it submersion\/} if $f_{rp}$ is a
submersion in the sense of Definitions \ref{kh2def4} and
\ref{kh2def7}, and $\hat f_{rp}:E_p\ra f_{rp}^*(F_r)$ is surjective.
Following Definition \ref{kh2def11}, we call $(f_{rp},\hat f_{rp})$
an {\it equivalence\/} if there exist open neighbourhoods $V_p'$ of
$\psi_p^{-1}(p)$ in $V_p$ and $W_r'$ of $\xi_r^{-1}(r)$ in $W_r$
such that $f_{rp}\vert_{V_p'}:V'_p \ra W'_r$ is a diffeomorphism,
and $\hat f_{rp}\vert_{V_p'}:E_p\vert_{V_p'}\ra
f_{rp}^*(F_r\vert_{W_r'})$ is an isomorphism of orbibundles.

A {\it strongly smooth map\/} $\bs f:X\ra Y$ inducing a given
continuous map $f:X\ra Y$ consists of, for all $p$ in $X$ with
$r=f(p)$ in $Y$, for all sufficiently small $(W_r,F_r,t_r,\xi_r)$ in
the germ of Kuranishi neighbourhoods at $r$, for all sufficiently
small $(V_p,E_p,s_p,\psi_p)$ in the germ at $p$, a smooth map
$(f_{rp},\hat f_{rp})$ from $(V_p,\ldots,\psi_p)$ to
$(W_r,\ldots,\xi_r)$, such that if $(\phi_{pq},\hat\phi_{pq})$ lies
in the germ of coordinate changes on $X$ then $f_{rp}\ci\phi_{pq}
=f_{rq}$ and $\hat f_{rp}\ci\hat\phi_{pq}=\hat f_{rq}$, and if
$(\phi_{rs},\hat\phi_{rs})$ lies in the germ of coordinate changes
on $Y$ then $\phi_{rs}\ci f_{sp}=f_{rp}$ and $\hat\phi_{rs}\ci\hat
f_{sp}=\hat f_{rp}$.

We call $\bs f$ a {\it strong submersion\/} if all the
$(f_{rp},\smash{\hat f_{rp}})$ are submersions. We call $\bs f$ a
{\it strong diffeomorphism\/} if $f:X\ra Y$ is a homeomorphism and
all the $(f_{rp},\hat f_{rp})$ are equivalences. We call two
strongly maps $\bs f:X\ra Y$, $\bs f':X\ra Y$ {\it equal\/} if they
induce the same continuous maps $f:X\ra Y$, $f':X\ra Y$, and for all
$p\in X$ and $r=f(p)$ in $Y$ and all sufficiently small
$(V_p,\ldots,\psi_p)$, $(W_r,\ldots,\xi_r)$ in the germs at $p,r$ in
$X,Y$, we have~$(f_{rp},\hat f_{rp})=(f'_{rp},\hat f'_{rp})$.
\label{kh2def17}
\end{dfn}

\begin{dfn} Let $X,Y,Z$ be Kuranishi spaces, and $\bs f:X\ra Y$,
$\bs g:Y\ra Z$ be strongly smooth maps inducing continuous maps
$f:X\ra Y$, $g:Y\ra Z$. Define the {\it composition\/} $\bs g\ci\bs
f:X\ra Z$ as follows: let $p\in X$, $q=f(p)\in Y$ and $r=g(q)\in Z$,
let $(V_p,\ldots,\psi_p)$, $(V'_q,\ldots,\psi'_q)$,
$(V''_r,\ldots,\psi''_r)$ be sufficiently small Kuranishi
neighbourhoods in the germs at $p,q,r$ in $X,Y,Z$, let $(f_{qp},\hat
f_{qp})$ be the smooth map from $(V_p,\ldots,\psi_p)$ to
$(V'_q,\ldots,\psi'_q)$ in the germ of $\bs f$, and $(g_{rq},\hat
g_{rq})$ the smooth map from $(V'_q,\ldots,\psi'_q)$ to
$(V''_r,\ldots,\psi''_r)$ in the germ of $\bs g$. Then the
composition $(g_{rq}\ci f_{qp},\hat g_{rq}\ci\hat f_{qp})$ is a
smooth map from $(V_p,\ldots,\psi_p)$ to $(V''_r,\ldots,\psi''_r)$.
It is easy to verify that the germ of all such compositions forms a
strongly smooth map $\bs g\ci\bs f:X\ra Z$, which induces the
continuous map $g\ci f:X\ra Z$. Clearly, if $\bs f,\bs g$ are strong
submersions then so is $\bs g\ci\bs f$, and composition is
associative,~$(\bs h\ci\bs g)\ci\bs f=\bs h\ci(\bs g\ci\bs f)$.

If $\bs f:X\ra Y$ is a strong diffeomorphism inducing the
homeomorphism $f:X\ra Y$, it is not difficult to show that there
exists a strong diffeomorphism $\bs f{}^{-1}:Y\ra X$ inducing
$f^{-1}:Y\ra X$, such that the composition $\bs f{}^{-1}\ci\bs
f:X\ra X$ equals the identity on $X$, and $\bs f\ci\bs f{}^{-1}:Y\ra
Y$ equals the identity on $Y$, using the notion of {\it equal\/}
strongly smooth maps in Definition~\ref{kh2def17}.
\label{kh2def18}
\end{dfn}

\subsection{Boundaries and fibre products}
\label{kh26}

We now define the {\it boundary\/} $\pd X$ of a Kuranishi space $X$,
which is itself a Kuranishi space of dimension $\vdim X-1$. Our
definition extends \cite[Lem.~17.8]{FuOn1} to the case of Kuranishi
spaces with (g-)corners, and is modelled on
Definition~\ref{kh2def3}.

\begin{dfn} Let $X$ be a Kuranishi space with g-corners. We shall
define a Kuranishi space $\pd X$ called the {\it boundary\/} of $X$.
The points of $\pd X$ are equivalence classes
$[p,(V_p,\ldots,\psi_p),B]$ of triples $(p,(V_p,\ldots,\psi_p),B)$,
where $p\in X$, $(V_p,\ldots,\psi_p)$ lies in the germ of Kuranishi
neighbourhoods at $p$, and $B$ is a local boundary component of
$V_p$ at $\psi_p^{-1}(p)$, in the sense of Definition~\ref{kh2def3}.

Two triples $(p,(V_p,\ldots,\psi_p),B),(q,(\ti
V_q,\ldots,\ti\psi_q),\ti B)$ are {\it equivalent\/} if $p=q$, and
the Kuranishi neighbourhoods $(V_p,\ldots,\psi_p),(\ti
V_q,\ldots,\ti\psi_q)$ are equivalent so that we are given an
isomorphism $(\al,\hat\al):(U_p,\ldots,\psi_p\vert_{U_p})\ra(\ti
U_q,\ldots,\ti\psi_q\vert_{\ti U_q})$ for open $\psi_p^{-1}(p)\in
U_p\subseteq V_p$ and $\ti\psi_q^{-1}(q)\in\ti U_q\subseteq\ti V_q$,
and $\al_*(B)=\ti B$ near $\ti\psi_q^{-1}(q)$.

We can define a unique natural {\it topology\/} and {\it Kuranishi
structure\/} on $\pd X$, such that $(\pd V_p,E_p\vert_{\pd
V_p},s_p\vert_{\pd V_p},\psi'_p)$ is a Kuranishi neighbourhood on
$\pd X$ for each Kuranishi neighbourhood $(V_p,\ldots,\psi_p)$ on
$X$, where $\psi'_p:(s_p\vert_{\pd V_p})^{-1}(0)\ra\pd X$ is given
by $\psi'_p:(q,B)\mapsto[\psi_p(q),(V_p,\ldots\psi_p),B]$ for
$(q,B)\in\pd V_p$ with $s_p(q)=0$. Then $\vdim(\pd X)=\vdim X-1$,
and $\pd X$ is compact if $X$ is compact.

Note that since $\pd V_p$ is not actually a subset of $V_p$, but a
set of points $(v,B)$ for $v\in V_p$, writing $E_p\vert_{\pd V_p}$
and $s_p\vert_{\pd V_p}$ here is an abuse of notation. A more
accurate way to say it is that we have an immersion $\io:\pd V_p\ra
V_p$ acting by $\io:(v,B)\mapsto v$, and by $E_p\vert_{\pd V_p}$ we
mean $\io^*(E_p)$, and by $s_p\vert_{\pd V_p}$ we mean~$s_p\ci\io$.

If $X$ is a manifold with g-corners then Definition \ref{kh2def5}
gave a natural smooth free involution $\si:\pd^2X\ra\pd^2X$, and
Definition \ref{kh2def6} noted this also works for orbifolds, though
$\si$ may not be free. In the same way, if $X$ is a Kuranishi space
there is a natural strong diffeomorphism $\bs\si:\pd^2X\ra\pd^2X$
with $\bs\si^2=\bs\id_X$, such that if $(V_p,E_p,s_p,\psi_p)$ is a
Kuranishi neighbourhood on $X$ and $(\pd^2V_p,\ldots,\psi_p
\vert_{\pd^2V_p})$ is the induced Kuranishi neighbourhood on
$\pd^2X$, then $\bs\si$ is represented by $(\si_p,\hat\si_p):
(\pd^2V_p,\ldots,\psi_p \vert_{\pd^2V_p})\ra(\pd^2V_p,\ldots,\psi_p
\vert_{\pd^2V_p})$, where $\si_p:\pd^2V_p\ra\pd^2V_p$ is the
involution for $V_p$, and $\hat\si_p$ its natural lift to
$E_p\vert_{\pd^2V_p}=\pd^2E_p$. If $X$ is oriented in the sense of
\S\ref{kh27} below, then $\bs\si:\pd^2X\ra\pd^2X$ is
orientation-reversing.
\label{kh2def19}
\end{dfn}

Fukaya et al.\ \cite[Def.~A1.37]{FOOO} define {\it fibre products\/}
of Kuranishi spaces when $\bs f,\bs f'$ are strong submersions to a
manifold $Y$. We generalize to orbifolds $Y$. One can also
generalize to Kuranishi spaces $Y$, but we will not need to.

\begin{dfn} Let $X,X'$ be Kuranishi spaces, $Y$ be a smooth
orbifold, and $\bs f:X\ra Y$, $\bs f':X'\ra Y$ be strongly smooth
maps inducing continuous maps $f:X\ra Y$ and $f':X'\ra Y$. Suppose
at least one of $\bs f,\bs f'$ is a strong submersion. We shall
define the {\it fibre product\/} $X\t_YX'$ or $X\t_{\bs f,Y,\bs
f'}X'$, a Kuranishi space. As a set, the underlying topological
space $X\t_YX'$ is given by \eq{kh2eq10}, as for fibre products of
orbifolds. The definition is modelled on Definition~\ref{kh2def9}.

Let $p\in X$, $p'\in X'$ and $q\in Y$ with $f(p)=q=f'(p')$. Let
$(V_p,E_p,s_p,\psi_p)$, $(V'_{\smash{p'}},E'_{\smash{p'}},
s'_{\smash{p'}},\psi'_{\smash{p'}})$ be sufficiently small Kuranishi
neighbourhoods in the germs at $p,p'$ in $X,X'$, and $f_p:V_p\ra Y$,
$f'_{\smash{p'}}:V'_{\smash{p'}}\ra Y$ be smooth maps in the germs
of $\bs f,\bs f'$ at $p,p'$ respectively. Define a Kuranishi
neighbourhood on $X\t_YX'$ by
\e
\begin{split}
\bigl(V_p\t_{f_p,Y,f'_{\smash{p'}}}&V'_{\smash{p'}},\pi_{V_p}^*(E_p)\op
\pi_{V'_{\smash{p'}}}^*(E'_{\smash{p'}}),\\
&(s_p\ci\pi_{V_p})\op(s'_{\smash{p'}}\ci\pi_{V'_{\smash{p'}}}),
(\psi_p\ci\pi_{V_p})\t(\psi'_{\smash{p'}}\ci\pi_{V'_{\smash{p'}}})\t
\chi_{pp'}\bigr).
\end{split}
\label{kh2eq14}
\e

Here $V_p\t_{f_p,Y,f'_{\smash{p'}}}V'_{\smash{p'}}$ is the fibre
product of orbifolds from Definition \ref{kh2def9}, defined as at
least one of $f_p,f'_{\smash{p'}}$ is a submersion, and
$\pi_{V_p},\pi_{V'_{\smash{p'}}}$ are the projections from
$V_p\t_YV'_{\smash{p'}}$ to $V_p,V'_{\smash{p'}}$. The final term
$\chi_{pp'}$ in \eq{kh2eq14} maps the biquotient terms in
\eq{kh2eq10} for $V_p\t_YV'_{\smash{p'}}$ to the same terms in
\eq{kh2eq10} for the set $X\t_YX'$, so that for $v\in V_p$ and
$v'\in V'_{\smash{p'}}$ with
$f_p(v)=f(\psi_p(v))=y=f'_{\smash{p'}}(v')=
f'(\psi'_{\smash{p'}}(v'))$ in $Y$ we define
\begin{align*}
\chi_{pp'}\vert_{(v,v')}:\,&(f_p)_*(\Stab(v))\backslash\Stab
(y)/(f'_{\smash{p'}})_*(\Stab(v')) \longra\\
&f_*(\Stab(\psi_p(v)))\backslash\Stab(y)/
f'_*(\Stab(\psi'_{\smash{p'}}(v'))
\end{align*}
to be the map induced by the isomorphisms
$(\psi_p)_*:\Stab(v)\ra\Stab(\psi_p(v))$
and~$(\psi'_{p'})_*:\Stab(v')\ra\Stab(\psi'_{\smash{p'}}(v'))$.

It is easy to verify that coordinate changes between Kuranishi
neighbourhoods in $X$ and $X'$ induce coordinate changes between
neighbourhoods \eq{kh2eq14}. So the systems of germs of Kuranishi
neighbourhoods and coordinate changes on $X,X'$ induce such systems
on $X\t_YX'$. This gives a {\it Kuranishi structure\/} on $X\t_YX'$,
making it into a {\it Kuranishi space}. Clearly
$\vdim(X\t_YX')=\vdim X+\vdim X'-\dim Y$, and $X\t_YX'$ is compact
if $X,X'$ are compact.
\label{kh2def20}
\end{dfn}

\subsection{Orientations and orientation conventions}
\label{kh27}

We can now discuss {\it orientations\/} of Kuranishi spaces $X$. Our
definition is basically equivalent to Fukaya et al.\
\cite[Def.~A1.17]{FOOO}, noting as in Remarks \ref{kh2rem5} and
\ref{kh2rem6} that our notion of Kuranishi structure is roughly
equivalent to that of Kuranishi structure with tangent bundles in
the 2008 version of~\cite{FOOO}.

\begin{dfn} Let $X$ be a Kuranishi space. An {\it orientation\/} on $X$
assigns, for all $p\in X$ and all sufficiently small Kuranishi
neighbourhoods $(V_p,E_p,s_p,\psi_p)$ in the germ at $p$,
orientations on the fibres of the orbibundle $TV_p\op E_p$ varying
continuously over $V_p$. These must be compatible with coordinate
changes, in the following sense. Let $q\in\Im\psi_p$,
$(V_q,\ldots,\psi_q)$ be sufficiently small in the germ at $q$, and
$(\phi_{pq},\hat\phi_{pq})$ be the coordinate change from
$(V_q,\ldots,\psi_q)$ to $(V_p,\ldots,\psi_p)$ in the germ. Define
$\d\hat s_p$ near $s_q^{-1}(0)\subseteq V_q$ as in~\eq{kh2eq12}.

Locally on $V_q$, choose any orientation for the fibres of
$\phi_{pq}^*(TV_p)/(\d\phi_{pq})(TV_q)$, and let
$\phi_{pq}^*(E_p)/\hat\phi_{pq}(E_q)$ have the orientation induced
from this by the isomorphism $\d\hat s_p$ in \eq{kh2eq12}. These
induce an orientation on $\frac{\phi_{pq}^*(TV_p)}{
(\d\phi_{pq})(TV_q)}\op\frac{\phi_{pq}^*(E_p)}{
\hat\phi_{pq}(E_q)}$, which is independent of the choice for
$\phi_{pq}^*(TV_p)/(\d\phi_{pq})(TV_q)$. Thus, these local choices
induce a natural orientation on the orbibundle
$\frac{\phi_{pq}^*(TV_p)}{(\d\phi_{pq})(TV_q)}\op
\frac{\phi_{pq}^*(E_p)}{\hat\phi_{pq}(E_q)}$ near $s_q^{-1}(0)$. We
require that in oriented orbibundles over $V_q$ near $s_q^{-1}(0)$,
we have
\e
\begin{split}
\phi_{pq}^*\bigl[TV_p\op E_p\bigr]\cong (-1)^{\dim V_q(\dim
V_p-\dim V_q)}\bigl[TV_q\op E_q\bigr]\op&\\
\bigl[\ts\frac{\phi_{pq}^*(TV_p)}{(\d\phi_{pq})(TV_q)}\op
\frac{\phi_{pq}^*(E_p)}{\hat\phi_{pq}(E_q)}\bigr]&,
\end{split}
\label{kh2eq15}
\e
where $TV_p\op E_p$ and $TV_q\op E_q$ have the orientations assigned
by the orientation on $X$. An {\it oriented Kuranishi space\/} is a
Kuranishi space with an orientation.
\label{kh2def21}
\end{dfn}

Suppose $X,X'$ are oriented Kuranishi spaces, $Y$ is an oriented
orbifold, and $\bs f:X\ra Y$, $\bs f':X'\ra Y$ are strong
submersions. Then by \S\ref{kh26} we have Kuranishi spaces $\pd X$
and $X\t_YX'$. These can also be given orientations in a natural
way. We follow the orientation conventions of Fukaya et al.~\cite[\S
45]{FOOO}.

\begin{conv} First, our conventions for manifolds:
\begin{itemize}
\setlength{\itemsep}{0pt}
\setlength{\parsep}{0pt}
\item[(a)] Let $X$ be an oriented manifold with boundary $\pd
X$. Then we define the orientation on $\pd X$ such that
$TX\vert_{\pd X}=\R_{\rm out}\op T(\pd X)$ is an isomorphism of
oriented vector bundles, where $\R_{\rm out}$ is oriented by an
outward-pointing normal vector to~$\pd X$.
\item[(b)] Let $X,X',Y$ be oriented manifolds, and $f:X\ra
Y$, $f':X'\ra Y$ be submersions. Then $\d f:TX\ra f^*(TY)$ and $\d
f':TX'\ra (f')^*(TY)$ are surjective maps of vector bundles over
$X,X'$. Choosing Riemannian metrics on $X,X'$ and identifying the
orthogonal complement of $\Ker\d f$ in $TX$ with the image $f^*(TY)$
of $\d f$, and similarly for $f'$, we have isomorphisms of vector
bundles over $X,X'$:
\e
TX\cong f^*(TY)\op\Ker\d f \quad\text{and}\quad TX'\cong
(f')^*(TY)\op\Ker\d f'.
\label{kh2eq16}
\e

Define orientations on the fibres of $\Ker\d f$, $\Ker\d f'$ over
$X,X'$ such that \eq{kh2eq16} are isomorphisms of oriented vector
bundles, where $TX,TX'$ are oriented by the orientations on $X,X'$,
and $f^*(TY),(f')^*(TY)$ by the orientation on $Y$. Then we define
the orientation on $X\t_YX'$ so that
\e
\begin{split}
T(X\t_YX')&\cong (f\ci\pi_X)^*(TY)\op\pi_X^*(\Ker\d f)
\op\pi_{X'}^*(\Ker\d f')\\
&\cong \pi_X^*(TX)\op\pi_{X'}^*(\Ker\d f')\\
&\cong (-1)^{\dim Y(\dim X-\dim Y)}\pi_X^*(\Ker\d
f)\op\pi_{X'}^*(TX')
\end{split}
\label{kh2eq17}
\e
are isomorphisms of oriented vector bundles. Here $\pi_X:X\t_YX'\ra
X$ and $\pi_{X'}:X\t_YX'\ra X'$ are the natural projections, and
$f\ci\pi_X\equiv f'\ci\pi_{X'}$.

Note that the second line of \eq{kh2eq17} makes sense if $f'$ is a
submersion but $f$ is only smooth, and the third line makes sense if
$f$ is a submersion but $f'$ is only smooth. Thus, our convention
extends to fibre products $X\t_{f,Yf'}X'$ in which only one of
$f,f'$ is a submersion.
\end{itemize}
These extend immediately to orbifolds. They also extend to the
Kuranishi space versions in Definitions \ref{kh2def19} and
\ref{kh2def20}; for Definition \ref{kh2def20} they are described in
\cite[Conv.~45.1(4)]{FOOO}. We will not give details, but here is an
algorithm to deduce Kuranishi space orientation conventions from
manifold ones.

Let $M$ be an oriented manifold or orbifold. Then $M$ is a Kuranishi
space with a single Kuranishi neighbourhood
$(V_q,E_q,s_q,\psi_q)=(M,M,0,\id_M)$. The orientation on $M$ induces
an orientation on the fibres of $TV_q\op E_q=TM$, so $M$ is oriented
in the Kuranishi space sense.

Let $E\ra M$ be an vector bundle or orbibundle. Write $\bar E$ for
the total space of $E$, as an orbifold, with projection $\pi:\bar
E\ra M$, and regard $M$ as a suborbifold of $\bar E$, the zero
section of $E$. Then we can define a second Kuranishi neighbourhood
$(V_p,E_p,s_p,\psi_p)=(\bar E,\pi^*(E),\id_E,\id_M)$. Here
$s_p=\id_E:\bar E\ra\pi^*(E)$ is just the identity map on $E$, so
$s_p^{-1}(0)$ is $M\subset\bar E$, the zero section of $E$. There is
an obvious coordinate change~$(\phi_{pq},\smash{\hat\phi_{pq}})
=(\id_M,0)$.

Now given any operation requiring an orientation convention which
makes sense for manifolds, such as boundaries or fibre products of
smooth submersions above, we can define it for manifolds, and then
use the example above and Definition \ref{kh2def21} to deduce what
the convention must be for Kuranishi neighbourhoods of the form
$(\bar E,\pi^*(E),\id_E,\id_M)$. This is enough to define the
convention for general Kuranishi neighbourhoods.
\label{kh2conv1}
\end{conv}

If $X$ is an oriented Kuranishi space, we often write $-X$ for the
same Kuranishi space with the opposite orientation. The next result
comes from \cite[Lem.~45.3]{FOOO}, except the second line of
\eq{kh2eq18}, which is elementary. In (a) we assume $\pd Y=\es$  to
avoid having to split $\pd X_1=\pd_+^{\bs f_1}X_1\amalg\pd_-^{\bs
f_1}X_1$, and so on, as in the manifolds case in Definition
\ref{kh2def4} and Proposition~\ref{kh2prop1}.

\begin{prop} Let\/ $X_1,X_2,\ldots$ be oriented Kuranishi spaces,
$Y,Y_1,\ldots$ be oriented orbifolds, and\/ $\bs f_1:X_1\ra
Y,\ldots$ be strongly smooth maps, with at least one strong
submersion in each fibre product below. Then the following hold, in
oriented Kuranishi spaces:
\begin{itemize}
\setlength{\itemsep}{0pt}
\setlength{\parsep}{0pt}
\item[{\rm(a)}] If\/ $\pd Y=\es,$ for\/ $\bs f_1:X_1\ra Y$ and\/
$\bs f_2:X_2\ra Y$ we have
\e
\!\!\!\!\!\!
\begin{gathered}
\pd(X_1\t_YX_2)\cong(\pd X_1)\t_YX_2\amalg (-1)^{\vdim X_1+\dim Y}
X_1\t_Y(\pd X_2)\\
\text{and}\quad X_1\t_YX_2\cong(-1)^{(\vdim X_1-\dim Y)(\vdim
X_2-\dim Y)}X_2\t_YX_1.
\end{gathered}
\label{kh2eq18}
\e
\item[{\rm(b)}] For\/ $\bs f_1:X_1\ra Y_1,$ $\bs f_2:X_2\ra Y_1\t
Y_2$ and\/ $\bs f_3:X_3\ra Y_2$, we have
\begin{equation*}
(X_1\t_{Y_1}X_2)\t_{Y_2}X_3\cong X_1\t_{Y_1}(X_2\t_{Y_2}X_3).
\end{equation*}
\item[{\rm(c)}] For\/ $\bs f_1\t\bs f_2:X_1\ra Y_1\t Y_2,$ $\bs
f_3:X_2\ra Y_1$ and\/ $\bs f_4:X_3\ra Y_2$, we have
\begin{equation*}
\!\!\!\!\!\!\!\! X_1\t_{Y_1\t Y_2}(X_2\t X_3)\cong (-1)^{\dim
Y_2(\dim Y_1+\vdim X_2)} (X_1\t_{Y_1}X_2)\t_{Y_2}X_3.
\end{equation*}
\end{itemize}
\label{kh2prop2}
\end{prop}

A useful special case of Proposition \ref{kh2prop2}(a) with
$X_1=[0,1]$, $X_2=X$ and $Y$ a point gives the decomposition in
oriented Kuranishi spaces:
\begin{equation*}
\pd\bigl([0,1]\t X\bigr)\cong\bigl(\{1\}\t X\bigr)\amalg
-\bigl(\{0\}\t X\bigr)\amalg - \bigl([0,1]\t\pd X\bigr).
\end{equation*}

\subsection{Coorientations}
\label{kh28}

We will also need a notion of {\it relative orientation} for a
strong submersion $\bs f:X\ra Y$. We call it a {\it coorientation},
as we use orientations in Kuranishi homology, and coorientations in
Kuranishi cohomology.

\begin{dfn} Let $X$ be a Kuranishi space, $Y$ an orbifold, and $\bs
f:X\ra Y$ a strong submersion. A {\it coorientation\/} for $(X,\bs
f)$ assigns, for all $p\in X$ and all sufficiently small Kuranishi
neighbourhoods $(V_p,E_p,s_p,\psi_p)$ in the germ at $p$ with
submersion $f_p:V_p\ra Y$ representing $\bs f$, orientations on the
fibres of the orbibundle $\Ker\d f_p\op E_p$ varying continuously
over $V_p$, where $\d f_p:TV_p\ra f_p^*(TY)$ is the (surjective)
derivative of $f_p$.

These must be compatible with coordinate changes, in the following
sense. Let $q\in\Im\psi_p$, $(V_q,\ldots,\psi_q)$ be sufficiently
small in the germ at $q$, let $f_q:V_q\ra Y$ represent $\bs f$, and
$(\phi_{pq},\hat\phi_{pq})$ be the coordinate change from
$(V_q,\ldots,\psi_q)$ to $(V_p,\ldots,\psi_p)$ in the germ. Then we
require that in oriented orbibundles over $V_q$ near $s_q^{-1}(0)$,
we have
\e
\begin{split}
\phi_{pq}^*\bigl[\Ker\d f_p\op E_p\bigr]\cong (-1)^{\dim V_q(\dim
V_p-\dim V_q)}\bigl[\Ker\d f_q\op E_q\bigr]\op&\\
\bigl[\ts\frac{\phi_{pq}^*(TV_p)}{(\d\phi_{pq})(TV_q)}\op
\frac{\phi_{pq}^*(E_p)}{\hat\phi_{pq}(E_q)}\bigr]&,
\end{split}
\label{kh2eq19}
\e
by analogy with \eq{kh2eq15}, where $\frac{\phi_{pq}^*(TV_p)}{
(\d\phi_{pq})(TV_q)}\op\frac{\phi_{pq}^*(E_p)}{\hat\phi_{pq}(E_q)}$
is oriented as in Definition~\ref{kh2def21}.

Suppose now that $Y$ is oriented. Then an orientation on $X$ is
equivalent to a coorientation for $(X,\bs f)$, since for all $p,
(V_p,\ldots,\psi_p),f_p$ as above, the isomorphism $TV_p\cong
f_p^*(TY)\op \Ker\d f_p$ induces isomorphisms of orbibundles
over~$V_p$:
\e
\bigl(TV_p\op E_p\bigr)\cong f_p^*(TY)\op\bigl(\Ker\d f_p\op
E_p\bigr).
\label{kh2eq20}
\e
There is a 1-1 correspondence between orientations on $X$ and
coorientations for $(X,\bs f)$ such that \eq{kh2eq20} holds in
oriented orbibundles, where $TV_p\op E_p$ is oriented by the
orientation on $X$, and $\Ker\d f_p\op E_p$ by the coorientation for
$(X,\bs f)$, and $f_p^*(TY)$ by the orientation on $Y$. Taking the
direct sum of $f_p^*(TY)$ with \eq{kh2eq19} and using \eq{kh2eq20}
yields \eq{kh2eq15}, so this is compatible with coordinate changes.
\label{kh2def22}
\end{dfn}

\begin{conv} Here are our conventions for manifolds:
\begin{itemize}
\setlength{\itemsep}{0pt}
\setlength{\parsep}{0pt}
\item[(a)] Let $X,Y$ be manifolds with $\pd Y=\es$ and $f:X\ra Y$ a
submersion. Suppose $(X,f)$ is cooriented. Then we define a
coorientation for $(\pd X,f\vert_{\pd X})$ such that $\Ker(\d
f)\vert_{\pd X}=(-1)^{\dim Y}\R_{\rm out}\op\Ker\d (f\vert_{\pd X})$
is an isomorphism of oriented vector bundles, where $(\d
f)\vert_{\pd X}$ maps $TX\vert_{\pd X}\ra f^*(TY)\vert_{\pd X}$ and
$\d (f\vert_{\pd X})$ maps $T(\pd X)\ra f^*(TY)\vert_{\pd X}$ and
$\R_{\rm out}$ is oriented by an outward-pointing normal vector to
$\pd X$ in~$\Ker(\d f)\vert_{\pd X}$.

Now drop the assumption that $\pd Y=\es$. Then as in Definition
\ref{kh2def4} we have a decomposition $\pd
X=\pd_+^fX\amalg\pd_-^fX$, and submersions $f_+:\pd_+^fX\ra Y$ and
$f_-:\pd_-^fX\ra\pd Y$. Define a coorientation for $(\pd_+^fX,f_+)$
such that $\Ker(\d f)\vert_{\pd_+^fX}=(-1)^{\dim Y}\R_{\rm
out}\op\Ker\d f_+$. Define a coorientation for $(\pd_-^fX,f_-)$ such
that~$\Ker(\d f)\vert_{\pd_-^fX}=\Ker\d f_-$.
\item[(b)] Let $X,X',Y$ be manifolds, $f:X\ra Y$, $f':X'\ra Y$ be
submersions, and $(X,f),(X',f')$ be cooriented. Define a
coorientation on $(X\t_YX',\pi_Y)$ such that
\begin{equation*}
\Ker\bigl(\d\pi_Y:T(X\t_YX')\ra\pi_Y^*(TY)\bigr)\cong\pi_X^*(\Ker\d
f)\op\pi_{X'}^*(\Ker\d f')
\end{equation*}
is an isomorphism of oriented vector bundles, where $\Ker\d f$ and
$\Ker\d f'$ are oriented by the coorientations for $(X,f),(X',f')$.
Note that we do not need $Y$ oriented, in contrast to
Convention~\ref{kh2conv1}(b).
\item[(c)] We can also combine orientations and coorientations in
fibre products, as follows. Let $X,X',Y$ be manifolds, $f:X\ra Y$ be
smooth, $f':X'\ra Y$ be a submersion, and suppose $X$ is oriented
and $(X',f')$ is cooriented. As in \eq{kh2eq17}, define orientations
on $X\t_YX'$ and $X'\t_YX$ such that
\begin{gather*}
T(X\t_YX')\cong \pi_X^*(TX)\op\pi_{X'}^*(\Ker\d f'),\\
T(X'\t_YX)\cong (-1)^{\dim Y(\dim X'-\dim Y)}\pi_{X'}^*(\Ker\d
f')\op\pi_X^*(TX)
\end{gather*}
are isomorphism of oriented vector bundles, where $TX,\Ker\d f'$ are
oriented by the orientation on $X$ and the coorientation on
$(X',f')$. Note that we do not need $Y$ oriented, nor $f$ a
submersion.
\end{itemize}
These extend to orbifolds and Kuranishi spaces as for
Convention~\ref{kh2conv1}.
\label{kh2conv2}
\end{conv}

As in Definition \ref{kh2def22}, if $\bs f:X\ra Y$ is a strong
submersion and $Y$ is oriented, there is a 1-1 correspondence
between orientations on $X$ and coorientations for $(X,\bs f)$. The
signs in Convention \ref{kh2conv2} then correspond with those in
Convention \ref{kh2conv1}. Therefore the analogue of Proposition
\ref{kh2prop2} holds for coorientations. In particular, for strong
submersions $\bs f_a:X_a\ra Y$ with $(X_a,\bs f_a)$ cooriented for
$a=1,2,3$, taking $\pd Y=\es$ in \eq{kh2eq21}, we have
\begin{gather}
\begin{split}
\bigl(\pd(X_1\t_YX_2),\bs\pi_Y\bigr)&\cong\bigl((\pd
X_1)\t_YX_2,\bs\pi_Y\bigr)\,\amalg \\
&\qquad\quad (-1)^{\vdim X_1+\dim Y} \bigl(X_1\t_Y(\pd
X_2),\bs\pi_Y\bigr),
\end{split}
\label{kh2eq21}\\
\bigl(X_1\t_YX_2,\bs\pi_Y\bigr)\cong(-1)^{(\vdim X_1-\dim Y)(\vdim
X_2-\dim Y)}\bigl(X_2\t_YX_1,\bs\pi_Y\bigr),
\label{kh2eq22}\\
\bigl((X_1\t_YX_2)\t_YX_3,\bs\pi_Y\bigr)\cong
\bigl(X_1\t_Y(X_2\t_YX_3),\bs\pi_Y\bigr).
\label{kh2eq23}
\end{gather}
Similarly, if $X_1$ is oriented, $\bs f_1:X_1\ra Y$ is strongly
smooth, $\bs f_2:X_2\ra Y$ is a cooriented strong submersion, and
$\pd Y=\es$ then
\e
\pd(X_1\t_YX_2)\cong \bigl((\pd X_1)\t_YX_2\bigr)\amalg (-1)^{\vdim
X_1+\dim Y}\bigl(X_1\t_Y(\pd X_2)\bigr)
\label{kh2eq24}
\e
in oriented Kuranishi spaces, combining orientations and
coorientations as in Convention~\ref{kh2conv2}(c).

Here is how to generalize equations \eq{kh2eq21} and \eq{kh2eq24} to
the case $\pd Y\ne\es$. If $X$ is a Kuranishi space and $\bs f:X\ra
Y$ a strong submersion then as for submersions of manifolds in
Definition \ref{kh2def4}, we can define a decomposition $\pd
X=\pd_+^{\bs f}X\amalg\pd_-^{\bs f}X$, and strong submersions $\bs
f_+:\pd_+^{\bs f}X\ra Y$ and $\bs f_-:\pd_-^{\bs f}X\ra\pd Y$. When
$\bs f_1:X_1\ra Y$ and $\bs f_2:X_2\ra Y$ are cooriented strong
submersions, and we do not assume $\pd Y=\es$, equation \eq{kh2eq21}
should be replaced by
\ea
\begin{split}
\bigl(\pd_+^{\bs\pi_Y}(X_1\t_YX_2),\bs\pi_Y\bigr)&\cong\bigl((\pd_+^{\bs
f_1}X_1)\t_YX_2,\bs\pi_Y\bigr)\,\amalg \\
&\qquad\quad (-1)^{\vdim X_1+\dim Y} \bigl(X_1\t_Y(\pd_+^{\bs f_2}
X_2),\bs\pi_Y\bigr),
\end{split}
\label{kh2eq25}\\
\bigl(\pd_-^{\bs\pi_Y}(X_1\t_YX_2),\bs\pi_{\pd
Y}\bigr)&\cong\bigl((\pd_-^{\bs f_1}X_1)\t_{\pd Y}(\pd_-^{\bs
f_2}X_2),\bs\pi_{\pd Y}\bigr).
\label{kh2eq26}
\ea

If $X_1$ is oriented, $\bs f_1:X_1\ra Y$ is strongly smooth, and
$\bs f_2:X_2\ra Y$ is a cooriented strong submersion, then equation
\eq{kh2eq24} should be replaced by
\e
\pd(X_1\t_YX_2)\cong \bigl((\pd X_1)\t_YX_2\bigr)\amalg (-1)^{\vdim
X_1+\dim Y}\bigl(X_1\t_Y(\pd_+^{\bs f_2}X_2)\bigr).
\label{kh2eq27}
\e

\subsection{Almost complex and almost CR structures}
\label{kh29}

We now define {\it almost complex structures} on Kuranishi spaces.
They are a substitute for Fukaya and Ono's notion of {\it stably
almost complex\/} Kuranishi spaces \cite[Def.~5.17]{FuOn1}. We will
use them to define almost complex Kuranishi bordism $KB_*^{\rm
ac}(Y;R)$ in Chapter \ref{kh5}, and in the author's approach to the
Gopakumar--Vafa Integrality Conjecture in Chapter~\ref{kh6}.

\begin{dfn} Let $V$ be an orbifold without boundary. An {\it almost
complex structure\/} $J$ on $V$ is a tensor $J=J_a^b$ in
$C^\iy(TV\ot T^*V)$ with $J^2=-1$, that is, $J_a^bJ_b^c=-\de_a^c$ in
index notation. Let $E\ra V$ be an orbibundle. An {\it almost
complex structure\/ $K$ on the fibres of\/} $E$ is $K\in
C^\iy(E^*\ot E)$ with~$K^2=-1$.

Let $X$ be a Kuranishi space without boundary. An {\it almost
complex structure\/} $(\bs J,\bs K)$ on $X$ assigns for all $p\in X$
and all sufficiently small $(V_p,E_p,s_p,\psi_p)$ in the germ at
$p$, a choice of almost complex structure $J_p$ on $V_p$, and a
choice of almost complex structure $K_p$ on the fibres of $E_p$.
These choices must satisfy the following conditions. For all $p\in
X$, for all $(V_p,\ldots,\psi_p)$ sufficiently small in the germ at
$p$ with almost complex structures $J_p,K_p$, for all
$q\in\Im\psi_p$, and for all sufficiently small
$(V_q,\ldots,\psi_q)$ in the germ at $q$ with almost complex
structures $J_q,K_q$, if $(\phi_{pq},\hat\phi_{pq})$ is the
coordinate change from $(V_q,\ldots,\psi_q)$ to
$(V_p,\ldots,\psi_p)$ in the germ of coordinate changes, then:
\begin{itemize}
\setlength{\itemsep}{0pt}
\setlength{\parsep}{0pt}
\item[(a)] $\d\phi_{pq}\ci J_q=\phi_{pq}^*(J_p)\ci\d\phi_{pq}$ as
morphisms of orbibundles $TV_q\ra\phi_{pq}^*(TV_p)$, that is,
$\phi_{pq}$ is a pseudoholomorphic map of almost complex orbifolds;
\item[(b)] $\hat\phi_{pq}\ci K_q=\phi_{pq}^*(K_p)\ci\hat\phi_{pq}$
as morphisms of orbibundles $E_q\ra\phi_{pq}^*(E_p)$; and
\item[(c)] Parts (a) and (b) imply that the orbibundles
$\phi_{pq}^*(TV_p)/(\d\phi_{pq})(TV_q)$ and $\phi_{pq}^*(E_p)/
\hat\phi_{pq}(E_q)$ over $V_q$ appearing in \eq{kh2eq12} have almost
complex structures $J_{pq},K_{pq}$ on their fibres, by projection
from $\phi_{pq}^*(J_p),\phi_{pq}^*(K_p)$. We require that
$K_{pq}\ci\d\hat s_p= \d\hat s_p\ci J_{pq}$ over $s_q^{-1}(0)$ in
$V_q$, for $\d\hat s_p$ as in~\eq{kh2eq12}.
\end{itemize}
\label{kh2def23}
\end{dfn}

If $(J_p,K_p)$ exists for $(V_p,E_p,s_p,\psi_p)$ as above then $\dim
V_p$ and $\rank E_p$ are even, so $\vdim X=\dim V_p-\rank E_p$ is
even. Thus $X$ can admit an almost complex structure only if $\vdim
X$ is even. Our next definition is an odd-dimensional analogue of
almost complex structures, in the same way that contact structures
are odd-dimensional analogues of symplectic structures.

\begin{dfn} Let $V$ be an orbifold with boundary, but without
(g-)corners. An {\it almost CR structure\/} $(D,J)$ on $V$ is an
orbifold vector subbundle $D$ of $TV$ with $\dim V-\rank D=1$, such
that $D\vert_{\pd V}=T(\pd V)\subset TV\vert_{\pd V}$, and an almost
complex structure $J$ on the fibres of $D$.

Let $X$ be a Kuranishi space with boundary, but without (g-)corners.
An {\it almost CR structure\/} $(\bs D,\bs J,\bs K)$ on $X$ assigns
for all $p\in X$ and all sufficiently small $(V_p,E_p,s_p,\psi_p)$
in the germ at $p$, a choice of almost CR structure $(D_p,J_p)$ on
$V_p$, and a choice of almost complex structure $K_p$ on the fibres
of $E_p$. These choices must satisfy the following conditions. For
all $p\in X$, for all $(V_p,\ldots,\psi_p)$ sufficiently small in
the germ at $p$ with structures $(D_p,J_p),K_p$, for all
$q\in\Im\psi_p$, and for all sufficiently small
$(V_q,\ldots,\psi_q)$ in the germ at $q$ with structures
$(D_q,J_q),K_q$, if $(\phi_{pq},\hat\phi_{pq})$ is the coordinate
change from $(V_q,\ldots,\psi_q)$ to $(V_p,\ldots,\psi_p)$ in the
germ of coordinate changes, then:
\begin{itemize}
\setlength{\itemsep}{0pt}
\setlength{\parsep}{0pt}
\item[(a)] $\d\phi_{pq}(D_q)$ is an orbisubbundle of
$\phi_{pq}^*(D_p)$ in $\phi_{pq}^*(TV_p)$, and $\d\phi_{pq}\ci
J_q=\phi_{pq}^*(J_p)\ci\d\phi_{pq}\vert_{D_q}$ as morphisms of
orbibundles $D_q\ra\phi_{pq}^*(D_p)$;
\item[(b)] $\hat\phi_{pq}\ci K_q=\phi_{pq}^*(K_p)\ci\hat\phi_{pq}$
as morphisms of orbibundles $E_q\ra\phi_{pq}^*(E_p)$; and
\item[(c)] Part (a) and the fact that $D_q,D_p$ both have codimension
1 yield an isomorphism $\phi_{pq}^*(D_p)/(\d\phi_{pq})(D_q)\cong
\phi_{pq}^*(TV_p)/(\d\phi_{pq})(TV_q)$ of orbibundles over $V_q$. So
projecting $\phi_{pq}^* (J_p)$ gives an almost complex structure
$J_{pq}$ on $\phi_{pq}^*(TV_p)/(\d\phi_{pq})(TV_q)$. Part (b)
implies that projecting $\phi_{pq}^*(K_p)$ gives an almost complex
structure $K_{pq}$ on $\phi_{pq}^*(E_p)/\hat\phi_{pq}(E_q)$. We
require that $K_{pq}\ci\d\hat s_p= \d\hat s_p\ci J_{pq}$ over
$s_q^{-1}(0)$ in $V_q$, for $\d\hat s_p$ as in~\eq{kh2eq12}.
\end{itemize}

If $V$ is an orbifold with boundary but without (g-)corners, and
$(D,J)$ is an almost CR structure on $V$, then $\pd V$ is an
orbifold without boundary, and $J\vert_{\pd V}$ is an almost complex
structure on $\pd V$, as $T(\pd V)=D\vert_{\pd V}$. In the same way,
if $X$ is a Kuranishi space with boundary but without (g-)corners
and $(\bs D,\bs J,\bs K)$ is an almost CR structure on $X$ then $\pd
X$ is a Kuranishi space without boundary, and $(\bs J,\bs
K)\vert_{\pd X}$ is an almost complex structure on~$\pd X$.
\label{kh2def24}
\end{dfn}

\subsection{Co-almost complex and co-almost CR structures}
\label{kh210}

Just as we can generalize orientations on $X$ in \S\ref{kh27} to
coorientations for $(X,\bs f)$ for $\bs f:X\ra Y$ a strong
submersion in \S\ref{kh28}, so we can generalize \S\ref{kh29} to
give notions of {\it co-almost complex structures} and {\it
co-almost CR structures} for~$(X,\bs f)$.

\begin{dfn} Let $V,Y$ be orbifolds without boundary, and $f:V\ra Y$
a submersion. A {\it co-almost complex structure} $L$ for $(V,f)$ is
an almost complex structure $L$ on the fibres of the orbibundle
$\Ker\bigl(\d f:TV\ra f^*(TY)\bigr)$ over~$V$.

Let $X$ be a Kuranishi space without boundary, and $\bs f:X\ra Y$ a
strong submersion. A {\it co-almost complex structure\/} $(\bs L,\bs
K)$ for $(X,\bs f)$ assigns for all $p\in X$ and all sufficiently
small $(V_p,E_p,s_p,\psi_p)$ in the germ at $p$ with submersion
$f_p:V_p\ra Y$, a co-almost complex structure $L_p$ for $(V_p,f_p)$,
and an almost complex structure $K_p$ on the fibres of $E_p$,
satisfying analogues of Definition~\ref{kh2def23}(a)--(c).

Let $V$ be an orbifold with boundary but without (g-)corners, and
$f:V\ra Y$ a submersion. A {\it co-almost CR structure\/} $(D,L)$
for $(V,f)$ is an orbifold vector subbundle $D$ of $\Ker\bigl(\d
f:TV\ra f^*(TY)\bigr)$ with $\rank D=\dim V-\dim Y-1$, such that
$D\vert_{\pd V}=\Ker\bigl(\d (f\vert_{\pd V}):T(\pd V)\ra
f\vert_{\pd V}^*(TY)\bigr)$, and an almost complex structure $L$ on
the fibres of~$D$.

Let $X$ be a Kuranishi space with boundary, but without (g-)corners.
A {\it co-almost CR structure\/} $(\bs D,\bs L,\bs K)$ for $(X,\bs
f)$ assigns a co-almost CR structure $(D_p,L_p)$ for $(V_p,f_p)$,
and an almost complex structure $K_p$ on the fibres of $E_p$, for
all $p,(V_p,\ldots,\psi_p),f_p$ as above, satisfying analogues of
Definition~\ref{kh2def24}(a)--(c).

As in Definition \ref{kh2def24}, if $(D,L)$ is a co-almost CR
structure for $(V,f)$ then $L\vert_{\pd V}$ is a co-almost complex
structure for $(\pd V,f\vert_{\pd V})$, and if $(\bs D,\bs L,\bs K)$
is a co-almost CR structure for $(X,\bs f)$ then $(\bs L,\bs
K)\vert_{\pd X}$ is a co-almost complex structure for~$(\pd X,\bs
f\vert_{\pd X})$.
\label{kh2def25}
\end{dfn}

We can {\it combine co-almost complex structures on fibre products}.
Let $X,X'$ be Kuranishi spaces, $Y$ an orbifold, $\bs f:X\ra Y$,
$\bs f':X'\ra Y$ strong submersions, and $(\bs L,\bs K),(\bs L',\bs
K')$ co-almost structures for $(X,\bs f),(X',\bs f')$. Then we can
define a natural co-almost complex structure $(\bs L,\bs K)\t_Y(\bs
L',\bs K')$ for $(X\t_YX',\bs\pi_Y)$. Suppose $p\in X$, $p'\in X'$
with $f(p)=f'(p')$ in $Y$, and $(V_p,\ldots,\psi_p),
(V_{\smash{p'}}',\ldots,\psi_{\smash{p'}}')$ are sufficiently small
Kuranishi neighbourhoods of $p,p'$ in $X,X'$, and $f_p:V_p\ra Y$,
$f_{\smash{p'}}':V_{\smash{p'}}'\ra Y$ are submersions representing
$\bs f,\bs f'$, and $(L_p,K_p),(L_{\smash{p'}}',K_{\smash{p'}}')$
represent $(\bs L,\bs K),(\bs L',\bs K')$ on
$(V_p,\ldots,\psi_p),(V_{\smash{p'}}',\ldots,\psi_{\smash{p'}}')$.

Then $(V_{pp'},\ldots,\psi_{pp'})$ is a Kuranishi neighbourhood on
$X\t_YX'$, where $V_{pp'}=V_p\t_YV_{\smash{p'}}'$, and
$E_{pp'}=\pi_{V_p}^*(E_p)\op\pi_{V_{\smash{p'}}'}^*(E_{\smash{p'}}')$,
and $\pi_{pp'}=\pi_Y:V_{pp'}\ra Y$ represents $\bs\pi_Y:X\t_YX'\ra
Y$. We have a natural isomorphism
\e
\begin{split}
\Ker\bigl(\d\pi_{pp'}:TV_{pp'}\ra \pi_{pp'}^*(TY)\bigr)\cong\,&
\pi_{V_p}^*\bigl(\Ker(\d f_p:TV_p\ra f_p^*(TY))\bigr)\op\\
&\pi_{V_{\smash{p'}}'}^*\bigl(\Ker(\d
f_{\smash{p'}}':TV_{\smash{p'}}'\ra (f_{\smash{p'}}')^*(TY))\bigr).
\end{split}
\label{kh2eq28}
\e
Define an almost complex structure $L_{pp'}$ on the fibres of the
l.h.s.\ of \eq{kh2eq28} to be $\pi_{V_p}^*(L_p)\op
\pi_{V_{\smash{p'}}'}^*(L_{\smash{p'}}')$ on the r.h.s.\ of
\eq{kh2eq28}. Define an almost complex structure $K_{pp'}$ on
$E_{pp'}$ to be $\pi_{V_p}^* (K_p)\op\pi_{V_{\smash{p'}}'}^*
(K_{\smash{p'}}')$ on $\pi_{V_p}^*(E_p)\op\pi_{V_{\smash{p'}}'}^*
(E_{\smash{p'}}')$. Then these $(L_{pp'},K_{pp'})$ for all $p,p'$
induce a co-almost complex structure $(\bs L,\bs K)\t_Y(\bs L',\bs
K')$ for~$(X\t_YX',\bs\pi_Y)$.

Similarly, we can naturally define the fibre product $(\bs L,\bs
K)\t_Y(\bs D',\bs L',\bs K')$ of a co-almost complex structure $(\bs
L,\bs K)$ for $(X,\bs f)$ and a co-almost CR structure $(\bs D',\bs
L',\bs K')$ for $(X',\bs f')$, which is a co-almost CR structure
for~$(X\t_YX',\bs\pi_Y)$.

Given $\bs f:X\ra Y$ strongly smooth and $\bs f':X'\ra Y$ a strong
submersion, it is also possible to define the fibre product of an
almost complex structure $(\bs J,\bs K)$ on $X$ and a co-almost
complex structure $(\bs L',\bs K')$ for $(X',\bs f')$, to get an
almost complex structure $(\bs J,\bs K)\t_Y(\bs L',\bs K')$ for
$X\t_{\bs f,Y,\bs f'}X'$. However, this construction is not quite
natural. The issue is that if $p\in X$ with $(J_p,K_p),f_p$
representing $(\bs J,\bs K),\bs f$ on $(V_p,\ldots,\psi_p)$, and
$p'\in X'$ with $(L'_{p'},K'_{p'}),f'_{p'}$ representing $(\bs
L',\bs K'),\bs f'$ on $(V'_{p'},\ldots,\psi'_{p'})$, then $J_p$ is
an almost complex structure on the fibres of $TV_p$, and $L'_{p'}$
is an almost complex structure on the fibres of $\Ker(\d
f'_{p'}:TV'_{p'}\ra (f')^*(TY))$. We can write
\e
T\bigl(V_p\t_{f_p,Y,f'_{p'}}V'_{p'}\bigr)\cong \pi_{V_p}^*
\bigl(TV_p\bigr)\op\pi_{V'_{p'}}^* \bigl(\Ker(\d f'_{p'}:TV'_{p'}\ra
(f')^*(TY))\bigr),
\label{kh2eq29}
\e
but this splitting is not natural; only the short exact sequence
\begin{equation*}
0\!\ra\!\pi_{V'_{p'}}^*\bigl(\Ker(\d f'_{p'}:TV'_{p'}\!\ra\!
(f')^*(TY))\bigr)\!\ra\! T\bigl(V_p\t_{f_p,Y,f'_{p'}}V'_{p'}\bigr)
\!\ra\!\pi_{V_p}^* \bigl(TV_p\bigr)\!\ra\!0
\end{equation*}
is natural. So to define $(\bs J,\bs K)\t_Y(\bs L',\bs K')$ we must
choose splittings \eq{kh2eq29} for all $p,p'$ with $f(p)=f'(p')$,
compatible with coordinate changes $(\phi_{pq},\hat\phi_{pq})$ and
$(\phi'_{p'q'},\ab\hat\phi'_{p'q'})$, and then define $J_{pp'}$ to
be $\pi_{V_p}^*(J_p)\op \pi_{V_{\smash{p'}}'}^*(L_{\smash{p'}}')$,
using~\eq{kh2eq29}.

In the same way, after choosing splittings \eq{kh2eq29}, we can
define fibre products $(\bs D,\bs J,\bs K)\t_Y(\bs L',\bs K')$ of an
almost CR structure $(\bs D,\bs J,\bs K)$ on $X$ and a co-almost
complex structure $(\bs L',\bs K')$ for $(X',\bs f')$, and $(\bs
J,\bs K)\t_Y(\bs D',\bs L',\bs K')$ of an almost complex structure
$(\bs J,\bs K)$ on $X$ and a co-almost CR structure $(\bs D',\bs
L',\bs K')$ for $(X',\bs f')$. Both are almost CR structures
on~$X\t_{\bs f,Y,\bs f'}X'$.

\section{Gauge-fixing data and co-gauge-fixing data}
\label{kh3}

The basic idea of this book is to define a homology theory
$KH_*(Y;R)$ for an orbifold $Y$ using as chains pairs $(X,\bs f)$,
where $X$ is a compact Kuranishi space and $\bs f:X\ra Y$ is a
strongly smooth map. However, we show in \S\ref{kh49} below that if
we allow chains $(X,\bs f)$ for which the {\it automorphism group\/}
$\Aut(X,\bs f)$ of strong diffeomorphisms $\bs a:X\ra X$ up to
equality such that $\bs f\ci\bs a=\bs f$ is infinite, the resulting
homology groups are always zero.

To deal with this problem, we add some {\it extra data\/} $\bs G$ to
$(X,\bs f)$ for which $\Aut(X,\bs f,\bs G)$ is finite, and use as
chains triples $(X,\bs f,\bs G)$. Exactly what this extra data is
does not matter very much, but here are the properties we would like
it to satisfy for this book and in \cite{AkJo,Joyc2,Joyc3,Joyc4}.
Parts (a),(e) are not essential for well-behaved (co)homology
theories, but will be needed in \cite{AkJo,Joyc2,Joyc3,Joyc4}. In
(e), we restrict to $X$ with corners (not g-corners) because
Principle \ref{kh2pri} fails for smooth extensions from boundaries
of general manifolds or orbifolds with g-corners.

\begin{property} We want to define some kind of {\it extra
data\/} $\bs G$ for pairs $(X,\bs f)$, where $X$ is a compact
Kuranishi space with g-corners, $Y$ an orbifold and $\bs f:X\ra Y$ a
strongly smooth map, such that:
\begin{itemize}
\setlength{\itemsep}{0pt}
\setlength{\parsep}{0pt}
\item[(a)] Every pair $(X,\bs f)$ should admit extra data $\bs G$.
If $\Ga$ is a finite subgroup of $\Aut(X,\bs f)$, then there should
exist $\Ga$-invariant extra data for~$(X,\bs f)$.
\item[(b)] The {\it automorphism group\/} $\Aut(X,\bs f,\bs G)$ of
isomorphisms $(\bs a,\bs b):(X,\bs f,\bs G)\ab\ra(X,\bs f,\bs G)$ is
finite, where $\bs b$ is a lift of $\bs a$ to~$\bs G$.
\item[(c)] Suppose $\Ga$ is a finite group acting on $(X,\bs f,\bs
G)$ by isomorphisms, that is, we are given a group morphism
$\rho:\Ga\ra\Aut(X,\bs f,\bs G)$, which need not be injective. Then
we can form the quotient $\ti X=X/\Ga$, a compact Kuranishi space,
with projection $\bs\pi:X\ra\ti X$, and $\bs f$ pushes down to
$\bs{\ti f}:\ti X\ra Y$ with $\bs f=\bs{\ti f}\ci\bs\pi$. We require
that $\bs G$ should also naturally push down to extra data $\bs{\ti
G}$ for~$(\ti X,\bs{\ti f})$.
\item[(d)] If $\bs G$ is extra data for $(X,\bs f)$, it has a {\it
restriction\/} $\bs G\vert_{\pd X}$ which is extra data for $(\pd
X,\bs f\vert_{\pd X})$.
\item[(e)] Let $\bs\si:\pd^2X\ra\pd^2X$ be the natural involution
of Definition \ref{kh2def19}. Suppose $X$ has corners (not
g-corners), and $\bs H$ is extra data for $(\pd X,\bs f\vert_{\pd
X})$. Then there should exist extra data $\bs G$ for $(X,\bs f)$
with $\bs G\vert_{\pd X}=\bs H$ if and only if $\bs H\vert_{\pd^2X}$
is invariant under~$\bs\si$.

If also $\Ga$ is a finite subgroup of $\Aut(X,\bs f)$, and $\bs H$
is invariant under $\Ga\vert_{\pd X}$, then we can choose $\bs G$ to
be $\Ga$-invariant.
\item[(f)] For the case of extra data for Kuranishi homology:
let $Y,Z$ be orbifolds, and $h:Y\ra Z$ a smooth map. Suppose $X$ is
a compact Kuranishi space and $\bs f:X\ra Y$ is strongly smooth. If
$\bs G$ is extra data for $(X,\bs f)$, then there should exist extra
data $h_*(\bs G)$ for $(X,h\ci\bs f)$. It should satisfy $(g\ci
h)_*(\bs G)=g_*(h_*(\bs G))$.
\item[(g)] For the case of extra data for Kuranishi cohomology:
let $Y,Z$ be orbifolds, and $h:Y\ra Z$ a smooth, proper map. Suppose
$X$ is a compact Kuranishi space and $\bs f:X\ra Z$ is a strong
submersion. If $\bs G$ is extra data for $(X,\bs f)$, then there
should exist extra data $h^*(\bs G)$ for $(Y\t_{h,Z,\bs
f}X,\bs\pi_Y)$. It should satisfy $(g\ci h)^*(\bs G)=h^*(g^*(\bs
G))$.
\item[(h)] (Mainly for Kuranishi cohomology). Let $X_1,X_2,X_3$ be
compact Kuranishi spaces, $Y$ an orbifold, $\bs f_e:X_e\ra Y$ be
strongly smooth for $e=1,2,3$ with at least two $\bs f_e$ strong
submersions, and $\bs G_e$ be extra data for $(X_e,\bs f_e)$ for
$e=1,2,3$. Then we should construct extra data $\bs G_1\t_Y\bs G_2$
for $(X_1\t_{\bs f_1,Y,\bs f_2}X_2,\bs\pi_Y)$ from~$\bs G_1,\bs
G_2$.

This construction should be {\it symmetric}, in that it yields
isomorphic extra data for $(X_1\t_{\bs f_1,Y,\bs f_2}X_2,\bs\pi_Y)$
and $(X_2\t_{\bs f_2,Y,\bs f_1}X_1,\bs\pi_Y)$ under the natural
isomorphism $X_1\t_{\bs f_1,Y,\bs f_2}X_2\cong X_2\t_{\bs f_2,Y,\bs
f_1}X_1$. It should also be {\it associative}, in that it yields
isomorphic extra data for $\bigl((X_1\!\t_Y\!X_2)\!\t_Y\!X_3,
\bs\pi_Y\bigr)$ and $\bigl(X_1\!\t_Y\!(X_2\!\t_Y\!X_3),\bs\pi_Y
\bigr)$ under $(X_1\!\t_Y\!X_2)\!\t_Y\!X_3\!\cong\!X_1\!\t_Y\!
(X_2\!\t_Y\!X_3)$. For Kuranishi cohomology, products should be {\it
functorial for pullbacks} in part (g), $h^*(\bs G_1\t_Z\bs G_2)\cong
\bs h^*(\bs G_1)\t_Yh^*(\bs G_2)$.
\end{itemize}
\label{kh3pr}
\end{property}

We shall define {\it gauge-fixing data\/} $\bs G$ (for Kuranishi
homology), and {\it co-gauge-fixing data\/} $\bs C$ (for Kuranishi
cohomology), which satisfy all these, and also have other technical
features which will be important in Appendices \ref{khB} and
\ref{khC}. The name comes from physics, where one often has an
infinite-dimensional symmetry group called the gauge group, and one
breaks the symmetry by `fixing a gauge'. It was not easy to arrange
for all of Property \ref{kh3pr}(a)--(h) to hold at once.

\subsection{Definition of (co-)gauge-fixing data}
\label{kh31}

{\it Good coordinate systems\/} are convenient choices of finite
coverings of $X$ by Kuranishi neighbourhoods, due to Fukaya and Ono
\cite[Def.~6.1]{FuOn1}, \cite[Lem.~A1.11]{FOOO}.

\begin{dfn} Let $X$ be a compact Kuranishi space. A {\it good
coordinate system\/ $\bs I$ for\/} $X$ consists of a finite indexing
set $I$, a total order $\le$ on $I$, a family
$\bigl\{(V^i,E^i,s^i,\psi^i):i\in I\bigr\}$ of compatible Kuranishi
neighbourhoods on $X$ with $V^i\ne\es$ for all $i\in I$ and
$X=\bigcup_{i\in I}\Im\psi^i$, and for all $i,j\in I$ with $j\le i$
a triple $(V^{ij},\phi^{ij},\hat\phi^{ij})$, where $V^{ij}$ is an
open neighbourhood of $(\psi^j)^{-1}(\Im\psi^i)$ in $V^j$, and
$(\phi^{ij},\hat\phi^{ij})$ is a coordinate change from
$(V^{ij},E^j\vert_{V^{ij}},s^j\vert_{V^{ij}},\psi^j\vert_{V^{ij}})$
to $(V^i,E^i,s^i,\psi^i)$. When $i=j$ these should satisfy
$V^{ij}=V^i$, with $(\phi^{ii},\hat\phi^{ii})$ the identity map on
$(V^i,\ldots,\psi^i)$. When $i,j,k\in I$ with $k\le j\le i$ these
should satisfy $\phi^{ij}\ci\phi^{jk}=\phi^{ik}$ and
$\hat\phi^{ij}\ci\hat\phi^{jk}=\hat\phi^{ik}$ over
$(\phi^{jk})^{-1}(V^{ij})\cap V^{ik}$, and
\e
\phi^{ij}(V^{ij})\cap\phi^{ik}(V^{ik})=
\phi^{ik}\bigl((\phi^{jk})^{-1}(V^{ij})\cap V^{ik}\bigr).
\label{kh3eq1}
\e

Now let $Y$ be an orbifold and $\bs f:X\ra Y$ a strongly smooth map.
A {\it good coordinate system $\bs I$ for\/} $(X,\bs f)$ consists of
a good coordinate system as above, together with smooth $f^i:V^i\ra
Y$ representing $\bs f$, such that if $j\le i$ in $I$ then
$f^j\vert_{V^{ij}}\equiv f^i\ci\phi^{ij}$. If $\bs f$ is a {\it
strong submersion}, we take the $f^i$ submersions.
\label{kh3def1}
\end{dfn}

The point of good coordinate systems is that one often needs to
choose some geometric data (for instance, a transverse multisection,
as in as \cite[Th.~6.4]{FuOn1}) on each of a system of Kuranishi
neighbourhoods satisfying compatibility conditions on the overlaps.
With a good coordinate system one can make these choices on each
$(V^i,\ldots,\psi^i)$ in the order $\le$ on $I$ by induction,
imposing compatibility conditions on the $V^{ij}$. Equation
\eq{kh3eq1} was added by the author; without it, in inductive
arguments when making the choice on $(V^i,\ldots,\psi^i)$ there may
be problems at points in $V^i$ lying in the left hand side of
\eq{kh3eq1} but not the right.

We shall now define three variations on the notion of good
coordinate system, which we call {\it very good}, {\it really
good\/} and {\it excellent coordinate systems}. A {\it very good
coordinate system\/} is a good coordinate system in which there is
at most one Kuranishi neighbourhood $(V^i,\ldots,\psi^i)$ of each
dimension $\dim V^i$, so we use the dimension $\dim V^i$ as the
index~$i$.

\begin{dfn} Let $X$ be a compact Kuranishi space, $Y$ an orbifold
and $\bs f:X\ra Y$ a strongly smooth map. A good coordinate system
$\bs I$ for $(X,\bs f)$ is called {\it very good\/} if
$I\subset\N=\{0,1,2,\ldots\}$, and the order $\le$ on $I$ is the
restriction of $\le$ on $\N$, and $\dim V^i=i$ for all $i\in I$.
\label{kh3def2}
\end{dfn}

Note that it is necessary that the order $\le$ on $I$ is the
restriction of $\le$ on $\N$, since coordinate changes
$(\phi^{ij},\hat\phi^{ij})$ from $V^{ij}\subseteq V^j$ to $V^i$ can
exist if $\dim V^j\le\dim V^i$, that is, if $j\le i$, but not if
$\dim V^j>\dim V^i$. A {\it really good coordinate system\/} is a
very good coordinate system $\bs I$ with some extra partition of
unity data $\bs\eta$. We will use $\bs\eta$ for defining really good
coordinate systems on fibre products in~\S\ref{kh38}.

\begin{dfn} Let $X$ be a compact Kuranishi space, $Y$ an orbifold,
and $\bs f:X\ra Y$ a strongly smooth map. A {\it really good
coordinate system\/} $(\bs I,\bs\eta)$ for $(X,\bs f)$ is a very
good coordinate system $\bs I=\bigl(I,(V^i,\ldots,\psi^i),f^i:i\in
I,\ldots\bigr)$ for $(X,\bs f)$ together with data
$\bs\eta=(\eta_i:i\in I$, $\eta_i^j:i,j\in I)$ satisfying:
\begin{itemize}
\setlength{\itemsep}{0pt}
\setlength{\parsep}{0pt}
\item[(i)] $\eta_i:X\ra[0,1]$ is a continuous function for each $i\in
I$, with $\sum_{i\in I}\eta_i\equiv 1$, and $\eta_i>0$ on
$\Im\psi^i$. If $j<i$ and $p\in X$ with $\eta_i(p)>0$,
$\eta_j(p)>0$ and $p\in\Im\psi^i$ then~$p\in\Im\psi^j$.
\item[(ii)] $\eta_i^j:V^j\ra[0,1]$ is a continuous function for
all $i,j\in I$, with $\sum_{i\in I}\eta_i^j\equiv 1$ on $V^j$,
and $\eta_i^i>0$ on $V^i$. If $j\le i$ then $\phi^{ij}(V^{ij})$
is a {\it closed\/} subset of $\{v\in V^i:\eta_j^i(v)>0\}$. If
$k<j<i$ and $v\in V^k$ with $\eta_i^k(v)>0$, $\eta_j^k(v)>0$ and
$v\in V^{ik}$ then~$v\in V^{jk}$.
\item[(iii)] $\eta_i^j\vert_{(s^j)^{-1}(0)}\equiv\eta_i\ci\psi^j$
for all~$i,j\in I$.
\item[(iv)] If $i,j,k\in I$ with $k\le j$ then
$\eta_i^k\vert_{V^{jk}}\equiv\eta_i^j\ci\phi^{jk}$.
\end{itemize}

We also impose two extra conditions on the~$V^i,f^i$:
\begin{itemize}
\setlength{\itemsep}{0pt}
\setlength{\parsep}{0pt}
\item[(v)] There should exist $M\ge 0$ such that for each $i\in I$,
each $v\in V^i$ lies in at most $M$ local boundary components,
in the sense of Definition \ref{kh2def3}. Here if the orbifold
with g-corners $V^i$ near $v$ is locally modelled on $U/\Ga$
near $\Ga u$ as in \S\ref{kh22}, where $\Ga$ is a finite group
acting linearly on $\R^n$ and $U$ is a $\Ga$-invariant region
with g-corners in $\R^n$, and $u\in U$, then we define the {\it
number of local boundary components} of $V^i$ at $v$ to be the
number of local boundary components of $U$ at $u$, even if some
of these local boundary components are identified by the action
of~$\Ga$.

Then the natural maps $\io:\pd V^i\ra V^i$ and
$\io:E^i\vert_{\pd V^i}=\pd E^i\ra E^i$ satisfy
$\md{\io^{-1}(v)}\le M$ and $\md{\io^{-1}(e)}\le M$ for all
$v\in V^i$ and~$e\in E^i$.
\item[(vi)] There should exist a compact subset $T\subseteq Y$ with
$f^i(V^i)\subseteq T$ for all $i\in I$.
\end{itemize}
If the $V^i$ have corners (not g-corners) then (v) holds
with~$M=\max\{i\in I\}$.
\label{kh3def3}
\end{dfn}

{\it Excellent coordinate systems\/} are really good coordinate
systems satisfying an extra condition, which we will need in
\S\ref{kh33} to ensure $\Aut(X,\bs f,\bs G)$ is finite.

\begin{dfn} Let $({\bs I},\bs\eta)$ be a really good coordinate
system for $(X,\bs f)$, and suppose $j\in I$, $l\ge 0$ and $W$ is a
connected component of $\pd^lV^j$. We call $W$ {\it unnecessary\/}
if $W\cap\pd^lV^{ij}=\es$ for all $j<i\in I$ and $\Im\psi^j\vert_W
\subseteq\bigcup_{j\ne i\in I}\Im\psi^i\vert_{\pd^lV^i}\subseteq
\pd^lX$. We call $({\bs I},\bs\eta)$ {\it excellent\/} if $\pd^lV^j$
has no unnecessary connected components for all $j\in I$ and~$l\ge
0$.
\label{kh3def4}
\end{dfn}

Excellent coordinate systems can be split into {\it connected
components}, a fact which will be important in~\S\ref{kh39}.

\begin{dfn} Let $Y$ be an orbifold, $X$ a compact Kuranishi space,
$\bs f:X\ra Y$ be strongly smooth, and $(\bs I,\bs\eta)$ an
excellent coordinate system for $(X,\bs f)$, with $\bs
I=\bigl(I,(V^i,\ldots,\psi^i):i\in I,\ldots\bigr)$. A {\it
splitting} of the quadruple $(X,\bs f,\bs I,\bs\eta)$ is a
decomposition $X=X_1\amalg\cdots\amalg X_n$ of $X$ as a disjoint
union of compact, nonempty Kuranishi spaces $X_1,\ldots,X_n$, such
that for each $i\in I$ we may write $V^i=V^i_1\amalg\cdots\amalg
V^i_n$ for open and closed subsets $V^i_1,\ldots,V^i_n$ of $V^i$,
such that $\Im\psi^i\vert_{V^i_a}\subseteq X_a$ for $a=1,\ldots,n$.
Note that as $(\bs I,\bs\eta)$ is excellent, every connected
component of $V^i$ intersects $\Im\psi^i$, so $V^i_1,\ldots,V^i_n$
are uniquely determined by~$X_1,\ldots,X_n$.

Write $\bs f_a=\bs f\vert_{X_a}$, so that $\bs f_a:X_a\ra Y$ is
strongly smooth. Define subsets $I_1,\ldots,I_n\subseteq I$ by
$I_a=\{i\in I:V^i_a\ne\es\}$ for $a=1,\ldots,n$. For $i\in I_a$ let
$E^i_a,s^i_a, \psi^i_a,f^i_a$ be the restrictions of
$E^i,s^i,\psi^i,f^i$ to $V^i_a$. For $i,j\in I_a$ set
$\eta_{i,a}=\eta_i\vert_{X_a}$ and
$\eta_{i,a}^j=\eta_i^j\vert_{V^j_a}$, and if $j\le i$ write
$\smash{V^{ij}_a,\phi^{ij}_a,\hat\phi^{ij}_a}$ for the restrictions
of $\smash{V^{ij},\phi^{ij},\hat\phi^{ij}}$ to
$\smash{V^j_a,E^j_a}$. It is easy to show this defines an {\it
excellent coordinate system\/} $(\bs I_a,\bs\eta_a)=(\bs
I,\bs\eta)\vert_{X_a}$ for~$(X_a,\bs f_a)$.

Call $(X,\bs f,\bs I,\bs\eta)$ {\it connected\/} if it does not
admit a splitting $X=X_1\amalg X_2$ as above. Since $X_1,X_2$ are
open and closed in $X$ and nonempty, if $X$ is connected as a
topological space then $(X,\bs f,\bs I,\bs\eta)$ is connected, but
the converse is false.
\label{kh3def5}
\end{dfn}

\begin{lem} Let\/ $Y$ be an orbifold, $X$ a compact Kuranishi space,
$\bs f:X\ra Y$ strongly smooth, and\/ $(\bs I,\bs\eta)$ an excellent
coordinate system for $(X,\bs f)$. Then there is a splitting
$X=X_1\amalg X_2\amalg\cdots\amalg X_n$ of\/ $(X,\bs f,\bs
I,\bs\eta),$ unique up to the order of\/ $X_1,\ldots,X_n,$ such
that\/ $\bigl(X_a,\bs f_a,\bs I_a,\bs\eta_a\bigr)$ is connected for
all\/ $a=1,\ldots,n$. We call $X_1,\ldots,X_n$ the
\begin{bfseries}connected components\end{bfseries} of\/~$(X,\bs
f,\bs I,\bs\eta)$.
\label{kh3lem}
\end{lem}

\begin{proof} Write $\bs I=\bigl(I,(V^i,\ldots,\psi^i):i\in
I,\ldots\bigr)$. Then for each $p\in X$, there exists $i\in I$ with
$p\in\Im\psi^i$. Let $\ti V^i_p$ be the connected component of $V^i$
containing $(\psi^i)^{-1}(p)$. Then $U_p^i=\Im\psi^i\vert_{\ti
V^i_p}$ is an open neighbourhood of $p$ in $X$. Such $U_p^i$ for all
$p\in X$ form an open cover for $X$, so as $X$ is compact we can
choose a finite subcover $U_{p_1}^{i_1},\ldots,U_{p_m}^{i_m}$
for~$X$.

Suppose $X=X_1\amalg\cdots\amalg X_n$ is a splitting of $(X,\bs
f,\bs I,\bs\eta)$. By Definition \ref{kh3def5}, if $p_b\in X_a$ with
$U_{p_b}^{i_b}=\Im\psi^{i_b}\vert_{\ti V^{i_b}_{p_b}}$ then $\ti
V^{i_b}_{p_b}\subseteq V^{i_b}_a$, so that $U_{p_b}^{i_b}\subseteq
X_a$. As $X_a\ne\es$ and $X=\bigcup_{b=1}^mU_{p_b}^{i_b}$, it
follows that for each $a=1,\ldots,n$ we have $p_b\in X_a$ for some
$b$. Therefore $n\le m$, and the number of pieces we can split
$(X,\bs f,\bs I,\bs\eta)$ into using Definition \ref{kh3def5} is
bounded.

Choose such a splitting $X=X_1\amalg\cdots\amalg X_n$ with $n$
largest. Then each $\bigl(X_a,\bs f_a,\bs I_a,\bs\eta_a\bigr)$ is
{\it connected\/}, since otherwise by splitting $\bigl(X_a,\bs
f_a,\bs I_a,\bs\eta_a\bigr)$ into 2 pieces we could split $X$ into
$n+1$ pieces, a contradiction. If $X=X_1'\amalg\cdots\amalg X_n'$ is
another such splitting, then $X_a=\coprod_{b=1,\ldots,n:X_a\cap
X_b'\ne\es}X_a\cap X_b'$ is a splitting of $\bigl(X_a,\bs f_a,\bs
I_a,\bs\eta_a\bigr)$. Since it is connected, we must have $X_a\cap
X_b'\ne\es$ for precisely one $b=1,\ldots,n$, with $X_a\cap
X_b'=X_a=X_b'$. Therefore $X_1',\ldots,X_n'$ are a reordering of
$X_1,\ldots,X_n$, and the decomposition is unique up to order.
\end{proof}

We can now define {\it gauge-fixing data} and {\it co-gauge-fixing
data}.

\begin{dfn} Let the symmetric group $S_k$ of permutations of
$\{1,\ldots,k\}$ act on $\R^k$ by $\si:(x_1,\ldots,x_k)\mapsto
(x_{\si(1)},\ldots,x_{\si(k)})$. Define $P=\coprod_{k=0}^\iy
\R^k/S_k$. For $n=0,1,2,\ldots$, define $P_n\subset P$ by $P_n=
\coprod_{k=0}^n\R^k/S_k$. Here $\R^0$ is a single point $(\es)$, the
ordered list of 0 elements of $\R$, and $S_0=\{1\}$ acts trivially
upon it, so $(\es)$ is a special element of $P$. Define a {\it
multiplication} $\mu:P\t P\ra P$ by
\e
\mu\bigl(S_k(x_1,\ldots,x_k),S_l(y_1,\ldots,y_l)\bigr)=S_{k+l}
(x_1,\ldots,x_k,y_1,\ldots,y_l),
\label{kh3eq2}
\e
for $k,l\ge 0$ and $(x_1,\ldots,x_k)\in\R^k$, $(y_1,\ldots,y_l)\in
\R^l$. Then $\mu$ is {\it commutative} and {\it associative}, with
{\it identity} $(\es)$, and also a {\it finite map}, since for any
point $S_n(x_1,\ldots,x_n)$ in $P$,
$\mu^{-1}\bigl(S_n(x_1,\ldots,x_n)\bigr)$ is at most $2^n$ points in
$P\t P$. We will use $\mu$ in defining fibre products of
(co-)gauge-fixing data in~\S\ref{kh38}.

Let $X$ be a compact Kuranishi space, $Y$ an orbifold and $\bs
f:X\ra Y$ strongly smooth. A set of {\it gauge-fixing data\/} $\bs
G$ for $(X,\bs f)$ consists of an excellent coordinate system $(\bs
I,\bs\eta)$ for $(X,\bs f)$, where $\bs I=\bigl(I,(V^i,
E^i,s^i,\psi^i):i\in I,\ldots\bigr)$, and maps $G^i:E^i\ra
P_n\subset P$ for some $n\gg 0$ and all $i\in I$, which are {\it
globally finite}, that is, there exists $N\ge 0$ such that
$\bmd{(G^i)^{-1}(p)}\le N$ for all~$p\in P$.

By calling $(X,\bs f,\bs G)$ a {\it triple}, we will usually mean
that $X$ is a compact Kuranishi space (possibly oriented), $\bs
f:X\ra Y$ is strongly smooth, and $\bs G$ is gauge-fixing data for
$(X,\bs f)$. A {\it splitting} of a triple $(X,\bs f,\bs G)$ is a
splitting $X=X_1\amalg\cdots\amalg X_n$ of $(X,\bs f,\bs I,\bs\eta)$
in the sense of Definition \ref{kh3def5}, where $\bs G=(\bs I,
\bs\eta,G^i:i\in I)$. Then Definition \ref{kh3def5} gives excellent
coordinate systems $(\bs I_a,\bs\eta_a)$ for $(X_a,\bs f_a)$, and
setting $G_a^i=G^i\vert_{E^i_a}$ for $i\in I_a$, it is easy to see
that $\bs G\vert_{X_a}=\bs G_a=(\bs I_a,\bs\eta_a,G^i_a:i\in I_a)$
is gauge-fixing data for $(X_a,\bs f_a)$.

Call $(X,\bs f,\bs G)$ {\it connected\/} if $(X,\bs f,\bs
I,\bs\eta)$ is connected in the sense of Definition \ref{kh3def5}.
Then Lemma \ref{kh3lem} implies that any triple $(X,\bs f,\bs G)$
has a splitting into {\it connected components} $(X_1,\bs f_1,\bs
G_1),\ldots,(X_n,\bs f_n,\bs G_n)$, unique up to order.

Let $X,\ti X$ be compact Kuranishi spaces, $\bs f:X\ra Y$, $\bs{\ti
f}:\ti X\ra Y$ be strongly smooth, and $\bs G,\bs{\ti G}$ be
gauge-fixing data for $(X,\bs f),(\ti X,\bs{\ti f})$. Suppose the
indexing sets $I,\ti I$ of $\bs I,\bs{\ti I}$ in $\bs G,\bs{\ti G}$
are equal, $I=\ti I$. An {\it isomorphism\/} $(\bs a,\bs b):(X,\bs
f,\bs G)\ra(\ti X,\bs{\ti f},\bs{\ti G})$ consists of a strong
diffeomorphism $\bs a:X\ra\ti X$ and a collection $\bs b$ of
isomorphisms $(b^i,\hat b^i):(V^i,\ldots,\psi^i)\ra(\ti
V^i,\ldots,\ti\psi^i)$ for $i\in I=\ti I$ lifting $\bs a:X\ra\ti X$,
satisfying:
\begin{itemize}
\setlength{\itemsep}{0pt}
\setlength{\parsep}{0pt}
\item[(a)] $\ti f^i\ci b^i\equiv f^i$ for all $i\in I$;
\item[(b)] $\ti G^i\ci\hat b^i\equiv G^i$ for all $i\in
I$; and
\item[(c)] If $j\le i$ in $I$ then $b^j(V^{ij})=\ti V^{ij}$, and
the following diagram commutes:
\begin{equation*}
\xymatrix@C=90pt@R=25pt{ (V^{ij},\ldots,\psi^j\vert_{V^{ij}})
\ar[r]_{(b^j,\hat b^j)\vert_{V^{ij}}}
\ar[d]^{(\phi^{ij},\hat\phi^{ij})} & (\ti
V^{ij},\ldots,\ti\psi^{j}\vert_{\ti V^{ij}})
\ar[d]_{(\ti\phi^{ij},\hat{\ti\phi}^{ij})}\\
(V^i,\ldots,\psi^i) \ar[r]^{(b^i,\hat b^i)} & (\ti
V^i,\ldots,\ti\psi^i). }
\end{equation*}
\end{itemize}
Isomorphisms compose, and have inverses which are isomorphisms, in
the obvious way. Write $\Aut(X,\bs f,\bs G)$ for the {\it
automorphism group\/} of $(X,\bs f,\bs G)$, that is, the group of
isomorphisms~$(\bs a,\bs b):(X,\bs f,\bs G)\ra(X,\bs f,\bs G)$.
\label{kh3def6}
\end{dfn}

\begin{dfn} Let $X$ be a compact Kuranishi space, $Y$ an orbifold
and $\bs f:X\ra Y$ a strong submersion. A set of {\it
co-gauge-fixing data\/} $\bs C$ for $(X,\bs f)$ consists of an
excellent coordinate system $(\bs I,\bs\eta)$ for $(X,\bs f)$ in the
sense of \S\ref{kh31}, where $\bs I=\bigl(I,(V^i,E^i,s^i,\psi^i):
i\in I,f^i:i\in I,\ldots\bigr)$, such that $f^i:V^i\ra Y$ is a
submersion for all $i\in I$, and maps $C^i:E^i\ra P_n\subset P$ for
some $n\gg 0$ and all $i\in I$, such that $C^i\t(f^i\ci\pi^i):E^i\ra
P\t Y$ is a {\it globally finite\/} map for all $i\in I$, as in
Definition \ref{kh3def6}. {\it Isomorphisms} $(\bs a,\bs b):(X,\bs
f,\bs C)\ra(\ti X,\bs{\ti f},\bs{\ti C})$ and $\Aut(X,\bs f,\bs C)$
are defined as in Definition \ref{kh3def6}, but with $\ti C^i\ci\hat
b^i\equiv C^i$ in~(b).

By calling $(X,\bs f,\bs C)$ a {\it triple}, we will usually mean
that $X$ is a compact Kuranishi space, $\bs f:X\ra Y$ a strong
submersion (possibly cooriented), and $\bs C$ co-gauge-fixing data
for $(X,\bs f)$. A {\it splitting} of $(X,\bs f,\bs C)$ is a
splitting $X=X_1\amalg\cdots\amalg X_n$ of $(X,\bs f,\bs
I,\bs\eta)$, where $\bs C=(\bs I, \bs\eta,C^i:i\in I)$. Then
Definition \ref{kh3def5} gives excellent $(\bs I_a,\bs\eta_a)$ for
$(X_a,\bs f_a)$, and $\bs C\vert_{X_a}=\bs C_a=(\bs
I_a,\bs\eta_a,C^i_a:i\in I_a)$ is co-gauge-fixing data for $(X_a,\bs
f_a)$, where $C_a^i=C^i\vert_{E^i_a}$ for~$i\in I_a$.

Call $(X,\bs f,\bs C)$ {\it connected\/} if $(X,\bs f,\bs
I,\bs\eta)$ is connected in the sense of Definition \ref{kh3def5}.
Then by Lemma \ref{kh3lem}, any $(X,\bs f,\bs C)$ has a splitting
into {\it connected components} $(X_1,\bs f_1,\bs
C_1),\ldots,(X_n,\bs f_n,\bs C_n)$, unique up to order.
\label{kh3def7}
\end{dfn}

\begin{rem} Note that we do {\it not\/} require the maps $G^i:E^i\ra
P$, $C^i:E^i\ra P$ to be continuous, or smooth, etc., they are just
arbitrary maps. Nor do we impose compatibility conditions between
$G^i,G^j$ under coordinate transformations
$(\phi^{ij},\hat\phi^{ij})$. Their only purpose is to ensure that
$\Aut(X,\bs f,\bs G),\Aut(X,\bs f,\bs C)$ are finite, and this will
follow from $G^i,C^i\t(f^i\ci\pi^i)$ globally finite.

To see why gauge-fixing data is natural for homology, and
co-gauge-fixing data for cohomology, read \S\ref{kh37} on
pushforwards and pullbacks.
\label{kh3rem1}
\end{rem}

For Property \ref{kh3pr}(e) we define (co-)gauge-fixing data {\it
with corners}.

\begin{dfn} Let $X$ be a compact Kuranishi space, and $\bs f:X\ra Y$
be strongly smooth. We say that a good or very good coordinate
system $\bs I$, or a really good or excellent coordinate system
$(\bs I,\bs\eta)$, or gauge-fixing data $\bs G=(\bs
I,\bs\eta,G^i:i\in I)$, or co-gauge-fixing data $\bs C=(\bs
I,\bs\eta,C^i:i\in I)$, for $(X,\bs f)$, are {\it without boundary},
or {\it with boundary}, or {\it with corners}, if the orbifolds
$V^i$ for all $i\in I$ are without boundary, or with boundary, or
with corners, respectively. These imply that the Kuranishi space $X$
is without boundary, or with boundary, or with corners,
respectively.
\label{kh3def8}
\end{dfn}

Sections \ref{kh32}--\ref{kh38} prove Property \ref{kh3pr}(a)--(h)
for (co-)gauge-fixing data. Since the two are so similar, we
generally give proofs only for gauge-fixing data.

\subsection{Existence of (co-)gauge-fixing data}
\label{kh32}

To show that (co-)gauge-fixing data exists for every pair $(X,\bs
f)$, we begin by showing that $(X,\bs f)$ admits good, very good,
really good and excellent coordinate systems. Fukaya and Ono prove
\cite[Lem.~6.3]{FuOn1}, \cite[Lem.~A1.11]{FOOO}:

\begin{prop} Let\/ $X$ be a compact Kuranishi space, $Y$ an
orbifold, $\bs f:X\ra Y$ a strongly smooth map, and\/
$\{U_\al:\al\in A\}$ an open cover for\/ $X$. Then there exist good
coordinate systems $\bs I$ for\/ $X$ and\/ $(X,\bs f)$ such that for
each\/ $i\in I$ we have\/ $\Im\psi^i\subseteq U_\al$ for
some\/~$\al\in A$.
\label{kh3prop1}
\end{prop}

\begin{prop} Let\/ $X$ be a compact Kuranishi space, $Y$ an
orbifold, and\/ $\bs f:X\ra Y$ a strongly smooth map. Then there
exists a very good coordinate system $\bs I$ for\/~$(X,\bs f)$.
\label{kh3prop2}
\end{prop}

\begin{proof} By Proposition \ref{kh3prop1} there exists a good
coordinate system $\bs I=\bigl(I,\ab\le,\ab(V^i,\ab\ldots,
\psi^i):i\in I,\ldots\bigr)$ for $(X,\bs f)$. Define $\check
I=\{\dim V^i:i\in I\}$. For each $k\in\check I$, we shall define a
Kuranishi neighbourhood $(\check V^k,\ldots,\check\psi^k)$ on $X$ by
gluing together the $(V^i,\ldots,\psi^i)$ for all $i\in I$ with
$\dim V^i=k$. We first explain how to glue together two such
Kuranishi neighbourhoods. Let $i,j\in I$ with $i\ne j$ and $\dim
V^i=\dim V^j=k$. Either $i\le j$ or $j\le i$ as $\le$ is a total
order; suppose $j\le i$. Then we have a triple $(V^{ij},\phi^{ij},
\hat\phi^{ij})$ by Definition~\ref{kh3def1}.

As $\dim V^i=\dim V^j$, it follows that the coordinate change
$(\phi^{ij},\hat\phi^{ij})$ is an {\it isomorphism\/} from $(V^{ij},
E^j\vert_{V^{ij}},s^j\vert_{V^{ij}},\psi^j\vert_{V^{ij}})$ to its
image in $(V^i,\ldots,\psi^i)$. In particular, $\phi^{ij}:V^{ij}\ra
V^i$ is a diffeomorphism $V^{ij}\ra\Im\phi^{ij}$. The obvious thing
to do now is to glue $V^i$ and $V^j$ together using $\phi^{ij}$,
that is, to form the topological space $(V^i\amalg V^j)/\!\sim$,
where $\sim$ is the equivalence relation on $V^i\amalg V^j$
generated by $\phi^{ij}$, with $v\sim v'$ if $v=v'$, or $v\in
V^{ij}$ and $v'=\phi^{ij}(v)\in V^i$, or vice versa. However, there
is a problem: since $V^{ij}$ and $\Im\phi^{ij}$ may not be closed in
$V^j,V^i$, this $(V^i\amalg V^j)/\!\sim$ may not be {\it Hausdorff},
and so not an orbifold.

Our solution is to make $V^i,V^j$ a little smaller, that is, we
replace them by open sets $\ti V^i\subset V^i$ and $\ti V^j\subset
V^j$, such that $(\ti V^i\amalg\ti V^j)/\!\sim$ is Hausdorff, and so
an orbifold. If $(V^i\amalg V^j)/\!\sim$ is not Hausdorff then there
exist sequences $(v_n)_{n=1}^\iy$ in $V^{ij}$ so that $v_n\ra v^j$
in $V^j\sm V^{ij}$ as $n\ra\iy$ and $\phi^{ij}(v_n)\ra v^i$ in
$V^i\sm\Im\phi^{ij}$ as $n\ra\iy$. Then $v^i,v^j$ are distinct
points in $(V^i\amalg V^j)/\!\sim$ which do not admit disjoint open
neighbourhoods, contradicting Hausdorffness, since any
neighbourhoods contain $v_n$ for~$n\gg 0$.

To make $(\ti V^i\amalg\ti V^j)/\!\sim$ Hausdorff, we have to ensure
that for all such $(v_n)_{n=1}^\iy$ and $v^i,v^j$, either
$v^i\notin\ti V^i$ or $v^j\notin\ti V^j$. Loosely speaking, the
`boundaries' of $V^{ij}$ in $\ti V^i$ and $\ti V^j$ may not
intersect. Note however that we cannot make $\ti V^i,\ti V^j$ too
small, as we need to preserve the condition $\bigcup_{i\in
I}\Im\psi^i=X$, that is, the choices of $\ti V^i$ for all $i\in I$
must satisfy $\bigcup_{i\in I}\Im\psi^i\vert_{\ti V^i}=X$. These
conditions are consistent, since gluing $(s^i)^{-1}(0)$ and
$(s^j)^{-1}(0)$ together using $\phi^{ij}$ yields an open subset
$\Im\psi^i\cup\Im\psi^j$ of $X$, which is Hausdorff, so the points
$v^i,v^j$ which need to be excluded cannot lie in $(s^i)^{-1}(0),
(s^j)^{-1}(0)$, the domains of~$\psi^i,\psi^j$.

We need to choose open $\ti V^i\subset V^i$ for all $i\in I$, which
are small enough that whenever $j\le i\in I$ with $\dim V^i=\dim
V^j$ then $(\ti V^i\amalg\ti V^j)/\!\sim$ is Hausdorff, and which
are large enough that $\bigcup_{i\in I}\Im\psi^i\vert_{\ti V^i}=X$.
This is possible. One way to do it is to use ideas on {\it really
good coordinate systems\/} below. Follow the proof of Proposition
\ref{kh3prop3} to choose auxiliary partitions of unity $\eta_i:i\in
I$ on $X$ and $\eta_i^j:i\in I$ on $V^j$ satisfying conditions
(a)--(d) of that proof, possibly making the $V^j$ smaller along the
way, and define $\ti V^i=\{v\in V^i:\eta_i^i(v)>0\}$. As in
Proposition \ref{kh3prop3} the resulting $\ti V^j$ have partitions
of unity $\eta_i^j\vert_{\ti V^j}:i\in I$ satisfying Definition
\ref{kh3def3} with $\bigcup_{i\in I}\Im\psi^i\vert_{\ti V^i}=X$, and
by Remark \ref{kh3rem2} this ensures that $(\ti V^i\amalg\ti
V^j)/\!\sim$ is Hausdorff.

This explains how to glue two Kuranishi neighbourhoods
$(V^i,\ldots,\psi^i)$, $(V^j,\ab\ldots,\ab\psi^j)$ for $i\ne j\in I$
with $\dim V^i=\dim V^j=k$. To make $\check V^k$ we repeat this
construction finitely many times to glue together $V^i$ for all
$i\in I$ with $\dim V^i=k$, so that $\check V^k=\bigl(\coprod_{i\in
I:\dim V^i=k}\ti V^i\bigr)\big/\!\sim$, where $\sim$ is the
equivalence relation generated by $\phi^{ij}$ for all $j\le i\in I$
with $\dim V^i=\dim V^j=k$, and the $\ti V^i$ are chosen small
enough to make $\check V^k$ Hausdorff, and thus an orbifold.

Since the $(\phi^{ij},\hat{\phi}{}^{ij})$ are isomorphisms with
their images, we can also glue the restrictions of $E^i,s^i,\psi^i$
to $\ti V^i$ together to make $\check E^k,\check s^k,\check\psi^k$
on $\check V^k$. That is, we have a Kuranishi neighbourhood $(\check
V^k,\check E^k,\check s^k,\check\psi^k)$ on $X$ with given
isomorphisms $(\check E^k,\check s^k,\check\psi^k)\vert_{\ti
V^i}\cong(E^i,s^i,\psi^i)\vert_{\ti V^i}$, for all $i\in I$ with
$\dim V^i=k$, regarding $\ti V^i$ as an open subset of both $\check
V^k$ and $V^i$. This gives the Kuranishi neighbourhoods $(\check
V^k,\check E^k,\check s^k,\check\psi^k)$ for $k\in\check I$. Since
$\bigcup_{i\in I}\Im\psi^i\vert_{\ti V^i}=X$, we have~$\bigcup_{k\in
\check I}\Im\check\psi^k=X$.

We have smooth $f^i: V^i\ra Y$ for $i\in I$ representing $\bs f$.
The compatibilities $f^i\ci\phi^{ij}\equiv f^j\vert_{V^{ij}}$ imply
that we may glue these to give smooth maps $\check f^k:\check V^k\ra
Y$ representing $\bs f$, with $\check f^k\vert_{\ti
V^i}=f^i\vert_{\ti V^i}$, for all $i\in I$ with~$\dim V^i=k$.

If $k,l\in\check I$ with $l<k$ then for all $i\in I$ with $\dim
V^i=k$ and $j\in I$ with $\dim V^j=l$ we have $j\le i$, so
$(V^{ij},\phi^{ij},\hat\phi^{ij})$ is given. Define an open $\check
V^{kl}\subseteq\check V^l$ by $\check V^{kl}=\ts\bigcup_{i,j\in
I:\dim V^i=k,\;\dim V^j=l}\bigl( V^{ij}\cap\ti
V^j\cap(\phi^{ij})^{-1}(\ti V^i)\bigr)$. There is then a unique
coordinate transformation
$(\check\phi^{kl},\smash{\hat{\check\phi}{}^{kl}}): \bigl(\check
V^{kl},\check E^l\vert_{\check V^{kl}},\check s^l\vert_{\check
V^{kl}},\ab\check \psi^l\vert_{V^{kl}}\bigr)\ra(\check
V^k,\ldots,\check\psi^k)$ which restricts to
$(\phi^{ij},\smash{\hat{\phi}{}^{ij}})$ on $V^{ij}\cap \ti
V^j\cap(\phi^{ij})^{-1}(\ti V^i)$ for all $i,j\in I$ with $\dim
V^i=k$ and $\dim V^j=l$. It is now easy to check that all this data
$\bs{\check I}=\bigl(\check I,(\check V^k,\ldots,\check
\psi^k),\check f^k:k\in\check I,\ldots\bigr)$ satisfies the
conditions for a very good coordinate system for $(X,\bs f)$. This
completes the proof.
\end{proof}

\begin{rem} In the proof of Proposition \ref{kh3prop2} we glued
topological spaces $V^j,V^i$ together on $V^{ij}\subseteq V^j$ and
$\phi^{ij}(V^{ij})\subseteq V^i$ using $\phi^{ij}$, and we had a
problem because the resulting topological space $V^i\amalg
V^j/\!\sim$ may not be {\it Hausdorff}. The existence of functions
$\eta_i^j$ satisfying Definition \ref{kh3def3}(ii),(iv) eliminates
this problem, so that $(\coprod_{i\in I}V^i)/\!\sim$ is a Hausdorff
topological space. This is useful, as it excludes some bad behaviour
in constructions where we choose data over each $V^i$ with
compatibilities on the $V^{ij}$, as we will do many times later.

To see this, suppose $(V^i\amalg V^j)/\!\sim$ is not Hausdorff. Then
there exists a sequence $(v_n)_{n=1}^\iy$ in $V^{ij}$ such that
$\phi^{ij}(v_n)\ra v^i$ in $V^i\sm\phi^{ij}(V^{ij})$ and $v_n\ra
v^j$ in $V^j\sm V^{ij}$ as $n\ra\iy$. Since $\eta_j^i,\eta_j^j$ are
continuous functions on $V^i,V^j$ we have $\eta_j^i\ci\phi^{ij}(v_n)
\ra\eta_j^i(v^i)$ and $\eta_j^j(v_n)\ra\eta_j^j(v^j)$ as $n\ra\iy$.
Definition \ref{kh3def3}(iv) implies that $\eta_j^i\ci\phi^{ij}
(v_n)=\eta_j^j(v_n)$ for all $n$, so $\eta_j^i(v^i)=\eta_j^j(v^j)$.
But Definition \ref{kh3def3}(ii) gives $\eta_j^j(v^j)>0$, and as
$\phi^{ij}(V^{ij})$ is a closed in $\{v\in V^i:\eta_j^i(v)>0\}$ and
$v^i\in\ov{\phi^{ij}(V^{ij})}\sm \phi^{ij}(V^{ij})$ we deduce that
$\eta_j^i(v^i)=0$, a contradiction.
\label{kh3rem2}
\end{rem}

\begin{prop} Let\/ $X$ be a compact Kuranishi space, $Y$ an
orbifold, and\/ $\bs f:X\ra Y$ a strongly smooth map. Then there
exists a really good coordinate system $(\bs I,\bs\eta)$
for\/~$(X,\bs f)$.
\label{kh3prop3}
\end{prop}

\begin{proof} By Proposition \ref{kh3prop2} we can choose a very
good coordinate system $\bs I$ for $(X,\bs f)$. Then
$\{\Im\psi^i:i\in I\}$ is an open cover for $X$. By general topology
we may choose a partition of unity $\{\eta_i:i\in I\}$ for $X$
subordinate to $\{\Im\psi^i:i\in I\}$, that is, $\eta_i:X\ra[0,1]$
is continuous, and $\sum_{i\in I}\eta_i\equiv 1$, and the closure of
the support of $\eta_i$ is contained in $\Im\psi^i$, that
is,~$\overline{\{p\in X:\eta_i(p)>0\}}\subseteq\Im\psi^i$.

Next, for each $j\in I$ we choose a partition of unity
$\{\eta_i^j:i\in I\}$ on $V^j$, by induction in the order $\le$ on
$I$, satisfying:
\begin{itemize}
\setlength{\itemsep}{0pt}
\setlength{\parsep}{0pt}
\item[(a)] $\eta_i^j:V^j\ra[0,1]$ is a continuous function for all
$i,j\in I$, with~$\sum_{i\in I}\eta_i^j\equiv 1$.
\item[(b)] If $j\le i$ in $I$ then $\supp(\eta_i^j)\subseteq V^{ij}$
and $\phi^{ij}\bigl(V^{ij}\cap \ov{\supp(\eta_j^j)}\bigr)$ is closed
in~$V^i$.
\item[(c)] $\eta_i^j\vert_{(s^j)^{-1}(0)}\equiv\eta_i\ci\psi^j$
for all~$i\in I$.
\item[(d)] If $k\in I$ with $k<j$ then~$\eta_i^j\ci\phi^{jk}\equiv
\eta_i^k\vert_{V^{jk}}$.
\end{itemize}
Here $\eta_i^k\vert_{V^{jk}}$ makes sense in (d) since when we are
choosing $\{\eta_i^j:i\in I\}$ we have already chosen $\{\eta_i^k:
i\in I\}$, as $k<j$. The second part of (b) means that if
$\phi^{ij}(V^{ij})$ has an `open edge' in $V^i$, then $\eta_j^j=0$
near this edge in $V^{ij}$, or equivalently by (d),
$\eta_j^i\vert_{\phi^{ij} (V^{ij})}=0$ near this edge
in~$\phi^{ij}(V^{ij})$.

Observe that conditions (c),(d) prescribe $\eta_i^j$ for $i\in I$ on
the subset
\e
(s^j)^{-1}(0)\cup\ts\bigcup_{k\in I:k<j}\phi^{jk}(V^{jk}).
\label{kh3eq3}
\e
These prescribed values are continuous on \eq{kh3eq3}, and
consistent on the overlaps between $(s^j)^{-1}(0)$ and
$\phi^{jk}(V^{jk})$ and $\phi^{jl}(V^{jl})$ for different $k,l<j$ in
$I$, using Definition \ref{kh3def1} and conditions (a)--(d) for
previous choices $\{\eta_i^k:i\in I\}$ for $k<j$. However, there is
a problem in extending the $\eta_i^j$ continuously from \eq{kh3eq3}
to all of $V^j$, since $\phi^{jk}(V^{jk})$ may not be closed in
$V^j$, and the prescribed values of $\eta_i^j$ may not extend
continuously to the closure of \eq{kh3eq3} in~$V^j$.

We can avoid this problem by making all the sets $V^k,V^{jk}$ {\it
slightly smaller}. That is, we replace $V^k$ and $V^{jk}$ by open
sets $\ti V^j\subset V^j$ and $\ti V^{jk}\subset V^{jk}$ such that
the previous conditions still hold, in particular, $X=\bigcup_{i\in
I}\Im\psi_i\vert_{\ti V^i}$, and in (b), we shrink the `open edges'
of $\phi^{jk}(V^{jk})$ in $V^k$ by little enough that $\eta_j^j=0$
near the new open edges, but we make $\ti V^{jk}$ small enough that
$\overline{\phi^{jk}(\ti V^{jk})} \subseteq\phi^{jk}(V^{jk})$,
taking closures in $V^j$, so the closure of the new set \eq{kh3eq3}
lies in the old set \eq{kh3eq3}. This is automatic if the closure of
$\ti V^{ij}$ in $V^{ij}$ is compact. With this replacement, the
$\eta_i^j$ have unique continuous extensions from their prescribed
values on \eq{kh3eq3} to the closure of \eq{kh3eq3}.

We can then choose arbitrary continuous extensions
$\eta_i^j:V^j\ra[0,1]$ from the values on this closed set to $V^j$,
satisfying (b) and $\smash{\sum_{i\in I}\eta_i^j\equiv 1}$.
Therefore by induction we can choose $\eta_i^j$ for all $i,j\in I$
satisfying (a)--(d), possibly with smaller sets $V^k,V^{jk}$. We
have now constructed a very good coordinate system $\bs I$ and
functions $\eta_i,\eta_i^j$ satisfying all of Definition
\ref{kh3def3} except the last two parts of (i), the last three parts
of (ii), and (v),(vi). But we do know that $\{p\in
X:\eta_i(p)>0\}\subseteq\Im\psi^i$, and $\{v\in
V^j:\eta_i^j(v)>0\}\subseteq V^{ij}$ if~$j\le i$.

For (v), as $X$ is compact and the map taking $p\in X$ to the number
of local boundary components of $X$ at $p$ is upper semicontinuous,
there exists $M\ge 0$ such that each $p\in X$ lies in at most $M$
local boundary components. If $v\in(\psi^i)^{-1} (0)\subseteq V^i$
with $\psi^i(v)=p$ in $X$ then the number of local boundary
components of $V^i$ containing $v$ is the same as the number of
local boundary components of $X$ containing $p$, and so at most $M$.
Conversely, if $v\in V^i$ and the number of local boundary
components of $V^i$ containing $v$ is greater than $M$, then
$v\notin(\psi^i)^{-1}(0)$. Thus, the subset of $v\in V^i$ lying in
more than $M$ local boundary components is a closed subset of $V^i$
not intersecting $(\psi^i)^{-1}(0)$. We will delete this subset,
which will make (v) hold.

For (vi), as $X$ is compact, $f(X)\subseteq Y$ is compact. Choose an
open neighbourhood $U$ of $f(X)$ in $Y$ whose closure $T=\bar U$ in
$Y$ is compact. To complete the proof we now replace $V^i$ by the
open subset $\ti V^i=\bigl\{v\in V^i:\eta_i^i(v)>0$, $v$ lies in at
most $M$ local boundary components, $f^i(y)\in U\bigr\}$ for each
$i\in I$, and replace $I$ by $\ti I=\{i\in I:\ti V^i\ne\es\}$, and
replace $V^{ij}$ by the open subset $\ti V^{ij}=V^{ij}\cap\ti
V^j\cap(\phi^{ij})^{-1}(\ti V^i)$ for all $i,j\in\ti I$. With these
replacements, we have $\Im\psi^i=\{p\in X:\eta_i(p)>0\}$ and
$V^{ij}=\{v\in V^j:\eta_i^j(v)>0\}$ for all $i,j$. These imply the
remaining parts of (i),(ii) except $\phi^{ij}(V^{ij})$ closed in
$\{v\in V^i:\eta_j^i(v)>0\}$ for $j\le i$, and this follows from (b)
above. Hence $(\bs I,\bs\eta)$ is a {\it really good coordinate
system\/} for~$(X,\bs f)$.
\end{proof}

Here is an algorithm for converting really good coordinate systems
into excellent coordinate systems, by throwing away all unnecessary
connected components of $\pd^kV^i$ for all $i\in I$ and~$k\ge 0$.

\begin{alg} Let $(\bs I,\bs\eta)$ be a really good coordinate
system for $(X,\bs f)$. We shall construct open sets $\check V^i$ in
$V^i$ and $\check V^{ij}\subseteq V^{ij}$ for all $j\le i\in I$,
such that defining $\check I=\{i\in I:\check V^i\ne\es\}$, and for
$i,j\in\check I$ setting $(\check V^i,\check E^i,\check
s^i,\check\psi^i)=(\check V^i,E^i\vert_{\check V^i},s^i\vert_{\check
V^i},\psi^i\vert_{\check V^i})$, $\check f^i=f^i\vert_{\check V^i}$,
$\check\eta_i=\eta_i$, $\check\eta_i^j=\eta_i^j\vert_{\check V^j}$,
and for $j\le i$ in $\check I$ defining $\check
V^{ij}=V^{ij}\cap\check V^j\cap(\phi^{ij})^{-1}(\check V^i)$ and
$(\check\phi^{ij},\smash{\hat{\check\phi}{}^{ij}})=
(\phi^{ij},\hat\phi^{ij})\vert_{\check V^{ij}}$, this data defines
an {\it excellent coordinate system\/} $(\bs{\check I},
\bs{\check\eta})$ for $(X,\bs f)$. We define $\check V^j$ by
induction on decreasing~$j\in I$.

In the inductive step, for fixed $j\in I$, suppose we have defined
$\check V^i$ for all $i\in I$ with $j<i\in I$, and $\check V^{ik}$
for all $j<i\le k$ in $I$. Define $\check V^{ij}\subseteq
V^{ij}\subseteq V^j$ for all $j<i$ in $I$ by $\check V^{ij}=
(\phi^{ij})^{-1}(\check V^i)$. Then define $\check V^j$ to be the
complement in $V^j$ of the union over all $l=0,\ldots,\dim V^j$ of
the images under the projection $\pd^lV^j\ra V^j$ of all connected
components $W$ of $\pd^lV^j$ such that $W\cap\pd^l\check V^{ij}=\es$
for all $j<i\in I$ and
\e
\Im\psi^j\vert_W\subseteq\bigcup\nolimits_{i\in I:i<j}\Im\psi^i
\vert_{\pd^l V^i}\cup\bigcup\nolimits_{i\in I:i>j}\Im\psi^i
\vert_{\pd^l\check V^i}.
\label{kh3eq4}
\e
For each $l$, the union of such $W$ is an open and closed set in
$\pd^lV^j$, so its image is closed in $V^j$. Thus $\check V^j$ is
the complement of a finite union of closed sets in $V^j$, and so is
open in $V^j$. Set $\check V^{jj}=\check V^j$. This completes the
inductive step.

The condition \eq{kh3eq4} on $W$ ensures that subtracting these
components $W$ does not change the fact that $X=\bigcup_{i\in
I}\Im\psi^i$, and therefore $X=\bigcup_{i\in I}\Im\check\psi^i$. One
can prove this by induction on $i\in I$, by checking that if $X=
\bigcup_{i\in I:i\le j}\Im\psi^i\cup\bigcup_{i\in I:i>j}
\Im\check\psi^i$ then $X=\bigcup_{i\in I:i<j}\Im\psi^i\cup
\bigcup_{i\in I:i\ge j}\Im\check\psi^i$. Also, subtracting
components $W$ in this way does not change the fact that the last
sentence of Definition \ref{kh3def3}(i) and the last two sentences
of Definition \ref{kh3def3}(ii) hold. For the middle sentence of
Definition \ref{kh3def3}(ii), note that for $j<i$ we only ever
delete some component $W$ of $\pd^lV^{ij}$ if in a previous
inductive step we have also deleted the component of $\pd^lV^i$
containing $\phi^{ij}(W)$. Doing both these deletions together does
not change the fact that $\phi^{ij}(V^{ij})$ is a closed subset of
$\{v\in V^i:\eta_j^i(v)>0\}$. So $(\bs{\check I},\bs{\check\eta})$
is a {\it really good coordinate system\/} for $(X,\bs f)$, since
$(\bs I,\bs\eta)$ is.

We claim that $(\bs{\check I},\bs{\check\eta})$ is an {\it excellent
coordinate system\/} for $(X,\bs f)$. Suppose for a contradiction
that $\check W$ is an unnecessary connected component of some
$\pd^l\check V^j$. Then $\check W$ is an open dense subset of a
unique connected component $W$ of $\pd^lV^j$. The difference
$W\sm\check W$ is a union of unnecessary components of $\pd^{l'}V^j$
for $l'>l$, and does not affect the following argument. As $\check
W\cap \check V^{ij}=\es$ for $j<i\in I$ we have $W\cap\check
V^{ij}=\es$ for $j<i$. Also $\Im\check\psi^j\vert_{\check W}
\subseteq\bigcup_{j\ne i\in I}\Im\check\psi^i\vert_{\pd^l\check
V^i}\subseteq\pd^lX$ implies that \eq{kh3eq4} holds. So $W$ is one
of the connected components of $\pd^lV^j$ which is deleted from
$\check V^j$, a contradiction. Thus $(\bs{\check
I},\bs{\check\eta})$ is excellent.
\label{kh3alg}
\end{alg}

Together with Proposition \ref{kh3prop3} this yields:

\begin{cor} Let\/ $X$ be a compact Kuranishi space, $Y$ an
orbifold, and\/ $\bs f:X\ra Y$ a strongly smooth map. Then there
exists an excellent coordinate system $(\bs I,\bs\eta)$
for\/~$(X,\bs f)$.
\label{kh3cor}
\end{cor}

We now prove Property \ref{kh3pr}(a):

\begin{thm} Let\/ $X$ be a compact Kuranishi space, $Y$ an orbifold,
and\/ $\bs f:X\ra Y$ a strongly smooth map. Then\/ $(X,\bs f)$
admits gauge-fixing data\/~$\bs G$.

If\/ $\Ga$ is a finite group of strong diffeomorphisms\/ $\bs a:X\ra
X$ with $\bs f\ci\bs a=\bs f$ then we can choose\/ $\bs G$
$\Ga$-invariant, i.e., there is a morphism\/ $\Ga\!\ra\!\Aut(X,\bs
f,\bs G)$ mapping\/ $\bs a\mapsto(\bs a,\bs b_{\bs a})$.

If\/ $X$ is without boundary, or with boundary, or with corners,
then we can choose $\bs G$ without boundary, or with boundary, or
with corners, respectively. All this also holds for co-gauge-fixing
data $\bs C,$ with $\bs f$ a strong submersion.
\label{kh3thm1}
\end{thm}

\begin{proof} By Corollary \ref{kh3cor} there exists an excellent
coordinate system $\bs I=\bigl(I,\ab(V^i,\ab E^i,\ab s^i,\psi^i),
f^i:i\in I),\ldots\bigr)$ for $(X,\bs f)$. For each $i\in I$ choose
an {\it injective} map $G^i:E^i\ra P_n$ for $n\gg 0$. This is
clearly possible, and we can even choose $G^i$ to be smooth, as a
generic smooth map $G^i:E^i\ra\R^n/S_n$ is injective provided
$n>2\dim E^i$. As $G^i$ is injective it is globally finite, with
$N=1$. Thus $\bs G=(\bs I,\bs\eta,G^i:i\in I)$ is gauge-fixing data
for~$(X,\bs f)$.

When $\bs f$ is a strong submersion, to construct co-gauge-fixing
data $\bs C$ we first note that in Proposition
\ref{kh3prop1}--Corollary \ref{kh3cor} we can work throughout with
good, very good, really good and excellent coordinate systems in
which $f^i:V^i\ra Y$ are submersions. We again choose injective
$C^i:E^i\ra P_n$ for $n\gg 0$, and $C^i$ injective implies
$C^i\t(f^i\ci\pi^i)$ is globally finite, with $N=1$, so $\bs C=(\bs
I,\bs\eta,C^i:i\in I)$ is co-gauge-fixing data for~$(X,\bs f)$.

Now let $\Ga\subseteq\Aut(X,\bs f)$ be a finite subgroup. We must
show we can choose $\bs G,\bs C$ to be $\Ga$-invariant. This is
similar to the problem of choosing (co-)gauge-fixing data for
$(X/\Ga,\bs f\ci\bs\pi)$ for $\bs\pi:X\ra X/\Ga$ the natural
projection, which we have already proved is possible in the first
part. We begin with the proof of Proposition \ref{kh3prop2}, which
constructs a very good coordinate system $\bs{\check I}$ for $(X,\bs
f)$ starting with a good coordinate system $\bs I$ for $(X,\bs f)$,
which exists by Proposition \ref{kh3prop1}. This $\bs I$ is built
from finitely many neighbourhoods $(V_p,\ldots,\psi_p)$ for~$p\in
X$.

To choose $\bs I$ to be $\Ga$-invariant, we take each $(V_p,\ldots,
\psi_p)$ to be invariant under $\{\bs a\in\Ga:a(p)=p\}$, with $V_p$
small enough that $\Im\psi_p\cap(a\cdot\Im\psi_p)=\es$ whenever $\bs
a\in\Ga$ with $a(p)\ne p$. Every $p\in X$ admits such a Kuranishi
neighbourhood. Given such $(V_p,E_p,s_p,\psi_p)$, for $\bs a\in\Ga$
define $(V_{a(p)},E_{a(p)},s_{a(p)},\psi_{a(p)})$ to be $\bs
a_*(V_p,E_p,s_p,\psi_p)$, a Kuranishi neighbourhood of~$a(p)\in X$.

As $(V_p,\ldots,\psi_p)$ is invariant under $\{\bs a\in\Ga:
a(p)=p\}$, this $(V_{a(p)},\ldots,\psi_{a(p)})$ depends only on
$a(p)$ rather than on $\bs a$, and the finite set of neighbourhoods
$(V_{a(p)},\ldots,\psi_{a(p)})$ for $\bs a\in\Ga$ is invariant under
$\Ga$, with $\Im\psi_{a(p)}\cap \Im\psi_{a'(p)}=\es$ when $a(p)\ne
a'(p)$. Thus the disjoint union $\coprod_{a(p):\bs a\in\Ga}
(V_{a(p)},\ldots,\psi_{a(p)})$ is a single $\Ga$-invariant Kuranishi
neighbourhood. We can now follow the proof of Proposition
\ref{kh3prop1} to make $\bs I$ using $\Ga$-{\it invariant\/}
Kuranishi neighbourhoods~$(V^i,\ldots,\psi^i)=\coprod_{a(p):\bs
a\in\Ga}(V_{a(p)},\ldots, \psi_{a(p)})$.

Then in Proposition \ref{kh3prop2} we construct a very good
coordinate system $\bs{\check I}$ for $(X,\bs f)$ by gluing together
the neighbourhoods $(V^i,\ldots,\psi^i)$ of each dimension in $\bs
I$. As each $(V^i,\ldots,\psi^i)$ is $\Ga$-invariant, we can do this
gluing in a $\Ga$-invariant way, to make $\bs{\check I}$
$\Ga$-invariant. That is, $\bs a\in\Ga$ lifts to $(\bs a,\bs b_{\bs
a})\in\Aut(X,\bs f,\bs{\check I})$, where $\bs b_{\bs a}=(b_{\bs
a}^i,\hat b_{\bs a}^i:i\in\check I)$, with each $(b_{\bs a}^i,\hat
b_{\bs a}^i):(\check V^i,\ldots,\check\psi^i)\ra (\check
V^i,\ldots,\check\psi^i)$ an isomorphism lifting~$\bs a$.

In Proposition \ref{kh3prop3} it is easy to extend this to a
$\Ga$-invariant really good coordinate system $(\bs I,\bs\eta)$ for
$(X,\bs f)$, we just choose $\eta_i,\eta_i^j$ to be $\Ga$-invariant.
Algorithm \ref{kh3alg} then automatically yields a $\Ga$-invariant
excellent coordinate system. Thus in Corollary \ref{kh3cor}, we can
take $(\bs I,\bs\eta)$ to be $\Ga$-invariant. Finally, in the first
part of the theorem, rather than choosing $G^i:E^i\ra P_n$,
$C^i:E^i\ra P_n$ injective, we take $G^i,C^i$ to be the pullbacks to
$E^i$ of injective maps $E^i/\Ga\ra P_n$ for $n\gg 0$. Then
$G^i,C^i$ are $\Ga$-invariant, and globally finite in the sense of
Definition \ref{kh3def6}, with $N=\md{\Ga}$. This implies
$C^i\t(f^i\ci\pi^i)$ is also globally finite, with $N=\md{\Ga}$.
Hence we can choose $\bs G,\bs C$ to be $\Ga$-invariant.

If $X$ is without boundary, or with boundary, or with corners, then
we can go through the whole of \S\ref{kh32} working only with $V^i$
without boundary, or with boundary, or with corners, respectively,
and so construct $\bs G,\bs C$ without boundary, or with boundary,
or with corners, as required.
\end{proof}

\subsection{Finiteness of automorphism groups}
\label{kh33}

We now prove Property \ref{kh3pr}(b) for (co-)gauge-fixing data.

\begin{prop} If\/ $(\bs I,\bs\eta)$ is an excellent coordinate
system for\/ $(X,\bs f),$ with\/ $X$ compact, then for each\/ $j\in
I$ there are only finitely many connected components\/ $W$ of\/
$V^j$ such that\/ $W\cap V^{ij}=\es$ for all\/ $j<i\in I$.
\label{kh3prop4}
\end{prop}

\begin{proof} Let $j\in I$ and $W^j_a$ for $a\in A^j$ be the connected
components of $V^j$ with\/ $W^j_a\cap V^{ij}=\es$ for all $j<i\in
I$. Define $U^j=V^j\sm\bigcup_{a\in A}W^j_a$. Then $V^j$ is the
disjoint union of $U^j$ and the $W^j_a$ for $a\in A^j$, and $U^j$
and each $W^j_a$ are open and closed in $V^j$. Note that openness
holds as $V^j$ is an orbifold; for general topological spaces
connected components are closed but need not be open. For instance,
the compact subset $\{0\}\cup\{1/n:n=1,2,\ldots\}$ in $\R$, has a
connected component $\{0\}$ which is closed but not open.

Now $C=X\sm\bigl(\Im\psi^j\vert_{U^j}\cup\bigcup_{j\ne i\in
I}\Im\psi^i\bigr)$ is a closed subset of $X$, since $\Im\psi^j
\vert_{U^j}$ and each $\Im\psi^i$ is open, so $C$ is compact as $X$
is compact and Hausdorff. Also $C\subseteq\bigcup_{a\in
A^j}\Im\psi^j\vert_{W^j_a}\subseteq X$ since $X=\bigcup_{j\in
I}\Im\psi^j$. And as $(\bs I,\bs\eta)$ is excellent, Definition
\ref{kh3def4} implies that $C\cap\Im\psi^j\vert_{W^j_a}\ne\es$ for
each $a\in A^j$, since otherwise $W^j_a$ would be an unnecessary
component of $V^j$. Therefore the sets $C\cap\Im\psi^j\vert_{W^j_a}$
for $a\in A^j$ form a cover of $C$ by nonempty, disjoint open sets.
As $C$ is compact there must be a finite subcover, but the only
possible subcover is the whole thing, since the sets are nonempty
and disjoint. Hence $A^j$ is finite.
\end{proof}

\begin{thm} Let\/ $X$ be a compact Kuranishi space, $Y$ an
orbifold, $\bs f:X\ra Y$ be strongly smooth, and\/ $\bs G$
gauge-fixing data for\/ $(X,\bs f)$. Then the automorphism group\/
$\Aut(X,\bs f,\bs G)$ is finite. The same holds for co-gauge-fixing
data.
\label{kh3thm2}
\end{thm}

\begin{proof} As in Definition \ref{kh3def6} we write elements of
$\Aut(X,\bs f,\bs G)$ as $(\bs a,\bs b)$ with $\bs b=((b^i,\hat
b^i):i\in I)$ for $b^i:V^i\ra V^i$, $\hat b^i:E^i\ra E^i$, and
$G^i\ci\hat b^i\equiv G^i:E^i\ra P$. As $G^i$ is globally bounded,
by Definition \ref{kh3def6} there exists $N\ge 1$ such that
$\md{(G^i)^{-1}(p)}\le N$ for all~$p\in P$.

For each such $\hat b^i$ and $p\in P$, $(G^i)^{-1}(p)$ is at most
$N$ points in $E^i$, and $\hat b^i$ must permute these points. But
the order of a permutation of a set of at most $N$ points must
divide $N!$, so $(\hat b^i)^{N!}$ is the identity on
$(G^i)^{-1}(p)\subset E^i$. As this is true for all $p\in P$, $(\hat
b^i)^{N!}$ is the identity on $E^i$. So, every such $\hat b^i$ is a
{\it diffeomorphism of finite order}. Thus the fixed point locus of
$\hat b^i$ in $E^i$ is a locally finite union of closed, embedded,
connected suborbifolds of $E^i$, possibly of varying dimensions.
Hence for each connected component $F$ of $E^i$, {\it either\/}
$\hat b^i$ is the identity on $F$, {\it or\/} $\hat b^i$ does not
fix generic points of~$F$.

Suppose for a contradiction that for some connected component $F$ of
$E^i$, we can find $(\bs a_1,\bs b_1),\ldots,(\bs a_{N+1},\bs
b_{N+1})$ in $\Aut(X,\bs f,\bs G)$ such that $\hat b^i_c\vert_F$ are
distinct maps $F\ra E^i$ for $c=1,\ldots,N+1$. Then for $1\le c<d\le
N+1$ we have finite order diffeomorphisms $\hat b^i_d\ci(\hat
b^i_c)^{-1}:E^i\ra E^i$ which are not the identity on $\hat
b^i_c(F)$, so $\hat b^i_d\ci(\hat b^i_c)^{-1}$ does not fix generic
points in $\hat b^i_c(F)$, and $\hat b^i_c,\hat b^i_d$ do not agree
on generic points in $F$. Therefore $\hat b^i_1(f),\ldots,\hat
b^i_{N+1}(f)$ are distinct points in $E^i$ for generic $f\in F$. But
$G^i\ci\hat b^i_c(f)=G^i(f)$ for $c=1,\ldots,N+1$, so $\hat
b^i_1(f),\ldots,\hat b^i_{N+1}(f)$ lie in $(G^i)^{-1}(G^i(f))$, and
$\bmd{(G^i)^{-1}(G^i(f))}\le N$, a contradiction. Therefore, for
each connected component $F$ of $E^i$, the maps $\hat b^i:E^i\ra
E^i$ coming from $(\bs a,\bs b)\in\Aut(X,\bs f,\bs G)$ can realize
at most $N$ distinct maps~$\hat b^i\vert_F:F\ra E^i$.

Now use the notation of the proof of Proposition \ref{kh3prop4}, so
that $V^j=U^j\amalg\coprod_{a\in A^j}W^j_a$ with each $W^j_a$ a
connected component of $V^j$, and $A^j$ finite. Then each
$E^j\vert_{W^j_a}$ is a connected component of $E^j$, so the $\hat
b^j\vert_{E^j\vert_{W^j_a}}$ realize at most $N$ distinct maps.
Hence, taken over all $j\in I$ and $a\in A^j$, the family of
restrictions $\hat b^j\vert_{E^j\vert_{W^j_a}}$ can realize at most
$N^{\sum_{j\in I}\md{A^j}}$ different possibilities.

Let $(\bs a,\bs b),(\bs a',\bs b')\in\Aut(X,\bs f,\bs G)$ with $\hat
b^j\vert_{E^j\vert_{W^j_a}}=\hat b^{\prime j}\vert_{E^j\vert_{
W^j_a}}$ for all $j\in I$ and $a\in A^j$. We will show that $(\bs
a,\bs b)=(\bs a',\bs b')$. Suppose for a contradiction that $\hat
b^j\ne\hat b^{\prime j}$ for some $j\in I$, and let $j$ be largest
for which this holds. From $V^j=U^j\amalg\coprod_{a\in A^j}W^j_a$
and $\hat b^j\vert_{E^j\vert_{W^j_a}}=\hat b^{\prime
j}\vert_{E^j\vert_{W^j_a}}$ we see that $\hat
b^j\vert_{E^j\vert_{U^j}}\ne\hat b^{\prime
j}\vert_{E^j\vert_{U^j}}$. Thus, there exists some connected
component $W$ of $U^j$ for which $\hat b^j\vert_{E^j\vert_W}\ne\hat
b^{\prime j}\vert_{E^j\vert_W}$. Then $W\cap V^{ij}\ne\es$ for some
$j<i\in I$, by definition of $W^j_a$ and~$U^j$.

On $E^j\vert_{W\cap V^{ij}}$ we have $\hat b^j=\hat
b^i\ci\hat\phi^{ij}=\hat b^{\prime i}\ci\hat\phi^{ij}=\hat b^{\prime
j}$ by Definition \ref{kh3def6}(c), since $j$ is largest with $\hat
b^j\ne\hat b^{\prime j}$ and $j<i$, so $\hat b^i=\hat b^{\prime i}$.
Hence, we know that $\hat b^j\vert_{E^j\vert_W}\ne\hat b^{\prime
j}\vert_{E^j\vert_W}$, so that $\hat b^j,\hat b^{\prime j}$ differ
at generic points of $E^j\vert_W$, but also $\hat b^j=\hat b^{\prime
j}$ on $E^j\vert_{W\cap V^{ij}}$, which is a nonempty open set in
$E^j\vert_W$. Choosing a generic point in $E^j\vert_{W\cap V^{ij}}$
gives a contradiction. Hence $\hat b^j=\hat b^{\prime j}$ for all
$j\in I$. But $\hat b^j:E^j\ra E^j$ is a lift of $b^j:V^j\ra V^j$,
so each $\hat b^j$ determines $b^j$. Thus $b^j=b^{\prime j}$ for
$j\in I$. Also each $(b^j,\hat b^j)$ determines $\bs a$ over
$\Im\psi^j$, and $X=\bigcup_{j\in I}\Im\psi^j$, so $\bs b$
determines $\bs a$. Hence $(\bs a,\bs b)=(\bs a',\bs b')$.

This proves that each $(\bs a,\bs b)\in\Aut(X,\bs f,\bs G)$ is
determined by the restrictions $\hat
b^j\vert_{E^j\vert_{\smash{W^j_a}}}$ for all $j\in I$ and $a\in
A^j$, which realize at most $N^{\sum_{j\in I}\md{A^j}}$ different
possibilities. So $\Aut(X,\bs f,\bs G)$ is finite, with
$\bmd{\Aut(X,\bs f,\bs G)}\le N^{\sum_{j\in I}\md{A^j}}$.
\end{proof}

\subsection{Quotients by finite groups}
\label{kh34}

(Co-)gauge-fixing data descends to finite quotients, as in
Property~\ref{kh3pr}(c).

\begin{dfn} Let $X$ be a compact Kuranishi space, $Y$ an orbifold,
$\bs f:X\ra Y$ a strongly smooth map, and $\bs G$ be gauge-fixing
data for $(X,\bs f)$. Suppose $\Ga$ is a finite group, and $\rho$ an
action of $\Ga$ on $(X,\bs f,\bs G)$ by isomorphisms. That is,
$\rho:\Ga\ra\Aut(X,\bs f,\bs G)$ is a group morphism. Note that we
do not require $\rho$ to be injective, so we cannot regard $\Ga$ as
a subgroup of $\Aut(X,\bs f,\bs G)$. Write $\bs G=(\bs I,\bs\eta,
G^i:i\in I)$, where $\bs I=\bigl(I,(V^i,E^i,s^i,\psi^i),f^i:i\in
I,\ldots\bigr)$. Then $\rho$ induces actions of $\Ga$ on $X$ by
strong diffeomorphisms and on $V^i,E^i$ for $i\in I$ by
diffeomorphisms, where if $\ga\in\Ga$ with $\rho(\ga)=(\bs a,\bs b)$
and $\bs b=\bigl((b^i,\hat b^i):i\in I\bigr)$, then $\ga$ acts on
$X$ by $\bs a$, on $V^i$ by $b^i$, and on $E^i$ by~$\hat b^i$.

Define a compact Kuranishi space $\ti X=X/\Ga$ and orbifolds $\ti
V^i=V^i/\Ga$ and $\ti E^i=E^i/\Ga$ for $i\in I$. As $\rho$ need not
be injective, for example we allow the case that $\Ga$ is a
nontrivial group but $\rho\equiv 1$, so that $\Ga$ acts trivially on
$X$ and $V^i,E^i$. In this case $\ti X,\ti V^i,\ti E^i$ coincide
with $X,V^i,E^i$ as topological spaces, but not as Kuranishi spaces
or orbifolds, since taking the quotient by $\Ga$ adds a factor of
$\Ga$ to the stabilizer groups of $X,V^i,E^i$ at every point, so
that $\ti V^i,\ti E^i$ are {\it noneffective\/} orbifolds.

Let $\pi^i:E^i\ra V^i$ be the projection. If $\ga\in\Ga$ and
$\rho(\ga)=(\bs a,\bs b)$ then $\pi^i\ci\hat b^i=b^i\ci\pi^i$ by
definition of $\Aut(X,\bs f,\bs G)$, so $\pi^i:E^i\ra V^i$ is
equivariant with respect to the actions of $\Ga$ on $V^i,E^i$, and
$\pi^i$ pushes down to a submersion $\ti\pi^i:\ti E^i\ra\ti V^i$.
The orbibundle structure on $\pi^i:E^i\ra V^i$ descends to an
orbibundle structure on $\ti\pi^i:\ti E^i\ra\ti V^i$. Since $\hat
b^i\ci s^i\equiv s^i\ci b^i$ for all $(\bs a,\bs b)=\rho(\ga)$,
$s^i$ pushes down to a smooth section $\ti s^i$ of $\ti E^i\ra\ti
V^i$, and $(\ti s^i)^{-1}(0)=(s^i)^{-1}(0)/\Ga$. As $\psi^i\ci
b^i\equiv a\ci\psi^i$, where $a:X\ra X$ is the homeomorphism
underlying $\bs a$ in $(\bs a,\bs b)=\rho(\ga)$, $\psi^i$ pushes
down to a continuous $\ti\psi^i:(\ti s^i)^{-1}(0)\ra\ti X=X/\Ga$.
Then $(\ti V^i,\ti E^i,\ti s^i,\ti\psi^i)$ is a Kuranishi
neighbourhood on~$\ti X$.

Write $\pi_{\ti V^i}:V^i\ra\ti V^i$ and $\pi_{\ti E^i}:E^i\ra\ti
E^i$ for the natural projections. Since $f^i:V^i\ra Y$ and
$G^i:E^i\ra P$ are $\Ga$-invariant they descend to $\ti f^i:\ti
V^i\ra Y$ and $\ti G^i:\ti E^i\ra P$ with $\ti f^i\ci\pi_{\ti
V^i}\equiv f^i$ and $\ti G^i\ci\pi_{\ti E^i}\equiv G^i$. Similarly
$\eta_i:X\ra[0,1]$ and $\eta_i^j:V^j\ra[0,1]$ descend to
$\ti\eta_i:\ti X\ra[0,1]$ and $\ti\eta_i^j:\ti V^j\ra[0,1]$. If
$j\le i$ in $I$ then $V^{ij}$ is a $\Ga$-invariant subset of $V^j$,
so $\ti V^{ij}=V^{ij}/\Ga$ is an open subset of $\ti V^j=V^j/\Ga$,
and $(\phi^{ij},\hat\phi^{ij})$ descends to a coordinate
transformation $(\ti\phi^{ij},\hat{\ti\phi}{}^{ij}):(\ti V^{ij},\ti
E^j\vert_{\ti V^{ij}},\ti s^j\vert_{\ti V^{ij}},\ti \psi^j\vert_{\ti
V^{ij}})\ra(\ti V^i,\ldots,\ti\psi^i)$. It is now easy to verify
that all this comprises {\it gauge-fixing data\/} $\bs{\ti G}$ for
$(\ti X,\bs{\ti f})$. We shall also write $\bs{\ti f}=\bs\pi_*(\bs
f)$ and $\bs{\ti G}=\bs\pi_*(\bs G)$, where $\bs\pi:X\ra\ti X$ is
the projection.

Exactly the same construction works for {\it co-gauge-fixing
data}~$\bs{\ti C}=\bs\pi_*(\bs C)$.
\label{kh3def9}
\end{dfn}

\subsection{Restriction of (co-)gauge-fixing data to boundaries}
\label{kh35}

Very good, really good and excellent coordinate systems on $X$
restrict to~$\pd X$.

\begin{dfn} Let $X$ be a compact Kuranishi space, $Y$ an orbifold,
$\bs f:X\ra Y$ be strongly smooth, and $\bs
I=\smash{\bigl(I,(V^i,\ldots, \psi^i),f^i:i\in I,\ldots\bigr)}$ be a
very good coordinate system for $(X,\bs f)$. Define a very good
coordinate system $\bs{\ti I}=\bs I\vert_{\pd X}$ for $(\pd X,\bs
f\vert_{\pd X})$ to have indexing set $\ti I=\{i:i+1\in I$, $\pd
V^{i+1}\ne\es\}$ and Kuranishi neighbourhoods
\begin{equation*}
\bigl(\ti V^i,\ti E^i,\ti s^i,\ti\psi^i\bigr)=\bigl(\pd
V^{i+1},E^{i+1}\vert_{\pd V^{i+1}},s^{i+1}\vert_{\pd
V^{i+1}},\psi^{i+1}\vert_{\pd V^{i+1}}\bigr)\quad\text{for $i\in\ti
I$.}
\end{equation*}
Set $\ti V^{ij}=\pd V^{(i+1)(j+1)}$ and
$(\ti\phi^{ij},\smash{\hat{\ti\phi}{}^{ij}})=
(\phi^{(i+1)(j+1)},\hat\phi^{(i+1)(j+1)})\vert_{\ti V^{ij}}$ for
$j\le i$ in $\ti I$. Set $\ti f^i=f^{i+1}\vert_{\pd V^{i+1}}$ for
$i\in\ti I$. It is now easy to verify that $\bs{\ti I}=\bigl(\ti
I,(\ti V^i,\ldots,\ti\psi^i):i\in\ti I,\ldots\bigr)$ is a {\it very
good coordinate system\/} for $(\pd X,\bs f\vert_{\pd X})$, which we
write as~$\bs I\vert_{\pd X}$.

If $(\bs I,\bs\eta)$ is a {\it really good coordinate system\/} for
$(X,\bs f)$ then we extend $\bs I\vert_{\pd X}$ to a really good
coordinate system $(\bs I\vert_{\pd X},\bs\eta\vert_{\pd X})$ in the
obvious way, so that $\bs\eta\vert_{\pd X}$ consists of the
functions $\ti\eta_i=\eta_{i+1}\vert_{\pd X}$ and $\ti\eta_i^j=
\eta_{i+1}^{j+1}\vert_{\pd V^{j+1}}$ for $i,j\in\ti I$. Definition
\ref{kh3def3}(v) holds with $M-1$ in place of $M$, or 0 when $M=0$,
as the local boundary components $\ti B$ of $\pd V^i$ at $(v,B)$
correspond to a subset of the local boundary components $B'$ of
$V^i$ at $v$ with~$B'\ne B$. If $(\bs I,\bs\eta)$ is an {\it
excellent coordinate system\/} for $(X,\bs f)$ then $(\bs
I\vert_{\pd X},\bs\eta\vert_{\pd X})$ is also an excellent
coordinate system, since the conditions on $V^i$ in Definition
\ref{kh3def4} are stable under taking boundaries.
\label{kh3def10}
\end{dfn}

We can restrict (co-)gauge-fixing data to boundaries, as in
Property~\ref{kh3pr}(d).

\begin{dfn} Let $X$ be a compact Kuranishi space, $Y$ an orbifold,
$\bs f:X\ra Y$ a strongly smooth map, and $\bs G=(\bs
I,\bs\eta,G^i:i\in I)$ be gauge-fixing data for $(X,\bs f)$.
Definition \ref{kh3def10} gives an excellent coordinate system $(\bs
I\vert_{\pd X},\bs\eta\vert_{\pd X})$ for $(\pd X,\bs f\vert_{\pd
X})$. For each $i$ in $\ti I=\{i:i+1\in I$, $\pd V^{i+1}\ne
\es\}$ we define $\ti G^i=G^{i+1}\vert_{(E^{i+1}\vert_{\pd
V^{i+1}})}$. Let $\io:E^{i+1}\vert_{\pd V^{i+1}}=\pd E^{i+1}\ra
E^{i+1}$ be the natural map, $M$ be as in Definition
\ref{kh3def3}(v) for $(\bs I,\bs\eta)$, and $N$ be as in Definition
\ref{kh3def6} for $G^{i+1}$. Then $\bmd{(G^{i+1})^{-1}(p)}\le N$ for
each $p\in P$, and $\bmd{\io^{-1}(e)}\le M$ for each $e\in
E^{i+1}\vert_{\pd V^{i+1}}$. Hence $\ti G^i=G^{i+1}\ci\io$ satisfies
$\bmd{(\ti G^i)^{-1}(p)}\le MN$ for each $p\in P$, that is, $\ti
G^i$ is {\it globally finite}, with constant $\ti N=MN$. Therefore
$\bs G\vert_{\pd X}=(\bs I\vert_{\pd X},\bs\eta\vert_{\pd X},\ti
G^i:i\in\ti I)$ is gauge-fixing data for~$(\pd X,\bs f\vert_{\pd
X})$.

Exactly the same construction works for co-gauge-fixing data.
\label{kh3def11}
\end{dfn}

\subsection{Extension of (co-)gauge-fixing data from boundaries}
\label{kh36}

Next we prove Property \ref{kh3pr}(e) for (co-)gauge-fixing data.

\begin{prop} Let\/ $X$ be a compact Kuranishi space with corners,
$Y$ an orbifold, and\/ $\bs f:X\ra Y$ be strongly smooth. Recall
from Definition {\rm\ref{kh2def19}} that there is a natural strong
diffeomorphism\/ $\bs\si:\pd^2X\ra\pd^2X$ with\/
$\bs\si^2=\bs\id_{\pd^2X}$. Let\/ $\bs J$ be a very good coordinate
system for\/ $(\pd X,\bs f\vert_{\pd X}),$ with corners. Then there
exists a very good coordinate system\/ $\bs I$ for\/ $(X,\bs f)$
with corners with\/ $\bs I\vert_{\pd X}=\bs J$ if and only if\/ $\bs
J\vert_{\pd^2X}$ is invariant under\/ $\bs\si$. The analogous result
holds for really good or excellent coordinate systems\/ $(\bs
J,\bs\ze),(\bs I,\bs\eta)$ for\/~$(\pd X,\bs f\vert_{\pd X}),(X,\bs
f)$.
\label{kh3prop5}
\end{prop}

\begin{proof} If $\bs I$ is a very good coordinate system with
corners for $(X,\bs f)$ then restricting to $\pd X$ and then to
$\pd^2X$ gives $\bs{\bar I}=\bs I\vert_{\pd^2X}$ for $(\pd^2X,\bs
f\vert_{\pd^2X})$, with indexing set $\bar I=\{i:i+2\in I$,
$\pd^2V^i\ne\es\}$ and Kuranishi neighbourhoods $(\bar
V^i,\ldots,\bar\psi^i)$ with corners for $i\in \bar I$ with $\bar
V^i=\pd^2V^{i+2}$ and $(\bar E^i,\bar s^i,\bar\psi^i)=
(E^{i+2},s^{i+2},\psi^{i+2})\vert_{\pd^2V^{i+2}}$. Clearly, $\bs
I\vert_{\pd^2X}$ is invariant under $\bs\si$, with $\bs\si$ acting
on $\bar V^i$ as the involution $\si:\pd^2V^{i+2}\ra \pd^2V^{i+2}$
of Definitions \ref{kh2def5} and \ref{kh2def6}. Thus, a necessary
condition for $\bs I\vert_{\pd X}=\bs J$ is that $\bs
J\vert_{\pd^2X}$ is invariant under $\bs\si$. This proves the `only
if' part.

For the `if' part, let $\bs J$ be given and $\bs\si$-invariant on
$\pd^2X$ as above. Let $(W^j,F^j,t^j,\xi^j)$ for $j\in J$ be one of
the Kuranishi neighbourhoods with corners in $\bs J$. Then $(\pd
W^j,F^j\vert_{\pd W^j},t^j\vert_{\pd W^j},\xi^j\vert_{\pd W^j})$ is
a Kuranishi neighbourhood with corners on $\pd^2X$, and $\bs\si$ on
$\pd^2X$ lifts to an action on $(\pd W^j,\ldots,\xi^j\vert_{\pd
W^j})$. This is the condition necessary to extend $(W^j,\ldots,
\xi^j)$ to a Kuranishi neighbourhood $(V^{j+1},\ldots,\psi^{j+1})$
with corners on $X$ with $(\pd V^{j+1},\ldots,\psi^{j+1}\vert_{\pd
V^{j+1}})=(W^j,\ldots,\xi^j)$, which is `small' in the sense that it
extends only a little way into the interior of $X$, and $V^{j+1}$
extends only a little way from~$\pd V^{j+1}=W^j$.

Note that to extend the Kuranishi map $s^{j+1}$ smoothly from its
prescribed values on $\pd V^{j+1}=W^j$ to $V^{j+1}$ requires the
{\it Extension Principle}, Principle \ref{kh2pri}(c), which holds
only for manifolds and orbifolds with corners, not g-corners, and
this is why we restrict to Kuranishi spaces and (co-)gauge-fixing
data with corners, not g-corners. The Kuranishi structure on $X$ is
defined in terms of germs, and so determines $(V^{j+1},\ldots,
\psi^{j+1})$ only on an arbitrarily small open neighbourhood of
$(s^{j+1})^{-1}(0)$. Away from $(s^{j+1})^{-1}(0)$ we have to choose
$s^{j+1}$, and near $\pd V^{j+1}$ away from $(t^j)^{-1}(0)$ we need
the Extension Principle to do this.

Modifying the proof of Proposition \ref{kh3prop1}, we can build a
good ({\it not\/} very good, at this point) coordinate system $\bs
I$ for $(X,\bs f)$ with corners, with indexing set $I=\{j+1:j\in
J\}\amalg K$, using Kuranishi neighbourhoods
$(V^{j+1},\ldots,\psi^{j+1})$ for $j\in J$ as above with $V^{j+1}$
sufficiently small with $W^j=\pd V^{j+1}$ fixed, together with other
Kuranishi neighbourhoods $(V^k,\ldots,\psi^k)$ for $k$ in some
finite set $K$ on the interior of $X$, with $\pd V^k=\es$. Then we
use Proposition \ref{kh3prop2} to convert $\bs I$ to a very good
coordinate system $\bs{\check I}$ with corners for $(X,\bs f)$ by
gluing together the Kuranishi neighbourhoods of each dimension
in~$\bs I$.

In doing this, we replace the $V^i$ by open subsets $\ti V^i\subset
V^i$. We choose the $\ti V^{j+1}$ for $j\in J$ so that $\pd\ti
V^{j+1}=\pd V^{j+1}=W^j$, that is, we do not make the boundaries of
the $\ti V^{j+1}$ smaller. This is possible as there is at most one
nonempty $\pd V^i$ for $i\in I$ in each dimension, so the boundaries
$\pd V^i$ do not need to be glued together, and do not need to be
made smaller to ensure Hausdorffness. It then follows that
$\bs{\check I}\vert_{\pd X}=\bs J$, proving the `if' part.

To prove the analogous `if' result for {\it really good coordinate
systems}, we first construct a very good coordinate system
$\bs{\check I}$ with corners with $\bs{\check I}\vert_{\pd X}=\bs J$
as above, and then extend to really good $(\bs{\check
I},\bs{\check\eta})$ with $(\bs{\check I},\bs{\check\eta})\vert_{\pd
X}=(\bs J,\bs\ze)$. Thus we must construct continuous
$\check\eta_i:X\ra[0,1]$ and $\check\eta_i^j:\check V^j\ra[0,1]$
with prescribed values $\ze_{j-1}$ on $\pd X$ and $\ze_{i-1}^{j-1}$
on $\pd\check V^j$, satisfying Definition \ref{kh3def3}(i)--(iv).
This is possible by extending Proposition~\ref{kh3prop3}.

Two new issues arise. Firstly, in constructing the $\eta_i^j$ in
Proposition \ref{kh3prop3}, we made the $V^j,V^{jk}$ smaller at each
inductive step, so that $\eta_i^j$ should extend continuously from
its prescribed values on \eq{kh3eq3} to the closure of \eq{kh3eq3}
in $V^j$. However, as we need $\pd\check V^j=W^{j-1}$, $\pd\check
V^{jk}=W^{(j-1)(k-1)}$, we cannot make $\pd\check V^j,\pd\check
V^{jk}$ smaller. This does not matter, since we have prescribed
values $\check\eta_i^j\vert_{\pd\check V^j}\equiv\ze_{i-1}^{j-1}$
for $\check\eta_i^j$ on $\pd\check V^j$ consistent with the
prescribed values on \eq{kh3eq3}, so there is no need to make
$\pd\check V^j$ smaller to ensure consistency.

Secondly, Proposition \ref{kh3prop3} constructed a really good
coordinate system satisfying the extra properties $\Im\psi^i=\{p\in
X:\eta_i(p)>0\}$ and $V^{ij}=\{v\in V^j:\eta_i^j(v)>0\}$ for all
$i,j$, so the final sentences of both Definition
\ref{kh3def3}(i),(ii) hold. However, the boundary data $(\bs
J,\bs\ze)$ may not satisfy $\Im\xi^j=\{p\in \pd X: \ze_j(p)>0\}$ and
$W^{jk}=\{w\in W^k:\ze_j^k(w)>0\}$ for $k<j\in J$. Therefore we must
relax the requirements that $\supp(\eta_i)\subseteq\Im\psi^i$ and
$\supp(\eta_i^j)\subseteq V^{ij}$ in the proof of Proposition
\ref{kh3prop3}. We do this by requiring that $\supp(\eta_i)$ should
be supported in the union of $\Im\psi^i$ and a thin neighbourhood in
$X$ of $\supp(\ze_{i-1})$ in $\pd X$, and $\supp(\eta_i^j)$ should
be supported in the union of $V^{ij}$ and a thin neighbourhood in
$V^j$ of $\supp(\ze_{i-1}^{j-1})$ in $\pd V^j$. We choose these thin
neighbourhoods of $\supp(\ze_{i-1}),\supp(\ze_{i-1}^{j-1})$ small
enough that the final sentences of Definition \ref{kh3def3}(i),(ii)
for $(\bs J,\bs\ze)$ imply the final sentences of Definition
\ref{kh3def3}(i),(ii) for $(\bs I,\bs\eta)$ in the thin
neighbourhoods. Then replacing $V^i,V^{ij}$ by $\ti V^i,\ti V^{ij}$
as in the proof of Proposition \ref{kh3prop3}, one can show that we
get a {\it really good coordinate system} $(\bs I,\bs\eta)$,
with~$(\bs I,\bs\eta)\vert_{\pd X}=(\bs J,\bs\ze)$.

For {\it excellent coordinate systems}, if $(\bs J,\bs\ze)$ is
excellent, we construct a really good coordinate system $(\bs
I,\bs\eta)$ for $(X,\bs f)$ with corners, with $(\bs
I,\bs\eta)\vert_{\pd X}=(\bs J,\bs\ze)$, as above, and then apply
Algorithm \ref{kh3alg} to get an excellent coordinate system
$(\bs{\check I},\bs{\check\eta})$ with corners. It is automatic that
$\pd\check V^i=\pd V^i=W^{i-1}$, since $(\bs J,\bs\ze)$ is an
excellent coordinate system, so every connected component of
$\pd^kW^{i-1}$ for all $k\ge 0$ is necessary, and Algorithm
\ref{kh3alg} does not discard any of them. Hence $(\bs{\check I},
\bs{\check\eta})\vert_{\pd X}=(\bs J,\bs\ze)$, as we want.
\end{proof}

\begin{thm} Suppose\/ $X$ is a compact Kuranishi space with
corners, $Y$ an orbifold, $\bs f:X\ra Y$ a strongly smooth map,
and\/ $\bs H$ is gauge-fixing data for\/ $(\pd X,\bs f\vert_{\pd
X})$ with corners. Let\/ $\bs\si:\pd^2X\ra\pd^2X$ be the involution
of Definition {\rm\ref{kh2def19}.} Then there exists gauge-fixing
data\/ $\bs G$ for\/ $(X,\bs f)$ with corners with\/ $\bs
G\vert_{\pd X}\!=\!\bs H$ if and only if\/ $\bs H\vert_{\pd^2X}$
is\/ $\bs\si$-invariant.

Suppose\/ $\Ga$ is a finite group of strong diffeomorphisms\/ $\bs
a:X\ra X$ with\/ $\bs f\ci\bs a=\bs f$ and\/ $\bs H$ is invariant
under\/ $\Ga\vert_{\pd X}$, that is, we are given a group morphism\/
$\Ga\ra\Aut(\pd X,\bs f\vert_{\pd X},\bs H)$ mapping\/ $\bs
a\mapsto(\bs a\vert_{\pd X},\bs c_{\bs a}),$ such that the action
of\/ $(\bs a\vert_{\pd^2X},\bs c_{\bs a}\vert_{\pd^2X})$ on\/ $\bs
H\vert_{\pd^2X}$ commutes with the action of\/ $\bs\si$ on\/ $\bs
H\vert_{\pd^2X}$. Then we can choose\/ $\bs G$ to be\/
$\Ga$-invariant, that is, there is a group morphism
$\Ga\ra\Aut(X,\bs f,\bs G)$ mapping $\bs a\mapsto(\bs a,\bs b_{\bs
a}),$ such that\/ $\bs b_{\bs a}\vert_{\pd X}=\bs c_{\bs a}$ for
all\/~$\bs a\in\Ga$.

If\/ $\pd Y=\es,$ all this also holds for co-gauge-fixing data $\bs
C,\bs D$ for\/ $(X,\bs f),$ $(\pd X,\bs f\vert_{\pd X}),$ with $\bs
f$ a strong submersion.
\label{kh3thm3}
\end{thm}

\begin{proof} If $\bs G$ is gauge-fixing data for $(X,\bs f)$ with
corners then restricting to $\pd X$ and then to $\pd^2X$ gives $\bs
G\vert_{\pd^2X}$. Clearly, $\bs G\vert_{\pd^2X}$ is invariant under
$\bs\si$. Thus, $\bs G\vert_{\pd X}=\bs H$ implies that $\bs
H\vert_{\pd^2X}$ is $\bs\si$-invariant. This proves the `only if'
part.

For the `if' part, suppose $\bs H\vert_{\pd^2X}$ is
$\bs\si$-invariant as above. Write $\bs H=(\bs J,\bs\ze,H^j:j\in
J)$, where $\bs J=\bigl(J,(W^j,F^j,t^j,\xi^j):j\in J,\ldots\bigr)$.
As $(\bs J,\bs\ze)$ is an excellent coordinate system with corners
for $(\pd X,\bs f\vert_{\pd X})$ and $(\bs J,\bs\ze)\vert_{\pd^2X}$
is $\bs\si$-invariant, Proposition \ref{kh3prop5} gives an excellent
coordinate system $(\bs I,\bs\eta)$ with corners for $(X,\bs f)$
with $(\bs I,\bs\eta)\vert_{\pd X}=(\bs J,\bs\ze)$, where~$\bs
I=\bigl(I,(V^i,E^i,s^i,\psi^i),i\in I,\ldots\bigr)$.

Then $\pd E^i=F^{i-1}$ for $i\in I$ with $i-1\in J$, and $\pd
E^i=\es$ for $i\in I$ with $i-1\notin J$. When $i\in I$ with $i-1\in
J$, we have $H^{i-1}:\pd E^i\ra P$ with $H^{i-1}\vert_{\pd^2E^i}$
invariant under $\si:\pd^2E^i\ra E^i$. Therefore by Principle
\ref{kh2pri}(a) we can choose $G^i:E^i\ra P_n\subset P$ with
$G^i\vert_{\pd E^i}\equiv H^i$. Since $G^i$ is arbitrary, not
continuous or smooth, the values of $G^i$ on $(E^i)^\ci$ need not be
related to those on~$\pd E^i$.

Choose $G^i\vert_{(E^i)^\ci}:(E^i)^\ci\ra P_n\subset P$ to be an
arbitrary {\it injective} map, for $n\gg 0$. As $H^{i-1}$ is {\it
globally finite} with constant $N\ge 0$ in the sense of Definition
\ref{kh3def6}, it follows that $G^i$ is globally finite with
constant $N+1$. If $i\in I$ but $i-1\notin J$ then there are no
boundary conditions on $G^i$, so we choose $G^i:E^i\ra P_n\subset P$
to be an arbitrary injective map for $n\gg 0$, and then $G^i$ is
globally finite with $N=1$. Thus $\bs G$ is gauge-fixing data with
corners for $(X,\bs f)$ with $\bs G\vert_{\pd X}=\bs H$, proving the
`if' part. For the second part on $\Ga$ actions, we combine the
proof above with the second part of the proof of Theorem
\ref{kh3thm1}. The analogues for co-gauge-fixing data are immediate.
\end{proof}

\subsection{Pushforwards and pullbacks of (co-)gauge-fixing data}
\label{kh37}

Our next definitions, which give Property \ref{kh3pr}(f),(g) for
gauge-fixing data and co-gauge-fixing data, encapsulate the reasons
for the differences between Definitions \ref{kh3def6} and
\ref{kh3def7}. We will use gauge-fixing data in Kuranishi homology,
so it should have a good notion of {\it pushforward\/} $h_*$ under
smooth $h:Y\ra Z$. We will use co-gauge-fixing data in
(compactly-supported) Kuranishi cohomology, so it should have a good
notion of {\it pullback\/} $h^*$ under smooth, proper~$h:Y\ra Z$.

\begin{dfn} Let $Y,Z$ be orbifolds and $h:Y\ra Z$ a smooth map.
Suppose $X$ is a compact Kuranishi space, $\bs f:X\ra Y$ is strongly
smooth, and $\bs G=(\bs I,\bs\eta,G^i:i\in I)$ is gauge-fixing data
for $(X,\bs f)$, where $\bs I=\bigl(I,(V^i,E^i,s^i,\psi^i),\ab
f^i:i\in I,\ldots\bigr)$. Then $h\ci\bs f:X\ra Z$ is strongly
smooth. Define {\it pushforward gauge-fixing data} $h_*(\bs G)$ for
$(X,h\ci\bs f)$ to be $h_*(\bs G)=(h_*(\bs I),\bs\eta,G^i:i\in I)$,
where $h_*(\bs I)=\bigl(I,(V^i,E^i,s^i,\psi^i),h\ci f^i:i\in
I,\ldots\bigr)$. That is, $h_*(\bs G)$ is the same as $\bs G$,
except that $f^i:V^i\ra Y$ for $i\in I$ are replaced by $h\ci
f^i:V^i\ra Z$.

Clearly $h_*(\bs G)$ is gauge-fixing data for $(X,h\ci\bs f)$. Note
however that the same construction does {\it not\/} work for
co-gauge-fixing data, as if $h:Y\ra Z$ is not globally finite, then
$C^i\t(f^i\ci\pi^i):E^i\ra P\t Y$ globally finite does not imply
$C^i\t(h\ci f^i\ci\pi^i):E^i\ra P\t Z$ globally finite, so
Definition \ref{kh3def7} fails for~$h_*(\bs C)$.
\label{kh3def12}
\end{dfn}

\begin{dfn} Let $Y,Z$ be orbifolds without boundary and $h:Y\ra Z$
a smooth, proper map. Suppose $X$ is a compact Kuranishi space, $\bs
f:X\ra Z$ is a strong submersion, and $\bs C=(\bs I,\bs\eta,C^i:i\in
I)$ is co-gauge-fixing data for $(X,\bs f)$, where $\bs
I=\bigl(I,(V^i,E^i,s^i,\psi^i),f^i:i\in I,\ldots\bigr)$. Then
$Y\t_{h,Z,\bs f}X$ is a compact Kuranishi space, and
$\bs\pi_Y:Y\t_ZX\ra Y$ is a strong submersion. We will define
co-gauge-fixing data $h^*(\bs C)$ for $(Y\t_ZX,\bs\pi_Y)$, where
$h^*(\bs C)=(\bs{\check I},\bs{\check\eta},\check C^i:i\in\check
I)$, with~$\bs{\check I}=\bigl(\check I,(\check V^i,\check
E^i,\check s^i,\check\psi^i),\check\pi_Y^i:i\in\check
I,\ldots\bigr)$.

Let $k=\dim Y-\dim Z$, and $\ti I=\{i+k:i\in I\}$. Then $\dim V^i=i$
for each $i\in I$, and as $f^i:V^i\ra Z$ is a submersion $i=\dim
V^i\ge\dim Z$, so $i+k\ge\dim Y\ge 0$. Thus $\ti I\subset\N$. For
$i\in I$, define $(\ti V^{i+k},\ti E^{i+k},\ti
s^{i+k},\ti\psi^{i+k})=\bigl(Y\t_{h,Z,f^i}V^i,\pi_{V^i}^*(E^i),
s^i\ci\pi_{V^i},h\t_Z\psi^i\bigr)$. This is a Kuranishi
neighbourhood on $Y\t_{h,Z,\bs f}X$, with $\dim\ti V^{i+k}=\dim
V^i+\dim Y-Z=i+k$, as we want. Let $\ti\pi^{i+k}_Y:\ti V^{i+k}\ra Y$
be the projection from the fibre product. This is a submersion, as
$f^i$ is, and represents~$\bs\pi_Y$.

For $i,j\in I$ with $j\le i$ we have a triple $(V^{ij},\phi^{ij},
\hat\phi^{ij})$. Define $\ti V^{(i+k)(j+k)}=\pi_{V^j}^{-1}(V^{ij})$,
where $\pi_{V^j}:\ti V^{j+k}\ra V^j$ is the projection from the
fibre product. Define $\ti\phi{}^{(i+k)(j+k)}:\ti V^{(i+k)(j+k)}\ra
\ti V^{i+k}$ by $\ti\phi{}^{(i+k)(j+k)}=h^*(\phi^{ij})$, mapping
$Y\t_ZV^{ij}\ra Y\t_ZV^i$. Define $\hat{\ti\phi}{}^{(i+k)(j+k)}: \ti
E^{j+k}\vert_{\ti V^{(i+k)(j+k)}}\ra(\ti\phi{}^{(i+k)(j+k)})^*(\ti
E^{i+k})$ by $\hat{\ti\phi}{}^{(i+k)(j+k)}=h^*(\hat\phi^{ij})$. Then
$\bs{\ti I}=\bigl(\ti I,(\ti V^i,\ti E^i,\ti
s^i,\ti\psi^i),\ti\pi_Y^i: i\in\ti I,(\ti V^{ij},\ti\phi^{ij},
\hat{\ti\phi}{}^{ij}):j\le i$ in $\ti I\bigr)$ is a {\it very good
coordinate system} for~$(Y\t_ZX,\bs\pi_Y)$.

Now for $i\in I$, define $\ti\eta_{i+k}:Y\t_ZX\ra[0,1]$ by
$\ti\eta_{i+k}=\eta_i\ci\pi_X$, and for $i,j\in I$ define
$\ti\eta_{i+k}^{j+k}:\ti V^{j+k}\ra[0,1]$ by
$\ti\eta_{i+k}^{j+k}=\eta_i^j\ci\pi_{V^j}$. Write
$\bs{\ti\eta}=(\ti\eta_i:i\in\ti I$, $\ti\eta_i^j:i,j\in\ti I)$. It
is not difficult to check that $(\bs{\ti I},\bs{\ti\eta})$ is a {\it
really good coordinate system}, as $(\bs I,\bs\eta)$ is. In
Definition \ref{kh3def3}(v) we can use the same $M$ for $(\bs{\ti
I}, \bs{\ti\eta})$ as  for $(\bs I,\bs\eta)$, and if $T\subseteq Z$
is as in Definition \ref{kh3def3}(vi) for $(\bs I,\bs\eta)$ then
$\ti T=h^{-1}(T)\subseteq Y$ is compact as $T$ is compact and $h$ is
proper, and satisfies Definition \ref{kh3def3}(vi) for~$(\bs{\ti
I},\bs{\ti\eta})$.

However, although $(\bs I,\bs\eta)$ is an {\it excellent coordinate
system}, $(\bs{\ti I},\bs{\ti\eta})$ need not be, as taking fibre
products $V^i\mapsto Y\t_ZV^i=\ti V^{i+k}$ may introduce unnecessary
connected components of $\ti V^{i+k}$ or $\pd^l\ti V^{i+k}$. Let
$(\bs{\check I},\bs{\check\eta})$ be the excellent coordinate system
obtained from $(\bs{\ti I},\bs{\ti\eta})$ by applying Algorithm
\ref{kh3alg}, where~$\bs{\check I}=\bigl(\check I,(\check V^i,\check
E^i,\check s^i,\check\psi^i),\check\pi_Y^i:i\in\check
I,\ldots\bigr)$.

Let $i\in\check I$, so that $i-k\in I$. Then $\check V^i$ is an open
subset of $\ti V^i=Y\t_{h,Z,f^{i-k}}V^{i-k}$, and $\check
E^i\ra\check V^i$ is the pullback of $E^{i-k}\ra V^{i-k}$ by the
projection $\pi_{V^{i-k}}:\check V^i\ra V^{i-k}$. Thus there is a
natural projection $\pi_{E^{i-k}}:\check E^i\ra E^{i-k}$ between the
total spaces of the bundles $\check E^i,E^{i-k}$. Define $\check
C^i:\check E^i\ra P$ by $\check C^i=C^{i-k}\ci \pi_{E^{i-k}}$. We
shall show that $\check C^i\t(\check\pi_Y^i\ci\check\pi^i):\check
E^i\ra P\t Y$ for $i\in\check I$ is {\it globally finite}.

Definition \ref{kh3def3}(vi) gives compact $T\subseteq Z$ with
$f^{i-k}(V^{i-k})\subseteq T$. The map $Z\ra\N$ taking
$z\mapsto\md{\Stab(z)}$ is upper semicontinuous, so as $T$ is
compact there exists $L\ge 1$ such that $\md{\Stab(z)}\le L$ for all
$z\in T$. As $C^{i-k}\t(f^{i-k}\ci \pi^{i-k}):E^{i-k}\ra P\t Z$ is
globally finite, there exists $N\ge 0$ such that
$\bmd{(C^{i-k}\t(f^{i-k}\ci\pi^{i-k}))^{-1}(p,z)}\le N$ for all
$(p,z)\in P\t Z$. Let $(p,y)\in P\t Y$, and $z=h(y)$. Now
$\pi_{E^{i-k}}:\check E^i\ra E^{i-k}$ restricts to a map
\e
\pi_{E^{i-k}}:(\check C^i\t(\check\pi_Y^i
\ci\check\pi^i))^{-1}(p,y)\longra
(C^{i-k}\t(f^{i-k}\ci\pi^{i-k}))^{-1}(p,z).
\label{kh3eq5}
\e

Using $\check E^i\subseteq Y\t_{h,Z,f^{i-k}\ci\pi^{i-k}}E^{i-k}$ and
the definition \eq{kh2eq10} of fibre products of orbifolds, we see
that each $e$ on the right hand side of \eq{kh3eq5} pulls back to a
subset of $h_*(\Stab(y))\backslash\Stab(z)/(f^{i-k}\ci\pi^{i-k})_*
(\Stab(e))$ on the left. Since $\md{\Stab(z)}\le L$, at most $L$
points on the left of \eq{kh3eq5} project to each point on the
right. But the right hand side of \eq{kh3eq5} has at most $N$
elements, so the left has at most $LN$ elements. Thus $\check
C^i\t(\check\pi_Y^i\ci\check\pi^i)$ is globally finite, with
constant $\check N=LN$. Therefore $h^*(\bs C)$ is co-gauge-fixing
data for~$(Y\t_ZX,\bs\pi_Y)$.

Note that the same construction does {\it not\/} work for
gauge-fixing data, as if $h$ is not globally finite then $G^i:E^i\ra
P$ globally finite does not imply $G^i\ci\pi_{E^i}:Y\t_ZE^i\ra P$ is
globally finite, so $\check G^i:\check E^i\!\ra\! P$ may not be
globally finite, and Definition \ref{kh3def6} fails for~$h^*(\bs
G)$.
\label{kh3def13}
\end{dfn}

\begin{rem} We suppose $Y,Z$ are {\it without boundary\/} in
Definition \ref{kh3def13} in order that $\pd(Y\t_ZX)\cong\pm
Y\t_Z\pd X$. This will be needed in \S\ref{kh44} to ensure that
pullbacks $h^*$ on Kuranishi cochains satisfy $\d\ci h^*=h^*\ci\d$.
If $\pd Y\ne\es$ then $\pd(Y\t_ZX)$ has an extra term $\pd Y\t_ZX$,
which would invalidate the proof of $\d\ci h^*=h^*\ci\d$. Actually
the assumption $\pd Z=\es$ is unnecessary, since our definitions of
(strong) submersions imply that $\bs f:X\ra Z$ and $f^i:V^i\ra Z$
for $i\in I$ map to $Z^\ci$, so $\pd Z$ does not affect
$\pd(Y\t_ZX)$ or $\pd(Y\t_{h,Z,f^i}V^i)$.
\label{kh3rem3}
\end{rem}

\subsection{Fibre products of (co-)gauge-fixing data}
\label{kh38}

When $(\bs I,\bs\eta)$ and $(\bs{\ti I},\bs{\ti\eta})$ are really
good coordinate systems for $(X,\bs f)$ and $(\ti X,\bs{\ti f})$,
where one of $\bs f:X\ra Y$ and $\bs{\ti f}:\ti X\ra Y$ is a strong
submersion, we shall define a really good coordinate system $(\bs
I,\bs\eta)\t_Y(\bs{\ti I},\bs{\ti\eta})$ for the fibre product
$(X\t_Y\ti X,\bs\pi_Y)$. The construction is long and complicated.

\begin{dfn} Let $X,\ti X$ be compact Kuranishi spaces, $Y$ be an
orbifold, $\bs f:X\ra Y$ and $\bs{\ti f}:\ti X\ra Y$ be strongly
smooth and $(\bs I,\bs\eta)$, $(\bs{\ti I},\bs{\ti\eta})$ be really
good coordinate systems for $(X,\bs f)$ and $(\ti X,\bs{\ti f})$.
Suppose one of $\bs f,\bs{\ti f}$ is a strong submersion, and the
corresponding $f^i:V^i\ra Y$ or $\ti f^{\ti\imath}:\ti
V^{\ti\imath}\ra Y$ in $\bs I,\bs{\ti I}$ are submersions. Then
$X\t_Y\ti X$ is a compact Kuranishi space, as in \S\ref{kh26}. As in
\eq{kh2eq14}, for $i\in I$ and $\ti\imath\in\ti I$ define a
Kuranishi neighbourhood on $X\t_Y\ti X$ by
\e
\begin{split}
\bigl(V^{(i,\ti\imath)},E^{(i,\ti\imath)},s^{(i,\ti\imath)}&,
\psi^{(i,\ti\imath)}\bigr)= \bigl(V^i\t_{f^i,Y,\ti f^{\ti\imath}}
\ti V^{\ti\imath},\pi_{V^i}^*(E^i)\op
\pi_{\ti V^{\ti\imath}}^*(\ti E^{\ti\imath}),\\
&(s^i\ci\pi_{V^i})\op(\ti s^{\ti\imath}\ci\pi_{\ti V^{\ti\imath}}),
(\psi^i\ci\pi_{V^i})\t(\ti\psi^{\ti\imath}\ci\pi_{\ti
V^{\ti\imath}})\t \chi^{(i,\ti\imath)}\bigr).
\end{split}
\label{kh3eq6}
\e
Let $\pi^{(i,\ti\imath)}_Y:V^{(i,\ti\imath)}\ra Y$ be the projection
from the fibre product. If $i,j\in I$ and $\ti\imath,\ti\jmath\in\ti
I$ with $j\le i$ and $\ti\jmath\le\ti\imath$ we have open sets
$V^{ij}\subseteq V^j$ and $\ti V^{\ti\imath\ti\jmath}\subseteq\ti
V^{\ti\jmath}$ and coordinate transformations
$(\phi^{ij},\hat\phi^{ij})$ and
$(\ti\phi^{\ti\imath\ti\jmath},\hat{\ti\phi}{}^{\ti\imath\ti\jmath})$.
Define $V^{(i,\ti\imath)(j,\ti\jmath)}\subset V^{(j,\ti\jmath)}$ to
be $V^{ij}\t_{f^j,Y,\ti f^{\ti\jmath}}\ti V^{\ti\imath\ti\jmath}$.
Then $(\phi^{ij},\hat\phi^{ij})$ and
$(\ti\phi^{\ti\imath\ti\jmath},\hat{\ti\phi}{}^{\ti\imath\ti\jmath})$
induce a coordinate transformation
\begin{align*}
\bigl(\phi^{(i,\ti\imath)(j,\ti\jmath)},\hat\phi^{(i,\ti\imath)
(j,\ti\jmath)}\bigr):\bigl(&V^{(i,\ti\imath)(j,\ti\jmath)},
E^{(j,\ti\jmath)}\vert_{V^{(i,\ti\imath)(j,\ti\jmath)}},
s^{(j,\ti\jmath)}\vert_{V^{(i,\ti\imath)(j,\ti\jmath)}},\\
&\psi^{(j,\ti\jmath)}\vert_{V^{(i,\ti\imath)(j,\ti\jmath)}}\bigr)
\longra\bigl(V^{(i,\ti\imath)},\ldots,\psi^{(i,\ti\imath)}\bigr).
\end{align*}

We want to build a really good coordinate system $(\bs{\dot
I},\bs{\dot\eta})=(\bs I,\bs\eta)\t_Y(\bs{\ti I},\bs{\ti\eta})$ for
$(X\t_Y\ti X,\bs\pi_Y)$ with Kuranishi neighbourhoods $(\dot
V^k,\ldots,\dot\psi^k)$ for $k\in\dot I$, using the Kuranishi
neighbourhoods \eq{kh3eq6} and maps $\pi^{(i,\ti\imath)}_Y$. Since
$\dim V^{(i,\ti\imath)}=\dim V^i+\dim\ti V^{\ti\imath}-\dim
Y=i+\ti\imath-\dim Y$, the obvious definition is $\dot
V^k=\coprod_{i\in I,\; \ti\imath\in\ti I:i+\ti\imath=k+\dim
Y}V^{(i,\ti\imath)}$, with $(\dot E^k,\dot
s^k,\dot\psi^k)\vert_{V^{(i,\ti\imath)}}=
(E^{(i,\ti\imath)},s^{(i,\ti\imath)},\psi^{(i,\ti\imath)})$ and
$\dot\pi_Y^k\vert_{V^{(i,\ti\imath)}}=\pi^{(i,\ti\imath)}_Y$. There
are two main problems with this:
\begin{itemize}
\setlength{\itemsep}{0pt}
\setlength{\parsep}{0pt}
\item[$(*)$] If $i\ne j\in I$ and $\ti\imath\ne \ti\jmath\in\ti I$ with
$i+\ti\imath=j+\ti\jmath=k+\dim Y$ and $\Im\psi^{(i,\ti\imath)}\cap
\Im\psi^{(j,\ti\jmath)}\ne\es$ then $\dot\psi^k$ is not injective,
so $(\dot V^k,\ldots,\dot\psi^k)$ is not a Kuranishi neighbourhood;
and
\item[$(**)$] Suppose $i<j\in I$ and $\ti\imath>\ti\jmath\in\ti I$ with
$i+\ti\imath=k+\dim Y$ and $j+\ti\jmath=l+\dim Y$ with $l<k$, and
that $\Im\psi^{(i,\ti\imath)}\cap\Im\psi^{(j,\ti\jmath)} \ne\es$.
Then we expect to define an open $\dot V^{kl}\subset\dot V^l$ and a
coordinate transformation $(\dot\phi^{kl},\hat{\dot\phi}{}^{kl})$
from $(\dot V^{kl},\dot E^l\vert_{\dot V^{kl}},\dot s^l\vert_{\dot
V^{kl}},\dot\psi^l \vert_{\dot V^{kl}})$ to $(\dot
V^k,\ldots,\dot\psi^k)$. This should include a transformation from a
nonempty subset of $V^{(j,\ti\jmath)}$ to $V^{(i,\ti\imath)}$.
However, since $i<j$ we have no transformation from $V^{ij}\subseteq
V^j$ to $V^i$ to build this from.
\end{itemize}

We shall deal with both these problems by making the
$V^{(i,\ti\imath)}$ smaller. That is, we shall define
\e
\dot V^k=\ts\coprod_{i\in I,\; \ti\imath\in\ti I:i+\ti\imath=k+\dim
Y}\dot V^{(i,\ti\imath)},
\label{kh3eq7}
\e
where $\dot V^{(i,\ti\imath)}$ is an open subset of
$V^{(i,\ti\imath)}$. We shall choose these $\dot V^{(i,\ti\imath)}$
so that in $(*)$ and $(**)$ above we always have
$\Im\psi^{(i,\ti\imath)}\vert_{\dot V^{(i,\ti\imath)}}\cap
\Im\psi^{(j,\ti\jmath)}\vert_{\dot V^{(j,\ti\jmath)}}=\es$. But we
cannot make the $\dot V^{(i,\ti\imath)}$ too small, as we need to
ensure that $X=\bigcup_{i\in I,\;\ti\imath\in\ti
I}\Im\psi^{(i,\ti\imath)}\vert_{\dot V^{(i,\ti\imath)}}$. To define
the $\dot V^{(i,\ti\imath)}$, we will construct partitions of unity
on $X\t_Y\ti X$ and on the $V^{(i,\ti\imath)}$ for all $i,\ti\imath$
using the data $\bs\eta,\bs{\ti\eta}$. We will also use these
partitions of unity to construct $\bs{\dot\eta}$ in the really good
coordinate system $(\bs{\dot I},\bs{\dot\eta})$ for~$(X\t_Y\ti
X,\bs\pi_Y)$.

For all $i,j\in I$ and $\ti\imath,\ti\jmath\in\ti I$ we will define
continuous functions $\eta_{(i,\ti\imath)}:X\t_Y\ti X\ra[0,1]$ and
$\eta_{(i,\ti\imath)}^{(j,\ti\jmath)}:V^{(j,\ti\jmath)}\ra[0,1]$
satisfying the conditions:
\begin{itemize}
\setlength{\itemsep}{0pt}
\setlength{\parsep}{0pt}
\item[(a)] $\sum_{i\in I,\;\ti\imath\in\ti I}\eta_{(i,\ti\imath)}
\!\equiv\!1$ on $X\!\t_Y\!\ti X$, and $\sum_{i\in I,\;\ti\imath
\in\ti I}\eta_{(i,\ti\imath)}^{(j,\ti\jmath)}\!\equiv\!1$ on
$V^{(j,\ti\jmath)}$ for all~$j,\ti\jmath$.
\item[(b)] $\eta_{(i,\ti\imath)}^{(j,\ti\jmath)}\vert_{(s^{(j,
\ti\jmath)})^{-1}(0)}\equiv\eta_{(i,\ti\imath)}\ci\psi^{(j,
\ti\jmath)}$ for all $i,j\in I$ and $\ti\imath,\ti\jmath\in\ti I$.
\item[(c)] Let $i,j,k\in I$ and $\ti\imath,\ti\jmath,\ti k\in\ti I$
with $k\le j$ and $\ti k\le\ti\jmath$. Then $\eta_{(i,\ti\imath)}^{
(k,\ti k)}\vert_{V^{(j,\ti\jmath)(k,\ti k)}}\equiv\eta_{(i,
\ti\imath)}^{(j,\ti\jmath)}\ci\phi^{(j,\ti\jmath)(k,\ti k)}$.
\item[(d)] If $j\le i$ in $I$, $\ti\jmath\le\ti\imath$ in $\ti I$
and $q\in X\t_Y\ti X$ with $\eta_{(i,\ti\imath)}(q)>0$,
$\eta_{(j,\ti\jmath)}(q)>0$ and $q\in\Im\psi^{(i,\ti\imath)}$
then~$q\in\Im\psi^{(j,\ti\jmath)}$.
\item[(e)] if $j\le i$ in $I$ and $\ti\jmath\le\ti\imath$ in $\ti I$
then $\phi^{(i,\ti\imath)(j,\ti\jmath)}(V^{(i,\ti\imath)(j,
\ti\jmath)})$ is a {\it closed\/} subset of $\{v\in
V^{(i,\ti\imath)}:\eta^{(i,\ti\imath)}_{(j,\ti\jmath)}(v)>0\bigr\}$.
If $k\le j\le i$ in $I$ and $\ti k\le\ti\jmath\le\ti\imath$ in $\ti
I$ and $v\in V^{(k,\ti k)}$ with $\eta_{(i,\ti\imath)}^{(k,\ti
k)}(v)>0$, $\eta_{(j,\ti\jmath)}^{(k,\ti k)}(v)>0$ and $v\in
V^{(i,\ti\imath)(k,\ti k)}$ then $v\in V^{(j,\ti\jmath)(k,\ti k)}$.
\item[(f)] Suppose $i,j,k\in I$ and $\ti\imath,\ti\jmath,\ti k\in\ti
I$ and either $i<j$ and $\ti\imath>\ti\jmath$, or $i>j$ and
$\ti\imath<\ti\jmath$. Then $\eta_{(i,\ti\imath)}
\eta_{(j,\ti\jmath)}\equiv 0$ on $X\t_Y\ti X$, and $\eta_{(i,
\ti\imath)}^{(k,\ti k)}\eta_{(j,\ti\jmath)}^{(k,\ti k)}\equiv 0$ on
$V^{(k,\ti k)}$. That is, the supports of $\eta_{(i,\ti\imath)}$ and
$\eta_{(j,\ti\jmath)}$ are disjoint, and the supports of
$\eta_{(i,\ti\imath)}^{(k,\ti k)}$ and $\eta_{(j,\ti\jmath)}^{(k,\ti
k)}$ are disjoint.
\item[(g)] Suppose $q\in X\t_YX$, set $p=\pi_X(q)$ and $\ti
p=\pi_{\ti X}q$, and let $i\in I$ and $\ti\imath\in\ti I$ be
least with $\eta_i(p)>0$ and $\ti\eta_{\ti\imath}(\ti p)>0$.
Then $\eta_{(i,\ti\imath)}(q)>0$
and~$q\in\Im\psi^{(i,\ti\imath)}$.
\end{itemize}
With the exception of (f),(g) these are analogues of Definition
\ref{kh3def3}(i)--(iv); note that (d),(e) are partial analogues of
the parts of Definition~\ref{kh3def3}(i),(ii).

Suppose for the moment we have constructed such
$\eta_{(i,\ti\imath)},\eta_{(i,\ti\imath)}^{(j,\ti\jmath)}$. Define
\e
\dot V^{(i,\ti\imath)}=\bigl\{v\in V^{(i,\ti\imath)}:
\eta_{(i,\ti\imath)}^{(i,\ti\imath)}(v)>0\bigr\}\subseteq
V^{(i,\ti\imath)}.
\label{kh3eq8}
\e
Then $\dot V^{(i,\ti\imath)}$ is open, as $\eta_{(i,\ti\imath)
}^{(i,\ti\imath)}$ is continuous, and (b) with $j=i$,
$\ti\jmath=\ti\imath$ implies that
\e
\begin{split}
\Im\psi^{(i,\ti\imath)}\vert_{\dot V^{(i,\ti\imath)}}&=\bigl\{q\in
\pi_X^{-1}(\Im\psi^i)\cap\pi_{\ti X}^{-1}(\Im\ti\psi^{\ti\imath}):
\eta_{(i,\ti\imath)}(q)>0\bigr\}\\
&\subseteq\bigl\{q\in X\t_Y\ti X:\eta_{(i,\ti\imath)}(q)>0\bigr\}.
\end{split}
\label{kh3eq9}
\e
Define $\dot I=\bigl\{i+\ti\imath-\dim Y:i\in I$, $\ti\imath\in\ti
I$, $\dot V^{(i,\ti\imath)}\ne\es\bigr\}$. Since one of $\bs
f,\bs{\ti f}$ is a strong submersion, it follows that either
$i\ge\dim Y$ for all $i\in I$ or $\ti\imath\ge\dim Y$ for all
$\ti\imath\in\ti I$, and so $\dot I$ is a finite subset of
$\N=\{0,1,2,\ldots\}$. For each $k\in\dot I$, define $\dot V^k$ by
\eq{kh3eq7}. Then $\dot V^k$ is an orbifold with $\dim\dot V^k=k$.
Define an orbibundle $\dot E^k\ra\dot V^k$, a smooth section $\dot
s^k$ of $\dot E^k$, a continuous map $\dot\psi^k:(\dot
s^k)^{-1}(0)\ra X\t_Y\ti X$ and a smooth map $\dot\pi_Y^k:\dot
V^k\ra Y$ by $(\dot E^k,\dot s^k,\dot\psi^k)\vert_{\dot
V^{(i,\ti\imath)}}=(E^{(i,\ti\imath)},s^{(i,\ti\imath)},
\psi^{(i,\ti\imath)})\vert_{\dot V^{(i,\ti\imath)}}$,
$\dot\pi_Y^k\vert_{\dot
V^{(i,\ti\imath)}}=\pi^{(i,\ti\imath)}_Y\vert_{\dot
V^{(i,\ti\imath)}}$ for all $i,\ti\imath$ in~\eq{kh3eq7}.

We claim that $(\dot V^k,\ldots,\dot\psi^k)$ is a {\it Kuranishi
neighbourhood\/} on $X\t_Y\ti X$. Since it is a disjoint union of
open subsets of Kuranishi neighbourhoods \eq{kh3eq6}, the only thing
to check is that $\dot\psi^k$ is injective, the problem identified
in $(*)$ above. It is enough to show that if $i,j\in I$ and
$\ti\imath,\ti\jmath\in\ti I$ with $(i,\ti\imath)\ne(j,\ti\jmath)$
and $i+\ti\imath=j+\ti\jmath=k+\dim Y$ then
$\Im\psi^{(i,\ti\imath)}\vert_{\dot V^{(i,\ti\imath)}}\cap
\Im\psi^{(j,\ti\jmath)}\vert_{\dot V^{(j,\ti\jmath)}}=\es$. But
these conditions imply that either $i<j$ and $\ti\imath>\ti\jmath$,
or $i>j$ and $\ti\imath<\ti\jmath$, so (f) implies that the supports
of $\eta_{(i,\ti\imath)}$ and $\eta_{(j,\ti\jmath)}$ are disjoint,
and thus $\Im\psi^{(i,\ti\imath)}\vert_{\dot V^{(i,\ti\imath)}},
\Im\psi^{(j,\ti\jmath)}\vert_{\dot V^{(j,\ti\jmath)}}$ are disjoint
by \eq{kh3eq9}. Hence $(\dot V^k,\ldots,\dot\psi^k)$ is a Kuranishi
neighbourhood on~$X\t_Y\ti X$.

Next, for $j\le i$ in $I$ and $\ti\jmath\le\ti\imath$ in $\ti I$
define
\e
\dot V^{(i,\ti\imath)(j,\ti\jmath)}=
V^{(i,\ti\imath)(j,\ti\jmath)}\cap \dot V^{(j,\ti\jmath)}\cap
(\phi^{(i,\ti\imath)(j,\ti\jmath)})^{-1}(\dot V^{(i,\ti\imath)}).
\label{kh3eq10}
\e
For $k,l\in\dot I$ with $l\le k$, define
\e
\dot V^{kl}=\ts\bigcup_{\begin{subarray}{l} i,j\in I,\;
\ti\imath,\ti\jmath\in\ti I:\; j\le i,\; \ti\jmath\le\ti\imath,\\
i+\ti\imath=k+\dim Y,\; j+\ti\jmath=l+\dim Y \end{subarray}} \dot
V^{(i,\ti\imath)(j,\ti\jmath)}\,\,\subseteq \dot V^l.
\label{kh3eq11}
\e
We claim that \eq{kh3eq11} is a {\it disjoint union}, that is,
distinct $\dot V^{(i,\ti\imath)(j,\ti\jmath)}$ in \eq{kh3eq11} do
not intersect. Since $\dot V^{(i,\ti\imath)(j,\ti\jmath)}\subseteq
\dot V^{(j,\ti\jmath)}$ and \eq{kh3eq7} is a disjoint union, it is
enough to show that if $i,i',j\in I$ and
$\ti\imath,\ti\imath',\ti\jmath\in\ti I$ with
$(i,\ti\imath)\ne(i',\ti\imath')$, $j\le i$, $j\le i'$,
$\ti\jmath\le\ti\imath$, $\ti\jmath\le\ti\imath'$,
$i+\ti\imath=i'+\ti\imath'=k+\dim Y$ and $j+\ti\jmath=l+\dim Y$ then
$V^{(i,\ti\imath)(j,\ti\jmath)}\cap
V^{(i',\ti\imath')(j,\ti\jmath)}=\es$. Now
$\eta_{(i,\ti\imath)}^{(i,\ti\imath)}>0$ on $\dot V^{(i,\ti\imath)}$
by \eq{kh3eq8}, and $\dot V^{(i,\ti\imath)(j,\ti\jmath)}\subseteq
(\phi^{(i,\ti\imath)(j,\ti \jmath)})^{-1}(\dot V^{(i,\ti\imath)})$
by \eq{kh3eq10}, and
$\eta_{(i,\ti\imath)}^{(j,\ti\jmath)}\vert_{V^{(i,\ti\imath)(j,\ti
\jmath)}}\equiv\eta_{(i,\ti\imath)}^{(i,\ti\imath)}
\ci\phi^{(i,\ti\imath)(j,\ti\jmath)}$ by (c), so
$\eta_{(i,\ti\imath) }^{(j,\ti\jmath)}>0$ on
$V^{(i,\ti\imath)(j,\ti\jmath)}$, and similarly
$\eta_{(i',\ti\imath')}^{(j,\ti\jmath)}>0$ on
$V^{(i',\ti\imath')(j,\ti\jmath)}$. Hence $\eta_{(i,\ti\imath)
}^{(j,\ti\jmath)}\eta_{(i',\ti\imath')}^{(j,\ti\jmath)}>0$ on
$V^{(i,\ti\imath)(j,\ti\jmath)}\cap V^{(i',\ti\imath')(j,
\ti\jmath)}$. But the conditions imply that either $i<i'$ and
$\ti\imath>\ti\imath'$, or $i>i'$ and $\ti\imath<\ti\imath'$, so (f)
gives $\eta_{(i,\ti\imath)}^{(j,\ti\jmath)}\eta_{(i',\ti\imath')
}^{(j,\ti\jmath)}\equiv 0$ on $V^{(j,\ti\jmath)}$. Thus
$V^{(i,\ti\imath)(j,\ti\jmath)}\cap
V^{(i',\ti\imath')(j,\ti\jmath)}=\es$, and \eq{kh3eq11} is a
disjoint union.

Define a coordinate transformation
\e
\bigl(\dot\phi^{kl},\hat{\dot\phi}{}^{kl}\bigr):\bigl(\dot
V^{kl},\dot E^l\vert_{\dot V^{kl}},\dot s^l\vert_{\dot
V^{kl}},\dot\psi^l\vert_{\dot V^{kl}}\bigr)\longra (\dot
V^k,\ldots,\dot\psi^k)
\label{kh3eq12}
\e
by $\bigl(\dot\phi^{kl},\hat{\dot\phi}{}^{kl}\bigr)\vert_{ \dot
V^{(i,\ti\imath)(j,\ti\jmath)}}=\bigl(\phi^{(i,\ti\imath)
(j,\ti\jmath)},\hat\phi^{(i,\ti\imath)(j,\ti\jmath)}\bigr)\vert_{
\dot V^{(i,\ti\imath)(j,\ti\jmath)}}$ for all
$i,j,\ti\imath,\ti\jmath$ in \eq{kh3eq11}. This is well-defined as
\eq{kh3eq11} is a disjoint union. Since each $\bigl(\phi^{
(i,\ti\imath) (j,\ti\jmath)},\hat\phi^{(i,\ti\imath)(j,\ti\jmath)}
\bigr)$ is a coordinate transformation, to check \eq{kh3eq12} is a
coordinate transformation it only remains to check that
$\dot\phi^{kl}$ in \eq{kh3eq12} is injective. It is enough to show
that if $i,j,j'\in I$ and $\ti\imath,\ti\jmath,\ti\jmath'\in\ti I$
with $(j,\jmath)\ne(j',\jmath')$, $j\le i$, $j'\le i$,
$\ti\jmath\le\ti\imath$, $\ti\jmath'\le\ti\imath$,
$i+\ti\imath=k+\dim Y$ and $j+\ti\jmath=j'+\ti\jmath'=l+\dim Y$ then
$\phi^{(i,\ti\imath)(j,\ti\jmath)}(\dot V^{(i,\ti\imath)
(j,\ti\jmath)})\cap \phi^{(i,\ti\imath)(j',\ti\jmath')}(\dot
V^{(i,\ti\imath)(j',\ti\jmath')})=\es$ in $\dot V^{(i,\ti\imath)}$.
This holds as by a similar argument
$\eta^{(i,\ti\imath)}_{(j,\ti\jmath)}>0$ on
$\phi^{(i,\ti\imath)(j,\ti\jmath)}(\dot V^{(i,\ti\imath)
(j,\ti\jmath)})$ and $\eta^{(i,\ti\imath)}_{(j',\ti\jmath')}>0$ on
$\phi^{(i,\ti\imath)(j',\ti\jmath')}(\dot V^{(i,\ti\imath)
(j',\ti\jmath')})$, but $\eta^{(i,\ti\imath)}_{(j,\ti\jmath)}
\eta^{(i,\ti\imath)}_{(j',\ti\jmath')}\equiv 0$ on
$\dot V^{(i,\ti\imath)}$.

We claim that $\bs{\dot I}=\bigl(\dot I,(\dot V^k,\ldots,
\dot\psi^k),\dot\pi_Y^k:k\in\dot I,(\dot V^{kl},\dot\phi^{kl},
\hat{\dot\phi}{}^{kl}):l\le k\in\dot I\,\bigr)$ is a {\it very good
coordinate system\/} for $(X\t_Y\ti X,\bs\pi_Y)$. Only two
nontrivial things remain to check. Firstly, (g) implies that every
$q\in X\t_Y\ti X$ lies in $\Im\psi^{(i,\ti\imath)}\vert_{\dot
V^{(i,\ti\imath)}}$ for some $i,\ti\imath$ and hence in
$\Im\dot\psi^k$ for some $k$, so $X\t_Y\ti X=\bigcup_{k\in\dot
I}\Im\dot\psi^k$, as we need for $\bs{\dot I}$ to be a good
coordinate system. Secondly, we must show that $\dot V^{kl}$ is an
open neighbourhood of $(\dot\psi^l)^{-1}(\Im\dot\psi^k)$ in $\dot
V^l$ whenever $l\le k$ in $\dot I$. It is enough to show that if
$i,j\in I$ and $\ti\imath,\ti\jmath\in\ti I$ with
$i+\ti\imath=k+\dim Y$ and $j+\ti\jmath=l+\dim Y$ then:
\begin{itemize}
\setlength{\itemsep}{0pt}
\setlength{\parsep}{0pt}
\item[(A)] if $j\!\le\!i$ and $\ti\jmath\!\le\!\ti\imath$ then $\dot
V^{(i,\ti\imath)(j,\ti\jmath)}$ is an open neighbourhood of
$(\psi^{(j,\ti\jmath)}\vert_{\dot V^{(j,\ti\jmath)}})^{-1}\ab
(\Im\psi^{(i,\ti\imath)}\vert_{\dot V^{(i,\ti\imath)}})$ in $\dot
V^{(j,\ti\jmath)}$, and
\item[(B)] if either $j>i$ and $\ti\jmath<\ti\imath$, or $j<i$ and
$\ti\jmath>\ti\imath$, then $(\Im\psi^{(i,\ti\imath)}\vert_{\dot
V^{(i,\ti\imath)}})\cap (\Im\psi^{(j,\ti\jmath)}\vert_{\dot
V^{(j,\ti\jmath)}})=\es$.
\end{itemize}

Part (A) is immediate from \eq{kh3eq10} and the fact that
$V^{(i,\ti\imath)(j,\ti\jmath)}$ is an open neighbourhood of
$(\psi^{(j,\ti\jmath)})^{-1} (\Im\psi^{(i,\ti\imath)})$ in
$V^{(j,\ti\jmath)}$, which follows from the corresponding statements
for $V^{ij},\psi^j,\psi^i$ in $X$ and $\ti V^{\ti\imath\ti\jmath},
\ti\psi^{\ti\jmath},\ti\psi^{\ti \imath}$ in $\ti X$. Part (B) is
the problem identified in $(**)$ above. To prove it, note that (f)
gives $\eta_{(i,\ti\imath)}\eta_{(j,\ti\jmath)}\equiv 0$ on
$X\t_Y\ti X$, but \eq{kh3eq9} yields $\eta_{(i,\ti\imath)}>0$ on
$\Im\psi^{(i,\ti\imath)} \vert_{\dot V^{(i,\ti\imath)}}$ and
$\eta_{(j,\ti\jmath)}>0$ on $\Im\psi^{(j,\ti\jmath)}\vert_{\dot
V^{(j,\ti\jmath)}}$. Therefore $\bs{\dot I}$ is a very good
coordinate system.

Next we define the partition of unity data $\bs{\dot\eta}$ for
$\bs{\dot I}$. For $k,l\in\dot I$, define continuous
$\dot\eta_k:X\t_Y\ti X\ra[0,1]$ and $\dot\eta_k^l:\dot V^l\ra[0,1]$
by
\e
\dot\eta_k=\ts\sum_{i\in I,\; \ti\imath\in\ti I:i+\ti\imath=k+\dim
Y}\eta_{(i,\ti\imath)},
\label{kh3eq13}
\e
and for all $j\in I$, $\ti\jmath\in\ti I$ with $j+\ti\jmath=l+\dim
Y$,
\e
\dot\eta_k^l\vert_{\dot V^{(j,\ti\jmath)}}=\ts\sum_{i\in I,\;
\ti\imath\in\ti I:i+\ti\imath=k+\dim
Y}\eta_{(i,\ti\imath)}^{(j,\ti\jmath)}.
\label{kh3eq14}
\e
All of Definition \ref{kh3def3}(i)--(iv) except the last sentence of
(i) and the last two sentences of (ii) follow easily from (a)--(c)
above, and these parts of (i),(ii) follow from (d),(e) respectively.
If $M,\ti M,T,\ti T$ are as in Definition \ref{kh3def3}(v),(vi) for
$(\bs I,\bs\eta),(\bs{\ti I},\bs{\ti\eta})$, it is easy to see that
Definition \ref{kh3def3}(v),(vi) hold for $(\bs{\dot
I},\bs{\dot\eta})$ with constant $\dot M=M+\ti M$ and compact set
$\dot T=T\cap\ti T$. Hence $(\bs{\dot I},\bs{\dot\eta})$ is a {\it
really good coordinate system}, which we will also write as~$(\bs
I,\bs\eta)\t_Y(\bs{\ti I},\bs{\ti\eta})$.

It remains to define the continuous functions $\eta_{(i,\ti\imath)}$
and $\eta_{(i,\ti\imath)}^{(j,\ti\jmath)}$ satisfying (a)--(g)
above. The obvious definitions for these are
$\eta_{(i,\ti\imath)}=(\eta_i\ci\pi_X)\cdot(\ti\eta_{\ti\imath}\ci
\pi_{\ti X})$ and $\eta_{(i,\ti\imath)}^{(j,\ti\jmath)}=(\eta_i^j\ci
\pi_{V^j})\cdot(\ti\eta_{\ti\imath}^{\ti\jmath}\ci\pi_{\ti V^{\ti
\jmath}})$, but these will not do as they satisfy (a)--(e) and (g)
but not (f), which was crucial in solving problems $(*)$ and $(**)$
above.

There may be many ways to define $\eta_{(i,\ti\imath)},
\eta_{(i,\ti\imath)}^{(j,\ti\jmath)}$ satisfying (a)--(g). The
definition we give was chosen very carefully to ensure that Theorem
\ref{kh3thm4} below holds, so that fibre products of really good
coordinate systems are {\it commutative\/} and {\it associative}. It
seems remarkable that such a definition is possible.

Before defining the $\eta_{(i,\ti\imath)},
\eta_{(i,\ti\imath)}^{(j,\ti\jmath)}$ we set up some notation. For
each $k\ge 0$, define the $k$-{\it simplex\/} to be
\begin{equation*}
\De_k=\bigl\{(x_0,\ldots,x_k)\in\R^{k+1}:x_i\ge 0,\;\>
x_0+\cdots+x_k=1\bigr\},
\end{equation*}
as in \eq{kh4eq1} below. Write $I=\{i_0,i_1,\ldots,i_k\}$ with $0\le
i_0<i_1<\cdots<i_k$ and $\ti I=\{\ti\imath_0,\ti\imath_1,
\ldots,\ti\imath_l\}$ with $0\le\ti\imath_0<\ti\imath_1<\cdots<
\ti\imath_l$, so that $\md{I}=k+1$ and $\md{\ti I}=l+1$. For
$a=0,\ldots,k$ and $b=0,\ldots,l$ we will define continuous
functions $\ze^{k,l}_{(a,b)}:\De_k\t\De_l\ra[0,1]$ satisfying the
conditions:
\begin{itemize}
\setlength{\itemsep}{0pt}
\setlength{\parsep}{0pt}
\item[(A)] $\sum_{a=0}^k\sum_{b=0}^l\ze^{k,l}_{(a,b)}\equiv 1$
on~$\De_k\t\De_l$.
\item[(B)] $\ze^{k,l}_{(a,b)}$ is supported on $\bigl\{\bigl((x_0,
\ldots,x_k),(y_0,\ldots,y_l)\bigr)\in\De_k\t\De_l:x_a>0$,
$y_b>0\bigr\}$, for all~$a,b$.
\item[(C)] Suppose $a,a'=0,\ldots,k$ and $b,b'=0,\ldots,l$ and
either $a<a'$ and $b>b'$ or $a>a'$ and $b<b'$. Then
$\ze^{k,l}_{(a,b)}\ze^{k,l}_{(a',b')}\equiv 0$ on~$\De_k\t\De_l$.
\item[(D)] If $x_0=\cdots=x_{a-1}=0$, $x_a>0$ and $y_0=\cdots=
y_{b-1}=0$, $y_b>0$ then $\ze^{k,l}_{(a,b)}\bigl((x_0,\ldots,x_k),
(y_0,\ldots,y_l)\bigr)>0$.
\end{itemize}
Then we define continuous $\eta_{(i_a,\ti\imath_b)}:X\t_Y\ti
X\!\ra\![0,1]$ and $\eta_{(i_a,\ti\imath_b)}^{(j,\ti\jmath)}:
V^{(j,\ti\jmath)}\!\ra\![0,1]$~by
\ea
\begin{split}
\eta_{(i_a,\ti\imath_b)}(q)=\ze^{k,l}_{(a,b)}\bigl[&\bigl(\eta_{i_0}\ci
\pi_X(q),\eta_{i_1}\ci\pi_X(q),\ldots,\eta_{i_k}\ci\pi_X(q)\bigr),\\
&\bigl(\ti\eta_{\ti\imath_0}\ci\pi_{\ti X}(q),\ti\eta_{\ti\imath_1}
\ci\pi_{\ti X}(q),\ldots,\ti\eta_{\ti\imath_l}\ci\pi_{\ti
X}(q)\bigr)\bigr],
\end{split}
\label{kh3eq16}\\
\begin{split}
\eta_{(i_a,\ti\imath_b)}^{(j,\ti\jmath)}(v)=
\ze^{k,l}_{(a,b)}\bigl[&\bigl(\eta_{i_0}^j\ci\pi_{V^j}(v),\eta_{i_1}^j
\ci\pi_{V^j}(v),\ldots,\eta_{i_k}^j\ci\pi_{V^j}(v)\bigr),\\
&\bigl(\ti\eta_{\ti\imath_0}^{\ti\jmath}\ci\pi_{\ti
V^{\ti\jmath}}(v),\ti\eta_{\ti\imath_1}^{\ti\jmath}\ci\pi_{\ti
V^{\ti\jmath}}(v),\ldots,\ti\eta_{\ti\imath_l}^{\ti\jmath}\ci
\pi_{\ti V^{\ti\jmath}}(v)\bigr)\bigr].
\end{split}
\label{kh3eq17}
\ea

Here $\bigl(\eta_{i_0}\ci\pi_X(p),\eta_{i_1}\ci\pi_X(p),
\ldots,\eta_{i_k}\ci\pi_X(p)\bigr)\in\De_k$ since
$\sum_{a=0}^l\eta_{i_a}=\sum_{i\in I}\eta_i\equiv 1$ on $X$ and
$\eta_{i_a}\ge 0$ for all $a=0,\ldots,l$ as $\eta_{i_a}$ maps
$X\ra[0,1]$. In the same way $\bigl(\ti\eta_{\ti\imath_0}\ci\pi_{\ti
X}(p),\ti\eta_{\ti\imath_1} \ci\pi_{\ti X}(p),\ldots,
\ti\eta_{\ti\imath_l}\ci\pi_{\ti X}(p)\bigr)\in\De_l$, so
\eq{kh3eq16} is well-defined, and similarly \eq{kh3eq17} is
well-defined. Part (a) above follows from (A) and
\eq{kh3eq16}--\eq{kh3eq17}. Parts (b),(c) are immediate from
\eq{kh3eq16}--\eq{kh3eq17}, Definition \ref{kh3def3}(iii),(iv) and
the identities $\pi_X\ci\psi^{(j,\ti\jmath)}\equiv\psi^j\ci
\pi_{V^j}:(s^{(j,\ti\jmath)})^{-1}(0)\ra X$, $\pi_{\ti
X}\ci\psi^{(j,\ti\jmath)}\equiv\ti\psi^{\ti\jmath}\ci \pi_{\ti
V^{\ti\jmath}}:(s^{(j,\ti\jmath)})^{-1}(0)\ra\ti X$,
$\pi_{V^j}\ci\phi^{(j,\ti\jmath)(k,\ti k)}\equiv \phi^{jk}\ci
\pi_{V^k}:V^{(j,\ti\jmath)(k,\ti k)}\ra V^j$, and $\pi_{\ti
V^{\ti\jmath}}\ci\phi^{(j,\ti\jmath)(k,\ti k)}\equiv\ti\phi^{\ti
\jmath\ti k}\ci\pi_{\ti V^{\ti k}}:V^{(j,\ti\jmath)(k,\ti k)}\ra\ti
V^{\ti\jmath}$.

Part (d) follows from (B), the last sentence of Definition
\ref{kh3def3}(i) for $\Im\psi^i$ and $\Im\ti\psi^{\ti \imath}$, and
$\Im\psi^{(i,\ti\imath)}=(\pi_X\t\pi_{\ti X})^{-1}
(\Im\psi^i\t\Im\ti\psi^{\ti\imath})$. Part (e) follows from (B), the
last two sentences of Definition \ref{kh3def3}(ii) for $V^{ij}$ and
$\ti V^{\ti\imath\ti\jmath}$ and $V^{(i,\ti\imath)(j,\ti\jmath)}=
V^{ij} \t_Y\ti V^{\ti\imath\ti\jmath}$. Part (f) follows from (C)
and \eq{kh3eq16}--\eq{kh3eq17}. For (g), suppose $q\in X\t_YX$, set
$p=\pi_X(q)$ and $\ti p=\pi_{\ti X}q$, and let $i_a\in I$ and
$\ti\imath_b\in\ti I$ be least with $\eta_{i_a}(p)>0$ and
$\ti\eta_{\ti\imath_b}(\ti p)>0$. Then by \eq{kh3eq16} we have
$\eta_{(i_a,\ti\imath_b)}(q)=\ze^{k,l}_{(a,b)}\bigl((x_0,\ldots,
x_k),(y_0,\ldots,y_l)\bigr)$ with $x_c=\eta_{i_c}(p)$ and
$y_d=\ti\eta_{\ti\imath_d}(\ti p)$, so $x_0=\cdots=x_{a-1}=0$,
$x_a>0$ and $y_0=\cdots=y_{b-1}=0$, $y_b>0$ by choice of $i_a,
\ti\imath_b$, and thus $\eta_{(i_a,\ti\imath_b)}(q)>0$ by (D), the
first part of (g). For the second part, note that as
$X=\bigcup_{j\in I}\Im\psi^j$ we have $p\in\Im\psi^j$ for some $j\in
I$, so $\eta_j(p)>0$ by Definition \ref{kh3def3}(i). Therefore $j\ge
i_a$, as $i_a$ is least with $\eta_{i_a}(p)>0$. The final part of
Definition \ref{kh3def3}(i) now implies that $p\in\Im\psi^{i_a}$.
Similarly $\ti p\in\Im\ti\psi^{\ti\imath_b}$. Thus $q\in \pi_X^{-1}
(\Im\psi^i)\cap\pi_{\ti X}^{-1}(\Im\ti\psi^{\ti\imath})$, so
equation \eq{kh3eq9} implies that~$q\in\Im\psi^{(i,\ti\imath)}$.

Finally we construct the functions $\ze^{k,l}_{(a,b)}:
\De_k\t\De_l\ra [0,1]$. To do this we first write down a {\it
triangulation\/} of $\De_k\t\De_l$ into $(k+l)$-simplices
$\De_{k+l}$, and then define the $\ze^{k,l}_{(a,b)}$ on each simplex
$\De_{k+l}$ to be the unique affine function with prescribed values
at each vertex of $\De_{k+l}$. For $0\le a<k$ and $0\le b< l$,
define the wall $W_{a,b}\subset \De_k\t\De_l$ to be the set of
$\bigl((x_0,\ldots,x_k),(y_0,\ldots,y_l)\bigr)$ for which
$x_0+\cdots+x_a=y_0+\cdots+y_b$. These walls $W_{a,b}$ for all $a,b$
cut $\De_k\t\De_l$ into $\binom{k+l}{k}$ regions affine-isomorphic
to the $(k+l)$-simplex~$\De_{k+l}$.

Here is one way to describe these regions. Suppose
$\{1,\ldots,k+l\}=A\amalg B$ with $\md{A}=k$ and $\md{B}=l$. Write
$A=\{r_0,r_1,\ldots,r_{k-1}\}$ with $r_0<\cdots<r_{k-1}$ and
$B=\{s_0,s_1,\ldots,s_{k-1}\}$ with $s_0<\cdots<s_{l-1}$. Define
\begin{gather*}
C_{A,B}=\bigl\{\bigl((x_0,\ldots,x_k),(y_0,\ldots,y_l)\bigr)\in
\De_k\t\De_l:t_1\le t_2\le\cdots\le t_{k+l}\bigr\},\\
\text{where $t_c=x_0+\cdots+x_a$ if $c=r_a$, and
$t_c=y_0+\cdots+y_a$ if $c=s_a$.}
\end{gather*}
Then $C_{A,B}$ for all such $A,B$ are the closures of the regions
into which the $W_{a,b}$ divide the $\De_k\t\De_l$. Write
$p_0,\ldots,p_k$ for the vertices of $\De_k$, where $p_a=(\de_{0a},
\de_{1a},\ldots,\de_{ka})$, and $q_0,\ldots,q_l$ for the vertices of
$\De_l$, where $q_b=(\de_{0b},\ldots,\de_{lb})$. Then each $C_{A,B}$
is a $(k+l)$-simplex, with $k+l+1$ vertices, all of the form
$(p_c,q_d)$ for some $c,d$ determined by~$A,B$.

For $a=0,\ldots,k$ and $b=0,\ldots,l$, define
$\ze^{k,l}_{(a,b)}:\De_k\t\De_l\ra\R$ to be affine linear on the
$(k+l)$-simplex $C_{A,B}$ for all $A,B$ as above, and to satisfy
$\ze^{k,l}_{(a,b)}(p_c,q_d)=\de_{ac}\de_{bd}$ for all $c=0,\ldots,k$
and $d=0,\ldots,k$. This determines $\ze^{k,l}_{(a,b)}$ uniquely on
each $C_{A,B}$, since we have prescribed $\ze^{k,l}_{(a,b)}$ on all
the vertices $(p_c,q_d)$ of the $(k+l)$-simplex $C_{A,B}$, and there
is a unique affine linear function taking any given values at the
vertices.

Each intersection $C_{A,B}\cap C_{A',B'}$ in $\De_k\t\De_l$ is empty
or a simplex $\De_m$ for $0\le m\le k+l$, whose vertices $(p_c,q_d)$
are the common vertices of $C_{A,B},C_{A',B'}$. The definitions of
$\ze^{k,l}_{(a,b)}$ on $C_{A,B}$ and $C_{A',B'}$ both restrict to
the unique affine linear function on $\De_m=C_{A,B}\cap C_{A',B'}$
with $\ze^{k,l}_{(a,b)}(p_c,q_d)=\de_{ac}\de_{bd}$ on the $m+1$
vertices $(p_c,q_d)$ of $\De_m$. Thus the definitions of
$\ze^{k,l}_{(a,b)}$ on $C_{A,B}$ and $C_{A',B'}$ agree on
$C_{A,B}\cap C_{A',B'}$, so $\ze^{k,l}_{(a,b)}$ is well-defined. The
maximum and minimum values of an affine linear function on a simplex
are achieved at the vertices, and we have prescribed the values 0 or
1 for each vertex. Thus, $\ze^{k,l}_{(a,b)}$ maps
$\De_k\t\De_l\ra[0,1]$. Also, $\ze^{k,l}_{(a,b)}\vert_{C_{A,B}}$ is
continuous as it is affine linear, so since
$\De_k\t\De_l=\bigcup_{A,B}C_{A,B}$ and the $C_{A,B}$ are closed,
$\ze^{k,l}_{(a,b)}$ is continuous.

We must verify (A)--(D) above. We have
$\sum_{a=0}^k\sum_{b=0}^l\ab\ze^{k,l}_{(a,b)}\ab(p_c,\ab
q_d)=\ab\sum_{a=0}^k\ab\sum_{b=0}^l\de_{ac}\de_{bd} =1$. Thus,
$\sum_{a=0}^k\sum_{b=0}^l \ze^{k,l}_{(a,b)}\vert_{C_{A,B}}$ is the
unique affine linear function which is 1 on each vertex $(p_c,q_d)$
of $C_{A,B}$, so $\sum_{a=0}^k\sum_{b=0}^l
\ze^{k,l}_{(a,b)}\vert_{C_{A,B}}\equiv 1$. This holds for all
$C_{A,B}$, giving (A). For (B), each $C_{A,B}$ has $k+l+1$ vertices
of the form $(p_c,q_d)$. If $(p_a,q_b)$ is not a vertex of $C_{A,B}$
then $\ze^{k,l}_{(a,b)}$ is zero at each vertex $(p_c,q_d)$ of
$C_{A,B}$, so $\ze^{k,l}_{(a,b)}\equiv 0$ on $C_{A,B}$, and (B)
holds trivially on $C_{A,B}$. If $(p_a,q_b)$ is a vertex of
$C_{A,B}$ then $\ze^{k,l}_{(a,b)}$ is 1 at this vertex and 0 at all
other vertices of $C_{A,B}$. Then the support of
$\ze^{k,l}_{(a,b)}\vert_{C_{A,B}}$, and the subset of $C_{A,B}$ on
which $x_a>0$ and $y_b>0$, are both the complement in $C_{A,B}$ of
the unique codimension one face of $C_{A,B}$ not containing
$(p_a,q_b)$. This proves~(B).

For  (C), suppose $0\le a<a'\le k$ and $0\le b'<b\le l$. Then
$(p_a,q_b)$ lies strictly on one side
$x_0+\cdots+x_a>y_0+\cdots+y_{b'}$ of the wall $W_{a,b'}$, and
$(p_{a'},q_{b'})$ lies strictly on the other side
$x_0+\cdots+x_a<y_0+\cdots+y_{b'}$ of the wall $W_{a,b'}$. Since
each $C_{A,B}$ lies on one or other side of $W_{a,b'}$, it follows
that for each $C_{A,B}$, at most one of $(p_a,p_b)$,
$(p_{a'},p_{b'})$ can be a vertex of $C_{A,B}$, so at most one of
$\ze^{k,l}_{(a,b)},\ze^{k,l}_{(a',b')}$ can be nonzero on $C_{A,B}$,
and $\ze^{k,l}_{(a,b)}\ze^{k,l}_{(a',b')}\equiv 0$ on $C_{A,B}$.
Therefore $\ze^{k,l}_{(a,b)}\ze^{k,l}_{(a',b')}\equiv 0$ on
$\De_k\t\De_l$, as this holds for all $A,B$. The case $0\le a'<a\le
k$ and $0\le b<b'\le l$ follows in the same way using $W_{a',b}$.
This proves~(C).

For (D), let $x_0=\cdots=x_{a-1}=0$, $x_a>0$ and $y_0=\cdots=
y_{b-1}=0$, $y_b>0$. For any $0\le a'<k$ and $0\le b'<l$, by
considering the cases $a'<a$, $a'\ge a$ and $b'<b$, $b'\ge b$
separately it is easy to see that $\bigl((x_0,\ldots,x_k),
(y_0,\ldots,y_l)\bigr)$ and $(p_a,q_b)$ cannot lie strictly on
opposite sides of the wall $W_{a',b'}$ in $\De_k\t\De_l$. Hence
there is some $(k+l)$-simplex $C_{A,B}$ containing both
$\bigl((x_0,\ldots,x_k),(y_0,\ldots,y_l)\bigr)$ and $(p_a,q_b)$. Now
$\ze^{k,l}_{(a,b)}$ is affine linear on $C_{A,B}$, and is 1 on the
vertex $(p_a,q_b)$, and 0 on the $\De_{k+l-1}$ face of $C_{A,B}$
opposite $(p_a,q_b)$, and positive otherwise. As $x_a>0$, $y_b>0$ we
see that $\bigl((x_0,\ldots,x_k),(y_0,\ldots,y_l)\bigr)$ cannot lie
in the $\De_{k+l-1}$ opposite $(p_a,q_b)$, so $\ze^{k,l}_{(a,b)}
\bigl((x_0,\ldots,x_k),(y_0,\ldots,y_l)\bigr)>0$. This proves (D),
and completes the construction of $(\bs I, \bs\eta)\t_Y(\bs{\ti
I},\bs{\ti\eta})$, and Definition~\ref{kh3def14}.
\label{kh3def14}
\end{dfn}

The construction of Definition \ref{kh3def14} is {\it commutative\/}
and {\it associative}. This will imply fibre products of
(co-)gauge-fixing data are commutative and associative.

\begin{thm} Let\/ $Y$ be an orbifold, and for\/ $i=1,2,3$ let\/
$X_i$ be a compact Kuranishi space, $\bs f_i:X_i\ra Y$ a strong
submersion, and\/ $(\bs I_i,\bs\eta_i)$ a really good coordinate
system for\/ $(X_i,\bs f_i)$. Then the really good coordinate
systems\/ $(\bs I_1,\bs\eta_1)\t_Y(\bs I_2,\bs\eta_2)$ for\/
$(X_1\t_YX_2,\bs\pi_Y)$ and\/ $(\bs I_2,\bs\eta_2)\t_Y(\bs
I_1,\bs\eta_1)$ for $(X_2\t_YX_1,\bs\pi_Y)$ are identified by the
natural strong diffeomorphism
\e
\smash{X_1\t_{\bs f_1,Y,\bs f_2}X_2\cong X_2\t_{\bs f_2,Y,\bs
f_1}X_1.}
\label{kh3eq18}
\e
Also, the really good coordinate systems\/ $\bigl((\bs
I_1,\bs\eta_1)\t_Y(\bs I_2,\bs\eta_2)\bigr)\t_Y(\bs I_3,\bs\eta_3)$
for\/ $\bigl((X_1\t_YX_2)\t_YX_3,\bs\pi_Y\bigr)$ and\/ $(\bs
I_1,\bs\eta_1)\t_Y\bigl((\bs I_2,\bs\eta_2)\t_Y(\bs
I_3,\bs\eta_3)\bigr)$ for\/ $\bigl(X_1\t_Y(X_2\t_YX_3),\bs\pi_Y
\bigr)$ are identified by the natural strong diffeomorphism
\e
\smash{\bigl(X_1\t_{\bs f_1,Y,\bs f_2}X_2\bigr)\t_{\bs\pi_Y,Y,\bs
f_3}X_3 \cong X_1\t_{\bs f_1,Y,\bs\pi_Y}\bigl(X_2\t_{\bs f_2,Y,\bs
f_3}X_3\bigr).}
\label{kh3eq19}
\e
\label{kh3thm4}
\end{thm}

\begin{proof} In Definition \ref{kh3def14}, adapting notation in
the obvious way, for $i_1\in I_1$ and $i_2\in I_2$ we have Kuranishi
neighbourhoods $(V^{(i_1,i_2)}_{12},\ldots,\psi^{(i_1,i_2)}_{12})$
on $X_1\t_YX_2$ and $(V^{(i_2,i_1)}_{21},\ldots,\psi^{(i_2,i_1)
}_{21})$ on $X_2\t_YX_1$ from \eq{kh3eq6}, with
$V^{(i_1,i_2)}_{12}=V_1^{i_1}\t_{f_1^{i_1},Y,f_2^{i_2}} V_2^{i_2}$
and $V^{(i_2,i_1)}_{21}=V_2^{i_2}\t_{f_2^{i_2},Y,f_1^{i_1}}
V_1^{i_1}$. Parallel to \eq{kh3eq18} we have a natural
diffeomorphism $V_1^{i_1} \t_{f_1^{i_1},Y,f_2^{i_2}} V_2^{i_2}\cong
V_2^{i_2}\t_{f_2^{i_2}, Y,f_1^{i_1}} V_1^{i_1}$, that is,
$V^{(i_1,i_2)}_{12}\cong V^{(i_2,i_1)}_{21}$. We also have functions
$\eta_{(i_1,i_2)}:X_1\t_YX_2\ra[0,1]$ and $\eta_{(i_2,i_1)}:X_2
\t_YX_1\ra[0,1]$, and $\eta_{(i_1,i_2)}^{(j_1,j_2)}:
V^{(j_1,j_2)}_{12}\ra[0,1]$ and $\eta_{(i_2,i_1)}^{(j_2,j_1)}:
V^{(j_2,j_1)}_{21}\ra[0,1]$ for $j_1\in I_1$, $j_2\in I_2$.

We claim that \eq{kh3eq18} identifies $\eta_{(i_1,i_2)}:X_1\t_YX_2
\ra[0,1]$ and $\eta_{(i_2,i_1)}:X_2\t_YX_1\ra[0,1]$, and the
diffeomorphism $V^{(j_1,j_2)}_{12}\cong V^{(j_2,j_1)}_{21}$
identifies $\eta_{(i_1,i_2)}^{(j_1,j_2)}$ and
$\eta_{(i_2,i_1)}^{(j_2,j_1)}$. Following the definitions through,
this is immediate since the $\smash{\ze^{k,l}_{(a,b)}}$
satisfy~$\ze^{k,l}_{
(a,b)}\bigl((x_0,\ldots,x_k),(y_0,\ldots,y_l)\bigr)\equiv
\ze^{l,k}_{(b,a)}\bigl((y_0,\ldots,y_l),(x_0,\ldots,x_k)\bigr)$.

Equation \eq{kh3eq8} now shows that the open subsets $\dot
V^{(i_1,i_2)}_{12}\subseteq V^{(i_1,i_2)}_{12}$ and $\dot
V^{(i_2,i_1)}_{21}\subseteq V^{(i_2,i_1)}_{21}$ are identified by
the diffeomorphism $V^{(i_1,i_2)}_{12}\cong V^{(i_2,i_1)}_{21}$. So
by \eq{kh3eq7}, for each $k$ we have a natural diffeomorphism
between $\dot V^k_{12}$ and $\dot V^k_{21}$, which easily extends to
an identification of the Kuranishi neighbourhoods $(\dot
V^k_{12},\ldots\dot\psi^k_{12})$ in $(\bs I_1,\bs\eta_1)\t_Y(\bs
I_2,\bs\eta_2)$ and $(\dot V^k_{21},\ldots\dot\psi^k_{21})$ in $(\bs
I_2,\bs\eta_2)\t_Y(\bs I_1,\bs\eta_1)$ under \eq{kh3eq18}. It is an
exercise to show that all the rest of the data in $(\bs
I_1,\bs\eta_1)\t_Y(\bs I_2,\bs\eta_2)$ and $(\bs
I_2,\bs\eta_2)\t_Y(\bs I_1,\bs\eta_1)$ is compatible with these
identifications for all $k$, which proves the first part of the
theorem.

For the second part, we follow a similar method, except that the
conditions on the $\ze^{k,l}_{(a,b)}$ are more complicated. Let
$i_e\in I_e$ for $e=1,2,3$, and define
\begin{align*}
V^{(i_1,i_2,i_3)}_{(12)3}&=(V_1^{i_1}\t_{f_1^{i_1},Y,f_2^{i_2}}
V_2^{i_2})\t_{\pi_Y,Y,f_3^{i_3}}V_3^{i_3}\quad\text{and}\\
V^{(i_1,i_2,i_3)}_{1(23)}&=V_1^{i_1}\t_{f_1^{i_1},Y,\pi_Y}
(V_2^{i_2}\t_{f_2^{i_2},Y,f_3^{i_3}}V_3^{i_3}).
\end{align*}
These naturally extend to Kuranishi neighbourhoods
$(V^{(i_1,i_2,i_3)}_{(12)3},\ldots,\psi^{(i_1,i_2,i_3)}_{(12)3})$ on
$(X_1\t_YX_2)\t_YX_3$ and $(V^{(i_1,i_2,i_3)}_{1(23)},\ldots,
\psi^{(i_1,i_2,i_3)}_{1(23)})$ on $X_1\t_Y(X_2\t_YX_3)$. As for
\eq{kh3eq19}, there is a natural diffeomorphism
$V^{(i_1,i_2,i_3)}_{(12)3}\cong V^{(i_1,i_2,i_3)}_{1(23)}$, which
extends to identifications of the Kuranishi neighbourhoods
$(V^{(i_1,i_2,i_3)}_{(12)3},\ldots,\psi^{(i_1,i_2,i_3)}_{(12)3})$
and $(V^{(i_1,i_2,i_3)}_{1(23)},\ldots, \psi^{(i_1,i_2,
i_3)}_{1(23)})$ under~\eq{kh3eq19}.

Set $I=\{i_1+i_2+i_3:i_1\in I_1$, $i_2\in I_2$, $i_3\in I_3\}$.
Applying Definition \ref{kh3def14} twice, we have functions
$\eta_k:(X_1\t_YX_2)\t_YX_3\ra[0,1]$ and $\eta_k^{(i_1,i_2,i_3)}:
V^{(i_1,i_2,i_3)}_{(12)3}\ra [0,1]$ for $k\in I$ satisfying
analogues of Definition \ref{kh3def3}(i)--(iv), and open subsets
$\dot V^{(i_1,i_2,i_3)}_{(12)3}=\bigl\{v\in
V^{(i_1,i_2,i_3)}_{(12)3}:
\eta_{i_1+i_2+i_3}^{(i_1,i_2,i_3)}(p)>0\bigr\}$ in
$V^{(i_1,i_2,i_3)}_{(12)3}$, such that the Kuranishi neighbourhoods
in $\bigl((\bs I_1,\bs\eta_1)\t_Y(\bs I_2,\bs\eta_2)\bigr)\t_Y(\bs
I_3,\bs\eta_3)$ are $(V^k_{(12)3},\ldots,\psi^k_{(12)3})$ with
\begin{equation*}
V^k_{(12)3}=\ts\coprod_{i_1\in I_1,\;i_2\in I_2,\;i_3\in
I_3:i_1+i_2+i_3=k}\dot V^{(i_1,i_2,i_3)}_{(12)3},
\end{equation*}
and $(E^k_{(12)3},s^k_{(12)3},\psi^k_{(12)3})\vert_{\dot
V^{(i_1,i_2,i_3)}_{(12)3}}=(E^{(i_1,i_2,i_3)}_{(12)3},
s^{(i_1,i_2,i_3)}_{(12)3},\psi^{(i_1,i_2,i_3)}_{(12)3})\vert_{\dot
V^{(i_1,i_2,i_3)}_{(12)3}}$. The partition of unity data in
$\bigl((\bs I_1,\bs\eta_1)\t_Y(\bs I_2,\bs\eta_2)\bigr)\t_Y(\bs
I_3,\bs\eta_3)$ is $\eta_k$ as above and $\eta_k^l$ given by
$\eta_k^l\vert_{\dot V^{(i_1,i_2,i_3)}_{(12)3}}\equiv
\eta_k^{(i_1,i_2,i_3)}$ when $i_1+i_2+i_3=l$. We also have the
analogous functions, sets and Kuranishi neighbourhoods for
$X_1\t_Y(X_2\t_YX_3)$; for clarity, we denote the corresponding
functions by $\bar\eta_k,\bar\eta_k^{(i_1, i_2,i_3)},\bar
\eta_k^l$.

Write $I_e=\{i_0^e,i_1^e,\ldots, i_{k_e}^e\}$ with
$i_0^e<i_1^e<\cdots<i_{k_e}^e$ for $e=1,2,3$, and write
$\{i_1+i_2:i_1\in I_1,$ $i_2\in I_2\}=\{\bar\imath_0,
\bar\imath_1,\ldots,\bar\imath_{\bar k}\}$ with
$\bar\imath_0<\bar\imath_1<\cdots<\bar\imath_{\bar k}$. Then
computing with Definition \ref{kh3def14} shows that for $k\in I$ we
have
\e
\begin{gathered}
\eta_k^{(i_1,i_2,i_3)}(v)=\!\!\!\sum_{\begin{subarray}{l}\bar
a=0,\ldots,\bar k,\\
a_3=0,\ldots,k_3:\;\bar\imath_{\bar a}+i^3_{a_3}=k
\!\!\!\!\!\!\!\!\!\!\!\!\!\!\!\!\!\!\!\!\!\!\!\!\!\!\!\!\end{subarray}}
\ze^{\bar k,k_3}_{(\bar a,a_3)}\bigl((\bar x_0,\ldots,\bar x_{\bar
k}),
\begin{aligned}[t]
(&\eta_{i_0^3}^{k_3}\ci\pi_{V_3^{i_3}}(v),\eta_{i_1^3}^{k_3}\ci\pi_{
V_3^{i_3}}(v),\\
&\ldots,\eta_{i_{k_3}^3}^{k_3}\ci\pi_{V_3^{i_3}}(v))\bigr),
\end{aligned}
\\
\text{where}\quad \bar x_c=\sum_{\begin{subarray}{l}\bar
a_1=0,\ldots,k_1,\\
a_2=0,\ldots,k_2:\;i^1_{a_1}+i^2_{a_2}=\bar\imath_c
\!\!\!\!\!\!\!\!\!\!\!\!\!\!\!\!\!\!\!\!\!\end{subarray}}
\ze^{k_1,k_2}_{(a_1,a_2)}\bigl(
\begin{aligned}[t]&(\eta_{i_0^1}^{k_1}\ci\pi_{
V_1^{i_1}}(v),\ldots,\eta_{i_{k_1}^1}^{k_1}\ci\pi_{V_1^{i_1}}(v)),\\
&(\eta_{i_0^2}^{k_2}\ci\pi_{V_2^{i_2}}(v),
\ldots,\eta_{i_{k_2}^2}^{k_2}\ci\pi_{V_2^{i_2}}(v))\bigr).
\end{aligned}
\end{gathered}
\label{kh3eq20}
\e
There is also an analogous expression
for~$\bar\eta_k^{(i_1,i_2,i_3)}:V^{(i_1,i_2,i_3)}_{1(23)}\ra[0,1]$.

It is now easy to see that $\bigl((\bs I_1,\bs\eta_1)\t_Y(\bs
I_2,\bs\eta_2)\bigr)\t_Y(\bs I_3,\bs\eta_3)$ and $(\bs
I_1,\bs\eta_1)\t_Y\bigl((\bs I_2,\bs\eta_2)\t_Y(\bs
I_3,\bs\eta_3)\bigr)$ are identified by \eq{kh3eq19} provided the
functions $\eta_k^{(i_1,i_2,i_3)}$ on $V^{(i_1,i_2,i_3)}_{(12)3}$
and $\bar\eta_k^{(i_1,i_2,i_3)}$ on $V^{(i_1,i_2,i_3)}_{1(23)}$ are
identified by the diffeomorphism $V^{(i_1,i_2,i_3)}_{(12)3}\cong
V^{(i_1,i_2,i_3)}_{1(23)}$, for all $i_1,i_2,i_3$ and $k$. To prove
this, we claim that
\e
\begin{gathered}
\eta_k^{(i_1,i_2,i_3)}(v)=\!\!\!
\sum_{\begin{subarray}{l}
a_1=0,\ldots,k_1,\\ a_2=0,\ldots,k_2,\\
a_3=0,\ldots,k_3:\\ i^1_{a_1}+i^2_{a_2}+i^3_{a_3}=k\end{subarray}
\!\!\!\!\!\!\!} \ze^{k_1,k_2,k_3}_{(a_1,a_2,a_3)}\bigl(
\begin{aligned}[t]
&(\eta_{i_0^1}^{k_1}\ci\pi_{V_1^{i_1}}(v),\ldots,
\eta_{i_{k_1}^1}^{k_1}\ci\pi_{V_1^{i_1}}(v)),\\
&(\eta_{i_0^2}^{k_2}\ci\pi_{V_2^{i_2}}(v),\ldots,
\eta_{i_{k_2}^2}^{k_2}\ci\pi_{V_2^{i_2}}(v)),\\
&(\eta_{i_0^3}^{k_3}\ci\pi_{V_3^{i_3}}(v),\ldots,
\eta_{i_{k_3}^3}^{k_3}\ci\pi_{V_3^{i_3}}(v))\bigr).
\end{aligned}
\end{gathered}
\label{kh3eq21}
\e

Here $\ze^{k_1,k_2,k_3}_{(a_1,a_2,a_3)}:\De_{k_1}\t\De_{k_2}
\t\De_{k_3}\ra[0,1]$ is defined by analogy with the functions
$\ze^{k,l}_{(a,b)}$ of Definition \ref{kh3def14}. Write points of
$\De_{k_1}\t\De_{k_2}\t\De_{k_3}$ as $\bigl((x_0,\ldots,x_{k_1}),
(y_0,\ldots,y_{k_2}),(z_0,\ldots,z_{k_3})\bigr)$. Then we introduce
three kinds of walls in $\De_{k_1}\t\De_{k_2}\t\De_{k_3}$, defined
by the equations $x_0+\cdots+x_a=y_0+\cdots+y_b$, or
$x_0+\cdots+x_a=z_0+\cdots+z_c$, or $y_0+\cdots+y_b=z_0+\cdots+z_c$,
for $0\le a< k_1$, $0\le b<k_2$, $0\le c<k_3$. These walls
triangulate $\De_{k_1}\t\De_{k_2}\t\De_{k_3}$ into
$\frac{(k_1+k_2+k_3)!}{k_1!k_2!k_3!}$ $(k_1+k_2+k_3)$-simplices
$C_{A,B,C}$ which have vertices of the form $(p_a,q_b,r_c)$, where
$p_a,q_b,r_c$ are vertices of $\De_{k_1},\De_{k_2},\De_{k_3}$
respectively. Then we define $\ze^{k_1,k_2,k_3}_{(a_1,a_2,a_3)}$ to
be the unique function which is affine linear on each $C_{A,B,C}$
and satisfies $\ze^{k_1,k_2,k_3}_{(a_1,a_2,a_3)}(p_a,q_b,r_c)
=\de_{a_1a}\de_{a_2b}\de_{a_3c}$ for all~$a,b,c$.

To prove \eq{kh3eq21} from \eq{kh3eq20}, we first show that on each
simplex $C_{A,B,C}$ on which the $\ze^{k_1,k_2,k_3}_{(a_1,a_2,a_3)}$
in \eq{kh3eq21} are affine linear maps in \eq{kh3eq20} to some
simplices $C_{A',B'}$ and $C_{A'',B''}$ in the domains of $\ze^{\bar
k,k_3}_{(\bar a,a_3)}$ and $\ze^{k_1,k_2}_{(a_1,a_2)}$ in which
$\ze^{\bar k,k_3}_{(\bar a,a_3)},\ze^{k_1,k_2}_{(a_1,a_2)}$ are
affine linear. Since compositions of affine linear functions are
affine linear, it follows that on each simplex $C_{A,B,C}$, the
functions in \eq{kh3eq20} and \eq{kh3eq21} are affine linear. But it
also quickly follows that \eq{kh3eq20} and \eq{kh3eq21} coincide on
the vertices $(p_a,q_b,r_c)$ of $\De_{k_1}\t\De_{k_2}\t\De_{k_3}$,
since both are 1 if $i^1_a+i^2_b+i^3_c=k$ and 0 otherwise. Therefore
\eq{kh3eq20} and \eq{kh3eq21} coincide on each $C_{A,B,C}$, and so
on the whole of~$\De_{k_1}\t\De_{k_2}\t\De_{k_3}$.

This proves \eq{kh3eq21}. In the same way we obtain an analogous
expression for $\bar\eta_k^{(i_1,i_2,i_3)}(v)$. Identifying
$V^{(j_1,j_2,j_3)}_{(12)3}$ and $V^{(j_1,j_2,j_3)}_{1(23)}$ using
the natural diffeomorphism, these expressions coincide, so
$\eta_k^{(i_1,i_2,i_3)}$ and $\bar\eta_k^{(i_1,i_2,i_3)}$ are
identified. This completes the proof of Theorem~\ref{kh3thm4}.
\end{proof}

Note that associativity in Theorem \ref{kh3thm4} worked because of
special properties of the $\ze^{k,l}_{(a,b)}$ in Definition
\ref{kh3def14}, implying that \eq{kh3eq20},\eq{kh3eq21} are
equivalent. We define fibre products of (co-)gauge-fixing data, and
prove Property \ref{kh3pr}(h). As we want products on cohomology,
not homology, we explain the co-gauge-fixing case.

\begin{dfn} Let $X,\ti X$ be compact Kuranishi spaces, $Y$ an
orbifold, and $\bs f:X\ra Y$, $\bs{\ti f}:\ti X\ra Y$ be strong
submersions. Suppose $\bs C=(\bs I,\bs\eta,C^i:i\in I)$ and $\bs{\ti
C}=(\bs{\ti I},\bs{\ti\eta},\ti C^{\ti\imath}: \ti\imath\in\ti I)$
are co-gauge-fixing data for $(X,\bs f)$ and $(\ti X,\bs{\ti f})$.
Let $(\bs I,\bs\eta)\t_Y(\bs{\ti I},\bs{\ti\eta})$ be the really
good coordinate system for $(X\t_Y\ti X,\bs\pi_Y)$ given in
Definition \ref{kh3def14}. Write $((\bs I,\bs\eta)\t_Y(\bs{\ti
I},\bs{\ti\eta}))\kern .1em\check{}\,$ for the excellent coordinate
system for $(X\t_Y\ti X,\bs\pi_Y)$ constructed from $(\bs
I,\bs\eta)\t_Y(\bs{\ti I},\bs{\ti\eta})$ by Algorithm~\ref{kh3alg}.

In the notation of Definition \ref{kh3def14} we have $(\bs
I,\bs\eta)\t_Y(\bs{\ti I},\bs{\ti\eta})=(\bs{\dot
I},\bs{\dot\eta})$, where $\bs{\dot I}=\bigl(\dot I,(\dot V^k,\dot
E^k,\dot s^k,\dot\psi^k),\dot\pi_Y^k:k\in\dot I,\ldots\bigr)$, and
$\dot V^k=\coprod_{i\in I,\; \ti\imath\in\ti I:i+\ti\imath=k+\dim
Y}\dot V^{(i,\ti\imath)}$ by \eq{kh3eq7}, with $\dot
V^{(i,\ti\imath)}$ an open subset in
$V^{(i,\ti\imath)}=V^i\t_{f^i,Y,\ti f^{\ti\imath}} \ti
V^{\ti\imath}$. Also $\dot E^k\vert_{\dot V^{(i,\ti\imath)}}=
E^{(i,\ti\imath)}\vert_{\dot V^{(i,\ti\imath)}}$, where
$E^{(i,\ti\imath)}=\pi_{V^i}^*(E^i)\op\pi_{\ti V^{\ti\imath}}^*(\ti
E^{\ti\imath})$ by \eq{kh3eq6}. Thus by Algorithm \ref{kh3alg} we
may write $((\bs I,\bs\eta)\t_Y(\bs{\ti I},\bs{\ti\eta}))\kern
.1em\check{}\,=(\bs{\check I},\bs{\check\eta})$, where $\bs{\check
I}=\bigl(\check I,(\check V^k,\check E^k,\check
s^k,\check\psi^k),\check\pi_Y^k: k\in\check I,\ldots\bigr)$, with
$\check V^k$ an open subset of $\dot V^k$ for each $k\in\dot I$, and
$\check I=\{k\in\dot I:\check V^k\ne\es\}$, and $(\check E^k,\check
s^k,\check\psi^k)=(\dot E^k,\dot s^k,\dot\psi^k)\vert_{\check V^k}$,
$\check\pi_Y^k=\dot\pi_Y^k \vert_{\check V^k}$ for all~$k\in\check
I$.

For $k\in\check I$, define $\check C^k:\check E^k\ra P$ by $\check
C^k(\check e)=\mu\bigl(C^i\ci\pi_{E^i}(\check e),\ti
C^{\ti\imath}\ci\pi_{\ti E^{\ti\imath}}(\check e)\bigr)$ whenever
$\check e$ lies in $\check E^k\cap(\dot E^k\vert_{\check V^k\cap\dot
V^{(i,\ti\imath)}})\subseteq E^{(i,\ti\imath)}$, for all $i\in I$
and $\ti\imath\in\ti I$ with $i+\ti\imath=k+\dim Y$, where
$\pi_{E^i}:E^{(i,\ti\imath)}\ra E^i$ and $\pi_{\ti
E^{\ti\imath}}:E^{(i,\ti\imath)}\ra\ti E^{\ti\imath}$ are the
projections from $E^{(i,\ti\imath)}=\pi_{V^i}^*(E^i)\op\pi_{\ti
V^{\ti\imath}}^*(\ti E^{\ti\imath})\cong E^i\t_Y\ti E^{\ti\imath}$,
and $\mu:P\t P\ra P$ is as in \eq{kh3eq2}. We shall show $\check
C^k\t(\check f^k\ci\check\pi^k):\check E^k\ra P\t Y$ is {\it
globally finite}.

Let $i\in I$ and $\ti\imath\in\ti I$ with $i+\ti\imath=k+\dim Y$.
Then $C^i:E^i\ra P_{n^i}\subset P$ and $\ti C^{\ti\imath}:\ti
E^{\ti\imath}\ra P_{\ti n^{\ti\imath}}\subset P$ for some $n^i,\ti
n^{\ti\imath}\gg 0$. Definition \ref{kh3def3}(vi) gives compact
$T,\ti T\subseteq Y$ with $f^i(V^i)\subseteq T$ and $\ti
f^{\ti\imath}(\ti V^{\ti\imath})\subseteq\ti T$. The map $Y\ra\N$
taking $y\mapsto\md{\Stab(y)}$ is upper semicontinuous, so as
$T\cap\ti T$ is compact there exists $L\ge 1$ such that
$\md{\Stab(y)}\le L$ for all $y\in T\cap\ti T$. As
$C^i\t(f^i\ci\pi^i):E^i\ra P\t Y$ and $\ti C^{\ti\imath}\t(\ti
f^{\ti\imath}\ci\ti\pi^{\ti\imath}):\ti E^{\ti\imath}\ra P\t Y$ are
globally finite, there exist $N^i,\ti N^{\ti\imath}\ge 0$ such that
$\bmd{(C^i\t(f^i\ci\pi^i))^{-1}(p,y)}\le N^i$ and $\bmd{(\ti
C^{\ti\imath}\t(\ti f^{\ti\imath}\ci\ti\pi^{\ti\imath}))^{-1}(\ti
p,y)}\le\ti N^{\ti\imath}$ for all $p,\ti p\in P$ and $y\in Y$. For
$p'\in P$ and $y\in Y$, consider the projection
\e
\begin{split}
\pi_{E^i}\t&\pi_{\ti E^{\ti\imath}}:\bigl((\check C^k\t(\check
f^k\ci\check\pi^k))^{-1}(p',y)\bigr)\cap E^{(i,\ti\imath)} \longra\\
&\coprod\nolimits_{\begin{subarray}{l}(p,\ti p)\in P\t P:\\
\mu(p,\ti p)=p'\end{subarray}} \bigl(C^i\t(f^i\ci\pi^i)
\bigr)^{-1}(p,y)\t\bigl(\ti C^{\ti\imath}\t(\ti
f^{\ti\imath}\ci\ti\pi^{\ti\imath})\bigr)^{-1}(\ti p,y).
\end{split}
\label{kh3eq22}
\e

Since $C^i$ maps $E^i\ra P_{n^i}$, $\ti C^{\ti\imath}$ maps $\ti
E^{\ti\imath}\ra P_{\ti n^{\ti\imath}}$, and $\mu$ maps $P_{n^i}\t
P_{\ti n^{\ti\imath}}\ra P_{n^i+\ti n^{\ti\imath}}$, the only
nonempty terms in \eq{kh3eq22} occur when $p\in P_{n^i}$, $\ti p\in
P_{\ti n^{\ti\imath}}$ and $p'\in P_{n^i+\ti n^{\ti\imath}}$. From
\eq{kh3eq2} we see that each $p'\in P_{n^i+\ti n^{\ti\imath}}$ may
be written as $\mu(p,\ti p)$ for at most $\binom{n^i+\ti
n^{\ti\imath}}{n^i}$ pairs $(p,\ti p)$. Thus on the second line of
\eq{kh3eq22} there are at most $\binom{n^i+\ti n^{\ti\imath}}{n^i}$
nonempty sets, each of which has at most $N^i \ti N^{\ti\imath}$
elements, so the second line of \eq{kh3eq22} is a finite set of size
at most~$\binom{n^i+\ti n^{\ti\imath}}{n^i}N^i \ti N^{\ti\imath}$.

But $\check E^k\cap E^{(i,\ti\imath)}$ is an open subset of
$E^i\t_Y\ti E^{\ti\imath}$, so the definition \eq{kh2eq10} of fibre
products of orbifolds implies that if $e\in E^i$ and $\ti e\in\ti
E^{\ti\imath}$ with $f^i\ci\pi^i(e)=y=\ti f^{\ti\imath}\ci\ti
\pi^{\ti\imath}(\ti e)$ then the set of $\check e\in\check E^k\cap
E^{(i,\ti\imath)}$ with $(\pi_{E^i}\t\pi_{\ti E^{\ti\imath}})(\check
e)=(e,\ti e)$ is a subset of $(f^i\ci\pi^i)_*(\Stab(e))\backslash
\Stab(y)/(\ti f^{\ti\imath}\ci\ti \pi^{\ti\imath})_*(\Stab(\ti e))$,
and so consists of at most $L$ points, as $\md{\Stab(y)}\le L$.
Hence each point on the second line of \eq{kh3eq22} pulls back to at
most $L$ points on the first line. Therefore for all $(p',y)\in P\t
Y$ we have
\begin{equation*}
\bmd{((\check C^k\t(\check f^k\ci\check\pi^k))^{-1}(p',y)\bigr)\cap
E^{(i,\ti\imath)}}\le L\ts \binom{n^i+\ti n^{\ti\imath}}{n^i}N^i \ti
N^{\ti\imath}.
\end{equation*}
Summing this over all $i,\ti\imath$ implies that $\check
C^k\t(\check f^k\ci\check\pi^k)$ is {\it globally finite}, with
constant $\check N^k=\sum_{i,\ti\imath:\; i+\ti\imath=k+\dim Y}L\ts
\binom{n^i+\ti n^{\ti\imath}}{n^i}N^i \ti N^{\ti\imath}$. Hence
$\bs{\check C}=(\bs{\check I},\bs{\check\eta},\check C^k:k\in\check
I)$ is {\it co-gauge-fixing data\/} for $(X\t_Y\ti X,\bs\pi_Y)$,
which we write as~$\bs C\t_Y\bs{\ti C}$.

Now suppose that $\bs f:X\ra Y$ is strongly smooth and $\bs G$ is
gauge-fixing data for $(X,\bs f)$, and $\bs{\ti f}:\ti X\ra Y$ is a
strong submersion and $\bs{\ti C}$ is co-gauge-fixing data for $(\ti
X,\bs{\ti f})$. Then the construction above, replacing $C^i,\check
C^k$ by $G^i,\check G^k$, yields {\it gauge-fixing data} $\bs{\check
G}=\bs G\t_Y\bs{\ti C}$ for $(X\t_Y\ti X,\bs\pi_Y)$. The proof that
the $\check G^k$ are globally finite is slightly different, but the
definitions are unchanged.
\label{kh3def15}
\end{dfn}

\begin{prop} Let\/ $Y$ be an orbifold without boundary, and for
$a=1,2,3$ let\/ $X_a$ be a compact Kuranishi space, $\bs f_a:X_a\ra
Y$ a strong submersion, and\/ $\bs C_a$ co-gauge-fixing data for
$(X_a,\bs f_a)$. Then:
\begin{itemize}
\setlength{\itemsep}{0pt}
\setlength{\parsep}{0pt}
\item[{\rm(a)}] The natural strong diffeomorphism\/ $\bs
a:X_1\t_YX_2\ra X_2\t_YX_1$ extends to an isomorphism in the sense
of Definition~{\rm\ref{kh3def7}}
\e
\!\!\!(\bs a,\bs b):(X_1\t_YX_2,\bs\pi_Y,\bs C_1\t_Y\bs
C_2)\longra(X_2 \t_YX_1,\bs\pi_Y,\bs C_2\t_Y\bs C_1).
\label{kh3eq23}
\e
\item[{\rm(b)}] The natural strong diffeomorphism\/ $\bs
a':\pd(X_1\t_YX_2)\ra \pd X_1\t_YX_2\amalg X_1\t_Y\pd X_2$ extends
to an isomorphism
\ea
\!\!\!\!&(\bs a',\bs b'):\bigl(\pd(X_1\t_YX_2),\bs\pi_Y,(\bs
C_1\t_Y\bs C_2)\vert_{\pd(X_1\t_YX_2)}\bigr)\longra
\label{kh3eq24}\\
\!\!\!\!&\bigl(\pd X_1\!\t_Y\!X_2,\bs\pi_Y,(\bs C_1\vert_{\pd
X_1})\!\t_Y\!\bs C_2\bigr)\!\amalg\!\bigl(X_1\!\t_Y\!\pd
X_2,\bs\pi_Y,\bs C_1\!\t_Y\!(\bs C_2\vert_{\pd X_2})\bigr).
\nonumber
\ea
\item[{\rm(c)}] The natural strong diffeomorphism $\bs
a'':(X_1\t_YX_2)\t_YX_3\ra X_1\t_Y(X_2\t_YX_3)$ extends to an
isomorphism
\e
\begin{split}
\!\!\!(\bs a'',\bs b''):\bigl((X_1&\t_YX_2)\t_YX_3,\bs\pi_Y,(\bs
C_1\t_Y
\bs C_2)\t_Y\bs C_3\bigr)\longra\\
&\bigl(X_1\t_Y(X_2\t_YX_3),\bs\pi_Y,\bs C_1\t_Y(\bs C_2\t_Y\bs
C_3)\bigr).
\end{split}
\label{kh3eq25}
\e
\item[{\rm(d)}] If instead\/ $\bs f_1:X_1\ra Y$ is strongly smooth
and\/ $\bs G_1$ is gauge-fixing data for $(X_1,\bs f_1),$ then in
{\rm(b),(c)} above $\bs a',\bs a''$ extend to isomorphisms
\ea
\!\!\!\!\!\!&(\bs a',\bs b'):\bigl(\pd(X_1\t_YX_2),\bs\pi_Y,(\bs
G_1\t_Y\bs C_2)\vert_{\pd(X_1\t_YX_2)}\bigr)\longra
\label{kh3eq26}\\
\!\!\!\!\!\!&\bigl(\pd X_1\!\t_Y\!X_2,\bs\pi_Y,(\bs G_1\vert_{\pd
X_1})\!\t_Y\!\bs C_2\bigr)\!\amalg\!\bigl(X_1\!\t_Y\!\pd
X_2,\bs\pi_Y,\bs G_1\!\t_Y\!(\bs C_2\vert_{\pd X_2})\bigr),
\nonumber\\
\begin{split}
\!\!\!\!\!\!&(\bs a'',\bs
b''):\bigl((X_1\t_YX_2)\t_YX_3,\bs\pi_Y,(\bs G_1\t_Y
\bs C_2)\t_Y\bs C_3\bigr)\longra\\
\!\!\!\!\!\!&\qquad\qquad\quad\bigl(X_1\t_Y(X_2\t_YX_3),\bs\pi_Y,\bs
G_1\t_Y(\bs C_2\t_Y\bs C_3)\bigr),
\end{split}
\label{kh3eq27}
\ea
\end{itemize}
\label{kh3prop6}
\end{prop}

\begin{proof} By Theorem \ref{kh3thm4}, the strong diffeomorphism
$\bs a$ identifies the really good coordinate systems $(\bs
I_1,\bs\eta_1)\t_Y(\bs I_2,\bs\eta_2)$ and $(\bs
I_2,\bs\eta_2)\t_Y(\bs I_1,\bs\eta_1)$. So applying Algorithm
\ref{kh3alg}, $\bs a$ also identifies $((\bs I_1,\bs\eta_1)\t_Y(\bs
I_2,\bs\eta_2))\kern .1em\check{}\,$ and $((\bs
I_2,\bs\eta_2)\t_Y(\bs I_1,\bs\eta_1))\kern .1em\check{}\,$. That
is, there is a natural isomorphism $\bs b:((\bs
I_1,\bs\eta_1)\t_Y(\bs I_2,\bs\eta_2))\kern .1em\check{}\,\ra((\bs
I_2,\bs\eta_2)\t_Y(\bs I_1,\bs\eta_1))\kern .1em\check{}\,$
compatible with $\bs a$. Write the really good coordinate system in
$((\bs I_1,\bs\eta_1)\t_Y(\bs I_2,\bs\eta_2))\kern .1em\check{}\,$
as $\bigl(\check I_{12},(\check
V_{12}^k,\ldots,\check\psi_{12}^k):k\in\check I_{12},\ldots\bigr)$
with functions $\check C_{12}^k:\check E_{12}^k\ra P$, and similarly
for $((\bs I_2,\bs\eta_2)\t_Y(\bs I_1,\bs\eta_1))\kern
.1em\check{}\,$ with~$\check I_{21},\ldots$.

To show that $(\bs a,\bs b)$ is an isomorphism in \eq{kh3eq23}, it
only remains to check that $\check C_{21}^k\ci\hat b^k\equiv\check
C_{12}^k$ for all $k\in\check I_{12}=\check I_{21}$. But this is
obvious since $C_{12}^k\vert_{\check E^k\cap E_{12}^{(i_1,i_2)}}
\equiv\mu\ci\bigl((C_1^{i_1}\ci\pi_{E_1^{i_1}})\t(C_2^{i_2}\ci
\pi_{E_2^{i_2}})\bigr)$ for $i_1\in I_1$, $i_2\in I_2$ with
$i_1+i_2= k+\dim Y$, and $C_{21}^k\vert_{\check E^k\cap
E_{21}^{(i_2,i_1)}}\equiv \mu\ci\bigl((C_2^{i_2}\ci\pi_{E_2^{i_2}})
\t(C_1^{i_1}\ci \pi_{E_1^{i_1}})\bigr)$, and we have
$(\pi_{E_1^{i_1}}\t\pi_{E_2^{i_2}})\ci\hat b^k\equiv
\pi_{E_1^{i_1}}\t\pi_{E_2^{i_2}}$ as maps $\smash{\check E^k\cap
E_{12}^{(i_1,i_2)}\ra E_1^{i_1}\t E_2^{i_2}}$, and $\mu$ in
\eq{kh3eq2} is {\it commutative}. This proves part~(a).

For (b)--(d) we follow a similar method, using the fact that $\mu$
in \eq{kh3eq2} is {\it associative} for \eq{kh3eq25} and
\eq{kh3eq27}, but there is one complication. The argument above
proves that the excellent coordinate systems $((\bs I_1\!\t_Y\!\bs
I_2)\t_Y\!\bs I_3)\kern .7pt\check{}\,$ and $(\bs I_1\!\t_Y(\bs
I_2\!\t_Y\!\bs I_3))\kern .7pt\check{}\,$ are isomorphic under $\bs
a''$. But we need to know that $((\bs I_1\!\t_Y\!\bs I_2)\kern
.7pt\check{}\t_Y\!\bs I_3)\kern .7pt\check{}\,$ and $(\bs
I_1\!\t_Y\!(\bs I_2\!\t_Y\!\bs I_3)\kern .7pt\check{}\,)\kern
.7pt\check{}\,$ are isomorphic. One can show that $((\bs I_1\t_Y\bs
I_2)\kern .7pt\check{}\,\t_Y\bs I_3)\kern .7pt\check{}\,=((\bs
I_1\t_Y\bs I_2)\t_Y\bs I_3)\kern .7pt\check{}\,$ and $(\bs
I_1\t_Y(\bs I_2\t_Y\bs I_3)\kern .7pt\check{}\,)\kern
.7pt\check{}\,=(\bs I_1\t_Y(\bs I_2\t_Y\bs I_3))\kern
.7pt\check{}\,$, that is, applying Algorithm \ref{kh3alg} twice
gives the same result as applying it once. Equations \eq{kh3eq25}
and \eq{kh3eq27} follow.
\end{proof}

Products are also functorial for proper pullbacks, and pushforwards.
The proof is similar to that of Proposition \ref{kh3prop6}, so we
leave it as an exercise.

\begin{prop} Suppose $h:Y\ra Z$ is a smooth, proper map of orbifolds.
\begin{itemize}
\setlength{\itemsep}{0pt}
\setlength{\parsep}{0pt}
\item[{\rm(a)}] Let\/ $X_i$ be a compact Kuranishi space, $\bs
f_i:X_i\ra Z$ a strong submersion, and\/ $\bs C_i$
co-gauge-fixing data for $(X_i,\bs f_i)$ for $i=1,2$. Then the
natural strong diffeomorphism\/ $\bs a:Y\t_Z(X_1\t_ZX_2)\ra
(Y\t_ZX_1)\t_Y(Y\t_ZX_2)$ extends to an isomorphism
\begin{align*}
(\bs a,\bs b):\,&\bigl(Y\t_Z(X_1\t_ZX_2),\bs\pi_Y,h^*(\bs C_1\t_Z\bs
C_2)\bigr)\longra\\
&\bigl((Y\t_ZX_1)\t_Y(Y\t_ZX_2),\bs\pi_Y,h^*(\bs C_1)\t_Yh^*(\bs
C_2)\bigr).
\end{align*}
\item[{\rm(b)}] Let\/ $X_1,X_2$ be compact Kuranishi spaces, $\bs
f_1:X_1\ra Y$ strongly smooth, $\bs G_1$ gauge-fixing data for
$(X_1,\bs f_1),$ $\bs f_2:X_2\ra Z$ a strong submersion, and\/ $\bs
C_2$ co-gauge-fixing data for $(X_1,\bs f_2)$. Then the natural
strong diffeomorphism\/ $\bs a':X_1\t_{h\ci\bs f_1,Z,\bs f_2}X_2\ra
X_1\t_{\bs f_1,Y,\bs\pi_Y}(Y\t_{h,Z,\bs f_2}X_2)$ extends to an
isomorphism
\begin{align*}
(\bs a',\bs b'):\,&\bigl(X_1\t_ZX_2,\bs\pi_Z,h_*(\bs
G_1)\t_Z\bs C_2)\bigr)\longra\\
&\bigl(X_1\t_Y(Y\t_ZX_2),h\ci\bs\pi_Y,h_*(\bs G_1\t_Yh^*(\bs
C_2))\bigr).
\end{align*}
\end{itemize}
\label{kh3prop7}
\end{prop}

\begin{rem}{\bf(a)} The construction of fibre products $\bs C\t_Y\bs{\ti
C}$ in Definition \ref{kh3def15} can actually be split into two
stages. Given compact Kuranishi spaces $X,\ti X$, orbifolds $Y,\ti
Y$, strongly smooth maps $\bs f:X\ra Y$ and $\bs{\ti f}:\ti X\ra\ti
Y$, and co-gauge-fixing data $\bs C,\bs{\ti C}$ for $(X,\bs f)$ and
$(\ti X,\bs{\ti f})$, we can define {\it product co-gauge-fixing
data} $\bs C\t\bs{\ti C}$ for $(X\t\ti X,\bs f\t\bs{\ti f})$, that
is, for the strongly smooth map $\bs f\t\bs{\ti f}:X\t\ti X\ra
Y\t\ti Y$. To do this we simply replace fibre products over $Y$ by
ordinary products (fibre products over a point) throughout the
construction of $\bs C\t_Y\bs{\ti C}$. Note that the construction of
$\bs C\t\bs{\ti C}$ is still extremely complex. Products of
gauge-fixing data $\bs C\t\bs{\ti C}$ are commutative and
associative, by an easy modification of Proposition~\ref{kh3prop6}.

Then to construct $\bs C\t_Y\bs{\ti C}$, take $\ti Y=Y$, and define
$\De:Y\ra Y\t Y$ to be the diagonal map, $\De:y\mapsto(y,y)$. Then
$\bs C\t_Y\bs{\ti C}=\De^*(\bs C\t\bs{\ti C})$, that is, $\bs
C\t_Y\bs{\ti C}$ is the pullback of $\bs C\t\bs{\ti C}$ under $\De$,
as in~\S\ref{kh37}.

If $\bs f:X\ra Y$, $\bs{\ti f}:\ti X\ra\ti Y$ are strongly smooth
and $\bs G,\bs{\ti G}$ are gauge-fixing data for $(X,\bs f),(\ti
X,\bs{\ti f})$, then the definition of $\bs C\t\bs{\ti C}$
generalizes immediately to give {\it product gauge-fixing data\/}
$\bs G\t\bs{\ti G}$ for $(X\t\ti X,\bs f\t\bs{\ti f})$, where $\bs
f\t\bs{\ti f}:X\t\ti X\ra Y\t\ti Y$ is strongly smooth.

\noindent{\bf(b)} We can use products $\bs C\t\bs{\ti C}$ in (a) and
pullbacks to define generalizations of fibre products of
co-gauge-fixing data to more than one orbifold $Y$, which we leave
as an exercise. For example, if $X_1,X_2$ are compact Kuranishi
spaces, $Y_3,Y_4,Y_5$ are orbifolds, $\bs f_{13}\t\bs f_{14}:X_1\ra
Y_3\t Y_4$, $\bs f_{24}\t\bs f_{25}:X_2\ra Y_4\t Y_5$ strong
submersions, and $\bs C_1,\bs C_2$ are co-gauge-fixing data for
$(X_1,\bs f_{13}\t\bs f_{14})$ and $(X_2,\bs f_{24}\t\bs f_{25})$,
then we can define co-gauge-fixing data $\bs C_1\t_{Y_4}\bs C_2$ for
$\bigl(X_1\t_{\bs f_{14},Y_4,\bs f_{24}}X_2,(\bs
f_{13}\ci\bs\pi_{X_1})\t\bs\pi_{Y_4}\t(\bs
f_{25}\ci\bs\pi_{X_2})\bigr)$ over~$Y_3\t Y_4\t Y_5$.

We will generalize these ideas on products of (co-)gauge-fixing data
in our discussion of Kuranishi (co)homology as a {\it bivariant
theory\/} in~\S\ref{kh48}.
\label{kh3rem4}
\end{rem}

\subsection{Effective gauge-fixing data and co-gauge-fixing data}
\label{kh39}

For reasons that will become clearer in Chapter \ref{kh4}, because
we allow triples $(X,\bs f,\bs G)$ or $(X,\bs f,\bs C)$ to have
finite automorphism groups, as in \S\ref{kh33}, we can only use them
to define Kuranishi (co)homology theories $KH_*,KH^*(Y;R)$ in which
the coefficient ring $R$ is a $\Q$-{\it algebra}, with the homology
theory isomorphic to singular homology $KH_*(Y;R)\cong
H^\rsi_*(Y;R)$. Essentially this is because $R$ must contain
$\md{\Aut(X,\bs f,\bs G)}^{-1}$, to enable us to average
over~$\Aut(X,\bs f,\bs G)$.

We will also define {\it effective Kuranishi (co)homology}
$KH_*^\ef,KH^*_\ec(Y;R)$, which work for $R$ any commutative ring,
with $KH_*^\ef(Y;R)\cong H^\rsi_*(Y;R)$ and (for $Y$ a manifold)
$KH^*_\ec(Y;R)\cong H^*_\cs(Y;R)$. To make this work, we need to
ensure that (connected) triples $(X,\bs f,\ubG)$ or $(X,\bs f,\ubC)$
have only {\it trivial automorphism groups} $\{1\}$. Our goal now is
to define notions of {\it effective\/} (co-)gauge-fixing data
$\ubG,\ubC$, and partially extend \S\ref{kh32}--\S\ref{kh38} to
them.

The basic idea is that we shall require the maps $G^i,C^i\t
(f^i\ci\pi^i)$ of Definitions \ref{kh3def6} and \ref{kh3def7} to be
{\it injective on} $(E^i)^\ci$ as well as globally finite, which
easily implies that $(X,\bs f,\ubG)$ and $(X,\bs f,\ubC)$ have only
trivial automorphisms. However, there are significant problems in
carrying this out. Here are some of them:
\begin{itemize}
\setlength{\itemsep}{0pt}
\setlength{\parsep}{0pt}
\item[(i)] Let $(X,\bs f,\ubG)$ be a triple, with $\ubG$
effective gauge-fixing data. Gauge-fixing data must restrict to
boundaries, so we can consider $(\pd^2X,\bs f\vert_{\pd^2X},
\ubG\vert_{\pd^2X})$. Let $\bs\si:\pd^2X\ra\pd^2X$ be the natural
involution of Definition \ref{kh2def19}. Then $\bs\si\in\Aut
(\pd^2X,\bs f\vert_{\pd^2X},\ubG\vert_{\pd^2X})$. This is
incompatible with $\ubG\vert_{\pd^2X}$ being {\it effective\/}
gauge-fixing data. Thus, restricting effective (co-)gauge-fixing
data to boundaries will not yield effective (co-)gauge-fixing data.
\item[(ii)] Even if we did define a homology theory
$KH_*^\ef(Y;R)$ using some notion of effective gauge-fixing data, we
could not prove $KH_*^\ef(Y;R)\cong H^\rsi_*(Y;R)$ for rings
$R\not\supseteq\Q$ unless we make some restriction on the {\it
stabilizer groups} of Kuranishi spaces $X$ in triples $(X,\bs
f,\ubG)$, because otherwise we cannot perturb $X$ to a manifold or
orbifold, and triangulate it by simplices.
\item[(iii)] If $f^i\ci\pi^i:E^i\ra Y$, $\ti f^{\ti\imath}\ci
\ti\pi^{\ti\imath}:\ti E^{\ti\imath}\ra Y$ are submersions of
orbifolds, the definition \eq{kh2eq10} of $E^i\t_Y\ti E^{\ti\imath}$
involving $f_*(\Stab(p))\backslash\Stab(f(p))/f'_*(\Stab(p'))$ means
that the projection $\pi_{E^i}\t\pi_{\ti E^{\ti\imath}}:E^i\t_Y\ti
E^{\ti\imath}\ra E^i\t\ti E^{\ti\imath}$ need not be injective, only
finite, if $Y$ is an orbifold, not a manifold. This can cause
$\check C^k\t(\check f^k\ci\check\pi^k):(\check E^k)^\ci\ra P\t Y$
in fibre products of effective co-gauge-fixing data $\ubC,\ubtC$ in
Definition \ref{kh3def15} not to be injective, so that
$\ubC\t_Y\ubtC$ is not effective co-gauge-fixing data.
\item[(iv)] Here is a second problem with taking fibre products of
effective co-gauge-fixing data. Let $(X,\bs f,\ubC)$ be a triple,
with $\ubC$ effective co-gauge-fixing data and $\bs f:X\ra Y$ a
strong submersion. Consider the triple
$(X\t_YX,\bs\pi_Y,\ubC\t_Y\ubC)$, with $\ubC\t_Y\ubC$ as in
\S\ref{kh38}. This has an obvious $\Z_2$-symmetry exchanging the two
factors of $X,\ubC$. So $\Aut(X\t_YX,\bs\pi_Y,\ubC\t_Y\ubC)$ may be
nontrivial, and $\ubC\t_Y\ubC$ is not effective.
\end{itemize}
To overcome each of these problems, we have to sacrifice one of the
good properties of our theory.

Problems (i)--(iii) force us to restrict the Kuranishi spaces $X$
and morphisms $\bs f$ allowed in triples $(X,\bs f,\ubG)$ or $(X,\bs
f,\ubC)$ in the theory. To solve (i) we require the injectivity
conditions to hold only for each {\it connected component\/} of
$(\pd^kX,\bs f\vert_{\pd^kX},\ubG\vert_{\pd^kX})$. Then we can
ensure that $\Aut(X,\bs f,\ubG)=\{1\}$ for {\it connected\/} triples
$(X,\bs f,\ubG)$, which is enough to make the theory work. If
$\bs\si:\pd^2X\ra\pd^2X$ acts freely on the connected components of
$(\pd^2X,\bs f\vert_{\pd^2X},\ubG\vert_{\pd^2X})$, there is then no
contradiction in (i). This constrains the codimension $k$ corners
$\pd^kX$ for $k=2,3,\ldots$. To deal with (ii) we restrict to {\it
effective} Kuranishi spaces $X$, defined below, which can always be
perturbed to effective orbifolds.

For (iii) we require the smooth maps $f^i\ci\pi^i:E^i\ra Y$ in
triples $(X,\bs f,\ubC)$ to be surjective on stabilizer groups,
which makes $\pi_{E^i}\t\pi_{\ti E^{\ti\imath}}:E^i\t_Y\ti
E^{\ti\imath}\ra E^i\t\ti E^{\ti\imath}$ injective. Thus, Property
\ref{kh3pr}(a) will not hold: only some, not all, pairs $(X,\bs f)$
admit effective (co-)gauge-fixing data. This will limit the
applications of the theory --- for instance, in \S\ref{kh62} we
cannot define Gromov--Witten type invariants in $KH_*^\ef(Y;\Z)$ by
taking $X$ to be a moduli space $\oM_{g,m}(J,\be)$ of
$J$-holomorphic curves in a symplectic manifold $(M,\om)$, since
such $\oM_{g,m}(J,\be)$ may not be effective.

Problem (iv) means that products $\ubC_1\t_Y\ubC_2$ of effective
co-gauge-fixing data {\it cannot be defined to be commutative}, and
therefore {\it the cup product\/ $\cup$ on effective Kuranishi
cochains $KC^*_\ec(Y;R)$ will not be supercommutative}. Thus,
symmetry in Property \ref{kh3pr}(h) fails. We will deal with this by
replacing $P=\coprod_{k=0}^\iy\R^k/S_k$ in Definition \ref{kh3def6}
by $\uP=\coprod_{k=0}^\iy\R^k$, and taking effective
(co-)gauge-fixing data $\ubG,\ubC$ to include maps
$\uG^i,\uC^i:E^i\ra\uP$. The product $\umu:\uP\t\uP\ra\uP$
corresponding to $\mu$ in \eq{kh3eq2} is associative but not
commutative, and $\umu:\R^k\t\R^l\ra\R^{k+l}$ is injective, which
will be enough to preserve the injectivity conditions for effective
co-gauge-fixing data. In Remark \ref{kh4rem6} we explain a reason
from algebraic topology why $\cup$ cannot be supercommutative for
cohomology theories over~$R=\Z$.

We define {\it effective} Kuranishi neighbourhoods and Kuranishi
spaces.

\begin{dfn} An orbifold $Y$ is {\it effective} if for all orbifold
charts $(U,\Ga,\phi)$ on $Y$ in Definition \ref{kh2def6} with
$U\subset\R^n$, the representation of $\Ga$ on $\R^n$ is {\it
effective}. Here a representation is {\it effective} if every
non-identity element of $\Ga$ acts nontrivially on $\R^n$.
Equivalently, an orbifold is effective if generic $y\in Y$
have~$\Stab_Y(y)=\{1\}$.

A Kuranishi neighbourhood $(V,E,s,\psi)$ on $X$ is called {\it
effective\/} if $V$ is an effective orbifold, and for each $v\in V$,
the stabilizer group $\Stab_V(v)$ acts trivially on the fibre
$E\vert_v$ of $E$ over $v$. That is, in the notation of Definition
\ref{kh2def8}, if $(U,\Ga,\phi)$ is an orbifold chart on $V$ and
$(E_U,\Ga,\hat\phi)$ the associated orbifold chart on $E$, then for
all $u\in U$, the action of $\{\ga\in\Ga:\ga\cdot u=u\}$ on
$E_U\vert_u$ is trivial.

If $(V_p,\ldots,\psi_p),(V_q,\ldots,\psi_q)$ are Kuranishi
neighbourhoods on $X$, $(\phi_{pq},\hat\phi_{pq}):
(V_q,\ldots,\psi_q)\ra(V_p,\ldots,\psi_p)$ is a coordinate change,
and $(V_p,\ldots,\psi_p)$ is effective, it is easy to show that
$(V_q,\ldots,\psi_q)$ is also effective.

A Kuranishi space $X$ is called {\it effective} if for all $p\in X$,
all sufficiently small Kuranishi neighbourhoods
$(V_p,\ldots,\psi_p)$ in the germ at $p$ are effective.
\label{kh3def16}
\end{dfn}

Here are the reasons for this definition:
\begin{itemize}
\setlength{\itemsep}{0pt}
\setlength{\parsep}{0pt}
\item For proving isomorphism of effective Kuranishi homology
$KH_*^\ef(Y;R)$ with singular homology $H^\rsi_*(Y;R)$, chains
involving {\it effective\/} orbifolds are as good as chains
involving manifolds. An effective orbifold can be triangulated by
simplices to define a singular chain.
\item However, if $Y$ is a {\it non-effective} orbifold, then
for generic $y\in Y$, $\Stab(y)$ is some nontrivial finite group
$\Ga$. When $Y$ is triangulated by simplices to define a
singular chain, they should be weighted by $\md{\Ga}^{-1}$. This
only makes sense if $\md{\Ga}^{-1}$ lies in the coefficient ring
$R$, that is, if $R$ is a $\Q$-algebra.
\item If $(V,E,s,\psi)$ is an {\it effective} Kuranishi
neighbourhood, and $\ti s$ is a small generic perturbation of
$s$, then $\ti s$ is transverse and $\ti s^{-1}(0)$ is an
effective orbifold. This is because the obstruction to deforming
$s$ to transverse is nontrivial actions of the stabilizer groups
$\Stab(v)$ on the fibres $E\vert_v$ of $E$, and $\ti s^{-1}(0)$
is effective as $V$ is effective.
\item If $X$ is a compact, {\it effective\/} Kuranishi space,
then by choosing a finite cover of $X$ by effective Kuranishi
neighbourhoods $(V^i,E^i,s^i,\psi^i)$, and perturbing the $s^i$,
we can deform $X$ to a compact, effective orbifold $\ti X$,
which we can then triangulate by simplices to get a singular
chain.
\end{itemize}

We can now define {\it effective\/} gauge-fixing data.

\begin{dfn} Modifying $P,\mu$ in Definition \ref{kh3def6}, define
$\uP=\coprod_{k=0}^\iy\R^k$ and $\umu:\uP\t\uP\ra\uP$ by
\e
\umu\bigl((x_1,\ldots,x_k),(y_1,\ldots,y_l)\bigr)=
(x_1,\ldots,x_k,y_1,\ldots,y_l).
\label{kh3eq28}
\e
Then $\umu$ is {\it associative}, with {\it identity} $(\es)$, but
is {\it not\/} commutative. The obvious projection $\Pi:\uP\ra P$
mapping $\Pi:(x_1,\ldots,x_k)\mapsto S_k(x_1,\ldots,x_k)$ satisfies
$\mu\bigl(\Pi(p),\Pi(q)\bigr)=\Pi\bigl(\umu(p,q)\bigr)$
for~$p,q\in\uP$.

Let $X$ be a compact, effective Kuranishi space, $Y$ an orbifold,
and $\bs f:X\ra Y$ strongly smooth. A set of {\it effective
gauge-fixing data\/} for $(X,\bs f)$ is $\ubG=(\bs
I,\bs\eta,\uG^i:i\in I)$, where $(\bs I,\bs\eta)$ is an excellent
coordinate system for $(X,\bs f)$ with $\bs
I=\bigl(I,(V^i,E^i,s^i,\psi^i):i\in I,\ldots\bigr)$, and
$\uG^i:E^i\ra\uP$ for $i\in I$, satisfying:
\begin{itemize}
\setlength{\itemsep}{0pt}
\setlength{\parsep}{0pt}
\item[(a)] $(V^i,E^i,s^i,\psi^i)$ is an effective Kuranishi
neighbourhood for all $i\in I$; and
\item[(b)] Write $\pd^lX=X_1^l\amalg \cdots\amalg
X_{\smash{n^l}}^l$ for the splitting of $\bigl(\pd^lX,\bs
f\vert_{\pd^lX},\bs I\vert_{\pd^lX},\bs\eta\vert_{\pd^lX}\bigr)$
into connected components in Lemma \ref{kh3lem} for $l\ge 0$, and
$(\bs I^l_a,\bs\eta^l_a)=(\bs I\vert_{X^l_a},\bs\eta
\vert_{X^l_a})$, with $\bs I_a^l=\bigl(I_a^l,
(V^{l,i}_a,\ldots,\psi^{l,i}_a):i\in I_a^l,\ldots\bigr)$. Then
$(E^{l,i}_a)^\ci\subseteq E^{l,i}_a\subseteq\pd^lE^{i+l}$, where
$(E^{l,i}_a)^\ci$ is the {\it interior\/} of $E^{l,i}_a$, so we can
consider $\uG^{i+l}\vert_{(E^{l,i}_a)^\ci}:(E^{l,i}_a)^\ci\ra\uP$.
We require that $\uG^{i+l}\vert_{(E^{l,i}_a)^\ci}$ should map
$(E^{l,i}_a)^\ci$ to $\R^{k^l_a}\subset\uP$ for some $k^l_a\ge 0$
which depends on $l,a$ but {\it not\/} on $i\in I_a^l$, for all
$l\ge 0$, $a=1,\ldots,n^l$ and $i\in I_a^l$. We also require that
$\coprod_{i\in I_a^l}\uG^{i+l}\vert_{(E^{l,i}_a)^\ci}:\coprod_{i\in
I_a^l}(E^{l,i}_a)^\ci\ra\R^{k^l_a}\subset\uP$ should be {\it
injective} for all $l\ge 0$ and~$a=1,\ldots,n^l$.
\end{itemize}
Here (a) implies $X$ is effective, as for sufficiently small
$(V_p,\ldots,\psi_p)$ in the germ at $p\in X$ we have
$(\phi_p^i,\hat\phi_p^i):(V_p,\ldots,\psi_p)\ra(V^i,\ldots,\psi^i)$
for some $i\in I$, and $(V^i,\ldots,\psi^i)$ effective implies
$(V_p,\ldots,\psi_p)$ effective, as in Definition~\ref{kh3def16}.

Define {\it triples} $(X,\bs f,\ubG)$, {\it splittings} $X=X_1\amalg
\cdots\amalg X_n$ of triples $(X,\bs f,\ubG)$, {\it connected\/}
triples $(X,\bs f,\ubG)$, {\it isomorphisms} $(\bs a,\bs b):(X,\bs
f,\ubG)\ra(\ti X,\bs{\ti f},\ubtG)$, and {\it automorphism groups}
$\Aut(X,\bs f,\ubG)$ as in Definition~\ref{kh3def6}.

If $\ubG$ is effective gauge-fixing data for $(X,\bs f)$, we define
the {\it restriction} $\ubG\vert_{\pd X}$ as in Definition
\ref{kh3def7}. Since the conditions in (a),(b) above are stable
under taking boundaries, $V^i,E^i\mapsto\pd V^i,\pd E^i$,
$\ubG\vert_{\pd X}$ is effective gauge-fixing data for~$(\pd X,\bs
f\vert_{\pd X})$.

Let $\ubG=(\bs I,\bs\eta,\uG^i:i\in I)$ be effective gauge-fixing
data for $(X,\bs f)$. We shall show that $\Pi\ci\uG^i:E^i\ra P$ is a
{\it globally finite} map, where $\Pi:\uP\ra P$ is as above. To see
this, note that $E^i=\bigcup_{l\ge 0,\; a=1,\ldots,n^a_l:
E^{l,i-l}_a\ne\es}\io\bigl((E^{l,i-l}_a)^\ci\bigr)$, where
$\io:\pd^lE^i_a\ra E^i_a$ is the natural map, and $\uG^i$ maps
$(E^{l,i-l}_a)^\ci\ra\R^{k^l_a}\subset\uP$, and $\Pi:\R^{k^l_a}\ra
\R^{k^l_a}/S_{k^l_a}$ pulls back one point to at most $k^l_a!$
points. Hence if $p\in P$ then $\bmd{(\Pi\ci\uG^i)^{-1}(p)}$ is at
most $N$ points, where $N=\sum_{l\ge 0,\;
a=1,\ldots,n^a_l:E^{l,i-l}_a\ne\es}k^l_a!$. As this is a finite sum,
$\Pi\ci\uG^i$ is globally finite, and so $\Pi(\ubG)=\bs G=(\bs
I,\bs\eta,\Pi\ci\uG^i:i\in I)$ is {\it gauge-fixing data} for
$(X,\bs f)$. This gives a natural map from effective gauge-fixing
data to gauge-fixing data.
\label{kh3def17}
\end{dfn}

Here are analogues of Definitions \ref{kh3def16} and \ref{kh3def17}
for co-gauge-fixing data.

\begin{dfn} Let $X$ be a topological space, $(V,E,s,\psi)$ a
Kuranishi neighbourhood on $X$, $Y$ an orbifold, and $f:V\ra Y$ a
submersion. We call the pair $(V,E,s,\psi),f$ {\it coeffective} if
for all $v$ in $V$ with $y=f(v)$ in $Y$, the morphism
$f_*:\Stab_V(v)\ra\Stab_Y(y)$ is surjective, and the subgroup $\Ker
f_*\subseteq\Stab_V(v)$ acts effectively on $T_vV$, and $\Stab_V(v)$
acts trivially on the fibre $E\vert_v$ of $E$ at~$v$.

If $(V,E,s,\psi),f$ is coeffective and $e\in E$ with $\pi(e)=v$ in
$V$ and $f(v)=y$ in $Y$ then $\pi_*:\Stab_E(e)\ra\Stab_V(v)$ is an
isomorphism as $\Stab_V(v)$ acts trivially on $E\vert_v$, so
$(f\ci\pi)_*=f_*\ci\pi_*:\Stab_E(e)\ra\Stab_Y(y)$ is surjective.
This will help us deal with problem (iii) above.

Suppose $Y$ is an {\it effective\/} orbifold. Then $(V,E,s,\psi),f$
is coeffective if and only if $(V,E,s,\psi)$ is effective and
$f:V\ra Y$ is surjective on stabilizer groups. The only nontrivial
thing to prove is that $(V,E,s,\psi),f$ coeffective implies $V$ is
effective. If $v$ is generic in $V$ then $y=f(v)$ is generic in $Y$
as $f$ is a submersion, so $\Stab_Y(y)=\{1\}$ as $Y$ is effective,
giving $\Ker f_*=\Stab_V(v)$. Thus $\Stab_V(v)$ acts effectively on
$T_vV$, which implies that generic $v'$ close to $v$ in $V$ have
$\Stab_V(v')=\{1\}$. Hence generic points in $V$ have trivial
stabilizers, and $V$ is effective. When $Y$ is a {\it manifold},
$f:V\ra Y$ is trivially surjective on stabilizer groups, so
$(V,E,s,\psi),f$ is coeffective if and only if $(V,E,s,\psi)$ is
effective.

If $(V_p,\ldots,\psi_p),(V_q,\ldots,\psi_q)$ are Kuranishi
neighbourhoods on $X$, $(\phi_{pq},\hat\phi_{pq}):
(V_q,\ldots,\psi_q)\ra(V_p,\ldots,\psi_p)$ is a coordinate change,
$f_p:V_p\ra Y$ and $f_q:V_q\ra Y$ are submersions with $f_q\equiv
f_p\ci\phi_{pq}$ and $(V_p,\ldots,\psi_p),f_p$ is coeffective, it is
easy to show that $(V_q,\ldots,\psi_q),f_q$ is also coeffective.

Let $X$ be a Kuranishi space, $Y$ an orbifold, and $\bs f:X\ra Y$ a
strong submersion. We call $(X,\bs f)$ {\it coeffective} if for all
$p\in X$ and all sufficiently small Kuranishi neighbourhoods
$(V_p,\ldots,\psi_p)$ in the germ at $p$ with $f_p:V_p\ra Y$
representing $\bs f$, the pair $(V_p,\ldots,\psi_p),f_p$ is
coeffective.
\label{kh3def18}
\end{dfn}

\begin{dfn} Now let $X$ be a compact Kuranishi space,
$Y$ an orbifold, and $\bs f:X\ra Y$ a strong submersion with $(X,\bs
f)$ coeffective. A set of {\it effective co-gauge-fixing data\/} for
$(X,\bs f)$ is $\ubC=(\bs I,\bs\eta,\uC^i:i\in I)$, where $(\bs
I,\bs\eta)$ is an excellent coordinate system for $(X,\bs f)$ with
$\bs I=\bigl(I,(V^i,E^i,s^i,\psi^i),f^i:i\in I,\ldots\bigr)$, and
$\uC^i:E^i\ra\uP$ are maps for $i\in I$, satisfying
\begin{itemize}
\setlength{\itemsep}{0pt}
\setlength{\parsep}{0pt}
\item[(a)] the pair $(V^i,E^i,s^i,\psi^i),f^i$ is coeffective for
all $i\in I$; and
\item[(b)] Write $\pd^lX=X_1^l\amalg \cdots\amalg
X_{\smash{n^l}}^l$ for the splitting of $\bigl(\pd^lX,\bs
f\vert_{\pd^lX},\bs I\vert_{\pd^lX},\bs\eta\vert_{\pd^lX}\bigr)$
into connected components in Lemma \ref{kh3lem} for $l\ge 0$, and
$(\bs I^l_a,\bs\eta^l_a)=(\bs I\vert_{X^l_a},\bs\eta\vert_{X^l_a})$,
with $\bs I_a^l=\bigl(I_a^l, (V^{l,i}_a, \ldots,\psi^{l,i}_a):i\in
I_a^l,\ldots\bigr)$. Then $(E^{l,i}_a)^\ci\subseteq
E^{l,i}_a\subseteq\pd^lE^{i+l}$, so we can consider
$\uC^{i+l}\vert_{(E^{l,i}_a)^\ci}:(E^{l,i}_a)^\ci\ra\uP$. We require
that $\uC^{i+l}\vert_{(E^{l,i}_a)^\ci}$ should map $(E^{l,i}_a)^\ci$
to $\R^{k^l_a}\subset\uP$ for some $k^l_a\ge 0$ which depends on
$l,a$ but {\it not\/} on $i\in I_a^l$, for all $l\ge 0$,
$a=1,\ldots,n^l$ and $i\in I_a^l$. We also require that
$\coprod_{i\in I_a^l}\bigl(\uC^{i+l}\t(f^{i+l}\ci\pi^{i+l})\bigr)
\vert_{(E^{l,i}_a )^\ci}:\coprod_{i\in I_a^l}(E^{l,i}_a)^\ci\ra\uP\t
Y$ should be {\it injective} for all $l\ge 0$ and~$a=1,\ldots,n^l$.
\end{itemize}
Here (a) implies $(X,\bs f)$ is coeffective, as for sufficiently
small $(V_p,\ldots,\psi_p),f_p$ in the germ at $p\in X$ we have
$(\phi_p^i,\hat\phi_p^i):(V_p,\ldots,\psi_p)\ra(V^i,\ldots, \psi^i)$
for some $i\in I$ with $f_p\equiv f^i\ci\phi_p^i$, and
$(V^i,\ldots,\psi^i),f^i$ coeffective implies
$(V_p,\ldots,\psi_p),f_p$ coeffective as in
Definition~\ref{kh3def18}.

As in Definition \ref{kh3def17} we define {\it triples} $(X,\bs
f,\ubC)$, {\it splittings} $X=X_1\amalg \cdots\amalg X_n$ of $(X,\bs
f,\ubC)$, {\it connected\/} $(X,\bs f,\ubC)$, {\it isomorphisms}
$(\bs a,\bs b):(X,\bs f,\ubC)\ra(\ti X,\bs{\ti f},\ubtC)$, {\it
automorphism groups} $\Aut(X,\bs f,\ubC)$, {\it restrictions\/}
$\ubC\vert_{\pd X}$, and a natural map $\ubC\mapsto\Pi(\ubC)$ from
effective co-gauge-fixing data to co-gauge-fixing data.
\label{kh3def19}
\end{dfn}

Here are partial analogues of Theorems \ref{kh3thm1}, \ref{kh3thm2}
and \ref{kh3thm3}, which show that versions of Property
\ref{kh3pr}(a),(b) and (e) hold for effective (co-)gauge-fixing
data.

\begin{thm}{\bf(a)} Let\/ $X$ be a compact, effective Kuranishi
space with\/ $\pd^2X=\es,$ $Y$ an orbifold, and\/ $\bs f:X\ra Y$ a
strongly smooth map. Then\/ $(X,\bs f)$ admits effective
gauge-fixing data\/~$\ubG$.

\noindent{\bf(b)} If\/ $\ubG$ is effective gauge-fixing data for
$(X,\bs f)$ and\/ $(X,\bs f,\ubG)$ is connected then $\Aut(X,\bs
f,\ubG)=\{1\}$.

\noindent{\bf(c)} Suppose\/ $X$ is a compact, effective Kuranishi
space with corners, $Y$ an orbifold, $\bs f:X\ra Y$ a strongly
smooth map, and\/ $\ubH$ is effective gauge-fixing data for\/ $(\pd
X,\bs f\vert_{\pd X})$ with corners. Let\/ $\bs\si:\pd^2X\ra\pd^2X$
be the involution of Definition {\rm\ref{kh2def19}.} Then there
exists effective gauge-fixing data\/ $\ubG$ for\/ $(X,\bs f)$ with
corners with\/ $\ubG\vert_{\pd X}\!=\!\ubH$ if and only if\/ $\ubH
\vert_{\pd^2X}$ is\/ $\bs\si$-invariant.

Parts\/ {\rm(a)--(c)} also hold for effective co-gauge-fixing data
$\ubC,\ubD,$ taking $\bs f$ to be a strong submersion and $(X,\bs
f)$ to be coeffective rather than $X$ effective.
\label{kh3thm5}
\end{thm}

\begin{proof} For (a), we modify the proof of Theorem \ref{kh3thm1}.
In the first part, since $X$ is effective and $\pd^2X=\es$, in
Proposition \ref{kh3prop2} we can choose the very good coordinate
system $\bs I=\bigl(I,(V^i,E^i,s^i,\psi^i),f^i:i\in I),\ldots\bigr)$
for $(X,\bs f)$ with $(V^i,\ldots,\psi^i)$ effective and
$\pd^2V^i=\es$ for all $i\in I$. Then Definition \ref{kh3def17}(a)
holds. In Theorem \ref{kh3thm1} we choose the $\uG^i:E^i\ra\uP$ such
that $\coprod_{i\in I}\uG^i:\coprod_{i\in I}E^i\ra\R^k\subset\uP$ is
injective, for some $k\gg 0$. As $\pd^2V^i=\es$ we have
$\pd^2E^i=\es$, so the immersion $\io:\pd E^i\ra E^i$ is injective,
rather than merely finite. Hence $\coprod_{i\in I}\uG^i\vert_{\pd
E^i}:\coprod_{i\in I}\pd E^i\ra\uP$ is injective. As $\pd^lE^i=\es$
for $l\ge 2$, Definition \ref{kh3def17}(b) follows.

Part (b) is a simplification of the proof of Theorem \ref{kh3thm2}:
for $(\bs a,\bs b)$ in $\Aut(X,\bs f,\ubG)$, as $\uG^i\vert_{
(E^i)^\ci} :(E^i)^\ci\!\ra\!\uP$ is injective, the diffeomorphisms
$\hat b^i: E^i\ra E^i$ with $\ubG^i\ci\hat b^i\equiv \ubG^i:E^i\ra
P$ must be the identity on $(E^i)^\ci$, so $\hat b^i$ is the
identity by continuity. Therefore $b^i:V^i\ra V^i$ is the identity,
as $\hat b^i$ lifts $b^i$. Since $(b^i,\hat b^i)$ determines $\bs a$
on $\Im\psi^i$ and $X=\bigcup_{i\in I}\Im\psi^i$ this gives $\bs
a=\bs\id_X$, so $\Aut(X,\bs f,\ubG)=\{1\}$. Part (c) is also a
simple modification of Theorem \ref{kh3thm3}: for the `if' part, if
$\ubH\vert_{\pd^2X}$ is $\bs\si$-invariant, then $\ubH$ prescribes
consistent values for $\uG^i:E^i\ra\uP$ on $E^i\sm(E^i)^\ci$ for all
$i\in I$, satisfying Definition \ref{kh3def17}(b) for $l\ge 1$. We
then choose $\uG^i\vert_{\smash{(E^i)^\ci}}:(E^i)^\ci\ra\uP$ for
$i\in I$ such that $\coprod_{i\in I}\uG^i\vert_{\smash{(E^i)^\ci}}:
\coprod_{i\in I}(E^i)^\ci\ra\R^k\subset\uP$ is injective for
some~$k\gg 0$.
\end{proof}

\begin{rem} In Theorem \ref{kh3thm5}(a) we omit the parts on finite
groups $\Ga$ in Property \ref{kh3pr}(a), since Theorem
\ref{kh3thm5}(b) shows that we cannot have nontrivial symmetry
groups $\Ga$ of $(X,\bs f,\ubG)$ for connected $X$. We also suppose
$\pd^2X=\es$. This can be weakened, but some assumption on $\pd^2X$
is necessary for effective gauge-fixing data to exist on $(X,\bs
f)$, as the following argument shows.

Suppose that $\ubG$ is effective gauge-fixing data for $(X,\bs f)$
with $X$ oriented, and that $\pd^2X$ is nonempty, and connected as a
topological space. Restrict $\ubG$ to $\pd X$ and $\pd^2X$ to get
effective gauge-fixing data $\ubG\vert_{\pd^2X}$ for $(\pd^2X,\bs
f\vert_{\pd^2X})$. Then $(\pd^2X,\bs f\vert_{\pd^2X},\ubG
\vert_{\pd^2X})$ is a connected triple, as $\pd^2X$ is a connected
topological space, so $\Aut(\pd^2X,\bs f\vert_{\pd^2X},\ubG
\vert_{\pd^2X})=\{1\}$ by Theorem \ref{kh3thm5}(b). However, the
natural involution $\bs\si:\pd^2X\ra\pd^2X$ clearly preserves
$\ubG\vert_{\pd^2X}$ and lies in $\Aut(\pd^2X,\bs
f\vert_{\pd^2X},\ubG \vert_{\pd^2X})=\{1\}$, a contradiction, as
$\bs\si$ is orientation-reversing.
\label{kh3rem5}
\end{rem}

Following Definitions \ref{kh3def12} and \ref{kh3def13}, we can
define {\it pushforwards\/} $h_*(\ubG)$ and {\it pullbacks\/}
$h^*(\ubC)$ of effective (co-)gauge-fixing data~$\ubG,\ubC$.

\begin{dfn} Let $Y,Z$ be orbifolds and $h:Y\ra Z$ a smooth map.
Suppose $X$ is a compact, effective Kuranishi space, $\bs f:X\ra Y$
is strongly smooth, and $\ubG=(\bs I,\bs\eta,\uG^i:i\in I)$ is
effective gauge-fixing data for $(X,\bs f)$, where $\bs
I=\bigl(I,(V^i,E^i,s^i,\psi^i),\ab f^i:i\in I,\ldots\bigr)$. Then
$h\ci\bs f:X\ra Z$ is strongly smooth. Define $h_*(\ubG)$ to be
$h_*(\ubG)=(h_*(\bs I),\bs\eta,\uG^i:i\in I)$, where $h_*(\bs
I)=\bigl(I,(V^i,E^i,s^i,\psi^i),h\ci f^i:i\in I,\ldots\bigr)$. Then
$h_*(\ubG)$ is effective gauge-fixing data for~$(X,h\ci\bs f)$.
\label{kh3def20}
\end{dfn}

\begin{dfn} Let $Y,Z$ be orbifolds without boundary and $h:Y\ra Z$
a smooth, proper map. Suppose $X$ is a compact Kuranishi space, $\bs
f:X\ra Y$ a strong submersion with $(X,\bs f)$ coeffective, and
$\ubC=(\bs I,\bs\eta,\uC^i:i\in I)$ is effective co-gauge-fixing
data for $(X,\bs f)$, where $\bs
I=\bigl(I,(V^i,E^i,s^i,\psi^i),f^i:i\in I,\ldots\bigr)$. Then
$Y\t_{h,Z,\bs f}X$ is a compact Kuranishi space, and
$\bs\pi_Y:Y\t_ZX\ra Y$ is a strong submersion.

Define an excellent coordinate system $(\bs{\check I},\bs{\check
\eta})$ for $(Y\t_ZX,\bs\pi_Y)$ as in Definition \ref{kh3def13}.
Then $\check I\subseteq\{i+k:i\in I\}$, where $k=\dim Y-\dim Z$, and
for each $i\in\check I$, $\check V^i$ is an open subset of
$Y\t_{h,Z,f^{i-k}}V^{i-k}$, and $\check E^i$ an open subset of
$Y\t_{h,Z,f^{i-k}\ci\pi^{i-k}}E^{i-k}$. Define $\ucC^i:\check
E^i\ra\uP$ by $\ucC^i=\uC^{i-k}\ci\pi_{E^{i-k}}$. We claim that
$h^*(\ubC)=(\bs{\check I},\bs{\check\eta},\ucC^i:i\in\check I)$ is
{\it effective co-gauge-fixing data\/} for~$(Y\t_ZX,\bs\pi_Y)$.

We must verify Definition \ref{kh3def19}(a),(b) hold. For (a), as
the pair $(V^{i-k},\ab\ldots,\ab\psi^{i-k}),f^{i-k}$ is coeffective
and $\check V^i\subseteq Y\t_{h,Z,f^{i-k}}V^{i-k}$, $\check
E^i\subseteq Y\t_{h,Z,f^{i-k}\ci\pi^{i-k}}E^{i-k}$, one can show
directly that the pair $(\check
V^i,\ldots,\check\psi^i),\check\pi_Y^i$ is coeffective. This then
implies that $(Y\t_ZX,\bs\pi_Y)$ is coeffective, as in
Definition~\ref{kh3def19}.

For (b), if $\pd^lX=\coprod_{a=1}^{m^l}X_a^l$ and
$\pd^l(Y\t_ZX)=\coprod_{b=1}^{n^l}(Y\t_ZX)_b^l$ are the splittings
of $\bigl(\pd^lX,\ab\bs f\vert_{\pd^lX},\bs I\vert_{\pd^lX},\bs\eta
\vert_{\pd^lX}\bigr)$ and $\bigl(\pd^l(Y\t_ZX),\bs\pi_Y
\vert_{\pd^l(Y\t_ZX)},\bs{\check I}\vert_{\pd^l(Y\t_ZX)},\ab
\bs{\check\eta}\vert_{\pd^l(Y\t_ZX)}\bigr)$ into connected
components in Lemma \ref{kh3lem}, then for each $b=1,\ldots,n^l$,
$(Y\t_ZX)_b^l$ is a connected component of $\bigl(Y\t_Z
(X_a^l),\ab\bs\pi_Y\vert_{Y\t_ZX_a^l},\ab\bs{\check I} \vert_{Y\t_Z
X_a^l},\ab\bs{\check\eta}\vert_{Y\t_ZX_a^l}\bigr)$ for some
$a=1,\ldots,m^l$. If $(\check V^{l,i}_b,\check E^{l,i}_b,\check
s^{l,i}_b,\check \psi^{l,i}_b)$ is a Kuranishi neighbourhood on
$(Y\t_ZX)_b^l$ in $\bs{\check I}{}_b^l$, then $(\check
E^{l,i}_b)^\ci\subseteq Y\t_Z(E^{l,i-k}_a)^\ci$, and
$\ucC^{i+l}\vert_{(\check E^{l,i}_b)^\ci}\equiv\uC^{i-k+l}\ci
\pi_{(E^{l,i-k}_a)^\ci}$. As $\uC^{i-k+l}$ maps $(E^{l,i-k}_a)^\ci
\ra\R^{k^l_a}\subset\uP$, it follows that $\ucC^{i+l}$ maps $(\check
E^{l,i}_b)^\ci\ra\R^{k^l_a}\subset\uP$, giving the first part of
Definition \ref{kh3def19}(b). For the second, since $\coprod_{i-k\in
I_a^l}\bigl(\uC^{i-k+l}\t
(f^{i-k+l}\ci\pi^{i-k+l})\bigr)\vert_{(E^{l,i-k}_a
)^\ci}:\coprod_{i-k\in I_a^l}(E^{l,i-k}_a)^\ci\ra\uP\t Z$ is
injective by Definition \ref{kh3def19}(b), and $\coprod_{i\in\check
I_b^l}(\check E^{l,i}_b)^\ci\subset Y\t_Z \bigl(\coprod_{i-k\in
I_a^l}(E^{l,i-k}_a)^\ci\bigr)$, it follows that $\coprod_{i\in\check
I_b^l}\bigl(\ucC^{i+l} \t (\check
f^{i+l}\ci\check\pi^{i+l})\bigr)\vert_{(\check E^{l,i}_b)^\ci}:
\coprod_{i\in\check I_b^l}(\check E^{l,i}_b)^\ci\ra\uP\t Y$ is
injective, as we want.

This last part is a little more subtle than it appears. The
important point is this: suppose $(V,E,s,\psi)$ is a Kuranishi
neighbourhood on $X$, $f:V\ra Z$ a submersion with $(V,E,s,\psi),f$
coeffective, and $\uC:E\ra\uP$ is a map with
$\uC\t(f\ci\pi):E\ra\uP\t Z$ injective. As $(V,E,s,\psi),f$ is
coeffective, $f\ci\pi:E\ra Z$ is surjective on stabilizer groups.
Therefore the biquotient terms $f_*(\Stab(p))\backslash\Stab(f(p))
/f'_*(\Stab(p'))$ in the definition \eq{kh2eq10} of the orbifold
fibre product $Y\t_{h,Z,f\ci\pi}E$ are all points, so that
$\pi_E\t\pi_Y:Y\t_ZE\ra E\t Y$ is injective. Combining this with
$\uC\t(f\ci\pi)$ injective implies that
$(\uC\ci\pi_E)\t\pi_Y:Y\t_ZE\ra\uP\t Y$ is injective. Thus, we need
the coeffective assumptions for $\ubC$ in Definition
\ref{kh3def19}(a) to prove the injectivity for $h^*(\ubC)$ in
Definition \ref{kh3def19}(b).
\label{kh3def21}
\end{dfn}

Pushforwards and pullbacks are {\it functorial}, that is, $(g\ci
h)_*(\ubG)=g_*\ci h_*(\ubG)$, and $(g\ci h)^*(\ubC)=h^*\ci
g^*(\ubC)$, as for (co-)gauge-fixing data. We modify Definition
\ref{kh3def15} to define {\it fibre products} $\ubC\t_Y\ubtC$,
$\ubG\t_Y\ubtC$ of effective (co-)gauge-fixing data.

\begin{dfn} Let $X,\ti X$ be compact Kuranishi spaces, $Y$ an
orbifold, and $\bs f:X\ra Y$, $\bs{\ti f}:\ti X\ra Y$ be strong
submersions with $(X,\bs f),(\ti X,\bs{\ti f})$ coeffective. Suppose
$\ubC=(\bs I,\bs\eta,\uC^i:i\in I)$ and $\ubtC=(\bs{\ti
I},\bs{\ti\eta}, \utC^{\ti\imath}:\ti\imath\in\ti I)$ are effective
co-gauge-fixing data for $(X,\bs f)$ and $(\ti X,\bs{\ti f})$. Let
$(\bs I,\bs\eta)\t_Y(\bs{\ti I},\bs{\ti\eta})$ be the really good
coordinate system for $(X\t_Y\ti X,\bs\pi_Y)$ given in Definition
\ref{kh3def14}, and $(\bs{\check I},\bs{\check\eta})=((\bs
I,\bs\eta)\t_Y(\bs{\ti I},\bs{\ti\eta}))\kern .1em\check{}\,$ the
excellent coordinate system for $(X\t_Y\ti X,\bs\pi_Y)$ constructed
from $(\bs I,\bs\eta)\t_Y(\bs{\ti I},\bs{\ti\eta})$ by Algorithm
\ref{kh3alg}. In the notation of Definitions \ref{kh3def14} and
\ref{kh3def15} we have $\bs{\check I}=\bigl(\check I,\ab(\check
V^k,\ab\check E^k,\ab\check s^k,\ab\check\psi^k), \check\pi_Y^k:
k\in\check I,\ldots\bigr)$, with $\check V^k,\check E^k$ open sets
in $\coprod_{i\in I,\; \ti\imath\in\ti I:i+\ti\imath=k+\dim
Y}V^{(i,\ti\imath)}$, $\coprod_{i\in I,\; \ti\imath\in\ti
I:i+\ti\imath=k+\dim Y}E^{(i,\ti\imath)}$, where $V^{(i,\ti\imath)}=
V^i\t_{f^i,Y,\ti f^{\ti\imath}}\ti V^{\ti\imath}$
and~$E^{(i,\ti\imath)}=E^i\t_{f^i\ci\pi^i,Y,\ti f^{\ti\imath}
\ci\ti\pi^{\ti\imath}}\ti E^{\ti\imath}$.

For $k\in\check I$, define $\ucC^k:\check E^k\ra\uP$ by
$\ucC^k(\check e)=\umu\bigl(\uC^i\ci\pi_{E^i}(\check e),
\utC^{\ti\imath}\ci\pi_{\ti E^{\ti\imath}}(\check e)\bigr)$ whenever
$\check e$ lies in $\check E^k\cap E^{(i,\ti\imath)}$, for all $i\in
I$ and $\ti\imath\in\ti I$ with $i+\ti\imath=k+\dim Y$, where
$\pi_{E^i}:E^{(i,\ti\imath)}\ra E^i$ and $\pi_{\ti
E^{\ti\imath}}:E^{(i,\ti\imath)}\ra\ti E^{\ti\imath}$ are the
projections from $E^{(i,\ti\imath)}=E^i\t_Y\ti E^{\ti\imath}$, and
$\umu:\uP\t\uP\ra\uP$ is as in \eq{kh3eq28}. We shall show that
$\ubC\t_Y\ubtC=\ubcC=(\bs{\check I},\bs{\check\eta},\check
\uC^i:i\in\check I)$ is {\it effective co-gauge-fixing data}
for~$(X\t_Y\ti X,\bs\pi_Y)$.

We must verify Definition \ref{kh3def19}(a),(b) hold. The proof is
similar to that for $h^*(\ubC)$ in Definition \ref{kh3def21}. For
(a), as $(V^i,\ldots,\psi^i),f^i$ and $(\ti
V^{\ti\imath},\ldots,\ti\psi^{\ti\imath}),\ti f^{\ti\imath}$ are
coeffective one can show directly that $(V^{(i,\ti\imath)},\ldots,
\psi^{(i,\ti\imath)}),\pi_Y^{(i,\ti\imath)}$ in \eq{kh3eq6} is
coeffective. Therefore $(\check V^k,\ldots,\check\psi^k),
\check\pi_Y^k$ is coeffective over the open set $\check V^k\cap
V^{(i,\ti\imath)}$ in $V^{(i,\ti\imath)}$, and as this holds for all
$i,\ti\imath$ with $i+\ti\imath=k+\dim Y$, $(\check
V^k,\ldots,\check\psi^k),\check\pi_Y^k$ is coeffective. This
proves~(a).

For Definition \ref{kh3def19}(b), first consider the case $l=0$, and
suppose $(X,\bs f,\ubC),\ab(\ti X,\ab\bs{\ti f},\ubtC)$ are
connected. Then Definition \ref{kh3def19}(b) implies that $\uC^i$
maps $(E^i)^\ci\ra\R^{k^0_1}\subset\uP$ for some $k^0_1\ge 0$
independent of $i$ and all $i\in I$, and $\utC^{\ti\imath}$ maps
$(\ti E^{\ti\imath})^\ci\ra\R^{\ti k^0_1}\subset\uP$ for some $\ti
k^0_1\ge 0$ independent of $\ti\imath$ and all $\ti\imath\in\ti I$.
As $(\check E^k)^\ci\subseteq\coprod_{i\in I,\; \ti\imath\in\ti
I:i+\ti\imath=k+\dim Y}(E^i)^\ci\t_Y(\ti E^{\ti\imath})^\ci$, it
follows from \eq{kh3eq28} that $\ucC^k$ maps $(\check E^k)^\ci
\ra\R^{k^0_1+\ti k^0_1}\subset\uP$ for all $k\in\check I$, where
$k^0_1+\ti k^0_1$ is independent of $k\in\check I$. This proves the
first part of Definition \ref{kh3def19}(b) in this case.

For the second part, observe that $\coprod_{k\in\check I}(\check
\uC^k \t(\check f^k\ci\check\pi^k))$ is the composition
\begin{equation*}
\xymatrix@C=12pt{\coprod\limits_{k\in\check I}(\check E^k)^\ci
\ar@{^{(}->}[r]
& {\begin{subarray}{l}\bigl(\coprod\limits_{i\in I}(E^i)^\ci\bigr)\t_Y\\
\bigl(\coprod\limits_{\ti\imath\in\ti I}(\ti E^{\ti\imath})^\ci\bigr)
\end{subarray}}
\ar[rrrr]^{\begin{subarray}{l}(\coprod_{i\in I}\uC^i\ci\pi_{E^i})\t\\
(\coprod_{\ti\imath\in\ti I}\utC^{\ti\imath}\ci \pi_{\ti
E^{\ti\imath}})\t\pi_Y\end{subarray}\!\!\!}
&&&&\R^{k^0_1}\!\t\!\R^{\ti k^0_1}\!\t\!Y \ar[rr]^{\umu\t\id_Y}
&&\R^{k^0_1+\ti k^0_1}\t Y.}
\end{equation*}
Here the first arrow `$\hookra$' is an inclusion of open sets, and
so injective. The second arrow is injective as $\coprod_{i\in
I}(\uC^i\t(f^i\ci\pi^i)): \coprod_{i\in I}(E^i)^\ci\ra\R^{k^0_1}\t
Y\subset\uP\t Y$ is injective, and $\coprod_{\ti\imath\in\ti I}(\ti
\uC^{\ti\imath}\t(\ti f^{\ti\imath}\ci\ti\pi^{\ti\imath})):
\coprod_{\ti\imath\in\ti I} (\ti E^{\ti\imath})^\ci\ra\R^{\ti
k^0_1}\t Y\subset\uP\t Y$ is injective, and $\pi_{E^i}\t\pi_{\ti
E^{\ti\imath}}:E^i\t_Y\ti E^{\ti\imath}\ra E^i\t\ti E^{\ti\imath}$
is injective as $f^i\ci\pi^i:E^i\ra Y$, $\ti f^{\ti\imath}\ci
\ti\pi^{\ti\imath}:\ti E^{\ti\imath}\ra Y$ are surjective on
stabilizer groups. The third arrow is injective as
$\umu:\R^{k^0_1}\t\R^{\ti k^0_1}\ra\R^{k^0_1+\ti k^0_1}$ in
\eq{kh3eq28} is injective.

Therefore the composition $\coprod_{k\in\check I}(\ucC^k \t(\check
f^k\ci\check\pi^k))$ is injective, and Definition \ref{kh3def19}(b)
holds in the case when $l=0$ and $(X,\bs f,\ubC),\ab(\ti
X,\ab\bs{\ti f},\ab\ubtC)$ are connected. For the general case, we
apply the argument above with connected components $\bigl(X_a^j,\bs
f\vert_{X_a^j},\ubC\vert_{X_a^j}\bigr)$ of $\bigl(\pd^jX,\bs
f\vert_{\pd^jX},\ubC\vert_{\pd^jX}\bigr)$ and $\bigl(\ti X_{\ti
a}^{\ti\jmath},\bs{\ti f}\vert_{\ti X_{\ti
a}^{\ti\jmath}},\ubtC\vert_{\ti X_{\ti a}^{\ti\jmath}}\bigr)$ of
$\bigl(\pd^{\ti\jmath}\ti X,\bs{\ti f} \vert_{\pd^{\ti\jmath}\ti
X},\ab\ubtC\vert_{\pd^{\ti\jmath}\ti X}\bigr)$ in place of $(X,\bs
f,\ubC),(\ti X,\bs{\ti f},\ubtC)$, with $j+\ti\jmath=l$. This proves
Definition \ref{kh3def19}(b), so $\ubC\t_Y\ubtC$ is {\it effective
co-gauge-fixing data}.

Now let $X$ be a compact, effective Kuranishi space, $Y$ an
orbifold, $\bs f:X\ra Y$ strongly smooth, $\ubG$ effective
gauge-fixing data for $(X,\bs f)$, $\ti X$ a compact Kuranishi
space, $\bs{\ti f}:\ti X\ra Y$ a strong submersion with $(\ti
X,\bs{\ti f})$ coeffective, and $\ubtC$ be effective co-gauge-fixing
data for $(\ti X,\bs{\ti f})$. Then we can modify Definition
\ref{kh3def15} in the same way to define effective gauge-fixing data
$\ubG\t_Y\ubtC$ for $(X\t_Y\ti X,\bs\pi_Y)$, defining $\ucG^k:\check
E^k\ra\uP$ by $\ucG^k(\check e)=\umu\bigl(\uG^i\ci\pi_{E^i}(\check
e), \utC^{\ti\imath}\ci\pi_{\ti E^{\ti\imath}}(\check e)\bigr)$. The
proof that $\ubG\t_Y\ubtC$ satisfies Definition \ref{kh3def17} is
similar to the above.

The maps $\ubG\mapsto\Pi(\ubG)$, $\ubC\mapsto\Pi(\ubC)$ from
effective (co-)gauge-fixing data to (co-)\ab gauge-fixing data in
Definitions \ref{kh3def17} and \ref{kh3def19} are compatible with
fibre products, that is,
$\Pi\bigl(\ubC\t_Y\ubtC)=\Pi(\ubC)\t_Y\Pi(\ubtC)$ and $\Pi\bigl(
\ubG\t_Y\ubtC)=\Pi(\ubG)\t_Y\Pi(\ubtC)$, because $\mu\bigl(\Pi(p),
\Pi(q)\bigr)=\Pi\bigl(\umu(p,q)\bigr)$ for~$p,q\in\uP$.
\label{kh3def22}
\end{dfn}

From Propositions \ref{kh3prop6} and \ref{kh3prop7}, with proofs
almost unchanged, we deduce:

\begin{prop} Effective (co-)gauge-fixing data satisfies analogues of
Propositions {\rm\ref{kh3prop6}(b),(c),(d)} and\/
{\rm\ref{kh3prop7},} taking $X_a$ effective or $(X_a,\bs f_a)$
coeffective where appropriate. However, the analogue of Proposition
{\rm\ref{kh3prop6}(a)} does not apply, that is, fibre products of
effective co-gauge-fixing data are not commutative, since $\umu$
in\/ \eq{kh3eq28} is not commutative.
\label{kh3prop8}
\end{prop}

Problem (iv) at the beginning of \S\ref{kh39} implies that it is
necessary for fibre products $\ubC\t_Y\ubtC$ of effective
co-gauge-fixing data not to be commutative. We give another reason,
involving {\it Steenrod squares}, in Remark~\ref{kh4rem6}.

\section{Kuranishi homology and cohomology}
\label{kh4}

We now define two slightly different homology theories of an
orbifold $Y$: {\it Kuranishi homology\/} $KH_*(Y;R)$, which is
defined when $R$ is a $\Q$-algebra, and {\it effective Kuranishi
homology} $KH_*^\ef(Y;R)$, which is defined when $R$ is a
commutative ring. Both of them are isomorphic to {\it singular
homology} $H^\rsi_*(Y;R)$. They are variations of the same idea,
which is to define a homology theory using as chains strongly smooth
maps $\bs f:X\ra Y$ from a compact Kuranishi space~$X$.

However, we show in \S\ref{kh49} that just using isomorphism classes
$[X,\bs f]$ of pairs $(X,\bs f)$ as chains would make $KH_*(Y;R)=0$
for all $Y,R$. The problems are caused by pairs $(X,\bs f)$ whose
automorphism groups are infinite. To avoid them, we include ({\it
effective\/}) {\it gauge-fixing data\/} $\bs G$ on $X$, as in
Chapter \ref{kh3}, and take chains to be generated by isomorphism
classes $[X,\bs f,\bs G]$ of triples~$(X,\bs f,\bs G)$.

We also define Poincar\'e dual cohomology theories, {\it Kuranishi
cohomology} $KH^*(Y;R)$ for $R$ a $\Q$-algebra, and {\it effective
Kuranishi cohomology} $KH^*_\ec(Y;R)$ for $R$ a commutative ring.
Here $KH^*(Y;R)$ is always isomorphic to {\it compactly supported
cohomology} $H^*_\cs(Y;R)$, and $KH^*_\ec(Y;R)\cong H^*_\cs(Y;R)$ if
$Y$ is a manifold. The main point of introducing them is as a tool
to use in areas where Kuranishi spaces naturally arise, such as
closed or open Gromov--Witten theory, or Lagrangian Floer
cohomology.

Kuranishi (co)homology $KH_*,KH^*(Y;R)$ is better behaved at the
(co)chain level than effective Kuranishi (co)homology
$KH_*^\ef,KH^*_\ec(Y;R)$. For example, the cup product $\cup$ on
Kuranishi cochains $KC^*(Y;R)$ is supercommutative, but the cup
product $\cup$ on effective Kuranishi cochains $KC^*_\ec(Y;R)$ is
not. Also, we can form chains $[X,\bs f,\bs G]$ in $KC_*(Y;R)$ for
$X$ an arbitrary compact oriented Kuranishi space, but for
$KC_*^\ef(Y;R)$, $X$ must satisfy restrictions on its stabilizer
groups and corners.

The advantage of effective Kuranishi (co)homology
$KH_*^\ef,KH^*_\ec(Y;R)$ is that they work over any commutative ring
$R$, such as $\Z$, not just over $\Q$-algebras. This will be
important in the author's approach to the integrality conjecture for
Gopakumar--Vafa invariants, to be discussed in Chapter \ref{kh6}. We
will also use effective Kuranishi homology in showing Kuranishi
homology is isomorphic to singular homology, as we will first prove
that $KH_*^\ef(Y;R)\cong H^\rsi_*(Y;R)$ for $R$ a commutative ring,
and then that $KH_*(Y;R)\cong KH_*^\ef(Y;R)$ for $R$ a $\Q$-algebra.

\subsection{Classical homology and cohomology}
\label{kh41}

Let $Y$ be a topological space and $R$ a commutative ring. Then
there are many ways of defining {\it homology groups} $H_k(Y;R)$ of
$Y$ with coefficients in $R$, for $k=0,1,2,\ldots$, for example,
singular homology, cellular homology (of a CW-complex), and
Borel--Moore homology. A general feature of all these theories is
that if $Y$ is sufficiently well-behaved, for instance, if $Y$ is a
manifold, then the homology groups obtained from any two homology
theories are canonically isomorphic. (Theorems \ref{kh4thm1} and
\ref{kh4thm2} below are results of this kind.) So for many purposes,
it does not matter which homology theory one uses.

There are also many ways of defining {\it cohomology groups}
$H^k(Y;R)$ and {\it compactly-supported cohomology groups}
$H^k_\cs(Y;R)$, for example, singular cohomology, \v Cech
cohomology, and de Rham cohomology (of a smooth manifold, over
$\R$). Again, if $Y$ is sufficiently well-behaved, then the
(compactly-supported) cohomology groups obtained from any two
cohomology theories are isomorphic.

We briefly recall the definition of {\it singular homology}, from
Bredon~\cite[\S IV]{Bred}.

\begin{dfn} For $k=0,1,\ldots$, the $k$-{\it simplex\/} $\De_k$ is
\e
\De_k=\bigl\{(x_0,\ldots,x_k)\in\R^{k+1}:x_i\ge 0,\;\>
x_0+\cdots+x_k=1\bigr\}.
\label{kh4eq1}
\e
It is a compact, oriented $k$-manifold, with boundary and corners.

Let $Y$ be a topological space, and $R$ a commutative ring. Define
$C_k^\rsi(Y;R)$ to be the $R$-module spanned by {\it singular
simplices}, which are continuous maps $\si:\De_k\ra Y$. Elements of
$C_k^\rsi(Y;R)$, which are called {\it singular chains}, are finite
sums $\sum_{a\in A}\rho_a\,\si_a$, where $A$ is a finite indexing
set, $\rho_a\in R$, and $\si_a:\De_k\ra Y$ is continuous for~$a\in
A$.

The {\it boundary operator} $\pd:C_k^\rsi(Y;R)\ra C_{k-1}^\rsi(Y;R)$
is \cite[\S IV.1]{Bred}:
\e
\pd:\ts\sum_{a\in A}\rho_a\,\si_a\longmapsto \ts\sum_{a\in
A}\sum_{j=0}^k(-1)^j\rho_a(\si_a\ci F_j^k),
\label{kh4eq2}
\e
where for $j=0,\ldots,k$ the map $F_j^k:\De_{k-1}\ra\De_k$ is given
by $F_j^k(x_0,\ldots,x_{k-1})=(x_0,\ldots,x_{j-1},0,x_j,\ldots,
x_{k-1})$. As a manifold with boundary, we have
$\pd\De_k=\coprod_{j=0}^k\De_j^k$, where $\De_j^k$ is the connected
component of $\pd\De_k$ on which $x_j\equiv 0$, so that
$F_j^k:\De_{k-1}\ra\De_j^k$ is a diffeomorphism. The orientation on
$\De_k$ induces one on $\pd\De_k$, and so induces orientations on
$\De_j^k$ for $j=0,\ldots,k$. It is easy to show that under
$F_j^k:\De_{k-1}\ra\De_j^k$, the orientations of $\De_{k-1}$ and
$\De_j^k$ differ by a factor $(-1)^j$, which is why this appears in
\eq{kh4eq2}. So $\pd$ in \eq{kh4eq2} basically restricts from
$\De_k$ to $\pd\De_k$, as oriented manifolds with boundary and
corners.

From \eq{kh4eq2} we find that $\pd\ci\pd=0$, since each codimension
2 face $\De_{k-2}$ of $\De_k$ contributes twice to $\pd\ci\pd$, once
with sign 1 and once with sign $-1$. Thus we may define the {\it
singular homology group}
\begin{equation*}
H_k^\rsi(Y;R)=\frac{\Ker\bigl(\pd:C_k^\rsi(Y;R)\ra
C_{k-1}^\rsi(Y;R)\bigr)}{\Im\bigl(\pd:C_{k+1}^\rsi(Y;R)\ra
C_k^\rsi(Y;R)\bigr)}\,.
\end{equation*}
If $Y$ is a smooth manifold or orbifold, we can instead define
$C_k^\rsi(Y;R),H_k^\rsi(Y;R)$ using {\it smooth\/} maps
$\si:\De_k\ra Y$, as in \cite[\S V.5]{Bred}, which we call {\it
smooth singular simplices}, and this gives the same homology groups.
We shall always take $C_k^\rsi(Y;R)$ and $H_k^\rsi(Y;R)$ to be
defined using smooth singular simplices.
\label{kh4def1}
\end{dfn}

We will not choose a particular cohomology theory. We will only be
interested in compactly-supported cohomology of manifolds or
orbifolds, and this can be characterized in terms of (singular)
homology by Poincar\'e duality isomorphisms, as in \eq{kh4eq3}
below. So when we prove that (effective) Kuranishi cohomology is
isomorphic to compactly-supported cohomology, we will go via
singular homology using Poincar\'e duality.

Here are some general properties of homology and cohomology, which
can mostly be found in Bredon \cite{Bred}.
\begin{itemize}
\setlength{\itemsep}{0pt}
\setlength{\parsep}{0pt}
\item Let $Y,Z$ be topological spaces, $h:Y\ra Z$ a continuous
map, and $R$ a commutative ring. Then there is a {\it pushforward
map} $h_*:H_*(Y;R)\ra H_*(Z;R)$. On singular homology, this is
induced by $h_*:C_*^\rsi(Y;R)\ra C_*^\rsi(Z;R)$ mapping $h_*:
\sum_{a\in A}\rho_a\,\si_a\mapsto\sum_{a\in A}\rho_a(h\ci\si_a)$.
Pushforwards are functorial, $(g\ci h)_*=g_*\ci h_*$.
\item Let $Y,Z$ be topological spaces, $h:Y\ra Z$ a continuous
map, and $R$ a commutative ring. Then there is a {\it pullback map}
$h^*:H^*(Z;R)\ra H^*(Y;R)$ on cohomology. Pullbacks are
functorial,~$(g\ci h)^*=h^*\ci g^*$.

For {\it compactly-supported\/} cohomology, we can only define
pullbacks $h^*:H^*_\cs(Z;R)\ra H^*_\cs(Y;R)$ if $h:Y\ra Z$ is {\it
proper}, that is, if $h^{-1}(S)\subseteq Y$ is compact in $Y$
whenever $S\subseteq Z$ is compact in $Z$. This is because $h^*$
must pull compactly-supported cochains in $Z$ back to
compactly-supported cochains in~$Y$.
\item There is a natural morphism $H^*_\cs(Y;R)\ra H^*(Y;R)$. If
$Y$ is compact this is an isomorphism, $H^*_\cs(Y;R)\cong H^*(Y;R)$.
\item There are associative, supercommutative, graded multiplications
on both cohomology and compactly-supported cohomology, the {\it cup
product\/} $\cup:H^k(Y;R)\t H^l(Y;R)\ra H^{k+l}(Y;R)$. There is an
{\it identity} $1\in H^0(Y;R)$, and if $Y$ is compact there is an
identity in $H^0_\cs(Y;R)$.

There are also {\it cap products} $\cap:H_k(Y;R)\t H^l(Y;R)\ra
H_{k-l}(Y;R)$ and $\cap:H_k(Y;R)\t H^l_\cs(Y;R)\ra H_{k-l}(Y;R)$,
which make $H_*(Y;R)$ into a {\it module} over $H^*(Y;R)$ and
$H^*_\cs(Y;R)$.
\item Suppose $Y$ is an oriented manifold, of dimension $n$,
without boundary, and not necessarily compact. Then there are {\it
Poincar\'e duality isomorphisms}
\e
\Pd:H^k_\cs(Y;R)\longra H_{n-k}(Y;R)
\label{kh4eq3}
\e
between compactly-supported cohomology, and homology.

If $Y$ is also {\it compact\/} then it has a {\it fundamental class}
$[Y]\in H_n(Y;R)$, and we can write the Poincar\'e duality map $\Pd$
of \eq{kh4eq3} in terms of the cap product by $\Pd(\al)=[Y]\cap\al$
for $\al\in H^k_\cs(Y;R)$. If $Y$ is noncompact and not too badly
behaved, a similar interpretation is possible, but the fundamental
class $[Y]$ exists not in $H_n(Y;R)$ but in the
`non-compactly-supported homology group' $H_n^{\rm nc}(Y;R)$ which
is $H_n(\bar Y,\{\iy\};R)$ in relative homology, where $\bar
Y=Y\cup\{\iy\}$ is the one-point compactification of~$Y$.

Poincar\'e duality for manifolds $Y$ with boundary and (g-)corners
will be discussed in~\S\ref{kh45}.
\item If $Y$ is an oriented $n$-manifold without boundary then
Poincar\'e duality and the cup product on $H^*_\cs(Y;R)$ induce an
associative, supercommutative {\it intersection product\/}
$\bu:H_k(Y;R)\t H_l(Y;R)\ra H_{k+l-n}(Y;R)$ on homology $H_*(Y;R)$.
This can be described geometrically at the chain level in singular
homology, as in Bredon \cite[VI.11]{Bred} and Lefschetz~\cite[\S
IV]{Lefs}.
\item Fulton and MacPherson \cite{FuMa} define the formalism of
{\it bivariant theories}, which simultaneously generalize homology
and cohomology. In \S\ref{kh48} we will explain this, and show how
(effective) Kuranishi (co)homology can be extended to a bivariant
theory, {\it at the level of (co)chains}.
\item Here is an important point: {\it by the homology or cohomology of an
orbifold $Y$, we will always mean the (co)homology of the underlying
topological space}. However, there are other ways to define
(co)homology of an orbifold, which also use the orbifold structure
of $Y$, and even additional data.

Behrend \cite{Behr2} defines versions of de Rham cohomology
$H_{\sst\rm dR}^*(Y;\R)$, singular homology $H_*^{\sst\rm Be}(Y;R)$
and singular cohomology $H^*_{\sst\rm Be}(Y;R)$ of an orbifold,
regarded as a stack. In Behrend's notation, the underlying
topological space $\bar Y$ of an orbifold $Y$ is called the {\it
coarse moduli space} of $Y$. Behrend shows \cite[Prop.~36]{Behr2}
that $H_*^{\sst\rm Be}(Y;R)\cong H_*^\rsi(\bar Y;R)$ and $H^*_{\sst
\rm Be}(Y;R)\cong H^*(\bar Y;R)$ when $R$ is a $\Q$-algebra, but for
general $R$ such as $R=\Z$ the (co)homology of $Y$ and $\bar Y$ may
differ.

Also, if $Y$ is an orbifold with an almost complex structure, Chen
and Ruan \cite{ChRu2} define the {\it orbifold cohomology}
$H^*_{\sst\rm CR}(Y;R)$, which appears to be natural in
Gromov--Witten theory and String Theory of orbifolds. In general
$H^*_{\sst\rm CR}(Y;R)$ is different from both $H^*(\bar Y;R)$
and~$H^*_{\sst\rm Be}(Y;R)$.
\end{itemize}

We shall be interested in how Poincar\'e duality for manifolds
extends to {\it orbifolds} $Y$. Satake \cite[Th.~3]{Sata} showed
that Poincar\'e duality isomorphisms \eq{kh4eq3} exist when $Y$ is
an oriented orbifold and $R$ is a $\Q$-{\it algebra}. However, {\it
Poincar\'e duality need not hold for orbifolds for general\/} $R$.
(Here, as in the last remark above, we are discussing Poincar\'e
duality for the (co)homology {\it of the underlying topological
space}.) Poincar\'e duality can fail at two levels.

Let $Y$ be an oriented $n$-orbifold, and suppose for simplicity that
$Y$ is compact. Firstly, if $Y$ is a {\it non-effective} orbifold,
then one can in general define the fundamental class $[Y]\in
H_n(Y;R)$ only if $R$ is a $\Q$-algebra. This is because generic
points in $Y$ have non-trivial stabilizer group $\Ga$, so working
say in singular homology, when we triangulate $Y$ by simplices to
define a cycle representing $[Y]$ we must weight each simplex by
$\pm\md{\Ga}^{-1}$, so we need $\md{\Ga}^{-1}\in R$, and
$\Q\subseteq R$. Since $\Pd$ in \eq{kh4eq3} is defined using $[Y]$,
for non-effective $Y$ and general $R$ the Poincar\'e duality map
$\Pd$ is undefined.

Secondly, if $Y$ is an {\it effective} orbifold, we can define $[Y]$
and $\Pd$ in \eq{kh4eq3}, but $\Pd$ may not be an isomorphism. Here
is an example of a compact, oriented, effective orbifold for which
$\Pd$ is not an isomorphism over~$\Z$.

\begin{ex} Let $Y=(\CP^1\t\CP^1)/\Z_2$, where the $\Z_2$-action is
generated by $\si:\CP^1\t\CP^1\ra\CP^1\t\CP^1$ mapping $\si:\bigl(
[x_0,x_1],[y_0,y_1]\bigr)\mapsto\bigl([x_0,-x_1],[y_0,-y_1]\bigr)$.
It is a compact 4-orbifold, with an orientation induced by the
complex structure, that has four orbifold points $\{[1,0],[0,1]\}\t
\{[1,0],[0,1]\}$ modelled on~$\R^4/\{\pm 1\}$.

Let $\al,\be$ in $H_2(Y;\Z)\subset H_2(Y;\Q)$ be the homology
classes of the suborbifolds $\bigl(\CP^1\t\{[1,0]\}\bigr)/\Z_2$ and
$\bigl(\{[1,0]\}\t\CP^1\bigr)/\Z_2$. Since Poincar\'e duality holds
over $\Q$, there is an intersection product $\bu:H_2(Y;\Q)\t
H_2(Y;\Q)\ra H_0(Y;\Q)\cong\Q$. As $\bigl(\CP^1\t\{[1,0]\}\bigr)
/\Z_2$, $\bigl(\{[1,0]\}\t \CP^1\bigr)/\Z_2$ intersect transversely
in one point $\bigl([1,0],[1,0]\bigr)$, which has stabilizer group
$\Z_2$, it follows that~$\al\bu\be=\ha\in\Q$.

If Poincar\'e duality \eq{kh4eq3} were an isomorphism for $R=\Z$
then there would exist an intersection product $\bu:H_2(Y;\Z)\t
H_2(Y;\Z)\ra H_0(Y;\Z)\cong\Z$ agreeing with that on $H_*(Y;\Q)$
under the morphism $H_*(Y;\Z)\ra H_*(Y;\Q)$ induced by
$\Z\hookra\Q$. However, this contradicts $\al,\be\in H_2(Y;\Z)$ but
$\al\bu\be=\ha\notin H_0(Y;\Z)\subset H_0(Y;\Q)$. This means that
$\al,\be$ do not lie in the image of $\Pd$ (though $2\al,2\be$ do).
Thus $\Pd$ in \eq{kh4eq3} is not surjective, and so not an
isomorphism.
\label{kh4ex1}
\end{ex}

\subsection{(Effective) Kuranishi homology}
\label{kh42}

We now define the (effective) Kuranishi homology groups of an
orbifold.

\begin{dfn} Let $Y$ be an orbifold. Consider triples $(X,\bs f,\bs
G)$, where $X$ is a compact oriented Kuranishi space, $\bs f:X\ra Y$
is strongly smooth, and $\bs G$ is gauge-fixing data for $(X,\bs
f)$. Write $[X,\bs f,\bs G]$ for the isomorphism class of $(X,\bs
f,\bs G)$ under isomorphisms $(\bs a,\bs b):(X,\bs f,\bs G)\ra(\ti
X,\bs{\ti f},\bs{\ti G})$ as in Definition \ref{kh3def6}, where in
addition we require $\bs a$ to identify the orientations of~$X,\ti
X$.

Let $R$ be a $\Q$-algebra, for instance $\Q,\R$ or $\C$. For each
$k\in\Z$, define $KC_k(Y;R)$ to be the $R$-module of finite
$R$-linear combinations of isomorphism classes $[X,\bs f,\bs G]$ for
which $\vdim X=k$, with the relations:
\begin{itemize}
\setlength{\itemsep}{0pt}
\setlength{\parsep}{0pt}
\item[(i)] Let $[X,\bs f,\bs G]$ be an isomorphism class,
and write $-X$ for $X$ with the opposite orientation. Then in
$KC_k(Y;R)$ we have
\e
[X,\bs f,\bs G]+[-X,\bs f,\bs G]=0.
\label{kh4eq4}
\e
\item[(ii)] Let $[X,\bs f,\bs G]$ be an isomorphism class, and
suppose there exists an isomorphism $(\bs a,\bs b):(X,\bs f,\bs
G)\ra(X,\bs f,\bs G)$ in the sense of Definition \ref{kh3def6}, such
that $\bs a$ reverses the orientation of $X$. Then
\e
[X,\bs f,\bs G]=0 \quad\text{in $KC_k(Y;R)$.}
\label{kh4eq5}
\e
\item[(iii)] Let $[X,\bs f,\bs G]$ be an isomorphism class, and
let $X=X_1\amalg X_2\amalg\cdots\amalg X_n$ be a {\it splitting} of
$(X,\bs f,\bs G)$, in the sense of Definition \ref{kh3def6}. Then
\e
[X,\bs f,\bs G]=\ts\sum_{a=1}^n[X_a,\bs f\vert_{X_a},\bs
G\vert_{X_a}]\quad\text{in $KC_k(Y;R)$.}
\label{kh4eq6}
\e
\item[(iv)] Let $[X,\bs f,\bs G]$ be an isomorphism class, $\Ga$ a
finite group, and $\rho$ an action of $\Ga$ on $(X,\bs f,\bs G)$ by
orientation-preserving automorphisms. That is,
$\rho:\Ga\ra\Aut(X,\bs f,\bs G)$ is a group morphism, and if
$\rho(\ga)=(\bs a,\bs b)$ for $\ga\in\Ga$ then $(\bs a,\bs b):(X,\bs
f,\bs G)\ra(X,\bs f,\bs G)$ is an isomorphism with $\bs a:X\ra X$
orientation-preserving. Note that we do not require $\rho$ to be
injective, so we cannot regard $\Ga$ as a subgroup of $\Aut(X,\bs
f,\bs G)$. Then $\Ga$ acts on $X$, and $\ti X=X/\Ga$ is a compact
oriented Kuranishi space, with a projection $\bs\pi:X\ra\ti X$. As
in Property \ref{kh3pr}(c), Definition \ref{kh3def9} gives strongly
smooth $\bs\pi_*(\bs f)=\bs{\ti f}:\ti X\ra Y$ and gauge-fixing data
$\bs\pi_*(\bs G)=\bs{\ti G}$ for $(\ti X,\bs{\ti f})$. Thus
$[X/\Ga,\bs\pi_*(\bs f),\bs\pi_*(\bs G)]$ is an isomorphism class.
We impose the relation
\e
\bigl[X/\Ga,\bs\pi_*(\bs f),\bs\pi_*(\bs G)\bigr]=\frac{1}{\md{\Ga}}
\,\bigl[X,\bs f,\bs G\bigr]\qquad\text{in $KC_k(Y;R)$.}
\label{kh4eq7}
\e
\end{itemize}
Elements of $KC_k(Y;R)$ will be called {\it Kuranishi chains}. We
require $R$ to be a $\Q$-{\it algebra\/} so that the factor
$\frac{1}{\md{\Ga}}$ in \eq{kh4eq7} makes sense.

Define the {\it boundary operator\/} $\pd:KC_k(Y;R)\ra
KC_{k-1}(Y;R)$ by
\e
\pd:\ts\sum_{a\in A}\rho_a[X_a,\bs f_a,\bs G_a]\longmapsto
\ts\sum_{a\in A}\rho_a[\pd X_a,\bs f_a\vert_{\pd X_a},\bs G_a
\vert_{\pd X_a}],
\label{kh4eq8}
\e
where $A$ is a finite indexing set and $\rho_a\in R$ for $a\in A$.
This is a morphism of $R$-modules. Clearly, $\pd$ takes each
relation \eq{kh4eq4}--\eq{kh4eq7} in $KC_k(Y;R)$ to the
corresponding relation in $KC_{k-1}(Y;R)$, and so $\pd$ is
well-defined.

Recall from Definition \ref{kh2def19} that if $X$ is an oriented
Kuranishi space then there is a natural orientation-reversing strong
diffeomorphism $\bs\si:\pd^2X\ra\pd^2X$, with $\bs\si^2=
\bs\id_{\pd^2X}$. If $[X,\bs f,\bs G]$ is an isomorphism class then
this $\bs\si$ extends to an isomorphism $(\bs\si,\bs\tau)$ of
$(\pd^2X,\bs f\vert_{\pd^2X},\bs G\vert_{\pd^2X})$, where $J$ is the
indexing set of $\bs I\vert_{\pd^2X}$, and for each $j\in J$,
$(\tau^j,\hat\tau^j)$ is the automorphism of
$(\pd^2V^{j+2},\ldots,\psi^{j+2}\vert_{\pd^2V^{j+2}})$ induced by
the orientation-reversing involution
$\si:\pd^2V^{j+2}\ra\pd^2V^{j+2}$ described in Definitions
\ref{kh2def5} and \ref{kh2def6}. So part (ii) in $KC_{k-2}(Y;R)$
yields
\begin{equation*}
[\pd^2X,\bs f\vert_{\pd^2X},\bs G\vert_{\pd^2X}]=0 \quad\text{in
$KC_{k-2}(Y;R)$.}
\end{equation*}
Clearly this implies that $\pd\ci\pd=0$ as a map~$KC_k(Y;R)\ra
KC_{k-2}(Y;R)$.

Define the {\it Kuranishi homology group\/} $KH_k(Y;R)$ of $Y$ for
$k\in\Z$ to be
\begin{equation*}
KH_k(Y;R)=\frac{\Ker\bigl(\pd:KC_k(Y;R)\ra KC_{k-1}(Y;R)\bigr)}{
\Im\bigl(\pd:KC_{k+1}(Y;R)\ra KC_k(Y;R)\bigr)}\,.
\end{equation*}
This is a well-defined $R$-module, as~$\pd\ci\pd=0$.
\label{kh4def2}
\end{dfn}

Note that if 2 is invertible in $R$, as it is when $R$ is a
$\Q$-algebra, then Definition \ref{kh4def2}(ii) follows from
Definition \ref{kh4def2}(i), since $(\bs a,\bs b)$ in (ii) implies
the isomorphism classes $[X,\bs f,\bs G]$ and $[-X,\bs f,\bs G]$ are
equal, so \eq{kh4eq4} gives $2[X,\bs f,\bs G]=0$, and multiplying by
$\frac{1}{2}$ gives \eq{kh4eq5}. But we state it as a separate axiom
as we shall re-use it in our next definition, in which 2 need not be
invertible in~$R$.

\begin{dfn} Let $Y$ be an orbifold. Consider triples $(X,\bs
f,\ubG)$, where $X$ is a compact, oriented, {\it effective\/}
Kuranishi space, $\bs f:X\ra Y$ is strongly smooth, and $\ubG$ is
{\it effective\/} gauge-fixing data for $(X,\bs f)$. Write $[X,\bs
f,\ubG]$ for the isomorphism class of $(X,\bs f,\ubG)$ under
isomorphisms $(\bs a,\bs b):(X,\bs f,\ubG)\ra(\ti X,\bs{\ti
f},\ubtG)$, where $\bs a$ identifies the orientations of~$X,\ti X$.

Let $R$ be a commutative ring, for instance $R=\Z,\Z_p,\Q,\R$ or
$\C$. For each $k\in\Z$, define $KC_k^\ef(Y;R)$ to be the $R$-module
of finite $R$-linear combinations of isomorphism classes $[X,\bs f,
\ubG]$ as above for which $\vdim X=k$, with the analogues of
relations Definition \ref{kh4def2}(i)--(iii), but {\it not\/}
Definition \ref{kh4def2}(iv). Elements of $KC_k^\ef(Y;R)$ will be
called {\it effective Kuranishi chains}.

Define the {\it boundary operator\/} $\pd:KC_k^\ef(Y;R)\ra
KC_{k-1}^\ef(Y;R)$ by
\begin{equation*}
\pd:\ts\sum_{a\in A}\rho_a[X_a,\bs f_a,\ubG_a]\longmapsto
\ts\sum_{a\in A}\rho_a[\pd X_a,\bs f_a\vert_{\pd X_a},\ubG_a
\vert_{\pd X_a}],
\end{equation*}
as in \eq{kh4eq8}. Then $\pd\ci\pd=0$. Define the {\it effective
Kuranishi homology group\/} $KH_k^\ef(Y;R)$ of $Y$ for $k\in\Z$ to
be
\begin{equation*}
KH_k^\ef(Y;R)=\frac{\Ker\bigl(\pd:KC_k^\ef(Y;R)\ra
KC_{k-1}^\ef(Y;R)\bigr)}{\Im\bigl(\pd:KC_{k+1}^\ef(Y;R)\ra
KC_k^\ef(Y;R)\bigr)}\,.
\end{equation*}
\label{kh4def3}
\end{dfn}

\begin{rem}{\bf(a)} We define $KC_k,KH_k,KC_k^\ef,KH_k^\ef(Y;R)$
for all $k\in\Z$, not just for $k\ge 0$. This is because nontrivial
Kuranishi spaces $X$ can have $\vdim X=k$ for any $k\in\Z$. Theorem
\ref{kh4thm1} and Corollary \ref{kh4cor1} will show that
$KH_k(Y;R)=KH_k^\ef(Y;R)=0$ when $k<0$, so we may as well restrict
to $k\ge 0$ when considering homology. However, at the chain level
we need to allow $k<0$, as $KC_k,KC_k^\ef(Y;R)$ can be nonzero for
any $k\in\Z$. To define $KH_0,KH_0^\ef(Y;R)$ we use
$\pd:KC_0(Y;R)\ra KC_{-1}(Y;R)$ and $\pd:KC_0^\ef(Y;R)\ra
KC_{-1}^\ef(Y;R)$, so we need to consider $-1$-chains. Also, the
inductive proofs in Appendices \ref{khB} and \ref{khC} can
implicitly involve Kuranishi chains of arbitrarily negative degree.

\noindent{\bf(b)} The relation Definition \ref{kh4def2}(iv) in
$KC_*(Y;R)$ is not necessary to get a well-behaved homology theory,
and we omit it in the definition of $KC_*^\ef(Y;R)$. We include it
because it will be useful in some of our applications, in
particular, the open Gromov--Witten invariants to be discussed in
\S\ref{kh67} and \cite{Joyc2}, for which we need it to prove
equation \eq{kh6eq28} of Corollary \ref{kh6cor1}. The author
believes that Theorem \ref{kh4thm2} below would still hold if we
defined $KC_*(Y;R)$ omitting Definition \ref{kh4def2}(iv), and in
Remark \ref{khCrem1}(a) we sketch how to modify the proof in this
case.

Note that if we omit Definition \ref{kh4def2}(iv) then the
definitions of $KC_*(Y;R)$ and $KH_*(Y;R)$ make sense for $R$ an
arbitrary commutative ring, rather than a $\Q$-algebra, since we no
longer need $R$ to contain the factors $\md{\Ga}^{-1}$ in
\eq{kh4eq7}. However, as discussed in Remark \ref{khCrem1}(b), we
would still need $R$ a $\Q$-algebra to prove Theorem \ref{kh4thm2}.
In fact the author has an outline proof that if we defined
$KC_*(Y;R),KH_*(Y;R)$ omitting Definition \ref{kh4def2}(iv) for $R$
a general commutative ring, then~$KH_*(Y;R)\cong
H_*^\rsi(Y;R\ot_\Z\Q)$.

\noindent{\bf(c)} One could probably weaken the condition that $X$
is {\it effective} in effective Kuranishi chains $[X,\bs f,\ubG]$,
and still get a homology theory with $KH_*^\ef(Y;R)\cong
H_*^\rsi(Y;R)$ for all commutative rings $R$, including $R=\Z$. In
\S\ref{kh56} we will define the {\it orbifold strata} $X^{\Ga,\rho}$
of a Kuranishi space $X$, and in Remark \ref{kh5rem5} we will show
that if $X$ is an oriented, {\it effective\/} Kuranishi space then
$\vdim X^{\Ga,\rho}\le\vdim X-2$ for all nonempty orbifold strata
$X^{\Ga,\rho}$ with~$\Ga\ne\{1\}$.

For a compact, oriented Kuranishi space $X$ without boundary, the
condition that $\vdim X^{\Ga,\rho}\le\vdim X-2$ for all
$X^{\Ga,\rho}\ne\es$ with $\Ga\ne\{1\}$ is enough to ensure that the
virtual class of $X$ can be defined in homology over $\Z$, not just
over $\Q$. This holds because if we take a generic (single-valued)
perturbation $\ti X$ of such a Kuranishi space $X$, although $\ti X$
may not be a manifold or orbifold, the singularities of $\ti X$ due
to the orbifold strata $X^{\Ga,\rho}$ occur in codimension at least
2, so $\ti X$ is a {\it pseudo-manifold}, and thus has an integral
fundamental class. This is the important fact underlying the
definition of Gromov--Witten invariants in {\it integral\/} homology
for {\it semi-positive} symplectic manifolds, as discussed briefly
in \S\ref{kh63}, and proved in McDuff and Salamon~\cite[\S 7]{McSa}.

It seems likely that one could define some modification of
$KH_*^\ef(Y;R)$ isomorphic to $H_*^\rsi(Y;R)$ for all commutative
rings $R$, with chains $[X,\bs f,\ubG]$ in which we require $\vdim
X^{\Ga,\rho}\le\vdim X-2$ for all $X^{\Ga,\rho}\ne\es$ with
$\Ga\ne\{1\}$, instead of requiring $X$ to be effective. This
modification might then have interesting applications, for instance
in defining Lagrangian Floer cohomology over $\Z$ for semi-positive
symplectic manifolds. But we will not attempt this here.
\label{kh4rem1}
\end{rem}

\subsection{Morphisms between singular and Kuranishi homology}
\label{kh43}

There is a natural projection $\Pi_\ef^\Kh$ from effective Kuranishi
chains and homology to Kuranishi chains and homology. Also,
(effective) Kuranishi chains and homology have functorial
pushforwards $h_*$ by smooth maps~$h:Y\ra Z$.

\begin{dfn} Let $Y$ be an orbifold, and $R$ a commutative ring. Define
\e
\begin{split}
\Pi_\ef^\Kh:KC_k^\ef(Y;R)&\longra KC_k(Y;R\ot_\Z\Q)\quad\text{by}\\
\Pi_\ef^\Kh:\ts\sum_{a\in A}\rho_a\bigl[X_a,\bs f_a,\ubG_a\bigr]
&\longmapsto \ts\sum_{a\in A}\pi(\rho_a)\bigl[X_a,\bs
f_a,\Pi(\ubG_a)\bigr],
\end{split}
\label{kh4eq9}
\e
for $k\in\Z$, where $\pi:R\ra R\ot_\Z\Q$ is the natural ring
morphism, and $\Pi(\ubG_a)$ is as in Definition \ref{kh3def17}. If
$R$ is a $\Q$-algebra then $R\ot_\Z\Q=R$ and $\pi$ is the identity.
The $\Pi_\ef^\Kh$ map relations in $KC_k^\ef(Y;R)$ to relations in
$KC_k(Y;R\ot_\Z\Q)$, and so are well-defined. They satisfy
$\Pi_\ef^\Kh\ci\pd=\pd\ci\Pi_\ef^\Kh$, so they induce morphisms
\e
\Pi_\ef^\Kh:KH_k^\ef(Y;R)\longra KH_k(Y;R\ot_\Z\Q).
\label{kh4eq10}
\e

Let $Y,Z$ be orbifolds, and $h:Y\ra Z$ a smooth map. Using the
notation of Definitions \ref{kh3def12} and \ref{kh3def20}, define
the {\it pushforwards}
\e
\begin{aligned}
h_*:KC_k(Y;R)\longra KC_k(Z;R),&\quad
h_*:KC_k^\ef(Y;R)\longra KC_k^\ef(Z;R)\\
\text{by}\quad h_*:\ts\sum_{a\in A}\rho_a\bigl[X_a,\bs f_a,\bs
G_a\bigr]&\longmapsto \ts\sum_{a\in A}\rho_a\bigl[X_a,h\ci\bs
f_a,h_*(\bs G_a)\bigr]\\
\text{and}\quad h_*:\ts\sum_{a\in A}\rho_a\bigl[X_a,\bs f_a,\ubG_a
\bigr] &\longmapsto \ts\sum_{a\in A}\rho_a\bigl[X_a,h\ci\bs
f_a,h_*(\ubG_a)\bigr].
\end{aligned}
\label{kh4eq11}
\e
These take relations (i)--(iv) or (i)--(iii) in $KC_k,KC_k^\ef(Y;R)$
to relations (i)--(iv) or (i)--(iii) in $KC_k,KC_k^\ef(Z;R)$, and so
are well-defined. They satisfy $h_*\ci\pd=\pd\ci h_*$, so they
induce morphisms of homology groups
\begin{equation*}
h_*:KH_k(Y;R)\longra KH_k(Z;R)\;\text{and}\;\>
h_*:KH_k^\ef(Y;R)\longra KH_k^\ef(Z;R).
\end{equation*}
Pushforward is functorial, that is, $(g\ci h)_*=g_*\ci h_*$, on
chains and homology.
\label{kh4def4}
\end{dfn}

\begin{dfn} Let $Y$ be an orbifold, $\De_k$ be the $k$-simplex of
\eq{kh4eq1}, and $\si:\De_k\ra Y$ be a smooth map. Then $\De_k$ is a
compact oriented manifold with corners, so we may regard it as an
oriented Kuranishi space, and $\si$ as a strongly smooth map.

Define an excellent coordinate system $(\bs I_{\De_k},
\bs\eta_{\De_k})$ for $(\De_k,\si)$ to have indexing set
$I_{\De_k}=\{k\}$, Kuranishi neighbourhood
$(V^k_{\De_k},E^k_{\De_k},s^k_{\De_k},\psi^k_{\De_k})
=(\De_k,\ab\De_k,\ab 0,\ab\id_{\De_k})$, and map
$\si^k:V^k_{\De_k}\ra Y$ representing $\si$ given by $\si^k=\si$.
Here $E^k_{\De_k}\ra V^k_{\De_k}$ is the zero vector bundle, so that
as manifolds $E^k_{\De_k}=V^k_{\De_k}=\De_k$. Set
$\bs\eta_{\De_k}=(\eta_{k,\De_k},\eta_{k,\De_k}^k)$, where
$\eta_{k,\De_k}\equiv 1$ on $\De_k$ and $\eta_{k,\De_k}^k\equiv 1$
on $V^k_{\De_k}$. Define $G^k_{\De_k}:E^k_{\De_k}\ra P$ and
$\uG^k_{\De_k}:E^k_{\De_k}\ra\uP$ by
\e
\begin{split}
G^k_{\De_k}\bigl((x_0,\ldots,x_k)\bigr)&=S_{l+1}\cdot
(y_0,y_0+y_1,\ldots,y_0+y_1+\cdots+y_l),\\
\uG^k_{\De_k}\bigl((x_0,\ldots,x_k)\bigr)&=
(y_0,y_0+y_1,\ldots,y_0+y_1+\cdots+y_l),
\end{split}
\label{kh4eq12}
\e
where $y_0,y_1,\ldots,y_l$ are those $x_0,x_1,\ldots,x_k$ which are
nonzero, in the same order, so that $(x_0,\ldots,x_k)$ is
$(y_0,y_1,\ldots,y_l)$ with $k-l$ zeros inserted in the list. Then
\e
\text{$G^k_{\De_k}\ci F_j^k\equiv G^{k-1}_{\De_{k-1}}$,\quad
$\uG^k_{\De_k}\ci F_j^k\equiv\uG^{k-1}_{\De_{k-1}}$ for
$j=0,\ldots,k$,}
\label{kh4eq13}
\e
with $F_j^k$ as in \S\ref{kh41}. Also, for
$(x_0,\ldots,x_k)\in\De_k^\ci$ we have $x_j>0$ for all $j$, so that
$(y_0,\ldots,y_l)=(x_0,\ldots,x_k)$, and
$G^k_{\De_k}:(x_0,\ldots,x_k)\mapsto S_{k+1}\cdot(x_0,\ldots,
x_0+\cdots+x_k)$, and $\uG^k_{\De_k}:(x_0,\ldots,x_k)\mapsto
(x_0,\ldots,x_0+\cdots+x_k)$. Since
$x_0,x_0+x_1,\ldots,x_0+\cdots+x_k$ are arranged in strictly
increasing order, $S_{k+1}\cdot(x_0,\ldots,x_0+\cdots+x_k)$
determines $x_0,x_0+x_1,\ldots,x_0+\cdots+x_k$ and hence
$x_0,\ldots,x_k$. Thus $G^k_{\De_k}\vert_{\De_k^\ci}$ is injective.
Similarly, $\uG^k_{\De_k}\vert_{\De_k^\ci}$ maps
$\De_k^\ci\ra\R^{k+1}\subset\uP$, and is injective. It is now easy
to check that Definitions \ref{kh3def6} and \ref{kh3def17} hold, so
that $\bs G_{\De_k}$ is {\it gauge-fixing data}, and $\ubG_{\De_k}$
is {\it effective gauge-fixing data}, for~$(\De_k,\si)$.

Taking $R$ to be a $\Q$-algebra in \eq{kh4eq14}, define maps
\ea
\Pi_\rsi^\Kh&:C_k^\rsi(Y;R)\ra KC_k(Y;R),\;\>
\Pi_\rsi^\Kh:\ts\sum\limits_{a\in A}\rho_a\si_a\mapsto
\ts\sum\limits_{a\in A}\rho_a\bigl[\De_k,\si_a,\bs G_{\De_k}\bigr],
\label{kh4eq14}\\
\Pi_\rsi^\ef&:C_k^\rsi(Y;R)\ra KC_k^\ef(Y;R),\;\>
\Pi_\rsi^\ef:\ts\sum\limits_{a\in A}\rho_a\si_a\mapsto
\ts\sum\limits_{a\in A}\rho_a\bigl[\De_k,\si_a,\ubG_{\De_k}\bigr].
\label{kh4eq15}
\ea
As $\pd\De_k\!=\!\coprod_{j=1}^k(-1)^j F_j^k(\De_{k-1})$ in oriented
$(k-1)$-manifolds, \eq{kh4eq13} gives
\e
\pd\bigl[\De_k,\si_a,\bs G_{\De_k}\bigr]=\ts\sum_{j=0}^k(-1)^j
\bigl[\De_{k-1},(\si_a\ci F_j^k),\bs G_{\De_{k-1}}\bigr]
\label{kh4eq16}
\e
in $KC_{k-1}(Y;R)$. Comparing \eq{kh4eq2}, \eq{kh4eq8} and
\eq{kh4eq16} shows that $\pd\ci\Pi_\rsi^\Kh=\Pi_\rsi^\Kh\ci\pd$.
Similarly $\pd\ci\Pi_\rsi^\ef=\Pi_\rsi^\ef\ci\pd$. Therefore
\eq{kh4eq14}--\eq{kh4eq15} induce projections
\ea
\Pi_\rsi^\Kh&:H_k^\rsi(Y;R)\longra KH_k(Y;R),
\label{kh4eq17}\\
\Pi_\rsi^\ef&:H_k^\rsi(Y;R)\longra KH_k^\ef(Y;R),
\label{kh4eq18}
\ea
requiring $R$ to be a $\Q$-algebra in \eq{kh4eq17}. The morphisms
\eq{kh4eq9}, \eq{kh4eq14}, \eq{kh4eq15} on chains, and \eq{kh4eq10},
\eq{kh4eq17}, \eq{kh4eq18} on homology, satisfy
\e
\Pi_\ef^\Kh\ci\Pi_\rsi^\ef=\Pi_\rsi^\Kh.
\label{kh4eq19}
\e
\label{kh4def5}
\end{dfn}

Our next two theorems, the most important of the book, will be
proved in Appendices \ref{khB} and \ref{khC} respectively, after
preparatory work in Appendix~\ref{khA}.

\begin{thm} Let\/ $Y$ be an orbifold and\/ $R$ a commutative ring.
Then the morphism $\Pi_\rsi^\ef:H_k^\rsi(Y;R)\ra KH_k^\ef(Y;R)$ in
\eq{kh4eq18} is an isomorphism for $k\ge 0,$ with\/
$KH_k^\ef(Y;R)=0$ when~$k<0$.
\label{kh4thm1}
\end{thm}

\begin{thm} Let\/ $Y$ be an orbifold and\/ $R$ a $\Q$-algebra. Then
the morphism $\Pi_\ef^\Kh:KH_k^\ef(Y;R)\ra KH_k(Y;R)$ in
\eq{kh4eq10} is an isomorphism for~$k\in\Z$.
\label{kh4thm2}
\end{thm}

These and equation \eq{kh4eq19} imply:

\begin{cor} Let\/ $Y$ be an orbifold and\/ $R$ a $\Q$-algebra. Then
the morphism $\Pi_\rsi^\Kh:H_k^\rsi(Y;R)\ra KH_k(Y;R)$ in
\eq{kh4eq17} is an isomorphism for $k\ge 0,$ with\/ $KH_k(Y;R)=0$
when~$k<0$.
\label{kh4cor1}
\end{cor}

\begin{rem}{\bf(a)} Theorems \ref{kh4thm1} and \ref{kh4thm2} are
our substitute for results of Fukaya and Ono on {\it
multisections\/} and {\it virtual cycles\/} \cite[\S 3, \S
6]{FuOn1}, \cite[\S A1]{FOOO}. If $X$ is a compact Kuranishi space,
$Y$ an orbifold, and $\bs f:X\ra Y$ is strongly smooth, then by
choosing an abstract, multivalued perturbation of $X$ called a
multisection which locally turns $X$ into a (non-Hausdorff) smooth
manifold $X^{\rm per}$, with weights in $\Q$, and then triangulating
$X^{\rm per}$ by simplices, Fukaya and Ono produce a singular chain
$C(X,\bs f)\in C_{\vdim X}^\rsi(Y;\Q)$ called the {\it virtual
chain}, depending on many choices. If $\pd X=\es$ then $\pd C(X,\bs
f)=0$, and $\bigl[C(X,\bs f)\bigr]\in H_{\vdim X}^\rsi(Y;\Q)$ is
independent of choices.

The relation of this to Theorems \ref{kh4thm1} and \ref{kh4thm2} and
Corollary \ref{kh4cor1} is that if $\pd X=\es$ then choosing any
gauge-fixing data $\bs G$ for $(X,\bs f)$, we have a class
$\bigl[[X,\bs f,\bs G]\bigr]\in KH_{\vdim X}(Y;\Q)$, and
$\Pi_\rsi^\Kh:\bigl[C(X,\bs f)\bigr]\mapsto \bigl[[X,\bs f,\bs
G]\bigr]$. That is, in proving $\Pi_\rsi^\Kh$ is invertible we
construct an inverse $(\Pi_\rsi^\Kh)^{-1}:KH_*(Y;\Q) \ra
H_*^\rsi(Y;\Q)$. Very roughly speaking, at the chain level
$(\Pi_\rsi^\Kh)^{-1}$ takes $[X,\bs f,\bs G]$ to $C(X,\bs f)$,
Fukaya and Ono's virtual chain for~$(X,\bs f)$.

In proving Theorems \ref{kh4thm1} and \ref{kh4thm2} in Appendices
\ref{khB} and \ref{khC}, we encounter most of the problems Fukaya
and Ono were tackling with their virtual cycle technology, and we
freely borrow their ideas for our proof, in particular, good
coordinate systems, and inductive choices of small perturbations.
But we do not use multisections: instead we use Definition
\ref{kh4def2}(iv) to represent each Kuranishi homology class by a
sum of chains $[X,\bs f,\bs G]$ in which $X$ has trivial stabilizer
groups, so that we can lift to effective Kuranishi homology, and
then we can perturb each $X$ to a manifold using single-valued
perturbations.

\noindent{\bf(b)} The proofs of Theorems \ref{kh4thm1} and
\ref{kh4thm2} use almost no properties of the target space $Y$. In
the proofs, given some $[X,\bs f,\bs G]$ or $[X,\bs f,\ubG]$, we
have a cover of $X$ by Kuranishi neighbourhoods
$(V^i,E^i,s^i,\psi^i)$ with maps $f^i:V^i\ra Y$ and $G^i:E^i\ra P$
or $\uG^i:E^i\ra\uP$. All we require of these maps $f^i$ is that
they satisfy $f^j\vert_{V^{ij}}\equiv f^i\ci\phi^{ij}$ when $j\le i$
in $I$. The proofs involve cutting $V^i$ into pieces, changing
$G^i,\uG^i$, and deforming $s^i$, which changes $X$ as
$X\cong(s^i)^{-1}(0)$, but we never change the maps~$f^i$.

Because of this, Theorems \ref{kh4thm1} and \ref{kh4thm2} and their
proofs would still hold, essentially unchanged, if we took $Y$ to be
an {\it arbitrary topological space}, rather than an orbifold, and
maps $\bs f:X\ra Y$ in chains $[X,\bs f,\bs G]$ or $[X,\bs f,\ubG]$
to be {\it strongly continuous}, rather than strongly smooth. We
will need $Y$ an orbifold and $\bs f$ strongly smooth or a strong
submersion in \S\ref{kh47} to define {\it products\/}
$\cup,\cap,\bu$ on Kuranishi (co)homology, but to establish
isomorphism with singular homology, this is unnecessary.

In the same way, Theorems \ref{kh4thm1} and \ref{kh4thm2} and their
proofs will still hold, essentially unchanged, if we take $Y$ to be
an {\it infinite-dimensional\/} manifold or orbifold, such as a {\it
loop space} ${\cal L}M$ of smooth maps $\la:{\cal S}^1\ra M$ for $M$
a finite-dimensional manifold, and $\bs f:X\ra Y$ to be strongly
smooth in the appropriate sense. (That is, we prove that
$KH_*,KH_*^\ef(Y;R)$ are isomorphic to singular homology
$H_*^\rsi(Y;R)$ defined using smooth singular simplices in $Y$.)
This will be important in \cite{Joyc3}, when we apply Kuranishi
(co)homology to the {\it String Topology} programme of Chas and
Sullivan~\cite{ChSu}.

We will define a homology theory which computes $H_*^\rsi({\cal
L}M,M;\Q)$, the singular homology of ${\cal L}M$ relative to
constant loops $M\subset{\cal L}M$, and on which Chas and Sullivan's
String Topology operations are defined {\it at the chain level},
without any transversality assumptions, and satisfy the desired
identities {\it exactly} at the chain level, not just up to
homotopy. We will use this to prove results sketched by Fukaya
\cite[\S 6]{Fuka2}, which use moduli spaces of $J$-holomorphic discs
in a symplectic manifold $(M,\om)$ with boundary in a Lagrangian $L$
to define chains in our homology theory of the loop space of $L$
satisfying identities involving String Topology operations.
\label{kh4rem2}
\end{rem}

\subsection{Kuranishi cohomology and Poincar\'e duality for $\pd Y=\es$}
\label{kh44}

Here is our dual notion of {\it Kuranishi cohomology} of an orbifold
$Y$. The differences with Kuranishi homology are that we replace
strongly smooth maps by strong submersions, and orientations by
coorientations, and gauge-fixing data by co-gauge-fixing data, and
we grade $[X,\bs f,\bs C]$ by codimension $\dim Y-\vdim X$ rather
than by dimension $\vdim X$. In this section we restrict to $Y$ {\it
without boundary}. The extension to $\pd Y\ne\es$ will be explained
in~\S\ref{kh45}.

\begin{dfn} Let $Y$ be an orbifold without boundary. Consider
triples $(X,\bs f,\bs C)$, where $X$ is a compact Kuranishi space,
$\bs f:X\ra Y$ is a strong submersion with $(X,\bs f)$ cooriented,
and $\bs C$ is co-gauge-fixing data for $(X,\bs f)$. Write $[X,\bs
f,\bs C]$ for the isomorphism class of $(X,\bs f,\bs C)$ under
isomorphisms $(\bs a,\bs b):(X,\bs f,\bs C)\ra(\ti X,\bs{\ti
f},\bs{\ti C})$ as in Definition \ref{kh3def7}, where $\bs a$ must
identify the coorientations of~$(X,\bs f),(\ti X,\bs{\ti f})$.

Let $R$ be a $\Q$-algebra. For $k\in\Z$, define $KC^k(Y;R)$ to be
the $R$-module of finite $R$-linear combinations of isomorphism
classes $[X,\bs f,\bs C]$ for which $\vdim X=\dim Y-k$, with the
analogues of relations Definition \ref{kh4def2}(i)--(iv), replacing
gauge-fixing data $\bs G$ by co-gauge-fixing data $\bs C$. Elements
of $KC^k(Y;R)$ are called {\it Kuranishi cochains}. Define
$\d:KC^k(Y;R)\ra KC^{k+1}(Y;R)$ by
\e
\d:\ts\sum_{a\in A}\rho_a[X_a,\bs f_a,\bs C_a]\longmapsto
\ts\sum_{a\in A}\rho_a[\pd X_a,\bs f_a\vert_{\pd X_a},\bs C_a
\vert_{\pd X_a}].
\label{kh4eq20}
\e
As in Definition \ref{kh4def2} we have $\d\ci\d=0$. Define the {\it
Kuranishi cohomology groups} $KH^k(Y;R)$ of $Y$ for $k\in\Z$ to be
\begin{equation*}
KH^k(Y;R)=\frac{\Ker\bigl(\d:KC^k(Y;R)\ra KC^{k+1}(Y;R)\bigr)}{
\Im\bigl(\d:KC^{k-1}(Y;R)\ra KC^k(Y;R)\bigr)}\,.
\end{equation*}

Let $Y,Z$ be orbifolds without boundary, and $h:Y\ra Z$ be a smooth,
proper map. Using the notation of Definition \ref{kh3def13}, define
the {\it pullback\/} $h^*:KC^k(Z;R)\ra KC^k(Y;R)$ by
\e
h^*:\ts\sum_{a\in A}\rho_a\bigl[X_a,\bs f_a,\bs C_a\bigr]\longmapsto
\ts\sum_{a\in A}\rho_a [Y\t_{h,Z,\bs f_a}X_a,\bs\pi_Y,h^*(\bs C_a)].
\label{kh4eq21}
\e
Here the coorientation for $(X_a,\bs f_a)$ pulls back to a
coorientation for $(Y\t_{h,Z,\bs f_a}X_a,\bs\pi_Y)$, as follows. Let
$p\in X_a$ and $(V_p,E_p,s_p,\psi_p)$ be a sufficiently small
Kuranishi neighbourhood in the germ at $p$, with submersion
$f_p:V_p\ra Z$ representing $\bs f$. Then the coorientation for
$(X,\bs f)$ gives orientations on the fibres of the orbibundle
$\bigl(\Ker\d f_p\op E_p\bigr)\ra V_p$. As in \eq{kh2eq14}, we lift
$(V_p,\ldots,\psi_p)$ to a Kuranishi neighbourhood
\begin{equation*}
\bigl(V_p^Y,E_p^Y,s_p^Y,\psi_p^Y\bigr)=
\bigl(Y\t_{h,Z,f_p}V_p,\pi_{V_p}^*(E_p),s_p\ci\pi_{V_p},
\pi_Y\t(\psi_p\ci\pi_{V_p})\t\chi\bigr)
\end{equation*}
on $Y\t_ZX_a$. Then the orbibundle $\bigl(\Ker\d\pi_Y\op
E_p^Y\bigr)\ra V_p^Y$ is naturally isomorphic to $\pi_Y^*\bigl(
(\Ker\d f_p\op E_p)\ra V_p\bigr)$, so the orientations on the fibres
of $\bigl(\Ker\d f_p\op E_p\bigr)\ra V_p$ lift to orientations on
the fibres of $\bigl(\Ker\d\pi_Y\op E_p^Y\bigr)\ra V_p^Y$, which
define a {\it coorientation\/} for~$(Y\t_ZX_a,\bs\pi_Y)$.

These $h^*:KC^k(Z;R)\ra KC^k(Y;R)$ satisfy $h^*\ci\d=\d\ci h^*$, so
they induce morphisms of cohomology groups $h^*:KH^k(Z;R)\ra
KH^k(Y;R)$. Pullbacks are functorial, that is, $(g\ci h)^*=h^*\ci
g^*$, on both cochains and cohomology.
\label{kh4def6}
\end{dfn}

\begin{rem} In Kuranishi cochains $[X,\bs f,\bs C]$ we assume $X$
is compact, so that $f(X)\subseteq Y$ is also compact. Thus, our
cochains are {\it compactly-supported}, and $KH^*(Y;R)$ is a form of
{\it compactly-supported cohomology}.

One could also try to define versions of Kuranishi cohomology
analogous to ordinary cohomology. The chains should be $[X,\bs f,\bs
C]$ for $X$ a {\it not necessarily compact\/} Kuranishi space, and
$\bs f:X\ra Y$ a {\it proper\/} strong submersion. The notion of
co-gauge-fixing data $\bs C$ for $(X,\bs f)$ when $X$ is noncompact
would certainly need revision. Working with noncompact $X$ raises
important issues to do with finiteness --- finiteness of indexing
sets for good coordinate systems, finiteness of $\Aut(X,\bs f,\bs
G)$ in Theorem \ref{kh3thm2} (this relies on Proposition
\ref{kh3prop4}, which uses the compactness of $X$), and so on. As we
do not need noncompact Kuranishi spaces in any of our applications
\cite{AkJo,Joyc2,Joyc3,Joyc4}, we will not develop this.
\label{kh4rem3}
\end{rem}

We shall prove {\it Poincar\'e duality} for Kuranishi (co)homology.
We do not assume $Y$ is compact.

\begin{dfn} Let $Y$ be an orbifold which is oriented, without
boundary, and of dimension $n$, and $R$ a $\Q$-algebra. Then for
$\bs f:X\ra Y$ a strong submersion, as in Definition \ref{kh2def22}
there is a 1-1 correspondence between coorientations for $(X,\bs f)$
and orientations for $X$. Choose an injective map $G_Y:Y\ra
P_m\subset P$, which is possible for $m\gg 0$. If $\bs C=(\bs
I,\bs\eta,C^i:i\in I)$ is co-gauge-fixing data for $(X,\bs f)$ with
$\bs I=\bigl(I,(V^i,\ldots,\psi^i),f^i:i\in I,\ldots\bigr)$, define
$\bs G_{\bs C}=(\bs I,\bs\eta,G^i:i\in I)$, where $G^i:E^i\ra P$ is
given by $G^i=\mu\ci\bigl((G_Y\ci f^i\ci\pi^i)\t C^i\bigr)$, for
$\mu$ in \eq{kh3eq2}. Definition \ref{kh3def6} for $\bs C$ says that
$C^i$ maps $E^i\ra P_n\subset P$ for some $n\gg 0$, and
$C^i\t(f^i\ci\pi^i):E^i\ra P\t Y$ is globally finite. Together with
$G_Y:Y\ra P_m$ injective and $\mu:P_m\t P_n\ra P_{m+n}$ globally
finite this implies that $G^i$ maps $E^i\ra P_{m+n}$ and is globally
finite, so that $\bs G_{\bs C}$ is {\it gauge-fixing data}
for~$(X,\bs f)$.

Define $R$-module morphisms $\Pi_\Kch^\Kh:KC^k(Y;R)\ra
KC_{n-k}(Y;R)$ for $k\in\Z$ by $\Pi_\Kch^\Kh:[X,\bs f,\bs C]\mapsto
[X,\bs f,\bs G_{\bs C}]$, using the coorientation for $(X,\bs f)$
from $[X,\bs f,\bs C]\in KC^k(Y;R)$ to determine the orientation on
$X$ for $[X,\bs f,\bs G_{\bs C}]\in KC_{n-k}(Y;R)$. Then
$\Pi_\Kch^\Kh$ takes Definitions \ref{kh4def2}(i)--(iv) to the same
relations, so they are well-defined, and $\pd\ci\Pi_\Kch^\Kh=
\Pi_\Kch^\Kh\ci\d$, so they induce morphisms of (co)homology groups
\e
\Pi_\Kch^\Kh:KH^k(Y;R)\longra KH_{n-k}(Y;R).
\label{kh4eq22}
\e

We will now construct an inverse $\Pi^\Kch_\Kh$ for \eq{kh4eq22}. On
the (co)chain level, we cannot just define $\Pi^\Kch_\Kh:[X,\bs
f,\bs G]\mapsto [X,\bs f,\bs C_{\bs G}]$, turning gauge-fixing data
$\bs G$ into co-gauge-fixing data $\bs C_{\bs G}$, since $\bs f$ may
not be a strong submersion. Instead, we {\it modify the Kuranishi
structure of\/} $X$ to make $\bs f$ a strong submersion.

Write $\pi:TY\ra Y$ for the tangent orbibundle of $Y$, and $z:Y\ra
TY$ for the zero section. Choose an orbifold Riemannian metric $g$
on $Y$. Then we can define the {\it exponential map\/} $\exp:V\ra
Y$, a smooth submersion from an open neighbourhood $V$ of $z(Y)$ in
$TY$, such that for $v\in V$ with $\pi(v)=y$ we have
$\exp(v)=\ga(1)$, where $\ga:[0,1]\ra Y$ is the unique geodesic
interval with $\ga(0)=y$, $\dot\ga(0)=v$ and
$\nabla_{\dot\ga(t)}\dot\ga(t)=0$ for $t\in[0,1]$. Choose $V$ such
that $\exp$ extends to the closure $\bar V$ of $V$ in $TY$, and
$\pi:\bar V\ra Y$ is a {\it proper} map. This is possible provided
$V$ is small enough.

Define a new Kuranishi structure $\ka_Y$ on $Y$ by the single
Kuranishi neighbourhood $(V,E,s,\psi)=(V,\pi^*(TY),
\id_{TY},z^{-1})$. Here $s=\id_{TY}$ is just the identity map on
$TY$, restricted to $V$ and regarded as a section of the bundle
$E=\pi^*(TY)\ra V$. Thus $s^{-1}(0)$ is $z(Y)$, so
$\psi=z^{-1}:z(Y)\ra Y$ makes sense. This induces $\ka_Y$ by taking
the germ of $\ka_Y$ at each $p\in Y$ to be the equivalence class of
$(V,E,s,\psi)$, regarded as a Kuranishi neighbourhood at $p$, and
the coordinate changes between such neighbourhoods to be the
identity.

Define strong submersions $\bs\pi:(Y,\ka_Y)\ra Y$ and
$\bs\exp:(Y,\ka_Y)\ra Y$ by $\pi:V\ra Y$ and $\exp:V\ra Y$ on the
Kuranishi neighbourhood $(V,E,s,\psi)$. The continuous maps induced
by $\bs\pi$ and ${\bf exp}$ are both the identity map $\id:Y\ra Y$,
but nonetheless $\bs\pi$ and $\bf exp$ are different as strongly
smooth maps, provided that $\dim Y>0$. Also, $z:Y\ra V$ induces an
embedding $\bs z:Y\ra(Y,\ka_Y)$ inducing the identity $\id:Y\ra Y$.
Choose an injective map $C_Y:E\ra P_{m'}\subset P$ for some $m'$,
which is possible for~$m'\gg 0$.

Suppose $X$ is a compact oriented Kuranishi space, $\bs f:X\ra Y$ is
strongly smooth, and $\bs G$ is gauge-fixing data for $(X,\bs f)$.
We will define modifications $X^Y,\bs f{}^Y,\bs C_{\bs G}^Y$ of
$X,\bs f,\bs G$ such that $X^Y$ is a compact Kuranishi space, $\bs
f{}^Y:X^Y\ra Y$ is a strong submersion with $(X^Y,\bs f{}^Y)$
cooriented, and $\bs C_{\bs G}^Y$ is co-gauge-fixing data for
$(X^Y,\bs f{}^Y)$. Define $X^Y$ to be the fibre product
$(Y,\ka_Y)\t_{\bs\pi,Y,\bs f}X$. The underlying topological space of
$X^Y$ is $X$, so $X^Y$ is really $X$ with a new Kuranishi structure.
Define a strong submersion $\bs f{}^Y:X^Y\ra Y$ by $\bs
f{}^Y=\bs\exp\ci\bs\pi_{(Y,\ka_Y)}$. There is a natural
coorientation on $\bigl((Y,\ka_Y),\bs\pi \bigr)$, and combining this
with the orientation on $X$ gives an orientation on $X^Y$, as in
Convention \ref{kh2conv2}(c). As $Y,X^Y$ are oriented and $\bs
f{}^Y:X^Y\ra Y$ is a strong submersion, as in Definition
\ref{kh2def22} we obtain a {\it coorientation} for~$(X^Y,\bs
f{}^Y)$.

Write $\bs G=(\bs I,\bs\eta,G^i:i\in I)$ with $\bs
I=\bigl(I,(V^i,\ldots,\psi^i),f^i:i\in I,\ldots\bigr)$. Define
$I^Y=\{i+n:i\in I\}$. Following \eq{kh2eq14}, for each $i\in I$
define a Kuranishi neighbourhood on $X^Y$ to be
\e
\begin{gathered}
(V^{i+n,Y},E^{i+n,Y},s^{i+n,Y},\psi^{i+n,Y})=
\bigl(V\t_{\pi,Y,f^i}V^i,\pi_V^*(E)\op
\pi_{V^i}^*(E^i),\\
(s\ci\pi_{V})\op(s^i\ci\pi_{V^i}),
(\psi\ci\pi_{V})\t(\psi^i\ci\pi_{V^i})\t\chi^i\bigr).
\end{gathered}
\label{kh4eq23}
\e
Define $f^{i+n,Y}:V^{i+n,Y}\ra Y$ by $f^{i+n,Y}=\exp\ci\pi_V$.
Define $\eta_{i+n}^Y:X^Y\ra[0,1]$ by $\eta_{i+n}^Y\equiv\eta_i$,
identifying $X^Y\cong X$. For $i,j\in I$, define
$\eta_{i+n}^{j+n,Y}:V^{j+n,Y}\ra[0,1]$ by $\eta_{i+n}^{j+n,Y}\equiv
\eta_i^j\ci\pi_{V^j}$. When $j\le i$ in $I$, define $V^{(i+n)(j+n),Y}
\subseteq V^{j+n,Y}$ to be $V\t_{\pi,Y,f^j}V^{ij}$, and let
$(\phi^{(i+n)(j+n),Y},\hat\phi^{(i+n)(j+n),Y}):\ab(V^{(i+n)(j+n),Y},\ab
\ldots,\ab\psi^{j+n,Y}\vert_{V^{(i+n)(j+n),Y}})\ab\ra
(V^{i+n,Y},\ldots,\psi^{i+n,Y})$ be induced by the identity on
$(V,E,s,\psi)$, and~$(\phi^{ij},\hat\phi^{ij})$.

As $(\bs I,\bs\eta)$ is a really good coordinate system for $(X,\bs
f)$, it is easy to show $\bigl(\bigl(I^Y,(V^{k,Y},
\ldots,\psi^{k,Y}),f^{k,Y}:k\in I^Y,\ldots\bigr),\eta_k^Y,
\eta_k^{l,Y}:k,l\in I^Y\bigr)$ is a really good coordinate system
for $(X^Y,\bs f{}^Y)$. To verify Definition \ref{kh3def3}(vi), let
$T\subseteq Y$ be the compact subset for $(\bs I,\bs\eta)$ in
Definition \ref{kh3def3}(vi). Then $\pi^{-1}(T)\cap\bar V$ is
compact as $\pi:\bar V\ra Y$ is proper, so
$T^Y=\exp\bigl(\pi^{-1}(T)\cap\bar V\bigr)\subseteq Y$ is
well-defined as $\exp$ extends to $\bar V$, and compact as
$\pi^{-1}(T)\cap\bar V$ is compact. Since $f^{k-n}(V^{k-n})\subseteq
T$ we find that $f^{k,Y}(V^{k,Y})\subseteq T^Y$ for all $k\in I^Y$,
proving Definition~\ref{kh3def3}(vi).

Let $(\bs{\check I}{}^Y,\bs{\check\eta}{}^Y)$ be the excellent
coordinate system for $(X^Y,\bs f{}^Y)$ constructed from this by
Algorithm \ref{kh3alg}, with $\bs{\check I}{}^Y=\bigl(\check
I^Y,(\check V^{k,Y},\ldots,\check\psi^{k,Y}): k\in\check
I^Y,\ldots\bigr)$. For each $k\in\check I^Y$, define $\check
C^{k,Y}:\check E^{k,Y}\ra P$ by $\check
C^{k,Y}=\mu\ci\bigl((C_Y\ci\pi_E)\t(G^{k-n}\ci
\pi_{E^{k-n}})\bigr)$. Since $C_Y$ maps $E\ra P_{m'}\subset P$, and
$G^{k-n}$ maps $E^{k-n}\ra P_{n'}\subset P$ for some $n'\ge 0$, and
$\mu$ maps $P_{m'}\t P_{n'}\ra P_{m'+n'}$, it follows that $\check
C^{k,Y}$ maps $\check E^{k,Y}\ra P_{m'+n'}$, one of the conditions
of Definition~\ref{kh3def7}.

We will prove $\check C^{k,Y}$ is {\it globally finite}. Observe
that $\check C^{k,Y}$ is the composition
\e
\smash{\xymatrix@C=13pt{\check E^{k,Y}\, \ar@{^{(}->}[r] &
E\t_YE^{k-n} \ar[rr]^{\pi_E\t\pi_{E^{k-n}}} && E\t E^{k-n}
\ar[rr]^{C_Y\t G^{k-n}} && P_{m'}\t P_{n'} \ar[r]^\mu & P_{m'+n'}.}}
\label{kh4eq24}
\e
We shall show that each morphism in \eq{kh4eq24} is globally finite,
in the sense of Definition \ref{kh3def6}. The first map in
\eq{kh4eq24} is an inclusion of open sets, and so globally finite
with $N=1$. For the second, let $T\subseteq Y$ be as above. The map
$Y\ra\N$ taking $y\mapsto\md{\Stab(y)}$ is upper semicontinuous, so
as $T$ is compact there exists $L\ge 1$ such that $\md{\Stab(y)}\le
L$ for all $y\in T$. As $f^{k-n}(V^{k-n})\subseteq T$, in the
formula \eq{kh2eq10} for the orbifold fibre product
$E\t_{\pi,Y,f^{k-n}\ci\pi^{k-n}}E^{k-n}$ the biquotient term
$f_*(\Stab(p))\backslash\Stab(f(p))/f'_*(\Stab(p'))$ is at most $L$
points, so $\pi_E\t\pi_{E^{k-n}}$ is globally finite with~$N=L$.

The third morphism $C_Y\t G^{k-n}$ in \eq{kh4eq24} is globally
finite as $C_Y$ is injective and $G^{k-n}$ is globally finite. For
the last, \eq{kh3eq2} implies that $\mu:P_{m'}\t P_{n'}\ra
P_{m'+n'}$ is globally finite with $N=\binom{m'+n'}{m'}$. Thus the
composition $\check C^{k,Y}$ is globally finite, and $\bs C_{\bs
G}^Y=(\bs{\check I}{}^Y,\bs{\check\eta}{}^Y,\check
C^{k,Y}:k\in\check I{}^Y)$ is {\it co-gauge-fixing data}
for~$(X^Y,\bs f{}^Y)$.

Define $R$-module morphisms $\Pi^\Kch_\Kh:KC_{n-k}(Y;R)\ra
KC^k(Y;R)$ by $\Pi^\Kch_\Kh:[X,\bs f,\bs G]\mapsto [X^Y,\bs
f{}^Y,\bs C_{\bs G}^Y]$. These are well-defined, and
$\d\ci\Pi^\Kch_\Kh= \Pi^\Kch_\Kh\ci\pd$, so they induce morphisms of
(co)homology groups
\e
\Pi^\Kch_\Kh:KH_{n-k}(Y;R)\longra KH^k(Y;R).
\label{kh4eq25}
\e
\label{kh4def7}
\end{dfn}

\begin{thm} Let\/ $Y$ be an oriented\/ $n$-orbifold without boundary.
Then $\Pi^\Kch_\Kh$ in \eq{kh4eq25} is the inverse of\/
$\Pi_\Kch^\Kh$ in {\rm\eq{kh4eq22},} so they are both isomorphisms.
Also {\rm\eq{kh4eq22}, \eq{kh4eq25}} are independent of choices in
their constructions.
\label{kh4thm3}
\end{thm}

\begin{proof} If $\dim Y=0$ then taking $G_Y\equiv(\es)$, $C_Y\equiv
(\es)$ we see that $\Pi^\Kch_\Kh$ is inverse to $\Pi_\Kch^\Kh$ at
the chain level, and the theorem is trivial. So suppose $n=\dim
Y>0$. Define a Kuranishi structure $\ka'$ on $[0,1]\t Y$ as follows.
On $[0,\ha)\t Y$ this comes from the orbifold structure on $[0,1]\t
Y$. On $[\ha,1]\t Y$ it comes from the product Kuranishi structure
on $[0,1]\t(Y,\ka_Y)$. Near $\{\ha\}\t Y$ we glue these two
Kuranishi structures together using $\bs z:Y\ra(Y,\ka_Y)$ in
Definition \ref{kh4def7}. Define strong submersions $\bs\pi,\bs\exp:
\bigl([0,1]\t Y,\ka'\bigr) \ra Y$ to be induced by $\id_Y:Y\ra Y$ on
$[0,\ha)\t Y$ and by $\bs\pi,\bs\exp:(Y,\ka_Y)\ra Y$ on $[\ha,1]\t
Y$. Using the orientation on $[0,1]$ we can construct a natural {\it
coorientation} for $\bigl(([0,1]\t Y,\ka'),\bs\pi\bigr)$,
independent of the orientation on~$Y$.

We shall define an excellent coordinate system $(\bs J_Y,\bs\ze_Y)$
for $\bigl(([0,1]\t Y,\ka',\bs\pi\bigr)$. The indexing set is
$J_Y=\{n+1,2n+1\}$, with Kuranishi neighbourhoods
\e
\begin{split}
(V^{n+1}_Y,\ldots,\psi^{n+1}_Y)&=\bigl([0,\ha)\t Y,[0,\ha)\t
Y,0,\id_{[0,\frac{1}{2})\t Y}\bigr),\\
(V^{2n\!+\!1}_Y,\ldots,\psi^{2n\!+\!1}_Y)&\!=\!
\bigl((\ts\frac{1}{3},1]\!\t\!V,(\frac{1}{3},1]\!\t\!
E,s\ci\pi_E,\id_{(\frac{1}{3},1]}\t\psi\bigr).
\end{split}
\label{kh4eq26}
\e
Here $E^{n+1}_Y\ra V^{n+1}_Y$ is the zero vector bundle, so
$E^{n+1}_Y\cong V^{n+1}_Y$ as manifolds. Define $V^{(2n+1)(n+1)}_Y=
(\frac{1}{3},\ha)\t Y\subset V^{n+1}_Y$, with coordinate change
$(\phi^{(2n+1)(n+1)}_Y,\hat\phi^{(2n+1)(n+1)}_Y):(V^{(2n+1)(n+1)}_Y,
\ldots,\psi^{n+1}_Y\vert_{V^{(2n+1)(n+1)}})\ra(V^{2n+1}_Y,\ab
\ldots,\ab\psi^{2n+1}_Y)$ given by $\phi^{(2n+1)(n+1)}_Y=
\smash{\id_{(\frac{1}{3},\frac{1}{2})}}\t z$, where $z:Y\ra V$ is
the zero section $z:Y\ra TY$, recalling that $V$ is a neighbourhood
of $z(Y)$ in $TY$, and $\hat\phi^{(2n+1)(n+1)}_Y=
\id_{(\frac{1}{3},\frac{1}{2})}\t(z_E\ci z)$, where $z_E:V\ra E$ is
the zero section of~$E$.

These Kuranishi neighbourhoods \eq{kh4eq26} are compatible with the
Kuranishi structure $\ka'$ on $[0,1]\t Y$, since for each point $p$
in $\Im\psi^{n+1}_Y$ or $\Im\psi^{2n+1}_Y$, there is a sufficiently
small $(V_p,\ldots,\psi_p)$ in the germ of $\ka'$ at $p$ and a
coordinate change from $(V_p,\ldots,\psi_p)$ to
$(V^{n+1}_Y,\ldots,\psi^{n+1}_Y)$ or
$(V^{2n+1}_Y,\ldots,\psi^{2n+1}_Y)$. Also $\bs\pi,\bs\exp$ are
represented by $\pi^{n+1}_Y=\exp^{n+1}_Y=\pi_Y:[0,\ha)\t Y \ra Y$ on
$(V^{n+1}_Y,\ldots,\psi^{n+1}_Y)$ and by $\pi^{2n+1}_Y=\pi\ci\pi_V,
\exp^{2n+1}_Y=\exp\ci\pi_V:(\ts\frac{1}{3},1]\t V\ra Y$
on~$(V^{2n+1}_Y,\ldots,\psi^{2n+1}_Y)$.

Let $\ze:[0,1]\ra[0,1]$ be continuous with $\ze\equiv 0$ on
$[0,\frac{1}{3}]$ and $\ze\equiv 1$ on $[\ha,1]$. Define
$\ze_{j,Y}:[0,1]\t Y\ra[0,1]$ for $j\in J_Y$ by $\ze_{n+1,Y}(t,y)=
1-\ze(t)$ and $\ze_{2n+1,Y}(t,y)=\ze(t)$. Define
$\ze_{i,Y}^j:V^j\ra[0,1]$ for $i,j\in J_Y$ by $\ze_{n+1,Y}^j\equiv
1-\ze\ci\pi_{[0,1]}$ and $\ze_{2n+1,Y}^j\equiv\ze\ci\pi_{[0,1]}$.
This data defines an {\it excellent coordinate system} $(\bs
J_Y,\bs\ze_Y)$ for~$\bigl(([0,1]\t Y,\ka'),\bs\pi\bigr)$.

Choose maps $C^j_Y:E^j_Y\ra P_m\subset P$ for $j\in J_Y$ and $m\gg
0$ such that $C^{n+1}_Y\vert_{\{0\}\t Y}\equiv(\es)$ and
$C^{2n+1}_Y\vert_{\{1\}\t E}\equiv\mu\ci\bigl(C_Y\t
(G_Y\ci\pi)\bigr)$, identifying $\{1\}\t E\cong E$, and
$C^j_Y\vert_{(E^j_Y)^\ci}:(E^j_Y)^\ci\ra P$ is injective. Then
$C^j_Y\t\pi_Y$ is globally finite, and $\bs C_Y=(\bs
J_Y,\bs\ze_Y,C_j^Y: j\in J_Y)$ is {\it co-gauge-fixing data}
for~$\bigl(([0,1]\t Y,\ka'),\bs\pi\bigr)$.

Define $R$-module morphisms $\Xi^k:KC^k(Y;R)\ra KC^{k-1}(Y;R)$~by
\e
\Xi^k:[X,\bs f,\bs C]\mapsto\bigl[([0,1]\t Y,\ka')\t_{\bs\pi,Y,\bs
f}X,\bs\exp\ci\bs\pi_{([0,1]\t Y,\ka')},\bs D_{\bs C}\bigr],
\label{kh4eq27}
\e
where $\bigl(([0,1]\t Y,\ka')\t_{\bs\pi,Y,\bs f}X,\bs\exp\ci
\bs\pi_{[0,1]\t Y}\bigr)$ has a coorientation constructed from that
of $(X,\bs f)$, and the co-gauge-fixing data $\bs D_{\bs C}$ is
defined as follows. We form co-gauge-fixing data $\bs C_Y\t_Y\bs C$
for $\bigl(([0,1]\t Y,\ka')\t_{\bs\pi,Y,\bs f}X,\bs\pi_Y\bigr)$ as
in Definition \ref{kh3def15}. Then we take $\bs D_{\bs C}$ to be the
same as $\bs C_Y\t_Y\bs C$, except that the maps $\check\pi^k_Y:
\check V^k\ra Y$ for $k\in\check I$, representing $\bs\pi_Y$, are
replaced by $\check\exp^k=\exp^j\ci\pi_{V^j_Y}$ on $\check V^k\cap
(V^j_Y\t_YV^i)$, representing~$\bs\exp\ci\bs\pi_{([0,1]\t Y,\ka')}$.

To verify $\bs D_{\bs C}$ is co-gauge-fixing data, there is one
nontrivial thing to check. The proof that $\bs C_Y\t_Y\bs C$ is
co-gauge-fixing data shows that $\check C^k\t(\check\pi^k_Y
\ci\check\pi^k):\check E^k\ra P\t Y$ is a globally finite map, but
we have replaced $\check\pi^k_Y$ by $\check\exp^k$, so we must prove
that $\check C^k\t(\check\exp^k\ci\check\pi^k):\check E^k\ra P\t Y$
is globally finite. Write $\check E^k=\check E^{k\prime}\amalg
\check E^{k\prime\prime}$, where $\check E^{k\prime}$ lies over
$V^{n+1}_Y$ and $\check E^{k\prime\prime}$ lies over $V^{2n+1}_Y$.
As $\pi^{n+1}_Y=\exp^{n+1}_Y$ on $V^{n+1}_Y$, $\check\pi^k_Y$ and
$\check\exp^k$ coincide on $\check E^{k\prime}$, so $\check
C^k\t(\check\exp^k\ci\check\pi^k)$ is globally finite on $\check
E^{k\prime}$ as $\check C^k\t(\check\pi^k_Y\ci\check\pi^k)$ is. On
$\check E^{k\prime\prime}$ we use the fact that
$C^{2n+1}_Y:E^{2n+1}_Y\ra P_m\subset P$ is globally finite to deduce
that $\check C^k:\check E^{k\prime\prime}\ra P$ is globally finite,
following the proof that $\bs G\t_Y\bs{\ti C}$ is gauge-fixing data,
and temporarily treating $\bs C_Y$ as gauge-fixing data. But $\check
C^k:\check E^{k\prime\prime}\ra P$ globally finite implies $\check
C^k\t(\check\exp^k\ci\check\pi^k):\check E^{k\prime\prime}\ra P\t Y$
is globally finite, as we want.

From \eq{kh2eq21} we find that in cooriented Kuranishi spaces we
have
\begin{align*}
\pd &\bigl(([0,1]\t Y,\ka')\t_{\bs\pi,Y,\bs
f}X,\bs\exp\ci\bs\pi_{[0,1]\t Y}\bigr)\cong -\bigl(\{0\}\t X,\bs
f\bigr)\,\amalg\\
&\amalg \bigl(\{1\}\t X^Y,\bs f{}^Y\bigr)-\bigl(([0,1]\t
Y,\ka')\t_{\bs\pi,Y,\bs f\vert_{\pd X}}\pd
X,\bs\exp\ci\bs\pi_{[0,1]\t Y}\bigr).
\end{align*}
Adding co-gauge-fixing data to this, in $KC^k(Y;R)$ we have
\e
\begin{split}
\d&\bigl[([0,1]\t Y,\ka')\t_{\bs\pi,Y,\bs
f}X,\bs\exp\ci\bs\pi_{[0,1]\t Y},\bs D_{\bs C}\bigr]=-[X,\bs f,\bs
C]\\
& [X^Y,\bs f{}^Y,\bs C^Y_{{\bs G}_{\bs C}}]-\bigl[([0,1]\t
Y,\ka')\t_{\bs\pi,Y,\bs f\vert_{\pd X}}\pd
X,\bs\exp\ci\bs\pi_{[0,1]\t Y},\bs D_{\bs C\vert_{\pd X}}\bigr],
\end{split}
\label{kh4eq28}
\e
where $\bs G_{\bs C}$ is the gauge-fixing data for $(X,\bs f)$
constructed from $\bs C$ in Definition \ref{kh4def7}, and $\bs
C^Y_{{\bs G}_{\bs C}}$ is the co-gauge-fixing data for $(X^Y,\bs
f{}^Y)$ constructed from $\bs G_{\bs C}$ in Definition
\ref{kh4def7}. The first term on the r.h.s.\ of \eq{kh4eq28} holds
because $C^{n+1}_Y$ restricts to $(\es)$ on $\{0\}\t Y$, where
$(\es)$ is the identity for $\mu$. The second term holds because
$C^{2n+1}_Y$ restricts on $\{1\}\t E$ to
$\mu\ci\bigl(C_Y\t(G_Y\ci\pi)\bigr)$, which combines the extra
factor $G_Y\ci\pi$ in going from $\bs C$ to ${\bs G}_{\bs C}$ with
the extra factor $C_Y$ in going from ${\bs G}_{\bs C}$ to $\bs
C^Y_{{\bs G}_{\bs C}}$, and using associativity of~$\mu$.

Combining \eq{kh4eq27}, \eq{kh4eq28} and Definition \ref{kh4def7}
shows that
\e
\d\ci\Xi^k=-\id+\Pi^\Kch_\Kh\ci\Pi_\Kch^\Kh-\Xi^{k+1}\ci\d:
KC^k(Y;R)\longra KC^k(Y;R).
\label{kh4eq29}
\e
Passing to cohomology, this implies $\Pi^\Kch_\Kh\ci\Pi_\Kch^\Kh:
KH^k(Y;R)\ra KH^k(Y;R)$ is the identity. An almost identical proof
yields $\Upsilon^{n-k}:KC_{n-k}(Y;R)\!\ra\!KC_{n-k+1}(Y;R)$
satisfying
\begin{equation*}
\pd\ci\Upsilon^{n-k}=-\id+\Pi_\Kch^\Kh\ci\Pi^\Kch_\Kh-
\Upsilon^{n-k-1}\ci\pd:KC_{n-k}(Y;R)\ra KC_{n-k}(Y;R).
\end{equation*}
In defining $\Upsilon^{n-k}$ we replace the condition on
$C^{2n+1}_Y$ above by $C^{2n+1}_Y\vert_{\{1\}\t
E}\equiv\mu\ci\bigl((G_Y\ci\exp\ci \pi_V)\t C_Y\bigr)$, identifying
$\{1\}\t E\cong E$. Passing to homology shows that
$\Pi_\Kch^\Kh\ci\Pi^\Kch_\Kh: KH_{n-k}(Y;R)\ra KH_{n-k}(Y;R)$ is the
identity. Therefore \eq{kh4eq25} is the inverse of~\eq{kh4eq22}.

For the final part of the theorem, by the first part it is enough to
show that \eq{kh4eq22} is independent of choices, and the only
choice in its construction is that of $G_Y:Y\ra P_m\subset P$. If
$G_Y,G_Y'$ are possible choices inducing operators
$\Pi_\Kch^\Kh,\Pi_\Kch^{\Kh\prime}$ on (co)chains, then choose
globally finite $G_{[0,1]\t Y}:[0,1]\t Y\ra P_{m''}\subset P$ for
$m''\gg 0$ with $G_{[0,1]\t Y}\vert_{\{0\}\t Y}=G_Y$, $G_{[0,1]\t
Y}\vert_{\{1\}\t Y}=G_Y'$.

Let $\bs f:X\ra Y$ be a cooriented strong submersion and $\bs C$ be
co-gauge-fixing data for $(X,\bs f)$, where $\bs C=(\bs
I,\bs\eta,C^i:i\in I)$ with $\bs I=\bigl(I,(V^i,\ldots,\psi^i),
f^i:i\in I,\ldots\bigr)$. Set $\ti I=\{i+1:i\in I\}$. For $i\in I$,
define $(\ti V^{i+1},\ti E^{i+1},\ti
s^{i+1},\ti\psi^{i+1})=\bigl([0,1]\t V^i,[0,1]\t E^i,
s^i\ci\pi_{V^i},\id_{[0,1]}\t\psi^i\bigr)$, $\ti
f^{i+1}=f^i\ci\pi_{V^i}:\ti V^{i+1}\ra Y$, and
$\ti\eta_{i+1}:[0,1]\t X\ra[0,1]$ by $\ti\eta_{i+1}=\eta_i\ci\pi_X$,
and $H^{i+1}:\ti E^{i+1}\ra P$ by $H^{i+1}(u,e)=\mu\bigl(G_{[0,1]\t
Y}(u,f^i\ci\pi_{V^i}(e)),C^i(e)\bigr)$.

For $i,j\in I$, define $\ti V^{(i+1)(j+1)}=[0,1]\t V^{ij}$,
$\ti\phi^{(i+1)(j+1)}= \id_{[0,1]}\t\phi^{ij}$,
$\hat{\ti\phi}{}^{(i+1)(j+1)} =\id_{[0,1]}\t\hat\phi^{ij}$ and
$\ti\eta_{i+1}^{j+1}:\ti V^{j+1}\ra[0,1]$ by
$\ti\eta_{i+1}^{j+1}=\eta_i^j\ci\pi_{V^j}$. Write $\bs{\ti
I}=\bigl(\ti I,(\ti V^i,\ldots,\ti\psi^i),\ti f^i:i\in\ti
I,\ldots\bigr)$, $\bs{\ti\eta}=(\ti\eta_i:i\in\ti I$,
$\ti\eta_i^j:i,j\in\ti I)$, and $\bs H_{\bs C}=(\bs{\ti I},
\bs{\ti\eta},H^i:i\in\ti I)$. Then $\bs H_{\bs C}$ is {\it
gauge-fixing data} for $\bigl([0,1]\t X,\bs f\ci\bs\pi_X\bigr)$.

Define $R$-module morphisms $\Th^k:KC^k(Y;R)\ra KC_{n-k-1}(Y;R)$ by
\begin{equation*}
\Th^k:[X,\bs f,\bs C]\longmapsto \bigl[[0,1]\t X,\bs
f\ci\bs\pi_X,\bs H_{\bs C}].
\end{equation*}
A similar proof to \eq{kh4eq29} using $G_{[0,1]\t Y}\vert_{\{0\}\t
Y}=G_Y$, $G_{[0,1]\t Y}\vert_{\{1\}\t Y}=G_Y'$ gives
\begin{equation*}
\pd\ci\Th^k=-\Pi_\Kch^\Kh+\Pi_\Kch^{\Kh\prime}-
\Th^{k+1}\ci\d:KC^k(Y;R)\ra KC_{n-k}(Y;R).
\end{equation*}
Passing to (co)homology yields $\Pi_\Kch^\Kh=\Pi_\Kch^{\Kh\prime}:
KH^k(Y;R)\ra KH_{n-k}(Y;R)$. This completes the proof.
\end{proof}

We show Kuranishi and compactly-supported cohomology are isomorphic.

\begin{dfn} Let $Y$ be an $n$-orbifold without boundary, and $R$ a
$\Q$-algebra. We will define $R$-module morphisms
$\Pi_\cs^\Kch:H^k_\cs(Y;R)\ra KH^k(Y;R)$ for all $k\ge 0$. First
suppose $Y$ is orientable, and choose an orientation. Define
$\Pi_\cs^\Kch$ to be the composition
\begin{equation*}
\smash{\xymatrix@C=30pt{ H^k_\cs(Y;R) \ar[r]^(0.4){\Pd} &
H_{n-k}^\rsi(Y;R) \ar[r]^{\Pi_\rsi^\Kh} & KH_{n-k}(Y;R)
\ar[r]^{\Pi^\Kch_\Kh} & KH^k(Y;R), }}
\end{equation*}
where $\Pd,\Pi_\rsi^\Kh,\Pi^\Kch_\Kh$ are as in \eq{kh4eq3},
\eq{kh4eq17}, \eq{kh4eq25} respectively. Changing the orientation of
$Y$ changes the signs of $\Pd,\Pi^\Kch_\Kh$. Thus $\Pi_\cs^\Kch$ is
independent of the choice of orientation on~$Y$.

Even if $Y$ is not orientable, it has an {\it orientation bundle}, a
principal $\Z_2$-bundle $O\ra Y$. The argument above still works if
we replace $H_*^\rsi(Y;R)$ and $KH_*(Y;R)$ by the {\it twisted
homology groups} $H_*^\rsi,KH_*(Y;O\t_{\{\pm 1\}}R)$. Geometrically,
$H_k^\rsi(Y;\ab O\t_{\{\pm 1\}}R)$ is defined using chains $(\si,o)$
where $\si:\De_k\ra Y$ is a smooth map and $o$ is an orientation of
the fibres of $f^*(TY)\ra\De_k$, with $(\si,-o)=-(\si,o)$, and
$KH_*(Y;O\t_{\{\pm 1\}}R)$ is defined using chains $[X,\bs f,\bs G]$
with $(X,\bs f)$ cooriented, rather than $X$ oriented. (One can
still make sense of coorientations for $(X,\bs f)$ even if $\bs f$
is only strongly smooth).
\label{kh4def8}
\end{dfn}

Poincar\'e duality for orbifolds $Y$ without boundary when $R$ is a
$\Q$-algebra as discussed in \S\ref{kh41}, Corollary \ref{kh4cor1},
and Theorem \ref{kh4thm3} imply:

\begin{cor} Let\/ $Y$ be an orbifold without boundary, and\/ $R$ a
$\Q$-algebra. Then $\Pi_\cs^\Kch:H^k_\cs(Y;R)\ra KH^k(Y;R)$ is an
isomorphism for $k\ge 0,$ with\/ $KH^k(Y;R)=0$ when~$k<0$.
\label{kh4cor2}
\end{cor}

\subsection{The case $\pd Y\ne\es$, and relative Kuranishi
(co)homology}
\label{kh45}

We now generalize \S\ref{kh44} to define Kuranishi cohomology
$KH^*(Y;R)$ for orbifolds $Y$ with $\pd Y\ne\es$. To show that
$KH^*(Y;R)\cong H^*_\cs(Y;R)$ using Poincar\'e duality involves {\it
relative\/} Kuranishi homology and cohomology, so we define these as
well. We will be brief, and give proofs only in the
case~$\pd^2Y=\es$.

We first explain the parts of classical (relative) (co)homology
theory we will need. In \S\ref{kh41} we discussed {\it Poincar\'e
duality isomorphisms} \eq{kh4eq3} for manifolds and orbifolds $Y$
without boundary. Here is how to extend this to manifolds with
boundary and g-corners. Suppose $Y$ is an oriented manifold with
g-corners, of dimension $n$, not necessarily compact. Write $\io:\pd
Y\ra Y$ for the natural immersion, and $\io(\pd Y)$ for its image as
a subset of $Y$. Then equation \eq{kh4eq3} extends to two families
of Poincar\'e duality isomorphisms:
\ea
&\Pd':H^k_\cs(Y,\io(\pd Y);R)\longra H_{n-k}\bigl(Y;R\bigr),
\label{kh4eq30}\\
&\Pd'':H^k_\cs(Y;R)\longra H_{n-k}\bigl(Y,\io(\pd Y);R\bigr),
\label{kh4eq31}
\ea
where $H^*_\cs(Y,\io(\pd Y);R)$ is the {\it compactly-supported
relative cohomology\/} and\break $H_*\bigl(Y,\io(\pd Y);R\bigr)$ the
{\it relative homology} of $(Y,\io(Y))$. The usual form in which
Poincar\'e duality is stated, as in Bredon \cite[\S VI.9]{Bred} for
instance, is \eq{kh4eq31} for {\it compact\/} $Y$, with $H^*(Y;R)$
in place of~$H^*_\cs(Y;R)$.

The interior $Y^\ci$ is a manifold without boundary, and there are
natural isomorphisms $H^*_\cs(Y,\io(\pd Y);R)\cong
H^*_\cs(Y^\ci;R)$, $H_*(Y;R)\cong H_*(Y^\ci;R)$ which identify
$\Pd'$ in \eq{kh4eq30} with $\Pd$ in \eq{kh4eq3} for $Y^\ci$. Thus,
\eq{kh4eq30} follows trivially from~\eq{kh4eq3}.

If $\pd^2Y\ne\es$ then $\io:\pd Y\ra Y$ is not injective, and its
image $\io(\pd Y)$ is not a manifold. However, if $\pd^2Y=\es$ then
$\io$ is injective and $\io(\pd Y)$ is a submanifold of $Y$. {\it
Relative singular homology\/} when $\pd^2Y=\es$ is the homology of
the chain complex $\bigl(C_*^\rsi(Y,\io(\pd Y);R),\pd\bigr)$, where
\begin{equation*}
C_k^\rsi\bigl(Y,\io(\pd Y);R\bigr)=C_k^\rsi(Y;R)/
\io_*\bigl(C_k^\rsi(\pd Y;R)\bigr),
\end{equation*}
and $\io_*:C_k^\rsi(\pd Y;R)\ra C_k^\rsi(Y;R)$ is as in
\S\ref{kh41}. This gives a short exact sequence of chain complexes:
\e
\xymatrix@C=15pt{ 0 \ar[r] & C_*^\rsi(\pd Y;R) \ar[r]^{\io_*} &
C_*^\rsi(Y;R) \ar[r]^(0.35){\pi} & C_*^\rsi\bigl(Y,\io(\pd
Y);R\bigr) \ar[r] & 0,}
\label{kh4eq32}
\e
where $\pi$ is the projection $C_*^\rsi(Y;R)\ra C_*^\rsi(Y;R)/
\io_*(C_k^\rsi(\pd Y;R))$. Equation \eq{kh4eq32} induces a long
exact sequence of (relative) homology groups:
\e
\text{\begin{small}$\displaystyle \xymatrix@C=11pt{ \!\cdots\!\ar[r]
& H_k^\rsi(\pd Y;R) \ar[r]^{\io_*} & H_k^\rsi(Y;R)
\ar[r]^(0.4){\pi_*} & H_k^\rsi\bigl(Y,\io(\pd Y);R\bigr)
\ar[r]^{\pd_*} & H_{k-1}^\rsi(\pd Y;R) \ar[r]
&\!\cdots.}$\end{small}}\!\!\!
\label{kh4eq33}
\e
Similarly, when $\pd^2Y=\es$ we have a long exact sequence of
cohomology groups:
\e
\text{\begin{small}$\displaystyle \xymatrix@C=11pt{ \!\cdots\!\ar[r]
& H^{k-1}_\cs(\pd Y;R) \ar[r]^(0.45){\pd^*} & H^k_\cs\bigl(Y,\io(\pd
Y);R\bigr) \ar[r]^(0.6){\pi^*} & H^k_\cs(Y;R) \ar[r]^{\io^*} &
H^k_\cs(\pd Y;R) \ar[r] &\!\cdots.}$\end{small}}\!\!\!
\label{kh4eq34}
\e

Now $\pd Y$ is an oriented $(n-1)$-manifold without boundary as
$\pd^2Y=\es$, so Poincar\'e duality isomorphisms \eq{kh4eq3} hold
for $\pd Y$. These combine with the isomorphisms
\eq{kh4eq30}--\eq{kh4eq31} and sequences \eq{kh4eq33}--\eq{kh4eq34}
in a commutative diagram
\begin{equation*}
\xymatrix@C=15pt{ \cdots\ar[r] & H^{k-1}_\cs(\pd Y;R)
\ar[d]^{\Pd}_\cong \ar[r]^(0.45){\pd^*} & H^k_\cs\bigl(Y,\io(\pd
Y);R\bigr) \ar[d]^{\Pd'}_\cong \ar[r]^(0.6){\pi^*} & H^k_\cs(Y;R)
\ar[d]^{\Pd''}_\cong \ar[r]^{\io^*} & \cdots\\
\cdots\ar[r] & H_{n\!-\!k}^\rsi(\pd Y;R) \ar[r]^{\io_*} &
H_{n\!-\!k}^\rsi(Y;R) \ar[r]^(0.4){\pi_*} &
H_{n\!-\!k}^\rsi\bigl(Y,\io(\pd Y);R\bigr) \ar[r]^(0.7){\pd_*}  &
\cdots. }
\end{equation*}
As in \S\ref{kh41}, the Poincar\'e duality isomorphisms
\eq{kh4eq30}--\eq{kh4eq31} also hold for orbifolds $Y$ provided $R$
is a $\Q$-algebra. We now define {\it relative Kuranishi homology}.

\begin{dfn} Let $Y,Z$ be orbifolds, $h:Y\ra Z$ a smooth map, and $R$
a $\Q$-algebra. Then Definition \ref{kh4def4} gives
$h_*:KC_k(Y;R)\ra KC_k(Z;R)$ with $\pd\ci h_*=h_*\ci\pd$. Define the
{\it relative Kuranishi chains} $KC_k(Z,h(Y);R)$ to be
\begin{equation*}
KC_k\bigl(Z,h(Y);R\bigr)=KC_k(Z;R)/h_*\bigl(KC_k(Y;R)\bigr).
\end{equation*}
Then $\pd:KC_k(Z;R)\ra KC_{k-1}(Z;R)$ descends to
$\pd:KC_k(Z,h(Y);R)\ra KC_{k-1}(Z,h(Y);R)$ with $\pd^2=0$. Define
{\it relative Kuranishi homology} $KH_*\bigl(Z,\ab h(Y);R\bigr)$ to
be the homology of $\bigl(KC_*(Z,h(Y);R),\pd\bigr)$.
\label{kh4def9}
\end{dfn}

Suppose now that $h:Z\ra Y$ is an {\it embedding}. Then
$h_*:KC_*(Y;R)\ra KC_*(Z;R)$ and $h_*:C_*^\rsi(Y;R)\ra
C_*^\rsi(Z;R)$ are injective, so we have a commutative diagram of
chain complexes with exact rows:
\e
\begin{gathered}
\xymatrix@C=15pt{ 0 \ar[r] & C_*^\rsi(Y;R) \ar[r]^{h_*}
\ar[d]^{\Pi_\rsi^\Kh} & C_*^\rsi(Z;R) \ar[r]^(0.35){\pi}
\ar[d]^{\Pi_\rsi^\Kh} & C_*^\rsi\bigl(Z,h(Y);R\bigr)
\ar[d]^{\Pi_\rrsi^\rKh}
\ar[r] & 0\\
0 \ar[r] & KC_*(Y;R) \ar[r]^{h_*} & KC_*(Z;R) \ar[r]^(0.35){\pi} &
KC_*\bigl(Z,h(Y);R\bigr) \ar[r] & 0,}
\end{gathered}
\label{kh4eq35}
\e
where $\Pi_\rrsi^\rKh$ is the projection from relative singular
chains to relative Kuranishi chains defined as for $\Pi_\rsi^\Kh$ in
Definition~\ref{kh4def4}.

Equation \eq{kh4eq35} induces a commutative diagram of homology
groups
\e
\begin{gathered}
\!\!\!\!\text{\begin{small}$\displaystyle \xymatrix@C=11pt{
\!\cdots\!\ar[r] & H_k^\rsi(Y;R) \ar[d]^{\Pi_\rsi^\Kh}_\cong
\ar[r]^{h_*} & H_k^\rsi(Z;R) \ar[d]^{\Pi_\rsi^\Kh}_\cong
\ar[r]^(0.4){\pi_*} & H_k^\rsi\bigl(Z,h(Y);R\bigr)
\ar[d]^{\Pi_\rrsi^\rKh} \ar[r]^(0.55){\pd_*}
& H_{k-1}^\rsi(Y;R) \ar[d]^{\Pi_\rsi^\Kh}_\cong \ar[r] &\cdots\\
\!\cdots\!\ar[r] & KH_k(Y;R) \ar[r]^{h_*} & KH_k(Z;R)
\ar[r]^(0.4){\pi_*} & KH_k\bigl(Z,h(Y);R\bigr) \ar[r]^(0.55){\pd_*}
& KH_{k-1}(Y;R) \ar[r] &\cdots,}$\end{small}}
\!\!\!\!\!\!\!\!\!\!\!\!\!\!\!\!\!\!
\end{gathered}
\label{kh4eq36}
\e
with long exact rows. Corollary \ref{kh4cor1} implies that the
columns $\Pi_\rsi^\Kh$ in \eq{kh4eq36} are isomorphisms. Since the
rows are exact, properties of long exact sequences then imply that
the remaining column $\Pi_\rrsi^\rKh$ is an isomorphism. This
implies:

\begin{cor} Suppose $h:Y\ra Z$ is a smooth embedding of
orbifolds, and\/ $R$ is a $\Q$-algebra. Then $\Pi_\rrsi^\rKh:
H_k^\rsi\bigl(Z,h(Y);R\bigr)\ra KH_k\bigl(Z,h(Y);R\bigr)$ is an
isomorphism for all\/ $k\ge 0,$ with\/ $KH_k\bigl(Z,h(Y);R\bigr)=0$
for~$k<0$.
\label{kh4cor3}
\end{cor}

The author expects that the assumption $h$ is an embedding is
unnecessary in Corollary \ref{kh4cor3}. Note that if $Y$ is an
orbifold with $\pd^2Y=\es$ then $\io:\pd Y\ra Y$ is an embedding, so
$\Pi_\rrsi^\rKh:H_k^\rsi\bigl(Y,\io(\pd Y);R\bigr)\ra
KH_k\bigl(Y,\io(\pd Y);R\bigr)$ is an isomorphism. Next we define
$KH^*(Y;R)$ and $KH^*(Y,\io(\pd Y);R)$ when $Y$ is an orbifold with
g-corners, allowing $\pd Y\ne\es$.

\begin{dfn} Let $X$ be a Kuranishi space, $Y$ an orbifold with
g-corners, and $\bs f:X\ra Y$ a strong submersion. Then if $\pd
Y\ne\es$, as for submersions of manifolds in Definition
\ref{kh2def4} we have a decomposition $\pd X=\pd_+^{\bs
f}X\amalg\pd_-^{\bs f}X$ and strong submersions $\bs f_+:\pd_+^{\bs
f}X\ra Y$ and $\bs f_-:\pd_-^{\bs f}X\ra\pd Y$, where $\bs f_+=\bs
f\vert_{\smash{\pd_+^{\bs f}X}}$. In particular, if $\pd_-^{\bs
f}X\ne\es$ then $\bs f\vert_{\pd X}$ is {\it not a strong
submersion}, but only strongly smooth.

Let $Y$ be an orbifold and $R$ a $\Q$-algebra. When $\pd Y=\es$,
Definition \ref{kh4def6} defined the Kuranishi cochains $KC^*(Y;R)$,
with $\d:KC^k(Y;R)\ra KC^{k+1}(Y;R)$ given by \eq{kh4eq20}. If $\pd
Y\ne\es$, equation \eq{kh4eq20} is no longer valid, since $\bs
f_a\vert_{\pd X_a}$ may not be a strong submersion. Instead, when
$\pd Y\ne\es$ we define $KC^*(Y;R)$ exactly as in Definition
\ref{kh4def6}, but we define $\d:KC^k(Y;R)\ra KC^{k+1}(Y;R)$ by
\begin{equation*}
\d:\ts\sum_{a\in A}\rho_a[X_a,\bs f_a,\bs C_a]\longmapsto
\ts\sum_{a\in A}\rho_a[\pd_+^{\bs f_a}X_a,\bs f_{a,+},\bs C_a
\vert_{\pd_+^{\bs f_a}X_a}].
\end{equation*}
Since the natural orientation-reversing involution
$\bs\si:\pd^2X_a\ra\pd^2X_a$ restricts to an involution of
$\pd_+^{\bs f_{a,+}}\bigl(\pd_+^{\bs f_a}X_a\bigr)$, relation
Definition \ref{kh4def2}(ii) still implies that $\d^2=0$. We then
define {\it Kuranishi cohomology\/} $KH^*(Y;R)$ to be the cohomology
of~$\bigl(KC^*(Y;R),\d\bigr)$.

Define a morphism $\io^*:KC^k(Y;R)\ra KC^k(\pd Y;R)$ by
\e
\io^*:\ts\sum_{a\in A}\rho_a[X_a,\bs f_a,\bs C_a]\longmapsto
(-1)^k\ts\sum_{a\in A}\rho_a[\pd_-^{\bs f_a}X_a,\bs f_{a,-},\bs C_a
\vert_{\pd_-^{\bs f_a}X_a}],
\label{kh4eq37}
\e
where $\bs f_{a,-}:\pd_-^{\bs f_a}X_a\ra\pd Y$ is a strong
submersion as above. Then $\io^*$ takes relations Definition
\ref{kh4def2}(i)--(iv) in $KC^k(Y;R)$ to the same relations in
$KC^k(\pd Y;R)$, and so is well-defined. We think of $\io^*$ as the
{\it pullback\/} under~$\io:\pd Y\ra Y$.

If $X$ is an oriented Kuranishi space and $\bs f:X\ra Y$ a strong
submersion, then the natural orientation-reversing involution
$\bs\si:\pd^2X\ra\pd^2X$ restricts to an orientation-reversing
isomorphism $\bs\si:\pd_-^{\bs f_+}(\pd_+^{\bs f}X)\ra \pd_+^{\bs
f_-}(\pd_-^{\bs f}X)$. Using this we see that
$\d\ci\io^*=\io^*\ci\d$ as maps $KC^k(Y;R)\ra KC^{k+1}(\pd Y;R)$.
The purpose of the factor $(-1)^k$ in \eq{kh4eq37} is to make
$\d\ci\io^*=\io^*\ci\d$, so that $\io^*$ is a chain map; without it
we would get $\d\ci\io^*=-\io^*\ci\d$, since $\bs\si$ is
orientation-reversing. Since $\io^*$ is a chain map it induces
morphisms~$\io^*:KH^k(Y;R)\ra KH^k(\pd Y;R)$.

Define the {\it relative Kuranishi cochains} $KC^k(Y,\io(Y);R)$ to
be the kernel of $\io^*:KC^k(Y;R)\ra KC^k(\pd Y;R)$. Then
$\d:KC^k(Y;R)\ra KC^{k+1}(Y;R)$ restricts to $\d:KC^k(Y,\io(Y);R)\ra
KC^{k+1}(Y,\io*(Y);R)$, since $\d\ci\io^*=\io^*\ci\d$. Define {\it
relative Kuranishi cohomology} $KH^*(Y,\io(Y);R)$ to be the
cohomology of~$\bigl(KC^k(Y,\io(Y);R),\d\bigr)$.
\label{kh4def10}
\end{dfn}

For the rest of this section we restrict to the case that
$\pd^2Y=\es$, so that $\io:\pd Y\ra Y$ is an embedding. Then one can
show that $\io^*:KC^k(Y;R)\ra KC^k(\pd Y;R)$ is {\it surjective};
one way to do this is to identify $Y\cong[0,\ep]\t\pd Y$ near $\pd
Y$ for small $\ep>0$, and then for each $[X,\bs f,\bs C]\in KC^k(\pd
Y;R)$, to construct explicit $\bigl[[0,\ep]\t X,\bs f',\bs
C'\bigr]\in KC^k(Y;R)$ with $\io^*:\bigl[[0,\ep]\t X,\bs f',\bs
C'\bigr]\mapsto[X,\bs f,\bs C]$. Thus we have a short exact sequence
of cochain complexes
\begin{equation*}
\xymatrix@C=18pt{ 0 \ar[r] & KC^*(Y,\io(Y);R) \ar[r]^(0.55){\rm inc}
& KC^*(Y;R) \ar[r]^{\io^*} &  KC^*(\pd Y;R) \ar[r] & 0,}
\end{equation*}
where `inc' is the inclusion. This induces a long exact sequence in
cohomology, the analogue of~\eq{kh4eq34}:
\begin{equation*}
\text{\begin{small}$\displaystyle\xymatrix@C=9.5pt{ \cdots\!\ar[r] &
KH^{k-1}(\pd Y;R) \ar[r]^(0.45){\pd^*} & KH^k\bigl(Y,\io(\pd
Y);R\bigr) \ar[r]^(0.58){{\rm inc}*} & KH^k(Y;R) \ar[r]^{\io^*} &
KH^k(\pd Y;R) \ar[r] &\!\cdots.}$\end{small}}
\end{equation*}

Our goal is to prove that $KH^*(Y;R)\cong H^*_\cs(Y;R)$ and
$KH^*(Y,\io(Y);R)\cong H^*_\cs(Y,\io(Y);R)$, generalizing Corollary
\ref{kh4cor2}. To do this we need analogues of the Poincar\'e
duality morphisms $\Pi_\Kch^\Kh,\Pi^\Kch_\Kh$ of Definition
\ref{kh4def7} for $KH^*(Y;R)$ and $KH^*(Y,\io(Y);R)$. That is, we
wish to define chain maps
\ea
&\xymatrix{\Pi_\rKch^\Kh:KC^k(Y,\io(\pd Y);R) \ar[r]&
KC_{n-k}(Y;R),}
\label{kh4eq38}\\
&\xymatrix{\Pi^\rKch_\Kh:KC_{n-k}(Y;R)\ar[r]& KC^k(Y,\io(\pd Y);R),}
\label{kh4eq39}\\
&\xymatrix{\Pi_\Kch^\rKh:KC^k(Y;R)\ar[r]& KC_{n-k}(Y,\io(\pd Y);R),}
\label{kh4eq40}\\
&\xymatrix{\Pi^\Kch_\rKh:KC_{n-k}(Y,\io(\pd Y);R)\ar@{-->}[r]&
KC^k(Y;R),}
\label{kh4eq41}
\ea
which should induce the analogues of $\Pd',\Pd''$ and their
inverses.

We define $\Pi_\rKch^\Kh,\Pi_\Kch^\rKh$ in \eq{kh4eq38},
\eq{kh4eq40} to map $[X,\bs f,\bs C]\mapsto[X,\bs f,\bs G_{\bs C}]$
exactly as for $\Pi_\rKch^\Kh$ in Definition \ref{kh4def7}. The
proof that $\pd\ci\Pi_\Kch^\rKh=\Pi_\Kch^\rKh\ci\d$ then becomes
nontrivial: we have
\begin{align*}
&\pd\ci\Pi_\Kch^\rKh:[X,\bs f,\bs C]\longmapsto\bigl[\pd X,\bs
f\vert_{\pd X},\bs G_{\bs C}\vert_{\pd X}\bigr],\qquad\text{and}\\
&\Pi_\Kch^\rKh\ci\d:[X,\bs f,\bs C]\longmapsto\bigl[\pd_+^{\bs
f}X,\bs f_+,\bs G_{\bs C\vert_{\smash{\pd_+^{\bs f}X}}}\bigr].
\end{align*}
The difference between these two is $\smash{\bigl[\pd_-^{\bs f}X,\bs
f\vert_{\pd_-^{\bs f}X},\bs G_{\bs C}\vert_{\pd_-^{\bs f}X}\bigr]}$,
but $\bs f\vert_{\pd_-^{\bs f}X}=\io\ci\bs f_-$, so this difference
lies in $\io_*\bigl(KC_k(\pd Y;R)\bigr)$, and is zero in
$KC_{n-k-1}(Y,\io(\pd Y);R)$, giving $\pd\ci\Pi_\Kch^\rKh=
\Pi_\Kch^\rKh\ci\d$. To prove that $\pd\ci\Pi_\rKch^\Kh=
\Pi_\rKch^\Kh\ci\d$, we show that $\pd\ci\Pi_\rKch^\Kh-
\Pi_\rKch^\Kh\ci\d=(-1)^k\io_*\ci\io^*$, and then note that
$KC^k(Y,\io(\pd Y);R)$ is the kernel of $\io^*$. Hence
$\Pi_\rKch^\Kh,\Pi_\Kch^\rKh$ are well-defined chain maps, and
induce morphisms $\Pi_\rKch^\Kh:KH^k(Y,\io(\pd Y);R)\ra
KH_{n-k}(Y;R)$ and~$\Pi_\Kch^\rKh:KH^k(Y;R)\ra KH_{n-k}(Y,\io(\pd
Y);R)$.

We can also define $\Pi^\rKch_\Kh$ in \eq{kh4eq39} generalizing
$\Pi^\Kch_\Kh$ in Definition \ref{kh4def7}, though this is a little
more complicated, since one must also add in a correction term over
$\pd Y$ to ensure that $\Pi^\Kch_\Kh$ maps to $\Ker\io^*$ in
$KC^k(Y;R)$. It is a chain map, and induces $\Pi^\rKch_\Kh:
KH_{n-k}(Y;R)\longra KH^k(Y,\io(\pd Y);R)$. We can then prove the
analogue of Theorem \ref{kh4thm3}, that $\Pi_\rKch^\Kh:
KH^k(Y,\io(\pd Y);R)\ra KH_{n-k}(Y;R)$ and $\Pi^\rKch_\Kh:
KH_{n-k}(Y;R)\longra KH^k(Y,\io(\pd Y);R)$ are inverse, and so are
both isomorphisms. This is not surprising, since we expect
$KH^*(Y,\io(\pd Y);R)\cong KH^*(Y^\ci;R)$ and $KH_*(Y;R)\cong
KH_*(Y^\ci;R)$ as for \eq{kh4eq30} above, and this analogue would
follow immediately from Theorem~\ref{kh4thm3}.

However, the author does not have a good definition for
$\Pi^\Kch_\rKh$ in \eq{kh4eq41}, on the (co)chain level. Instead, we
consider the commutative diagram
\e
\begin{gathered}
\!\!\!\!\text{\begin{small}$\displaystyle \xymatrix@C=13pt@R=15pt{
\cdots \ar[r] & KH^{k-1}(\pd Y;R)
\ar[d]^(0.4){\Pi_\Kch^\Kh}_(0.4)\cong \ar[r] & KH^k(Y,\pd Y;R)
\ar[r] \ar[d]^(0.4){\Pi_\rKch^\Kh} & KH^k(Y;R) \ar[r]
\ar[d]^(0.4){\Pi_\Kch^\rKh}_(0.4)\cong & \cdots \\
\cdots \ar[r] & KH_{n-k}(\pd Y;R) \ar[r] & KH_{n-k}(Y;R) \ar[r] &
KH_{n-k}(Y,\pd Y;R) \ar[r] & \cdots,}$\end{small}}
\!\!\!\!\!\!\!\!\!\!\!\!\!\!\!\!\!\!
\end{gathered}
\label{kh4eq42}
\e
with exact rows. Here the left hand column is an isomorphism by
Theorem \ref{kh4thm3}, as $\pd Y$ is an oriented $(n-1)$-manifold
without boundary, and the right hand column is an isomorphism by the
analogue of Theorem \ref{kh4thm3} discussed above. Thus the central
column is also an isomorphism, by properties of long exact
sequences. Therefore we have proved:

\begin{thm} Suppose $Y$ is an oriented\/ $n$-orbifold with\/
$\pd^2Y=\es,$ and\/ $R$ is a $\Q$-algebra. Then we have isomorphisms
$\Pi_\rKch^\Kh:KH^k(Y,\io(\pd Y);R)\ra KH_{n-k}(Y;R)$ and\/
$\Pi_\Kch^\rKh:KH^k(Y;R)\ra KC_{n-k}(Y,\io(\pd Y);R)$.
\label{kh4thm4}
\end{thm}

The author expects the theorem holds without assuming $\pd^2Y=\es$.
The proof for $\Pi_\rKch^\Kh$ looks straightforward, but for
$\Pi_\Kch^\rKh$ the long exact sequence argument in \eq{kh4eq42}
will need modification. Here is the analogue of
Definition~\ref{kh4def8}:

\begin{dfn} Let $Y$ be an $n$-orbifold, and $R$ a
$\Q$-algebra. We will define $R$-module morphisms
$\Pi_\rcs^\rKch:H^k_\cs(Y,\io(\pd Y);R)\ra KH^k(Y,\io(\pd Y);R)$ and
$\Pi_\cs^\Kch:H^k_\cs(Y;R)\ra KH^k(Y;R)$ for all $k\ge 0$. First
suppose $Y$ is orientable, and choose an orientation. Define
$\Pi_\rcs^\rKch$ and $\Pi_\cs^\rKch$ to be the compositions
\begin{align*}
\xymatrix@C=15pt{ H^k_\cs(Y,\io(\pd Y);R) \ar[r]^(0.55){\Pd'} &
H_{n-k}^\rsi(Y;R) \ar[r]^{\Pi_\rsi^\Kh} & KH_{n-k}(Y;R)
\ar[r]^{\Pi^\rKch_\Kh} & KH^k(Y,\io(\pd Y);R),} \\
\xymatrix@C=15pt{H^k_\cs(Y;R) \ar[r]^(0.4){\Pd''} &
H_{n-k}^\rsi(Y,\io(\pd Y);R) \ar[r]^{\Pi_\rrsi^\rKh} &
KH_{n-k}(Y,\io(\pd Y);R) \ar[r]^(0.58){\Pi^\Kch_\rKh} & KH^k(Y;R). }
\end{align*}
Changing the orientation of $Y$ changes the signs of
$\Pd'\Pd'',\Pi^\rKch_\Kh,\Pi^\Kch_\rKh$. Thus $\Pi_\rcs^\rKch$ and
$\Pi_\cs^\rKch$ are independent of the choice of orientation on $Y$.
If $Y$ is not orientable, we define $\Pi_\rcs^\rKch,\Pi_\cs^\rKch$
using {\it twisted\/} (relative) homology groups, by the method of
Definition~\ref{kh4def8}.
\label{kh4def11}
\end{dfn}

Poincar\'e duality, Corollaries \ref{kh4cor1} and \ref{kh4cor3}, and
Theorem \ref{kh4thm4} imply:

\begin{cor} Let\/ $Y$ be an orbifold with\/ $\pd^2Y=\es,$ and\/ $R$ a
$\Q$-algebra. Then $\Pi_\rcs^\rKch:H^k_\cs(Y,\io(\pd Y);R)\ra
KH^k(Y,\io(\pd Y);R)$ and\/ $\Pi_\cs^\Kch:H^*_\cs(Y;R)\ra KH^*(Y;R)$
are isomorphisms.
\label{kh4cor4}
\end{cor}

The author expects the corollary holds without assuming
$\pd^2Y=\es$.

The other material we give on Kuranishi cochains $KC^*(Y;R)$ and
Kuranishi cohomology $KH^*(Y;R)$ when $\pd Y=\es$, such as pullbacks
$h^*$ in \S\ref{kh44}, effective Kuranishi cohomology in
\S\ref{kh46}, and cup and cap products in \S\ref{kh47}, all has
generalizations to the case $\pd Y\ne\es$ using the ideas of this
section. We leave this as an exercise for the reader.

\subsection{Effective Kuranishi cohomology}
\label{kh46}

We extend \S\ref{kh44} to {\it effective\/} Kuranishi cohomology.

\begin{dfn} Let $Y$ be an orbifold without boundary. Consider
triples $(X,\bs f,\ubC)$, where $X$ is a compact Kuranishi space,
$\bs f:X\ra Y$ is a strong submersion with $(X,\bs f)$ coeffective
and cooriented, and $\ubC$ is effective co-gauge-fixing data for
$(X,\bs f)$. Write $[X,\bs f,\ubC]$ for the isomorphism class of
$(X,\bs f,\ubC)$ under coorientation-preserving isomorphisms $(\bs
a,\bs b):(X,\bs f,\ubC)\ra(\ti X,\bs{\ti f},\ubtC)$.

Let $R$ be a commutative ring. For each $k\in\Z$, define
$KC^k_\ec(Y;R)$ to be the $R$-module of finite $R$-linear
combinations of isomorphism classes $[X,\bs f,\ubC]$ as above for
which $\vdim X=\dim Y-k$, with the analogues of relations Definition
\ref{kh4def2}(i)--(iii), replacing gauge-fixing data $\bs G$ by
effective co-gauge-fixing data $\ubC$. Elements of $KC^k_\ec(Y;R)$
are called {\it effective Kuranishi cochains}. Define
$\d:KC^k_\ec(Y;R)\ra KC^{k+1}_\ec(Y;R)$ by \eq{kh4eq20}, replacing
$\bs C_a$ by $\ubC_a$. Then $\d\ci\d=0$. Define the {\it effective
Kuranishi cohomology groups} $KH^k_\ec(Y;R)$ of $Y$ to be
\begin{equation*}
KH_\ec^k(Y;R)=\frac{\Ker\bigl(\d:KC_\ec^k(Y;R)\ra
KC_\ec^{k+1}(Y;R)\bigr)}{\Im\bigl(\d:KC_\ec^{k-1}(Y;R)\ra
KC_\ec^k(Y;R)\bigr)}
\end{equation*}
for $k\in\Z$. Arguably one should call this `coeffective Kuranishi
cohomology', as it is defined using $[X,\bs f,\ubC]$ with $(X,\bs
f)$ coeffective, but we choose not to.

As for $\Pi_\ef^\Kh$ in Definition \ref{kh4def4}, define
\e
\begin{split}
\Pi_\ec^\Kch:KC^k_\ec(Y;R)&\longra KC^k(Y;R\ot_\Z\Q)\quad\text{by}\\
\Pi_\ec^\Kch:\ts\sum_{a\in A}\rho_a\bigl[X_a,\bs f_a,\ubC_a\bigr]
&\longmapsto \ts\sum_{a\in A}\pi(\rho_a)\bigl[X_a,\bs
f_a,\Pi(\ubC_a) \bigr].
\end{split}
\label{kh4eq43}
\e
Then $\Pi_\ec^\Kch\ci\d=\d\ci\Pi_\ec^\Kch$, so they induce morphisms
of cohomology groups
\e
\Pi_\ec^\Kch:KH^k_\ec(Y;R)\longra KH^k(Y;R\ot_\Z\Q).
\label{kh4eq44}
\e
Let $Y,Z$ be orbifolds without boundary, and $h:Y\ra Z$ be a smooth,
proper map. As in Definition \ref{kh4def6}, define the {\it
pullback\/} $h^*:KC^k_\ec(Z;R)\ra KC^k_\ec(Y;R)$ by \eq{kh4eq21},
with $\ubC_a$ in place of $\bs C_a$. These satisfy $h^*\ci\d=\d\ci
h^*$, and so induce morphisms $h^*:KH^k_\ec(Z;R)\ra KH^k_\ec(Y;R)$.
We have $(g\ci h)^*=h^*\ci g^*$ and $\Pi_\ec^\Kch\ci
h^*=h^*\ci\Pi_\ec^\Kch$, on both cochains and cohomology.
\label{kh4def12}
\end{dfn}

We can also try to extend Definition \ref{kh4def7} and Theorem
\ref{kh4thm3} to effective Kuranishi (co)homology. However, as in
\S\ref{kh41} Poincar\'e duality fails for general orbifolds $Y$ and
commutative rings $R$, so something must go wrong. The problem is in
generalizing the morphisms $\Pi^\Kh_\Kch:KC^k(Y;R)\ra KC_{n-k}(Y;R)$
and $\Pi^\Kch_\Kh:KC_{n-k}(Y;R)\ra KC^k(Y;R)$ to the effective case.
We can define $\Pi^\ef_\ec:KC^k_\ec(Y;R)\ra KC_{n-k}^\ef(Y;R)$ only
if $Y$ is an {\it effective} orbifold, as otherwise
$(V^i,\ldots,\psi^i), f^i$ coeffective does not imply
$(V^i,\ldots,\psi^i)$ effective. This corresponds to the observation
in \S\ref{kh41} that if $Y$ is not effective then we cannot define
the Poincar\'e duality map $\Pd$ in \eq{kh4eq3}, which is the
analogue of~$\Pi^\ef_\ec$.

Worse still, we can define $\Pi^\ec_\ef:KC_{n-k}^\ef(Y;R)\ra
KC^k_\ec(Y;R)$ only if $Y$ is a {\it manifold}, since otherwise the
Kuranishi neighbourhoods \eq{kh4eq23} may not be coeffective, as the
stabilizer groups of $V$ can act nontrivially on the fibres of $E$.
This corresponds to the observation in \S\ref{kh41} that if $Y$ is
an effective orbifold, but not a manifold, then $\Pd$ in \eq{kh4eq3}
is defined but may not be invertible, so we cannot define its
inverse $\Pd^{-1}$, which is the analogue of~$\Pi^\ec_\ef$.

\begin{dfn} Let $Y$ be an {\it effective\/} orbifold without
boundary, which is oriented of dimension $n$, and $R$ a commutative
ring. Choose an injective map $\uG_Y:Y\ra\R^m\subset\uP$ for some
$m\gg 0$. Suppose $X$ is a compact, oriented Kuranishi space, $\bs
f:X\ra Y$ is a strong submersion, $(X,\bs f)$ is coeffective and
cooriented, and $\ubC$ is effective co-gauge-fixing data for $(X,\bs
f)$. Define effective gauge-fixing data $\ubG_{\ubC}$ for $(X,\bs
f)$ as for $\bs G_{\bs C}$ in Definition \ref{kh4def7}, but with
$\uG^i=\umu\ci\bigl((\uG_Y\ci f^i\ci\pi^i)\t\uC^i\bigr)$, for $\umu$
as in \eq{kh3eq28}. A similar proof to those for $h^*(\ubC)$ and
$\ubC\t_Y\ubtC$ in Definitions \ref{kh3def21} and \ref{kh3def22}
shows that $\ubG_{\ubC}$ is {\it effective gauge-fixing data}. To
verify Definition \ref{kh3def17}(a), note that as
$(V^i,\ldots,\psi^i),f^i$ in $\ubC$ is coeffective and $Y$ is
effective, $(V^i,\ldots,\psi^i)$ is effective as in
Definition~\ref{kh3def18}.

Define $R$-module morphisms $\Pi_\ec^\ef:KC^k_\ec(Y;R)\ra
KC_{n-k}^\ef(Y;R)$ for $k\in\Z$ by $\Pi_\ec^\ef:[X,\bs
f,\ubC]\mapsto [X,\bs f,\ubG_{\ubC}]$, using the coorientation for
$(X,\bs f)$ and orientation on $Y$ to determine the orientation on
$X$. Then $\pd\ci\Pi_\ec^\ef=\Pi_\ec^\ef\ci\d$, so they induce
morphisms of (co)homology groups
\e
\Pi_\ec^\ef:KH^k_\ec(Y;R)\longra KH_{n-k}^\ef(Y;R).
\label{kh4eq45}
\e

Now suppose $Y$ is a {\it manifold}, without boundary. Let $X$ be a
compact, effective, oriented Kuranishi space, $\bs f:X\ra Y$ a
strongly smooth map, and $\ubG$ be effective gauge-fixing data for
$(X,\bs f)$. Define $X^Y,\bs f{}^Y$ as in Definition \ref{kh4def7}.
As $Y$ is a manifold, one can show from \eq{kh4eq23} that $(X^Y,\bs
f{}^Y)$ is coeffective. Define $\ubC_\ubG^Y$ for $(X^Y,\bs f{}^Y)$
as for $\bs C_{\bs G}^Y$ in Definition \ref{kh4def7}, except that we
set $\ucC^{k,Y}=\umu\ci\bigl((\uC_Y\ci\pi_E)\t(\uG^{k-n}
\ci\pi_{E^{k-n}})\bigr)$, where $\uC_Y:E\ra\R^{m'}\subset\uP$ is
injective for some $m'\gg 0$. Again, we can check that $\ubC_\ubG^Y$
is {\it effective co-gauge-fixing data}.

Define $R$-module morphisms $\Pi^\ec_\ef:KC_{n-k}^\ef(Y;R)\ra
KC^k_\ec(Y;R)$ by $\Pi^\ec_\ef:[X,\bs f,\ubG]\mapsto [X^Y,\bs
f{}^Y,\ubC_\ubG^Y]$. Then $\d\ci\Pi^\ec_\ef=\Pi^\ec_\ef\ci\pd$, so
they induce morphisms of (co)homology groups
\e
\smash{\Pi^\ec_\ef:KH_{n-k}^\ef(Y;R)\longra KH^k_\ec(Y;R).}
\label{kh4eq46}
\e
\label{kh4def13}
\end{dfn}

Theorem \ref{kh4thm3}, Definition \ref{kh4def8} and Corollary
\ref{kh4cor2} then generalize as follows.

\begin{thm} Let\/ $Y$ be an oriented\/ $n$-manifold without
boundary, and\/ $R$ a commutative ring. Then $\Pi^\ec_\ef$ in
\eq{kh4eq46} is the inverse of\/ $\Pi_\ec^\ef$ in {\rm\eq{kh4eq45},}
so they are both isomorphisms. Also \eq{kh4eq45}--\eq{kh4eq46} are
independent of choices.
\label{kh4thm5}
\end{thm}

\begin{dfn} Let $Y$ be an $n$-manifold without boundary, and $R$ a
commutative ring. We will define $R$-module morphisms
$\Pi_\cs^\ec:H^k_\cs(Y;R)\ra KH^k_\ec(Y;R)$ for all $k\ge 0$. First
suppose $Y$ is oriented. Define $\Pi_\cs^\ec$ to be the composition
\begin{equation*}
\smash{\xymatrix@C=30pt{ H^k_\cs(Y;R) \ar[r]^(0.4){\Pd} &
H_{n-k}^\rsi(Y;R) \ar[r]^{\Pi_\rsi^\ef} & KH_{n-k}^\ef(Y;R)
\ar[r]^{\Pi^\ec_\ef} & KH^k_\ec(Y;R), }}
\end{equation*}
where $\Pd,\Pi_\rsi^\ef,\Pi^\ec_\ef$ are as in \eq{kh4eq3},
\eq{kh4eq18}, \eq{kh4eq46}. Independence of the orientation of $Y$,
and extension to non-orientable $Y$, are as in
Definition~\ref{kh4def8}.
\label{kh4def14}
\end{dfn}

\begin{cor} Let\/ $Y$ be a manifold without boundary, and\/ $R$ a
commutative ring. Then $\Pi_\cs^\ec:H^*_\cs(Y;R)\ra KH^*_\ec(Y;R)$
is an isomorphism.
\label{kh4cor5}
\end{cor}

This still leaves the question of whether $H^*_\cs(Y;R)\cong
KH^*_\ec(Y;R)$ for $Y$ a general orbifold without boundary, and $R$
a commutative ring. The author expects this is true, even though our
proof via homology and Poincar\'e duality fails in this case. When
$R$ is a $\Q$-algebra, one approach might be to generalize the proof
of Theorem \ref{kh4thm2} in Appendix \ref{khC} to show that
$\Pi_\ec^\Kch:KH^k_\ec(Y;R)\ra KH^k(Y;R)$ in \eq{kh4eq44} is an
isomorphism, and then use Corollary~\ref{kh4cor2}.

\subsection{Products on (effective) Kuranishi (co)homology}
\label{kh47}

There are {\it cup products} $\cup$, {\it cap products} $\cap$ and
{\it intersection products} $\bu$ on (effective) Kuranishi
(co)chains and (co)homology.

\begin{dfn} Let $Y$ be an orbifold without boundary, and $R$ a
$\Q$-algebra. Define the {\it cup product\/} $\cup:KC^k(Y;R)\t
KC^l(Y;R)\ra KC^{k+l}(Y;R)$ by
\e
[X,\bs f,\bs C]\cup[\ti X,\bs{\ti f},\bs{\ti C}]=
\smash{\bigl[X\t_{\bs f,Y,\bs{\ti f}}\ti X,\bs\pi_Y,\bs C\t_Y
\bs{\ti C}\bigr]},
\label{kh4eq47}
\e
extended $R$-bilinearly. Here $\bs\pi_Y:X\t_{\smash{\bs f,Y,\bs{\ti
f}}}\ti X\ra Y$ is the projection from the fibre product, which is a
strong submersion as both $\bs f,\bs{\ti f}$ are strong submersions,
and $\bs C\t_Y\bs{\ti C}$ is as in Definition \ref{kh3def15}. As for
$\d$, it is easy to show $\cup$ takes relations (i)--(iv) in both
$KC^k(Y;R)$ and $KC^l(Y;R)$ to the same relations in
$KC^{k+l}(Y;R)$. Thus $\cup$ is well-defined.

Proposition \ref{kh3prop6}(a)--(c) and \eq{kh2eq21}--\eq{kh2eq23}
imply that for $\ga\in KC^k(Y;R)$, $\de\in KC^l(Y;R)$ and $\ep\in
KC^m(Y;R)$ we have
\begin{gather}
\ga\cup\de=(-1)^{kl}\de\cup\ga,
\label{kh4eq48}\\
\d(\ga\cup\de)=(\d\ga)\cup\de+(-1)^k\ga\cup(\d\de)\;\>\text{and}\;\>
(\ga\cup\de)\cup\ep=\ga\cup(\de\cup\ep).
\label{kh4eq49}
\end{gather}
Therefore in the usual way $\cup$ induces an associative,
supercommutative product $\cup:KH^k(Y;R)\t KH^l(Y;R)\ra
KH^{k+l}(Y;R)$ given for $\ga\in KC^k(Y;R)$ and $\de\in KC^l(Y;R)$
with $\d\ga=\d\de=0$ by
\begin{equation*}
(\ga+\Im\d_{k-1})\cup(\de+\Im\d_{l-1})=(\ga\cup\de)+\Im\d_{k+l-1}.
\end{equation*}

Now suppose $Y$ is compact, and of dimension $n$. Then $\id_Y:Y\ra
Y$ is a (strong) submersion, with a trivial coorientation giving the
positive sign to the zero vector bundle over $Y$. Define
co-gauge-fixing data $\bs C_Y$ for $(Y,\id_Y)$ by $I=\{n\}$,
$(V^n,E^n,s^n,\psi^n)=(Y,Y,0,\id_Y)$, $\id_Y^n=\id_Y:V^n\ra Y$, and
$\eta_n:Y\ra[0,1]$, $\eta_n^n:V^n\ra[0,1]$ and $C^n:E^n\ra P$ given
by $\eta_n\equiv 1$, $\eta_n^n\equiv 1$ and $C^n\equiv(\es)$, where
$(\es)\in P$ is the identity for $\mu$, as in Definition
\ref{kh3def6}. Then $[Y,\id_Y,\bs C_Y]\in KC^0(Y;R)$, with
$\d[Y,\id_Y,\bs C_Y]=0$. Following the definitions through shows
that for $[X,\bs f,\bs C]\in KC^k(Y;R)$ we have
\begin{equation*}
[Y,\id_Y,\bs C_Y]\cup [X,\bs f,\bs C]=[X,\bs f,\bs
C]\cup[Y,\id_Y,\bs C_Y]=[X,\bs f,\bs C].
\end{equation*}
Thus $[Y,\id_Y,\bs C_Y]$ is the identity for $\cup$, at the cochain
level. Passing to cohomology, $\bigl[[Y,\id_Y,\bs C_Y]\bigr]\in
KH^0(Y;R)$ is the identity for $\cup$ in $KH^*(Y;R)$. We call
$\bigl[[Y,\id_Y,\bs C_Y]\bigr]$ the {\it fundamental class} of~$Y$.

Define the {\it cap product\/} $\cap:KC_k(Y;R)\t KC^l(Y;R)\ra
KC_{k-l}(Y;R)$ by
\e
\smash{[X,\bs f,\bs G]\cap[\ti X,\bs{\ti f},\bs{\ti C}]=\bigl[X
\t_{\bs f,Y,\bs{\ti f}}\ti X,\bs\pi_Y,\bs G\t_Y\bs{\ti C}\bigr]},
\label{kh4eq50}
\e
extended $R$-bilinearly, with $\bs G\t_Y\bs{\ti C}$ as in Definition
\ref{kh3def15}. For $\ga\in KC_k(Y;R)$ and $\de,\ep\in KC^*(Y;R)$,
the analogue of \eq{kh4eq49}, using Proposition \ref{kh3prop6}(d),
is
\e
\pd(\ga\cap\de)=(\pd\ga)\cap\de+(-1)^{\dim Y-k}\ga\cup(\d\de),
\quad(\ga\cap\de)\cap\ep=\ga\cap(\de\cup\ep).
\label{kh4eq51}
\e
If also $Y$ is compact then $\ga\cap[Y,\id_Y,\bs C_Y]=\ga$. Thus
$\cap$ induces a {\it cap product\/} $\cap:KH_k(Y;R)\t KH^l(Y;R)\ra
KH_{k-l}(Y;R)$. These products $\cap$ make Kuranishi chains and
homology into {\it modules} over Kuranishi cochains and cohomology.

Let $Y,Z$ be orbifolds without boundary, and $h:Y\ra Z$ a smooth,
proper map. Then Proposition \ref{kh3prop7} implies that pullbacks
$h^*$ and pushforwards $h_*$ are compatible with $\cup,\cap$ on
(co)chains, in the sense that if $\al\in KC_*(Y;R)$ and $\be,\ga\in
KC^*(Z;R)$ then
\e
h^*(\be\cup\ga)=h^*(\be)\cup h^*(\ga) \quad\text{and}\quad
h_*(\al\cap h^*(\be))=h_*(\al)\cap\be.
\label{kh4eq52}
\e
Since $\cup,\cap,h^*,h_*$ are compatible with $\d,\pd$, passing to
(co)homology shows that \eq{kh4eq52} also holds for $\al\in
KH_*(Y;R)$ and $\be,\ga\in KH^*(Z;R)$. If $Z$ is compact then $Y$
is, with $h^*\bigl([Z,\id_Z,\bs C_Z]\bigr)=[Y,\id_Y,\bs C_Y]$ in
$KC^0(Y;R)$ and $h^*\bigl(\bigl[[Z,\id_Z,\bs C_Z]\bigr]\bigr)\ab=
\bigl[[Y,\id_Y,\bs C_Y]\bigr]$ in~$KH^0(Y;R)$.

To summarize: Kuranishi cochains $KC^*(Y;R)$ form a {\it
supercommutative, associative, differential graded\/ $R$-algebra},
and Kuranishi cohomology $KH^*(Y;R)$ is a {\it supercommutative,
associative, graded\/ $R$-algebra}. These algebras are {\it with
identity} if $Y$ is compact without boundary, and {\it without
identity} otherwise. Pullbacks $h^*$ induce {\it algebra morphisms}
on both cochains and cohomology. Kuranishi chains $KC_*(Y;R)$ are a
{\it graded module} over $KC^*(Y;R)$, and Kuranishi homology
$KH_*(Y;R)$ is a {\it graded module} over~$KH^*(Y;R)$.

Let $Y$ be oriented, of dimension $n$. Then Poincar\'e duality,
Theorem \ref{kh4thm3}, and the cup product $\cup$ on Kuranishi
cohomology induce an {\it intersection product\/} $\bu:KH_k(Y;R)\t
KH_l(Y;R)\ra KH_{k+l-n}(Y;R)$ on Kuranishi homology. It is
supercommutative (with degrees shifted by $n$) and associative, with
identity in $KH_n(Y;R)$ if $Y$ is compact. Intersection products are
{\it not\/} preserved by pushforwards $h_*:KH_*(Y;R)\ra KH_*(Z;R)$,
the degrees of $h_*(\al)\bu h_*(\be)$ and $h_*(\al\bu\be)$ differ
by~$\dim Y-\dim Z$.

On chains we can define $\bu:KC_k(Y;R)\t KC_l(Y;R)\ra
KC_{k+l-n}(Y;R)$ by
\e
\al\bu\be=\Pi_\Kch^\Kh\bigl(\Pi_\Kh^\Kch(\al)\cup\Pi_\Kh^\Kch
(\be)\bigr),
\label{kh4eq53}
\e
where $\Pi_\Kch^\Kh,\Pi_\Kh^\Kch$ are as in Definition
\ref{kh4def7}. Then $\bu$ on chains induces $\bu$ on homology, and
$\bu$ on chains is supercommutative, but {\it not\/} associative,
and does not have an identity. Thus, the intersection product $\bu$
is not as well-behaved on Kuranishi chains as the cup and cap
products~$\cup,\cap$.
\label{kh4def15}
\end{dfn}

\begin{rem} It is an attractive and unusual feature of our theory of
Kuranishi (co)homology that the products $\cup,\cap$ are everywhere
defined, and supercommutative, and associative, {\it at the
(co)chain level}, as well as on (co)homology. Furthermore,
pushforwards and pullbacks are functorial and compatible with
$\cup,\cap$ at the (co)chain level, and for compact $Y$ there is an
identity $[Y,\id_Y,\bs C_Y]$ for $\cup,\cap$ at the cochain level.
There are few homology or cohomology theories for which this holds
--- it works for de Rham cohomology, for instance, but fails for
singular and simplicial (co)homology.

The fact that products $\cup,\cap$ are {\it everywhere defined\/} on
(co)chains, rather than only under some transversality assumptions
such as transverse intersection of chains in $Y$, is a feature of
working with Kuranishi spaces rather than manifolds or orbifolds. If
$X$ is a manifold or orbifold then submersions $f:X\ra Y$ only exist
if $\dim X\ge\dim Y$. So for defining (co)homology using pairs
$(X,f)$ as (co)chains with $X$ a manifold or orbifold (for example,
singular homology, with $X=\De_k$), making $f$ a submersion is too
strong an assumption, and this leads to transversality problems when
defining products. But for Kuranishi spaces $X$ we can have strong
submersions $\bs f:X\ra Y$ with $\vdim X<\dim Y$, since the
obstruction bundles on $X$ decrease $\vdim X$. Thus we can use pairs
$(X,\bs f)$ with $\bs f$ a strong submersion as cochains, and
$\cup,\cap$ are everywhere defined.

A lot of work in Chapter \ref{kh3} went into defining
(co-)gauge-fixing data with associative, supercommutative products.
This good behaviour of Kuranishi cochains will be very important,
and lead to huge simplifications, in the sequels
\cite{AkJo,Joyc2,Joyc3}, where we have to work at the cochain level.
\label{kh4rem5}
\end{rem}

\begin{dfn} Let $Y$ be an orbifold without boundary, and $R$ a
commutative ring. For cup and cap products on {\it effective\/}
Kuranishi (co)chains and (co)homology, we follow Definition
\ref{kh4def15}, except that we use the fibre products
$\ubC\t_Y\ubtC$, $\ubG\t_Y\ubtC$ of effective (co-)gauge-fixing data
from Definition \ref{kh3def22}. Thus we define $\cup:KC^k_\ec(Y;R)\t
KC^l_\ec(Y;R)\ra KC^{k+l}_\ec(Y;R)$ by
\begin{equation*}
\smash{[X,\bs f,\ubC]\cup[\ti X,\bs{\ti f},\ubtC]=\bigl[X \t_{\bs
f,Y,\bs{\ti f}}\ti X,\bs\pi_Y,\ubC\t_Y\ubtC\bigr]}.
\end{equation*}
Then $\cup$ satisfies \eq{kh4eq49} by Propositions
\ref{kh3prop6}(b),(c) and \ref{kh3prop8}, but it does not satisfy
\eq{kh4eq48}, since fibre products of effective co-gauge-fixing data
are not commutative. Hence $\cup$ induces an associative cup product
\e
\cup:KH^k_\ec(Y;R)\t KH^l_\ec(Y;R)\longra KH^{k+l}_\ec(Y;R).
\label{kh4eq54}
\e

Suppose $Y$ is compact of dimension $n$. As for $\bs C_Y$ in
Definition \ref{kh4def15}, define effective co-gauge-fixing data
$\ubC_Y$ for $(Y,\id_Y)$ by $I=\{n\}$,
$(V^n,E^n,s^n,\psi^n)=(Y,Y,0,\id_Y)$, $\id_Y^n=\id_Y:V^n\ra Y$, and
$\eta_n:Y\ra[0,1]$, $\eta_n^n:V^n\ra[0,1]$ and $\uC^n:E^n\ra\uP$
given by $\eta_n\equiv 1$, $\eta_n^n\equiv 1$ and
$\uC^n\equiv(\es)$, where $(\es)\in\uP$ is the identity for $\umu$,
as in Definition \ref{kh3def17}. Then $[Y,\id_Y,\ubC_Y]\in
KC^0_\ec(Y;R)$, with $\d[Y,\id_Y,\ubC_Y]=0$, and $[Y,\id_Y,\ubC_Y]$
is the identity for $\cup$, at the cochain level. Passing to
cohomology, $\bigl[[Y,\id_Y,\ubC_Y]\bigr]\in KH^0_\ec(Y;R)$ is the
identity for $\cup$ in $KH^*_\ec(Y;R)$. We call
$\bigl[[Y,\id_Y,\ubC_Y]\bigr]$ the {\it fundamental class\/} of~$Y$.

Define the {\it cap product\/} $\cap:KC_k^\ef(Y;R)\t
KC^l_\ec(Y;R)\ra KC_{k-l}^\ef(Y;R)$ by
\begin{equation*}
\smash{[X,\bs f,\ubG]\cap[\ti X,\bs{\ti f},\ubtC]=\bigl[X \t_{\bs
f,Y,\bs{\ti f}}\ti X,\bs\pi_Y,\ubG\t_Y\ubtC\bigr]}.
\end{equation*}
It satisfies \eq{kh4eq51}, by Propositions \ref{kh3prop6}(d) and
\ref{kh3prop8}. Thus $\cap$ induces a {\it cap product\/}
\e
\cap:KH_k^\ef(Y;R)\t KH^l_\ec(Y;R)\longra KH_{k-l}^\ef(Y;R).
\label{kh4eq55}
\e
If $Y$ is compact then $\al\cap[Y,\id_Y,\ubC_Y]=\al$ for all $\al$
in $KC_*^\ef(Y;R)$, and $\be\cap\bigl[[Y,\id_Y,\ubC_Y]\bigr]=\be$
for all $\be$ in~$KH_*^\ef(Y;R)$.

Let $Y,Z$ be orbifolds without boundary, and $h:Y\ra Z$ a smooth,
proper map. Then Propositions \ref{kh3prop7} and \ref{kh3prop8}
imply \eq{kh4eq52} holds for $\al\in KC_*^\ef(Y;R)$ and $\be,\ga\in
KC^*_\ec(Z;R)$, and hence for $\al\in KH_*^\ef(Y;R)$ and $\be,\ga\in
KH^*_\ec(Z;R)$. If $Z$ is compact then $Y$ is, with
$h^*\bigl([Z,\id_Z,\ubC_Z]\bigr)=[Y,\id_Y,\ubC_Y]$ in
$KC^0_\ec(Y;R)$ and $h^*\bigl(\bigl[[Z,\id_Z,\ubC_Z]\bigr]\bigr)\ab=
\bigl[[Y,\id_Y,\ubC_Y]\bigr]$ in~$KH^0_\ec(Y;R)$.

Let $Y$ be an oriented $n$-manifold without boundary. Then
Poincar\'e duality, Theorem \ref{kh4thm5}, and $\cup$ in
\eq{kh4eq54} induce an associative {\it intersection product\/}
$\bu:KH^\ef_k(Y;R)\t KH^\ef_l(Y;R)\ra KH^\ef_{k+l-n}(Y;R)$ on
effective Kuranishi homology. On chains we can define
$\bu:KC_k^\ef(Y;R)\t KC_l^\ef(Y;R)\ra KC_{k+l-n}^\ef(Y;R)$ by
\e
\al\bu\be=\Pi_\ec^\ef\bigl(\Pi_\ef^\ec(\al)\cup\Pi_\ef^\ec
(\be)\bigr).
\label{kh4eq56}
\e
This is neither associative nor supercommutative, at the chain
level.

Since $\mu:P\t P\ra P$, $\umu:\uP\t\uP\ra\uP$ and $\Pi:\uP\ra P$ in
Definitions \ref{kh3def6} and \ref{kh3def17} satisfy
$\Pi\bigl(\umu(p,q)\bigr)=\mu\bigl(\Pi(p),\Pi(q)\bigr)$ for all
$p,q\in\uP$, it easily follows that the morphisms
$\Pi_\ef^\Kh,\Pi_\ec^\Kch$ from effective Kuranishi (co)chains to
Kuranishi (co)chains in \eq{kh4eq9} and \eq{kh4eq43} satisfy
\e
\begin{gathered}
\Pi_\ec^\Kch(\ga)\cup\Pi_\ec^\Kch(\de)=\Pi_\ec^\Kch(\ga\cup\de),
\qquad
\Pi_\ef^\Kh(\al)\cap\Pi_\ec^\Kch(\ga)=\Pi_\ef^\Kh(\al\cap\ga)\\
\text{and}\qquad
\Pi_\ef^\Kh(\al)\bu\Pi_\ef^\Kh(\be)=\Pi_\ef^\Kh(\al\bu\be),
\end{gathered}
\label{kh4eq57}
\e
for all $\al,\be\in KC_*^\ef(Y;R)$ and $\ga,\de\in KC^*_\ec(Y;R)$,
supposing $Y$ oriented for the third equation. Passing to
(co)homology, the morphisms $\Pi_\ef^\Kh,\Pi_\ec^\Kch$ in
\eq{kh4eq10} and \eq{kh4eq44} also satisfy \eq{kh4eq57} for all
$\al,\be\in KH_*^\ef(Y;R)$ and~$\ga,\de\in KH^*_\ec(Y;R)$.
\label{kh4def16}
\end{dfn}

Although $\cup$ on $KC^*_\ec(Y;R)$ is not supercommutative, it is up
to homotopy, and so $\cup$ on $KH^*_\ec(Y;R)$ is supercommutative.

\begin{prop} Let\/ $Y$ be an orbifold without boundary and\/ $R$ a
commutative ring. Then $\cup$ on $KH^*_\ec(Y;R)$ in \eq{kh4eq54}
satisfies {\rm\eq{kh4eq48},} that is, $\cup$ is supercommutative.
If\/ $Y$ is an oriented manifold without boundary then $\bu$ is
supercommutative on~$KH_*^\ef(Y;R)$.
\label{kh4prop}
\end{prop}

\begin{proof} Let $X,\ti X$ be compact Kuranishi spaces,
$\bs f:X\ra Y$, $\bs{\ti f}:\ti X\ra Y$ coeffective, cooriented
strong submersions, and $\ubC,\ubtC$ effective co-gauge-fixing data
for $(X,\bs f),(\ti X,\bs{\ti f})$. Then Definition \ref{kh3def22}
defines effective co-gauge-fixing data $\ubC\t_Y\ubtC$,
$\ubtC\t_Y\ubC$ for $(X\t_Y\ti X,\bs\pi_Y)$ and $(\ti
X\t_YX,\bs\pi_Y)$. Our problem is that the natural isomorphism
$X\t_Y\ti X\cong\ti X\t_Y X$ does not identify $\ubC\t_Y\ubtC$ with
$\ubtC\t_Y\ubC$, because in the definition $\ucC^k(\check
e)=\umu\bigl(\uC^i\ci\pi_{E^i}(\check e),\ti
\uC{}^{\ti\imath}\ci\pi_{\ti E^{\ti\imath}}(\check e)\bigr)$ of
$\ucC^k:\check E^k\ra\uP$ in $\ubC\t_Y\ubtC$, the function $\umu$ in
\eq{kh3eq28} is not commutative.

As in Definition \ref{kh3def15}, write $(\bs{\check
I},\bs{\check\eta})$ for the excellent coordinate system for
$(X\t_Y\ti X,\bs\pi_Y)$ in $\ubC\t_Y\ubtC$, where $\bs{\check
I}=\bigl(\check I,(\check V^k,\check E^k,\check
s^k,\check\psi^k),\check\pi_Y^k: k\in\check I,\ldots\bigr)$. Define
an excellent coordinate system $(\bs J,\bs\ze)$ for
$\bigl([0,1]\t(X\t_Y\ti X),\bs\pi_Y\bigr)$ to have indexing set
$J=\{k+1:k\in\check I\}$, Kuranishi neighbourhoods
\begin{equation*}
(\ti V^{k+1},\ti E^{k+1},\ti s^{k+1},\ti\psi^{k+1})= \bigl([0,1]\t
\check V^k,[0,1]\t\check E^k,\id_{[0,1]}\t\check
s^k,\id_{[0,1]}\t\check\psi^k\bigr)
\end{equation*}
for $k\in\check I$, maps $\ti\pi_Y^{k+1}:\ti V^{k+1}\ra Y$,
$\ze_{k+1}:[0,1]\t(X\t_Y\ti X)\ra [0,1]$ and $\ze_{k+1}^{l+1}:\ti
V^{l+1}\ra[0,1]$ given by
$\ti\pi_Y^{k+1}=\check\pi_Y^k\ci\pi_{\check V^k}$,
$\ze_{k+1}=\check\eta_k\ci\pi_{X\t_Y\ti X}$ and
$\ze_{k+1}^{l+1}=\check\eta_k^l\ci\pi_{\check V^l}$, and triples
$(\ti V^{(k+1)(l+1)},\ti\phi^{(k+1)(l+1)},\hat{\ti\phi}{}^{(k+1)
(l+1)})$ given by $\ti V^{(k+1)(l+1)}=[0,1]\t\check V^{kl}$,
$\ti\phi^{(k+1)(l+1)}=\id_{[0,1]}\t\check\phi^{kl}$
and~$\hat{\ti\phi}{}^{(k+1)(l+1)}=\id_{[0,1]}\t\hat{\check\phi}{}^{kl}$.

Define $\ul\nu:[0,1]\t\uP\t\uP\ra\uP$ by
\e
\ul\nu\bigl(t,(x_1,\ldots,x_k),(y_1,\ldots,y_l)\bigr)=
\begin{cases} (x_1,\ldots,x_k,y_1,\ldots,y_l), & t=0, \\
(t,x_1,\ldots,x_k,y_1,\ldots,y_l), & t\in(0,1), \\
(y_1,\ldots,y_l,x_1,\ldots,x_k), & t=1.
\end{cases}
\label{kh4eq58}
\e
Then $\ul\nu(0,p,q)=\umu(p,q)$ and $\ul\nu(1,p,q)=\umu(q,p)$. Define
$\uD^{k+1}:\ti E^{k+1}\ra P$ for $k\in\check I$ and $u\in[0,1]$,
$\check e\in\check E^k\cap E^{(i,\ti\imath)}$ by
\e
\begin{split}
\uD^{k+1}(u,\check e)=\ul\nu\bigl(u,\uC^i\ci\pi_{E^i}(\check e),
\utC^{\ti\imath}\ci\pi_{\ti E^{\ti\imath}}(\check e)\bigr).
\end{split}
\label{kh4eq59}
\e
A similar proof to those for $h^*(\ubC)$ and $\ubC\t_Y\ubtC$ in
Definitions \ref{kh3def21} and \ref{kh3def22} shows that
$\ubD_{\ubC, \ubtC}=\bigl((\bs J,\bs\ze),\uD^j:j\in J\bigr)$ is {\it
effective co-gauge-fixing data\/} for~$\bigl([0,1]\t(X\t_Y\ti
X),\bs\pi_Y \bigr)$.

Using \eq{kh2eq21} we see that in cooriented Kuranishi spaces we
have
\e
\begin{split}
\pd\bigl([0,1]&\t(X\t_Y\ti X),\bs\pi_Y\bigr)\cong
-\bigl(\{0\}\t(X\t_Y\ti X),\bs\pi_Y\bigr)\\
&\amalg\bigl(\{1\}\t(X\t_Y\ti X),\bs\pi_Y\bigr)\amalg
-\bigl([0,1]\t(\pd X\t_Y\ti X),\bs\pi_Y\bigr)\\
&\amalg -(-1)^{\vdim X-\dim Y}\bigl([0,1]\t(X\t_Y\pd\ti
X),\bs\pi_Y\bigr).
\end{split}
\label{kh4eq60}
\e
Let us add co-gauge-fixing data to \eq{kh4eq60}. The restriction of
$\ubD_{\ubC,\ubtC}$ to the first term on the r.h.s.\ is isomorphic
to $\ubC\t_Y\ubtC$ under $\{0\}\t(X\t_Y\ti X)\cong X\t_Y\ti X$,
since $\uD^{k+1}:[0,1]\t\check E^k\ra P$ in \eq{kh4eq59} restricts
on $\{0\}\t\check E^k\cong\check E^k$ to $\ucC^k:\check E^k\ra P$ in
the definition of $\ubC\t_Y\ubtC$, by the case $t=0$ in
\eq{kh4eq58}. The restriction of $\bs D_{\ubC,\ubtC}$ to the second
term is isomorphic to $\ubtC\t_Y\ubC$ under $\{1\}\t(X\t_Y\ti
X)\cong\ti X\t_YX$, since $\uD^{k+1}$ restricts on $\{1\}\t\check
E^k\cong\check E^k$ to $\ucC^k$ with the r\^oles of $\ubC,\ubtC$
reversed, by the case $t=1$ in \eq{kh4eq58}. In cooriented Kuranishi
spaces we have $\bigl(\{1\}\t(X\t_Y\ti X),\bs\pi_Y\bigr)\cong
(-1)^{(\vdim X-\dim Y)(\vdim\ti X-\dim Y)}(\ti X\t_YX,\bs\pi_Y)$ by
\eq{kh2eq22}. The restrictions of $\ubD_{\ubC,\ubtC}$ to the third
and fourth terms are $\ubD_{\ubC\vert_{\pd X},\ubtC}$ and
$\ubD_{\ubC,\ubtC\vert_{\pd\ti X}}$. Putting all this together, in
$KC^*_\ec(Y;R)$ we have
\e
\begin{split}
\d\bigl[[0,1]&\t(X\t_Y\ti X),\bs\pi_Y,\ubD_{\ubC,\ubtC}\bigr]=
-\bigl[X\t_Y\ti X,\bs\pi_Y,\ubC\t_Y\ubtC\bigr]\\
&+(-1)^{(\vdim X-\dim Y)(\vdim\ti X-\dim Y)}\bigl[\ti
X\t_YX,\bs\pi_Y,\ubtC\t_Y\ubC\bigr]\\
&-\bigl[[0,1]\t(\pd X\t_Y\ti X),\bs\pi_Y,\ubD_{\ubC\vert_{\pd
X},\ubtC}\bigr]\\
&-(-1)^{\vdim X-\dim Y}\bigl[[0,1]\t(X\t_Y\pd\ti X),\bs\pi_Y,\bs
D_{\ubC,\ubtC\vert_{\pd\ti X}}\bigr].
\end{split}
\label{kh4eq61}
\e

Define $R$-bilinear $\Xi^{k,l}:KC^k_\ec(Y;R)\t KC^l_\ec(Y;R)\ra
KC^{k+l-1}_\ec(Y;R)$ by
\begin{equation*}
\Xi^{k,l}:\bigl([X,\bs f,\ubC],[\ti X,\bs{\ti f},\ubtC]\bigr)
\longmapsto\bigl[[0,1]\t(X\t_Y\ti
X),\bs\pi_Y,\ubD_{\ubC,\ubtC}\bigr].
\end{equation*}
Using \eq{kh4eq61} and putting $\vdim X=\dim Y-k$, $\vdim\ti X=\dim
Y-l$ yields
\begin{equation*}
\d\ci\Xi^{k,l}(\al,\be)\!=\!-\al\cup\be\!+\!(-1)^{kl}\be\cup\al
\!-\!\Xi^{k+1,l}(\d\al,\be)\!-\!(-1)^{k}\Xi^{k,l+1}(\d\al,\d\be)
\end{equation*}
for $\al\in KC_\ec^k(Y;R)$ and $\be\in KC_\ec^l(Y;R)$. Passing to
cohomology, this gives $\ga\cup\de =(-1)^{kl}\de\cup\ga$ for $\ga\in
KH^k_\ec(Y;R)$, $\de\in KH^l_\ec(Y;R)$, that is, $\cup$ is
supercommutative on $KH^*_\ec(Y;R)$. The last part is immediate as
$\cup,\bu$ are identified by $\Pi_\ec^\ef$ in~\eq{kh4eq45}.
\end{proof}

\begin{rem} We have now defined Kuranishi cohomology $KH^*(Y;R)$,
which has the advantage of a cup product $\cup$ which is
supercommutative and associative on cochains, but the disadvantage
that it works only for $R$ a $\Q$-algebra. We have also defined
effective Kuranishi cohomology $KH^*_\ec(Y;R)$, which has the
advantage that it works for $R$ any commutative ring, but the
disadvantage that the cup product $\cup$ is not supercommutative on
cochains.

In fact this dichotomy is necessary: it should not be possible to
define a cohomology theory which both works when $R=\Z$, and has a
cup product $\cup$ on cochains that is everywhere defined and
supercommutative. This because of {\it Steenrod squares\/} \cite[\S
VI.15--\S VI.16]{Bred}, which are cohomology operations defined
using the failure of $\cup$ to be supercommutative at the cochain
level. Since Steenrod squares are known to be nontrivial when $R=\Z$
(though they are trivial for $R=\Q$), they are an obstruction to
defining such a cohomology theory.

To make our theories of effective Kuranishi (co)homology work over
$R=\Z$ and compute $H_*^\rsi(Y;R)$ and $H^*_\cs(Y;R)$, we had to
sacrifice some good properties of Kuranishi (co)homology, including
supercommutativity of $\cup$ on cochains. The argument above
suggests these sacrifices were unavoidable.
\label{kh4rem6}
\end{rem}

The next theorem, proved in Appendix \ref{khD}, shows that (at least
for manifolds $Y$) the isomorphisms $\Pi_\cs^\Kch,\Pi_\cs^\ec,
\Pi_\rsi^\Kh,\Pi_\rsi^\ef$ identify the products $\cup,\cap,\bu$ on
$H^*_\cs,H^\rsi_*(Y;R)$ with those on $KH^*,KH_*(Y;R)$
and~$KH_\ec^*,KH^\ef_*(Y;R)$.

\begin{thm} Let\/ $Y$ be a manifold without boundary, and\/ $R$ a
commutative ring {\rm(\kern.5pt}for $KH_*^\ef,KH^*_\ec(Y;R))$ or
$\Q$-algebra {\rm(\kern.5pt}for $KH_*,KH^*(Y;R))$. Then
\begin{itemize}
\setlength{\itemsep}{0pt}
\setlength{\parsep}{0pt}
\item[{\rm(a)}] The isomorphisms $\Pi_\cs^\Kch:H^*_\cs(Y;R)\ra
KH^*(Y;R)$ and\/ $\Pi_\cs^\ec:H^*_\cs(Y;R)\ab\ra KH^*_\ec(Y;R)$
identify the cup product\/ $\cup$ on $H^*_\cs(Y;R)$ with those on
$KH^*,\ab KH_\ec^*\ab(Y;R)$. That is, for all\/ $\al,\be\in
H^*_\cs(Y;R)$ we have
\begin{equation*}
\smash{\Pi_\cs^\Kch(\al\!\cup\!\be)\!=\!\bigl(\Pi_\cs^\Kch(\al)\bigr)\!
\cup\!\bigl(\Pi_\cs^\Kch(\be)\bigr), \;\>
\Pi_\cs^\ec(\al\!\cup\!\be)\!=\!\bigl(\Pi_\cs^\ec(\al)\bigr)\!
\cup\!\bigl(\Pi_\cs^\ec(\be)\bigr).}
\end{equation*}
\item[{\rm(b)}] $\Pi_\rsi^\Kh,\Pi_\rsi^\ef,\Pi_\cs^\Kch,\Pi_\cs^\ec$
identify the cap product\/ $\cap:H_*^\rsi(Y;R)\t H^*_\cs(Y;R)\ra
H_*^\rsi(Y;R)$ with those on {\rm(}effective\/{\rm)} Kuranishi
{\rm(}co\/{\rm)}homology.
\item[{\rm(c)}] If\/ $Y$ is oriented then $\Pi_\rsi^\Kh,
\Pi_\rsi^\ef$ identify the intersection product\/ $\bu$ on\/
$H_*^\rsi(Y;R)$ with those on\/ $KH_*,KH_*^\ef(Y;R)$.
\item[{\rm(d)}] If\/ $Y$ is compact then $\Pi_\cs^\Kch,\Pi_\cs^\ec$
take the identity $1$ in $H^0_\cs(Y;R)$ to the identities
$\bigl[[Y,\id_Y,\bs C_Y]\bigr],\bigl[[Y,\id_Y,\ubC_Y]\bigr]$
in~$KH^0,KH_\ec^0(Y;R)$.
\end{itemize}
\label{kh4thm6}
\end{thm}

\subsection[Kuranishi (co)homology as a bivariant theory on
(co)chains]{Kuranishi (co)homology as a bivariant theory, at the
(co)chain level}
\label{kh48}

For manifolds, because of Poincar\'e duality, homology and
cohomology can be identified into a single theory with {\it
products\/} (cup, cap and intersection products), {\it
pushforwards\/} (covariant functoriality), and {\it pullbacks}
(contravariant functoriality). Fulton and MacPherson \cite{FuMa}
defined {\it bivariant theories}, a way of naturally integrating
pairs of a homology theory and a cohomology theory into a single
larger theory, which works for singular spaces.

We now summarize the notion of a bivariant theory, following
\cite[\S 1.1]{FuMa}. Let ${\cal C}$ be a category. A bivariant
theory from $\cal C$ to abelian groups (or graded vector spaces,
etc.) associates a group $T(f:X\ra Y)$ to each morphism $f:X\ra Y$
in $\cal C$. It also has three operations:
\begin{itemize}
\setlength{\itemsep}{0pt}
\setlength{\parsep}{0pt}
\item {\bf Products.} If $f:X\ra Y$ and $g:Y\ra Z$ are morphisms in
$\cal C$, and $\al\in T(f:X\ra Y)$, $\be\in T(Y\ra Z)$, we can form
the product~$\be\cdot\al\in T(g\ci f:X\ra Z)$.
\item {\bf Pushforwards.} For a certain class of morphisms $f:X\ra
Y$ in $\cal C$ called {\it confined morphisms}, if $g:Y\ra Z$ is a
morphism in $\cal C$ and $\al\in T(g\ci f:X\ra Z)$ we can form the
pushforward $f_*(\al)\in T(g:Y\ra Z)$.
\item {\bf Pullbacks.} For a certain class of commutative squares in
$\cal C$ called {\it independent squares}, if
\e
\begin{gathered}
\xymatrix@C=40pt@R=10pt{
X' \ar[r]_{g'} \ar[d]_{f'} & X \ar[d]^f \\
Y' \ar[r]^g & Y, }
\end{gathered}
\label{kh4eq62}
\e
is an independent square then for each $\al\in T(f:X\ra Y)$, we can
form the pullback $g^*(\al)\in T(f':X'\ra Y')$.
\end{itemize}
All this data must satisfy a number of compatibility axioms.

The following example, taken from \cite[\S 3.1]{FuMa} modified as in
\cite[\S 3.3.2]{FuMa}, combines homology and compactly-supported
cohomology into a bivariant theory.

\begin{ex} Let {\bf Top} be the category of topological spaces $X$
which admit an embedding as a closed subspace of $\R^n$ for some
$n\gg 0$ (in particular, $X$ must be Hausdorff and paracompact), and
morphisms are continuous maps. Let $R$ be a commutative ring.
Suppose $f:X\ra Y$ is a morphism in {\bf Top} with $X$ {\it
compact}. Choose an embedding $\phi:X\ra\R^n$ for some $n\gg 0$, and
define $H^k(f:X\ra Y;R)=H^{k+n}\bigl(Y\t\R^n,Y\t\R^n\sm(f\t\phi)(X);
R\bigr)$, using relative cohomology. Then $H^*(f:X\ra Y;R)$ is a
graded $R$-module, which turns out to be independent of the choice
of $\phi$ up to canonical isomorphism. For general morphisms in {\bf
Top}, we define $H^*(f:X\ra Y;R)$ to be the direct limit $\varinjlim
H^*(f\ci p:K\ra Y;R)$ over all $p:K\ra X$ in {\bf Top} with $K$
compact.

If $X$ is a topological space in {\bf Top}, we can associate two
natural morphisms in {\bf Top} to $X$: the projection to a point
$\pi:X\ra\{0\}$, and the identity map $\id_X:X\ra X$. By standard
algebraic topology one can show that $H^k(\pi:X\ra\{0\};R)\cong
H_{-k}(X;R)$ and $H^k(\id_X:X\ra X;R)\cong H^k_\cs(X;R)$. Thus, the
bivariant theory $H^*$ specializes to homology and
compactly-supported cohomology over~$R$.

{\it Products\/} are defined as in \cite[\S 3.1.7]{FuMa}; the
product $H^k(\id_X:X\ra X;R)\t H^l(\id_X:X\ra X;R)\ra
H^{k+l}(\id_X:X\ra X;R)$ is $\cup:H^k_\cs(X;R)\t H^l_\cs(X;R)\ra
H^{k+l}_\cs(X;R)$, and the product $H^k(\pi:X\ra\{0\};R)\t
H^l(\id_X:X\ra X;R)\ra H^{k+l}(\pi:X\ra \{0\};R)$ is
$\cap:H_{-k}(X;R)\t H^l_\cs(X;R)\ra H_{-k-l}(X;R)$.

All morphisms in {\bf Top} are confined, and {\it pushforwards} are
defined as in \cite[\S 3.1.8]{FuMa}. For $f:X\ra Y$ a morphism in
{\bf Top}, the pushforward $f_*:H^k(\pi:X\ra\{0\};R)\ra H^k(\pi:Y\ra
\{0\};R)$ is $f_*:H_{-k}(X;R)\ra H_{-k}(Y;R)$ on homology.

We define {\it independent squares\/} to be squares \eq{kh4eq62}
which are Cartesian, so that $X'\cong X\t_{f,Y',g}Y$, and with $g$
proper, so that $g'$ is also proper. {\it Pullbacks\/} are defined
as in \cite[\S 3.1.6]{FuMa}. If $g:Y'\ra Y$ is a proper morphism in
{\bf Top} then $g^*:H^k(\id_Y:Y\ra Y;R)\ra H^k(\id_{Y'}:Y'\ra Y';R)$
is $g^*:H_\cs^k(Y;R)\ra H_\cs^k(Y';R)$ on compactly-supported
cohomology.
\label{kh4ex2}
\end{ex}

Let {\bf Man} be the category of manifolds, with morphisms smooth
maps. Applying the natural functor ${\bf Man}\ra{\bf Top}$ taking a
manifold to its underlying topological space, the bivariant theory
$H^*$ of Example \ref{kh4ex2} induces a bivariant theory on {\bf
Man}. However, because of Poincar\'e duality the picture simplifies,
and we can write $H^*$ on {\bf Man} more directly. The following is
an equivalent way of defining the bivariant theory of Example
\ref{kh4ex2} on the category {\bf Man}.

\begin{ex} Let {\bf Man} be the category of manifolds, with
morphisms smooth maps, and $R$ be a commutative ring. Suppose
$f:X\ra Y$ is a morphism in {\bf Man}, with $\dim X=m$ and $\dim
Y=n$. Define
\e
\begin{split}
H^k(f:X\ra Y;R)&=H^{k+m-n}_\cs(X;O_X\ot f^*(O_Y)\t_{\{\pm 1\}}R)\\
&\cong H_{n-k}(X;f^*(O_Y)\t_{\{\pm 1\}}R),
\end{split}
\label{kh4eq63}
\e
where as in Definition \ref{kh4def8}, $O_X,O_Y$ are the {\it
orientation bundles\/} of $X,Y$, which are principal $\Z_2$-bundles
over $X,Y$, so that $O_X\ot f^*(O_Y)$ is a principal $\Z_2$-bundle
over $X$, and $H^*_\cs(X;O_X\ot f^*(O_Y)\t_{\{\pm 1\}}R)$ is
compactly-supported cohomology of $X$ with coefficients in $R$ {\it
twisted by\/} $O_X\ot f^*(O_Y)$, and $H_*(X;f^*(O_Y)\t_{\{\pm
1\}}R)$ is homology of $X$ with coefficients in $R$ {\it twisted
by\/} $f^*(O_Y)$, and the groups on the r.h.s.\ of \eq{kh4eq63} are
isomorphic by Poincar\'e duality.

If $f:X\ra Y$ is $\pi:X\ra\{0\}$ then $O_Y$ and $f^*(O_Y)$ are
trivial and $n=0$, so $H^k(\pi:X\ra\{0\};R)\cong H_{-k}(X;R)$. If
$f:X\ra Y$ is $\id_X:X\ra X$ then $f^*(O_Y)\cong O_X$ and $O_X\ot
f^*(O_Y)$ is the trivial $\Z_2$-bundle over $X$, so that
$H^k(\id_X:X\ra X;R)=H^k_\cs(X;R)$. Thus, \eq{kh4eq63} specializes
to homology and compactly-supported cohomology.

To define products, note that the cup product generalizes to twisted
cohomology as follows: if $P\ra X$ and $Q\ra X$ are principal
$\Z_2$-bundles then there is a cup product $\cup:H_\cs^k(X;
P\t_{\{\pm 1\}}R)\t H_\cs^l(X;Q\t_{\{\pm 1\}}R)\ra
H_\cs^{k+l}(X;(P\ot Q)\t_{\{\pm 1\}}R)$. Suppose $f:X\ra Y$ and
$g:Y\ra Z$ are smooth with $\dim X=m$, $\dim Y=n$ and $\dim Z=o$.
For $\al\in H^k(f:X\ra Y;R)$ and $\be\in H^l(g:Y\ra Z;R)$, define
$\be\cdot\al=f^*(\be)\cup\al$. Here $\al\in H^{k+m-n}_\cs(X;O_X\ot
f^*(O_Y)\t_{\{\pm 1\}}R)$ and $f^*(\be)\in
H^{l+n-o}_\cs(X;f^*(O_Y\ot g^*(O_Z))\t_{\{\pm 1\}}R)$, and
$\bigl(O_X\ot f^*(O_Y)\bigr)\ot f^*\bigl(O_Y\ot g^*(O_Z)\bigr)\cong
O_X\ot(g\ci f)^*(O_Z)$ as principal $\Z_2$-bundles over $X$, so that
$f^*(\be)\cup\al$ lies in $H^{k+l+m-o}_\cs(X;O_X\ot(g\ci
f)^*(O_Z)\t_{\{\pm 1\}}R)=H^{k+l}(g\ci f:X\ra Z;R)$.

All morphisms in {\bf Man} are confined. {\it Pushforwards\/}
$f_*:H^k(g\ci f:X\ra Z;R)\ra H^k(g:Y\ra Z;R)$ are given by
pushforwards $f_*:H_{o-k}(X;(g\ci f)^*(O_Z)\t_{\{\pm 1\}}R)\ra
H_{o-k}(Y;g^*(O_Z)\t_{\{\pm 1\}}R)$ in twisted homology. {\it
Independent squares\/} are Cartesian squares \eq{kh4eq62} in {\bf
Man} with $g$ proper, and hence $g'$ proper. Then with $\dim
X\!=\!m$, $\dim Y\!=\!n$, $\dim X'\!=\!m'$, $\dim Y'\!=\!n'$, the
Cartesian property gives $(g')^*\bigl(O_X\!\ot\!
f^*(O_Y)\bigr)\!\cong\! O_{X'}\!\ot\! (f')^*(O_{Y'})$ and
$m-n\!=\!m'-n'$. So we define $g^*:H^k(f:X\ra Y;R)\!\ra\!
H^k(f':X'\ra Y';R)$ to be
\begin{align*}
(g')^*:\,&H^{k+m-n}_\cs(X;O_X\ot f^*(O_Y)\t_{\{\pm 1\}}R)\longra\\
&H^{k+m-n}_\cs(X';(g')^*(O_X\ot f^*(O_Y))\t_{\{\pm 1\}}R)\cong\\
&H^{k+m'-n'}_\cs(X';O_{X'}\ot (f')^*(O_{Y'})\t_{\{\pm 1\}}R),
\end{align*}
using proper pullbacks on twisted compactly-supported cohomology.
\label{kh4ex3}
\end{ex}

For {\it oriented\/} manifolds $X,Y$, with $O_X,O_Y$ trivial,
equation \eq{kh4eq63} becomes
\begin{equation*}
H^k(f:X\ra Y;R)=H^{k+m-n}_\cs(X;R)\cong H_{n-k}(X;R),
\end{equation*}
and we lose the dependence on $f,Y$ entirely. Thus, for oriented
manifolds and smooth maps, the bivariant theories of Examples
\ref{kh4ex2} and \ref{kh4ex3} reduce to homology,
compactly-supported cohomology, and Poincar\'e duality. That is,
{\it bivariant (co)homology theories really give us nothing new for
manifolds}, they are of interest {\it only for singular spaces}.

As our Kuranishi (co)homology theories are defined only for
manifolds and orbifolds, and bivariant (co)homology theories give
nothing new for manifolds, one might expect that there is nothing to
be gained from recasting Kuranishi (co)homology as a bivariant
theory. However, as we have already emphasized, Kuranishi
(co)homology is particularly well behaved at the (co)chain level,
with pushforwards, pullbacks, and cup and cap products on (co)chains
all well-defined and with all the properties one would wish.

We will now explain how to formulate (effective) Kuranishi
(co)homology as a bivariant theory {\it at the (co)chain level}. The
author does not know of any other (co)homology theories for which
this can be done. This bivariant framework is likely to be valuable
in situations in which one is forced to work with (co)chains rather
than just with (co)homology, for example, in defining
$A_\iy$-algebras in Lagrangian Floer cohomology, as in \cite{AkJo}
and sketched in~\S\ref{kh66}.

\begin{dfn} Let {\bf Orb} be the category of orbifolds without
boundary, with morphisms smooth maps, and $R$ be a commutative ring.
Roughly speaking, we will define a bivariant theory from {\bf Orb}
to the category of differential graded $R$-modules, generalizing
Kuranishi (co)chains. However, we will define $KC^*(g:Y\ra Z;R)$
only when $g:Y\ra Z$ is a {\it submersion\/} of orbifolds~$Y,Z$.

Let $Y,Z$ be orbifolds without boundary and $g:Y\ra Z$ be a
submersion. Consider triples $(X,\bs f,\bs B)$, where $X$ is a
compact Kuranishi space, $\bs f:X\ra Y$ is a strongly smooth map
such that $g\ci\bs f:X\ra Z$ is a strong submersion, and we are
given a coorientation for $(X,g\ci\bs f)$, and $\bs B=(\bs
I,\bs\eta,B^i:i\in I)$, where $(\bs I,\bs\eta)$ is an excellent
coordinate system for $(X,\bs f)$ with $\bs
I=\bigl(I,(V^i,\ldots,\psi^i),f^i:i\in I,\ldots\bigr)$, and $g\ci
f^i:V^i\ra Z$ are submersions representing $g\ci\bs f$ for all $i\in
I$, and $B^i:E^i\ra P$ for $i\in I$ are maps such that each $B^i$
maps $E^i\ra P_l\subset P$ for $l\gg 0$, and $B^i\t(g\ci
f^i\ci\pi^i):E^i\ra P\t Y$ is a globally finite map for all $i\in
I$. Such $\bs B$ generalize the notions of (co-)gauge-fixing data in
\S\ref{kh31}, since if $g:Y\ra Z$ is $\pi:Y\ra\{0\}$ then $\bs B$ is
gauge-fixing data for $(X,\bs f)$, and if $g:Y\ra Z$ is $\id_Y:Y\ra
Y$ then $\bs B$ is co-gauge-fixing data for $(X,\bs f)$. We call
$\bs B$ {\it bi-gauge-fixing data}.

Define {\it isomorphisms\/} $(\bs a,\bs b)$ of such triples $(X,\bs
f,\bs B)$ as for (co-)gauge-fixing data in \S\ref{kh31}. Write
$[X,\bs f,\bs B]$ for the isomorphism class of $(X,\bs f,\bs B)$
under isomorphisms $(\bs a,\bs b):(X,\bs f,\bs B)\ra(\ti X,\bs{\ti
f},\bs{\ti B})$ which identify the coorientations of $(X,g\ci\bs f)$
and $(\ti X,g\ci\bs{\ti f})$.

For each $k\in\Z$, define $KC^k(g:Y\ra Z;R)$ to be the $R$-module of
finite $R$-linear combinations of isomorphism classes $[X,\bs f,\bs
G]$ for which $\vdim X=\dim Z-k$, with the analogues of relations
Definition \ref{kh4def2}(i)--(iv). Then $KC^k(\pi:Y\ra\{0\};R)=
KC_{-k}(Y;R)$ and $KC^k(\id_Y:Y\ra Y;R)=KC^k(Y;R)$, so this
generalizes the definitions of Kuranishi (co)chains in \S\ref{kh42}
and \S\ref{kh44}. Elements of $KC^*(g:Y\ra Z;R)$ will be called {\it
Kuranishi bichains}.

Define $\d:KC^k(g:Y\ra Z;R)\ra KC^{k+1}(g:Y\ra Z;R)$ by
\begin{equation*}
\d:\ts\sum_{a\in A}\rho_a[X_a,\bs f_a,\bs B_a]\longmapsto
\ts\sum_{a\in A}\rho_a[\pd X_a,\bs f_a\vert_{\pd X_a},\bs B_a
\vert_{\pd X_a}],
\end{equation*}
as in \eq{kh4eq8} and \eq{kh4eq20}. Then $\d\ci\d=0$. Define
$KH^*(g:Y\ra Z;R)$ to be the cohomology of $\bigl(KC^*(g:Y\ra
Z;R),\d\bigr)$. Then $KH^k(\pi:Y\ra\{0\};R)=KH_{-k}(Y;R)$ and
$KH^k(\id_Y:Y\ra Y;R)=KH^k(Y;R)$. Combining the proofs of Corollary
\ref{kh4cor1} and Theorem \ref{kh4thm3} shows that $KH^k(g:Y\ra
Z;R)\cong  H_{\dim Z-k}^\rsi(Y;g^*(O_Z)\t_{\{\pm 1\}}R)$, as
in~\eq{kh4eq63}.

We define the three bivariant operations on Kuranishi bichains:
\begin{itemize}
\setlength{\itemsep}{0pt}
\setlength{\parsep}{0pt}
\item {\bf Products.} Let $X,Y,Z$ be orbifolds without boundary, and
$f:X\ra Y$ and $g:Y\ra Z$ be submersions. Define a product
$KC^k(g:Y\ra Z;R)\t KC^l(f:X\ra Y;R)\ra KC^{k+l}(g\ci f:X\ra
Z;R)$~by
\e
[W,\bs e,\bs B]\cdot[\ti W,\bs{\ti e},\bs{\ti B}]=[W\t_{\bs
e,Y,f\ci\bs{\ti e}}\ti W,\bs{\ti e}\ci\bs\pi_{\ti W},\bs
B\t_Y\bs{\ti B}],
\label{kh4eq64}
\e
extended $R$-bilinearly. Here $[W,\bs e,\bs B]\in KC^k(g:Y\ra Z;R)$,
so that $\bs e:W\ra Y$ is strongly smooth, and $[\ti W,\bs{\ti
e},\bs{\ti B}]\in KC^l(f:X\ra Y;R)$ so that $f\ci\bs{\ti e}:\ti W\ra
Y$ is a strong submersion. Thus $W\t_{\bs e,Y,f\ci\bs{\ti e}}\ti W$
is a compact Kuranishi space, and $\bs{\ti e}\ci\bs\pi_{\ti
W}:W\t_{\bs e,Y,f\ci\bs{\ti e}}\ti W$ is strongly smooth.

We have $f\ci(\bs{\ti e}\ci\bs\pi_{\ti W})=(f\ci\bs{\ti
e})\ci\bs\pi_{\ti W}=\bs\pi_Y=\bs e\ci\bs\pi_W$. Composing with $g$
then gives $(g\ci f)\ci(\bs{\ti e}\ci\bs\pi_{\ti W})=(g\ci\bs
e)\ci\bs\pi_W$. But $g\ci\bs e$ is a strong submersion as $[W,\bs
e,\bs B]\in KC^k(g:Y\ra Z;R)$, and $\bs\pi_W$ is a strong submersion
as $f\ci\bs{\ti e}$ is, by properties of fibre products. Hence
$(g\ci f)\ci(\bs{\ti e}\ci\bs\pi_{\ti W})$ is a strong submersion,
as we need for the r.h.s.\ of \eq{kh4eq64} to lie in $KC^{k+l}(g\ci
f:X\ra Z;R)$. The coorientations for $(W,g\ci\bs e)$ and $(\ti
W,f\ci\bs{\ti e})$ combine to give a coorientation for
$\bigl(W\t_{\bs e,Y,f\ci\bs{\ti e}}\ti W,(g\ci f)\ci(\bs{\ti
e}\ci\bs\pi_{\ti W})\bigr)$, using $(g\ci f)\ci(\bs{\ti
e}\ci\bs\pi_{\ti W})=(g\ci\bs e)\ci\bs\pi_W$. The fibre product of
bi-gauge-fixing data $\bs B\t_Y\bs{\ti B}$ is defined as for $\bs
C\t_Y\bs{\ti C}$ in~\S\ref{kh38}.

When $X=Y=Z$ and $f=g=\id_Y$, equation \eq{kh4eq64} agrees with
$\cup$ in \eq{kh4eq47}, and when $X=Y$, $f=\id_Y$, $Z=\{0\}$ and
$g=\pi:Y\ra\{0\}$, equation \eq{kh4eq64} agrees with $\cap$ in
\eq{kh4eq50}. Thus, the product \eq{kh4eq64} generalizes $\cup,\cap$
on Kuranishi (co)chains.
\item {\bf Pushforwards.} Let $X,Y,Z$ be orbifolds without
boundary, and $f:X\ra Y$ and $g:Y\ra Z$ be smooth with $g,g\ci f$
submersions. Define $f_*:KC^k(g\ci f:X\ra Z;R)\ra KC^k(g:Y\ra Z;R)$
by
\begin{equation*}
f_*:\ts\sum_{a\in A}\rho_a\bigl[W_a,\bs e_a,\bs B_a \bigr]
\longmapsto \ts\sum_{a\in A}\rho_a\bigl[W_a,f\ci\bs e_a,f_*(\bs
B_a)\bigr],
\end{equation*}
as in \eq{kh4eq11}, where $f_*(\bs B_a)$ is defined as for $h_*(\bs
G)$ in \S\ref{kh37}. Note that as $(g\ci f)\ci\bs e_a$ is a strong
submersion with $\bigl(W_a,(g\ci f)\ci\bs e_a\bigr)$ cooriented,
trivially $g\ci(f\ci\bs e_a)$ is a strong submersion with
$\bigl(W_a,g\ci(f\ci\bs e_a)\bigr)$ cooriented.
\item {\bf Pullbacks.} Let \eq{kh4eq62} be a Cartesian square of
orbifolds without boundary, with $f$ (and hence $f'$) a submersion,
and $g$ (and hence $g'$) proper. Define $g^*:KC^k(f:X\ra Y;R)\ra
KC^k(f':X'\ra Y';R)$ by
\begin{equation*}
f^*:\ts\sum_{a\in A}\rho_a\bigl[W_a,\bs e_a,\bs B_a\bigr]\longmapsto
\ts\sum_{a\in A}\rho_a [X'\t_{g',X,\bs
e_a}W_a,\bs\pi_{X'},(g')^*(\bs B_a)],
\end{equation*}
as in \eq{kh4eq21}, where $(g')^*(\bs B_a)$ is defined as for
$h^*(\bs C)$ in \S\ref{kh48}. Since $f\ci\bs e_a:W_a\ra Y$ is a
strong submersion and \eq{kh4eq62} is Cartesian, $f'\ci\bs\pi_{X'}:
X'\t_{g',X,\bs e_a}W_a\ra Y'$ is a strong submersion. Also, using
the coorientation for $(W_a,f\ci\bs e_a)$ we construct a
coorientation for~$(X'\t_{g',X,\bs e_a}W_a,f'\ci\bs\pi_{X'})$.
\end{itemize}
It is easy to show that these operations satisfy the axioms for a
bivariant theory in \cite[\S 2.2]{FuMa}. Also, the three operations
are compatible with differentials $\d$ on $KC^*(g:Y\ra Z;R)$, and so
induce operations on~$KH^*(g:Y\ra Z;R)$.
\label{kh4def17}
\end{dfn}

It is easy to generalize Definition \ref{kh4def17} to {\it
effective} Kuranishi (co)homology. The differences are that we
restrict to submersions $g:Y\ra Z$ for which
$g_*:\Stab_Y(y)\ra\Stab_Z(g(y))$ is {\it surjective\/} for all $y\in
Y$, and in defining $KC^k_\ef(g:Y\ra Z;R)$, we use triples $(X,\bs
f,\bs{\ul B})$ in which $(X,g\ci\bs f)$ is also {\it coeffective},
and $\bs{\ul B}=(\bs I,\bs\eta,\ul B^i:i\in I)$, where
$(V^i,\ldots,\psi^i),g\ci f^i$ is {\it coeffective\/} for all $i\in
I$, and $\ul B^i:E^i\ra\uP$ are maps such that $\ul B^i\t(g\ci
f^i\ci\pi^i):E^i\ra\uP\t Z$ satisfy the conditions of
Definition~\ref{kh3def19}.

\subsection{Why we need (co-)gauge-fixing data}
\label{kh49}

To show that we really do need some form of (co-)gauge-fixing data
to define sensible (co)homology theories using Kuranishi spaces as
chains, we will now prove that if we omit (co-)gauge-fixing data
from the definitions above, then Kuranishi (co)homology is always
zero.

Let $Y$ be an orbifold without boundary, and $R$ a commutative ring.
Define {\it na\"\i ve Kuranishi chains} $KC^\na_k(Y;R)$ to be the
$R$-module spanned by isomorphism classes $[X,\bs f]$ of pairs
$(X,\bs f)$ with $X$ a compact, oriented Kuranishi space with $\vdim
X=k$ and $\bs f:X\ra Y$ strongly smooth, with the analogues of
relations Definition \ref{kh4def2}(i)--(iii). The boundary operator
is $\pd:[X,\bs f]\mapsto[\pd X,\bs f\vert_{\pd X}]$. Let {\it na\"\i
ve Kuranishi homology} $KH^\na_*(Y;R)$ be the homology of~$\bigl(
KC^\na_*(Y;R),\pd\bigr)$.

Similarly, define {\it na\"\i ve Kuranishi cochains} $KC_\na^k(Y;R)$
to be spanned by isomorphism classes $[X,\bs f]$ with $X$ a compact
Kuranishi space with $\vdim X=\dim Y-k$ and $\bs f:X\ra Y$ a
cooriented strong submersion, and the same relations. Define
$\d:[X,\bs f]\mapsto[\pd X,\bs f\vert_{\pd X}]$, and let {\it na\"\i
ve Kuranishi cohomology} $KH_\na^*(Y;R)$ be the cohomology of
$\bigl(KC_\na^*(Y;R),\d\bigr)$. One can then define cup and cap
products $\cup,\cap$ on $KC_\na^*,KC^\na_*(Y;R)$ and
$KH_\na^*,KH^\na_*(Y;R)$, as in \S\ref{kh47}, but omitting all
(co-)gauge-fixing data.

Take $Y$ to be a point $\{0\}$, and for any Kuranishi space $X$
write $\bs\pi:X\ra\{0\}$ for the trivial projection. Let $R$ be a
commutative ring. In the next example we show that the identity
$\bigl[[\{0\},\id_{\{0\}}]\bigr]$ in $KH_\na^*(\{0\};R)$ is zero.
From this we will deduce that $KH_\na^*(Y;R)=KH^\na_*(Y;R)=0$ for
all~$Y$.

\begin{ex} Let $L\ra\CP^1$ be the complex line bundle ${\mathcal
O}(1)$. Define a Kuranishi space $X_k$ for $k\in\Z$ to be the
topological space $\CP^1$ with the Kuranishi structure induced by
the Kuranishi neighbourhood $(\CP^1,L^k,0,\id_{\CP^1})$, with
obstruction bundle $L^k\ra\CP^1$. Give $(X_k,\bs\pi)$ the obvious
coorientation. Then $\vdim X_k=0$, and $\pd X_k=\es$, so
$[X_k,\bs\pi]$ defines a class $\bigl[[X_k,\bs\pi]\bigr]$
in~$KH^0_\na(\{0\};R)$.

When $k\ge 0$ we can choose a generic smooth section $s$ of
$L^k\ra\CP^1$ which has exactly $k$ zeroes $x_1,\ldots,x_k$, each of
multiplicity 1. Let $t$ be the coordinate on $[0,1]$. Then $ts$ is a
section of $L^k\ra[0,1]\t\CP^1$, with $(ts)^{-1}(0)=\{0\}\t\CP^1\cup
[0,1]\t\{x_1,\ldots,x_k\}$, and $\bigl([0,1]\t\CP^1,L^k,ts,
\id_{(ts)^{-1}(0)}\bigr)$ is a Kuranishi neighbourhood on
$(ts)^{-1}(0)$, making it into a Kuranishi space of virtual
dimension 1. By taking the boundary of this we see that
\e
\bigl[[X_k,\bs\pi]\bigr]=\bigl[[\{x_1,\ldots,x_k\},\bs\pi]\bigr]=
k\bigl[[\{0\},\id_{\{0\}}]\bigr]
\label{kh4eq65}
\e
in $KH^0_\na(\{0\};R)$, that is, $\bigl[[X_k,\bs\pi]\bigr]$ is $k$
times the fundamental class of~$Y=\{0\}$.

Write $[z_0,z_1]$ for the homogeneous coordinates on $\CP^1$, and
define
\begin{equation*}
V=\bigl\{(t,[z_0,z_1])\in\R\t\CP^1:\min(\md{z_0}^2,\md{z_1}^2)
\max(\md{z_0}^2,\md{z_1}^2)^{-1}\le t\le 2\bigr\}.
\end{equation*}
Then $V$ is a compact oriented 3-manifold with corners, and $\pd V$
is the disjoint union of three pieces, a copy of $\CP^1$ with $t=2$,
the hemisphere $H_+=\bigl\{[z_0,z_1]\in\CP^1:\md{z_0}\le
\md{z_1}\bigr\}$ with $t=\md{z_0}^2/\md{z_1}^2$, and the hemisphere
$H_-=\bigl\{[z_0,z_1]\in\CP^1:\md{z_0}\ge\md{z_1}\bigr\}$
with~$t=\md{z_1}^2/\md{z_0}^2$.

Define the Kuranishi space $W_k$ for $k\in\Z$ to be the topological
space $V$ with the Kuranishi structure induced by the neighbourhood
$(V,\pi^*(L^k),0,\id_V)$, where $\pi:V\ra\CP^1$ is the projection.
Define Kuranishi spaces $X_+,X_-$ to be the topological spaces
$H_+,H_-$ with the Kuranishi structure induced by the Kuranishi
neighbourhoods $(H_\pm,L^0,0,\id_{H_\pm})$. Now the line bundles
$L^k\ra H^\pm$ are for $k\in\Z$ are isomorphic to $L^0\ra H^\pm$.
Thus there is an isomorphism of cooriented Kuranishi spaces $(\pd
W_k,\bs\pi)\cong(X_k,\bs\pi)\amalg -(X_+,\bs\pi)\amalg
-(X_-,\bs\pi)$, so in $KH^0_\na(\{0\};R)$ we have
$\bigl[[X_k,\bs\pi]\bigr]=\bigl[[X_+,\bs\pi]+[X_-,\bs\pi]\bigr]$,
and the class $\bigl[[X_k,\bs\pi]\bigr]$ is independent of $k\in\Z$.
Combining this with \eq{kh4eq65} gives $\bigl[[\{0\},\id_{\{0\}}]
\bigr]=0$.
\label{kh4ex4}
\end{ex}

If $Y$ is any compact orbifold without boundary, and $\pi:Y\ra\{0\}$
is the projection, then the pullback $\pi^*:KH_\na^*(\{0\};R)\ra
KH_\na^*(Y;R)$ takes $\pi^*:\bigl[[\{0\},\id_{\{0\}}]\bigr]\mapsto
\bigl[[Y,\id_Y]\bigr]$, so $\bigl[[\{0\},\id_{\{0\}}]\bigr]=0$
implies $\bigl[[Y,\id_Y]\bigr]=0$. But $\bigl[[Y,\id_Y]\bigr]$ is
the identity for $\cup$ on $KH_\na^*(Y;R)$ and $\cap$ on
$KH^\na_*(Y;R)$, so this implies $KH_\na^*(Y;R)=KH^\na_*(Y;R)=0$.
Using a related argument on (co)chains, we can also show
$KH_\na^*(Y;R)= KH^\na_*(Y;R)=0$ for any orbifold $Y$. Hence,
Kuranishi (co)homology without (co-)gauge-fixing data is vacuous.

The reason this example goes so badly wrong is that the Kuranishi
space $\pd X^+=-\pd X^-$, the circle $C$ with Kuranishi
neighbourhood $(C,\R^2\t C,0,\id_C)$, has a {\it large automorphism
group}, including topologically nontrivial automorphisms which fix
$C$ but rotate the fibres of the obstruction bundle $\R^2\t C$ by a
degree $l$ smooth map $C\ra{\rm SO}(2)$. Thus, by cutting $X_k$ into
two pieces $X^+,X^-$ and then gluing them together again twisted by
such a nontrivial automorphism, we can get $X_{k+l}$ for any
$l\in\Z$, which forces $\bigl[[X_k,\bs\pi]\bigr]=\bigl[[X_{k+l},
\bs\pi]\bigr]$ in $KH^0_\na(\{0\};R)$. Including (co-)gauge-fixing
data avoids this problem, as by Theorem \ref{kh3thm2} it allows only
finite automorphism groups.

\section{Kuranishi bordism and cobordism}
\label{kh5}

We now discuss {\it bordism\/} and {\it cobordism\/} theories of an
orbifold $Y$ and a commutative ring $R$, defined in a similar way to
Kuranishi (co)homology. In fact we will define two different
classical and five different Kuranishi bordism theories, and five
corresponding Kuranishi cobordism theories. But these are all
variations on the same idea: {\it Kuranishi bordism groups\/}
$KB_*(Y;R)$ are spanned over $R$ by isomorphism classes $[X,\bs f]$
for $X$ a compact oriented Kuranishi space without boundary and $\bs
f:X\ra Y$ strongly smooth, up to bordism in Kuranishi spaces, and
{\it Kuranishi cobordism groups\/} $KB^*(Y;R)$ are spanned over $R$
by $[X,\bs f]$ for $X$ a compact Kuranishi space without boundary
and $\bs f:X\ra Y$ a strong submersion with $(X,\bs f)$ cooriented,
up to bordism in Kuranishi spaces.

\subsection{Classical bordism and cobordism groups}
\label{kh51}

Bordism groups were introduced by Atiyah \cite{Atiy}, and Connor
\cite[\S I]{Conn} gives a good introduction. Other useful references
on bordism and cobordism are Stong \cite{Ston} and Connor and Floyd
\cite{CoFl1,CoFl2}. Our definition is not standard, but fits in with
our later work.

\begin{dfn} Let $Y$ be an orbifold. Consider pairs $(X,f)$, where
$X$ is a compact, oriented manifold without boundary or corners, not
necessarily connected, and $f:X\ra Y$ is a smooth map. An {\it
isomorphism\/} between two such pairs $(X,f),(\ti X,\ti f)$ is an
orientation-preserving diffeomorphism $i:X\ra\ti X$ with $f=\ti f\ci
i$. Write $[X,f]$ for the isomorphism class of~$(X,f)$.

Let $R$ be a commutative ring, for instance $\Z,\Q,\R$ or $\C$. For
each $k\ge 0$, define the $k^{\it th}$ {\it bordism group\/}
$B_k(Y;R)$ of $Y$ with coefficients in $R$ to be the $R$-module of
finite $R$-linear combinations of isomorphism classes $[X,f]$ for
which $\dim X=k$, with the relations:
\begin{itemize}
\setlength{\itemsep}{0pt}
\setlength{\parsep}{0pt}
\item[(i)] $[X,f]+[X',f']=[X\amalg X',f\amalg f']$ for all
classes $[X,f],[X',f']$; and
\item[(ii)] Suppose $W$ is a compact, oriented $(k\!+\!1)$-manifold
with boundary but without (g-)corners, and $e:W\ra Y$ is smooth.
Then~$[\pd W,e\vert_{\pd W}]=0$.
\end{itemize}

If $h:Y\ra Z$ is a smooth map of orbifolds, the {\it pushforward\/}
$h_*:B_k(Y;R)\ra B_k(Z;R)$ acts by $h_*:\sum_{a\in A}\rho_a[X_a,f_a]
\mapsto\sum_{a\in A}\rho_a[X_a,h\ci f_a]$. This takes relations
(i),(ii) in $B_k(Y;R)$ to (i),(ii) in $B_k(Z;R)$, and so is
well-defined.

Define a projection $\Pi_\bo^\rsi:B_k(Y;R)\!\ra\!H^\rsi_k(Y;R)$ by
$\Pi_\bo^\rsi:\sum_{a\in A}\rho_a[X_a,f_a]\!\mapsto\!\sum_{a\in
A}\rho_a (f_a)_*([X_a])$, where $[X_a]\in H_k^\rsi(X_a;R)$ is the
fundamental class of~$X_a$.
\label{kh5def1}
\end{dfn}

The following consequence of relations Definition
\ref{kh5def1}(i),(ii) will be useful, and also works for all the
other (co)bordism theories below. Let $[X,f]\in KB_k(Y;R)$, and
write $-X$ for $X$ with the opposite orientation. Taking $W=[0,1]\t
X$ and $e:W\ra Y$ to be $f\ci\pi_X$, we have $\pd W=-\{0\}\t X\amalg
\{1\}\t X\cong -X\amalg X$ in oriented manifolds. Using (i),(ii) we
see that in $B_k(Y;R)$ we have
\e
[-X,f]=-[X,f].
\label{kh5eq1}
\e

Here is how this definition relates to those in \cite{Atiy,Conn}.
When $Y$ is a manifold and $R=\Z$, our $B_k(Y;\Z)$ is equivalent to
Connor's {\it differential bordism group\/} $D_k(Y)$, \cite[\S
I.9]{Conn}. Atiyah \cite[\S 2]{Atiy} and Connor \cite[\S I.4]{Conn}
also define {\it bordism groups\/} $MSO_k(Y)$ as for $B_k(Y;\Z)$
above, but only requiring $f:X\ra Y$ to be continuous, not smooth.
Connor \cite[Th.~I.9.1]{Conn} shows that when $Y$ is a manifold, the
natural projection $D_k(Y)\ra MSO_k(Y)$ is an isomorphism.

As in \cite[\S I.5]{Conn}, bordism is a {\it generalized homology
theory\/} in the sense of Whitehead \cite{Whit}, that is, it
satisfies all the Eilenberg--Steenrod axioms for a homology theory
except the dimension axiom. The bordism groups of a point
$MSO_*({\rm pt})$ are well understood \cite[\S I.2]{Conn}: Thom
showed that $MSO_*({\rm pt})\ot_\Z\Q$ is the free commutative
algebra over $\Q$ generated by $\ze_{4k}=[\CP^{2k}]\in MSO_{4k}({\rm
pt})$ for $k=1,2,\ldots$, and Wall and others determined the torsion
of $MSO_*({\rm pt})$, which is all of order 2. This gives some
information on bordism groups of general spaces $Y$: for any
generalized homology theory $GH_*(Y)$, there is a spectral sequence
from the singular homology $H^\rsi_*\bigl(Y;GH_*({\rm pt})\bigr)$ of
$Y$ with coefficients in $GH_*({\rm pt})$ converging to $GH_*(Y)$,
so that $GH_*({\cal S}^n)\cong H^\rsi_*\bigl({\cal S}^n;GH_*({\rm
pt})\bigr)$, for instance.

Atiyah \cite{Atiy} and Connor \cite[\S 13]{Conn} also define {\it
cobordism groups} $MSO^k(Y)$ for $k\in\Z$, which are a {\it
generalized cohomology theory} dual to bordism $MSO_k(Y)$. There is
a natural product $\cup$ on $MSO^*(Y)$, making it into a
supercommutative ring. If $Y$ is a compact oriented $n$-manifold
without boundary then \cite[Th.~3.6]{Atiy}, \cite[Th.~13.4]{Conn}
there are canonical Poincar\'e duality isomorphisms
\e
MSO^k(Y)\cong MSO_{n-k}(Y)\quad\text{for $k\in\Z$.}
\label{kh5eq2}
\e
The definition of $MSO^*(Y)$ uses homotopy theory, direct limits of
$k$-fold suspensions, and classifying spaces. There does not seem to
be a satisfactory differential-geometric definition of cobordism
groups parallel to Definition~\ref{kh5def1}.

If $Y$ is compact and oriented then from the isomorphisms
\eq{kh5eq2} and the product $\cup$ on cobordism we obtain an
associative, supercommutative {\it intersection product\/} $\bu$ on
bordism groups. (Actually, $Y$ does not need to be compact here.) We
can write this product geometrically for the groups $B_k(Y;R)$ as
follows. Suppose $[X,f]$ and $[\ti X,\ti f]$ are isomorphism classes
in $B_*(Y;R)$ such that if $f(x)=\ti f(\ti x)= y\in Y$ then $T_yY=\d
f\vert_x(T_xX)+\d\ti f\vert_{\ti x} (T_{\ti x}\ti X)$; essentially
this says that $f(X)$ and $\ti f(\ti X)$ intersect transversely in
$Y$. Then $X\t_{\smash{f,Y,\ti f}}\ti X$ is a compact, oriented
manifold with smooth $\pi:X\t_{\smash{f,Y,\ti f}}\ti X\ra Y$, and
one can show that
\e
[X,f]\bu[\ti X,\ti f]=[X\t_{\smash{f,Y,\ti f}}\ti X,\pi].
\label{kh5eq3}
\e

The bordism and cobordism theories $MSO_*(Y),MSO^*(Y)$ are
generalized (co)homology theories defined using oriented
$n$-manifolds $X$, that is, manifolds $X$ whose tangent bundle has
an $SO(n)$-structure. In a similar way, for each of the series of
Lie groups $O(n),U(n),SU(n),Sp(n)$, one can define bordism and
cobordism theories built using manifolds $X$ whose (stable) tangent
bundle has a $G$-structure for $G=O(n),U(n),SU(n)$ or $Sp(n)$. So we
can define {\it unoriented bordism} $MO_*(Y)$, {\it unoriented
cobordism} $MO^*(Y)$, {\it complex bordism} or {\it unitary bordism}
$MU_*(Y)$, {\it complex cobordism} or {\it unitary cobordism}
$MU^*(Y)$, and so on. This is explained by Connor and Floyd
\cite{CoFl1,CoFl2} and Stong \cite{Ston}. We will not discuss them.
However, note that $MU_*(Y),MU^*(Y)$ should be related to the almost
complex Kuranishi bordism groups $KB_*^\ac,KB_*^\eac(Y;R)$ and
$KB^*_\ca,KB^*_\eca(Y;R)$ defined below.

One can also repeat Definition \ref{kh5def1} using orbifolds rather
than manifolds. We choose to do this using {\it effective\/}
orbifolds, as the resulting groups will turn out to be closely
related to effective Kuranishi bordism and effective Kuranishi
homology. As far as the author knows, the next definition is new.

\begin{dfn} Let $Y$ be an orbifold, and $R$ be a commutative ring.
Define the {\it effective orbifold bordism group} $B_k^\eo(Y;R)$ of
$Y$ as in Definition \ref{kh5def1}, except that we allow $X,W$ to be
effective orbifolds rather than manifolds. {\it Pushforwards\/}
$h_*:B_k^\eo(Y;R)\ra B_k^\eo(Z;R)$ for $h:Y\ra Z$ smooth, and
projections $\Pi_\eo^\rsi:B_k^\eo(Y;R)\ra H^\rsi_k(Y;R)$, are
defined as in Definition \ref{kh5def1}. For $\Pi_\eo^\rsi$, when $R$
is not a $\Q$-algebra, it is essential that the $X_a$ are {\it
effective} orbifolds, as the fundamental class $[X]\in
H_k^\rsi(X;R)$ of a compact oriented orbifold $X$ without boundary
is defined for all $R$ when $X$ is effective, but only when
$\Q\subseteq R$ otherwise.
\label{kh5def2}
\end{dfn}

\subsection{Kuranishi bordism groups}
\label{kh52}

Motivated by \S\ref{kh42} and \S\ref{kh51}, we define:

\begin{dfn} Let $Y$ be an orbifold. Consider pairs $(X,\bs f)$,
where $X$ is a compact oriented Kuranishi space without boundary,
and $\bs f:X\ra Y$ is strongly smooth. An {\it isomorphism\/}
between two pairs $(X,\bs f),(\ti X,\bs{\ti f})$ is an
orientation-preserving strong diffeomorphism $\bs i:X\ra\ti X$ with
$\bs f=\bs{\ti f}\ci\bs i$. Write $[X,\bs f]$ for the isomorphism
class of~$(X,\bs f)$.

Let $R$ be a commutative ring. For each $k\in\Z$, define the $k^{\it
th}$ {\it Kuranishi bordism group\/} $KB_k(Y;R)$ of $Y$ with
coefficients in $R$ to be the $R$-module of finite $R$-linear
combinations of isomorphism classes $[X,\bs f]$ for which $\vdim
X=k$, with the relations:
\begin{itemize}
\setlength{\itemsep}{0pt}
\setlength{\parsep}{0pt}
\item[(i)] $[X,\bs f]+[X',\bs f']=[X\amalg X',\bs f\amalg
\bs f']$ for all classes $[X,\bs f],[X',\bs f']$; and
\item[(ii)] Suppose $W$ is a compact oriented Kuranishi space with
boundary but without (g-)corners, with $\vdim W=k+1$, and $\bs
e:W\ra Y$ is strongly smooth. Then~$[\pd W,\bs e\vert_{\pd W}]=0$.
\end{itemize}
Elements of $KB_k(Y;R)$ will be called {\it Kuranishi bordism
classes}.

Let $h:Y\ra Z$ be a smooth map of orbifolds. Define the {\it
pushforward\/} $h_*:KB_k(Y;R)\ra KB_k(Z;R)$ by $h_*:\sum_{a\in
A}\rho_a[X_a,\bs f_a]\mapsto\sum_{a\in A}\rho_a[X_a,h\ci\bs f_a]$.
This takes relations (i),(ii) in $KB_k(Y;R)$ to (i),(ii) in
$KB_k(Z;R)$, and so is well-defined. Pushforward is functorial, that
is,~$(g\ci h)_*=g_*\ci h_*$.

We also define two minor variations on~$KB_*(Y;R)$:
\begin{itemize}
\setlength{\itemsep}{0pt}
\setlength{\parsep}{0pt}
\item Define {\it effective Kuranishi bordism} $KB_*^\eb(Y;R)$ as
above, except that we require $X,W$ to be {\it effective\/}
Kuranishi spaces, as in Definition~\ref{kh3def16}.
\item Similarly, define {\it trivial stabilizers Kuranishi bordism}
$KB_*^\tr(Y;R)$ as above, except that we require $X,W$ to be
Kuranishi spaces with {\it trivial stabilizers}, that is,
$\Stab_X(p)=\{1\}$ for all~$p\in X$.
\end{itemize}
Pushforwards $h_*$ also work on $KB_*^\eb(Y;R)$ and~$KB_*^\tr(Y;R)$.
\label{kh5def3}
\end{dfn}

Our reason for introducing effective Kuranishi bordism
$KB_*^\eb(Y;R)$ is that it maps naturally to effective Kuranishi
homology $KH_*^\ef(Y;R)$, and hence to singular homology
$H^\rsi_*(Y;R)$, as we will see in \S\ref{kh53}. Since this works
over $R=\Z$, this makes $KB_*^\eb(Y;R)$ a useful tool for
considering questions of integrality of homology classes, such as
the Gopakumar--Vafa Integrality Conjecture for Gromov--Witten
invariants, to be discussed in~\S\ref{kh63}.

Our reason for defining trivial stabilizers Kuranishi bordism
$KB_*^\tr(Y;R)$ is that in \S\ref{kh53} we will show that
$KB_*^\tr(Y;R)\cong B_*(Y;R)$, and so $KB_*^\tr(Y;R)$ helps to
clarify the relationship between classical bordism theories and our
new Kuranishi bordism theories. We will also prove in \S\ref{kh57}
that $KB_*(Y;R)$ and $KB_*^\eb(Y;R)$ are generally very large. These
are nontrivial facts, as it is not obvious from Definition
\ref{kh5def3} that $KB_*,KB_*^\eb,KB_*^\tr(Y;R)$ can ever be
nonzero.

For example, if we omitted orientations in Definition \ref{kh5def3},
so that $KB_k(Y;R)$ were generated by $[X,\bs f]$ for $X$ a compact,
{\it unoriented\/} Kuranishi space, then we would have
$KB_k(Y;R)=\{0\}$ for all $Y,k$. This is because for any $[X,\bs f]$
we may take $W=([-1,1]/\{\pm 1\})\t X$, a non-orientable Kuranishi
space, and $\bs e:W\ra Y$ to be $\bs e=\bs f\ci\bs\pi_X$, and then
$[X,\bs f]=[\pd W,\bs e\vert_{\pd W}]=0$ by Definition
\ref{kh5def3}(ii) with orientations omitted. Here $[-1,1]/\{\pm 1\}$
is a compact, non-orientable 1-orbifold with boundary one point.

Using the notation of \S\ref{kh29}, we add almost complex structures
to Definition \ref{kh5def3} to define {\it almost complex Kuranishi
bordism\/} $KB^\ac_{2l}(Y;R)$. We restrict to {\it even dimensions}
$\vdim X=2l$, as almost complex structures only exist then.

\begin{dfn} Let $Y$ be an orbifold. Consider triples $\bigl(X,(\bs
J,\bs K),\bs f\bigl)$, where $X$ is a compact oriented Kuranishi
space without boundary, $(\bs J,\bs K)$ an almost complex structure
on $X$, and $\bs f:X\ra Y$ is strongly smooth. The orientation on
$X$ need not match that induced by $(\bs J,\bs K)$. An {\it
isomorphism\/} $\bs i:\bigl(X,(\bs J,\bs K),\bs f\bigr)\ra\bigl(\ti
X,(\bs{\ti J},\bs{\ti K}),\bs{\ti f}\bigr)$ is an
orientation-preserving strong diffeomorphism $\bs i:X\ra\ti X$ with
$\bs i_*(\bs J,\bs K)=(\bs{\ti J},\bs{\ti K})$ and $\bs f=\bs{\ti
f}\ci\bs i$. Write $\bigl[X,(\bs J,\bs K),\bs f\bigr]$ for the
isomorphism class of~$\bigl(X,(\bs J,\bs K),\bs f\bigl)$.

Let $R$ be a commutative ring. For each $l\in\Z$, define the {\it
almost complex Kuranishi bordism group\/} $KB^\ac_{2l}(Y;R)$ of $Y$
with coefficients in $R$ to be the $R$-module of finite $R$-linear
combinations of isomorphism classes $\bigl[X,(\bs J,\bs K),\bs
f\bigr]$ for which $\vdim X=2l$, with the relations:
\begin{itemize}
\setlength{\itemsep}{0pt}
\setlength{\parsep}{0pt}
\item[(i)] $\bigl[X,(\bs J,\bs K),\bs f\bigr]+\bigl[X',(\bs J',\bs
K'),\bs f'\bigr]=\bigl[X\amalg X',(\bs J,\bs K)\amalg(\bs J',\bs
K'),\bs f\amalg \bs f'\bigr]$ for all $\bigl[X,(\bs J,\bs K),\bs
f\bigr],\bigl[X',(\bs J',\bs K'),\bs f'\bigr]$; and
\item[(ii)] Let $W$ be a compact oriented Kuranishi space with
boundary but without (g-)corners, with $\vdim W=2l+1$, $(\bs D,\bs
J,\bs K)$ an almost CR structure on $W$, and $\bs e:W\ra Y$ a strong
submersion. Then $(\bs J,\bs K)\vert_{\pd W}$ is an almost complex
structure on $\pd W$, as in \S\ref{kh29}. We set $\bigl[\pd W,(\bs
J,\bs K)\vert_{\pd W},\bs e\vert_{\pd W}\bigr]=0$
in~$KB_{2l}^\ac(Y;R)$.
\end{itemize}

Following Definition \ref{kh5def3}, define {\it effective almost
complex Kuranishi bordism} $KB_{2l}^\eac(Y;R)$ as above, but
requiring $X,W$ to be {\it effective\/} Kuranishi spaces.

Let $h:Y\ra Z$ be a smooth map of orbifolds. Define the {\it
pushforward\/} $h_*:KB^\ac_{2l}(Y;R)\ra KB^\ac_{2l}(Z;R)$ and
$h_*:KB^\eac_{2l}(Y;R)\ra KB^\eac_{2l}(Z;R)$ by $h_*:\sum_{a\in
A}\rho_a[X_a,(\bs J_a,\bs K_a),\bs f_a]\mapsto\sum_{a\in
A}\rho_a[X_a,(\bs J_a,\bs K_a),h\ci\bs f_a]$. Pushforwards are
functorial,~$(g\ci h)_*=g_*\ci h_*$.
\label{kh5def4}
\end{dfn}

Definition \ref{kh5def4} and its cobordism analogue Definition
\ref{kh5def7} are motivated by applications to Gromov--Witten
theory, to be discussed in~\S\ref{kh62}--\S\ref{kh63}.

\subsection{Morphisms of Kuranishi bordism groups}
\label{kh53}

There are many natural projections between the bordism groups of
\S\ref{kh51}--\S\ref{kh52}, and from them to (effective) Kuranishi
homology or singular homology.

\begin{dfn} We adopt the following uniform notation for projections:
$\Pi_A^B$ represents a projection from the bordism/homology group of
type $A$ to the bordism/homology group of type $B$, where $A,B$ can
be `bo' for $B_*(Y;R)$, `eo' for $B_*^\eo(Y;R)$, `Kb' for
$KB_*(Y;R)$, `eb' for $KB_*^\eb(Y;R)$, `tr' for $KB_*^\tr(Y;R)$,
`ac' for $KB_*^\ac(Y;R)$, `eac' for $KB_*^\eac(Y;R)$, `Kh' for
$KH_*(Y;R\ot_\Z\Q)$, `ef' for $KH_*^\ef(Y;R)$, and `si' for
$H_*^\rsi(Y;R)$ or $H_*^\rsi(Y;R\ot_\Z\Q)$. With this notation:
\begin{itemize}
\setlength{\itemsep}{0pt}
\setlength{\parsep}{0pt}
\item[(a)] Define $\Pi_\bo^\eo,\Pi_\bo^\tr,\Pi_\bo^\eb,\Pi_\bo^\Kb,
\Pi_\eo^\eb,\Pi_\eo^\Kb,\Pi_\tr^\eb,\Pi_\tr^\Kb,\Pi_\eb^\Kb$ by
$\Pi_A^B:\sum_{a\in A}\rho_a\ab[X_a,\ab f_a]\ab\mapsto\sum_{a\in
A}\rho_a[X_a,f_a]$ or $\Pi_A^B:\sum_{a\in A}\rho_a[X_a,\bs
f_a]\mapsto\ab\sum_{a\in A}\rho_a[X_a,\bs f_a]$.
\item[(b)] Define $\Pi_\eac^\ac:\sum_{a\in A}\rho_a[X_a,(\bs
J_a,\bs K_a),\bs f_a]\mapsto\sum_{a\in A}\rho_a[X_a,\ab(\bs
J_a,\ab\bs K_a),\ab\bs f_a]$.
\item[(c)] Define $\Pi_\eac^\eb,\Pi_\eac^\Kb,\Pi_\ac^\Kb:\sum_{a\in
A}\rho_a[X_a,(\bs J_a,\bs K_a),\bs f_a]\mapsto\ab\sum_{a\in
A}\rho_a[X_a,\bs f_a]$.
\item[(d)] Define $\Pi_\tr^\Kh,\Pi_\eb^\Kh,\Pi_\Kb^\Kh,\Pi_\eac^\Kh,
\Pi_\ac^\Kh$ as maps into $KH_k(Y;R\ot_\Z\Q)$ by
$\Pi_A^\Kh:\sum_{a\in A}\rho_a\bigl[X_a,\bs f_a\bigr]\mapsto
\bigl[\ts\sum_{a\in A}\pi(\rho_a)[X_a,\ab\bs f_a,\bs G_a]\bigr]$ or
$\Pi_A^\Kh:\sum_{a\in A}\ab\rho_a\ab\bigl[X_a,\ab(\bs J_a,\bs
K_a),\bs f_a\bigr]\mapsto \bigl[\ts\sum_{a\in A}\pi(\rho_a)[X_a,\bs
f_a,\bs G_a]\bigr]$, where $\bs G_a$ is some choice of gauge-fixing
data for $(X_a,\bs f_a)$, which exists by Theorem \ref{kh3thm1}, and
$\pi:R\ra R\ot_\Z\Q$ is the natural morphism.
\item[(e)] Define $\Pi_\tr^\ef,\Pi_\eb^\ef,\Pi_\eac^\ef$ by
$\Pi_A^\ef:\sum_{a\in A}\rho_a\bigl[X_a,\bs f_a\bigr]\mapsto
\bigl[\ts\sum_{a\in A}\rho_a[X_a,\ab\bs f_a,\ab\ubG_a]\bigr]$ or
$\Pi_\eac^\ef:\sum_{a\in A}\rho_a\bigl[X_a,(\bs J_a,\bs K_a),\bs
f_a\bigr]\mapsto \bigl[\ts\sum_{a\in A}\rho_a[X_a,\bs
f_a,\ubG_a]\bigr]$, where $\ubG_a$ is some choice of effective
gauge-fixing data for $(X_a,\bs f_a)$, which exists by
Theorem~\ref{kh3thm5}(a).
\item[(f)] Define $\Pi_\Kb^\rsi,\Pi_\ac^\rsi$ as maps into
$H_k^\rsi(Y;R\ot_\Z\Q)$ by $\Pi_A^\rsi=(\Pi^\Kh_\rsi)^{-1}\ci
\Pi_A^\Kh$, where $(\Pi^\Kh_\rsi)^{-1}$ over $R\ot_\Z\Q$ exists by
Corollary \ref{kh4cor1}, and $\Pi_A^\Kh$ is as in~(d).

Define $\Pi_\tr^\rsi,\Pi_\eb^\rsi,\Pi_\eac^\rsi$ as maps into
$H_k^\rsi(Y;R)$ by $\Pi_A^\rsi=(\Pi^\ef_\rsi)^{-1}\ci \Pi_A^\ef$,
where $(\Pi^\ef_\rsi)^{-1}$ exists by Theorem \ref{kh4thm1} and
$\Pi_A^\ef$ is as in~(e).
\end{itemize}
Here (a)--(c) follow from the obvious implications manifold $\Ra$
trivial stabilizers Kuranishi space $\Ra$ effective Kuranishi space
$\Ra$ Kuranishi space, and manifold $\Ra$ effective orbifold $\Ra$
effective Kuranishi space.

We will show $\Pi_\Kb^\Kh$ in (d) is well-defined, the other
$\Pi_A^B$ in (d),(e) are also well-defined by the same argument. We
must prove three things:
\begin{itemize}
\setlength{\itemsep}{0pt}
\setlength{\parsep}{0pt}
\item[(i)] that $\sum_{a\in A}\rho_a[X_a,\bs f_a,\bs G_a]$ is a
cycle in~$KC_k(Y;R)$;
\item[(ii)] that $\bigl[\ts\sum_{a\in A}\rho_a[X_a,\bs f_a,\bs
G_a]\bigr]$ is independent of the choice of $\bs G_a$; and
\item[(iii)] that $\Pi_\Kb^\Kh$ takes relations in $KB_k(Y;R)$ to
equations in $KH_k(Y;R)$.
\end{itemize}
For (i), $\pd X_a=\emptyset$ by definition of $KB_k(Y;R)$, so
$\pd\bigl(\sum_{a\in A}\rho_a[X_a,\bs f_a,\bs G_a]\bigr)=0$ is
immediate from \eq{kh4eq8}. For (ii), suppose $\bs G_a,\bs G_a'$ are
alternative choices. Then by Theorem \ref{kh3thm3}, there exists
gauge-fixing data $\bs H_a$ for $\bigl([0,1]\t X_a,\bs f_a\ci\bs
\pi_{X_a}\bigr)$ with $\bs H_a\vert_{X_a\t\{0\}}=\bs G_a$ and $\bs
H\vert_{X_a\t\{1\}}=\bs G_a'$. Hence
\begin{equation*}
\pd\bigl[[0,1]\t X_a,\bs f_a\ci\bs\pi_{X_a},\bs H_a\bigr]= [X_a,\bs
f_a,\bs G_a']-[X_a,\bs f_a,\bs G_a],
\end{equation*}
by Definition \ref{kh4def2}(i),(iii), so $\sum_{a\in
A}\rho_a[X_a,\bs f_a,\bs G_a]$ and $\sum_{a\in A}\rho_a[X_a,\bs
f_a,\bs G_a']$ are homologous and define the same class
in~$KH_k(Y;R)$.

For (iii), let $[X,\bs f],[X',\bs f']$ be as in Definition
\ref{kh5def3}(i), and choose gauge-fixing data $\bs G,\bs G'$ for
$(X,\bs f),(X',\bs f')$ by Theorem \ref{kh3thm1}. Then
\begin{equation*}
[X,\bs f,\bs G]+[X',\bs f',\bs G']=[X\amalg X',\bs f\amalg \bs
f',\bs G\amalg\bs G']\quad\text{in $KC_k(Y;R)$}
\end{equation*}
by Definition \ref{kh4def2}(iii). Let $Z,\bs g$ be as in Definition
\ref{kh5def3}(ii), and choose gauge-fixing data $\bs H$ for $(Z,\bs
g)$. Then $\pd[Z,\bs g,\bs H]\!=\![\pd Z,\bs g\vert_{\pd Z},\bs
H\vert_{\pd Z}]$, so $\bigl[[\pd Z,\bs g\vert_{\pd Z},\bs
H\vert_{\pd Z}]\bigr]\!=\!0$ in $KH_k(Y;R)$. But $\bigl[[\pd Z,\bs
g\vert_{\pd Z},\bs H\vert_{\pd Z}]\bigr]\!=\!\Pi_\Kb^\Kh\bigl([\pd
Z,\bs g\vert_{\pd Z}]\bigr)$. Thus $\Pi_\Kb^\Kh$ maps relations
Definition \ref{kh5def3}(i),(ii) to equations in $KH_k(Y;R)$, and
$\Pi_\Kb^\Kh$ is well-defined.

Whenever $\Pi_A^B,\Pi_B^C$ are defined, $\Pi_A^C$ is also defined,
with $\Pi_A^C=\Pi_B^C\ci\Pi_A^B$. Also, these projections all
commute with pushforwards,~$h_*\ci\Pi_A^B=\Pi_A^B\ci h_*$.
\label{kh5def5}
\end{dfn}

\begin{figure}[htb]
\begin{footnotesize}
$\displaystyle \!\!\!\begin{gathered} \xymatrix@R=15pt@!C=24pt{
B_k(Y;R) \ar[dr]_(0.4){\Pi_\bo^\tr}^(0.55)\cong
\ar[rr]_{\Pi_\bo^\eo} && B_k^\eo(Y;R)
\ar[dr]^(0.55){\Pi_\eo^\eb}_(0.4)\cong && KB^\eac_k(Y;R)
\ar[dl]_(0.55){\Pi_\eac^\eb} \ar[rr]_{\Pi_\eac^\ac} &&
KB^\ac_k(Y;R) \ar[dl]^(0.4){\Pi_\ac^\Kb} \\
& KB_k^\tr(Y;R) \ar[rr]^(0.45){\Pi_\tr^\eb}
\ar[dr]^(0.65){\Pi_\tr^\ef} && KB_k^\eb(Y;R) \ar[rr]^{\Pi_\eb^\Kb}
\ar[dl]_(0.65){\Pi_\eb^\ef} \ar[dr]^(0.65){\Pi_\eb^\Kh} &&
KB_k(Y;R) \ar[dl]_(0.65){\Pi_\Kb^\Kh} \\
{}\,\,\,H^\rsi_k(Y;R) \ar[rr]^{\Pi_\rsi^\ef}_\cong && KH_k^\ef(Y;R)
\ar[rr]^{\Pi_\ef^\Kh} && KH_k(Y;R\ot_\Z\Q) &&
H^\rsi_k(Y;R\ot_\Z\Q).\,\,\,\,\,\,\,\,{}
\ar[ll]_{\Pi_\rsi^\Kh}^\cong }
\end{gathered}\!\!\!\!\!
$
\medskip

\centerline{{\bf NB:} $KB_k^\eac\!,KB_k^\ac(Y;R)$ defined only for $k$ even.
Also $B_k(Y;\Z)\!\cong\!MSO_k(Y)$ for manifolds.}
\end{footnotesize}
\caption{Morphisms between bordism and homology groups}
\label{kh5fig1}
\end{figure}

These morphisms are illustrated in the commutative diagram Figure
\ref{kh5fig1}. The isomorphisms in the figure are proved in Theorems
\ref{kh4thm1}, \ref{kh5thm1} and Corollary \ref{kh4cor1}. We show
that $\Pi_\bo^\tr,\Pi_\eo^\eb$ in Definition \ref{kh5def5} are
isomorphisms. The proof follows part of that of Theorem
\ref{kh4thm1} in Appendix~\ref{khB}.

\begin{thm} The morphisms $\Pi_\bo^\tr:B_k(Y;R)\ra KB_k^\tr(Y;R)$
and\/ $\Pi_\eo^\eb:B_k^\eo(Y;R)\ra KB_k^\eb(Y;R)$ are isomorphisms
for $k\ge 0$. Also $KB_k^\tr(Y;R)=KB_k^\eb(Y;R)=\{0\}$ for~$k<0$.
\label{kh5thm1}
\end{thm}

\begin{proof} First we show $\Pi_\bo^\tr$ is surjective. Let $[X,\bs
f]$ be a generator of $KB_k^\tr(Y;R)$. Choose an excellent
coordinate system $(\bs I,\bs\eta)$ for $(X,\bs f)$, which is
possible by Corollary \ref{kh3cor}, with $\bs I=\bigl(I,
(V^i,\ldots,\psi^i),f^i:i\in I,\ldots\bigr)$, where the $V^i$ can be
taken to be {\it manifolds without boundary}, since $X$ has trivial
stabilizers and is without boundary.

Follow the first part of \S\ref{khB2} to choose $C^1$ small, smooth
transverse perturbations $\ti s^i$ of $s^i$ for $i\in I$, define a
compact, oriented $k$-manifold $\ti X$, where $\ti X=\es$ if $k<0$,
and a smooth map $\ti f:\ti X\ra Y$, a compact, oriented Kuranishi
space $Z$ of virtual dimension $k+1$, and a strongly smooth map $\bs
g:Z\ra Y$. As $\pd X=\es$, we have $\pd\ti X=\es$, and the component
$Z^\pd$ of $\pd Z$ lying over $\pd X$ is also empty. So equation
\eq{khBeq16} gives $\pd Z\cong \ti X\amalg -X$ in oriented Kuranishi
spaces. Also $\bs g\vert_{\pd Z}$ agrees with $\ti f\amalg\bs f$
under this identification. Thus Definition \ref{kh5def3}(i),(ii)
give $[X,\bs f]=[\ti X,\ti f]$ in $KB_k^\tr(Y;R)$ if $k\ge 0$, and
$[X,\bs f]=0$ in $KB_k^\tr(Y;R)$ if $k<0$ so that $\ti X=\es$. When
$k<0$ this implies $KB_k^\tr(Y;R)=\{0\}$. When $k\ge 0$, $[\ti X,\ti
f]$ is also a generator of $B_k(Y;R)$, so as an element of
$KB_k^\tr(Y;R)$ it lies in the image of $\Pi_\bo^\tr$. Therefore
$\Pi_\bo^\tr$ is {\it surjective}.

To prove that $\Pi_\bo^\tr$ is injective, it is enough to show that
if $[X,f]\in B_k(Y;R)$ with $\Pi_\bo^\tr([X,f])=0$ in
$KB_k^\tr(Y;R)$, then $[X,f]=0\in B_k(Y;R)$. Since $[X,f]=0\in
KB_k^\tr(Y;R)$, there exists a compact oriented Kuranishi space $W$
of virtual dimension $k+1$ and a strongly smooth $\bs e:W\ra Y$ with
$\pd W=X$ and $\bs e\vert_{\pd W}=f$. Since $\pd W=X$ is a
$k$-manifold, and this is an open condition, $W$ is a
$(k+1)$-manifold near $\pd W$. As above, choose an excellent
coordinate system $(\bs I,\bs\eta)$ for $(W,\bs e)$, such that near
$\pd W$ the only Kuranishi neighbourhood in $\bs I$ is
$(V^{k+1},E^{k+1},s^{k+1},\psi^{k+1})$, where $V^{k+1}$ is a
$(k+1)$-manifold with $\pd V^{k+1}=X$, and $E^{k+1}\ra V^{k+1}$ is
the zero vector bundle, so that $s^{k+1}\equiv 0$. This is possible
as $W$ is a $(k+1)$-manifold near $\pd W$.

Again, follow the first part of \S\ref{khB2} to choose $C^1$ small,
smooth transverse perturbations $\ti s^i$ of $s^i$ for $i\in I$, and
define a compact, oriented $(k+1)$-manifold $\ti W$ and a smooth map
$\ti e:\ti W\ra Y$. Since $E^{k+1}\ra V^{k+1}$ is the zero vector
bundle we have $\ti s^{k+1}\equiv 0\equiv s^{k+1}$. Thus, near $\pd
W$ we do not perturb $W$ at all, so $\ti W,W$ coincide near $\pd W$,
and $\pd\ti W=\pd W=X$, and $\ti e\vert_{\pd\ti W}=e\vert_{\pd
W}=f$. Hence applying Definition \ref{kh5def1}(ii) to $(\ti W,\ti
e)$ gives $[X,f]=0$ in $B_k(Y;R)$. Thus $\Pi_\bo^\tr$ is {\it
injective}, so it is an isomorphism.

For $\Pi_\eo^\eb$ we follow the same proof, but we take the
$(V^i,\ldots,\psi^i)$ to be {\it effective Kuranishi neighbourhoods}
in the sense of Definition \ref{kh3def16}, and $\ti X,\ti W$ to be
{\it effective orbifolds} rather than manifolds. Since the
$(V^i,\ldots,\psi^i)$ are effective Kuranishi neighbourhoods, the
stabilizers of $V^i$ act trivially on the fibres of $E^i$, and so
the proof in \S\ref{khB2} that we can choose the $\ti s^i$ to be
{\it transverse} perturbations is still valid. This would not work
if we used non-effective Kuranishi neighbourhoods, and the
corresponding result for $KB_k(Y;R)$ is false.
\end{proof}

As $B_*(Y;\Z)\cong MSO_*(Y)$ for manifolds, this implies that
$KB_*^\tr(Y;\Z)\cong MSO_*(Y)$ for manifolds $Y$.

\subsection{Kuranishi cobordism and Poincar\'e duality}
\label{kh54}

(Effective) Kuranishi bordism $KB_*,KB_*^\eb(Y;R)$ is basically
(effective) Kuranishi homology $KH_*,KH_*^\ef(Y;R)$ with corners and
gauge-fixing data omitted. In \S\ref{kh44} and \S\ref{kh46} we
defined dual theories of (effective) Kuranishi cohomology. So the
obvious thing to do is to omit corners and co-gauge-fixing data from
\S\ref{kh44} and \S\ref{kh46} to define dual theories of (effective)
Kuranishi cobordism. For simplicity we restrict to orbifolds $Y$
{\it without boundary}.

\begin{dfn} Let $Y$ be an orbifold without boundary. Consider pairs
$(X,\bs f)$, where $X$ is a compact Kuranishi space without
boundary, $\bs f:X\ra Y$ is a strong submersion, and $(X,\bs f)$ is
cooriented. An {\it isomorphism\/} between two pairs $(X,\bs f),(\ti
X,\bs{\ti f})$ is a coorientation-preserving strong diffeomorphism
$\bs i:X\ra\ti X$ with $\bs f=\bs{\ti f}\ci\bs i$. Write $[X,\bs f]$
for the isomorphism class of~$(X,\bs f)$.

Let $R$ be a commutative ring. For each $k\in\Z$, define the $k^{\it
th}$ {\it Kuranishi cobordism group\/} $KB^k(Y;R)$ of $Y$ with
coefficients in $R$ to be the $R$-module of finite $R$-linear
combinations of isomorphism classes $[X,\bs f]$ for which $\vdim
X=\dim Y-k$, with the relations:
\begin{itemize}
\setlength{\itemsep}{0pt}
\setlength{\parsep}{0pt}
\item[(i)] $[X,\bs f]+[X',\bs f']=[X\amalg X',\bs f\amalg
\bs f']$ for all classes $[X,\bs f],[X',\bs f']$; and
\item[(ii)] Suppose $W$ is a compact Kuranishi space with
boundary but without \hbox{(g-)}\ab corners, with $\vdim W=\dim
Y-k+1$, and $\bs e:W\ra Y$ is a strong submersion, with $(W,\bs e)$
cooriented. Then $\bs e\vert_{\pd W}:\pd W\ra Y$ is a strong
submersion which is cooriented as in Convention \ref{kh2conv2}(a),
and we impose the relation $[\pd W,\bs e\vert_{\pd W}]=0$ in
$KB^k(Y;R)$.
\end{itemize}
Elements of $KB^k(Y;R)$ will be called {\it Kuranishi cobordism
classes}.

Let $Y,Z$ be orbifolds without boundary, and $h:Y\ra Z$ a smooth,
proper map. Motivated by pullback in Kuranishi cohomology
\eq{kh4eq21}, define the {\it pullback\/} $h^*:KB^k(Z;R)\ra
KB^k(Y;R)$ by $h^*:\sum_{a\in A}\rho_a[X_a,\bs f_a]\mapsto\sum_{a\in
A}\rho_a [Y\t_{h,Z,\bs f_a}X_a,\bs\pi_Y]$. Here the coorientation
for $(X_a,\bs f_a)$ pulls back to a coorientation for $(Y\t_{h,Z,\bs
f_a}X_a,\bs\pi_Y)$ as in Definition \ref{kh4def6}. This takes
relations (i),(ii) in $KB^k(Z;R)$ to (i),(ii) in $KB^k(Y;R)$, and so
is well-defined. Pullbacks are functorial, $(g\ci h)^*=h^*\ci g^*$.

Following Definition \ref{kh5def1} we define two minor variations
on~$KB^k(Y;R)$:
\begin{itemize}
\setlength{\itemsep}{0pt}
\setlength{\parsep}{0pt}
\item Define {\it effective Kuranishi cobordism} $KB^k_\ecb(Y;R)$ as
above, except that we require $(X,\bs f)$ and $(W,\bs e)$ to be {\it
coeffective}, as in Definition~\ref{kh3def18}.
\item Similarly, define {\it trivial stabilizers Kuranishi cobordism}
$KB^k_\trc(Y;R)$ as above, except that we require $\bs f:X\ra Y$ and
$\bs e:W\ra Y$ to be {\it coeffective} and {\it isomorphisms on
stabilizer groups}. That is, for all $p\in X$, $\bs
f_*:\Stab_X(p)\ra \Stab_Y(f(p))$ is an isomorphism, and $\Stab_X(p)$
acts trivially on the fibre of the obstruction bundle of $X$ at~$p$.
\end{itemize}
Pullbacks $h^*$ also work on $KB^*_\ecb(Y;R)$ and~$KB^*_\trc(Y;R)$.
\label{kh5def6}
\end{dfn}

\begin{dfn} Let $Y$ be an orbifold without boundary, and $R$ a
commutative ring. Consider triples $\bigl(X,(\bs L,\bs K),\bs
f\bigl)$, where $X$ is a compact Kuranishi space without boundary,
$\bs f:X\ra Y$ a strong submersion with $(X,\bs f)$ cooriented, and
$(\bs L,\bs K)$ a co-almost complex structure for $(X,\bs f)$, as in
\S\ref{kh210}. An {\it isomorphism\/} $\bs i:\bigl(X,(\bs L,\bs
K),\bs f\bigr)\ra\bigl(\ti X,(\bs{\ti L},\bs{\ti K}),\bs{\ti
f}\bigr)$ is a strong diffeomorphism $\bs i:X\ra\ti X$ with $\bs
f=\bs{\ti f}\ci\bs i$ and $\bs i_*(\bs L,\bs K)=(\bs{\ti L},\bs{\ti
K})$, which identifies the coorientations of $(X,\bs f)$ and $(\ti
X,\bs{\ti f})$. Write $\bigl[X,(\bs L,\bs K),\bs f\bigr]$ for the
isomorphism class of~$\bigl(X,(\bs L,\bs K),\bs f\bigl)$.

Define the {\it almost complex Kuranishi cobordism group\/}
$KB_\ca^{2l}(Y;R)$ of $Y$ with coefficients in $R$ for $l\in\Z$ to
be the $R$-module of finite $R$-linear combinations of isomorphism
classes $\bigl[X,(\bs L,\bs K),\bs f\bigr]$ for which $\vdim X=\dim
Y-2l$, with relations (i),(ii) combining those of Definitions
\ref{kh5def4}, \ref{kh5def6}, using co-almost CR structures in~(ii).

Define {\it effective almost complex Kuranishi cobordism}
$KB^{2l}_\eca(Y;R)$ as above, but requiring $(X,\bs f)$ and $(W,\bs
e)$ to be {\it coeffective}.

Let $h:Y\ra Z$ be a smooth, proper map of orbifolds without
boundary. Define the {\it pullback\/} $h^*:KB_\ca^{2l}(Z;R)\ra
KB_\ca^{2l}(Y;R)$ and $h^*:KB_\eca^{2l}(Z;R)\ra KB_\eca^{2l}(Y;R)$
by $h^*:\sum_{a\in A}\rho_a[X_a,(\bs L_a,\bs K_a),\bs
f_a]\mapsto\sum_{a\in A}\rho_a [Y\t_{h,Z,\bs
f_a}X_a,\ab\bs\pi_{X_a}^*(\bs L_a,\bs K_a),\bs\pi_Y]$, where
$\bs\pi_{X_a}^*(\bs L_a,\bs K_a)$ is the pullback of the co-almost
complex structure to $Y\t_ZX_a$ defined in the obvious way.
Pullbacks are functorial, $(g\ci h)^*=h^*\ci g^*$.
\label{kh5def7}
\end{dfn}

\begin{figure}[htb]
\begin{footnotesize}
$\displaystyle \!\!\!\begin{gathered} \xymatrix@R=15pt@!C=24pt{ &&&&
KB_\eca^k(Y;R) \ar[dl]_(0.55){\Pi_\eca^\ecb} \ar[rr]_{\Pi_\eca^\ca}
&&
KB_\ca^k(Y;R) \ar[dl]^(0.4){\Pi_\ca^\Kcb} \\
& KB^k_\trc(Y;R) \ar[rr]^(0.45){\Pi_\trc^\ecb}
\ar[dr]^(0.65){\Pi_\trc^\ec} && KB^k_\ecb(Y;R)
\ar[rr]^{\Pi_\ecb^\Kcb} \ar[dl]_(0.65){\Pi_\ecb^\ec}
\ar[dr]^(0.65){\Pi_\ecb^\Kch} &&
KB^k(Y;R) \ar[dl]_(0.65){\Pi_\Kcb^\Kch} \\
{}\,\,\,H_\cs^k(Y;R) \ar[rr]^{\Pi_\cs^\ec}_\cong && KH^k_\ec(Y;R)
\ar[rr]^{\Pi_\ec^\Kch} && KH^k(Y;R\ot_\Z\Q) &&
H_\cs^k(Y;R\ot_\Z\Q).\,\,\,\,\,\,\,\,{} \ar[ll]_{\Pi_\cs^\Kch}^\cong
}
\end{gathered}\!\!\!\!\!
$
\medskip

\centerline{{\bf NB:} $KB^k_\eca,KB^k_\ca(Y;R)$ defined only for
$k$ even; $\Pi_\cs^\ec$ defined only for $Y$ a manifold.}
\end{footnotesize}
\caption{Morphisms between cobordism and cohomology groups}
\label{kh5fig2}
\end{figure}

As in \S\ref{kh53}, there are many natural projections between these
cobordism groups, (effective) Kuranishi cohomology, and
compactly-supported cohomology. The definitions are the obvious
modifications of Definition \ref{kh5def5}, and we leave them to the
reader. These morphisms are illustrated in Figure \ref{kh5fig2}. As
we have not defined cobordism analogues of $B_*(Y;R)$ and
$B_*^\eo(Y;R)$, the corresponding entries in Figure \ref{kh5fig1}
are missing. Corollaries \ref{kh4cor2} and \ref{kh4cor5} say that
$\Pi_\cs^\Kch,\Pi_\cs^\ec$ are isomorphisms, supposing $Y$ a
manifold for~$\Pi_\cs^\ec$.

\begin{rem}{\rm(a)} The cobordism theories defined above are
forms of {\it compactly-supported cobordism}. One can also define
versions of Kuranishi cobordism analogous to ordinary cobordism, as
discussed for the cohomology case in Remark \ref{kh4rem3}. Since we
do not have to worry about co-gauge-fixing data for noncompact
Kuranishi spaces, this is straightforward. We just omit the
assumption that the Kuranishi space $X$ is compact in Definitions
\ref{kh5def6} and \ref{kh5def7}, and require instead that the strong
submersion $\bs f:X\ra Y$ should be {\it proper}. Pullbacks $h^*$
are then defined for all smooth $h:Y\ra Z$, not just when $h$ is
proper.

\noindent{\bf(b)} The definitions above for orbifolds $Y$ without
boundary can be extended to general orbifolds $Y$ using the ideas of
\S\ref{kh45} in the cohomology case. An important new point is that
if $\pd Y\ne\es$ then we should define $KB^k(Y;R)$ in terms of
isomorphisms $[X,\bs f]$ for which $X$ is a compact Kuranishi space
with $\vdim X=\dim Y-k$, $\bs f:X\ra Y$ is a cooriented strong
submersion, and $\pd_+^{\bs f}X=\es$, {\it without assuming\/} $\pd
X=\es$. Then we can define $\io^*:KB^k(Y;R)\ra KB^k(\pd Y;R)$
mapping $\io^*:[X,\bs f]\mapsto[\pd_-^{\bs f}X,\bs f_-]$.
\label{kh5rem1}
\end{rem}

Motivated by \S\ref{kh44}, we define {\it Poincar\'e duality maps\/}
between our Kuranishi (co)bordism groups. Basically, we just omit
(co-)gauge-fixing data from Definitions \ref{kh4def7}
and~\ref{kh4def13}.

\begin{dfn} Let $Y$ be an orbifold which is oriented, without
boundary, and of dimension $n$, and $R$ a commutative ring. Then for
$\bs f:X\ra Y$ a strong submersion, as in Definition \ref{kh2def22}
there is a 1-1 correspondence between coorientations for $(X,\bs f)$
and orientations for $X$. Define $R$-module morphisms $\Pi_\Kcb^\Kb:
KB^k(Y;R)\ra KB_{n-k}(Y;R)$ for $k\in\Z$ by $\Pi_\Kcb^\Kb:\sum_{a\in
A}\rho_a[X_a,\bs f_a]\mapsto\sum_{a\in A}\rho_a[X_a,\bs f_a]$, using
the coorientation for $\bs f_a$ from $[X_a,\bs f_a]\in KB^k(Y;R)$
and the orientation on $Y$ to determine the orientation on $X_a$ for
$[X_a,\bs f_a]\in KB_{n-k}(Y;R)$. Then $\Pi_\Kcb^\Kb$ takes
Definition \ref{kh5def3}(i),(ii) to Definition
\ref{kh5def6}(i),(ii), so it is well-defined.

Define $\Pi^\Kcb_\Kb:KB_{n-k}(Y;R)\ra KB^k(Y;R)$ by
$\Pi^\Kcb_\Kb:\sum_{a\in A}\rho_a[X_a,\ab\bs f_a]\mapsto
\sum_{a\in A}\rho_a[X_a^Y,\bs f{}_a^Y]$, where in the notation of
Definition \ref{kh4def7}, $X_a^Y$ is the fibre product
$(Y,\ka_Y)\t_{\bs\pi,Y,\bs f}X$, $\bs f{}^Y:X^Y\ra Y$ is the strong
submersion $\bs\exp\ci\bs\pi_{(Y,\ka_Y)}$, and the orientations on
$X_a$ and $Y$ determine a coorientation for $(X_a^Y,\bs f{}_a^Y)$.
Then $\Pi_\Kb^\Kcb$ takes Definition \ref{kh5def6}(i),(ii) to
Definition \ref{kh5def3}(i),(ii), so it is well-defined.

For effective Kuranishi (co)bordism, following Definition
\ref{kh4def13}, when $Y$ is an {\it effective\/} orbifold define
$\Pi_\ecb^\eb:KB^k_\ecb(Y;R)\ra KB_{n-k}^\eb(Y;R)$ as for
$\Pi_\Kcb^\Kb$, and when $Y$ is a {\it manifold}, define
$\Pi^\ecb_\eb:KB_{n-k}^\eb(Y;R)\ra KB^k_\ecb(Y;R)$ as for
$\Pi^\Kcb_\Kb$. Here for $\Pi_\ecb^\eb$ we need $Y$ to be an
effective orbifold in order that $(X_a,\bs f_a)$ coeffective implies
$X_a$ effective, and for $\Pi_\eb^\ecb$ we need $Y$ to be a manifold
so that $X_a$ effective implies $(X_a^Y,\bs f{}_a^Y)$ coeffective.

For trivial stabilizers Kuranishi (co)bordism, when $Y$ is a {\it
manifold\/} define $\Pi_\trc^\tr:KB^k_\trc(Y;R)\ra
KB_{n-k}^\tr(Y;R)$ as for $\Pi_\Kcb^\Kb,\Pi_\ecb^\eb$, and define
$\Pi^\trc_\tr:KB_{n-k}^\tr(Y;R)\ra KB^k_\trc(Y;R)$ as
for~$\Pi^\Kcb_\Kb,\Pi^\ecb_\eb$.
\label{kh5def8}
\end{dfn}

The proof of the next theorem follows those of Theorems
\ref{kh4thm3} and \ref{kh4thm5}, omitting (co-)gauge-fixing data,
and we leave it as an exercise.

\begin{thm} Let\/ $Y$ be an oriented\/ $n$-orbifold without
boundary. Then $\Pi^\Kcb_\Kb$ is the inverse of\/ $\Pi^\Kb_\Kcb$ in
Definition {\rm\ref{kh5def8},} so they are both isomorphisms. If
also $Y$ is a manifold, then $\Pi^\ecb_\eb$ is the inverse of\/
$\Pi^\eb_\ecb$ and\/ $\Pi^\trc_\tr$ is the inverse of\/
$\Pi^\tr_\trc$ in Definition {\rm\ref{kh5def8},} so they are all
isomorphisms.
\label{kh5thm2}
\end{thm}

Let $Y$ be a compact, oriented $n$-manifold without boundary. Then
we have isomorphisms $MSO^k(Y)\cong MSO_{n-k}(Y)$ by \eq{kh5eq2},
and $MSO_{n-k}(Y)\cong B_{n-k}(Y;\Z)$ from \S\ref{kh51}, and
$\Pi_\bo^\tr:B_{n-k}(Y;\Z)\ra KB_{n-k}^\tr(Y;\Z)$ by Theorem
\ref{kh5thm1}, and $\Pi^\trc_\tr:KB_{n-k}^\tr(Y;\Z)\ra
KB^k_\trc(Y;\Z)$ by Theorem \ref{kh5thm2}. Combining all these gives
an isomorphism $MSO^k(Y)\cong KB^k_\trc(Y;\Z)$, in a similar way to
Corollary \ref{kh4cor2}. We can extend this to the case $Y$ not
oriented using orientation bundles, as in Definition \ref{kh4def8},
giving:

\begin{cor} If\/ $Y$ is a compact manifold without boundary, there
is a natural isomorphism~$MSO^*(Y)\cong KB^*_\trc(Y;\Z)$.
\label{kh5cor1}
\end{cor}

This shows that $KB^*_\trc(Y;\Z)$ is a {\it differential-geometric
realization\/} of cobordism, very different to the conventional
algebraic topological definition \cite{Atiy}. We can also drop the
assumption that $Y$ is compact, if we replace ordinary cobordism
$MSO^*(Y)$ by {\it compactly-supported cobordism\/}~$MSO^*_\cs(Y)$.

\begin{rem} We can also consider Poincar\'e duality for {\it almost
complex\/} Kuranishi (co)bordism $KH^\ac_{2*},KB_\ca^{2*}(Y;R)$, and
{\it effective almost complex\/} Kuranishi (co)bordism
$KH^\eac_{2*},KB_\eca^{2*}(Y;R)$. But to define the corresponding
Poincar\'e duality maps, we need some extra data: we must choose an
{\it almost complex structure\/} $J$ on $Y$. (So in particular,
$\dim Y$ must be even, $\dim Y=2m$ say.)

To define $\Pi_\ca^\ac$, for $[X,(\bs L,\bs K),\bs f]\in
KB^{2l}_\ca(Y;R)$ we must define an almost complex structure on $X$
using $J$ and the co-almost complex structure $(\bs L,\bs K)$ for
$(X,\bs f)$. Similarly, to define $\Pi_\ac^\ca$, for $[X,(\bs J,\bs
K),\bs f]\in KB_{2m-2l}^\ac(Y;R)$ we must define a co-almost complex
structure for $(X^Y,\bs f{}^Y)$ using $J$ and the almost complex
structure $(\bs J,\bs K)$ for $X$. Both of these can be done, though
there are some arbitrary choices involved, as for the construction
of $(\bs J,\bs K)\t_Y(\bs L',\bs K')$ in \S\ref{kh210}. One can then
extend Theorem \ref{kh5thm2} to show that $KB^{2l}_\ca(Y;R)\cong
KB_{2m-2l}^\ac(Y;R)$, and $KB^{2l}_\eca(Y;R)\cong
KB_{2m-2l}^\eac(Y;R)$ for $Y$ a manifold.
\label{kh5rem2}
\end{rem}

\subsection{Products on Kuranishi (co)bordism}
\label{kh55}

Motivated by \S\ref{kh47}, we define cup, cap and intersection
products on Kuranishi (co)bordism. Let $Y$ be an orbifold without
boundary, not necessarily oriented. Define the {\it cup product\/}
$\cup:KB^k(Y;R)\t KB^l(Y;R)\ra KB^{k+l}(Y;R)$ by
\e
\raisebox{-4pt}{\begin{Large}$\displaystyle\Bigl[$\end{Large}}
\sum_{a\in A}\rho_a\bigl[X_a,\bs f_a\bigr]
\raisebox{-4pt}{\begin{Large}$\displaystyle\Bigr]$\end{Large}}\!\cup\!
\raisebox{-4pt}{\begin{Large}$\displaystyle\Bigl[$\end{Large}}
\sum_{b\in B}\si_b\bigl[\ti X_b,\bs{\ti f}_b\bigr]
\raisebox{-4pt}{\begin{Large}$\displaystyle\Bigr]$\end{Large}}\!=\!\!
\sum_{a\in A,\; b\in B\!\!\!\!\!\!}\!\!\rho_a\si_b
\bigl[X_a\t_{\bs f_a,Y,\bs{\ti f}_b}\ti X_b,\bs\pi_Y\bigr],
\label{kh5eq5}
\e
for $A,B$ finite and $\rho_a,\si_b\in R$. Here $X_a\t_{\bs
f_a,Y,\bs{\ti f}_b}\ti X_b$ is the fibre product of \S\ref{kh26},
which is a compact Kuranishi space without boundary as $X_a,\ti X_b$
are, and $\bs\pi_Y: X_a\t_{\smash{\bs f_a,Y,\bs{\ti f}_b}}\ti X_b\ra
Y$ is the projection from the fibre product, which is a strong
submersion as $\bs f_a,\bs{\ti f}_b$ are, and the coorientations on
$(X_a,\bs f_a)$ and $(\ti X_b,\bs{\ti f}_b)$ induce a coorientation
on $(X_a\t_Y\ti X_b,\bs\pi_Y)$ as in Convention \ref{kh2conv2}(b).
Since $\vdim X_a=\dim Y-k$, $\vdim\ti X_b=\dim Y-l$ we have $\vdim
X_a\t_Y\ti X_b=\dim Y-(k+l)$, so $\cup$ does map~$KB^k(Y;R)\t
KB^l(Y;R)\ra KB^{k+l}(Y;R)$.

To show $\cup$ is well-defined we must show that \eq{kh5eq5} takes
relations (i),(ii) in $KB^k(Y;R),KB^l(Y;R)$ to (i),(ii) in
$KB^{k+l}(Y;R)$. For (i) this holds as
\begin{align*}
\bigl([X,\bs f]&+[X',\bs f']\bigr)\cup[\ti X,\bs{\ti f}]=
\bigl[X\t_{\bs f,Y,\bs{\ti f}}\ti X,\bs\pi_Y\bigr]+\bigl[X'\t_{\bs
f',Y,\bs{\ti f}}\ti X,\bs\pi_Y\bigr]\\
&=\bigl[(X\t_{\bs f,Y,\bs{\ti f}}\ti X)\amalg (X'\t_{\bs
f',Y,\bs{\ti f}}\ti X),\bs\pi_Y\amalg\bs\pi_Y\bigr]\\
&=\bigl[(X\amalg X')\t_{\bs f\amalg \bs f',Y,\bs{\ti f}}\ti X,\bs
\pi_Y \bigr]=[X\amalg X',\bs f\amalg \bs f']\cup[\ti X,\bs{\ti f}],
\end{align*}
and similarly for $[\ti X,\bs{\ti f}]\cup\bigl([X,\bs f]+[X',\bs
f']\bigr)$. For (ii) it holds as
\begin{align*}
[\pd Z,\bs g\vert_{\pd Z}]\cup[\ti X,\bs{\ti f}]=\bigl[\pd Z\t_{\bs
g\vert_{\pd Z},Y,\bs{\ti f}}\ti X,\bs\pi_Y\bigr]=\bigl[\pd(Z\t_{\bs
g,Y,\bs{\ti f}}\ti X),\bs\pi_Y\vert_{\pd(\cdots)}\bigr]=0,
\end{align*}
using Proposition \ref{kh2prop2}(a) and $\pd\ti X=\emptyset$, and
similarly for~$[\ti X,\bs{\ti f}]\cup[\pd Z,\bs g\vert_{\pd Z}]$.

If $[X,\bs f]\in KB^k(Y;R)$, $[\ti X,\bs{\ti f}]\in KB^l(Y;R)$ then
equation \eq{kh2eq22} gives $[X,\bs f]\cup[\ti X,\bs{\ti
f}]=(-1)^{kl}[\ti X,\bs{\ti f}]\cup[X,\bs f]$, that is, $\cup$ is
{\it supercommutative}. Also \eq{kh2eq23} implies that $\cup$ is
{\it associative}. If $Y$ is compact then using the trivial
coorientation for $\id_Y:Y\ra Y$, we have $[Y,\id_Y]\in KB^0(Y;R)$,
which is the {\it identity\/} for $\cup$. Thus, $KB^*(Y;R)$ is a
{\it graded, supercommutative, associative $R$-algebra, with
identity\/} if $Y$ is compact, and {\it without identity\/}
otherwise.

We can also mix Kuranishi bordism and cobordism, and define a {\it
cap product\/} $\cap:KB_k(Y;R)\t KB^l(Y;R)\ra KB_{k-l}(Y;R)$ by
\eq{kh5eq5}, where now $[X_a,\bs f_a]\in KB_k(Y;R)$ so that $\bs
f_a$ is strongly smooth and $X_a$ oriented, and $[\ti X_b,\bs{\ti
f}_b]\in KB^l(Y;R)$ so that $\bs{\ti f}_b$ is a strong submersion
and $(\ti X_b,\bs{\ti f}_b)$ cooriented, and $X_a\t_{\smash{\bs
f_a,Y,\bs{\ti f}_b}}\ti X_b$ is well-defined as $\bs{\ti f}_b$ is a
strong submersion, and the orientation for $X_a$ and coorientation
for $(\ti X_b,\bs{\ti f}_b)$ combine to give an orientation for
$X_a\t_Y\ti X_b$ as in Convention \ref{kh2conv2}(c). This has the
associativity property $\al\cap(\be\cup\ga)=(\al\cap\be)\cap\ga$ for
$\al\in KB_*(Y;R)$ and $\be,\ga\in KB^*(Y;R)$. Thus, {\it Kuranishi
bordism $KB_*(Y;R)$ is a module over Kuranishi
cobordism}~$KB_*(Y;R)$.

As for Kuranishi (co)homology in \S\ref{kh47}, cup and cap products
are compatible with pullbacks and pushforwards. That is, if $Y,Z$
are orbifolds without boundary, and $h:Y\ra Z$ a smooth, proper map,
and $\al\in KB_*(Y;R)$ and $\be,\ga\in KB^*(Z;R)$ then as in
\eq{kh4eq52} we have
\e
h^*(\be\cup\ga)=h^*(\be)\cup h^*(\ga) \quad\text{and}\quad
h_*(\al\cap h^*(\be))=h_*(\al)\cap\be.
\label{kh5eq6}
\e
If $Z$ is compact then $Y$ is compact as $h$ is proper, and
$h^*([Z,\id_Z])=[Y,\id_Y]$. Thus, $h^*$ is an $R$-{\it algebra
morphism}.

If $Y$ is an oriented $n$-orbifold without boundary, then using
$\cup$ and the isomorphism $KB^k(Y;R)\cong KB_{n-k}(Y;R)$ of Theorem
\ref{kh5thm2} we obtain an {\it intersection product\/}
$\bu:KB_k(Y;R)\t KB_l(Y;R)\ra KB_{k+l-n}(Y,R)$, which is also
supercommutative (with degrees shifted by $n$) and associative. From
\eq{kh5eq3} we see that $\Pi_\bo^\Kb:B_*(Y;R)\ra KB_*(Y;R)$
intertwines intersection products $\bu$ on~$B_*(Y;R),KB_*(Y;R)$.

All of the above also holds without change for cup, cap and
intersection products on {\it effective Kuranishi (co)bordism\/}
$KB^*_\ecb,KB_*^\eb(Y;R)$, and on {\it trivial stabilizers Kuranishi
(co)bordism\/} $KB^*_\trc,KB_*^\tr(Y;R)$, except that we must
suppose $Y$ is a manifold to define $\bu$ on
$KB_*^\eb,KB_*^\tr(Y;R)$.

We can also generalize these ideas to (effective) almost complex
Kuranishi (co)bordism. Define $R$-bilinear maps
$\cup:KB^{2k}_\ca(Y;R)\t KB^{2l}_\ca(Y;R)\ra KB^{2(k+l)}_\ca(Y;R)$
and $\cup:KB^{2k}_\eca(Y;R)\t KB^{2l}_\eca(Y;R)\ra
KB^{2(k+l)}_\eca(Y;R)$ by
\begin{equation*}
\bigl[X,(\bs L,\bs K),\bs f\bigr]\cup\bigl[\ti X,(\bs{\ti L},\bs{\ti
K}),\bs{\ti f}\bigr]=\bigl[\smash{X\t_{\bs f,Y,\bs{\ti f}}\ti X,(\bs
L,\bs K)\t_Y(\bs{\ti L},\bs{\ti K}),\bs\pi_Y}\bigr],
\end{equation*}
generalizing \eq{kh5eq5}, where $\smash{(\bs L,\bs K)\t_Y(\bs{\ti
L},\bs{\ti K})}$ is the fibre product of co-almost complex
structures from \S\ref{kh210}. It is associative and
supercommutative. If $Y$ is compact then $[Y,(\bs 0,\bs 0),\id_Y]$
is the identity for $\cup$, where $(\bs 0,\bs 0)$ is the trivial
co-almost complex structure for $(Y,\id_Y)$, involving only almost
complex structures on zero vector bundles. Thus $KB^*_\ca(Y;R)$ is a
{\it supercommutative ring}.

Similarly, define $\cap:KB_{2k}^\ac(Y;R)\t KB^{2l}_\ca(Y;R)\ra
KB_{2(k-l)}^\ac(Y;R)$ and $\cap:KB_{2k}^\eac(Y;R)\t
KB^{2l}_\eca(Y;R)\ra KB_{2(k-l)}^\eac(Y;R)$ by
\begin{equation*}
\bigl[X,(\bs J,\bs K),\bs f\bigr]\cap\bigl[\ti X,(\bs{\ti L},\bs{\ti
K}),\bs{\ti f}\bigr]=\bigl[\smash{X\t_{\bs f,Y,\bs{\ti f}}\ti X,(\bs
J,\bs K)\t_Y(\bs{\ti L},\bs{\ti K}),\bs\pi_Y}\bigr].
\end{equation*}
They satisfy $\al\cap(\be\cup\ga)=(\al\cap\be)\cap\ga$ and
$h_*(\al\cap h^*(\be))=h_*(\al)\cap\be$ as above.

All of the morphisms between bordism and homology groups defined in
\S\ref{kh53}, and their analogues for cobordism and cohomology
groups, intertwine cup, cap and intersection products in the obvious
way. So, for instance, $\Pi_\Kcb^\Kch:KB^*(Y;R)\ra
KH^*(Y;R\ot_\Z\Q)$ satisfies $\Pi_\Kcb^\Kch(\al\cup\be)=
\Pi_\Kcb^\Kch(\al)\cup\Pi_\Kcb^\Kch(\be)$ for all $\al,\be\in
KB^*(Y;R)$. In each case this is either obvious, or follows from
Theorem~\ref{kh4thm6}.

\begin{rem} {\bf(a)} As in \S\ref{kh48} for Kuranishi (co)homology,
each of the matching pairs of (co)bordism theories $KB_*,KB^*(Y;R)$,
$KB_*^\eb,KB^*_\ecb(Y;R)$, $KB_*^\tr,\ab KB^*_\trc(Y;R)$,
$KB_{2*}^\ac,KB^{2*}_\ca(Y;R)$ and $KB_{2*}^\eac,KB^{2*}_\eca(Y;R)$
can be combined into a single {\it bivariant theory}. In some sense
this works `at the (co)chain level'.
\smallskip

\noindent{\bf(b)} The author expects that one can generalize the
proof in the classical case to show that each of our various
different Kuranishi bordism theories is a {\it generalized homology
theory\/} in the sense of Whitehead \cite{Whit} in the category of
smooth orbifolds, that is, it satisfies all the Eilenberg--Steenrod
axioms for homology except the dimension axiom, and that the
Kuranishi cobordism theories are the corresponding {\it generalized
compactly-supported cohomology theory}. To do this, for the case of
Kuranishi bordism, we should introduce a notion of {\it relative
Kuranishi bordism\/} $KB_*(Y,Z;R)$ when $Y$ is an orbifold and
$Z\subseteq Y$ is an open set, and prove that it lives in an exact
sequence
\begin{equation*}
\cdots \ra KB_k(Z;R)\ra KB_k(Y;R)\ra KB_k(Y,Z;R)\,{\buildrel\pd\over
\longra}\, KB_{k-1}(Z;R)\ra \cdots,
\end{equation*}
and satisfies the excision axiom, and various other natural
properties. We leave this as an exercise for the interested reader.
\label{kh5rem3}
\end{rem}

\subsection{Orbifold strata and operators $\Pi^{\Ga,\rho}$}
\label{kh56}

Any orbifold $V$ may be written as $V=\coprod_\Ga V^\Ga$, where the
disjoint union is over isomorphism classes of finite groups $\Ga$,
and $V^\Ga$ is the subset of $v\in V$ with stabilizer group
$\Stab_V(v)$ isomorphic to $\Ga$. This is the rough idea behind our
definition of the {\it orbifold strata\/} of $V$. However, for our
purposes we must modify this idea in three ways:
\begin{itemize}
\setlength{\itemsep}{0pt}
\setlength{\parsep}{0pt}
\item[(i)] $V^\Ga$ above can be a union of components of different
dimensions, so we pass to a finer stratification
$V=\coprod_{\Ga,\rho}V^{\Ga,\rho}$, where $\rho$ is an {\it
isomorphism class of nontrivial representations of\/} $\Ga$, and the
normal bundle of $V^{\Ga,\rho}$ in $V$ has fibre the representation
$\rho$, so that~$\dim V^{\Ga,\rho}=\dim V-\dim\rho$.
\item[(ii)] $V^{\Ga,\rho}$ above is not a {\it closed\/} subset of
$V$. To make it closed, we replace $V^{\Ga,\rho}$ by a subset of
points $v\in V$ such that $\Stab_V(v)$ {\it has a subgroup
isomorphic to} $\Ga$. Once we do this, the disjoint union
$V=\coprod_{\Ga,\rho}V^{\Ga,\rho}$ no longer holds.
\item[(iii)] Making the replacement (ii), but regarding
$V^{\Ga,\rho}$ as a subset of $V$, we find that $V^{\Ga,\rho}$ is
{\it not an orbifold}, since several pieces of $V^{\Ga,\rho}$ may
intersect in $V$ locally to give singularities. To make
$V^{\Ga,\rho}$ into an orbifold, we must lift to an immersion
$\io^{\Ga,\rho}:V^{\Ga,\rho}\ra V$, which may take finitely many
points to one point. Then $V^{\Ga,\rho}$ being closed is replaced by
$\io^{\Ga,\rho}$ being {\it proper}. (This is very similar to the
idea in \S\ref{kh21} that if $X$ is a manifold with (g-)corners,
then to make $\pd X$ a manifold, we cannot regard $\pd X$ as a
subset of $X$, but instead we have a finite immersion $\io:\pd X\ra
X$.)
\end{itemize}
Here then is our definition of orbifold strata:

\begin{dfn} Let $\Ga$ be a finite group, and consider
(finite-dimensional) real representations $(W,\om)$ of $\Ga$, that
is, $W$ is a finite-dimensional real vector space and
$\om:\Ga\ra\Aut(W)$ is a group morphism. Call $(W,\om)$ a {\it
trivial representation} if $\om\equiv\id_W$, that is, if
$\Fix(\om(\Ga))=W$. Call $(W,\om)$ a {\it nontrivial representation}
if it has no positive-dimensional trivial subrepresentation, that
is, if $\Fix(\om(\Ga))=\{0\}$. Then every $\Ga$-representation
$(W,\om)$ has a unique decomposition $W=W^{\rm t}\op W^{\rm nt}$ as
the direct sum of a trivial representation $(W^{\rm t},\om^{\rm t})$
and a nontrivial representation $(W^{\rm nt},\om^{\rm nt})$,
where~$W^{\rm t}=\Fix(\om(\Ga))$.

Let $(I_0,\io_0),\ldots,(I_N,\io_N)$ be representatives of the
isomorphism classes of {\it irreducible\/} $\Ga$-representations,
where $(I_0,\io_0)=(\R,1)$ is the trivial irreducible
representation, and $(I_i,\io_i)$ is nontrivial for $i=1,\ldots,N$.
Then every $\Ga$-representation $(W,\om)$ is isomorphic to a direct
sum $\bigop_{i=0}^Na_i(I_i,\io_i)$ with $a_0,\ldots,a_N\in\N$, where
$(W,\om)$ is trivial if $a_1=\cdots=a_N=0$, and $(W,\om)$ is
nontrivial if $a_0=0$. Two $(W,\om),(W,\om')$ are isomorphic if and
only if~$(a_0,\ldots,a_N)=(a_0',\ldots,a_N')$.

Hence, {\it isomorphism classes of nontrivial\/
$\Ga$-representations} are in 1-1 correspondence with $N$-tuples
$(a_1,\ldots,a_N)$, with each $a_i\in\Z$. We will often use $\rho$
to denote such an isomorphism class. The {\it dimension\/}
$\dim\rho$ is $\sum_{i=1}^Na_i\dim I_i$. If $(W^{\rm nt},\om^{\rm
nt})$ is a nontrivial $\Ga$-representation, we will write $[(W^{\rm
nt},\om^{\rm nt})]$ for the isomorphism class of nontrivial
$\Ga$-representations containing $(W^{\rm nt},\om^{\rm nt})$, so
that $[(W^{\rm nt},\om^{\rm nt})]=\rho$ means $(W^{\rm nt},\om^{\rm
nt})$ lies in the isomorphism class~$\rho$.

Now let $V$ be an orbifold, of dimension $n$. Each $v\in V$ has a
{\it stabilizer group\/} $\Stab_V(v)$. In this section we consider
the tangent space $T_vV$ to be an $n$-dimensional vector space with
a representation $\tau_v$ of $\Stab_V(v)$. Let
$\la:\Ga\ra\Stab_V(v)$ be an injective group morphism, so that
$\la(\Ga)$ is a subgroup of $\Stab_V(v)$ isomorphic to $\Ga$. Hence
$\tau_v\ci\la:\Ga\ra\Aut(T_vV)$ is a $\Ga$-representation, and we
can split $T_vV=(T_vV)^{\rm t}\op(T_vV)^{\rm nt}$ into trivial and
nontrivial $\Ga$-representations, and form the isomorphism class
$\bigl[(T_vV)^{\rm nt},(\tau_v\ci\la)^{\rm nt}\bigr]$.

As a set, define the {\it orbifold stratum} $V^{\Ga,\rho}$ to be
\e
\begin{split}
V^{\Ga,\rho}=\bigl\{\Stab_V(v)\cdot(v,\la):\text{$v\in V$,
$\la:\Ga\ra\Stab_V(v)$ is an injective}& \\
\text{group morphism, $\bigl[(T_vV)^{\rm nt},(\tau_v\ci\la)^{\rm
nt}\bigr]=\rho$}&\bigr\},
\end{split}
\label{kh5eq7}
\e
where $\Stab_V(v)$ acts on pairs $(v,\la)$ by $\si:(v,\la)\mapsto
(v,\la^\si)$, where $\la^\si:\Ga\ra\Stab_V(v)$ is given by
$\la^\si(\ga)=\si\la(\ga)\si^{-1}$. Define a map $\io^{\Ga,\rho}:
V^{\Ga,\rho}\ra V$ by~$\io^{\Ga,\rho}:\Stab_V(v)\cdot(v,\la)\mapsto
v$.
\label{kh5def9}
\end{dfn}

\begin{prop} In Definition {\rm\ref{kh5def9},} $V^{\Ga,\rho}$
naturally has the structure of an orbifold of dimension
$n-\dim\rho,$ and\/ $\io^{\Ga,\rho}$ lifts to a proper, finite
immersion.
\label{kh5prop1}
\end{prop}

\begin{proof} Let $(U,\De,\psi)$ be an orbifold chart on $V$. We
will construct a corresponding orbifold chart on $V^{\Ga,\rho}$. Let
$\la:\Ga\ra\De$ be an injective group morphism. As $\De$ acts on $U$
by diffeomorphisms, composing with $\la$ gives an action of $\Ga$ on
$U$ by diffeomorphisms. Write $\Fix(\la(\Ga))$ for the fixed point
set of this action in $U$. It is closed, and a disjoint union of
embedded submanifolds of $U$, possibly with different dimensions. If
$u\in\Fix(\la(\Ga))$ then $\Ga$ has a representation on $T_uU$, so
we have a splitting $T_uU=(T_uU)^{\rm t}\op(T_uU)^{\rm nt}$ into
trivial and nontrivial $\Ga$-representations, with $(T_uU)^{\rm
t}=T_u\Fix(\la(\Ga))$. For each isomorphism class of nontrivial
$\Ga$-representations $\rho$, write $\Fix(\la(\Ga))^\rho$ for the
subset of $u\in\Fix(\la(\Ga))$ with $\bigl[(T_uU)^{\rm
nt}\bigr]=\rho$. Then $\Fix(\la(\Ga))^\rho$ is a closed, embedded
submanifold of $U$, of dimension $n-\dim\rho$, since~$\dim U=\dim
V=n$.

Define a manifold $U^{\Ga,\rho}$ of dimension $n-\dim\rho$ to be
\e
U^{\Ga,\rho}=\ts\coprod_{\text{injective group morphisms
$\la:\Ga\ra\De$}}\Fix(\la(\Ga))^\rho.
\label{kh5eq8}
\e
Note that this is a {\it disjoint union}, not a union:
$\Fix(\la(\Ga))^\rho$ and $\Fix(\la'(\Ga))^\rho$ may intersect in
$U$ for distinct $\la,\la'$, but we do not identify such
intersections in $U^{\Ga,\rho}$. Define $i^{\Ga,\rho}:U^{\Ga,\rho}
\ra U$ to be the inclusion $\Fix(\la(\Ga))^\rho\hookra U$ on each
component $\Fix(\la(\Ga))^\rho$ in \eq{kh5eq8}. As each
$\Fix(\la(\Ga))^\rho$ is a closed, embedded submanifold, and the
disjoint union \eq{kh5eq8} is finite, we see that $i^{\Ga,\rho}$ is
a {\it proper, finite immersion}.

Define an action of $\De$ on $U^{\Ga,\rho}$ by $\de\in\De$ takes $u$
in $\Fix(\la(\Ga))^\rho$ to $\de\cdot u$ in
$\Fix(\la^\de(\Ga))^\rho$, where $\la^\de:\Ga\ra\De$ is defined by
$\la^\de(\ep)=\de\la(\ep)\de^{-1}$, as for $\la^\si$ in Definition
\ref{kh5def9}. Then $\De$ acts on $U^{\Ga,\rho}$ by diffeomorphisms,
and $i^{\Ga,\rho}$ is $\De$-equivariant. Define
$\psi^{\Ga,\rho}:U^{\Ga,\rho}/\De\ra V^{\Ga,\rho}$ as follows: let
$\la:\Ga\ra\De$ be an injective group morphism and $u\in
\Fix(\la(\Ga))^\rho\subseteq U^{\Ga,\rho}$. Set $v=\psi(\De\cdot u)$
in $V$. Then we have a natural isomorphism $\Stab_V(v)\cong\{\de\in
\De:\de\cdot u=u\}$, and $\la$ maps $\Ga$ to the subgroup
$\{\de\in\De:\de\cdot u=u\}$ of $\De$ as $u\in\Fix(\la(\Ga))$. Thus
we can regard $\la$ as an injective group morphism
$\la:\Ga\ra\Stab_V(v)$. With this identification,
define~$\psi^{\Ga,\rho}\bigl(\De u)=\Stab_V(v)\cdot(v,\la)$.

It is elementary to show that $\psi^{\Ga,\rho}$ is well-defined, and
$(U^{\Ga,\rho},\De,\psi^{\Ga,\rho})$ is an orbifold chart on
$V^{\Ga,\rho}$ with image $(\io^{\Ga,\rho})^{-1}(\Im\psi)$. Also, if
$(U,\De,\psi)$ and $(U',\De',\psi')$ are orbifold charts on $V$
compatible on their overlaps, then $(U^{\Ga,\rho},\De,\ab
\psi^{\Ga,\rho}),\ab(U^{\Ga,\rho\prime},\De',\psi^{\Ga,\rho
\prime})$ are compatible on their overlaps. Hence this system of
compatible orbifold charts on $V^{\Ga,\rho}$ constructed from the
system of compatible orbifold charts on $V$ gives $V^{\Ga,\rho}$ the
structure of an {\it orbifold}, of dimension $n-\dim\rho$. Also
$\io^{\Ga,\rho}$ is represented in orbifold charts by the
$i^{\Ga,\rho}$ above, which are proper, finite immersions, so
$\io^{\Ga,\rho}$ lifts to a proper, finite immersion.
\end{proof}

\begin{rem}{\bf(a)} If $v^{\Ga,\rho}=\Stab_V(v)\cdot(v,\la)\in
V^{\Ga,\rho}$ then the stabilizer group $\Stab_{V^{\Ga,\rho}}
(v^{\Ga,\rho})$ of $v^{\Ga,\rho}$ in $V^{\Ga,\rho}$ is
\begin{equation*}
\{\si\in\Stab_V(v):\la^\si\!=\!\la\}\!=\!\{\si\!\in\!\Stab_V(v):
\la(\ga)\si\!=\!\si\la(\ga)\;
\,\,\forall\ga\!\in\!\Ga\}\!=\!C(\la(\Ga)),
\end{equation*}
the centralizer of $\la(\Ga)$ in $\Stab_V(v)$. The induced morphism
$\io^{\Ga,\rho}_*:\Stab_{V^{\Ga,\rho}}
(v^{\Ga,\rho})\ab\ra\Stab_V(v)$ is just the inclusion
$C(\la(\Ga))\hookra\Stab_V(v)$. Thus $\io^{\Ga,\rho}$ is {\it
injective on stabilizer groups}.
\smallskip

\noindent{\bf(b)} Write $\Aut(\Ga)$ for the group of isomorphisms
$\ze:\Ga\ra\Ga$. If $\ze\in\Aut(\Ga)$ and $\la:\Ga\ra\Stab_V(v)$ is
an injective group morphism then $\la\ci\ze^{-1}:\Ga\ra\Stab_V(v)$
is an injective group morphism. Also $\Aut(\Ga)$ acts naturally on
the set of isomorphism classes of nontrivial $\Ga$-representations,
$\rho\mapsto\ze\cdot\rho$. If $\ze\in\Aut(\Ga)$ then there is a
natural diffeomorphism $\ze_*:V^{\Ga,\rho}\ra V^{\Ga,\ze\cdot\rho}$
acting by $\Stab_V(v)\cdot(v,\la)\mapsto\Stab_V(v)\cdot(v,\la\ci
\ze^{-1})$, with $\io^{\Ga,\rho}\equiv\io^{\Ga,\ze\cdot\rho}
\ci\ze_*$.
\smallskip

\noindent{\bf(c)} If $V$ is an {\it effective\/} orbifold then
$\Stab_V(v)$ acts effectively on $T_vV$ for all $v\in V$. It is then
easy to see that $V^{\Ga,\rho}=\es$ unless $\rho$ is the isomorphism
class of an {\it effective} nontrivial $\Ga$-representation.
Conversely, if $V^{\Ga,\rho}=\es$ for all finite groups $\Ga$ and
isomorphism classes of non-effective nontrivial
$\Ga$-representations $\rho$, then $V$ is an effective orbifold.
\smallskip

\noindent{\bf(d)} If $V$ is an {\it oriented\/} orbifold, the
orbifold strata $V^{\Ga,\rho}$ {\it need not be oriented, or even
orientable}. Here is an example. Take $V$ to be the 4-orbifold
$(\RP^3\t\R)/\Z_2$, where the generator $\si$ of $\Z_2$ acts on
$\RP^3\t\R$ by
\e
\si:\bigl([x_0,x_1,x_2,x_3],y\bigr)\mapsto
\bigl([x_0,x_1,x_2,-x_3],-y\bigr).
\label{kh5eq9}
\e
Then $\Fix(\si)$ is the disjoint union of $\bigl\{\bigl([x_0,x_1,
x_2,0],0\bigr):[x_0,x_1,x_2]\in\RP^2\bigr\}$, a copy of $\RP^2$, and
the single point $\bigl([0,0,0,1],0\bigr)$. Hence $V$ has two
nontrivial orbifold strata, $V^{\Z_2,\rho_2}$, a copy of $\RP^2$,
and $V^{\Z_2,\rho_4}$, a point, where $\rho_2,\rho_4$ are the
isomorphism classes of nontrivial representations of $\Z_2$ of
dimensions 2,4. Now $\RP^k$ is orientable if and only if $k$ is odd.
Using this, it is easy to show that $V$ is orientable, but
$V^{\Z_2,\rho_2}$ is not.

Here is a sufficient condition on $\rho$ for $V^{\Ga,\rho}$ to be
oriented whenever $V$ is oriented. Let $(W,\om)$ be a nontrivial
$\Ga$-representation with $[(W,\om)]=\rho$, and suppose that $W$
{\it does not admit any orientation-reversing automorphisms as a
$\Ga$-representation}. Then choosing an orientation for $W$
determines an orientation for the fibres of the normal bundle of the
immersed submanifold $\io^{\Ga,\rho}(V^{\Ga,\rho})$ in $V$, since
these fibres are isomorphic to $(W,\om)$ as $\Ga$-representations.
Hence an orientation for $V$ determines one for $V^{\Ga,\rho}$. This
condition fails for any nonzero representation of $\Ga=\Z_2$. But if
$\md{\Ga}$ is odd, the condition holds for all $\rho$, as in
Definition \ref{kh5def11} below.
\smallskip

\noindent{\bf(e)} Suppose $V$ is a $2m$-orbifold and $J$ is an {\it
almost complex structure} on $V$, as in Definition \ref{kh2def23}.
Then for each $v\in V$, $J\vert_v$ gives $T_vV$ the structure of a
{\it complex vector space}, and as $J\vert_v$ is invariant under
$\Stab_V(v)$ we can regard $(T_vV,\tau_v)$ as a {\it complex
representation} of $\Ga$. We can therefore repeat Definition
\ref{kh5def9} and Proposition \ref{kh5prop1} using isomorphism
classes $\rho$ of {\it complex\/} nontrivial $\Ga$-representations,
rather than real representations, to define $V^{\Ga,\rho}$ and
$\io^{\Ga,\rho}$. Also $J$ on $V$ restricts to an almost complex
structure $J^{\Ga,\rho}$ on~$V^{\Ga,\rho}$.

Observe that $V$ and $V^{\Ga,\rho}$ are automatically orientable,
with the orientations induced by $J,J^{\Ga,\rho}$. Thus, by working
with almost complex structures, we avoid the problem with
non-oriented orbifold strata discussed in (d). This is an important
reason for working with almost complex Kuranishi bordism.

Perhaps surprisingly, the same thing works for {\it almost CR
structures}. Suppose $V$ is a $(2m+1)$-orbifold and $(D,J)$ is an
almost CR structure on $V$, as in Definition \ref{kh2def24}. Suppose
that $V$ is {\it orientable}. Then for each $v\in V$, we have a
codimension one subspace $D\vert_v$ of $T_vV$, an almost complex
structure $J\vert_v$ on $D\vert_v$, and we can choose an orientation
on $T_vV$. All three structures are invariant under the action of
$\Stab_V(v)$. As $J\vert_v$ induces an orientation on $D\vert_v$ we
have an orientation on $T_vV/D\vert_v\cong\R$ which is fixed by
$\Stab_V(v)$. Therefore $\Stab_V(v)$ acts trivially on the
quotient~$(T_vV)/(D\vert_v)$.

It follows that if $\la:\Ga\ra\Stab_V(v)$ is an injective group
morphism, then in the induced splitting $T_vV=(T_vV)^{\rm
t}\op(T_vV)^{\rm nt}$, $(T_vV)^{\rm nt}$ is a vector subspace of
$D\vert_v$, which must be invariant under $J\vert_v$. So we can
regard $(T_vV)^{\rm nt}$ as a {\it complex} representation of $\Ga$.
This is what we need to repeat Definition \ref{kh5def9} and
Proposition \ref{kh5prop1} using isomorphism classes $\rho$ of
complex nontrivial $\Ga$-representations. Also $(D,J)$ on $V$
restricts to an almost CR structure $(D^{\Ga,\rho},J^{\Ga,\rho})$
on~$V^{\Ga,\rho}$.

When working with complex representations, we will use $\dim\rho$ to
mean the {\it complex\/} dimension of the representations, not the
real dimension.
\smallskip

\noindent{\bf(f)} Suppose $V$ is an orbifold with boundary, but {\it
without (g-)corners}. Then $\io:\pd V\ra V$ is an embedding, and
induces isomorphisms $\io_*:\Stab_{\pd V}(v,B)\ra\Stab_V(v)$ on
stabilizer groups. It is easy to show that restricting to orbifold
strata commutes with taking boundaries, that
is,~$\pd(V^{\Ga,\rho})=(\pd V)^{\Ga,\rho}$.

However, as we explained in Remark \ref{kh2rem2}(b), if $V$ has
(g-)corners then $\io:\pd V\ra V$ is an immersion, and the induced
maps $\io_*:\Stab_{\pd V}(v,B)\ra\Stab_V(v)$ are injective, but may
not be surjective. Because of this, {\it for orbifolds with corners
or g-corners, restricting to orbifold strata does not commute with
taking boundaries}, that is, $\pd(V^{\Ga,\rho})$ may not agree
with~$(\pd V)^{\Ga,\rho}$.

This leads to an important difference between Kuranishi (co)bordism
and (effective) Kuranishi (co)homology. In Kuranishi (co)bordism,
information from orbifold strata $X^{\Ga,\rho}$ of generators
$[X,\bs f]$ survives in $KB_*,KB^*(Y;R)$, as we show below, since we
do not allow Kuranishi spaces with (g-)corners in defining
$KB_*,KB^*(Y;R)$. But in (effective) Kuranishi (co)homology, where
we do allow Kuranishi spaces with (g-)corners, all information from
orbifold strata is lost. (In fact relation Definition
\ref{kh4def2}(iv) would also kill information from orbifold strata,
but this does not apply to {\it effective} Kuranishi (co)homology.)
\label{kh5rem4}
\end{rem}

We generalize Definition \ref{kh5def9} to Kuranishi spaces $X$.
Since we must consider the representation of $\Ga$ on the
obstruction bundle of $X$, we must replace $\rho$ by an isomorphism
class of {\it virtual\/} nontrivial $\Ga$-representations.

\begin{dfn} Let $\Ga$ be a finite group. A {\it virtual
representation} of $\Ga$ is a formal difference
$(W,\om)\ominus(W',\om')$ of two $\Ga$-representations
$(W,\om),(W',\om')$. The {\it dimension} of
$(W,\om)\ominus(W',\om')$ is $\dim W-\dim W'$. We call
$(W,\om)\ominus(W',\om')$ {\it nontrivial\/} if $(W,\om),(W',\om')$
are nontrivial $\Ga$-representations. We say that two virtual
$\Ga$-representations $(W_1,\om_1)\ominus(W_1',\om_1')$ and
$(W_2,\om_2)\ominus(W_2',\om_2')$ are {\it isomorphic} if there
exist $\Ga$-representations $(W_3,\om_3),(W_4,\om_4)$ with
$(W_1,\om_1)\op(W_3,\om_3)\cong(W_2,\om_2)\op(W_4,\om_4)$ and
$(W_1',\om_1')\op(W_3,\om_3)\cong(W_2',\om_2')\op(W_4,\om_4)$.

Let $(I_0,\io_0),\ldots,(I_N,\io_N)$ be representatives of the
isomorphism classes of {\it irreducible\/} $\Ga$-representations,
where $(I_0,\io_0)=(\R,1)$ is the trivial irreducible
representation, and $(I_i,\io_i)$ is nontrivial for $i=1,\ldots,N$.
Then every virtual $\Ga$-representation $(W,\om)\ominus(W',\om')$ is
isomorphic to $\bigop_{i=0}^Na_i(I_i,\io_i)$ with
$a_0,\ldots,a_N\in\Z$, where $\bigop_{i=0}^Na_i(I_i,\io_i)$ is a
shorthand for $\bigl(\bigop_{i:a_i>0}a_i(I_i,\io_i)\bigr)\ominus
\bigl(\bigop_{i:a_i<0}(-a_i)(I_i,\io_i)\bigr)$, and two virtual
representations are isomorphic if and only if they yield the same
$a_0,\ldots,a_N$.

Hence, {\it isomorphism classes of virtual\/ $\Ga$-representations}
are in 1-1 correspondence with $(N+1)$-tuples $(a_0,a_1,\ldots,a_N)$
with $a_i\in\Z$, and isomorphism classes of virtual nontrivial
$\Ga$-representations are in 1-1 correspondence with $N$-tuples
$(a_1,\ldots,a_N)$ with $a_i\in\Z$. We often write $\rho$ for such
an isomorphism class. The {\it dimension\/} of $\rho$ is
$\dim\rho=\sum_{i=0}^Na_i\dim I_i$. Equivalently, an isomorphism
class of virtual $\Ga$-representatives is an element of
$K_0(\text{mod-}\Ga)$, the Grothendieck group of the abelian
category mod-$\Ga$ of finite-dimensional real $\Ga$-representations.

Now let $X$ be a Kuranishi space, $p\in X$ and $(V_p,\ldots,\psi_p)$
be a Kuranishi neighbourhood in the germ at $p$. Set
$v=\psi_p^{-1}(p)$ in $V$. The stabilizer group
$\Stab_X(p)\cong\Stab_{V_p}(v)$ has natural representations on the
finite-dimensional vector spaces $T_vV_p$ and $E_p\vert_v$. Let
$\la:\Ga\ra\Stab_{V_p}(v)$ be an injective group morphism. Then
$\Ga$ acts on $T_vV_p$ and $E_p\vert_v$ via $\la$, so we may split
into trivial and nontrivial $\Ga$-representations
$T_vV_p=(T_vV_p)^{\rm t}\op(T_vV_p)^{\rm nt}$ and
$E_p\vert_v=(E_p\vert_v)^{\rm t}\op(E_p\vert_v)^{\rm nt}$. Hence
$(T_vV_p)^{\rm nt}\ominus(E_p\vert_v)^{\rm nt}$ is a virtual
nontrivial $\Ga$-representation, with isomorphism
class~$\bigl[(T_vV_p)^{\rm nt}\ominus(E_p\vert_v)^{\rm nt}\bigr]$.

As a set, define the {\it orbifold stratum} $X^{\Ga,\rho}$ to be
\e
\begin{split}
X^{\Ga,\rho}=\bigl\{\Stab_X(p)\cdot(p,\la):\text{$p\in X$,
$\la:\Ga\ra\Stab_X(p)$ is an injective}& \\
\text{group morphism, $\bigl[(T_vV_p)^{\rm
nt}\ominus(E_p\vert_v)^{\rm nt}\bigr]=\rho$}\bigr\},
\end{split}
\label{kh5eq10}
\e
where to interpret the condition $\bigl[(T_vV_p)^{\rm
nt}\ominus(E_p\vert_v)^{\rm nt}\bigr]=\rho$ we choose
$(V_p,\ldots,\psi_p)$ in the germ at $p$, set $v=\psi_p^{-1}(p)$ and
identify $\Stab_{V_p}(v)$ with $\Stab_X(p)$, so that $\la$ makes
$T_vV_p$ and $E_p\vert_v$ into $\Ga$-representations. The condition
is independent of the choice of $(V_p,\ldots,\psi_p)$, by definition
of equivalence in the germ. Define a map $\io^{\Ga,\rho}:
X^{\Ga,\rho}\ra X$ by~$\io^{\Ga,\rho}:\Stab_X(p)\cdot(p,\la)\mapsto
p$.
\label{kh5def10}
\end{dfn}

We have the following analogue of Proposition \ref{kh5prop1}. The
proof generalizes that of Proposition \ref{kh5prop1}, and we leave
it as an exercise. The most important new point to note is that
coordinate transformations $(\phi_{pq},\hat\phi_{pq})$ preserve the
condition on virtual nontrivial $\Ga$-representations. Suppose
$v_q\in V_{pq}\subseteq V_q$ with $s_q(v_q)=0$,
$v_p=\phi_{pq}(v_q)$, and $\la_q:\Ga\ra\Stab_{V_q}(v)$ is an
injective group morphism. Then $(\phi_{pq})_*:\Stab_{V_q}(v_q)
\ra\Stab_{V_p}(v_p)$ is an isomorphism, as $\phi_{pq}$ is an
embedding, and Definition \ref{kh2def12}(e) implies that
$(T_{v_p}V_p)^{\rm nt}\ominus(E_p\vert_{v_p})^{\rm nt}$ and
$(T_{v_q}V_q)^{\rm nt}\ominus(E_q\vert_{v_q})^{\rm nt}$ are
isomorphic as virtual $\Ga$-representations. Hence
$\bigl[(T_{v_p}V_p)^{\rm nt}\ominus(E_p\vert_{v_p})^{\rm
nt}\bigr]=\rho$ if and only if~$\bigl[(T_{v_q}V_q)^{\rm
nt}\ominus(E_q\vert_{v_q})^{\rm nt}\bigr]=\rho$.

\begin{prop} In Definition {\rm\ref{kh5def10},} $X^{\Ga,\rho}$
naturally has the structure of a Kuranishi space of virtual
dimension $\vdim X-\dim\rho,$ and\/ $\io^{\Ga,\rho}$ lifts to a
proper, finite, strongly smooth map
$\bs\io^{\Ga,\rho}:X^{\Ga,\rho}\ra X$.
\label{kh5prop2}
\end{prop}

\begin{rem} The analogues of Remark \ref{kh5rem4}(a)--(f) hold for
Kuranishi spaces. In particular, for (c), if $X$ is an {\it
effective} Kuranishi space then in Definition \ref{kh5def10}
$T_vV_p$ is an {\it effective\/} $\Ga$-representation, so
$(T_vV_p)^{\rm nt}$ is effective, and $E_p\vert_v$ is a {\it
trivial\/} $\Ga$-representation, so $(E_p\vert_v)^{\rm nt}=0$. Hence
$X^{\Ga,\rho}=\es$ unless $\rho$ is the isomorphism class of an {\it
effective nontrivial\/ $\Ga$-representation}, not just a virtual
one. For $\Ga\ne\{1\}$ this gives $\dim\rho>0$, so that $\vdim
X^{\Ga,\rho}<\vdim X$. If $X$ is an {\it orientable}, effective
Kuranishi space, we can also exclude the case $\dim\rho=1$, since
$\Ga$ cannot act on $T_vV_p$ preserving orientations such that
$\dim(T_vV_p)^{\rm nt}=1$, because no finite group has an
orientation-preserving, nontrivial representation on $\R$. Hence, if
$X$ is an orientable, effective Kuranishi space then
$X^{\Ga,\rho}=\es$ for $\Ga\ne\{1\}$ unless $\dim\rho\ge 2$, so
that~$\vdim X^{\Ga,\rho}\le\vdim X-2$.

Also, for (e), if $X$ is a Kuranishi space without boundary, $(\bs
J,\bs K)$ is an almost complex structure on $X$, $\Ga$ a finite
group, and $\rho$ an isomorphism class of virtual nontrivial {\it
complex\/} $\Ga$-representations, then we can define the orbifold
stratum $X^{\Ga,\rho}$, strongly smooth $\bs\io^{\Ga,\rho}:
X^{\Ga,\rho}\ra X$, and an almost complex structure $(\bs
J^{\Ga,\rho},\bs K^{\Ga,\rho})$ on $X^{\Ga,\rho}$. Similarly, if $X$
is an {\it orientable\/} Kuranishi space with boundary but without
(g-)corners and $(\bs D,\bs J,\bs K)$ an almost CR structure on $X$,
we define $X^{\Ga,\rho},\bs\io^{\Ga,\rho}$, and an almost CR
structure $(\bs D^{\Ga,\rho},\bs J^{\Ga,\rho},\bs K^{\Ga,\rho})$
on~$X^{\Ga,\rho}$.
\label{kh5rem5}
\end{rem}

When $\md{\Ga}$ is odd, given an orientation on $X$, we define
orientations on~$X^{\Ga,\rho}$.

\begin{dfn} Let $\Ga$ be a finite group. Consider the following
condition on $\Ga$: {\it no nontrivial representation $(W,\om)$ of\/
$\Ga$ should admit an orientation-reversing automorphism}, as in
Remark \ref{kh5rem4}(d). It is sufficient to apply this condition to
the nontrivial irreducible representations $(I_1,\io_1),\ldots,
(I_N,\io_N)$ of $\Ga$. Using Schur's lemma we can show that
$(I_i,\io_i)$ admits an orientation-reversing automorphism if and
only if $\dim I_i$ is odd, and then $-1$ is such an automorphism.
One can prove using the character theory of finite groups that $\dim
I_i$ is even for $i=1,\ldots,N$ if and only if $\md{\Ga}$ is odd, so
the condition is equivalent to $\md{\Ga}$ odd.

Suppose $\md{\Ga}$ is odd. Choose an orientation on $I_i$ for
$i=1,\ldots,N$. Let $(W,\om)$ be a nontrivial $\Ga$-representation.
Then $(W,\om)$ is isomorphic to $\bigop_{i=1}^Na_i(I_i,\io_i)$ for
$a_1,\ldots,a_N\in\N$. Thus the orientations on $I_i$ induce an
orientation on $W$. This orientation on $W$ is independent of the
choice of isomorphism $W\cong\bigop_{i=1}^Na_iI_i$, since by the
condition above $(W,\om)$ admits no orientation-reversing
automorphisms. Thus, we have chosen orientations on all nontrivial
$\Ga$-representations, which are consistent with direct sums.

Now let $X$ be an oriented Kuranishi space, and $\rho$ an
isomorphism class of virtual nontrivial $\Ga$-representations. We
shall define an {\it orientation\/} on the orbifold stratum
$X^{\Ga,\rho}$. As for orbifold charts in the proof of Proposition
\ref{kh5prop1}, if $(V,E,s,\psi)$ is a compatible Kuranishi
neighbourhood on $X$ then we construct a compatible Kuranishi
neighbourhood $(V^{\Ga,\rho},E^{\Ga,\rho},s^{\Ga,\rho},
\psi^{\Ga,\rho})$ on $X^{\Ga,\rho}$, and the strongly smooth map
$\bs\io^{\Ga,\rho}$ is represented by a smooth map
$(\io^{\Ga,\rho},\hat\io^{\Ga,\rho}):(V^{\Ga,\rho},\ldots,
\psi^{\Ga,\rho})\ra(V,\ldots,\psi)$. Here points of $V^{\Ga,\rho}$
are of the form $v^{\Ga,\rho}=\Stab_V(v)\cdot(v,\la)$, where
$v=\io^{\Ga,\rho}\bigl(\Stab_V(v)\cdot(v,\la)\bigr)$ in $V$,
$\la:\Ga\ra\Stab_V(v)$ is an injective group morphism.

Now $\la$ makes $T_vV$ and $E\vert_v$ into $\Ga$-representations,
with $T_{v^{\Ga,\rho}}V^{\Ga,\rho}=(T_vV)^{\rm t}$ and
$E^{\Ga,\rho}\vert_{v^{\Ga,\rho}}=(E\vert_v)^{\rm t}$. Also
$T_vV=(T_vV)^{\rm t}\op(T_vV)^{\rm nt}$ and
$E\vert_v=(E\vert_v)^{\rm t}\op(E\vert_v)^{\rm nt}$, and
$(T_vV)^{\rm nt},(E\vert_v)^{\rm nt}$ both have orientations as they
are nontrivial $\Ga$-representations, and the orientation on $X$
gives an orientation on $T_vV\op E\vert_v$. Thus, there is a unique
orientation on $T_{v^{\Ga,\rho}}V^{\Ga,\rho}\op
E^{\Ga,\rho}\vert_{v^{\Ga,\rho}}$ such that
\begin{equation*}
\bigl(T_vV\op E\vert_v\bigr)\cong\bigl(T_{v^{\Ga,\rho}} V^{\Ga,\rho}
\op E^{\Ga,\rho}\vert_{v^{\Ga,\rho}}\bigr)\op(T_vV)^{\rm nt}\op(E
\vert_v)^{\rm nt}
\end{equation*}
holds in oriented vector spaces.

Using the fact that nontrivial $\Ga$-representations have even
dimension, so that changing their order in a direct sum preserves
orientations, we find that the compatibility condition in
\S\ref{kh27} for orientations under coordinate changes
$(\phi_{pq},\hat\phi_{pq})$ on $X$ implies the corresponding
compatibility on $X^{\Ga,\rho}$. Hence the orientations on
$T_{v^{\Ga,\rho}}V^{\Ga,\rho}\op E^{\Ga,\rho}\vert_{v^{\Ga,\rho}}$
above induce an orientation on $X^{\Ga,\rho}$, as we want. This
orientation does depend on the choice of orientations for
$I_1,\ldots,I_N$; changing the orientations of $I_i$ by $\ep_i=\pm
1$ changes the orientation of $X^{\Ga,\rho}$ by
$\prod_{i=1}^N\ep_i^{a_i}$, where $\rho=\bigl[\bigop_{i=1}^N
a_i(I_i,\io_i)\bigr]$ for~$a_1,\ldots,a_N\in\Z$.
\label{kh5def11}
\end{dfn}

We define operators $\Pi^{\Ga,\rho},\Pi^{\Ga,\rho}_\ac$ on (almost
complex) Kuranishi bordism.

\begin{dfn} Let $\Ga$ be a finite group with $\md{\Ga}$ odd. Let
$(I_1,\io_1),\ldots,(I_N,\io_N)$ be the nontrivial irreducible
representations of $\Ga$, and choose orientations on $I_i$ for
$i=1,\ldots,N$. Let $Y$ be an orbifold and $R$ a commutative ring.
For each isomorphism class $\rho$ of virtual nontrivial
$\Ga$-representations and each $k\in\Z$, define an operator
$\Pi^{\Ga,\rho}:KB_k(Y;R)\ra KB_{k-\dim\rho}(Y;R)$ by
\e
\ts\Pi^{\Ga,\rho}:\sum_{a\in A}\rho_a\bigl[X_a,\bs
f_a\bigr]\longmapsto\sum_{a\in A}\rho_a\bigl[X_a^{\Ga,\rho},\bs
f_a\ci\bs\io_a^{\Ga,\rho}\bigr].
\label{kh5eq11}
\e
Here $X_a^{\Ga,\rho}$ and $\bs\io_a^{\Ga,\rho}:X_a^{\Ga,\rho}\ra
X_a$ are as in Definition \ref{kh5def10} for $X_a$. Since $X_a$ is
compact and $\bs\io_a^{\Ga,\rho}$ is proper, $X_a^{\Ga,\rho}$ is a
{\it compact\/} Kuranishi space. We define the orientation on
$X_a^{\Ga,\rho}$ as in Definition~\ref{kh5def11}.

Equation \eq{kh5eq11} takes relations Definition
\ref{kh5def3}(i),(ii) in $KB_k(Y;R)$ to Definition
\ref{kh5def3}(i),(ii) in $KB_{k-\dim\rho}(Y;R)$, and so is
well-defined. For (ii), this depends on the Kuranishi space analogue
of Remark \ref{kh5rem4}(f), that is, as $W$ is a Kuranishi space
with boundary but {\it without (g-)corners}, we have
$\pd(W^{\Ga,\rho})=(\pd W)^{\Ga,\rho}$. Note that since $\md{\Ga}$
is odd, $\dim\rho$ is even, as in Definition \ref{kh5def11},
so~$\dim\rho\in 2\Z$.

If $Y,Z$ are orbifolds and $h:Y\ra Z$ is smooth then the
$\Pi^{\Ga,\rho}$ clearly commute with pushforwards, that is,
$h_*\ci\Pi^{\Ga,\rho}=\Pi^{\Ga,\rho}\ci h_*:KB_k(Y;R)\ra
KB_{k-\dim\rho}(Z;R)$.
\label{kh5def12}
\end{dfn}

\begin{dfn} Let $Y$ be an orbifold, $R$ a commutative ring,
$l$ an integer, $\Ga$ a finite group, and $\rho$ an isomorphism
class of virtual nontrivial {\it complex\/} $\Ga$-representations.
Define $\Pi^{\Ga,\rho}_\ac:KB_{2l}^\ac(Y;R)\ra
KB_{2l-2\dim\rho}^\ac(Y;R)$ by
\e
\ts\Pi^{\Ga,\rho}:\sum_{a\in A}\rho_a\bigl[X_a,(\bs J_a,\bs K_a),\bs
f_a\bigr]\longmapsto\sum_{a\in A}\rho_a\bigl[X_a^{\Ga,\rho},(\bs
J_a^{\Ga,\rho},\bs K_a^{\Ga,\rho}),\bs f_a\ci
\bs\io_a^{\Ga,\rho}\bigr].
\label{kh5eq12}
\e
Here $X_a^{\Ga,\rho}$, $\bs\io_a^{\Ga,\rho}:X_a^{\Ga,\rho}\ra X_a$
and $(\bs J_a^{\Ga,\rho},\bs K_a^{\Ga,\rho})$ are as in Remark
\ref{kh5rem5} for $X_a$. Since $X_a$ is compact and
$\bs\io_a^{\Ga,\rho}$ is proper, $X_a^{\Ga,\rho}$ is a {\it
compact\/} Kuranishi space. We take $\dim\rho$ to be the {\it
complex} dimension of the virtual representation, so $2\dim\rho$ is
the real dimension.

To define the orientation on $X_a^{\Ga,\rho}$, note that $X_a$ is
oriented by Definition \ref{kh5def4} (this orientation need not
coincide with the orientation induces by $(\bs J,\bs K)$). We
construct the orientation on $X_a^{\Ga,\rho}$ by the method of
Definition \ref{kh5def11}, but instead of choosing orientations on
$I_1,\ldots,I_N$ in order to define orientations on $(T_vV)^{\rm
nt},(E\vert_v)^{\rm nt}$, we note that $(T_vV)^{\rm nt},
(E\vert_v)^{\rm nt}$ are {\it complex\/} vector spaces, and so have
natural orientations.

Equation \eq{kh5eq12} takes relations Definition
\ref{kh5def4}(i),(ii) in $KB_{2l}^\ac(Y;R)$ to Definition
\ref{kh5def4}(i),(ii) in $KB_{2l-2\dim\rho}^\ac(Y;R)$, and so is
well-defined. For (ii), given an oriented Kuranishi space $W$ with
boundary and an almost CR structure $(\bs D,\bs J,\bs K)$ on $W$, we
define $W^{\Ga,\rho}$, $\bs\io^{\Ga,\rho}:W^{\Ga,\rho}\ra W$ and
$(\bs D^{\Ga,\rho},\bs J^{\Ga,\rho},\bs K^{\Ga,\rho})$ as in Remark
\ref{kh5rem5} for $W$, and we construct the orientation on
$W^{\Ga,\rho}$ as for $X_a^{\Ga,\rho}$ above, again noting that
$(T_vV)^{\rm nt},(E\vert_v)^{\rm nt}$ are {\it complex\/} vector
spaces by the argument of Remark \ref{kh4rem5}(e).

The $\Pi^{\Ga,\rho}_\ac$ commute with pushforwards
$h_*:KB_{2*}^\ac(Y;R)\ra KB_{2*}^\ac(Z;R)$.
\label{kh5def13}
\end{dfn}

\begin{rem} We have not defined operators $\Pi^{\Ga,\rho}$ on {\it
effective} Kuranishi bordism $\Pi^\eb_*(Y;R)$, as there does not
seem to be a good way to make $X_a^{\Ga,\rho}$ into an effective
Kuranishi space. In particular, using Remark \ref{kh5rem4}(a) one
can show that if $(V_p,\ldots,\psi_p)$ is a sufficiently small
Kuranishi neighbourhood in the germ at $p$ in $X_a^{\Ga,\rho}$ then
for all $v\in V_p$ the stabilizer group $\Stab_{V_p}(v)$ contains
the centre $C(\Ga)$ of $\Ga$, so if $C(\Ga)\ne\{1\}$ then $V_p$ is
not an effective orbifold, and $X_a^{\Ga,\rho}$ not an effective
Kuranishi space. There is no point in defining $\Pi^{\Ga,\rho}$ on
{\it trivial stabilizers} Kuranishi bordism $\Pi^\tr_*(Y;R)$, as
they would be zero for all~$\Ga\ne\{1\}$.

However, we can consider the compositions $\Pi^{\Ga,\rho}\ci
\Pi_\ef^\Kb:\Pi^\eb_k(Y;R)\ra KB_{k-\dim\rho}(Y;R)$ and
$\Pi^{\Ga,\rho}\ci\Pi_\tr^\Kb:\Pi^\tr_k(Y;R)\ra KB_{k-\dim\rho}
(Y;R)$. By Remark \ref{kh5rem5} we have $\Pi^{\Ga,\rho}\ci
\Pi_\ef^\Kb=0$ unless $\rho$ is the isomorphism class of an {\it
effective nontrivial\/ $\Ga$-representation}, not just a virtual
one. For $\Ga\ne\{1\}$ this implies that $\Pi^{\Ga,\rho}\ci
\Pi_\ef^\Kb=0$ unless $\dim\rho>0$. We can also exclude the case
$\dim\rho=1$ as in Remark \ref{kh5rem5}, since the Kuranishi spaces
are oriented. Hence $\Pi^{\Ga,\rho}\ci \Pi_\ef^\Kb=0$ unless
$\dim\rho\ge 2$. For {\it trivial stabilizers} Kuranishi bordism we
clearly have $\Pi^{\Ga,\rho}\ci\Pi_\tr^\Kb=0$ unless $\Ga=\{1\}$.
These imply that $\Pi_\ef^\Kb: \Pi^\eb_*(Y;R)\ra KB_*(Y;R)$ and
$\Pi_\tr^\Kb:\Pi^\tr_*(Y;R)\ra KB_*(Y;R)$ are far from surjective,
their images lie in small subspaces of~$KB_*(Y;R)$.

In the same way, $\Pi^{\Ga,\rho}\ci\Pi_\eac^\ac:\Pi^\eac_{2l}
(Y;R)\ra KB_{2l-2\dim\rho}^\ac(Y;R)$ is zero unless $\rho$ is the
isomorphism class of an {\it effective nontrivial complex\/
$\Ga$-representation}, not just a virtual one, so for $\Ga\ne\{1\}$,
$\Pi^{\Ga,\rho}_\ac\ci\Pi_\eac^\ac=0$ unless $\dim\rho>0$.
\label{kh5rem6}
\end{rem}

All of the above has analogues for the corresponding Kuranishi
cobordism groups, which we now briefly explain. The starting point
is the notion of orbifold strata for a submersion $f:V\ra Y$ of
orbifolds~$V,Y$.

\begin{dfn} Let $V,Y$ be orbifolds, and $f:V\ra Y$ a submersion.
Use the notation of Definition \ref{kh5def9}. Generalizing equation
\eq{kh5eq7}, as a set, define the {\it orbifold stratum}
$V^{\Ga,\rho,f}$ of $(V,f)$ to be
\ea
V^{\Ga,\rho,f}&\!=\!\bigl\{\Stab_V(v)\cdot(v,\la):v\in V,\;
\la:\Ga\!\ra\!\Ker\bigl(f_*:\Stab_V(v)\!\ra\!\Stab_Y(f(v))\bigr)
\nonumber\\
&\quad\,\,\text{is an injective group morphism, $\bigl[(T_vV)^{\rm
nt},(\tau_v\ci\la)^{\rm nt}\bigr]=\rho$}\bigr\}.
\label{kh5eq13}
\ea
That is, $V^{\Ga,\rho,f}$ is the subset of $\Stab_V(v)\cdot(v,\la)$
in $V^{\Ga,\rho}$ with $\la(\Ga)\subseteq\Ker f_*\subseteq
\Stab_V(v)$. Define $\io^{\Ga,\rho,f}:V^{\Ga,\rho,f}\ra V$
by~$\io^{\Ga,\rho,f}=\io^{\Ga,\rho}\vert_{V^{\Ga,\rho,f}}$.

Then $V^{\Ga,\rho,f},\io^{\Ga,\rho,f}$ have two important
properties, easily verified from the definition of the orbifold
structure on $V^{\Ga,\rho}$ in the proof of Proposition
\ref{kh5prop1}. Firstly, $V^{\Ga,\rho,f}$ is an open and closed
subset of $V^{\Ga,\rho}$. Hence $V^{\Ga,\rho,f}$ is an {\it
orbifold}, by Proposition \ref{kh5prop1}, and
$\io^{\Ga,\rho,f}:V^{\Ga,\rho,f}\ra V$ lifts to a {\it proper,
finite immersion}.

Secondly, $f\ci\io^{\Ga,\rho,f}:V^{\Ga,\rho,f}\ra Y$ is a {\it
submersion}. This is not immediately obvious, as the restriction of
a submersion to a submanifold is generally no longer a submersion.
To see it, let $v^{\Ga,\rho,f}=\Stab_V(v)\cdot(v,\la)$ in
$V^{\Ga,\rho,f}$, and set $y=f(v)$. Then $\la$ makes $T_vV$ into a
$\Ga$-representation, so we may split into trivial and nontrivial
representations $T_vV=(T_vV)^{\rm t}\op(T_vV)^{\rm nt}$, and
$T_{v^{\Ga,\rho,f}}V^{\Ga,\rho,f}\cong(T_vV)^{\rm t}$. We have
linear maps $\d f\vert_v:T_vV\ra T_yY$, which is {\it surjective} as
$f$ is a submersion, and $\d(f\ci\io^{\Ga,\rho,f})\vert_{v^{\Ga,
\rho,f}}:T_{v^{\Ga,\rho,f}}V^{\Ga,\rho,f}\ra T_yY$, and under the
identification $T_{v^{\Ga,\rho,f}}V^{\Ga,\rho,f}\cong(T_vV)^{\rm
t}$, we have $\d(f\ci\io^{\Ga, \rho,f})\vert_{v^{\Ga,\rho,f}}\cong
\bigl(\d f\vert_v\bigr)\vert_{(T_vV)^{\rm t}}$.

Now $\d f\vert_v:T_vV\ra T_yY$ is equivariant under the actions of
$\Stab_V(v)$ on $T_vV$ and $\Stab_Y(y)$ on $T_yY$ and the group
morphism $f_*:\Stab_V(v)\ra\Stab_Y(v)$. Since $\la$ maps
$\Ga\ra\Ker(f_*)$, it follows that $\Ga$ acts trivially on $T_y$,
and so $(T_vV)^{\rm nt}$ lies in the kernel $\Ker\bigl(\d
f\vert_v:T_vV\ra T_yY\bigr)$. As $\d f\vert_v:T_vV\ra T_yY$ is
surjective and $T_vV=(T_vV)^{\rm t}\op(T_vV)^{\rm nt}$, this implies
that $\bigl(\d f\vert_v\bigr)\vert_{(T_vV)^{\rm t}}:(T_vV)^{\rm
nt}\ra T_yY$ is surjective. Therefore $\d(f\ci\io^{\Ga,\rho,f})
\vert_{v^{\Ga,\rho,f}}:T_{v^{\Ga,\rho,f}}V^{\Ga,\rho,f}\ra T_yY$ is
{\it surjective}, for all $v^{\Ga,\rho,f}\in V^{\Ga,\rho,f}$. Hence
$f\ci\io^{\Ga,\rho,f}$ is a {\it submersion}.

Now let $X$ be a Kuranishi space and $\bs f:X\ra Y$ a strong
submersion. Combining \eq{kh5eq10} and \eq{kh5eq14}, define the {\it
orbifold stratum} $X^{\Ga,\rho,\bs f}$ as a set to be
\ea
X^{\Ga,\rho,\bs f}&\!=\!\bigl\{\Stab_X(p)\cdot(p,\la):p\in X,\;
\la:\Ga\!\ra\!\Ker\bigl(\bs f_*:\Stab_X(p)\!\ra\!\Stab_Y(f(p))\bigr)
\nonumber\\
&\quad\,\,\,\text{is an injective group morphism,
$\bigl[(T_vV_p)^{\rm nt}\ominus(E_p\vert_v)^{\rm
nt}\bigr]=\rho$}\bigr\}.\!\!\!
\label{kh5eq14}
\ea
Then $X^{\Ga,\rho,\bs f}$ is an open and closed subset of
$X^{\Ga,\rho}$, so it is a Kuranishi space by Proposition
\ref{kh5prop2}, and $\bs\io^{\Ga,\rho,\bs f}=\bs\io^{\Ga,\rho}
\vert_{X^{\Ga,\rho,\bs f}}:X^{\Ga,\rho,\bs f}\ra X$ is a proper,
finite, strongly smooth map. Also $\bs f\ci\bs\io^{\Ga,\rho,\bs f}:
X^{\Ga,\rho,\bs f}\ra Y$ is a {\it strong submersion}.

Let $\Ga$ be a finite group with $\md{\Ga}$ odd. Choose orientations
on $I_1,\ldots,I_N$ as in Definition \ref{kh5def12}. Let $Y$ be an
orbifold without boundary, and $R$ a commutative ring. Following
\eq{kh5eq11}, for each isomorphism class $\rho$ of virtual
nontrivial $\Ga$-representations and each $k\in\Z$, define an
operator $\Pi^{\Ga,\rho}:KB^k(Y;R)\ra KB^{k+\dim\rho}(Y;R)$ by
\e
\ts\Pi^{\Ga,\rho}:\sum_{a\in A}\rho_a\bigl[X_a,\bs
f_a\bigr]\longmapsto\sum_{a\in A}\rho_a\bigl[X_a^{\Ga,\rho,\bs
f_a},\bs f_a\ci\bs\io_a^{\Ga,\rho,\bs f_a}\bigr].
\label{kh5eq15}
\e
Here $X_a^{\Ga,\rho,\bs f_a}$ is a compact Kuranishi space as $X_a$
is compact and $\bs\io_a^{\Ga,\rho,\bs f_a}$ is proper, and $\bs
f_a\ci\bs\io_a^{\Ga,\rho,\bs f_a}$ is a strong submersion as above.
We define the coorientation for $(X_a^{\Ga,\rho,\bs f_a},\bs
f_a\ci\bs\io_a^{\Ga,\rho,\bs f_a})$ from the coorientation for
$(X_a,\bs f_a)$ by the analogue of Definition \ref{kh5def11}. We
also define $\Pi^{\Ga,\rho}_\ca:KB^{2l}_\ca(Y;R)\ra
KB^{2l+2\dim\rho}_\ca(Y;R)$ by including co-almost complex
structures in the obvious way.

Suppose $Y,Z$ are orbifolds, $X$ is a Kuranishi space, $\bs f:X\ra
Z$ is a strong submersion, and $h:Y\ra Z$ is a smooth map. Then
$Y\t_{h,Z,\bs f}X$ is a Kuranishi space, and $\bs\pi_Y:Y\t_ZX\ra Y$
a strong submersion. One can show that there is a natural
isomorphism of Kuranishi spaces $(Y\t_{h,Z,\bs f}X)^{\Ga,\rho,\bs
\pi_Y}\cong Y\t_{h,Z,\bs f\ci\io^{\Ga,\rho,\bs f}}X^{\Ga,\rho,\bs
f}$. Using this, we find that the $\Pi^{\Ga,\rho}$ commute with
pullbacks, $h^*\!\ci\!\Pi^{\Ga,\rho}\!=\!\Pi^{\Ga,\rho}\!\ci\!
h^*:KB^k(Z;R)\!\ra\! KB^{k+\dim\rho}(Y;R)$, and similarly~$h^*
\!\ci\!\Pi^{\Ga,\rho}_\ca\!=\!\Pi^{\Ga,\rho}_\ca\!\ci\!h^*$.
\label{kh5def14}
\end{dfn}

\subsection{How large are Kuranishi (co)bordism groups?}
\label{kh57}

We will now use the operators $\Pi^{\Ga,\rho},\Pi^{\Ga,\rho}_\ac$ of
\S\ref{kh56} composed with projection to singular homology
$\Pi_\Kb^\rsi,\Pi_\ac^\rsi$ to show that the bordism groups
$KB_*,KB_*^\ac(Y;R)$ are {\it very large}, infinitely generated at
least in every even degree, and that $KB_*^\ef,KB_*^\eac(Y;R)$ are
{\it large}, infinitely generated at least in every positive even
degree. Poincar\'e duality then implies the corresponding results
for the cobordism groups $KB^*,KB^*_\ca,KB^*_\ecb,KB^*_\eac(Y;R)$.

\begin{ex} Let $Y$ be a connected, nonempty orbifold, and fix $y\in
Y$, so that $H_0(Y;\Q)\cong\Q$ is generated by $[y]\in H_0(Y;\Q)$.
Let $\De$ be a finite group. Then $\{0\}/\De$ is a compact
0-orbifold without boundary, consisting of a single point with
stabilizer group $\De$. Give $\{0\}/\De$ the positive orientation.
Let $f:\{0\}/\De\ra Y$ be the smooth map of orbifolds mapping
$0\mapsto y$, with induced map $f_*:\Stab_{\{0\}/\De}(0)
\ra\Stab_Y(y)$ given by $f_*\equiv 1$. Regarding $\{0\}/\De$ as a
Kuranishi space and $f$ as strongly smooth, $[\{0\}/\De,f]$ is
well-defined in $KB_0(Y;\Z)$. Write $(\bs 0,\bs 0)$ for the trivial
almost complex structure on $\{0\}/\De$, giving almost complex
structures on zero vector spaces. Then $[\{0\}/\De,(\bs 0,\bs 0),f]$
is well-defined in~$KB_0^\ac(Y;\Z)$.

Let $\Ga$ be a finite group, and take $\rho=0$ to be the isomorphism
class of the zero (virtual) $\Ga$-representation. Then \eq{kh5eq7}
and \eq{kh5eq10} imply that $(\{0\}/\De)^{\Ga,0}\cong {\rm
InjHom}(\Ga,\De)/\De$ as an orbifold, where ${\rm InjHom}(\Ga,\De)$
is the set of injective group morphisms $\mu:\Ga\ra\De$, and $\De$
acts on ${\rm InjHom}(\Ga,\De)$ by $\de:\mu\ra\mu^\de$ for
$\de\in\De$, where $\mu^\de(\ga)=\de\mu(\ga)\de^{-1}$
for~$\ga\in\Ga$.

Suppose $\md{\Ga}$ is odd. Then Definition \ref{kh5def11} gives an
orientation on $(\{0\}/\De)^{\Ga,0}$, which is just the positive
orientation on ${\rm InjHom}(\Ga,\De)/\De$, and Definition
\ref{kh5def12} defines $\Pi^{\Ga,0}: KB_0(Y;\Z)\ra KB_0(Y;\Z)$,
where
\e
\Pi^{\Ga,0}\bigl([\{0\}/\De,f]\bigr)=\bigl[{\rm
InjHom}(\Ga,\De)/\De,f^{\Ga,0}\bigr],
\label{kh5eq16}
\e
where $f^{\Ga,0}$ maps ${\rm InjHom}(\Ga,\De)/\De\ra y$ with
$f^{\Ga,0}_*\equiv 1$ on stabilizers.

Now apply the morphism $\Pi_\Kb^\rsi:KB_0(Y;\Z)\ra H_0(Y;\Q)$ of
Definition \ref{kh5def5}(f). (Note that the definition of
$\Pi_\Kb^\rsi$ uses the fact that $\Pi_\Kh^\rsi:KH_0(Y;\Q)\ra
H_0(Y;\Q)$ is invertible, a highly nontrivial fact from Theorems
\ref{kh4thm1} and \ref{kh4thm2} which took all of Appendices
\ref{khA}--\ref{khC} to prove.) From \eq{kh5eq16} we see that
\e
\Pi_\Kb^\rsi\ci\Pi^{\Ga,0}\bigl([\{0\}/\De,f]\bigr)=\frac{\md{{\rm
InjHom}(\Ga,\De)}}{\md{\De}}\,[y]\in H_0(Y;\Q).
\label{kh5eq17}
\e

Similarly, for all finite groups $\Ga$ we find that
\e
\Pi_\ac^\rsi\ci\Pi^{\Ga,0}\bigl([\{0\}/\De,(\bs 0,\bs 0),f]\bigr)=
\frac{\md{{\rm InjHom}(\Ga,\De)}}{\md{\De}}\,[y]\in H_0(Y;\Q).
\label{kh5eq18}
\e
\label{kh5ex1}
\end{ex}

Using Example \ref{kh5ex1} we show that $KB_0(Y;\Z)$ and
$KB_0^\ac(Y;\Z)$ contain a copy of $\Z^\iy$, so they are very large.
The same proof shows that if $R$ is a commutative ring without
torsion, so that $\pi:R\ra R\ot_\Z\Q$ is injective, then $KB_0(Y;R)$
and $KB_0^\ac(Y;R)$ contain a copy of~$R^\iy$.

\begin{prop} In Example {\rm\ref{kh5ex1},} the elements
$[\{0\}/\De,f]$ taken over all isomorphism classes of finite groups
$\De$ with $\md{\De}$ odd are linearly independent over $\Z$ in
$KB_0(Y;\Z)$. Similarly, the elements $[\{0\}/\De,(\bs 0,\bs 0),f]$
taken over all isomorphism classes of finite groups $\De$ odd are
linearly independent over $\Z$ in both $KB_0^\ac(Y;\Z)$. Hence
$KB_0(Y;\Z)$ and\/ $KB_0^\ac(Y;\Z)$ are infinitely generated
over~$\Z$.
\label{kh5prop3}
\end{prop}

\begin{proof} Suppose for a contradiction that the
$[\{0\}/\De,f]$ are not linearly independent over $\Z$ in
$KB_0(Y;\Z)$. Then there exist pairwise non-isomorphic finite groups
$\De_1,\ldots,\De_N$ with $\md{\De_i}$ odd and nonzero integers
$a_1,\ldots,a_N$ such that $\sum_{i=1}^Na_i[\{0\}/\De_i,f_i]=0$ in
$KB_0(Y;\Z)$. Pick $i=1,\ldots,N$ with $\md{\De_i}$ largest, and set
$\Ga=\De_i$. If $j\ne i$ then either $\md{\De_j}<\md{\Ga}$, or
$\md{\De_j}=\md{\Ga}$ but $\De_j\not\cong\Ga$, and in both cases
${\rm InjHom}(\Ga,\De_j)=\es$. Thus $\Pi_\Kb^\rsi\ci\Pi^{\Ga,0}
([\{0\}/\De_j,f])=0$ for $j\ne i$ by \eq{kh5eq17}. So applying
$\Pi_\Kb^\rsi\ci\Pi^{\Ga,0}$ to $\sum_{i=1}^Na_i[\{0\}/\De_i,f_i]=0$
and using \eq{kh5eq17} gives $a_i\md{\Aut(\Ga)}/\md{\Ga}[y]=0$ in
$H_0(Y;\Q)\cong\Q$, a contradiction as $a_i,\md{\Aut(\Ga)},\md{\Ga}$
and $[y]$ are all nonzero. The proof for $KB_0^\ac(Y;\Z)$ is the
same.
\end{proof}

We can prove similar results for $\Pi^{\Ga,\rho}$ with $\rho\ne 0$,
and so show that $KB_{2k},\ab KB_{2k}^\ac(Y;\Z)$ are infinitely
generated for all $k\in\Z$. To do this, we need to replace
$\{0\}/\De$ by a more complicated Kuranishi space~$X$.

\begin{ex} Let $Y,y$ be as above. Suppose $\De$ is a finite group,
and let $(W,\om),(W',\om')$ be finite-dimensional nontrivial complex
representations of $\De$. Consider the complex projective space
$P(\C\op W)=W\amalg P(W)$. The representation $1\op\om$ of $\De$ on
$\C\op W$ induces an action of $\De$ on $P(\C\op W)$, so that
$P(\C\op W)/\De$ is a compact complex orbifold. The point $[1,0]\in
P(\C\op W)$ is fixed by $\De$, so its image $[1,0]\in P(\C\op
W)/\De$ has stabilizer group~$\De$.

Define $X$ to be $P(\C\op W)/\De$ as a topological space, and define
a Kuranishi structure on $X$ by the single Kuranishi neighbourhood
\begin{equation*}
(V,E,s,\psi)=\bigl(P(\C\op W)/\De, (P(\C\op W)\t W')/\De,0,\id_X).
\end{equation*}
Here $E=(P(\C\op W)\t W')/\De\ra P(\C\op W)/\De$ is an orbifold
vector bundle with fibre $W'$, and $\De$ acts on $P(\C\op W)\t W'$
by the given action on $P(\C\op W)$, and the representation $\om'$
on $W'$. Then $X$ is a compact Kuranishi space without boundary, of
virtual dimension $2(\dim W -\dim W')$. The complex structures on
$P(\C\op W)$ and $W'$ induce an {\it almost complex structure} $(\bs
J,\bs K)$ on $X$. This in turn induces an {\it orientation} on $X$.
Define a strongly smooth map $\bs f:X\ra Y$ to be represented by the
projection $f:V\ra Y$ mapping $f:v\mapsto y$ for all $v\in V$, with
$f_*\equiv 1$ on stabilizer groups. Then $[X,\bs f]$ is well-defined
in $KB_{2(\dim W-\dim W')}(Y;R)$, and $[X,(\bs J,\bs K),\bs f]$ is
well-defined in~$KB_{2(\dim W-\dim W')}^\ac(Y;R)$.

For each $p\in X$, the stabilizer group $\Stab_X(p)$ is a subgroup
of $\De$. Thus if $\md{\Ga}>\md{\De}$, or if $\md{\Ga}=\md{\De}$ but
$\Ga\not\cong\De$, then there can exist no injective morphisms
$\la:\Ga\ra\Stab_X(p)$, so \eq{kh5eq10} gives $X^{\Ga,\rho}=\es$.
Hence
\e
\Pi^{\Ga,\rho}\bigl([X,\bs f]\bigr)\!=\!\Pi^{\Ga,\rho}_\ac
\bigl([X,(\bs J,\bs K),\bs f]\bigr)\!=\!0\;\>\text{if
$\md{\Ga}\!>\!\md{\De}$ or $\md{\Ga}\!=\!\md{\De}$,
$\Ga\!\not\cong\!\De$.}\!\!
\label{kh5eq19}
\e

Now set $\Ga=\De$. The points in $X$ with stabilizer group $\De$ are
the fixed point set of $\De$ in $P(\C\op W)$. Writing $P(\C\op
W)=W\amalg P(W)$, the only fixed point of $\De$ in $W$ is 0, since
$W$ is a {\it nontrivial\/} $\De$-representation. The fixed points
of $\De$ in $P(W)$ correspond to $\De$-subrepresentations of $W$
isomorphic to $\C$. We want to ensure that {\it none of these fixed
points in $P(W)$ is isolated}. We therefore impose the condition:
\begin{itemize}
\setlength{\itemsep}{0pt}
\setlength{\parsep}{0pt}
\item[$(*)$] No one-dimensional irreducible complex representation
of $\De$ occurs with multiplicity exactly one in $(W,\om)$.
\end{itemize}
This implies that the fixed points of $\De$ in $P(W)$ are a
(possibly empty) finite disjoint union of projective subspaces of
$P(W)$ of {\it positive dimension}.

Thus, if $X^{\De,\rho}\ne\es$, there are two possibilities:
\begin{itemize}
\setlength{\itemsep}{0pt}
\setlength{\parsep}{0pt}
\item[(a)] $\rho=\bigl[(W,\om)\ominus(W',\om')\bigr]$, in which case
$X^{\De,\rho}$ comes from $0\in W$, and
$X^{\De,\rho}\cong\Aut(\De)/\De$ as an orbifold, and $\vdim
X^{\De,\rho}=0$; or
\item[(b)] $X^{\De,\rho}$ is one or more positive-dimensional
projective subspaces of $P(W)$, an orbifold, and $\vdim
X^{\De,\rho}>0$.
\end{itemize}
Thus as for \eq{kh5eq18} we see that
\e
\begin{split}
\Pi_\Kb^\rsi&\ci\Pi^{\De,\rho}\bigl([X,\bs f]\bigr)
=\Pi_\ac^\rsi\ci\Pi^{\De,\rho}\bigl([X,(\bs J,\bs K),\bs f]\bigr)\\
&=\begin{cases} \frac{\md{\Aut(\De)}}{\md{\De}}\,[y], &
\rho=\bigl[(W,\om)\ominus(W',\om')\bigr], \\
0, & \rho\ne\bigl[(W,\om)\ominus(W',\om')\bigr],\; \dim_\R\rho=\vdim
X. \end{cases}
\end{split}
\label{kh5eq20}
\e
\label{kh5ex2}
\end{ex}

In the almost complex case, for all isomorphism classes $\rho$ of
complex virtual nontrivial $\De$-representations, we can choose
complex nontrivial representations $(W,\om),(W',\om')$ satisfying
$(*)$ with $\rho=\bigl[(W,\om)\ominus(W',\om')\bigr]$. In the real
case, if $\md{\De}$ is odd then every real nontrivial representation
of $\De$ has an underlying complex representation, so again, for all
isomorphism classes $\rho$ of real virtual nontrivial
$\De$-representations, we can choose complex nontrivial
representations $(W,\om),(W',\om')$ satisfying $(*)$ with
$\rho=\bigl[(W,\om)\ominus(W',\om')\bigr]$ in real virtual
representations. Thus, using \eq{kh5eq19}--\eq{kh5eq20} and a
similar argument to the proof of Proposition \ref{kh5prop3} gives:

\begin{prop} In Example {\rm\ref{kh5ex2},} the elements $[X,\bs f]$ taken
over all isomorphism classes of finite groups $\De$ with\/
$\md{\De}$ odd and all isomorphism classes $\rho$ of real virtual
nontrivial\/ $\De$-representations are linearly independent over
$\Z$ in $KB_*(Y;\Z)$. There are infinitely many such generators
$[X,\bs f]$ in $KB_{2l}(Y;\Z)$ for all $l\in\Z$. Hence
$KB_{2l}(Y;\Z)$ is infinitely generated over $\Z$ for
all\/~$l\in\Z$.

Similarly, the elements $[X,(\bs J,\bs K),\bs f]$ taken over all
isomorphism classes of finite groups $\De$ odd and all isomorphism
classes $\rho$ of complex virtual nontrivial\/ $\De$-representations
are linearly independent over $\Z$ in $KB_{2*}^\ac(Y;\Z)$. Hence
$KB_{2l}^\ac(Y;\Z)$ is infinitely generated over $\Z$ for
all\/~$l\in\Z$.
\label{kh5prop4}
\end{prop}

We can also use these ideas to study {\it effective} Kuranishi
bordism. It is easy to see that the Kuranishi space $X$ of Example
\ref{kh5ex2} is effective if and only if $W'=0$ and $W$ is an {\it
effective} representation, and then $X$ is actually an {\it
effective orbifold}, so that $[X,\bs f]$ is well-defined in
$B_*^\eo(Y;\Z)$. Thus as for Proposition \ref{kh5prop4} we obtain:

\begin{prop} In Example {\rm\ref{kh5ex2},} taking $W'=0$ and\/
$(W,\om)$ to be an effective nontrivial\/ $\De$-representation, the
elements $[X,\bs f]$ taken over all isomorphism classes of finite
groups $\De$ with\/ $\md{\De}$ odd and all isomorphism classes of\/
$(W,\om)$ satisfying $(*)$ are linearly independent over $\Z$ in
both\/ $KB_*^\ef(Y;\Z)$ and\/ $B_*^\eo(Y;\Z)$. There are infinitely
many such generators $[X,\bs f]$ in $KB_{2l}(Y;\Z)$ and\/
$B_{2l}^\eo(Y;\Z)$ for all\/ $l=2,3,4,\ldots$. Hence
$KB_{2l}(Y;\Z),$ $B_{2l}^\eo(Y;\Z),$ and similarly
$KB_{2l}^\eac(Y;\Z),$ are infinitely generated over $\Z$ for
all\/~$l\ge 2$.
\label{kh5prop5}
\end{prop}

We exclude the case $l=1$ here as condition $(*)$ in Example
\ref{kh5ex2} would fail. But using different examples of Kuranishi
spaces $X$ one can also show that $KB_2,B_2^\eo,KB_2^\eac(Y;\Z)$ are
infinitely generated over $\Z$. Using the Poincar\'e duality ideas
of \S\ref{kh54} we can deduce:

\begin{cor} Let\/ $Y$ be an oriented\/ $n$-orbifold without
boundary. Then $KB^{n-2l}(Y;\Z)$ is infinitely generated over $\Z$
for all\/ $l\in\Z,$ and $KB^{n-2l}_\ec(Y;\Z)$ is infinitely
generated over $\Z$ for all\/~$l\ge 2$.
\label{kh5cor2}
\end{cor}

Similar results hold for $KB^*_\ac(Y;\Z)$ and $KB^*_\eac(Y;\Z)$. The
fact that our Kuranishi (co)bordism groups are very large has both
advantages and disadvantages. One advantage is that any invariant
defined in a Kuranishi (co)bordism group {\it contains a lot of
information}, as it lies in a very large group. We will see this
with Gromov--Witten cobordism invariants in \S\ref{kh6}. One
disadvantage is that it is not really feasible to compute the groups
and write them down.

\section{Applications in Symplectic Geometry}
\label{kh6}

We now discuss some applications of Kuranishi (co)homology and
Kuranishi (co)bordism in areas of Symplectic Geometry concerned with
moduli spaces of $J$-holomorphic curves. Our aim is not to give a
complete treatment --- we postpone this to the sequels
\cite{Joyc2,Joyc4} --- but rather to demonstrate the potential of
our theories as tools both for proving new theorems, and for
reframing existing theories in a much cleaner and more economical
way.

Sections \ref{kh61}--\ref{kh63} concern closed Gromov--Witten
invariants. They explain how to define new Gromov--Witten type
invariants in (almost complex) Kuranishi bordism groups, and
consider how these might be used to study integrality properties of
Gromov--Witten invariants, and prove the Gopakumar--Vafa Integrality
Conjecture. Sections \ref{kh64}--\ref{kh67} discuss Lagrangian Floer
cohomology and open Gromov--Witten invariants, and \S\ref{kh68}
draws some conclusions.

\subsection{Moduli spaces of closed Riemann surfaces}
\label{kh61}

We now discuss (closed) Riemann surfaces and stable maps,
following~\cite{FuOn1}.

\begin{dfn} A {\it Riemann surface\/} $\Si$ is a compact connected
complex 1-manifold. A {\it prestable\/} or {\it nodal Riemann
surface\/} $\Si$ is a compact connected complex variety whose only
singularities are {\it nodes\/} ({\it ordinary double points\/}),
that is, $\Si$ has finitely many singularities, each modelled on the
{\it node\/} $(0,0)$ in $\bigl\{(x,y)\in\C^2:xy=0\bigr\}$. The {\it
genus\/} of a prestable Riemann surface $\Si$ is the genus of its
smoothing $\Si'$, that is, we smooth nodes modelled on
$\bigl\{(x,y)\in\C^2:xy=0\bigr\}$ to a nonsingular curve modelled on
$\bigl\{(x,y)\in\C^2:xy=\ep\bigr\}$. A ({\it prestable\/}) {\it
Riemann surface with $m$ marked points\/} $(\Si,\vec z)$ is a
(prestable) Riemann surface $\Si$ with $\vec z=(z_1,\ldots,z_m)$,
where $z_1,\ldots,z_m$ are distinct points in $\Si$, none of which
is a node. A {\it stable\/} Riemann surface with marked points is a
prestable $(\Si,\vec z)$ whose automorphism group is finite.

Write $\oM_{g,m}$ for the moduli space of isomorphism classes
$[\Si,\vec z\,]$ of stable Riemann surfaces $(\Si,\vec z)$ of genus
$g$ with $m$ marked points, and $\M_{g,m}$ for the subset of
$[\Si,\vec z\,]$ in $\oM_{g,m}$ with $\Si$ nonsingular. Then
$\oM_{g,m}$ is known as the {\it Deligne--Mumford
compactification\/} of $\M_{g,m}$. It is easy to show that when
$(g,m)=(0,0),(0,1),(0,2)$ or (1,0) we have
$\M_{g,m}=\oM_{g,m}=\emptyset$. So we generally restrict to~$2g+m\ge
3$.

Let $(\Si,\vec z)$ be a prestable Riemann surface with genus $g$ and
$m$ marked points, and suppose $2g+m\ge 3$. If $(\Si,\vec z)$ is not
stable then $\Si$ has a $\CP^1$ component with less than 3 nodes or
marked points. Collapse this $\CP^1$ component to a point. After
repeating this process finitely many times we obtain a unique stable
Riemann surface $(\ti\Si,\vec{\ti z}\,)$ with genus $g$ and $m$
marked points, called the {\it stabilization\/} of~$(\Si,\vec z)$.
\label{kh6def1}
\end{dfn}

A great deal is known about the moduli spaces $\M_{g,m}$ and
$\oM_{g,m}$, see for instance Harris and Morrison \cite{HaMo}. The
following theorem is proved in~\cite[\S 9]{FuOn1}.

\begin{thm} In Definition {\rm\ref{kh6def1},} $\oM_{g,m}$ is a compact
complex orbifold of real dimension $2(3g+m-3),$ without boundary or
corners.
\label{kh6thm1}
\end{thm}

Observe that $\oM_{1,1}$ is a {\it non-effective orbifold}, whose
generic points have stabilizer group~$\Z_2$.

\begin{dfn} Let $(M,\om)$ be a compact symplectic manifold and $J$
an almost complex structure on $M$ compatible with $\om$. A {\it
stable map to $M$ from a Riemann surface with\/ $m$ marked points\/}
is a triple $(\Si,\vec z,w)$, where $(\Si,\vec z)$ is a prestable
Riemann surface with $m$ marked points, and $w:\Si\ra M$ a
$J$-pseudoholomorphic map, such that the group $\Aut(\Si,\vec z,w)$
of biholomorphisms $\ga:\Si\ra\Si$ with $\ga(z_j)=z_j$ for all $j$
and $w\ci\ga=w$ is finite.

For $\be\in H^2(M,\Z)$ and $g,m\ge 0$, write $\oM_{g,m}(M,J,\be)$
for the moduli space of $J$-pseudoholomorphic stable maps $(\Si,\vec
z,w)$ to $M$ from a Riemann surface $\Si$ with genus $g$ and $m$
marked points, with $[w(\Si)]=\be\in H^2(M,\Z)$. Points of
$\oM_{g,m}(M,J,\be)$ are isomorphism classes $[\Si,\vec z,w]$, where
$(\Si,\vec z,w)$, $(\Si',\vec z',w')$ are isomorphic if there exists
a biholomorphism $\ga:\Si\ra\Si'$ with $\ga(z_j)=z_j'$ for all $j$
and $w'\ci\ga\equiv w$. Define {\it evaluation maps\/}
$\ev_j:\oM_{g,l}(M,J,\be)\ra M$ for $j=1,\ldots,m$ by
$\ev_j:[\Si,\vec z,w]\mapsto w(z_j)$. When $2g+m\ge 3$, define
$\pi_{g,m}:\oM_{g,m}(M,J,\be)\ra\oM_{g,m}$ by $\ev_j:[\Si,\vec
z,w]\mapsto[\ti\Si,\vec{\ti z}\,]$, where $(\ti\Si,\vec{\ti z}\,)$
is the stabilization of~$(\Si,\vec z)$.
\label{kh6def2}
\end{dfn}

There is a natural topology on $\oM_{g,m}(M,J,\be)$ called the
$C^\iy$ {\it topology}, due to Gromov \cite{Grom} and defined in
\cite[\S 10]{FuOn1}. The next result, due in principle to Gromov, is
proved by Fukaya and Ono \cite[\S 10-\S 11]{FuOn1}.

\begin{thm} In Definition {\rm\ref{kh6def2},} $\oM_{g,m}(M,J,\be)$ is
compact and Hausdorff in the\/ $C^\iy$ topology.
\label{kh6thm2}
\end{thm}

\begin{rem} In the following, I will label Theorems \ref{kh6thm3},
\ref{kh6thm4}, \ref{kh6thm5}, \ref{kh6thm9}, \ref{kh6thm10} and
\ref{kh6thm11} as ``Theorem'' rather than Theorem. This is intended
to indicate that {\it I am not confident that a complete, detailed
proof has yet been given for these results}; that they should be
regarded as {\it partially proved}, as somewhere between a Theorem
and a Conjecture.

I state these as ``Theorems'' because actually proving them properly
would be a massive undertaking, beyond the scope of this book and
beyond my abilities and patience, but I cannot proceed to
applications in Symplectic Geometry without them, and I need some
applications to persuade the reader that the effort of defining my
Kuranishi (co)homology machinery was worthwhile.

The results of \S\ref{kh62}--\S\ref{kh63} will then assume the
``Theorems'' of \S\ref{kh61}, and the results of
\S\ref{kh65}--\S\ref{kh67} will assume the ``Theorems'' of
\S\ref{kh64}. For comparison, note that whole areas of mathematics
have the Riemann hypothesis as a standing assumption, although this
has not yet been proved. In this case, we have good reasons to hope
that complete proofs of our ``Theorems'' will be available soon.

The worries I have about the proofs of the ``Theorems'' below mostly
concern the definition of the Kuranishi structure on the moduli
spaces at and near singular $J$-holomorphic curves, that is, curves
with interior or boundary nodes. Especially, I am bothered about
{\it smoothness issues}: whether one can define Kuranishi
neighbourhoods $(V_p,E_p,s_p,\psi_p)$ on the moduli spaces near such
a singular curve for which $s_p$ is a smooth section of $E_p$, and
what arbitrary choices are involved; whether coordinate changes
$(\phi_{pq},\hat\phi_{pq})$ between such neighbourhoods can be made
smooth, and whether they satisfy Definition \ref{kh2def12}(e).

These questions are very difficult, concerning as they do
singularities of solutions of nonlinear elliptic p.d.e.s, with
bubbling, neck stretching, and so on --- the `Analytical Chamber of
Horrors', as Hofer describes it \cite[p.~2]{Hofe1}. Fukaya et al.\
discuss smoothness in \cite[\S A1.4]{FOOO}. Their underlying idea is
related to the notion of {\it gluing profile\/} in Hofer \cite[\S
4.1]{Hofe1}, which is how these smoothness issues are resolved in
the theory of polyfolds.

Part of the problem is the mismatch between the definitions of
Kuranishi spaces used in this book, and in our references. For
example, Fukaya and Ono give a proof of ``Theorem'' \ref{kh6thm3}
below using a different and rather weaker notion of Kuranishi space:
for them, the Kuranishi maps $s_p$ only have to be continuous, so
they avoid most of the smoothness issues, and Definition
\ref{kh2def12}(e) is replaced by a weaker notion involving
isomorphisms $\chi_{pq}$ as in Remark \ref{kh2rem5}. But I want
``Theorem'' \ref{kh6thm3} to be true with our stronger notion of
Kuranishi space.

Part of the problem too is that I am too lazy to work through all
the proofs that have been written down, and satisfy myself that they
are correct, and extend to our stronger definition of Kuranishi
spaces. But also, there are parts of some proofs in
\cite{FuOn1,FOOO} which I have read, which appeared to me to be
either lacking in detail at some important points, or to have
problems which needed fixing. Fortunately, as we will discuss in
\S\ref{kh68}, there may soon be an alternative proof of ``Theorem''
\ref{kh6thm3} and related results, via the theory of polyfolds of
Hofer, Wysocki and Zehnder~\cite{Hofe1,Hofe2,HWZ1,HWZ2,HWZ3}.
\label{kh6rem1}
\end{rem}

The following ``Theorem'' is proved in \cite[Th.s 7.10--7.11 \&
Prop.~16.1]{FuOn1}, but with a weaker definition of Kuranishi
structure. The claims about strong submersions are evident from the
construction.

\begin{qthm} In Definition {\rm\ref{kh6def2},} there exists an oriented
Kuranishi structure $\ka$ on\/ $\oM_{g,m}(M,J,\be),$ without
boundary or (g-)corners, depending on choices, with
\e
\vdim\oM_{g,m}(M,J,\be)=2\bigl(c_1(M)\cdot\be+(n-3)(1-g)+m\bigr),
\label{kh6eq1}
\e
where $\dim M=2n$. The maps\/ $\ev_j$ in Definition
{\rm\ref{kh6def2}} are continuous and extend for\/ $j=1,\ldots,m$ to
strong submersions $\bev_j:\bigl(\oM_{g,m}(M,J,\be),\ka\bigr)\ra M,$
and\/ $\bev_1\t\cdots\t\bev_m$ is a strong submersion. When $2g+m\ge
3,$ the map\/ $\pi_{g,m}$ in Definition {\rm\ref{kh6def2}} is
continuous and extends to a strong submersion
$\bs\pi_{g,m}:\bigl(\oM_{g,m}(M,J,\be),\ka\bigr)\!\ra\!\oM_{g,m},$
and\/ $\bev_1\!\t\!\cdots\!\t\!\bev_m\!\t\bs\pi_{g,m}$ is a strong
submersion.
\label{kh6thm3}
\end{qthm}

Our next ``Theorem'' describes the dependence of
$\oM_{g,m}(M,J,\be)$ on the almost complex structure $J$. It is
discussed briefly in the proof of \cite[Th.~17.11]{FuOn1}, and
proved by the same method as ``Theorem''~\ref{kh6thm3}.

\begin{qthm} Let\/ $(M,\om)$ be a compact symplectic manifold and\/
$J_t$ for\/ $t\in[0,1]$ a smooth\/ $1$-parameter family of almost
complex structures on\/ $M$ compatible with\/ $\om$. Write
\begin{align*}
&\oM_{g,m}(M,J_t:t\in[0,1],\be)=\\
&\bigl\{\bigl(t,[\Si,\vec z,w]\bigr): t\in[0,1],\;\>[\Si,\vec
z,w]\in\oM_{g,m}(M,J_t,\be)\bigr\}.
\end{align*}
There is a natural\/ $C^\iy$ topology on\/ $\oM_{g,m}(M,J_t:t\in
[0,1],\be)$ making it into a compact Hausdorff topological space.
Let\/ $\ka_0,\ka_1$ be possible oriented Kuranishi structures on\/
$\oM_{g,m}(M,J_t,\be)$ for\/ $t=0,1$ given by ``Theorem''
{\rm\ref{kh6thm3}}. Then there exists an oriented Kuranishi
structure\/ $\ti\ka$ on\/ $\oM_{g,m}(M,J_t:t\in[0,1],\be),$ with
\begin{equation*}
\vdim\bigl(\oM_{g,m}(M,J_t:t\in[0,1],\be),\ti\ka\bigr)
=2\bigl(c_1(M)\cdot\be+(n-3)(1-g)+m\bigr)+1,
\end{equation*}
without (g-)corners, and with boundary
\e
\begin{split}
\pd\bigl(\oM_{g,m}&(M,J_t:t\in[0,1],\be),\ti\ka\bigr)\\
&=\bigl(\oM_{g,m}(M,J_1,\be),\ka_1\bigr)\amalg
-\bigl(\oM_{g,m}(M,J_0,\be),\ka_0\bigr).
\end{split}
\label{kh6eq2}
\e

There are strong submersions\/ $\widetilde\bev_j:\bigl(\oM_{g,m}
(M,J_t:t\in[0,1],\be),\ti\ka\bigr)\ra M$ for\/ $j=1,\ldots,m$ and\/
$\widetilde{\bs\pi}_{g,m}:\bigl(\oM_{g,m}(M,J_t:t\in[0,1],\be),
\ti\ka\bigr)\ra\oM_{g,m}$ when\/ $2g+m\ge 3,$ which restrict to\/
$\bev_j,\bs\pi_{g,m}$ in ``Theorem'' {\rm\ref{kh6thm3}} on\/
$\oM_{g,m}(M,J_t,\be)$ for\/ $t=0,1$ under\/ {\rm\eq{kh6eq2},} and\/
$\widetilde\bev_1\t\cdots\t\widetilde\bev_m,$ $\widetilde
\bev_1\t\cdots\t\widetilde\bev_m\t\widetilde{\bs\pi}_{g,m}$ are
strong submersions. Note that we can take\/ $J_t\equiv J$ for
$t\in[0,1],$ so this result relates possible choices of\/ $\ka$ in
``Theorem'' {\rm\ref{kh6thm3}} for a single\/ $J$ on\/~$M$.
\label{kh6thm4}
\end{qthm}

We can put an {\it almost complex structure} on
$\oM_{g,m}(M,J,\be)$, in the sense of \S\ref{kh29}. The (incomplete)
proof is based on Fukaya and Ono \cite[Prop.~16.5]{FuOn1}, who prove
that $\oM_{g,m}(M,J,\be)$ is {\it stably almost complex}, that is,
that the virtual tangent bundle of $\oM_{g,m}(M,J,\be)$ is
equivalent to a virtual complex vector bundle in a version of
K-theory for Kuranishi spaces.

\begin{qthm} In the situation of ``Theorem'' {\rm\ref{kh6thm3},} one can
construct an almost complex structure\/ $(\bs J,\bs K)$ on
$\bigl(\oM_{g,m}(M,J,\be),\ka\bigr)$. This construction depends on
some choices, and $\ka$ must also be chosen appropriately.

These $\ka,(\bs J,\bs K)$ are unique up to isotopy in the following
sense. Let\/ $J_t$ for $t\in[0,1]$ be a smooth family of almost
complex structures on $M$ compatible with\/ $\om$. Let\/ $\ka_0,(\bs
J_0,\bs K_0)$ and\/ $\ka_1,(\bs J_1,\bs K_1)$ be outcomes of the
construction above for $\oM_{g,m}(M,J_0,\be),\oM_{g,m}(M,J_1,\be)$.
Then in ``Theorem'' {\rm\ref{kh6thm4}} we can choose\/ $\ti\ka$ and
an almost CR structure\/ $(\bs D,\bs J,\bs K)$ on
$\bigl(\oM_{g,m}(M,J_t:t\in[0,1],\be),\ti\ka\bigr)$ with
\e
\begin{split}
&\pd\bigl((\oM_{g,m}(M,J_t:t\in[0,1],\be),\ti\ka),(\bs D,\bs
J,\bs K)\bigr)=\\
&\bigl((\oM_{g,m}(M,J_1,\be),\ka_1),(\bs J_1,\bs K_1)\bigr)\amalg
-\bigl((\oM_{g,m}(M,J_0,\be),\ka_0),(\bs J_0,\bs K_0)\bigr).
\end{split}
\label{kh6eq3}
\e

Analogues of the above results also hold for co-almost complex
structures and co-almost CR structures, using the ideas
of\/~{\rm\S\ref{kh210}}.
\label{kh6thm5}
\end{qthm}

\begin{proof}[Partial proof] In \cite[Prop.~16.5]{FuOn1}, Fukaya
and Ono observe that for a Kuranishi neighbourhood
$(V_p,E_p,s_p,\psi_p)$ on $\oM_{g,m}(M,J,\be)$, the fibres of the
orbibundles $TV_p$ and $E_p$ over $V_p$ are roughly the kernels and
cokernels of a family of linear elliptic operators $D_v:B_1\ra B_2$
between {\it complex\/} Banach spaces $B_1,B_2$, parametrized by
$v\in V_p$. Now the $D_v$ are not complex linear, but their symbols
$D_v'$ are. So writing $D_{(v,t)}=(1-t)D_v+tD_v'$, we obtain a
smooth family of elliptic operators $D_{(v,t)}:B_1\ra B_2$
parametrized by $(v,t)\in V_p\t[0,1]$, such that $D_{(v,0)}=D_v$ and
$D_{(v,1)}$ is complex linear. Using this and ideas in the proofs of
\cite[Prop.~16.5]{FuOn1} and in \S\ref{kh32}, we can construct the
following data:
\begin{itemize}
\setlength{\itemsep}{0pt}
\setlength{\parsep}{0pt}
\item[(i)] A {\it good coordinate system\/} $\bs I=\bigl(I,\le,(V^i,E^i,
s^i,\psi^i):i\in I,\ldots\bigr)$ for $\oM_{g,m}\ab(M,\ab J,\be)$, as
in Definition~\ref{kh3def1};
\item[(ii)] Real orbibundles $F_1^i,F_2^i$ over $V^i\t[0,1]$ for all~
$i\in I$;
\item[(iii)] Isomorphisms $\io_1^i:F_1^i\vert_{V^i\t\{0\}}\ra TV^i$,
$\io_2^i:F_2^i\vert_{V^i\t\{0\}}\ra E^i$ for all~$i\in I$;
\item[(iv)] Almost complex structures $J_1^i,J_2^i$ on the fibres of
$F_1^i\vert_{V^i\t\{1\}}$ and $F_2^i\vert_{V^i\t\{1\}}$ for all
$i\in I$; and
\item[(v)] For all $i,j\in I$ with $j\le i$ and
$\Im\psi^i\cap\Im\psi^j\ne\emptyset$, embeddings of orbibundles
$\Phi^{ij}_a:F_a^j\vert_{V^{ij}}\ra(\phi^{ij})^*(F_a^i)$ over
$V^{ij}$ for $a=1,2$, and an isomorphism:
\e
\Psi^{ij}:\frac{(\phi^{ij})^*(F_1^i)}{\Phi^{ij}_1(F_1^j\vert_{V^{ij}})}
\longra\frac{(\phi^{ij})^*(F_2^i)}{\Phi^{ij}_2(F_2^j\vert_{V^{ij}})}\,.
\label{kh6eq4}
\e
\end{itemize}

This data satisfies the following conditions:
\begin{itemize}
\setlength{\itemsep}{0pt}
\setlength{\parsep}{0pt}
\item[(vi)] For all $i,j\in I$ with $j\le i$ and
$\Im\psi^i\cap\Im\psi^j\ne\emptyset$ we have
$(\phi^{ij})^*(\io_1^i)\ci\Phi_1^{ij}=\d\phi^{ij}\ci\io_1^j$ and
$(\phi^{ij})^*(\io_2^i)\ci\Phi_2^{ij}=\hat\phi^{ij}\ci\io_2^j$ as
morphisms $F_1^j\vert_{V^{ij}\t\{0\}}\ra(\phi^{ij})^*(TV^i)$ and
$F_2^j\vert_{V^{ij}\t\{0\}}\ra(\phi^{ij})^*(E^i)$ over
$V^{ij}\t\{0\}\cong V^{ij}$, and the following commutes near
$(s^j)^{-1}(0)\cap V^{ij}$ in morphisms over $V^{ij}\t\{0\}\cong
V^{ij}$, where $\d\hat s^i$ is as in~\eq{kh2eq12}:
\e
\begin{gathered}
\xymatrix@C=70pt@R=20pt{
\displaystyle\frac{(\phi^{ij})^*(F_1^i)}{\Phi^{ij}_1(F_1^j
\vert_{V^{ij}})}\Big\vert_{V^{ij}\t\{0\}} \ar[r]^{\Psi^{ij}}
\ar[d]^{(\io^i_1)_*} &
\displaystyle\frac{(\phi^{ij})^*(F_2^i)}{\Phi^{ij}_2(F_2^j
\vert_{V^{ij}})}\Big\vert_{V^{ij}\t\{0\}} \ar[d]^{(\io^i_2)_*}
\\
\displaystyle\frac{(\phi^{ij})^*(TV^i)}{(\d\phi^{ij})(TV^j)}
\ar[r]^{\d\hat s^i} &
\displaystyle\frac{(\phi^{ij})^*(E^i)}{\hat\phi^{ij}(E^j)}\,, }
\end{gathered}
\label{kh6eq5}
\e
where the columns are well-defined because of the previous
conditions.
\item[(vii)] For all $i,j\in I$ with $j\le i$ and $\Im\psi^i\cap\Im
\psi^j\ne\emptyset$ we have $(\phi^{ij})^*(J_a^i)\ci\Phi^{ij}_a=
\Phi^{ij}_a\ci J_a^j$ as morphisms $F_a^j\vert_{V^{ij}}\ra
(\phi^{ij})^*(F_a^i)$ for $a=1,2$. This implies that
$(\phi_a^{ij})^*(J^i_a)$ for $a=1,2$ push down to almost complex
structures $J^{ij}_a$ for $a=1,2$ on the orbibundles over $V^{ij}$
appearing in \eq{kh6eq4}, and we require
that~$J^{ij}_2\ci\Psi^{ij}=\Psi^{ij}\ci J^{ij}_1$.
\item[(viii)] Whenever $i,j,k\in I$ with $k\le j\le i$ and $\Im\psi^i
\cap\Im\psi^j\cap\Im\psi^k\ne\emptyset$, we have
$\Phi^{ik}_a=(\phi^{jk})^*(\Phi^{ij}_a)\ci\Phi^{jk}$ over
$(\phi^{jk})^{-1}(V^{ij})\cap V^{jk}\cap V^{ik}$ for $a=1,2$, and
the following diagram of orbibundles over $(\phi^{jk})^{-1}
(V^{ij})\cap V^{jk}\cap V^{ik}$ commutes:
\end{itemize}
\begin{equation*}
\xymatrix@C=10pt@R=14pt{ 0 \ar[r] &
\displaystyle\frac{(\phi^{jk})^*(F^j_1)}{\Phi^{jk}_1(F^k_1)}
\ar[rrr]^{(\phi^{jk})^*(\Phi^{ij}_1)} \ar[d]^{\Psi^{jk}} &&&
\displaystyle\frac{(\phi^{ik})^*(F^i_1)}{\Phi^{ik}_1(F^k_1)}
\ar[rrr]^(0.45){\text{\tiny projection}} \ar[d]^{\Psi^{ik}} &&&
\displaystyle
\frac{(\phi^{ik})^*(F^i_1)}{(\phi^{jk})^*(\Phi^{ij}_1(F^j_1))}
\ar[r] \ar[d]^{(\phi^{jk})^*(\Psi^{ij})} & 0 \\
0 \ar[r] & \displaystyle
\frac{(\phi^{jk})^*(F^j_2)}{\Phi^{jk}_2(F^k_2)}
\ar[rrr]^{(\phi^{jk})^*(\Phi^{ij}_2)} &&&
\displaystyle\frac{(\phi^{ik})^*(F^i_2)}{\Phi^{ik}_2(F^k_2)}
\ar[rrr]^(0.45){\text{\tiny projection}} &&& \displaystyle
\frac{(\phi^{ik})^*(F^i_2)}{(\phi^{jk})^*(\Phi^{ij}_2(F^j_2))}
\ar[r] & 0. }
\end{equation*}

Roughly speaking, the data $F^i_a$ and overlap maps $\Phi^{ij}_a$
and $\Psi^{ij}$ define a bundle system $(\bs F_1,\bs F_2)$ on
$\oM_{g,m}(M,J,\be)\t[0,1]$, and (iii), (vi) say that $(\bs F_1,\bs
F_2)$ restricts to the tangent bundle system on $\oM_{g,m}(M,J,\be)
\t\{0\}$, and (iv), (vii) say that $(\bs F_1,\bs F_2)$ restricts to
a complex bundle system on $\oM_{g,m}(M,J,\be)\t\{1\}$, as in the
discussion of \cite[Prop.~16.5]{FuOn1} above. But working with
neighbourhoods of the form $(V^i\t[0,1],E^i,s^i,\psi^i)$ on
$\oM_{g,m}(M,J,\be)\t[0,1]$, rather than covering
$\oM_{g,m}(M,J,\be)\t[0,1]$ with small Kuranishi neighbourhoods,
gives an advantage.

Now any orbibundle over $V^i\t[0,1]$ is (non-canonically) isomorphic
to the pullback of an orbibundle over $V^i$, by the orbifold version
of well-known facts for manifolds. Thus by (iii) there exist
isomorphisms $\ti\io_1^i:F_1^i\ra(\pi^i)^*(TV^i)$ and
$\ti\io_2^i:F_2^i\ra(\pi^i)^*(E^i)$ with
$\ti\io_a^i\vert_{V^i\t\{0\}}=\io_a^i$ for $a=1,2$ and all $i\in I$,
where $\pi^i:V^i\t[0,1]\ra V^i$ is the projection. We claim we can
choose these isomorphisms $\ti\io_a^i$ to satisfy $(\phi^{ij})^*
(\ti\io_1^i)\ci\Phi_1^{ij}=\d\phi^{ij}\ci\ti\io_1^j$ and
$(\phi^{ij})^*(\ti\io_2^i)\ci\Phi_2^{ij}=\hat\phi^{ij}\ci\ti\io_2^j$
over $V^{ij}\t[0,1]$, as in (vi), and the obvious extension of
\eq{kh6eq5} to~$V^{ij}\t[0,1]$.

To make these choices we work by induction on $I$ in the order
$\le$. Suppose $i\in I$ and we have already chosen
$\ti\io_1^j,\ti\io_2^j$ for all $j\le i$, $j\ne i$. Then we have to
choose $\ti\io_1^i,\ti\io_2^i$ satisfying some conditions over
$\phi^{ij}(V^{ij})$ for all $j\le i$, $j\ne i$. Equation \eq{kh3eq1}
and the previously chosen compatibility between
$\ti\io_1^j,\ti\io_2^j$, $\ti\io_1^k,\ti\io_2^k$ for $k\le j\le i$
ensure that the conditions are consistent over
$\phi^{ij}(V^{ij})\cap \phi^{ik}(V^{ik})$. Thus we can choose
$\ti\io_1^i,\ti\io_2^i$ satisfying these conditions, and the claim
is true by induction.

Now define almost complex structures $J^i,K^i$ on the fibres of
$TV^i$ and $E^i$ over $V^i$ by
$J^i=(\ti\io_1^i\vert_{V^i\t\{1\}})^*(J^i_1)$ and
$K^i=(\ti\io_2^i\vert_{V^i\t\{1\}})^*(J^i_2)$. Part (vii) and the
conditions on the $\ti\io_a^i$ ensure that for all $i,j\in I$ with
$j\le i$, $j\ne i$ and $\Im\psi^i\cap\Im \psi^j\ne\emptyset$, the
following analogues of Definition \ref{kh2def23}(a)--(c) hold:
\begin{itemize}
\setlength{\itemsep}{0pt}
\setlength{\parsep}{0pt}
\item[(a)] $\d\phi^{ij}\ci J^j=(\phi^{ij})^*(J^i)\ci\d\phi^{ij}$ as
morphisms $TV^j\ra(\phi^{ij})^*(TV^i)$;
\item[(b)] $\hat\phi^{ij}\ci K^j=(\phi^{ij})^*(K^i)\ci\hat\phi^{ij}$
as morphisms $E^j\ra(\phi^{ij})^*(E^i)$; and
\item[(c)] Parts (a) and (b) imply that the orbibundles
$(\phi^{ij})^*(TV^i)/(\d\phi^{ij})(TV^j)$ and $(\phi^{ij})^*(E^i)/
\hat\phi^{ij}(E^j)$ over $V^{ij}$ have almost complex structures
$J^{ij},K^{ij}$ on their fibres, by projection from
$(\phi^{ij})^*(J^i),(\phi^{ij})^*(K^i)$. We require that
$K^{ij}\ci\d\hat s^i=\d\hat s^i\ci J^{ij}$ over $(s^j)^{-1}(0)$ in
$V^{ij}$, for $\d\hat s^i$ as in~\eq{kh2eq12}.
\end{itemize}

We can now define $\ka$ and $(\bs J,\bs K)$. For each
$p\in\oM_{g,m}(M,J,\be)$, let $i\in I$ be least in the order $\le$
such that $p\in\Im\psi^i$. Define the germ of Kuranishi
neighbourhoods of $p$ in $\ka$ to be the equivalence class of
$(V^i,E^i,s^i,\psi^i)$, regarded as a Kuranishi neighbourhood of
$p$. To define the germ of coordinate changes in $\ka$, let
$(V_p,\ldots,\psi_p)$ be sufficiently small in the germ of $\ka$ at
$p$, let $q\in\Im\psi_p$, and let $(V_q,\ldots,\psi_q)$ be
sufficiently small in the germ of $\ka$ at $q$. Let $i,j\in I$ be
least in the order $\le$ such that $p\in\Im\psi^i$
and~$q\in\Im\psi^j$.

Then by definition there are open neighbourhoods $U_p$ and $U_q$ of
$(\psi^i)^{-1}(p)$, $(\psi^j)^{-1}(q)$ in $V^i,V^{ij}\subseteq V^j$
such that $(V_p,\ldots,\psi_p)$ and $(V_q,\ldots,\psi_q)$ are
isomorphic to $(U_p,E^i\vert_{U_p}, s^i\vert_{U_p},
\psi^i\vert_{U_p})$ and $(U_q,E^j\vert_{U_q},s^j
\vert_{U_q},\psi^j\vert_{U_q})$ respectively. Define a coordinate
change $(\phi_{pq},\hat\phi_{pq})$ from $(V_q,\ldots,\psi_q)$ to
$(V_p,\ldots,\psi_p)$ to be that identified with
$(\phi^{ij}\vert_{U_q},\hat\phi^{ij}\vert_{U_q})$ by these
isomorphisms. Using Definition \ref{kh3def1}, it is easy to see that
$\ka$ is a {\it Kuranishi structure}. It is one of those allowed by
``Theorem'' \ref{kh6thm3}, which uses the same construction without
the almost complex structures.

Let $p\in\oM_{g,m}(M,J,\be)$ and $(V_p,\ldots,\psi_p)$ be
sufficiently small in the germ of $\ka$ at $p$. Define complex
structures $J_p$ on $V_p$ and $K_p$ on the fibres of $E_p$ to be
those identified with $J^i\vert_{U_p}$ and $K^i\vert_{U_p}$ by the
isomorphisms above. Then (a)--(c) above imply that Definition
\ref{kh2def23}(a)--(c) hold, so these $J_p,K_p$ define an {\it
almost complex structure\/} $(\bs J,\bs K)$ on $\bigl(\oM_{g,m}
(M,J,\be),\ka\bigr)$. This completes the first part of the theorem.

The second part works by a similar argument for
$\oM_{g,m}(M,J_t:t\in[0,1],\be)$, so we will be brief. To define
$\ka_0,(\bs J_0,\bs K_0)$ and $\ka_1,(\bs J_1,\bs K_1)$ using the
first part, we must choose good coordinate systems on
$\oM_{g,m}(M,J_0,\be),\oM_{g,m}(M,J_1,\be)$. In the second part we
choose a good coordinate system on $\oM_{g,m}(M,J_t:t\in[0,1],\be)$
which restricts to these choices for $J_0,J_1$ at $t=0,1$. Then we
go through the whole proof, choosing data
$F_a^i,\io_a^i,J^i_a,\Phi^{ij}_a,\Psi^{ij}$ and $\ti\io_a^i$ over
$\oM_{g,m}(M,J_t:t\in[0,1],\be)\t[0,1]$ for this new good coordinate
system, ensuring at each stage that the choices made agree at
$t=0,1$ with those made in defining $\ka_0,(\bs J_0, \bs K_0)$ and
$\ka_1,(\bs J_1,\bs K_1)$. This is always possible.

Note that (iii) above has to be modified as follows. Let
$(V^i,E^i,s^i,\psi^i)$ be a Kuranishi neighbourhood in the good
coordinate system for $\oM_{g,m}(M,J_t:t\in[0,1],\be)$, and let
$\bs\pi$ in ``Theorem'' \ref{kh6thm4} be represented by the
submersion $\pi^i:V^i\ra[0,1]$. Define $D^i$ to be the codimension 1
orbisubbundle $\Ker(\d\pi^i)$ in $TV^i$. In (iii) we should take
$\io^i_1$ to be an isomorphism $\io_1^i:F_1^i\vert_{V^i\t\{0\}}\ra
D^i$, mapping to $D^i$ rather than $TV^i$. Then $J^i$ becomes an
almost complex structure on the fibres of $D^i$ rather then $TV^i$,
so that $(D^i,J^i)$ is an almost CR structure on $V^i$, of
codimension 1. This is why, in the last part of the theorem, we
construct an almost CR structure $(\bs D,\bs J,\bs K)$ of
codimension 1 on $\bigl(\oM_{g,m}(M,J_t:t\in[0,1],\be),\ti\ka
\bigr)$, rather than an almost complex structure.

Finally, the modifications for the co-almost complex and co-almost
CR structures case are straightforward: we replace $TV^i$ by the
kernel of the projection to $T(M^m)$ or $T(\oM_{g,m}\t M^m)$
throughout, in particular in (iii) above. Note that these
projections are actually complex linear, with respect to the almost
complex structure $J_1^i$ on $TV^i$ defined using $D_{(v,1)}=D_v'$
above, the almost complex structure $J\t\cdots\t J$ on $M^m$, and
the natural complex structure on $\oM_{g,m}$. Hence $J_1^i$
restricts to an almost complex structure on the kernel of the
projection $TV^i\ra T(M^m)$ or $T(\oM_{g,m}\t M^m)$, as we need for
the modified (iv) above.
\end{proof}

\subsection{Gromov--Witten invariants}
\label{kh62}

In this section and \S\ref{kh63}, we assume the ``Theorems'' of
\S\ref{kh61} throughout. For our attitude to these, see
Remark~\ref{kh6rem1}.

{\it Gromov--Witten invariants} are invariants `counting'
$J$-holomorphic curves $\Si$ of genus $g$ with marked points in a
complex or symplectic manifold. They are important in String Theory,
and attracted widespread attention through their r\^ole in the
Mirror Symmetry story for Calabi--Yau 3-folds. They can be defined
either in algebraic geometry or symplectic geometry. A
Gromov--Witten invariant is basically a {\it virtual cycle} for the
moduli space $\oM_{g,m}(M,J,\be)$ of Definition \ref{kh6def2}, where
$\oM_{g,m}(M,J,\be)$ is regarded as a proper Deligne--Mumford
$\C$-stack with an obstruction theory in the algebraic case, and as
a Kuranishi space (or something similar) in the symplectic case. For
a good introduction, see Cox and Katz~\cite[\S 7]{CoKa}.

In Algebraic Geometry, Behrend \cite{Behr1,BeFa} and Li and Tian
\cite{LiTi1} defined Gromov--Witten invariants for a smooth
projective variety $M$ over $\C$. In Symplectic Geometry, Fukaya and
Ono \cite{FuOn1}, Li and Tian \cite{LiTi2}, Ruan \cite{Ruan}, and
Siebert \cite{Sieb1,Sieb2} defined Gromov--Witten invariants for a
general compact symplectic manifold $(M,\om)$, building on earlier
work of Ruan and Tian for semi-positive symplectic manifolds. It is
generally believed, and partly proved, that these definitions agree
with each other on their common domain of validity. Li and Tian
\cite{LiTi3} show that their algebraic \cite{LiTi1} and symplectic
\cite{LiTi2} invariants coincide, and Siebert \cite{Sieb3} show that
his symplectic invariants \cite{Sieb1,Sieb2} agree with Behrend's
algebraic ones \cite{Behr1}. Kontsevich and Manin \cite{KoMa} wrote
down a system of {\it axioms\/} for Gromov--Witten invariants, which
all these definitions satisfy.

One can define Gromov--Witten invariants in either homology or
cohomology; because of Poincar\'e duality there is not much
difference. We choose to define them in homology, and also in
bordism, as this simplifies issues in~\S\ref{kh63}.

\begin{dfn} Let $(M,\om)$ be a compact symplectic manifold of
dimension $2n$, let $\be\in H_2(M;\Z)$, and $g,m\ge 0$. Choose an
almost complex structure $J$ on $M$ compatible with $\om$. Set
$k=c_1(M)\cdot\be+(n-3)(1-g)+m$. Then Theorem \ref{kh6thm2} and
``Theorem'' \ref{kh6thm3} yield a compact oriented Kuranishi space
without boundary $\bigl(\oM_{g,m}(M,J,\be),\ka\bigr)$, with virtual
dimension $2k$, and strong submersions $\bev_1\t\cdots\t\bev_m:
\bigl(\oM_{g,m}(M,J,\be), \ka\bigr)\ra M^m$, and
$\bev_1\!\t\!\cdots\!\t\!\bev_m\!\t\bs\pi_{g,m}:\bigl(\oM_{g,m}(M,J,\be),
\ka\bigr)\!\ra\!M^m\t\oM_{g,m}$ if $2g\!+\!m\!\ge\!3$. ``Theorem''
\ref{kh6thm5} gives an almost complex structure $(\bs J,\bs K)$ on
$\bigl(\oM_{g,m}(M,J,\be),\ka\bigr)$, for suitably chosen $\ka$.
When $2g+m<3$, define
\begin{align*}
&GW^\Kb_{g,m}(M,\om,\be)\!\in\!KB_{2k}(M^m;\Z)\;\text{and}\;
GW^\Kb_{g,m}(M,\om,\be)\!\in\!KB_{2k}^\ac(M^m;\Z)\;\text{by}\\
&GW^\Kb_{g,m}(M,\om,\be)=\bigl[[(\oM_{g,m}(M,J,\be),\ka),
\bev_1\t\cdots\t\bev_m]\bigr],\\
&GW^\ac_{g,m}(M,\om,\be)=\bigl[[(\oM_{g,m}(M,J,\be),\ka),(\bs J,\bs
K),\bev_1\t\cdots\t\bev_m]\bigr].
\end{align*}
When $2g+m\ge 3$, define
\begin{align*}
&GW^\Kb_{g,m}(M,\om,\be)\in KB_{2k}(M^m\t\oM_{g,m};\Z)\;\>\text{and}\\
&GW^\ac_{g,m}(M,\om,\be)\in
KB_{2k}^\ac(M^m\t\oM_{g,m};\Z)\;\>\text{by}\\
&GW^\Kb_{g,m}(M,\om,\be)=\bigl[[(\oM_{g,m}(M,J,\be),\ka),
\bev_1\t\cdots\t\bev_m\t\bs\pi_{g,m}]\bigr],\\
&GW^\ac_{g,m}(M,\om,\be)\!=\!\bigl[[(\oM_{g,m}(M,J,\be),\ka),(\bs
J,\bs K),\bev_1\!\t\!\cdots\!\t\!\bev_m\!\t\bs\pi_{g,m}]\bigr].
\nonumber
\end{align*}
We call $GW^\Kb_{g,m},GW^\ac_{g,m}(M,\om,\be)$ ({\it almost
complex\/}) {\it Gromov--Witten bordism invariants}, as they lie in
(almost complex) Kuranishi bordism groups.

Define {\it Kuranishi Gromov--Witten invariants\/} by
$GW^\Kh_{g,m}(M,\om,\be)=\ab\Pi_\Kb^\Kh\ab\bigl(GW^\Kb_{g,m}
(M,\om,\be)\bigr)$, and {\it ordinary Gromov--Witten invariants\/}
$GW^\rsi_{g,m}(M,\om,\be)=\Pi_\Kb^\rsi\bigl(
GW^\Kb_{g,m}(M,\om,\be)\bigr)$, for $\Pi_\Kb^\Kh,\Pi_\Kb^\rsi$ as in
Definition \ref{kh5def5}. Then
\begin{align*}
&GW^\Kh_{g,m}(M,\om,\be)\in \begin{cases} KH_{2k}(M^m;\Q),& 2g+m<3,\\
KH_{2k}(M^m\t\oM_{g,m};\Q), & 2g+m\ge 3, \end{cases}\\
&GW^\rsi_{g,m}(M,\om,\be)\in \begin{cases} H_{2k}^\rsi(M^m;\Q), & 2g+m<3,\\
H_{2k}^\rsi(M^m\t\oM_{g,m};\Q),\,\,\, & 2g+m\ge 3.
\end{cases}
\end{align*}
Note that as $GW^\Kb_{g,m}(M,\om,\be)\!=\!\Pi_\ac^\Kb\bigl(
GW^\ac_{g,m}(M,\om,\be)\bigr)$, we also have
$GW^\Kh_{g,m}\ab(M,\ab\om,\ab\be)\!=\!\Pi_\ac^\Kh\bigl(GW^\ac_{g,m}
(M,\om,\be)\bigr)$ and~$GW^\rsi_{g,m}(M,\om,\be)\!=\!\Pi_\ac^\rsi
\bigl(GW^\ac_{g,m}(M,\om,\be)\bigr)$.
\label{kh6def3}
\end{dfn}

Our next result is analogous to Fukaya and
Ono~\cite[Th.~17.11]{FuOn1}.

\begin{thm} $GW^\Kb_{g,m},GW^\ac_{g,m},GW^\Kh_{g,m},GW^\rsi_{g,m}
(M,\om,\be)$ depend only on\/ $g,\ab m,\ab M,\ab\om,\be,$ and are
independent of choices of\/~$J,\ka,(\bs J,\bs K)$.
\label{kh6thm6}
\end{thm}

\begin{proof} Suppose $J_0,\ka_0,(\bs J_0,\bs K_0)$ and
$J_1,\ka_1,(\bs J_1,\bs K_1)$ are two possible sets of choices of
$J,\ka,(\bs J,\bs K)$ in defining $GW^\Kb_{g,m},GW^\ac_{g,m}
(M,\om,\be)$. Choose a smooth 1-parameter family of almost complex
structures $J_t$ for $t\in[0,1]$ interpolating between $J_0$ and
$J_1$. ``Theorems'' \ref{kh6thm4} and \ref{kh6thm5} then give a
Kuranishi structure $\ti\ka$ on $\oM_{g,m}(M,J_t:t\in[0,1],\be)$,
strong submersions $\widetilde\bev_1,\ldots,
\widetilde\bev_m,\ab\widetilde{\bs\pi}_{g,m}$, and an almost CR
structure $(\bs D,\bs J,\bs K)$ on~$\bigl(\oM_{g,m}(M,J_t:
t\in[0,1],\be),\ti\ka\bigr)$.

When $2g+m<3$, applying Definition \ref{kh5def6}(ii) with
$W=\bigl(\oM_{g,m}(M,J_t:t\in[0,1],\be),\ti\ka\bigr)$ and $\bs
e=\widetilde\bev_1\t\cdots\t\widetilde\bev_m$ and using \eq{kh6eq2}
shows that
\begin{align*}
\bigl[&\bigl(\oM_{g,m}(M,J_1,\be),\ka_1\bigr)\amalg
-\bigl(\oM_{g,m}(M,J_0,\be),\ka_0\bigr),\\
&\bev_1\t\cdots\t\bev_m\amalg\bev_1\t\cdots\t\bev_m\bigr]=0\quad
\text{in $KB^{2k}(M^m;\Z)$.}
\end{align*}
Using Definition \ref{kh5def6}(i) and $[-X,\bs f]=-[X,\bs f]$ as in
\eq{kh5eq1} gives
\begin{equation*}
\ts\bigl[\bigl(\oM_{g,m}(M,J_1,\be),\ka_1\bigr),\prod_{i=1}^m\bev_i
\bigr]=\bigl[\bigl(\oM_{g,m}(M,J_0,\be),\ka_0\bigr),
\prod_{i=1}^m\bev_i\bigr],
\end{equation*}
so $GW^\Kb_{g,m}(M,\om,\be)$ is independent of choices $J,\ka$. The
proofs for the other bordism cases are the same, including
$\bs\pi_{g,m},\widetilde{\bs\pi}_{g,m},(\bs J_i,\bs K_i)$ or $(\bs
D,\bs J,\bs K)$ as appropriate. The result for $GW^\Kh_{g,m},
GW^\rsi_{g,m}(M,\om,\be)$ then follows.
\end{proof}

The claim in Remark \ref{kh4rem2}(a) that the construction of the
inverse $(\Pi_\rsi^\Kh)^{-1}:KH_k(Y;\Q)\ra H_k^\rsi(Y;\Q)$ in
Theorems \ref{kh4thm1} and \ref{kh4thm2} is basically equivalent to
Fukaya and Ono's construction of {\it virtual cycles\/} for compact
oriented Kuranishi spaces without boundary in \cite[\S 6]{FuOn1}
also implies:

\begin{thm} The invariants\/ $GW^\rsi_{g,m}(M,\om,\be)$ of
Definition {\rm\ref{kh6def3}} agree with the symplectic
Gromov--Witten invariants of Fukaya and Ono\/ {\rm\cite{FuOn1}}. In
the notation of\/ {\rm\cite[\S 17]{FuOn1},} $GW^\rsi_{g,m}
(M,\om,\be)$ coincides with\/ $ev_*\bigl([C\M_{g,m}(M,J,\be)]\bigl)$
when\/ $2g+m<3,$ and with\/ $\Pi_*(C\M_{g,m}(M,J,\be))$
when\/~$2g+m\ge 3$.
\label{kh6thm7}
\end{thm}

\begin{rem}{\bf(a)} Since $\bev_1\t\cdots\t\bev_m$ and $\bev_1\t
\cdots\t\bev_m\t\bs\pi_{g,m}$ in ``Theorem'' \ref{kh6thm3} are
strong submersions, and ``Theorem'' \ref{kh6thm5} has an analogue
for co-almost complex structures, we could equally well have defined
Gromov--Witten invariants $GW^\Kcb_{g,m},GW^\ca_{g,m},GW^\Kch_{g,m},
GW^\cs_{g,m}(M,\om,\be)$ in (almost complex) Kuranishi cobordism,
Kuranishi cohomology, and compactly-supported cohomology.

In fact, as in Kontsevich and Manin \cite{KoMa}, it is best to
regard Gromov--Witten invariants as maps $H^*(M;\Q)^{\ot^m}\ra
H^*(\oM_{g,m};\Q)$, and so as elements of $H_*(M^m;\Q)\ot
H^*(\oM_{g,m};\Q)$, mixing homology and cohomology. This suggests
that in the {\it bivariant theories\/} discussed in \S\ref{kh48}, we
should regard Gromov--Witten invariants as elements
of~$H^*(\pi:M^m\t\oM_{g,m}\ra\oM_{g,m};\Q)$.

We have opted for defining Gromov--Witten invariants in bordism and
homology, rather than cobordism and cohomology or a mixture of the
two, because of our applications in \S\ref{kh63} to {\it integrality
properties\/} of Gromov--Witten invariants. Essentially, we will ask
the question: does $GW^\rsi_{g,m}(M,\om,\be)$ lie in the image of
the projection $H_{2k}^\rsi(M^m\t\oM_{g,m};\Z)\ra H_{2k}^\rsi
(M^m\t\oM_{g,m};\Q)$, and if not, can we modify $GW^\rsi_{g,m}
(M,\om,\be)$ so that it does lie in the image?

Since $\oM_{g,m}$ is an orbifold, working in cobordism and
cohomology has two disadvantages compared with bordism and homology.
Firstly, as in \S\ref{kh46}, we have not even defined maps
$\Pi_\cs^\ec: H^*_\cs(Y;\Z)\ra KH^*_\ec(Y;\Z)$ for orbifolds $Y$,
nor proved them to be isomorphisms, so we would have a problem going
from $KH^*_\ec(M^m\t\oM_{g,m};\Z)$ to $H^*_\cs(M^m\t\oM_{g,m};\Z)$.
Secondly, the `blow up functor' $B:KB_*^\ac(Y;\Z)\ra
KB_*^\eac(Y;\Z)$ discussed in \S\ref{kh63} does not seem to have a
good cobordism analogue $B:KB^*_\ca(Y;\Z)\ra KB^*_\eca(Y;\Z)$ when
$Y$ is an orbifold.
\smallskip

\noindent{\bf(b)} We claim that the whole of symplectic closed
Gromov--Witten theory --- Gromov--Witten invariants, Quantum
Cohomology, and so on --- will become simpler and more streamlined
if it is rewritten using Kuranishi (co)homology rather than
(co)homology; and simpler still if we lift to Kuranishi (co)bordism.

That is, in Kuranishi (co)homology and Kuranishi (co)bordism we do
not have to perturb moduli spaces to obtain a virtual cycle; the
moduli spaces $\oM_{g,m}(M,J,\be)$ are their own virtual classes.
For Kuranishi (co)homology we must choose (co-)gauge-fixing data for
$\oM_{g,m}(M,J,\be)$, but this is a much milder process than
perturbing moduli spaces; in Kuranishi (co)bordism there are no
extra choices to make. So for instance, the proof that
Gromov--Witten invariants satisfy the Kontsevich--Manin axioms
\cite{KoMa} would be simplified, as one does not have to worry about
the effects of perturbing the moduli spaces.
\smallskip

\noindent{\bf(c)} As in \S\ref{kh57} Kuranishi bordism groups are
{\it huge}, much larger than homology groups. Therefore we expect
$GW^\Kb_{g,m},GW^\ac_{g,m}(M,\om,\be)$ in these groups {\it to
contain more information\/} than conventional Gromov--Witten
invariants. For example, we can apply the operators $\Pi^{\Ga,\rho}$
of \S\ref{kh56} to $GW^\Kb_{g,m},GW^\ac_{g,m}(M,\om,\be)$ and then
project to singular homology, to obtain Gromov--Witten type
invariants counting $\Ga$-invariant stable maps.
\smallskip

\noindent{\bf(d)} Our Gromov--Witten bordism invariants
$GW^\Kb_{g,m},GW^\ac_{g,m}(M,\om,\be)$ are defined {\it in groups
over\/} $\Z$, not over $\Q$. Therefore they are good tools for
studying {\it integrality properties\/} of Gromov--Witten
invariants. We discuss this in~\S\ref{kh63}.
\label{kh6rem3}
\end{rem}

\subsection{Integrality properties of Gromov--Witten invariants}
\label{kh63}

Gromov--Witten invariants lie in {\it rational\/} homology
$H_{2k}^\rsi(M^m;\Q)$ or $H_{2k}^\rsi(M^m\t\oM_{g,m};\Q)$, rather
than integral homology. Essentially this is because if we think of
Gromov--Witten invariants as `counting' $J$-holomorphic maps
$w:\Si\ra Y$ with the images of the marked points $w(z_i)$ lying on
cycles $C_i$ in $Y$ for $i=1,\ldots,m$, then $[\Si,\vec z,w]$ should
be counted with rational weight $\md{\Aut(\Si,\vec z,w)}^{-1}$. Thus
the reason why Gromov--Witten invariants are defined over $\Q$
rather than $\Z$ is because of points $[\Si,\vec z,w]\in
\oM_{g,m}(M,J,\be)$ with nontrivial finite automorphism groups
$\Aut(\Si,\vec z,w)\ne\{1\}$. Equivalently, it is because of {\it
nontrivial orbifold strata\/} $\oM_{g,m}(M,J,\be)^{\Ga,\rho}$ in
$\oM_{g,m}(M,J,\be)$, in the sense of~\S\ref{kh56}.

To better understand these orbifold strata, we define some notation.

\begin{dfn} We work in the situation of \S\ref{kh61}. Suppose
$[\Si,\vec z]\in\oM_{g,m}$ and $\rho:\Ga\ra\Aut(\Si,\vec z)$ is an
injective morphism. Define $\Si'=\Si/\Ga$, with projection
$\pi:\Si\ra\Si'$. Then there is a natural way to give $\Si'$ the
structure of a prestable Riemann surface, with genus $g'\le g$, such
that $\pi$ is a holomorphic branched cover. Define $z_i'=\pi(z_i)$
for $i=1,\ldots,m$, and $\vec z'=(z_1',\ldots,z_m')$. Then
$[\Si',\vec z']\in\oM_{g',m}$. Define $\oM_{g,m}^\Ga$ to be the
moduli space of isomorphism classes $[\Si,\vec z,\rho]$ of triples
$(\Si,\vec z,\rho)$ as above, and define a smooth map
$Q^\Ga:\oM_{g,m}^\Ga \ra\coprod_{g'\le g}\oM_{g',m}$
by~$Q^\Ga:[\Si,\vec z,\rho]\mapsto[\Si',\vec z']$.

Now consider a Kuranishi space $\oM_{g,m}(M,J,\be)$, and one of its
orbifold strata $\oM_{g,m}(M,J,\be)^{\Ga,\rho}$. Points of
$\oM_{g,m}(M,J,\be)^{\Ga,\rho}$ consist of points $[\Si,\vec z,w]\in
\oM_{g,m}(M,J,\be)$ with an injective group morphism
$\rho:\Ga\ra\Aut(\Si,\vec z,w)$. Define $\Si'=\Si/\Ga$ and $\vec z'$
as above. Then $w:\Si\ra M$ factors as $w=w'\ci\pi$ for some
$J$-holomorphic $w':\Si'\ra M'$. Write $\be'=w_*([\Si'])\in
H_2(M;\Z)$. If the action of $\Ga$ on $\Si$ is locally effective
then $\be'=\md{\Ga}^{-1}\be$, but in general $\be'$ is not simply
related to $\be$, although we do have $[\om]\cdot\be'\le[\om]\cdot
\be$. Thus we can form $[\Si',\vec z',w']\in\oM_{g',m}(M,J,\be')$.
Define a map $Q^{\Ga,\rho}:\oM_{g,m}(M,J,\be)^{\Ga,
\rho}\!\ra\!\coprod_{g',\be'}\oM_{g',m}(M,J,\be')$
by~$Q^{\Ga,\rho}:[\Si,\vec z,w,\rho]\!\mapsto\![\Si',\vec z',w']$.
\label{kh6def4}
\end{dfn}

We now make the following, slightly tentative claim. It may need
modification over points $(\Si,\vec z,w),\rho$ where $\Aut(\Si,\vec
z,w)$ is finite but $\Aut(\Si,\vec z)$ is not, so we have to pass to
the {\it stabilizations\/} of $(\Si,\vec z),(\Si',\vec z')$, as in
Definition~\ref{kh6def1}.

\begin{cla} In Definition {\rm\ref{kh6def4},} $Q^{\Ga,\rho}$
lifts to a strongly smooth map $\bs Q^{\Ga,\rho}:\oM_{g,m}(M,\ab
J,\ab\be)^{\Ga,\rho}\ra \coprod_{g',\be'}\oM_{g',m}(M,J,\be')$.
Furthermore, $\oM_{g,m}(M,J,\be)^{\Ga,\rho}$ is naturally isomorphic
as a Kuranishi space to an open and closed set in the fibre product
of Kuranishi spaces
\begin{equation*}
\ts\oM_{g,m}^\Ga\t_{Q,\coprod_{g'\le
g}\oM_{g',m},\coprod_{g',\be'}\bs\pi_{g',m}} \coprod_{g',
\be'}\oM_{g',m}(M,J,\be'),
\end{equation*}
such that\/ $\bs Q^{\bs\Ga,\rho}$ is identified with projection to
the second factor in the fibre product. Hence the orbifold strata
$\oM_{g,m}(M,J,\be)^{\Ga,\rho}$ of\/ $\oM_{g,m}(M,J,\be)$ can be
completely described, up to isomorphism as Kuranishi spaces, in
terms of other moduli spaces $\oM_{g',m}(M,J,\be')$ with\/ $g'\le g$
and\/~$[\om]\cdot\be'\le[\om]\cdot\be$.
\label{kh6cla1}
\end{cla}

Assuming Claim \ref{kh6cla1}, we see that applying the operators
$\Pi^{\Ga,\rho}$ of \S\ref{kh56} to the Gromov--Witten bordism
invariants $GW^\Kb_{g,m},GW^\ac_{g,m}(M,\om,\be)$ should give
answers which can be expressed in terms of other Gromov--Witten
bordism invariants $GW^\Kb_{g',m},GW^\ac_{g',m}(M,\om,\be')$.
Projecting to homology, we expect that $\Pi_\Kb^\rsi\ci
\Pi^{\Ga,\rho}\bigl(GW^\Kb_{g,m}(M,\om,\be)\bigr)$ and
$\Pi_\ac^\rsi\ci\Pi^{\Ga,\rho}\bigl(GW^\ac_{g,m}(M,\om,\be)\bigr)$
should be expressed in terms of other Gromov--Witten
invariants~$GW^\rsi_{g',m}(M,\om,\be')$.

These ideas suggest an approach to the question of whether one can
define Gromov--Witten type invariants in {\it integral\/} homology.
Since Gromov--Witten invariants fail to lie in integral homology
because of nontrivial orbifold strata $\oM_{g,m}(M,J,\be)^{\Ga,
\rho}$, and (assuming Claim \ref{kh6cla1}) these orbifold strata can
be described in terms of $\oM_{g',m}(M,J,\be')$ with $g'\le g$ and
$[\om]\cdot\be'\le[\om]\cdot\be$, then perhaps one could modify the
Gromov--Witten invariant $GW^\rsi_{g,m}(M,\om,\be)$ by adding on a
series of `corrections' involving `lower' Gromov--Witten invariants
$GW^\rsi_{g',m}(M,\om,\be')$, to end up with a Gromov--Witten type
invariant which lies in the image of integral homology in rational
homology.

Much progress has been made on these issues, which we briefly
discuss.
\smallskip

\noindent{\bf Semi-positive symplectic manifolds and genus zero
invariants.} A homology class $\be\in H_2(M;\Z)$ is called {\it
spherical\/} if it lies in the image of $\pi_2(M)\ra H_2(M;\Z)$. A
symplectic $2n$-manifold $(M,\om)$ is called {\it semi-positive} if
for all spherical $\be\in H_2(M;\Z)$, we never have
$[\om]\cdot\be>0$ and $3-n\le c_1(M)\cdot\be<0$. This applies in
particular whenever $n\le 3$. Suppose $(M,\om)$ is semi-positive,
and $\be\in H_2(M;\Z)$ may not be written as $\be=l\ga$ for
spherical $\ga\in H_2(M;\Z)$ and $l\ge 2$ an integer. Then the genus
zero Gromov--Witten invariants $GW^\rsi_{0,m}(M,\om,\be)$ can be
defined in {\it integral\/} homology $H_{2k}^\rsi(M^m;\Z)$
or~$H_{2k}^\rsi(M^m\t\oM_{0,m};\Z)$.

This is proved in detail in McDuff and Salamon \cite[\S 7]{McSa}. In
our notation, the reason is that if $J$ is a generic almost complex
structure, then one can show that all nonempty orbifold strata
$\oM_{0,m}(M,J,\be)^{\Ga,\rho}$ of $\oM_{0,m}(M,J,\be)$ with
$\Ga\ne\{1\}$ satisfy $\vdim\oM_{0,m}(M,J,\be)^{\Ga,\rho}\le
\vdim\oM_{0,m}(M,J,\be)-2$. In a very similar way to our use of {\it
effective\/} Kuranishi spaces to define a virtual cycle map
$\Pi_\eb^\rsi:KB_*(Y;\Z)\ra H_*^\rsi(Y;\Z)$, this condition on the
orbifold strata of $\oM_{0,m}(M,J,\be)$ is enough to ensure that its
virtual cycle can be defined over $\Z$ rather than $\Q$. In fact, in
Remark \ref{kh5rem5} we showed that if $X$ is an oriented, effective
Kuranishi space then $\vdim X^{\Ga,\rho}\le\vdim X-2$ for all
nonempty orbifold strata $X^{\Ga,\rho}$ of $X$ with~$\Ga\ne\{1\}$.
\smallskip

\noindent{\bf Gopakumar--Vafa invariants and their Integrality
Conjecture.} Let $(M,\om)$ be a compact symplectic 6-manifold with
$c_1(M)=0$, that is, a {\it symplectic Calabi--Yau $3$-fold}. Take
the number $m$ of marked points to be zero. Then \eq{kh6eq1} gives
$\vdim\oM_{g,0}(M,J,\be)=0$ for all $g,\be$. Hence the
Gromov--Witten invariants $GW^\rsi_{g,0}(M,\om,\be)$ lie in
$H_0^\rsi(M^m;\Q)\cong\Q$ or $H_0^\rsi(M^m\t\oM_{g,0};\Q)\cong\Q$.
So we regard the $GW^\rsi_{g,0}(M,\om,\be)$ as {\it rational
numbers\/} for all~$g,\be$.

Using physical reasoning about counting BPS states, the String
Theorists Gopakumar and Vafa \cite{GoVa1,GoVa2} conjectured the
existence of {\it Gopakumar--Vafa invariants\/}
$GV_g(M,\om,\be)\in\Z$ for Calabi--Yau 3-folds $M$, which morally
speaking `count' {\it embedded\/} $J$-holomorphic curves of genus
$g$ in class $\be$ in $M$. By expressing any $J$-holomorphic curve
in $M$ as the branched cover of an embedded curve, they
conjecturally wrote Gromov--Witten in terms of Gopakumar--Vafa
invariants, and vice versa, by the equation in formal power series
\e
\begin{split}
&\ts\sum_{g,\be}GW^\rsi_{g,0}(M,\om,\be)\, t\sp{2g-2}q^{\be}=\\
&\ts\sum_{k>0,\,g,\be}GV_g(M,\om,\be)\,\frac{1}{k}(2\sin(k
t/2))^{2g-2}q^{k\be}.
\end{split}
\label{kh6eq6}
\e
The {\it Gopakumar--Vafa Integrality Conjecture} says Gromov--Witten
invariants of Calabi--Yau 3-folds satisfy \eq{kh6eq6} for some
integers $GV_g(M,\om,\be)$. Gopakumar and Vafa also conjecture that
$GV_g(M,\om,\be)=0$ for all fixed $\be$ and~$g\gg 0$.

There are two approaches to proving this conjecture. The first would
be to define some curve-counting invariants $GV_g(M,\om,\be)$ which
are automatically integers, and then prove that they satisfy
\eq{kh6eq6}. Pandharipande and Thomas \cite{PaTh1,PaTh2} have taken
an important step in this direction, in the context of algebraic
rather than symplectic geometry, by defining integer-valued
invariants which count `stable pairs' $(F,s)$ of a coherent sheaf
$F$ supported on a curve in $M$ and a section $s\in H^0(F)$. Proving
equivalence of Pandharipande--Thomas invariants with Gromov--Witten
invariants, the analogue of \eq{kh6eq6}, still seems difficult.

The second is to regard \eq{kh6eq6} as a {\it definition\/} of
numbers $GV_g(M,\om,\be)\in\Q$, and then to prove that these
$GV_g(M,\om,\be)$ actually lie in $\Z$. We will sketch a method of
attack following this approach below.
\smallskip

\noindent{\bf More general integrality conjectures.} Pandharipande
\cite{Pand1,Pand2} extends Go\-pa\-ku\-mar--Vafa invariants and
their integrality conjecture to all smooth projective complex
algebraic 3-folds. Klemm and Pandharipande \cite{KlPa} and
Pandharipande and Zinger \cite{PaZi} study integrality of
Gromov--Witten invariants of Calabi--Yau $m$-folds for $m>3$. In
these cases, for dimension reasons, only genus $g=0$ and $g=1$ can
yield nonzero invariants. An integrality conjecture for $g=0$ and
all $m>3$ is given in \cite{KlPa}, and conjectures for $g=1$, $m=4$
and $g=1$, $m=5$ are given in \cite{KlPa} and \cite{PaZi},
respectively.
\medskip

The rest of the section outlines an approach to proving integrality
results for Gromov--Witten invariants, including the Gopakumar--Vafa
Integrality Conjecture, and the conjectures of Pandharipande et al.\
\cite{KlPa,Pand1,Pand2,PaZi}. We stress that this approach is still
very incomplete, and we are not yet claiming a proof of any of these
conjectures. The author hopes to fill in some of the gaps in the
sketch in \cite{Joyc4}. We divide our method into four steps. Steps
1--3 are very general, and work in all dimensions.

\subsubsection{Step 1: a `blow-up functor' $B:KB_*^\ac(Y;R)\ra
KB_*^\eac(Y;R)$}
\label{kh631}

Let $X$ be a compact Kuranishi space without boundary, with an
almost complex structure $(\bs J,\bs K)$. We would like to find a
functorial procedure to modify $X,(\bs J,\bs K)$ to yield an {\it
effective\/} compact Kuranishi space $\ti X$ without boundary, with
an almost complex structure $(\bs{\ti J},\bs{\ti K})$, and with a
strongly smooth map $\bs\pi_X:\ti X\ra X$. Now $X$ fails to be
effective if its {\it orbifold strata\/} $X^{\Ga,\rho}$ have the
wrong behaviour, in particular, if $\rho$ is not the isomorphism
class of an effective $\Ga$-representation and $X^{\Ga,\rho}\ne\es$
then $X$ is not effective. Thus, to make $X$ effective we must
modify it around its orbifold strata~$X^{\Ga,\rho}$.

Suppose $\md{\Ga}$ is largest such that $X$ has an orbifold stratum
$X^{\Ga,\rho}$ along which $X$ is not effective. Then $X^{\Ga,\rho}$
is in some sense `nonsingular' in $X$. In this case the author has a
procedure for `blowing up' an almost complex Kuranishi space $X,(\bs
J,\bs K)$ at $X^{\Ga,\rho}$ with respect to $(\bs J,\bs K)$, to get
a new almost complex Kuranishi space $\hat X,(\bs{\hat J},\bs{\hat
K})$ with strongly smooth $\bs\pi_X:\hat X\ra X$. We aim to do the
blow up in such a way that $\hat X$ is `closer to being effective'
than~$X$.

If $(V_p,E_p,s_p,\psi_p)$ is a sufficiently small Kuranishi
neighbourhood in the germ at $p\in X$ with $(J_p,K_p)$ representing
$K_p$ then we define a Kuranishi neighbourhood $(\hat V_p,\hat
E_p,\hat s_p,\hat\psi_p)$ on $\hat X$, with $(\hat J_p,\hat K_p)$
representing $(\bs{\hat J},\bs{\hat K})$, where $(\hat V_p,\hat
J_p)$ is the blow up of $(V_p,J_p)$ along $V_p^{\Ga,\rho}$, in the
usual sense of (almost) complex geometry. (Note a potential
confusion here: we are {\it not\/} treating the orbifold strata of
$V_p$ as singularities, and resolving them to get a manifold $\hat
V_p$. In fact $\hat V_p$ is an orbifold, and its orbifold strata may
be no simpler than those of $V_p$.) The crucial part of the
construction is to define $\hat E_p,\hat s_p$ in such a way that
these blow ups are compatible with $(\phi_{pq},\hat\phi_{pq})$ in
the germ of coordinate changes on~$X$.

By applying this blow up procedure iteratively by reverse induction
on $\md{\Ga}$ to modify all orbifold strata of $X$ to make them
effective, we make the:

\begin{cla} There is a procedure to modify $X,(\bs J,\bs K)$ to
get\/ $\ti X,(\bs{\ti J},\bs{\ti K}),\bs\pi_X$ as above, with\/ $\ti
X$ effective. This induces a functor $B:KB^\ac_{2k}(Y;R)\!\ra\!
KB^\eac_{2k}(Y;R)$ mapping $B:[X,(\bs J,\bs K),\bs f]\!\mapsto\![\ti
X,(\bs{\ti J},\bs{\ti K}),\bs f\ci\bs\pi_X],$ such that if\/ $X$ has
trivial stabilizers then~$B\bigl([X,(\bs J,\bs K),\bs
f]\bigr)=[X,(\bs J,\bs K),\bs f]$.
\label{kh6cla2}
\end{cla}

Applying the functor $B$ to the invariants $GW^\ac_{g,m}
(M,\om,\be)$ of \S\ref{kh62} and composing with $\Pi_\eac^\rsi$
yields invariants $\Pi_\eac^\rsi\ci B\bigl(GW^\ac_{g,m}
(M,\om,\be)\bigr)$ in {\it integral\/} homology $H_{2k}^\rsi
(M^m;\Z)$ or $H_{2k}^\rsi(M^m\t\oM_{g,m};\Z)$, since $\Pi_\eac^\rsi$
maps $KB^\eac_*(Y;\Z)\ra H^\rsi_*(Y;\Z)$, although $\Pi_\ac^\rsi$
maps $KB^\ac_*(Y;\Z)\ra H^\rsi_*(Y;\Q)$. So we can think of the
$\Pi_\eac^\rsi\ci B\bigl(GW^\ac_{g,m}(M,\om,\be)\bigr)$ as
Gopakumar--Vafa type invariants, curve counting invariants that are
automatically integer-valued.

This is similar to an idea sketched by Fukaya and Ono \cite{FuOn2}.
Translated into our notation, roughly speaking, Fukaya and Ono claim
that given a compact Kuranishi space $X$ without boundary and an
almost complex structure $(\bs J,\bs K)$, then one should be able to
carry out their virtual cycle construction for $X$ using
single-valued sections rather than multisections, by using $(\bs
J,\bs K)$ to determine a canonical form for the perturbation near
each orbifold stratum $X^{\Ga,\rho}$. As the sections are
single-valued, the virtual cycle lies in homology over~$\Z$.

\subsubsection{Step 2: an expression for
$\Pi_\eac^\Kb\ci B:KB_*^\ac(Y;\Z)\ra KB_*(Y;\Z)$}
\label{kh632}

Let $[X,(\bs J,\bs K),\bs f]\in KB^\ac_{2k}(Y;\Z)$. If as in
\S\ref{kh631} we blow up $X,(\bs J,\bs K)$ along $X^{\Ga,\rho}$ to
get $\hat X,(\bs{\hat J},\bs{\hat K}),\bs\pi_X$, we can form $[\hat
X,(\bs{\hat J},\bs{\hat K}),\bs f\ci\bs\pi_X]\in KB^\ac_{2l}(Y;\Z)$.
There is a notion of virtual vector bundle over a Kuranishi space
due to Fukaya and Ono \cite[Def.~5.10]{FuOn1}, called a {\it bundle
system}. The {\it virtual normal bundle\/} $\bs\nu_{X^{\Ga,\rho}}$
of $X^{\Ga,\rho}$ in $X$ is a complex bundle system. We expect the
difference $[\hat X,(\bs{\hat J},\bs{\hat K}),\bs
f\ci\bs\pi_X]-[X,(\bs J,\bs K),\bs f]$ should depend only on
$X^{\Ga,\rho}$ and $\bs\nu_{X^{\Ga,\rho}}$.

This is more likely to work if we project to Kuranishi bordism
$KB_{2k}(Y;\Z)$, as constructing the relevant bordism $W$ between
Kuranishi spaces will be easier if $W$ need not carry an almost CR
structure. Thus we expect that
\begin{equation*}
\Pi_\ac^\Kb\bigl([\hat X,(\bs{\hat J},\bs{\hat K}),\bs
f\!\ci\!\bs\pi_X]\!-\![X,(\bs J,\bs K),\bs f]\bigr)\!=\!\bigl[
C^{\Ga,\rho}(X^{\Ga,\rho},\bs\nu^{\Ga,\rho}),\bs
f\!\ci\!\bs\io^{\Ga,\rho}\!\ci\!\bs\ga^{\Ga,\rho}\bigr],
\end{equation*}
where $C^{\Ga,\rho}(X^{\Ga,\rho},\bs\nu^{\Ga,\rho})$ is some compact
Kuranishi space without boundary constructed from $X^{\Ga,\rho}$ and
$\bs\nu^{\Ga,\rho}$ in some functorial way, with $\vdim
C^{\Ga,\rho}(X^{\Ga,\rho},\bs\nu^{\Ga,\rho})\ab=\vdim X$, and
$\bs\ga^{\Ga,\rho}: C^{\Ga,\rho}(X^{\Ga,\rho},\bs\nu^{\Ga,\rho})\ra
X^{\Ga,\rho}$ is a cooriented strong submersion. Define
$\Pi^{\Ga,\rho,C^{\Ga,\rho}}:KB_{2k}^\ac(Y;\Z)\ra KB_{2k}(Y;\Z)$ by
\begin{equation*}
\Pi^{\Ga,\rho,C^{\Ga,\rho}}:[X,(\bs J,\bs K),\bs f]\mapsto\bigl[
C^{\Ga,\rho}(X^{\Ga,\rho},\bs\nu^{\Ga,\rho}),\bs f\ci
\bs\io^{\Ga,\rho}\ci\bs\ga^{\Ga,\rho}\bigr].
\end{equation*}

The functor $B$ is the result of composing many such blowups along
$X^{\Ga,\rho}$ for different $\Ga,\rho$. This suggests the:

\begin{cla} In operators $KB_{2l}^\ac(Y;\Z)\ra KB_{2l}(Y;\Z)$ we have
\e
\Pi_\eac^\Kb\ci B=\Pi_\ac^\Kb+
\ts\sum_{\Ga,\rho}\Pi^{\Ga,\rho,C^{\Ga,\rho}},
\label{kh6eq7}
\e
where the sum is over isomorphisms classes of finite groups $\Ga$
with\/ $\Ga\ne\{1\}$ and isomorphism classes of virtual nontrivial\/
$\Ga$-representations $\rho,$ and each\/ $\Pi^{\Ga,\rho,C^{\Ga,
\rho}}$ is an operator of the form above.
\label{kh6cla3}
\end{cla}

\subsubsection{Step 3: projecting to singular homology}
\label{kh633}

Apply equation \eq{kh6eq7} to the Gromov--Witten type invariant
$GW^\ac_{g,m}(M,\om,\be)$, and compose with projection
$\Pi_\Kb^\rsi$ to rational singular homology. This yields
\e
\begin{split}
\Pi_\eac^\rsi\ci B\bigl(GW^\ac_{g,m}(M,\om,\be)\bigr)=\,&
GW^\rsi_{g,m}(M,\om,\be)+\\
&\ts\sum_{\Ga,\rho}\Pi_\Kb^\rsi\ci\Pi^{\Ga,\rho,C^{\Ga,\rho}}
\bigl(GW^\ac_{g,m}(M,\om,\be)\bigr).
\end{split}
\label{kh6eq8}
\e
Here the left hand side of \eq{kh6eq8} lies in the image of integral
homology in rational homology, $\pi_*:H_{2k}^\rsi(M^m;\Z)\ra
H_{2k}^\rsi(M^m;\Q)$ or $\pi_*:H_{2k}^\rsi(M^m\t\oM_{g,m};\Z)\ra
H_{2k}^\rsi(M^m\t\oM_{g,m};\Q)$. The first term on the right hand
side is an ordinary Gromov--Witten invariant, and the remaining
terms are corrections from the nontrivial orbifold strata
$\oM_{g,m}(M,J,\be)^{\Ga,\rho}$ of~$\oM_{g,m}(M,J,\be)$.

Now Claim \ref{kh6cla1} says that we can write the
$\oM_{g,m}(M,J,\be)^{\Ga,\rho}$ in terms of other `lower' moduli
spaces $\oM_{g',m}(M,J,\be')$. Thus it seems reasonable to hope that
we can write the corrections $\Pi_\Kb^\rsi\ci\Pi^{\Ga,\rho,
C^{\Ga,\rho}}\bigl(GW^\ac_{g,m}(M,\om,\be)\bigr)$ in terms of other
`lower' ordinary Gromov--Witten invariants $GW^\rsi_{g',m}(M,\om,
\be')$. However, there is a problem. The $\Pi^{\Ga,\rho,C^{\Ga,
\rho}}\bigl(GW^\ac_{g,m}(M,\om,\be)\bigr)$ depend not only on the
orbifold strata $\oM_{g,m}(M,J,\be)^{\Ga,\rho}$, but also on its
{\it virtual normal bundle\/} $\bs\nu^{\Ga,\rho}$ in $\oM_{g,m}
(M,J,\be)$. The author does {\it not\/} expect $\bs\nu^{\Ga,\rho}$
to depend only on $\oM_{g',m}(M,J,\be')$, it is extra information.
But we claim that if $\dim\rho\ge 0$ then $\Pi_\Kb^\rsi\ci\Pi^{\Ga,
\rho,C^{\Ga,\rho}}\bigl(GW^\ac_{g,m}(M,\om,\be)\bigr)$ does not
depend on~$\bs\nu^{\Ga,\rho}$.

\begin{cla} In operators $KB^\ac_{2k}(Y;\Z)\!\ra\!H^\rsi_{2k}(Y;\Q)$
we have $\Pi_\Kb^\rsi\!\ci\!\Pi^{\Ga,\rho,C^{\Ga,\rho}}\!=\!0$ if\/
$\dim\rho>0,$ and\/ $\Pi_\Kb^\rsi\ci\Pi^{\Ga,\rho,C^{\Ga,\rho}}\!=\!
q^{\Ga,\rho}\Pi_\ac^\rsi\!\ci\!\Pi^{\Ga,\rho}_\ac$ if\/
$\dim\rho=0,$ where $q^{\Ga,\rho}\in\Q$ and\/
$\Pi^{\Ga,\rho}_\ac:KB^\ac_{2k}(Y;\Z)\!\ra\!KB^\ac_{2k}(Y;\Z)$ is as
in~{\rm\S\ref{kh56}}.
\label{kh6cla4}
\end{cla}

To justify this, note that $\Pi_\Kb^\rsi\ci\Pi^{\Ga,\rho,
C^{\Ga,\rho}}\bigl([X,(\bs J,\bs K),\bs f]\bigr)$ is represented by
a virtual cycle $VC\bigl(C^{\Ga,\rho}(X^{\Ga,\rho},\bs\nu^{\Ga,
\rho})\bigr)$ for $\bigl(C^{\Ga,\rho}(X^{\Ga,\rho},\bs\nu^{\Ga,
\rho}),\bs f\ci\bs\io^{\Ga,\rho}\ci\bs\ga^{\Ga,\rho}\bigr)$, and
$\bs f\ci\bs\io^{\Ga,\rho}\ci\bs\ga^{\Ga,\rho}$ factorizes through
$\bs\ga^{\Ga,\rho}:C^{\Ga,\rho}(X^{\Ga,\rho},\bs\nu^{\Ga,\rho})\ra
X^{\Ga,\rho}$, where $\vdim C^{\Ga,\rho}(X^{\Ga,\rho},\bs\nu^{\Ga,
\rho})=\vdim X$ and $\vdim X^{\Ga,\rho}=\vdim X-\dim\rho$.

If $VC(X^{\Ga,\rho})$ is a virtual cycle for $(X^{\Ga,\rho},\bs
f\ci\bs\io^{\Ga,\rho})$, then we can choose $VC\bigl(C^{\Ga,\rho}
(X^{\Ga,\rho},\bs\nu^{\Ga,\rho})\bigr)$ to be supported in an
arbitrarily small open neighbourhood of the support of
$VC(X^{\Ga,\rho})$. If $\dim\rho>0$ then $\dim VC\bigl(C^{\Ga,\rho}
(X^{\Ga,\rho},\bs\nu^{\Ga,\rho})\bigr)\ab>\dim VC(X^{\Ga,\rho})$, so
$\bigl[VC\bigl(C^{\Ga,\rho} (X^{\Ga,\rho},\bs\nu^{\Ga,\rho})\bigr)
\bigr]=0$ in $H^\rsi_{2k}(Y;\Q)$ for dimensional reasons. If
$\dim\rho=0$ then $\bigl[VC\bigl(C^{\Ga,\rho} (X^{\Ga,\rho},
\bs\nu^{\Ga,\rho})\bigr)\bigr]=q^{\Ga,\rho}\bigl[VC(X^{\Ga,\rho})
\bigr]$, where $q^{\Ga,\rho}\in\Q$ is the `multiplicity'
of~$\bs\ga^{\Ga,\rho}$.

Applying Claim \ref{kh6cla4} to \eq{kh6eq8} yields
\e
\begin{split}
\Pi_\eac^\rsi\ci B\bigl(GW^\ac_{g,m}&(M,\om,\be)\bigr)=
GW^\rsi_{g,m}(M,\om,\be)+\\
&\ts\sum_{\Ga,\rho:\,\dim\rho=0}q^{\Ga,\rho}\Pi_\ac^\rsi
\ci\Pi^{\Ga,\rho}_\ac\bigl(GW^\ac_{g,m}(M,\om,\be)\bigr)+\\
&\ts\sum_{\Ga,\rho:\,\dim\rho<0}\Pi_\Kb^\rsi\ci\Pi^{\Ga,\rho,
C^{\Ga,\rho}}\bigl(GW^\ac_{g,m}(M,\om,\be)\bigr).
\end{split}
\label{kh6eq9}
\e
When $\dim\rho=0$, the terms $\Pi_\ac^\rsi\ci\Pi^{\Ga,\rho}_\ac
\bigl(GW^\ac_{g,m}(M,\om,\be)\bigr)$ on the second line of
\eq{kh6eq9} no longer depend on the virtual normal bundle
$\bs\nu^{\Ga,\rho}$, so by Claim \ref{kh6cla1}, we expect we can
write them solely in terms of `lower' Gromov--Witten
invariants~$GW^\rsi_{g',m}(M,\om,\be')$.

The terms on the third line of \eq{kh6eq9}, when $\dim\rho<0$, come
from orbifold strata $\oM_{g,m}(M,J,\be)^{\Ga,\rho}$ with
$\vdim\oM_{g,m}(M,J,\be)^{\Ga,\rho}>\vdim\oM_{g,m}(M,J,\be)$. The
author expects that these terms really do depend on
$\bs\nu^{\Ga,\rho}$, and so {\it cannot\/} be written solely in
terms of $GW^\rsi_{g',m}(M,\om,\be')$, except in special
circumstances. In effect, these terms use `characteristic classes'
of $\bs\nu^{\Ga,\rho}$ --- some kind of analogue of Chern classes,
perhaps --- to lower the dimension by~$-\dim\rho$.

\subsubsection{Step 4: the Gopakumar--Vafa Integrality Conjecture}
\label{kh634}

We now specialize to the case when $(M,\om)$ is a symplectic
Calabi--Yau 3-fold, that is, a compact symplectic 6-manifold with
$c_1(M)=0$. Let $J$ be a {\it generic} almost complex structure on
$M$ compatible with $\om$. Using genericness of $J$ and $\vdim
\oM_{g,0}(M,J,\be)=0$ for all $g,\be$, one can show that for all
$[\Si,\vec z,w]\in \oM_{g,m}(M,J,\be)$, if $\be=0$ then $w\equiv p$
for some $p\in M$, and if $\be\ne 0$ then $w$ factors uniquely as
$w=\ti w\ci\pi$, where $\ti\Si$ is a {\it nonsingular} Riemann
surface of genus $\ti g\le g$, and $\ti w:\ti\Si\ra M$ is a
$J$-holomorphic {\it embedding}, and $\pi:\Si\ra\ti\Si$ is
surjective and holomorphic. Furthermore, distinct embedded
nonsingular $J$-holomorphic curves do not intersect in~$M$.

Setting $m=0$, and regarding the $GW^\rsi_{g,0}(M,\om,\be)$ as {\it
rational numbers\/} for all $g,\be$, we expect \eq{kh6eq9} to reduce
to an equation of the form
\e
\Pi_\eac^\rsi\ci B\bigl(GW^\ac_{g,0}(M,\om,\be)\bigr)=
\sum_{0\le g'\le g,\; k\ge 1:k\vert\be}c^{g',g,k}
GW^\rsi_{g',0}(M,\om,\be/k),
\label{kh6eq10}
\e
where $c^{g',g,k}\in\Q$ are some universal family of rational
numbers with $c^{g,g,1}=1$, and \eq{kh6eq10} is an equation in $\Q$,
with the left hand side in~$\Z$.

Now the Gopakumar--Vafa equation \eq{kh6eq6} may be rewritten in the
form
\e
GV_g(M,\om,\be)=
\sum_{0\le g'\le g,\; k\ge 1:k\vert\be}d^{g',g,k}
GW^\rsi_{g',0}(M,\om,\be/k),
\label{kh6eq11}
\e
where $d^{g',g,k}\in\Q$ are some universal family of explicitly
computable rational numbers with $d^{g,g,1}=1$, and the
Gopakumar--Vafa Integrality Conjecture says that the left hand side
of \eq{kh6eq11} lies in $\Z$ for all~$g,\be$.

We do not expect it will be feasible to compute the coefficients
$c^{g',g,k}$, as they are the sum over all $\Ga,\rho$ of many
complicated contributions. Nor do we claim that
$c^{g',g,k}=d^{g',g,k}$. However, we do hope to prove that the
integrality properties of the family of $GW^\rsi_{g,0}(M,\om,\be)$
for all $g,\be$ implied by \eq{kh6eq10} are equivalent to the
integrality properties implied by \eq{kh6eq11}, and this will be
enough to prove the Gopakumar--Vafa Integrality Conjecture.

By induction on $g',g,k$ using the fact that $c^{g,g,1}=d^{g,g,1}
=1$, we can show there exist unique $e^{g',g,k},f^{g',g,k}\in\Q$ for
$0\le g'\le g$ and $k\ge 1$, with $e^{g,g,1}=f^{g,g,1}=1$,
satisfying
\e
c^{g'',g,k''}=\sum_{g',k,k':\; g''\le g'\le g,\; k''=kk'
\!\!\!\!\!\!\!\!\!\!\!\!\!\!\!\!\!\!\!\!\!\!\!\!\!\!\!\!\!\!\!\!
\!\!\!\!\!\!\!\!\!\!\!\!} e^{g',g,k}d^{g'',g',k}\;\>\text{and}\;\>
d^{g'',g,k''}=\sum_{g',k,k':\; g''\le g'\le g,\; k''=kk'
\!\!\!\!\!\!\!\!\!\!\!\!\!\!\!\!\!\!\!\!\!\!\!\!\!\!\!\!\!\!\!\!
\!\!\!\!\!\!\!\!\!\!\!\!} f^{g',g,k}c^{g'',g',k}
\label{kh6eq12}
\e
for all $0\le g''\le g$ and $k''\ge 1$. Combining
\eq{kh6eq10}--\eq{kh6eq12} yields
\begin{gather}
\Pi_\eac^\rsi\ci B\bigl(GW^\ac_{g,0}(M,\om,\be)\bigr)=
\sum_{0\le g'\le g,\; k\ge 1:k\vert\be\!\!\!\!}e^{g',g,k}\,
GV_{g'}(M,\om,\be/k),
\label{kh6eq13}\\
GV_g(M,\om,\be)=
\sum_{0\le g'\le g,\; k\ge 1:k\vert\be\!\!\!\!}f^{g',g,k}\,
\Pi_\eac^\rsi\ci B\bigl(GW^\ac_{g',0}(M,\om,\be/k)\bigr).
\label{kh6eq14}
\end{gather}
Hence, if all $e^{g',g,k}$ lie in $\Z$ then the integrality of
\eq{kh6eq11} implies that of \eq{kh6eq10}, and if all $f^{g',g,k}$
lie in $\Z$ then the integrality of \eq{kh6eq10} implies that of
\eq{kh6eq11}, and if all $e^{g',g,k},f^{g',g,k}$ lie in $\Z$ then
the integrality of \eq{kh6eq10} and \eq{kh6eq11} is equivalent.

The Gopakumar--Vafa Integrality Conjecture has been verified in many
examples, for instance, for local curves by Bryan and Pandharipande
\cite{BrPa}. In each of these examples, it has already been proved
that the right hand side of \eq{kh6eq11} is an integer, and if we
can carry through our programme above, then we will have shown that
the right hand side of \eq{kh6eq10} is an integer. The author hopes
that by applying this to a large enough family of examples, one will
be able to prove by an induction on $g,g',k$ that all
$e^{g',g,k},f^{g',g,k}$ lie in $\Z$, so that integrality of
\eq{kh6eq10} implies the Gopakumar--Vafa Integrality Conjecture.

\subsection{Moduli spaces of stable holomorphic discs}
\label{kh64}

In \S\ref{kh64}--\S\ref{kh67} we discuss $J$-holomorphic curves in a
symplectic manifold $(M,\om)$ {\it with boundary in a Lagrangian
submanifold\/} $L$ in $M$. As our aim is just to give an idea of
possible applications in this area, not to develop a complete
theory, we restrict our attention to {\it stable holomorphic discs},
rather than curves of higher genus, and we consider only {\it
boundary marked points}, not interior marked points. The
generalization of this section to higher genus curves and both kinds
of marked points will be explained in \cite{Joyc2}. Most of this
section follows Fukaya et al.\ \cite[\S 2]{FOOO}, see also Katz and
Liu \cite[\S 3]{KaLi} and Liu~\cite[\S 2]{Liu}.

\begin{dfn} A {\it bordered Riemann surface\/} $\Si$ is a compact
connected complex 1-manifold with nonempty boundary $\pd\Si$. Given
such a $\Si$ we can construct the {\it double\/}
$\ti\Si=\Si\cup_{\pd\Si}\bar\Si$, that is, the union of $\Si$ and
its complex conjugate $\bar\Si$, glued along their common boundary
$\pd\Si$. Then $\ti\Si$ is a Riemann surface without boundary, and
there is an antiholomorphic involution $\si:\ti\Si\ra\ti\Si$
exchanging $\Si$ and $\bar\Si$, such that the fixed points of $\si$
are $\pd\Si$, and $\ti\Si$ is the disjoint union
$\Si^\ci\amalg\pd\Si\amalg\si(\Si^\ci)$, where $\Si^\ci$ is the
interior of~$\Si$.

A {\it prestable\/} or {\it nodal bordered Riemann surface\/} $\Si$
is a compact connected singular complex 1-manifold with nonempty
boundary $\pd\Si$, whose only singularities are nodes. Again, $\Si$
has a {\it double} $\ti\Si=\Si\cup_{\pd\Si}\bar\Si$, with an
antiholomorphic involution $\si:\ti\Si\ra\ti\Si$. We allow one
interior and two boundary kinds of nodes:
\begin{itemize}
\setlength{\itemsep}{0pt}
\setlength{\parsep}{0pt}
\item[(a)] $\Si$ has a node $z$ in $\Si^\ci$ locally modelled on
$(0,0)$ in~$\bigl\{(x,y)\in\C^2:xy=0\bigr\}$.

Such a node can occur as the limit of a family of nonsingular
bordered Riemann surfaces $\Si_\ep$ modelled on
$\bigl\{(x,y)\in\C^2:xy=\ep\bigr\}$, for~$\ep\in\C\sm\{0\}$.
\item[(b)] $\ti\Si$ has a node $z$ in $\pd\Si$ locally modelled
on $(0,0)$ in $\bigl\{(x,y)\in\C^2:xy=0\bigr\}$, where $\si$ acts by
$(x,y)\mapsto(\bar x,\bar y)$. Then $\pd\Si$ is locally modelled on
$\bigl\{(x,y)\in\R^2:xy=0\bigr\}$, the fixed point set of $\si$, and
we take $\Si$ to be locally modelled on~$\bigl\{(x,y)\in\C^2:xy=0$,
$\Im x\ge 0$, $\Im y\ge 0\bigr\}$.

Such a node can occur as the limit of a family of nonsingular
$\Si_\ep$ modelled on $\bigl\{(x,y)\in\C^2:xy=-\ep$, $\Im x\ge 0$,
$\Im y\ge 0\bigr\}$, for $\ep>0$, with doubles $\ti\Si_\ep$ modelled
on $\bigl\{(x,y)\in\C^2:xy=-\ep\bigr\}$.
\item[(c)] $\ti\Si$ has a node $z$ in $\pd\Si$ locally modelled
on $(0,0)$ in $\bigl\{(x,y)\in\C^2:xy=0\bigr\}$, where $\si$ acts by
$(x,y)\mapsto(\bar y,\bar x)$. Then $\pd\Si$ is locally modelled on
$(0,0)$, the fixed point set of $\si$, and we take $\Si$ to be
locally modelled on~$\bigl\{(x,0)\in\C^2\bigr\}$.

Such a node can occur as the limit of a family of nonsingular
$\Si_\ep$ modelled on $\bigl\{(x,y)\in\C^2:xy=\ep$,
$\md{x}\ge\md{y}\bigr\}$, for $\ep>0$, with doubles $\ti\Si_\ep$
modelled on $\bigl\{(x,y)\in\C^2:xy=\ep\bigr\}$, where $\si$ acts by
$(x,y)\mapsto(\bar y,\bar x)$, and boundary the circle $\pd\Si_\ep=
\bigl\{(\ep^{1/2}{\rm e}^{i\th},\ep^{1/2}{\rm e}^{-i\th}):
\th\in[0,2\pi)\bigr\}$.
\end{itemize}
Note that near a node of type (c), $\Si$ is actually a {\it
nonsingular\/} Riemann surface {\it without boundary}, save that one
point $z$ in $\Si$ is designated as a `boundary node'. Also $\pd\Si$
near $z$ is just one point $\{z\}$, not a real curve. We allow as
prestable bordered Riemann surfaces curves whose only boundary is
nodes of type (c), so for example $\Si$ could be $\CP^1$ with one
point designated as a boundary node.

If $\Si$ is a prestable bordered Riemann surface, the {\it
smoothing\/} $\Si'$ of $\Si$ is the compact, connected, oriented
2-manifold with boundary obtained by smoothing the nodes of $\Si$.
Note that there are two topologically distinct ways of
desingularizing nodes of type (a)--(c), and we choose those
described as families of surfaces $\Si_\ep$ above. The rule is that
smoothing each node decreases the Euler characteristic by 1, but the
other ways of desingularizing nodes of type (a)--(c) increase Euler
characteristic by 1,1,0 respectively. We call $\Si$ a {\it prestable
holomorphic disc\/} if its smoothing $\Si'$ is a holomorphic disc.
\label{kh6def5}
\end{dfn}

\begin{dfn} Suppose $(M,\om)$ is a compact symplectic manifold, $L$
a compact, embedded Lagrangian submanifold in $M$, and $J$ an almost
complex structure on $M$ compatible with $\om$. Let $\be\in
H^2(M;\Z)$ and $l\ge 0$. Consider triples $(\Si,\vec y,w)$ where
$\Si$ is a {\it prestable holomorphic disc\/} in the sense of
Definition \ref{kh6def5}, $\vec y=(y_1,\ldots,y_l)$, where
$y_1,\ldots,y_l$ are distinct points of $\pd\Si$, none of which are
nodes, and which occur in this cyclic order in $\pd\Si$, and
$w:\Si\ra M$ is a $J$-holomorphic map, with $w(\pd\Si)\subset L$ and
$w_*([\Si])=\be$ in $H_2(M,L;\Z)$. Call two such triples $(\Si,\vec
y,w)$ and $(\Si',\vec y',w')$ {\it isomorphic\/} if there exists a
biholomorphism $f:\Si\ra\Si'$ with $f(y_i)=y_i'$ for $i=1,\ldots,l$
and $w\equiv w'\ci f$. The {\it automorphism group\/} $\Aut(\Si,\vec
y,w)$ is the group of isomorphisms $f:(\Si,\vec y,w)\ra(\Si,\vec
y,w)$. Call a triple $(\Si,\vec y,w)$ {\it stable\/} if
$\Aut(\Si,\vec y,w)$ is finite.

Write $\oM_l^\ma(M,L,J,\be)$ for the moduli space of isomorphism
classes $[\Si,\vec y,w]$ of stable triples $(\Si,\vec y,w)$
satisfying the above conditions. Define {\it boundary evaluation
maps\/} $\ev_i:\oM_l(M,L,J,\be)\ra L$ for $i=1,\ldots,l$ by
$\ev_i:[\Si,\vec y,w]\mapsto w(y_i)$. Write $\M_l^\ma(M,L,J,\be)$
for the subset of $[\Si,\vec y,w]\in\oM_l^\ma(M,L,J,\be)$ with $\Si$
{\it nonsingular}, that is, $\Si$ has no nodes. The `ma' in
$\oM_l^\ma(M,L,J,\be)$ follows Fukaya et al.\ \cite[\S 2]{FOOO}, who
first define moduli spaces $\oM_l(M,L,J,\be)$ in which
$y_1,\ldots,y_l$ need not be in this cyclic order in $\pd\Si$, and
then define $\oM_l^\ma(M,L,J,\be)$ to be the `main component' of
$[\Si,\vec y,w]$ in $\oM_l(M,L,J,\be)$ in which $y_1,\ldots,y_l$ are
cyclically ordered.

There is a natural action of $\Z_l$ on $\oM_l(M,L,J,\be)$, by cyclic
permutation of the boundary marked points~$y_1,\ldots,y_l$.
\label{kh6def6}
\end{dfn}

\begin{rem} Fukaya et al.\ \cite{FOOO} do not consider boundary
nodes of type (c) in Definition \ref{kh6def5}, though they appear to
fall within \cite[Def.~2.16]{FOOO}. Because of this they consider
the moduli spaces $\oM_0(M,L,J,\be)$ when $l=0$ to be possibly {\it
noncompact}. For discussion of this, see remarks preceding
\cite[Fig.~2.4, Th.~2.29 \& Prop.~13.27]{FOOO}, and~\cite[\S
32.1]{FOOO}. They deal with their noncompactness problem by allowing
certain unstable maps in a rather ad hoc way. We prefer our approach
including type (c) nodes, based on Liu~\cite[Def.~3.4]{Liu}.
\label{kh6rem4}
\end{rem}

These moduli spaces have natural topologies called the $C^\iy$ {\it
topology}, due to Gromov \cite{Grom} and defined in \cite[\S
29]{FOOO}. The next result, the analogue of Theorem \ref{kh6thm2},
is proved by Fukaya et al.~\cite[Th.~29.43]{FOOO}, with a slightly
different notion of Kuranishi space. See Remark \ref{kh6rem1} for
what we mean by ``Theorem''.

\begin{thm} In Definition {\rm\ref{kh6def6}} the\/ $\oM_l^\ma
(M,L,J,\be)$ are compact and Hausdorff in the\/ $C^\iy$ topology,
and the\/ $\ev_i$ are continuous.
\label{kh6thm8}
\end{thm}

We define the {\it Maslov index\/} of~$\be\in H_2(M,L;\Z)$.

\begin{dfn} Let $(M,\om)$ be a symplectic manifold of dimension
$2n$, and choose a compatible almost complex structure $J$. Then the
canonical bundle $K_M=\La^{n,0}M$ of $(M,J)$ is a complex line
bundle over $M$, with Hermitian metrics on the fibres $\C$. Let $L$
be an embedded Lagrangian submanifold in $M$. Then there is a
natural isomorphism $K_M\vert_L\cong\La^nT^*L\ot_\R\C$. Choose a
connection $A$ on $K_M$ compatible with the metrics on the fibres,
which respects the isomorphism $K_M\vert_L\cong\La^nT^*L\ot_\R\C$.
Let $F_A$ be the curvature of this connection. Then $A$ is a $\U(1)$
connection, so identifying ${\mathfrak u}(1)\cong\R$, $F_A$ is a
closed 2-form. On $L$, this $A$ reduces to an ${\rm
O}(1)$-connection on $\La^nT^*L$, so $A$ is flat,
and~$F_A\vert_L\equiv 0$.

Thus we may define the {\it relative first Chern class\/} $c_1(M,L)$
of $(M,L)$ to be the relative de Rham cohomology class
$\frac{1}{2\pi}[F_A]$ in $H^2(M,L;\R)$, recalling that $H^l(M,L;\R)$
is defined using closed $l$-forms $\al$ on $M$ with
$\al\vert_L\equiv 0$. One can show that $c_1(M,L)$ lies in the image
of $H^2(M,L;\ha\Z)$ in $H^2(M,L;\R)$, and it is independent of the
choices of $J,A$, and if $L$ is orientable then $c_1(M,L)$ lies in
the image of $H^2(M,L;\Z)$ in~$H^2(M,L;\R)$.

For $\be\in H_2(M,L;\Z)$ we define the {\it Maslov index\/}
$\mu_L(\be)=2c_1(M,L)\cdot\be$ in $\Z$. If $\be$ is represented by
an oriented 2-submanifold (or 2-cycle) $\Si$ in $M$ with
$\pd\Si\subset L$ then $\mu_L(\be)=\frac{1}{\pi}\int_\Si F_A$. If
$L$ is orientable then~$\mu_L(\be)\in 2\Z$.
\label{kh6def7}
\end{dfn}

This is different to the usual definition of Maslov index \cite[\S
2]{FOOO}, but is equivalent to it. Our definition makes it obvious
that $\mu_L(\be)$ depends only on the relative homology class $\be$;
Katz and Liu \cite[Rem.~4.2.2]{KaLi} claim incorrectly that the
Maslov index is a homotopy but not a homology invariant. The next
result is proved by Fukaya et al.~\cite[\S 29]{FOOO}.

\begin{qthm} In Definition {\rm\ref{kh6def6}} we can construct
Kuranishi structures\/ $\ka$ on\/ $\oM_l^\ma(M,L,J,\be)$ with
corners (not g-corners), with
\e
\vdim\bigl(\oM_l^\ma(M,L,J,\be),\ka\bigr)=\mu_L(\be)+l+n-3,
\label{kh6eq15}
\e
where\/ $\dim M=2n$. The maps\/ $\ev_i$ extend to strong submersions
\e
\bev_i:\bigl(\oM_l^\ma(M,L,J,\be),\ka\bigr)\longra L.\\
\label{kh6eq16}
\e
Also $\bev_1\t\cdots\t\bev_l:\bigl(\oM_l^\ma(M,L,J,\be),
\ka\bigr)\ra L^l$ is a strong submersion.
\label{kh6thm9}
\end{qthm}

We can choose $\ka$ to be invariant under the action of $\Z_l$ on
$\oM_l(M,L,J,\be)$ by cyclic permutation of $y_1,\ldots,y_l$, that
is, $\Z_l$ acts on $\oM_l(M,L,J,\be)$ by strong diffeomorphisms,
with the obvious compatibility with $\bev_1,\ldots,\bev_l$. Next we
discuss {\it orientations\/} on $\oM_l^\ma(M,L,J,\be)$. For this we
need the notion of a {\it relative spin structure\/} on $(M,L)$,
\cite[Def.~44.2]{FOOO}:

\begin{dfn} Let $L$ be a submanifold of a manifold $M$. Fix a
triangulation of $M$ such that $L$ is a subcomplex. A {\it relative
spin structure\/} on $(M,L)$ consists of an orientation on $L$; a
class ${\rm st}\in H^2(M;\Z_2)$ such that ${\rm st}\vert_L=w_2(L)\in
H^2(L;\Z_2)$, the second Stiefel--Whitney class of $L$; an oriented
vector bundle $V$ on the 3-skeleton $M_{[3]}$ of $M$ with $w_2(V)=
{\rm st}$; and a spin structure on~$(TL\op V)\vert_{L_{[2]}}$.

Here $L_{[2]}$ is the 2-skeleton of $L$, and as $w_2(V
\vert_{L_{[2]}})={\rm st}\vert_{L_{[2]}}= w_2(L)\vert_{L_{[2]}}$ we
have $w_2\bigl((TL\op V)\vert_{L_{[2]}}\bigr)=0$, so $(TL\op
V)\vert_{L_{[2]}}$ admits a spin structure. If $L$ is spin then
$w_2(L)=0$, so we can take ${\rm st}=0$ and $V=0$ and the spin
structure on $TL\vert_{L_{[2]}}$ to be the restriction of that on
$TL$. Hence, an orientation and spin structure on $L$ induce a
relative spin structure on~$(M,L)$.
\label{kh6def8}
\end{dfn}

Then Fukaya et al.\ \cite[\S 44.1 \& \S 46]{FOOO} prove:

\begin{qthm} A relative spin structure on\/ $(M,L)$ determines
orientations on the Kuranishi spaces\/ $\bigl(\oM_l^\ma(M,L,J,\be),
\ka\bigr)$.
\label{kh6thm10}
\end{qthm}

From now on we suppress the Kuranishi structures $\ka$, and regard
$\oM_l^\ma(M,\ab L,\ab J,\ab\be)$ as a Kuranishi space. Finally we
state formulae for the boundaries of our moduli spaces, in oriented
Kuranishi spaces. Our next theorem can be deduced from results of
Fukaya et al.\ \cite{FOOO}. The closest result to \eq{kh6eq17} is
\cite[Prop.~46.2]{FOOO}, but here in effect the authors omit the sum
over $j$ in \eq{kh6eq17}, and fix $j=1$. The difficult thing to get
right is the orientations (signs). These depend on the order in
which the remaining boundary marked points of the two sides of
\eq{kh6eq17} are identified, which is why we have been careful to
specify this order through the isomorphism between \eq{kh6eq18} and
\eq{kh6eq19}. The correct sign $(-1)^{(l-k)(k-j) +j+n}$ in
\eq{kh6eq17} when $j=1$ is given by Fukaya et al.\ in
\cite[Prop.~46.2 \& Rem.~46.3]{FOOO}, and we can deduce the sign in
the general case by permuting the $l$ marked points.

In equation \eq{kh6eq20}, the second line is in effect the term
$j=k=l=0$ from \eq{kh6eq17}, and the sign $(-1)^n$ is the same as
$(-1)^{(l-k)(k-j)+j+n}$ in \eq{kh6eq17} with $j=k=l=0$. However,
unlike in \eq{kh6eq17}, we must quotient by an action of $\Z_2$. The
point is that at the boundary of $\oM_l^\ma(M,L,J,\be)$, the
prestable holomorphic disc $\Si$ breaks into the union of two
prestable holomorphic discs $\Si_1,\Si_2$ joined by a node. The
marked points $y_1,\ldots,y_l$ in $\pd\Si$ each end up in either
$\pd\Si_1$ or $\pd\Si_2$. In \eq{kh6eq17}, when $l\ge 1$, we
distinguish between $\Si_1$ and $\Si_2$ by decreeing that the marked
point $y_1$ should lie on $\pd\Si_1$; this is why we specify $j\ge
1$ rather than $j\ge 0$ in \eq{kh6eq17}. However in \eq{kh6eq20},
when $l=0$, we cannot use marked points to distinguish $\Si_1$ and
$\Si_2$, so the symmetry of exchanging $\Si_1$ and $\Si_2$ means
that we must divide by~$\Z_2$.

The third line of \eq{kh6eq20} comes from boundary nodes of type (c)
in Definition \ref{kh6def5}. The sign $1$ on the third line is
computed by Fukaya et al.\ in \cite[Prop.~53.6]{FOOO}, though they
include an extra interior marked point. As we explained in Remark
\ref{kh6rem4}, Fukaya et al.\ do not consider nodes of type (c), and
so regard $\oM_0^\ma(M,L,J,\be)$ as possibly noncompact. See also
\cite[Prop.~13.27]{FOOO} and \cite[\S 32.1]{FOOO} for results
analogous to the third line of~\eq{kh6eq20}.

\begin{qthm} In the situation above, choose a relative spin
structure on\/ $(M,L),$ so that\/ $\oM_l^\ma(M,L,J,\be)$ for\/ $l\ge
0$ are compact, oriented Kuranishi spaces with corners by Theorem
{\rm\ref{kh6thm8}} and ``Theorems'' {\rm\ref{kh6thm9}} and\/
{\rm\ref{kh6thm10}}. Then for $l\ge 1$ there is an isomorphism of
oriented Kuranishi spaces:
\e
\begin{split}
&\pd\oM_l^\ma(M,L,J,\be)\cong\\
&\coprod_{1\le j\le k\le l}\,
\coprod_{\begin{subarray}{l}\ga,\de\in H_2(M,L;\Z):\\
\be=\ga+\de\end{subarray}}
\begin{aligned}[t]
&(-1)^{(l-k)(k-j)+j+n}\cdot\oM_{l+j-k+1}^\ma(M,L,J,\ga)\\
&\quad\t_{\bev_{j+1},L,\bev_1}\oM_{k-j+1}^\ma(M,L,J,\de).
\end{aligned}
\end{split}
\label{kh6eq17}
\e
Here on the\/ $j,k,\ga,\de$ term\/ \eq{kh6eq17} identifies the
strong submersions
\ea
\begin{split}
&\bigl(\bev_1\t\bev_1\t\cdots\t\bev_j\t\\
&\bev_{j+1}\t\bev_{j+2}\t\cdots\t\bev_k\t\\
&\bev_{k+1}\t\bev_{k+2}\t\cdots\t\bev_l\bigr)
\vert_{\pd\oM_l^\ma(M,L,J,\be)}:\\
&\pd\oM_l^\ma(M,L,J,\be)\longra L^j \t L^{k-j} \t
L^{l-k}=L^l\qquad\text{and}
\end{split}
\label{kh6eq18}\\
\begin{split}
&(\bev_1\ci\bs\pi_1)\t(\bev_1\ci\bs\pi_1)\t\cdots\t
(\bev_j\ci\bs\pi_1)\t\\
&(\bev_2\ci\bs\pi_2)\t(\bev_3\ci\bs\pi_2)\t\cdots\t
(\bev_{k-j+1}\ci\bs\pi_2)\t\\
&(\bev_{j+2}\ci\bs\pi_1)\t(\bev_{j+3}\ci\bs\pi_1)\t\cdots\t
(\bev_{l+j-k+1}\ci\bs\pi_1):\\
&\oM_{l+j-k+1}^\ma(M,L,J,\ga)\t_{\bev_{j+1},L,\bev_1}
\oM_{k-j+1}^\ma(M,L,J,\de)\\
&\qquad\qquad\qquad\longra L^j \t L^{k-j} \t L^{l-k}=L^l,
\end{split}
\label{kh6eq19}
\ea
where\/ $\bs\pi_1,\bs\pi_2$ project\/ $\oM_{l+j-k+1}^\ma
(M,L,J,\ga)\t_{\bev_{j+1},L,\bev_1}\oM_{k-j+1}^\ma(M,\ab L,\ab
J,\ab\de)$ to $\oM_{l+j-k+1}^\ma(M,L,J,\ga)$ and\/
$\oM_{k-j+1}^\ma(M,L,J,\de)$ respectively.

For\/ $l=0$ there is an isomorphism of oriented Kuranishi spaces:
\e
\begin{split}
&\pd\oM_0^\ma(M,L,J,\be)\cong\\
&\raisebox{-10pt}{\begin{Large}$\displaystyle\biggl[$\end{Large}}
\coprod_{\begin{subarray}{l}\ga,\de\in H_2(M,L;\Z):\; \be=\ga+\de,\\
[\om]\cdot\ga>0,\; [\om]\cdot\de>0\end{subarray}
\!\!\!\!\!\!\!\!\!\!\!\!\!\!\!\!\!\!\!\!\!\!\!\!\!\!\!\!\!\!\!\!\!
\!\!\!} (-1)^n\cdot\oM_1^\ma(M,L,J,\ga)\!\t_{\bev_1,L,\bev_1\!}
\oM_1^\ma(M,L,J,\de)
\raisebox{-10pt}{\begin{Large}$\displaystyle\biggr]\Bigr/$
\end{Large}$\!\!\!\!\!\Z_2$}\\
&\amalg\coprod_{\be'\in H_2(M;\Z):\,\imath(\be')=\be}
\oM_{0,1}(M,J,\be')\t_{\bev_1,M,{\rm inc}}L.
\end{split}
\label{kh6eq20}
\e
Here the\/ $\Z_2$-action exchanges\/ $\ga,\de$ and\/
$\oM_1^\ma(M,L,J,\ga),\oM_1^\ma(M,L,J,\de)$. In the third line,
${\rm inc}:L\ra M$ is the inclusion and\/ $\imath$ comes from the
exact sequence\/ $\smash{H_2(M;\Z)\,{\buildrel\imath\over\longra}
\,H_2(M,L;\Z)\,{\buildrel\pd\over\longra}\,H_1(L;\Z)}$. The third
line is zero unless\/~$\pd\be=0$.
\label{kh6thm11}
\end{qthm}

\subsection{Choosing co-gauge-fixing data for the
${\,\,\overline{\!\!\mathcal M\!}\,}{}_l^{\rm ma}(M,L,J,\be)$}
\label{kh65}

In \S\ref{kh65}--\S\ref{kh67}, we assume the ``Theorems'' of
\S\ref{kh64} throughout. For our attitude to these, see
Remark~\ref{kh6rem1}.

In the situation of \S\ref{kh64}, $\oM_l^\ma(M,L,J,\be)$ is a
compact, oriented Kuranishi space, $\bev_1\t\cdots\t\bev_l:
\oM_l^\ma(M,L,J,\be)\ra L^l$ is a strong submersion, and $L^l$ is
oriented, so $\bigl(\oM_l^\ma(M,L,J,\be),\bev_1\t\cdots\t\bev_l
\bigr)$ is cooriented. Thus by Theorem \ref{kh3thm1} we can choose
co-gauge-fixing data $\bs C_l(M,L,J,\be)$ for $\bigl(\oM_l^\ma
(M,L,J,\be),\ab\bev_1\t\cdots\t\bev_l\bigr)$, so that
$\bigl[\oM_l^\ma(M,L,J,\be),\bev_1\t\cdots\t\bev_l,\bs
C_l(M,L,J,\be)\bigr]$ lies in $KC^*(L^l;R)$, for $R$ a $\Q$-algebra.
We will now show that we can choose the $\bs C_l(M,L,J,\be)$ for all
$l,\be$ so that equations \eq{kh6eq17} and \eq{kh6eq20} of
``Theorem'' \ref{kh6thm11} lift to equations in Kuranishi cochains.

However, thinking carefully about \eq{kh6eq17} and \eq{kh6eq20}, we
see that there is a puzzle. The fibre products in \eq{kh6eq17} and
\eq{kh6eq20} do not quite correspond to any operations we have
defined on co-gauge-fixing data, even in Remark \ref{kh3rem4}(b). In
the $j,k,\ga,\de$ term in \eq{kh6eq17}, it is natural to regard the
term $\oM_{l+j-k+1}^\ma(M,L,J,\ga)$ as giving a cochain in
$KC^*(L^{l+j-k}\t L;R)$, and the term $\oM_{k-j+1}^\ma(M,L,J,\de)$
as giving a cochain in $KC^*(L^{k-j}\t L;R)$. As in Remark
\ref{kh3rem4}(b) we can take the fibre product of these to get a
cochain in $KC^*(L^{l+j-k}\t L^{k-j}\t L;R)= KC^*(L^{l+1};R)$. But
we actually want a cochain in $KC^*(L^l;R)$. We want to forget one
factor of $L$, that is, we want to push forward by $L^{l+1}\ra L^l$.
But pushforwards are defined for Kuranishi chains, not Kuranishi
cochains.

This suggests that we need to mix chains and cochains, so that we
should be working in the {\it bivariant theory\/} of \S\ref{kh48}.
This may be the most natural point of view for Lagrangian Floer
cohomology. In $\oM_l^\ma(M,L,J,\be)$ for $l\ge 1$, we can
distinguish the first marked point $y_1$ from the rest, and regard
$\bigl(\oM_l^\ma(M,L,J,\be),\ab\bev_1\t\cdots\t\bev_l\bigr)$ as
defining a {\it Kuranishi bichain} in $KB^*(\pi_1:L^l\ra L;R)$ where
$\pi_1:(y_1,\ldots,y_l)\mapsto y_1$. Then the fibre products in
\eq{kh6eq17} do correspond to natural operations on Kuranishi
bichains, made out of products, pushforwards, and pullbacks, but
those in \eq{kh6eq20} do not.

Unfortunately, using the bivariant theory in this way is not good
for defining open Gromov--Witten invariants. As in \S\ref{kh65},
there is an action of $\Z_l$ on $\oM_l^\ma(M,L,J,\be)$ by cyclic
permutation of the marked points. For open Gromov--Witten theory we
need to use this cyclic symmetry, but distinguishing the first
marked point $y_1$ destroys it. So we adopt a different solution. We
choose some extra data, in effect gauge-fixing data $\bs G_L$ for
$(L,\id_L)$. Using this we define a {\it pushforward\/}
$KC^*(L^{l+1};R)\ra KC^*(L^l;R)$, which is in effect a cap product
with $(L,\id_L,\bs G_L)$. Then we can define products $*$ of
co-gauge-fixing data and of Kuranishi cochains $KC^*(L^a;R)\t
KC^*(L^b;R)\ra KC^*(L^{a+b-2};R)$ that we will need to make sense of
adding co-gauge-fixing data to \eq{kh6eq17} and~\eq{kh6eq20}.

\begin{dfn} Let $P_m,P,\mu$ be as in Definition \ref{kh3def6}.
Suppose $L$ is an oriented $n$-manifold, not necessarily compact.
Choose an injective map $G_L:L\ra P_m\subset P$ for some $m\gg 0$.
Let $a,b\ge 1$, $X,\ti X$ be compact Kuranishi spaces, $\bs f_1\t
\cdots\t\bs f_a:X\ra L^a$, $\bs{\ti f}_1\t\cdots\t\bs {\ti f}_b:\ti
X\ra L^b$ be strong submersions, $\bs C,\bs{\ti C}$ be
co-gauge-fixing data for $(X,\bs f_1\t\cdots\t\bs f_a)$ and $(\ti X,
\bs{\ti f}_1\t\cdots\t\bs{\ti f}_b)$. Let $1\le i\le a$ and $1\le
j\le b$. Then $X\t_{\bs f_i,L,\bs{\ti f}_j}\ti X$ is a compact
Kuranishi space. Define a strong submersion $\bs g_1\t\cdots\t\bs
g_{a+b-2}:X\t_{\bs f_i,L,\bs{\ti f}_j}\ti X\ra L^{a+b-2}$ to be
\e
\begin{split}
&(\bs f_1\!\ci\!\bs\pi_X)\!\t\!\cdots\!\t\!(\bs
f_{i-1}\!\ci\!\bs\pi_X)\!\t\!(\bs f_{i+1}\!\ci\!\bs\pi_X)\!\t
\!\cdots\!\t\!(\bs f_a\!\ci\!\bs\pi_X)\t\\
&(\bs{\ti f}_1\!\ci\!\bs\pi_{\ti X})\!\t\!\cdots\!\t\!(\bs{\ti
f}_{j-1}\!\ci\!\bs\pi_{\ti X})\!\t\!(\bs{\ti
f}_{j+1}\!\ci\!\bs\pi_{\ti X})\!\t\!\cdots\!\t\!(\bs{\ti
f}_b\!\ci\!\bs\pi_{\ti X}).
\end{split}
\label{kh6eq21}
\e

We will define co-gauge-fixing data $\bs C*_{\bs f_i,L,\bs{\ti
f}_j}\bs{\ti C}$ for $\bigl(X\t_{\bs f_i,L,\bs{\ti f}_j}\ti X,\bs
g_1\t\cdots\t\bs g_{a+b-2}\bigr)$. Remark \ref{kh3rem4}(b) explains
how to modify the construction of \S\ref{kh38} to define
co-gauge-fixing data $\bs C\t_{\bs f_i,L,\bs{\ti f}_j}\bs{\ti C}$
which we can consider as co-gauge-fixing data for $\bigl(X\t_{\bs
f_i,L,\bs{\ti f}_j}\ti X,\bs g_1\t\cdots\t\bs g_{a+b-2}\t\bs\pi_L
\bigr)$, where $\bs\pi_L:X\t_{\bs f_i,L,\bs{\ti f}_j}\ti X\ra L$ is
the projection from the fibre product. Write $\bs C\t_{\bs
f_i,L,\bs{\ti f}_j}\bs{\ti C}$ as $(\bs{\check
I},\bs{\check\eta},\check C^k:k\in\check I)$, where $\bs{\check
I}=\bigl(\check I,(\check V^k,\ldots,\check\psi^k),\check
g_1^k\t\cdots\t\check g_{a+b-2}^k\t\check\pi_L^k:k\in\check
I,\ldots\bigr)$, and $\check C^k$ maps $\check E^k\ra P$ for $k\in
I$.

Define $\bs{\grave I}=\bigl(\check I,(\check V^k,\ldots,\check
\psi^k),\check g_1^k\t\cdots\t\check g_{a+b-2}^k:k\in\check
I,\ldots\bigr)$, and $\bs{\grave\eta}=\bs{\check\eta}$, and for each
$k\in\check I$ define $\grave C^k:\check E^k\ra P$ by $\grave
C^k=\mu\ci\bigl((G_L\ci\check\pi_L^k\ci\check\pi^k)\t \check
C^k\bigr)$, where $\check\pi^k:\check E^k\ra\check V^k$ is the
natural projection. Define $\bs C*_{\bs f_i,L,\bs{\ti f}_j}\bs{\ti
C}=(\bs{\grave I},\bs{\grave\eta},\grave C^k:k\in\check I)$. That
is, $\bs C*_{\bs f_i,L,\bs{\ti f}_j}\bs{\ti C}$ is the same as $\bs
C\t_{\bs f_i,L,\bs{\ti f}_j}\bs{\ti C}$ except that we omit the maps
$\check\pi_L^k$, and multiply the $\check C^k$ by the pullback of
$G_L$ to $\check E^k$ using $\check\pi_L^k\ci\check\pi^k$. This is
similar to the construction of $\bs G_{\bs C}$ in Definition
\ref{kh4def7}, and a similar proof shows that $\bs C*_{\bs
f_i,L,\bs{\ti f}_j}\bs{\ti C}$ is {\it co-gauge-fixing data}
for~$\bigl(X\t_{\bs f_i,L,\bs{\ti f}_j}\ti X,\bs g_1\t\cdots\t\bs
g_{a+b-2}\bigr)$.

Note that using the commutativity and associativity of $\mu$, one
can show that these products $*_{\bs f_i,L,\bs{\ti f}_j}$ have the
commutativity and associativity properties one would expect, in a
similar way to Proposition~\ref{kh3prop6}.

If $R$ is a $\Q$-algebra, $L$ is oriented, and $1\!\le\!i\!\le\! a$,
$1\!\le\!j\!\le\! b$, define $R$-module morphisms
$*_{i,j}:KC^k(L^a;R)\!\t\!
KC^l(L^b;R)\!\ra\!KC^{k+l-n}(L^{a+b-2};R)$~by
\e
\begin{split}
\bigl[X,\bs f_1&\t\cdots\t\bs f_a,\bs C\bigr]*_{i,j}\bigl[\ti X,
\bs{\ti f}_1\t\cdots\t\bs{\ti f}_b,\bs{\ti C}\bigr]=\\
&\bigl[X\t_{\bs f_i,L,\bs{\ti f}_j}\ti X,\bs g_1\t\cdots\t\bs
g_{a+b-2},\bs C*_{\bs f_i,L,\bs{\ti f}_j}\bs{\ti C}\bigr],
\end{split}
\label{kh6eq22}
\e
where we use the orientation on $L$ to convert the product
coorientation for $\bigl(X\t_{\bs f_i,L,\bs{\ti f}_j}\ti X,\bs
g_1\t\cdots\t\bs g_{a+b-2}\t\bs\pi_L\bigr)$ into a coorientation for
$\bigl(X\t_{\bs f_i,L,\bs{\ti f}_j}\ti X,\bs g_1\t\cdots\t\bs
g_{a+b-2}\bigr)$. Equation \eq{kh6eq22} takes relations in
$KC^k(L^a;R),KC^l(L^b;R)$ to relations in $KC^{k+l-n}(L^{a+b-2};R)$,
and so is well-defined. As for cup products in \eq{kh4eq49} we find
that $\d(\ga*_{i,j}\de)=(\d\ga)*_{i,j}\de+(-1)^k\ga*_{i,j}(\d\de)$,
so the $*_{i,j}$ induce morphisms  on Kuranishi cohomology
\e
*_{i,j}:KH^k(L^a;R)\t KH^l(L^b;R)\longra KH^{k+l-n}(L^{a+b-2};R).
\label{kh6eq23}
\e
\label{kh6def9}
\end{dfn}

We can now choose $\bs C_l(M,L,J,\be)$ compatible with \eq{kh6eq17},
\eq{kh6eq20} for all~$l,\be$.

\begin{thm} In the situation of\/ {\rm\S\ref{kh64}} we can choose
co-gauge-fixing data $\bs C_l(M,L,J,\be)$ for $\bigl(\oM_l^\ma
(M,L,J,\be),\ab\bev_1\t\cdots\t\bev_l\bigr)$ for all\/ $l\ge 1$
and\/ $\be$ in $H_2(M,L;\Z)$ such that under {\rm\eq{kh6eq17},} the
restriction $\bs C_l(M,L,\ab J,\ab\be)\vert_{\pd
\oM_l^\ma(M,L,J,\be)}$ is identified on the $j,k,\ga,\de$ term
with\/ $\bs C_{l+j-k+1}(M,L,J,\ga)*_{\bev_{j+1}, L,\bev_1}\ab\bs
C_{k-j+1}\ab(M,\ab L,J,\de),$ using the notation of Definition {\rm
\ref{kh6def9}}. Moreover, we can choose $\bs C_l(M,L,J,\be)$ to be
invariant under the action of\/ $\Z_l$ on $\oM_l^\ma(M,L,J,\be)$ by
cyclic permutation of the marked points.

For the case\/ $l=0,$ choose $\Z_l$-invariant\/ $\bs C_l(M,L,J,\be)$
for all\/ $l\ge 1$ and\/ $\be$ as above, and choose co-gauge-fixing
data $\bs C_{0,1}(M,J,\be')$ for $\bigl(\oM_{0,1}(M,J,\be'),
\bev_1\bigr)$ for all\/ $\be'$ in $H_2(M;\Z)$. Write
$\bs\pi:\oM_0^\ma(M,L,J,\be)\ra \{0\}$ for the strong submersion
projecting to one point\/ $\{0\}$. Then we can choose
co-gauge-fixing data $\bs C_0(M,L,J,\be)$ for
$\bigl(\oM_0^\ma(M,L,J,\be),\bs\pi\bigr)$ for all\/ $\be$ in
$H_2(M,L;\Z)$ such that under {\rm\eq{kh6eq20},} the restriction
$\bs C_0(M,L,J,\be)\vert_{\pd\oM_0^\ma(M,L,J,\be)}$ is identified on
the $\ga,\de$ term with\/ $\bs C_1(M,L,J,\ga)*_{\bev_1, L,\bev_1}\bs
C_1(M,L,J,\de),$ where we take the quotient of this co-gauge-fixing
data by the action of\/ $\Z_2$ as in {\rm\S\ref{kh34},} and on the
$\be'$ term with\/ ${\rm inc}^*\bigl(\bs C_{0,1}(M,J,\be')\bigr)
*_{\bev_1, L,\id_L}[L,\id_L,\bs C_L]$ using pullbacks of
co-gauge-fixing data as in\/ {\rm\S\ref{kh38},} defining
$[L,\id_L,\bs C_L]$ as in {\rm\S\ref{kh47},} and
identifying~$L^0\cong\{0\}$.
\label{kh6thm12}
\end{thm}

\begin{proof} We only need choose $\bs C_l(M,L,J,\be)$ when
$\oM_l^\ma(M,L,J,\be)\ne\es$. Using compactness results from
Geometric Measure Theory, we find that for any $A\ge 0$, there are
only finitely many $\be\in H_2(M,L;\Z)$ with $[\om]\cdot\be\le A$
for which there exist $J$-holomorphic curves in class $\be$. Thus we
can arrange all such $\be$ in an ordered list
$0=\be_0,\be_1,\be_2,\ldots$ such that $m\le n$ implies that
$[\om]\cdot\be_m\le[\om]\cdot\be_n$.

For the first part of the theorem, we will choose $\bs C_l(M,L,J,
\be_m)$ by a double induction on $l,m$. The outer induction is over
$m=0,1,2,\ldots$, and the inductive hypothesis at the $m^{\rm th}$
step is that we have chosen $\bs C_{l'}(M,L,J,\be_{m'})$ satisfying
the conditions for all $1\le l'$ and $0\le m'<m$. The inner
induction, for fixed $m$, is over $l=1,2,\ldots$, and the inductive
hypothesis is that we have chosen $\bs C_{l'}(M,L,J,\be_{m'})$
satisfying the conditions for all $1\le l'$ and $0\le m'<m$, and
also that we have chosen $\bs C_{l'}(M,L,J,\be_m)$ satisfying the
conditions for all $1\le l'<l$. In the inductive step of the inner
induction we must choose $\bs C_l(M,L,J,\be_m)$ satisfying the
conditions.

We suppose these inductive hypotheses, and prove the inner inductive
step. Consider which $j,k,\ga,\de$ terms in \eq{kh6eq17} can be
nonempty, where $\be=\be_m$. For $\oM_{l+j-k+1}^\ma(M,L,
J,\ga)\ne\es$ and $\oM_{k-j+1}^\ma(M,L,J,\de)\ne\es$ we must have
$\ga=\be_{m'}$ and $\de=\be_{m''}$ for some $m',m''$. First suppose
$\ga,\de\ne 0$. Then $[\om]\cdot\ga,[\om]\cdot\de>0$, so
$[\om]\cdot\ga,[\om]\cdot\de<[\om]\cdot\be_m$ as
$[\om]\cdot\ga+[\om]\cdot\de=[\om]\cdot\be_m$, and so $m',m''<m$.
Thus in this case we have already chosen $\bs C_{l+j-k+1}(M,L,J,
\ga)$ and $\bs C_{k-j+1}\ab(M,L,J,\de)$ in a previous outer
inductive step. If $\ga=0$ then $\de=\be_m$, so $m''=m$, and we must
have $l+j-k+1\ge 3$, since $\oM_i^\ma(M,L, J,0)=\es$ unless $i\ge
3$, so $k-j+1<l$. Thus, if $\ga=0$ then $\oM_{k-j+1}^\ma(M,L,J,\de)$
is of the form $\oM_{l'}^\ma(M,L,J,\be_m)$ for $l'<l$, and we have
already chosen $\bs C_{k-j+1}(M,L,J,\de)$ in a previous inner
inductive step. Similarly, if $\de=0$ then we have already chosen
$\bs C_{l+j-k+1}(M,L,J,\ga)$ in a previous inner inductive step.

Thus in all cases, if the $j,k,\ga,\de$ term in \eq{kh6eq17} is
nonempty then we have chosen $\bs C_{l+j-k+1} (M,L,J,\ga)$ and $\bs
C_{k-j+1}(M,L,J,\de)$ in a previous inductive step. So the
conditions in the theorem uniquely determine co-gauge-fixing data
$\bs D_l(M,L,J,\be_m)$ for $\bigl(\pd\oM_l^\ma
(M,L,J,\be_m),(\bev_1\t\cdots\t\bev_l)\vert_{\pd \oM_l^\ma(M,L,
J,\be_m)}\bigr)$, and we must choose $\bs C_l(M,L,J,\be_m)$ with
$\bs C_l(M,L,J,\be_m)\vert_{\pd\oM_l^\ma(M,L,J,\be_m)}\!=\!\bs
D_l(M,\ab L,\ab J,\ab\be_m)$. As in ``Theorem'' \ref{kh6thm9},
$\oM_l^\ma(M,L,J,\be_m)$ is a Kuranishi space {\it with corners, not
g-corners}. So we can apply Theorem \ref{kh3thm3}. We must check
that the restriction $\bs D_l(M,L,J,\be_m)\vert_{\pd^2\oM_l^\ma
(M,L,J,\be_m)}$ is invariant under the natural
involution~$\bs\si:\pd^2\oM_l^\ma(M,L,J,\be_m) \ra\pd^2\oM_l^\ma
(M,L,J,\be_m)$.

Now taking the boundary of \eq{kh6eq17}, using Proposition
\ref{kh2prop2}(a) to rewrite the r.h.s.\ in terms of
$\pd\oM_{l+j-k+1}^\ma(M,L,J,\ga),\pd\oM_{k-j+1}^\ma(M,L,J,\de)$, and
then substituting for these using \eq{kh6eq17}, gives an expression
for $\pd^2\oM_l^\ma (M,L,J,\be_m)$ in terms of double fibre products
of three factors $\oM_i^\ma(M,L,J,\ep)$. By previous inductive
steps, using this we can write $\bs D_l(M,L,J,\be_m)\vert_{\pd^2
\oM_l^\ma(M,L,J,\be_m)}$ as a disjoint union of double products
$*_{\bev_i,L,\bev_j}$ of three factors~$\bs C_i(M,L,J,\ep)$.

When written this way, $\bs\si$-invariance of $\bs D_l(M,L,J,\be_m)
\vert_{\pd^2\oM_l^\ma (M,L,J,\be_m)}$ turns out to follow from
associativity of $*_{\bev_i,L,\bev_j}$ on co-gauge-fixing data in
Definition \ref{kh6def9}. That is, $\bs\si$ acts by exchanging terms
$X_1\t_L(X_2\t_LX_3)$ and $(X_1\t_LX_2)\t_LX_3$, and this must
identify $(\bs C_1*_L\bs C_2)*_L\bs C_3$ and $\bs C_1*_L(\bs
C_2*_L\bs C_3)$, where $X_1,X_2,X_3$ are of the form $\bs
C_i(M,L,J,\ep)$, and $\bs C_1,\bs C_2,\bs C_3$ of the form $\bs
C_i(M,L,J,\ep)$. Therefore $\bs
D_l(M,L,J,\be_m)\vert_{\pd^2\oM_l^\ma(M,L,J,\be_m)}$ is
$\bs\si$-invariant, and we can choose $\bs C_l(M,L,J,\be_m)$
satisfying the conditions of the theorem by Theorem \ref{kh3thm3}.
This completes the inner inductive step. The outer inductive step
the initial steps are trivial, so the first part of the theorem
follows by induction.

To choose $\bs C_l(M,L,J,\be_m)$ to be $\Z_l$-invariant, we use the
$\Ga$-invariance part of Theorem \ref{kh3thm3} with $\Ga=\Z_l$. If
the $\bs C_{l'}(M,L,J,\be_{m'})$ in previous steps have been chosen
$\Z_{l'}$-invariant then $\bs D_l(M,L,J,\be_m)$ turns out to be
$\Z_l$-invariant, with $\Z_l$-action on $\bs
D_l(M,L,J,\be_m)\vert_{\pd^2\oM_l^\ma (M,L,J,\be_m)}$ commuting with
$\bs\si$, so Theorem \ref{kh3thm3} tells us we can choose $\bs
C_l(M,L,J,\be_m)$ to $\Z_l$-invariant.

For the final part, when $l=0$, a similar proof works, though no
induction is needed. The conditions in the theorem uniquely
determine co-gauge-fixing data $\bs D_0(M,L,J,\be)$ for
$\bigl(\pd\oM_0^\ma(M,L,J,\be),\bs\pi\vert_{\pd\oM_0^\ma(M,L,
J,\be)}\bigr)$, and we must choose $\bs C_0(M,L,J,\be)$ with $\bs
C_0(M,L,J,\be)\vert_{\pd\oM_0^\ma(M,L,J,\be)}\!=\!\bs
D_0(M,L,J,\be)$ using Theorem \ref{kh3thm3}. Note that
$\pd^2\oM_0^\ma(M,L,J,\be)$ involves terms in
$\oM_2^\ma(M,L,J,\ep)$, and we need $\bs C_2(M,L,J,\ep)$ to be
$\Z_2$-invariant to make $\bs D_0(M,L,J,\be) \vert_{\pd^2\oM_0^\ma
(M,L,J,\be)}$ $\bs\si$-invariant. So we do need the $\bs
C_l(M,L,J,\be)$ for $l\ge 1$ to be $\Z_l$-invariant.
\end{proof}

We define some notation.

\begin{dfn} In Theorem \ref{kh6thm12}, choose $\Z_l$-invariant
co-gauge-fixing data $\bs C_l(M,L,J,\be)$ for all $l\ge 1$, $\be$,
and $\bs C_{0,1}(M,J,\be')$ for all $\be'$, and $\bs C_0(M,L,J,\be)$
for all $\be$. Let $R$ be a $\Q$-algebra. Define elements
${\mathfrak M}_l(\be)\in KC(L^l;R)$ for all $l\ge 1$ and $\be\in
H_2(M,L;\Z)$, and ${\mathfrak M}_0(\be)\in KC(\{0\};R)$ for all
$\be\in H_2(M,L;\Z)$, and ${\mathfrak N}_0(\be')\in KC(\{0\};R)$ for
all $\be'\in H_2(M;\Z)$, by
\ea
{\mathfrak M}_l(\be)&=\bigl[\oM_l^\ma(M,L,J,\be),\bev_1\t\cdots\t
\bev_l,\bs C_l(M,L,J,\be)\bigr],
\label{kh6eq24}\\
{\mathfrak M}_0(\be)&=\bigl[\oM_0^\ma(M,L,J,\be),\bs\pi,\bs
C_0(M,L,J,\be)\bigr],
\label{kh6eq25}\\
\begin{split}
{\mathfrak N}_0(\be')&=\bigl[\oM_{0,1}(M,J,\be')\t_{\bev_1,M,{\rm
inc}}L,\bs\pi,\\
&\qquad\qquad\qquad {\rm inc}^*\bigl(\bs
C_{0,1}(M,J,\be')\bigr)*_{\bev_1, L,\id_L}[L,\id_L,\bs C_L]\bigr],
\end{split}
\label{kh6eq26}
\ea
where $\bs\pi$ is the projection to $\{0\}$ in both cases. We
identify $L^0=\{0\}$, so that~${\mathfrak M}_0(\be),{\mathfrak
N}_0(\be')\in KC(L^0;R)$.
\label{kh6def10}
\end{dfn}

Combining ``Theorem'' \ref{kh6thm11}, Theorem \ref{kh6thm12} and
Definitions \ref{kh6def9} and \ref{kh6def10} gives:

\begin{cor} In the situation above, ${\mathfrak M}_l(\be)\in
KC(L^l;R)$ is unchanged by cyclic permutation of the $l$ factors
of\/ $L$ in $L^l,$ for all\/ $l\ge 1$ and\/ $\be$. We have
\begin{gather}
\d\bigl({\mathfrak M}_l(\be)\bigr)=\sum_{1\le j\le k\le l}\,
\sum_{\begin{subarray}{l} \ga,\de\in H_2(M,L;\Z):\\
\be=\ga+\de\end{subarray}}
\begin{aligned}[t]
&(-1)^{(l-k)(k-j)+j+n}\cdot\\[-3pt]
&{}\,\,\,{\mathfrak M}_{l+j-k+1}(\ga)*_{j+1,1} {\mathfrak
M}_{k-j+1}(\de),\end{aligned}
\label{kh6eq27}\\
\begin{split}
&\d\bigl({\mathfrak M}_0(\be)\bigr)=\ha\!\!
\sum_{\begin{subarray}{l}\ga,\de\in H_2(M,L;\Z):\; \be=\ga+\de,\\
[\om]\cdot\ga>0,\; [\om]\cdot\de>0\end{subarray}
\!\!\!\!\!\!\!\!\!\!\!\!\!\!\!\!\!\!\!\!\!\!\!\!\!\!\!\!\!\!\!\!\!
\!\!\!} (-1)^n\cdot{\mathfrak M}_1(\ga)*_{1,1}{\mathfrak M}_1(\de)
+\!\!\sum_{\be'\in H_2(M;\Z):\,\imath(\be')=\be
\!\!\!\!\!\!\!\!\!\!\!\!\!\!\!\!\!\!\!\!\!\!} {\mathfrak N}_0(\be'),
\end{split}
\label{kh6eq28}
\end{gather}
for all\/ $l\ge 1$ and\/ $\be\in H_2(M,L;\Z),$ where in the term
${\mathfrak M}_{l+j-k+1} (\ga)*_{j+1,1}{\mathfrak M}_{k-j+1}(\de)$
in {\rm\eq{kh6eq27},} we reorder the factors of\/ $L$ as in
\eq{kh6eq18} and\/ \eq{kh6eq19}. Also $\d\bigl({\mathfrak
N}_0(\be')\bigr)=0$ for all\/ $\be'\in H_2(M;\Z),$ and\/ ${\mathfrak
M}_l(0)=0$ for\/~$l<3$.
\label{kh6cor1}
\end{cor}

The last line holds as $\pd\oM_{0,1}(M,J,\be')=\pd L=\es$ and
$\oM_l^\ma(M,L,J,0)=\es$ for $l<3$. Observe that we have replaced
the quotient by $\Z_2$ in \eq{kh6eq20} by the factor $\ha$ in
\eq{kh6eq28}, using relation Definition \ref{kh4def3}(iv) in
$KC^*(\{0\};R)$. This was an important reason for including
Definition \ref{kh4def3}(iv) amongst the relations for Kuranishi
(co)chains; as discussed in Remarks \ref{kh4rem1}(b) and
\ref{khCrem1}(a), we do not need Definition \ref{kh4def3}(iv) to
define a good theory of Kuranishi (co)homology, and we do not use it
in effective Kuranishi (co)homology.

Here is an important point. As we will see in \S\ref{kh66}, most of
Lagrangian Floer cohomology for one Lagrangian is a {\it purely
algebraic consequence\/} of equation \eq{kh6eq27}, together with an
analogue of \eq{kh6eq27} for families of almost complex structures
$J_t:t\in[0,1]$, to express the dependence of the ${\mathfrak
M}_l(\be)$ on the choice of almost complex structure $J$. Also, as
in \S\ref{kh66}, much of the theory of genus zero open
Gromov--Witten invariants is a purely algebraic consequence of
equations \eq{kh6eq27}--\eq{kh6eq28}, and analogues for families
$J_t:t\in[0,1]$. So, Corollary \ref{kh6cor1} and results of the same
kind are {\it powerful tools for turning geometry into algebra}.

It is a natural question whether any parts of Theorem \ref{kh6thm12}
and Corollary \ref{kh6cor1} have analogues for {\it effective\/}
co-gauge-fixing data $\ubC_l(M,L,J,\be)$, so that we can work in
$KC^*_\ef(L^l;R)$ for $R$ a commutative ring, rather than a
$\Q$-algebra. When $l=0$ the answer is definitely not: in the
$\ga,\de$ term in \eq{kh6eq20}, to make the co-gauge-fixing data
descend to the quotient $\Z_2$, we need fibre products of
co-gauge-fixing data to be commutative, and this does not hold for
effective co-gauge-fixing data. For similar reasons, the part in
Theorem \ref{kh6thm12} on choosing $\bs C_l(M,L,J,\be)$
$\Z_l$-invariant for $l\ge 2$ cannot extend to
effective~$\ubC_l(M,L,J,\be)$.

However, the first part of Theorem \ref{kh6thm12}, on choosing $\bs
C_l(M,L,J,\be)$ for $l\ge 1$, works in the effective case provided
the Kuranishi spaces $\oM_l^\ma(M,L,J,\be)$ are {\it effective}.
This is because the proof relies on the {\it associativity\/} of
products $\bs C*_{\bs f_i,L,\bs{\ti f}_j}\bs{\ti C}$, but not their
(super)commutativity. We can generalize Definition \ref{kh6def9} to
give an associative product $\ubC*_{\bs f_i,L,\bs{\ti f}_j}\ubtC$ of
effective co-gauge-fixing data. We choose an injective map
$\uG_L:L\ra\R^m\subset\uP$ for some $m\gg 0$. In defining the maps
$\ul{\grave C}^k:\check E^k\ra\uP$ we must multiply together three
factors using the associative, non-commutative multiplication $\umu$
on $\uP$. These factors come from $\uC^i$ in $\ubC$, from $\uG_L$,
and from $\utC^{\ti\imath}$ in $\utC$, and we must multiply them in
this order, so that $\uG_L$ is in the middle, to get the
associativity we need.

Now if $[\Si,\vec y,w]\in\oM_l^\ma(M,L,J,\be)$ for $l\ge 1$, then
the only way in which we can have $\Aut(\Si,\vec y,w)\ne\{1\}$ is if
$\Si$ has a nonconstant $\CP^1$ component (or several $\CP^1$
components) joined to the rest of $\Si$ by a node, and
$\Aut(\Si,\vec y,w)$ acts nontrivially on this $\CP^1$, but
trivially on the disc part of $\Si$. If we assume that there are no
nonconstant $J$-holomorphic $\CP^1$'s in $M$, then this cannot
happen, so $\Aut(\Si,\vec y,w)=\{1\}$ for all $[\Si,\vec y,w]$, so
$\oM_l^\ma(M,L,J,\be)$ for $l\ge 1$ has {\it trivial stabilizers},
and is an {\it effective} Kuranishi space. (Note that this does not
work when $l=0$, as we need at least one marked point to prevent the
disc part of $\Si$ from having automorphisms.) Then we can prove:

\begin{thm} In the situation above, suppose there are no
nonconstant\/ $J$-holomorphic {\rm$\CP^1$'}s in $M$. This happens
for instance if there exist no spherical $\be'\in H_2(M;\Z)$ with
$[\om]\cdot\be'>0,$ or if\/ $J$ is generic and for all spherical
$\be'\in H_2(M;\Z)$ with\/ $[\om]\cdot\be'>0$ we have
$c_1(M)\cdot\be'<3-n,$ where\/~$\dim M=2n$.

Extend Definitions {\rm\ref{kh6def9}} and\/ {\rm\ref{kh6def10}} to
effective co-gauge-fixing data in the obvious way. Then we can
choose effective co-gauge-fixing data $\ubC_l (M,L,J,\be)$ for
$\bigl(\oM_l^\ma\ab(M,\ab L,\ab J,\be),\bev_1\t\cdots\t
\bev_l\bigr)$ for all\/ $l\ge 1$ and\/ $\be$ in $H_2(M,L;\Z)$ such
that under {\rm\eq{kh6eq17},} the restriction $\ubC_l(M,L,\ab
J,\ab\be)\vert_{\pd \oM_l^\ma(M,L,J,\be)}$ is identified on the
$j,k,\ga,\de$ term with\/ $\ubC_{l+j-k+1}(M,L,J,\ga)*_{\bev_{j+1},
L,\bev_1}\ubC_{k-j+1}\ab(M,\ab L,J,\de)$. The analogue of\/
\eq{kh6eq27} is
\e
\d\bigl(\ul{{\mathfrak M}\kern -.1em}\kern
.1em{}_l(\be)\bigr)=\sum_{1\le j\le k\le l}\,
\sum_{\begin{subarray}{l} \ga,\de\in H_2(M,L;\Z):\\
\be=\ga+\de\end{subarray}}
\begin{aligned}[t]
&(-1)^{(l-k)(k-j)+j+n}\cdot\\[-3pt]
&{}\,\,\, \ul{{\mathfrak M}\kern -.1em}\kern
.1em{}_{l+j-k+1}(\ga)*_{j+1,1} \ul{{\mathfrak M}\kern -.1em}\kern
.1em{}_{k-j+1}(\de),\end{aligned}
\label{kh6eq29}
\e
which holds in $KC^*_\ec(L^l;R)$ for all\/ $l\ge 1$ and\/ $\be,$
with\/ $R$ a commutative ring.
\label{kh6thm13}
\end{thm}

\begin{rem} In \S\ref{kh631} we sketched a blow-up procedure for
turning an almost complex Kuranishi space $X,(\bs J,\bs K)$ into an
{\it effective\/} almost complex Kuranishi space $\ti X,(\bs{\ti
J},\bs{\ti K})$. This was central to our programme for understanding
integrality properties of Gromov--Witten invariants. We compared our
idea to a scheme of Fukaya and Ono \cite{FuOn2} for defining
$\Z$-valued Gromov--Witten type invariants.

Now in \cite{FuOn2}, Fukaya and Ono also apply the same ideas to
indicate how one could define Lagrangian Floer cohomology over the
integers, for general symplectic manifolds. Parallel to this, one
should be able to adapt our ideas in \S\ref{kh63} to define
Lagrangian Floer cohomology for general symplectic manifolds using
{\it effective} Kuranishi cochains, which would then work over
arbitrary rings $R$, not just $\Q$-algebras. That is, the
Gopakumar--Vafa integrality story for closed Gromov--Witten
invariants should have an open version, in which Lagrangian Floer
cohomology works over $\Z$ not~$\Q$.

We would start by using the blow-up procedure of \S\ref{kh631} to
modify the moduli spaces $\oM_l^\ma(M,L,J,\be)$ to make new, {\it
effective\/} spaces $\ti{\mathcal M}_l^\ma(M,L,J,\be)$, satisfying
the analogue of \eq{kh6eq17}. Replacing $\oM_l^\ma(M,L,J,\be)$ by
$\ti{\mathcal M}_l^\ma(M,L,J,\be)$, we could then prove the analogue
of Theorem \ref{kh6thm13} {\it without\/} assuming there are no
nonconstant $J$-holomorphic $\CP^1$'s in~$M$.

Now the $\oM_l^\ma(M,L,J,\be)$ do not have almost complex
structures, so apparently we cannot use the blow-up procedure on
them. However, as above, all the nontrivial orbifold strata
$\oM_l^\ma(M,L,J,\be)^{\Ga,\rho}$ come from $\CP^1$ bubbles with a
$\Ga$-action. These $\CP^1$ bubbles come lie in moduli spaces
$\oM_{0,1}(M,J,\be')$ which do have almost complex structures $(\bs
J,\bs K)$. From these $(\bs J,\bs K)$ we can construct almost
complex structures on the virtual normal bundle of
$\oM_l^\ma(M,L,J,\be)^{\Ga,\rho}$ in $\oM_l^\ma(M,L,J,\be)$, and
this is enough to perform the blow up procedure.
\label{kh6rem5}
\end{rem}

\subsection{Lagrangian Floer cohomology for one Lagrangian}
\label{kh66}

{\it Lagrangian Floer cohomology} is defined and studied by Fukaya,
Oh, Ohta and Ono in the monumental \cite{FOOO}, which is nearly 1400
pages long. Much of the length and technical complexity of
\cite{FOOO} is due to problems caused by perturbing moduli spaces
$\oM_{k+1}^\ma(M,L,J,\be)$ to choose virtual moduli chains in
singular homology, and ensuring these choices of perturbations and
virtual chains are compatible at the boundaries with choices for
other moduli spaces under the isomorphism~\eq{kh6eq17}.

One difficulty Fukaya et al.\ face is making their singular chains
suitably transverse; a singular chain of positive codimension is
never transverse to itself, so they can never arrange for all the
singular chains they consider to be pairwise transverse. A second
difficulty is that one can in effect make only finitely many choices
of perturbations and virtual chains compatible at once, so Fukaya et
al.\ build their $A_\iy$ algebras through an algebraic limiting
process from a series of finite geometric approximations called
$A_{N,K}$ algebras.

If we develop the theory of Lagrangian Floer cohomology using
Kuranishi cochains instead of singular chains, there is no need to
perturb moduli spaces, and no problems with transverseness. Instead
we choose co-gauge-fixing data, which is a much less disruptive
thing to do, and there is no problem in making infinitely many
choices of co-gauge-fixing data for the $\oM_{k+1}^\ma(M,L,J,\be)$
compatible at their boundaries, as we have shown in Theorem
\ref{kh6thm12}. Because of this, the treatment of \cite{FOOO} can be
hugely simplified and shortened.

Manabu Akaho and the author will define Lagrangian Floer cohomology
using Kuranishi cochains in \cite{AkJo}. Our goal here is only to
whet the reader's appetite, to show that Kuranishi cochains are an
elegant and economical setting for Lagrangian Floer cohomology, and
so to justify all the effort we have put into developing Kuranishi
(co)homology, and perhaps persuade the reader that it is worth
taking the time to learn. So here we will define Lagrangian Floer
cohomology only for one Lagrangian, not two, and we will work with a
fixed almost complex structure $J$ on $M$, without considering the
effect of changing~$J$.

In defining Lagrangian Floer cohomology, we have to consider sums
involving infinitely many terms, coming from $J$-holomorphic discs
of larger and larger area. To ensure these sums converge, we work
over a ring of formal power series known as a {\it Novikov ring}, as
in Fukaya et al.\ \cite[Def.~6.2]{FOOO}.

\begin{dfn} Let $T$ and $e$ be formal variables, graded of degree 0
and 2, respectively. Define two {\it universal Novikov rings} (over
$\Q$) by
\begin{align*}
\La_\nov&=\bigl\{\ts\sum_{i=0}^\iy a_iT^{\la_i}e^{\mu_i}:
\text{$a_i\in\Q$, $\la_i\in\R$, $\mu_i\in\Z$, $\lim_{i\ra\iy}\la_i
=\iy$}\bigr\},\\
\La^0_\nov&=\bigl\{\ts\sum_{i=0}^\iy a_iT^{\la_i}e^{\mu_i}:
\text{$a_i\in\Q$, $\la_i\in[0,\iy)$, $\mu_i\in\Z$,
$\lim_{i\ra\iy}\la_i =\iy$}\bigr\}.
\end{align*}
Then $\La_\nov^0\subset\La_\nov$ are $\Q$-vector spaces. Define
multiplications `$\,\cdot\,$' on $\La_\nov,\La_\nov^0$ by
$\bigl(\sum_{i=0}^\iy a_iT^{\la_i}e^{\mu_i}\bigr)\cdot
\bigl(\sum_{j=0}^\iy b_jT^{\nu_j}e^{\xi_j}\bigr)=\sum_{i,j=0}^\iy
a_ib_jT^{\la_i+\nu_j}e^{\mu_i+\xi_j}$. Then $\La_\nov$, $\La_\nov^0$
are {\it commutative $\Q$-algebras} with identity~$1=1T^0e^0$.

Write $\La_\nov^*$ for either $\La_\nov$ or $\La_\nov^0$. Define
{\it filtrations\/} of $\La_\nov^*$ by
\begin{align*}
F^\la\La_\nov^*&=\bigl\{\ts\sum_{i=0}^\iy
a_iT^{\la_i}e^{\mu_i}\in\La_\nov^*:\text{$\la_i\ge\la$ for all
$i=0,1,\ldots$}\bigr\},
\end{align*}
for $\la\in\R$. These filtrations induce {\it topologies} on
$\La_\nov^*$, and notions of {\it convergence\/} for sequences and
series. An infinite sum $\sum_{k=0}^\iy\al_k$ in $\La_\nov^*$
converges in $\La_\nov^*$ if and only if for all $\la\in\R$ we have
$\al_k\in F^\la\La_\nov^*$ for all except finitely many~$k$.

As $T,e$ are graded of degrees 0,2, we can regard
$\La_\nov,\La_\nov^0$ as {\it graded\/} rings. Write
$\La^{(k)}_\nov,\La_\nov^{0\,(k)}$ for the degree $k$ parts of
$\La_\nov,\La_\nov^0$, for $k\in\Z$. Then
\e
\La_\nov^{(2k)}=\bigl\{\ts\sum\limits_{i=0}^\iy a_iT^{\la_i}e^k:
\text{$a_i\in\Q$, $\la_i\in\R$,
$\lim\limits_{i\ra\iy}\la_i=\iy$}\bigr\}, \quad \La^{(2k+1)}_\nov=0.
\label{kh6eq30}
\e
We can also consider {\it modules\/} over $\La_\nov^*$. Let
$V=\bigop_{k\in\Z}V^k$ be a graded $\Q$-vector space. Then
$V\ot\La_\nov^*$ is a $\La_\nov^*$-module, with grading
$(V\ot\La_\nov^*)^l=\bigop_{j+k=l}V^j\ot\La_\nov^{*\,(k)}$, and
filtration $F^\la(V\ot\La^*_\nov)=V\ot F^\la\La^*_\nov$ for
$\la\in\R$. Write $V\hat\ot\La_\nov^*$ for the {\it completion} of
$V\ot\La_\nov^*$ with respect to this filtration. Then
$V\hat\ot\La_\nov^*$ has a topology, and an infinite sum
$\sum_{k=0}^\iy v_k$ converges in $V\hat\ot\La_\nov^*$ if for all
$\la\in\R$ we have $v_k\in F^\la\bigl(V\hat\ot\La_\nov^*\bigr)$ for
all except finitely many~$k$.

If $L$ is a manifold $KC^*(L;\La_\nov^*)=KC^*(L;\Q)\ot \La_\nov^*$.
Write $\widehat{KC}{}^*(L;\La_\nov^*)$ for the completion
$KC^*(L;\Q)\hat\ot\La_\nov^*$. We consider the grading on
$\widehat{KC}{}^*(L;\La_\nov^*)$ to be the combination of those on
$KC^*(L;\Q)$ and~$\La_\nov^*$.
\label{kh6def11}
\end{dfn}

Next we define {\it gapped filtered\/ $A_\iy$ algebras}, following
Fukaya et al.~\cite[\S 7.2]{FOOO}.

\begin{dfn} A {\it gapped filtered\/ $A_\iy$ algebra}
$(A\ot\La_\nov^*,\m)$ consists of:
\begin{itemize}
\setlength{\itemsep}{0pt}
\setlength{\parsep}{0pt}
\item[(a)] A graded $\Q$-vector space $A=\bigop_{d\in\Z}A^d$;
\item[(b)] Graded $\La_\nov^*$-multilinear maps
$\smash{\m_k:{\buildrel {\!\ulcorner\,\,\,\text{$k$ copies }
\,\,\,\urcorner\!} \over
{\!\vphantom{m}\smash{(A\hat\ot\La_\nov^*)\!
\t\cdots\t\!(A\hat\ot\La_\nov^*)}}}\!\ra}A\hat\ot\La_\nov^*$ for
$k=0,1,2,\ldots$, of degree $+1$. Write~$\m=(\m_k)_{k\ge 0}$.
\end{itemize}
These must satisfy the following conditions:
\begin{itemize}
\setlength{\itemsep}{0pt}
\setlength{\parsep}{0pt}
\item[(i)] there exists a subset ${\cal G}\subset[0,\iy)\t\Z$, closed
under addition, such that ${\cal G}\cap(\{0\}\t\Z)=\{(0,0)\}$ and
${\cal G}\cap([0,C]\t\Z)$ is finite for any $C\ge 0$, and the maps
$\m_k$ for $k\ge 0$ may be written
$\smash{\m_k=\sum_{(\la,\mu)\in{\cal G}}T^\la
e^\mu\m_k^{\smash{\la,\mu}}}$, for unique $\Q$-multilinear maps
$\m_k^{\smash{\la,\mu}}:{\buildrel {\ulcorner\,\,\,\text{$k$ copies
} \,\,\,\urcorner} \over {\vphantom{m}\smash{A\t\cdots\t A}}}\ra A$
graded of degree $1-2\mu$. When $k=0$, we take
$\m_0\in(A\ot\La_\nov^0)^{(1)}$ and $\m_0^{\smash{\la,\mu}}\in
A^{1-2\mu}$;
\item[(ii)] $\m_0^{\smash{0,0}}=0$, in the notation of (i); and
\item[(iii)] call $a\in A\ot\La_\nov^0$ {\it pure} if $a\in
(A\ot\La_\nov^0)^{(d)}\sm\{0\}$ for some $d\in\Z$, and then define
the {\it degree} of $a$ to be $\deg a=d$. Then we require that for
all $k\ge 0$ and all pure $a_1,\ldots,a_k$ in $A\ot\La_\nov^0$, we
have
\end{itemize}
\begin{equation}
\begin{gathered}
\sum_{\begin{subarray}{l}i,k_1,k_2:1\le i\le k_1,\\
k_2\ge 0,\; k_1+k_2=k+1\end{subarray}}
\begin{aligned}[t]
(-1)^{\sum_{l=1}^{i-1}\deg a_l}\m_{k_1}(a_1,\ldots,a_{i-1},\m_{k_2}
(a_i,\ldots,a_{i+k_2-1}), \\
a_{i+k_2}\ldots,a_k)=0.
\end{aligned}
\end{gathered}
\label{kh6eq31}
\end{equation}
Equation \eq{kh6eq31} may be rewritten in terms of the
$\m_k^{\smash{\la,\mu}}$ as follows: for all $k\ge 0$, all
$(\la,\mu)\in{\cal G}$ and all pure $a_1,\ldots,a_k$ in $A$, we have
\begin{equation}
\begin{gathered}
\sum_{\begin{subarray}{l}i,k_1,k_2,\la_1,\la_2,\mu_1,\mu_2:\;
1\le i\le k_1, k_2\ge 0,\\
 k_1+k_2=k+1,\;\la_1+\la_2=\la,\; \mu_1+\mu_2=\mu\end{subarray}
\!\!\!\!\!\!\!\!\!\!\!\!\!\!\!\!\!\!\!\!\!\!\!\!\!\!\!\!\!\!\!
\!\!\! } \!\!\!\!\!\!\!\!\!\!\!\!\!\!\!\!\!\!\!
\begin{aligned}[t]
(-1)^{\sum_{l=1}^{i-1}\deg
a_l}\m_{k_1}^{\la_1,\mu_1}(a_1,\ldots,a_{i-1},\m_{k_2}^{\la_2,\mu_2}
(a_i,\ldots,a_{i+k_2-1}), \\
a_{i+k_2}\ldots,a_k)=0.
\end{aligned}
\end{gathered}
\label{kh6eq32}
\end{equation}

A gapped filtered $A_\iy$ algebra $(A\hat\ot\La_\nov^*,\m)$ is
called {\it strict\/} if $\m_0=0$. Then (iii) implies that
$\m_1\ci\m_1=0$, so $(A\hat\ot\La_\nov^*,\m_1)$ is a complex of
$\La_\nov^*$-modules, and we can form its {\it cohomology}
$H^*(A\hat\ot\La_\nov^*,\m_1)$.
\label{kh6def12}
\end{dfn}

To define Lagrangian Floer cohomology we will need {\it strict\/}
gapped filtered $A_\iy$ algebras. {\it Bounding cochains} are a
method of modifying gapped filtered $A_\iy$ algebras to make them
strict, introduced by Fukaya et al.~\cite[\S 5.7, \S 11]{FOOO}.

\begin{dfn} Let $(A\hat\ot\La_\nov^*,\m)$ be a gapped filtered
$A_\iy$ algebra, and suppose $b\in F^\la(A\hat\ot\La_\nov^*)^{(0)}$
for some $\la>0$. Define graded $\La_\nov^*$-multilinear maps
$\m_k^b:{\buildrel {\ulcorner\,\,\,\text{$k$ copies }
\,\,\,\urcorner} \over {\vphantom{m}\smash{(A\hat\ot\La_\nov^*)
\!\t\!\cdots\!\t\!(A\hat\ot\La_\nov^*)}}}\!\ra\! A\hat\ot\La_\nov^*$
for $k=0,1,\ldots$, of degree $+1$,~by
\begin{gather*}
\m_k^b(a_1,\ldots,a_k)=\!\!\sum_{n_0,\ldots,n_k\ge
0}\!\!\!\m_{k+n_0+\cdots+n_k}\begin{aligned}[t]
\bigl({\buildrel{\ulcorner\,\,\text{$n_0$} \,\,\urcorner}
\over{\vphantom{m}\smash{b,\ldots,b}}},a_1, {\buildrel
{\ulcorner\,\,\text{$n_1$}\,\,\urcorner}
\over{\vphantom{m}\smash{b,\ldots,b}}},a_2,{\buildrel
{\ulcorner\,\,\text{$n_2$}\,\,\urcorner}
\over{\vphantom{m}\smash{b,\ldots,b}}}&,\\[-6pt]
\ldots, {\buildrel {\ulcorner\,\,\text{$n_{k-1}$}\,\,\urcorner}
\over{\vphantom{m}\smash{b,\ldots,b}}},a_k,{\buildrel
{\ulcorner\,\,\text{$n_k$} \,\,\urcorner}
\over{\vphantom{m}\smash{b,\ldots,b}}}&\bigr).
\end{aligned}
\end{gather*}
This converges as $b\in F^\la(A\hat\ot\La_\nov^*)$ for $\la>0$.
Write $\m^b=(\m_k^b)_{k\ge 0}$. We call $b$ a {\it bounding
cochain\/} for $(A\ot\La_\nov^0,\m)$ if $\m_0^b=0$, that is, if
\e
\ts\sum_{k\ge 0}\m_k(b,\ldots,b)=0.
\label{kh6eq33}
\e
This is the {\it Maurer--Cartan equation}, or {\it
Batalin--Vilkovisky master equation}.
\label{kh6def13}
\end{dfn}

It is then easy to prove \cite[Prop.~11.10]{FOOO}:

\begin{lem} In Definition {\rm\ref{kh6def13},} $(A\hat\ot
\La_\nov^*,\m^b)$ is a gapped filtered\/ $A_\iy$ algebra, which is
strict if and only if\/ $b$ is a bounding cochain.
\label{kh6lem}
\end{lem}

Thus, if $b$ is a bounding cochain then $(A\ot\La_\nov^0,\m_1^b)$ is
a complex, and we may form its cohomology $H^*(A\hat\ot\La_\nov^*,
\m_1^b)$. We will now show that in the situation of
\S\ref{kh64}--\S\ref{kh65}, we can define a {\it gapped filtered\/
$A_\iy$ algebra structure} on the shifted Kuranishi cochains
$\widehat{KC}{}^*(L;\La_\nov^*)[1]$ of~$L$.

\begin{dfn} Suppose $(M,\om)$ is a compact symplectic manifold, $L$
a compact, embedded Lagrangian submanifold in $M$, and $J$ an almost
complex structure on $M$ compatible with $\om$. Fix a relative spin
structure on $(M,L)$, so that $\oM_l^\ma(M,L,J,\be)$ for $l\ge 1$
are compact, oriented Kuranishi spaces with corners by Theorem
\ref{kh6thm8} and ``Theorems'' \ref{kh6thm9} and \ref{kh6thm10}.
Note that the relative spin structure includes an orientation on
$L$. Choose co-gauge-fixing data $\bs C_l(M,L,J,\be)$ for
$\bigl(\oM_l^\ma (M,L,J,\be),\bev_1\t\cdots\t\bev_l\bigr)$ for all
$l\ge 1$ and $\be$ in $H_2(M,L;\Z)$ by the first part of Theorem
\ref{kh6thm12}. We do not require $\bs C_l(M,L,J,\be)$ to be
$\Z_l$-invariant, and we do not choose $\bs C_0(M,L,J,\be)$. Define
${\mathfrak M}_l(\be)\in KC(L^l;R)$ for all $l\ge 1$ and $\be\in
H_2(M,L;\Z)$ as in Definition \ref{kh6def10}. Then Corollary
\ref{kh6cor1} shows that equation \eq{kh6eq27} holds for all $l\ge
1$ and~$\be$.

Write $KC^*(L;\Q)[1]$ for the Kuranishi cochains {\it with grading
shifted by one}, that is, $\bigl(KC^*(L;\Q)[1]\bigr){}^k=
KC^{k+1}(L;\Q)$, and similarly for $\widehat{KC}{}^*(L;
\La_\nov^*)[1]$. We will define a {\it gapped filtered\/ $A_\iy$
algebra} $\bigl(\widehat{KC}{}^*(L;\La_\nov^*)[1],\m\bigr)$ on the
shifted cochains $\widehat{KC}{}^*(L;\La_\nov^*)[1]
\!=\!KC^*(L;\Q)[1]\hat\ot \La_\nov^*$. Define ${\cal
G}\!\subset\![0,\iy)\!\t\!\Z$~by
\begin{align*}
{\cal G}=\bigl\{(\la,\mu)\in[0,\iy)&\t\Z:\text{there exist
$\be_1,\ldots,\be_k\in H_2(M,L;\Z)$ with}\\
&\text{$\oM_3^\ma(M,L,J,\be_i)\ne\es$ for $i=1,\ldots,k$, and}\\
&\text{$\la=[\om]\cdot(\be_1+\cdots+\be_k)$,
$\mu=\ha\mu_L(\be_1+\cdots+\be_k)$}\bigr\}.
\end{align*}

Define $\m_k^{\smash{\la,\mu}}:{\buildrel {\ulcorner\,\,\,\text{$k$
copies } \,\,\,\urcorner} \over
{\vphantom{m}\smash{KC^*(L;\Q)[1]\t\cdots\t KC^*(L;\Q)[1]}}}\ra
KC^*(L;\Q)[1]$ for all $k\ge 0$ and $(\la,\mu)\in{\cal G}$  by
\e
\begin{split}
&\m_1^{\smash{0,0}}(f_1)=(-1)^n\d f_1,\qquad\qquad
\m_k^{\la,\mu}(f_1,\ldots,f_k)=\\
&\sum_{\begin{subarray}{l}
\be\in H_2(M,L;\Z):\\
[\om]\cdot\be=\la,\;\> \mu_L(\be)=2\mu,\\
\oM_{k+1}^\ma(M,L,J,\be)\ne\es\end{subarray}
\!\!\!\!\!\!\!\!\!\!\!\!\!\!\!\!\!\!\!\!\!\!\!\!\!\!\!\!\!\!\!\!\!\!}
\!\!\!
\begin{aligned}[t]
(\cdots(({\mathfrak M}_{k+1}(\be)*_{2,L,1}f_1)*_{3,L,1}f_2)\cdots)
*_{k+1,L,1}f_k)&,\\
(k,\la,\mu)\ne(1,0,0)&.
\end{aligned}
\end{split}
\label{kh6eq34}
\e
Following the definitions through shows that for
$(k,\la,\mu)\ne(1,0,0)$ on generators $[X_i,\bs f_i,\bs C_i]$ of
$KC^*(L;\Q)$ this gives
\e
\begin{split}
&\m_k^{\la,\mu}\bigl([X_1,\bs f_1,\bs C_1],\ldots,[X_k,\bs f_k,\bs
C_k]\bigr)=\\
&\sum_{\begin{subarray}{l}
\be\in H_2(M,L;\Z):\\
[\om]\cdot\be=\la,\;\> \mu_L(\be)=2\mu,\\
\oM_{k+1}^\ma(M,L,J,\be)\ne\es\end{subarray}
\!\!\!\!\!\!\!\!\!\!\!\!\!\!\!\!\!\!\!\!\!\!\!\!\!\!\!\!\!\!\!\!\!\!}
\begin{aligned}[t]
\bigl[\oM_{k+1}^\ma(M,L,J,\be)\t_{\bev_2\t\cdots\bev_{k+1},L^k,\bs
f_1\t\cdots\t\bs f_k}(X_1\t\cdots\t X_k)&,\\
\bev_1\ci\bs\pi_{\oM_{k+1}^\ma(M,L,J,\be)},\bs C_{\bs
C_l(M,L,J,\be),\bs C_1,\ldots,\bs C_k}\bigr]&,
\end{aligned}
\end{split}
\label{kh6eq35}
\e
where $\bs C_{\bs C_l(M,L,J,\be),\bs C_1,\ldots,\bs C_k}$ is
co-gauge-fixing data for $\bigl(\oM_{k+1}^\ma(M,L,J,\be)
\t_{L^k}(X_1\t\cdots\t X_k),\bev_1\ci\bs\pi_{\oM_{k+1}^\ma(M,L,
J,\be)}\bigr)$ obtained as a fibre product of the co-gauge-fixing
data $\bs C_l(M,L,J,\be),\bs C_1,\ldots,\bs C_k$ on the $k+1$
factors by $k$ applications of the product $*_{\bs f_i,L,\bs{\ti
f}_j}$ of Definition \ref{kh6def9}. Suppose $[X_i,\bs f_i,\bs C_i]$
has degree $a_i$ in $KC^*(L;\Q)[1]$, so that $[X_i,\bs f_i,\bs
C_i]\in KC^{a_i+1}(L;\Q)$, and $\vdim X_i=n-a_i-1$, where $n=\dim
L$. Then using \eq{kh6eq15} we find that the fibre product in
\eq{kh6eq35} has dimension $2\mu+n-2-a_1-\cdots -a_k$, so
\eq{kh6eq35} has degree $1-2\mu+a_1+\cdots+a_k$ in $KC^*(L;\Q)[1]$,
and $\m_k^{\la,\mu}$ is graded of degree~$1-2\mu$.

Now define $\smash{\m_k:{\buildrel {\!\ulcorner\,\,\,\text{$k$
copies } \,\,\,\urcorner\!} \over
{\!\vphantom{m}\smash{\widehat{KC}{}^*(L; \La_\nov^*)[1]\!
\t\cdots\t\!\widehat{KC}{}^*(L;
\La_\nov^*)[1]}}}\!\ra}\widehat{KC}{}^*(L; \La_\nov^*)[1]$ for
$k=0,1,2,\ldots$ by $\smash{\m_k=\sum_{(\la,\mu)\in{\cal G}}T^\la
e^\mu\m_k^{\smash{\la,\mu}}}$. As the $\m_k^{\la,\mu}$ are graded of
degree $1-2\mu$, $\m_k$ is graded of degree $+1$. Using \eq{kh6eq27}
and \eq{kh6eq34} one can show that the $\m_k^{\la,\mu}$ satisfy
\eq{kh6eq32}, and thus the $\m_k$ satisfy \eq{kh6eq31}. Write
$\m=(\m_k)_{k\ge 0}$. Then $\bigl(\widehat{KC}{}^*(L;
\La_\nov^*)[1],\m\bigr)$ is a {\it gapped filtered\/ $A_\iy$
algebra}.

Let $b$ be a {\it bounding cochain} for $\bigl(\widehat{KC}{}^*(L;
\La_\nov^*)[1],\m\bigr)$, and define $\m^b$ as in Definition
\ref{kh6def13}. Then $\m^b_1:\widehat{KC}{}^*(L;\La_\nov^*)\ra
\widehat{KC}{}^*(L;\La_\nov^*)$ is graded of degree $+1$ with $
\m^b_1\ci\m^b_1=0$. Define the {\it Lagrangian Floer cohomology\/}
$HF^*(L,b;\La_\nov^*)$ to be the cohomology of $\bigl(
\widehat{KC}{}^*(L;\La_\nov^*),\m^b_1\bigr)$. Note that we drop the
shift in gradings, so that the grading of $HF^*(L,b;\La_\nov^*)$
agrees with that of ordinary cohomology. We can think of
$HF^*(L,b;\La_\nov^*)$ as a `quantum deformation' of Kuranishi
cohomology~$KH^*(L;\La_\nov^*)$.
\label{kh6def14}
\end{dfn}

The important thing is that $HF^*(L,b;\La_\nov^*)$ is {\it
independent of all the arbitrary choices involved up to canonical
isomorphism}, including the choice of almost complex structure $J$.
(It does depend on the bounding cochain $b$, so to make sense of
this we must introduce a notion of equivalence between bounding
cochains $b$ for $J$, and bounding cochains $b'$ for another almost
complex structure $J'$.) Thus, Lagrangian Floer cohomology
$HF^*(L,b;\La_\nov^*)$ basically depends only on $(M,\om)$ and $L$,
as for Gromov--Witten invariants in Theorem \ref{kh6thm6}. Also, the
gapped filtered $A_\iy$ algebra above is independent of choices up
to $A_\iy$ quasi-isomorphism.

Now compare Definition \ref{kh6def14} with the construction of
gapped filtered $A_\iy$ algebras for Lagrangian submanifolds in
Fukaya et al.\ \cite[Chapter 3]{FOOO}. Briefly, they construct a
countably generated subcomplex $\mathcal X$ of the complex of smooth
singular chains $C^\rsi_*(L;\La_\nov^*)$ on $L$ with values in
$\La_\nov^*$, and multilinear operators $\m_k:{\mathcal
X}^{\ot^k}\ra{\mathcal X}$ satisfying the $A_\infty$ relations
\eq{kh6eq31}. Roughly speaking, $\m_k$ is a sum over all $\be$ of
intersections with virtual cycles for moduli
spaces~$\oM_{k+1}^\ma(M,L,J,\be)$.

The details of the construction, however, are horrible. For each
Kuranishi space $\oM_{k+1}^\ma(M,L,J,\be)$ we must choose a virtual
moduli chain, an abstract perturbation, which must be compatible at
the boundary $\pd\oM_{k+1}^\ma(M,L,J,\be)$ with the choices of
virtual chains for other $\oM_{l+1}^\ma(M,L,J,\ga)$ under
\eq{kh6eq17}. Worse, one virtual moduli chain for
$\oM_{k+1}^\ma(M,L,J,\be)$ is not enough, because of transversality
problems: we must choose one for each $k$-tuple $(P_1,\ldots,P_k)$
in ${\mathcal X}^k$. Worse still, it is not actually possible to
make all these choices consistently for all $k\ge 0$ and
$(P_1,\ldots,P_k)$, but only for finitely generated $\mathcal X$ and
bounded~$k$.

Therefore Fukaya et al.\ construct $A_{N,K}$-algebras ${\mathcal
X}_{N,K}$ for increasingly large $N,K$, prove homotopy equivalences
between them, and finally use an abstract algebraic construction to
prove existence of the (now non-geometric) $A_\iy$-algebra they
want. With Kuranishi cohomology we can avoid all this, and go
directly to a geometric $A_\iy$-algebra structure on the full
Kuranishi cochains~$\widehat{KC}{}^*(L;\La_\nov^*)$.

If there are no nonconstant $J$-holomorphic $\CP^1$'s in $M$ then we
can use Theorem \ref{kh6thm13} to define Lagrangian Floer cohomology
using {\it effective\/} Kuranishi cochains, and the construction
then works with Novikov rings {\it over the integers\/} $\Z$, or
over an arbitrary commutative ring $R$, rather than just over~$\Q$.

These ideas are also relevant to the construction of the {\it Fukaya
category\/} ${\rm Fuk}(M,\om)$ of a general symplectic manifold
$(M,\om)$, as in \cite{Fuka1,FOOO,Seid1,Seid3}. The {\it derived
Fukaya category\/} $D^b{\rm Fuk}(M,\om)$ is an essential part of
Kontsevich's Homological Mirror Symmetry programme \cite{Kont}, and
so is a very important piece of mathematics. Roughly speaking, the
objects of ${\rm Fuk}(M,\om)$ should be pairs $(L,b)$ of a
Lagrangian $L$ and a bounding cochain $b$ for $L$, and the morphisms
between such pairs in $D^b{\rm Fuk}(M,\om)$ are given by Lagrangian
Floer cohomology. The author believes that by using Kuranishi
cochains as above, it should be possible to write down a reasonably
clean and simple model for~${\rm Fuk}(M,\om)$.

\subsection{Open Gromov--Witten invariants}
\label{kh67}

As in \S\ref{kh62}, (closed) Gromov--Witten invariants `count'
$J$-holomorphic curves without boundary in a symplectic manifold
$(M,\om)$, and are independent of the choice of almost complex
structure $J$. In a similar way, if $L$ is Lagrangian submanifold of
$(M,\om)$ we can attempt to define {\it open Gromov--Witten
invariants} which `count' $J$-holomorphic curves in $M$ with
boundary in $L$, and hope to make them independent of $J$. This is a
more difficult problem than the closed case, since as in
\S\ref{kh61} and \S\ref{kh64} the moduli spaces have no boundary in
the closed case, but have boundary and corners in the open case.

Quite a lot of work has already been done on open Gromov--Witten
invariants. In String Theory, it is expected that open
Gromov--Witten invariants should exist and play a r\^ole, and String
Theorists can even compute them in examples; see for example Ooguri
and Vafa \cite{OoVa}, Aganagic and Vafa \cite{AgVa}, Aganagic, Klemm
and Vafa \cite{AKV}, and Graber and Zaslow \cite{GrZa}, and others.
Some mathematics papers on open Gromov--Witten invariants are Katz
and Liu \cite{KaLi}, Liu \cite{Liu}, Li and Song \cite{LiSo},
Welschinger \cite{Wels1,Wels2}, Cho \cite{Cho1,Cho2}, Solomon
\cite{Solo}, and Pandharipande, Solomon and Walcher~\cite{PSW}.

Despite these, there does not yet seem to be a satisfactory
definition of open Gromov--Witten invariants (by `satisfactory' I
really mean that the invariants should be defined using any almost
complex structure $J$ on $M$ compatible with $\om$, and give an
answer independent of $J$). It is clear from these references (see
Aganagic et al.\ \cite[\S 5]{AKV}, and computations by Cho
\cite{Cho2}, for example), that just fixing $(M,\om),L$ is not
enough to get open Gromov--Witten invariants independent of $J$, one
needs some extra assumptions or extra data.

Liu \cite{Liu} defines open Gromov--Witten invariants for the case
in which there is a $\U(1)$-action on $(M,\om),L$, and $J$ is
$\U(1)$-invariant. Welschinger \cite{Wels1,Wels2}, Cho \cite{Cho1}
and Solomon \cite{Solo} define open Gromov--Witten invariants in the
case in which $\dim M$ is 4 or 6, and $L$ is the fixed point set of
an involution $\si:M\ra M$, which is antiholomorphic with respect to
$J$; these invariants are calculated when $M$ is the quintic
Calabi--Yau 3-fold by Pandharipande et al.~\cite{PSW}.

We now present a part of the author's solution to the problem of
defining open Gromov--Witten invariants; the complete picture will
be explained in \cite{Joyc2}. For genus zero open Gromov--Witten
invariants, the key idea is this: we need to fix some extra data to
define open Gromov--Witten invariants independent of $J$, and {\it
the right kind of extra data is a choice of bounding cochain $b$ for
Lagrangian Floer cohomology of} $L$, in the sense of \S\ref{kh66}.
The open Gromov--Witten invariants will then depend on $(L,b)$, and
they will be independent of the choice of $J$ in the same sense that
Lagrangian Floer cohomology is. To correct for the
$\oM_{0,1}(M,J,\be')$ terms in \eq{kh6eq20} we need to assume
$[L]=0$ in $H_n(M;\Q)$, and choose an $(n+1)$-chain $B$ in $M$ with
$\pd B=L$. In \cite{Joyc2} we will also explain how to define
higher-genus open Gromov--Witten invariants under stronger
assumptions, including requiring $L$ to be a $\Q$-homology sphere.
\smallskip

\noindent{\bf Warning:} {\it most of the signs below are probably
wrong, but I'll fix them in the final version.}
\smallskip

We work in the situation of \S\ref{kh64}--\S\ref{kh66}. Let $b$ be a
bounding cochain for $\bigl(\widehat{KC}{}^*(L;
\La_\nov^0)[1],\m\bigr)$. Then $b$ satisfies \eq{kh6eq33}, which by
\eq{kh6eq34} means that
\e
\begin{split}
&(-1)^n\d b+\\
&\sum_{k\ge 0,\; \be\in H_2(M,L;\Z):\;
\oM_{k+1}^\ma(M,L,J,\be)\ne\es
\!\!\!\!\!\!\!\!\!\!\!\!\!\!\!\!\!\!\!\!\!\!\!\!\!\!\!\!\!\!\!\!
\!\!\!\!\!\!\!\!\!\!\!\!\!\!\!\!\!\!\!\!\!\!\!\!\!\!\!\!\!\!\!\!
\!\!\!\!\!\!\!\!\!\!\!\!\!\!\!\!\!\!\!\!\!\!\!\!\!\!\!\!\!\!\!\! }
\begin{aligned}[t]
T^{[\om]\cdot\be}e^{\mu_L(\be)/2}(\cdots({\mathfrak
M}_{k+1}(\be)*_{2,L,1}b)*_{3,L,1}b)\cdots)
*_{k+1,L,1}b)=0.
\end{aligned}
\end{split}
\label{kh6eq36}
\e

Let ${\rm inc}:L\ra M$ be the inclusion. Define gauge-fixing data
$\bs G_L$ for $(L,{\rm inc})$ to be $\bs G_L=(\bs
I_L,\bs\eta_L,G_L)$, where $\bs I_L$ has indexing set $I=\{n\}$, one
Kuranishi neighbourhood $(L,L,0,\id_L)$, and map ${\rm inc}:L\ra M$,
and $\bs\eta=(\eta_n,\eta_n^n)$ with $\eta_n\equiv 1\equiv\eta_n^n$,
and $G_L:L\ra P$ is the function used in the definition of the
products $*_{\bs f_i,L,\bs{\ti f}_j}$ in Definition \ref{kh6def9}.
As $L$ is compact, oriented and without boundary, we have $[L,{\rm
inc},\bs G_L]\in KC_n(M;\Q)$ with $\pd[L,{\rm inc},\bs G_L]=0$, so
the homology class $\bigl[[L,{\rm inc},\bs G_L]\bigr]$ is defined in
$KH_n(M;\Q)$. We have $\Pi_\rsi^\Kh\bigl([L]\bigr)=\bigl[[L,{\rm
inc},\bs G_L]\bigr]$, where $[L]\in H_n^\rsi(M;\Q)$ is the homology
class of $L$ in $M$, and $\Pi_\rsi^\Kh:H_n^\rsi(M;\Q)\ra KH_n(M;\Q)$
is an isomorphism by Corollary~\ref{kh4cor1}.

Suppose that $[L]=0$ in $H_n^\rsi(M;\Q)$. Note that we don't require
$L$ to be connected, it can be the union of several compact
lagrangian submanifolds whose homology classes sum to zero. Then
$\bigl[[L,{\rm inc},\bs G_L]\bigr]=0$ in $KH_n(M;\Q)$, so we can
choose a Kuranishi chain $B\in KC_{n+1}(M;\Q)$ with $\pd B=[L,{\rm
inc},\bs G_L]$. Let $\be'\in H_2(M;\Z)$, and form the cap product
\begin{equation*}
C\cap\bigl[\oM_{0,1}(M,J,\be'),\bev_1,C_{0,1}(M,J,\be')\bigr]\in
KC_*(M;\Q),
\end{equation*}
where $\bs C_{0,1}(M,J,\be')$ is the co-gauge-fixing data for
$\bigl(\oM_{0,1}(M,J,\be'),\bev_1\bigr)$ chosen in Theorem
\ref{kh6thm12}. Writing $\pi:M\ra\{0\}$ for the projection, we apply
the pushforward $\pi_*:KC_*(M;\Q)\ra KC_*(\{0\};\Q)$ to define
\e
\pi_*\bigl(C\cap\bigl[\oM_{0,1}(M,J,\be'),\bev_1,C_{0,1}(M,J,\be')
\bigr]\bigr)\in KC_*(\{0\};\Q).
\label{kh6eq37}
\e

Now over a single point $\{0\}$, strongly smooth $\bs f:X\ra\{0\}$
coincide with strong submersions, and orientations for $X$ coincide
with coorientations for $(X,\bs f)$, and gauge-fixing data and
co-gauge-fixing for $(X,\bs f)$ coincide. Thus there is a natural
identification $KC_k(\{0\};\Q)=KC^{-k}(\{0\};\Q)$, which identifies
$\pd$ on $KC_*(\{0\};\Q)$ with $\d$ on $KC^*(\{0\};\Q)$. Using this
we regard \eq{kh6eq37} as an element of $KC^*(\{0\};\Q)$. Then using
$\pd B=[L,{\rm inc},\bs G_L]$ and comparing the way $G_L:L\ra P$
appears in $\bs G_L$ and the product $*_{\bev_1, L,\id_L}$, we can
show that
\e
\begin{split}
\d\bigl(\pi_*&(C\cap[\oM_{0,1}(M,J,\be'),\bev_1,C_{0,1}(M,J,\be')])
\bigr)= \\
\bigl[&\oM_{0,1}(M,J,\be')\t_{\bev_1,M,{\rm
inc}}L,\bs\pi,\\
& {\rm inc}^*\bigl(\bs C_{0,1}(M,J,\be')\bigr)*_{\bev_1,
L,\id_L}[L,\id_L,\bs C_L]\bigr]={\mathfrak N}_0(\be')
\end{split}
\label{kh6eq38}
\e
in $KC_*(\{0\};\Q)$, by \eq{kh6eq26}. We will use this to cancel the
${\mathfrak N}_0(\be')$ terms in~\eq{kh6eq28}.

Now define an expression $OGW(M,L)$ in
$\widehat{KC}{}^*(\{0\};\La_\nov^0)$ by
\ea
&OGW(M,L)=\sum_{\be\in
H_2(M,L;\Z)\!\!\!\!}T^{[\om]\cdot\be}e^{\mu_L(\be)/2} {\mathfrak
M}_0(\be) +
\label{kh6eq39}\\
&\sum_{\begin{subarray}{l}l\ge 1,\; \be\in H_2(M,L;\Z):\\
\oM_l^\ma(M,L,J,\be)\ne\es\end{subarray}
\!\!\!\!\!\!\!\!\!\!\!\!\!\!\!\!\!\!\!\!\!\!\!\!\!\!\!\!}
T^{[\om]\cdot\be}e^{\mu_L(\be)/2}\frac{2-l}{l}(\cdots({\mathfrak
M}_l(\be)*_{1,L,1}b)*_{2,L,1}b)\cdots)*_{l,L,1}b)
\nonumber\\
&-\sum_{\begin{subarray}{l}\be'\in H_2(L;\Z):\\
\oM_{0,1}(M,J,\be')\ne\es\end{subarray}
\!\!\!\!\!\!\!\!\!\!\!\!\!\!\!\!\!\!\!\!\!\!\!\!\!\!\!\!}
T^{[\om]\cdot\be'}e^{c_1(M)\cdot\be'}
\pi_*\bigl(C\cap\bigl[\oM_{0,1}(M,J,\be'),\bev_1,C_{0,1}
(M,J,\be')\bigr]\bigr).
\nonumber
\ea
Using \eq{kh6eq27}, \eq{kh6eq28}, \eq{kh6eq36} and \eq{kh6eq38}, we
compute that $\d\bigl(OGW(M,L)\bigr)=0$. Thus we may form
$[OGW(M,L)]\in \widehat{KH}{}^*(\{0\};\La_\nov^0)=
KH^*(\{0\};\Q)\hat\ot\La_\nov^0$.

Identifying $KH^*(\{0\};\Q)\cong H_\cs^*(\{0\};\Q)\cong\Q$ we have
$\widehat{KH}{}^*(\{0\}; \La_\nov^0)\cong\La_\nov^0$. Computing
dimensions shows that $OGW(M,L)$ is graded of degree $3-n$, so
$[OGW(M,L)]\in\La_\nov^{0\,(3-n)}$. When $n$ is even we have
$\La_\nov^{0\,(3-n)}$ by \eq{kh6eq30}, so $[OGW(M,L)]=0$. When $n$
is odd, we may write
\e
[OGW(M,L)]=\ts\sum_{\la\ge 0}OGW(M,L,\la)T^\la e^{(3-n)/2},
\label{kh6eq40}
\e
with $OGW(M,L,\la)\in\Q$.

We interpret these rational numbers $OGW(M,L,\la)$ as {\it open
Gromov--Witten invariants\/} which `count' all stable
$J$-holomorphic discs in $M$ with boundary in $L$ in classes $\be\in
H_2(M,L;\Z)$ with $[\om]\cdot\be=\la$ and $\mu_L(\be)=3-n$. The
meaning of the restriction $\mu_L(\be)=3-n$ is that by \eq{kh6eq15}
this is the condition for $\vdim\oM_0^\ma(M,L,J,\be)=0$, and
$[OGW(M,L)]$ `counts' only moduli spaces $\oM_0^\ma(M,L,J,\be)$ of
virtual dimension zero. When $n$ is even all $\oM_0^\ma(M,L,J,\be)$
have odd virtual dimension, so there are no $\oM_0^\ma(M,L,J,\be)$
with virtual dimension zero. Probably the most interesting case is
that of a graded Lagrangian $L$ in a Calabi--Yau 3-fold $(M,\om)$,
when $\mu_L(\be)=3-n$ holds automatically.

We will show in \cite{Joyc2} that the $OGW(M,L,\la)$ are {\it
independent of the almost complex structure\/} $J$, in the same
sense that Lagrangian Floer cohomology is. They do depend on the
choice of bounding cochain $b$. They depend on the chain $B$ with
$\pd B=[L,{\rm inc},\bs G_L]$ only through the relative homology
class of $B$ in $H_{n+1}(M,L;\Q)$. That is, if we replace $B$ by
$B'=B+\pd C$ for $C\in KC_{n+2}(M;\Q)$, then the $OGW(M,L,\la)$ are
unchanged.

The fact that we can only count discs with fixed area and Maslov
index, rather than counting discs in a fixed $\be\in H_2(M,L;\Z)$,
is because we are using a bounding cochain $b$ defined using a
Novikov ring $\La_\nov^0$ which only keeps track of the area and
Maslov index of curves. As we explain in \cite{Joyc2}, we can refine
the invariants if we either work over a refined Novikov ring, or
make stronger assumptions on the topology of~$M,L$.

\subsection{Conclusions}
\label{kh68}

There is a large area of symplectic geometry which uses moduli
spaces of $J$-holomorphic curves in a symplectic manifold $(M,\om)$
to define some homological structure, and then shows that this
structure is independent of the choice of almost complex structure
$J$, up to some kind of isomorphism, so that the structure basically
depends only on $(M,\om)$. This includes Gromov--Witten invariants
and Lagrangian Floer cohomology, as discussed above, and also
contact homology \cite{EES}, Symplectic Field Theory \cite{EGH},
Fukaya categories \cite{Fuka1,Seid1,Seid3}, and symplectic
cohomology~\cite{Seid2}.

To carry out such a programme, one has to solve four problems:
\begin{itemize}
\setlength{\itemsep}{0pt}
\setlength{\parsep}{0pt}
\item[(a)] Define a geometric structure which is the model
structure for the moduli spaces of $J$-holomorphic curves we are
interested in -- for example, Fukaya--Ono Kuranishi spaces. This
structure should be powerful enough that analogues of
differential-geometric ideas for manifolds such as orientations,
smooth maps, submersions, fibre products work for it.
\item[(b)] Show that the moduli spaces of $J$-holomorphic curves we
wish to study carry the structure in (a), and that natural maps of
these moduli spaces are `smooth'. Prove identities between these
structures, such as the boundary formulae \eq{kh6eq17}
and~\eq{kh6eq20}.
\item[(c)] Use this geometric structure to construct {\it virtual
cycles\/} or {\it virtual chains\/} for the moduli spaces, in
some (co)chain complex. Translate relationships between moduli
spaces into algebraic identities upon their virtual chains.
\item[(d)] Draw some interesting geometrical conclusions --- for
instance, define Gromov--Witten invariants, prove the Arnold
Conjecture, and so on.
\end{itemize}

In this book we have developed new tools for solving problem (c).
Previous virtual cycle or chain constructions for solving problem
(c) in symplectic geometry have been given by Fukaya and Ono
\cite{FuOn1,FOOO}, Li and Tian \cite{LiTi2}, Liu and Tian
\cite{LiuTi}, Lu and Tian \cite{LuTi}, McDuff \cite{McDu}, Ruan
\cite{Ruan}, Siebert \cite{Sieb1}, and in the context of their
theory of polyfolds Hofer, Wysocki and Zehnder \cite{HWZ3}. All
these other constructions involve {\it perturbing the moduli
space\/} in some way.

Defining virtual chains and cycles using Kuranishi (co)chains is
very simple, and it eliminates perturbation of moduli spaces, and
problems with transversality. Relationships between moduli spaces
can often be translated into algebraic identities between Kuranishi
(co)chains in a very direct way, as we saw in passing from
``Theorem'' \ref{kh6thm11} to Corollary \ref{kh6cor1}, and these
identities are likely to hold exactly, rather than in some weak
homotopy sense owing to errors introduced during the virtual cycle
construction.

Also, for moduli spaces without boundary, Kuranishi (co)bordism
classes are in effect a kind of virtual cycle, and so a new kind of
solution to problem (c). Kuranishi (co)bordism groups are huge, much
larger than the corresponding (co)homology groups. This has the
disadvantage that they it is probably not feasible to compute them,
but the advantage that they contain more information than
(co)homology, and are also useful tools for studying integrality
questions.

We have based our theory upon Fukaya and Ono's solution to problem
(a), that is, their notion of Kuranishi spaces, with some
modifications of our own. For our applications above we have simply
assumed the solution to problem (b), when in \S\ref{kh61} and
\S\ref{kh64} we quoted results of Fukaya et al.\ \cite{FuOn1,FOOO}
on existence of Kuranishi structures on the moduli spaces
$\oM_{g,m}(M,J,\be),\oM^\ma_l(M,L,J,\be)$.

Both of these might be seen as a weakness. There are other solutions
to problem (a) in the literature, notably the theory of {\it
polyfolds} due to Hofer, Wysocki and Zehnder
\cite{Hofe1,Hofe2,HWZ1,HWZ2,HWZ3}. If the theory of polyfolds (or
some other alternative) becomes widely adopted, then our Kuranishi
space framework might come to seem obsolete. Also, as discussed in
Remark \ref{kh6rem1}, it is not clear to the author that the proofs
of the solution to problem (b) we assume are complete, and rigorous
in every detail, and some of the proofs are based on weaker
definitions of Kuranishi space than the one we use. In their
polyfold programme, Hofer et al.\ aim to give a fully rigorous
treatment of problem~(b).

In fact this is not a weakness, because a polyfold structure can be
truncated to a Kuranishi structure in a simple way: Hofer
\cite[p.~2-3]{Hofe2} says:
\begin{quotation}
\noindent `The work to build a Kuranishi structure or to build a
polyfold structure (which is a finer structure) is approximately the
same and seems to have more or less the same ingredients. The
polyfold structure captures not only more information but comes with
a strikingly easier formalism. There is in fact a trivial functor
from the polyfold Fredholm theory into the Kuranishi structures.'
\end{quotation}

Kuranishi structures are probably more-or-less the weakest feasible
solution to problem (a), because they are built out of
finite-dimensional data, that is, Kuranishi neighbourhoods
$(V_p,E_p,s_p,\psi_p)$. Other solutions such as polyfolds are
generally built from infinite-dimensional ingredients, and are much
richer structures, with more information; one passes to Kuranishi
structures by forgetting much of this information. Thus, one can use
the theory of polyfolds to solve problems (a) and (b), and then
truncate to Kuranishi structures, and use our Kuranishi (co)homology
or Kuranishi (co)bordism to solve problems (c) and (d). This gives
an alternative to the virtual cycles for polyfolds defined by Hofer
et al.\ \cite{HWZ3}, which work by perturbing the polyfold structure
using $\Q$-weighted multisections in a similar way to Fukaya and
Ono~\cite{FuOn1}.

\appendix

\section{Tent functions}
\label{khA}

We will now define and study {\it tent functions}, which are the
principal tools we will use in the proofs of Theorems \ref{kh4thm1}
and \ref{kh4thm2} in Appendices \ref{khB} and \ref{khC}. Here is the
rough idea. For $S$ a set, write $F(S)$ for the set of finite,
nonempty subsets of $S$. Let $X$ be a compact, oriented $n$-manifold
with g-corners. A {\it tent function} on $X$ is a function $T:X\ra
F([1,\iy))$ which locally near $p\in X$ may be written
$T(p)=\{f_1(p),\ldots,f_k(p)\}$ for smooth functions
$f_1,\ldots,f_k:X\ra[1,\iy)$ called the {\it branches} of $T$, such
that $f_i-f_j:T\ra\R$ is a submersion near $(f_i-f_j)^{-1}(0)$ for
$1\le i<j\le k$, and more complicated transversality conditions hold
when several $f_i$ are equal in~$X$.

Then $\min T:X\ra[1,\iy)$ is a continuous, piecewise smooth
function, and
\e
Z_{X,T}=\bigl\{(t,p)\in[0,\iy)\t X:t\le\min T(p)\bigr\}
\label{khAeq1}
\e
is a compact, oriented $(n+1)$-manifold with g-corners. Also
\e
\pd Z_{X,T}=-Z_{\pd X,T\vert_{\pd X}}\amalg -\{0\}\t X
\amalg\ts\coprod_{c\in C}X_c
\label{khAeq2}
\e
in oriented $n$-manifolds with g-corners, where $C$ is a finite
indexing set, and $X_c$ for $c\in C$ are the connected components of
the portion of $\pd Z_{X,T}$ meeting $(0,\iy)\t X^\ci$. Locally the
$X_c$ correspond to branches $f_c:T\ra[1,\iy)$ of $T$, and $X_c\cong
\bigl\{p\in X:f_c(p)=\min T(p)\bigr\}$. The reason for the name
`tent function' is that when we sketch $Z_{X,T}$ in \eq{khAeq1}, it
looks like a tent, as in Figure \ref{khAfig1} on
page~\pageref{khAfig1}.

We will also define tent functions on {\it orbifolds}, and on {\it
(effective) Kuranishi chains} $[X,\bs f,\bs G]$ or $[X,\bs f,\ubG]$.
The Kuranishi chain version of \eq{khAeq2} is
\e
\begin{split}
\pd\bigl[Z_{X,T},\bs f\ci\bs\pi,\bs H_{X,\bs T}\bigr]=\,&-
\bigl[Z_{\pd X,\bs T\vert_{\pd X}},\bs f\vert_{\pd X}\ci\bs\pi,\bs
H_{\pd X,\bs T\vert_{\pd X}}\bigr]\\
&-[X,\bs f,\bs G]+\ts\sum_{c\in C}[X_c,\bs f_c,\bs G_c].
\end{split}
\label{khAeq3}
\e
This says that $[X,\bs f,\bs G]$ {\it is homologous to} $\sum_{c\in
C}[X_c,\bs f_c,\bs G_c\bigr]$, modulo terms over $\pd X$. Much of
Appendices \ref{khB} and \ref{khC} will involve showing that some
Kuranishi cycle which is a combination of $[X,\bs f,\bs G]$ is
homologous to a second Kuranishi cycle which is a combination of
$\sum_{c\in C}[X_c,\bs f_c,\bs G_c]$, and that {\it this second
cycle has some better property than the first cycle}. It is helpful
to think of \eq{khAeq2} and \eq{khAeq3} as meaning that {\it we have
cut\/ $X$ into finitely many pieces $X_c$ for\/}~$c\in C$.

In fact we will use tent functions for four different purposes:
\begin{itemize}
\setlength{\itemsep}{0pt}
\setlength{\parsep}{0pt}
\item[(a)] To cut a Kuranishi space $X$ into {\it arbitrarily small
pieces\/} $X_c$ for $c\in C$. This will be needed in Step 1 of the
proof of Theorem \ref{kh4thm2} in Appendix~\ref{khC}.
\item[(b)] To cut an effective orbifold $X$ into pieces $X_c$ for
$c\in C$ which are {\it manifolds}, and more generally, to cut an
effective Kuranishi space $X$ into pieces $X_c$ for $c\in C$ which
are Kuranishi spaces with {\it trivial stabilizers}. This will be
needed in Step 1 of the proof of Theorem \ref{kh4thm1} in
Appendix~\ref{khB}.
\item[(c)] To cut a compact $n$-manifold $X$ with g-corners into pieces
$X_c$ for $c\in C$ which are $n$-{\it simplices\/} $\De_n$. This
will be needed in Step 3 of the proof of Theorem \ref{kh4thm1} in
Appendix \ref{khB}, to construct a cycle in singular homology.
\item[(d)] To get round the failure of the smooth Extension
Principle for manifolds with g-corners, Principle \ref{kh2pri}(c),
by defining a class of piecewise smooth functions on manifolds with
g-corners for which the Extension Principle holds. This is important
in Step 2 of the proof of Theorem \ref{kh4thm1} in Appendix
\ref{khB}, when we perturb Kuranishi spaces to get manifolds; in
effect we are choosing piecewise smooth perturbations with
prescribed boundary values, and this would not be possible working
only with smooth perturbations.
\end{itemize}
These different goals are responsible for the length of this
appendix, and also for a certain lack of coherence, as we jump
between (a)--(d) or mix them up. We start by developing our material
on compact manifolds in \S\ref{khA1}, then generalize to orbifolds
in \S\ref{khA2}, and finally to Kuranishi chains in~\S\ref{khA3}.

\subsection{Tent functions on compact manifolds}
\label{khA1}

\subsubsection{The definition of tent functions}
\label{khA11}

\begin{dfn} For $S$ a set, write $F(S)$ for the set of finite,
nonempty subsets of $S$. Let $X$ be an $n$-manifold with g-corners.
We call $T:X\ra F\bigl([1,\iy)\bigr)$ a {\it tent function\/} if for
each $x\in X$ there exists an open neighbourhood $U$ of $x$ in $X$,
open subsets $U_1,\ldots,U_N\subseteq U$ and smooth functions
$t_i:U_i\ra[1,\iy)$ for $i=1,\ldots,N$ such that
$T(u)=\bigl\{t_i(u):i=1,\ldots,N$, $u\in U_i\bigr\}$ for all $u\in
U$, satisfying the following conditions:
\begin{itemize}
\setlength{\itemsep}{0pt}
\setlength{\parsep}{0pt}
\item[(a)] For every subset $\{i_1,\ldots,i_l\}$ of
$\{1,\ldots,N\}$ with $l\ge 2$, the subset
$S_{\{i_1,\ldots,i_l\}}\ab=\bigl\{u\in U_{i_1}\cap\cdots\cap
U_{i_l}:t_{i_1}(u)=\cdots=t_{i_l}(u)\bigr\}$ is an embedded
submanifold of $U_{i_1}\cap\cdots\cap U_{i_l}\subseteq X$. We
require that $S_{\{i_1,\ldots,i_l\}}$ should intersect the
codimension $k$ stratum $S_k(X)$ in Definition \ref{kh2def3}
transversely for all $k\ge 0$. This implies that
$S_{\{i_1,\ldots,i_l\}}$ has boundary or (g-)corners where it
intersects the boundary or (g-)corners of $X$, and we require that
these should be the only boundary or (g-)corners of
$S_{\{i_1,\ldots,i_l\}}$. That is, we have $\pd^kS_{\{i_1,\ldots,
i_l\}}=\bigl\{u\in\pd^k(U_{i_1}\cap\cdots\cap U_{i_l})\subseteq
\pd^kX:t_{i_1}(u)=\cdots=t_{i_l}(u)\bigr\}$.
\item[(b)] In (a), for all $u\in S_{\{i_1,\ldots,i_l\}}$, as vector
subspaces of $T_u^*X$ we have
\begin{equation*}
\!\!\!\!\!\!\!\!\!\!\!\!\!
\big\langle\d(t_{i_2}\!-\!t_{i_1})\vert_u,\d(t_{i_3}\!-\!t_{i_1})
\vert_u,\ldots,\d(t_{i_l}\!-\!t_{i_1})\vert_u\big\rangle\!=\!
\bigl\{\al\!\in\!T_u^*X:\al\vert_{T_uS_{\{i_1,
\ldots,i_l\}}}\!\equiv\!0\bigr\}.
\end{equation*}
More generally, for all $k\ge 0$ and $u$ in
$\pd^kS_{\{i_1,\ldots,i_l\}}$, we have
\begin{gather*}
\big\langle\d\bigl(t_{i_2}-t_{i_1}\vert_{\pd^k(U_{i_1}\cap\cdots\cap
U_{i_l})}\bigr)\vert_u,\ldots,\d\bigl(t_{i_l}-t_{i_1}\vert_{\pd^k
(U_{i_1}\cap\cdots\cap U_{i_l})}\bigr)\vert_u\big\rangle=\\
\bigl\{\al\in T_u^*\pd^kX:\al\vert_{T_u\pd^kS_{\{i_1,
\ldots,i_l\}}}\equiv 0\bigr\}.
\end{gather*}
\item[(c)] For each $i=1,\ldots,N$, the set $\bigl\{u\in
U_i:t_i(u)=\min T(u)\bigr\}$ is closed in $U$.
\end{itemize}
We call the functions $t_i:U_i\ra[1,\iy)$ {\it branches} of $T$ on
$U$, and the sets $U_1,\ldots,U_N$ {\it branch domains}. If $T:X\ra
F([1,\iy))$ is a tent function, then the {\it minimum\/} of $T$ is
$\min T:X\ra[1,\iy)$, the function taking $x\mapsto\min T(x)$, which
is well-defined as $T(x)$ is a finite, nonempty subset of $[1,\iy)$.
Part (c) implies that $\min T$ is {\it continuous}.

The point here is that locally the branches of $T$ do not have to be
defined on all of $X$, but only on open subsets. Thus, as $x$ moves
about in $X$, the number of elements $T(x)$ can change as branches
of $T$ pop in and out of existence. However, (c) implies that a
branch $t_i$ of $T$ cannot disappear near $x$ if $t_i(x)=\min T(x)$.
Thus, disappearing branches do not cause $\min T$ to become
discontinuous, and $\min T$ near $x$ is locally the minimum of
smooth functions $t_i:X\ra[1,\iy)$ defined near $x$, so $\min T$ is
continuous, and piecewise-smooth.

As (a)--(c) are stable under taking boundaries, if $T:X\ra
F([1,\iy))$ is a tent function then $T\vert_{\pd X}:\pd X\ra
F([1,\iy))$ is also a tent function.
\label{khAdef1}
\end{dfn}

The next proposition shows how we will use tent functions.
Definition \ref{khAdef1} was carefully designed just to make
$Z_{X,T}$ a manifold with g-corners. Think of the proposition this
way: we are cutting $X$ into finitely many pieces $X_c$ for $c\in
C$, which as we will see later can be chosen to be arbitrarily
small, to be simplices $\De_n$, or to have other good properties.
The $(n+1)$-manifold $Z_{X,T}$ is a kind of (solid) graph of $\min
T$, and we will use it to define a homology between a chain
involving $X$, and a sum over $c\in C$ of chains involving~$X_c$.

\begin{prop} Let\/ $X$ be a compact, oriented\/ $n$-manifold with
g-corners, and\/ $T:X\ra F\bigl([1,\iy)\bigr)$ a tent function. Then
\e
Z_{X,T}=\bigl\{(t,x)\in [0,\iy)\t X:t\le\min T(x)\bigr\}
\label{khAeq4}
\e
is a compact, oriented\/ $(n+1)$-manifold with g-corners, and
\e
\pd Z_{X,T}=-Z_{\pd X,T\vert_{\pd X}}\amalg -\{0\}\t X
\amalg\ts\coprod_{c\in C}X_c
\label{khAeq5}
\e
in oriented\/ $n$-manifolds with g-corners, where\/ $C$ is a finite
indexing set, each\/ $X_c$ is connected, and if\/ $\pi_X:[0,\iy)\t
X\ra X$ is the projection, then\/ $\pi=\pi_X\vert_{X_c}:X_c\ra X$ is
an orientation-preserving immersion of\/ $X_c$ as an\/
$n$-submanifold of\/ $X$ for each $c\in C,$ which is an embedding on
$X_c^\ci,$ with\/ $X=\bigcup_{c\in C}\pi(X_c),$ and\/
$\pi(X_c^\ci)\cap \pi(X_{c'}^\ci)=\emptyset$ for all\/ $c\ne c'$
in\/~$C$.
\label{khAprop1}
\end{prop}

\begin{proof} Since $\min T$ is continuous and $X$ is compact, $\min
T$ is bounded on $X$, and $Z_{X,T}$ is compact. As $X$ is an
oriented $n$-manifold with g-corners, $[0,\iy)\t X$ is an oriented
$(n+1)$-manifold with g-corners. We claim $Z_{X,T}$ is an embedded
$(n+1)$-submanifold of $[0,\iy)\t X$. It is enough to verify that
$Z_{X,T}$ is locally modelled on a region with g-corners in
$\R^{n+1}$ near each point $(t,x)\in Z_{X,T}$ with $t=\min T(x)$,
since if $(t,x)\in Z_{X,T}$ with $t<\min T(x)$ then $Z_{X,T}$
coincides with $[0,\iy)\t X$ near~$(t,x)$.

For such $(t,x)$, let $U,U_1,\ldots,U_N,f_1,\ldots,f_N$ be as in
Definition \ref{khAdef1}, and let $\{j_1,\ldots,j_k\}$ be the subset
of $j$ in $\{1,\ldots,N\}$ with $x\in U_j$ and $t_j(x)=\min T(x)$.
Choose a small open neighbourhood $\ti U$ of $x$ in $U$ such that
$\ti U\subseteq U_{j_a}$ for $a=1,\ldots,k$, and if $i\in
\{1,\ldots,N\}\sm\{j_1,\ldots,j_k\}$ and $u\in U_i\cap\ti U$ then
$t_i(u)>\min T(u)$. This last part is possible by Definition
\ref{khAdef1}(c), as for each $i\in \{1,\ldots,N\}\sm
\{j_1,\ldots,j_k\}$, we have to choose $\ti U$ not to intersect a
{\it closed\/} set in $U$ not containing $x$. Then
$t_{j_1}\vert_{\ti U},\ldots,t_{j_k}\vert_{\ti U}:\ti U\ra\R$ are
well-defined, and $\min T(u)=\min\bigl\{t_{j_1}(u),\ldots,
t_{j_k}(u)\bigr\}$ for all~$u\in\ti U$.

Now let $(U',\phi')$ be a chart with g-corners on $X$ with
$x\in\phi'(U')\subseteq\ti U$. Then $U'$ is a region with g-corners
in $\R^n$, defined as in Definition \ref{kh2def1} using some
$W'\subseteq\R^n$ and $f_1',\ldots,f'_{N'}$. As in Definition
\ref{kh2def2}, the smooth maps $t_{j_a}\ci\phi':U'\ra\R$ extend
smoothly to some open subset of $\R^n$ containing $U'$. Make $W'$
smaller if necessary so that $t_{j_a}\ci\phi'$ extends to $W'$ for
$a=1,\ldots,k$, and choose such extensions so that Definition
\ref{khAdef1}(a),(b) still hold.

Define $\dot W\subset\R^{n+1}$ by $\dot W=(0,\iy)\t W'$ and $\dot
f_1,\ldots,\dot f_{k+N'}:\dot W\ra\R$ by
\begin{equation*}
\dot f_a(t,x_1,\ldots,x_n)=\begin{cases}
t_{j_a}\ci\phi'(x_1,\ldots,x_n)-t, & a=1,\ldots,k, \\
f_{a-k}'(x_1,\ldots,x_n), & a=k+1,\ldots,k+N'.
\end{cases}
\end{equation*}
Then Definitions \ref{kh2def1}(a),(b) and \ref{khAdef1}(a),(b) for
$f_1',\ldots,f'_{N'}$ and $t_{j_1},\ldots,t_{j_k}$ imply that $\dot
W,\dot f_1,\ldots,\dot f_{k+N'}$ satisfy Definition
\ref{kh2def1}(a),(b), and the corresponding subset
\begin{equation*}
\dot U=\bigl\{(s,x_1,\ldots,x_n):(x_1,\ldots,x_n)\in U',\;\> 0<s\le
\min T\ci\phi'(x_1,\ldots,x_n)\bigr\}
\end{equation*}
is a region with g-corners in $\R^{n+1}$. Define $\dot\phi:\dot U\ra
Z_{X,T}$ by $\dot\phi:(s,x_1,\ldots,x_n)\mapsto
\bigl(s,\phi'(x_1,\ldots,x_n)\bigr)$. Then $(\dot U,\dot\phi)$ is a
chart with g-corners on $Z_{X,T}$ with $(t,x)$ in $\dot\phi(\dot
U)$. Compatibility between such charts $(U',\phi')$ on $X$ implies
compatibility between the corresponding charts $(\dot U,\dot\phi)$
on $Z_{X,T}$. Therefore $Z_{X,T}$ is a {\it compact, oriented\/
$(n+1)$-manifold with g-corners}.

In \eq{khBeq2}, the first two terms $-Z_{\pd X,T\vert_{\pd X}}\amalg
-\{0\}\t X$ on the right hand side are the portion of $\pd Z_{X,T}$
coming from the boundary of $[0,\iy)\t X$, and their signs come from
$\pd\bigl([0,\iy)\t X\bigr)=-[0,\iy)\t\pd X\amalg -\{0\}\t X$. Apart
from these two terms, $\pd Z_{X,T}$ also has other contributions
from points $(t,x)$ with $t=\min T(x)$. As $\pd Z_{X,T}$ is a
compact manifold, these other contributions have only finitely many
connected components. Define these connected components to be $X_c$
for $c$ in $C$, some finite indexing set. Then \eq{khAeq5} holds, by
definition. If $(t,x),\ldots$ are as above with $\ti U^\ci$
connected and $((t,x),B)\in X_c$ for some $c\in C$ then the local
boundary component $B$ is of the form $\bigl\{(t,u)\in
Z_{X,T}:u\in\ti U^\ci$, $t=\phi_{j_a}(u)\bigr\}$ for some
$a=1,\ldots,k$, and then $\bigl\{(t,u)\in Z_{X,T}:u\in\ti U$,
$t=\phi_{j_a}(u)\bigr\}$ is a local model for $X_c$. The remaining
claims about the $X_c$ follow.
\end{proof}

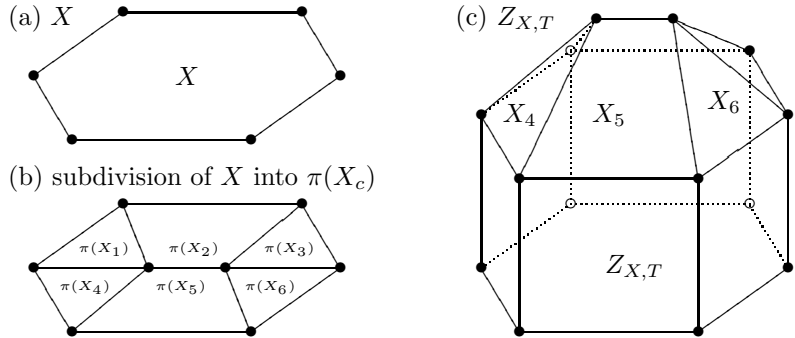
\begin{figure}[htb]
\centerline{
$\begin{xy} 0;<1.7mm,0mm>: ,(9,20)*{\bullet} ,(5,10)*{\bullet}
,(12,15)*{\bullet} ,(-5,20)*{\bullet} ,(-9,10)*{\bullet}
,(-12,15)*{\bullet} ,(0,15)*{X} ,(-14,19.5)*!L{\text{(a) $X$}}
,(3,0)*{\bullet} ,(9,5)*{\bullet} ,(5,-5)*{\bullet}
,(12,0)*{\bullet} ,(-3,0)*{\bullet} ,(-5,5)*{\bullet}
,(-9,-5)*{\bullet} ,(-12,0)*{\bullet} ,(-6.5,1.5)*{\sst\pi(X_1)}
,(.5,1.5)*{\sst\pi(X_2)} ,(8,1.5)*{\sst\pi(X_3)}
,(-8,-1.5)*{\sst\pi(X_4)} ,(-.5,-1.5)*{\sst\pi(X_5)}
,(6.5,-1.5)*{\sst\pi(X_6)} ,(-14,7)*!L{\text{(b) subdivision of $X$
into $\pi(X_c)$}} ,(44,5)*{\circ} ,(40,-5)*{\bullet}
,(47,0)*{\bullet} ,(30,5)*{\circ} ,(26,-5)*{\bullet}
,(23,0)*{\bullet} ,(44,17)*{\bullet} ,(40,7)*{\bullet}
,(47,12)*{\bullet} ,(30,17)*{\circ} ,(26,7)*{\bullet}
,(23,12)*{\bullet} ,(38,19.5)*{\bullet} ,(32,19.5)*{\bullet}
,(35,0)*{Z_{X,T}} ,(21,19.5)*!L{\text{(c) $Z_{X,T}$}} ,(26,12)*{X_4}
,(33,12)*{X_5} ,(42,13)*{X_6} ,\ar@{-}(9,20);(12,15)
,\ar@{-}(5,10);(12,15) ,\ar@{-}(-5,20);(-12,15)
,\ar@{-}(-9,10);(-12,15) ,\ar@{-}(-5,20);(9,20)
,\ar@{-}(-9,10);(5,10) ,\ar@{-}(-3,0);(3,0) ,\ar@{-}(3,0);(9,5)
,\ar@{-}(3,0);(12,0) ,\ar@{-}(3,0);(5,-5) ,\ar@{-}(9,5);(12,0)
,\ar@{-}(5,-5);(12,0) ,\ar@{-}(-3,0);(-5,5) ,\ar@{-}(-3,0);(-12,0)
,\ar@{-}(-3,0);(-9,-5) ,\ar@{-}(-5,5);(-12,0)
,\ar@{-}(-9,-5);(-12,0) ,\ar@{-}(-5,5);(9,5) ,\ar@{-}(-9,-5);(5,-5)
,\ar@{.}(44,5);(47,0) ,\ar@{-}(40,-5);(47,0) ,\ar@{.}(30,5);(23,0)
,\ar@{-}(26,-5);(23,0) ,\ar@{.}(30,5);(44,5) ,\ar@{-}(26,-5);(40,-5)
,\ar@{-}(44,17);(47,12) ,\ar@{-}(40,7);(47,12)
,\ar@{.}(30,17);(23,12) ,\ar@{-}(26,7);(23,12)
,\ar@{.}(30,17);(44,17) ,\ar@{-}(26,7);(40,7) ,\ar@{.}(44,5);(44,17)
,\ar@{-}(47,0);(47,12) ,\ar@{.}(30,5);(30,17) ,\ar@{-}(23,0);(23,12)
,\ar@{-}(26,-5);(26,7) ,\ar@{-}(40,-5);(40,7)
,\ar@{-}(44,17);(38,19.5) ,\ar@{-}(47,12);(38,19.5)
,\ar@{.}(30,17);(32,19.5) ,\ar@{-}(23,12);(32,19.5)
,\ar@{-}(26,7);(32,19.5) ,\ar@{-}(40,7);(38,19.5)
,\ar@{-}(32,19.5);(38,19.5)
\end{xy}$}
\caption{Example of a tent function used to subdivide a polygon}
\label{khAfig1}
\end{figure}

The reason for the name {\it tent function\/} is that when we sketch
$Z_{X,T}$ in \eq{khAeq4}, it looks rather like a tent. We illustrate
this in Figure~\ref{khAfig1}.

\subsubsection{Cutting a compact manifold into arbitrarily small
pieces}
\label{khA12}

Here is our first method for constructing tent functions.

\begin{dfn} Let $X$ be a compact $n$-manifold with g-corners. Let
$g$ be a Riemannian metric on $X$, and choose $\ep>0$ with
$2\ep<\de(g)$, where $\de(g)$ is the {\it injectivity radius\/} of
$g$. Write $d(\,,\,)$ for the metric on $X$ determined by $g$. For
$p\in X$ and $r>0$ write $B_r(p)$ for the open ball of radius $r$
about $p$, that is, $B_r(p)=\bigl\{x\in X:d(x,p)<r\bigr\}$. Then
using geodesic normal coordinates on $X$ near $p$, we see that for
any $p\in X$, the map $X\ra[0,\iy)$ given by $x\mapsto d(p,x)^2$ is
continuous, and smooth wherever $d(p,x)<\de(g)$. Hence $x\mapsto
d(p,x)^2$ is smooth on $B_{2\ep}(p)$, as~$2\ep<\de(g)$.

The open sets $B_\ep(p)$ for $p\in X$ cover $X$, so as $X$ is
compact we can choose a finite subset $\{p_1,\ldots,p_N\}$ of $X$
with $X=\bigcup_{c=1}^NB_\ep(p_c)$. Choose $p_1,\ldots,p_N$ to be
{\it generic\/} amongst all such $N$-tuples $p_1',\ldots,p_N'$; this
can be done by first choosing any $p_1',\ldots,p_N'$ with
$X=\bigcup_{c=1}^NB_\ep(p_c')$, and then taking $p_c\in X$ to be
generic with $d(p_c,p_c')\ll\ep$. Define $T:X\ra
F\bigl([1,\iy)\bigr)$ by
\e
T(x)=\bigl\{1+d(p_c,x)^2:c=1,\ldots,N,\;\> d(p_c,x)<2\ep\bigr\}.
\label{khAeq6}
\e
\label{khAdef2}
\end{dfn}

\begin{prop} In Definition {\rm\ref{khAdef2},} $T$ is a tent
function on\/ $X$. In {\rm\eq{khAeq5},} we can take the indexing
set\/ $C$ to be $\{1,\ldots,N\}$ with\/ $p_c\in\pi(X_c^\ci)\subseteq
\pi(X_c)\subseteq B_\ep(p_c)$ for\/ $c=1,\ldots,N$. We can choose\/
$T$ so that the pieces\/ $X_c$ are `arbitrarily small', in the sense
that if\/ $\{V_i:i\in I\}$ is any open cover for $X$ then we can
choose $T$ so that for all\/ $c=1,\ldots,N$ there exists\/ $i\in I$
with\/~$\pi(X_c)\subseteq V_i$.
\label{khAprop2}
\end{prop}

\begin{proof} Since $X=\bigcup_{c=1}^NB_\ep(p_c)$, if $x\in X$ there
exists at least one $c$ with $d(p_c,x)<2\ep$, and then
$1+d(p_c,x)^2\in T(x)$. So $T(x)$ is nonempty, and $T$ does map
$X\ra F\bigl([1,\iy)\bigr)$. In Definition \ref{khAdef1}, for any
$x\in X$ we can take $U=X$ and $U_c=B_{2\ep}(p_c)$ and define
$t_c:U_c\ra [1,\iy)$ by $t_c(u)=1+d(p_c,u)^2$. Then $U_c$ is open in
$U$, $t_c$ is smooth, and $T(u)=\bigl\{t_i(u):i=1,\ldots,N$, $u\in
U_i\bigr\}$.

Definition \ref{khAdef1}(a),(b) hold because $p_1,\ldots,p_N$ are
{\it generic}, and we also have the extra property that $\dim
S_{\{i_1,\ldots,i_l\}}=n-l+1$ for all
$\{i_1,\ldots,i_l\}\subseteq\{1,\ldots,N\}$ with $l\ge 2$, so that
$\d(t_{i_2}-t_{i_1}),\ldots,\d(t_{i_l}-t_{i_1})$ are linearly
independent along $S_{\{i_1,\ldots,i_l\}}$. To see this, note that
$S_{\{i_1,\ldots,i_l\}}=\bigl\{x\in X: d(p_{i_1},x)=d(p_{i_2},x)=
\cdots=d(p_{i_l},x)<2\ep\bigr\}$, and as $p_{i_1},\ldots,p_{i_l}$
are generic, the $l-1$ equations $d(p_{i_1},x)=d(p_{i_2},x)=
\cdots=d(p_{i_l},x)$ are transverse.

For $i=1,\ldots,N$, we have $\bigl\{u\in U_i:t_i(u)=\min
T(u)\bigr\}=\bigl\{u\in B_{2\ep}(p_i):d(p_i,u)=\min
\{d(p_c,u):c=1,\ldots,N\}\bigr\}$. Since
$X=\bigcup_{c=1}^NB_\ep(p_c)$, we have $\min
\{d(p_c,u):c=1,\ldots,N\}<\ep$, so $\bigl\{u\in U_i:t_i(u)=\min
T(u)\bigr\}\subseteq B_{\ep}(p_i)$. Thus $\bigl\{u\in
U_i:t_i(u)=\min T(u)\bigr\}$ is a closed subset of
$U_i=B_{2\ep}(p_i)$, and is contained in
$\,\overline{\!B}_\ep(p_i)$. Since $\,\overline{\!B}_\ep(p_i)$ is
closed in both $U_i$ and $X$, it follows that $\bigl\{u\in
U_i:t_i(u)=\min T(u)\bigr\}$ is closed in $X$, and Definition
\ref{khAdef1}(c) holds. Thus $T$ is a {\it tent function}.

For the second part, it is easy to see that the remaining boundary
components $\coprod_{c\in C}X_c$ in \eq{khAeq5} may be written
\e
\ts\coprod_{c=1}^N\bigl\{\bigl(1\!+\!d(p_c,u)^2,u\bigr):u\in X,\;
d(p_c,u)\!=\!\min\{d(p_a,u):a\!=\!1,\ldots,N\}\bigr\},
\label{khAeq7}
\e
omitting local boundary components $B$. Now each $u$ in the set
$\bigl\{u\in X:d(p_c,u)=\min\{d(p_a,u):a=1,\ldots,N\}\bigr\}$ is
joined to $p_c$ by a unique geodesic segment of length $d(p_c,u)$
which also lies in the set. Hence the $c$ term in \eq{khAeq7} is
nonempty and connected. But by definition the connected components
of \eq{khAeq7} are $X_c$ for $c\in C$. Therefore we can take
$C=\{1,\ldots,N\}$, and $X_c$ to be the $c$ term in \eq{khAeq7}.
This gives $p_c\in\pi(X_c^\ci)\subseteq \pi(X_c)\subseteq
B_\ep(p_c)$, since $\min\{d(p_a,u):a=1,\ldots,N\}<\ep$ for all~$u\in
X$.

Given an open cover $\{V_i:i\in I\}$ for $X$, we fix $g$ in
Definition \ref{khAdef2} and choose $\ep$ sufficiently small that
for all $x\in X$ there exists $i\in I$ with $B_\ep(x)\subseteq V_i$,
which is possible by compactness of $X$. The final part of the
proposition then follows from $\pi(X_c)\subseteq B_\ep(p_c)$ in the
second part.
\end{proof}

\subsubsection{Cutting a compact manifold into simplices}
\label{khA13}

In Definition \ref{khAdef2}, regard $X$ and $g$ as fixed, and choose
$\ep>0$ very small compared to all the natural length-scales of
$(X,g)$. That is, choose $\ep$ with $\ep\ll\de(g)$, and $\ep\ll
\de(g\vert_{S_k(X)})$ for $0<k<n$, and $\Vert R(g)\Vert_{C^0}
\ll\ep^{-2}$, and $\Vert R(g\vert_{S_k(X)})\Vert_{C^0} \ll\ep^{-2}$
for $0<k<n$, where $\de(g)$ is the injectivity radius and $R(g)$ the
Riemann curvature of $g$, and the second fundamental form of
$S_k(X)$ in $X$ is $\ll\ep^{-1}$ for $0<k<n$. Then using {\it
geodesic normal coordinates}, we see that balls $B_{2\ep}(x)$ of
radius $2\ep$ in the interior of $(X,g)$ are approximately isometric
to balls of radius $2\ep$ in $\R^n$ with its Euclidean metric, and
balls of radius $2\ep$ close to the boundary are approximately
isometric to balls of radius $2\ep$ in a {\it polyhedral cone} in
$\R^n$ with its Euclidean metric, since as in Remark
\ref{kh2rem1}(c),(d), manifolds with g-corners are approximately
locally modelled on polyhedral cones in~$\R^n$.

Therefore geodesic normal coordinates at $p_c$ approximately
identify $\pi(X_c)\subset X$ with a subset of $\R^n$, or a
polyhedral cone in $\R^n$, defined by inequalities $d_{\rm
Eu}(p_c,u)\le d_{\rm Eu}(p_a,u)$ for all $a=1,\ldots,N$, for some
distinct points $p_1,\ldots,p_N$ in $\R^n$, where $d_{\rm Eu}$ is
the Euclidean metric on $\R^n$. Such inequalities define a {\it
convex polyhedron} in $\R^n$, which is compact as $X_c$ is compact.
So we deduce:

\begin{lem} In Definition {\rm\ref{khAdef2},} if\/ $\ep$ is
small compared to the natural length-scales of\/ $(X,g),$ then
geodesic normal coordinates at\/ $p_c$ in $X$ approximately identify
$\pi(X_c)\subset X$ with a compact convex polyhedron
in\/~$T_{p_c}X\cong\R^n$.
\label{khAlem1}
\end{lem}

Note that Remark \ref{kh2rem1}(e) implies that $\pi(X_c)$ need not
be diffeomorphic with a convex polyhedron in $\R^n$, so our notion
of `approximately identifies' does not imply diffeomorphism, just
that the defining equations of $\pi(X_c)$ can be slightly perturbed,
without changing the discrete structure of $\pi(X_c)$ as a manifold
with g-corners, to give a compact convex polyhedron.

There is a standard way of dividing a compact convex polyhedron in
$\R^n$ into $n$-simplices $\De_n$, called {\it barycentric
subdivision}, as in Bredon~\cite[\S IV.17]{Bred}.

\begin{dfn} Let $K$ be a compact convex polyhedron in $\R^n$. For
each $k$-dimensional face $F$ of $K$ for $k=0,\ldots,n$, define the
{\it barycentre\/} $b_F$ of $F$ to be the point in the interior of
$F$ whose position vector is the average of the position vectors of
the vertices of $F$. Let $\bs F=F_0\subset F_1\subset\cdots\subset
F_n=K$ be a chain of faces of $K$ with $\dim F_j=j$ for
$j=0,\ldots,n$; we call $\bs F$ a {\it flag} for $K$. Write
$p_0,\ldots,p_n$ for the vertices of $\De_n$,
where~$p_j=(\de_{j0},\de_{j1},\ldots,\de_{jn})$.

Let $\si_{\bs F}:\De_n\ra\R^n$ be the unique affine map with
$\si_{\bs F}(p_j)=b_{F_j}$ for $j=0,\ldots,n$. Then $K$ is the union
of $\si_{\bs F}(\De_n)$ over all flags $\bs F$ for $K$, and
$\si_{\bs F}(\De_n^\ci)\cap \si_{\bs F'}(\De_n^\ci)=\emptyset$ for
$\bs F\ne\bs F'$. This division of $K$ into $n$-simplices $\si_{\bs
F}(\De_n)$ is called the {\it barycentric subdivision} of $K$. It is
compatible with boundaries in two ways:
\begin{itemize}
\setlength{\itemsep}{0pt}
\setlength{\parsep}{0pt}
\item[(i)] Suppose $\bs F=F_0\subset F_1\subset\cdots\subset
F_n=K$ and $\bs F'=F_0'\subset F_1'\subset\cdots\subset F_n'=K$
are two flags for $K$, and $F_j=F_j'$ for $j\in
\{j_0,j_1,\ldots,j_l\}\subseteq\{0,\ldots,k\}$, with
$j_0<j_1<\cdots<j_l$. Let $\De_l^{\bs j}\subseteq\De_n$ be the
$l$-dimensional face $\De_l$ of $\De_n$ with vertices
$p_{j_0},\ldots,p_{j_l}$. Then $\si_{\bs F}\vert_{\De_l^{\bs
j}}\equiv\si_{\bs F'}\vert_{\De_l^{\bs j}}$.
\item[(ii)] The barycentric subdivision of $K$ restricts on $\pd
K$ to the disjoint union of the barycentric subdivisions of each
boundary face of~$K$.
\end{itemize}

We can construct the barycentric subdivision of $K$ using a {\it
tent function}. Choose constants $0=c_0 < c_1\ll c_2\ll\cdots\ll
c_n$. For each flag $\bs F=F_0\subset F_1\subset \cdots\subset
F_n=K$ for $K$, define $A_{\bs F}:\R^n\ra\R$ to be the unique affine
function with $A_{\bs F}(b_{F_j})=c_j$ for $j=0,\ldots,n$. Then
define $T:K\ra F\bigl([1,\iy)\bigr)$ by
\begin{equation*}
T(x)=\bigl\{2+\de A_{\bs F}(x):\text{$\bs F$ is a flag for
$K$}\bigr\},
\end{equation*}
where $\de>0$ is sufficiently small that $\de A_{\bs F}\ge -1$ on
$K$ for all $\bs F$. One can show that if the ratios
$c_2/c_1,\ldots,c_n/c_{n-1}$ are large enough then $T$ is a tent
function, where in Definition \ref{khAdef1} we take
$U=U_1=\ldots=U_N=K$ and $t_j=2+\de A_{\bs F_j}$, where $\bs
F_1,\ldots,\bs F_N$ are the possible flags for $K$, and
\begin{equation*}
\si_{\bs F}(\De_n)=\bigl\{x\in K:\text{$A_{\bs F}(x)\le A_{\bs
F'}(x)$ for all flags $\bs F'$ for $K$}\bigr\},
\end{equation*}
so in \eq{khAeq5}, the pieces $\pi(X_c)$ for all $c\in C$ are
$\si_{\bs F}(\De_n)$ for all flags $\bs F$ for $K$.
\label{khAdef3}
\end{dfn}

We can combine the ideas of the last four definitions and results to
construct a tent function $T$ for $X$ such that each $X_c$ for $c\in
C$ is a simplex $\De_n$. That is, we use a tent function to
construct a {\it triangulation} of $X$ by $n$-simplices, in the
sense of Bredon \cite[p.~246]{Bred}. This will enable us to show
that chains defined using compact manifolds $X$ with g-corners are
homologous to chains defined using simplices $\De_n$, that is, to
chains in singular homology.

\begin{thm} Let\/ $X$ be a compact\/ $n$-manifold with g-corners.
Then we can construct a tent function $T:X\ra F\bigl([1,\iy)\bigr)$
such that the components $X_c,$ $c\in C$ of\/ $\pd Z_{X,T}$ in
\eq{khAeq5} are all diffeomorphic to the $n$-simplex $\De_n,$ with
diffeomorphisms $\si_c:\De_n\ra X_c$ for $c\in C$. Thus
\e
\ts\pd^2Z_{X,T}\supseteq \coprod_{c\in C}\pd X_c=\coprod_{c\in
C}\si_c(\pd\De_n)=\coprod_{c\in C}\coprod_{j=0}^n\si_c\ci
F_j^n(\De_{n-1}).
\label{khAeq8}
\e
Definition \ref{kh2def5} gives a free, orientation-reversing
involution $\si:\pd^2Z_{X,T}\ra\ab \pd^2Z_{X,T}$.  We can choose the
$\si_c$ to have the following boundary compatibility: suppose that\/
$\si$ exchanges two $(n-1)$-simplices $\si_c\ci F_j^n(\De_{n-1})$
and\/ $\si_{c'}\ci F_{j'}^n(\De_{n-1})$ in \eq{khAeq8}. Then $j=j',$
$c\ne c'$ and\/ $\si\ci\si_c\ci F_j^n\equiv\si_{c'}\ci
F_j^n:\De_{n-1}\ra\pd^2Z_{X,T}$.
\label{khAthm1}
\end{thm}

\begin{proof} We first sketch the construction, and then fill in
some details. As in Lemma \ref{khAlem1}, use Definition
\ref{khAdef2} with $\ep$ small compared to the length-scales of
$(X,g)$ to construct a tent function we will write as $\ti T:X\ra
F\bigl([1,\iy)\bigr)$, using points $p_1,\ldots,p_N$, and functions
$\ti t_i:B_{2\ep}(p_i)\ra[1,\iy)$ given by $\ti
t_i(u)=1+d(p_i,u)^2$, for $i=1,\ldots,N$. Use the notation $\ti X_i$
for $i=1,\ldots,N$ for the $X_c$ in \eq{khAeq5} for $Z_{X,\ti T}$,
where $p_i\in\pi(\ti X_i)$. Then as in Lemma \ref{khAlem1}, using
geodesic normal coordinates at $p_i$, $\pi(\ti X_i)$ approximates a
compact convex polyhedron $K_i$ in $T_{p_i}X\cong\R^n$
for~$i=1,\ldots,N$.

Define $C=\bigl\{(i,\bs F):i=1,\ldots,N$, $\bs F$ is a flag for
$K_i\bigr\}$. We shall first choose embeddings $\ti\si_{(i,\bs
F)}:\De_n\ra\ti X_i$ for all $(i,\bs F)\in C$ such that for each
$i=1,\ldots,N$, the $\ti\si_{(i,\bs F)}$ for all faces $\bs F$ of
$K_i$ form a triangulation of $\ti X_i$, where $\ti\si_{(i,\bs F)}$
approximates $\si_{\bs F}$ in Definition \ref{khAdef3} for $K_i$
under the approximation $K_i\cong\ti X_i$ in geodesic normal
coordinates. These $\ti\si_{(i,\bs F)}$ will satisfy boundary
compatibilities analogous to Definition~\ref{khAdef3}(i),(ii).

Choose constants $0=c_0<c_1\ll c_2\ll\cdots\ll c_n$, such that the
ratios $c_2/c_1,\ldots,c_n/c_{n-1}$ are large enough for the tent
function construction in Definition \ref{khAdef3} to work for $K_i$
for all $i=1,\ldots,N$. Choose smooth functions $A_{(i,\bs
F)}:B_{2\ep}(p_i)\ra\R$ for all $(i,\bs F)\in C$, such that
$A_{(i,\bs F)}\ci\ti\pi\ci\ti\si_{(i,\bs F)}:\De_n\ra\R$ is the
unique affine function with $A_{(i,\bs F)}\ci\ti\pi\ci\ti\si_{(i,\bs
F)}(p_j)=c_j$ for $j=0,\ldots,N$, where $\ti\pi:Z_{X,\ti T}\ra X$ is
the projection so that $\ti\pi\ci\ti\si_{i,\bs F}:\De_n\ra
B_{2\ep}(p_i)\subseteq X$, and using geodesic normal coordinates to
identify $B_{2\ep}(p_i)$ with the ball of radius $2\ep$ in
$T_{p_i}X$, $A_{(i,\bs F)}$ approximates an affine function on the
whole of~$B_{2\ep}(p_i)$.

Choose $\de>0$ sufficiently small that $\de A_{(i,\bs F)}\ge -1$ and
$\bmd{\de A_{(i,\bs F)}}\ll \ep^2$ on $B_{2\ep}(p_i)$ for all
$(i,\bs F)\in C$. Generalizing \eq{khAeq6}, define $T:X\ra
F\bigl([1,\iy)\bigr)$ by
\e
T(x)=\bigl\{2+d(p_i,x)^2+\de A_{(i,\bs F)}:(i,\bs F)\in C,\;\>
d(p_i,x)<2\ep\bigr\}.
\label{khAeq9}
\e
Then similar proofs to Proposition \ref{khAprop1} and Definition
\ref{khAdef3} show that $T$ is a tent function, and the components
$X_c$ for $c\in C$ of $Z_{X,T}$ in \eq{khAeq5} can be written with
the indexing set $C$ above such that $\pi(X_{(i,\bs F)})=\ti\pi\ci
\ti\si_{(i,\bs F)}(\De_n)$, so that $\ti\si_{(i,\bs F)}$ lifts to a
unique diffeomorphism $\si_{(i,\bs F)}:\De_n\ra X_{(i,\bs F)}$
with~$\pi\ci\si_{(i,\bs F)}\equiv\ti\pi\ci\ti\si_{(i,\bs F)}$.

In carrying out this programme, two points need special care:
\begin{itemize}
\setlength{\itemsep}{0pt}
\setlength{\parsep}{0pt}
\item Because of the local variation of polyhedral cone models along
codimension $l$ strata of manifolds with g-corners for $l\ge 3$
discussed in Remark \ref{kh2rem1}(e), we cannot assume that $\ti
X_i$ is diffeomorphic to $K_i$ for~$i=1,\ldots,N$.
\item The $\ti X_i$ for $i=1,\ldots,N$ are in general manifolds
with g-corners, not corners, and the smooth Extension Principle,
Principle \ref{kh2pri}(c), fails for manifolds with g-corners. So in
the proof we must avoid choosing data such as $A_{(i,\bs F)}$ on
$\pd\ti X_i$ and trying to extend it smoothly over~$\ti X_i$.
\end{itemize}

We now explain how to choose the embeddings $\ti\si_{(i,\bs
F)}:\De_n\ra\ti X_i$ for all $(i,\bs F)\in C$. For the
$\ti\si_{(i,\bs F)}$ to lift to $\si_{(i,\bs F)}$ satisfying the
boundary compatibility $\si\ci\si_c\ci F_j^n\equiv\si_{c'}\ci
F_j^n:\De_{n-1}\ra\pd^2Z_{X,T}$ in the theorem, the $\ti\si_{(i,\bs
F)}$ must satisfy boundary compatibilities on $\pd\De_n$, and these
in turn imply compatibilities between the $\ti\si_{(i,\bs F)}$ on
$\pd^m\De_n$ for $m=1,2,\ldots,n$. In order to satisfy all these
compatibilities, we must choose the $\ti\si_{(i,\bs F)}
\vert_{\pd^m\De_n}$ {\it by induction on decreasing\/}
$m=n,n-1,\ldots,1,0$. Similar proofs where we choose data on
$\pd^mX_a$ by induction on decreasing codimension $m$ will appear in
\S\ref{khB1}--\S\ref{khB3} and~\S\ref{khC1}.

We introduce some notation. For $m=0,\ldots,n$, let $S^m_n$ be the
set of $(n-m)$-dimensional faces of $\De_n$. Each $\De$ in $S^m_n$
is a simplex $\De_{n-m}$, the convex hull of its $n-m+1$ vertices,
which are a subset of the vertices $p_0,\ldots,p_n$ of $\De_n$, and
any $n-m+1$ out of $p_0,\ldots,p_n$ determine some $\De$ in $S^m_n$.
Hence $\md{S^m_n}=\binom{n+1}{n-m+1}$. Note that it is {\it not\/}
true that $\pd^m\De_n=\coprod_{\De\in S^m_n}\De$, but instead, each
$\De$ in $S^m_n$ appears $m!$ times in $\pd^m\De_n$, indexed by a
choice of order of the $m$ codimension 1 faces of $\De_n$ that meet
at $\De$. We shall describe the compatibilities between the
$\ti\si_{(i,\bs F)}\vert_{\De}$ for all $(i,\bs F)\in C$, $\De\in
S^m_n$ using an equivalence relation $\sim$ on~$C\t S^m_n$.

Let $(i,\bs F,\De)\in C\t S^m_n$. The flag $\bs F$ for $K_i$
determines a simplex $\si_{\bs F}:\De_n\ra K_i$ in the barycentric
subdivision of $K_i$ in Definition \ref{khAdef3}, so $\si_{\bs
F}(\De)$ is one of the $(n-m)$-dimensional simplices in the
associated triangulation of $K_i$. The interior $\si_{\bs
F}(\De^\ci)$ is a subset of the codimension $k$ stratum $S_k(K_i)$
for some unique $k=0,\ldots,m$; in fact $n-k=\max\{j:p_j$ is a
vertex of $\De\}$, so $k$ is determined by $\De$. Let $K_i^{\bs
F,\De}$ be the unique $(n-k)$-dimensional face of $K_i$ containing
$\si_{\bs F}(\De)$. Then $K_i^{\bs F,\De}$ is an $(n-k)$-dimensional
compact convex polyhedron in an affine subspace of $T_{p_i}X$
isomorphic to $\R^{n-k}$, and the barycentric subdivision of $K_i$
restricts to the barycentric subdivision of $K_i^{\bs F,\De}$, with
$\si_{\bs F}(\De)$ one of the simplices of the barycentric
subdivision of~$K_i^{\bs F,\De}$.

Now $\ti X_i$ approximates $K_i$, and this approximation identifies
the face structure of $\ti X_i$ and $K_i$ as manifolds with
g-corners, the connected components of $S_k(\ti X_i)$ and
$S_k(K_i)$, and so on. As $(K_i^{\bs F,\De})^\ci$ is a connected
component of $S_k(K_i)$, there is a unique corresponding connected
component of $S_k(\ti X_i)$. Let $\ti X_i^{\bs F,\De}$ be its
closure. Then $\ti X_i^{\bs F,\De}\subseteq\ti X_i$ is an
$(n-k)$-dimensional face of $\ti X_i$, which approximates the convex
polyhedron $K_i^{\bs F,\De}$. As $\ti X_i\subset\pd Z_{X,\ti T}$,
the natural immersion $\io:\pd Z_{X,\ti T}\ra Z_{X,\ti T}$ restricts
to $\io:\ti X_i\ra Z_{X,\ti T}$, and so to $\io:\ti X_i^{\bs
F,\De}\ra Z_{X,\ti T}$. Then $\io\bigl((\ti X_i^{\bs
F,\De})^\ci\bigr)$ is a connected component of $S_{k+1}(Z_{X,\ti
T})$, and $\io(\ti X_i^{\bs F,\De})$ is an $(n-k)$-dimensional face
of~$Z_{X,\ti T}$.

Let $(i,\bs F),(i',\bs F')$ lie in $C$. We use the same $\De$ and
$k$ for both. Then we get $(n-k)$-dimensional faces $\io(\ti
X_i^{\bs F,\De})$, $\io(\ti X_{i'}^{\bs F',\De})$ of $Z_{X,\ti T}$.
Suppose $\io(\ti X_i^{\bs F,\De})=\io(\ti X_{i'}^{\bs F',\De})$. As
$\io(\ti X_i^{\bs F,\De}),\io(\ti X_{i'}^{\bs F',\De})$ are
approximately identified with convex polyhedra $K_i^{\bs F,\De},
K_{i'}^{\bs F',\De}$, there is an approximate identification between
$K_i^{\bs F,\De}$ and $K_{i'}^{\bs F',\De}$, which identifies the
faces of $K_i^{\bs F,\De}$ with the faces of $K_{i'}^{\bs F',\De}$,
and the simplices of their barycentric subdivisions. Now $\si_{\bs
F}(\De)$ and $\si_{\bs F'}(\De)$ are $(n-m)$-simplices in the
barycentric subdivisions of $K_i^{\bs F,\De}$ and $K_{i'}^{\bs
F',\De}$.

Define a binary relation $\sim$ on $C\t S^m_n$ by $(i,\bs F,\De)\sim
(i',\bs F',\De)$ if $\io(\ti X_i^{\bs F,\De})=\io(\ti X_{i'}^{\bs
F',\De})$, and $\si_{\bs F}(\De)$ and $\si_{\bs F'}(\De)$ correspond
under the identification of simplices of the barycentric
subdivisions of $K_i^{\bs F,\De}$ and $K_{i'}^{\bs F',\De}$. Define
$(i,\bs F,\De)\not\sim (i',\bs F',\De')$ if $\De\ne\De'$. It is easy
to see that $\sim$ is an {\it equivalence relation} on $C\t S^m_n$.
Note that we allow $i=i'$ but $\bs F\ne\bs F'$, and then $\si_{\bs
F}(\De_n)$, $\si_{\bs F'}(\De_n)$ are two $n$-simplices in the
barycentric subdivision of $K_i$ which share an $(n-m)$-dimensional
face~$\si_{\bs F}(\De)=\si_{\bs F'}(\De)$.

By reverse induction on $m=n,n-1,\ldots,0$, we will choose maps
$\ti\si_{(i,\bs F,\De)}:\De\ra\ti X_i^{\bs F,\De}\subseteq\ti X_i$
for all $(i,\bs F)\in C$ and $\De\in S^m_n$ satisfying the
conditions:
\begin{itemize}
\setlength{\itemsep}{0pt}
\setlength{\parsep}{0pt}
\item[(i)] Under the approximate identification of $\ti X_i$
and $K_i$, $\ti\si_{(i,\bs F,\De)}$ approximates $\si_{\bs
F}\vert_\De$.
\item[(ii)] If $\De\in S^m_n$ for $m<n$ then there is a natural
identification $\pd\De\cong\De^0\amalg \cdots\amalg\De^{n-m}$, where
$\De^0,\ldots,\De^{n-m}\in S^{m+1}_n$. Under this identification we
require that $\si_{(i,\bs F,\De)}\vert_{\De^j}\equiv\si_{(i,\bs
F,\De^j)}$ for all $(i,\bs F)\in C$ and $j=0,\ldots,n-m$, where
$\si_{(i,\bs F,\De^j)}$ was chosen in the previous inductive step.
\item[(iii)] If $(i,\bs F),(i',\bs F')\in C$ and $\De\in S^m_n$
with $(i,\bs F,\De)\sim(i',\bs F',\De)$ then $\io\ci\ti\si_{(i, \bs
F,\De)}\equiv\io\ci\ti \si_{(i',\bs F',\De)}:\De\ra Z_{X,\ti T}$.
\end{itemize}
The idea here is that $\ti\si_{(i,\bs F,\De)}=\ti\si_{(i,\bs F)}
\vert_\De$. Then (i),(ii) are obvious, and (iii) is the boundary
compatibilities that the $\ti\si_{(i,\bs F)}$ must satisfy over
$\pd^m\De$. For the initial step $m=n$, each $\De$ in $S^n_n$ is a
point, one of the vertices $p_0,\ldots,p_n$ of $\De_n$, and $\ti
X_i^{\bs F,\De}$ is a point, one of the vertices of $\ti X_i$. So
the maps $\ti\si_{(i,\bs F,\De)}$ are uniquely determined. Parts
(i),(iii) are immediate, and (ii) is vacuous.

For the inductive step, suppose we have chosen $\ti\si_{(i,\bs
F,\De)}$ satisfying (i)--(iii) for all $(i,\bs F)\in C$ and $\De\in
S^m_n$ for $m=n,n-1,\ldots,k+1$, where $0\le k<n$. We shall choose
$\ti\si_{(i,\bs F,\De)}$ satisfying (i)--(iii) for all $(i,\bs F)\in
C$ and $\De\in S^k_n$. Part (ii) determines $\ti\si_{(i,\bs
F,\De)}\vert_{\pd\De}$ uniquely. Part (ii) for $m=k+1$ implies that
these values for $\ti\si_{(i,\bs F,\De)}\vert_{\pd\De}$, when
restricted to $\pd^2\De$, are invariant under the natural involution
$\si:\pd^2\De\ra\pd^2\De$ from \S\ref{kh21}. Since $\De$ is a
manifold with corners (not g-corners), the Extension Principle,
Principle \ref{kh2pri}(c) holds, and it is possible to smoothly
extend the prescribed values for $\ti\si_{(i,\bs F,\De)}$ on
$\pd\De$ smoothly to $\De$, at least near $\pd\De$. As in (i) we
know approximately what $\ti\si_{(i,\bs F,\De)}$ should be, and (i)
when $m=k+1$ implies that the prescribed values for $\ti\si_{(i,\bs
F,\De)}\vert_{\pd\De}$ are approximately right.

Thus there are no global topological problems in extending
$\ti\si_{(i,\bs F,\De)}$ from its prescribed boundary values to a
global embedding satisfying (i),(ii), and for any given $(i,\bs
F,\De)$ in $C\t S^k_n$ we can choose $\ti\si_{(i,\bs F,\De)}$
satisfying (i),(ii). To make our choices satisfy (iii), we choose
$\ti\si_{(i,\bs F,\De)}$ as above for one representative in each
equivalence class of $\sim$ in $C\t S^k_n$, and this choice and
(iii) determine $\ti\si_{(i',\bs F',\De)}$ uniquely for all $(i',\bs
F',\De)$ in $C\t S^k_n$ with $(i',\bs F',\De)\sim(i,\bs F,\De)$.
These choices then satisfy (i)--(iii) for all $(i,\bs F,\De)$ in
$C\t S^k_n$, completing the inductive step. When $m=0$ we have
$S^0_n=\{\De_n\}$. Define embeddings $\ti\si_{(i,\bs
F)}=\ti\si_{(i,\bs F,\De_n)}:\De_n\ra\ti X_i$ for all~$(i,\bs F)\in
C$.

Most of the rest of the proof follows the sketch above with few
complications. Our choice of $\ti\si_{(i,\bs F)}$ determines
$A_{(i,\bs F)}$ on the simplex $\ti\pi\ci\ti\si_{(i,\bs F)}(\De_n)$
in $B_{2\ep}(p_i)$. Our definition of smooth maps between subsets of
$\R^n$ in Definition \ref{kh2def2} implies that $A_{(i,\bs F)}$
extends smoothly to an open neighbourhood of $\ti\pi\ci\ti
\si_{(i,\bs F)}(\De_n)$ in $B_{2\ep}(p_i)$, so there is no problem
in extending $A_{(i,\bs F)}$ from its prescribed values on
$\ti\pi\ci\ti\si_{(i,\bs F)}(\De_n)$ to $B_{2\ep}(p_i)$, and in
getting an approximately affine function. Then provided
$c_0,\ldots,c_n$ and $\de$ are chosen as above, we get a tent
function $T$ with the properties in the theorem, and unique
diffeomorphisms $\si_{(i,\bs F)}:\De_n\ra X_{(i,\bs F)}$ with
$\pi\ci\si_{(i,\bs F)}\equiv\ti\pi\ci\ti\si_{(i,\bs F)}$.

For the final part, writing $c=(i,\bs F)$ and $c'=(i',\bs F')$,
suppose that $\si$ exchanges two $(n-1)$-simplices $\si_{(i,\bs
F)}\ci F_j^n(\De_{n-1})$ and $\si_{(i',\bs F')}\ci
F_{j'}^n(\De_{n-1})$ in \eq{khAeq8}. Applying $\io:\pd^2Z_{X,T}\ra
Z_{X,T}$ gives $\io\ci\si_{(i,\bs F)}\ci
F_j^n(\De_{n-1})=\io\ci\si_{(i',\bs F')}\ci F_{j'}^n(\De_{n-1})$.
Applying $\pi:Z_{X,T}\ra X$ gives $\pi\ci\si_{(i,\bs F)}\ci
F_j^n(\De_{n-1})=\pi\ci\si_{(i',\bs F')}\ci F_{j'}^n(\De_{n-1})$,
since $\pi\ci\io=\pi$. But $\ti\pi\ci\ti\si_{(i,\bs
F)}=\pi\ci\si_{(i,\bs F)}$ and $\ti\pi\ci\ti\si_{(i',\bs
F')}=\pi\ci\si_{(i',\bs F')}$, so $\ti\pi\ci\ti\si_{(i,\bs F)}\ci
F_j^n(\De_{n-1})=\ti\pi\ci\ti\si_{(i',\bs F')}\ci F_{j'}^n
(\De_{n-1})$. This then lifts to $\io\ci\ti\si_{(i,\bs F)}\ci
F_j^n(\De_{n-1})=\io\ci\ti\si_{(i',\bs F')}\ci F_{j'}^n
(\De_{n-1})$.

Now $F_j^n(\De_{n-1}),F_{j'}^n(\De_{n-1})$ lie in $S^1_n$, so
$\ti\si_{(i,\bs F)}=\ti\si_{(i,\bs F,\De_n)}$ and (ii) above for
$m=0$ give $\ti\si_{(i,\bs F)}\vert_{F_j^n
(\De_{n-1})}=\ti\si_{(i,\bs F,F_j^n(\De_{n-1}))}$, and similarly
$\ti\si_{(i',\bs F')}\vert_{F_{j'}^n(\De_{n-1})}=\ti\si_{(i',\bs
F',F_{j'}^n(\De_{n-1}))}$. Therefore $\io\ci\ti\si_{(i,\bs
F,F_j^n(\De_{n-1}))}\bigl(F_j^n(\De_{n-1})\bigr)=\io\ci\ti\si_{(i',\bs
F',F_{j'}^n(\De_{n-1}))}\ab\bigl(F_{j'}^n(\De_{n-1})\bigr)$, as in
part (iii) above. This implies that $j=j'$, since the image of
$\io\ci\ti\si_{(i,\bs F,\De)}$ in $Z_{X,\ti T}$ determines $\De$
uniquely; this is why we put $(i,\bs F,\De)\not\sim (i',\bs
F',\De')$ if $\De\ne\De'$.

By considering the faces of $K_i$, $K_{i'}$ corresponding to
$\ti\si_{(i,\bs F)}\ci F_j^n(\De_{n-1})$ and $\ti\si_{(i',\bs
F')}\ci F_j^n(\De_{n-1})$, we see from the definition of $\sim$ that
$(i,\bs F,F_j^n(\De_{n-1}))\sim(i',\bs F',F_j^n(\De_{n-1}))$, and so
(iii) above with $m=1$ gives $\io\ci\ti\si_{(i,\bs F,F_j^n
(\De_{n-1}))}\equiv\io\ci\ti\si_{(i',\bs F',F_j^n(\De_{n-1}))}$,
yielding $\io\ci\ti\si_{(i,\bs F)}\ci
F_j^n\equiv\io\ci\ti\si_{(i',\bs F')}\ci F_j^n$. Applying $\ti\pi$
gives $\ti\pi\ci\ti\si_{(i,\bs F)}\ci F_j^n\equiv\ti\pi\ci
\ti\si_{(i',\bs F')}\ci F_j^n$, as $\ti\pi\ci\io=\ti\pi$, so
$\pi\ci\si_{(i,\bs F)}\ci F_j^n\equiv\pi\ci\si_{(i',\bs F')}\ci
F_j^n$. Since $\si$ exchanges $\si_{(i,\bs F)}\ci F_j^n(\De_{n-1})$
and $\si_{(i',\bs F')}\ci F_{j'}^n(\De_{n-1})$, and $\pi\ci\si=\pi$,
this lifts to $\si\ci\si_{(i,\bs F)}\ci F_j^n\equiv\si_{(i',\bs
F')}\ci F_j^n$, as we have to prove.
\end{proof}

\subsubsection{Tent functions and piecewise smooth extensions}
\label{khA14}

As in \S\ref{kh21}, the smooth Extension Principle, Principle
\ref{kh2pri}(c), fails for manifolds with g-corners. We will now
introduce a notion of {\it piecewise-smooth\/} function associated
to a tent function, and show that the Extension Principle does hold
for such piecewise smooth functions.

\begin{dfn} Let $X$ be a compact manifold with g-corners, $T$ a tent
function on $X$, and $\pi:X_c\ra X$ for $c\in C$ be as in
Proposition \ref{khAprop1}. A function $f:X\ra\R$ is called {\it
piecewise smooth subordinate to\/} $T$ if $f$ is continuous and
$f\vert_{\pi(X_c)}$ is smooth for each $c\in C$, that is,
$f\ci\pi:X_c\ra\R$ is smooth for each~$c\in C$.

Similarly, if $E\ra X$ is a vector bundle then a section $s$ of $E$
is called {\it piecewise smooth subordinate to\/} $T$ if $s$ is
continuous with $s\vert_{\pi(X_c)}$ smooth for each~$c\!\in\! C$.
\label{khAdef4}
\end{dfn}

We shall define a tent function $T$ on $X$ for which a piecewise
version of the Extension Principle holds.

\begin{dfn} Let $X$ be a compact $n$-manifold with g-corners. Choose a
smooth vector field $v$ on $X$ such that $v$ is nonzero and {\it
inward pointing} at every point of $\pd X$. That is, we require that
whenever $(p,B)\in\pd X$ and $(x_1,\ldots,x_n)$ are coordinates on
$X$ near $p$ such that $x_1\ge 0$ on $X$ and $x_1=0$ on $B$, then
$\d x_1(v\vert_p)>0$. Then for all $t\ge 0$, the exponential map
$\exp(tv):X\ra X$ is a well-defined smooth embedding, which is the
identity for $t=0$, and maps $X$ to a subset of $X^\ci$ for $t>0$.
Essentially, the flow $\exp(tv)$ for $t\ge 0$ {\it shrinks $X$ away
from its boundary}.

Write $X_1=\exp(v)X$ and $X_2=\exp(2v)X$. Then $X_1,X_2$ are
embedded submanifolds of $X^\ci$, and $\exp(v):X\ra X_1$,
$\exp(2v):X\ra X_2$ are diffeomorphisms. The inclusion $\io:\pd
X_2\ra X$ is an immersion. We wish to extend $\pd X_2$ to a
$(n-1)$-manifold $Y$ without boundary, not necessarily compact, with
$\pd X_2\subset Y$ an embedded $(n-1)$-submanifold, and extend
$\io:\pd X_2\ra X$ to an immersion $\io:Y\ra X$ with the properties:
\begin{itemize}
\setlength{\itemsep}{0pt}
\setlength{\parsep}{0pt}
\item[(i)] $\io(Y\sm\pd X_2)\subset X_1\sm X_2$.
\item[(ii)] $v$ is transverse to $\io(Y)$ at every point of $Y$.
\item[(iii)] there exists $M\ge 1$ such that for each $x\in X$
there are at most $M$ points $(t,y)\in[0,2]\t Y$
with~$\exp(tv)x=\io(y)$.
\item[(iv)] $S_l=\bigl\{(y_1,\ldots,y_l)\in Y^l:\io(y_1)=\cdots=
\io(y_l)\bigr\}$ is a disjoint union of embedded submanifolds of
$Y^l$, for each~$l\ge 1$.
\item[(v)] In (iv), the map $i_l:S_l\ra X$,
$i_l:(y_1,\ldots,y_l)\mapsto\io(y_1)$ is an immersion, and $\d
i_l\bigl(T_{(y_1,\ldots,y_l)}S_l\bigr)=\bigcap_{i=1}^l\d\io(T_{y_i}Y)$
in $T_{\io(y_1)}X$ for all $(y_1,\ldots,y_l)\in S_l$.
\end{itemize}

The basic idea is that we are adding a small `collar' to $\pd X_2$
near $\pd^2X_2$. If $\pd^3X_2=\emptyset$ then we can think of $Y$ as
$\pd X_2\amalg (0,\ep)\t\pd^2X_2$ for $\ep>0$ small. If we choose
any $Y,\io$ extending $\pd X_2,\io$ then $\io(Y\sm\pd X_2)\subset
X_1\sm X_2$ near $\pd^2X_2$, so making $Y$ smaller if necessary (i)
holds. Since $v$ is $\exp(2v)$-invariant and $v$ is nonzero and
inward-pointing at each point of $\pd X$, it follows that $v$ is
transverse to $\io(\pd X_2)$ at every point of $\pd X_2$. This is an
open condition, so for any $Y,\io$ extending $\pd X_2,\io$, $v$ is
transverse to $\io(Y)$ on an open neighbourhood of $\pd X_2$ in $Y$.
Thus making $Y$ smaller if necessary, (ii) holds.

For (iii), using (ii) we find that the set of $(t,y)\in[0,2]\t Y$
with $\exp(tv)x=\io(y)$ is discrete, that is, it has no limit
points. So if $Z$ is any compact subset of $Y$, then there are only
finitely many $(t,z)\in[0,2]\t Z$ with $\exp(tv)x=\io(z)$. The
number of such $(t,z)$ is an upper-semi-continuous function of $x\in
X$, and so bounded by $M\ge 1$ as $X$ is compact. Thus, given any
$Y,\io$ satisfying (ii), replacing $Y$ by an open neighbourhood $Y'$
of $\pd X_2$ in $Y$ with $Z=\overline{Y'}$ compact in $Y$, there
exists $M\ge 1$ with at most $M$ pairs $(t,y)\in[0,2]\t Y'$ with
$\exp(tv)x=\io(y)$ for any $x\in X$. So making $Y$ smaller if
necessary, (iii) holds.

To arrange for (iv) and (v) to hold, let $x\in\io(\pd X_2)\subset
X$, and choose a chart $(U,\phi)$ for $X$ near $x$, so that
$U\subset\R^n$ is open and $\phi:U\ra X$ is a diffeomorphism with
its image, with $x\in\phi(U)$. As $X_2$ is an embedded submanifold
of $X$, $\phi^{-1}(X_2)$ is a {\it region with g-corners} in $\R^n$.
Thus there exist smooth $f_1,\ldots,f_N:U\ra\R$ satisfying
Definition \ref{kh2def1}(a),(b), such that
$\phi^{-1}(X_2)=\bigl\{u\in U:f_i(u)\ge 0$, $i=1,\ldots,N\bigr\}$.
Choosing such $f_1,\ldots,f_N$ with $N$ least, and making $U$
smaller if necessary, we see that $\pd X_2$ is locally modelled on
$\bigl\{(u,i):u\in U$, $i=1,\ldots,N$, $f_i(u)=0$, $f_j(u)\ge 0$,
$j=1,\ldots,N\bigr\}$ over $\phi(U)$. We choose $Y$ to be locally
modelled on $\bigl\{(u,i):u\in U$, $i=1,\ldots,N$, $f_i(u)=0\bigr\}$
over $\phi(U)$. Then Definition \ref{kh2def1}(a),(b) for
$f_1,\ldots,f_N$ imply (iv),(v) above. Hence we can choose $Y,\io$
satisfying (i)--(v) near any point of $\io(\pd X_2)$. Patching
together such local choices, we can define $Y,\io$ globally. Now
define $T:X\ra F\bigl([1,\iy)\bigr)$ by
\e
\ts T(x)=\bigl\{2\bigr\}\cup\bigl\{3-t:t\in(0,2],\;\>
\exp(tv)x\in\io(Y)\bigr\}.
\label{khAeq10}
\e
\label{khAdef5}
\end{dfn}

\begin{prop} In Definition {\rm\ref{khAdef5},} $T$ is a tent
function on $X$. Furthermore, if\/ $\pi:X_c\ra X$ and\/ $C$ are as
in Proposition {\rm\ref{khAprop1}} for $T,$ then there is a unique
diffeomorphism $\Phi:X_1\amalg [0,1]\t\pd X\ra\coprod_{c\in C}X_c$
such that\/ $\pi\ci\Phi:X_1\amalg [0,1]\t\pd X\ra X$ maps $x\mapsto
x$ for $x\in X_1$ and\/ $\bigl(t,(x,B)\bigr)\mapsto
\exp(tv)x$ for $\bigl(t,(x,B)\bigr)\in [0,1]\t\pd X$. If\/ $X$ is
oriented then $\Phi$ is orientation-preserving on $X_1$ and
orientation-reversing on $[0,1]\!\t\!\pd X,$ so we can regard\/
$\Phi$ as an orientation-preserving diffeomorphism
$\Phi:X_1\!\amalg\! -[0,1]\!\t\!\pd X\!\ra\!\coprod_{c\in C}X_c$.
Also~$\min T\vert_{\pd X}\!\equiv\! 1$.
\label{khAprop3}
\end{prop}

\begin{proof} Define $Z=\bigl\{(x,t,y)\in X\t(0,3)\t Y:\exp(tv)x
=\io(y)\bigr\}$. Regarding $Z$ as a fibre product
$\bigl(X\t(0,3)\bigr)\t_XY$ and noting that $(x,t)\mapsto\exp(tv)x$
is a submersion $X\t(0,3)\ra X$, we see that $Z$ is an $n$-manifold.
Using Definition \ref{khAdef5}(ii) and the fact that $(0,3)\t Y$ has
no boundary, we can show that the projection $\pi_X:Z\ra X$,
$\pi_X:(x,t,y)\mapsto x$ is a {\it local diffeomorphism}. Definition
\ref{khAdef5}(ii) implies that $\io(Y)\cap X_2^\ci=\emptyset$. But
if $x\in X$ and $t>2$ then $\exp(tv)x\in X_2^\ci$, so there exists
no $y\in Y$ with $\exp(tv)x=\io(y)$. Hence $Z\subseteq X\t(0,2]\t
Y$. Therefore Definition \ref{khAdef5}(iii) implies that
$\bmd{\pi_X^{-1}(x)}\le M$ for all $x\in X$, and $\pi_X$ is a {\it
globally finite} map. Note too that~$T(x)=\{2\}\cup\{3-t:(x,t,y)\in
Z\}$.

Let $x\in X$. Since $\pi_X:Z\ra X$ is a globally finite local
diffeomorphism, if $U$ is a sufficiently small open neighbourhood of
$x$ in $X$, then we may split $\pi_X^{-1}(U)=W_2\amalg\cdots\amalg
W_N$ for some open subsets $W_2,\ldots,W_N$ of $Z$ for which
$\pi_X\vert_{W_i}:W_i\ra X$ is injective. Set $U_1=U$ and
$U_i=\pi_X(W_i)$ for $i=2,\ldots,N$. Then $U_i$ are open subsets of
$U$, and $\pi_X\vert_{W_i}:W_i\ra U_i$ is a diffeomorphism for
$i=2,\ldots,N$. Define $t_i:U_i\ra[1,\iy)$ for $i=1,\ldots,N$ by
$t_1\equiv 2$ and $t_i(u)=3-t$ when $(u,t,y)\in W_i$ for
$i=2,\ldots,N$. Since $\pi_X\vert_{W_i}:W_i\ra U_i$ is a
diffeomorphism there is a unique $(u,t,y)\in W_i$ for each $u\in
U_i$, which depends smoothly on $u$, so $t_i$ is well-defined and
smooth.

We now see that $T(u)=\bigl\{t_i(u):i=1,\ldots,N$, $u\in U_i\bigr\}$
for all $u\in U$, as in Definition \ref{khAdef1}. We must verify
Definition \ref{khAdef1}(a)--(c). For subsets $\{i_1,\ldots,i_l\}$
of $\{1,\ldots,N\}$ with $1\notin\{i_1,\ldots,i_l\}$, Definition
\ref{khAdef1}(a),(b) follow directly from Definition
\ref{khAdef5}(iv),(v), with $S_{\{i_1,\ldots,i_l\}}\subseteq
\bigl\{x\in X:\exp(tv)x\in i_l(S_l)$, $t\in(0,3)\bigr\}$. When
$1\in\{i_1,\ldots,i_l\}$, Definition \ref{khAdef1}(a),(b) follow
from Definition \ref{khAdef5}(iv),(v), with
$S_{\{i_1,\ldots,i_l\}}\subseteq \bigl\{x\in X:\exp(v)x\in
i_{l-1}(S_{l-1})\bigr\}$.

To prove Definition \ref{khAdef1}(c), we shall show using
\eq{khAeq10} that
\e
\min T(x)=\begin{cases} 2, & x\in X_1, \\
1+\inf\bigl\{t\in[0,1]:x\in\exp(tv)X\bigr\}, & x\in X\sm X_1.
\end{cases}
\label{khAeq11}
\e
For the first line of \eq{khAeq11}, note that if $x\in X_1$ and
$t\in(1,2]$ then $\exp(tv)x\in X_2^\ci$. But Definition
\ref{khAdef5}(i) implies that $\io(Y)\cap X_2^\ci=\emptyset$. So
$T(x)$ cannot contain $3-t$ for $t\in(1,2]$, implying that $\min
T(x)=2$. For the second line, let $x\in X\sm X_1$, and set
$s=\inf\bigl\{t\in[0,1]:x\in\exp(tv)X\bigr\}$. (This set contains 1,
so the infimum is well-defined.) By compactness of $X$ we see that
there exists $x'\in X$ with $\exp(sv)x'=x$. As $s$ is smallest we
must have $x'\in\io(\pd X)$, so we can lift to $(x',B')\in\pd X$. We
have a diffeomorphism $\exp(2v):\pd X\ra \pd X_2\subset Y$. Define
$y=\exp(2v)(x',B')$, and set $t=2-s$. Then
$\io(y)=\exp(2v)x'=\exp((t+s)v)x'=\exp(tv)\exp(sv)x'=\exp(tv)x$.
Hence $3-t=1+s\in T(x)$, and $\min T\le 1+s$. But using Definition
\ref{khAdef5}(i) we can show that $\min T\ge 1+s$,
proving~\eq{khAeq11}.

Equation \eq{khAeq11} implies that $\min T(x)$ occurs either at the
extra branch 2, or at $3-t$ when $\exp(tv)x=\io(y)$ for $y\in\pd
X_2$. The minimum can never occur at $3-t$ when $\exp(tv)x=\io(y)$
for $y\in Y\sm\pd X_2$. Definition \ref{khAdef1}(c) now follows
using compactness of $\pd X_2$. Therefore $T$ is a {\it tent
function} on $X$. The remainder of the proposition is a consequence
of \eq{khAeq11} and its proof. We showed that either $x\in X_1$ and
$\min T(x)=2$, or $\min T(x)=1+t$ for $t\in[0,1]$ and there exists
$(x',B')\in\pd X$ with $\exp(tv)x'=x$. These possibilities
correspond to the fibres of $\pi:\coprod_{c\in C}X_c\ra X$ over $x$,
and show that we can write $\coprod_{c\in C}X_c\cong X_1\amalg
[0,1]\t\pd X$ as in the proposition.
\end{proof}

Here is our first version of Principle \ref{kh2pri}(c) for piecewise
smooth functions. The proof is easy, we just have to check
\eq{khAeq12} is well-defined.

\begin{prop} Let\/ $X$ be a compact manifold with g-corners, and
suppose $\ze:\pd X\ra\R$ is smooth, with\/ $\ze\vert_{\pd^2X}$
invariant under the natural involution $\si:\pd^2X\ra\pd^2X$ of\/
\S{\rm\ref{kh21}}. Let\/ $v,T,X_1$ be as in Definition
{\rm\ref{khAdef5},} and\/ $\Phi$ be as in Proposition
{\rm\ref{khAprop3}}. Choose any smooth function $\eta_1:X_1\ra\R$.
Then there is a piecewise smooth function $\eta:X\ra\R$ subordinate
to the tent function $T$ with\/ $\eta\vert_{\pd X}\equiv\ze,$ given
by
\e
\eta(x)=\begin{cases} \eta_1(x), & x\in X_1, \\
(1-t)\ze(x',B')+t\eta_1(\exp(v)x'), &
\begin{subarray}{l} \ts x=\pi_X\ci\Phi\bigl(t,(x',B')\bigr), \\
\ts\bigl(t,(x',B')\bigr)\in[0,1]\t\pd X.\end{subarray}
\end{cases}
\label{khAeq12}
\e

Similarly, if\/ $E\ra X$ is a vector bundle and\/ $f$ is a smooth
section of\/ $E\vert_{\pd X}$ with\/ $f\vert_{\pd^2X}$
$\si$-invariant, then by choosing a smooth section $e_1$ of\/
$E\vert_{X_1}$ we can write down a piecewise smooth section $e$ of\/
$E$ subordinate to $T$ with\/~$e\vert_{\pd X}\equiv f$.
\label{khAprop4}
\end{prop}

Proposition \ref{khAprop4} is not adequate for our later
applications as we cannot apply it iteratively, by induction on
dimension, since the inputs $\ze,f$ are smooth but the outputs
$\eta,e$ are piecewise smooth. We shall generalize it to a
construction in which both input and output are piecewise smooth.
First we generalize Definition \ref{khAdef5} to extend a tent
function $T'$ on $\pd X$ to a tent function $T$ on~$X$.

\begin{dfn} Let $X$ be a compact $n$-manifold with g-corners, and let
$T':\pd X\ra F\bigl([1,\iy)\bigr)$ be a tent function on $\pd X$,
such that $\min T'\vert_{\pd^2X}$ is invariant under the natural
involution $\si:\pd^2X\ra\pd^2X$. We will construct a tent function
$T:X\ra F\bigl([1,\iy)\bigr)$ with $T\vert_{\pd X}\equiv T'$, with
useful properties we will explain in Propositions \ref{khAprop5}
and~\ref{khAprop6}.

Let $v,X_1,X_2$ and $\io:Y\ra X$ be as in Definition \ref{khAdef5}.
Then $\exp(2v)$ induces a diffeomorphism $\exp(2v):\pd X\ra\pd X_2$,
so that $T'\ci\exp(2v)$ is a tent function on $\pd X_2$, and $\pd
X_2\subset Y$. Extend $T'\ci\exp(2v)$ to a tent function $T''$ on
$Y$. Making $Y$ smaller if necessary, this is possible, except that
we may need to relax the condition that $\min T''\ge 1$. So we take
$T'':Y\ra F\bigl([\ha,\iy)\bigr)$ to satisfy Definition
\ref{khAdef1} replacing $[1,\iy)$ by $[\ha,\iy)$, with
$T''\vert_{\pd X_2}\equiv T'\ci\exp(2v)$.

As $\min T':\pd X\ra[1,\iy)$ is continuous and $\pd X$ is compact,
$\min T'$ is bounded above. Choose $D>1$ such that $\min T'<D$ on
$\pd X$. So $\min T''<D$ on $\pd X_2$, as $T''\vert_{\pd X_2}\equiv
T'\ci\exp(2v)$. Making $Y$ smaller if necessary, we can also assume
that $\min T''<D$ on $Y$. Define $T:X\ra F\bigl([1,\iy)\bigr)$ by
\e
\ts T(x)\!=\!\bigl\{D\bigr\}\cup\bigl\{u+\th(2-t):
\text{$t\!\in\!(0,2]$, $y\!\in\! Y$, $u\!\in\! T''(y)$,
$\exp(tv)x=\io(y)$}\bigr\},
\label{khAeq13}
\e
where $\th>D-1$ is chosen large enough to satisfy conditions in
Proposition~\ref{khAprop5}.
\label{khAdef6}
\end{dfn}

\begin{prop} In Definition {\rm\ref{khAdef6},} making $Y$ smaller
if necessary, if\/ $\th>D-1$ is sufficiently large then\/ $T$ is a
tent function on $X,$ and\/ $\min T\vert_{\pd X}\equiv\min T',$ so
that\/ $Z_{\pd X,T\vert_{\pd X}}=Z_{\pd X,T'},$ and equation\/
\eq{khAeq5} becomes
\e
\pd Z_{X,T}=-Z_{\pd X,T'}\amalg -\{0\}\t X \amalg\ts\coprod_{c\in
C}X_c.
\label{khAeq14}
\e

Let\/ $X'_{c'}$ for $c'\in C'$ be as in \eq{khAeq5} for $Z_{\pd
X,T'}$, and\/ $X_c$ for $c\in C$ be as in \eq{khAeq5} for $Z_{X,T}$.
Define $\ti X=\bigl\{x\in X:\min T(x)=D\bigr\}$. Then $\ti X$ is a
compact, embedded submanifold of $X$ with g-corners, with
$X_1\subset \ti X\subset X^\ci$. There are unique diffeomorphisms
\e
\ts\Phi:\ti X\amalg\bigl([0,1]\t\coprod_{c'\in
C'}X'_{c'}\bigr)\longra\coprod _{c\in C}X_c,\quad
\Psi:\coprod_{c'\in C'}X'_{c'}\longra\pd\ti X
\label{khAeq15}
\e
which satisfy $\pi_X\ci\Phi(x)=x$ for $x\in\ti X,$ and for
$t\in[0,1],$ $x'\in\coprod_{c'\in C'}X'_{c'},$
\ea
\pi_X\ci\Phi:(t,x')&\longmapsto
\exp\Bigl(\frac{t}{\th(D-\min T'\ci\pi_{\pd X}(x'))}\,v\Bigr)\,\pi_X(x'),
\label{khAeq16}\\
\pi_X\ci\Psi:x'&\longmapsto
\exp\Bigl(\frac{1}{\th(D-\min T'\ci\pi_{\pd X}(x'))}\,v\Bigr)\,\pi_X(x'),
\label{khAeq17}\\
\text{and\/}\quad \min T\ci\pi_X\ci\Phi(t,x')&=(1-t)\min
T'\ci\pi_{\pd X}(x')+t D,
\label{khAeq18}
\ea
writing $\pi_{\pd X}:\coprod_{c'\in C'}X'_{c'}\ra\pd X,$ $\pi_X:
\coprod_{c'\in C'}X'_{c'}$ or $\coprod_{c\in C}X_c\ra X$ for the
natural projections. If\/ $X$ is oriented then $\Phi$ is
orientation-preserving on $\ti X$ and orientation-reversing on
$[0,1]\!\t\!\coprod_{c'\in C'}X'_{c'},$ and\/ $\Psi$ is
orientation-preserving.
\label{khAprop5}
\end{prop}

\begin{proof} Our first goal is to prove that making $Y$ smaller if
necessary and taking $\th>D-1$ sufficiently large, if $x\in X$ and
$\min T(x)= u+\th(2-t)$ for $t\in(0,2]$, $y\in Y$ and $u\in T''(y)$
with $\exp(tv)x=\io(y)$, then $y\in\pd X_2$. That is, the choice of
`collar' $Y\sm\pd X_2$ and the values of $T''$ upon it do not affect
$\min T$, since the branches of $T$ coming from $T''$ on $Y\sm\pd
X_2$ are never minimal.

Let $y\in Y\sm\pd X_2$. Then there exists unique $s\in(0,1)$ with
$\exp(sv)\io(y)\in\io(\pd X_2)$. Thus $\exp(sv)\io(y)=\io(z)$ for
some (not necessarily unique) $z\in\pd X_2$. We claim that if $Y$ is
small enough and $\th$ is large enough, then $\min T''(y)+\th s>\min
T''(z)$ for all such $y,s,z$. This implies what we want, for if
$x\in X$ and $t\in(0,2]$ with $\exp(tv)x=\io(y)$, then
$\exp((s+t)v)x=\io(z)$, so $\min T''(z)+\th(2-s-t)\in T(x)$, and
$\min T''(z)+\th(2-s-t)<\min T''(y)+\th(2-t)\le u+\th(2-t)$ for
all~$u\in T''(y)$.

To prove the claim, choose a Riemannian metric $g$ on $X$, let $g_Y$
be the pullback metric on $Y$, and $d_Y$ the induced metric on $Y$.
As $\min T''$ is continuous and locally the minimum of
differentiable functions, it is Lipschitz. Thus by compactness of
$\pd X_2$, making $Y$ smaller if necessary, there exists $E\ge 0$
such that $\bmd{\min T''(y_1)-\min T''(y_2)}\le E d_Y(y_1,y_2)$ for
all $y_1,y_2\in Y$. Now let $y,s,z$ be as above. Then $y,z$ are on
different sheets of $Y$ in $X$ locally, but these sheets should
intersect in $X$, so we can choose $y',z'\in\pd X_2$ such that $y'$
is close to $y$ in $Y$, and $z'$ is close to $z$ in $Y$,
and~$\io(y')=\io(z')$.

Since $\min T'\vert_{\pd^2X}$ is $\si$-invariant, it follows that
$\min T'$ pushes down to a function on $\io(\pd X)\subset X$. Hence
$\min T''\vert_{\pd X_2}$ pushes down to a function on $\io(\pd
X_2)\subset X$, and $\io(y')=\io(z')$ implies $\min T''(y')=\min
T''(z')$. Think of $P=\io(y)$, $Q=\io(z)$ and $R=\io(y')=\io(z')$ as
the vertices of a small triangle in $X$. The side $PR$ has `length'
$d_Y(y,y')$. The side $QR$ has `length' $d_Y(z,z')$. The side $PQ$
has `length' at most $s\Vert v\Vert_{C^0}$, where $\Vert
v\Vert_{C^0}$ is computed using $g$. The angle in $PQR$ at $R$ is
the angle at $\io(y')$ between the two sheets of $\io(Y)$ containing
$\io(y),\io(z)$. By compactness of $\pd^2X_2$, and using Definition
\ref{khAdef5}(iv,)(v), there is a global, positive lower bound for
all angles of different sheets of $Y$ meeting at $\io(\pd^2X_2)$.
Hence there exists small $\ep>0$ independent of $y,s,z$ such that
$\md{PQ}\ge\ep\bigl(\md{PR}+\md{QR}\bigr)$. So we have
\begin{align*}
s\Vert v\Vert_{C^0}&\ge \ep\bigl(d_Y(y,y')+d_Y(z,z')\bigr) \\
&\ge\ep E^{-1}\bigl(\bmd{\min T''(y)-\min T''(y')}+ \bmd{\min
T''(z)-\min T''(z')}\bigr)\\
&\ge \ep E^{-1}\bigl(\bmd{\min T''(y)-\min T''(z)}\bigr),
\end{align*}
since $\min T''(y')=\min T''(z')$. Hence if $\th>\Vert v\Vert_{C^0}
E\ep^{-1}$ then $\min T''(y)+\th s>\min T''(z)$ for all such
$y,s,z$, as we have to prove.

One can now show that $T$ is a {\it tent function} by a modification
of the proof of Proposition \ref{khAprop3}, using the fact that
$T''$ is a tent function. We can then give an expression for $\min
T$: if $(x,B)\in\pd X$ and $s\in[0,2)$ then
\e
\min T\bigl(\exp(s v)x\bigr)=\min\bigl(\min T'(x,B)+\th s,D\bigr).
\label{khAeq19}
\e
This holds as putting $t=2-s$, we have $t\in(0,2]$ and
$\exp(tv)\exp(s v)x=\exp(2v)x=\io\bigl(\exp(2v)(x,B)\bigr)$, where
$\exp(2v)(x,B)\in\pd X_2$, and $T''\bigl(\exp(2v)(x,\ab B)\bigr)=
T'(x,B)$. Hence by \eq{khAeq13}, $T\bigl(\exp(s v)x\bigr)$ is the
union of $\{D\}$, and $\{u+\th s:u\in T'(x,B)\bigr\}$, and
contributions from $y\in Y\sm\pd X_2$ which from above cannot be
minimum in $T\bigl(\exp(s v)x\bigr)$. Equation \eq{khAeq19} follows.
When $s=0$, it implies that $\min T\vert_{\pd X}\equiv\min T'$, as
we want.

Since $\min T'\ge 1$ as $T'$ is a tent function, and $\th>D-1$, we
see from \eq{khAeq19} that $\min T\bigl(\exp(s v)x\bigr)=D$ if
$(x,B)\in\pd X$ and $s\in[1,2)$. Hence $\min T\equiv D$ on $X_1\sm
X_2$. But \eq{khAeq13} gives $T(x)=\{D\}$ for $x\in X_2$, so $\min
T\equiv D$ on $X_2$, and so on $X_1$. Thus $X_1\subset\ti X$. Since
$\min T'<D$ on $\pd X$ and $\min T\vert_{\pd X}\equiv\min T'$ we
have $\min T<D$ on $\io(\pd X_1)$. Therefore $X_1\subset\ti X\subset
X^\ci$, as we have to prove.

We can now describe the pieces $X_c$ for $c\in C$ in \eq{khAeq5} for
$Z_{X,T}$. Each $X_c$ for $c\in C$ is the graph of a smooth branch
of $T$ on a connected component of the subset where it minimizes
$T$. Therefore $\{D\}\t\ti X$ is a union of some of the $X_c$ in
$C$, those associated to the constant branch $x\mapsto D$ of $T$.
Hence $\ti X$ is a compact, embedded submanifold of $X$ with
g-corners. We define $\Phi$ on $\ti X$ by $\Phi:x\mapsto(D,x)$, and
this is a diffeomorphism of $\ti X$ with these $X_c$ with
$\pi_X\ci\Phi(x)=x$ for $x\in\ti X$, as we want.

For the remaining $X_c$ for $c\in C$ on which $\min T\not\equiv D$,
we can see from \eq{khAeq19} that these correspond to $X'_{c'}$ for
$c'\in C'$, since the smooth branches of $T$ on their minimal
subsets in $X$ correspond to the smooth branches of $T'$ on their
minimal subsets in $\pd X$. In fact, if $c'\in C'$ then
\e
X_c=\bigl\{\bigl(u+\th s,\exp(s v)x\bigr):\bigl(u,(x,B)\bigr)\in
X'_{c'},\;\> 0\le \th s\le D-u\bigr\}
\label{khAeq20}
\e
is one of the $X_c$. We define a diffeomorphism $\Phi:[0,1]\t
X'_{c'}\ra X_c$ by
\begin{equation*}
\Phi:\bigl(t,(x,B)\bigr)\longmapsto\Bigl((1-t)\min T'(x,B)+ tD,
\exp\Bigl(\frac{t}{\th(D-\min T'(x,B))}\,v\Bigr)\,x\Bigr).
\end{equation*}
Equations \eq{khAeq16} and \eq{khAeq17} follow. For $\Psi$ in
\eq{khAeq15} and \eq{khAeq17}, we note that $\Phi\vert_{\{1\}\t
\coprod_{c'\in C'}X'_{c'}}$ maps $\coprod_{c'\in C'}X'_{c'}$ to the
faces of $X_c$ in \eq{khAeq20} on which $\min T\equiv D$, and these
correspond to faces of $\ti X$. The remainder of the proposition
follows as for Proposition~\ref{khAprop3}.
\end{proof}

Here is a generalization of Proposition \ref{khAprop4}. The proof is
easy.

\begin{prop} Suppose $X$ is a compact manifold with g-corners, and\/
$T':\pd X\ra F\bigl([1,\iy)\bigr)$ is a tent function, and\/
$\ze:\pd X\ra\R$ is a piecewise smooth function subordinate to $T',$
such that\/ $\min T'\vert_{\pd^2X}$ and\/ $\ze\vert_{\pd^2X}$ are
invariant under the natural involution $\si:\pd^2X\ra\pd^2X$.

Let\/ $T$ be as in Definition {\rm\ref{khAdef6},} and\/ $\ti
X,\Phi,\Psi$ be as in Proposition {\rm\ref{khAprop5}}. Choose any
smooth function $\ti\eta:\ti X\ra\R$. Then there is a piecewise
smooth function $\eta:X\ra\R$ subordinate to the tent function $T$
with\/ $\eta\vert_{\pd X}\equiv\ze,$ given by
\e
\eta(x)\!=\!\begin{cases} \ti\eta(x), & x\in\ti X, \\
(1\!-\!t)\ze\!\ci\!\pi'(x')\!+\!t\ti\eta\vert_{\pd\ti
X}\!\ci\!\Psi(x'), &
\begin{subarray}{l} \ts x=\pi_X\ci\Phi(t,x'), \\
\ts(t,x')\!\in\![0,1]\!\t\!(\coprod_{c'\in
C'}\!X'_{c'}).\!\!\end{subarray}
\end{cases}
\label{khAeq21}
\e

Similarly, if\/ $E\ra X$ is a vector bundle and\/ $f$ is a piecewise
smooth section of\/ $E\vert_{\pd X}$ subordinate to $T'$ with\/
$\min T'\vert_{\pd^2X}$ and\/ $f\vert_{\pd^2X}$ $\si$-invariant,
then by choosing a smooth section $\ti e$ of\/ $E\vert_{\ti X}$ we
can write down a piecewise smooth section $e$ of\/ $E$ subordinate
to $T$ with\/~$e\vert_{\pd X}\equiv f$.
\label{khAprop6}
\end{prop}

\subsubsection[Cutting into small pieces with boundary
conditions]{Cutting a manifold into small pieces with boundary
conditions}
\label{khA15}

We now combine the ideas of \S\ref{khA12} and \S\ref{khA14}. Given a
compact manifold $X$ and a tent function $T'$ on $\pd X$ with
$T'\vert_{\pd^2X}$ $\si$-invariant, we construct a tent function $T$
on $X$ with $\min T\vert_{\pd X}\equiv\min T'$, which cuts $X$ into
`arbitrarily small pieces' $X_c$ for $c\in C$, as far as this is
possible given the boundary conditions.

\begin{dfn} Let $X$ be a compact $n$-manifold with g-corners and
$T':\pd X\ra F\bigl([1,\iy)\bigr)$ a tent function such that $\min
T'\vert_{\pd^2X}$ is invariant under the involution
$\si:\pd^2X\ra\pd^2X$. Let $v,\io:Y\ra X,T'',D,\th$ be as in
Definition~\ref{khAdef6}.

Choose a Riemannian metric $g$ on $X$. Note that $g$ is not required
to take any particular values on $\pd X$ --- it is unrelated to any
metric used to construct $T'$. Choose $\ep>0$ with $2\ep<\de(g)$,
the injectivity radius of $g$. Write $d(\,,\,)$ for the metric on
$X$ determined by $g$. For $p\in X$ and $r>0$ write $B_r(p)$ for the
open ball of radius $r$ about $p$. Then $x\mapsto d(p,x)^2$ is
smooth on $B_{2\ep}(p)$, as~$2\ep<\de(g)$.

The open sets $B_\ep(p)$ for $p\in X^\ci$ cover $X$, so as $X$ is
compact we can choose a finite subset $\{p_1,\ldots,p_N\}$ of
$X^\ci$ with $X=\bigcup_{i=1}^NB_\ep(p_i)$. As in Definition
\ref{khAdef2}, choose $p_1,\ldots,p_N$ to be {\it generic\/} amongst
all such $N$-tuples $p_1',\ldots,p_N'$. Since $p_i\in X^\ci$ and
$\pd X$ is compact, there exists $\ze\in(0,\ep)$ such that
$d(p_i,x)\ge\ze$ for all $i=1,\ldots,N$ and $(x,B)\in\pd X$. Make
$\th$ larger if necessary so that $\th>D+\ep^2-1$ and $1+\th\ze\Vert
v\Vert_{C^0}>D$, where $\Vert v\Vert_{C^0}$ is computed using $g$.
Define $T:X\ra F\bigl([1,\iy)\bigr)$ by
\e
\begin{split}
\ts T(x)=\,&\bigl\{D+d(p_i,x)^2:i=1,\ldots,N,\;\>
d(p_i,x)<2\ep\bigr\}\cup\\
&\bigl\{u+\th(2-t):\text{$t\in(0,2]$, $y\in Y$, $u\in T''(y)$,
$\exp(tv)x=\io(y)$}\bigr\}.
\label{khAeq22}
\end{split}
\e
\label{khAdef7}
\end{dfn}

\begin{prop} In Definition {\rm\ref{khAdef7},} $T$ is a tent
function on\/ $X$, and\/ $\min T'\ab\equiv\min T\vert_{\pd X}$.
Let\/ $X'_{c'}$ for $c'\in C'$ be as in \eq{khAeq5} for $Z_{\pd
X,T'},$ and\/ $X_c$ for $c\in C$ be as in \eq{khAeq5} for $Z_{X,T}$.
Then we can take $C=C'\amalg\{1,\ldots,N\},$ where for $c'\in C'$
there is a natural identification between $X'_{c'}$ and a component
of\/ $\pd X_{c'},$ and for $i=1,\ldots,N$ we
have~$p_i\in\pi(X_i^\ci)\subseteq\pi(X_i)\subseteq B_\ep(p_i)$.

We can choose\/ $T$ so that the pieces\/ $X_c$ are `arbitrarily
small', given the $X'_{c'}$ which are already fixed. That is, if\/
$\{V_i:i\in I\}$ is any open cover for\/ $X$ such that for all\/
$c'\in C'$ there exists\/ $i\in I$ with\/ $\pi'(X'_{c'})\subseteq\pd
V_i,$ then we can choose\/ $T$ so that for all\/ $c\in C$ there
exists\/ $i\in I$ with\/~$\pi(X_c)\subseteq V_i$.
\label{khAprop7}
\end{prop}

\begin{proof} The proof is a straightforward generalization of those
of Propositions \ref{khAprop2} and \ref{khAprop5}. The first and
second lines of \eq{khAeq22} are analogues of \eq{khAeq6} and
\eq{khAeq13} respectively. Since $X=\bigcup_{i=1}^NB_\ep(p_i)$, the
minimum of the terms $D+d(p_i,x)^2$ in the first line of
\eq{khAeq22} is a continuous function on $X$ with values between $D$
and $D+\ep^2$. The condition $\th>D+\ep^2-1$ strengthening $\th>D-1$
in Definition \ref{khAdef6} ensures that in the second line of
\eq{khAeq22}, $u+\th(2-t)$ is never minimal in $T(x)$ for
$t\in[1,2)$, which helps to ensure Definition \ref{khAdef1}(c) holds
and $\min T$ is continuous. It follows as in Propositions
\ref{khAprop2} and \ref{khAprop5} that $T$ is a {\it tent function},
and~$\min T'\ab\equiv\min T\vert_{\pd X}$.

For $i=1,\ldots,N$, we claim that $\min T(p_i)=D$, and
$D+d(p_i,x)^2$ is the only branch of $T$ achieves this minimum. For
$j\ne i$ we have $d(p_j,p_i)>0$, so $D+d(p_i,x)^2>D$ at $x=p_i$.
Suppose $t\in(0,2]$, $y\in Y$, $u\in T''(y)$ and
$\exp(tv)x=\io(y)=p_i$ such that $u+\th(2-t)$ is minimal amongst the
branches of the second line of \eq{khAeq22} at $p_i$. Then as in the
proof of Proposition \ref{khAprop5}, we have $y\in\pd X_2$, so
$y=\exp(2v)(x',B')$ for $(x',B')\in\pd X$. Hence
$p_i=\exp((2-t)v)x'$. Therefore $(2-t)\Vert v\Vert_{C^0}\ge\ze$,
since $d(p_i,x')\ge\ze$. Also $u=\min T'(x',B')\ge 1$. The condition
$1+\th\ze\Vert v\Vert_{C^0}>D$ thus gives $u+\th(2-t)>D$, proving
the claim.

The pieces $X_c$ for $c\in C$ are thus of two kinds: the branches
$D+d(p_i,x)^2$ of $T$ for $i=1,\ldots,N$ from the first line of
\eq{khAeq22} are minimal near $p_i$ and so yield components $X_i$
with $p_i\in\pi(X_i^\ci)\subseteq\pi(X_i)\subseteq B_\ep(p_i)$, as
in Proposition \ref{khAprop2}. The remaining $X_c$ come from minimal
branches $u+\th(2-t)$ from the second line of \eq{khAeq22}, and
these correspond to minimal branches of $T'$ on $\pd X$ and so to
components $X'_{c'}$, with $\pi(X_{c'})\subset X$ the union of
segments of flow lines of $v$ starting in $\pi(X'_{c'})\subset
\io(\pd X)$. This proves the next part of the proposition. The final
part holds provided $\ep>0$ is taken sufficiently small with fixed
$g$ and $\th$ is taken sufficiently large, since then the
$\pi(X_{c'})$ for $c'\in C'$ become very thin neighbourhoods of the
$\pi(X'_{c'})$, and the $\pi(X_i)$ for $i=1,\ldots,N$ are contained
in small balls~$B_\ep(p_i)$.
\end{proof}

\subsubsection[Cutting into simplices with boundary
conditions]{Cutting a manifold into simplices with boundary
conditions}
\label{khA16}

We would also like to extend Lemma \ref{khAlem1} and Theorem
\ref{khAthm1} to the case of prescribed boundary data. However, we
cannot do this with arbitrary $T':\pd X\ra F\bigl([1,\iy)\bigr)$,
since Lemma \ref{khAlem1} relies on taking $\ep$ very small so that
the pieces $\pi(X_c)$ approximate small convex polyhedra in $\R^n$,
but as in Proposition \ref{khAprop7} the $X'_{c'}$ for $c'\in C'$
appear as components of $\pd X_c$, and the $X'_{c'}$ depend on $T'$
and are of a fixed size, so this prevents us making the $X_c$
arbitrarily small.

To get round this, we assume that $T'$ has been constructed using
Lemma \ref{khAlem1} or Theorem \ref{khAthm1} (or, recursively, using
Lemma \ref{khAlem2} or Theorem \ref{khAthm2} below with $\pd X$ in
place of $X$), using a sufficiently small constant $\ep'>0$, {\it
where `sufficiently small' depends not just on $\pd X$ and choices
made during the construction of\/ $T',$ but also on $X$ and choices
to be made during the construction of\/} $T$. Thus, $T'$ depends not
only on $\pd X$, it also has a minimal dependence on $X$ as well.
The argument here appears circular: we are trying to choose $T'$
depending on $T$, and $T$ depending on $T'$. However, it is not
circular, as the only data on $X$ that $\ep'$ must depend on are the
choices of $g,v,D$ in Definition \ref{khAdef7}, and these are
independent of $T'$, that is, they are chosen without reference
to~$T'$.

Here is our generalization of Lemma \ref{khAlem1}. There is one
subtle point: at the junction between a region $\pi(X_i)$ in $X$
coming from a branch $D+d(p_i,x)^2$ of $T$ in the first line of
\eq{khAeq22}, and a region $\pi(X_{c'})$ in $X$ coming from a branch
$u+\th(2-t)$ in the second line of \eq{khAeq22}, the common boundary
may not approximate a hyperplane $\R^{n-1}$ in $\R^n$, but rather a
piece of a quadric in $\R^n$. The condition $\th\inf_{\pd
X}\md{v}\gg\ep^2$ ensures this piece is approximately flat and
parallel to $\pd X$, so we can still approximate $\pi(X_i)$ and
$\pi(X_{c'})$ by compact convex polyhedra.

\begin{lem} Suppose $X$ is a compact manifold with g-corners. Let\/
$T':\pd X\ra F\bigl([1,\iy)\bigr)$ be constructed in Definition
{\rm\ref{khAdef2} (}or Definition {\rm\ref{khAdef7})} using some
metric $g'$ on $\pd X$ and constant\/ $\ep'>0,$ and suppose\/ $\min
T'\vert_{\pd^2X}$ is invariant under the natural involution
$\si:\pd^2X\ra\pd^2X$. Let\/ $T:X\ra F\bigl([1,\iy)\bigr)$ be
constructed in Definition {\rm\ref{khAdef7}} for this $T',$ using
some metric $g$ on $X,$ vector field\/ $v$ and constants
$D,\th,\ep$. We regard\/ $g',g,v,D$ as chosen in advance and fixed,
and\/ $\ep',\ep,\th$ as satisfying inequalities which may depend
on\/ $g',g,v,D$. In particular, $\ep'$ can depend on~$g,v,D$.

Suppose that\/ $\ep$ is small compared to the natural length-scales
of\/ $(X,g)$ and\/ $v,$ that\/ $\ep'$ is small compared to the
natural length-scales of\/ $(X,g),v,\ab(\pd X,\ab g')$ and constants
comparing $g\vert_{\pd X}$ and\/ $g',$ and that\/ $\th$ is chosen
large enough that\/ $\th\inf_{\pd X}\md{v}\gg\ep^2$. Then geodesic
normal coordinates in $X$ approximately identify the
$\pi(X_c)\subset X$ for $c\in C$ with compact convex polyhedra
in~$\R^n$.

For $c'\in C'$ Lemma {\rm\ref{khAlem1}} approximately identifies
$\pi'(X'_{c'})\subset\pd X$ with a compact convex polyhedron in
$\R^{n-1}$. Also, Proposition {\rm\ref{khAprop7}} identifies
$X'_{c'}$ with a component of\/ $\pd X_c,$ and we have approximately
identified $X_c$ with a compact convex polyhedron in $\R^n$.
Combining these three gives an approximate identification of a
compact convex polyhedron in $\R^{n-1}$ with a codimension one face
of a compact convex polyhedron in $\R^n$. This is approximately
affine.
\label{khAlem2}
\end{lem}

Here is our generalization of Theorem \ref{khAthm1}. The proof is an
easy modification, as we have already done the hard work in
Definition \ref{khAdef7} and Lemma \ref{khAlem2}. The main points
are to take $T'$ and $\si'_{(i',\bs F')}$ to be built in Theorem
\ref{khAthm1} using data on $\pd X$ in Lemma \ref{khAlem1} chosen
with $\ep'$ sufficiently small as in Lemma \ref{khAlem2}. Then when
we choose the $\ti\si_{(i,\bs F,\De)}$ inductively in Theorem
\ref{khAthm1}, those lying over $\pd X$ rather than $X^\ci$ are
equal to $\si'_{(i',\bs F')}\vert_{\De'}$ for some natural choice
of~$i',\bs F',\De'$.

In the last part of Lemma \ref{khAlem2} we comment that
identifications of $X'_{c'}$ and a component of $\pd X_c$ with
convex polyhedra in $\R^{n-1}$ and $\R^n$ are related by an
approximately affine transformation. The point here is that we have
two metrics $g'$ and $g\vert_{\pd X}$ on $\pd X$, with geodesic
normal coordinates on $X'_{c'}$ are defined using $g'$, and on the
component of $\pd X_c$ using $g\vert_{\pd X}$ (approximately). Since
$\ep'$ is chosen small compared to the length scales of both $g'$
and $g\vert_{\pd X}$ and any comparison between them, the geodesic
normal coordinates on $\pd X$ using $g'$ and $g\vert_{\pd X}$ differ
by approximately affine transformations at length scale $\ep'$. This
is important in Theorem \ref{khAlem2}, since the barycentric
subdivision of a compact convex polyhedron is invariant under affine
transformations, so the barycentric subdivision steps in Theorem
\ref{khAthm1} for $T'$ and Theorem \ref{khAthm2} for $T$ are
compatible.

In (ii) below, only the last sentence is an assumption, the rest
follows from Theorem \ref{khAthm1} for $T'$ and from~(i).

\begin{thm} Suppose $X$ is a compact manifold with g-corners. Let\/
$T':\pd X\ra F\bigl([1,\iy)\bigr),$ $X'_{c'},$ $c'\in C'$ and
diffeomorphisms\/ $\si'_{c'}:\De_{n-1}\ra X'_{c'}$ for\/ $c'\in C'$
be constructed in Theorem {\rm\ref{khAthm1} (}or Theorem
{\rm\ref{khAthm2})} for $\pd X$. Suppose:
\begin{itemize}
\setlength{\itemsep}{0pt}
\setlength{\parsep}{0pt}
\item[{\rm(i)}] Assume $\min T'\vert_{\pd^2X}$ is invariant under the
involution\/ $\si:\pd^2X\ra\pd^2X$.
\item[{\rm(ii)}] The maps\/ $\pi'\ci\si'_{c'}:\De_{n-1}\ra\pd X$ for
$c'\in C'$ are a triangulation of\/ $\pd X$ by\/
$(n-1)$-simplices. They induce a triangulation of\/ $\pd^2X$
by\/ $(n-2)$-simplices, of the following form: let\/
$F_{n-1}^{n-1}:\De_{n-2}\ra\De_{n-1}$ be as in
{\rm\S\ref{kh41}}. Let\/ $\check C'$ be the subset of\/ $c'\in
C'$ for which\/ $\pi'\ci\si'_{c'}\ci
F_{n-1}^{n-1}:\De_{n-2}\ra\pd X$ maps\/ $\De_{n-2}$ to\/
$\io(\pd^2X)$. For each\/ $c'\in\check C'$ there is a unique
embedding\/ $\check\si'_{c'}:\De_{n-2}\ra\pd^2X$ with\/ $\io\ci
\check\si'_{c'}\equiv\pi'\ci\si'_{c'}\ci F_{n-1}^{n-1},$ where\/
$\io:\pd^2X\ra\pd X$ is the natural immersion. Then\/
$\check\si'_{c'}$ for $c'\in\check C'$ are a triangulation of\/
$\pd^2X$ by\/ $(n-2)$-simplices. Part {\rm(i)} implies that the
set of images $\check\si'_{c'}(\De_{n-2})$ in $\pd^2X$ for
$c'\in\check C'$ is $\si$-invariant, since
$\check\si'_{c'}(\De_{n-2})$ are the $\pi(X_c)$ for
$T'\vert_{\pd^2X}$.

Assume the $\check\si'_{c'}$ for $c'\in\check C'$ are invariant
under $\si:\pd^2X\ra\pd^2X$, in the sense that if\/
$\si\bigl(\check\si'_{c'}(\De_{n-2})\bigr)=\check\si'_{\ti
c'}(\De_{n-2})$ for $c',\ti c'\in\check C'$ then
$\si\ci\check\si'_{c'}\equiv\check\si'_{\ti c'}$.
\item[{\rm(iii)}] Let the initial application of Definition
{\rm\ref{khAdef2}} in the construction of\/ $T'$ use some metric
$g'$ on $\pd X$ and constant\/ $\ep'>0$. Assume $\ep'$ is chosen
small enough, compared to geometry on\/ $X,$ that Lemma
{\rm\ref{khAlem2}} applies.
\end{itemize}
Then we can construct a tent function\/ $T:X\ra
F\bigl([1,\iy)\bigr),$ with\/ $\min T\vert_{\pd X}\equiv T',$ such
that the components $X_c,$ $c\in C$ of\/ $\pd Z_{X,T}$ in
\eq{khAeq5} are all diffeomorphic to the $n$-simplex $\De_n,$ with
diffeomorphisms $\si_c:\De_n\ra X_c$ for $c\in C,$ and for each
$c'\in C'$ there exists a unique $c\in C$ such that
$\io\ci\pi'\ci\si'_{c'}\equiv\pi\ci\si_c\ci F_n^n$ as maps
$\De_{n-1}\ra X,$ where $\io:\pd X\ra X$ is the natural immersion.
We can also choose the $\si_c$ to have the boundary compatibility
described in Theorem \ref{khAthm1} under $\si:\pd^2Z_{X,T}\ra\ab
\pd^2Z_{X,T}$.
\label{khAthm2}
\end{thm}

The reason why we use $F_{n-1}^{n-1}:\De_{n-2}\ra\De_{n-1}$ in (ii),
and $F_n^n:\De_{n-1}\ra\De_n$ in the last part, rather than just
referring to some codimension 1 faces of $\De_{n-1}$ and $\De_n$, is
that the barycentric subdivision of a convex $n$-polyhedron $K$ in
Definition \ref{khAdef3} automatically triangulates $\pd K$ by
$(n-1)$-simplices $\si_{\bs F}\ci F_n^n(\De_{n-1})$, so we know that
the $(n-2)$-simplices in the triangulation of $\pd^2X$ induced by
$\si_{c'}$, $c'\in C'$ are all of the form $\pi'\ci\si'_{c'}\ci
F_{n-1}^{n-1}(\De_{n-2})$ for some $c'\in C'$, and similarly the
$(n-1)$-simplices in the triangulation of $\pd X$ induced by
$\si_c$, $c\in C$ are all of the form $\pi\ci\si_c\ci
F_n^n(\De_{n-1})$ for some~$c\in C$.

\subsection{Tent functions on orbifolds}
\label{khA2}

Next we extend \S\ref{khA1} to orbifolds $X$. Here is the analogue
of Definition~\ref{khAdef1}.

\begin{dfn} Let $X$ be an $n$-orbifold with g-corners. Then as in
\S\ref{kh22}, $X$ is covered by orbifold charts $(U,\Ga,\phi)$,
where $\Ga$ acts linearly on $\R^n$, $U\subseteq\R^n$ is a
$\Ga$-invariant region with g-corners, and $\phi:U/\Ga\ra X$ is a
homeomorphism with an open set $\phi(U/\Ga)\subseteq X$. Write
$\pi:U\ra U/\Ga$ for the natural projection. Then $\phi\ci\pi$ maps
$U\ra X$. A function $T:X\ra F\bigl([1,\iy)\bigr)$ is called a {\it
tent function\/} if for all orbifold charts $(U,\Ga,\phi)$ on $X$,
the map $T\ci\phi\ci\pi:U\ra F\bigl([1,\iy)\bigr)$ is a tent
function on the $n$-manifold $U$, in the sense of
Definition~\ref{khAdef1}.

Equivalently, $T$ is a tent function if for each $x\in X$ there
exists a small orbifold chart $(U,\Ga,\phi)$ with $x\in\phi(U/\Ga)$,
open subsets $U_1,\ldots,U_N\subseteq U$ and smooth functions
$t_i:U_i\ra[1,\iy)$ for $i=1,\ldots,N$ satisfying Definition
\ref{khAdef1}(a)--(c), such that $T\ci\phi\ci\pi(u)=\bigl\{t_i(u):
i=1,\ldots,N$, $u\in U_i\bigr\}$ for all $u\in U$. We can take the
$t_i$ to be locally distinct, that is, $t_i\not\equiv t_j$ in any
nonempty subset of $U_i\cap U_j$, and the $U_j$ to be connected.
Then $U$ and $T\ci\phi\ci\pi$ determine the $U_i,t_i$ uniquely up to
permutations of~$\{1,\ldots,N\}$.

Since $T\ci\phi\ci\pi$ is invariant under the action of $\Ga$ on
$U$, it follows that $\Ga$ acts on $\{1,\ldots,N\}$ such that
$\ga(U_i)=U_{\ga\cdot i}$ and $t_{\ga\cdot i}\ci\ga\equiv t_i$ for
all $\ga\in\Ga$ and $i=1,\ldots,N$. Note that we do {\it not\/}
require this action of $\Ga$ on $\{1,\ldots,N\}$ to be trivial. This
means that tent functions on an orbifold $X$ in general {\it
cannot\/} be written locally as the union of finitely many smooth,
single-valued functions on $X$ near orbifold strata of $X$, it
really is necessary to lift to orbifold charts.
\label{khAdef8}
\end{dfn}

Much of \S\ref{khA1} extends from manifolds to orbifolds with almost
no change; this holds for \S\ref{khA14}, for instance. We will
comment only on new issues in the orbifold case. There is one
significant idea in Definition \ref{khAdef2} and Proposition
\ref{khAprop2} which does not work for orbifolds, and needs
revision. It is this: in Definition \ref{khAdef2}, with $X$ a
compact manifold, $g$ a Riemannian metric on $X$, and $d(\,,\,)$ the
metric induced by $g$, we used the notion of {\it injectivity
radius} $\de(g)>0$, and the fact that for any $p\in X$ and
$0<r<\de(g)$ the map $B_r(p)\ra[0,r^2)$ taking $q\mapsto d(p,q)^2$
is smooth.

For orbifolds this is no longer true: if $X$ is an orbifold and
$g,d$ are as above, then for any $r>0$, if $p\in X$ is sufficiently
close to an orbifold stratum of $X$, then the map $B_r(p)\ra[0,r^2)$
taking $q\mapsto d(p,q)^2$ is only piecewise smooth, not smooth. We
get round this by replacing $q\mapsto d(p,q)^2$ by a smooth, {\it
multivalued\/} function measuring the squared lengths of all
geodesic segments of length less than $r$ joining $p$ and~$q$.

Here are our generalizations of Definition \ref{khAdef2} and
Proposition~\ref{khAprop2}.

\begin{dfn} Let $X$ be a compact $n$-orbifold with g-corners.
Choose a Riemannian metric $g$ on $X$, and write $d(\,,\,)$ for the
metric on $X$ determined by $g$. As $X$ is compact we can choose a
finite system of orbifold charts $(U^j,\Ga^j,\phi^j)$ for $j\in J$
with $X=\bigcup_{j\in J}\phi^j(U^j/\Ga^j)$. Write $\pi^j:U^j\ra
U^j/\Ga^j$, so that $\phi^j\ci\pi^j$ maps $U^j\ra X$. Define
$g^j=(\phi^j\ci\pi^j)^* (g)$, a Riemannian metric on $U^j$. Write
$d^j(\,,\,)$ for the metric on $U^j$ induced by $g^j$. Even though
$U^j$ may not be compact, one can use the compactness of $X$ to show
that $g^j$ has a positive injectivity radius $\de(g^j)>0$ on $U^j$,
so that $u\mapsto d^j(p,u)^2$ is a smooth map $B_r^j(p)\ra[0,r^2)$
for all $p\in U$ and $0<r\le\de(g^j)$, where $B_r^j(p)$ is the open
ball of radius $r$ about $p$ in $U$ defined using~$d^j$.

Choose $\ep>0$ satisfying two conditions: (a) $2\ep<\de(g^j)$ for
all $j\in J$, and (b) for all $p\in X$, there exists $j\in J$ such
that $B_{3\ep}(p)\subseteq\phi^j(U^j/\Ga^j)$. This holds for all
sufficiently small $\ep$ as $J$ is finite, $\de(g^j)>0$ and $X$ is
compact. The open sets $B_\ep(p)$ for $p\in X$ cover $X$, so as $X$
is compact we can choose a finite subset $\{p_1,\ldots,p_N\}$ of $X$
with $X=\bigcup_{c=1}^NB_\ep(p_c)$. Choose $p_1,\ldots,p_N$ to be
{\it generic\/} amongst all such $N$-tuples $p_1',\ldots,p_N'$, as
in Definition \ref{khAdef2}. Define $T:X\ra F\bigl([1,\iy)\bigr)$ by
\e
\begin{split}
T(x)=\bigl\{&1+l^2:\text{$l\in[0,2\ep)$, there exists a geodesic
segment of length $l$}\\
&\text{in $(X,g)$ with end points $x$ and $p_c$ for some
$c=1,\ldots,N$}\bigr\}.
\end{split}
\label{khAeq23}
\e
Note that $T$ does not depend on the choice of $(U^j,\Ga^j,\phi^j)$,
$j\in J$, except through the smallness conditions (a),(b) on $\ep$.
These are a substitute for the condition $2\ep<\de(g)$ in Definition
\ref{khAdef2}, and are needed because the notion of injectivity
radius does not behave well on orbifolds.
\label{khAdef9}
\end{dfn}

\begin{prop} In Definition {\rm\ref{khAdef9},} $T$ is a tent
function on\/ $X$. In {\rm\eq{khAeq5},} we can take the indexing
set\/ $C$ to be $\{1,\ldots,N\}$ with\/ $p_c\in\pi(X_c^\ci)\subseteq
\pi(X_c)\subseteq B_\ep(p_c)$ for\/ $c=1,\ldots,N$. We can choose\/
$T$ so that the pieces\/ $X_c$ are `arbitrarily small', in the sense
that if\/ $\{V_i:i\in I\}$ is any open cover for $X$ then we can
choose $T$ so that for all\/ $c=1,\ldots,N$ there exists\/ $i\in I$
with\/~$\pi(X_c)\subseteq V_i$.

If\/ $X$ is an \begin{bfseries}effective\end{bfseries} orbifold then
the $X_c$ for $c\in C$ in \eq{khAeq5} are compact
\begin{bfseries}manifolds\end{bfseries} with g-corners.
\label{khAprop8}
\end{prop}

\begin{proof} Let $x\in X$. Then by choice of $\ep>0$ in Definition
\ref{khAdef9}, there exists $j\in J$ such that $B_{3\ep}(p)\subseteq
\phi^j(U^j/\Ga^j)$. Choose $y\in U^j$ with $\phi\ci\pi(y)=x$. Let
$\Ga=\{\ga\in\Ga^j:\ga\cdot y=\ga\}$, the stabilizer group of $y$ in
$\Ga$. Then if $\ga\in\Ga^j\sm\Ga$ we have $y\ne\ga\cdot y$, so
$d^j(y,\ga\cdot y)>0$. Choose $\de\in(0,\ep]$ such that
$d^j(y,\ga\cdot y)\ge 2\de$ for all $\ga\in\Ga^j\sm\Ga$, and define
$U=B_\de^j(y)$, an open neighbourhood of $y$ in $U$. Then $U$ is a
$\Ga$-invariant subset of $U^j$, and $U\cap\ga\cdot U=\emptyset$ for
all $\ga\in\Ga^j\sm\Ga$. Hence the projection $U/\Ga\ra U^j/\Ga^j$
given by $\Ga u\mapsto\Ga^ju$ is a homeomorphism with its image.
Define $\phi:U/\Ga\ra X$ by $\phi:\Ga u\mapsto\phi^j(\Ga^ju)$. Then
$(U,\Ga,\phi)$ is an orbifold chart on $X$ with~$x\in\phi(U/\Ga)$.

Write $q_1,\ldots,q_M$ for the points of $B_{2\ep+\de}^j(y)\cap
(\phi^j\ci\pi^j)^{-1}(\{p_1,\ldots,p_N\})$. That is,
$q_1,\ldots,q_M$ are those preimages of $p_1,\ldots,p_N$ in $U^j$
with distance less than $2\ep+\de$ from $y$. Since
$\bmd{(\phi^j\ci\pi^j)^{-1}\bigl(\{p_1,\ldots,p_N\})}\le
N\md{\Ga^j}$ this is a finite set, and $M\le\md{\Ga^j}N$. The action
of $\Ga$ on $B_{2\ep+\de}^j(y)\subseteq U^j$ permutes
$q_1,\ldots,q_M$, and $\ga\in\Ga^j\sm\Ga$ may also identify distinct
$q_a,q_b$. For $a=1,\ldots,M$, define an open set $U_a\subseteq U$
by $U_a=U\cap B_{2\ep}(q_j)$, and define a smooth function
$t_a:U_a\ra[1,\iy)$ by $t_a(u)=1+d^j(q_a,u)^2$. It is now easy to
check from \eq{khAeq23} that
\e
T\ci\phi\ci\pi:u\longmapsto\bigl\{t_a(u):a=1,\ldots,M,\;\> u\in
U_a\bigr\}.
\label{khAeq24}
\e

We claim that Definition \ref{khAdef1}(a)--(c) hold for these
$U_1,\ldots,U_M,t_1,\ldots,t_M$. The proof is a little different to
that in the manifold case in Proposition \ref{khAprop2}, and the
extra property $\dim S_{\{i_1,\ldots,i_l\}}=n-l+1$ for all
$\{i_1,\ldots,i_l\}\subseteq\{1,\ldots,N\}$ with $l\ge 2$ in
Proposition \ref{khAprop2} in general will not hold here. The point
is that although $p_1,\ldots,p_N$ are generic in $X$, as each $p_i$
in $X$ lifts to finitely many $q_a$ in $U$, the $M$-tuple
$q_1,\ldots,q_M$ may {\it not\/} be generic in~$U$.

However, $(\phi^j\ci\pi^j)^{-1}\bigl(\{p_1,\ldots,p_N\})$ is generic
amongst $\Ga^j$-invariant finite subsets of $U^j$, and this is
enough to imply Definition \ref{khAdef1}(a),(b). The conditions
$d^j(q_{i_1},u)=\cdots=d^j(q_{i_l},u)$ may intersect
non-transversely, so that $\dim S_{\{i_1,\ldots,i_l\}}>n-l+1$, but
this can happen only if a nontrivial subset of
$\{q_{i_1},\ldots,q_{i_l}\}$ is preserved by a nontrivial subgroup
$G$ of $\Ga^j$, and $S_{\{i_1,\ldots,i_l\}}$ is contained in the
fixed-point set Fix$(G)$ of $G$ in $U^j$. Since Fix$(G)$ is a
submanifold of $U^j$, one can show $S_{\{i_1,\ldots,i_l\}}$ is a
submanifold of $U$. Therefore $T\ci\phi\ci\pi$ in \eq{khAeq24} is a
tent function on the manifold $U$. As every $x\in X$ has an orbifold
chart $(U,\Ga,\phi)$ with $x\in\phi(U/\Ga)$ and $T\ci\phi\ci\pi$ a
tent function on $U$, equation \eq{khAeq23} defines a tent function
on the orbifold $X$, as we have to prove.

The second and third parts of the proposition, that we can take
$C=\{1,\ldots,N\}$ and choose the $X_c$ arbitrarily small, follow as
in Proposition \ref{khAprop2}. For the final part, suppose $X$ is an
{\it effective} orbifold. Then $\Ga^j$ acts effectively on $U^j$ for
each $j\in J$. Let $x\in X$, and define $(U,\Ga,\phi),U_1,
\ldots,U_M,t_1,\ldots,t_M$ as above. Write $\pi:Z_{X,T}\ra X$ for
the projection. Then there is a diffeomorphism
\begin{equation*}
Z_{X,T}\!\supseteq\! \pi^{-1}\bigl(\phi(U/\Ga)\bigr)\cong\bigl\{
(t,u): u\!\in\! U,\; 0\!\le\! t\!\le\!\min
T\!\ci\!\phi\!\ci\!\pi(u)\bigr\}/\Ga\!=\!Z_{U,T\ci\phi \ci\pi}/\Ga.
\end{equation*}
Points of $\coprod_{c\in C}X_c\subset\pd Z_{X,T}$ are of the form
$\bigl((t,x'),B\bigr)$ for $(t,x')\in Z_{X,T}$ with $t=\min T(x')$,
and $B$ a local boundary component of $Z_{X,T}$ containing $(t,x')$.
Lifting from $Z_{X,T}$ up to $Z_{U,T\ci\phi\ci\pi}$, choices of $B$
correspond to choices of branch $t_b$ of $T\ci\phi\ci\pi$ inducing
the value $\min T(x')$, for $b=1,\ldots,M$. Therefore
\e
\begin{split}
\smash{\ts\coprod\limits_{c\in C}X_c}\supseteq
(\pi\ci\io)^{-1}\bigl(\phi(U/\Ga)\bigr)
&\cong\bigl\{(t,u,b):u\in U,\;\> t=\min T\ci\phi\ci\pi(u),\\
&b=1,\ldots,M,\;\min T\ci\phi\ci\pi(u)=t_b(u)\bigr\}/\Ga,
\end{split}
\label{khAeq25}
\e
where $\io:\coprod_{c\in C}X_c\ra Z_{X,T}$ is the natural immersion.

Here in \eq{khAeq25}, $\Ga$ acts trivially on $t\in[0,\iy)$, in the
usual way on $u\in U$, and on $b\in\{1,\ldots,M\}$, the action is
induced from the action of $\Ga$ on $\{q_1,\ldots,q_M\}\subset U$.
Now $\Ga$ acts effectively on $U$, and $p_1,\ldots,p_N$ are generic
in $X$, so that each $q_b$ is generic in $U$. Together these imply
that the stabilizer group of each $q_b$ in $\Ga$ is trivial, so
$\Ga$ acts freely on $\{1,\ldots,M\}$. Thus the action of $\Ga$ in
\eq{khAeq25} is free, so the quotient $(\pi\ci\io)^{-1}(\phi
(U/\Ga))$ is a manifold. Since $\coprod_{c\in C}X_c$ is covered by
such open sets $(\pi\ci\io)^{-1}(\phi(U/\Ga))$, we see that
$\coprod_{c\in C}X_c$ and each $X_c$ for $c\in C$ are manifolds.
This completes the proof.
\end{proof}

\begin{rem} In Remark \ref{kh2rem2}(b), we explained that for
orbifolds $X$ with corners or g-corners, restricting to {\it
orbifold strata\/} $X^{\Ga,\rho}$ in the sense of \S\ref{kh56} does
not commute with taking boundaries, and therefore no information
from orbifold strata $X^{\Ga,\rho}$ survives in Kuranishi homology.
Proposition \ref{khAprop8} is an illustration of this. Given some
Kuranishi cycle $K$ defined using effective orbifolds $X$,
Propositions \ref{khAprop1} and \ref{khAprop8} enable us to
construct a homologous Kuranishi cycle $K'$ defined using manifolds
$X_c$ for $c\in C$, where the homology between $K$ and $K'$ is
defined using the orbifolds $Z_{X,T}$. This works because
$\pd(Z_{X,T}^{\Ga,\rho})\not\cong(\pd Z_{X,T})^{\Ga,\rho}$. So, a
homology theory defined using effective orbifolds with (g-)corners
as chains is equivalent to a homology theory defined using manifolds
as chains. This will be important in making the $X_{ac}$ have
trivial stabilizers in the proof of Step 1 in~\S\ref{khB1}.
\label{khArem1}
\end{rem}

The rest of \S\ref{khA1} extends to orbifolds $X$ without any
further significant problems or changes, except those already
explained in Definitions \ref{khAdef8}--\ref{khAdef9} and
Proposition \ref{khAprop8}. In fact we will not need analogues of
\S\ref{khA13} and \S\ref{khA16} for orbifolds, but we will need
versions of \S\ref{khA14}--\S\ref{khA15}. Section \ref{khA14} holds
for orbifolds with essentially no change. The choice of vector field
$v$ and flows $\exp(tv)$ work for orbifolds $X$ as well as
manifolds; note that $v$ is automatically tangent to each orbifold
stratum $X^{\Ga,\rho}$ of $X$, for $\Ga$ a finite group, and
$\exp(tv)$ takes $X^{\Ga,\rho}$ to $X^{\Ga,\rho}$. In \S\ref{khA15}
we restrict to {\it effective\/} orbifolds, so that we can use the
last part of Proposition \ref{khAprop8}, and we must modify the
first line of \eq{khAeq22} as in \eq{khAeq23}, and the proof of
Proposition \ref{khAprop7} as in Proposition~\ref{khAprop8}.

\begin{thm} The results of\/ {\rm\S\ref{khA14}--\S\ref{khA15}}
all hold when $X$ is a compact\/ $n$-orbifold with g-corners,
restricting to effective orbifolds in {\rm\S\ref{khA15}}. In
Proposition {\rm\ref{khAprop7},} if\/ $X$ is an effective orbifold
and the $X'_{c'}$ are manifolds for all\/ $c'\in C'$, then the $X_c$
are manifolds for all\/~$c\in C$.
\label{khAthm3}
\end{thm}

\subsection[Tent functions on Kuranishi chains {$[X,\bs f,\bs G],
[X,\bs f,$\underline{$\boldsymbol G\!$}$\,]$}]{Tent functions on
Kuranishi chains $[X,\bs f,\bs G],[X,\bs f,\ubG]$}
\label{khA3}

In \S\ref{khA31} we discuss some of the new problems in generalizing
\S\ref{khA1}--\S\ref{khA2} from manifolds and orbifolds to Kuranishi
chains $[X,\bs f,\bs G],[X,\bs f,\ubG]$, and how we solve them.
Section \ref{khA32} defines tent functions on Kuranishi chains, and
\S\ref{khA33}--\S\ref{khA35} give analogues of sections \ref{khA12},
\ref{khA14} and \ref{khA15} respectively. We will not need Kuranishi
chain analogues of \S\ref{khA13} and~\S\ref{khA16}.

\subsubsection{Preliminary discussion}
\label{khA31}

Let $X$ be a compact Kuranishi space, $Y$ an orbifold, $\bs f:X\ra
Y$ be strongly smooth, and $\bs G$ (or $\ubG$) be (effective)
gauge-fixing data for $(X,\bs f)$. Then $\bs G,\ubG$ include an
excellent coordinate system $(\bs I,\bs\eta)$ for $(X,\bs f)$, where
$\bs I=(I,(V^i,E^i,\ab s^i,\ab\psi^i):i\in I,\ldots)$. Roughly
speaking, a {\it tent function\/} for $(X,\bs f,\bs G)$ or $(X,\bs
f,\ubG)$ is $\bs T=(T^i:i\in I)$, where $T^i:V^i\ra
F\bigl([1,\iy)\bigr)$ is a tent function on $V^i$ for each $i\in I$,
such that if $j\le i$ in $I$ then~$\min T^j\vert_{V^{ij}}\equiv\min
T^i\ci\phi^{ij}$.

However, we need to modify this definition a little, to make it
easier to construct tent functions $\bs T$. Here are the two main
new problems we will meet in trying to extend the results of
\S\ref{khA1}--\S\ref{khA2} to Kuranishi chains:
\begin{itemize}
\setlength{\itemsep}{0pt}
\setlength{\parsep}{0pt}
\item[(a)] One thing we will need to do very often below is to
choose some smooth data $\de^i$ on $V^i$ for each $i\in I$, (such as
a Riemannian metric $g$, or a transverse perturbation $\ti s^i$ of
the Kuranishi map $s^i$), with compatibility conditions
$\de^j\vert_{V^{ij}}\equiv(\phi^{ij})^*(\de^i)$ when $j\le i$ in
$I$. Our basic method for doing this is to choose $\de^i$ by
induction on increasing $i\in I$, where
$\de^i\vert_{\phi^{ij}(V^{ij})}$ for all $j<i$ in $I$ is (partially)
prescribed by $\de^j\vert_{V^{ij}}\equiv(\phi^{ij})^*(\de^i)$ and
the choice of $\de^j$ in a previous inductive step.

Since $\phi^{ij}(V^{ij})$ need not be closed in $V^i$, we may not be
able to extend smooth prescribed values for $\de^i$ on
$\phi^{ij}(V^{ij})$ smoothly to $V^i$. In fact, if $\phi^{ij}
(V^{ij})$ is badly placed in $V^i$, then smooth prescribed values
for $\de^i$ on $\phi^{ij}(V^{ij})$ may not even extend continuously
to the closure $\overline{\phi^{ij}(V^{ij})}$ of $\phi^{ij}(V^{ij})$
in $V^i$. Thus, the inductive step may not be possible.
\item[(b)] All our results in \S\ref{khA1}--\S\ref{khA2} involve
{\it compact\/} manifolds and orbifolds. In particular, compactness
is used to choose finite covers by balls $B_\ep(p_c)$. However, the
$V^i$ are in general noncompact orbifolds.
\end{itemize}

Here is an example of problem (a).

\begin{ex} Let $X$ be the compact topological space
$\bigl\{(x,y)\in[-1,1]^2:xy=0\bigr\}$. Define Kuranishi
neighbourhoods $(V^1,E^1,s^1, \psi^1)$, $(V^2,E^2,s^2,\psi^2)$ on
$X$ by $V^1=\bigl\{(x,0):0<\md{x}\le
1\bigr\}\cup\bigl\{(0,y):0<\md{y}\le 1\bigr\}$, and $E^1=V^1$,
regarded as the zero vector bundle over $V^1$, $s^1\equiv 0$ and
$\psi^1=\id_{V^1}$, and $V^2=[-1,1]^2$, $E^2=\R\t V^2$,
$s^2:V^2\ra\R$ is $s^2(x,y)=xy$, and $\psi^2=\id_X$. Then $X$ has a
unique Kuranishi structure whose Kuranishi neighbourhoods
$(V_p,\ldots,\psi_p)$ are equivalent to $(V^1,\ldots,\psi^1)$ for
$p\ne(0,0)$ and to $(V^2,\ldots,\psi^2)$ for~$p=(0,0)$.

We can complete $(V^1,\ldots,\psi^1),(V^2,\ldots,\psi^2)$ to a very
good coordinate system for $X$, with indexing set $I=\{1,2\}$, and
$V^{21}=V^1$. Suppose we want to choose some continuous or smooth
data for each $(V^i,\ab\ldots,\ab\psi^i)$ by induction over $i$ in
$I=\{1,2\}$ with compatibilities on $V^{ij}$. Having chosen data for
$(V^1,\ldots,\psi^1)$, we have to choose data for
$(V^2,\ldots,\psi^2)$ taking prescribed values over the subset $V^1$
of $V^2$. In general this is {\it not possible}, even continuously,
near $(0,0)$ in $\bar V^1\sm V^1$, since the four different ends of
$V^1$ at $(0,0)$ can prescribe four different limits at $(0,0)$, if
the limits exist at all.

This example also illustrates a different problem with choosing {\it
transverse perturbations}, which appears to have been overlooked in
the proofs of Fukaya and Ono \cite[Th.~6.4]{FuOn1},
\cite[Th.~A1.23]{FOOO}. To carry out Step 2 in \S\ref{khB2}, we
would like to choose small, smooth, transverse perturbations $\ti
s^1,\ti s^2$ of the Kuranishi maps $s^1,s^2$, satisfying
$\hat\phi^{21}\ci \ti s^1\equiv \ti s^2\ci\phi^{21}$ on $V^{21}$.
This is not possible near $(0,0)$, not even using multisections, and
not because bad choices of $\ti s^1$ may not extend smoothly, but
for a different reason. As $E^1$ is the zero vector bundle, the only
choice for $\ti s^1$ is $\ti s^1\equiv 0$, and thus we must have
$\ti s^2\equiv 0$ on $V^1$, and hence on $\bar V^1$ in $V^2$. But
this implies $\ti s^2\vert_{(0,0)}=\d\ti s^2\vert_{(0,0)}=0$, and
$\ti s^2$ is not transverse at~$(0,0)$.
\label{khAex1}
\end{ex}

Our solution to problem (a) above, which we used in the proof of
Proposition \ref{kh3prop3} above, is to replace the $V^i$ by smaller
open subsets $\dot V^i$. We now define these $\dot V^i$ using the
partition of unity data $\bs\eta$ in $(\bs I,\bs\eta)$, and complete
the $\dot V^i$ to a modified really good coordinate system
$(\bs{\dot I},\bs{\dot\eta})$. The advantage of doing it this way,
rather than choosing some arbitrary smaller subsets $\dot V^i\subset
V^i$, is that passing from $(\bs I,\bs\eta)$ to $(\bs{\dot
I},\bs{\dot\eta)}$ is functorial in chains $[X,\bs f,\bs G]$ or
$[X,\bs f,\ubG]$, and compatible with relations in
$KC_*,KC_*^\ef(Y;R)$ that we need to preserve.

\begin{dfn} Let $X$ be a compact Kuranishi space, $Y$ an orbifold,
$\bs f:X\ra Y$ be strongly smooth, and $(\bs I,\bs\eta)$ an
excellent coordinate system for $(X,\bs f)$, where $\bs
I=(I,(V^i,E^i,s^i,\psi^i):i\in I,\ldots)$ and $\bs\eta=(\eta_i:i\in
I,\; \eta_i^j:i,j\in I)$.

Let $\vartheta>0$ be small. For each $j\in I$, define
$\dot\eta_i^j:V^j\ra[0,1]$ for $i\in I$ by
\e
\dot\eta_i^j(v)={\rm
mid}\bigl(0,(1+\vartheta)\eta_i^j(v)+c(v),1\bigr),
\label{khAeq26}
\e
where ${\rm mid}(x,y,z)$ is the middle one of $x,y,z$, that is,
${\rm mid}(x,y,z)=y$ when $x\le y\le z$ or $z\le y\le x$, and so on,
and $c(v)\in\R$ is chosen so that $\sum_{i\in I}\dot\eta_i^j(v)=1$.
From \eq{khAeq26} we see that $\sum_{i\in I}\dot\eta_i^j(v)$ is a
continuous, monotone increasing function of $c(v)$ with minimum 0
and maximum $\md{I}\ge 1$, so there exists $c(v)$ with $\sum_{i\in
I}\dot\eta_i^j(v)=1$ by the Intermediate Value Theorem, and this
defines the $\dot\eta_i^j(v)$ (though not necessarily $c(v)$)
uniquely. It is easy to see that the $\dot\eta_i^j$ are continuous.
Define continuous $\dot\eta_i:X\ra[0,1]$ from the $\eta_i$ in the
same way.

By construction $\sum_{i\in I}\dot\eta_i^j\equiv 1$ and $\sum_{i\in
I}\dot\eta_i\equiv 1$, as in Definition \ref{kh3def3}(i),(ii),
Definition \ref{kh3def3}(iii) implies that
$\dot\eta_i^j\vert_{(s^j)^{-1}(0)}\equiv\dot\eta_i\ci\psi^j$ for all
$i,j\in I$, and Definition \ref{kh3def3}(iv) implies that if
$i,j,k\in I$ with $k\le j$ then $\dot\eta_i^k\vert_{V^{jk}}
\equiv\dot\eta_i^j\ci\phi^{jk}$. Define open subsets $\dot
V^i\subset V^i$ and $\dot V^{ij}\subset V^{ij}$ by $\dot V^i=\{v\in
V^i:\dot\eta_i^i(v)>0\}$ and $\dot V^{ij}=\{v\in
V^{ij}:\dot\eta_i^j(v)>0\}$. Let $\dot E^i,\dot s^i,\dot\psi^i$ be
the restrictions of $E^i,s^i,\psi^i$ to $\dot V^i$ for all $i\in I$,
and $\dot\phi{}^{ij},\smash{\hat{\dot\phi}}{}^{ij}$ the restrictions
of $\phi^{ij},\hat\phi^{ij}$ to $\dot V^{ij}$. Set~$\dot I=\{i\in
I:\dot V^i\ne\es\}$.

When $\vartheta>0$, the effect of \eq{khAeq26} is to move the points
$\eta_i^j(v)$ in $[0,1]$ away from the middle and towards the end
points 0,1. Thus, if $\eta_i^j(v)>0$ is small then $\dot
\eta_i^j(v)$ is smaller, and may be zero. So we have $\supp(\dot
\eta_i^j)\subseteq\supp(\eta_i^j)$, and $\dot V^i\subseteq V^i$.
When $\vartheta=0$ we have $\dot\eta_i^j=\eta_i^j$ and $\dot V^i=
V^i$. Also $\dot\eta^i_j$ varies continuously with $\vartheta$, and
$\supp(\dot\eta^i_j)$ is a monotone decreasing, lower semicontinuous
function of $\vartheta$, so $\dot V^i$ and $\Im\dot\psi^i$ are also
monotone decreasing and lower semicontinuous in $\vartheta$. Hence,
when $\vartheta>0$, is small, $\dot V^i,\Im\dot\psi^i$ are just a
little smaller than~$V^i,\Im\psi^i$.

We have $X=\bigcup_{i\in I}\Im\psi^i$, where $X$ is compact, and
$\{\Im\psi^i:i\in I\}$ is a finite open cover for $X$. Being a
finite open cover of a compact space is not changed by making all
the sets a little smaller. Therefore, for sufficiently small
$\vartheta>0$, we have $X=\bigcup_{i\in I}\Im\dot\psi^i$. Let
$\vartheta>0$ be chosen small enough to achieve this.

It is not difficult to show that apart from the condition
$X=\bigcup_{i\in I}\Im\dot\psi^i$, all the conditions for this data
$\dot I,(\dot V^i,\ldots,\dot\psi^i),\ldots$ to define a really good
coordinate system follow from the corresponding conditions on
$I,(V^i,\ldots,\psi^i),\ldots$. Thus we have defined a {\it really
good coordinate system} $(\bs{\dot I},\bs{\dot\eta})$ for $(X,\bs
f)$. We can apply Algorithm \ref{kh3alg} to $(\bs{\dot
I},\bs{\dot\eta})$ to get an excellent coordinate system
$(\bs{\check I},\bs{\check\eta})$ if we wish.

We shall show that if $j\le i$ in $I$ then $\ov{\phi^{ij}(\dot
V^{ij})}\subseteq\phi^{ij}(V^{ij})$, where $\ov{\phi^{ij}(\dot
V^{ij})}$ is the closure of $\phi^{ij}(\dot V^{ij})$ in $V^i$. This
property will be important in problems involving choosing smooth
data on the $\dot V^i$ with compatibilities over the $\dot V^{ij}$.
Set $\ep=\smash{\frac{\vartheta}{(1+\vartheta)\md{I}}}>0$. Suppose
$\eta_j^j(v)\le\ep$ for some $v\in V^j$. We will prove that
$\dot\eta_j^j(v)=0$. Since $\sum_{i'\in I}\eta_{i'}^j\equiv 1$, the
average value of the $\eta_{i'}^j(v)$ for $i'\in\md{I}$ is
$\frac{1}{\md{I}}>\ep$, so $\eta_j^j(v)$ is not the maximum of the
$\eta_{i'}^j(v)$. If any $\dot\eta_{i''}^j(v)=1$ then
$\eta_{i''}^j(v)$ must be maximum amongst the $\eta_{i'}^j(v)$, so
$i''\ne i$, and thus $\dot\eta_j^j(v)=0$. So suppose
$\dot\eta_{i'}^j(v)\ne 1$ for any $i'\in I$. Then \eq{khAeq26} gives
$\dot\eta_{i'}^j(v)\ge (1+\vartheta)\eta_{i'}^j(v)+c(v)$, so summing
over $i'\in\md{I}$ and using $\sum_{i'\in I}\eta_{i'}^j\equiv
\sum_{i'\in I}\dot\eta_{i'}^j\equiv 1$ implies that $c(v)\le
-\frac{\vartheta}{\md{I}}$. But then $(1+\vartheta)\eta_j^j(v)+
c(v)\le 0$, so \eq{khAeq26} gives $\dot\eta_j^j(v)=0$, as we want.
Thus $\dot\eta_j^j=0$ wherever~$\eta_j^j\le\ep$.

Since $\dot\eta_j^j>0$ on $\dot V^{ij}$ we have $\eta_j^j>\ep$ on
$\dot V^{ij}$, so $\eta_j^i>\ep$ on $\phi^{ij}(\dot V^{ij})$ by
Definition \ref{kh3def3}(iv), and thus $\eta_j^i\ge\ep>0$ on
$\ov{\phi^{ij}(\dot V^{ij})}$. But $\phi^{ij}(V^{ij})$ is closed in
$\{v\in V^i:\eta_j^i(v)>0\}$ by Definition \ref{kh3def3}(ii), and
$\ov{\phi^{ij}(\dot V^{ij})}$ lies in the intersection of the
closure of $\phi^{ij}(V^{ij})$ and $\{v\in V^i:\eta_j^i(v)>0\}$,
so~$\ov{\phi^{ij}(\dot V^{ij})}\subseteq\phi^{ij}(V^{ij})$.
\label{khAdef10}
\end{dfn}

We can now explain our solutions to problems (a) and (b) above.

\begin{rem} Here is how to use Definition \ref{khAdef10} to avoid
problem (a) above. Suppose we want to choose some smooth data
$\de^i$ on $\dot V^i$ for $i\in I$ with $\de^j\vert_{\dot V^{ij}}
\equiv\de^i\ci\dot\phi^{ij}$ when $j\le i$ in $I$. Then, what we
actually do is to choose $\de^i$ on $\dot V^i$ by induction on
increasing $i\in I$, with $\de^j\vert_{\dot V^{ij}}
\equiv\de^i\ci\dot\phi^{ij}$ for all $j<i$ in $I$, such that $\de^i$
{\it extends smoothly to an open neighbourhood\/ $U^i$ of the
closure $\overline{\dot V^i}$ of\/ $\dot V^i$ in} $V^i$. Thus, in
the inductive step where we choose $\de^i$, we have prescribed
values for $\de^i$ on the subsets $\phi^{ij}(\dot V^{ij})$ of $V^i$
for $j<i$ in $I$. The prescribed values are consistent on the
overlaps $\phi^{ij}(\dot V^{ij})\cap\phi^{ik}(\dot V^{ik})$ for
$k<j<i$ because~$\de^k\vert_{\dot V^{jk}}
\equiv\de^j\ci\dot\phi^{jk}$.

Since $\ov{\phi^{ij}(\dot V^{ij})}\subseteq\phi^{ij}(V^{ij})$ as in
Definition \ref{khAdef10}, we see that $\ov{\phi^{ij}(\dot
V^{ij})}=\phi^{ij}(\ov{\dot V^{ij}})\subseteq\phi^{ij}(V^{ij}\cap
U^j)$. As $\de^j$ extends smoothly to $U^j$, it follows that the
prescribed values for $\de^i$ on $\phi^{ij}(\dot V^{ij})$ extend
smoothly to $\ov{\phi^{ij}(\dot V^{ij})}$ in $V^i$, and thus, to a
neighbourhood of $\phi^{ij}(\dot V^{ij})$ in $V^i$. This overcomes
problem~(a).

Here is how we deal with problem (b) above. In choosing data such as
a tent function $T^i$ on $\dot V^i$, the values of $T^i$ only really
matter near $\dot V^i\cap(s^i)^{-1}(0)$, since this is the part of
$\dot V^i$ identified with $X$. So we can use compactness of $X$ to
control the $T^i$ near the $\dot V^i\cap(s^i)^{-1}(0)$, and then
take $T^i$ to be constant away from $\dot V^i\cap (s^i)^{-1}(0)$,
for instance.

Note that when we apply Definition \ref{khAdef10} we must choose
some $\vartheta>0$, which must be chosen small enough that
$X=\bigcup_{i\in I}\Im\dot\psi^i$. When we do this in Appendices
\ref{khB} and \ref{khC}, we will generally be working not with one
$X,\bs f,(\bs I,\bs\eta)$ but with a finite collection of $X_a,\bs
f_a,(\bs I_a,\bs\eta_a)$ for $a\in A$. We choose the same
$\vartheta$ for all $a\in A$, and make it small enough that
$X_a=\bigcup_{i\in I_a}\Im\dot\psi^i_a$ for all $a\in A$, which is
possible as $A$ is finite.
\label{khArem2}
\end{rem}

\subsubsection{The definition of tent functions on Kuranishi chains}
\label{khA32}

\begin{dfn} Let $X$ be a compact Kuranishi space, $Y$ an orbifold,
$\bs f:X\ra Y$ be strongly smooth, and either $\bs G$ be
gauge-fixing data for $(X,\bs f)$, or $\ubG$ be effective
gauge-fixing data for $(X,\bs f)$. Then $\bs G$ or $\ubG$ includes
an excellent coordinate system $(\bs I,\bs\eta)$ for $(X,\bs f)$,
where $\bs I=(I,(V^i,E^i,s^i,\psi^i):i\in I,\ldots)$. Let $(\bs{\dot
I},\bs{\dot\eta})$ be as in Definition \ref{khAdef10}, where
$\bs{\dot I}=(I,(\dot V^i,\dot E^i,\dot s^i,\dot\psi^i):i\in
I,\ldots)$, depending on the choice of some small~$\vartheta>0$.

A {\it tent function\/} for $(X,\bs f,\bs G)$ or $(X,\bs f,\ubG)$ is
$\bs T=(T^i:i\in I)$, where $T^i:\dot V^i\ra F\bigl([1,\iy)\bigr)$
is a tent function on the orbifold $\dot V^i$ for each $i\in I$,
such that $T^i$ extends to a tent function on an open neighbourhood
of the closure $\overline{\dot V^i}$ of $\dot V^i$ in $V^i$, and if
$j\le i$ in $I$ then $\min T^j\vert_{\dot V^{ij}}\equiv \min
T^i\ci\dot\phi^{ij}$, and the submanifolds $S_{\{i_1,\ldots,i_l\}}$
of Definition \ref{khAdef1} for $T^i$ intersect $\dot\phi^{ij}(\dot
V^{ij})$ transversely in $\dot V^i$ wherever $t_{i_1}(u)=\min
T^i(u)$. When $\dot V^i$ is an orbifold, this last condition must be
expressed in orbifold charts, as follows. Suppose $(U,\Ga,\phi)$ is
an orbifold chart for $\dot V^i$. Since $\dot\phi^{ij}(\dot V^{ij})$
is an embedded suborbifold of $\dot V^i$, $V=(\phi\ci\pi)^{-1}\bigl(
\dot\phi^{ij} (\dot V^{ij})\bigr)$ is a submanifold of $U$. By
Definition \ref{khAdef8}, $T^i\ci\phi\ci\pi:U\ra
F\bigl([1,\iy)\bigr)$ is a tent function on the manifold $U$. We
require that each of the submanifolds $S_{\{i_1,\ldots,i_l\}}$ of
$U$ given in Definition \ref{khAdef1} for $T=T^i\ci\phi\ci\pi$
should intersect $V$ transversely wherever~$t_{i_1}(u)=\min
T^i\ci\phi\ci\pi(u)$.

If $\bs T,\bs{\ti T}$ are tent functions for $(X,\bs f,\bs G)$ or
$(X,\bs f,\ubG)$ then we write $\min\bs T=\min\bs{\ti T}$ if $\min
T^i\equiv\min\ti T^i$ on $\dot V^i$ for all $i\in I$. If $\bs
T=(T^i:i\in I)$ is a tent function for $(X,\bs f,\bs G)$ or $(X,\bs
f,\ubG)$ then $\bs T\vert_{\pd X}=\bigl(T_\pd^{i-1}=
T^i\vert_{\pd\dot V^i}:i\in I$, $\pd V^i\ne\emptyset\bigr)$ is a
tent function for $(\pd X,\bs f\vert_{\pd X},\bs G\vert_{\pd X})$ or
$(\pd X,\bs f\vert_{\pd X},\ubG\vert_{\pd X})$. Often, given some
tent function $\bs T'$ for $(\pd X,\bs f\vert_{\pd X},\bs
G\vert_{\pd X})$, we will wish to construct a tent function $\bs T$
for $(X,\bs f,\bs G)$ with~$\min\bigl(\bs T\vert_{\pd
X}\bigr)=\min\bs T'$.
\label{khAdef11}
\end{dfn}

For tent functions $T$ on an oriented manifold or orbifold $X$, we
defined an oriented manifold or orbifold $Z_{X,T}$ in \eq{khAeq4},
with $\pd Z_{X,T}$ given by \eq{khAeq5}. Here are analogues of these
for tent functions on (effective) Kuranishi chains.

\begin{dfn} Suppose $X$ is a compact Kuranishi space, $Y$ an orbifold,
$\bs f:X\ra Y$ strongly smooth, and $\bs G$ (or $\ubG$) is
(effective) gauge-fixing data for $(X,\bs f)$. Let $\bs T$ be a tent
function for $(X,\bs f,\bs G)$ or $(X,\bs f,\ubG)$. Since the
continuous functions $\min T^i:\dot V^i\ra[1,\iy)$ for $i\in I$
satisfy $\min T^i\ci\dot\phi^{ij}\equiv\min T^j\vert_{\dot V^{ij}}$,
they induce a unique continuous function $\min T:X\ra[1,\iy)$ with
$\min T\ci\dot\psi^i\equiv\min T^i\vert_{(\dot s^i)^{-1}(0)}$ for
$i\in I$. As in \eq{khAeq4}, define $Z_{X,\bs T}$ as a compact
topological space by
\begin{equation*}
Z_{X,\bs T}=\bigl\{(t,x)\in [0,\iy)\t X:t\le\min T(x)\bigr\}.
\end{equation*}

We will define an oriented Kuranishi structure on $Z_{X,\bs T}$ and
(effective) gauge-fixing data $\bs H_{X,\bs T}$ or $\ubH_{X,\bs T}$
for $(Z_{X,\bs T},\bs f\ci\bs\pi)$, where $\bs\pi:Z_{X,\bs T}\ra X$
is the projection. Set $J=\{i+1:i\in I\}$. For each $i\in I$, define
\e
W^{i+1}=[0,\ha)\t V^i\amalg \bigl\{(t,v)\in [\ha,\iy)\t\dot V^i:
t\le\min T^i(v)\bigr\}\subset V^i\t[0,\iy).
\label{khAeq27}
\e
Then $W^{i+1}$ is an {\it orbifold with g-corners}, as $T^i:\dot
V^i\ra F\bigl([1,\iy)\bigr)$ is a tent function, by the orbifold
version of Proposition \ref{khAprop1}. Write $\pi_{V^i}:W^{i+1}\ra
V^i$ for the projection. Define an orbifold vector bundle
$F^{i+1}\ra W^{i+1}$ by $F^{i+1}=\pi_{V^i}^*(E^i)$, a section
$t^{i+1}$ of $F^{i+1}$ by $t^{i+1}=s^i\ci\pi_{V^i}$, and a
continuous map $\xi^{i+1}:(t^{i+1})^{-1}(0)\ra Z_{X,\bs T}$ by
$\xi^{i+1}(t,v)=(t,\psi^i(v))$. Define a smooth map
$g^{i+1}:W^{i+1}\ra Y$ by $g^{i+1}(t,v)=f^i(v)$. Then
$(W^j,F^j,t^j,\xi^j)$ for $j\in J$ are {\it Kuranishi
neighbourhoods} on $Z_{X,\bs T}$.

If $j\le i$ in $I$, define an open subset $W^{(i+1)(j+1)}$ in
$W^{j+1}$ by $W^{(i+1)(j+1)}=W^{j+1}\cap\bigl([0,\iy)\t
V^{ij}\bigr)$, and define a {\it coordinate change}
$(\psi^{(i+1)(j+1)},\hat\psi^{(i+1)(j+1)})$ from
$(W^{(i+1)(j+1)},\ldots,\xi^{j+1}\vert_{W^{(i+1)(j+1)}})$ to
$(W^{i+1},\ldots,\xi^{i+1})$ by $\psi^{(i+1)(j+1)}(t,v)\ab=\ab(t,\ab
\phi^{ij}(v))$ and $\hat\psi^{(i+1)(j+1)}(t,e)=(t,\hat\phi^{ij}
(e))$. The condition in Definition \ref{khAdef11} that the
$S_{\{i_1,\ldots,i_l\}}$ for $T^i$ should intersect
$\dot\phi^{ij}(\dot V^{ij})$ transversely in $\dot V^i$ ensures that
$\psi^{(i+1)(j+1)}$ is compatible with boundaries and g-corners, as
in Definition~\ref{kh2def12}(a).

Since $X$ is an oriented Kuranishi space, $[0,\iy)\t X$ is an
oriented Kuranishi space. It is now easy to check that there is a
unique oriented Kuranishi structure on $Z_{X,\bs T}$ such that
$Z_{X,\bs T}$ is a Kuranishi subspace of $[0,\iy)\t X$, with natural
strongly smooth map $\bs\pi:Z_{X,\bs T}\ra X$ induced by the
projection $[0,\iy)\t X\ra X$, and $\bs J=\bigl(J,(W^j,F^j,t^j,
\xi^j),g^j:j\in J$, $(W^{jk},\psi^{jk},\hat\psi^{jk}):k\le j\in
J\bigr)$ is a {\it very good coordinate system} for~$(Z_{X,\bs
T},\bs f\ci\bs\pi)$.

For $i,j\in I$, define continuous $\ze_{i+1}:Z_{X,\bs T}\ra[0,1]$,
$\ze_{i+1}^{j+1}:W^{j+1}\ra[0,1]$~by
\e
\begin{split}
\ze_{i+1}(t,x)=\begin{cases} (1-2t)\eta_i(x)+2t\dot\eta_i(x), &
t\in[0,\ha), \\ \dot\eta_i(x), & t\ge\ha,\end{cases} \\
\ze_{i+1}^{j+1}(t,v)=\begin{cases}
(1-2t)\eta_i^j(v)+2t\dot\eta_i^j(v), & t\in[0,\ha), \\
\dot\eta_i^j(v), & t\ge\ha.\end{cases}
\end{split}
\label{khAeq28}
\e
Write $\bs\ze=(\ze_j:j\in J$, $\ze_j^k:j,k\in J)$. Since $(\bs
I,\bs\eta)$ and $(\bs{\dot I},\bs{\dot\eta})$ are both really good
coordinate systems for $(X,\bs f)$, one can show that $(\bs
J,\bs\ze)$ is a {\it really good coordinate system} for~$(Z_{X,\bs
T},\bs f\ci\bs\pi)$.

Let $(\bs{\check J},\bs{\check\ze})$ be the excellent coordinate
system for $(Z_{X,\bs T},\bs f\ci\bs\pi)$ constructed from $(\bs
J,\bs\ze)$ by Algorithm \ref{kh3alg}. Then $\check W^j$ is an open
subset of $W^j$ for each $j\in J$. As $(\bs I,\bs\eta)$ is
excellent, no pieces of $[0,\ha)\t V^i\subseteq W^{i+1}$ are
discarded, so $[0,\ha)\t V^i\subseteq \check W^{i+1}$ for all $i\in
I$. This implies~$\check J=J=\{i+1:i\in I\}$.

If we have gauge-fixing data $\bs G$ for $(X,\bs f)$, define
$H^{i+1}:\check W^{i+1}\ra P$ for each $i\in I$ by
$H^{i+1}(t,v)=G^i(v)$ if $t=0$, and
$H^{i+1}(t,v)=S_{n+1}(t,x_1,\ldots,x_n)$ if $t>0$ and
$G^i(v)=S_n(x_1,\ldots,x_n)$, and write $\bs H_{Z,\bs T}=\bigl(
(\bs{\check J},\bs{\check\ze}),H^j:j\in\check J\bigr)$. If we have
effective gauge-fixing data $\ubG$ for $(X,\bs f)$, define
$\uH^{i+1}:\check W^{i+1}\ra\uP$ for each $i\in I$ by
$\uH^{i+1}(t,v)=\uG^i(v)$ if $t=0$, and
$\uH^{i+1}(t,v)=(t,x_1,\ldots,x_n)$ if $t>0$ and
$\uG^i(v)=(x_1,\ldots,x_n)$, and write $\ubH_{Z,\bs T}=\bigl(
(\bs{\check J},\bs{\check\ze}),\uH^j:j\in\check J\bigr)$. Then $\bs
H_{Z,\bs T}$ or $\ubH_{Z,\bs T}$ is (effective) gauge-fixing data
for~$(Z_{X,\bs T},\bs f\ci\bs\pi)$.

If $\bs T,\bs{\ti T}$ are tent functions for $(X,\bs f,\bs G)$ or
$(X,\bs f,\ubG)$ and $\min\bs T=\min\bs{\ti T}$ then $Z_{X,\bs
T},\bs\pi$ and $\bs H_{X,\bs T}$ or $\ubH_{X,\bs T}$ coincide with
$Z_{X,\bs{\ti T}},\bs\pi,\bs H_{X,\bs{\ti T}}$ or~$\ubH_{X,\bs{\ti
T}}$.
\label{khAdef12}
\end{dfn}

Here is the analogue of Proposition \ref{khAprop1}. The proof is
straightforward.

\begin{prop} In Definition {\rm\ref{khAdef12}} we have
\ea
\begin{split}
\pd\bigl[Z_{X,\bs T},\bs f\ci\bs\pi,\bs H_{X,\bs T}\bigr]=\,&-
\bigl[Z_{\pd X,\bs T\vert_{\pd X}},\bs f\vert_{\pd X}\ci\bs\pi,\bs
H_{\pd X,\bs T\vert_{\pd X}}\bigr]\\
&-[X,\bs f,\bs G]+\ts\sum_{c\in C}[X_c,\bs f_c,\bs G_c],
\end{split}
\label{khAeq29}\\
\begin{split}
\pd\bigl[Z_{X,\bs T},\bs f\ci\bs\pi,\ubH_{X,\bs
T}\bigr]=\,&-\bigl[Z_{\pd X,\bs T\vert_{\pd X}},\bs f\vert_{\pd
X}\ci\bs\pi,\ubH_{\pd X,\bs T\vert_{\pd X}}\bigr]\\
&-[X,\bs f,\ubG]+\ts\sum_{c\in C}[X_c,\bs f_c,\ubG_c],
\end{split}
\label{khAeq30}
\ea
in $KC_*(Y;R)$ and\/ $KC_*^\ef(Y;R)$. Here $C$ is a finite indexing
set and\/ $(X_c,\bs f_c,\bs G_c)$ for $c$ in $C$ in \eq{khAeq29} are
the connected components {\rm(}in the sense of Lemma
{\rm\ref{kh3lem})} of\/ $\bigl(\pd Z_{X,\bs T},\bs
f\ci\bs\pi\vert_{\pd Z_{X,\bs T}},\bs H_{X,\bs T}\vert_{\pd Z_{X,\bs
T}}\bigr)$ which meet\/ $(0,\iy)\t X^\ci,$ and similarly for
$(X_c,\bs f_c,\ubG_c)$ in \eq{khAeq30}. For each\/ $c\in C,$
$\bs\pi\vert_{X_c}:X_c\ra X$ is an orientation-preserving immersion
of Kuranishi spaces which is an embedding on $X_c^\ci,$ with\/
$X=\bigcup_{c\in C}\pi(X_c),$ and\/
$\pi(X_c^\ci)\cap\pi(X_{c'}^\ci)=\emptyset$ for all\/ $c\ne c'$
in\/~$C$.
\label{khAprop9}
\end{prop}

We think of this as cutting $X$ into finitely many pieces $\pi(X_c)$
for~$c\in C$.

\subsubsection{Cutting a Kuranishi chain into arbitrarily small pieces}
\label{khA33}

We now generalize \S\ref{khA12} to Kuranishi chains. Here are
analogues of Definitions \ref{khAdef2} and \ref{khAdef9} and
Propositions \ref{khAprop2} and~\ref{khAprop8}.

\begin{dfn} Suppose $X$ is a compact Kuranishi space, $Y$ an
orbifold, $\bs f:X\ra Y$ strongly smooth, and $\bs G$ (or $\ubG$) is
(effective) gauge-fixing data for $(X,\bs f)$. Let $\bs I=(I,(V^i,
E^i,s^i,\psi^i):i\in I,\ldots)$ and $\bs{\dot I}=(I,(\dot V^i,\dot
E^i,\dot s^i,\dot\psi^i):i\in I,\ldots)$ be as in
Definition~\ref{khAdef11}.

In the following, when we choose some data on $\dot V^i$ and require
that it should extend to an open neighbourhood $U^i$ of
$\overline{\dot V^i}$ in $V^i$, we do not mean a fixed neighbourhood
$U^i$, but a different $U^i$ each time. We should think of $U^i$ as
getting smaller each time it is used in an inductive proof, so that
in effect we are choosing a series of neighbourhoods
$U^i_1,U^i_2,U^i_3,\ldots$ with
\begin{equation*}
V^i\supseteq \overline{U^i_1} \supseteq
U^i_1\supseteq\overline{U^i_2} \supseteq U^i_2\supseteq
\overline{U^i_3}\supseteq U^i_3 \supseteq \cdots \supseteq
\overline{\dot V^i}\supseteq \dot V^i.
\end{equation*}
This convention eliminates problems about extending data from $U^i$
to~$\overline{U^i}$.

By induction on increasing $i\in I$, choose a Riemannian metric
$g^i$ on $\dot V^i$ which extends to a neighbourhood $U^i$ of
$\overline{\dot V^i}$ in $V^i$, such that whenever $j\le i$ in $I$
we have $(\phi^{ij})^*(g^i)\equiv g^j$ on $V^{ij}\cap
U^j\cap(\phi^{ij})^{-1}(U^i)$, and $\phi^{ij}(V^{ij})$ is
geodesically closed in $(U^i,g^i)$. This is possible as in Remark
\ref{khArem2}. These conditions imply that the metrics $d^i,d^j$ on
$\dot V^i,\dot V^j$ induced by $g^i,g^j$ satisfy
$(\phi^{ij})^*(d^i)\equiv d^j$ on $\dot V^{ij}$, for small
distances, and $\phi^{ij}$ maps $g^j$-geodesic segments in $\dot
V^{ij}$ to $g^i$-geodesics in~$\dot V^i$.

In Definition \ref{khAdef9}, to construct tent functions $T$ on an
orbifold $X$, we introduced a finite system of orbifold charts on
$X$, but the only way in which we used them was to determine
injectivity-radius-type smallness conditions on $\ep>0$. We need to
do the same here, but for simplicity we will not go into detail. In
brief, for each $i\in I$ we choose a finite system of orbifold
charts on the orbifold $U^i\subseteq V^i$ which cover the compact
subset $\overline{\dot V^i}\cap(s^i)^{-1}(0)$, but need not cover
the possibly noncompact orbifold $U^i$. Then we choose $\ep>0$
sufficiently small such that the conditions of Definition
\ref{khAdef9} hold for these charts on $U^i$ with metric $g^i$, for
all~$i\in I$.

Now choose finite subsets $\{p_1^i,\ldots,p_{N^i}^i\}$ in $\dot V^i$
for each $i\in I$ satisfying:
\begin{itemize}
\setlength{\itemsep}{0pt}
\setlength{\parsep}{0pt}
\item[(a)] $\overline{\dot V^i}\cap(s^i)^{-1}(0)\subseteq
\bigcup_{b=1}^{N^i}B^i_\ep(p_b^i)$, where $B_r^i(p)$ is the open
ball of radius $r$ about $p$ in $(U^i,g^i)$.
\item[(b)] $\phi^{ij}\bigl(\{p_1^j,\ldots,p_{N^j}^j\}\cap \dot
V^{ij}\bigr)=\{p_1^i,\ldots,p_{N^i}^i\}\cap\phi^{ij}(\dot V^{ij})$
for all $j<i$ in~$I$.
\item[(c)] $\{p_1^i,\ldots,p_{N^i}^i\}\setminus\bigcup_{j\in J:j<i}
\phi^{ij}(\dot V^{ij})$ is generic amongst finite subsets of $\dot
V^i$, for all~$i\in I$.
\item[(d)] if $j<i$ in $I$ and $v\in\dot V^{ij}\subseteq\dot V^j$,
then either $d^j(p_b^j,v)>\ep$ for all $b=1,\ldots,N^j$ with
$p_b^j\notin\dot V^{ij}$, or there exists $p_c^j\in\dot V^{ij}$ for
some $c=1,\ldots,N^j$ such that $d^j(p_c^j,v)<d^j(p_b^j,v)$ for all
$b=1,\ldots,N^j$ with $p_b^j\notin\dot V^{ij}$.
\item[(e)] if $j<i$ in $I$ and $v\in\phi^{ij}(\dot V^{ij})\subseteq
\dot V^i$, then either $d^i(p_b^i,v)>\ep$ for all $b=1,\ldots,N^i$
with $p_b^j\notin\phi^{ij}(\dot V^{ij})$, or there exists
$p_c^i\in\phi^{ij}(\dot V^{ij})$ for some $c=1,\ldots,N^i$ such that
$d^i(p_c^i,v)<d^i(p_b^i,v)$ for all $b=1,\ldots,N^i$
with~$p_b^j\notin\phi^{ij}(\dot V^{ij})$.
\end{itemize}

Here is one way to do this. For each subset $\emptyset\ne J\subseteq
I$, let $j=\min J$, and choose a finite subset
\e
S_J\subset \bigcap_{i\in J}\dot V^{ij} \setminus\bigcup_{k\in
I:k<j}\phi^{jk}(\dot V^{jk}) \subseteq \dot V^j
\label{khAeq31}
\e
such that
\e
(s^j)^{-1}(0)\cap \bigcap_{i\in J}\overline{\dot V^{ij}}\subseteq
\bigcup_{p\in S_J}B^j_\ep(p),
\label{khAeq32}
\e
where in \eq{khAeq32}, the closure is taken in $U^j$, and balls
$B^j_\ep(p)$ are in $(U^j,g^j)$. It is possible to choose such a
finite set as the left hand side of \eq{khAeq32} is compact, and is
also contained in the closure of the right hand side of
\eq{khAeq31}, since the deletions $\phi^{jk}(\dot V^{jk})$ have
strictly smaller dimension than $\dot V^j$, so they do not change
the closure. Then for each $i\in I$, define
\begin{equation*}
\{p_1^i,\ldots,p_{N^i}^i\}=\bigcup\nolimits_{\begin{subarray}{l}
\text{sets $J$ with $i\in J\subseteq I$:}\\
\text{set $j=\min J$}\end{subarray}}\phi^{ij}(S_J).
\end{equation*}

Then conditions (a),(b) above hold automatically. For (c), it is
enough to take each $S_J$ to be generic amongst such subsets. For
(d),(e) we have to work a bit harder. Basically we should choose the
$S_J$ by induction on $\md{J}$, and when we choose $S_K$ with
$k=\min K$, we should ensure that for all $\emptyset\ne J\subset K$
with $J\ne K$ and $j=\min J$, each point of
$\phi^{jk}\bigl(\bigcap_{i\in K}\dot V^{ik}\bigr)$ is either
distance at least $\ep$ from each point of $S_J$ (which was chosen
in a previous inductive step), or is closer to some point of
$\phi^{jk}(S_K)$ than to any point of $S_J$. This is possible if the
points in $S_K$ are sufficiently close together. But to deal with
some noncompactness issues we may need to make $\ep$ smaller, and to
require the points in $S_J$ to be close to $(s^j)^{-1}(0)$, so that
we are really working in a compact subset of~$U^j$.

Following \eq{khAeq23}, define $T^i:\dot V^i\ra
F\bigl([1,\iy)\bigr)$ for each $i\in I$ by
\e
\begin{split}
&T^i(v)\!=\!\bigl\{1\!+\!l^2:\text{$l\in[0,2\ep)$, there exists a
geodesic segment of length $l$}\\
&\;\>\text{in $(\dot V^i,g^i)$ with end points $v$ and $p_b^i$, some
$b=1,\ldots,N^i$}\bigr\}\!\cup\!\{1\!+\!\ep^2\}.
\end{split}
\label{khAeq33}
\e
Then by a very similar proof to Proposition \ref{khAprop8}, $T^i$ is
a {\it tent function\/} on $\dot V^i$. We added the extra branch
$1+\ep^2$ in \eq{khAeq33} because the balls $B_\ep^i(p_b^i)$ for
$b=1,\ldots,N^i$ need not cover $\dot V^i$, but only $\dot V^i\cap
(s^i)^{-1}(0)$, so without including $1+\ep^2$ we might have
$T^i(v)=\emptyset$, which is not allowed.

Let $j\le i$ in $I$. Then (d) above implies that $\min
T^j\vert_{\dot V^{ij}}$ is the same whether we compute it using all
of $p_1^j,\ldots,p_{N^j}^j$, or only those $p_b^j$ lying in $\dot
V^{ij}$. Part (e) implies $\min T^i\vert_{\dot\phi^{ij}(\dot
V^{ij})}$ is the same whether we compute it using all of
$p_1^i,\ldots,p_{N^i}^i$, or only those $p_b^i$ lying in
$\dot\phi^{ij}(\dot V^{ij})$. Part (a), $(\phi^{ij})^*(g^i)\equiv
g^j$ on $\dot V^{ij}$, and $\phi^{ij}(\dot V^{ij})$ geodesically
closed in $(\dot V^i,g^i)$, imply that $\min T^j\vert_{\dot V^{ij}}$
computed using only $p_b^i$ in $\dot\phi^{ij}(\dot V^{ij})$ agrees
with $\min T^i\vert_{\dot\phi^{ij}(\dot V^{ij})}$ computed using
only $p_b^i$ in $\dot\phi^{ij}(\dot V^{ij})$ under $\dot\phi^{ij}$.
Therefore $\min T^i\ci\dot\phi^{ij}\equiv\min T^j\vert_{\dot
V^{ij}}$, as we have to prove.

It remains only to verify the condition in Definition \ref{khAdef11}
that the submanifolds $S_{\{i_1,\ldots,i_l\}}$ of Definition
\ref{khAdef1} for $T^i$ intersect $\dot\phi^{ij}(\dot V^{ij})$
transversely in $\dot V^i$ wherever $t_{i_1}(u)=\min T^i(u)$. Part
(e) above implies that any branch of $t_k$ of $T^i$ on
$\dot\phi^{ij}(\dot V^{ij})$ with $t_k(u)=\min T^i(u)$ comes either
from a point $p_b^i$ in $\dot\phi^{ij}(\dot V^{ij})$, or the extra
branch $1+\ep^2$ in \eq{khAeq33}; we cannot get such branches $t_k$
from $p_c^i$ in $\dot V^i\sm\dot\phi^{ij}(\dot V^{ij})$. Using this
and the genericness condition (c), one can show that the
$S_{\{i_1,\ldots,i_l\}}$ intersect $\dot\phi^{ij}(\dot V^{ij})$
transversely. This proves that $\bs T=(T^i:i\in I)$ is a {\it tent
function} for $(X,\bs f,\bs G)$ or~$(X,\bs f,\ubG)$.
\label{khAdef13}
\end{dfn}

\begin{prop} In Definition {\rm\ref{khAdef13},} if\/ $\bs
I_c=\bigl(I_c,(V^i_c,\ldots,\psi^i_c):i\in I_c,\ldots\bigr)$ is the
very good coordinate system for $(X_c,\bs f_c)$ in $\bs G_c$ or
$\ubG_c$ in Proposition {\rm\ref{khAprop9},} then by construction
$I_c\subseteq I$ and $\bs\pi:X_c\ra X$ lifts to natural embeddings
$\pi:V^i_c\ra\dot V^i$ for $i\in I_c$. For each\/ $i\in I_c$ there
is a unique $b^{i,c}=1,\ldots,N^i$ with\/
$p^i_{\smash{b^{i,c}}}\in\pi(V^i_c),$ and\/ $\pi(V^i_c) \subseteq
B^i_\ep(p^i_{\smash{b^{i,c}}})$. If\/ $j=\min I_c$ then for each
$i\in I_c$ we have $p^j_{\smash{b^{j,c}}}\in\dot V^{ij}$
and\/~$p^i_{\smash{b^{i,c}}}=\dot\phi^{ij}(p^j_{\smash{b^{j,c}}})$.

By taking $\ep>0$ sufficiently small we can choose $\bs T$ to make
the $X_c$ and\/ $V^i_c$ `arbitrarily small', in the sense that given
open covers $\{U_j:j\in J\}$ for $X$ and\/ $\{U_j^i:j\in J^i\}$ for
$\smash{\overline{\dot V^i}}\subseteq V^i$ for each\/ $i\in I,$ we
can choose $\bs T$ such that for all\/ $c\in C$ there exists\/ $j\in
J$ with\/ $\pi(X_c)\subseteq U_j,$ and for all\/ $i\in I_c$ there
exists\/ $j\in J^i$ with\/~$\pi(V^i_c)\subseteq U_j^i$.

In the case of \begin{bfseries}effective\end{bfseries} gauge-fixing
data $\ubG,\ubH,\ubG_c$ the Kuranishi space $X_c$ has trivial
stabilizers for all\/ $c\in C,$ and the $V^i_c$ are manifolds for
all\/~$i\in I_c$.
\label{khAprop10}
\end{prop}

\begin{proof} The first part follows quickly from the construction.
To make sense of it, note that $W^{i+1}$ in \eq{khAeq27} has
boundary composed of $\{0\}\t V^i$, and a portion lying over $\pd
V^i$, and a connected component for each $p^i_b$, $b=1,\ldots,N^i$,
and possibly also a piece from the branch $1+\ep^2$ in \eq{khAeq33}.
However, when we apply Algorithm \ref{kh3alg} and pass to $\check
W^{i+1}$, we discard boundary components not intersecting
$(t^{i+1})^{-1}(0)$. This includes the whole of the piece from the
branch $1+\ep^2$, since Definition \ref{khAdef13}(a) implies that
$\min T^i<1+\ep^2$ on $(s^i)^{-1}(0)\cap\dot V^i$. It may also
include some components from the $p^i_b$. Hence $\pd\check W^{i+1}$
is the union of $\{0\}\t V^i$, a portion lying over $\pd V^i$, and a
component associated to $p^i_b$ for a subset of
$b\in\{1,\ldots,N^i\}$. Definition \ref{khAdef13}(b) then implies
that the $X_c$ for $c\in C$ must correspond to points
$p_{\smash{b^{i,c}}}^i$ as in the proposition. In particular, no
$X_c$ corresponds to the branch $1+\ep^2$ in~\eq{khAeq33}.

For the second part, note that each $\pi(V^i_c)$ is contained in
$B_\ep^i(p^i_{\smash{b^{i,c}}})$ and intersects $(s^i)^{-1}(0)\cap
\dot V^i$, and $(s^i)^{-1}(0)\cap\overline{\dot V^i}$ is compact.
Given any open cover $\{U_j^i:j\in J^i\}$ for $\overline{\dot V^i}$,
as $(s^i)^{-1}(0)\cap\overline{\dot V^i}$ is compact, if $\ep>0$ is
sufficiently small then every ball of radius $\ep$ intersecting
$(s^i)^{-1}(0)\cap\overline{\dot V^i}$ must lie in $U_j^i$ for some
$j\in J^i$. Hence, if $\ep>0$ is sufficiently small then each
$\pi(V^i_c)$ lies in some $U_j^i$. For $X$ we use a similar
argument, noting that $\pi(X_c)$ has diameter at most $2\ep$ in the
metric on $X$ obtained by gluing together the metrics
$\psi^i_*(d^i)$ on $\Im\psi^i\subseteq X$ for $i\in I$. For the last
part, as $\ubG$ is effective gauge-fixing data the $V^i$ are
effective orbifolds, so the $V^i_c$ are manifolds as in Proposition
\ref{khAprop8}, and thus $X_c$ has trivial stabilizers as
$\Stab_{X_c}(p)=\Stab_{V^i_c}(v)=\{1\}$ for $p$ in $X_c$ and $v$ in
$(s^i_c)^{-1}(0)\subseteq V^i_c$ with~$p=\psi^i_c(v)$.
\end{proof}

\subsubsection{Piecewise smooth extensions on Kuranishi chains}
\label{khA34}

We extend \S\ref{khA14} to Kuranishi chains. Here is our analogue of
Definition~\ref{khAdef5}.

\begin{dfn} Suppose $X$ is a compact Kuranishi space, $Y$ an
orbifold, $\bs f:X\ra Y$ strongly smooth, and $\bs G$ (or $\ubG$) is
(effective) gauge-fixing data for $(X,\bs f)$. Let $\bs
I=(I,(V^i,E^i, s^i,\psi^i):i\in I,\ldots)$ and $\bs{\dot I}=(I,(\dot
V^i,\dot E^i,\dot s^i,\dot\psi^i):i\in I,\ldots)$ be as in
Definition~\ref{khAdef11}.

We will choose a smooth vector field $v^i$ on $V^i$ for each $i\in
I$ such that:
\begin{itemize}
\setlength{\itemsep}{0pt}
\setlength{\parsep}{0pt}
\item[(a)] $v^i$ is nonzero and {\it inward pointing} at every point
of~$\pd V^i$.
\item[(b)] For any $x_0\in V^i$, there exists a smooth path
$x:[0,\iy)\ra V^i$ such that $x(0)=x_0$ and $\frac{\d x}{\d
t}(t)=v^i\vert_{x(t)}$ for all $t\in[0,\iy)$. That is, $v^i$ is a
{\it complete} vector field on $V^i$. This implies that
$\exp(tv^i):V^i\ra V^i$ is a well-defined smooth embedding for
all~$t\ge 0$.
\item[(c)] There exist open neighbourhoods $U^i$ of
$\overline{\dot V^i}$ in $V^i$ for each $i\in I$, such that whenever
$j\le i$ in $I$ we have $\phi^{ij}_*(v^j)\equiv v^i$
on~$\phi^{ij}(V^{ij}\cap U^j)\cap U^i$.
\item[(d)] The $U^i$ in (c) have the following property: if $x\in
V^i$ and $\exp(tv^i)\in\overline{\dot V^i}$ for some $t\in[0,2]$,
then~$x\in U^i$.
\end{itemize}

Before doing this, we pause to explain the issues involved. One of
our problems is ensuring that the exponential maps
$\exp(tv^i):V^i\ra V^i$ are well-defined, at least where we need
them. If $X$ is a compact orbifold without boundary and $v$ is a
smooth vector field on $X$ then $\exp(tv):X\ra X$ is a well-defined
diffeomorphism for all $t\in\R$. If $X$ is a compact orbifold with
boundary or (g-)corners and $v$ is a vector field on $X$ then
$\exp(tv):X\ra X$ is a well-defined embedding for all $t\ge 0$
provided $v$ is inward-pointing along~$\pd X$.

However, if $X$ is a {\it noncompact\/} orbifold and $v$ a vector
field on $X$ then conditions for $\exp(tv):X\ra X$ to exist for
$t\ge 0$ are more complicated. In effect, $X$ has two kinds of
boundary: the usual boundary $\pd X$, along which $v$ must be
inward-pointing for $\exp(tv)$ to exist, plus some noncompact ends,
which we think of as being `at infinity' in $X$. If some flow-line
of $v$ in $X$ goes to infinity in finite time, then $\exp(tv)$ does
not exist for $t>0$. For instance, if $X=\R$ then $v=\frac{\pd}{\pd
x}$ is complete, with $\exp(tv):x\mapsto x+t$, but
$v=x^2\frac{\pd}{\pd x}$ is not complete, as the flow-line $t\mapsto
(1-t)^{-1}$ of $v$ for $t\in[0,1)$ goes to infinity in finite time.
To avoid this kind of behaviour we have to ensure that $\md{v}$ does
not grow too quickly near infinity in $X$, so that flow-lines of $v$
can only go to infinity in infinite time.

Unfortunately, this issue of completeness of vector fields causes
difficulties with compatibility of $v^i,v^j$ under $\phi^{ij}$. In
general, if $v^i,v^j$ are to be compatible on $U^i,U^j$ as in (c),
we cannot also choose $v^i,v^j$ to be complete on $U^i,U^j$, because
as $v^j$ must grow slowly near the noncompact boundary of $U^j$, we
find that the limit of $\phi^{ij}_*(v^j)$ at any point of
$\overline{\phi^{ij}(V^{ij}\cap U^j)}\sm \phi^{ij}(V^{ij}\cap U^j)$
must be zero, so $v^i$ must be zero there by continuity, but this
could contradict $v^i$ being nonzero at boundary points in (a).
Therefore we make $v^i$ complete on $V^i$ but not on $U^i$. Part (d)
will be used later to ensure that the tent functions $T^i$ on $\dot
V^i$ we construct using $v^i$ satisfy $\min T^j\vert_{\dot
V^{ij}}\equiv \min T^i\ci\dot\phi^{ij}$, as in
Definition~\ref{khAdef11}.

Next we explain how to choose $v^i,U^i$ satisfying (a)--(d).
Firstly, for each fixed $i\in I$, here is a way of constructing
$v^i$ satisfying (a),(b). Choose some Riemannian metric $g^i$ on
$V^i$. For each $x\in V^i$, define $f^i(x)=\sup\bigl\{\ep\in(0,1]:
\overline{B}_\ep^i(x)$ is a compact subset of $V^i\bigr\}$, where
$\overline{B}_\ep^i(x)$ is the closed ball of radius $\ep$ about $x$
in $(V^i,g^i)$. Then $f^i(x)$ is the minimum of 1 and the `distance
$d(x,\iy)$ from $x$ to infinity in $V^i$', measured using $g^i$, and
$f^i:V^i\ra(0,1]$ is continuous.

Let $v^i$ be a vector field on $V^i$ which is nonzero and inward
pointing at every point of $\pd V^i$. Then a sufficient condition
for $v^i$ to satisfy (b) is that $\md{v^i(x)}\le Cf^i(x)$ for some
$C>0$ and all $x\in V^i$, where $\md{v^i(x)}$ is computed using
$g^i$. This is because if $d(x,\iy)=\iy$ then as flow of $v$ moves
at speed at most $C$ then it must take infinite time to reach
infinity, and if $d(x,\iy)$ is finite then $\log d(x,\iy)$ decreases
with speed at most $C$, so $\log d(x,\iy)\ra -\iy$ in infinite time,
and the flow of $v$ takes infinite time to reach infinity in~$V^i$.

Thus, given any $v^i$ satisfying (a), which are easy to make by
joining local choices together by a partition of unity, we can make
it satisfy (b) as well by multiplying it by a smooth function
$F:V^i\ra(0,1]$ which approaches zero sufficiently fast at infinity
that $\md{F(x)v^i(x)}\le Cf^i(x)$ for all $x\in V^i$ and some $C>0$.
Note too that if $v^i$ satisfies (a),(b) and $F:V^i\ra(0,1]$ is
smooth then $Fv^i$ also satisfies (a),(b), since the flow-lines of
$Fv^i$ are just the flow-lines of $v^i$ reparametrized to move more
slowly, so they still must take infinite time to reach infinity
in~$V^i$.

Secondly, by induction on increasing $i\in I$, we choose $v^i,U^i$
for each $i\in I$ satisfying (a)--(c). In the inductive step, having
chosen $v^j,U^j$ for all $j<i$, choosing smooth $v^i,U^i$ satisfying
(c) is possible as in Remark \ref{khArem2}, making the $U^j$ for
$j<i$ smaller if necessary as in Definition \ref{khAdef13}. To
choose $v^i$ to satisfy (a) is also possible, since the prescribed
values $\phi^{ij}_*(v^j)$ for $\phi^{ij}(V^{ij}\cap U^j)\cap U^i$
satisfy (a) by (a) for $v^j$, and (a) is an open convex condition.
To make $v^i$ satisfy (b) as well, we make some initial choice not
necessarily satisfying (b), and then multiply it by some smooth
function $F:V^i\ra(0,1]$ which is 1 on $\phi^{ij}(V^{ij}\cap
U^j)\cap U^i$ for each $j<i$ in $I$ and decreases sufficiently fast
near infinity in $V^i$. For portions of $\phi^{ij}(V^{ij}\cap
U^j)\cap U^i$ which go to infinity in $V^i$, the prescribed values
for $v^i$ are already complete by (b) for $v^j$, so it is not
necessary to make $v^i$ smaller there.

Finally, having chosen $v^i,U^i$ for $i\in I$ satisfying (a)--(c)
but not necessarily (d), we fix these $U^i$, and choose some smooth
functions $F^i:V^i\ra(0,1]$ for $i\in I$ with $F^i\ci\phi^{ij}\equiv
F^j$ on $V^{ij}\cap U^j\cap(\phi^{ij})^{-1}(U^i)$, such that $\ti
v^i=F^iv^i$ for $i\in I$ satisfy (d) with the given $U^i$. This
holds provided $F^i$ is sufficiently small; what we need is that the
flow of $-\ti v^i=-F^iv^i$ cannot start in $\overline{\dot V^i}$ and
reach $V^i\sm U^i$ within time 2. Hence we can choose $v^i,U^i$ for
$i\in I$ satisfying (a)--(d) above.

We now continue the argument of Definition \ref{khAdef5}. Write
$V_1^i=\exp(v^i)V^i$ and $V_2^i=\exp(2v^i)V^i$. Then $V_1^i,V_2^i$
are embedded submanifolds of $(V^i)^\ci$, and $\exp(v^i):V^i\ra
V_1^i$, $\exp(2v^i):V^i\ra V_2^i$ are diffeomorphisms. The inclusion
$\io^i:\pd V_2^i\ra V^i$ is an immersion, with closed image. As for
$Y$ in Definition \ref{khAdef5}, for each $i\in I$ we extend $\pd
V_2^i$ to an $(i-1)$-orbifold $Y^i$ without boundary, with $\pd
V_2^i\subset Y$ a closed embedded $(i-1)$-submanifold, and extend
$\io^i:\pd V_2^i\ra V^i$ to an immersion $\io^i:Y^i\ra V^i$
satisfying Definition \ref{khAdef5}(i)--(v) with
$v^i,V^i,V^i_1,V^i_2$ in place of $v,X,X_1,X_2$, and where we only
need (iii) to hold locally in $V^i$, that is, such a bound $M$
should exist for each compact subset of $V^i$. In addition, we
require the $Y^i,\io^i$ to be compatible with coordinate changes
$\phi^{ij}$, in the sense that if $j<i$ in $I$ then defining
$Y^{ij}=(\io^j)^{-1}\bigl(V^{ij}\cap U^j\cap
(\phi^{ij})^{-1}(U^i)\bigr)$, and open set in $Y^j$, there should
exist an embedding $\psi^{ij}:Y^{ij}\ra Y^i$ such that
$\io^i\ci\psi^{ij}\equiv \phi^{ij}\ci\io^j$ on $Y^{ij}$. By Remark
\ref{khArem2} and the proof in Definition \ref{khAdef5}, we can
choose such $Y^i,\io^i$ for all $i\in I$. As in \eq{khAeq10}, define
$T^i:\dot V^i\ra F\bigl([1,\iy)\bigr)$ by
\e
\ts T^i(x)=\bigl\{2\bigr\}\cup\bigl\{3-t:t\in(0,2],\;\>
\exp(tv^i)x\in\io^i(Y^i)\bigr\}
\label{khAeq34}
\e
for each $i\in I$. Write~$\bs T=(T^i:i\in I)$.

The proof of Proposition \ref{khAprop3} shows that $T^i$ is a {\it
tent function} on $\dot V^i$ for each $i\in I$. It extends to a tent
function on $V^i$ also given by \eq{khAeq34}, so $T^i$ certainly
extends to a tent function on an open neighbourhood of
$\overline{\dot V^i}$ in $V^i$. The construction implies that if
$j\le i$ in $I$ then $T^j\vert_{\dot V^{ij}}\equiv
T^i\ci\dot\phi^{ij}$, and thus $\min T^j\vert_{\dot V^{ij}}\equiv
\min T^i\ci\dot\phi^{ij}$. Note that we use (d) above to prove this:
\eq{khAeq34} and (d) imply that for each $x\in\dot V^{ij}$, $T^j(x)$
depends only on the values of $v^j$ on $U^j$, rather than on the
whole of $V^j$, and then (c) shows that computing $T^j(x)$ in $V^j$
is equivalent to computing $T^i\ci\dot\phi^{ij}(x)$ in $V^i$.

The suborbifolds $S_{\{i_1,\ldots,i_l\}}$ of Definition
\ref{khAdef1} for $T^i$ are all locally either the images of some
codimension $k$ corner of $V^i$ under the flow of $v^i$, or the
intersection of such an image with $\io^i(Y^i)$. Since
$\dot\phi^{ij}:\dot V^{ij}\ra\dot V^i$ is an embedding which takes
$v^j$ to $v^i$, it follows that the $S_{\{i_1,\ldots,i_l\}}$ for
$T^i$ intersect $\dot\phi^{ij}(\dot V^{ij})$ transversely. Hence
$\bs T$ is a {\it tent function} for $(X,\bs f,\bs G)$ or $(X,\bs
f,\ubG)$ by Definition~\ref{khAdef11}.
\label{khAdef14}
\end{dfn}

We can also combine Definitions \ref{khAdef6} and \ref{khAdef14} to
extend a tent function for $(\pd X,\bs f\vert_{\pd X},\bs
G\vert_{\pd X})$ to one for~$(X,\bs f,\bs G)$.

\begin{dfn} Let $X$ be a compact Kuranishi space, $Y$ an orbifold,
$\bs f:X\ra Y$ strongly smooth, and $\bs G$ or $\ubG$ (effective)
gauge-fixing data for $(X,\bs f)$. Suppose $\bs T'$ is a tent
function for $(\pd X,\bs f\vert_{\pd X},\bs G\vert_{\pd X})$ or
$(\pd X,\bs f\vert_{\pd X},\ubG\vert_{\pd X})$ such that $\min(\bs
T')\vert_{\pd^2X}$ is invariant under the natural involution
$\bs\si:\pd^2X\ra\pd^2X$. That is, for $\bs I=(I,(V^i,E^i,s^i,
\psi^i):i\in I,\ldots)$ in $\bs G$ or $\ubG$, for each $i\in I$ the
tent function $\bs T'$ gives a tent function $T^{\prime (i-1)}$ on
$\pd\dot V^i$, and we require that $\min T^{\prime (i-1)}
\vert_{\pd^2\dot V^i}$ should be invariant under the involution
$\si:\pd^2\dot V^i\ra\pd^2\dot V^i$. For simplicity, suppose also
that $\min T^{\prime (i-1)}$ is bounded above on $\pd\dot V^i$ for
all $i\in I$. This holds automatically in all our constructions.

We will construct a tent function $\bs T$ for $(X,\bs f,\bs G$ or
$\ubG)$ with $\min\bigl(\bs T\vert_{\pd X}\bigr)=\min\bs T'$. Choose
$v^i,U^i,Y^i,\io^i$ for $i\in I$ as in Definition \ref{khAdef14}.
Then we have a diffeomorphism $\exp(2v^i):\pd V^i\ra\pd V^i_2$, so
that $T^{\prime(i-1)}\ci\exp(2v^i)$ is a tent function on
$\exp(2v^i)(\pd\dot V^i)\subseteq\pd V^i_2\subseteq Y^i$. Extend
$T^{\prime(i-1)}\ci\exp(2v^i)$ to a tent function
$T^{\prime\prime(i-1)}$ on an open neighbourhood $Y^{\prime i}$ of
$\exp(2v^i)(\overline{\pd\dot V^i})$ in $Y^i$, relaxing the
condition $\min T^{\prime\prime(i-1)}\ge 1$. We require too that
$\min T^{\prime\prime (i-1)}$ should be bounded above on $Y^{\prime
i}$. This is possible, as $\min T^{\prime(i-1)}$ is bounded above.

Choose $D>1$ such that $\min T^{\prime\prime (i-1)}<D$ on $Y^{\prime
i}$ for all $i\in I$. As in \eq{khAeq13}, define $T^i:\dot V^i\ra
F\bigl([1,\iy)\bigr)$ by
\e
\begin{split}
\ts T^i(x)=\bigl\{D\bigr\}\cup\bigl\{u+\th(2-t): \text{$t\in(0,2]$,
$y\in Y^{\prime i}$, $u\in T^{\prime\prime(i-1)}(y)$,}&\\
\exp(tv^i)x=\io^i(y)\bigr\},&
\end{split}
\label{khAeq35}
\e
where $\th>D-1$ is chosen large enough to satisfy conditions in
Proposition \ref{khAprop5} for each $i\in I$. Actually, because of
the noncompactness of $\dot V^i$ it may not be possible to choose
$\th$ sufficiently large everywhere. We can deal with this by
replacing $\th$ by smooth functions $\th^i:\dot V^i\ra (D-1,\iy)$
for $i\in I$ with $\th^i\ci\dot\phi^{ij}\equiv\th^j\vert_{\dot
V^{ij}}$ for $j<i$ in $I$, which are sufficiently large locally
in~$\dot V^i$.

The proof of Proposition \ref{khAprop5} now shows that $T^i$ is a
tent function on $\dot V^i$ for each $i\in I$, with
$\min\bigl(T^i\vert_{\pd\dot V^i}\bigr)\equiv\min T^{\prime(i-1)}$.
It extends to a tent function on an open neighbourhood of
$\overline{\dot V^i}$ in $V^i$, given by the same formula
\eq{khAeq35}. For $j<i$ in $I$ we find as in Definition
\ref{khAdef14} that $\min T^j\vert_{\dot V^{ij}}\equiv \min
T^i\ci\dot\phi^{ij}$, and the $S_{\{i_1,\ldots,i_l\}}$ for $T^i$
intersect $\dot\phi^{ij}(\dot V^{ij})$ transversely. Therefore $\bs
T=(T^i:i\in I)$ is a {\it tent function} for $(X,\bs f,\bs G)$ or
$(X,\bs f,\ubG)$, with $\min\bigl(\bs T\vert_{\pd X}\bigr)=\min\bs
T'$.
\label{khAdef15}
\end{dfn}

Analogues of Definition \ref{khAdef4} and Propositions
\ref{khAprop3}, \ref{khAprop4}, \ref{khAprop5} and \ref{khAprop6}
hold for $\bs T$ in Definitions \ref{khAdef14} and \ref{khAdef15},
so that we have a well-defined notion of {\it piecewise strongly
smooth functions subordinate to a tent function\/} $\bs T$ on
$(X,\bs f,\bs G)$, for which an Extension Principle holds. We will
not write these out, as they are cumbersome to state. However, Step
2 of the proof of Theorem \ref{kh4thm1} in \S\ref{khB2} depends on
this piecewise smooth Extension Principle for effective Kuranishi
chains, so we will discuss some of the ideas there.

\subsubsection{Cutting into small pieces with boundary conditions}
\label{khA35}

Section \ref{khA15} combined \S\ref{khA12} and \S\ref{khA14} to
define tent functions $T$ cutting a compact manifold $X$ into
arbitrarily small pieces, with boundary conditions $\min T\vert_{\pd
X}\equiv\min T'$ over $\pd X$. We now combine \S\ref{khA33} and
\S\ref{khA34} to do the same thing for Kuranishi chains. Here is our
analogue of Definition~\ref{khAdef7}.

\begin{dfn} Let $X,Y,\bs f,\bs G$ or $\ubG$ and $\bs T'$ be as in
Definition \ref{khAdef15}. We will construct a tent function $\bs T$
for $(X,\bs f,\bs G$ or $\ubG)$ with $\min\bigl(\bs T\vert_{\pd
X}\bigr)=\min\bs T'$, such that the pieces $X_c$ in
\eq{khAeq29}--\eq{khAeq30} can be made `arbitrarily small', subject
to the boundary conditions~$\bs T'$.

Choose metrics $g^i$ on $U^i\supset\overline{\dot V^i}$ for $i\in
I$, $\ep>0$, finite subsets $\{p_1^i,\ldots,p_{N^i}^i\}$ in $\dot
V^i$ for $i\in I$, as in Definition \ref{khAdef13}, with the extra
condition that $p_c^i\in(\dot V^i)^\ci$ for all $i,c$. Choose
$v^i,U^i,Y^i,\io^i$ for $i\in I$ as in Definition \ref{khAdef14}.
Choose $Y^{\prime i},T^{\prime\prime(i-1)}$ for $i\in I$ and $D>1$,
$\th>D-1$ or $\th^i:\dot V^i\ra (D-1,\iy)$ as in Definition
\ref{khAdef15}. Make $\th$, or the functions $\th^i$ for $i\in I$,
larger if necessary so that $\th>D+\ep^2-1$, and an analogue of the
condition $1+\th\ze\Vert v\Vert_{C^0}>D$ in Definition \ref{khAdef7}
holds locally in $\dot V^i$. Combining \eq{khAeq22}, \eq{khAeq33}
and \eq{khAeq35}, define $T^i:\dot V^i\ra F\bigl([1,\iy)\bigr)$ by
\e
\begin{split}
&T^i(v)\!=\!\bigl\{D\!+\!l^2:\text{$l\in[0,2\ep)$, there exists a
geodesic segment of length $l$}\\
&\;\>\text{in $(\dot V^i,g^i)$ with end points $v$ and $p_b^i$, some
$b=1,\ldots,N^i$}\bigr\}\!\cup\!\{D\!+\!\ep^2\}\cup\\
&\;\>\bigl\{u+\th(2-t): \text{$t\in(0,2]$, $y\in Y^{\prime i}$,
$u\in T^{\prime\prime(i-1)}(y)$, $\exp(tv^i)x=\io^i(y)$}\bigr\},
\end{split}
\label{khAeq36}
\e
for each $i\in I$. Write~$\bs T=(T^i:i\in I)$.
\label{khAdef16}
\end{dfn}

The proofs of Propositions \ref{khAprop7}, \ref{khAprop10} and
Definitions \ref{khAdef13}, \ref{khAdef15} then yield:

\begin{thm} In Definition {\rm\ref{khAdef16},} $\bs T$ is a tent
function for $(X,\bs f,\bs G)$ or $(X,\bs f,\ubG),$
with\/~$\min\bigl(\bs T\vert_{\pd X}\bigr)=\min\bs T'$.

Let\/ $X_c,\bs f_c,\bs G_c$ for $c\in C$ be as in Proposition
{\rm\ref{khAprop9}} for $Z_{X,\bs T},$ and\/ $X'_{c'},\bs
f'_{c'},\ab\bs G'_{c'}$ for $c'\in C'$ be as in Proposition
{\rm\ref{khAprop9}} for $Z_{\pd X,\bs T'}$. Then we can take
$C'\subseteq C,$ such that for each\/ $c'\in C'$ there is a natural
identification between $(X'_{c'},\bs f'_{c'},\bs G'_{c'})$ and a
component of\/ $(\pd X_{c'},\bs f_{c'}\vert_{\pd X_{c'}},\bs
G_{c'}\vert_{\pd X_{c'}})$. For $c\in C\sm C',$ the $X_c$ correspond
to points $p^i_{\smash{b^{i,c}}}$ as in Proposition
{\rm\ref{khAprop10}}. The same holds for~$\ubG_c,\ubG'_{c'}$.

By taking $\ep>0$ sufficiently small and\/ $\th^i$ for $i\in I$
sufficiently large, we can choose $\bs T$ to make the $X_c$ and\/
$V^i_c$ `arbitrarily small', given the $X'_{c'},V^{\prime i}_{c'}$
which are already fixed. That is, given an open cover $\{U_j:j\in
J\}$ for $X$ such that for all\/ $c'\in C'$ there exists\/ $j\in J$
with\/ $\pi'(X'_{c'})\subseteq U_j,$ then we can choose\/ $\bs T$ so
that for all\/ $c\in C$ there exists\/ $j\in J$ with\/
$\pi(X_c)\subseteq U_j,$ and similarly for the~$V^i_c,V^{\prime
i}_{c'}$.

In the case of \begin{bfseries}effective\end{bfseries} gauge-fixing
data $\ubG,\ubH,\ubG_c,$ if the $X'_{c'}$ have trivial stabilizers
and the $V^{\prime i}_{c'}$ are manifolds for all\/ $c'\in C'$ and\/
$i\in I'_{c'},$ then the $X_c$ have trivial stabilizers and the
$V^i_c$ are manifolds for all\/ $c\in C$ and\/~$i\in I_c$.
\label{khAthm5}
\end{thm}

The tent functions $\bs T$ of Definition \ref{khAdef16} are also
suitable for piecewise smooth extensions, as in \S\ref{khA34}, and
we will use them for this in Step 2 of the proof of Theorem
\ref{kh4thm1} in~\S\ref{khB2}.

\begin{rem} It is instructive to compare our notion of tent
functions on Kuranishi chains $(X,\bs f,\bs G)$ with Fukaya and
Ono's notion of {\it multisection} on a Kuranishi space $X$ \cite[\S
3, \S 6]{FuOn1}, \cite[\S A1.1]{FOOO}. There are many similarities,
as in writing this book I started with Fukaya and Ono's definitions
of Kuranishi spaces and virtual cycles and tried to improve them.
Here are some:
\begin{itemize}
\setlength{\itemsep}{0pt}
\setlength{\parsep}{0pt}
\item In both cases we fix a good coordinate system
$\bs I=\bigl(I,(V^i,\ldots,\psi^i):i\in I,\ldots\bigr)$ on $X$, and
then consider a {\it smooth, multivalued\/} function $T^i$ on $\dot
V^i\subset V^i$, or section $\bs s^i$ of $E^i\ra V^i$ for each $i\in
I$, which must be compatible with coordinate changes $\phi^{ij}$ for
$j<i$ in~$I$.
\item Fukaya and Ono use their multisections to perturb Kuranishi
spaces to (non-Hausdorff, $\Q$-weighted) manifolds, and then to
triangulate these manifolds by simplices. We use tent functions for
both these purposes.
\end{itemize}

However, there are also important differences:
\begin{itemize}
\setlength{\itemsep}{0pt}
\setlength{\parsep}{0pt}
\item For a multisection with $n$-branches, each branch has equal weight
$1/n$. Therefore by its nature the multisection method produces
singular chains over $\Q$, not over $\Z$. In contrast, though a tent
function $T$ may have many smooth branches, only the minimal branch
$\min T$ really counts, and it has weight 1. Because of this, we are
able to make effective Kuranishi (co)homology work over $\Z$ as well
as~$\Q$.
\item Fukaya and Ono use {\it multivalued\/} sections of {\it orbifold
vector bundles} so that they can perturb to transverse sections,
which is generally not possible for single-valued sections of
orbifold vector bundles.

For {\it effective} Kuranishi neighbourhoods $(V,E,s,\psi)$,
single-valued sections of $E$ can always be perturbed to transverse,
so this problem does not arise in Appendix \ref{khB}, where we
actually do our perturbations.

However, in Appendix \ref{khC} we use a tent function to cut a
general Kuranishi space $X$ into small pieces $X_c$ for $c\in C$
which are quotients $X_c'/\Ga_c$ for $X_c'$ an effective Kuranishi
space and $\Ga_c$ a finite group, and then we replace $X_c$ by
$X'_c$ with weight $1/\md{\Ga_c}$ in $\Q$. So we are forced to work
over $\Q$ in Appendix \ref{khC}, for essentially the same reason as
Fukaya and Ono.
\end{itemize}
\label{khArem3}
\end{rem}

\section{Proof that $H^{\rm si}_*(Y;R)\cong KH_*^{\rm ef}(Y;R)$}
\label{khB}

We now prove Theorem \ref{kh4thm1}. Let $Y$ be an orbifold and $R$ a
commutative ring. We must show that $\Pi_\rsi^\ef:H_k^\rsi(Y;R)\ra
KH_k^\ef(Y;R)$ in \eq{kh4eq18} is an isomorphism. We do this in four
steps, which we describe briefly below and then cover in detail in
\S\ref{khB1}--\S\ref{khB4}. In Steps 1--3 we start with an arbitrary
cycle $\sum_{a\in A}\rho_a[X_a,\bs f_a,\ubG_a]$ in $KC_k^\ef(Y;R)$
representing some $\al\in KH_k^\ef(Y;R)$, and progressively modify
it to get other homologous cycles representing $\al$ with better and
better properties, until we represent $\al$ by the image under
$\Pi_\rsi^\ef$ in \eq{kh4eq15} of a cycle in $C_k^\rsi(Y;R)$,
proving that $\Pi_\rsi^\ef$ in \eq{kh4eq18} is {\it surjective}.
Step 4 proves $\Pi_\rsi^\ef$ is {\it injective}, using the same
method but for chains with boundary rather than cycles. We allow
$k<0$ in Steps 1 and 2, and prove the last part of Theorem
\ref{kh4thm1} at the end of Step~2.
\medskip

\noindent{\bf Step 1.} Let $k\in\Z$, $\al\in KH_k^\ef(Y;R)$, and
$\sum_{a\in A}\rho_a[X_a,\bs f_a,\ubG_a]\in KC_k^\ef(Y;R)$ represent
$\al$. Then each $X_a$ is an {\it effective Kuranishi space}, and in
the excellent coordinate system $\bs I_a=\bigl(I_a,(V^i_a,\ldots,
\psi^i_a):i\in I_a,\ldots\bigr)$ for $(X_a,\bs f_a)$ in $\ubG_a$,
$V^i_a$ is an {\it effective orbifold}. As in \S\ref{khA3}, we use a
{\it tent function} $\bs T_a$ for $(X_a,\bs f_a,\ubG_a)$ to `cut'
each $X_a$ into finitely many smaller pieces $X_{ac}$ for $c\in
C_a$, in triples $(X_{ac},\ab\bs f_{ac},\ab\ubG_{ac})$, where
$X_{ac}$ is a Kuranishi space with {\it trivial stabilizers}, and in
the excellent coordinate system $\bs I_{ac}=\bigl(I_{ac},(V^i_{ac},
\ldots,\psi^i_{ac}):i\in I_{ac},\ldots\bigr)$ for $(X_{ac},\bs
f_{ac})$ in $\ubG_{ac}$, each $V^i_{ac}$ is a {\it manifold}. We
show $\sum_{a\in A}\sum_{c\in C_a}\rho_a[X_{ac},\ab\bs
f_{ac},\ab\ubG_{ac}]$ is homologous to $\sum_{a\in A}\rho_a[X_a,\bs
f_a,\ubG_a]$ in $KC_k^\ef(Y;R)$, and so represents~$\al$.

As $\pd\bigl(\sum_{a\in A}\rho_a[X_a,\bs f_a,\ubG_a]\bigr)=0$, the
$[\pd X_a,\bs f_a\vert_{\pd X_a},\ubG_a\vert_{\pd X_a}]$ must
satisfy relations in $KC_{k-1}^\ef(Y;R)$. The main problem is to
choose the $\bs T_a$ for $a\in A$ in a way compatible with these
relations, so that $\pd\bigl(\sum_{a\in A}\sum_{c\in
C_a}\ab\rho_a[X_{ac},\ab\bs f_{ac},\ab\ubG_{ac}]\bigr)=0$. That is,
the choices of $\bs T_a,C_a$ and $X_{ac},\bs f_{ac},\ubG_{ac}$ for
$c\in C_a$ cannot be made independently for each $a\in A$, but
rather, they must satisfy complicated boundary conditions over $\pd
X_a$ relating the choices for different~$a\in A$.

To make these choices in a way compatible with the boundary
conditions, we introduce some more notation. By Lemma \ref{kh3lem},
for each $a\in A$ and $m\ge 0$ the triple $(\pd^mX_a,\bs
f_a\vert_{\pd^mX_a},\ab\ubG_a\vert_{\pd^mX_a})$ has a splitting
$\pd^mX_a=X^m_{a1}\amalg\cdots\amalg X^m_{an_a^m}$ into {\it
connected\/} triples $(X_{ab}^m,\bs f_a\vert_{X_{ab}^m},
\ubG_a\vert_{X_{ab}^m})$, which is unique up to order of the
$X^m_{ab}$ for $b=1,\ldots,n_a^m$. Write $\bs f_{ab}^m=\bs
f_a\vert_{X_{ab}^m}$ and $\ubG_{ab}^m=\ubG_a\vert_{X_{ab}^m}$. Then
$\Aut(X_{ab}^m,\bs f_{ab}^m,\ubG_{ab}^m)=\{1\}$ by
Theorem~\ref{kh3thm5}(b).

Define $P^m=\bigl\{(a,b):a\in A$, $b=1,\ldots,n_a^m\bigr\}$. Define
an equivalence relation $\sim$ on $P^m$ by $(a,b)\sim(a',b')$ if
there exists an isomorphism $(\bs a,\bs b):(X_{ab}^m,\bs
f_{ab}^m,\ubG_{ab}^m)\ra(X_{a'b'}^m,\bs f_{a'b'}^m,\ubG_{a'b'}^m)$
as in Definition \ref{kh3def17}, where $\bs a$ need not identify
orientations. Note that $(\bs a,\bs b)$ is unique as
$\Aut(X_{ab}^m,\bs f_{ab}^m,\ubG_{ab}^m)=\{1\}$. Let $Q^m=P^m/\sim$
be the set of equivalence classes of $(a,b)$ in $P^m$. We shall
treat $Q^m$ as an {\it indexing set}. For each $q\in Q^m$, choose
$(a^q,b^q)\in P^m$ representing the equivalence class $q$. Write
$R^m=\bigl\{(a^q,b^q):q\in Q^m\bigr\}$. Define $\phi^m:P^m\ra R^m$
by $\phi^m(a,b)=(a^q,b^q)$ when the $\sim$-equivalence class of
$(a,b)$ in $Q^m$ is $q$. Write $(\bs a,\bs b)_{ab}^m$ for the unique
isomorphism $(\bs a,\bs b):(X_{ab}^m,\bs f_{ab}^m,\ubG_{ab}^m)
\ra(X_{a^qb^q}^m,\bs f_{a^qb^q}^m, \ubG_{a^qb^q}^m)$.

Since $A$ is finite and $\pd^mX_a=\es$ for $m\gg 0$, $P^m,Q^m,R^m$
are finite with $P^m=Q^m=R^m=\es$ for $m\gg 0$. Let $M\ge 0$ be
largest with $P^M,Q^M,R^M\ne\es$. Write $\ep_{ab}^m=1$ if $\bs a$ in
$(\bs a,\bs b)_{ab}^m$ is orientation-preserving, and
$\ep_{ab}^m=-1$ if $\bs a$ is orientation-reversing. As
$X_{ab}^m\subseteq\pd^mX_a$, we have $\pd X_{ab}^m\subseteq
\pd^{m+1}X_a=\coprod_{b'=1}^{n_a^{m+1}}X_{ab'}^{m+1}$. Clearly, $\pd
X_{ab}^m$ is a disjoint union of some subset of the $X_{ab'}^{m+1}$.
Define $B_{ab}^m$ to be the set of $b'=1,\ldots,n_a^{m+1}$ with
$X_{ab'}^{m+1}\subseteq\pd X_{ab}^m$. Then~$\pd
X_{ab}^m=\coprod_{b'\in B_{ab}^m} X_{ab'}^{m+1}$.

By induction on decreasing $m=M,M-1,\ldots,1,0$, we use Definition
\ref{khAdef16} and Theorem \ref{khAthm5} to choose a tent function
$\bs T_{ab}^m$ for $(X_{ab}^m,\bs f_{ab}^m,\ubG_{ab}^m)$ for all
$(a,b)$ in $R^m$, satisfying the following boundary condition: we
have $\pd X_{ab}^m=\coprod_{b'\in B_{ab}^m} X_{ab'}^{m+1}$, and for
each such $X_{ab'}^{m+1}$ we have $(\bar a,\bar
b)=\phi^{m+1}(a,b')\in R^{m+1}$, an isomorphism $(\bs a,\bs
b)^{m+1}_{ab'}:(X_{ab'}^{m+1},\ab \bs
f_{ab'}^{m+1},\ab\ubG_{ab'}^{m+1})\ab\ra(X_{\bar a\bar
b}^{m+1},\ab\bs f_{\bar a\bar b}^{m+1},\ab\ubG_{\bar a\bar
b}^{m+1})$, and a tent function $\bs T_{\bar a\bar b}^{m+1}$ for
$(X_{\bar a\bar b}^{m+1},\ab\bs f_{\bar a\bar b}^{m+1},\ab
\ubG_{\bar a\bar b}^{m+1})$ chosen in the previous inductive step.
We need $\min\bigl(\bs T_{ab}^m\vert_{\smash{X_{ab'}^{m+1}}} \bigr)$
to agree with $\min\bs T_{\bar a\bar b}^{m+1}$ under the isomorphism
$(\bs a,\bs b)^{m+1}_{ab'}$, for all $b'$ in $B_{ab}^m$. This is
possible as in~\S\ref{khA35}.

After completing the induction, for $a\in A$ we have
$X_a=\coprod_{b=1}^{n_a^0}X^0_{ab}$, and for each such $(a,b)$ we
have $(\bar a,\bar b)=\phi^0(a,b)\in R^0$, an isomorphism $(\bs
a,\bs b)^0_{ab}:(X_{ab}^0,\ab\bs f_{ab}^0,\ab\ubG_{ab}^0)\ra(X_{\bar
a\bar b}^0,\ab\bs f_{\bar a\bar b}^0,\ab\ubG_{\bar a\bar b}^0)$, and
a tent function $\bs T_{\bar a\bar b}^0$ for $(X_{\bar a\bar
b}^0,\bs f_{\bar a\bar b}^0,\ubG_{\bar a\bar b}^0)$. We define the
tent function $\bs T_a$ for $(X_a,f_a,\ubG_a)$ to be the pullback of
$\bs T_{\bar a\bar b}^0$ under $(\bs a,\bs b)^0_{ab}$ on $X_{ab}^0$,
for each $b=1,\ldots,n_a^0$. Then by construction, the $\bs
T_a\vert_{\pd X_a}$ for $a\in A$ satisfy relations corresponding
to~$\pd\bigl(\sum_{a\in A}\rho_a[X_a,\bs f_a,\ubG_a]\bigr)=0$.

Let $(X_{ac},\ab\bs f_{ac},\ab\ubG_{ac})$ for $c\in C_a$ be the
decomposition of $(X_a,\bs f_a,\ubG_a)$ into pieces defined using
$\bs T_a$ in \S\ref{khA32}. Theorem \ref{khAthm5} and induction on
$m$ implies that the $X_{ac}$ have trivial stabilizers, and the
$V^i_{ac}$ are manifolds. We show using Proposition \ref{khAprop9}
that $\sum_{a\in A}\sum_{c\in C_a}\rho_a[X_{ac},\ab\bs
f_{ac},\ab\ubG_{ac}]$ is homologous to $\sum_{a\in A}\rho_a[X_a,\bs
f_a,\ubG_a]$ in~$KC_k^\ef(Y;R)$.
\medskip

\noindent{\bf Step 2.} We shall make a small deformation of $X_{ac}$
into a {\it compact manifold\/ $\ti X_{ac}$ with g-corners}, for
each $a\in A$ and $c\in C_a$. This yields a new cycle $\sum_{a\in
A}\sum_{c\in C_a}\rho_a[\ti X_{ac},\ti f_{ac},\ab\ubtG_{ac}]\in
KC_k^\ef(Y;R)$ which is homologous to $\sum_{a\in A}\ab\sum_{c\in
C_a}\rho_a[X_{ac},\ab\bs f_{ac},\ab\ubG_{ac}]$, and so represents
$\al$. We write $\ti f_{ac}:\ti X_{ac}\ra Y$ rather than $\bs{\ti
f}_{ac}:\ti X_{ac}\ra Y$ since $\ti f_{ac}$ is a smooth map of
orbifolds, rather than a strongly smooth map of Kuranishi spaces.

The equation $\pd\bigl(\sum_{a\in A}\sum_{c\in
C_a}\ab\rho_a[X_{ac},\ab\bs f_{ac},\ab\ubG_{ac}]\bigr)=0$ yields
relations on the $[\pd X_{ac},\bs f_{ac}\vert_{\pd
X_{ac}},\ab\ubG_{ac}\vert_{\pd X_{ac}}]$ in $KC_{k-1}^\ef(Y;R)$. The
main problem is to choose the deformations $\ti X_{ac}$ in a way
compatible with these relations, so that $\pd\bigl(\sum_{a\in
A}\ab\sum_{c\in C_a}\rho_a[\ti X_{ac},\ti
f_{ac},\ubtG_{ac}]\bigr)=0$. That is, the perturbations $\ti X_{ac}$
must satisfy complicated boundary conditions. Before explaining how
we do this, we discuss how to perturb a single $(X_{ac},\bs
f_{ac},\ubG_{ac})$ without boundary conditions.

Let $(\bs I_{ac},\bs\eta_{ac})$ be the excellent coordinate system
for $(X_{ac},\bs f_{ac})$ in $\ubG_{ac}$, where $\bs I_{ac}=\bigl(
I_{ac},(V^i_{ac},E^i_{ac},s^i_{ac},\psi^i_{ac}),f^i_{ac}:i\in
I_{ac},\ldots\bigr)$, the $V^i_{ac}$ are manifolds, and each
$E^i_{ac}\ra V^i_{ac}$ is a vector bundle.

Working by induction on $i\in I_{ac}$, we will choose $C^1$ small
perturbations $\ti s^i_{ac}$ of $s^i_{ac}$ on $V^i_{ac}$ for $i\in
I_{ac}$ such that $\ti s^i_{ac}$ is {\it transverse\/} along $(\ti
s^i_{ac})^{-1}(0)$, that is, $\d\ti s^i_{ac}:TV^i_{ac}\ra E^i_{ac}$
is surjective on $(\ti s^i_{ac})^{-1}(0)$, and
$\hat\phi^{ij}_{ac}\ci\ti s^j_{ac}\equiv\ti
s^i_{ac}\ci\phi^{ij}_{ac}$ on $V^{ij}_{ac}$ whenever $j\le i$
in~$I_{ac}$.

By gluing together the sets $(\ti s^i_{ac})^{-1}(0)$ we make a {\it
compact, oriented\/ $k$-manifold with g-corners\/} $\ti X_{ac}$
which is a small transverse perturbation of $X_{ac}$, with a smooth
map $\ti f_{ac}:\ti X_{ac}\ra Y$, and an excellent coordinate system
$(\bs{\ti I}_{ac},\bs{\ti\eta}_{ac})$ for $(\ti X_{ac},\ti f_{ac})$
with indexing set $\ti I_{ac}=\{k\}$, one Kuranishi neighbourhood
$(\ti V_{ac}^k,\ab\ti E_{ac}^k,\ab\ti s_{ac}^k,\ti\psi_{ac}^k)=(\ti
X_{ac},\ti X_{ac},0,\id_{\ti X_{ac}})$ and $\ti\eta_{ac,k}\equiv
1\equiv\ti\eta_{ac,k}^k$. We choose $\utG^k_{ac}:\ti X_{ac}\ra \uP$
such that $\ubtG_{ac}\!=(\bs{\ti I}_{ac},\bs{\ti\eta}_{ac},
\utG^k_{ac})$ is effective gauge-fixing data for $(\ti X_{ac},\ti
f_{ac})$. We also construct a homology $[Z_{ac},\bs g_{ac},
\ubH_{ac}]$ between $[X_{ac},\ab\bs f_{ac},\ab\ubG_{ac}]$ and $[\ti
X_{ac},\ab\ti f_{ac},\ab\ubtG_{ac}]$ in $KC_{k+1}^\ef(Y;R)$, modulo
terms over~$\pd X_{ac}$.

Next we explain how to make these choices of perturbations $\ti
s^i_{ac}$ for $a\in A$ and $c\in C_a$ compatible with the boundary
conditions. Given how we chose the $\bs T_a$ compatible with
boundary conditions in Step 1, one obvious approach would be to use
the same method, and to work by induction on decreasing
$m=M,M-1,\ldots,0$ choosing the $\ti s^i_{ac}\vert_{\pd^mV^i_{ac}}$,
by splitting $\pd^m[X_{ac},\ab\bs f_{ac},\ab\ubG_{ac}]$ into
connected components, and ensuring our choices depend only on the
isomorphism class of the component.

However, there is a problem with this: as the $V^i_{ac}$ are
manifolds {\it with g-corners}, not corners, the Extension Principle
for smooth sections, Principle \ref{kh2pri}(c), does not hold. Thus,
given prescribed values for $\ti s^i_{ac}$ over $\pd V^i_{ac}$ which
are $\si$-invariant over $\pd^2V^i_{ac}$, we may {\it not\/} be able
to extend these prescribed values smoothly over $V^i_{ac}$ to choose
$\ti s^i_{ac}$ satisfying the necessary boundary conditions.

Our solution to this was explained in \S\ref{khA14}. There are
notions of {\it continuous, piecewise smooth functions and sections
subordinate to a tent function}, and an Extension Principle holds
for such functions and sections on manifolds and orbifolds with
g-corners. So, we could solve the problem by using continuous,
piecewise smooth sections $\ti s^i_{ac}$ on $V^i_{ac}$. That is, we
could choose tent functions $\bs T_{ac}$ for $(X_{ac},\ab\bs
f_{ac},\ab\ubG_{ac})$ which cut $X_{ac}$ into still smaller pieces
$X_{acd}$ for $d\in D_{ac}$, and continuous, piecewise smooth
sections $\ti s^i_{ac}$ over $V^i_{ac}$ subordinate to the tent
function $T^i_{ac}$ in $\bs T^i_{ac}$, satisfying the necessary
boundary conditions. The restriction of $\ti s^i_{ac}$ to
$V^i_{acd}$ would then be smooth for each~$d\in D_{ac}$.

In fact there is a simpler method. We do not need to choose more
tent functions $\bs T_{ac}$ and divide the $X_{ac}$ into smaller
pieces $X_{acd}$. Instead, we treat the smooth sections $\ti
s^i_{ac}$ on $V^i_{ac}$ for $c\in C_a$ as coming from a piecewise
smooth section $\ti s^i_a$ on $\dot V^i_a$ subordinate to $T^i_a$.
So, we can re-use the tent functions $\bs T_a$ chosen in Step 1 to
inductively construct piecewise smooth sections over $X_a$ and $\dot
V^i_a$ satisfying boundary conditions, which restrict to smooth
sections over $X_{ac}$ and $V^i_{ac}$ satisfying the necessary
boundary conditions.

The $\bs T_a$ are suitable for this purpose, as we were careful to
define them this way in \S\ref{khA35}. Also, despite the discussion
after Example \ref{khAex1} in \S\ref{khA3}, in this case we do not
need to define the $\ti s^i_{ac}$ on smaller subsets $\dot
V_{ac}^i\subset V_{ac}^i$, since the $V_{ac}^i$ are constructed from
subsets $\dot V_a^i$ in $V_a^i$, and by choosing the $\ti s^i_a$ to
extend to open neighbourhoods of $\overline{\dot V_a^i}$ in $V_a^i$
for $i\in I_a$, we can define $\ti s_a^i$ on all of $\dot V_a^i$ and
so define $\ti s^i_{ac}$ on all of~$V_{ac}^i$.

Here, then, is how we actually choose the $\ti s^i_{ac}$. Use the
notation of Step 1, and let the effective gauge-fixing data
$\ubG_{ab}^m$ for $(X_{ab}^m,\bs f_{ab}^m)$ contain $\bs I_{ab}^m=
\bigl((V^{i,m}_{ab},E^{i,m}_{ab},\ab s^{i,m}_{ab},\ab
\psi^{i,m}_{ab}):i\in I_{ab}^m,\ldots\bigr)$, and $\bs
T^m_{ab}=(T^{i,m}_{ab}:i\in I_{ab}^m)$, where $T^{i,m}_{ab}$ is a
tent function on $\dot V^{i,m}_{ab}\subset V^{i,m}_{ab}$. By
induction on decreasing $m=M,M-1,\ldots,1,0$, we choose a piecewise
smooth section $\ti s^{i,m}_{ab}$ of $E^{i,m}_{ab}$ over $\dot
V^{i,m}_{ab}$ subordinate to $T^{i,m}_{ab}$, which is a piecewise
$C^1$ small perturbation of $s^{i,m}_{ab}$, with $\ti s^{i,m}_{ab}$
prescribed over $\pd\dot V^{i,m}_{ab}$ in the same way as for the
boundary conditions for $\min T^{i,m}_{ab}$ in Step 1. This is
possible, in a similar way to~\S\ref{khA14}.

After completing the induction, we define piecewise smooth sections
$\ti s_a^i$ of $E^i_a$ over $\dot V^i_a$ from the $\ti s^{i,0}_{ab}$
subordinate to $T^i_a$, as for the construction of $\bs T_a$ from
$\bs T^0_{ab}$ in Step 1. These $\ti s_a^i$ pull back to smooth
sections $\ti s^i_{ac}$ on each $V^i_{ac}$, under the natural
embeddings $\pi_a:V^i_{ac}\ra\dot V^i_a$. The $\ti s^i_{ac}$ then
satisfy all the necessary boundary conditions.

We then show that $\sum_{a\in A}\sum_{c\in C_a}\rho_a[\ti X_{ac},\ti
f_{ac},\ab\ubtG_{ac}]\in KC_k^\ef(Y;R)$ is homologous to $\sum_{a\in
A}\sum_{c\in C_a}\rho_a[X_{ac},\ab\bs f_{ac},\ab\ubG_{ac}]$, and so
represents $\al$. If $k<0$ then $\ti X_a=\es$ as $\dim\ti X_a<0$, so
$\al=0$. Thus $KH_k^\ef(Y;R)=0$ for $k<0$, proving the last part of
Theorem \ref{kh4thm1}. Therefore we suppose $k\ge 0$ from now on.
\medskip

\noindent{\bf Step 3.} For simplicity we now {\it change notation\/}
from $\ti X_{ac}$ for $a\in A$ and $c\in C_a$ to $\ti X_a$ for $a\in
A$. That is, Steps 1 and 2 have shown that we can represent $\al$ by
a cycle $\sum_{a\in A}\rho_a[\ti X_a,\ti f_a,\ubtG_a]\in
KC_k^\ef(Y;R)$, with each $\ti X_a$ a compact, oriented $k$-manifold
with g-corners, $\ti f_a:\ti X_a\ra Y$ a smooth map, and
$\ubtG_a=(\bs{\ti I}_a,\bs{\ti\eta}_a,\utG^k_a)$ effective
gauge-fixing data for $(\ti X_a,\ti f_a)$, where $\bs{\ti I}_a$ has
indexing set $\ti I_a=\{k\}$, one Kuranishi neighbourhood $(\ti
V_a^k,\ti E_a^k,\ti s_a^k,\ti\psi_a^k)=(\ti X_a,\ti X_a,0,\id_{\ti
X_a})$, and $\bs{\ti\eta}_a=(\ti\eta_{a,k},\ti\eta_{a,k}^k)$ with
$\ti\eta_{a,k}\equiv 1\equiv\ti\eta_{a,k}^k$, and $\utG^k_a:\ti
X_a\ra\uP$ is some map satisfying injectivity conditions as
in~\S\ref{kh39}.

As in \S\ref{khA13} and \S\ref{khA16}, we can choose a tent function
$\ti T_a:\ti X_a\ra F\bigl([1,\iy)\bigr)$ for each $a\in A$ such
that the components $\ti X_{ac}$, $c\in C_a$ of $\pd Z_{\ti X_a,\ti
T_a}$ in \eq{khAeq5} are all diffeomorphic to the $k$-simplex
$\De_k$, with diffeomorphisms $\si_{ac}:\De_k\ra\ti X_{ac}$ for
$c\in C_a$. In order to make these choices in a way compatible with
the boundary relations involved in $\pd\bigl(\sum_{a\in A}\rho_a[\ti
X_a,\ti f_a,\ubtG_a]\bigr)=0$, we follow the strategy and use the
notation of Step 1, and by reverse induction on $m=k,k-1,\ldots,1,0$
we use Theorem \ref{khAthm2} to choose tent functions $\ti T_{ab}^m$
on $\ti X_{ab}^m$ for $(a,b)\in R^m$ and diffeomorphisms
$\si_{abc}^m:\De_{k-m}\ra\ti X_{abc}^m$ for $(a,b)\in R^m$ and $c\in
C_{ab}^m$, such that $\min\ti T_{ab}^m\vert_{\pd\ti X_{ab}^m}$ and
$\si_{abc}^m\vert_{\pd\ti X_{ab}^m}$ are identified on
$X_{ab'}^{m+1}$ for each $b'\in B_{ab}^m$ with $\min\ti T_{\bar
a\bar b}^{m+1}$ and $\si_{\bar a\bar b\bar c}^{m+1}$ under $(\bs
a,\bs b)_{ab'}^{m+1}$, for $(\bar a,\bar b)=\phi^{m+1}(a,b')$ in
$R^{m+1}$, and $\bar c\in C_{\bar a\bar b}^{m+1}$. At the end of the
induction we define $\ti T_a,C_a$ and $\si_{ac}$ using the $\ti
T_{ab}^0,C_{ab}^0$ and $\si_{abc}$ for $b=1,\ldots,n_a^0$ and~$c\in
C_{ab}^0$.

For $a\in A$ and $c\in C_a$, define $\ep_{ac}=1$ if the
diffeomorphism $\si_{ac}:\De_k\ra\ti X_{ac}$ is
orientation-preserving, and $\ep_{ac}=-1$ otherwise. We then show
that $\sum_{a\in A}\rho_a[\ti X_a,\ti f_a,\ubtG_a]$ is homologous in
$KC_k^\ef(Y;R)$ to the cycle
\e
\begin{split}
\ts\sum_{a\in A}&\ts\sum_{c\in C_a}\rho_a\ep_{ac}\bigl[\De_k,\ti f_a
\ci\pi\ci\si_{ac},\ubG_{\De_k}\bigr]\\
&=\Pi_\rsi^\ef\bigl(\ts\sum_{a\in A}\sum_{c\in C_a}
(\rho_a\ep_{ac})\,(\ti f_a \ci\pi\ci\si_{ac})\bigr).
\end{split}
\label{khBeq1}
\e
Since $\Pi_\rsi^\ef:C_*^\rsi(Y;R)\ra KC_*^\ef(Y;R)$ is injective at
the chain level, with $\Pi_\rsi^\ef\ci\pd=\pd\ci\Pi_\rsi^\ef$, as
\eq{khBeq1} is a cycle in $KC_k^\ef(Y;R)$ we have
\e
\pd\bigl(\ts\sum_{a\in A}\sum_{c\in C_a} (\rho_a\ep_{ac})\,(\ti f_a
\ci\pi\ci\si_{ac})\bigr)=0 \quad\text{in $C_{k-1}^\rsi(Y;R)$.}
\label{khBeq2}
\e
Thus $\be=\bigl[\sum_{a\in A}\sum_{c\in C_a}(\rho_a\ep_{ac})(\ti f_a
\!\ci\!\pi\!\ci\!\si_{ac})\bigr]$ is well-defined in
$H_k^\rsi(Y;R)$, and \eq{khBeq1} implies that
$\Pi_\rsi^\ef(\be)=\al$. Hence $\Pi_\rsi^\ef$ in \eq{kh4eq18} is
{\it surjective}.
\medskip

\noindent{\bf Step 4.} Suppose $\be\in H_k^\rsi(Y;R)$ with
$\Pi_\rsi^\ef(\be)=0$. Represent $\be$ by $\sum_{d\in
D}\eta_d\,\tau_d$ in $C^\rsi_k(Y;R)$, for $D$ a finite indexing set,
$\eta_d\in R$ and $\tau_d:\De_k\ra Y$ smooth. Then
$\Pi_\rsi^\ef\bigl(\sum_{d\in D}\eta_d\,\tau_d\bigr)$ is exact in
$KC_k^\ef(Y;R)$, so there exists $\sum_{a\in A}\rho_a[X_a,\bs
f_a,\ubG_a]\in KC_{k+1}^\ef(Y;R)$ with
\e
\ts\pd\bigl(\sum_{a\in A}\rho_a[X_a,\bs f_a,\ubG_a]\bigr)
=\sum_{d\in D}\eta_d\bigl[\De_k,\tau_d,\ubG_{\De_k}\bigr].
\label{khBeq3}
\e

We now apply Steps 1--3 to this $\sum_{a\in A}\rho_a[X_a,\bs
f_a,\ubG_a]$, replacing $k$-cycles by $(k+1)$-chains. This
eventually yields a singular $(k+1)$-chain $\ts\sum_{a\in
A}\sum_{c\in C_a}\ab(\rho_a\ep_{ac})\ab(\ti
f_a\!\ci\!\pi\!\ci\!\si_{ac})$. We make choices in Steps 1--3
so~that
\e
\pd\raisebox{-2pt}{$\displaystyle\Bigl($}\ts\sum\limits_{a\in
A}\sum\limits_{c\in C_a}(\rho_a\ep_{ac})\,(\ti f_a
\ci\pi\ci\si_{ac})\raisebox{-2pt}{$\displaystyle\Bigr)$}
=\ts\sum\limits_{d\in D}\sum\limits_{e\in
E}(\eta_d\ze_e)(\tau_d\ci\up_e),
\label{khBeq4}
\e
where $\{\up_e:e\in E\}$ are the simplices in the $N^{\it th}$ {\it
barycentric subdivision} of $\De_k$ for some $N\gg 0$, so that
$\up_e:\De_k\ra\De_k$ is an affine embedding for each $e\in E$, with
$\ze_e=1$ if $\up_e$ is orientation-preserving, and $\ze_e=-1$
otherwise. But $\sum_{d\in D}\sum_{e\in E}(\eta_d\ze_e)
(\tau_d\ci\up_e)$ is homologous to $\sum_{d\in D}\eta_d\,\tau_d$ in
$C^\rsi_k(Y;R)$ by Bredon \cite[\S IV.17]{Bred}, so $\sum_{d\in
D}\eta_d\tau_d$ is also a boundary in $C_k^\rsi(Y;R)$. Thus $\be=0$,
and $\Pi_\rsi^\ef$ in \eq{kh4eq18} is {\it injective}.
\medskip

\begin{rem}{\bf(a)} Steps 1--3 above are related to Fukaya and
Ono's construction \cite[\S 6]{FuOn1}, \cite[\S A1.1]{FOOO} of
virtual moduli chains for Kuranishi spaces by using {\it transverse
multisections}, small multi-valued perturbations ${\mathfrak
s}_{ac}^i$ of the obstruction maps $s_{ac}^i$. In our approach, we
first pass from Kuranishi spaces and $KH_*(Y;R)$ to {\it effective}
Kuranishi spaces and $KH_*^\ef(Y;R)$, as in Theorem \ref{kh4thm2}
and its proof in Appendix \ref{khC}. Then, Steps 1--3 above in
effect construct virtual moduli chains for effective Kuranishi
spaces using {\it single-valued\/} transverse perturbations $\ti
s_{ac}^i$. Since we use single-valued perturbations rather than
multisections, we are able to work over $\Z$, or an arbitrary
commutative ring $R$, rather than~$\Q$.

\noindent{\bf(b)} We require our perturbations $\ti s_{ac}^i$ in
Step 2 to be close to the Kuranishi maps $s_{ac}^i$ in $C^1$,
whereas Fukaya and Ono only require their multisections ${\mathfrak
s}_{ac}^i$ to be close to $s_{ac}^i$ in $C^0$. The difference comes
from our new Definition \ref{kh2def12}(e), which replaces Fukaya and
Ono's notion of `having a tangent bundle', and alters the notion of
compatibility of perturbations $\ti s_{ac}^i$ under coordinate
changes. Under a coordinate change $(\phi^{ij}_{ac},\hat
\phi^{ij}_{ac})$ from $(V^j_{ac},\ldots,\psi^j_{ac})$ to
$(V^i_{ac},\ldots,\psi^i_{ac})$, in our set-up $\ti s_{ac}^j$
determines $\ti s_{ac}^i$ only on $\phi^{ij}_{ac}(\ti V^{ij}_{ac})$,
but for Fukaya and Ono ${\mathfrak s}_{ac}^j$ determines ${\mathfrak
s}_{ac}^i$ on a neighbourhood of $\phi^{ij}_{ac}(V^{ij}_{ac})$. We
need $C^1$ closeness to ensure Definition \ref{kh2def12}(e) holds
for the coordinate changes \eq{khBeq15} below. Remark \ref{kh2rem5}
discussed ways of doing without Definition~\ref{kh2def12}(e).

\noindent{\bf(c)} In fact the notion that `$\ti s^i_{ac}$ is $C^1$
close to $s^i_{ac}$' in Step 2 does not really make sense, as we
have not said how close $\ti s^i_{ac},s^i_{ac}$ need to be. As in
Fukaya and Ono \cite[Th.~6.4]{FuOn1}, \cite[Th A1.23]{FOOO}, what we
should probably have done throughout Step 2 is not to choose a
single perturbation $\ti s^i_{ac}$ of $s^i_{ac}$, but a series $\ti
s^i_{ac,n}$ for $n=1,2,\ldots$, such that $\ti s^i_{ac,n}\ra
s^i_{ac}$ in $C^1$ as $n\ra\iy$, compatibilities between $\ti
s^i_{ac},\ti s^j_{ac}$ are replaced by compatibilities between $\ti
s^i_{ac,n},\ti s^j_{ac,n}$ for each $n=1,2,\ldots$, and so on. Then
we take $n\gg 0$ to ensure $C^1$-closeness wherever we need it. But
we suppress this point, as there are too many indices in the proof
already.
\label{khBrem1}
\end{rem}

\subsection{Step 1: cutting $X_a$ into $X_{ac}$ with trivial
stabilizers}
\label{khB1}

We work in the situation and with the notation of Step 1 above, so
that $\sum_{a\in A}\rho_a[X_a,\bs f_a,\ubG_a]\in KC_k^\ef(Y;R)$ is a
cycle representing $\al\in KH_k^\ef(Y;R)$. We will use Theorem
\ref{khAthm5} to choose a {\it tent function} $\bs T_a$ for
$(X_a,\bs f_a,\ubG_a)$ in the sense of \S\ref{khA32}, for each $a\in
A$. Then Definition \ref{khAdef12} defines $\bigl[Z_{X_a,\bs
T_a},\bs f_a\ci\bs\pi_a,\ubH_{X_a,\bs T_a}\bigr]$ in
$KC_{k+1}^\ef(Y;R)$, and equation \eq{khAeq30} of Proposition
\ref{khAprop9} yields
\e
\begin{split}
\pd\bigl[Z_{X_a,\bs T_a},\bs f_a\ci\bs\pi_a,\ubH_{X_a,\bs T_a}\bigr]=\,
&\ts\sum_{c\in C_a}[X_{ac},\bs f_{ac},\ubG_{ac}]-[X_a,\bs f_a,\ubG_a]\\
-\bigl[Z_{\pd X_a,\bs T_a\vert_{\pd X_a}},&\bs f_a \vert_{\pd X_a}
\ci\bs\pi_a,\ubH_{\pd X_a,\bs T_a\vert_{\pd X_a}}\bigr],
\end{split}
\label{khBeq5}
\e
where $C_a$ is some finite indexing set, and Theorem \ref{khAthm5}
says that each $X_{ac}$ has {\it trivial stabilizers}, and the
$V^i_{ac}$ in $\ubG_{ac}$ are {\it manifolds}. As $\sum_{a\in
A}\rho_a[X_a,\bs f_a,\ubG_a]$ is a cycle we have
\e
\ts\sum_{a\in A}\rho_a\bigl[\pd X_a,\bs f_a\vert_{\pd
X_a},\ubG_a\vert_{\pd X_a}\bigr]=0.
\label{khBeq6}
\e
We want to choose the $\bs T_a$ for $a\in A$ such that
\e
\ts\sum_{a\in A}\rho_a\bigl[Z_{\pd X_a,\bs T_a\vert_{\pd X_a}},\bs
f_a \vert_{\pd X_a} \ci\bs\pi_a,\ubH_{\pd X_a,\bs T_a\vert_{\pd
X_a}}\bigr]=0.
\label{khBeq7}
\e
Note that if the $\bs T_a$ depend functorially on $[X_a,\bs
f_a,\ubG_a]$ in some sense compatible with boundaries, then
\eq{khBeq7} would follow from \eq{khBeq6}.

Assuming \eq{khBeq7}, multiplying \eq{khBeq5} by $\rho_a$ and
summing over $a\in A$ yields
\e
\begin{split}
\pd\Bigl(\ts\sum_{a\in A}&\rho_a\bigl[Z_{X_a,\bs T_a},\bs
f_a\ci\bs\pi_a,\ubH_{X_a,\bs T_a}\bigr]\Bigr)\\
&=\ts\sum_{a\in A}\sum_{c\in C_a}\rho_a[X_{ac},\bs f_{ac},
\ubG_{ac}]-\sum_{a\in A}\rho_a[X_a,\bs f_a,\ubG_a].
\end{split}
\label{khBeq8}
\e
Thus $\sum_{a\in A}\sum_{c\in C_a}\rho_a[X_{ac},\bs f_{ac},
\ubG_{ac}]$ is a cycle homologous to $\sum_{a\in A}\rho_a
\bigl[X_a,\ab \bs f_a,\ab\ubG_a\bigr]$, and so represents $\al$, as
we have to prove.

The obvious way to do this is to first choose tent functions $\bs
T_a\vert_{\pd X_a}$ for $(\pd X_a,\bs f_a\vert_{\pd
X_a},\ubG_a\vert_{\pd X_a})$ satisfying the necessary
compatibilities that will ensure \eq{khBeq7} holds, and then extend
these choices from $\pd X_a$ to $X_a$. However, there is a problem.
Let $\bs\si_{X_a}:\pd^2X_a\ra \pd^2X_a$ be the orientation-reversing
involution from Definition \ref{kh2def19}. It is a general
principle, as in Principle \ref{kh2pri} and Property \ref{kh3pr}(e)
for instance, that a choice of data over $\pd X_a$ should extend to
$X_a$ only if the restriction of the choice to $\pd^2X_a$ is
invariant under $\bs\si_{X_a}$. So our choice of $\bs T_a$ over $\pd
X_a$ should satisfy extra conditions over~$\pd^2X_a$.

We can achieve this by first choosing $\bs T_a$ over $\pd^2X_a$
satisfying the compatibilities and these extra conditions, and then
extending the choice from $\pd^2X_a$ to $\pd X_a$. But this
extension is only possible provided the restriction of the choice to
$\pd^3X_a$ is invariant under $\bs\si_{\pd X_a}:\pd^2(\pd
X_a)\ra\pd^2(\pd X_a)$. Continuing this argument indefinitely, we
see that the right approach is to choose data $\bs T_a$ over
$\pd^mX_a$ for $a\in A$ by {\it reverse induction on\/} $m$, where
in the inductive step we choose data over $\pd^mX_a$ which restricts
on $\pd^{m+1}X_a$ to the data already chosen, and has the
compatibilities necessary to prove \eq{khBeq7}, and is invariant
under the strong diffeomorphisms $\bs\si_l:\pd^mX_a\ra\pd^mX_a$
induced by $\bs\si_{\pd^l X_a}:\pd^2(\pd^l X_a)\ra\pd^2(\pd^l X_a)$
for all~$0\le l\le m-2$.

Use the notation $n_a^m,(X_{ab}^m,\bs f_{ab}^m,
\ubG_{ab}^m),P^m,Q^m,R^m,M,\ab(a^q,b^q),\ab\phi^m,\ab(\bs a,\bs
b)_{ab}^m,\ab\ep_{ab}^m,\ab B_{ab}^m$ of Step 1 above. By induction
on decreasing $m=M,M-1,\ldots,1,0$, for all $(a,b)$ in $R^m$, we
choose tent functions $\bs T_{ab}^m$ for $(X_{ab}^m,\bs
f_{ab}^m,\ubG_{ab}^m)$ using Definition \ref{khAdef16} and Theorem
\ref{khAthm5}, satisfying the inductive hypothesis:
\begin{itemize}
\item[$(*)_{ab}^m$] By definition of $B_{ab}^m$ we have
$\pd X_{ab}^m=\coprod_{b'\in B_{ab}^m}X_{ab'}^{m+1}$, and if
$\phi^{m+1}(a,b')=(\bar a,\bar b)$ in $R^{m+1}$, we have an
isomorphism $(\bs a,\bs b)^{m+1}_{ab'}:(X_{ab'}^{m+1},\ab\bs
f_{ab'}^{m+1},\ab\ubG_{ab'}^{m+1})\ab \ra(X_{\bar a\bar
b}^{m+1},\ab\bs f_{\bar a\bar b}^{m+1},\ab\ubG_{\bar a\bar
b}^{m+1})$. So there is an isomorphism
\e
\!\!\!\!\!\!\!\!
\coprod_{b'\in B_{ab}^m}(\bs a,\bs
b)^{m+1}_{ab'}:\pd(X_{ab}^m,\bs f_{ab}^m, \ubG_{ab}^m)\longra
\coprod_{\!\!\!\!\text{$b'\in B_{ab}^m$: define $(\bar a,\bar
b)=\phi^{m+1}(a,b')$}
\!\!\!\!\!\!\!\!\!\!\!\!\!\!\!\!\!\!\!\!\!\!\!\!\!\!\!\!\!
\!\!\!\!\!\!\!\!\!\!\!\!\!\!\!\!\!\!\!\!\!\!\!\!\!\!\!\!\!
\!\!\!\!\!\!\!\!\!\!\!\!\!\!\!\!\!\!\!\!\! } (X_{\bar a\bar
b}^{m+1},\ab\bs f_{\bar a\bar b}^{m+1},\ab\ubG_{\bar a\bar
b}^{m+1}).
\label{khBeq9}
\e
On the left hand side $\pd(X_{ab}^m,\bs f_{ab}^m,\ubG_{ab}^m)$ of
\eq{khBeq9} we have a tent function $\bs T_{ab}^m\vert_{\pd
X_{ab}^m}$. On the right hand side we have a tent function
$\coprod_{b'}\bs T_{\bar a\bar b}^{m+1}$. We require that
\eq{khBeq9} should identify $\min\bigl(\bs T_{ab}^m\vert_{\pd
X_{ab}^m}\bigr)$ with~$\min\bigl(\coprod_{b'}\bs T_{\bar a\bar
b}^{m+1}\bigr)$.
\end{itemize}
For the first step, when $m=M$, as $P^{M+1}=\es$, if $(a,b)\in R^M$
then $\pd X_{ab}^M=\es$, so the boundary conditions in $(*)_{ab}^M$
are trivial. Thus for all $(a,b)\in R^M$ we choose $\bs T_{ab}^M$
for $(X_{ab}^M,\bs f_{ab}^M,\ubG_{ab}^M)$ using Definition
\ref{khAdef13} and Proposition~\ref{khAprop10}.

For the inductive step, suppose that for some $l=0,\ldots,M-1$ we
have chosen $\bs T_{ab}^m$ satisfying $(*)_{ab}^m$ for all $l<m\le
M$ and $(a,b)\in R^m$. We will choose $\bs T_{ab}^l$ satisfying
$(*)_{ab}^l$ for all $(a,b)\in R^l$. The important point is this.
Let $\bs T_{ab}^{\prime l}$ be the tent function for
$\pd(X_{ab}^l,\bs f_{ab}^l,\ubG_{ab}^l)$ obtained by pulling back
$\coprod_{b'}\bs T_{\bar a\bar b}^{l+1}$. By $(*)_{ab}^l$, we have
to choose a tent function $\bs T_{ab}^l$ for $(X_{ab}^l,\bs
f_{ab}^l,\ubG_{ab}^l)$ such that~$\min\bigl(\bs
T_{ab}^l\vert_{\smash{\pd X_{ab}^l}}\bigr)=\min\bs T_{ab}^{\prime
l}$.

As in Definition \ref{khAdef16}, we can do this provided $\min\bs
T_{ab}^{\prime l}\vert_{\pd^2X_{ab}^l}$ is invariant under the
natural involution $\bs\si:\pd^2X_{ab}^l\ra\pd^2X_{ab}^l$. To see
this, note that $\pd X_{ab}^l=\coprod_{b'\in
B_{ab}^l}X_{ab'}^{l+1}$, and so $\pd^2X_{ab}^l=\coprod_{b'\in
B_{ab}^l}\coprod_{b''\in B_{ab'}^{l+1}}X_{ab''}^{l+2}$. Thus, the
splitting of $\bigl( \pd^2X_{ab}^l,\bs
f_{ab}^l\vert_{\pd^2X_{ab}^l},\ubG_{ab}^l
\vert_{\pd^2X_{ab}^l}\bigr)$ into connected components in Lemma
\ref{kh3lem} is
\e
\bigl(\pd^2X_{ab}^l,\bs f_{ab}^l\vert_{\pd^2X_{ab}^l},\ubG_{ab}^l
\vert_{\pd^2X_{ab}^l}\bigr)\cong\coprod_{b'\in
B_{ab}^l}\coprod_{b''\in B_{ab'}^{l+1}}\bigl(X_{ab''}^{l+2},\bs
f_{ab''}^{l+2},\ubG_{ab''}^{l+2}\bigr).
\label{khBeq10}
\e

The natural orientation-reversing involution $\bs\si:\pd^2X_{ab}^l
\ra\pd^2X_{ab}^l$ extends to an involution $(\bs\si,\bs\tau)$ of
$\bigl(\pd^2X_{ab}^l,\bs f_{ab}^l\vert_{\pd^2
X_{ab}^l},\ubG_{ab}^l\vert_{\pd^2X_{ab}^l}\bigr)$. Thus, $\bs\si$
must permute the components $X_{ab''}^{l+2}$ in \eq{khBeq10}. Now
$\Aut\bigl(X_{ab''}^{l+2},\bs
f_{ab''}^{l+2},\ubG_{ab''}^{l+2}\bigr)=\{1\}$ by Theorem
\ref{kh3thm5}(b), as $\ubG_{ab''}^{l+2}$ is effective. Since
$\bs\si$ is orientation-reversing, it cannot act as the identity.
Therefore $\bs\si$ must act {\it freely} on the set of
$X_{ab''}^{l+2}$ appearing in \eq{khBeq10}, as if $\bs\si$ fixed any
component, that component would have a nontrivial automorphism.

Let $X_{ab''}^{l+2}$ be a component on the r.h.s.\ of \eq{khBeq10}.
Then $\bs\si$ takes $X_{ab''}^{l+2}$ to some unique
$X_{ab'''}^{l+2}$ with $b''\ne b'''$. The isomorphism
$(\bs\si,\bs\tau)$ between these two components implies that
$(a,b'')$ and $(a,b''')$ lie in the same equivalence class in
$P^{l+2}$, so that $\phi^{l+2}(a,b'')=\phi^{l+2}(a,b''')=(\bar
a,\bar b)\in R^{l+2}$. The associated isomorphisms $(\bs a,\bs
b)_{ab''}^{l+2}:\bigl(X_{ab''}^{l+2},\bs f_{ab''}^{l+2},
\ubG_{ab''}^{l+2}\bigr)\ra\bigl(X_{\bar a\bar b}^{l+2},\bs f_{\bar
a\bar b}^{l+2},\ubG_{\bar a\bar b}^{l+2}\bigr)$ and $(\bs a,\bs
b)_{ab'''}^{l+2}:\bigl(X_{ab'''}^{l+2},\bs f_{ab'''}^{l+2},
\ubG_{ab'''}^{l+2}\bigr)\!\ra\!\bigl(X_{\bar a\bar b}^{l+2},\bs
f_{\bar a\bar b}^{l+2},\ubG_{\bar a\bar b}^{l+2}\bigr)$
satisfy~$(\bs a,\bs b)_{ab''}^{l+2}\!=\!(\bs a,\bs
b)_{ab'''}^{l+2}\ci(\bs\si,\bs\tau)$.

The restrictions of $\min\bs T_{ab}^{\prime l}$ to $X_{ab''}^{l+2}$
and $X_{ab'''}^{l+2}$ are the pullbacks of $\min\bs T_{\bar a\bar
b}^{l+2}$ by $(\bs a,\bs b)_{ab''}^{l+2}$ and $(\bs a,\bs
b)_{ab'''}^{l+2}$ respectively, and so $(\bs a,\bs
b)_{ab''}^{l+2}\!=\!(\bs a,\bs b)_{ab'''}^{l+2}\ci(\bs\si,\bs\tau)$
implies that the restrictions of $\min\bs T_{ab}^{\prime l}$ to
$X_{ab''}^{l+2}$ and to $X_{ab'''}^{l+2}$ are identified by
$\bs\si$. Hence $\min\bs T_{ab}^{\prime l}\vert_{\pd^2X_{ab}^l}$ is
$\bs\si$-invariant, and Definition \ref{khAdef16} and Theorem
\ref{khAthm5} apply. This completes the inductive step.

Hence by induction we can choose $\bs T_{ab}^m$ satisfying
$(*)_{ab}^m$ for all $m,a,b$. Now $X_a=\coprod_{b=1}^{
n_a^0}X_{ab}^0$, and for each $b=1,\ldots,n_a^0$, if
$\phi^0(a,b)=(\bar a,\bar b)$ in $R^0$, we have an isomorphism $(\bs
a,\bs b)^0_{ab}:(X_{ab}^0,\ab\bs f_{ab}^0,\ab\ubG_{ab}^0)\ab
\ra(X_{\bar a\bar b}^0,\ab\bs f_{\bar a\bar b}^0,\ab\ubG_{\bar a\bar
b}^0)$, and a tent function $\bs T_{\bar a\bar b}^0$ for $(X_{\bar
a\bar b}^0,\ab\bs f_{\bar a\bar b}^0,\ab\ubG_{\bar a\bar b}^0)$.
Define a tent function $\bs T_a$ for $(X_a,\bs f_a,\ubG_a)$ by $\bs
T_a\vert_{X_{ab}^0}=((\bs a,\bs b)^0_{ab})^*(\bs T_{\bar a\bar
b}^0)$ for $b=1,\ldots,n_a^0$, for each $a\in A$. We have
\ea
0&=\pd\Bigl[\sum_{a\in A}\rho_a[X_a,\bs f_a,\ubG_a]\Bigr]
=\sum_{a\in A}\rho_a\sum\limits_{b=1}^{n_a^1}
\bigl[X_{ab}^1,\bs f_{ab}^1, \ubG_{ab}^1\bigr]
\label{khBeq11}\\
&=\sum\limits_{\begin{subarray}{l}\text{$a\in A$,
$b=1,\ldots,n_a^1$.}\\ \text{Set $(\bar a,\bar
b)=\phi^1(a,b)$}\end{subarray}\!\!\!\!\!\!\!\!\!\!\!\!\!\!\!\!
\!\!\!\!\!\!\!\!\!\!\!\!\!\! } \rho_a\ep_{ab}^1\bigl[X_{\bar a\bar
b}^1,\bs f_{\bar a\bar b}^1,\ubG_{\bar a\bar b}^1\bigr]=\!\!\!
\sum_{(\bar a,\bar b)\in R^1}
\raisebox{-9pt}{\begin{Large}$\displaystyle\biggl\{$\end{Large}}
\sum\limits_{\begin{subarray}{l}a\in A,\;
b=1,\ldots,n_a^1:\\ \phi^1(a,b)=(\bar a,\bar
b)\end{subarray}\!\!\!\!\!\!\!\!\!\!\!\!\!}\rho_a\ep_{ab}^1
\raisebox{-9pt}{\begin{Large}$\displaystyle\biggr\}$\end{Large}}
\bigl[X_{\bar a\bar b}^1,\bs f_{\bar a\bar b}^1,\ubG_{\bar a\bar
b}^1\bigr].
\nonumber
\ea

By construction the triples $\bigl(X_{\bar a\bar b}^1,\bs f_{\bar
a\bar b}^1,\ubG_{\bar a\bar b}^1\bigr)$ for $(\bar a,\bar b)\in R^1$
are connected, have no non-trivial automorphisms, and are mutually
non-isomorphic, even through orientation-reversing isomorphisms. It
follows that $\bigl[X_{\bar a\bar b}^1,\bs f_{\bar a\bar
b}^1,\ubG_{\bar a\bar b}^1\bigr]$ for $(\bar a,\bar b)\in R^1$ are
linearly independent over $R$ in $KC_{k-1}^\ef(Y;R)$, as Definition
\ref{kh4def2}(i)--(iii) impose no nontrivial relations between them.
Thus \eq{khBeq11} implies that the coefficient $\{\cdots\}$ on the
last line is zero, for all $(\bar a,\bar b)\in R^1$. The same
argument as for \eq{khBeq11} shows that
\e
\begin{split}
&\ts\sum_{a\in A}\rho_a\bigl[Z_{\pd X_a,\bs T_a\vert_{\pd X_a}},\bs
f_a \vert_{\pd X_a} \ci\bs\pi_a,\ubH_{\pd X_a,\bs T_a\vert_{\pd
X_a}}\bigr]=\\
&\sum\limits_{\!\!\!\!\!\!\!(\bar a,\bar b)\in R^1\!\!\!\!\!\!\!}
\,\,\,\,
\raisebox{-9pt}{\begin{Large}$\displaystyle\bigg\{$\end{Large}}
\sum\limits_{\begin{subarray}{l}a\in A,\;
b=1,\ldots,n_a^1:\\ \phi^1(a,b)=(\bar a,\bar
b)\end{subarray}\!\!\!\!\!\!\!\!\!\!\!\!\!\!\!}\rho_a\ep_{ab}^1
\raisebox{-9pt}{\begin{Large}$\displaystyle\biggr\}$\end{Large}}
\bigl[Z_{X_{\bar a\bar b}^1,\bs T_{\bar a\bar b}^1},\bs f_{\bar
a\bar b}^1\ci\bs\pi_{\bar a\bar b}^1,\ubH_{X_{\bar a\bar b}^1,\bs
T_{\bar a\bar b}^1}\bigr]=0,
\end{split}
\label{khBeq12}
\e
since the coefficient $\{\cdots\}$ on the second line is zero. This
proves \eq{khBeq7}, and \eq{khBeq8} follows, completing Step~1.

\subsection{Step 2: deforming $X_{ac}$ to compact manifolds $\ti
X_{ac}$}
\label{khB2}

We work in the situation and with the notation of Steps 1 and 2
above, so that $\sum_{a\in A}\sum_{c\in C_a}\rho_a[X_{ac},\ab\bs
f_{ac},\ab\ubG_{ac}]$ is a cycle representing $\al\in
KH_k^\ef(Y;R)$, where $X_{ac}$ has trivial stabilizers and the
$V^i_{ac}$ are manifolds, and our goal is to construct a homologous
cycle $\sum_{a\in A}\sum_{c\in C_a}\rho_a[\ti X_{ac},\ti
f_{ac},\ab\ubtG_{ac}]$ with $\ti X_{ac}$ a compact $k$-manifold with
g-corners, constructed from $X_{ac}$ by choosing $C^1$ small,
smooth, transverse perturbations $\ti s^i_{ac}$ of the Kuranishi
maps $s^i_{ac}$ of the Kuranishi neighbourhoods
$(V^i_{ac},E^i_{ac},s^i_{ac},\psi^i_{ac})$ on $X_{ac}$ in $\bs
I_{ac}$ in $\ubtG_{ac}$. Our argument is based on Fukaya and Ono
\cite[Th.~6.4]{FuOn1}, \cite[Th.~A1.23]{FOOO}, but using
single-valued sections rather than multisections.

To prove that $\sum_{a\in A}\sum_{c\in C_a}\rho_a[X_{ac},\ab\bs
f_{ac},\ab\ubG_{ac}]$ and $\sum_{a\in A}\sum_{c\in C_a}\rho_a[\ti
X_{ac},\ti f_{ac},\ab\ubtG_{ac}]$ are homologous, we will also
construct a homology $[Z_{ac},\bs g_{ac},\ubH_{ac}]$ between
$[X_{ac},\ab\bs f_{ac},\ab\ubG_{ac}]$ and $[\ti X_{ac},\ti
f_{ac},\ab\ubtG_{ac}]$, modulo terms over $\pd X_{ac}$. We first
explain how to choose the $\ti s^i_{ac}$ and define $\ti X_{ac},\ti
f_{ac},\ab\ubtG_{ac},Z_{ac},\bs g_{ac},\ubH_{ac}$ for a single pair
$a,c$, ignoring compatibility conditions between other pairs
$a',c'$. Then later we will explain how to choose the $\ti s^i_{ac}$
for all $a\in A$ and $c\in C_a$ so that the necessary compatibility
conditions hold. For simplicity, in this first part we omit the
subscripts `$ac$', so we begin with a triple $(X,\bs f,\ubG)$,
choose perturbations $\ti s^i$, and define~$\ti X,\ti f,\ubtG,Z,\bs
g,\ubH$.

Let $(X,\bs f,\ubG)$ be a triple with $\ubG=(\bs I,\bs\eta,\uG^i:
i\in I)$ and $\bs I=\bigl(I,(V^i,\ldots,\ab\psi^i),\ab f^i:i\in
I,\ldots \bigr)$, where $X$ has trivial stabilizers and the $V^i$
are manifolds. Define $\dot V^i,\dot V^{ij},(\bs{\dot I},\bs{\dot
\eta})$ as in Definition \ref{khAdef10}, so that $\dot V^i\subset
V^i$, $\dot V^{ij}\subset V^{ij}$ are open sets for $j\le i$ with
$\ov{\phi^{ij}(\dot V^{ij})}\subseteq\phi^{ij}(V^{ij})$. We shall
choose perturbations $\ti s^i$ of $s^i\vert_{\dot V^i}$ for $i\in
I$. The necessity to use smaller open subsets $\dot V^i\subset V^i$
on which to define the perturbations $\ti s^i$ is explained after
Example \ref{khAex1}: without shrinking the domains $V^i,V^{ij}$ we
may be unable to choose the $\ti s^i$ to be transverse over points
in $V^i$ lying in the boundary of $\phi^{ij}(V^{ij})$ for $j\le i$
in~$I$.

We will choose $C^1$ small perturbations $\ti s^i$ of
$s^i\vert_{\smash{\dot V^i}}$ with $\ti s^i$ {\it transverse\/}
along $(\ti s^i)^{-1}(0)$, such that if $j<i$ in $I$ then
$\hat\phi^{ij}\ci\ti s^j\equiv\ti s^i\ci\phi^{ij}$ on $\dot V^{ij}$,
by induction on $i\in I$. Here is our inductive hypothesis. Suppose
that $i\in I$, and that for $j<i$ in $I$ we have already chosen a
smooth section $\ti s^j$ of $E^j\vert_{\smash{\dot V^j}}$, such that
$\ti s^j$ extends smoothly and transversely from $\dot V^j$ to some
open neighbourhood $U^j$ of the closure $\overline{\dot V^j}$ of
$\dot V^j$ in $V^j$, and $\ti s^j$ is $C^1$ close to $s^j$, and $\ti
s^j$ is {\it transverse\/} along $(\ti s^j)^{-1}(0)$, and if $k\le
j<i$ in $I$ then $\hat\phi^{jk}\ci\ti s^k\equiv\ti s^j\ci\phi^{jk}$
on $\dot V^{jk}$. The first step $i=\min I$ is vacuous.

To do the inductive step we have to choose a smooth transverse
section $\ti s^i$ of $E^i\vert_{\smash{\dot V^i}}$ which extends
smoothly from $\dot V^i$ to some open neighbourhood $U^i$ of
$\overline{\dot V^i}$ in $V^i$, such that if $j<i$ in $I$ then
$\hat\phi^{ij}\ci\ti s^j\equiv\ti s^i\ci\phi^{ij}$ on $\dot V^{ij}$.
That is, $\ti s^i$ is prescribed on the subset $\bigcup_{\text{$j<i$
in $I$}}\phi^{ij}(\dot V^{ij})$ of $\dot V^i$, and we must extend it
smoothly and transversely to~$\dot V^i$.

As in the proof of Theorem \ref{kh3thm1}, the prescribed values for
$\ti s^i$ from different $j,k<i$ are consistent, and the subset
$\bigcup_{\text{$j<i$ in $I$}}\phi^{ij}(\dot V^{ij})$ is a locally
closed submanifold of $\dot V^i$ except possibly at points of the
{\it boundaries\/} $\pd\bigl(\phi^{ij}(\dot V^{ij})\bigr)=
\ov{\phi^{ij}(\dot V^{ij})}\sm \phi^{ij}(\dot V^{ij})$. Thus, except
possibly at these boundaries, locally we can extend $\ti s^i$
smoothly to $\dot V^i$, and since the $\ti s^j$ for $j<i$ are {\it
transverse}, locally we can choose $\ti s^i$ to be transverse too.
Note that as $V^i,\dot V^i$ are manifolds and $E^i$ is a vector
bundle, there are no problems with orbifold groups, and generic
sections are transverse.

We have $\overline{\phi^{ij}(\dot V^{ij})}\subseteq
\phi^{ij}(V^{ij})$ as in Definition \ref{khAdef10}. But we already
know that $\ti s^j$ extends smoothly and transversely to some open
neighbourhood $U^j$ of $\overline{\dot V^j}$ in $V^j$, and therefore
the prescribed values of $\ti s^i$ on $\phi^{ij}(\dot V^{ij})$
extend smoothly and transversely to $\phi^{ij}(U^j\cap V^{ij})$,
which contains $\pd\smash{\bigl(\phi^{ij}(\dot V^{ij})\bigr)}$. Thus
$\ti s^i$ can locally be extended smoothly and transversely over the
boundaries $\pd\bigl(\phi^{ij}(\dot V^{ij})\bigr)$, so the extension
is locally possible everywhere in $\dot V^i$. Patching the local
choices together with a partition of unity, we can choose $\ti s^i$
with all the properties we need. This completes the inductive step.

We shall now construct $\ti X,\ti f,\ubtG$. As in Remark
\ref{khBrem1}(a), this is similar to parts of Fukaya and Ono's
construction of {\it virtual cycles\/} \cite[\S 6]{FuOn1}, \cite[\S
A1]{FOOO}. Define
\e
\ti X=\ts\bigl[\,\coprod_{i\in I}(\ti
s^i)^{-1}(0)\,\bigr]\big/\!\sim,
\label{khBeq13}
\e
where $\sim$ is the equivalence relation induced by the maps
$\phi^{ij}:(\ti s^j)^{-1}(0)\cap\dot V^{ij}\ra(\ti s^i)^{-1}(0)$ for
all $i,j\in I$ with $j\le i$. Here $\phi^{ij}$ maps $(\ti
s^j)^{-1}(0)$ to $(\ti s^i)^{-1}(0)$ since $\hat\phi^{ij}\ci\ti
s^j\equiv\ti s^i\ci\phi^{ij}$, and $\hat\phi^{ij}$ is an embedding
and so injective.

Since the $\dot V^i$ are manifolds and the $\ti s^i$ are transverse,
each subset $(\ti s^i)^{-1}(0)$ is actually a manifold with
g-corners, of dimension $k$. As the gluing maps $\phi^{ij}$ are
embeddings $V^{ij}\ra V^i$ they restrict to local diffeomorphisms
from $(\ti s^j)^{-1}(0)$ to $(\ti s^i)^{-1}(0)$. Therefore $\ti X$
has the structure of a $k$-manifold with g-corners, which is {\it
Hausdorff\/} by the argument of Remark \ref{kh3rem2}, and {\it
compact\/} provided the $\ti s^i$ are sufficiently close to $s^i$ in
$C^0$, so that $(\ti s^i)^{-1}(0)$ is close to~$(s^i)^{-1}(0)$.

Define $\ti f:\ti X\ra Y$ by $\ti f\vert_{(\ti s^i)^{-1}(0)}
=f^i\vert_{(\ti s^i)^{-1}(0)}$ for $i\in I$. Since
$f^i\ci\phi^{ij}\equiv f^j$ on $V^{ij}$, this is compatible with the
gluing maps $\phi^{ij}$ defining $\sim$, so $\ti f$ is well-defined,
and is clearly a smooth map. Define an {\it excellent coordinate
system\/} $(\bs{\ti I}, \bs{\ti\eta})$ for $(\ti X,\ti f)$ to have
indexing set $\ti I=\{k\}$ and Kuranishi neighbourhood $(\ti V^k,\ti
E^k,\ti s^k,\ti\psi^k)=(\ti X,\ti X,0,\smash{\id_{\ti X}})$, where
$\ti E^k\ra\ti V^k=\ti X$ is the zero vector bundle, so that $\ti
E^k=\ti X$ as a manifold, and define $\ti f^k:\ti V^k\ra Y$ by $\ti
f^k=\ti f$, so that $\ti f^k$ represents $\ti f$. Define
$\ti\eta_k:\ti X\ra[0,1]$ and $\ti\eta_k^k:\ti V^k\ra[0,1]$
by~$\ti\eta_k\equiv 1\equiv \ti\eta_k^k$.

Choose some map $\utG^k:\ti X\ra\uP$ such that $\ubtG=(\bs{\ti
I},\bs{\ti\eta},\utG^k)$ is {\it effective gauge-fixing data} for
$(\ti X,\ti f)$. Ignoring boundary conditions, it is enough to take
$\utG^k:\ti X\ra\R^n\subset\uP$ injective for some $n\gg 0$, but
when we discuss how to choose $\ti s^i_{ac},\ti X_{ac},\ti
f_{ac},\ab\ubtG_{ac},Z_{ac},\bs g_{ac},\ubH_{ac}$ with compatibility
between different $a\in A$, $c\in C_a$, we will explain in more
detail how to choose these maps $\utG^k$. Then $[\ti X,\ti f,\ubtG]$
is well-defined in~$KC_k^\ef(Y;R)$.

Next we define $(Z,\bs g,\ubH)$ such that $Z$ is a perturbation of
$[0,1]\t X$, and $\bigl[Z,\bs g,\ubH\bigr]$ is a homology between
$[X,\bs f,\ubG]$ and $[\ti X,\ti f,\ubtG{}]$, modulo terms over a
perturbation of $[0,1]\t\pd X$. In $\ubH$ we will have a very good
coordinate system $\bs J$ with indexing set $J=\{i+1:i\in
I\}\cup\{k+1\}$ and Kuranishi neighbourhoods $(W^j,F^j,t^j,\xi^j)$
for $j\in J$, and we define these first. For $i\in I\sm\{k\}$,
define $W^{i+1}=\bigl([0,\frac{1}{3})\t V^i\bigr)\cup
\bigl([\frac{1}{3},\frac{3}{4})\t\dot V^i\bigr)$. This is an open
set in $[0,\frac{3}{4})\t V^i$, and so is a manifold with g-corners.
Define $F^{i+1}=\pi_{\smash{V^i}}^*(E^i)$, so that $F^{i+1}\ra
W^{i+1}$ is a vector bundle.

If $k\notin I$, define $W^{k+1}=(\frac{2}{3},1]\t\ti X$. If $k\in
I$, define $W^{k+1}=\bigl([0,\frac{1}{3})\t V^k\bigr)\cup
\bigl([\frac{1}{3},\frac{3}{4})\t\ti V^k\bigr)\cup
(\frac{2}{3},1]\t\ti X$. Here as $\vdim X=k$ we have $\rank E^k=0$,
so $(\ti s^k)^{-1}(0)=\ti V^k$, and we regard $\ti V^k$ as an open
subset of $\ti X$, so that $[\frac{1}{3},\frac{3}{4})\t\ti V^k$ and
$(\frac{2}{3},1]\t\ti X$ intersect in
$(\frac{2}{3},\frac{3}{4})\t\ti V^k$. In both cases define
$F^{k+1}\ra W^{k+1}$ to be the zero vector bundle, so that
$F^{k+1}=W^{k+1}$ as a manifold.

Fix a smooth map $\eta:[0,1]\ra[0,1]$ independent of $a\in A$ such
that $\eta(x)=0$ for $x\in[0,\frac{1}{3}]$ and $\eta(x)=1$ for
$x\in[\frac{2}{3},1]$. For $i\in I\sm\{k\}$ define a smooth section
$t^{i+1}$ of $F^{i+1}$ by $t^{i+1}(u,v)=s^i(v)$ for
$u\in[0,\frac{1}{3})$ and $t^{i+1}(u,v)=
(1-\eta(u))s^i(v)+\eta(u)\ti s^i(v)$ for
$u\in[\frac{1}{3},\frac{3}{4})$. Define $t^{k+1}=0$. Following
\eq{khBeq13}, define
\e
Z=\bigl[\,\ts\coprod_{j\in J}(t^j)^{-1}(0)\,\bigr]/\approx,
\label{khBeq14}
\e
where $\approx$ is the equivalence relation induced by
$(u,v^i)\approx (u,v^j)$ if $u\in[0,\frac{3}{4})$, $j\le i$ in $I$,
$(u,v^i)\in(t^{i+1})^{-1}(0)$, $(u,v^j)\in(t^{j+1})^{-1} (0)$, and
$v^j\in V^{ij}$ with $\phi^{ij}(v^j)=v^i$, and $(u,v^i)\approx
(u,v^k)$ if $u\in(\frac{2}{3},\frac{3}{4})$, $i\in I$,
$(u,v^i)\in(t^{i+1})^{-1}(0)$, $(u,v^k)\in(\frac{2}{3},1] \t\ti X$,
and $v^i=v^k$ in $\ti X$, where $\eta(u)=1$ and
$(u,v^i)\in(t^i)^{-1}(0)$ imply that $v^i\in(\ti s^i)^{-1}(0)$, so
we can regard $v^i$ as an element of $\ti X$. As the gluing maps are
local homeomorphisms, $Z$ has the structure of a topological space,
which is compact and Hausdorff provided $\ti s^i$ is close to $s^i$
in~$C^0$.

Define $\xi^j:(t^j)^{-1}(0)\ra Z$ by $\xi^j=\id_{(t^j )^{-1}(0)}$.
Then $(W^j,F^j,t^j,\xi^j)$ is a {\it Kuranishi neighbourhood\/} on
$Z$ for $j\in J$. When $k<j\le i$ in $I$, define an open subset
$W^{(i+1)(j+1)}$ in $W^{j+1}$ to be $\bigl([0,\frac{1}{3})\t
V^{ij}\bigr)\cup\bigl([\frac{1}{3}, \frac{3}{4})\t\dot
V^{ij}\bigr)$. Define
\e
\begin{split}
\!\!\!(\psi^{(i+1)(j+1)},\hat\psi^{(i+1)(j+1)})\!:\!
\bigl(&W^{(i+1)(j+1)},F^{j+1}\vert_{W^{(i+1)
(j+1)}},t^{j+1}\vert_{W^{(i+1)(j+1)}},\\
&\xi^{j+1}\vert_{W^{(i+1)(j+1)}}\bigr)\longra
\bigl(W^{i+1},F^{i+1},t^{i+1},\xi^{i+1}\bigr)
\end{split}
\label{khBeq15}
\e
by $\psi^{(i\!+\!1)(j\!+\!1)}=(\id_{[0,\frac{3}{4})}\t\phi^{ij})
\vert_{W^{(i\!+\!1)(j\!+\!1)}}$, $\hat\psi^{(i\!+\!1)(j\!+\!1)}
=(\id_{[0,\frac{3}{4})}\t\hat \phi^{ij})
\vert_{\pi^*(W^{(i\!+\!1)(j\!+\!1)})}$. We must check \eq{khBeq15}
satisfies Definition \ref{kh2def12}. Parts (a),(b) follow as
$(\phi^{ij},\smash{\hat\phi^{ij}})$ is a coordinate change, (c) as
$\hat\phi^{ij}\ci s^j\equiv s^i\ci\phi^{ij}$ on $V^{ij}$ and
$\hat\phi^{ij}\ci\ti s^j\equiv\ti s^i\ci\phi^{ij}$ on $\dot V^{ij}$,
and (d) by the definition of $Z$ using gluing maps
$\id_{[0,\frac{3}{4})}\t\phi^{ij}$ in $\approx$. The only subtle
part is (e). This holds as $\ti s^i,\ti s^j$ are $C^1$ close to
$s^i,s^j$, so $t^{i+1},t^{j+1}$ are $C^1$-close to
$s^i\ci\pi_{V^i},s^i\ci\pi_{V^i}$. We need $\d\hat t^{i+1}$ to be an
isomorphism over $(t^{j+1})^{-1}(0)$. This would hold if
$t^{i+1},t^{j+1}$ were replaced by
$s^i\ci\pi_{V^i},s^j\ci\pi_{V^j}$, and being an isomorphism is an
open condition, so $\d\hat t^{i+1}$ is an isomorphism over
$(t^{j+1})^{-1}(0)$ provided $\ti s^i,\ti s^j$ are sufficiently
close to $s^i,s^j$ in~$C^1$.

If $k\notin I$ and $k<i\in I$, define $W^{(i+1)(k+1)}$ in
$W^{k+1}=(\frac{2}{3},1]\t\ti X$ by
$W^{(i+1)(k+1)}=(\frac{2}{3},\frac{3}{4})\t(\ti s^i)^{-1}(0)$.
Define $\psi^{(i+1)(k+1)}: W^{(i+1)(k+1)}\ra W^{i+1}$ to be
inclusion, regarding $W^{(i+1)(k+1)}$ as an open subset of
$(t^{i+1})^{-1}(0)\subseteq W^{i+1}$. Define
$\hat\psi^{(i+1)(k+1)}:F^{k+1}\vert_{W^{ (i+1)(k+1)}}\ra F^{i+1}$ by
$\hat\psi^{(i+1)(k+1)}=z\ci\psi^{(i+1)(k+1)}$, where $z:W^{i+1}\ra
F^{i+1}$ is the zero section. This is valid as $F^{k+1}=W^{k+1}$, so
$F^{k+1} \vert_{W^{(i+1)(k+1)}}=W^{(i+1)(k+1)}$. If $k\in I$ and
$k<i\in I$, we combine the previous two definitions for
$W^{(i+1)(k+1)},\psi^{(i+1)(k+1)}, \hat\psi^{(i+1)(k+1)}$ in the
obvious way. We also define $W^{(k+1)(k+1)}=W^{k+1}$, and
$\psi^{(k+1) (k+1)},\hat\psi^{(k+1)(k+1)}$ to be the identities.
Then \eq{khBeq15} is a coordinate change also for $j=k$ or~$i=j=k$.

Define a Kuranishi structure on $Z$ as follows. There is a natural
projection $\pi_{[0,1]}:Z\ra[0,1]$, with homeomorphisms
$\pi_{[0,1]}^{-1}\bigl([0,\frac{1}{3}]\bigr)\cong[0,\frac{1}{3}]\t
X$ and $\pi_{[0,1]}^{-1}\bigl([\frac{2}{3},1]\bigr)\cong
[\frac{2}{3},1]\t\ti X$. On $\pi_{[0,1]}^{-1}\bigl([0,
\frac{1}{3})\bigr)$, let $Z$ have the Kuranishi structure of
$[0,\frac{1}{3})\t X$. On $\pi_{[0,1]}^{-1}\bigl((\frac{2}{3},1]
\bigr)$, let $Z$ have the Kuranishi structure (manifold structure)
of $(\frac{2}{3},1]\t\ti X$. For each $p\in\pi_{[0,1]}^{-1}\bigl(
[\frac{1}{3},\frac{2}{3}]\bigr)$, let $j\in J$ be least such that
$p\in(t^j)^{-1}(0)$, regarding $(t^j)^{-1}(0)$ as a subset of $Z$ by
\eq{khBeq14}. Then $(W^j,\ldots,\xi^j)$ is a Kuranishi neighbourhood
of $p$. Define the germ of Kuranishi neighbourhoods on $Z$ at $p$ to
be the equivalence class of $(W^j,\ldots, \xi^j)$.

Since we use the Kuranishi structures from $[0,1]\t X$ and
$[0,1]\t\ti X$ on open subsets of $Z$, and choose $j\in J$ least for
each $p\in\pi_{[0,1]}^{-1}\bigl([\frac{1}{3}, \frac{2}{3}]\bigr)$ so
that the function $p\mapsto j$ is upper semicontinuous, to define
the germ of coordinate changes in the Kuranishi structure on $Z$, we
need coordinate changes as follows:
\begin{itemize}
\setlength{\itemsep}{0pt}
\setlength{\parsep}{0pt}
\item[(i)] from $(V_q,\ldots,\psi_q)$ to $(W^j,\ldots,\xi^j)$,
where $(V_q,\ldots,\psi_q)$ lies in the germ of Kuranishi
neighbourhoods on $[0,1]\!\t\! X$ at $q\in[0,\frac{1}{3})\!\t\!\ti
X$, and $(W^j,\ldots,\xi^j)$ is in the germ of Kuranishi
neighbourhoods of $p\in\pi_{[0,1]}^{-1}\bigl(\{\frac{1}{3}\}\bigr)$;
\item[(ii)] from $(V_q,V_q,0,\id_{V_q})$ to $(W^j,\ldots,\xi^j)$,
where $V_q$ is a small open neighbourhood of $q\in
(\frac{2}{3},1]\t\ti X$, so that $(V_q,V_q,0,\id_{V_q})$ lies in the
germ of Kuranishi neighbourhoods on $[0,1]\t\ti X$, and
$(W^j,\ldots,\xi^j)$ is in the germ of Kuranishi neighbourhoods of
$p\in\pi_{[0,1]}^{-1}\bigl(\{\frac{2}{3}\}\bigr)$; and
\item[(iii)] If $p,q\in\pi_{[0,1]}^{-1}\bigl([\frac{1}{3},
\frac{2}{3}]\bigr)$ with $p\in(t^j)^{-1}(0)$, $j\in J$ least, and
$q\in(t^j)^{-1}(0)\cap(t^{j'})^{-1}(0)$ for $j'<j\in J$ least, a
coordinate change from $(W^{jj'},\ab \ldots,\ab\xi^{j'}
\vert_{W^{jj'}})$ to $(W^j,\ldots,\xi^j)$ for $j'<j$ in $J$, where
$(W^{jj'},\ldots,\xi^{j'}\vert_{W^{jj'}})$ is equivalent to
$(W^{j'},F^{j'},t^{j'},\xi^{j'})$ as a Kuranishi neighbourhood of
$q$ in $Z$, and so lies in the germ at~$q$.
\end{itemize}

For (i) the coordinate changes exist as $(W^j,\ldots,\xi^j)$ is a
compatible Kuranishi neighbourhood on $[0,1]\t X$ over
$[0,\frac{1}{3}]\t X$, since it is constructed from a compatible
Kuranishi neighbourhood $(V^{j-1},\ldots,\psi^{j-1})$ on $X$. For
(ii) the coordinate changes are easily constructed from the
inclusions $V_q\subseteq(t^j)^{-1}(0)\subseteq W^j$. For (iii), the
coordinate change is $\smash{(\psi^{jj'},\hat\psi{}^{jj'} )}$ in
\eq{khBeq15}. It is now not difficult to show this defines a {\it
Kuranishi structure} on $Z$, for which the $(W^j,\ldots,\xi^j)$ for
$j\in J$ are compatible Kuranishi neighbourhoods.

Then $Z$ is a compact Kuranishi space with $\vdim Z=k+1$. The
product orientation on $[0,1]\t X$ induces orientations on the
$(W^i,\ldots,\xi^i)$, yielding an orientation on $Z$. Then
\e
\pd Z\cong\bigl(\{1\}\t\ti X\bigr)\amalg-\bigl(\{0\}\t X\bigr)
\amalg -Z^\pd
\label{khBeq16}
\e
in oriented Kuranishi spaces, where $Z^\pd$ is a perturbation of
$[0,1]\t\pd X$, obtained by applying the construction above to $\pd
X$ rather than~$X$.

Define $g^{i+1}:W^{i+1}\ra Y$ by $g^i=f^i\ci\pi_{V^i}$ for $i\in
I\sm\{k\}$, and $g^{k+1}:W^{k+1}\ra Y$ by $g^{k+1}=\ti f\ci\pi_{\ti
X}$ when $k\notin I$ or on $(\frac{2}{3},1]\t\ti X$, and
$g^{k+1}=\ti f\ci\pi_{V^i}$ when $k\in I$ on $\bigl(
[0,\frac{1}{3})\t V^k\bigr)\cup\bigl([\frac{1}{3},\frac{3}{4})\t \ti
V^k\bigr)$. These represent a strongly smooth map $\bs g:Z\ra Y$,
which restricts to $\bs f$ on $\{0\}\t X$ and $\ti f$ on $\{1\}\t\ti
X$, under the identification \eq{khBeq16}. Then $\bs
J=\bigl(J,(W^j,F^j,t^j,\xi^j), g^j:j\in J,\ldots\bigr)$ is a {\it
very good coordinate system\/} for~$(Z,\bs g)$.

Choose continuous functions $\al_1,\al_2,\al_3:[0,1]\ra[0,1]$ with
$\al_1+\al_2+\al_3\equiv 1$ such that $\al_1\equiv 1$ on
$[0,\frac{1}{4}]$, $\al_1\equiv 0$ on $[\frac{1}{3},1]$,
$\al_2\equiv 0$ on $[0,\frac{1}{4}]$, $\al_2\equiv 1$ on
$[\frac{1}{3},\frac{2}{3}]$, $\al_2\equiv 0$ on $[\frac{3}{4},1]$,
$\al_3\equiv 0$ on $[0,\frac{2}{3}]$, and $\al_3\equiv 1$ on
$[\frac{3}{4},1]$. We shall define functions
$\ze_j^{\ti\jmath}:W^{\ti\jmath}\ra[0,1]$ for $j,\ti\jmath\in J$.
For $i,i'\in I\sm\{k\}$, set
\ea
\ze_{i+1}^{i'+1}(w)&\!=\!\al_1\!\ci\!\pi_{[0,1]}(w)\cdot
\eta_i^{i'}\!\ci\!\pi_{V^{i'}}(w)+\al_2\!\ci\!
\pi_{[0,1]}(w)\cdot\dot\eta_i^{i'}\!\ci\!\pi_{V^{i'}}(w),
\label{khBeq17}\\
\begin{split}
\ze_{k+1}^{i'+1}(w)&\!=\!\al_1\!\ci\!\pi_{[0,1]}(w)\cdot
\eta_k^{i'}\!\ci\!\pi_{V^{i'}}(w)\!+\!\al_2\!\ci\!\pi_{[0,1]}(w)
\cdot\dot\eta_k^{i'}\!\ci\!\pi_{V^{i'}}(w)\\
&\qquad\qquad\qquad\qquad\qquad +\al_3\ci\pi_{[0,1]}(w).
\end{split}
\label{khBeq18}
\ea
When $k\in I$, define $\ze_{i+1}^{k+1}$ on $\bigl([0,\frac{1}{3})\t
V^k\bigr)\cup \bigl([\frac{1}{3},\frac{3}{4})\t\ti V^k\bigr)$ by
\eq{khBeq17} and \eq{khBeq18} with $i'=k$. Regarding $\ti X$ as a
union of subsets $(\ti s^{i'})^{-1}(0)$, define $\ze_{i+1}^{k+1}$ on
$(\frac{2}{3},1]\t(\ti s^{i'})^{-1}(0)$ for each $i'\in I$ by
\ea
\ze_{i+1}^{k+1}(w)&=\al_2\ci\pi_{[0,1]}(w)\cdot
\dot\eta_i^{i'}\ci\pi_{(\ti s^{i'})^{-1}(0)},\qquad\qquad i\ne k,
\label{khBeq19}\\
\ze_{k+1}^{k+1}(w)&=\al_2\ci
\pi_{[0,1]}(w)\cdot\dot\eta_k^{i'}\ci\pi_{(\ti
s^{i'})^{-1}(0)}+\al_3\ci\pi_{[0,1]}(w).
\label{khBeq20}
\ea

One can show that \eq{khBeq17}--\eq{khBeq20} are consistent on their
overlaps, so the $\ze_j^{\ti\jmath}$ are well-defined, and that the
$\ze_j^{\ti\jmath}$ are continuous and satisfy $\sum_{j\in
J}\ze_j^{\ti\jmath}\equiv 1$ on $W^{\ti\jmath}$ and
$\ze_j^{\ti\jmath}\ci\psi^{\ti\jmath j'}\equiv
\ze_j^{j'}\vert_{W^{\ti\jmath j'}}$ for $j,\ti\jmath,j'\in J$ with
$j'\le\ti\jmath$. Define $\ze_j:Z\ra [0,1]$ for $j\in J$ by
$\ze_j\vert_{(t^{\ti\jmath})^{-1}(0)}\equiv
\ze_j^{\ti\jmath}\vert_{(t^{\ti\jmath})^{-1}(0)}$ for all
$\ti\jmath\in J$. Then the $\ze_j$ are well-defined and continuous
with $\sum_{j\in J}\ze_j\equiv 1$. Write $\bs\ze=
(\ze_j,\ze_j^{\ti\jmath}:j,\ti\jmath\in J)$. Then $(\bs J,\bs\ze)$
is a {\it really good coordinate system} for~$(Z,\bs g)$.

Furthermore, since $(\bs I,\bs\eta)$ is an excellent coordinate
system, $\pd^lV^i$ has no unnecessary connected components for all
$i\in I$ and $l\ge 0$, so $\pd^lW^{i+1}$ has no unnecessary
connected components. (To see this, note that for $i\ne k$,
$\pd^lW^{i+1}$ is the disjoint union of $l$ copies of
$\pd^{l-1}V^i$, and a piece which retracts onto $\pd^lV^i$. When
$i=k$, the whole of $W^{i+1}$ is necessary as $t^{k+1}\equiv 0$.)
Hence $(\bs J,\bs\ze)$ is an {\it excellent good coordinate system}
for~$(Z,\bs g)$.

Choose some maps $\uH^j:F^j\ra\uP$ for $j\in J$ such that $\ubH=
\bigl(\bs J,\bs\ze,\uH^j:j\in J\bigr)$ is {\it effective
gauge-fixing data} for $(Z,\bs g)$, satisfying the boundary
conditions
\e
\text{$\uH^{i+1}(0,e)=\uG^i(e)$, $i\in I$, $e\in E^i$, and
$\uH^{k+1}(0,e)=\utG^k(e)$, $e\in\ti X$.}
\label{khBeq21}
\e
Ignoring boundary conditions, it is enough to define
$\uH^{i+1}(u,e)$ for $u=0$ or 1 by \eq{khBeq21}, and to define
$\uH^{i+1}(u,e)$ for $u\in(0,1)$ by some injective map
$F^j\cap\pi_{[0,1]}^{-1}\bigl((0,1)\bigr)\ra\R^n\subset\uP$ for some
$n\gg 0$. When we discuss how to choose $Z_{ac},\bs g_{ac},
\ubH_{ac}$ with compatibility between different $a\in A$, $c\in C_a$
below, we will explain in more detail how to choose these maps
$\uH^j$. Then $[Z,\bs g,\ubH]$ is well-defined in
$KC_{k+1}^\ef(Y;R)$. Using \eq{khBeq16} we find that in
$KC_k^\ef(Y;R)$ we have
\e
\begin{split}
\pd\bigl[Z,\bs g,\ubH\bigr]=\bigl[\ti X,\ti
f,\ubtG\bigr]-\bigl[X,\bs f,\ubG\bigr]-\bigl[Z^\pd,\bs
g\vert_{Z^\pd},\ubH\vert_{Z ^\pd}\bigr].
\end{split}
\label{khBeq22}
\e
This concludes our discussion of perturbing a single~$(X,\bs
f,\ubG)$.

Next we explain how to choose such $(\ti X_{ac},\ti
f_{ac},\ubtG_{ac})$ and $(Z_{ac},\bs g_{ac},\ubH_{ac})$ for all
$a\in A$ and $c\in C_a$, with boundary compatibilities. The equation
we need to hold, similar to \eq{khBeq7}, is
\e
\ts\sum_{a\in A}\sum_{c\in C_a}\rho_a\bigl[Z_{ac}^\pd,\bs
g_{ac}\vert_{Z_{ac}^\pd},\ubH_{ac}\vert_{Z_{ac}^\pd}\bigr]=0.
\label{khBeq23}
\e
We would like \eq{khBeq23} to follow from $\sum_{a\in A}\sum_{c\in
C_a}\rho_a\bigl[\pd X_{ac},\bs f_{ac}\vert_{\pd X_{ac}},
\ubG_{ac}\vert_{\pd X_{ac}}\bigr]\!=\!0$, in the same way as
\eq{khBeq7} follows from \eq{khBeq6}. Assuming \eq{khBeq23},
inserting subscripts `$ac$' in \eq{khBeq22}, multiplying it by
$\rho_a$, and summing over all $a\in A$ and $c\in C_a$ yields
\e
\begin{split}
&\pd\Bigl(\ts\sum_{a\in A}\sum_{c\in C_a}\rho_a\bigl[Z_{ac},\bs
g_{ac},\ubH_{ac}\bigr]\Bigr)\\
&=\ts\sum_{a\in A}\sum_{c\in C_a}\rho_{ac}\bigl[\ti X_{ac},\ti
f_{ac},\ubtG_{ac}\bigr]-\sum_{a\in A}\sum_{c\in C_a}\rho_a
\bigl[X_{ac},\bs f_{ac},\ubG_{ac}\bigr].
\end{split}
\label{khBeq24}
\e
Thus $\sum_{a\in A}\sum_{c\in C_a}\rho_a\bigl[\ti X_{ac},\ti
f_{ac},\ubtG_{ac}\bigr]$ is a cycle homologous to $\sum_{a\in
A}\sum_{c\in C_a}\ab\rho_a\ab\bigl[X_{ac},\bs
f_{ac},\ubG_{ac}\bigr]$, and so represents $\al$, by Step~1.

The only arbitrary choices made in the construction of $\ti X,\ti
f,\ubtG,Z,\bs g,\ubH$ above were the sections $\ti s^i$ of
$E^i\vert_{\dot V^i}$ for $i\in I$, and the maps $\utG^k:\ti
X\ra\uP$ and $\uH^j:F^j\ra\uP$ for $j\in J$. Inserting subscripts
`$ac$', we will give a method for choosing $\ti s^i_{ac}$ for $i\in
I_{ac}$, $\utG^k_{ac}$ and $\uH^j_{ac}$ for $j\in J_{ac}$ for all
$a\in A$ and $c\in C_a$ such that \eq{khBeq23} holds. As explained
in Step 2, we will actually construct the smooth $\ti s^i_{ac}$ for
all $c\in C_a$ from piecewise smooth $\ti s^i_a$ on $X_a$, and
similarly, we will define $\utG^k_{ac},\uH^j_{ac}$ for $c\in C_a$
from maps $\utG^k_a$ and $\uH^j_a$ over~$X_a$.

We choose the $\ti s^i_a$ by an inductive procedure on decreasing
codimension like that in Step 1. Use the notation of Step 1 above.
For $m=0,\ldots,M$ and $(a,b)$ in $R^m$, here is our inductive
hypothesis:
\begin{itemize}
\item[$(\star)_{ab}^m$] In the triple $(X_{ab}^m,\bs f_{ab}^m,
\ubG_{ab}^m)$ of Step 1, write $(\bs I_{ab}^m,\bs\eta_{ab}^m)$ for
the excellent coordinate system in $\ubG_{ab}^m$, with $\bs
I_{ab}^m=\bigl(I_{ab}^m,(V_{ab}^{m,i},\ldots,\psi_{ab}^{m,i}),
f_{ab}^{mi}:i\in I_{ab}^m,\ldots\bigr)$, and construct $(\bs{\dot
I}{}_{ab}^m,\bs{\dot \eta}{}_{ab}^m),\dot V^{m,i}_{ab},\dot
V^{m,ij}_{ab},\ldots$ as in Definition \ref{khAdef10}. Let $\bs
T_{ab}^m=\bigl(T_{ab}^{m,i}:i\in I_{ab}^m\bigr)$ be the tent
function for $(X_{ab}^m,\bs f_{ab}^m,\ubG_{ab}^m)$ in $(*)_{ab}^m$
in Step 1, so that $T_{ab}^{m,i}$ is a tent function on~$\dot
V^{m,i}_{ab}$.

Then $\ti s^{m,i}_{ab}$ is a continuous, piecewise smooth section of
$E^{m,i}_{ab}\vert_{\dot V^{m,i}_{ab}}$ subordinate to the tent
function $T^{m,i}_{ab}$ for each $i\in I_{ab}^m$, satisfying:
\begin{itemize}
\setlength{\itemsep}{0pt}
\setlength{\parsep}{0pt}
\item[(a)] $\ti s^{m,i}_{ab}$ is $C^1$ close to $s^{m,i}_{ab}$, and
transverse (both in a piecewise smooth sense).
\item[(b)] If $j<i$ in $I_{ab}^m$ then
$\hat\phi_{ab}^{m,ij}\ci\ti s^{m,j}_{ab}\equiv\ti
s_{ab}^{m,i}\ci\phi_{ab}^{m,ij}$ on $\dot V_{ab}^{m,ij}$.
\item[(c)] As in $(*)_{ab}^m$ we have an isomorphism \eq{khBeq9}
between $\pd(X_{ab}^m,\bs f_{ab}^m,\ubG_{ab}^m)$ and $\coprod_{b'\in
B_{ab}^m}(X_{\bar a\bar b}^{m+1},\ab\bs f_{\bar a\bar
b}^{m+1},\ab\ubG_{\bar a\bar b}^{m+1})$. We require that this
isomorphism should identify $\ti s^{m,i}_{ab}\vert_{\pd\dot
V_{ab}^{m,i}}$ with $\coprod_{b'\in B_{ab}^m}\ti s_{\bar a\bar
b}^{m+1,i-1}$ for all $i\in I_{ab}^m$.
\end{itemize}
\end{itemize}

We now choose $\ti s^{m,i}_{ab}$ satisfying $(\star)_{ab}^m$ by
induction on decreasing $m=M,M-1,\ldots,1,0$, for all $(a,b)$ in
$R^m$ and $i\in I_{ab}^m$. For the first step, when $m=M$, as
$P^{M+1}=\es$, if $(a,b)\in R^M$ then $\pd X_{ab}^M=\es$, so the
boundary conditions (c) in $(\star)_{ab}^M$ are trivial. By
induction on increasing $i\in I_{ab}^M$, we choose smooth ({\it
not\/} piecewise smooth) transverse sections $\ti s^{M,i}_{ab}$ of
$E^{M,i}_{ab}\vert_{\dot V^{M,i}_{ab}}$ satisfying (a),(b) above.
This is possible as in the construction of $\ti X$ without boundary
conditions above, with one extra remark about orbifolds.

As $\ubG^M_{ab}$ is effective gauge-fixing data, $(V_{ab}^{M,i},
E_{ab}^{M,i},s_{ab}^{M,i},\psi_{ab}^{M,i})$ is an {\it effective}
Kuranishi neighbourhood in the sense of Definition \ref{kh3def16}.
Thus, even though $\dot V_{ab}^{M,i}$ may be an orbifold, its
orbifold groups act trivially on the fibres of $E_{ab}^{M,i}$. This
implies that generic smooth sections of $E_{ab}^{M,i}$ are
transverse, so we can indeed choose $\ti s^{M,i}_{ab}$ transverse as
in (a). However, if the orbifold groups of $\dot V_{ab}^{M,i}$ were
allowed to act nontrivially on the fibres of $E_{ab}^{M,i}$, it
could happen that $E_{ab}^{M,i}$ admits no transverse sections $\ti
s^{M,i}_{ab}$ at all, since being invariant under the orbifold
groups might force $\ti s^{M,i}_{ab}$ to be zero on an orbifold
stratum of $\dot V_{ab}^{M,i}$ in a non-transverse way. This
finishes the first step.

For the inductive step, suppose that for some $l=0,\ldots,M-1$ we
have chosen data $\ti s_{ab}^{m,i}$ satisfying $(\star)_{ab}^m$ for
all $l<m\le M$, $(a,b)\in R^m$ and $i\in I_{ab}^m$. We will
construct $\ti s_{ab}^{l,i}$ satisfying $(\star)_{ab}^l$ for all
$(a,b)\in R^l$ and $i\in I_{ab}^l$. First, by induction on
increasing $i\in I_{ab}^l$, choose generic smooth ({\it not\/}
piecewise smooth) sections $\hat s^{l,i}_{ab}$ of
$E^{l,i}_{ab}\vert_{\dot V^{l,i}_{ab}}$ satisfying (a),(b) above
with $\hat s^{l,i}_{ab}$ in place of $\ti s^{m,i}_{ab}$, but {\it
without\/} boundary conditions such as (c). This is possible as in
the first step.

To define $\ti s_{ab}^{l,i}$ we use a version of the piecewise
smooth extension idea of \S\ref{khA14}, in Definition \ref{khAdef6}
and Proposition \ref{khAprop6}. That is, $(\star)_{ab}^l$(c)
prescribes piecewise smooth values for $\ti s_{ab}^{l,i}$ over
$\pd\dot V_{ab}^{l,i}$, and we want to extend these to piecewise
smooth values over $\dot V_{ab}^{l,i}$. We do this by choosing an
arbitrary smooth section $\hat s^{l,i}_{ab}$ on the interior of
$\dot V_{ab}^{l,i}$, corresponding to $\ti\eta$ in Proposition
\ref{khAprop6}, and then interpolating between $\hat s^{l,i}_{ab}$
and the boundary prescribed values for $\ti s_{ab}^{l,i}$ near
$\pd\dot V_{ab}^{l,i}$ using the flow of $v^{l,i}_{ab}$. However, we
have to modify the construction of Proposition \ref{khAprop6}
slightly, as $T_{ab}^{l,i}$ was made not using Definition
\ref{khAdef6}, but with Definition \ref{khAdef14}, a noncompact,
orbifold version of Definition~\ref{khAdef7}.

Define subsets $W_1,W_2$ of $\dot V_{ab}^{l,i}$ by
\begin{align*}
W_1&=\bigl\{w\in\dot V_{ab}^{l,i}:\min T_{ab}^{l,i}(w)\ge
D\bigr\},\\
W_2&=\bigl\{w\in\dot V_{ab}^{l,i}:\text{$w=\exp(tv_{ab}^{l,i})x$,
$t\in[0,2)$, $(x,B)\in\pd U_{ab}^{l,i}$}\bigr\}.
\end{align*}
Then $W_1$ is closed in $\dot V_{ab}^{l,i}$ as $\min T_{ab}^{l,i}$
is continuous, and $W_2$ is open in $\dot V_{ab}^{l,i}$, and
$W_1\cup W_2=\dot V_{ab}^{l,i}$, since the branches $u+\th(2-t)$ in
\eq{khAeq36} are supported on $W_2$, so for $v\notin W_2$ we see
that $\min T_{ab}^{l,i}(v)$ is of the form $D+l^2$ or $D+\ep^2$, so
$\min T_{ab}^{l,i}(v)\ge D$ and~$v\in W_1$.

Write $\ti s_{ab}^{\prime l,i}$ for the piecewise smooth section of
$E_{ab}^{l,i}\vert_{\pd\dot V_{ab}^{l,i}}$ prescribed as boundary
conditions for $\ti s_{ab}^{l,i}$ in $(\star)_{ab}^l$(c). Then as in
\S\ref{khA34}--\S\ref{khA35}, $\ti s_{ab}^{\prime l,i}$ extends to
$\pd U^{l,i}_{ab}$, for some open neighbourhood $U^{l,i}_{ab}$ of
$\overline{\dot V_{ab}^{l,i}}$ in $V_{ab}^{l,i}$, such that
Definition \ref{khAdef14}(a)--(d) hold for $U^{l,i}_{ab}$ and the
vector field $v^{l,i}_{ab}$ on $V_{ab}^{l,i}$ used to construct
$T_{ab}^{l,i}$. From $(\star)_{\bar a\bar b}^{l+1}$(c) for $(\bar
a,\bar b)=\phi^{l+1}(a,b')$, all $b'\in B_{ab}^l$, we can deduce
that $\ti s_{ab}^{\prime l,i}\vert_{\pd^2U_{ab}^{l,i}}$ is invariant
under the natural involution $\si:\pd^2U_{ab}^{l,i}\ra\pd^2\dot
V_{ab}^{l,i}$, in the same way as we showed that $\min\bs
T_{ab}^{\prime l}\vert_{\pd^2X_{ab}^l}$ is invariant under the
natural involution $\bs\si:\pd^2X_{ab}^l\ra\pd^2X_{ab}^l$ in the
proof of $(*)_{ab}^m$ in~\S\ref{khB1}.

We will extend $\ti s_{ab}^{\prime l,i}$ to a continuous, piecewise
smooth section $\check s^{l,i}_{ab}$ of $E_{ab}^{l,i}$ over $W_2$
subordinate to the tent function $T_{ab}^{l,i}\vert_{W_2}$. Let
$\nabla_{ab}^{l,i}$ be a connection on the orbifold vector bundle
$E_{ab}^{l,i}\ra V_{ab}^{l,i}$. If $w\in W_2$ then
$w=\exp(tv_{ab}^{l,i})x$ for some $t\in[0,2)$ and $(x,B)\in\pd
U_{ab}^{l,i}$. Here $t,x$ are uniquely determined by $w$, as the
flow-line of $-v_{ab}^{l,i}$ first hits $\io(\pd U_{ab}^{l,i})$ at
$x$ in time $t$. The choice of local boundary component $B$ of
$U_{ab}^{l,i}$ at $x$ is not uniquely determined.

However, since $\ti s_{ab}^{\prime l,i}\vert_{\pd^2U_{ab}^{l,i}}$ is
invariant under the natural involution $\si:\pd^2U_{ab}^{l,i}\ab
\ra\pd^2\dot V_{ab}^{l,i}$, it follows that $\ti s_{ab}^{\prime
l,i}(x,B)\in E_{ab}^{l,i}\vert_x$ depends only on $x$, not on the
choice of $B$. Thus $\ti s_{ab}^{\prime l,i}(x,B)$ is also uniquely
determined by $w$. Define
\e
\check
s^{l,i}_{ab}(w)=s^{l,i}_{ab}(w)+P_{\nabla_{ab}^{l,i}}\bigl(\ti
s_{ab}^{\prime l,i}(x,B)-s^{l,i}_{ab}(x)\bigr),
\label{khBeq25}
\e
where $P_{\nabla_{ab}^{l,i}}:E_{ab}^{l,i}\vert_x\ra
E_{ab}^{l,i}\vert_w$ is the {\it parallel transport map} along the
flow-line of $v_{ab}^{l,i}$ from $x$ to $w$ using the connection
$\nabla_{ab}^{l,i}$. Considering the construction of $T_{ab}^{l,i}$,
we see that $\check s^{l,i}_{ab}$ is continuous and piecewise smooth
subordinate to $T_{ab}^{l,i}\vert_{W_2}$, with~$\check
s^{l,i}_{ab}\vert_{\pd\dot V_{ab}^{l,i}}\equiv\ti s_{ab}^{\prime
l,i}$.

By the construction in \S\ref{khA34}--\S\ref{khA35}, $\min
T_{ab}^{l,i}\vert_{\pd\dot V_{ab}^{l,i}}$ is bounded above, and $D$
is strictly greater than this bound. Thus there exists $D'<D$ with
$\min T_{ab}^{l,i}\vert_{\pd\dot V_{ab}^{l,i}}\ab\le\ab D'$. Choose
a smooth function $\eta:\R\ra[0,1]$ with $\eta(t)=0$ for $t\le D'$
and $\eta(t)=1$ for $t\ge D$. Define a section $\ti s^{l,i}_{ab}$ of
$E_{ab}^{l,i}$ over $\dot V_{ab}^{l,i}$ by
\e
\ti s^{l,i}_{ab}(w)\!=\!\begin{cases} \hat s^{l,i}_{ab}(w), & w\!\in\! W_1,\\
\bigl(\eta\!\ci\!\min T_{ab}^{l,i}(w)\bigr)\hat s^{l,i}_{ab}(w)\!+\!
\bigl(1\!-\!\eta\!\ci\!\min T_{ab}^{l,i}(w)\bigr)\check
s^{l,i}_{ab}(w), & w\!\in\!W_2.
\end{cases}\!
\label{khBeq26}
\e
On the overlap $W_1\cap W_2$ we have $\min T_{ab}^{l,i}(w)\ge D$ and
$\eta\ci\min T_{ab}^{l,i}(w)=1$, so $\ti s^{l,i}_{ab}$ is
well-defined. As $\hat s^{l,i}_{ab},\eta$ are smooth and $\min
T_{ab}^{l,i},\check s^{l,i}_{ab}$ are continuous and piecewise
smooth subordinate to $T_{ab}^{l,i}$, we see that $\ti s^{l,i}_{ab}$
is {\it continuous} and {\it piecewise smooth subordinate
to}~$T_{ab}^{l,i}$.

On $\pd\dot V_{ab}^{l,i}$ we have $\min T_{ab}^{l,i}\le D'$,
$\eta\ci\min T_{ab}^{l,i}\equiv 0$ and $\check s^{l,i}_{ab}
\vert_{\pd\dot V_{ab}^{l,i}}\equiv\ti s_{ab}^{\prime l,i}$, so $\ti
s^{l,i}_{ab}\vert_{\pd\dot V_{ab}^{l,i}}\equiv\ti s_{ab}^{\prime
l,i}$. Hence $\ti s^{l,i}_{ab}$ satisfies $(\star)_{ab}^l$(c). For
part (a), from $(\star)_{\bar a\bar b}^{l+1}$(a) for $(\bar a,\bar
b)=\phi^{l+1}(a,b')$, all $b'\in B_{ab}^l$, it follows that $\ti
s_{ab}^{\prime l,i}$ is $C^1$-close to $s_{ab}^{l,i}\vert_{\pd\dot
V^i}$ and transverse, in a piecewise smooth sense. Thus in
\eq{khBeq25}, $\ti s_{ab}^{\prime l,i}(x,B)-s^{l,i}_{ab}(x)$ is
$C^1$-small as a function of $(x,B)\in\pd\dot V_{ab}^{l,i}$, which
implies that $\check s^{l,i}_{ab}(w)-s^{l,i}_{ab}(w)$ is $C^1$-small
as a function of $w\in W_2$. So $\check s^{l,i}_{ab}$ is $C^1$-close
to $s^{l,i}_{ab}$. Also $\check s^{l,i}_{ab}$ was chosen above to be
$C^1$-close to $s^{l,i}_{ab}$. Therefore $\ti s^{l,i}_{ab}$ in
\eq{khBeq26} is also $C^1$-close to $s^{l,i}_{ab}$.

As $\ti s_{ab}^{\prime l,i}$ and $\hat s^{l,i}_{ab}$ were chosen
generic and transverse, we see that $\ti s^{l,i}_{ab}$ is
transverse, and $(\star)_{ab}^l$(a) holds. For part (b), the $\hat
s^{l,i}_{ab}$ for $i\in I_{ab}^l$ were chosen to satisfy
$\hat\phi_{ab}^{l,ij}\ci\hat s^{l,j}_{ab}\equiv\hat
s_{ab}^{l,i}\ci\phi_{ab}^{l,ij}$ if $j<i$ in $I_{ab}^l$. Also
$(\star)_{\bar a\bar b}^{l+1}$(b) for $(\bar a,\bar
b)=\phi^{l+1}(a,b')$, all $b'\in B_{ab}^l$, imply that the $\ti
s_{ab}^{\prime l,i}$ satisfy $\hat\phi_{ab}^{l,ij}\ci\ti s^{\prime
l,j}_{ab}\equiv\ti s_{ab}^{\prime l,i}\ci\phi_{ab}^{l,ij}$ if $j<i$
in $I_{ab}^l$. Now $\check s^{l,i}_{ab}$ was constructed using $\ti
s_{ab}^{\prime l,i}$, the vector field $v_{ab}^{l,i}$, and the
connection $\nabla_{ab}^{l,i}$ on $E_{ab}^{l,i}$. We have
$(\phi^{l,ij}_{ab})_* (v^{l,j}_{ab})\equiv v_{ab}^{l,i}$ if $j<i$ in
$I_{ab}^l$ by Definition \ref{khAdef14}(c). Thus, if the connections
$\nabla_{ab}^{l,i},\nabla_{ab}^{l,j}$ on $E_{ab}^{l,i},E_{ab}^{l,j}$
are chosen compatible under $(\phi_{ab}^{l,ij},\hat
\phi_{ab}^{l,ij})$, which we can do, then
$\hat\phi_{ab}^{l,ij}\ci\check s^{l,j}_{ab}\equiv\check
s_{ab}^{l,i}\ci\phi_{ab}^{l,ij}$ if $j<i$ in $I_{ab}^l$. So
$(\star)_{ab}^l$(b) follows, by \eq{khBeq26}. This completes the
inductive step. Hence by induction we can choose $\ti s^{m,i}_{ab}$
satisfying $(\star)_{ab}^m$ for all~$m,a,b$.

Now $X_a=\coprod_{b=1}^{n_a^0}X_{ab}^0$, and for each
$b=1,\ldots,n_a^0$, if $\phi^0(a,b)=(\bar a,\bar b)$ in $R^0$, we
have an isomorphism $(\bs a,\bs b)^0_{ab}:(X_{ab}^0,\ab\bs
f_{ab}^0,\ab\ubG_{ab}^0)\ab\ra(X_{\bar a\bar b}^0,\ab\bs f_{\bar
a\bar b}^0,\ab\ubG_{\bar a\bar b}^0)$, and piecewise smooth sections
$\ti s_{\bar a\bar b}^{0,i}$ of $E_{\bar a\bar b}^{0,i}\ra\dot
V_{\bar a\bar b}^{0,i}$ for all $i\in I_{\bar a\bar b}^{0,i}$.
Define sections $\ti s_a^i$ of $E_a^i\ra\dot V_a^i$ by $\hat
b_{ab}^{0,i}\ci\ti s_a^i\vert_{\dot V_{ab}^{0,i}}=\ti s_{\bar a\bar
b}^{0,i}\ci b_{ab}^{0,i}$ for $b=1,\ldots,n_a^0$, where
$(b_{ab}^{0,i},\hat b_{ab}^{0,i}):\bigl(V_{ab}^{0,i},\ldots,
\psi_{ab}^{0,i}\bigr)\ra\bigl(V_{\bar a\bar b}^{0,i},\ldots,
\psi_{\bar a\bar b}^{0,i}\bigr)$ is the coordinate change in $(\bs
a,\bs b)^0_{ab}$. Then $\ti s_a^i$ is continuous, and piecewise
smooth subordinate to~$T^i_a$.

For each $c\in C_a$ we have projections $\pi_a:V_{ac}^i\ra\dot
V_a^i$ for all $i\in I_{ac}$, with $E_{ac}^i=\pi_a^*(E_a^i)$. Define
a section $\ti s_{ac}^i$ of $E_{ac}^i\ra V_{ac}^i$ by $\ti
s_{ac}^i=\ti s_a^i\ci\pi_a$. Since $\ti s_a^i$ is piecewise smooth
subordinate to $T^i_a$ and $V_{ac}^i$ for $i\in I_{ac}$ are the
pieces in the decomposition of $\dot V_a^i$ induced by $T^i_a$, it
follows that $\ti s_{ac}^i$ is a {\it smooth\/} section, not just
piecewise smooth. From $(\star)_{\bar a\bar b}^0$(a),(b) we see that
$\ti s_{ac}^i$ is $C^1$ close to $s_{ac}^i$ and transverse, and
$j<i$ in $I_{ac}$ then $\hat\phi_{ac}^{ij}\ci\ti s^j_{ac}\equiv\ti
s_{ac}^i\ci\phi_{ac}^{ij}$ on~$V_{ac}^{ij}$.

We have now constructed perturbations $\ti s_{ac}^i$ of $s_{ac}^i$
for all $a\in A$, $c\in C_a$ and $i\in I_{ac}$ satisfying all the
conditions we need. Thus we can apply the construction of $\ti X,\ti
f,\ubtG,Z,\bs g,\ubH$ from $(X,\bs f,\ubG)$ above to $(X_{ac},\bs
f_{ac},\ubG_{ac})$ for $a\in A$ and $c\in C_a$ using these $\ti
s_{ac}^i$. This gives a compact, oriented $k$-manifold $\ti X_{ac}$,
a smooth map $\ti f_{ac}:\ti X_{ac}\ra Y$, a compact oriented
Kuranishi space $Z_{ac}$ and a strongly smooth map $\bs
g_{ac}:Z_{ac}\ra Y$.

To define the effective gauge-fixing data $\ubtG_{ac},\ubH_{ac}$ for
$(\ti X_{ac},\ti f_{ac}),(Z_{ac},\bs g_{ac})$ we need in addition to
choose maps $\utG^k_{ac}:\ti X_{ac}\ra\uP$ and
$\uH^j_{ac}:F^j_{ac}\ra\uP$ for $j\in J_{ac}$. To ensure the
necessary boundary compatibilities hold, these must also be chosen
by reverse induction on codimension by a procedure similar to
$(\star)_{ab}^m$ above. In effect, from the piecewise smooth
perturbations $\ti s_{ab}^{m,i}$ for $i\in I_{ab}^m$ in
$(\star)_{ab}^m$ above, we construct piecewise smooth $\ti
X_{ab}^m,\ti f_{ab}^m,Z_{ab}^m, \bs g_{ab}^m$ for $m=M,M-1,\ldots,0$
and $(a,b)\in R^m$, and then in an analogue of $(\star)_{ab}^m$ we
choose $\utG^{m,k-m}_{ab}:\ti X_{ab}^m\ra\uP$ and
$\uH^{m,j}_{ab}:F^{m,j}_{ab}\ra\uP$ for $j\in J_{ab}^m$ for all
$m=M,M-1,\ldots,0$, $(a,b)\in R^m$ and $j\in J_{ab}^m$ satisfying
boundary conditions similar to $(\star)_{ab}^m$(c). This is much
simpler than the inductive proof above --- at each step, the values
of $\utG^{m,k-m}_{ab}\vert_{\pd\ti X_{ab}^m}$ and
$\uH^{m,j}_{ab}\vert_{\pd F^{m,j}_{ab}}$ are prescribed, and we
choose arbitrary injective extensions $\utG^{m,k-m}_{ab}:(\ti
X_{ab}^m)^\ci\ra\R^n\subset\uP$ and
$\uH^{m,j}_{ab}:(F^{m,j}_{ab})^\ci\ra\R^{n'}\subset\uP$ for $n,n'\gg
0$ and all $j\in J_{ab}^m$, satisfying no compatibility conditions
at all.

This completes the construction of $(\ti X_{ac},\ti
f_{ac},\ubtG_{ac})$ and $(Z_{ac},\bs g_{ac},\ubH_{ac})$ for all
$a\in A$ and $c\in C_a$. We claim that the boundary compatibilities
we imposed on the data $\ti s_{ac}^i,\utG^k_{ac},\uH^j_{ac}$ during
the inductive construction $(\star)_{ab}^m$ for the $\ti s_{ac}^i$
and the analogues for $\utG^k_{ac},\uH^j_{ac}$ imply that
\eq{khBeq23} holds. To see this, note that $Z_{ac}^\pd$ is a
perturbation of $[0,1]\t\pd X_{ac}$. Now the $X_{ac}$ for $c\in C_a$
are the result of cutting $X_a$ into pieces using $\bs T_a$.
Therefore we can split $\pd X_{ac}=(\pd X_{ac})^{\rm int}\amalg(\pd
X_{ac})^{\rm ext}$, where $(\pd X_{ac})^{\rm int}$ is the {\it
internal\/} boundary, the new boundaries creating in the interior of
$X_a$ by cutting $X_a$ into pieces, and $(\pd X_{ac})^{\rm ext}$ is
the {\it external\/} boundary, part of $\pd X_a$. Thus there is a
corresponding splitting~$Z_{ac}^\pd= (Z_{ac}^\pd)^{\rm
int}\amalg(Z_{ac}^\pd)^{\rm ext}$.

Each component of $\bigl[(\pd X_{ac})^{\rm int},\bs
f_{ac}\vert_{(\pd X_{ac})^{\rm int}},\ab \ubG_{ac}\vert_{(\pd
X_{ac})^{\rm int}}\bigr]$ is cancelled by a component of $\bigl[(\pd
X_{ac'})^{\rm int},\bs f_{ac'}\vert_{(\pd X_{ac'})^{\rm
int}},\ubG_{ac'}\vert_{(\pd X_{ac'})^{\rm int}}\bigr]$ for some
other $c'\in C_a$, the `other side of the cut'. Therefore for each
$a\in A$ we have $\sum_{c\in C_a}\bigl[(\pd X_{ac})^{\rm int},\ab
\bs f_{ac}\vert_{(\pd X_{ac})^{\rm int}},\ubG_{ac}\vert_{(\pd
X_{ac})^{\rm int}}\bigr]=0$. Because the perturbations $\ti s_a^i$
for $i\in I_a$ were chosen {\it continuous} and piecewise smooth,
and the $\utG^k_{ac},\uH^j_{ac}$ are also restrictions of functions
independent of $c\in C_a$, the construction of $Z_{ac}^\pd,\bs
g_{ac}\vert_{Z_{ac}^\pd},\ubH_{ac}\vert_{Z_{ac}^\pd}$ over
$[0,1]\t\pd X_{ac}$ is compatible with this cancellation of internal
boundary components. Hence for each $a\in A$ we have
\e
\ts\sum_{c\in C_a}\bigl[(Z_{ac}^\pd)^{\rm int},\ab \bs
g_{ac}\vert_{(Z_{ac}^\pd)^{\rm
int}},\ubH_{ac}\vert_{(Z_{ac}^\pd)^{\rm int}}\bigr]=0.
\label{khBeq27}
\e

For the external boundaries, using the argument of \eq{khBeq11} and
\eq{khBeq12} we have
\e
\begin{split}
&\ts\sum_{a\in A}\rho_a\bigl[(Z_{ac}^\pd)^{\rm ext},\ab \bs
g_{ac}\vert_{(Z_{ac}^\pd)^{\rm ext}},\ubH_{ac}\vert_{(Z_{ac}^\pd
)^{\rm ext}}\bigr]=\\
&\sum\limits_{\!\!\!\!\!\!\!(\bar a,\bar b)\in R^1\!\!\!\!\!\!\!}
\,\,\,\,
\raisebox{-9pt}{\begin{Large}$\displaystyle\biggl\{$\end{Large}}
\sum\limits_{\begin{subarray}{l}a\in A,\;
b=1,\ldots,n_a^1:\\ \phi^1(a,b)=(\bar a,\bar
b)\end{subarray}\!\!\!\!\!\!\!\!\!\!\!}\rho_a\ep_{ab}^1
\raisebox{-9pt}{\begin{Large}$\displaystyle\biggr\}$\end{Large}}
\sum_{\bar c\in C_{\bar a\bar b}^1}\bigl[Z_{\bar a\bar b\bar c}^1,
\bs g_{\bar a\bar b\bar c}^1,\ubH_{\bar a\bar b\bar c}^1\bigr]=0,
\end{split}
\label{khBeq28}
\e
where $C_{\bar a\bar b}^1$ is an indexing set for the decomposition
induced by the tent function $\bs T_{\bar a\bar b}^1$ for $(X_{\bar
a\bar b}^1,\bs f_{\bar a\bar b}^1,\ubG_{\bar a\bar b}^1)$ chosen in
$(*)_{\bar a\bar b}^1$ in \S\ref{khB1}, and $(Z_{\bar a\bar b\bar
c}^1,\bs g_{\bar a\bar b\bar c}^1,\ubH_{\bar a\bar b\bar c}^1)$ for
$C_{\bar a\bar b}^1$ are the corresponding chains constructed above
using the piecewise smooth sections $\ti s_{\bar a\bar b}^1$ in
$(\star)_{\bar a\bar b}^1$. The coefficient $\{\cdots\}$ on the
second line of \eq{khBeq28} is zero as in the proof of \eq{khBeq12}
in \S\ref{khB1}. Combining \eq{khBeq27}, \eq{khBeq28} and
$Z_{ac}^\pd=(Z_{ac}^\pd)^{\rm int}\amalg(Z_{ac}^\pd)^{\rm ext}$
proves \eq{khBeq23}, and completes Step~2.

\subsection{Step 3: triangulating the $\ti X_a$ with simplices
$\si_{ac}(\De_k)$}
\label{khB3}

As in Step 3 above we now {\it change notation} from $\ti X_{ac}$
for $a\in A$ and $c\in C_a$ to $\ti X_a$ for $a\in A$. So, we
represent $\al$ by a cycle $\sum_{a\in A}\rho_a[\ti X_a,\ti
f_a,\ubtG_a]\in KC_k^\ef(Y;R)$, with each $\ti X_a$ a compact,
oriented $k$-manifold with g-corners, $\ti f_a:\ti X_a\ra Y$ a
smooth map, and $\ubtG_a=(\bs{\ti I}_a,\bs{\ti\eta}_a,\utG^k_a)$
effective gauge-fixing data for $(\ti X_a,\ti f_a)$, where $\bs{\ti
I}_a$ has indexing set $\ti I_a=\{k\}$, one Kuranishi neighbourhood
$(\ti V_a^k,\ti E_a^k,\ti s_a^k,\ti\psi_a^k)=(\ti X_a,\ti
X_a,0,\id_{\ti X_a})$, and
$\bs{\ti\eta}_a=(\ti\eta_{a,k},\ti\eta_{a,k}^k)$
with~$\ti\eta_{a,k}\equiv 1\equiv\ti\eta_{a,k}^k$.

We will use Theorem \ref{khAthm2} to choose a tent function $\ti
T_a:\ti X_a\ra F\bigl([1,\iy)\bigr)$ for each $a\in A$ such that the
components $\ti X_{ac}$, $c\in C_a$ of $\pd Z_{\ti X_a,T_a}$ in
\eq{khAeq5} are all diffeomorphic to the $k$-simplex $\De_k$, with
chosen diffeomorphisms $\si_{ac}:\De_k\ra\ti X_{ac}$ for $c\in C_a$.
(Note that these $C_a$ and $\ti X_{ac}$ are different to those in
Steps 1 and 2.) To make these choices of $\ti T_a,\si_{ac}$
compatible with the boundary relations in $\pd\bigl(\sum_{a\in
A}\rho_a[\ti X_a,\ti f_a,\ubtG_a]\bigr)=0$, we use reverse induction
on codimension in a similar way to Steps 1 and 2.

Define notation $n_a^m$, $(\ti X_{ab}^m,\ti f_{ab}^m,
\ubtG_{ab}^m)$, $P^m,Q^m,R^m,\ab(a^q,b^q),\ab\phi^m,\ab(\bs a,\bs
b)_{ab}^m,$ $\ep_{ab}^m,\ab B_{ab}^m$ as in Step 1 above, but using
$(\ti X_a,\ti f_a,\ubtG_a)$ for $a\in A$ rather than $(X_a,\ab\bs
f_a,\ab\ubG_a)$ for $a\in A$, so that $(\pd^m\ti X_a,\ti f_a\vert_{
\pd^m\ti X_a},\ubtG_a\vert_{\pd^m\ti X_a})=\coprod_{b=1,\ldots,
n_a^m}(\ti X_{ab}^m,\ti f_{ab}^m,\ubtG_{ab}^m)$, and so on. As $\ti
X_{ab}^m$ is a manifold of dimension $k-m$ we have $\ti
X_{ab}^m=\emptyset$ for $m>k$, so we fix $M=k$. By induction on
decreasing $m=k,k-1,\ldots,1,0$, for all $(a,b)$ in $R^m$, we choose
tent functions $\ti T_{ab}^m$ on $\ti X_{ab}^m$ using Theorem
\ref{khAthm2}, such that the expression \eq{khAeq5} for $\pd Z_{\ti
X_{ab}^m,\ti T_{ab}^m}$ involves components $\ti X_{abc}^m$ for
$c\in C_{ab}^m$, with $C_{ab}^m$ a finite indexing set, and for all
$(a,b)$ in $R^m$ and $c\in C_{ab}^m$ we choose diffeomorphisms
$\si_{abc}^m:\De_{k-m}\ra\ti X_{abc}^m$, satisfying the inductive
hypothesis:
\begin{itemize}
\item[$(\dag)_{ab}^m$] By definition of $B_{ab}^m$ we have
$\pd\ti X_{ab}^m=\coprod_{b'\in B_{ab}^m}\ti X_{ab'}^{m+1}$, and if
$\phi^{m+1}(a,b')=(\bar a,\bar b)$ in $R^{m+1}$, we have an
isomorphism $(\bs a,\bs b)^{m+1}_{ab'}:(\ti X_{ab'}^{m+1},\ti
f_{ab'}^{m+1},\ab\ubtG_{ab'}^{m+1})\ab\ra(\ti X_{\bar a\bar
b}^{m+1},\ab\ti f_{\bar a\bar b}^{m+1},\ab\ubtG_{\bar a\bar
b}^{m+1})$. So as for \eq{khBeq9} there is a diffeomorphism
\e
\ts\coprod_{b'\in B_{ab}^m}\bs a^{m+1}_{ab'}:\pd\ti X_{ab}^m\longra
\coprod_{\text{$b'\in B_{ab}^m$: define $(\bar a,\bar
b)=\phi^{m+1}(a,b')$}}\ti X_{\bar a\bar b}^{m+1}.
\label{khBeq29}
\e
On the left hand side $\pd\ti X_{ab}^m$ of \eq{khBeq29} we have a
tent function $\ti T_{ab}^m\vert_{\pd\ti X_{ab}^m}$. On the right
hand side we have a tent function $\coprod_{b'}\ti T_{\bar a\bar
b}^{m+1}$. We require:
\begin{itemize}
\setlength{\itemsep}{0pt}
\setlength{\parsep}{0pt}
\item[(a)] Equation \eq{khBeq29} must identify $\min\bigl(\ti
T_{ab}^m\vert_{\pd\ti X_{ab}^m}\bigr)$ with~$\min\bigl(
\coprod_{b'}\ti T_{\bar a\bar b}^{m+1}\bigr)$.
\item[(b)] Suppose the natural involution $\si:\pd^2Z_{\ti
X_{ab}^m, \ti T_{ab}^m}\ra\ab\pd^2Z_{\ti X_{ab}^m,\ti T_{ab}^m}$
exchanges two $(k-m-1)$-simplices $\si_{abc}^m\ci
F_j^{k-m}(\De_{k-m-1})$ and\/ $\si_{abc'}^m\ci
F_j^{k-m}(\De_{k-m-1})$ for $c,c'\in C_{ab}^m$. Then $j=j',$ $c\ne
c'$ and\/ $\si\ci\si_{abc}^m\ci F_j^{k-m}\equiv\si_{abc'}^m\ci
F_j^{k-m}$ as maps~$\De_{k-m-1}\ra\pd^2Z_{\ti X_{ab}^m,\ti
T_{ab}^m}$.
\item[(c)] If $b'\in B_{ab}^m$, $(\bar a,\bar b)=\phi^{m+1}(a,b')$
and $\bar c\in C_{\bar a\bar b}^{m+1}$ then $\ti X_{\bar a\bar b\bar
c}^{m+1}=\si_{\bar a\bar b\bar c}^{m+1}(\De_{k-m-1})$ is a component
of $\pd Z_{\ti X_{\bar a\bar b}^{m+1},\ti T_{\bar a\bar b}^{m+1}}$.
Applying the projection $\pi:\pd Z_{\ti X_{\bar a\bar b}^{m+1},\ti
T_{\bar a\bar b}^{m+1}}\ra\ti X_{\bar a\bar b}^{m+1}$ gives an
embedding $\pi\ci\si_{\bar a\bar b\bar c}^{m+1}:\De_{k-m-1}\ra\ti
X_{\bar a\bar b}^{m+1}$, so $\bigl((\bs a,\bs b)^{m+1}_{ab'}\bigr)
{}^{-1}\ci\pi\ci\si_{\bar a\bar b\bar c}^{m+1}$ is an embedding
$\De_{k-m-1}\ra\ti X_{ab'}^{m+1}\subset\pd\ti X_{ab}^m$. Composing
with $\io:\pd\ti X_{ab}^m\ra\ti X_{ab}^m$ yields an embedding
$\io\ci\bigl((\bs a,\bs b)^{m+1}_{ab'}\bigr) {}^{-1}\ci\pi\ci
\si_{\bar a\bar b\bar c}^{m+1}:\De_{k-m-1}\ra\ti X_{ab}^m$. There
should exist some unique $c\in C_{ab}^m$ such that
$\pi\ci\si_{abc}^m\ci F_{k-m}^{k-m}\equiv\io\ci\bigl((\bs a,\bs
b)^{m+1}_{ab'}\bigr) {}^{-1}\ci\pi\ci \si_{\bar a\bar b\bar
c}^{m+1}$ as embeddings $\De_{k-m-1}\ra\ti X_{ab}^m$.
\end{itemize}
\end{itemize}
Here part (a) is almost the same as $(*)_{ab}^m$. Part (b) gives
compatibility between the `internal' codimension one faces of
$\si_{abc}^m(\De_{k-m})$ within a fixed $X_{ab}^m$, and corresponds
to the boundary compatibility between the $\si_c$ in Theorem
\ref{khAthm1}. Part (c) gives compatibility of the `external'
codimension one faces of $\si_{abc}^m(\De_{k-m})$ between two
different $X_{ab}^m,X_{a'b'}^m$, and corresponds to the boundary
compatibility between $\si'_{c'},\si_c$ in Theorem~\ref{khAthm2}.

The first step, when $m=k$, is trivial. Each $\ti X_{ab}^k$ is a
compact, connected 0-manifold, that is, a single point. We define
$\ti T_{ab}^k\equiv\{1\}$, $C_{ab}^k=\{1\}$ for all $(a,b)\in R^k$,
and take $\si_{ab1}:\De_0\ra\ti X_{ab1}^k=\{1\}\t\ti X_{ab}^k$ to be
the unique diffeomorphism between the two points. Parts (a)--(c) are
vacuous as~$\pd\ti X_{ab}^k=\emptyset$.

For the inductive step, suppose that for some $l=0,\ldots,k-1$ we
have chosen $\ti T_{ab}^m,C_{ab}^m$ and $\si_{abc}^m$ satisfying
$(\dag)_{ab}^m$ for all $l<m\le k$, $(a,b)\in R^m$ and $c\in
C_{ab}^m$. For each $(a,b)\in R^l$, we will apply Theorem
\ref{khAthm2} to choose $\ti T_{ab}^l,C_{ab}^l$ and $\si_{abc}^l$
for $c\in C_{ab}^l$ satisfying $(\dag)_{ab}^l$. We have a tent
function $T'=\coprod_{b'\in B_{ab}^m}\ti T_{\bar a\bar
b}^{m+1}\ci\bs a^{m+1}_{ab'}$ on $\pd\ti T_{ab}^m$, and part (a) of
$(\dag)_{ab}^m$ says that $\min\ti T_{ab}^m\vert_{\pd
X_{ab}^m}\equiv\min T'$, which corresponds to $\min T\vert_{\pd
X}\equiv\min T'$ in Theorem~\ref{khAthm2}.

We must verify the hypotheses Theorem \ref{khAthm2}(i)--(iii). Parts
(i),(ii) hold by the argument showing that $\min\bs T_{ab}^{\prime
l}\vert_{\pd^2X_{ab}^l}$ is invariant under the natural involution
$\bs\si:\pd^2X_{ab}^l\ra\pd^2X_{ab}^l$ in the proof of the inductive
step of $(*)_{ab}^m$ in \S\ref{khB1}. Part (iii) is a little more
subtle. To satisfy it, before starting the induction, we first
choose Riemannian metrics $g_{ab}^m$ on $\ti X_{ab}^m$ for all
$m=0,\ldots,k$ and $(a,b)\in R^m$, and then we choose $\ep>0$ very
small compared to length scales in $(\ti X_{ab}^m,g_{ab}^m)$ and
constants comparing $g_{ab}^m\vert_{\ti X_{ab'}^{m+1}}$ and $g_{\bar
a\bar b}^{m+1}$ under the inclusion $\ti X_{ab'}^{m+1}
\subseteq\pd\ti X_{ab}^m$ and diffeomorphism $\bs a^{m+1}_{ab'}:\ti
X_{ab'}^{m+1}\ra\ti X_{\bar a\bar b}^{m+1}$, for all
$m,(a,b),b',(\bar a,\bar b)$. Then we use these $g_{ab}^m$ and $\ep$
to construct $\ti T_{ab}^l$ in Theorem~\ref{khAthm2}.

The point is that as the $\ep$ used to construct $T'$ and $\si_{c'}$
on $\pd\ti X_{ab}^m$ has been chosen small compared not just to
geometry on $\pd\ti X_{ab}^m$, but also to geometry on $\ti
X_{ab}^m$ and constants comparing the two, Theorem
\ref{khAthm2}(iii) holds. So Theorem \ref{khAthm2} gives $\ti
T_{ab}^m,C_{ab}^m$ and $\si_{abc}^m$ for $c\in C_{ab}^m$. Part (a)
of $(\dag)_{ab}^m$ follows from $\min T\vert_{\pd X}\equiv\min T'$
in Theorem \ref{khAthm2}, part (b) from the boundary compatibility
between the $\si_c$ in Theorem \ref{khAthm1} referred to in Theorem
\ref{khAthm2}, and (c) from the compatibility between the
$\si'_{c'}$ and $\si_c$ in Theorem \ref{khAthm2}. This proves the
inductive step. Hence by induction $(\dag)_{ab}^m$ holds for
all~$m,a,b$.

Now $\ti X_a=\coprod_{b=1}^{n_a^0}\ti X_{ab}^0$, and if
$b=1,\ldots,n_a^0$ and $\phi^0(a,b)=(\bar a,\bar b)$ in $R^0$, we
have a diffeomorphism $\bs a^0_{ab}:\ti X_{ab}^0\ra \ti X_{\bar
a\bar b}^0$ and a tent function $\ti T_{\bar a\bar b}^0$ on $\ti
X_{\bar a\bar b}^0$. Define a tent function $\ti T_a$ on $\ti X_a$
by $\ti T_a\vert_{\ti X_{ab}^0}=\ti T_{\bar a\bar b}^0\ci\bs
a^0_{ab}$ for $b=1,\ldots,n_a^0$, for each $a\in A$. Then we may
form a compact $(k\!+\!1)$-manifold $Z_{\ti X_a,\ti T_a}$ as in
\eq{khAeq4}. Let $\ti X_{ac}$ for $c\in C_a$ be the extra pieces in
the decomposition of $\pd Z_{\ti X_a,\ti T_a}$ in \eq{khAeq5}. Then
from above we may take $C_a=\coprod_{b=1}^{n_a^0}C_{\bar a\bar
b}^0$, where $(\bar a,\bar b)=\phi^0(a,b)$, and the diffeomorphisms
$\si_{\bar a\bar b\bar c}^0:\De_k\ra\ti X_{\bar a\bar b\bar c}^0$
for $b=1,\ldots,n_a^0$, $(\bar a,\bar b)=\phi^0(a,b)$ and $\bar c\in
C_{\bar a\bar b}^0$ induce diffeomorphisms $\si_{ac}:\De_k\ra\ti
X_{ac}$ for all~$c\in C_a$.

Write $\pi:Z_{\ti X_a,\ti T_a}\ra\ti X_a$ for the natural
projection. We will define effective gauge-fixing data
$\ubtH_a=(\bs{\ti J}_a,\bs{\ti\ze}_a,\utH^{k+1}_a)$ for $(Z_{\ti
X_a,\ti T_a},\ti f_a\ci\pi)$. Let $\bs{\ti J}_a$ have indexing set
$\ti J_a=\{k+1\}$, one Kuranishi neighbourhood $(\ti W_a^{k+1},\ti
F_a^{k+1},\ti t_a^{k+1},\ti\chi_a^{k+1})=(Z_{\ti X_a,\ti T_a},
Z_{\ti X_a,\ti T_a},0,\id_{Z_{\ti X_a,\ti T_a}})$, map $\ti
f_a\ci\pi:\ti W_a^{k+1}\ra Y$ representing $\ti f_a\ci\pi$, and let
$\bs{\ti\ze}_a=(\ti\ze_{a,k+1},\ti\ze_{a,k+1}^{k+1})$ with
$\ti\ze_{a,k+1}\equiv 1\equiv\ti\ze_{a,k+1}^{k+1}$. Then $(\bs{\ti
J}_a,\bs{\ti\ze}_a)$ is an excellent coordinate system for $(Z_{\ti
X_a,\ti T_a},\ti f_a\ci\pi)$. Define $\utH^{k+1}_a:Z_{\ti X_a,\ti
T_a}\ra\uP$ by
\e
\begin{split}
&\utH^{k+1}_a(u,x)=\\
&\begin{cases} \utG^k_a(x), & u=0, \\
(u,x_1,\ldots,x_n), & \text{$0\!<\!u\!<\!\min\ti T_a(x)$,
$\utG^k_a(x)\!=\!(x_1,\ldots,x_n)$,} \\
\uG_{\De_k}^k\bigl((y_0,\ldots,y_k)\bigr), & \begin{subarray}{l}\ts
\text{$u\!=\!\min\ti T_a(x)$, $(u,x)\!=\!\io\!\ci\!\si_{ac}
(y_0,\ldots,y_k)$,}\\
\ts\text{$c\in C_a$, $(y_0,\ldots,y_k)\in\De_k$,}\end{subarray}
\end{cases}
\end{split}
\label{khBeq30}
\e
where points of $Z_{\ti X_a,\ti T_a}$ are $(u,x)$ for $x\in\ti X_a$
and $0\le u\le\min\ti T_a(x)$, and $\uG_{\De_k}^k$ is as in
\eq{kh4eq12}, and $\io:\pd Z_{\ti X_a,\ti T_a}\ra Z_{\ti X_a,\ti
T_a}$ is the natural immersion, restricted to~$\ti X_{ac}\subset\pd
Z_{\ti X_a,\ti T_a}$.

Note that if $(u,x)\in Z_{\ti X_a,\ti T_a}$ with $u=\min\ti T_a(x)$
then $(u,x)\in\io(\ti X_{ac})$ for some $c\in C_a$, so that
$(u,x)=\io\ci\si_{ac}(y_0,\ldots,y_k)$ for some
$(y_0,\ldots,y_k)\in\De_k$ as $\si_{ac}:\De_k\ra\ti X_{ac}$ is a
diffeomorphism, and the third line of \eq{khBeq30} makes sense.
However, $c\in C_a$ may not be unique, and we may also have
$(u,x)\in\io(\ti X_{ac'})$ for some other $c'\in C_a$, with
$(u,x)=\io\ci\si_{ac'}(y_0',\ldots,y_k')$ for $(y_0',\ldots,y_k')
\in\De_k$. Thus, to show \eq{khBeq30} is well-defined we must prove
that $\uG_{\De_k}^k\bigl((y_0,\ldots,y_k)\bigr)=\uG_{\De_k}^k\bigl(
(y_0',\ldots,y_k')\bigr)$. Let $(z_1,\ldots,z_l)$ and
$(z'_1,\ldots,z'_{l'})$ be the result of deleting all zero entries
in $(y_0,\ldots,y_k)$ and $(y_0',\ldots,y_k')$ respectively. Then
using $(\dag)_{ab}^m$(b) one can show that $l=l'$ and
$(z_1,\ldots,z_l)=(z'_1,\ldots,z'_{l'})$, and \eq{kh4eq12} then
implies that~$\uG_{\De_k}^k\bigl((y_0,\ldots,y_k)\bigr)=
\uG_{\De_k}^k\bigl((y_0',\ldots,y_k')\bigr)$.

Hence $\utH^{k+1}_a$ is well-defined, and it is easy to check that
$\ubtH_a=(\bs{\ti J}_a,\bs{\ti\ze}_a,\utH^{k+1}_a)$ is {\it
effective gauge-fixing data} for $(Z_{\ti X_a,\ti T_a},\ti
f_a\ci\pi)$. As in \eq{khAeq30} we see that
\e
\begin{split}
\pd&\bigl[Z_{\ti X_a,\ti T_a},\ti f_a\ci\pi,\ubtH_a\bigr]=
\ts\sum_{c\in C_a}\bigl[\ti X_{ac},\ti f_a\ci\pi\vert_{\ti
X_{ac}},\ubtH_a\vert_{\ti X_{ac}}\bigr]\\
&-[\ti X_a,\ti f_a,\ubtG_a]-\bigl[Z_{\pd\ti X_a,\ti T_a\vert_{\pd\ti
X_a}},\ti f_a\vert_{\pd \ti
X_a}\ci\pi,\ubtH_a\vert_{\smash{Z_{\pd\ti X_a,\ti T_a\vert_{\pd\ti
X_a}}}}\bigr].
\end{split}
\label{khBeq31}
\e
For $a\in A$ and $c\in C_a$, define $\ep_{ac}=1$ if the
diffeomorphism $\si_{ac}:\De_k\ra\ti X_{ac}$ is
orientation-preserving, and $\ep_{ac}=-1$ otherwise. Then from
\eq{kh4eq15} and the last line of \eq{khBeq30} we see that
\e
\begin{split}
\bigl[\ti X_{ac},\ti f_a\ci\pi\vert_{\ti X_{ac}},\ubtH_a\vert_{\ti
X_{ac}}\bigr]&=\ep_{ac}\bigl[\De_k,\ti f_a\ci\pi\ci\si_{ac},\bs
f_{ac},\ubG_{\De_k}\bigr]\\
&=\ep_{ac}\Pi_\rsi^\ef(\ti f_a\ci\pi\ci\si_{ac}).
\end{split}
\label{khBeq32}
\e
A similar proof to \eq{khBeq11}--\eq{khBeq12} and \eq{khBeq28} shows
that
\e
\ts\sum_{a\in A}\rho_a\bigl[Z_{\pd\ti X_a,\ti T_a\vert_{\pd\ti
X_a}},\ti f_a\vert_{\pd \ti X_a}\ci\pi,\ubtH_a\vert_{Z_{\pd\ti
X_a,\ti T_a\vert_{\pd\ti X_a}}}\bigr]=0.
\label{khBeq33}
\e
Multiplying \eq{khBeq31} by $\rho_a$, summing over $a\in A$, and
using \eq{khBeq32} and \eq{khBeq33} yields
\e
\begin{split}
\pd\Bigl(&\ts\sum_{a\in A}\rho_a \bigl[Z_{\ti X_a,\ti T_a},\ti
f_a\ci\pi,\ubtH_a\bigr]\Bigr)\\
&=\Pi_\rsi^\ef\bigl(\ts\sum_{a\in A}\sum_{c\in
C_a}(\rho_a\ep_{ac})(\ti f_a\ci\pi\ci\si_{ac})\bigr)-\sum_{a\in
A}\rho_a[\ti X_a,\ti f_a,\ubtG_a].
\end{split}
\label{khBeq34}
\e
Thus $\Pi_\rsi^\ef\bigl(\ts\sum_{a\in A}\sum_{c\in
C_a}(\rho_a\ep_{ac})(\ti f_a\ci\pi\ci\si_{ac})\bigr)$ is homologous
to $\sum_{a\in A}\rho_a[\ti X_a,\ti f_a,\ab\ubtG_a]$, and is a cycle
representing~$\al$.

Now Definition \ref{kh4def2}(i)--(iii) impose no nontrivial
relations upon chains of the form $[\De_k,\si, \ubG_{\De_k}]$ for
smooth $\si:\De_k\ra Y$, and thus $\Pi_\rsi^\ef$ in \eq{kh4eq15} is
injective at the chain level. As $\Pi_\rsi^\ef\ci\pd=\pd\ci
\Pi_\rsi^\ef$ and $\Pi_\rsi^\ef\bigl(\ts\sum_{a\in A}\sum_{c\in
C_a}(\rho_a\ep_{ac})(\ti f_a\ci\pi\ci\si_{ac})\bigr)$ is a cycle in
$KC_k^\ef(Y;R)$, \eq{khBeq2} holds, and $\be=\bigl[\sum_{a\in
A}\sum_{c\in C_a}(\rho_a\ep_{ac}) \,(\ti
f_a\ci\pi\ci\si_{ac})\bigr]$ is well-defined in $H_k^\rsi(Y;R)$,
with $\Pi_\rsi^\ef(\be)=\al$. This completes Step~3.

\subsection{Step 4: $\Pi_\rsi^\ef:H_k^\rsi(Y;R)\ra KH_k^\ef(Y;R)$
is injective}
\label{khB4}

Let $\be$ in $H_k^\rsi(Y;R)$ with $\Pi_\rsi^\ef(\be)=0$, the
singular $k$-chain $\sum_{d\in D}\eta_d\tau_d$, and $\sum_{a\in
A}\rho_a[X_a,\bs f_a,\ubG_a]$ in $KC_{k+1}^\ef(Y;R)$ be as in Step
4, so that \eq{khBeq3} holds. We now explain how to apply Steps 1--3
to this $\sum_{a\in A}\rho_a[X_a,\bs f_a,\ubG_a]$, replacing
$k$-chains by $(k\!+\!1)$-chains, to construct a chain in
$C^\rsi_{k+1}(Y;R)$ with boundary $\sum_{d\in D}\eta_d\tau_d$,
proving that $\be=0$. Write $\pd X_a=X_{a1}^1\amalg\cdots\amalg
X^1_{an_a^1}$ for the decomposition of $[\pd X_a,\bs f_a\vert_{\pd
X_a},\ubG_a\vert_{\pd X_a}]$ into connected pieces given by
Lemma~\ref{kh3lem}.

From \eq{khBeq3} and the proof that $\Pi_\rsi^\ef:C_k^\rsi(Y;R)\ra
KC_k^\ef(Y;R)$ is injective, we see that we can divide the pieces
$X_{ab}^1$ for $a\in A$ and $b=1,\ldots,n_a^1$ into two kinds: (i)
those with $[X_{ab}^1,\bs f_{ab}^1,\bs G_{ab}^1]=\pm\bigl[
\De_k,\tau_d,\ubG_{\De_k}\bigr]$ for some $d\in D$, and (ii) the
rest. Moreover, the sum of $\rho_a[X_{ab}^1,\bs
f_{ab}^1,\ubG_{ab}^1]$ over all $(a,b)$ of type (i) is $\sum_{d\in
D}\eta_d\bigl[\De_k,\tau_p,\ubG_{\De_k}\bigr]$, and the sum over all
$(a,b)$ of type (ii) is zero. Let us rewrite the chain $\sum_{d\in
D}\eta_d\tau_d$ if necessary, allowing repeats among the $\tau_d$,
so that there is a 1-1 correspondence between $d\in D$ and those
$(a,b)$ of type (i), so that if $d$ corresponds to $(a,b)$ then
$[X_{ab}^1,\bs f_{ab}^1,\ubG_{ab}^1]=\pm\bigl[\De_k,\tau_d,\bs
G_{\De_k}\bigr]$ and $\rho_a=\pm\eta_d$. Clearly, this is possible.

If $X_{ab}^1$ is a boundary component of type (i) then as $X_{ab}^1$
is diffeomorphic to the smooth $k$-manifold $\De_k$, it follows that
$X_a$ is near $X_{ab}^1$ a smooth $(k+1)$-manifold. So in Step 1
there is no need to cut $X_a$ into pieces near $X_{ab}^1$ to get
trivial stabilizers, and in Step 2 there is no need to deform $X_a$
near $X_{ab}^1$ to make it a manifold; we can make choices in Steps
1 and 2 so that the type (i) boundary components $X_{ab}^1$ are not
changed at all.

In Step 1 we do this by choosing the tent functions $\bs T_a$ for
$a\in A$ so that $\min\bs T_a\vert_{X_{ab}^1}\equiv 1$ for each type
(i) boundary component $X_{ab}^1$. That is, in $(*)_{ab}^m$ we
impose the extra condition that $T_{\bar a\bar
b}^{m,k-m}\equiv\{1\}$ for any triple $(X_{\bar a\bar b}^m,\bs
f_{\bar a\bar b}^m,\ubG_{\bar a\bar b}^m)$ which is isomorphic to a
component of $\pd^{m-1}(X_{ab}^1,\bs f_{ab}^1,\ubG_{ab}^1)$ for some
type (i) boundary component $X_{ab}^1$ and all $m=k,k-1,\ldots,1$.
Then at the end of Step 1, to each $d\in D$ there corresponds a
unique type (i) boundary component $X_{ab}^1\cong\De_k$ with
$[X_{ab}^1,\bs f_{ab}^1,\ubG_{ab}^1]=\pm\bigl[
\De_k,\tau_d,\ubG_{\De_k}\bigr]$, and $X_{ab}^1$ is naturally
diffeomorphic to a unique component of $\pd X_{ac}$ for some
unique~$c\in C_a$.

In Step 2, on $X_{ac}$ near a component of $\pd X_{ac}$ naturally
diffeomorphic to a type (i) boundary component $X_{ab}^1$, since the
excellent coordinate system $(\bs I_{ab}^1,\bs\eta_{ab}^1)$ for
$(X_{ab}^1,\bs f_{ab}^1)$ in $\ubG_{ab}^1$ has indexing set
$I_{ab}^1=\{k\}$ and Kuranishi neighbourhood
$(V_{ab}^{1,k},E_{ab}^{1,k},s_{ab}^{1,k},
\psi_{ab}^{1,k})=(X_{ab}^1,X_{ab}^1,0,\id_{X_{ab}^1})$, the
excellent coordinate system $(\bs I_{ac},\bs\eta_{ac})$ for
$(X_{ac},\bs f_{ab})$ in $\ubG_{ac}$ includes a Kuranishi
neighbourhood $(V_{ac}^{k+1},E_{ac}^{k+1},s_{ac}^{k+1},
\psi_{ac}^{k+1})$ covering the component of $\pd X_{ac}$
diffeomorphic to~$X_{ab}^1$.

But $E_{ac}^{k+1}\equiv 0$, so $\ti s_{ac}^{k+1}\equiv 0$, the $\ti
s_{ac}^i$ for $i\in I_{ac}$ are zero near this component of $\pd
X_{ac}$. Thus $\ti X_{ac}$ coincides with $X_{ac}$ near this
component of $\pd X_{ac}$. We can also choose the map
$\utG_{ac}^{k+1}:\ti X_{ac}\ra\uP$ to agree with
$\uG_{ab}^1:X_{ab}^1\ra\uP$, or equivalently with
$\uG_{\De_k}^k:\De_k\ra\uP$, on the component of $\pd\ti X_{ac}$
diffeomorphic to $X_{ab}^1\cong\De_k$. Therefore at the end of Step
2, to each $d\in D$ there corresponds a type (i) boundary component
$X_{ab}^1\cong\De_k$ which is naturally diffeomorphic to a unique
component of $\pd\ti X_{ac}$ for unique~$c\in C_a$.

As in Step 3, we change notation from $\ti X_{ac}$ to $\ti X_a$.
Thus, at the beginning of Step 3 we have a chain $\sum_{a\in
A}\rho_a[\ti X_a,\ti f_a,\ab\ubtG_a]$ in $KC_{k+1}^\ef(Y;R)$ with
\e
\ts\pd\bigl(\sum_{a\in A}\rho_a[\ti X_a,\ti f_a,\ubtG_a]\bigr)
=\sum_{d\in D}\eta_d\bigl[\De_k,\tau_d,\ubG_{\De_k}\bigr],
\label{khBeq35}
\e
where for each $d\in D$ we are given $a\in A$ and a component $\ti
X_{ab}^1$ of $\pd(\ti X_a,\ti f_a,\ubtG_a)$ with a diffeomorphism
$\ti X_{ab}^1\cong\De_k$ such that $[\ti X_{ab}^1,\ti f_{ab}^1,\ti
G_{ab}^1]=\pm\bigl[\De_k,\tau_d,\ubG_{\De_k}\bigr]$ and
$\eta_d=\pm\rho_a$. Step 3 then chooses tent functions $\ti T_a$ on
$\ti X_a$ for $a\in A$ and diffeomorphisms $\si_{ac}:\De_{k+1}\ra\ti
X_{ac}$. If $d\in D$ and $\ti X_{ab}^1\cong\De_k$ is the given
diffeomorphism, then $\ti T_a\vert_{\ti X_{ab}^1}$ is a tent
function on $\De_k$, which cuts $\De_k$ up into pieces, each of
which has a diffeomorphism with $\De_k$ of the form
$\pi\ci\si_{ac}\ci F_{k+1}^{k+1}:\De_k\ra\ti X_a$ for some~$c\in
C_a$.

Unlike in Steps 1 and 2, we cannot choose the $\ti T_a$ in Step 3 so
that the type (i) boundary components $X_{ab}^1$, that is, the
chains $\bigl[\De_k,\tau_d,\ubG_{\De_k}\bigr]$ for $d\in D$, are
unchanged, because to apply the results of \S\ref{khA16} during Step
3 we need to cut the $\ti X_a$ into {\it small pieces}, which
includes the boundary components $\ti X_{ab}^1\cong\De_k$. This is
necessary in the proof of the inductive step in $(\dag)_{ab}^m$,
when we have to ensure that the smallness condition Theorem
\ref{khAthm2}(iii) holds, and means that we must cut the $\ti
X_{ab}^1\cong\De_k$ into pieces small compared to the length scales
involved in $\ti X_{ab}^m$ for all $m=k,k-1,\ldots,0$ and~$(a,b)\in
R^m$.

However, we need not cut the boundary components $\ti X_{ab}^1\cong
\De_k$ into small pieces in an arbitrary way. In Definition
\ref{khAdef3} we defined the {\it barycentric subdivision} of a
polyhedron, and showed how to realize it using a tent function. The
barycentric subdivision of $\De_k$ divides it into $(k+1)!$
simplices $\De_k$, each of which is smaller than the original
$\De_k$ by a factor of at most $\frac{k}{k+1}$ in every direction,
for $k\ge 1$. Barycentric subdivision can be iterated, so that for
$N=1,2,\ldots$ the $N^{\rm th}$ barycentric subdivision of $\De_k$
divides it up into $\bigl((k+1)!\bigr){}^N$ smaller simplices
$\De_k$, each of which is smaller than the original $\De_k$ by a
factor of at most $(\frac{k}{k+1})^N$ in every direction. Since
$(\frac{k}{k+1})^N\ra 0$ as $N\ra\iy$, by taking $N$ arbitrarily
large we make the simplices of the $N^{\rm th}$ barycentric
subdivision of $\De_k$ arbitrarily small. The $N^{\rm th}$
barycentric subdivision of $\De_k$ can also be realized using a tent
function, combining the constructions of Definition \ref{khAdef3}
and Theorem~\ref{khAthm1}.

Therefore we choose $N\gg 0$ such that the $N^{\rm th}$ barycentric
subdivision of each boundary component $\ti X_{ab}^1\cong\De_k$
corresponding to some $d\in D$ cuts $\ti X_{ab}^1$, and all its
boundary faces of dimension $k-1,k-2\ldots,1$, into sufficiently
small pieces that the smallness condition Theorem \ref{khAthm2}(iii)
holds for these faces if we use the $N^{\rm th}$ barycentric
subdivision and its associated tent function for $\ti
X_{ab}^1\cong\De_k$, rather than just some tent function chosen
arbitrarily by iteration using Theorem \ref{khAthm2}. Then in
$(\dag)_{ab}^m$ we impose the extra condition that for any triple
$(\ti X_{\bar a\bar b}^m,\ti f_{\bar a\bar b}^m,\ubtG_{\bar a\bar
b}^m)$ which is isomorphic to a component of $\pd^{m-1}(\ti
X_{ab}^1,\ti f_{ab}^1,\ubtG_{ab}^1)$ for some $\ti X_{ab}^1\cong
\De_k$ with $[\ti X_{ab}^1,\ti f_{ab}^1,\ti G_{ab}^1]=\pm
\bigl[\De_k,\tau_d,\ubG_{\De_k}\bigr]$ for some $d\in D$, so that we
have a natural isomorphism $\ti X_{\bar a\bar b}^m\cong\De_{k-m}$,
then $\ti T_{\bar a\bar b}^m$ and $\si_{\bar a\bar b\bar c}^m$ for
$c\in C_{\bar a\bar b}^m$ should come from the $N^{\rm th}$
barycentric subdivision of~$\De_{k-m}$.

Then writing $\{\up_e:e\in E\}$ for the simplices in the $N^{\rm
th}$ barycentric subdivision of $\De_k$, and setting $\ze_e=1$ if
$\up_e$ is orientation-preserving and $\ze_e=-1$ otherwise, the
proof of \eq{khBeq2} in Step 3 generalizes to give \eq{khBeq4}. Now
barycentric subdivision is used by Bredon \cite[\S IV.17]{Bred} to
prove the Excision Axiom for singular homology. Bredon defines
$\Up:C_*^\rsi(Y;R)\ra C_*^\rsi(Y;R)$ which replaces each singular
simplex by its barycentric subdivision, and $T:C_*^\rsi(Y;R)\ra
C_{*+1}^\rsi(Y;R)$ with $\pd\ci T+T\ci\pd=\Up-\id$. Thus $\Up$
induces the identity on $H^\rsi_*(Y;R)$. In this notation,
$\sum_{d\in D}\sum_{e\in E}(\eta_d \ze_e)\,(\tau_d\ci\up_e)=
\Up^N\bigl(\sum_{d\in D}\eta_d\,\tau_d\bigr)$, so
\begin{equation*}
\be=\bigl[\ts\sum_{d\in D}\eta_d\,\tau_d \bigr]=\bigl[\ts\sum_{d\in
D}\sum_{e\in E}(\eta_d \ze_e)\,(\tau_d\ci \up_e)\bigr]=0\quad
\text{in $H_k^\rsi(Y;R)$,}
\end{equation*}
as $\Up$ is the identity on $H^\rsi_*(Y;R)$, and using \eq{khBeq4}.
Therefore $\Pi_\rsi^\ef$ in \eq{kh4eq18} is {\it injective}, as we
want. This completes Step 4, and the proof of Theorem~\ref{kh4thm1}.

\section{Proof that $KH_*(Y;R)\cong KH_*^{\rm ef}(Y;R)$}
\label{khC}

We now prove Theorem \ref{kh4thm2}. Let $Y$ be an orbifold, $R$ a
$\Q$-algebra and $k\in\Z$. We must show that $\Pi_\ef^\Kh:
KH_k^\ef(Y;R)\ra KH_k(Y;R)$ in \eq{kh4eq10} is an isomorphism. In a
similar way to the proof of Theorem \ref{kh4thm1} in Appendix
\ref{khB}, we do this in five steps, which we describe briefly below
and then cover in detail in \S\ref{khC1}--\S\ref{khC5}. Steps 1--4
show that $\Pi_\ef^\Kh$ in \eq{kh4eq10} is surjective, and Step 5
that it is injective.
\medskip

\noindent{\bf Step 1.} Let $k\in\Z$, $\al\in KH_k(Y;R)$, and
$\sum_{a\in A}\rho_a[X_a,\bs f_a,\bs G_a]\in KC_k(Y;R)$ represent
$\al$. The goal of this step is to choose a {\it tent function} $\bs
T_a$ for $(X_a,\bs f_a,\bs G_a)$ for each $a\in A$ to `cut' each
$X_a$ into finitely many {\it arbitrarily small pieces} $X_{ac}$ for
$c\in C_a$, in triples $(X_{ac},\ab\bs f_{ac},\ab\bs G_{ac})$, using
the results of \S\ref{khA35}. This construction will be used for two
different purposes in Steps 2 and 3, and each time we apply it, we
will explain what we mean by `arbitrarily small'.

As in Step 1 of Appendix \ref{khB}, we wish to choose the $\bs T_a$
such that $\sum_{a\in A}\sum_{c\in C_a}\ab\rho_a[X_{ac},\ab\bs
f_{ac},\ab\bs G_{ac}]$ is homologous to $\sum_{a\in A}\rho_a[X_a,\bs
f_a,\bs G_a]$ in $KC_k(Y;R)$. Since $\pd\bigl(\sum_{a\in
A}\rho_a[X_a,\bs f_a,\bs G_a]\bigr)=0$, the $[\pd X_a,\bs
f_a\vert_{\pd X_a},\ab\bs G_a\vert_{\pd X_a}]$ must satisfy
relations in $KC_{k-1}(Y;R)$. The issue is to choose the $\bs T_a$
for $a\in A$ in a way compatible with these boundary relations, so
that $\pd\bigl(\sum_{a\in A}\sum_{c\in C_a}\ab\rho_a[X_{ac},\ab\bs
f_{ac},\ab\bs G_{ac}]\bigr)=0$. In Step 1 of Appendix \ref{khB} we
solved a similar problem by introducing notation $n_a^m,(X_{ab}^m,
\bs f_{ab}^m,\ubG_{ab}^m),\ab P^m,Q^m,R^m,\ldots$, and then choosing
data on $X_{ab}^m$ for $(a,b)\in R^m$ with boundary conditions by
reverse induction on codimension $m=M,M-1,\ldots,0$. Here we adopt a
similar strategy, but the problem is more difficult because of the
extra relation Definition~\ref{kh4def2}(iv).

Define $\widetilde{KC}_k(Y;R)$ to be the $R$-module of finite
$R$-linear combinations of isomorphism classes $[X,\bs f,\bs G]$ as
in Definition \ref{kh4def2} for which $\vdim X=k$, with relations
Definition \ref{kh4def2}(i)--(iii), but {\it not\/} (iv). Write
$\Pi:\widetilde{KC}_k(Y;R)\ab\ra KC_k(Y;R)$ for the natural,
surjective projection $\Pi:\sum_{a\in A}\rho_a[X_a,\ab\bs f_a,\ab\bs
G_a]\mapsto\sum_{a\in A}\rho_a[X_a,\ab\bs f_a,\bs G_a]$, whose
kernel is the $R$-submodule of $\widetilde{KC}_k(Y;R)$ set to zero
by equation \eq{kh4eq7} in Definition \ref{kh4def2}(iv). That is,
$\Ker\Pi$ is the subspace of elements of $\widetilde{KC}_k(Y;R)$ of
the form
\e
\ts\sum_{d\in D}\eta_d\bigl([X_d/\Ga_d,\bs\pi_*(\bs
f_d),\bs\pi_*(\bs G_d)]-\frac{1}{\md{\Ga_d}}\,[X_d,\bs f_d,\bs
G_d]\bigr),
\label{khCeq1}
\e
where $D$ is a finite indexing set, and for each $d\in D$ we have
$\eta_d\in R$ and $X_d,\bs f_d,\bs G_d,\Ga_d$ are as in Definition
\ref{kh4def2}(iv). Using Definition \ref{kh4def2}(i)--(iii) we can
also suppose that each $(X_d,\bs f_d,\bs G_d)$ is connected.

The proofs in Appendix \ref{khB} relied on the fact that triples
$[X,\bs f,\ubG]$ with $(X,\bs f,\ubG)$ connected form a basis for
$KC_k^\ef(Y;R)$ over $R$, modulo only the relation $[-X,\bs
f,\ubG]=-[X,\bs f,\ubG]$; this was how we showed the coefficient
$\{\cdots\}$ on the last line of \eq{khBeq11} is zero. Since
$\widetilde{KC}_k(Y;R)$ is defined using the same relations as
$KC_k^\ef(Y;R)$, it is easy to see that triples $[X,\bs f,\bs G]$
such that $(X,\bs f,\bs G)$ is connected and does not admit an
orientation-reversing isomorphism form a basis for
$\widetilde{KC}_k(Y;R)$ over $R$, modulo only the relation $[-X,\bs
f,\bs G]=-[X,\bs f,\bs G]$. However, we do {\it not\/} have a good
basis for $KC_k(Y;R)$. This makes choosing the $\bs T_a$ for $a\in
A$ compatible with relations in $KC_{k-1}(Y;R)$ more complicated.

Our strategy is to lift from $KC_*(Y;R)$ to $\widetilde{KC}_*(Y;R)$,
where we have a good basis, and deal with the extra relation
Definition \ref{kh4def2}(iv) `by hand'. As $\pd\bigl(\sum_{a\in
A}\rho_a[X_a,\bs f_a,\bs G_a]\bigr)=0$ in $KC_{k-1}(Y;R)$, it lifts
to $\Ker\Pi$ in $\widetilde{KC}_{k-1}(Y;R)$, and so is of the form
\eq{khCeq1}. Thus we may write
\ea
&\ts\sum_{a\in A}\rho_a[\pd X_a,\bs f_a\vert_{\pd X_a},\bs
G_a\vert_{\pd X_a}]=
\label{khCeq2}\\
&\ts\sum_{d\in D}\eta_d\bigl([\ti X_d/\Ga_d,\bs\pi_*(\bs{\ti
f}_d),\bs\pi_*(\bs{\ti G}_d)]\!-\!\frac{1}{\md{\Ga_d}}[\ti
X_d,\bs{\ti f}_d,\bs{\ti G}_d]\bigr) \;\>\text{in
$\widetilde{KC}_{k-1}(Y;R)$,} \nonumber
\ea
with each $(\ti X_d,\bs{\ti f}_d,\bs{\ti G}_d)$ connected.

We shall choose tent functions $\bs T_a$ for $(X_a,\bs f_a,\bs G_a)$
for $a\in A$, and $\bs{\ti T}_d$ for $(\ti X_d,\bs{\ti f}_d,\bs{\ti
G}_d)$ for $d\in D$, satisfying the conditions: \begin{itemize}
\setlength{\itemsep}{0pt}
\setlength{\parsep}{0pt}
\item[(a)] $\bs T_a\vert_{\pd X_a}$ is invariant under $\Aut(\pd
X_a,\bs f_a\vert_{\pd X_a},\bs G_a\vert_{\pd X_a})$ for each $a\in
A$, and $\bs{\ti T}_d$ is invariant under $\Aut(\ti X_d,\bs{\ti
f}_d,\bs{\ti G}_d)$, and $\bs\pi_*(\bs{\ti T}_d)$ is invariant under
$\Aut(\ti X_d/\Ga_d,\bs\pi_*(\bs{\ti f}_d),\bs\pi_*(\bs{\ti G}_d))$
for each~$d\in D$.
\item[(b)] If two connected components of any two of
$(\pd X_a,\bs f_a\vert_{\pd X_a},\bs G_a\vert_{\pd X_a})$ for $a\in
A$ and $(\ti X_d,\bs{\ti f}_d,\bs{\ti G}_d)$, $(\ti X_d/\Ga_d,
\bs\pi_*(\bs{\ti f}_d),\bs\pi_*(\bs{\ti G}_d))$ for $d\in D$ are
isomorphic, the restrictions of $\bs T_a,\bs{\ti T}_d$ to these
components are identified by the isomorphism.
\item[(c)] The $\bs T_a,\bs{\ti T}_d$ cut $(X_a,\bs f_a,\bs G_a),
(\ti X_d,\bs{\ti f}_d,\bs{\ti G}_d)$ into `arbitrarily small
pieces', where we will explain how small the pieces need be whenever
we apply the result.
\end{itemize}
To do this consistently, we must actually work by induction on
decreasing codimension $m$ in a similar way to Appendix \ref{khB},
and choose data on $\pd^m(X_a,\bs f_a,\bs G_a)$ and $\pd^{m-1}(\ti
X_d,\bs{\ti f}_d,\bs{\ti G}_d)$ for $m=M,M-1,\ldots,1,0$ satisfying
analogues of (a)--(c), with boundary values chosen in the previous
inductive step. The details will be explained in~\S\ref{khC1}.

Haven chosen $\bs T_a,\bs{\ti T}_d$ for all $a,d$ satisfying
(a)--(c), equation \eq{khAeq29} of \S\ref{khA3} gives, in both
$KC_k(Y;R)$ and $\widetilde{KC}_k(Y;R)$,
\ea
\pd\bigl[Z_{X_a,\bs T_a},&\bs f_a\ci\bs\pi,\bs H_{X_a,\bs
T_a}\bigr]=-\bigl[Z_{\pd X_a,\bs T_a\vert_{\pd X_a}},\bs
f_a\vert_{\pd X_a}\ci\bs\pi,\bs H_{\pd X_a,\bs T_a\vert_{\pd
X_a}}\bigr]
\nonumber\\
&-[X_a,\bs f_a,\bs G_a]+\ts\sum_{c\in C_a}[X_{ac},\bs f_{ac},\bs
G_{ac}].
\label{khCeq3}
\ea
From equation \eq{khCeq2} and (a),(b) above we deduce that
\ea
&\sum_{a\in A}\rho_a \bigl[Z_{\pd X_a,\bs T_a\vert_{\pd X_a}},\bs
f_a\vert_{\pd X_a}\ci\bs\pi,\bs H_{\pd X_a,\bs T_a\vert_{\pd
X_a}}\bigr]=
\label{khCeq4}\\
&\sum_{d\in D}\begin{aligned}[t]&\eta_d\bigl([Z_{\ti
X_d/\Ga_d,\bs\pi_*(\bs{\ti T}_d)},\bs\pi_*(\bs{\ti
f}_d)\ci\bs\pi,\bs H_{\ti X_d/\Ga_d,\bs\pi_*(\bs{\ti
T}_d)}]\\
&\ts-\!\frac{1}{\md{\Ga_d}}[Z_{\ti X_d,\bs{\ti T}_d},\bs{\ti
f}_d\ci\bs\pi,\bs H_{\ti X_d,\bs{\ti T}_d}]\bigr)
\end{aligned}\quad\text{in
$\widetilde{KC}_k(Y;R)$.} \nonumber
\ea
The action of $\Ga_d$ on $(\ti X_d,\bs{\ti f}_d,\bs{\ti G}_d)$
preserves $\bs{\ti T}_d$ and so lifts to an action of $\Ga_d$ on
$(Z_{\ti X_d,\bs{\ti T}_d},\bs{\ti f}_d\ci\bs\pi,\bs H_{\ti
X_d,\bs{\ti T}_d})$, and the quotient is isomorphic to $(Z_{\ti
X_d/\Ga_d,\bs\pi_*(\bs{\ti T}_d)},\ab\bs\pi_*(\bs{\ti
f}_d)\ci\bs\pi,\bs H_{\ti X_d/\Ga_d,\bs\pi_*(\bs{\ti T}_d)})$. Thus
the right hand side of \eq{khCeq4} lies in $\Ker\Pi$, and projecting
to $KC_k(Y;R)$ yields
\e
\ts\sum_{a\in A}\rho_a\bigl[Z_{\pd X_a,\bs T_a\vert_{\pd X_a}},\bs
f_a\vert_{\pd X_a}\ci\bs\pi,\bs H_{\pd X_a,\bs T_a\vert_{\pd
X_a}}\bigr]=0\;\>\text{in $KC_k(Y;R)$.}
\label{khCeq5}
\e

Multiplying \eq{khCeq3} by $\rho_a$, summing over $a\in A$ and using
\eq{khCeq5} yields
\e
\begin{split}
\pd\Bigl(&\ts\sum_{a\in A}\rho_a\bigl[Z_{X_a,\bs T_a},\bs
f_a\ci\bs\pi_a,\bs H_{X_a,\bs T_a}\bigr]\Bigr)\\
&=\ts\sum_{a\in A}\sum_{c\in C_a}\rho_a[X_{ac},\bs f_{ac}, \bs
G_{ac}]-\sum_{a\in A}\rho_a[X_a,\bs f_a,\bs G_a]
\end{split}
\label{khCeq6}
\e
in $KC_k(Y;R)$, as for \eq{khBeq8}. Thus $\sum_{a\in A}\sum_{c\in
C_a}\rho_a[X_{ac},\bs f_{ac},\bs G_{ac}]$ is a cycle homologous to
$\sum_{a\in A}\rho_a \bigl[X_a,\ab\bs f_a,\ab\bs G_a\bigr]$, as we
want.
\medskip

\noindent{\bf Step 2.} Let $k\in\Z$, $\al\in KH_k(Y;R)$, and
$\sum_{a\in A}\rho_a[X_a,\bs f_a,\bs G_a]\in KC_k(Y;R)$ represent
$\al$. Using Step 1, we choose a tent function $\bs T_a$ for
$(X_a,\bs f_a,\bs G_a)$ for all $a\in A$ which `cuts' each $X_a$
into small pieces $X_{ac}$ for $c\in C_a$, in triples
$(X_{ac},\ab\bs f_{ac},\ab\bs G_{ac})$, such that for each $a\in A$
and $c\in C_a$ there exists a triple $(\acute X_{ac},\ab\bs{\acute
f}_{ac},\ab \bs{\acute G}_{ac})$ with $\acute X_{ac}$ an oriented
Kuranishi space with {\it trivial stabilizers}, and writing
$\bs{\acute G}_{ac}=(\bs{\acute
I}_{ac},\bs{\acute\eta}_{ac},G^i_{ac}:i\in I_{ac})$ with $\bs{\acute
I}_{ac}=\bigl(\acute I_{ac},(\acute V^i_{ac},\ldots,
\acute\psi^i_{ac}):i\in\acute I_{ac},\ldots\bigr)$, each $\acute
V^i_{ac}$ is a {\it manifold}, and there exists a finite group
$\Ga_{ac}$ of automorphisms $(\bs a,\bs b):(\acute
X_{ac},\ab\bs{\acute f}_{ac},\ab \bs{\acute G}_{ac})\ra(\acute
X_{ac},\ab\bs{\acute f}_{ac},\ab \bs{\acute G}_{ac})$ with $\bs a:
\acute X_{ac}\ra\acute X_{ac}$ orientation-preserving, such that the
triple $\bigl(\acute X_{ac}/\Ga_{ac},\bs\pi_*(\bs{\acute f}_{ac}),
\bs\pi_*(\bs{\acute G}_{ac})\bigr)$ is isomorphic to $(X_{ac},\ab\bs
f_{ac},\ab\bs G_{ac})$. Relation Definition \ref{kh4def2}(iv) thus
gives
\e
\bigl[X_{ac},\ab\bs f_{ac},\ab\bs G_{ac}\bigr]=
\frac{1}{\md{\Ga_{ac}}} \,\bigl[\acute X_{ac},\ab\bs{\acute
f}_{ac},\ab \bs{\acute G}_{ac}\bigr]\qquad\text{in $KC_k(Y;R)$.}
\label{khCeq7}
\e
Such $(\acute X_{ac},\ab\bs{\acute f}_{ac},\ab \bs{\acute G}_{ac})$
exist provided the $X_{ac}$ and $V_{ac}^i$ in $\bs I_{ac}$ in $\bs
G_{ac}$ are `sufficiently small', just as in an orbifold $V$, any
sufficiently small open set in $V$ may be written as $U/\Ga$ for $U$
a manifold and $\Ga$ a finite group.

Combining equation \eq{khCeq7} and Step 1 shows that
\e
\ts\sum_{a\in A}\sum_{c\in C_a}\rho_a\md{\Ga_{ac}}^{-1}\ab[\acute
X_{ac},\ab\bs{\acute f}_{ac},\ab\bs{\acute G}_{ac}]=\sum_{a\in
A}\sum_{c\in C_a}\rho_a[X_{ac},\ab\bs f_{ac},\ab\bs G_{ac}]
\label{khCeq8}
\e
is homologous to $\sum_{a\in A}\rho_a[X_a,\bs f_a,\bs G_a]$ in
$KC_k(Y;R)$, and so is a cycle representing $\al$. Thus we can
represent $\al$ by a cycle \eq{khCeq8} in which the Kuranishi spaces
$\acute X_{ac}$ have {\it trivial stabilizers}, and in the Kuranishi
neighbourhoods $(\acute V^i_{ac},\ldots, \acute\psi^i_{ac})$ in
$\bs{\acute G}_{ac}$, the $\acute V^i_{ac}$ are {\it manifolds}.
Hence $\acute X_{ac}$ is an {\it effective} Kuranishi space, and the
$(\acute V^i_{ac},\ldots, \acute\psi^i_{ac})$ are {\it effective}
Kuranishi neighbourhoods. We will need this when we lift to
effective Kuranishi chains in Step~4.
\medskip

\noindent{\bf Step 3.} For simplicity we now {\it change notation\/}
from $\acute X_{ac}$ for $a\in A$ and $c\in C_a$ back to $X_a$ for
$a\in A$, and absorb the factors $\md{\Ga_{ac}}^{-1}$ in \eq{khCeq8}
into new constants $\rho_a\in R$. That is, Steps 1 and 2 have shown
that we can represent $\al$ by a cycle $\sum_{a\in A}\rho_a[X_a,\bs
f_a,\bs G_a]\in KC_k(Y;R)$, where each $X_a$ is a compact, oriented
Kuranishi space with {\it trivial stabilizers}, and writing $\bs
G_a=(\bs I_a,\bs\eta_a,G^i_a:i\in I_a)$ with $\bs I_a=\bigl(I_a,
(V^i_a,\ldots,\psi^i_a):i\in I_a,\ldots\bigr)$, each $V^i_a$ is a
{\it manifold}. Equation \eq{khCeq2} holds for connected $(\ti
X_d,\bs{\ti f}_d,\bs{\ti G}_d)$ for $d\in D$. We stress that these
$A,X_a,\bs f_a,\bs G_a,\rho_a,D,\ti X_d,\bs{\ti f}_d,\bs{\ti
G}_d,\Ga_d,\eta_d$ are different from those in Steps 1 and~2.

We now repeat Step 1, choosing tent functions $\bs T_a$ for
$(X_a,\bs f_a,\bs G_a)$ for $a\in A$, and $\bs{\ti T}_d$ for $(\ti
X_d,\bs{\ti f}_d,\bs{\ti G}_d)$ for $d\in D$, satisfying conditions
(a)--(c), such that $\bs T_a$ cuts $(X_a,\bs f_a,\bs G_a)$ into
finitely many small pieces $(X_{ac},\ab\bs f_{ac},\ab\bs G_{ac})$
for $c\in C_a$, and $\bs{\ti T}_d$ cuts $(\ti X_d,\bs{\ti
f}_d,\bs{\ti G}_d)$ into finitely many small pieces $(\ti
X_{df},\bs{\ti f}_{df},\bs{\ti G}_{df})$ for~$f\in F_d$.

As in Step 1, equations \eq{khCeq2}--\eq{khCeq4} hold in
$\widetilde{KC}_*(Y;R)$. Taking $\pd$ of \eq{khCeq3} and using
$\pd^2=0$ gives
\e
\begin{aligned}
&\pd\bigl[Z_{\pd X_a,\bs T_a\vert_{\pd X_a}},\bs f_a\vert_{\pd
X_a}\ci\bs\pi,\bs H_{\pd X_a,\bs T_a\vert_{\pd X_a}}\bigr]=\\
&-\pd[X_a,\bs f_a,\bs G_a]+\ts\sum_{c\in C_a}\pd[X_{ac},\bs
f_{ac},\bs G_{ac}]
\end{aligned}\quad\text{in $\widetilde{KC}_{k-1}(Y;R)$.}
\label{khCeq9}
\e
The analogue of \eq{khCeq3} for $(\ti X_d,\bs{\ti f}_d,\bs{\ti
G}_d)$ in $\widetilde{KC}_{k-1}(Y;R)$ is
\ea
\pd\bigl[Z_{\ti X_d,\bs{\ti T}_d},&\bs{\ti f}_d\ci\bs\pi,\bs H_{\ti
X_d,\bs{\ti T}_d}\bigr]=-\bigl[Z_{\pd\ti X_d,\bs{\ti
T}_d\vert_{\pd\ti X_d}},\bs{\ti f}_d\vert_{\pd\ti X_d}\ci\bs\pi,\bs
H_{\pd\ti X_d,\bs{\ti T}_d\vert_{\pd\ti X_d}}\bigr]
\nonumber\\
&-[\ti X_d,\bs{\ti f}_d,\bs{\ti G}_d]+\ts\sum_{f\in F_d}[\ti
X_{df},\bs{\ti f}_{df},\bs{\ti G}_{df}].
\label{khCeq10}
\ea
The action of $\Ga_d$ on $(\ti X_d,\bs{\ti f}_d,\bs{\ti G}_d)$ lifts
to $\coprod_{f\in F_d}(\ti X_{df},\bs{\ti f}_{df},\bs{\ti G}_{df})$,
and
\begin{equation*}
\raisebox{-2pt}{$\displaystyle\Bigl($}\coprod_{f\in F_d\!\!\!\!}
(\ti X_{df},\bs{\ti f}_{df},\bs{\ti G}_{df})\!
\raisebox{-2pt}{$\displaystyle\Bigr)\Big/$}\Ga_d\!\cong
\!\!\!\!\coprod_{f\Ga_d\in F_d/\Ga_d}\!\!\!\!\!\! \bigl(\ti
X_{df}/\Stab_{\Ga_d}(f),\bs\pi_*(\bs{\ti f}_{df}),\bs\pi_*(\bs{\ti
G}_{df})\bigr),
\end{equation*}
where to interpret the right hand side, for each orbit $f\Ga_d$ of
$\Ga_d$ in $F_d$ we pick a representative $f$, and then divide $(\ti
X_{df},\bs{\ti f}_{df},\bs{\ti G}_{df})$ by the stabilizer group
$\Stab_{\Ga_d}(f)$ of $f$ in $\Ga_d$. From this we deduce an
analogue of \eq{khCeq10} for $(\ti X_d/\Ga_d,\ab\bs\pi_*(\bs{\ti
f}_d),\ab\bs\pi_*(\bs{\ti G}_d))$ in $\widetilde{KC}_{k-1}(Y;R)$:
\e
\begin{split}
\pd&\bigl[Z_{\ti X_d/\Ga_d,\bs\pi_*(\bs{\ti T}_d)},\bs\pi_*(\bs{\ti
f}_d)\ci\bs\pi,\bs H_{\ti X_d/\Ga_d,\bs\pi_*(\bs{\ti
T}_d)}\bigr]=\\
&-\bigl[Z_{\pd(\ti X_d/\Ga_d),\bs\pi_*(\bs{\ti T}_d)\vert_{\pd(\ti
X_d/\Ga_d)}},\bs\pi_*(\bs{\ti f}_d)\vert_{\pd(\ti
X_d/\Ga_d)}\ci\bs\pi,\\
&\quad\bs H_{\pd(\ti X_d/\Ga_d),\bs\pi_*(\bs{\ti T}_d)\vert_{\pd(\ti
X_d/\Ga_d)}}\bigr]\!-\!\bigl[\ti X_d/\Ga_d,\bs\pi_*(\bs{\ti
f}_d),\bs\pi_*(\bs{\ti G}_d)\bigr]\\
&+\ts\sum\limits_{f\Ga_d\in F_d/\Ga_d}\bigl[\ti
X_{df}/\Stab_{\Ga_d}(f),\bs\pi_*(\bs{\ti f}_{df}),\bs\pi_*(\bs{\ti
G}_{df})\bigr].
\end{split}
\label{khCeq11}
\e

Taking the boundary of \eq{khCeq2} and using $\pd^2=0$ gives
\e
\begin{gathered}
0=\sum_{d\in D}\begin{aligned}[t]&\eta_d\bigl(\bigl[\pd(\ti
X_d/\Ga_d),\bs\pi_*(\bs{\ti f}_d)\vert_{\pd(\ti
X_d/\Ga_d)},\bs\pi_*(\bs{\ti G}_d)\vert_{\pd(\ti X_d/\Ga_d)}\bigr]\\
&\ts-\md{\Ga_d}^{-1}\bigl[\pd\ti X_d,\bs{\ti f}_d\vert_{\pd\ti
X_d},\bs{\ti G}_d\vert_{\pd\ti X_d}\bigr]\bigr)\quad \text{in
$\widetilde{KC}_{k-2}(Y;R)$.}
\end{aligned}
\end{gathered}
\label{khCeq12}
\e
As the choice of tent functions $\bs{\ti T}_d,\bs\pi_*(\bs{\ti
T}_d)$ in Step 1 actually depends only on the isomorphism class of
each connected component of $(\ti X_d,\bs{\ti f}_d,\ab\bs{\ti
G}_d)$, $(\ti X_d/\Ga_d,\ab\bs\pi_*(\bs{\ti f}_d),\ab
\bs\pi_*(\bs{\ti G}_d))$ in $\widetilde{KC}_{k-1}(Y;R)$, and by the
inductive construction this also holds for their boundaries in
$\widetilde{KC}_{k-2}(Y;R)$, we can lift \eq{khCeq12} to an equation
in the $Z_{\ldots}$ chains in $\widetilde{KC}_{k-1}(Y;R)$:
\e
\begin{gathered}
0=\sum_{d\in D}\begin{aligned}[t]&\eta_d\Bigl( \bigl[Z_{\pd(\ti
X_d/\Ga_d),\bs\pi_*(\bs{\ti T}_d)\vert_{\pd(\ti
X_d/\Ga_d)}},\bs\pi_*(\bs{\ti f}_d)\vert_{\pd(\ti
X_d/\Ga_d)}\ci\bs\pi,\\[-4pt]
&\qquad\qquad \bs H_{\pd(\ti X_d/\Ga_d),\bs\pi_*(\bs{\ti
T}_d)\vert_{\pd(\ti X_d/\Ga_d)}}\bigr]\\[-4pt]
&-\md{\Ga_d}^{-1}\bigl[Z_{\pd\ti X_d,\bs{\ti T}_d\vert_{\pd\ti
X_d}},\bs{\ti f}_d\vert_{\pd\ti X_d}\ci\bs\pi,\bs H_{\pd\ti
X_d,\bs{\ti T}_d\vert_{\pd\ti X_d}}\bigr]\Bigr).
\end{aligned}
\end{gathered}
\label{khCeq13}
\e

Take the boundary $\pd$ of \eq{khCeq4} to get an equation in
$\widetilde{KC}_{k-1}(Y;R)$. Use \eq{khCeq9}, \eq{khCeq11} and
\eq{khCeq10} to substitute for the boundaries of the first, second
and third lines of \eq{khCeq4} respectively. Add on \eq{khCeq2},
which cancels three sets of terms, and then add on \eq{khCeq13},
which cancels two more. This yields an equation
in~$\widetilde{KC}_{k-1}(Y;R)$:
\e
\begin{split}
&\ts\pd\bigl(\sum_{a\in A}\sum_{c\in C_a}\rho_a[X_{ac},\bs
f_{ac},\bs
G_{ac}]\bigr)=\\
&\sum_{d\in D}\eta_d\Bigl(\begin{aligned}[t]&\ts\sum_{f\Ga_d\in
F_d/\Ga_d}\bigl[\ti X_{df}/\Stab_{\Ga_d}(f),\bs\pi_*(\bs{\ti
f}_{df}),\bs\pi_*(\bs{\ti
G}_{df})\bigr]\\
&-\md{\Ga_d}^{-1}\ts\sum_{f\in F_d}[\ti X_{df},\bs{\ti
f}_{df},\bs{\ti G}_{df}]\Bigr).
\end{aligned}
\end{split}
\label{khCeq14}
\e

We will prove that provided the $\bs T_a,\bs{\ti T}_d$ are chosen to
make the $X_{ac},\ti X_{df}$ sufficiently small, the terms on the
right hand side of \eq{khCeq14} all cancel, that is, the right hand
side is zero. To get a rough idea of why, note that since the $X_a$
have trivial stabilizers, if $\ti X_d/\Ga_d$ is needed to cancel
part of any $\pd X_a$ in \eq{khCeq2}, then $\ti X_d/\Ga_d$ must have
trivial stabilizers, and so $\Ga_d$ acts freely on $\ti X_d$. But if
$\Ga_d$ acts freely on $\ti X_d$, and the pieces $\ti X_{df}$ for
$f\in F_d$ are sufficiently small, then $\Ga_d$ must also act freely
on $F_d$. And if $\Ga_d$ acts freely on $F_d$ then $\sum_{f\Ga_d\in
F_d/\Ga_d}[\ti X_{df}/\Stab_{\Ga_d}(f),
\ldots]=\md{\Ga_d}^{-1}\sum_{f\in F_d}[\ti X_{df},\ldots]$, so the
second and third lines of \eq{khCeq14} cancel.

We have now constructed $\sum_{a\in A}\sum_{c\in C_a}\rho_a[X_{ac},
\bs f_{ac},\bs G_{ac}]$, where $X_{ac}$ has trivial stabilizers and
the $V_{ac}^i$ are manifolds, and shown that regarded as an element
of $KC_k(Y;R)$ this is a cycle representing $\al$, and regarded as
an element of $\widetilde{KC}_k(Y;R)$ satisfies
\e
\ts\pd\bigl(\sum_{a\in A}\sum_{c\in C_a}\rho_a[X_{ac},\bs f_{ac},\bs
G_{ac}]\bigr)=0\quad\text{in $\widetilde{KC}_{k-1}(Y;R)$.}
\label{khCeq15}
\e
Thus, we have lifted from a cycle in $\bigl(KC_*(Y;R),\pd\bigr)$ to
one in $\bigl(\widetilde{KC}_*(Y;R),\pd\bigr)$. For the rest of the
proof we will work in the complex $\bigl(\widetilde{KC}_*(Y;R),
\pd\bigr)$, which as in Step 1 has an easily understood basis, and
we no longer need to worry about the relation
Definition~\ref{kh4def2}(iv).

In this step we also prove we can take the $(X_{ac},\bs f_{ac},\bs
G_{ac})$ to have an additional property that we will need later. Fix
$a\in A$ and $c\in C_a$, and for $m\ge 0$ write $(X_{acd}^m,\bs
f_{acd}^m,\bs G_{acd}^m)$ for $d=1,\ldots,n_{ac}^m$ for the
connected components of $(\pd^mX_{ac},\bs f_{ac}
\vert_{\pd^mX_{ac}},\bs G_{ac}\vert_{\pd^mX_{ac}})$. Then in the
usual notation the Kuranishi neighbourhoods $(V_{acd}^{m,i-m},
\ldots,\psi_{acd}^{m,i-m})$ for $i-m\in I_{ac}^m$ in $\bs G_{acd}^m$
have $V_{acd}^{m,i-m}$ an open subset of $\pd^mV_{ac}^i$, where
$V_{ac}^i$ is a {\it manifold}.

If $V$ is a manifold with corners and $m\ge 2$ then points of
$\pd^mV$ may be written as $(v,B_1,\ldots,B_m)$, where
$B_1,\ldots,B_m$ are distinct local boundary components of $V$.
Hence $(v,B_1),\ldots,(v,B_m)$ lie in $\pd V$. Thus, a point of
$V_{acd}^{m,i-m}$ is of the form $(v,B_1,\ldots,B_m)$, where
$(v,B_l)\in\pd V_{ac}^i$ for $l=1,\ldots,m$. It is necessary that
$(v,B_l)\in V_{acd_l}^{1,i-1}$ for some unique $d_l=1,\ldots,
n_{ac}^1$. That is, $(v,B_l)$ cannot lie in the part of $\pd
V_{ac}^i$ discarded by Algorithm \ref{kh3alg}, since then
$(v,B_1,\ldots,B_m)$ would also be discarded, so $(v,B_l)$ lies in
$\coprod_{d'=1}^{n_{ac}^1}V_{acd'}^{1,i-1}$, and so in
$V_{acd_l}^{1,i-1}$ for unique~$d_l$.

Hence, each $i\in I_{acd}^m$ and $(v,B_1,\ldots,B_m)$ in
$V_{acd}^{m,i-m}$ determine an $l$-tuple $(d_1,\ldots,d_m)$ with
$d_l=1,\ldots,n_{ac}^1$. We claim that this $(d_1,\ldots,d_m)$ is
independent of $i$ and $(v,B_1,\ldots,B_m)$, and depends only on
$d=1,\ldots,n_{ac}^m$. To see this, note that $(d_1,\ldots,d_m)$ is
locally constant on $V_{acd}^{m,i-m}$, and is also preserved under
coordinate changes $(\phi_{acd}^{m,(i-m)(j-m)},
\hat\phi_{acd}^{m,(i-m)(j-m)})$ in $\bs G_{acd}^m$, so the part of
$(\pd^mX_{ac},\bs f_{ac} \vert_{\pd^mX_{ac}},\bs G_{ac}
\vert_{\pd^mX_{ac}})$ with fixed $(d_1,\ldots,d_m)$ must be a union
of connected components. Thus, each component of $(\pd^mX_{ac},\bs
f_{ac} \vert_{\pd^mX_{ac}},\bs G_{ac}\vert_{\pd^mX_{ac}})$
determines an ordered $m$-tuple of components of~$(\pd X_{ac},\bs
f_{ac}\vert_{\pd X_{ac}},\bs G_{ac}\vert_{\pd X_{ac}})$.

We will prove that provided the $X_{ac}$ are chosen sufficiently
small, for all $a\in A$, $c\in C_a$, $m\ge 2$ and
$d=1,\ldots,n_{ac}^m$, the numbers $d_1,\ldots,d_m$ above are {\it
distinct}. That is, we show that each component of $(\pd^mX_{ac},\bs
f_{ac}\vert_{\pd^mX_{ac}},\bs G_{ac}\vert_{\pd^mX_{ac}})$ occurs as
the intersection of $m$ {\it distinct\/} components of~$(\pd X_{ac},
\bs f_{ac}\vert_{\pd X_{ac}},\bs G_{ac}\vert_{\pd X_{ac}})$.

Here is why we need this. Suppose $(X,\bs f,\ubG)$ is a triple in
effective Kuranishi homology. Then there is a natural action of the
symmetric group $S_m$ on $(\pd^mX,\bs f\vert_{\pd^mX},\ubG
\vert_{\pd^mX})$ by isomorphisms, which acts on points
$(v,B_1,\ldots,B_m)$ in $V^{m,i-m}$ in $\ubG\vert_{\pd^mX}$ by
permutation of $B_1,\ldots,B_m$. As connected components of
$(\pd^mX,\bs f\vert_{\pd^mX},\ubG\vert_{\pd^mX})$ have no
automorphisms by Theorem \ref{kh3thm5}(b), $S_m$ must act {\it
freely\/} on the set of connected components of~$(\pd^mX,\bs
f\vert_{\pd^mX},\ubG\vert_{\pd^mX})$.

In Step 4 we will modify $\bs G_{ac}$ to effective gauge-fixing
data. A necessary condition to be able to do this is that $S_m$ acts
freely on the connected components of $(\pd^mX_{ac},\bs f_{ac}
\vert_{\pd^mX_{ac}},\bs G_{ac}\vert_{\pd^mX_{ac}})$. Since each
component determines some $(d_1,\ldots,d_m)$ and $S_m$ permutes
$d_1,\ldots,d_m$, if $d_1,\ldots,d_m$ are {\it distinct\/} then the
action of $S_m$ on the set of components is automatically free, as
we want.
\medskip

\noindent{\bf Step 4.} For simplicity we again {\it change
notation\/} from $X_{ac}$ for $a\in A$ and $c\in C_a$ back to $X_a$
for $a\in A$. That is, Step 3 has shown that we can represent $\al$
by a cycle $\sum_{a\in A}\rho_a[X_a,\bs f_a,\bs G_a]$ in $KC_k(Y;R)$
or $\widetilde{KC}_k(Y;R)$, where each $X_a$ is a compact, oriented
Kuranishi space with {\it trivial stabilizers}, and writing $\bs
G_a=(\bs I_a,\bs\eta_a,G^i_a:i\in I_a)$ with $\bs I_a=\bigl(I_a,
(V^i_a,\ldots,\psi^i_a):i\in I_a,\ldots\bigr)$, each $V^i_a$ is a
{\it manifold}, each $(X_a,\bs f_a,\bs G_a)$ is {\it connected}, and
\e
\ts\pd\bigl(\sum_{a\in A}\rho_a[X_a,\bs f_a,\bs G_a]\bigr)=0
\quad\text{in $\widetilde{KC}_{k-1}(Y;R)$.}
\label{khCeq16}
\e
Moreover, for $a\in A$ and $m\ge 2$ each component of $(\pd^mX_a,\bs
f_a\vert_{\pd^mX_a},\bs G_a\vert_{\pd^mX_a})$ occurs as the
intersection of $m$ {\it distinct\/} components of~$(\pd X_a, \bs
f_a\vert_{\pd X_a},\bs G_a\vert_{\pd X_a})$.

Define notation $n_a^m,\kern-0.35pt(X_{ab}^m,\bs f_{ab}^m,\bs
G_{ab}^m),P^m,\sim,Q^m,R^m,\phi^m,\kern-0.35pt(\bs a,\bs
b)_{ab}^m,\ab M,\ab\ep_{ab}^m,\ab B_{ab}^m$ exactly as in Step 1 of
Appendix \ref{khB}, but with $\bs G_a,\bs G_{ab}^m$ in place of
$\ubG_a,\ubG_{ab}^m$. Note one significant difference: in Appendix
\ref{khB} we have $\Aut(X_{ab}^m,\bs f_{ab}^m,\ubG_{ab}^m)=\{1\}$ by
Theorem \ref{kh3thm5}(b), so $(\bs a,\bs b)_{ab}^m$ is unique.
However, in the noneffective case we know only that
$\Aut(X_{ab}^m,\bs f_{ab}^m,\bs G_{ab}^m)$ is finite by Theorem
\ref{kh3thm2}, so when $\phi^m(a,b)=(\bar a,\bar b)$ the isomorphism
$(\bs a,\bs b)_{ab}^m:(X_{ab}^m,\bs f_{ab}^m,\bs
G_{ab}^m)\ra(X_{\bar a\bar b}^m,\bs f_{\bar a\bar b}^m,\bs G_{\bar
a\bar b}^m)$ may not be unique, we have to make an arbitrary choice
of $(\bs a,\bs b)_{ab}^m$ out of finitely many possibilities.

Then as for equation \eq{khBeq11}, from \eq{khCeq16} in
$\widetilde{KC}_{k-1}(Y;R)$ we have
\ea
0&=\pd\Bigl[\sum_{a\in A}\rho_a[X_a,\bs f_a,\bs G_a]\Bigr]
=\sum\limits_{a\in
A}\rho_a\sum\limits_{b=1}^{n_a^1}\bigl[X_{ab}^1,\bs f_{ab}^1, \bs
G_{ab}^1\bigr]
\label{khCeq17}\\
&=\sum\limits_{\begin{subarray}{l}\text{$a\in A$,
$b=1,\ldots,n_a^1$.}\\ \text{Set $(\bar a,\bar
b)=\phi^1(a,b)$}\end{subarray}\!\!\!\!\!\!\!\!\!\!\!\!\!\!\!\!
\!\!\!\!\!\!\!\!\!\!\!\!\!\!\!\!\!\!\!\!\!\!\!\!\!\!\!\!\!\! }
\rho_a\ep_{ab}^1\bigl[X_{\bar a\bar b}^1,\bs f_{\bar a\bar b}^1,\bs
G_{\bar a\bar b}^1\bigr]=\!\!\!\sum\limits_{(\bar a,\bar b)\in R^1}
\raisebox{-9pt}{\begin{Large}$\displaystyle\biggl\{$\end{Large}}
\sum\limits_{\begin{subarray}{l}\text{$a\in A$,
$b=1,\ldots,n_a^1$:}\\ \phi^1(a,b)=(\bar a,\bar
b)\end{subarray}\!\!\!\!\!\!\!\!\!\!\!\!\!\!\!\!}\rho_a\ep_{ab}^1
\raisebox{-9pt}{\begin{Large}$\displaystyle\biggr\}$\end{Large}}
\bigl[X_{\bar a\bar b}^1,\bs f_{\bar a\bar b}^1,\bs G_{\bar a\bar
b}^1\bigr].
\nonumber
\ea
By construction the triples $\bigl(X_{\bar a\bar b}^1,\bs f_{\bar
a\bar b}^1,\bs G_{\bar a\bar b}^1\bigr)$ for $(\bar a,\bar b)\in
R^1$ are connected and mutually non-isomorphic, even through
orientation-reversing isomorphisms. As in Step 1, triples $[X,\bs
f,\bs G]$ such that $(X,\bs f,\bs G)$ is connected and does not
admit an orientation-reversing automorphism form a basis for
$\widetilde{KC}_{k-1}(Y;R)$ over $R$, modulo~$[-X,\bs f,\bs
G]=-[X,\bs f,\bs G]$.

We separate $(\bar a,\bar b)\in R^1$ into two cases: (A) $(X_{\bar
a\bar b}^1,\bs f_{\bar a\bar b}^1,\bs G_{\bar a\bar b}^1)$ admits an
orientation-reversing automorphism, and (B) it does not. In case (A)
we have $[X_{\bar a\bar b}^1,\bs f_{\bar a\bar b}^1,\bs G_{\bar
a\bar b}^1]=0$ by Definition \ref{kh4def2}(ii). The $[X_{\bar a\bar
b}^1,\bs f_{\bar a\bar b}^1,\bs G_{\bar a\bar b}^1]$ in case (B)
form part of a basis for $\widetilde{KC}_{k-1}(Y;R)$ over $R$, and
so are linearly independent over $R$. Thus, in the same way that we
argued in \S\ref{khB1} that the coefficients $\{\cdots\}$ in the
last line of \eq{khBeq11} are zero, from \eq{khCeq17} we see that
for each $(\bar a,\bar b)\in R^1$, either
\begin{itemize}
\setlength{\itemsep}{0pt}
\setlength{\parsep}{0pt}
\item[(A)] $(X_{\bar a\bar b}^1,\bs f_{\bar a\bar b}^1,\bs G_{\bar
a\bar b}^1)$ admits an orientation-reversing automorphism, so that
$[X_{\bar a\bar b}^1,\bs f_{\bar a\bar b}^1,\bs G_{\bar a\bar
b}^1]=0$ in $KC_{k-1}(Y;R)$ and $\widetilde{KC}_{k-1}(Y;R)$; or
\item[(B)] $\sum_{a\in A,\; b=1,\ldots,n_a^1:
\phi^1(a,b)=(\bar a,\bar b)}\rho_a\ep_{ab}^1=0$ in~$R$.
\end{itemize}

The rough idea of this step is to choose maps $\uG_a^i:E^i_a\ra\uP$
for all $a\in A$ and $i\in I_a$ such that $\ubG_a=(\bs I_a,
\bs{\eta}_a,\uG^i_a:i\in I_a)$ is effective gauge-fixing data for
$(X_a,\bs f_a)$. (Note that as $X_a$ has trivial stabilizers, it is
an {\it effective} Kuranishi space, and as $V^i_a$ is a manifold,
$(V^i_a,\ldots,\psi^i_a)$ is an {\it effective} Kuranishi
neighbourhood, so we have already satisfied some of the conditions
in Definition \ref{kh3def17} for $\ubG_a$.) Thus we can form a chain
$\sum_{a\in A}\rho_a[X_a,\bs f_a,\ubG_a]\in KC_k^\ef(Y;R)$, and
project to $\sum_{a\in A}\rho_a[X_a,\bs
f_a,\ab\Pi(\ubG_a)]=\Pi_\ef^\Kh \bigl(\sum_{a\in A}\rho_a[X_a,\bs
f_a,\ubG_a]\bigr)$ in~$KC_k(Y;R)$.

If we could choose the $\uG_a^i$ for all $a,i$ compatible with the
relations between the $[\pd X_a,\bs f_a\vert_{\pd X_a},\bs G_a
\vert_{\pd X_a}]$ implied by \eq{khCeq16}--\eq{khCeq17}, we would
expect to prove that $\pd\bigl(\sum_{a\in A}\rho_a[X_a,\bs
f_a,\ubG_a]\bigr) =0$ in $KC_{k-1}^\ef(Y;R)$, giving a homology
class $\be=\bigl[\sum_{a\in A}\rho_a[X_a,\bs f_a,\ubG_a]\bigr]$ in
$KH_k^\ef(Y;R)$, and to show $\sum_{a\in A}\rho_a[X_a,\bs f_a,
\Pi(\ubG_a)]$ is homologous to $\sum_{a\in A}\rho_a[X_a,\bs f_a,\bs
G_a]$ in $KC_k(Y;R)$, so that $\Pi_\ef^\Kh(\be)=\al$. As $\al\in
KH_k(Y;R)$ was arbitrary, this would show that $\Pi_\ef^\Kh:
KH_k^\ef(Y;R)\ra KH_k(Y;R)$ is {\it surjective}.

However, there is a problem. In the $\Aut(X_a,\bs f_a,\bs G_a)$ are
finite by Theorem \ref{kh3thm2}, but for connected $(X_a,\bs
f_a,\ubG_a)$ the $\Aut(X_a,\bs f_a,\ubG_a)$ are trivial by Theorem
\ref{kh3thm5}(b). Thus, when we lift from $\bs G_a$ to $\ubG_a$ we
must {\it break the symmetries\/} of $(X_a,\bs f_a,\bs G_a)$. But
breaking these symmetries is unnatural, and may cause
$\pd\bigl(\sum_{a\in A}\rho_a[X_a,\ab\bs f_a,\ubG_a]\bigr)$ to be
nonzero. To see this, note that \eq{khCeq17} implicitly uses an
isomorphism
\e
\bigl(\pd X_a,\bs f_a\vert_{\pd X_a},\bs G_a\vert_{\pd
X_a}\bigr)\cong\coprod_{\!\!\!\! \text{$b=1,\ldots,n_a^1$. Set
$(\bar a,\bar b)=\phi^1(a,b)$}
\!\!\!\!\!\!\!\!\!\!\!\!\!\!\!\!\!\!\!\!\!\!\!\!\!\!\!\!\!\!\!\!\!\!
\!\!\!\!\!\!\!\!\!\!\!\!\!\!\!} \bigl( X_{\bar a\bar b}^1,\bs{
f}_{\bar a\bar b}^1,\bs G_{\bar a\bar b}^1\bigr).
\label{khCeq18}
\e

The set of isomorphisms \eq{khCeq18} has a free, transitive action
of $\prod_b\Aut\bigl(X_{\bar a\bar b}^1,\ab\smash{\bs f_{\bar a\bar
b}^1,\bs G_{\bar a\bar b}^1\bigr)}$, and in general there is no
natural choice. To choose $\ubG_a$ for $a\in A$ satisfying boundary
relations, we expect (as part of an induction on reverse
codimension) to choose $\ubG_{\bar a\bar b}^1$ for all $(\bar a,\bar
b)\in R^1$, and to use these and \eq{khCeq18} to determine boundary
values $\ubG_a\vert_{\pd X_a}$ for $\ubG_a$, and to extend these
over $ X_a$ to get $\ubG_a$. But to do this we must make an
arbitrary choice of isomorphism \eq{khCeq18}, and similarly in
previous inductive steps. There is no reason for these arbitrary
choices to be compatible with involutions $\bs\si:\pd^2X_a\ra\pd^2
X_a$. Thus the prescribed values for $\ubG_a\vert_{\pd X_a}$, when
restricted to $\pd^2 X_a$, may not be $\bs\si$-invariant, and so by
Principle \ref{kh2pri}, we cannot extend them over~$\pd X_a$.

To get round this, we choose not one set of effective gauge-fixing
data $\ubG_a$ for $(X_a,\bs f_a)$ for each $a\in A$, but a finite
family $\ubG_{ao}$ for $o\in O_a$, such that $\Aut(X_a,\bs f_a,\bs
G_a)$ acts freely on $O_a$, and if $(\bs a,\bs b)\in\Aut(X_a,\bs
f_a,\bs G_a)$ and $d\in O_a$ with $o'=(\bs a,\bs b)\cdot o$ then
$(\bs a,\bs b)$ lifts to an isomorphism $(\bs a',\bs b'):(X_a,\bs
f_a,\ubG_{ao})\ra(X_a,\bs f_a,\ubG_{ao'})$. Similarly, for all
$m=M,M-1,\ldots,0$ and $(a,b)\in P^m$ we will choose a finite family
$\ubG_{abo}^m$ for $o\in O_{ab}^m$ of effective gauge-fixing data
for $(X_{ab}^m,\bs f_{ab}^m)$, where $\Aut(X_{ab}^m,\bs f_{ab}^m,\bs
G_{ab}^m)$ acts freely on $O_{ab}^m$, and lifts to isomorphisms
between the $(X_{ab}^m,\bs f_{ab}^m,\ubG_{abo}^m)$, and
$\ubG_{abo}^m\vert_{\pd X_{ab}^m}$ is determined in terms of
the~$\ubG_{\bar a\bar b\bar o}^{m+1}$.

The construction involves a complex and difficult double induction
on decreasing $m$. The underlying idea is that the families
$(X_{ab}^m,\bs f_{ab}^m,\ubG_{abo}^m)$ for $o\in O_{ab}^m$ have all
the symmetries of $(X_a,\bs f_a,\bs G_a)$, so the symmetries are not
broken, and we need not make arbitrary choices of isomorphisms
\eq{khCeq18}. Using \eq{khCeq17} and alternatives (A),(B) above we
prove that
\e
\ts\pd\bigl(\sum_{a\in A}\sum_{o\in O_a}\rho_a\md{O_a}^{-1}[X_a,\bs
f_a,\ubG_{ao}]\bigr)=0\quad\text{in $KC_{k-1}^\ef(Y;R)$.}
\label{khCeq19}
\e
We then set $\be=\bigl[\sum_{a\in A}\sum_{o\in O_a}\rho_a
\md{O_a}^{-1}[X_a,\bs f_a,\ubG_{ao}]\bigr]$ in $KH_k^\ef(Y;R)$, and
show that $\Pi_\ef^\Kh(\be)=\al$, so $\Pi_\ef^\Kh$ is {\it
surjective}, as we want.
\medskip

\noindent{\bf Step 5.} Suppose $\be\in KH_k^\ef(Y;R)$ with
$\Pi_\ef^\Kh(\be)=0$. We choose a cycle in $KC_k^\ef(Y;R)$
representing $\be$. As $KH_k^\ef(Y;R)\cong H^\rsi_k(Y;R)$ by Theorem
\ref{kh4thm1} we can take this to be the image of a cycle
$\sum_{s\in S}\ze_s\,\tau_s$ in $C^\rsi_k(Y;R)$, for $S$ a finite
indexing set, $\ze_s\in R$ and $\tau_s:\De_k\ra Y$ smooth. Thus
$\be$ is represented by $\Pi_\rsi^\ef(\sum_{s\in S}\ze_s\,\tau_s)$
in $KC_k^\ef(Y;R)$. The image in $KC_k(Y;R)$ under $\Pi_\ef^\Kh$ is
$\Pi_\rsi^\Kh(\sum_{s\in S}\ze_s\,\tau_s)$, which represents
$\Pi_\ef^\Kh(\be)=0$, and so is exact. Therefore as for \eq{khBeq3}
there exists $\sum_{a\in A}\rho_a[X_a,\bs f_a,\bs G_a]\in
KC_{k+1}(Y;R)$ with
\e
\ts\pd\bigl(\sum_{a\in A}\rho_a[X_a,\bs f_a,\bs G_a]\bigr)=
\Pi_\rsi^\Kh(\sum_{s\in S}\ze_s\,\tau_s)=\sum_{s\in
S}\ze_s\bigl[\De_k,\tau_s,\bs G_{\De_k}\bigr].
\label{khCeq20}
\e

We now apply Steps 1--4 to this $\sum_{a\in A}\rho_a[X_a,\bs f_a,\bs
G_a]$, replacing $k$-cycles by $(k+1)$-chains. In Step 1 we use tent
functions $\bs T_a$ to `cut' each $(X_a,\bs f_a,\bs G_a)$ into small
pieces $(X_{ac},\bs f_{ac},\bs G_{ac})$ for $c\in C_a$. When this
process is applied to \eq{khCeq20}, the $(\De_k,\tau_s,\bs
G_{\De_k})$ must also be cut into small pieces. To keep track of
these pieces, it is convenient to use the method of Step 4 of
\S\ref{khB} in \S\ref{khB4}: we choose the $\bs T_a$ such that each
$(\De_k,\tau_s,\bs G_{\De_k})$ is cut into the $N^{\it th}$ {\it
barycentric subdivision\/} of $\De_k$, for some $N\gg 0$ independent
of $s\in S$.

After slightly modifying the definitions of the functions $G^i_{ac}$
in $\bs G_{ac}$, in Step 1 we construct a new chain $\sum_{a\in
A}\sum_{c\in C_a}\rho_a[X_{ac},\bs f_{ac},\bs G_{ac}]$ with all
$X_{ac}$ `arbitrarily small', satisfying
\e
\ts\pd\bigl(\sum_{a\in A}\sum_{c\in C_a}\rho_a[X_{ac},\bs f_{ac},\bs
G_{ac}]\bigr)=\Pi_\rsi^\Kh\ci\Up^N(\sum_{s\in S}\ze_s\,\tau_s)
\label{khCeq21}
\e
for some $N\gg 0$, where $\Up:C^\rsi_k(Y;R)\ra C^\rsi_k(Y;R)$ is the
barycentric subdivision operator of \S\ref{khB4}, so that
$\Up^N=\Up\ci\cdots\ci\Up$ is the $N^{\rm th}$ barycentric
subdivision operator. Note that we will not need the second part of
Step 1, which constructs a homology $\sum_{a\in A}\rho_a\bigl[
Z_{X_a,\bs T_a},\bs f_a\ci\bs\pi_a,\bs H_{X_a,\bs T_a}\bigr]$
between $\sum_{a\in A}\rho_a[X_a,\bs f_a,\bs G_a]$ and $\sum_{a\in
A}\sum_{c\in C_a}\rho_a[X_{ac},\bs f_{ac},\bs G_{ac}]$, we only have
to prove \eq{khCeq21}.

Going through Steps 2--4 with minor modifications then yields a
effective Kuranishi chain $\sum_{a\in A}\sum_{o\in
O_a}\rho_a\md{O_a}^{-1}\ab[ X_a,\bs f_a,\ubG_{ao}]$ in
$KC_{k+1}^\ef(Y;R)$ satisfying
\e
\pd\raisebox{-2pt}{$\displaystyle\Bigl($}\ts\sum\limits_{a\in
A}\sum\limits_{o\in O_a}\rho_a\md{O_a}^{-1}[X_a,\bs f_a,
\ubG_{ao}]\raisebox{-2pt}{$\displaystyle\Bigr)$}=
\Pi_\rsi^\ef\ci\Up^N(\sum_{s\in S}\ze_s\,\tau_s)
\label{khCeq22}
\e
for some $N\gg 0$, where we change notation so that $A,\rho_a,X_a,
\bs f_a,\bs G_a,N$ are different from the first versions. Since
$\sum_{s\in S}\ze_s\,\tau_s$ is a cycle in $C^\rsi_k(Y;R)$, as in
\S\ref{khB4}, $\Up^N(\sum_{s\in S}\ze_s\,\tau_s)$ is homologous to
$\sum_{s\in S}\ze_s\,\tau_s$, and thus $\Pi_\rsi^\ef\ci\Up^N
(\sum_{s\in S}\ze_s\,\tau_s)$ is also a cycle in $KH_k^\ef(Y;R)$
representing $\be$. Therefore \eq{khCeq22} implies that $\be=0$, so
$\Pi_\ef^\Kh:KH_k^\ef(Y;R)\ra KH_k(Y;R)$ is {\it injective}.
\medskip

\begin{rem}{\bf(a)} We used relation Definition \ref{kh4def2}(iv)
in Step 2 to lift to Kuranishi chains with trivial stabilizers, as
in \eq{khCeq7}. But also, much of the work in Steps 1--3 is coping
with problems caused by Definition \ref{kh4def2}(iv), and the
equivalences between chains that it introduces.

We said in Remark \ref{kh4rem1}(b) that Kuranishi homology
$KH_*(Y;R)$ would still work, and Theorem \ref{kh4thm2} would still
hold, if relation Definition \ref{kh4def2}(iv) were omitted in the
definition of $KC_*(Y;R)$, that is, if we replaced $KC_*(Y;R)$ by
$\widetilde{KC}_*(Y;R)$ in the notation of Step 1 above. We now
briefly sketch how to prove Theorem \ref{kh4thm2} with this modified
definition. This time there are four steps:
\begin{list}{}{\setlength{\leftmargin}{35pt}
\setlength{\labelwidth}{35pt}}
\item[{\bf Step 1.}] Let $\al\in KH_*(Y;R)$. Represent $\al$ by a
cycle $\sum_{a\in A}\rho_a[X_a,\bs f_a,\bs G_a]$. Show that $\al$
can be represented by a cycle $\sum_{a\in A}\rho_a[X_a,\bs
f_a,\bs{\ti G_a}]$, in which we fix $X_a,\bs f_a$ but change to new
gauge-fixing data $\bs{\ti G}_a=(\bs{\ti I}_a,\bs{\ti\eta}_a,\ab\ti
G^i_a:i\in\ti I_a)$ with $\bs{\ti I}_a=\bigl(\ti I_a,(\ti
V^i_a,\ldots,\ti\psi^i_a),\ti f^i_a:i\in\ti I_a,\ldots\bigr)$ in
which each $\ti V^i_a$ is an {\it effective} orbifold, although we
do {\it not\/} require $(\ti V^i_a,\ldots,\ti\psi^i_a)$ to be an
effective Kuranishi neighbourhood, that is, the stabilizer groups of
$\ti V^i_a$ need not act trivially on the fibres of~$\ti E^i_a$.

\item[{\bf Step 2.}] As in \S\ref{khA2}, we can use a tent function
to cut an effective orbifold up into manifolds. Following Step 1 of
the proof of Theorem \ref{kh4thm1} in Appendix \ref{khB}, we use the
results of \S\ref{khA35} to choose tent functions $\bs T_a$ for
$(X_a,\bs f_a,\bs{\ti G}_a)$ for $a\in A$ which cut the $\ti V_a^i$
up into manifolds $\ti V_{ac}^i$, and the $X_a$ up into Kuranishi
spaces $X_{ac}$ with trivial stabilizers, for $c\in C_a$. Thus we
represent $\al$ by a cycle $\sum_{a\in A}\sum_{c\in
C_a}\rho_a[X_{ac},\ab\bs f_{ac},\ab\bs{\ti G}_{ac}]$ in which the
$X_{ac}$ have trivial stabilizers, and the $\ti V_{ac}^i$ in
$\bs{\ti G}_{ac}$ are manifolds.

\item[{\bf Steps 3 and 4.}] At the end of the new Step 2, we are in
exactly the same position as at the end of the old Step 3 above,
noting that the new $KC_*(Y;R)$ is the old $\widetilde{KC}_*(Y;R)$.
We can now follow the old Steps 4 and 5 to complete the proof.
\end{list}

\noindent{\bf(b)} In the proof of Theorem \ref{kh4thm1} we use the
fact that $R$ is a $\Q$-algebra in two different ways: firstly, in
\eq{khCeq7} of Step 2 the $\md{\Ga_{ac}}^{-1}$ must lie in $R$, and
secondly, in \eq{khCeq19} of Step 4 the $\md{O_a}^{-1}$ must lie in
$R$, so for both we need~$\Q\subseteq R$.

If as in (a) we omit relation Definition \ref{kh4def2}(iv) then we
no longer need $R$ to be a $\Q$-algebra to define
$KC_*(Y;R),KH_*(Y;R)$, and in proving Theorem \ref{kh4thm2} by the
method sketched in (a) we would no longer need $R$ to be a
$\Q$-algebra in the first way. However, we would still need $R$ to
be a $\Q$-algebra in the second way, in Step 4. If $R$ were not a
$\Q$-algebra, the proof would fail, and in fact the author can show
that $KH_*(Y;R)\cong H_*^\rsi(Y;R\ot_\Z\Q)$ in this case.
\label{khCrem1}
\end{rem}

In the rest of this appendix we go through Steps 1--5.

\subsection{Step 1: cutting the $X_a,\ti X_d$ into small pieces
$X_{ac},\ti X_{df}$}
\label{khC1}

Let $k\in\Z$, $\al\in KH_k(Y;R)$, $\sum_{a\in A}\rho_a[X_a,\bs
f_a,\bs G_a]\in KC_k(Y;R)$ representing $\al$, and $(\ti X_d,\bs{\ti
f}_d,\bs{\ti G}_d),\Ga_d$ for $d\in D$ satisfying \eq{khCeq2} be as
in Step 1. Let the indexing sets $A,D$ be disjoint, $A\cap
D=\emptyset$. We now explain how to choose tent functions $\bs T_a$
for $(X_a,\bs f_a,\bs G_a)$ for $a\in A$, and $\bs{\ti T}_d$ for
$(\ti X_d,\bs{\ti f}_d,\bs{\ti G}_d)$ for $d\in D$, satisfying
conditions (a)--(c) of Step 1, such that $\bs T_a$ cuts $X_a$ into
arbitrarily small pieces $X_{ac}$ for $c\in C_a$, and $\bs{\ti T}_d$
cuts $\ti X_d$ into arbitrarily small pieces $\ti X_{df}$ for $f\in
F_d$. This construction will be used in Steps 2 and 3. Choosing a
tent function $\bs T$ to cut one $(X,\bs f,\bs G)$ into arbitrarily
small pieces was explained in \S\ref{khA33} and \S\ref{khA35}. Our
problem is to make these choices satisfy (a)--(c).

We introduce notation similar to $P^m,Q^m,\ldots$ in Step 1 of
Appendix \ref{khB}. By Lemma \ref{kh3lem}, for each $a\in A$ and
$m\ge 0$ the triple $(\pd^mX_a,\bs f_a\vert_{\pd^mX_a},\ab\bs
G_a\vert_{\pd^mX_a})$ splits as $\pd^mX_a=X^m_{a1}\amalg\cdots
\amalg X^m_{an_a^m}$ into {\it connected\/} triples $(X_{ab}^m,\bs
f_{ab}^m,\bs G_{ab}^m)=(X_{ab}^m,\bs f_a\vert_{X_{ab}^m},\bs
G_a\vert_{X_{ab}^m})$ for $b=1,\ldots,n_a^m$. Similarly, for $d\in
D$ and $m\ge 1$ the triple $(\pd^{m-1}\ti X_d,\bs{\ti
f}_d\vert_{\pd^{m-1}\ti X_d},\ab\bs{\ti G}_d\vert_{\pd^{m-1}\ti
X_d})$ splits as $\pd^{m-1}\ti X_d=X^m_{d1}\amalg\cdots\amalg
X^m_{d\ti n_d^m}$ into connected triples $(X_{de}^m,\bs f_{de}^m,\bs
G_{de}^m)=(X_{de}^m,\bs{\ti f}_d\vert_{X_{de}^m},\bs{\ti
G}_d\vert_{X_{de}^m})$ for $e=1,\ldots,\ti n_d^m$. Also, for $d\in
D$ and $m\ge 1$ the triple $(\pd^{m-1}(\ti X_d/\Ga_d),
\bs\pi_*(\bs{\ti f}_d)\vert_{\pd^{m-1}(\ti X_d/\Ga_d)},\ab
\bs\pi_*(\bs{\ti G}_d)\vert_{\pd^{m-1}(\ti X_d/\Ga_d)})$ splits as
$\pd^{m-1}(\ti X_d/\Ga_d)=X^m_{d,-1}\amalg X^m_{d,-2}\amalg
\cdots\amalg X^m_{d,-\hat n_d^m}$ into connected $(X_{de}^m,\bs
f_{de}^m,\bs G_{de}^m)=(X_{de}^m,\bs{\ti
f}_d\vert_{X_{de}^m},\bs{\ti G}_d\vert_{X_{de}^m})$
for~$e=-1,-2,\ldots,-\hat n_d^m$.

Note that this choice of notation may at first seem odd: in
$X^m_{de}\subseteq\pd^{m-1}\ti X_d$ and $X^m_{de}\subseteq
\pd^{m-1}(\ti X_d/\Ga_d)$ we are changing from $m-1$ to $m$, and
dropping accents `$\,\tilde{\,\,}\,$'. The three sets of
$(X_{ab}^m,\bs f_{ab}^m,\bs G_{ab}^m)$ and $(X_{de}^m,\bs
f_{de}^m,\bs G_{de}^m)$ are distinguished as $A\cap D=\emptyset$, so
the subscripts `$ab$' and `$de$' never coincide, and
$X_{de}\subseteq\pd^{m-1}\ti X_d$ have $e>0$, and $X^m_{de}\subseteq
\pd^{m-1}(\ti X_d/\Ga_d)$ have~$e<0$.

Define $\ddot P^m=\bigl\{(a,b):a\in A$, $b=1,\ldots,n_a^m\bigr\}\cup
\bigl\{(d,e):d\in D$, $e=1,\ldots,\ti n_d^m\bigr\}\cup
\bigl\{(d,e):d\in D$, $e=-1,-2,\ldots,-\hat n_d^m\bigr\}$, omitting
the second and third sets when $m=0$. We write a general element of
$\ddot P^m$ as $(s,t)$. Thus $(X_{st}^m,\bs f_{st}^m,\bs G_{st}^m)$
can mean either $(X_{ab}^m,\bs f_{ab}^m,\bs G_{ab}^m)$, a component
of $\pd^m(X_a,\bs f_a,\bs G_a)$ for $a\in A$, or $(X_{de}^m,\bs
f_{de}^m,\bs G_{de}^m)$, a component of $\pd^{m-1}(\ti X_d,\bs{\ti
f}_d,\bs{\ti G}_d)$ for $d\in D$ and $e>0$ when $m\ge 1$, or
$(X_{de}^m,\bs f_{de}^m,\bs G_{de}^m)$, a component of
$\pd^{m-1}(\ti X_d/\Ga_d,\bs\pi_*(\bs{\ti f}_d),\ab\bs\pi_*(\bs{\ti
G}_d))$ for $d\in D$ and $e<0$ when~$m\ge 1$.

Since $A,D$ are finite and $\pd^mX_a=\pd^{m-1}\ti X_d=\es$ for $m\gg
0$, $\ddot P^m$ is finite with $\ddot P^m=\es$ for $m\gg 0$. Let
$M\ge 0$ be largest with $\ddot P^M\ne\es$. If $(s,t)\in\ddot P^m$
then $(\pd X_{st}^m,\bs f_{st}^m\vert_{\pd X_{st}^m},\bs
G_{st}^m\vert_{\pd X_{st}^m})$ is a union of connected components
each of which is of the form $(X_{st'}^{m+1},\bs f_{st'}^{m+1},\bs
G_{st'}^{m+1})$ for $(s,t')\in\ddot P^{m+1}$. Write $\ddot B_{st}^m$
for the set of such $t'$. Then~$\pd X_{st}^m=\coprod_{t'\in\ddot
B_{st}^m}X_{st'}^{m+1}$.

By induction on decreasing $m=M,M-1,\ldots,1,0$ we will choose tent
functions $\bs T_{st}^m$ for $(X_{st}^m,\bs f_{st}^m,\bs G_{st}^m)$
for all $(s,t)\in\ddot P^m$, satisfying the conditions:
\begin{itemize}
\setlength{\itemsep}{0pt}
\setlength{\parsep}{0pt}
\item[(i)] If $(s,t),(s',t')\in\ddot P^m$ and $(\bs a,\bs b):
(X_{st}^m,\bs f_{st}^m,\bs G_{st}^m)\ra(X_{s't'}^m,\bs
f_{s't'}^m,\bs G_{s't'}^m)$ is an isomorphism, not necessarily
orientation-preserving, then $\bs T_{st}^m=\bs b^*(\bs T_{s't'}^m)$.
That is, $T_{st}^{m,i}=T_{s't'}^{m,i}\ci b^i$ for all $i\in
I_{st}^m$, where $b^i:V^{m,i}_{st}\ra V^{m,i}_{s't'}$ is in $\bs b$
and $T_{st}^{m,i},T_{s't'}^{m,i}$ are the tent functions on
$V^{m,i}_{st},V^{m,i}_{s't'}$ in~$\bs T_{st}^m,\bs T_{s't'}^m$.
\item[(ii)] Suppose $m\ge 1$, $d\in D$, $e=1,\ldots,n_d^m$ and
$e=-1,-2,\ldots,-\hat n_d^m$ with $X^m_{de'}=\pi^{m-1}(X_{de}^m)$.
Then $\Ga_d$ acts on the indexing set $\{1,\ldots,n_d^m\}$ of
components of $\pd^{m-1}(\ti X_d,\bs{\ti f}_d,\bs{\ti G}_d)$. Write
$\Stab_{\Ga_d}(e)$ for the stabilizer subgroup of
$e\in\{1,\ldots,n_d^m\}$ in $\Ga_d$. Then $\Stab_{\Ga_d}(e)$ acts on
$(X_{de}^m,\bs f_{de}^m,\bs G_{de}^m)$ by automorphisms, and the
quotient $\bigl(X_{de}^m/\Stab_{\Ga_d}(e),\bs\pi^m_*(\bs f_{de}^m),
\bs\pi^m_*(\bs G_{de}^m)\bigr)$ is $(X_{de'}^m,\bs f_{de'}^m,\bs
G_{de'}^m)$. Write $(\bs\pi^m,\bs p^m):(X_{de}^m,\bs f_{de}^m,\bs
G_{de}^m)\ra(X_{de'}^m,\bs f_{de'}^m,\bs G_{de'}^m)$ for the
corresponding projection, so that $\bs\pi^m:X_{de}^m\ra X_{de'}^m$
is a strong submersion, and $p^{m,i}:V^i_{de}\ra V^i_{de}/
\Stab_{\Ga_d}(e)=V^i_{de'}$ is the natural projection. Then $\bs
T_{de}^m=(\bs p^m)^*(\bs T_{de'}^m)$, that is,
$T_{de}^{m,i}=T_{de'}^{m,i}\ci p^{m,i}$ for all $i\in I_{de}^m$.
\item[(iii)] If $(s,t)\in\ddot P^m$ then $(\pd X_{st}^m,\bs
f_{st}^m\vert_{\pd X_{st}^m},\bs G_{st}^m\vert_{\pd X_{st}^m})=
\coprod_{t'\in\ddot B_{st}^m}(X_{st'}^{m+1},\bs f_{st'}^{m+1},\ab\bs
G_{st'}^{m+1})$. We require $\min\bs T_{st}^m\vert_{\pd X_{st}^m}=
\coprod_{t'\in\ddot B_{st}^m}\min\bs T_{st'}^{m+1}$ under this.
\item[(iv)] Each $\bs T_{st}^m$ cuts $(X_{st}^m,\bs f_{st}^m,\bs
G_{st}^m)$ into `arbitrarily small pieces'.
\end{itemize}
Note that taking $(s,t)=(s',t')$ in (i) implies that $\bs T_{st}^m$
is invariant under $\Aut(X_{st}^m,\bs f_{st}^m,\bs G_{st}^m)$ for
all $(s,t)\in\ddot P^m$. Using this we will see later that (i),(ii)
imply (a),(b) in Step 1, and clearly (iv) implies~(c).

We will carry out the induction to choose the $\bs T_{st}^m$ after
Remark \ref{khCrem2}. First we have more notation to define, and
some work to do. Define an equivalence relation $\sim^m$ on $\ddot
P^m$ by $(s,t)\sim^m(s',t')$ if there exists an isomorphism $(\bs
a,\bs b):(X_{st}^m,\bs f_{st}^m,\bs G_{st}^m)\ra(X_{s't'}^m,\bs
f_{s't'}^m,\bs G_{s't'}^m)$, where $\bs a$ need not identify
orientations.

Define a second, weaker equivalence relation $\approx^m$ on $\ddot
P^m$ to be the equivalence relation generated by $\sim^m$ and the
additional relations when $m\ge 1$ that $(d,e)\approx^m(d,e')$ if
$e>0$ and $e'<0$, so that $X_{de}\subseteq\pd^{m-1}\ti X_d$ and
$X^m_{de'}\subseteq \pd^{m-1}(\ti X_d/\Ga_d)$, and
$X^m_{de'}=\pi^m(X_{de}^m)$, where $\pi^m:\pd^{m-1}\ti
X_d\ra\pd^{m-1}(\ti X_d/\Ga_d)$ is the continuous map associated to
$\bs\pi^m:\pd^{m-1}\ti X_d\ra \pd^{m-1}(\ti X_d/\Ga_d)$ induced by
the projection $\bs\pi:\ti X_d\ra\ti X_d/\Ga_d$. For each $(d,e)$
the image of $(X_{de}^m,\bs f_{de}^m,\bs G_{de}^m)$ under $\bs\pi$
is a connected component of $\pd^{m-1}(\ti
X_d/\Ga_d,\bs\pi_*(\bs{\ti f}_d),\ab\bs\pi_*(\bs{\ti G}_d))$, so for
each $d\in D$ and $e=1,\ldots,\ti n_d^m$ there is a unique
$e'=-1,\ldots,-\hat n_d^m$ with~$X^m_{de'}=\pi^m(X_{de}^m)$.

The relevance of $\sim^m,\approx^m$ is that if $(s,t)\sim^m(s',t')$
then (i) above implies that $\bs T_{st}^m=\bs b^*(\bs T_{s't'}^m)$.
More generally, if $(s,t)\approx^m(s',t')$ then by a finite number
of applications of (i),(ii) above we can determine $\bs T_{st}^m$ in
terms of $\bs T_{s't'}^m$, and vice versa. Thus, the equivalence
relations $\sim^m$ and $\approx^m$ correspond to the equivalences
imposed on the $\bs T_{st}^m$ by (i) and by (i),(ii) respectively,
{\it on the level of the indexing set\/} $\ddot P^m$. We can also
define such equivalence relations on the levels of {\it topological
spaces} $\coprod_{(s,t)\in\ddot P^m}X_{st}^m$, and {\it
orbifolds}~$\coprod_{(s,t)\in\ddot P^m:i\in I_{st}^m}V^{m,i}_{st}$.

Let $\sim_{\sst X}^m$ be the equivalence relation on
$\coprod_{(s,t)\in\ddot P^m}X_{st}^m$, as a set, given by
$x\sim_{\sst X}^m x'$ if $x\in X_{st}^m$, $x'\in X_{s't'}^m$, and
there exists an isomorphism $(\bs a,\bs b):(X_{st}^m,\bs
f_{st}^m,\bs G_{st}^m)\ra(X_{s't'}^m,\bs f_{s't'}^m,\bs G_{s't'}^m)$
with $a(x)=x'$, where $a:X_{st}^m\ra X_{s't'}^m$ is the
homeomorphism associated to $\bs a$. Let $\approx_{\sst X}^m$ be the
equivalence relation on $\coprod_{(s,t)\in\ddot P^m}X_{st}^m$
generated by $\sim_{\sst X}^m$ and the additional relations when
$m\ge 1$ and $d,e,e'$ and $(\bs\pi^m,\bs p^m):(X_{de}^m,\bs
f_{de}^m,\bs G_{de}^m)\ra(X_{de'}^m,\bs f_{de'}^m,\bs G_{de'}^m)$
are as in (ii) above and $x\in X_{de}^m$, $x'\in X_{de'}^m$ with
$\pi^m(x)=x'$, where $\pi^m:X_{de}^m\ra X_{de'}^m$ is the continuous
map associated to $\bs\pi^m$, then $x\approx_{\sst X}^mx'$. Define
$\ddot X^m$ to be $\bigl(\coprod_{(s,t)\in\ddot
P^m}X_{st}^m\bigr)/\approx_{\sst X}^m$, the set of equivalence
classes of $\approx_{\sst X}^m$, with the quotient topology, that
is, $\ddot U\subseteq\ddot X^m$ is open if its preimage
$U\subseteq\coprod_{(s,t)\in\ddot P^m}X_{st}^m$ is open. Write
$\pi_{\ddot X^m}:\coprod_{(s,t)\in\ddot P^m}X_{st}^m\ra\ddot X^m$
for the natural, continuous projection.

Similarly, let $\sim^{m,i}_{\sst V}$ be the equivalence relation on
$\coprod_{(s,t)\in\ddot P^m:i\in I_{st}^m}V^{m,i}_{st}$ for each
$i\in\N$ given by $v\sim^{m,i}_{\sst V}v'$ if $v\in V_{st}^{m,i}$,
$v'\in V_{s't'}^{m,i}$, and there exists an isomorphism $(\bs a,\bs
b):(X_{st}^m,\bs f_{st}^m,\bs G_{st}^m)\ra(X_{s't'}^m,\bs
f_{s't'}^m,\bs G_{s't'}^m)$ with $b^i(v)=v'$. Let
$\approx^{m,i}_{\sst V}$ be the equivalence relation on
$\coprod_{(s,t)\in\ddot P^m:i\in I_{st}^m}V^{m,i}_{st}$ generated by
$\sim^{m,i}_{\sst V}$ and the additional relations when $m\ge 1$ and
$d,e,e'$ and $(\bs\pi^m,\bs p^m):(X_{de}^m,\bs f_{de}^m,\bs
G_{de}^m)\ra(X_{de'}^m,\bs f_{de'}^m,\bs G_{de'}^m)$ are as in (ii)
above and $v\in V_{de}^{m,i}$, $v'\in V_{de'}^{m,i}$ with
$p^{m,i}(v)=v'$, then $v\approx^{m,i}_{\sst V}v'$. Define $\ddot
V^{m,i}=\bigl(\coprod_{(s,t)\in\ddot P^m:i\in
I_{st}^m}V^{m,i}_{st}\bigr)/\approx^{m,i}_{\sst V}$, as a
topological space. Write $\pi_{\ddot V^{m,i}}:\coprod_{
(s,t)\in\ddot P^m:i\in I_{st}^m}V^{m,i}_{st}\ra\ddot V^{m,i}$ for
the projection.

In the same way, we can define $\approx^{m,i}_{\sst E}$ on
$\coprod_{(s,t)\in \ddot P^m:i\in I_{st}^m}E^{m,i}_{st}$, and $\ddot
E^{m,i}=\bigl(\coprod_{(s,t)\in \ddot P^m:i\in I_{st}^m}E^{m,i}_{st}
\bigr)/\approx^{m,i}_{\sst E}$. Since the isomorphisms $(\bs a,\bs
b)$ and projections $(\bs\pi^m,\bs p^m)$ used to generate
$\approx_{\sst X}^m,\approx_{\sst V}^{m,i},\approx_{\sst E}^{m,i}$
are compatible with $s^{m,i}_{st},\psi^{m,i}_{st}$ in Kuranishi
neighbourhoods $(V_{st}^{m,i},\ldots,\psi^{m,i}_{st})$ and
coordinate changes between these neighbourhoods, and the data
$\eta_i^{m,j}$ used to define open sets $\dot V^{m,i},\dot
V^{m,ij}$, all these descend to the topological spaces~$\ddot
X^m,\ddot V^{m,i},\ddot E^{m,i}$.

That is, $\ddot E^{m,i}\ra\ddot V^{m,i}$ is a `bundle' in some weak
sense, with a well-defined `zero section', and we have continuous
section $\ddot s^{m,i}:\ddot V^{m,i}\ra\ddot E^{m,i}$ which lifts to
$s^{m,i}_{st}:V^{m,i}_{st}\ra E^{m,i}_{st}$ for each $(s,t)\in\ddot
P^m$, and a continuous map $\ddot\psi^{m,i}:(\ddot
s^{m,i})^{-1}(0)\ra\ddot X^m$ which is a homeomorphism with its
image. There is an open set $\dddot V{}^{m,i}\subset\ddot V^{m,i}$
which lifts to $\coprod_{(s,t)\in\ddot P^m:i\in I_{st}^m}\dot
V^{m,i}_{st}\subset\coprod_{ (s,t)\in\ddot P^m:i\in
I_{st}^m}V^{m,i}_{st}$.

When $j\le i$ we have an open set $\ddot V^{m,ij}\subset\ddot
V^{m,j}$ which lifts to $V^{m,ij}_{st}$ in each $V^{m,j}_{st}$, and
embeddings of topological spaces $\ddot\phi^{m,ij}:\ddot
V^{m,ij}\ra\ddot V^{m,i}$, $\hat{\ddot\phi}{}^{m,ij}:\ddot
E^{m,ij}\ra\ddot E^{m,i}$, which lift to
$(V^{m,ij}_{st},\phi^{m,ij}_{st},\hat\phi{}^{m,ij}_{st})$ on
$V^{m,j}_{st}$. Also we have an open set $\dddot
V{}^{m,ij}\subset\ddot V^{m,ij}$ which lifts to $\dot V^{m,ij}_{st}$
in each~$V^{m,ij}_{st}$.

Having got this far, it seems obvious (though perhaps optimistic) to
hope that $\ddot X^m$ is a compact Kuranishi space, $\ddot
V^{m,i},\ddot E^{m,i}$ are orbifolds, $(\ddot V^{m,i},\ab\ddot
E^{m,i},\ab\ddot s^{m,i},\ab\ddot\psi^{m,i})$ is a Kuranishi
neighbourhood on $\ddot X^m$, $(\ddot\phi^{m,ij},
\smash{\hat{\ddot\phi}{}^{m,ij}})$ is a coordinate change between
Kuranishi neighbourhoods, and all the other data in $\bs
f_{st}^m,\bs G_{st}^m$ descends to $\ddot X^m,\ddot V^{m,i},\ddot
E^{m,i}$ to define a triple $(\ddot X^m,\bs{\ddot f}{}^m,\bs{\ddot
G}{}^m)$ in Kuranishi homology. Then to choose data $\bs T_{st}^m$
satisfying (i),(ii) above for all $(s,t)\in\ddot P^m$, we could
choose a tent function $\bs{\ddot T}{}^m$ for $(\ddot X^m,\bs{\ddot
f}{}^m,\bs{\ddot G}{}^m)$ using the results of \S\ref{khA3}, and
lift it to each $V_{st}^{m,i}$ using the projections
$V_{st}^{m,i}\ra\ddot V^{m,i}$. The construction of $\approx_{\sst
V}^{m,i}$ would then imply~(i),(ii).

It turns out that this does not quite work, for rather subtle
reasons to do with the distinction between effective and
non-effective orbifolds, that will be explained in Remark
\ref{khCrem2}. We can give $\ddot V^{m,i},\ddot E^{m,i}$ the
structure of effective orbifolds in a canonical way, but
$\ddot\pi^{m,i}:\ddot E^{m,i}\ra\ddot V^{m,i}$ may fail to be an
orbibundle, and $\ddot\phi^{m,ij}:\ddot V^{m,ij}\ra\ddot V^{m,i}$
may not be an isomorphism on stabilizer groups. So we cannot make
$\ddot X^m$ into a Kuranishi space. However, these problems will not
matter when we apply the constructions of \S\ref{khA3}, so we can
still choose $\bs{\ddot T}{}^m$ for $(\ddot X^m,\bs{\ddot
f}{}^m,\bs{\ddot G}{}^m)$ and lift to $\bs T_{st}^m$ for
each~$(X_{st}^m,\bs f_{st}^m,\bs G_{st}^m)$.

We will need notions of ({\it weak\/}) {\it finite covers} for
manifolds and orbifolds. Covering maps, fundamental group, and
universal covers for topological spaces are discussed by Bredon
\cite[\S III]{Bred}, and for orbifolds by Adem at al.~\cite[\S
2.2]{ALR}.

\begin{dfn} Let $V,W$ be manifolds of the same dimension, which may
have boundary, corners, or g-corners, and $c:V\ra W$ a smooth map.
We call $c$ a {\it finite cover} if it is a surjective, proper local
diffeomorphism. Here $c$ is {\it proper} means that if $S\subseteq
W$ is compact in $W$ then $c^{-1}(S)$ is compact in $V$, and $c$ is
a {\it local diffeomorphism} means that every $v\in V$ admits an
open neighbourhood $U$ in $V$ such that $c(U)$ is open in $W$ and
$c\vert_U:U\ra c(U)$ is a diffeomorphism.

Let $c:V\ra W$ be a finite cover of {\it connected\/} manifolds
$V,W$. Pick $v_0\in V$, and set $w_0=c(v_0)\in W$. We can form the
{\it fundamental groups} $\pi_1(V,v_0),\pi_1(W,w_0)$ and {\it
universal covers} $\ti V,\ti W$ of $V,W$ with base-points $v_0,w_0$.
Then $c:(V,v_0)\ra(W,w_0)$ induces morphisms of fundamental groups
$c_*:\pi_1(V,v_0)\ra \pi_1(W,w_0)$, and $\ti c:\ti V\ra\ti W$ of
universal covers. The important fact is that $c_*:\pi_1(V,v_0)\ra
\pi_1(W,w_0)$ is {\it injective}, and $\ti c:\ti V\ra\ti W$ is a
{\it diffeomorphism}. Thus we may identify $\pi_1(V,v_0)$ with a
subgroup of $\pi_1(W,w_0)$, and $\ti V$ with $\ti W$. So we have
diffeomorphisms $W\cong \ti W/\pi_1(W,w_0)$ and $V\cong\ti
W/\pi_1(V,v_0)$. One class of examples of finite covers of manifolds
are {\it free finite group quotients}. Let $V$ be a manifold and
$\Ga$ a finite group acting freely on $V$ by diffeomorphisms. Then
$W=V/\Ga$ is a manifold, and the projection $\pi:V\ra V/\Ga$ is a
finite cover.
\label{khCdef1}
\end{dfn}

\begin{dfn} Let $V,W$ be orbifolds of the same dimension $n$, which may
have boundary, corners, or g-corners, and $c:V\ra W$ a smooth map.
We call $c$ a {\it finite cover} if the following condition holds.
Suppose $(U,\Ga,\phi)$ is an orbifold chart on $W$. Write $\pi:U\ra
U/\Ga$ for the projection. Then $\phi\ci\pi:U\ra W$ is a submersion,
and $c:V\ra W$ is a smooth map, so we can form the orbifold fibre
product $U\t_{\phi\ci\pi,W,c}V$, which is an $n$-orbifold with
smooth projection $\pi_U:U\t_{\phi\ci\pi,W,c}V\ra U$. We require
that for all such $(U,\Ga,\phi)$, $U\t_{\phi\ci\pi,W,c}V$ is an
$n$-manifold, and $\pi_U:U\t_{\phi\ci\pi,W,c}V\ra U$ is a finite
cover of manifolds.

An important class of examples of finite covers of orbifolds are
{\it finite group quotients}. Let $V$ be an orbifold, and $\Ga$ a
finite group acting on $V$ by diffeomorphisms; note that $\Ga$ need
not act freely, nor even effectively. Then $W=V/\Ga$ is an orbifold,
and the natural projection $\pi:V\ra V/\Ga$ is a finite cover of
orbifolds.

As for finite covers of manifolds, one can also describe finite
covers of connected orbifolds in terms of fundamental groups and
universal covers. However, these are not the usual fundamental
groups and universal covers of topological spaces \cite[\S
III]{Bred}, but special orbifold versions, the {\it orbifold
fundamental group} $\pi_1^\orb(V,v_0)$ and {\it orbifold universal
cover} $\ti V^\orb$ of a connected orbifold $V$. These are defined
by Adem at al.\ \cite[\S 2.2]{ALR} using their language of
groupoids, and the notion of classifying space of an orbifold.

One can also define $\pi_1^\orb(V,v_0)$ in a similar way to the
usual definition of $\pi_1(X)$ for a topological space $X$, in terms
of isotopy classes of {\it smooth\/} paths $\ga:[0,1]\ra V$ with
$\ga(0)=\ga(1)=v_0$. However, at this point the discussion of Remark
\ref{kh2rem3} becomes important. We must interpret $\ga:[0,1]\ra V$,
a smooth map of orbifolds, not in the weak sense of Satake, but as a
{\it strong map} \cite{MoPr}, {\it good map} \cite{ChRu1} or {\it
orbifold morphisms} \cite{ALR}. Thus, $\ga$ carries with it some
discrete extra data, which enables us to make natural choices of
maps $\si_{pq}$ in Definition \ref{kh2def7}. This discrete extra
data makes the difference between $\pi_1^\orb(V,v_0)$
and~$\pi_1(V,v_0)$.

For example, let $\Ga$ be a finite subgroup of GL$(n,\R)$. Then
$\R^n/\Ga$ is a connected orbifold. Choose $0\Ga$ as the basepoint.
As a topological space $\R^n/\Ga$ is contractible, so
$\pi_1(\R^n/\Ga,0\Ga)=\{1\}$, and the universal cover is $\R^n/\Ga$.
But $\pi_1^\orb(\R^n/\Ga,0\Ga)=\Ga$, with orbifold universal
cover~$\widetilde{(\R^n/\Ga)}{}^\orb=\R^n$.

As for manifolds, $\ti V^\orb$ is a connected orbifold which is
unique up to diffeomorphism with $\pi_1^\orb(\ti V^\orb,\ti
v_0)=\{1\}$, and $\pi_1^\orb(V,v_0)$ acts on $\ti V^\orb$ by
diffeomorphisms (though not necessarily freely), with $V$ naturally
diffeomorphic to $\ti V^\orb/\pi_1^\orb(V,v_0)$. If $c:V\ra W$ is a
finite cover of orbifolds with $V,W$ connected (in the usual sense
of topological spaces), $v_0\in V$, and $w_0=c(v_0)\in W$, then
$c:(V,v_0)\ra(W,w_0)$ induces an injective morphism of orbifold
fundamental groups $c_*:\pi_1^\orb(V,v_0)\ra\pi_1^\orb(W,w_0)$, and
a diffeomorphism $\ti c:\ti V^\orb\ra\ti W^\orb$ of orbifold
universal covers. Thus we may identify $\pi_1^\orb(V,v_0)$ with a
subgroup of $\pi_1^\orb(W,w_0)$, and $\ti V^\orb$ with $\ti W^\orb$,
and we have diffeomorphisms $W\cong\ti W^\orb/\pi_1^\orb(W,w_0)$
and~$V\cong\ti W^\orb/\pi_1^\orb(V,v_0)$.

The orbifold universal cover $\ti V^\orb$ of $V$ is always an {\it
effective\/} orbifold. Write $\Diff(\ti V^\orb)$ for the group of
diffeomorphisms of $\ti V^\orb$. Then $\pi_1^\orb(V,v_0)$ acts on
$\ti V^\orb$ by a group morphism $\rho:\pi_1^\orb(V,v_0)\ra\Diff(\ti
V^\orb)$. The kernel $\Ker\rho$ is a finite group, isomorphic to the
stabilizer group of a generic point of $V$. Now $V$ is an {\it
effective} orbifold if and only if generic points in $V$ have
trivial stabilizers. Thus, $V$ is effective if and only if $\rho$ is
injective. Therefore, {\it if a connected orbifold\/ $V$ is
effective, we can regard\/ $\pi_1^\orb(V,v_0)$ as a subgroup of\/
$\Diff(\ti V^\orb),$ but if\/ $V$ is not effective, we cannot}. This
will be important in Proposition~\ref{khCprop1}.
\label{khCdef2}
\end{dfn}

\begin{dfn} Let $V,W$ be orbifolds, and $c:V\ra W$ a smooth map.
Then as in Definition \ref{kh2def7}, for all $v\in V$ and $w=c(w)$
in $W$ there exist orbifold charts $(U,\Ga,\phi)$ on $V$ and
$(U',\Ga',\phi')$ on $W$ with $v\in\Im\phi$ and $w\in\Im\phi'$, a
morphism of finite groups $\rho:\Ga\ra\Ga'$, and a
$\rho$-equivariant smooth map $\si:U\ra U'$, such that $c\ci\phi(\Ga
u)=\phi'(\Ga'\si(u))$ for all $u\in U$. We can also make $U'$
smaller if necessary so that~$c(\Im\phi)=\Im\phi'$.

We call $c$ a {\it weak finite cover} if $c$ is surjective and
proper, and for all $v\in V$ we can choose
$(U,\Ga,\phi),(U',\Ga',\phi'),\rho,\si$ above such that $\si:U\ra
U'$ is a finite cover of manifolds. If we also required
$\rho:\Ga\ra\Ga'$ to be injective then this would be equivalent to
the definition of finite covers of orbifolds in Definition
\ref{khCdef2}. For weak finite covers we allow $\rho$ not to be
injective. As $\si$ is $\rho$-equivariant, $\Ker\rho$ must act
trivially on $U$. Thus, if $\rho$ is not injective then the
nontrivial group $\Ker\rho$ is contained in the stabilizer group of
every point in the nonempty open set $\phi(U/\Ga)$ in $V$, so $V$ is
a {\it non-effective} orbifold. Therefore a weak finite cover
$c:V\ra W$ is a finite cover if $V$ is effective.

If $c:V\ra W$ is a weak finite cover with $V,W$ connected then
choosing $v_0\in V$ and setting $w_0=c(v_0)$, we can form orbifold
fundamental groups $\pi_1^\orb(V,v_0),\pi_1^\orb(W,w_0)$ and
orbifold universal covers $\ti V^\orb,\ti W^\orb$, and we get a
group morphism $c_*:\pi_1^\orb(V,v_0)\ra\pi_1^\orb(W,w_0)$ and a
smooth map $\ti c:\ti V^\orb\ra\ti W^\orb$. As for finite covers of
orbifolds, it turns out that $\ti c$ is a {\it diffeomorphism}.
However, in the weak case, $c_*:\pi_1^\orb(V,v_0)\ra\pi_1^\orb(W,
w_0)$ {\it may not be injective}. Instead, if $v\in V$ is a generic
point and $w=c(v)$ (which is generic in $W$) then $\Ker c_*$ is
isomorphic to the kernel of the induced morphism of stabilizer
groups $\Stab_V(v)\ra\Stab_W(w)$. Note that $V$ is effective if and
only if $\Stab_V(v)=\{1\}$, and then $\Ker c_*$ is trivial, and $c$
is a finite cover.
\label{khCdef3}
\end{dfn}

We will show that the $\ddot V^{m,i}$ can be given the structure of
effective orbifolds. The proof uses the finiteness properties of the
gauge-fixing data~$\bs G_{st}^m$.

\begin{prop} In the situation above, there is a unique way to give
$\ddot V^{m,i}$ the structure of an effective $i$-orbifold, such
that\/ $\pi_{\ddot V^{m,i}}:\coprod_{(s,t)\in\ddot P^m:i\in
I_{st}^m} V^{m,i}_{st}\ra\ddot V^{m,i}$ lifts to a weak finite cover
of orbifolds for each\/ $i\in\N$. Furthermore:
\begin{itemize}
\setlength{\itemsep}{0pt}
\setlength{\parsep}{0pt}
\item[{\rm(a)}] $\ddot\phi^{m,ij}:\ddot V^{m,ij}\ra\ddot V^{m,i}$ lifts to
a smooth, injective map of orbifolds for $j\le i$ in $\N,$ and\/
$\d\ddot\phi^{m,ij}:T\ddot V^{m,ij}\ra(\ddot\phi^{m,ij})^*(T\ddot
V^{m,i})$ is injective. However, $\ddot\phi^{m,ij}$ need not be an
isomorphism on stabilizer groups, and so may not be an embedding in
the sense of Definition {\rm\ref{kh2def7}}. For $k\le j\le i$ we
have $\ddot\phi^{m,ij}\ci \ddot\phi^{m,jk}=\ddot\phi^{m,ik}$
on~$(\ddot\phi^{m,jk})^{-1}(\ddot V^{m,ij})\cap\ddot V^{m,ik}$.
\item[{\rm(b)}] We can also relate $\pd\ddot V^{m,i}$ to $\ddot
V^{m+1,i-1}$. By definition we have
\e
\ts\pd\bigl[\coprod_{(s,t)\in\ddot P^m:i\in I_{st}^m}
V^{m,i}_{st}\bigr]=\coprod_{(s',t')\in\ddot P^{m+1}:i-1\in
I_{st}^{m+1}}V^{m+1,i-1}_{st}.
\label{khCeq23}
\e
The equivalence relation $\approx_{\sst V}^{m,i}$ on $\coprod_{(s,t)
\in\ddot P^m:i\in I_{st}^m}V^{m,i}_{st}$ induces $\approx_{\sst\pd
V}^{m,i}$ on the l.h.s.\ of\/ {\rm\eq{khCeq23},} with homeomorphism
$(\pd[\,\coprod_{(s,t)\in\ddot P^m:i\in I_{st}^m} V^{m,i}_{st}])/
\approx_{\sst\pd V}^{m,i}\cong\pd\ddot V^{m,i},$ such that\/
$\pi_{\ddot V^{m,i}}\vert_{\pd(\cdots)}:\pd(\,\coprod_{(s,t)\in\ddot
P^m:i\in I_{st}^m} V^{m,i}_{st})\ra\pd\ddot V^{m,i}$ is the
projection from \eq{khCeq23} to its quotient by~$\approx_{\sst\pd
V}^{m,i}$.

The equivalence relation $\approx_{\sst V}^{m+1,i-1}$ on the right
hand side of\/ \eq{khCeq23} is stronger than $\approx_{\sst\pd
V}^{m,i},$ that is, $v\approx_{\sst\pd V}^{m,i}v'$ implies
$v\approx_{\sst V}^{m+1,i-1}v'$. Thus, there is a natural,
continuous, surjective projection
\begin{equation*}
\ts\bigl(\pd\bigl[\coprod_{\begin{subarray}{l}(s,t)\in\ddot P^m:\\
i\in I_{st}^m\end{subarray}}
V^{m,i}_{st}\bigr]\bigr)\big/\approx_{\sst\pd V}^{m,i}
\longra\bigl(\coprod_{\begin{subarray}{l}(s',t')\in\ddot P^{m+1}:\\
i-1\in I_{s't'}^{m+1}\end{subarray}}V^{m+1,i-1}_{s't'}\bigr)
\big/\approx_{\sst V}^{m+1,i-1},
\end{equation*}
that is, a projection $\pd\ddot V^{m,i}\ra\ddot V^{m+1,i-1}$. This
lifts to a weak finite cover of orbifolds $\ddot\pi_\pd^m:\pd\ddot
V^{m,i}\ra\ddot V^{m+1,i-1}$ with\/~$\ddot\pi_\pd^m\ci \pi_{\ddot
V^{m,i}} \vert_{\pd(\cdots)}\equiv\pi_{\ddot V^{m+1,i-1}}$.
\end{itemize}
\label{khCprop1}
\end{prop}

\begin{proof} For each $(s,t)\in\ddot P^m$ and $i\in I_{st}^m$, the
gauge-fixing data $\bs G_{st}^m$ includes a map $G_{st}^{m,i}:
E_{st}^{m,i}\ra P$, which is {\it globally finite} with constant
$N_{st}^{m,i}$, say. Identifying $V_{st}^{m,i}$ with the zero
section in $E_{st}^{m,i}$, we can regard $V_{st}^{m,i}$ as a
suborbifold of $E_{st}^{m,i}$, and so restrict $G_{st}^{m,i}$ to
$V_{st}^{m,i}$. Then $G_{st}^{m,i}\vert_{V_{st}^{m,i}}$ is also
globally finite with constant $N_{st}^{m,i}$. Thus we can consider
the map
\e
\ts\coprod_{(s,t)\in\ddot P^m:i\in I_{st}^m}
G_{st}^{m,i}\vert_{V_{st}^{m,i}}:\coprod_{(s,t)\in\ddot P^m:i\in
I_{st}^m}V^{m,i}_{st}\longra P.
\label{khCeq24}
\e
It is a globally finite map, with constant~$N=\sum_{(s,t)\in\ddot
P^m:i\in I_{st}^m}N_{st}^{m,i}$.

Now the equivalence relation $\approx_{\sst V}^{m,i}$ on
$\coprod_{(s,t)\in\ddot P^m:i\in I_{st}^m}V^{m,i}_{st}$ is generated
by maps $b^i:V_{st}^{m,i}\ra V_{s't'}^{m,i}$ from isomorphisms $(\bs
a,\bs b)$, and maps $p^{m,i}:V_{de}^{m,i}\ra V_{de'}^{m,i}$ from
quotient maps $(\bs\pi^m,\bs p^m)$. The definitions of $(\bs a,\bs
b)$ and $(\bs\pi^m,\bs p^m)$ imply that
$G_{s't'}^{m,i}\vert_{V_{s't'}^{m,i}}\ci b^i\equiv
G_{st}^{m,i}\vert_{G_{st}^{m,i}}$ and
$G_{de'}^{m,i}\vert_{V_{de'}^{m,i}}\ci p^{m,i}\equiv
G_{de}^{m,i}\vert_{V_{de}^{m,i}}$. Hence the maps
$G_{st}^{m,i}\vert_{G_{st}^{m,i}}$ are preserved by $\approx_{\sst
V}^{m,i}$, and descend to the quotient $\ddot V^{m,i}$. Thus, there
exists $\ddot G^{m,i}:\ddot V^{m,i}\ra P$ such that \eq{khCeq24} is
$\ddot G^{m,i}\ci\pi_{\ddot V^{m,i}}$. But as \eq{khCeq24} is
globally finite with constant $N$, it follows that $\pi_{\ddot
V^{m,i}}$ is also globally finite with constant $N$. That is, each
$\approx_{\sst V}^{m,i}$ equivalence class in
$\coprod_{(s,t)\in\ddot P^m:i\in I_{st}^m}V^{m,i}_{st}$ is of size
at most~$N$.

The maps $b^i,p^{m,i}$ generating $\approx_{\sst V}^{m,i}$ behave
well with respect to connected components: $b^i$ is a diffeomorphism
so identifies connected components, and $p^{m,i}:V_{de}^{m,i}\ra
V_{de'}^{m,i}$ is a finite quotient map and so takes each connected
component of $V_{de}^{m,i}$ to a connected component of
$V_{de'}^{m,i}$. It follows that $\pi_{\ddot V^{m,i}}$ behaves well
with respect to connected components. That is, the image of each
connected component of $V_{st}^{m,i}$ under $\pi_{\ddot V^{m,i}}$ is
a connected component of $\ddot V^{m,i}$, and conversely, the
inverse image of any connected component of $\ddot V^{m,i}$ is the
disjoint union of at most $N$ connected components
of~$\coprod_{(s,t)\in\ddot P^m:i\in I_{st}^m}V^{m,i}_{st}$.

Fix $\ddot v\in\ddot V^{m,i}$. Let $\ddot V_{\ddot v}$ be the
connected component of $\ddot V^{m,i}$ containing $\ddot v$. The
preimage $\pi_{\ddot V^{m,i}}^{-1}(\{\ddot v\})$ is an
$\approx_{\sst V}^{m,i}$ equivalence class, and so is of size at
most $N$. Write $\pi_{\ddot V^{m,i}}^{-1}(\{\ddot
v\})=\{v_1,v_2,\ldots,v_k\}$, with $k\le N$. Define $V_a$ to be the
connected component of $\coprod_{(s,t)\in\ddot P^m:i\in
I_{st}^m}V^{m,i}_{st}$ containing $v_a$ for $a=1,\ldots,k$. (Note
that although $v_1,\ldots,v_k$ are distinct, $V_1,\ldots,V_k$ may
not be distinct, there could be two or more $v_a$ in the same
connected component.) Then each $v_a$ lies in $V_{st}^{m,i}$ for
some $(s,t)\in\ddot P^m$, and $V_a$ is the connected component of
$V_{st}^{m,i}$ containing $v_a$. Write $\pi_1^\orb(V_a,v_a)$ for the
orbifold fundamental group, and $\ti V_a^\orb$ for the orbifold
universal cover, of $(V_a,v_a)$. Then~$V_a\cong\ti
V_a^\orb/\pi_1^\orb(V_a,v_a)$.

Now $\approx_{\sst V}^{m,i}$ is generated by diffeomorphisms
$b^i:V_{st}^{m,i}\ra V_{s't'}^{m,i}$, and finite group quotients
$p^{m,i}:V_{de}^{m,i}\ra V_{de'}^{m,i}$. In both cases, these are
{\it finite covers} of orbifolds. For each pair $a,b=1,\ldots,k$,
define $C_{ab}$ to be the set of smooth maps $c:V_a\ra V_b$ with
$c(v_a)=v_b$ that are the restrictions of $b^i$ or $p^{m,i}$ above
to connected components $V_a$. Then each $C_{ab}$ is finite, as
there are only finitely many possibilities for $(\bs a,\bs b)$ and
$(\bs\pi^m,\bs p^m)$. If $c:V_a\ra V_b$ lies in $C_{ab}$ then $c$ is
a finite cover of connected orbifolds $V_a,V_b$, and we have chosen
basepoints $v_a,v_b$ with $c(v_a)=v_b$. Hence $c$ induces an {\it
injective group morphism} $c_*:\pi_1^\orb(V_a,v_a)\ra
\pi_1^\orb(V_b,v_b)$, and a {\it diffeomorphism} $\ti c:\ti
V_a^\orb\ra \ti V_b^\orb$, as in Definition~\ref{khCdef2}.

It may be helpful to think of this data $(V_a,v_a)$ for
$a=1,\ldots,k$ and $C_{ab}$ for $a,b=1,\ldots,k$ as defining a {\it
quiver} $Q$, with vertices labelled by $1,\ldots,k$, and directed
edges $C_{ab}$ from vertex $a$ to vertex $b$ for $a,b=1,\ldots,k$.
Then $Q$ is {\it connected}, since $\{v_1,\ldots,v_k\}$ is an
$\approx_{\sst V}^{m,i}$ equivalence class, and equivalences are
induced by the maps $c:V_a\ra V_b$ in $C_{ab}$, so we can get
between any two vertices $v_d,v_e$ by a finite series of edges
$c:V_a\ra V_b$ or their inverses. However, $Q$ need not be
simply-connected. Choose a subquiver $Q'$ of $Q$ with the same set
of vertices $\{1,\ldots,k\}$ but a possibly smaller set of edges
$C'_{ab}\subseteq C_{ab}$ for $a,b=1,\ldots,k$, such that $Q'$ is
connected and simply-connected.

As above, each $c\in C'_{ab}$ induces a diffeomorphism $\ti c:\ti
V_a^\orb\ra\ti V_b^\orb$ of the universal covers. We use these
diffeomorphisms for all $a,b=1,\ldots,k$ and $c\in C'_{ab}$, that
is, all edges of $Q'$, to {\it identify\/} all the $\ti V_a^\orb$,
so we write $\ti V_a^\orb\cong\ti V^\orb$ for all $a=1,\ldots,k$.
This is well-defined as $Q'$ is connected and simply-connected.

The action of $\pi_1^\orb(V_a,v_a)$ on $\ti V_a^\orb\cong\ti V^\orb$
induces a group morphism $\rho_a:\pi_1^\orb(V_a,v_a)\ra\Diff(\ti
V^\orb)$ for $a=1,\ldots,k$, which is injective if and only if $V_a$
is effective. For all $a,b=1,\ldots,k$ and $c\in C_{ab}$, we have a
diffeomorphism $\ti c:\ti V_a^\orb\ra\ti V_b^\orb$, which under the
identifications $\ti V_a^\orb\cong\ti V^\orb$, $\ti V_b^\orb\cong\ti
V^\orb$ becomes a diffeomorphism $\ti c:\ti V^\orb\ra\ti V^\orb$, so
that $\ti c\in\Diff(\ti V^\orb)$, which is the identity if $c\in
C'_{ab}$. Define $G$ to be the subgroup of $\Diff(\ti V^\orb)$
generated by the images $\rho_a\bigl(\pi_1^\orb(V_a,v_a)\bigr)$ for
$a=1,\ldots,k$ and $\ti c\in\Diff(\ti V^\orb)$ corresponding to
$c\in C_{ab}$ for all~$a,b=1,\ldots,k$.

For each $a=1,\ldots,k$, consider the composition
\e
\ti V^\orb=\ti V^\orb_a\,{\buildrel\pi\over\longra}\,\ti
V^\orb_a/\pi_1^\orb(V_a,v_a)\cong V_a\,{\buildrel\pi_{\ddot
V^{m,i}}\over\longra}\,\ddot V_{\ddot v}.
\label{khCeq25}
\e
The projections $\pi_{\ddot V^{m,i}}:V_a\ra\ddot V_{\ddot v}$ are
quotients by the equivalence relation on $\coprod_{a=1}^kV_a$
induced by $c:V_a\ra V_b$ for $c\in C_{ab}$ and $a,b=1,\ldots,k$. It
is not difficult to see that the composition \eq{khCeq25} is
independent of $a=1,\ldots,k$, and factors as $\ti
V^\orb\,{\buildrel\pi\over\longra}\,\ti V^\orb/G\,{\buildrel
i\over\longra}\,\ddot V_{\ddot v}$, where $\pi$ is the projection to
the quotient, and $i$ is a homeomorphism.

Thus, $\ddot V_{\ddot v}$ is naturally homeomorphic to $\ti
V^\orb/G$. We claim that $\ti V^\orb/G$ is an {\it effective
orbifold}. As $G$ is a subgroup of $\Diff(\ti V^\orb)$, it acts
effectively on $\ti V^\orb$ by definition, and as $\ti V^\orb$ is
connected, the action is also locally effective. So it is enough to
show that stabilizer group $\Stab_G(\ti v)$ in $G$ of any $\ti v$ in
$\ti V^\orb$ is finite. By factoring the projection $\ti
V^\orb\ra\ti V^\orb/G$ as
\begin{equation*}
\ti V^\orb\,{\buildrel\pi\over\longra}\,\ti V^\orb/\rho_1
(\pi_1^\orb(V_1,v_1))\cong V_1'\,{\buildrel\pi\over\longra}\,\ti
V^\orb/G\cong\ddot V_{\ddot v},
\end{equation*}
where $V_1'$ is the effective orbifold underlying $V_1$, noting that
the projection $\pi_{\ddot V^{m,i}}:V_1'\ra\ddot V_{\ddot v}$ is
globally finite with constant $N$ as $V_1',V_1$ coincide as sets,
considering a generic point $\ti v'$ in $\ti V^\orb$, and using the
local effectiveness of $G$, we see that $\bmd{G/\rho_1
(\pi_1^\orb(V_1,v_1))}\le N$. Hence,
\begin{equation*}
\bmd{\Stab_G(\ti v)}\!\le\!N\bmd{\Stab_{\rho_1(\pi_1^\orb
(V_1,v_1))}(\ti v)}\!\le\!N\bmd{\Stab_{\pi_1^\orb(V_1,v_1)}(\ti
v)}\!\le\! N\Stab_{V_1}(\pi(\ti v)),
\end{equation*}
which is finite. Therefore $\ti V^\orb/G$ is an effective orbifold,
so $\ddot V_{\ddot v}$ is also an effective orbifold. As we can do
this for each connected component $\ddot V_{\ddot v}$ of $\ddot
V^{m,i}$, we have given $\ddot V^{m,i}$ the structure of an {\it
effective orbifold}, as we have to prove.

By construction, if $V$ is a connected component of $\coprod_{
(s,t)\in\ddot P^m:i\in I_{st}^m}V^{m,i}_{st}$ and $v\in V$, then the
restriction of $\pi_{\ddot V^{m,i}}$ to $V$ may be identified with
the projection $\ti V^\orb/\pi_1^\orb(V,v)\ra\ti V^\orb/G$, where
$G$ is a subgroup of $\Diff(\ti V^\orb)$ containing
$\rho(\pi_1^\orb(V,v))$ as a subgroup of finite index. By Definition
\ref{khCdef3} it follows that $\ti V^\orb/\pi_1^\orb(V,v)\ra\ti
V^\orb/G$ is a weak finite cover. Thus, the restriction of
$\pi_{\ddot V^{m,i}}:\coprod_{(s,t)\in\ddot P^m:i\in
I_{st}^m}V^{m,i}_{st}\ra\ddot V^{m,i}$ to each connected component
of $V^{m,i}_{st}$ is a weak finite cover of its image. Since
$\pi_{\ddot V^{m,i}}$ is surjective and globally finite, it follows
that $\pi_{\ddot V^{m,i}}:\coprod_{(s,t) \in\ddot P^m:i\in
I_{st}^m}V^{m,i}_{st}\ra\ddot V^{m,i}$ is a {\it weak finite cover},
as we want.

This completes the first paragraph of the proposition. Parts (a),(b)
now follow from the construction. To define $\ddot\phi^{m,ij}$ as a
map of orbifolds, for $(s,t)\in\ddot P^m$ consider
$\phi^{m,ij}_{st}:V^{m,ij}_{st}\ra V^{m,i}_{st}$, restrict to a
connected component $V^{ij}$ of $V^{m,ij}_{st}$ and the connected
component $V^i$ of $V^{m,i}_{st}$ containing
$\phi^{m,ij}_{st}(V^{ij})$, write $\ti V^{ij,\orb},\ti V^{i,\orb}$
for the orbifold universal covers of $V^{ij},V^i$. Then
$\phi^{m,ij}_{st}$ lifts to an embedding~$\ti\phi^{m,ij}_{st}:\ti
V^{ij,\orb}\ra\ti V^{i,\orb}$.

Now the image of $V^{ij}$ in $\ddot V^{m,ij}$ is isomorphic as an
orbifold to $\ti V^{ij,\orb}/G^{ij}$ for some
$G^{ij}\subset\Diff(\ti V^{ij,\orb})$, and the image of $V^i$ in
$\ddot V^{m,i}$ is isomorphic to $\ti V^{i,\orb}/G^i$ for some
$G^i\subset\Diff(\ti V^{i,\orb})$. Since $\phi^{m,ij}_{st}$ is
compatible with the equivalences defining $\approx_{\sst
V}^{m,i},\approx_{\sst V}^{m,j}$, we see that $\ti\phi^{m,ij}_{st}$
is equivariant under some morphism $\rho^{ij}:G^{ij}\ra G^i$, and so
descends to a morphism of orbifolds $\ddot\phi^{m,ij}_{st}:\ti
V^{ij,\orb}/G^{ij}\ra\ti V^{i,\orb}/G^i$. This lifts
$\ddot\phi^{m,ij}$ from a continuous map to a map of orbifolds over
the connected component $\pi_{\ddot V^{m,j}}(V^{ij})$ of $\ddot
V^{m,ij}$, and proves~(a).

For (b), to see that $\approx_{\sst V}^{m+1,i-1}$ is stronger than
$\approx_{\sst\pd V}^{m,i}$, note that the morphisms $(\bs a,\bs
b),\ab(\bs\pi^m,\bs p^m)$ on $(X_{st}^m,\bs f_{st}^m,\bs G_{st}^m)$
used to define $\approx_{\sst V}^{m,i}$ restrict on each connected
component $(X_{st'}^{m+1},\bs f_{st'}^{m+1},\bs G_{st'}^{m+1})$ of
$(\pd X_{st}^m,\bs f_{st}^m\vert_{\pd X_{st}^m},\bs
G_{st}^m\vert_{\pd X_{st}^m})$ to morphisms $(\bs a',\bs
b'),\ab(\bs\pi^{m+1},\bs p^{m+1})$ used to define $\approx_{\sst
V}^{m+1,i-1}$. But there may also be extra morphisms $(\bs a',\bs
b'),\ab(\bs\pi^{m+1},\bs p^{m+1})$ on $(X_{st'}^{m+1},\bs
f_{st'}^{m+1},\bs G_{st'}^{m+1})$ defining $\approx_{\sst
V}^{m+1,i-1}$, which are not the restriction of any $(\bs a,\bs
b),(\bs\pi^m,\bs p^m)$ to the boundary of $(X_{st}^m,\bs
f_{st}^m,\bs G_{st}^m)$.

Thus, the set of morphisms $b^{i-1},p^{m+1,i-1}$ used to define
$\approx_{\sst V}^{m+1,i-1}$ on the right hand side of \eq{khCeq23}
includes all the morphisms $b^i\vert_{\pd(\cdots)},p^{m,i}
\vert_{\pd(\cdots)}$ used to define $\approx_{\sst\pd V}^{m,i}$ on
the left hand side of \eq{khCeq23}, so $\approx_{\sst V}^{m+1,i-1}$
is stronger than $\approx_{\sst\pd V}^{m,i}$. The rest of (b) is
straightforward.
\end{proof}

\begin{rem} The proof of Proposition \ref{khCprop1} constructed an
{\it effective} orbifold structure on $\ddot V^{m,i}$ from the
orbifold structures of $V_{st}^{m,i}$. In doing so, it {\it
discarded information} encoded in non-effective orbifold structures.
If a connected component $V_a$ of some $V_{st}^{m,i}$ is
non-effective, then $\rho_a:\pi_1^\orb(V_a,v_a)\ra\Diff(\ti
V_a^\orb)$ is not injective. We made the group $G$ using only the
image $\rho_a\bigl(\pi_1^\orb(V_a,v_a)\bigr)$, so all information
about the kernel $\Ker\rho_a$ is lost.

This loss of information means that in general $\ddot E^{m,i}\ra
\ddot V^{m,i}$ may not be an orbibundle, so that $(\ddot
V^{m,i},\ddot E^{m,i},\ddot s^{m,i},\ddot\psi^{m,i})$ may not be a
Kuranishi neighbourhood; and $(\ddot\phi^{m,ij},
\smash{\hat{\ddot\phi}{}^{m,ij}})$ may not be a coordinate change
between Kuranishi neighbourhoods; and $\ddot X^m$ may not be a
Kuranishi space. It is easy to find examples of orbibundles $E\ra V$
in which $E$ is an effective orbifold, but $V$ is a non-effective
orbifold; for instance, $V=\{0\}/G$ and $E=\R^n/G$ for
$G\subseteq{\mathop{\rm GL}}(n,\R)$ a nontrivial finite subgroup.
The construction of Proposition \ref{khCprop1} would fix $E$, but
replace $V$ by its underlying effective orbifold, $\{0\}$ in our
example. Then the projection $E\ra V$, in our example
$\R^n/G\ra\{0\}$, is no longer an orbibundle.

Similarly, one can find examples of coordinate changes
$(\phi^{ij},\hat\phi^{ij})$ in which $V^{ij}$ is a non-effective
orbifold, but $V^i$ is effective, and $\phi^{ij}:V^{ij}\ra V^i$ maps
$V^{ij}$ to some orbifold stratum of $V^i$. The construction of
Proposition \ref{khCprop1} would fix $V^i$, but replace $V^{ij}$ by
its underlying effective orbifold. Then $\phi^{ij}:V^{ij}\ra V^i$
would no longer be an {\it embedding}, as required in Definition
\ref{kh2def12}(a), since as in Definition \ref{kh2def7} an embedding
induces isomorphisms on stabilizer groups. So $(\phi^{ij},
\hat\phi^{ij})$ would no longer be a coordinate change.

One might expect that this is a solvable problem, and that there
should be some way to modify the construction of Proposition
\ref{khCprop1} to include the information in non-effective
orbifolds, so that we can make $\ddot X^m$ into a compact Kuranishi
space, $\bs{\ddot f}{}^m:\ddot X^m\ra Y$ strongly smooth, and
$\bs{\ddot G}{}^m$ gauge-fixing data including the $(\ddot
V^{m,i},\ddot E^{m,i},\ddot s^{m,i},\ddot\psi^{m,i})$. However,
examples considered by the author show that {\it this really is not
possible}. If we try to construct $G$ not as a subgroup of
$\Diff(\ti V^\orb)$, but as an abstract group with morphism
$\rho:G\ra\Diff(\ti V^\orb)$ containing each $\pi_1^\orb(V_a,v_a)$
as a subgroup, there is no natural, functorial way to construct a
candidate finite group $\Ker\rho$ containing $\Ker\rho_a$
for~$a=1,\ldots,k$.
\label{khCrem2}
\end{rem}

We now return to the problem of inductively constructing tent
functions $\bs T_{st}^m$ for $(X_{st}^m,\bs f_{st}^m,\bs G_{st}^m)$
satisfying (i)--(iv) at the beginning of \S\ref{khC1}. Firstly, by
reverse induction on $m=M,M-1,M-2,\ldots,1,0$ we choose $\bs{\ddot
T}{}^m=(\ddot T^{m,i}:i\in\N)$, where $\ddot T^{m,i}$ is a tent
function on $\dddot V{}^{m,i}$ which extends to a tent function on
an open neighbourhood of the closure $\overline{\dddot V{}^{m,i}}$
of $\dddot V{}^{m,i}$ in $\ddot V^{m,i}$, satisfying:
\begin{itemize}
\setlength{\itemsep}{0pt}
\setlength{\parsep}{0pt}
\item[(A)] whenever $j\le i$ we have $\min\ddot T^{m,j}\vert_{\dddot
V{}^{m,ij}}\equiv\min\ddot T^{m,i}\ci\dot\phi^{m,ij}\vert_{\dddot
V{}^{m,ij}}$, and the submanifolds $S_{\{i_1,\ldots,i_l\}}$ of
Definition \ref{khAdef1} for $\ddot T^{m,i}$ intersect
$\ddot\phi^{ij}(\dddot V{}^{m,ij})$ transversely in $\dddot
V{}^{m,i}$ wherever $t_{i_1}(u)=\min\ddot T^{m,i}(u)$.
\item[(B)] for all $i$ we have $\min\ddot T^{m,i}\vert_{\pd\dddot
V{}^{m,i}}\equiv \min\ddot T^{m+1,i-1}\ci\ddot\pi_\pd^m$, where
$\ddot\pi_\pd^m$ is as in Proposition~\ref{khCprop1}(b).
\item[(C)] $\ddot T^{m,i}$ and its extension cut $\dddot V{}^{m,i},
\overline{\dddot V{}^{m,i}}$ into `arbitrarily small pieces'.
\end{itemize}
Here (A) is taken directly from Definition \ref{khAdef11}, the
definition of tent functions $\bs T$ for triples $(X,\bs f,\bs G)$.
Part (B) will imply (iii) at the beginning of \S\ref{khC1}, and also
corresponds to the condition $\min\bigl(\bs T\vert_{\pd
X}\bigr)=\min\bs T'$ in Theorem \ref{khAthm5}. Part (C) will imply
(iv) at the beginning of~\S\ref{khC1}.

Suppose for the moment that we could treat $(\ddot X^m,\bs{\ddot
f}{}^m,\bs{\ddot G}{}^m)$ as a triple in Kuranishi homology, and
that the morphisms $\ddot\pi_\pd^m:\pd\ddot V^{m,i}\ra\ddot
V^{m+1,i-1}$ lift to $(\bs{\ddot p}{}_\pd^m,
\bs{\ddot\pi}{}_\pd^m):(\pd\ddot X^m,\bs{\ddot f}{}^m\vert_{\pd\ddot
X^m},\bs{\ddot G}{}^m\vert_{\pd\ddot X^m})\ra(\ddot X^m,\bs{\ddot
f}{}^m,\bs{\ddot G}{}^m)$. Then $\bs{\ddot T}{}^m$ is a tent
function for $(\ddot X^m,\bs{\ddot f}{}^m,\bs{\ddot G}{}^m)$ which
cuts $(\ddot X^m,\bs{\ddot f}{}^m,\bs{\ddot G}{}^m)$ into
arbitrarily small pieces, and (B) says that $\min\bigl(\bs{\ddot
T}{}^m\vert_{\pd\ddot X^m}\bigr)=\min\bs{\ddot
T}{}^{m+1}\ci\bs{\ddot\pi}{}_\pd^m$. We could choose such $\bs{\ddot
T}{}^m$ by reverse induction on $m$ using the results
of~\S\ref{khA35}.

Here is the important point. As in Remark \ref{khCrem2}, we cannot
make $(\ddot X^m,\bs{\ddot f}{}^m,\ab\bs{\ddot G}{}^m)$ into a
triple in Kuranishi homology. {\it But to apply the results of\/}
\S\ref{khA35}, {\it this does not matter at all}. The ways in which
$(\ddot X^m,\bs{\ddot f}{}^m,\bs{\ddot G}{}^m)$ fails to be a triple
in Kuranishi homology --- that $\ddot\phi^{m,ij}:\ddot
V^{m,ij}\ra\ddot V^{m,i}$ may not be an embedding, since it may not
be an isomorphism on stabilizer groups as in Proposition
\ref{khCprop1}(a), for instance --- do not affect any of the proofs
in \S\ref{khA3}. So we can still use the arguments of \S\ref{khA3}
to choose $\bs{\ddot T}{}^m$ inductively, just as if $(\ddot
X^m,\bs{\ddot f}{}^m,\bs{\ddot G}{}^m)$ were a triple in Kuranishi
homology.

Therefore, by reverse induction on $m=M,M-1,\ldots,0$ we use
Definition \ref{khAdef16} and Theorem \ref{khAthm5} to choose
$\bs{\ddot T}{}^m$. To apply these we must verify the assumption in
Definition \ref{khAdef15} that the prescribed values $\min\bs{\ddot
T}{}^{m+1}\ci\bs{\ddot\pi}{}_\pd^m$ for $\min\bigl(\bs{\ddot
T}{}^m\vert_{\pd\ddot X^m}\bigr)$ restrict to $\bs\si$-invariant
data on $\pd^2\ddot X^m$. That is, to be able to choose $\ddot
T^{m,i}$, we must check that $\bigl(\min\ddot T^{m+1,i-1}\ci\ddot
\pi_\pd^m\bigr)\vert_{\pd^2\dddot V{}^{m,i}}$ is invariant under the
natural involution $\si:\pd^2\dddot V{}^{m,i}\ra\pd^2\dddot
V{}^{m,i}$. But (B) for $(m+1,i-1)$ gives
\e
\bigl(\min\ddot T^{m+1,i-1}\ci\ddot
\pi_\pd^m\bigr)\vert_{\pd^2\dddot V{}^{m,i}}\equiv\min\ddot
T^{m+2,i-2}\ci\ddot\pi_\pd^{m+1}\ci\bigl(\ddot\pi_\pd^m
\vert_{\pd^2\dddot V{}^{m,i}}\bigr).
\label{khCeq26}
\e

We claim that $\ddot\pi_\pd^{m+1}\ci\bigl(\ddot\pi_\pd^m
\vert_{\pd^2\dddot V{}^{m,i}}\bigr)\ci\si\equiv\ddot\pi_\pd^{m+1}
\ci\bigl(\ddot\pi_\pd^m\vert_{\pd^2\dddot V{}^{m,i}}\bigr)$ as maps
$\pd^2\dddot V{}^{m,i}\ra\dddot V{}^{m+2,i-2}$. Combining this with
\eq{khCeq26} shows that $\bigl(\min\ddot T^{m+1,i-1}\ci\ddot
\pi_\pd^m\bigr)\vert_{\pd^2\dddot V{}^{m,i}}$ is $\si$-invariant, as
we have to show. To prove the claim, note that for each
$(s,t)\in\ddot P^m$, we have a natural orientation-reversing
involution $\bs\si:\pd^2(X_{st}^m,\bs f_{st}^m,\bs G_{st}^m)
\ra\pd^2(X_{st}^m,\bs f_{st}^m,\bs G_{st}^m)$. The restriction of
$\bs\si$ to each connected component $(X_{st'}^{m+2},\bs
f_{st'}^{m+2},\bs G_{st'}^{m+2})$ of $\pd^2(X_{st}^m,\bs
f_{st}^m,\bs G_{st}^m)$ is one of the isomorphisms $(\bs a,\bs b)$
in the definition of~$\approx^{m+2}, \approx_{\sst
X}^{m+2},\ab\approx_{\sst V}^{m+2,i-2}$.

Hence, $\si:\pd^2V_{st}^{m,i}\ra\pd^2V_{st}^{m,i}$ is one of the
diffeomorphisms in the definition of $\approx_{\sst V}^{m+2,i-2}$.
As $\si:\pd^2\dddot V{}^{m,i}\ra\pd^2\dddot V{}^{m,i}$ lifts locally
to $\si:\pd^2V_{st}^{m,i}\ra\pd^2V_{st}^{m,i}$, which then projects
to the identity on $\dddot V{}^{m+2,i-2}$ by definition of
$\approx_{\sst V}^{m+2,i-2}$, the claim follows. Therefore by
induction, we can use the results of \S\ref{khA35} to choose
$\bs{\ddot T}{}^m$ satisfying (A)--(C) above for~$m=M,M-1,\ldots,0$.

Now, for all $m=0,\ldots,M$, $(s,t)\in\ddot P^m$ and $i\in I^m_{st}$
define $T^{m,i}_{st}:\dot V^{m,i}_{st}\ra F\bigl([1,\iy)\bigr)$ by
$T^{m,i}_{st}=\ddot T^{m,i}\ci\pi_{\ddot V^{m,i}}\vert_{\dot
V^{m,i}_{st}}$. Since $\ddot T^{m,i}:\dddot V{}^{m,i}\ra F\bigl(
[1,\iy)\bigr)$ is a tent function and $\pi_{\ddot V^{m,i}}$ is a
weak finite cover of orbifolds by Proposition \ref{khCprop1}, it
follows that $T^{m,i}_{st}$ is a {\it tent function} on $\dot
V^{m,i}_{st}$, which extends to a tent function on an open
neighbourhood of the closure $\ov{\dot V^{m,i}_{st}}$ of $\dot
V^{m,i}_{st}$ in $V^{m,i}_{st}$, since $\ddot T^{m,i}$ extends to an
open neighbourhood of $\ov{\dddot V{}^{m,i}}$ in $\ddot V{}^{m,i}$.
Define $\bs T^m_{st}=(T^{m,i}_{st}:i\in I^m_{st})$. Part (A) above
implies that the $T^{m,i}_{st}$ for $i\in I^m_{st}$ satisfy the
conditions of Definition \ref{khAdef11}, and therefore $\bs
T^m_{st}$ is a {\it tent function} for $(X_{st}^m,\bs f_{st}^m,\bs
G_{st}^m)$. We claim that these $\bs T^m_{st}$ for all
$m=0,\ldots,M$ and $(s,t)\in\ddot P^m$ satisfy conditions (i)--(iv)
at the beginning of~\S\ref{khC1}.

Parts (i),(ii) hold by construction of $\ddot V^{m,i}$ using
$\approx_{\sst V}^{m,i}$. For (i), if $(s,t),(s',t')$ lie in $\ddot
P^m$, $(\bs a,\bs b): (X_{st}^m,\bs f_{st}^m,\bs
G_{st}^m)\ra(X_{s't'}^m,\bs f_{s't'}^m,\bs G_{s't'}^m)$ is an
isomorphism, and $i\in I_{st}^m$, then $b^i:V_{st}^{m,i}\ra
V_{s't'}^{m,i}$ is one of the maps generating the equivalence
relation $\approx_{\sst V}^{m,i}$ on $\coprod_{(s,t)\in\ddot
P^m:i\in I_{st}^m}V^{m,i}_{st}$. Since $\ddot V^{m,i}=\bigl(
\coprod_{(s,t)\in\ddot P^m:i\in I_{st}^m}V^{m,i}_{st}\bigr)/
\approx^{m,i}_{\sst V}$, and $\pi_{\ddot V^{m,i}}:\coprod_{(s,t)
\in\ddot P^m:i\in I_{st}^m}V^{m,i}_{st}\ra\ddot V^{m,i}$ is the
natural projection, we see that $\pi_{\ddot V^{m,i}}\vert_{\dot
V^{m,i}_{st}}\equiv \pi_{\ddot V^{m,i}}\vert_{\dot V^{m,i}_{s't'}}
\ci b^i$, and composing with $\ddot T^{m,i}$ gives $T_{st}^{m,i}
=T_{s't'}^{m,i}\ci b^i$. Part (ii) is the same, as the $p^{m,i}$ in
(ii) are also maps generating $\smash{\approx_{\sst V}^{m,i}}$.
Parts (iii),(iv) easily follow from (B),(C) above, respectively.
This completes the construction of $\bs T^m_{st}$ for all $m,s,t$
satisfying (i)--(iv) at the beginning of~\S\ref{khC1}.

In the notation of the beginning of \S\ref{khC1}, for $a\in A$ and
$d\in D$ we have
\begin{gather*}
(X_a,\bs f_a,\bs G_a)=\ts\coprod_{b=1}^{n_a^0}(X_{ab}^0,\bs
f_{ab}^0,\bs G_{ab}^0),\quad (\ti X_d,\bs{\ti f}_d,\bs{\ti
G}_d)=\coprod_{e=1}^{\ti n_d^1} (X_{de}^1,\bs
f_{de}^1,\bs G_{de}^1)\\
\text{and}\qquad (\ti X_d/\Ga_d,\bs\pi_*(\bs{\ti
f}_d),\bs\pi_*(\bs{\ti G}_d))=\ts\coprod_{e=-1}^{-\hat
n_d^1}(X_{de}^1,\bs f_{de}^1,\bs G_{de}^1).
\end{gather*}
Using these equations, define tent functions $\bs T_a$ for $(X_a,\bs
f_a,\bs G_a)$ and $\bs{\ti T}_d$ for $(\ti X_d,\bs{\ti f}_d,\bs{\ti
G}_d)$ and $\bs{\hat T}_d$ for $(\ti X_d/\Ga_d,\bs\pi_*(\bs{\ti
f}_d),\bs\pi_*(\bs{\ti G}_d))$ by $\bs T_a=\coprod_{b=1}^{n_a^0}\bs
T_{ab}^0$, and $\bs{\ti T}_d=\coprod_{e=1}^{\ti n_d^1}\bs T_{de}^1$,
and $\bs{\hat T}_d=\coprod_{e=-1}^{-\hat n_d^1}\bs T_{de}^1$. We
claim that these $\bs T_a,\bs{\ti T}_d$ satisfy conditions (a)--(c)
of Step 1 of \S\ref{khC}, shortly after equation \eq{khCeq2}.

We prove (a)--(c) using (i)--(iv) at the beginning of \S\ref{khC1},
which the $\bs T_{st}^m$ satisfy. For part (a), we have $\bs
T_a\vert_{\pd X_a}=\coprod_{b=1}^{n_a^1}(X_{ab}^1,\bs f_{ab}^1,\bs
G_{ab}^1)$, so it follows from (i) that $\bs T_a\vert_{\pd X_a}$ is
invariant under $\Aut(\pd X_a,\bs f_a\vert_{\pd X_a},\bs
G_a\vert_{\pd X_a})$, and similarly $\bs{\ti T}_d$ is invariant
under $\Aut(\ti X_d,\bs{\ti f}_d,\bs{\ti G}_d)$. Also (ii) implies
that $\bs\pi_*(\bs{\ti T}_d)=\bs{\hat T}_d$, so again using (i) we
see that $\bs\pi_*(\bs{\ti T}_d)$ is invariant under $\Aut(\ti
X_d/\Ga_d,\bs\pi_*(\bs{\ti f}_d),\bs\pi_*(\bs{\ti G}_d))$. Part (b)
follows from (i) with $m=1$, and (c) from (iv). This completes
Step~1.

\subsection{Step 2: lifting to $\acute X_{ac}$ with trivial stabilizers}
\label{khC2}

We continue to use the notation of Step 1. For each $a\in A$, write
$(X_{ac},\ab\bs f_{ac},\ab\bs G_{ac})$ for $c\in C_a$ for the
arbitrarily small pieces into which $(X_a,\bs f_a,\bs G_a)$ is `cut'
using the tent function $\bs T_a$, in the notation of \eq{khAeq29}.
We must show that we can choose the $\bs T_a$ so that
$(X_{ac},\ab\bs f_{ac},\ab\bs G_{ac})\cong\bigl(\acute
X_{ac}/\Ga_{ac},\bs\pi_*(\bs{\acute f}_{ac}),\bs\pi_*(\bs{\acute
G}_{ac})\bigr)$ for each $c\in C_a$, where $(\acute
X_{ac},\bs{\acute f}_{ac},\bs{\acute G}_{ac})$ is a triple and
$\acute X_{ac}$ has {\it trivial stabilizers}.

Let $a\in A$ and $p\in X_a$, and write $\Ga_p$ for the stabilizer
group $\Stab(p)$. We first show that a small neighbourhood of $p$ in
$ X_a$ may be written, as a Kuranishi space, as a quotient $\acute
X_p/\Ga_p$, where $\acute X_p$ is a Kuranishi space with trivial
stabilizers.

Let $(V_p,\ldots,\psi_p)$ be a Kuranishi neighbourhood in the germ
of $ X_a$ at $p$, and set $v=\psi_p^{-1}(p)$. As $V_p$ is an
orbifold, by \S\ref{kh21} there exists an orbifold chart $(\acute
V_p,\Ga_p,\xi_p)$ for $V_p$ near $v$, with group $\Ga_p$. That is,
$\acute V_p$ is a smooth manifold acted on by $\Ga_p$, and
$\xi_p:\acute V_p/\Ga_p\ra V_p$ is a diffeomorphism with an open
neighbourhood of $v$ in $V_p$. There is a unique $\acute p\in\acute
V_p$ with $\xi_p(\Ga_p\acute p)=v$, which is fixed by $\Ga_p$.
Making $V_p$ smaller if necessary, we can assume $V_p=\xi_p(\acute
V_p/\Ga_p)$. Write $\pi_p:\acute V_p\ra V_p$ for the corresponding
projection, which is a submersion. Define $\acute E_p=\pi^*(E_p)$.
Then $\acute E_p\ra\acute V_p$ is a vector bundle, in the usual
manifold sense. Define $\acute s_p=s_p\ci\pi_p$, so that $\acute
s_p$ is a smooth section of~$\acute E_p$.

Define a topological space $\acute X_p=\acute s_p^{-1}(0)$, and let
$\acute\psi_p$ be the identity map on $\acute X_p$. Then $(\acute
V_p,\acute E_p,\acute s_p,\acute\psi_p)$ is a Kuranishi
neighbourhood covering $\acute X_p$. The action of $\Ga_p$ on
$\acute V_p$ lifts naturally to an action on $\acute E_p$, and
$\acute s_p$ is $\Ga_p$-equivariant, so the subset $\acute X_p$ is
preserved by $\Ga_p$, and $\Ga_p$ acts on $\acute X_p$, with
$\psi_p\ci\xi_p:\acute X_p/\Ga_p\ra X_a$ a homeomorphism with an
open neighbourhood $\Im\psi_p$ of $p$ in $ X_a$. Also
$\psi_p\ci\pi_p:\acute X_p\ra\Im\psi_p\subset X_a$ is continuous.

We shall define a Kuranishi structure $\acute\ka_p$ on $\acute X_p$.
Let $\acute q\in\acute X_p$, and set $q=\psi_p\ci\pi_p(\acute q)$.
Then $q\in\Im\psi_p$, so for all sufficiently small
$(V_q,\ldots,\psi_q)$ in the germ of $ X_a$ at $q$ there is a
coordinate change $(\phi_{pq},\hat\phi_{pq}):(V_q,\ldots,\psi_q)\ra
(V_p,\ldots,\psi_p)$. Define a Kuranishi neighbourhood $(\acute
V_{\acute q},\ldots,\acute\psi_{\acute q})$ of $\acute q$ in $\acute
X_p$ to be
\e
(\acute V_{\acute q},\acute E_{\acute q},\acute s_{\acute
q},\acute\psi_{\acute q})= \bigl(\acute
V_p\t_{\pi_p,V_p,\phi_{pq}}V_q,\pi_{V_q}^*(E_q),
s_q\ci\pi_{V_q},\pi_{\acute V_p}\bigr).
\label{khCeq27}
\e
Since $\acute V_p$ is a manifold, $\pi_p$ is a submersion, and
$\phi_{pq}$ is an embedding, $\acute V_{\acute q}$ is a manifold,
and $\acute E_{\acute q}\ra\acute V_{\acute q}$ a vector bundle.
Since $s_p\ci\phi_{pq}\equiv\hat\phi_{pq}\ci s_q$, $\phi_{pq}$ takes
$s_q^{-1}(0)$ to $s_p^{-1}(0)$, and so $\acute\psi_{\acute
q}=\pi_{\acute V_p}$ takes $\acute s_{\acute q}^{-1}(0)$ to $\acute
s_{\acute p}^{-1} (0)=\acute X_p$, as we need.

Define the germ of Kuranishi neighbourhoods at $\acute q$ in
$\acute\ka_p$ to be the equivalence class of $(\acute V_{\acute
q},\ldots,\acute\psi_{\acute q})$. This is independent of the choice
of $(V_q,\ldots,\psi_q)$, and so well-defined. It is now easy to
show that the germ of coordinate changes between Kuranishi
neighbourhoods on $\Im\psi_p\subset X_a$ induces a germ of
coordinate changes between Kuranishi neighbourhoods on $\acute X_p$,
giving a Kuranishi structure $\acute\ka_p$ on $\acute X_p$. Since
the $\acute V_{\acute q}$ are all manifolds, $(\acute
X_p,\acute\ka_p)$ has {\it trivial stabilizers}. Furthermore,
$\acute\ka_p$ is invariant under the action of $\Ga_p$ on $\acute
X_p$, and $\psi_p\ci\pi_p:\acute X_p\ra\Im\psi_p\subset X_a$ lifts
to an isomorphism of Kuranishi spaces between $(\acute
X_p,\acute\ka_p)/\Ga_p$ and the restriction of the Kuranishi
structure on $ X_a$ to~$\Im\psi_p$.

Thus, for all $a\in A$ and $p\in X_a$, we can find an open
neighbourhood $\Im\psi_p$ of $p$ in $X_a$, and an isomorphism of
Kuranishi spaces $\Im\psi_p\cong\acute X_p/\Ga_p$, where $\acute
X_p$ has trivial stabilizers, and $\Ga_p$ is a finite group. Such
$\Im\psi_p$ cover $X_a$, so by Step 1 we can cut the $X_a$ into
pieces $X_{ac}$ such that each $X_{ac}$ is contained in some such
subset $\Im\psi_p$ in~$X_a$.

Then we have $X_{ac}\cong\acute X_{ac}/\Ga_{ac}$, where $\acute
X_{ac}= \acute X_p\t_{\bs\pi, X_a,\bs\io} X_{ac}$ is a compact
Kuranishi space with trivial stabilizers, and $\Ga_{ac}=\Ga_p$ is a
finite group acting on $\acute X_{ac}$ by strong diffeomorphisms.
Here $\bs\pi:\acute X_p\ra\Im\psi_p\subset X_a$ is the obvious
projection, which is a strong submersion, and $\bs\io: X_{ac}\ra
X_a$ is the inclusion. Write $\bs\pi_{ac}:\acute X_{ac}\ra X_{ac}$
for the projection, which is a strong submersion, and $\bs{\acute
f}_{ac}=\bs{ f}_{ac}\ci\bs\pi_{ac}:\acute X_{ac}\ra Y$, which is
strongly smooth.

It remains to show that we can choose $\bs T_a$ for $a\in A$ and
$X_{ac},\acute X_{ac}$ for $c\in C_a$ so that $\bs G_{ac}$ on
$X_{ac}$ lifts to $\bs{\acute G}_{ac}$ on $\acute X_{ac}$. In doing
this, it is important that we can not only make the $X_{ac}$
arbitrarily small, we can also make the subsets $V^i_{ac}$ of
$V^i_a$ for $i\in I_a$ arbitrarily small. This is part (C) in
\S\ref{khC1}, and can be seen from the construction in \S\ref{khA33}
used to define~$\bs T_a$.

In the situation above, let $q\in\Im\psi_p$, and suppose $i\in I_a$
and $v\in(s^i_a)^{-1}(0)\subseteq V^i_a$ with $\psi^i_a(v)=q$. Let
$(V_q,\ldots,\psi_q)$ be a sufficiently small Kuranishi
neighbourhood in the germ of $ X_a$ at $q$. Then as
$(V^i_a,\ldots,\psi^i_a)$ is a compatible Kuranishi neighbourhood on
$ X_a$, we are given a coordinate change
$(\phi_q^i,\hat\phi_q^i):(V_q,\ldots,\psi_q)\ra(V^i_a,\ldots,
\psi^i_a)$. From \eq{khCeq27} we have a Kuranishi neighbourhood
$(\acute V_{\acute q},\acute E_{\acute q},\acute s_{\acute
q},\acute\psi_{\acute q})$ on $\acute X_p$ with $\acute V_{\acute
q}$ a manifold, a diffeomorphism $\xi_q:\acute V_{\acute q}/\Ga_p\ra
V_q$, and a submersion~$\pi_q=\pi_{V_q}:\acute V_{\acute q}\ra V_q$.

Let $W$ be an open neighbourhood of $v$ in $V^i_a$. Since
$\phi_q^i:V_q\ra V^i_a$ is an embedding with
$\phi_q^i(\psi_q^{-1}(q))=v$, if $W$ is sufficiently small we may
extend the diffeomorphism $\xi_q:\acute V_{\acute q}/\Ga_p\ra V_q$
to a diffeomorphism $\ti\xi_q:\acute W/\Ga_p\ra W$, where $\acute W$
is a manifold with a $\Ga_p$ action, such that $\phi_q^i$ lifts over
$(\phi_q^i)^{-1}(W)$ to a smooth $\Ga_p$-equivariant map
$\ti\phi_q^i:(\phi_q^i\ci\pi_q)^{-1}(W)\ra\acute W$ from the open
subset $(\phi_q^i\ci\pi_q)^{-1}(W)$ in $\acute V_{\acute q}$.
Furthermore, if $W$ is connected these $\acute
W,\ti\xi_q,\ti\phi_q^i$ are unique up to canonical isomorphism.

As we can do this for any $q\in\Im\psi_p\cap\Im\psi_a^i$, and the
diffeomorphisms $\xi_q:\acute V_{\acute q}/\Ga_p\ra V_q$ all come
from an isomorphism $\acute X_p/\Ga_p\cong\Im\psi_p\subset X_a$ and
so are compatible under coordinate changes, and the $\acute
W,\ti\xi_q,\ti\phi_q^i$ are unique up to {\it canonical\/}
isomorphism above, we can patch all these local lifts together. This
proves that if $W$ is any sufficiently small open neighbourhood of
some subset in $(\psi^i_a)^{-1}(\Im\psi_p)$ in $V^i_a$, and every
connected component of $W$ meets $(\psi^i_a)^{-1}(\Im\psi_p)$, then
we can write $W\cong\acute W/\Ga_p$ for some manifold $\acute W$, so
that all embeddings $\phi_q^i:V_q\ra W\subset V^i_a$ for
$q\in\Im\psi_p\cap\psi_a^i(W)$ lift canonically to
embeddings~$\ti\phi_q^i:\acute V_{\acute q}\ra\acute W$.

Let us choose the $\bs T_a$ for $a\in A$ so that all $V^i_{ac}$ are
sufficiently small in this sense, which is possible by Step 1, and
define $\acute X_{ac},\Ga_{ac}=\Ga_p$ as above with
$X_{ac}\cong\acute X_{ac}/\Ga_{ac}$. Then there exist a manifold
$\acute V_{ac}^i$, an action of $\Ga_{ac}$ on $\acute V_{ac}^i$, and
a diffeomorphism $\acute V_{ac}^i/\Ga_{ac}\cong V_{ac}^i$, which are
unique up to canonical isomorphism. To prove this, note that as
$(V_{ac}^i,\ldots, \psi_{ac}^i)$ is part of an excellent coordinate
system, every connected component of $V_{ac}^i$
intersects~$(s_{ac}^i)^{-1}(0)$.

Define a Kuranishi neighbourhood $(\acute V_{ac}^i,\acute
E_{ac}^i,\acute s_{ac}^i,\acute\psi_{ac}^i)$ on $\acute X_{ac}$,
where $\acute E_{ac}^i=\pi_{V_{ac}^i}^*(E_{ac}^i)$ and $\acute
s_{ac}^i\equiv s_{ac}^i\ci\pi_{V_{ac}^i}$, for
$\pi_{V_{ac}^i}:\acute V_{ac}^i\ra V_{ac}^i$ the projection, and
$\acute\psi_{ac}^i$ is the natural lift of
$\psi_{ac}^i\ci\pi_{V_{ac}^i}:(\acute s_{ac}^i)^{-1}(0)\ra X_{ac}$
to $\acute X_{ac}$. These extend to an {\it excellent coordinate
system\/} $(\bs{\acute I}_{ac},\bs{\acute\eta}_{ac})=\bigl(\bigl(
I_{ac},(\acute V_{ac}^i,\ab\ldots,\ab\acute\psi_{ac}^i),\acute
f^i_{ac}:i\in
I_{ac},\ldots\bigr),(\acute\eta_{i,ab},\acute\eta_{i,ab}^j:i,j\in
I_{ac})\bigr)$ for $(\acute X_{ac},\bs{\acute f}_{ac})$, where
$\acute f^i_{ac}=f^i_{ac}\ci\pi_{V^i_{ac}}$, $\acute\eta_{i,ab}=
\eta_{i,ab}\ci\pi_{X_{ac}}$, and $\acute\eta_{i,ab}^j=
\eta_{i,ab}^j\ci\pi_{V^j_{ac}}$. We also define $\bs{\acute
G}_{ac}=\bigl(\bs{\acute I}_{ac},\bs{\acute\eta}_{ac},\acute
G^i_{ac}:i\in\acute I_{ac}\bigr)$, where $\acute G^i_{ac}:\acute
E^i_{ac}\ra P$ is given by $\acute
G^i_{ac}=G^i_{ac}\ci\pi_{E^i_{ac}}$ for $\pi_{E^i_{ac}}:\acute
E_{ac}^i\ra E_{ac}^i$ the natural projection.

Then $(\acute X_{ac},\bs{\acute f}_{ac},\bs{\acute G}_{ac})$ is a
triple in Kuranishi homology, and $\Ga_{ac}$ acts on $(\acute
X_{ac},\bs{\acute f}_{ac},\bs{\acute G}_{ac})$ by automorphisms with
$(\acute X_{ac}/\Ga_{ac},\bs\pi_*(\bs{\acute
f}_{ac}),\bs\pi_*(\bs{\acute G}_{ac}))\cong\ab(X_{ac},\ab\bs
f_{ac},\ab\bs G_{ac})$, so Definition \ref{kh4def2}(iv) gives
\eq{khCeq7}, and Step 2 follows.

\subsection{Step 3: lifting to a cycle in $\widetilde{KC}_k(Y;R)$}
\label{khC3}

As in Step 3 we now {\it change notation\/} from $\acute X_{ac}$ for
$a\in A$ and $c\in C_a$ back to $X_a$ for $a\in A$. That is, by Step
2 we can represent $\al$ by a cycle $\sum_{a\in A}\rho_a[X_a,\bs
f_a,\bs G_a]$ in $KC_k(Y;R)$, where each $X_a$ is a compact,
oriented Kuranishi space with {\it trivial stabilizers}, and each
$V^i_a$ in $\bs G_a$ is a {\it manifold}. Equation \eq{khCeq2} holds
for $(\ti X_d,\bs{\ti f}_d,\bs{\ti G}_d)$ for $d\in D$. These
$A,X_a,\bs f_a,\bs G_a,\rho_a,D,\ab\ti X_d,\ab\bs{\ti
f}_d,\ab\bs{\ti G}_d,\Ga_d,\eta_d$ are different from those in Steps
1 and~2.

We now repeat Step 1, choosing tent functions $\bs T_a$ for
$(X_a,\bs f_a,\bs G_a)$ for $a\in A$, and $\bs{\ti T}_d$ for $(\ti
X_d,\bs{\ti f}_d,\bs{\ti G}_d)$ for $d\in D$, satisfying conditions
(a)--(c) of Step 1, such that $\bs T_a$ cuts $(X_a,\bs f_a,\bs G_a)$
into small pieces $(X_{ac},\ab\bs f_{ac},\ab\bs G_{ac})$ for $c\in
C_a$, and $\bs{\ti T}_d$ cuts $(\ti X_d,\bs{\ti f}_d,\bs{\ti G}_d)$
into small pieces $(\ti X_{df},\bs{\ti f}_{df},\bs{\ti G}_{df})$ for
$f\in F_d$. Equations \eq{khCeq9}--\eq{khCeq14} follow as in Step~3.

We must prove that provided the $X_{ac},\ti X_{df}$ are sufficiently
small, the right hand side of \eq{khCeq14} is zero. To do this, we
recall the notation of \S\ref{khC1}, fixing $m=1$. For $a\in A$ and
$d\in D$ we have
\e
\begin{aligned}
(\pd X_a,\bs f_a\vert_{\pd X_a},\bs G_a\vert_{\pd
X_a})&=\ts\coprod_{b=1}^{n_a^1}(X_{ab}^1,\bs f_{ab}^1,\bs
G_{ab}^1),\\
(\ti X_d,\bs{\ti f}_d,\bs{\ti G}_d)&=\ts\coprod_{e=1}^{\ti
n_d^1}(X_{de}^1,\bs f_{de}^1,\bs G_{de}^1),\qquad\text{and}\\
(\ti X_d/\Ga_d,\bs\pi_*(\bs{\ti f}_d),\bs\pi_*(\bs{\ti
G}_d))&=\ts\coprod_{e=-1}^{-\hat n_d^1}(X_{de}^1,\bs f_{de}^1,\bs
G_{de}^1).
\end{aligned}
\label{khCeq28}
\e
Write $(X^1,\bs f^1,\bs G^1)$ for the disjoint union of the triples
on the left hand sides of equations \eq{khCeq28} over all $a\in A$
and $d\in D$, and $(V^{1,i},\ldots,\psi^{1,i})$ for $i\in I^1$ for
the Kuranishi neighbourhoods in~$\bs G^1$.

Write $\bs T^1$ for the tent function for $(X^1,\bs f^1,\bs G^1)$
which is the disjoint union of the tent functions $\bs T_a\vert_{\pd
X_a}$, $\bs{\ti T}_d$ and $\bs\pi_*(\bs{\ti T}_d)$ on the triples in
\eq{khCeq28} chosen in Step 1. Write $(X^1_c,\bs f^1_c,\bs G^1_c)$
for $c\in C^1$ for the pieces into which $(X^1,\bs f^1,\bs G^1)$ is
`cut' by $\bs T^1$, as in \eq{khAeq29}. These $(X^1_c,\bs f^1_c,\bs
G^1_c)$ consist of triples $(\ti X_{df},\bs{\ti f}_{df},\bs{\ti
G}_{df})$ and $(\ti X_{df}/\Stab_{\Ga_d}(f),\bs\pi_*(\bs{\ti
f}_{df}),\bs\pi_*(\bs{\ti G}_{df}))$ which occur in \eq{khCeq14},
and also the decomposition of $(\pd X_a,\bs f_a\vert_{\pd X_a},\bs
G_a\vert_{\pd X_a})$ using $\bs T_a\vert_{\pd X_a}$ for $a\in A$.
Write $\bs\pi:X^1_c\ra X^1$ for the natural projection.

Now $\ddot P^1$ is the disjoint union of all $(a,b),(d,e)$ occurring
in \eq{khCeq28}. Hence
\e
\ts\coprod_{(s,t)\in\ddot P^1}X_{st}^1=\coprod_{a\in A}\pd X_a
\amalg \coprod_{d\in D}\ti X_d\amalg\coprod_{d\in D}\ti
X_d/\Ga_d=X^1.
\label{khCeq29}
\e
Also $\approx_{\sst X}^1$ is an equivalence relation on \eq{khCeq29}
with $\ddot X^1=X^1/\approx_{\sst X}^1$, and projection $\pi_{\ddot
X^1}:X^1\ra\ddot X^1$, which is continuous, surjective, and globally
finite, with $\ddot X^1$ a compact topological space. We shall prove
Step 3 by working locally on $\ddot X{}^1$. Fix some $\ddot
p\in\ddot X{}^1$. As $\pi_{\ddot X{}^1}$ is surjective and globally
finite, $\pi_{\ddot X{}^1}^{-1}(\ddot p)$ is finitely many points
$p_{(1)},\ldots,p_{(n)}$ in~$X^1$.

We now use the ideas of Step 2. For each $j=1,\ldots,n$, there
exists an open neighbourhood $U_{(j)}$ of $p_{(j)}$ in $X^1$ such
that as a Kuranishi space we may write $U_{(j)}\cong\acute
U_{(j)}/\Ga_{(j)}$, where $\acute U_{(j)}$ has {\it trivial
stabilizers}, and $\Ga_{(j)}=\Stab_{X^1}(p_{(j)})$ is a finite group
acting on $\acute U_{(j)}$ by strong diffeomorphisms, fixing the
pullback $\acute p_{(j)}$ of $p_{(j)}$ in $\acute U_{(j)}$. Making
the $U_{(j)}$ smaller if necessary, we suppose
$U_{(1)},\ldots,U_{(n)}$ are {\it disjoint}. Similarly, if $i\in\N$
with $p_{(j)}\in\Im\psi^{1,i}$, so that
$p_{(j)}=\psi^{1,i}(v^i_{(j)})$ for some unique $v^i_{(j)}\in
V^{1,i}$, we choose disjoint open neighbourhoods $W_{(j)}^i$ of
$\ddot v{}^i$ in $V^{1,i}$ for $j=1,\ldots,n$ with
$W_{(j)}^i\cong\acute W_{(j)}^i/\Ga_{(j)}$, where $\acute W_{(j)}^i$
is a {\it manifold}, and $\Ga_{(j)}$ is as before,
since~$\Stab_{X^1}(p_{(j)})\cong \Stab_{V^{1,i}}(\ddot
v{}^i_{(j)})$.

For $j=1,\ldots,n$, define $(X_{(j)},\bs f_{(j)},\bs G_{(j)})$ to be
the disjoint union of those components $(X^1_c,\bs f^1_c,\bs G^1_c)$
for $c\in C^1$ with $p_{(j)}\in\bs\pi(X^1_c)$. There is exactly one
such component if $p_{(j)}$ lies in the interior of
$\bs\pi(X_{(j)})$, but if $p_{(j)}$ lies in a codimension $k$ corner
for $k\ge 1$, there may be more than one. Write
$(V_{(j)}^i,\ldots,\psi_{(j)}^i)$ for $i\in I_{(j)}$ for the
Kuranishi neighbourhoods in $\bs G_{(j)}$, and $\pi:V_{(j)}^i\ra
V^{1,i}$ for the natural projections.

Suppose now that the $\bs T_a,\bs{\ti T}_d$ in Step 1 have been
chosen such that $\bs\pi(X_{(j)})\subseteq U_{(j)}$ and
$\pi(V^i_{(j)})\subseteq W_{(j)}^i$ for all $j=1,\ldots,n$ and $i\in
I_{(j)}$. This is possible as Step 1 allows us to cut $X^1,V^{1,i}$
into `arbitrarily small pieces', and so to make $\bs\pi(X_{(j)})$
and $\pi(V^i_{(j)})$ arbitrarily small. Then the proof of Step 2
gives a triple $(\acute X_{(j)},\bs{\acute f}_{(j)},\bs{\acute
G}_{(j)})$, an action of $\Ga_{(j)}$ on $(\acute X_{(j)},\bs{\acute
f}_{(j)},\bs{\acute G}_{(j)})$ by isomorphisms, and an isomorphism
$(\acute X_{(j)}/\Ga_{(j)}, \bs\pi_*(\bs{\acute
f}_{(j)}),\bs\pi_*(\bs{\acute G}_{(j)}))\cong \ab(X_{(j)},\ab\bs
f_{(j)},\ab\bs G_{(j)})$, for each~$j$.

Making the pieces smaller if necessary, we see that $(\acute
X_{(j)},\bs{\acute f}_{(j)},\bs{\acute G}_{(j)})$ is {\it unique up
to isomorphism}, that is, there are no arbitrary choices in the
construction of $(\acute X_{(j)},\bs{\acute f}_{(j)},\bs{\acute
G}_{(j)})$ from $(X_{(j)},\bs f_{(j)},\bs G_{(j)})$ and $p_{(j)}$.
We can think of $(\acute X_{(j)},\bs{\acute f}_{(j)},\ab\bs{\acute
G}_{(j)})$ as like a kind of `Kuranishi space universal cover' of
$(X_{(j)},\bs f_{(j)},\bs G_{(j)})$, by analogy with orbifold
universal covers in~\S\ref{khC1}.

We now come to one of the crucial steps in the proof. We claim that
{\it the triples $(\acute X_{(j)},\bs{\acute f}_{(j)},\bs{\acute
G}_{(j)})$ for $j=1,\ldots,n$ are all isomorphic}, through either
orientation-preserving or orientation-reversing isomorphisms. To
prove this, note that the equivalence relation $\approx_{\sst X}^1$
on $X^1$ is generated by (possibly orientation-reversing)
isomorphisms of components of $(X^1,\bs f^1,\bs G^1)$, and finite
group quotients between components of $(X^1,\bs f^1,\bs G^1)$. Since
$\{p_{(1)},\ldots,p_{(n)}\}$ is an $\approx_{\sst X}^1$-equivalence
class, any two $p_{(j)},p_{(k)}$ are connected by a finite chain of
such isomorphisms or finite group quotients. Thus, to show the
$(\acute X_{(j)},\bs{\acute f}_{(j)},\bs{\acute G}_{(j)})$ are all
isomorphic for $j=1,\ldots,n$, it is enough to prove that if
$p_{(j)}$ and $p_{(k)}$ are identified by an isomorphism or finite
group quotient of their components, then $(\acute X_{(j)},\bs{\acute
f}_{(j)},\bs{\acute G}_{(j)})$ and $(\acute X_{(k)},\bs{\acute
f}_{(k)},\bs{\acute G}_{(k)})$ are isomorphic.

The case of isomorphisms is easy: if $(\bs a,\bs b):(X^1_{(j)},\bs
f^1_{(j)},\bs G^1_{(j)})\ra(X^1_{(k)},\bs f^1_{(k)},\ab\bs
G^1_{(k)})$ is a (possibly orientation-reversing) isomorphism
between the components of $(X^1,\bs f^1,\bs G^1)$ containing
$p_{(j)},p_{(k)}$ with $a(p_{(j)})=p_{(k)}$, then condition (i) in
\S\ref{khC1} implies that $\smash{\bs T^1\vert_{X^1_{(j)}}=\bs
b^*(\bs T^1\vert_{X^1_{(k)}})}$, so $(\bs a,\bs b)$ lifts to a
(possibly orientation-reversing) isomorphism between $(\acute
X_{(j)},\bs{\acute f}_{(j)},\bs{\acute G}_{(j)})$ and~$(\acute
X_{(k)},\bs{\acute f}_{(k)},\bs{\acute G}_{(k)})$.

Suppose $p_{(j)}$ and $p_{(k)}$ are related by one of the finite
group quotients used to define $\approx_{\sst X}^1$. Then in
\S\ref{khC1}(ii) with $m=1$, we have $(X^1_{(j)},\bs f^1_{(j)},\bs
G^1_{(j)})=(X_{de}^1,\ab\bs f_{de}^1,\ab\bs G_{de}^1)$, and
$\Ga=\Stab_{\Ga_d}(e)$, and $(X^1_{(k)},\bs f^1_{(k)},\bs
G^1_{(k)})= (X^1_{(j)}/\Ga,\ab\bs\pi_*(\bs f^1_{(j)}),\ab
\bs\pi_*(\bs G^1_{(j)}))$ with $p_{(k)}=p_{(j)}\Ga$, where
$(\bs\pi,\bs p):(X^1_{(j)},\bs f^1_{(j)},\bs G^1_{(j)})\ra
(X^1_{(k)},\bs f^1_{(k)},\bs G^1_{(k)})$ is the projection to the
quotient. Part (ii) of \S\ref{khC1} gives~$\bs
T^1\vert_{X^1_{(j)}}=\bs p^*(\bs T^1\vert_{X^1_{(k)}})$.

It is now easy to see that the action of $\Ga$ on $(X^1_{(j)},\bs
f^1_{(j)},\bs G^1_{(j)})$ lifts to $(X_{(j)},\ab\bs f_{(j)},\ab\bs
G_{(j)})$, and $(X_{(k)},\bs f_{(k)},\bs G_{(k)})\cong
(X_{(j)}/\Ga,\ab\bs\pi_*(\bs f_{(j)}),\ab \bs\pi_*(\bs G_{(j)}))$.
Using the obvious shorthand notation, we now have isomorphisms
\e
\begin{gathered}
(X_{(j)},\bs f_{(j)},\bs G_{(j)})\cong(\acute X_{(j)},\bs{\acute
f}_{(j)},\bs{\acute G}_{(j)})/\Ga_{(j)},\qquad
\Ga_{(j)}=\Stab(p_{(j)}),\\
(X_{(k)},\bs f_{(k)},\bs G_{(k)})\cong(\acute X_{(k)},\bs{\acute
f}_{(k)},\bs{\acute G}_{(k)})/\Ga_{(k)},\qquad
\Ga_{(k)}=\Stab(p_{(k)}), \\
\text{and}\qquad (X_{(k)},\bs f_{(k)},\bs G_{(k)})\cong (X_{(j)},\bs
f_{(j)},\bs G_{(j)})/\Ga.
\end{gathered}
\label{khCeq30}
\e

As the $U_{(j)}$ are disjoint, $p_{(j)}$ is the only point in its
$\approx_{\sst X}^1$-equivalence class in $X_{(j)}$, so the action
of $\Ga$ on $X_{(j)}$ must fix $p_{(j)}$. The equalities
$X_{(k)}=X_{(j)}/\Ga$, $p_{(k)}=p_{(j)}\Ga$ now imply that
$\Stab_{X^1}(p_{(j)})$ is a normal subgroup of
$\Stab_{X^1}(p_{(k)})$ with $\Ga\cong\Stab_{X^1}(p_{(k)})/
\Stab_{X^1}(p_{(j)})$. So by \eq{khCeq30} we see that $\Ga_{(j)}$ is
normal in $\Ga_{(k)}$ with $\Ga\cong\Ga_{(k)}/\Ga_{(j)}$.
Considering the group quotients in \eq{khCeq30}, and using
uniqueness up to isomorphism of $(\acute X_{(j)},\bs{\acute
f}_{(j)},\bs{\acute G}_{(j)})$, $(\acute X_{(k)},\bs{\acute
f}_{(k)},\bs{\acute G}_{(k)})$, we see that $(\acute
X_{(k)},\bs{\acute f}_{(k)},\bs{\acute
G}_{(k)})/\Ga_{(j)}\cong(X_{(j)},\bs f_{(j)},\bs G_{(j)})$, and
deduce that $(\acute X_{(j)},\bs{\acute f}_{(j)},\bs{\acute
G}_{(j)})\cong (\acute X_{(k)},\bs{\acute f}_{(k)},\bs{\acute
G}_{(k)})$, in this case with an isomorphism-preserving isomorphism.

This proves that the $(\acute X_{(j)},\bs{\acute f}_{(j)},\bs{\acute
G}_{(j)})$ are all isomorphic for all $j$, possibly with orientation
changes. Choose one and write it as $(\acute X,\bs{\acute
f},\bs{\acute G})$, and then identify $(\acute X_{(j)},\bs{\acute
f}_{(j)},\bs{\acute G}_{(j)})$ with $(\ep_{(j)}\acute X,\bs{\acute
f},\bs{\acute G})$ for $j=1,\ldots,n$, where $\ep_{(j)}=\pm 1$ is to
change orientations if necessary.

We are now ready to prove Step 3. We first derive a modified version
of equation \eq{khCeq14}. Fix $a\in A$, and consider the term
$\sum_{c\in C_a}\pd[X_{ac},\bs f_{ac},\bs G_{ac}]$ on the left hand
side of \eq{khCeq14}. Since $X_a$ is `cut into pieces' $X_{ac}$ for
$c\in C_a$, we can divide the boundaries $\pd X_{ac}$ into
`external' components, which are part of $\pd X_a$, and `internal'
components, which are new boundary components created by cutting the
interior of $X_a$ along some hypersurface. When we cut $X_a$ along a
hypersurface, the two sides of the cut yield identical boundary
components with opposite orientations. Therefore the contributions
of the `internal' components to $\sum_{c\in C_a}\pd[X_{ac},\bs
f_{ac},\bs G_{ac}]$ cancel out.

This leaves only the `external' components, which are the pieces
into which $(\pd X_a,\bs f_a\vert_{\pd X_a},\bs G_a\vert_{\pd X_a})$
is `cut' using the tent function $\bs T_a\vert_{\pd X_a}$. Write
these as $(X_{ac'}^\pd,\bs f_{ac'}^\pd,\bs G_{ac'}^\pd)$ for $c'\in
C_a^\pd$. Then in $\widetilde{KC}_{k-1}(Y;R)$ we have
\e
\ts\sum_{c\in C_a}\pd[X_{ac},\bs f_{ac},\bs G_{ac}]=
\sum_{c'\in C_a^\pd}[X_{ac'}^\pd,\bs f_{ac'}^\pd,\bs G_{ac'}^\pd].
\label{khCeq31}
\e
Substituting \eq{khCeq31} into the left hand of \eq{khCeq14} yields
\e
\begin{split}
&\ts\sum_{a\in A}\sum_{c'\in C_a^\pd}\rho_a[X_{ac'}^\pd,\bs
f_{ac'}^\pd,\bs G_{ac'}^\pd]=\\
&\sum_{d\in D}\eta_d\Bigl(\begin{aligned}[t]&\ts\sum_{f\Ga_d\in
F_d/\Ga_d}\bigl[\ti X_{df}/\Stab_{\Ga_d}(f),\bs\pi_*(\bs{\ti
f}_{df}),\bs\pi_*(\bs{\ti
G}_{df})\bigr]\\
&-\md{\Ga_d}^{-1}\ts\sum_{f\in F_d}[\ti X_{df},\bs{\ti
f}_{df},\bs{\ti G}_{df}]\Bigr).
\end{aligned}
\end{split}
\label{khCeq32}
\e
We must prove that both sides of \eq{khCeq32} are zero. Note that as
$X_a$ for $a\in A$ have trivial stabilizers, so does $\pd X_a$, and
so all $X_{ac'}^\pd$ on the left hand side of \eq{khCeq32} have {\it
trivial stabilizers}.

Our approach is to work locally in $\ddot X^1$. Now each term
$X_{ac'}^\pd,\ti X_{df}/\Stab_{\Ga_d}(f)$, $\ti X_{df}$ in
\eq{khCeq32} has a natural projection $\bs\pi$ to $X^1$, and so to
$\ddot X^1$ via $\pi_{\ddot X^1}:X^1\ra\ddot X^1$. (This did not
hold for \eq{khCeq14}, as the `internal' boundary components of $\pd
X_a$ do not project to $X^1$.) Going through the proof of
\eq{khCeq32}, one can show that it also holds locally in $\ddot
X^1$. Thus, if we fix $\ddot p\in\ddot X^1$ and restrict to those
terms in \eq{khCeq32} whose image in $\ddot X^1$ contains $\ddot p$,
the resulting equation still holds. We will show that each side of
this resulting equation is zero.

For $a\in A$, the sum of terms in $\sum_{c'\in C_a^\pd}[X_{ac'}^\pd,
\bs f_{ac'}^\pd,\bs G_{ac'}^\pd]$ whose image in $\ddot X^1$ contain
$\ddot p$ is
\e
\ts\sum_{j=1,\ldots,n:\;p_{(j)}\in\pd X_a}[X_{(j)},\bs f_{(j)},\bs
G_{(j)}].
\label{khCeq33}
\e
Now if $p_{(j)}\in\pd X_a$ then $\Stab(p_{(j)})=\{1\}$ as $\pd X_a$
has trivial stabilizers, so $\Ga_{(j)}=\{1\}$ by \eq{khCeq30},
giving $(X_{(j)},\bs f_{(j)},\bs G_{(j)})\cong\ep_{(j)}(\acute
X_{(j)},\bs{\acute f}_{(j)},\bs{\acute G}_{(j)})$. So another way to
write \eq{khCeq33} is
\e
\ts\sum_{j=1,\ldots,n:\;p_{(j)}\in\pd X_a}\ep_{(j)}[\acute
X,\bs{\acute f},\bs{\acute G}].
\label{khCeq34}
\e

Similarly, for $d\in D$ the sum of terms in $\sum_{f\in F_d}[\ti
X_{df},\bs{\ti f}_{df},\bs{\ti G}_{df}]$ whose image in $\ddot X^1$
contains $\ddot p$ is
\e
\ts\sum_{j=1,\ldots,n:\;p_{(j)}\in\ti X_d}[X_{(j)},\bs f_{(j)},\bs
G_{(j)}].
\label{khCeq35}
\e
The sum of terms in $\sum_{f\Ga_d\in F_d/\Ga_d}\bigl[\ti
X_{df}/\Stab_{\Ga_d}(f),\bs\pi_*(\bs{\ti f}_{df}),\bs\pi_*(\bs{\ti
G}_{df})\bigr]$ whose images contain $\ddot p$ is the quotient of
\eq{khCeq35} by $\Ga_d$ acting on $\coprod_{j=1,\ldots,n:
\;p_{(j)}\in\ti X_d}X_{(j)}$. This action of $\Ga_d$ factors through
an action on $\{j=1,\ldots,n: p_{(j)}\in\ti X_d\}$. Hence we may
write the sum of terms in $\bigl[\ti X_{df}/\Stab_{\Ga_d}(f),
\bs\pi_*(\bs{\ti f}_{df}),\bs\pi_*(\bs{\ti G}_{df})\bigr]$ whose
images contain $\ddot p$ as
\e
\ts\sum_{j\Ga_d\in\{j=1,\ldots,n:\;p_{(j)}\in\ti X_d\}/\Ga_d}
\bigl[X_{(j)}/\Stab_{\Ga_d}(j),\bs\pi_*(\bs f_{(j)}),\bs\pi_*(\bs
G_{(j)})\bigr].
\label{khCeq36}
\e
Here $\bigl(X_{(j)}/\Stab_{\Ga_d}(j),\bs\pi_*(\bs
f_{(j)}),\bs\pi_*(\bs G_{(j)})\bigr)$ is one of the quotients
$(X^1_{(k)},\bs f^1_{(k)},\ab\bs G^1_{(k)})=
(X^1_{(j)}/\Ga,\ab\bs\pi_*(\bs f^1_{(j)}),\ab \bs\pi_*(\bs
G^1_{(j)}))$ above.

Combining \eq{khCeq32}--\eq{khCeq36}, replacing $[X_{(j)},\bs
f_{(j)},\bs G_{(j)}]$ by $\ep_{(j)}\bigl[\acute X/\Ga_{(j)},
\bs\pi_*(\bs{\acute f}),\ab\bs\pi_*(\bs{\acute G})\bigr]$, and
restricting to terms in \eq{khCeq32} whose images contain $\ddot p$
yields
\e
\begin{split}
&\sum_{a\in A}\sum_{j=1,\ldots,n:\;p_{(j)}\in\pd
X_a}\rho_a\ep_{(j)}[\acute X,\bs{\acute f},\bs{\acute G}]=\\
&\sum_{d\in D}\eta_d\Bigl(\begin{aligned}[t]&\ts
\sum\limits_{j\Ga_d\in\{j=1,\ldots,n:\;p_{(j)}\in\ti X_d\}/\Ga_d
\!\!\!\!\!\!\!\!\!\!\!\!\!\!\!\! }
\ep_{(j)}\bigl[\acute X/\Stab_{\Ga_d}(j)\ltimes\Ga_{(j)},\bs\pi_*(\bs{\acute
f}),\bs\pi_*(\bs{\acute G})\bigr] \\
&-\md{\Ga_d}^{-1}\ts\sum\limits_{j=1,\ldots,n:\;p_{(j)}\in\ti
X_d\!\!\!\!\!\!\!\!}\ep_{(j)}\bigl[\acute
X/\Ga_{(j)},\bs\pi_*(\bs{\acute f}),\bs\pi_*(\bs{\acute
G})\bigr]\Bigr).
\end{aligned}
\end{split}
\label{khCeq37}
\e
Equation \eq{khCeq37} holds in $\widetilde{KC}_{k-1}(Y;R)$.
Projecting to $KC_{k-1}(Y;R)$, each term $(\cdots)$ on the right
hand side vanishes by relation Definition \ref{kh4def2}(iv), which
holds in $KC_{k-1}(Y;R)$ but not in $\widetilde{KC}_{k-1}(Y;R)$.
Thus
\e
\ts\sum_{a\in A}\sum_{j=1,\ldots,n:\;p_{(j)}\in\pd
X_a}\rho_a\ep_{(j)}[\acute X,\bs{\acute f},\bs{\acute G}]=0
\quad\text{in $KC_{k-1}(Y;R)$.}
\label{khCeq38}
\e

There are now two possibilities: either there exists an
orientation-reversing isomorphism $(\bs a,\bs b)$ of $(\acute
X,\bs{\acute f},\bs{\acute G})$, so that $[\acute X,\bs{\acute
f},\bs{\acute G}]=0$ in both $KC_{k-1}(Y;R)$ and
$\widetilde{KC}_{k-1}(Y;R)$ by Definition \ref{kh4def2}(ii), or the
coefficient of $[\acute X,\bs{\acute f},\bs{\acute G}]$ is zero in
\eq{khCeq38}. (More complicated scenarios of having to use several
applications of relations in $KC_{k-1}(Y;R)$ to prove $[\acute
X,\bs{\acute f},\bs{\acute G}]=0$ can be excluded, because of the
general form of \eq{khCeq37}.) In both cases, \eq{khCeq38} also
holds in $\widetilde{KC}_{k-1}(Y;R)$. Hence both sides of
\eq{khCeq37} are zero in~$\widetilde{KC}_{k-1}(Y;R)$.

We have shown that for a given $\ddot p\in\ddot X^1$, if the $\bs
T_a,\bs{\ti T}_d$ are chosen to make the pieces $X_{ac},\ti X_{df}$
sufficiently small near the preimage of $\ddot p$, then on each side
of \eq{khCeq32}, the sum of terms whose image in $\ddot X^1$
contains $\ddot p$ is zero. On $\ddot X^1$, we can write the
`sufficiently small' condition as requiring that
$\bs\pi(X_{ac}),\bs\pi(\ti X_{df})\subseteq\ddot U_{\ddot p}$ for
any $X_{ac},\ti X_{df}$ whose image contains $\ddot p$, where $\ddot
U_{\ddot p}$ is the intersection in $\ddot X^1$ of the images of
$U_{(j)}$ for $j=1,\ldots,n$, so that $\ddot U_{\ddot p}$ is an open
neighbourhood of $\ddot p$ in $\ddot X^1$. Such open sets $\ddot
U_{\ddot p}$ for all $\ddot p$ in $\ddot X^1$ form an open cover of
$\ddot X^1$, which is compact. So we can choose a finite subset
$\ddot p^1,\ldots,\ddot p^l$ with~$\ddot X^1=\ddot U_{\ddot
p^1}\cup\cdots\cup\ddot U_{\ddot p^l}$.

We can then choose the $\bs T_a,\bs{\ti T}_d$ to make the
$X_{ac'}^\pd,\ti X_{df}$ small enough that each
$\bs\pi(X_{ac'}^\pd)$ or $\bs\pi(\ti X_{df})$ lies in $\ddot
U_{\ddot p^i}$ for some $i=1,\ldots,l$, and with the corresponding
smallness conditions in $V^{1,i}$, so that each $(X_{ac'}^\pd,\bs
f_{ac'}^\pd,\bs G_{ac'}^\pd)$ and $(\ti X_{df},\bs{\ti f}_{df},
\bs{\ti G}_{df})$ can be written naturally as a quotient of a triple
with trivial stabilizers by a finite group, which is the stabilizer
group of a preimage of some $\ddot p^i$ in~$X^1$.

This method does not necessarily ensure that every
$\bs\pi(X_{ac'}^\pd)$ or $\bs\pi(\ti X_{df})$ contains some $\ddot
p^i$, only that it lies in some $\ddot U_{\ddot p^i}$. But this does
not matter, as we can rewrite the argument above, replacing
restricting to terms in \eq{khCeq32} whose image in $\ddot X^1$
contains $\ddot p^i$ for some $i=1,\ldots,l$ by restricting to terms
in \eq{khCeq32} whose image in $\ddot X^1$ lies in $\ddot U_{\ddot
p^i}$ for some~$i=1,\ldots,l$.

Thus, with the $\bs T_a,\bs{\ti T}_d$ chosen in this way, equation
\eq{khCeq32} has the following properties. Firstly, for every term
in \eq{khCeq32}, the projection of $X_{ac'}^\pd$ or $\ti X_{df}/
\Stab_{\Ga_d}(f)$ or $\ti X_{df}$ to $\ddot X^1$ lies in $\ddot
U_{\ddot p^i}$ for at least one $i=1,\ldots,l$. Secondly, for each
$i=1,\ldots,l$, if we restrict to the terms in \eq{khCeq32} whose
projection to $\ddot X^1$ lies in $\ddot U_{\ddot p^i}$, the
resulting equation is true, and also, both sides are zero. Thirdly,
an easy generalization of the second property shows that for any
subset $S$ of $\{1,\ldots,l\}$, if we restrict to the terms in
\eq{khCeq32} whose projection to $\ddot X^1$ lies in $\ddot U_{\ddot
p^i}$ for all $i\in S$, the resulting equation is true, and also,
both sides are zero.

By the first property, we see that summing over all $\emptyset\ne
S\subseteq\{1,\ldots,l\}$ of the restriction to the terms in
\eq{khCeq32} whose projection to $\ddot X^1$ lies in $\ddot U_{\ddot
p^i}$ for all $i\in S$, multiplied by $(-1)^{\md{S}-1}$, yields
\eq{khCeq32}. But the third property says that each of these
restrictions has both sides zero, so \eq{khCeq32} has both sides
zero. Therefore both sides of \eq{khCeq14} are zero. This proves
\eq{khCeq15}, and the first part of Step~3.

It remains to show that if the $X_{ac}$ are chosen sufficiently
small, for all $a\in A$, $c\in C_a$ and $m\ge 2$, each component of
$(\pd^mX_{ac},\bs f_{ac} \vert_{\pd^mX_{ac}},\bs
G_{ac}\vert_{\pd^mX_{ac}})$ occurs as the intersection of $m$ {\it
distinct\/} components of $(\pd X_{ac}, \bs f_{ac}\vert_{\pd
X_{ac}},\bs G_{ac}\vert_{\pd X_{ac}})$. We use notation
$(X_{acd}^m,\bs f_{acd}^m,\bs G_{acd}^m)$ for $d=1,\ldots,n_{ac}^m$,
$V_{acd}^{m,i-m}\subseteq\pd^mV_{ac}^i$ and $(v,B_1,\ldots,B_m)$ for
points in $\pd^mV_{ac}^i$ as in Step~3.

Recall that $(X_{ac},\bs f_{ac},\bs G_{ac})$ is a piece in the
decomposition of $(X_a,\bs f_a,\bs G_a)$ using $\bs T_a=(T_a^i:i\in
I_a)$, where $T_a^i:\dot V_a^i\ra F\bigl([1,\iy)\bigr)$ is a tent
function on the {\it manifold\/} $\dot V_a^i$. By Definition
\ref{khAdef1}, each $v\in\dot V_a^i$ has an open neighbourhood $U$
in $\dot V_a^i$ such that $T_a^i\vert_U$ may be written
$T_a^i(u)=\bigl\{t_p(u):p=1,\ldots,N$, $u\in U_p\bigr\}$, for some
open sets $U_1,\ldots,U_N$ in $U$. Considering the constructions of
tent functions in \S\ref{khA}, we see that how small this open set
$U$ needs to be can be made to depend only on $\dot V_a^i$ and other
data in $(X_a,\bs f_a,\bs G_a)$, not on the choices in $\bs T_a$
which make the $(X_{ac},\bs f_{ac},\bs G_{ac})$ small.

Therefore by choosing $\bs T_a$ to make the pieces $(X_{ac},\bs
f_{ac},\bs G_{ac})$ sufficiently small, we can arrange that
$\pi(V_{ac}^i)$ for each $c\in C_a$ is contained in some such open
set $U$ in $\dot V_a^i$. Then $V_{ac}^i$ is a piece in the
decomposition of $U$ using the tent function $T_a^i\vert_U$ which
has the {\it global\/} definition $T_a^i(u)=\bigl\{t_p(u):
p=1,\ldots,N$, $u\in U_p\bigr\}$. Consider the boundary $\pd
V_{ac}^i$. It is a disjoint union of two types of components: (a)
the lift of part of $\pd V_a^i$, a piece in the decomposition of
$\pd\dot V_a^i$ using $T_a^i\vert_{\pd\dot V_a^i}$, and (b) part of
the hypersurface $t_p(u)=t_q(u)$ in $U$, for $p,q=1,\ldots,N$, $p\ne
q$. To specify the order of $p,q$ we require that $t_p-t_q$ should
be increasing in the direction of the outward normal to the boundary
component of~$X_{ac}$.

Let $d=1,\ldots,n_{ac}^1$ and $i-1\in I_{acd}^1$, so that
$\emptyset\ne V_{abd}^{1,i-1}\subseteq\pd V_{ac}^i$. Then for each
$(v,B)\in V_{abd}^{1,i-1}$, the classification of $B$ into either
(a) or (b) with $(p,q)$ depends only on $d$. That is, if
$(v,B),(v',B')$ lie in $V_{abd}^{1,i-1}$ then $B$ is a hypersurface
$t_p=t_q$ if and only if $B'$ is a hypersurface $t_p=t_q$, with the
same $(p,q)$. This holds even though $V_{abd}^{1,i-1}$ may not be
connected, because of the global definition of the $T_a^j$ near
$\pi(V_{ac}^j)$, as $(X_{acd}^1,\bs f_{acd}^1,\bs G_{acd}^1)$ is
connected, and since the way in which we make $T_a^i,T_a^j$
compatible under coordinate changes $(\phi_a^{ij},\hat\phi_a^{ij})$
ensures that the hypersurface $t_p=t_q$ in $U^j\cap\dot
V^{ij}_a\subseteq\dot V^j_a$ maps to a hypersurface $t_{p'}=t_{q'}$
in $U^i\subseteq\dot V^i_a$ under $\phi_a^{ij}$, for $p',q'$
depending only on~$p,q,i,j$.

Now let $m\ge 2$ and $d=1,\ldots,n_{ac}^m$, so that $(X_{acd}^m,\bs
f_{acd}^m,\bs G_{acd}^m)$ is a component of $(\pd^mX_{ac},\bs
f_{ac}\vert_{\pd^mX_{ac}},\bs G_{ac}\vert_{\pd^mX_{ac}})$. Suppose
$i-m\in I_{acd}^m$, so that $\emptyset\ne
V_{acd}^{m,i-m}\subseteq\pd^mV_{ac}^i$, let $(v,B_1,\ldots,B_m)\in
V_{acd}^{m,i-m}$, and let $(v,B_l)\in V_{acd_l}^{1,i-1}$ for
$l=1,\ldots,m$ and $d_l=1,\ldots, n_{ac}^1$. We must show that
$d_1,\ldots,d_m$ are {\it distinct}. Suppose $d_l=d_{l'}$ for $l\ne
l'$. As $V_a^i$ is a manifold, $B_1,\ldots,B_m$ are distinct, so
$B_l\ne B_{l'}$. As above, $B_l,B_{l'}$ are classified into either
(a) or (b) with $(p,q)$, and this classification depends only on
$d_l=d_{l'}$. Thus, if $B_l$ is the hypersurface $t_p=t_q$ near $v$
then so is $B_{l'}$, contradicting $B_l\ne B_{l'}$. So the only
possibility is that $B_l,B_{l'}$ are of type (a), that is, they are
local boundary components of~$V_a^i$.

Since $V_a^i$ is a manifold, we can exclude this possibility by
making the $X_{ac}$ sufficiently small, as it would mean that at
$v\in V_a^i$ we have two distinct local boundary components
$B_l,B_{l'}$ which are connected via some path in $\pd
V_a^i\cap\pi(V_{ac}^i)$, or by some more complicated path involving
coordinate changes $(\phi_a^{ij},\hat\phi_a^{ij})$. This can only
happen if some nontrivial global topology in $V_a^i$ happens in the
image $\pi(V_{ac}^i)$, which is ruled out by making the $V_{ac}^i$
small. Therefore $d_1,\ldots,d_m$ are distinct, as we have to prove.
This completes Step~3.

\subsection{Step 4: lifting to a cycle in $KC_k^\ef(Y;R)$}
\label{khC4}

As in Step 4, we again change notation from $X_{ac}$ back to $X_a$,
and introduce new notation $n_a^m,\ldots,B_{ab}^m$. Also, for $a\in
A$ and $0\le m\le m+l\le M$ define
$\th^{ml}_a:\{1,\ldots,n_a^{m+l}\}\ra \{1,\ldots,n_a^m\}$ by
$\th^{ml}_a(\bar b)=b$ if $X^{m+l}_{a\bar b}\subseteq\pd^l
X_{ab}^m$. Then $\th^{m0}_a$ is the identity, and
\e
\bigl(\pd^lX_{ab}^m,\bs f_{ab}^m\vert_{\pd^lX_{ab}^m},\bs
G_{ab}^m\vert_{\pd^lX_{ab}^m}\bigr)=\ts\coprod\limits_{\bar
b=1,\ldots,n_a^{m+l}:\th^{ml}_a(\bar b)=b\!\!\!\!}\bigl(X_{a\bar
b}^{m+l},\bs f_{a\bar b}^{m+l},\bs G_{a\bar b}^{m+l}\bigr)
\label{khCeq39}
\e
is the decomposition of $\pd^lX_{ab}^m$ given by Lemma~\ref{kh3lem}.

We represent $\al$ by a cycle $\sum_{a\in A}\rho_a[X_a,\bs f_a,\bs
G_a]$ in $\widetilde{KC}_k(Y;R)$ satisfying \eq{khCeq16}, where each
$X_a$ has trivial stabilizers, each $V^i_a$ is a manifold, each
component of $(\pd^mX_a,\bs f_a\vert_{\pd^mX_a},\bs
G_a\vert_{\pd^mX_a})$ occurs as the intersection of $m$ distinct
components of $(\pd X_a, \bs f_a\vert_{\pd X_a},\bs G_a\vert_{\pd
X_a})$, and \eq{khCeq17} holds with each $(\bar a,\bar b)\in R^1$
satisfying (A) or (B) in Step 4. Our goal is to choose finite
families $\ubG_{ao}$ for $o\in O_a$ of effective gauge-fixing data
for $(X_a,\bs f_a)$, such that $\sum_{a\in A}\sum_{o\in
O_a}\rho_a\md{O_a}^{-1}[X_a,\bs f_a,\ab\ubG_{ao}]$ is a effective
Kuranishi cycle in class $\be\in KH_k^\ef(Y;R)$
with~$\Pi_\ef^\Kh(\be)=\al$.

Here is our inductive hypothesis:

\begin{hypo} Let $N\ge 0$ be given. Suppose that for all $m=N,N+1,\ldots,M$
and all $(a,b)\in R^m$, we have chosen the following data:
\begin{itemize}
\setlength{\itemsep}{0pt}
\setlength{\parsep}{0pt}
\item[(i)] A nonempty finite set $O_{ab}^m$;
\item[(ii)] A free left action of the finite group $\Aut(X_{ab}^m,\bs
f_{ab}^m,\bs G_{ab}^m)$ on $O_{ab}^m$;
\item[(iii)] For all $o\in O^m_{ab}$ and $l\ge 0$, a map
\e
\!\!\! f_{abo}^{ml}:\bigl\{\bar b=1,\ldots,n_a^{m+l}:\th^{ml}_a(\bar
b)=b\bigr\}\longra\coprod\nolimits_{\substack{\bar b=1,\ldots,n_a^{m+l}: \\
\th^{ml}_a(\bar b)=b}}O_{\phi^{m+l}(a,\bar b)}^{m+l}
\label{khCeq40}
\e
with $f_{abo}^{ml}(\bar b)\in O_{\phi^{m+l}(a,\bar b)}^{m+l}$ for
all $\bar b$; and
\item[(iv)] For each $o\in O_{ab}^m$, a family of maps
$\ubG_{abo}^{m,i}:E_{ab}^{m,i}\ra\uP$ for $i\in I_{ab}^m$ such that
$\ubG_{abo}^m=\bigl(\bs I_{ab}^m,\bs\eta_{ab}^m,\ubG_{abo}^{m,i}:
i\in I_{ab}^m\bigr)$ is effective gauge-fixing data for
$(X_{ab}^m,\bs f_{ab}^m)$. Thus $\ubG_{abo}^{m,i}$ must satisfy the
injectivity condition Definition~\ref{kh3def17}(b).
\end{itemize}
These should satisfy the following conditions:
\begin{itemize}
\setlength{\itemsep}{0pt}
\setlength{\parsep}{0pt}
\item[(a)] If $o,o'\in O_{ab}^m$ and $\ubG_{abo}^{m,i}\equiv
\ubG_{abo'}^{m,i}$ for all $i\in I_{ab}^m$, then~$o=o'$.
\item[(b)] Let $o\in O_{ab}^m$ and $(\bs u,\bs v)\in\Aut(X_{ab}^m,
\bs f_{ab}^m,\bs G_{ab}^m)$, and set $o'=(\bs u,\bs v)\cdot o$ using
the action in (ii). Then for all $i\in I_{ab}^m$ we have
diffeomorphisms $\hat v^i:E^{m,i}_{ab}\ra E^{m,i}_{ab}$, and we
require that~$\ubG_{abo}^{m,i}\equiv\ubG_{abo'}^{m,i}\ci\hat v^i$.
\item[(c)] Let $(a,b)\in R^m$, $i\in I_{ab}^m$, $o\in O^m_{ab}$ and
$l\ge 0$. Then \eq{khCeq39} implies that
\e
\begin{split}
\pd^lE_{ab}^{m,i}&=\ts\coprod_{\bar
b=1,\ldots,n_a^{m+l}:\th^{ml}_a(\bar b)=b} E^{m+l,i-l}_{a\bar b}\\
&\!\!\!\!\!\!\!\!\!{\buildrel\coprod_{\bar b}\hat b_{a\bar
b}^{m+l,i-l}\over \cong}\ts\coprod_{\bar
b=1,\ldots,n_a^{m+l}:\th^{ml}_a(\bar b)=b}
E^{m+l,i-l}_{\phi^{m+l}(a,\bar b)},
\end{split}
\label{khCeq41}
\e
where in the notation of Step 4 we have an isomorphism $(\bs a,\bs
b)_{a\bar b}^{m+l}:(X_{a\bar b}^{m+l},\bs f_{a\bar b}^{m+l},\ab\bs
G_{a\bar b}^{m+l})\ab\ra(X_{\phi^{m+l}(a,\bar b)}^{m+l},\bs
f_{\phi^{m+l}(a,\bar b)}^{m+l},\bs G_{\phi^{m+l}(a,\bar b)}^{m+l})$
which includes a diffeomorphism $\hat b_{a\bar
b}^{m+l,i-l}:E^{m+l,i-l}_{a\bar b}\ra E^{m+l,i-l}_{\phi^{m+l}(a,\bar
b)}$, and we use these diffeomorphisms to make the identification in
the second line of~\eq{khCeq41}.

We have functions $\ubG_{abo}^{m,i}:E_{ab}^{m,i}\ra\uP$ and
$\ubG^{m+l,i-l}_{\phi^{m+l}(a,\bar b)f_{abo}^{ml}(\bar
b)}:E^{m+l,i-l}_{ \phi^{m+l}(a,\bar b)}\ra\uP$ for each $\bar b$
with $\th^{ml}_a(\bar b)=b$. We require that
\e
\ubG_{abo}^{m,i}\big\vert{}_{\pd^lE_{ab}^{m,i}}\equiv
\ts\coprod_{\bar b=1,\ldots,n_a^{m+l}:\th^{ml}_a(\bar
b)=b}\ubG^{m+l,i-l}_{\phi^{m+l}(a,\bar b)f_{abo}^{ml}(\bar b)}
\label{khCeq42}
\e
under the identification~\eq{khCeq41}.

When $l=0$, $\th^{m0}_a$ is the identity, so the domain of
$f_{abo}^{m0}$ is $\{b\}$, and if $f_{abo}^{m0}(b)=o'$ then
\eq{khCeq41} reduces to $E_{ab}^{m,i}=E_{ab}^{m,i}$ and \eq{khCeq42}
to $\ubG_{abo}^{m,i}\equiv\ubG_{abo'}^{m,i}$. Thus (a) implies
that~$f_{abo}^{m0}(b)=o$.
\end{itemize}

Before giving the two remaining conditions, we define some notation.
Fix $m\ge 0$, and suppose that we have chosen data $O_{ab}^{m'},
f_{abo}^{m'l},\ubG_{abo}^{m',i}$ for all $m'>m$ and all $a,b,o,i,l$
satisfying (i)--(iv) and (a)--(c) above. Then in an inductive step
we will wish to choose $O_{ab}^m,f_{abo}^{ml},\ubG_{abo}^{m,i}$ for
all $a,b,o,i,l$. For $l>0$, consider which maps $f_{ab}^{ml}:\{\bar
b=1,\ldots,n_a^{m+l}:\th^{ml}_a(\bar b)=b\}\ra \coprod_{\bar
b}O_{\phi^{m+l}(a,\bar b)}^{m+l}$ could be $f_{abo}^{ml}$ in
\eq{khCeq40} for some $o\in O_{ab}^m$. We will derive a necessary
condition on such~$f_{ab}^{ml}$.

Suppose that for some $o\in O_{ab}^m$ we have chosen maps
$f_{abo}^{ml}$ for $l\ge 0$ and $\ubG_{abo}^{m,i}$ for $i\in
I_{ab}^m$ satisfying (iii),(iv),(c) above. Let $l>0$. Then by (c)
for $a,b,m$ and (a),(c) for $m'=m+l>m$, we see that $f_{abo}^{ml}$
uniquely determines and is determined by the restrictions
$\ubG_{abo}^{m,i}\big\vert{}_{\pd^lE_{ab}^{m,i}}$ for $i\in
I_{ab}^m$. Since these determine $\ubG_{abo}^{m,i}\big\vert
{}_{\pd^{l'}E_{ab}^{m,i}}$ for $l'>l$, it follows that
$f_{abo}^{ml}$ determines $f_{abo}^{ml'}$ for all~$l'>l$.

Observe that there is a natural action of the symmetric group $S_l$
on $\pd^lE_{ab}^{m,i}$ generalizing the involution
$\si:\pd^2E_{ab}^{m,i}\ra\pd^2E_{ab}^{m,i}$ discussed in Definitions
\ref{kh2def5} and \ref{kh2def6}. The restriction (pullback)
$\ubG_{abo}^{m,i}\big\vert{}_{\pd^lE_{ab}^{m,i}}$ is invariant under
this $S_l$ action on $\pd^lE_{ab}^{m,i}$. This gives rise to a
consistency condition on $f_{abo}^{ml}$. Suppose $f_{ab}^{ml}$ is
some map as in \eq{khCeq40}, with $f_{ab}^{ml}(\bar b)\in
O_{\phi^{m+l}(a,\bar b)}^{m+l}$ for all $\bar b$. Then (c)
determines functions $\pd^lE_{ab}^{m,i}\ra\uP$ for $i\in I_{ab}^m$,
by the r.h.s.\ of \eq{khCeq42} with $f_{ab}^{ml}$ in place of
$f_{abo}^{ml}$. To be of the form $\ubG_{abo}^{m,i}\big\vert
{}_{\pd^lE_{ab}^{m,i}}$, these functions must be $S_l$-invariant.

We say $f_{ab}^{ml}$ is $S_l$-{\it invariant\/} if the functions it
determines are $S_l$-invariant in this way. We say $f_{ab}^{ml}$ is
{\it allowable\/} if it is $S_l$-invariant, and also each map
$f_{ab}^{ml'}$ for $l'>l$ determined by $f_{ab}^{ml}$ is
$S_{l'}$-invariant. Write $D_{ab}^{ml}$ for the set of allowable
maps $f_{ab}^{ml}$. When $l=0$, each $f_{ab}^{m0}$ maps $\{b\}\ra
O_{ab}^m$, and for each $o\in O_{ab}^m$ we have $f_{abo}^{m0}:b\ra
o$, which is allowable. Thus $D_{ab}^{m0}$ is the set of maps
$\{b\}\ra O_{ab}^m$, and $f_{ab}^{m0}\mapsto f_{ab}^{m0}(b)$ induces
a bijection $D_{ab}^{m0}\ra O_{ab}^m$.

When $l'>l\ge 0$, write $\Pi^{mll'}_{ab}:D_{ab}^{ml}\ra
D_{ab}^{ml'}$ for the natural projection taking
$\Pi^{mll'}_{ab}:f_{abo}^{ml}\mapsto f_{abo}^{ml'}$, defined by
restriction from $\pd^lE_{ab}^{m,i}$ to $\pd^{l'}E_{ab}^{m,i}$ as
above.

Suppose $\ti m>m\ge 0$ and $l>\ti l\ge 0$ with $m+l=\ti m+\ti l$,
and that $\hat b=1,\ldots,n_a^{\ti m}$ with $\th^{m(\ti m-m)}_a(\hat
b)=b$. Set $(\ti a,\ti b)=\phi^{\ti m}(a,\hat b)$. Then we have an
isomorphism $(\bs a,\bs b)_{a\hat b}^{\ti m}:(X_{a\hat b}^{\ti
m},\bs f_{a\hat b}^{\ti m},\bs G_{a\hat b}^{\ti m})\ra (X_{\ti a\ti
b}^{\ti m},\bs f_{\ti a\ti b}^{\ti m},\bs G_{\ti a\ti b}^{\ti m})$.
Composing this with the inclusion $X^{\ti m}_{a\hat b}\subseteq
\pd^{\ti m-m}X^m_{ab}$ and applying $\pd^{\ti l}$ to both sides
induces a morphism $(\bs a_{a\hat b}^{\ti m})^{-1}:\pd^{\ti l}X_{\ti
a\ti b}^{\ti m}\ra\pd^lX^m_{ab}$ which is an isomorphism with a
union of components of~$\pd^lX^m_{ab}$.

Since the components of $\pd^{\ti l}X_{\ti a\ti b}^{\ti m}$ are
$X_{\ti a\bar b}^{\ti m+\ti l}$ for $\bar b=1,\ldots,n_a^{\ti m+\ti
l}$ with $\th^{\ti m\ti l}_{\ti a}(\bar b)=\ti b$, and the
components of $\pd^lX^m_{ab}$ are $X_{ab'}^{m+l}$ for
$b'=1,\ldots,n_a^{m+l}$ with $\th^{ml}_a(b')=b$, there is a natural,
injective map $I_{ab\hat b}^{ml\ti m\ti l}:\bigl\{\bar b=1,\ldots,
n_{\ti a}^{\ti m+\ti l}:\th^{\ti m\ti l}_{\ti a}(\bar b)=\ti
b\bigr\}\ra\bigl\{b'=1,\ldots,n_a^{m+l}:\th^{ml}_a(b')=b\bigr\}$
such that $(\bs a_{a\hat b}^{\ti m})^{-1}\bigl(X_{\ti a\bar b}^{\ti
m+\ti l}\bigr)=X_{ab'}^{m+l}$ when $I_{ab\hat b}^{ml\ti m\ti l}(\bar
b)=b'$. Let $f_{ab}^{ml}\in D_{ab}^{ml}$, and define $\ti f_{\ti
a\ti b}^{\ti m\ti l}=f_{ab}^{ml}\ci I_{ab\hat b}^{ml\ti m\ti l}$.
Then the conditions for $f_{ab}^{ml}$ to be allowable restrict to
imply that $\ti f_{\ti a\ti b}^{\ti m\ti l}$ is allowable, so $\ti
f_{\ti a\ti b}^{\ti m\ti l}\in D_{\ti a\ti b}^{\ti m\ti l}$. Thus we
define a projection $P_{ab\hat b}^{ml\ti m\ti l}:D_{ab}^{ml}\ra
D_{\ti a\ti b}^{\ti m\ti l}$ mapping $P_{ab\hat b}^{ml\ti m\ti
l}:f_{ab}^{ml}\mapsto f_{ab}^{ml}\ci I_{ab\hat b}^{ml\ti m\ti l}$.

Here are the final conditions:
\begin{itemize}
\setlength{\itemsep}{0pt}
\setlength{\parsep}{0pt}
\item[(d)] When $l'>l\ge 0$ the projection
$\Pi^{mll'}_{ab}:D_{ab}^{ml}\ra D_{ab}^{ml'}$ is surjective, and the
preimage of each point is the same number $\bmd{D_{ab}^{ml}}/
\bmd{D_{ab}^{ml'}}$ of points.
\item[(e)] For all $m,\ti m,l,\ti l,a,b,\hat b$ as above,
$P_{ab\hat b}^{ml\ti m\ti l}:D_{ab}^{ml}\ra D_{\ti a\ti b}^{\ti m\ti
l}$ is surjective, and the preimage of each point is
$\bmd{D_{ab}^{ml}}/\bmd{D_{\ti a\ti b}^{\ti m\ti l}}$ points.
\end{itemize}
\label{khChyp}
\end{hypo}

Here in (ii), $\Aut(X_{ab}^m,\bs f_{ab}^m,\bs G_{ab}^m)$ is finite
by Theorem \ref{kh3thm2}. We shall prove by reverse induction on $N$
that Hypothesis \ref{khChyp} holds for all $N\ge 0$. For $m>M$ we
have $R^m=\es$, so Hypothesis \ref{khChyp} is vacuous for $N>M$, the
initial step. For the inductive step, we shall suppose Hypothesis
\ref{khChyp} holds for $N=n+1$, and prove it holds for~$N=n\ge 0$.

To do this, we will first fix $(a,b)\in R^n$ and construct the sets
$D_{ab}^{nl}$ of allowable maps $f_{ab}^{nl}$ for $l>0$, and prove
by reverse induction on $l$ that (d),(e) hold. Then setting $l=1$,
the maps $f_{ab}^{n1}$ in $D_{ab}^{n1}$ uniquely determine the
possible values of $\ubG_{abo}^{n,i}\big\vert{}_{\pd E_{ab}^{n,i}}$
for $i\in I_{ab}^n$. We shall choose $O_{ab}^n$ and
$\ubG_{abo}^{n,i}$ for $o\in O_{ab}^n$ and $i\in I_{ab}^n$ with
these boundary values, and prove that (a)--(e) hold.

Since the sets $D_{ab}^{nl}$ of allowable maps $f_{ab}^{nl}$ for
$l>0$ depend only on data $O_{ab}^m$ and $(\ubG_{abo}^{m,i})_{i\in
I_{ab}^m}$ for $m>n$ which we have already chosen, they are already
well-defined without making further choices. When $l>M-n$ we have
$\pd^lX_{ab}^n=\es$, so that $\bigl\{\bar b=1,\ldots,n_a^{n+l}
:\th^{nl}_a(\bar b)=b\bigr\}=\es$, and the only possible map
\eq{khCeq40} is the trivial map $\es\ra\es$, which we write $\es$.
It is allowable, so $D_{ab}^{nl}=\{\es\}$, that is, the one point
set with element $\es$. Thus for $l'>l>M-n$ part (d) is immediate,
and for $l>M-n$ part (e) is trivial as there are no such~$\hat b$.

Suppose by induction that for some $p>0$, (d) holds for $m=n$ and
all $l'>l>p$, and (e) holds for $m=n$ and all $l>p$. This holds for
$p\ge M-n$ from above. We shall prove that (d) holds for all $m=n$
and $l'>l=p$, and (e) holds for $m=n$ and $l=p$. Then by induction
(d),(e) hold for all $l',l>0$. For (d) it is sufficient to prove the
case $l=p$ and $l'=p+1$, as the case $l'>p+1$ follows by composing
projections $\smash{D_{ab}^{np}\ra D_{ab}^{n(p+1)}\ra D_{ab}^{nl'}}$
and using the inductive hypothesis.

For the first part of (d) we must show that the projection
$\Pi^{np(p+1)}_{ab}:D_{ab}^{np}\ra D_{ab}^{n(p+1)}$ is surjective.
Fix $f_{ab}^{n(p+1)}$ in $D_{ab}^{n(p+1)}$. We will construct
$f_{ab}^{np}$ in $D_{ab}^{np}$ which restricts to $f_{ab}^{n(p+1)}$.
Let $\hat b$ lie in the domain of $f_{ab}^{np}$, that is, $\hat
b=1,\ldots,n_a^{n+p}$ with $\th_a^{np}(\hat b)=b$. Set
$(a',b')=\phi^{n+p}(a,\hat b)$. We shall construct the set of
possible values of $f_{ab}^{np}(\hat b)$ in $O^{n+p}_{a'b'}$.
Consider the diagram:
\e
\begin{gathered}
\xymatrix@C=85pt@R=15pt{ D^{n(p+1)}_{ab} \ar[d]^(0.55){P_{ab\hat
b}^{n(p+1)(n+p)1}} & D^{np}_{ab} \ar[l]^(0.4){\Pi^{np(p+1)}_{ab}}
\ar[d]_(0.55){P_{ab\hat b}^{np(n+p)0}}
\ar[dr]^{{}\qquad f^{np}_{ab}\mapsto f^{np}_{ab}(\hat b)} \\
D^{(n+p)1}_{a'b'} & D^{(n+p)0}_{a'b'}
\ar[r]^{\cong}_{f^{(n+p)0}_{a'b'}\mapsto f^{(n+p)0}_{a'b'}(b')}
\ar[l]^{\Pi^{(n+p)01}_{a'b'}} & O^{n+p}_{a'b'}.}
\end{gathered}
\label{khCeq43}
\e
It is easy to show from the definitions that \eq{khCeq43} commutes.
Thus a necessary condition for $f_{ab}^{np}\in D_{ab}^{np}$ to
satisfy $\Pi^{np(p+1)}_{ab}(f_{ab}^{np})=f_{ab}^{n(p+1)}$, as in
(d), is that if $f_{ab}^{np}(\hat b)=o'\in O^{n+1}_{a'b'}$ then
$\Pi_{a'b'}^{(n+p)01}\bigl(\{b'\}\mapsto o'\bigr)=P_{ab\hat
b}^{n(p+1)(n+p)1}\bigl(f_{ab}^{n(p+1)}\bigr)$, where $\{b'\}\mapsto
o'$ in $D^{(n+p)0}_{a'b'}$ is the map $\{b'\}\ra O^{n+1}_{a'b'}$
taking~$b'\mapsto o'$.

Since $p>0$, Hypothesis \ref{khChyp}(d) with $m=n+p>n$ implies that
$\Pi_{a'b'}^{(n+p)01}$ is surjective, and the preimage of each point
is $\bmd{D^{(n+p)0}_{a'b'}}/\bmd{D^{(n+p)1}_{a'b'}}$ points. Thus,
for each $\hat b$ in the domain of $f_{ab}^{np}$ there is a nonempty
set $S_{\hat b}\subseteq O^{n+1}_{a'b'}$ of possible values $o'$ for
$f^{np}_{ab}(\hat b)$ making \eq{khCeq43} commute, with
$\bmd{S_{\hat b}}=\bmd{D^{(n+p)0}_{a'b'}}/ \bmd{D^{(n+p)1}_{a'b'}}$.
For $f^{np}_{ab}$ to lie in $D^{np}_{ab}$ it is necessary that
$f^{np}_{ab}(\hat b)\in S_{\hat b}$ for all $\hat b$ in the domain.
Any such map $f^{np}_{ab}$ induces the fixed map $f^{n(p+1)}_{ab}$
in $D^{np}_{ab}$, by construction. As $f^{n(p+1)}_{ab}$ is
allowable, any $f^{nl'}_{ab}$ induced from $f^{np}_{ab}$ for $l'>p$
is induced from $f^{n(p+1)}_{ab}$, and so is $S_{l'}$-invariant.
Therefore such an $f^{np}_{ab}$ is allowable if and only if it is
$S_p$-invariant.

There are natural $S_p$-actions on the domain $\bigl\{\hat
b=1,\ldots,n_a^{n+p}:\th^{np}_a(\hat b)=b\bigr\}$ and range
$\coprod_{\hat b:\th^{np}_a(\hat b)=b}O_{\phi^{m+l}(a,\hat
b)}^{m+l}$ of $f^{np}_{ab}$, and $f^{np}_{ab}$ is $S_p$-invariant if
and only if it is equivariant w.r.t.\ these actions. Here is the
crux of this part of the argument: by the last part of Step 3, each
point of $\pd^pX^n_{ab}\subseteq\pd^{n+p}X_a$ occurs at an ordered
intersection of $n+p$ {\it distinct\/} pieces $X_{ad}^1$ for
$d=1,\ldots,n_a^1$. Thus, to each point in $\pd^pX^n_{ab}$ we can
associate an ordered $(n+p)$-tuple $(d_1,\ldots,d_{n+p})$, where
$d_1,\ldots,d_{n+p}$ are distinct elements of $\{1,\ldots,n_a^1\}$.
On a connected component $X_{a\hat b}^{n+p}\subset\pd^pX^n_{ab}$
this $n+p$-tuple is constant. Hence, $(d_1,\ldots,d_{n+p})$ depends
only on $\hat b$ in $\bigl\{\hat b=1,\ldots,n_a^{n+p}:
\th^{np}_a(\hat b)=b\bigr\}$, where $X_{a\hat b}^{n+p}$ is the
component of $\pd^pX^n_{ab}$ containing the chosen point.

Now $S_p$ acts on $\pd^pX^n_{ab}$ by permuting $p$ out of the $n+p$
pieces $X_{ad}^1$ meeting at the point. Thus, the $S_p$-action is
compatible with an action of $S_p$ on $n+p$-tuples $(d_1,\ldots,
d_{n+p})$, by permuting $p$ of $d_i$. Since $d_1,\ldots,d_{n+p}$ are
{\it distinct}, this action of $S_p$ on $n+p$-tuples
$(d_1,\ldots,d_{n+p})$ is free. Since $(d_1,\ldots,d_{n+p})$ depends
only on $\hat b$, and the action of $S_p$ on the
$(d_1,\ldots,d_{n+p})$ is free, it follows that the action of $S_p$
on the domain $\bigl\{\hat b=1,\ldots,n_a^{n+p}: \th^{np}_a(\hat
b)=b\bigr\}$ of $f^{np}_{ab}$ is {\it free}.

We can now show that $\Pi^{np(p+1)}_{ab}:D_{ab}^{np}\ra
D_{ab}^{n(p+1)}$ is {\it surjective}. From above, for fixed
$f_{ab}^{n(p+1)}\in D_{ab}^{n(p+1)}$, the preimage of
$f_{ab}^{n(p+1)}$ is the set of $S_p$-invariant $f^{np}_{ab}$ with
$f^{np}_{ab}(\hat b)\in S_{\hat b}$ for all $\hat b$. To construct
such $f^{np}_{ab}$, we choose one point $\hat b$ in each $S_p$-orbit
in the domain, and choose $f^{np}_{ab}(\hat b)\in S_{\hat b}$
arbitrarily at each such $\hat b$. This is possible as $S_{\hat
b}\ne\es$. Then we extend $f^{np}_{ab}$ in an $S_p$-invariant way.
This is possible, uniquely, as $S_p$ acts freely on the domain. So
there exist such $f^{np}_{ab}$, and $\Pi^{np(p+1)}_{ab}$ is
surjective.

This proof also gives us the size of $\bmd{(\Pi^{np(p+1)}_{ab})^{-1}
(f_{ab}^{n(p+1)})}$: it is the product of $\md{S_{\hat b}}$ over
representatives $\hat b$ of each $S_p$-orbit. But $\bmd{S_{\hat
b}}=\bmd{D^{(n+p)0}_{a'b'}}/\bmd{D^{(n+p)1}_{a'b'}}$, which is
independent of $f_{ab}^{n(p+1)}$, so the preimage of each
$f_{ab}^{n(p+1)}$ is the same number of points, which must be
$\bmd{D_{ab}^{np}}/\bmd{D_{ab}^{n(p+1)}}$. If $l'>p$ then either
$l'=p+1$, when $\Pi^{npl'}_{ab}$ is surjective from above, or
$l'>p+1$, when $\Pi^{npl'}_{ab}=\Pi^{n(p+1)l'}_{ab}
\ci\Pi^{np(p+1)}_{ab}$. But $\Pi^{n(p+1)l'}_{ab}$ is surjective by
induction, and $\Pi^{np(p+1)}_{ab}$ is surjective from above. Hence
$\Pi^{npl'}_{ab}$ is surjective for all $l'>p$, which proves the
first part of (d). The second part also follows, since
$\Pi^{n(p+1)l'}_{ab}$ and $\Pi^{np(p+1)}_{ab}$ both have the same
number of points in the preimage of each point.

Next we prove (e) holds. Let $m,l,\ti m,\ti l,a,b,\hat b,\ti a,\ti
b$ be as in (e) with $m=n$ and $l=p$. We must show $P_{ab\hat
b}^{mp\ti m\ti l}:D_{ab}^{np}\ra D_{\ti a\ti b}^{\ti m\ti l}$ is
surjective. Consider the diagram
\e
\begin{gathered}
\xymatrix@C=100pt@R=10pt{ D^{np}_{ab} \ar[d]^{\Pi^{np(p+1)}_{ab}}
\ar[r]_{P_{ab\hat b}^{np\ti m\ti l}} & D_{\ti a\ti b}^{\ti m\ti l}
\ar[d]_{\Pi_{\ti a\ti b}^{\ti m\ti l(\ti l+1)}} \\
D^{n(p+1)}_{ab} \ar[r]^{P_{ab\hat b}^{n(p+1)\ti m(\ti l+1)}} &
D_{\ti a\ti b}^{\ti m(\ti l+1)}.}
\end{gathered}
\label{khCeq44}
\e
It is easy to show from the definitions that \eq{khCeq44} commutes.
Fix $f_{\ti a\ti b}^{\ti m\ti l}$ in $D_{\ti a\ti b}^{\ti m\ti l}$.
By (e) and induction, $P_{ab\hat b}^{n(p+1)\ti m(\ti l+1)}$ is
surjective, so we may choose $f^{n(p+1)}_{ab}$ in $D^{n(p+1)}_{ab}$
with~$P_{ab\hat b}^{n(p+1)\ti m(\ti l+1)}(f^{n(p+1)}_{ab})=\Pi_{\ti
a\ti b}^{\ti m\ti l(\ti l+1)}(f_{\ti a\ti b}^{\ti m\ti l})$.

We now claim that there exists $f^{np}_{ab}$ in $D^{np}_{ab}$ with
$\Pi^{np(p+1)}_{ab}(f^{np}_{ab})=f^{n(p+1)}_{ab}$ and $P_{ab\hat
b}^{np\ti m\ti l}(f^{np}_{ab})=f_{\ti a\ti b}^{\ti m\ti l}$. From
above, the set of $f^{np}_{ab}$ with $\Pi^{np(p+1)}_{ab}
(f^{np}_{ab})=f^{n(p+1)}_{ab}$ is the set of $S_p$-invariant maps
$f^{np}_{ab}$ with $f^{np}_{ab}(\hat b)\in S_{\hat b}$ for all $\hat
b$. By definition, $P_{ab\hat b}^{np\ti m\ti l}(f^{np}_{ab})=f_{\ti
a\ti b}^{\ti m\ti l}$ holds if and only if $f^{np}_{ab}$ takes
prescribed values on $\Im I_{ab\hat b}^{np\ti m\ti l}$. Since
$f_{\ti a\ti b}^{\ti m\ti l}$ is allowable, $\Im I_{ab\hat b}^{np\ti
m\ti l}$ is an $S_p$-invariant subset of the domain of
$f^{np}_{ab}$, and the prescribed values are $S_p$-equivariant. As
$P_{ab\hat b}^{n(p+1)\ti m(\ti l+1)}(f^{n(p+1)}_{ab})=\Pi_{\ti a\ti
b}^{\ti m\ti l(\ti l+1)}(f_{\ti a\ti b}^{\ti m\ti l})$, the
prescribed values on $\Im I_{ab\hat b}^{np\ti m\ti l}$ lie in
$S_{\hat b}$ for each~$\hat b$.

Thus, we take $f^{np}_{ab}$ to have the prescribed values on $\Im
I_{ab\hat b}^{np\ti m\ti l}$, and then for each $S_p$-orbit in the
complement of $\Im I_{ab\hat b}^{np\ti m\ti l}$ in the domain, we
fix a representative $\hat b$, choose $f^{np}_{ab}(\hat b)\in
S_{\hat b}\ne\es$ arbitrarily, and extend $f^{np}_{ab}$ uniquely to
be $S_p$-invariant. This $f^{np}_{ab}$ lies in $D^{np}_{ab}$ with
$P_{ab\hat b}^{np\ti m\ti l}(f^{np}_{ab})=f_{\ti a\ti b}^{\ti m\ti
l}$, so $P_{ab\hat b}^{np\ti m\ti l}$ is {\it surjective}. The
previous argument shows that the preimage of each $f_{\ti a\ti
b}^{\ti m\ti l}$ is the same number of points, proving~(e).

This concludes the inductive step. Thus, we have shown by induction
that (d) holds for $m=n$ and all $l'>l>0$, and (e) holds for $m=n$
and all $l>0$. Since $D_{ab}^{nl}$ is the one point set $\{\es\}$
for $l>M-n$, and $\Pi^{n1l}_{ab}:D_{ab}^{n1}\ra D_{ab}^{nl}$ is
surjective by (d), we see that when $l=1$, the set $D_{ab}^{n1}$ is
nonempty. But functions $f_{ab}^{n1}$ in $D_{ab}^{n1}$ prescribe
allowable boundary conditions for the data $(\ubG_{abo}^{n,i})_{i\in
I_{ab}^n}$ for $o\in O_{ab}^n$ that we must construct. Therefore, in
the proof so far we have shown that the set of allowable, consistent
boundary conditions for $(\ubG_{abo}^{n,i})_{i\in I_{ab}^n}$ is
nonempty.

We shall now choose the data of Hypothesis \ref{khChyp}(i)--(iv)
when $m=n$. For part (i), define $O_{ab}^n=\Aut(X_{ab}^n,\bs
f_{ab}^n,\bs G_{ab}^n)\t D_{ab}^{n1}$ for each $(a,b)\in R^n$. Since
$\Aut(X_{ab}^n,\bs f_{ab}^n,\bs G_{ab}^n)$ is a finite group by
Theorem \ref{kh3thm2} and $D_{ab}^{n1}$ is a nonempty finite set
from above, $O_{ab}^n$ is a nonempty finite set, so Hypothesis
\ref{khChyp}(i) holds.

For (ii), define a left action of $\Aut(X_{ab}^n,\bs f_{ab}^n, \bs
G_{ab}^n)$ on $O_{ab}^n$ to be the product of the left action of
$\Aut(X_{ab}^n,\bs f_{ab}^n,\bs G_{ab}^n)$ on itself with the
natural action of $\Aut(X_{ab}^n,\bs f_{ab}^n,\bs G_{ab}^n)$ on
$D_{ab}^{n1}$, combining the morphism $\Aut(X_{ab}^n,\bs
f_{ab}^n,\bs G_{ab}^n)\ra\Aut(\pd X_{ab}^n,\bs f_{ab}^n\vert_{\pd
X_{ab}^n},\bs G_{ab}^n \vert_{\pd X_{ab}^n})$ with the
interpretation of $D_{ab}^{n1}$ as a set of choices of data on $(\pd
X_{ab}^n,\bs f_{ab}^n\vert_{\pd X_{ab}^n},\bs G_{ab}^n \vert_{\pd
X_{ab}^n})$. This action of $\Aut(X_{ab}^n,\bs f_{ab}^n, \bs
G_{ab}^n)$ on $O_{ab}^n$ is {\it free}, since $\Aut(X_{ab}^n,\bs
f_{ab}^n,\bs G_{ab}^n)$ acts freely on itself, so part (ii) holds.

For (iii), we must define maps $f_{abo}^{nl}$ in \eq{khCeq40} for
all $o\in O^n_{ab}$ and $l\ge 0$. When $l=0$, as in the last part of
(c) we define $f_{abo}^{n0}(b)=o$ for all $o\in O^n_{ab}$. When
$l=1$, if $o\in O_{ab}^n$ then $o=\bigl((\bs u,\bs
v),f_{ab}^{n1}\bigr)$ for $(\bs u,\bs v)\in\Aut(X_{ab}^n,\bs
f_{ab}^n,\bs G_{ab}^n)$ and $f_{ab}^{n1}\in D_{ab}^{n1}$. As above,
$\Aut(X_{ab}^n,\bs f_{ab}^n,\bs G_{ab}^n)$ acts on $D_{ab}^{n1}$.
Define $\ti f_{ab}^{n1}=(\bs u,\bs v)\cdot f_{ab}^{n1}$ using this
action, so that $\ti f_{ab}^{n1}\in D_{ab}^{n1}$. Now define
$f_{abo}^{n1}=\ti f_{ab}^{n1}$. This is valid as by construction,
$D_{ab}^{n1}$ is the set of possibilities for $f_{abo}^{n1}$. When
$l>1$ we have a map $\Pi_{ab}^{n1l}:D^{n1}_{ab}\ra D^{nl}_{ab}$ and
we define $f_{abo}^{nl}=\Pi_{ab}^{n1l}(f_{abo}^{n1})$. The
definitions imply that $f_{abo}^{nl}(\bar b)\in O_{\phi^{n+l}(a,\bar
b)}^{n+l}$ for all $\bar b$, completing
Hypothesis~\ref{khChyp}(iii).

For (iv), write $(\bs 1,\bs 1)$ for the identity in
$\Aut(X_{ab}^n,\bs f_{ab}^n,\bs G_{ab}^n)$. We first choose the
$(\ubG_{abo}^{n,i})_{i\in I_{ab}^n}$ for each $(a,b)\in R^n$ and
$o\in O_{ab}^n$ of the form $\bigl((\bs 1,\bs 1),f_{ab}^{n1}\bigr)$
for $f_{ab}^{n1}\in D_{ab}^{n1}$. Via Hypothesis \ref{khChyp}(c),
this $f_{ab}^{n1}$ determines the restrictions
$\ubG_{abo}^{n,i}\big\vert{}_{\pd E_{ab}^{n,i}}$. That is,
$f_{ab}^{n1}$ determines maps $\ubH_{abo}^{n,i}:\pd
E_{ab}^{n,i}\ra\uP$ for all $i\in I_{ab}^n$, and we seek
$\ubG_{abo}^{n,i}:E_{ab}^{n,i}\ra\uP$ for $i\in I_{ab}^n$ with
$\ubG_{abo}^{n,i}\big\vert{}_{\pd E_{ab}^{n,i}}\equiv
\ubH_{abo}^{n,i}$, satisfying the conditions of Definition
\ref{kh3def17}, so that $\ubG_{abo}^n=\bigl(\bs I_{ab}^n,
\bs\eta_{ab}^n,\ubG_{abo}^{n,i}: i\in I_{ab}^n\bigr)$ is {\it
effective gauge-fixing data} for $(X_{ab}^n,\bs f_{ab}^n)$, as in
Hypothesis~\ref{khChyp}(iv).

Since $f_{ab}^{n1}$ is allowable, the induced $f_{ab}^{n2}$ is
invariant under $S_2\cong\Z_2$. Now $f_{ab}^{n2}$ is the restriction
of $f_{ab}^{n1}$ to $\pd^2X_{ab}^n$, and $S_2$-invariance means
$\ubH_{abo}^{n,i}\big\vert{}_{\pd^2E_{ab}^{n,i}}$ is invariant under
the natural involution $\si:\pd^2E_{ab}^{n,i}\ra\pd^2E_{ab}^{n,i}$,
for all $i\in I_{ab}^n$. But this is exactly the necessary and
sufficient condition for $\ubH_{abo}^{n,i}$ to extend to a map
$\ubG_{abo}^{n,i}:E_{ab}^{n,i}\ra\uP$ for $i\in I_{ab}^n$ with
$\ubG_{abo}^{n,i}\big\vert{}_{\pd E_{ab}^{n,i}}\equiv
\ubH_{abo}^{n,i}$, by Principle~\ref{kh2pri}(a).

Thus, $\ubH_{abo}^{n,i}$ prescribes consistent values for
$\ubG_{abo}^{n,i}$ on $E_{ab}^{n,i}\sm(E_{ab}^{n,i})^\ci$, and it
remains to specify $\ubG_{abo}^{n,i}$ on $(E_{ab}^{n,i})^\ci$.
Choose some $k_o\ge 0$, and maps
$\ubG_{abo}^{n,i}\vert_{(E_{ab}^{n,i})^\ci}:
(E_{ab}^{n,i})^\ci\ra\R^{k_o}\subset\uP$ for $i\in I_{ab}^n$ such
that $\coprod_{i\in I_{ab}^n}\ubG_{abo}^{n,i}
\vert_{(E_{ab}^{n,i})^\ci}:\coprod_{i\in
I_{ab}^n}(E_{ab}^{n,i})^\ci\ra\uP$ is {\it injective}. This is
possible if $k_o\gg 0$, and we can take $\ubG_{abo}^{n,i}
\vert_{(E_{ab}^{n,i})^\ci}$ smooth. In addition, we require that the
choices of integers $k_o$ for all $o=\bigl((\bs 1,\bs
1),f_{ab}^{n1}\bigr)$ for $f_{ab}^{n1}\in D_{ab}^{n1}$ should be
{\it distinct}.

We now claim that $\ubG_{abo}^n$ satisfies Definition
\ref{kh3def17}. Part (a) holds as the $V_{ab}^{n,i}$ are manifolds,
so $(V_{ab}^{n,i},\ldots,\psi_{ab}^{n,i})$ is automatically an
effective Kuranishi neighbourhood for all $i\in I_{ab}^n$. Part (b)
when $l=0$ holds as $\coprod_{i\in I_{ab}^n}\ubG_{abo}^{n,i}
\vert_{(E_{ab}^{n,i})^\ci}$ is injective, and when $l>0$ holds by
induction and Hypothesis \ref{khChyp}(iv) for $m=n+l$, since by
construction $\ubG_{abo}^n\vert_{\pd^lX_{ab}^n}$ is a disjoint union
of effective gauge-fixing data $\ubG_{a'b'o'}^{n+l}$, so Definition
\ref{kh3def17}(b) with $l=0$ for these $\ubG_{a'b'o'}^{n+l}$ implies
Definition \ref{kh3def17}(b) for $\ubG_{abo}^n$ with $l>0$. Hence
$\ubG_{abo}^n$ is effective gauge-fixing data.

This defines the data $(\ubG_{abo}^{n,i})_{i\in I_{ab}^n}$ for each
$(a,b)\in R^n$ and $o\in O_{ab}^n$ of the form $\bigl((\bs 1,\bs
1),f_{ab}^{n1}\bigr)$. For general $o\in O_{ab}^n$, write
$o=\bigl((\bs u,\bs v),f_{ab}^{n1}\bigr)$ and $\dot o=(\bs u,\bs
v)^{-1}\cdot o=\bigl((\bs 1,\bs 1),\dot f_{ab}^{n1}\bigr)$, where
$\dot f_{ab}^{n1}=(\bs u,\bs v)^{-1}\cdot f_{ab}^{n1}$, so that we
have already defined $\ubG_{ab\dot o}^n$, and set $\ubG_{abo}^n=(\bs
u,\bs v)_*(\ubG_{ab\dot o}^n)$. That is, we define
$\ubG_{abo}^{n,i}\equiv(\hat v_{ab}^{n,i})_*(\ubG_{ab\dot
o}^{n,i})\equiv\ubG_{ab\dot o}^{n,i}\ci(\hat v_{ab}^{n,i})^{-1}$ for
all $i\in I_{ab}^n$, where $\hat v_{ab}^{n,i}:E_{ab}^{n,i}\ra
E_{ab}^{n,i}$ is the diffeomorphism in $\bs v$. Then $\ubG_{abo}^n$
is effective gauge-fixing data for all $o\in O_{ab}^n$, giving
Hypothesis~\ref{khChyp}(iv).

We must now prove Hypothesis \ref{khChyp}(a)--(e) hold when $m=n$.
For (a), suppose $o=\bigl((\bs u,\bs v),f_{ab}^{n1}\bigr)$ and $\ti
o=\bigl((\bs{\ti u},\bs{\ti v}),\ti f_{ab}^{n1}\bigr)$ are distinct
elements of $O_{ab}^n$. We must prove that $\ubG_{abo}^n\ne
\ubG_{ab\ti o}^n$. Write $\dot f_{ab}^{n1}=(\bs u,\bs v)^{-1}\cdot
f_{ab}^{n1}$, $\dot o=\bigl((\bs 1,\bs 1),\dot f_{ab}^{n1}\bigr)$
and $\ddot f_{ab}^{n1}=(\bs{\ti u},\bs{\ti v})^{-1}\cdot \ti
f_{ab}^{n1}$, $\ddot o=\bigl((\bs 1,\bs 1),\ddot f_{ab}^{n1}\bigr)$.
Then $\ubG_{abo}^n=(\bs u,\bs v)_*(\ubG_{ab\dot o}^n)$ and
$\ubG_{ab\ti o}^n=(\bs{\ti u},\bs{\ti v})_*(\ubG_{ab\ddot o}^n)$.
Since the $\ubG_{ab\dot o}^{n,i}$ map $(E_{ab}^{n,i})^\ci\ra
\R^{k_{\dot o}}$, it follows that the $\ubG_{abo}^{n,i}$ map
$(E_{ab}^{n,i})^\ci\ra\R^{k_{\dot o}}$, and similarly the
$\ubG_{ab\ti o}^{n,i}$ map~$(E_{ab}^{n,i})^\ci \ra\R^{k_{\ddot o}}$.

If $\dot o\ne\ddot o$ then $k_{\dot o}\ne k_{\ddot o}$, as from
above the $k_o$ for all $o=\bigl((\bs 1,\bs 1),f_{ab}^{n1}\bigr)$
are distinct, so $\ubG_{abo}^{n,i}\vert_{(E_{ab}^{n,i})^\ci}$ and
$\ubG_{ab\ti o}^{n,i}\vert_{(E_{ab}^{n,i})^\ci}$ map to different
spaces, giving $\ubG_{abo}^{n,i}\ne\ubG_{ab\ti o}^{n,i}$ and
$\ubG_{abo}^n\ne \ubG_{ab\ti o}^n$ as we want. So suppose $\dot
o=\ddot o$. Then $\ubG_{abo}^n=(\bs u,\bs v)_*(\ubG_{ab\dot o}^n)$
and $\ubG_{ab\ti o}^n=(\bs{\ti u},\bs{\ti v})_*(\ubG_{ab\dot o}^n)$,
so that $\ubG_{ab\ti o}^n=(\bs{\bar u},\bs{\bar
v})_*(\ubG_{abo}^n)$, where $(\bs{\bar u},\bs{\bar v})=(\bs{\ti
u},\bs{\ti v})\ci(\bs u,\bs v)^{-1}$ in $\Aut(X_{ab}^n,\bs
f_{ab}^n,\bs G_{ab}^n)$. As $o\ne\ti o$ we have $(\bs{\bar
u},\bs{\bar v})\ne(\bs 1,\bs 1)$. Since $\ubG_{abo}^n$ is effective
gauge-fixing data and $(X_{ab}^n,\bs f_{ab}^n,\ubG_{abo}^n)$ is
connected we have $\Aut(X_{ab}^n,\bs f_{ab}^n,\ubG_{abo}^n)=\{1\}$
by Theorem \ref{kh3thm5}(b). Hence $(\bs{\bar u},\bs{\bar
v})\notin\Aut(X_{ab}^n,\bs f_{ab}^n, \ubG_{abo}^n)$, so
$\ubG_{abo}^n\ne(\bs{\bar u},\bs{\bar v})_*(\ubG_{abo}^n)$,
and~$\ubG_{abo}^n\ne\ubG_{ab\ti o}^n$.

This proves (a). Part (b) is immediate from $\ubG_{abo}^n=(\bs u,\bs
v)_*(\ubG_{ab\dot o}^n)$ when $o=\bigl((\bs u,\bs v),
f_{ab}^{n1}\bigr)$ and $\dot o=(\bs u,\bs v)^{-1}\cdot o$. Part (c)
when $l=0$ is trivial, as we defined $f_{abo}^{n0}(b)=o$. For
$o=\bigl((\bs 1,\bs 1),f_{ab}^{n1}\bigr)$ we have
$f_{abo}^{n1}=f_{ab}^{n1}$, and part (c) for $m=n$, $l=1$ and this
$o$ holds by construction, since we used (c) for $l=1$ to choose the
values of $\ubG_{abo}^{n,i}\big\vert{}_{\pd E_{ab}^{n,i}}$ for $i\in
I_{ab}^n$. Part (c) for $m=n$ and $l>1$ follows from (c) for $m=n$,
$l=1$, since the corresponding $f^{nl}_{ab}$ are induced by
$f^{n1}_{ab}$. This proves (c) for all $o$ of the form $o=\bigl((\bs
1,\bs 1),f_{ab}^{n1}\bigr)$. Part (c) for general $o\in O_{ab}^n$
follows by equivariance of the construction under $\Aut(X_{ab}^n,\bs
f_{ab}^n,\ubG_{abo}^n)$.

It remains to prove (d),(e). We have already done this for $m=n$ and
$l>0$ above, so only the case $l=0$ remains. In fact it is enough to
prove (d) in the case $m=n$, $l=0$ and $l'=1$, since we can then
deduce (d) with $l=0$ and $l'>1$ and (e) when $l=0$ as in the
inductive proof above. Thus, we have to prove that
$\Pi^{n01}_{ab}:D_{ab}^{n0}\ra D_{ab}^{n1}$ is surjective, and the
preimage of each point is the same number of points. But as above
each $f_{ab}^{n0}$ in $D_{ab}^{n0}$ maps $\{b\}\ra O_{ab}^n$, and
$f_{ab}^{n0}\mapsto f_{ab}^{n0}(b)$ induces a bijection
$D_{ab}^{n0}\ra O_{ab}^n$. Identifying $D_{ab}^{n0}$ with
$O_{ab}^n=\Aut(X_{ab}^n,\bs f_{ab}^n,\bs G_{ab}^n)\t D_{ab}^{n1}$,
we see that $\Pi^{n01}_{ab}$ maps $\Aut(X_{ab}^n,\bs f_{ab}^n,\bs
G_{ab}^n)\t D_{ab}^{n1}\ra D_{ab}^{n1}$, and is given explicitly by
$\Pi^{n01}_{ab}:\bigl((\bs u,\bs v),f_{ab}^{n1}\bigr)\mapsto
f_{ab}^{n1}$. Hence $\Pi^{n01}_{ab}$ is surjective, and the preimage
of each point is the same number of points $\bmd{\Aut(X_{ab}^n,\bs
f_{ab}^n,\bs G_{ab}^n)}$. This completes (d),(e), and the inductive
step in the proof of Hypothesis \ref{khChyp}. Therefore by
induction, Hypothesis \ref{khChyp} holds for all~$N\ge 0$.

We can now at last complete Step 4. Consider the data of Hypothesis
\ref{khChyp} when $m=0$. Since each $(X_a,\bs f_a,\bs G_a)$ is
connected by Step 3, we have $n_a^0=1$, with $X_{a1}^0=X_a$. Setting
$m=0$ and $b=1$, define $O_a=O_{a1}^0$ and $\ubG_{ao}=\ubG_{a1o}^0$
for all $a\in A$ and $o\in O_a$. Using the notation introduced in
Step 4 we have
\e
\begin{split}
\pd[X_a,\bs f_a,\bs G_a]&=\sum_{b=1}^{n_a^1}[X_{ab}^1,\bs
f_{ab}^1,\bs G_{ab}^1]=\sum_{b=1}^{n_a^1}\ep_{ab}^1
[X_{\phi^1(a,b)}^1,\bs f_{\phi^1(a,b)}^1,\bs G_{\phi^1(a,b)}^1]\\
&=\sum_{(\bar a,\bar b)\in R^1}
\raisebox{-10pt}{\begin{Large}$\displaystyle\biggl\{$\end{Large}}
\sum_{\begin{subarray}{l} b=1,\ldots,n_a^1:\\
\phi^1((a,b))=(\bar a,\bar b)\end{subarray}}\!\!\!\!\!\ep_{ab}^1
\raisebox{-10pt}{\begin{Large}$\displaystyle\biggr\}$\end{Large}}
[X_{\bar a\bar b}^1,\bs f_{\bar a\bar b}^1,\bs G_{\bar a\bar b}^1]
\end{split}
\label{khCeq45}
\e
in $\widetilde{KC}_{k-1}(Y;R)$, where in the second step we use the
isomorphism $(\bs a,\bs b)_{ab}^1:(X_{ab}^1,\bs f_{ab}^1,\bs
G_{ab}^1)\ra(X_{\phi^1(a,b)}^1,\bs f_{\phi^1(a,b)}^1,\bs
G_{\phi^1(a,b)}^1)$, which multiplies orientation by $\ep_{ab}^1$.
We need an analogue of \eq{khCeq45} with $\bs G_a$ replaced by
$\ubG_{ao}$ for $o\in O_a$. Now Hypothesis \ref{khChyp}(c) with
$m=0$, $b=1$ and $l=1$ implies that $(\bs a,\bs b)_{ab}^1$
identifies $(X_{ab}^1,\bs f_{ab}^1,\bs G_{ao}\vert_{X_{ab}^1})$ with
$(X_{\phi^1(a,b)}^1,\bs f_{\phi^1(a,b)}^1,\bs G_{\phi^1(a,b)\bar
o}^1)$, where $\bar o\in O_{\phi^1(a,b)}^1$ is given by $\bar
o=f_{a1o}^{01}(b)$. Thus following \eq{khCeq45}, in
$KC^\ef_{k-1}(Y;R)$ we have
\e
\pd[X_a,\bs f_a,\ubG_{ao}]=\sum_{\begin{subarray}{l}(\bar a,\bar
b)\in R^1,\\ \bar o\in O_{\bar a\bar b}^1\end{subarray}}
\raisebox{-10pt}{\begin{Large}$\displaystyle\biggl\{$\end{Large}}
\sum_{\begin{subarray}{l} b=1,\ldots,n_a^1:\\
\phi^1((a,b))=(\bar a,\bar b),\;\bar o=f_{a1o}^{01}(b)
\end{subarray}}\!\!\!\!\!\!\!\!\!\!\!\!\ep_{ab}^1\,\,\,
\raisebox{-10pt}{\begin{Large}$\displaystyle\biggr\}$\end{Large}}
[X_{\bar a\bar b}^1,\bs f_{\bar a\bar b}^1,\ubG_{\bar a\bar b\bar
o}^1].
\label{khCeq46}
\e

Consider the composition $O_a=O_{a1}^0\cong D_{a1}^{00}
\,{\buildrel\Pi_{a1}^{001}\over\longra}\, D_{a1}^{01} \,{\buildrel
P_{a1b}^{0110}\over\longra}\,D_{ab}^{10}\cong O_{ab}^1$ for $a\in A$
and $b=1,\ldots,n_a^1$. This maps $o$ to $\bar o=f_{a1o}^{01}(b)$.
Now Hypothesis \ref{khChyp}(d),(e) say that $\Pi_{a1}^{001}$ and
$P_{a1b}^{0110}$ are surjective and the preimages of each point is
the same size. Hence for fixed $a\in A$, $b=1,\ldots,n_a^1$ and
$(\bar a,\bar b)=\phi^1(a,b)$, the map $O_a\ra O_{\bar a\bar b}^1$
taking $o\mapsto \bar o=f_{a1o}^{01}(b)$ is surjective, and the
preimage of each point is $\md{O_a}/\md{O_{\bar a\bar b}^1}$ points.
Therefore, if we average equation \eq{khCeq46} over all $o\in O_a$,
on the right hand side each $O_{\bar a\bar b}^1$ occurs equally
often, yielding
\e
\begin{split}
&\ts\pd\bigl(\sum_{o\in O_a}\md{O_a}^{-1}[X_a,\bs
f_a,\ubG_{ao}]\bigr)=\\
&\sum_{(\bar a,\bar b)\in R^1}
\raisebox{-10pt}{\begin{Large}$\displaystyle\biggl\{$\end{Large}}
\sum_{\begin{subarray}{l} b=1,\ldots,n_a^1:\\
\phi^1((a,b))=(\bar a,\bar b)\end{subarray}}\!\!\!\!\!\ep_{ab}^1
\raisebox{-10pt}{\begin{Large}$\displaystyle\biggr\}$\end{Large}}
\sum_{\bar o\in O_{\bar a\bar b}^1}\md{O_{\bar a\bar b}^1}^{-1}[X_{\bar
a\bar b}^1,\bs f_{\bar a\bar b}^1,\ubG_{\bar a\bar b\bar o}^1].
\end{split}
\label{khCeq47}
\e
Multiplying \eq{khCeq47} by $\rho_a$ and summing over $a\in A$ gives
\e
\begin{split}
&\ts\pd\bigl(\sum_{a\in A}\sum_{o\in O_a}\rho_a\md{O_a}^{-1}[X_a,\bs
f_a,\ubG_{ao}]\bigr)=\\
&\sum_{(\bar a,\bar b)\in R^1}
\raisebox{-10pt}{\begin{Large}$\displaystyle\biggl\{$\end{Large}}
\sum_{\begin{subarray}{l} a\in A,\; b=1,\ldots,n_a^1:\\
\phi^1((a,b))=(\bar a,\bar
b)\end{subarray}}\!\!\!\!\!\rho_a\ep_{ab}^1
\raisebox{-10pt}{\begin{Large}$\displaystyle\biggr\}$\end{Large}}
\sum_{\bar o\in O_{\bar a\bar b}^1}
\md{O_{\bar a\bar b}^1}^{-1}[X_{\bar a\bar b}^1,\bs f_{\bar a\bar
b}^1,\ubG_{\bar a\bar b\bar o}^1].
\end{split}
\label{khCeq48}
\e

Now in Step 4 of \S\ref{khC} we showed that for each $(\bar a,\bar
b)\in R^1$, either (A) $(X_{\bar a\bar b}^1,\bs f_{\bar a\bar
b}^1,\bs G_{\bar a\bar b}^1)$ admits an orientation-reversing
automorphism, or (B) the coefficient $\{\cdots\}$ on the second line
of \eq{khCeq48} is zero in $R$. We claim that in case (A), the final
sum $\sum_{\bar o\in O_{\bar a\bar b}^1} \md{O_{\bar a\bar
b}^1}^{-1}[X_{\bar a\bar b}^1,\bs f_{\bar a\bar b}^1,\ubG_{\bar
a\bar b\bar o}^1]$ in \eq{khCeq48} is zero. Thus in both cases (A)
and (B), the product on the second line of \eq{khCeq48} is zero, so
both sides of \eq{khCeq48} are zero, proving equation~\eq{khCeq19}.

To prove the claim, suppose $(X_{\bar a\bar b}^1,\bs f_{\bar a\bar
b}^1,\bs G_{\bar a\bar b}^1)$ admits an orientation-reversing
automorphism. There is a natural group morphism $\ga:\Aut(X_{\bar
a\bar b}^1,\bs f_{\bar a\bar b}^1,\bs G_{\bar a\bar b}^1)\ra\{\pm
1\}$ with $\ga(\bs u,\bs v)=1$ if $\bs u$ is orientation-preserving,
and $\ga(\bs u,\bs v)=-1$ if $\bs u$ is orientation-reversing. As
there exists an orientation-reversing automorphism of $(X_{\bar
a\bar b}^1,\bs f_{\bar a\bar b}^1,\bs G_{\bar a\bar b}^1)$, this
$\ga$ is surjective.

By Hypothesis \ref{khChyp}(ii), $\Aut(X_{\bar a\bar b}^1,\bs f_{\bar
a\bar b}^1,\bs G_{\bar a\bar b}^1)$ acts freely on $O_{\bar a\bar
b}^1$. Choose a subset $\ti O_{\bar a\bar b}^1$ in $O_{\bar a\bar
b}^1$ containing one point from each orbit of $\Aut(X_{\bar a\bar
b}^1,\bs f_{\bar a\bar b}^1,\bs G_{\bar a\bar b}^1)$. Then there is
a bijection $\Aut(X_{\bar a\bar b}^1,\bs f_{\bar a\bar b}^1,\bs
G_{\bar a\bar b}^1)\t\ti O_{\bar a\bar b}^1\ra O_{\bar a\bar b}^1$
mapping $\bigl((\bs u,\bs v),o\bigr)\mapsto(\bs u,\bs v)\cdot o$.
Thus we have
\e
\begin{split}
\sum_{\bar o\in O_{\bar a\bar b}^1}[X_{\bar a\bar b}^1, \bs f_{\bar
a\bar b}^1,\ubG_{\bar a\bar b\bar o}^1] &=\sum_{\ti o\in\ti O_{\bar
a\bar b}^1}\sum_{(\bs u,\bs v)\in \Aut(X_{\bar a\bar b}^1,\bs
f_{\bar a\bar b}^1,\bs G_{\bar a\bar b}^1)\!\!\!\!\!\!\!\!\!\!}
[X_{\bar a\bar b}^1,\bs f_{\bar a\bar b}^1,\ubG_{\bar a\bar b
(\bs u,\bs v)\cdot\ti o}^1]\\
&=\sum_{\ti o\in\ti O_{\bar a\bar b}^1}\sum_{(\bs u,\bs v)\in
\Aut(X_{\bar a\bar b}^1,\bs f_{\bar a\bar b}^1,\bs G_{\bar a\bar
b}^1)\!\!\!\!\!\!\!\!\!\!\!\!\!\!\!\!\!\!\!\!}\ga(\bs u,\bs
v)[X_{\bar a\bar b}^1,\bs f_{\bar a\bar b}^1,\ubG_{\bar a\bar b\ti
o}^1]=0,
\end{split}
\label{khCeq49}
\e
where in the second step we use the isomorphism $(\bs u,\bs
v):(X_{\bar a\bar b}^1,\bs f_{\bar a\bar b}^1,\ubG_{\bar a\bar b\ti
o}^1)\ra (X_{\bar a\bar b}^1,\bs f_{\bar a\bar b}^1,\ubG_{\bar a\bar
b (\bs u,\bs v)\cdot\ti o}^1)$, which follows from Hypothesis
\ref{khChyp}(b) and multiplies orientations by $\ga(\bs u,\bs v)$,
and in the third step we note that $\sum_{(\bs u,\bs v)\in
\Aut(X_{\bar a\bar b}^1,\bs f_{\bar a\bar b}^1,\bs G_{\bar a\bar
b}^1)}\ab\ga(\bs u,\bs v)=0$ as $\ga$ is surjective, so exactly half
of $(\bs u,\bs v)$ have $\ga(\bs u,\bs v)=1$, and half have $\ga(\bs
u,\bs v)=-1$. Equation \eq{khCeq49} proves the claim, and
\eq{khCeq19} follows.

Thus $\sum_{a\in A}\sum_{o\in O_a}\rho_a \md{O_a}^{-1}[X_a,\bs
f_a,\ubG_{ao}]$ is a cycle in $KC_k^\ef(Y;R)$, so
$\be=\bigl[\sum_{a\in A}\sum_{o\in O_a}\rho_a \md{O_a}^{-1}[X_a,\bs
f_a,\ubG_{ao}]\bigr]$ is well-defined in $KH_k^\ef(Y;R)$. Finally,
we must prove that $\Pi_\ef^\Kh(\be)=\al$. By definition of
$\Pi_\ef^\Kh$, $\sum_{a\in A}\sum_{o\in
O_a}\rho_a\md{O_a}^{-1}\ab[X_a,\bs f_a,\Pi(\ubG_{ao})]$ is a cycle
in $KC_k(Y;R)$ representing $\Pi_\ef^\Kh(\be)$. Hence it is
sufficient to construct a homology between $\sum_{a\in A}\sum_{o\in
O_a}\rho_a\ab\md{O_a}^{-1}[X_a,\bs f_a,\Pi(\ubG_{ao})]$
and~$\sum_{a\in A}\rho_a[X_a,\bs f_a,\bs G_a]$.

For each $a\in A$, consider $[0,1]\t X_a$ as a compact oriented
Kuranishi space with the product orientation, and define $\bs
g_a=\bs f_a\ci\bs\pi_{X_a}:[0,1]\t X_a\ra Y$, a strongly smooth map.
Set $J_a=\{i+1:i\in I_a\}$. For $i\in I_a$, define a Kuranishi
neighbourhood on $[0,1]\t X_a$ by $(W^{i+1}_a,F^{i+1}_a,
t^{i+1}_a,\chi^{i+1}_a)=\bigl([0,1]\t V^i_a,[0,1]\t
E^i_a,s^i_a\ci\pi_{V^i_a},\id_{[0,1]}\t\psi^i\bigr)$, define
$g^{i+1}_a:W^{i+1}_a\ra Y$ by $g^{i+1}_a=f^i_a\ci\pi_{V^i_a}$, and
$\ze_{i+1,a}:[0,1]\t X_a\ra[0,1]$ by $\ze_{i+1,a}=\eta_{i,a}\ci
\pi_{X_a}$. For $i,j\in I_a$, define $W^{(i+1)(j+1)}_a=[0,1]\t
V^{ij}_a$, $\psi^{(i+1)(j+1)}_a=\id_{[0,1]}\t\phi^{ij}_a$,
$\hat\psi{}^{(i+1)(j+1)}_a=\id_{[0,1]}\t\hat\phi^{ij}_a$ and
$\ze_{i+1,a}^{j+1}:W^{j+1}_a\ra[0,1]$ by
$\ze_{i+1,a}^{j+1}=\eta_{i,a}^j\ci\pi_{V^j_a}$. Write $\bs
J_a=\bigl(J_a,(W^j_a,\ldots,\chi^j_a),g^j_a:j\in J,\ldots\bigr)$,
and $\bs\ze_a=(\ze_{j,a}:j\in J$, $\ze_{j,a}^k:j,k\in J)$. It is
easy to check that $(\bs J_a,\bs\ze_a)$ is an {\it excellent
coordinate system} for~$\bigl([0,1]\t X_a,\bs g_a\bigr)$.

For $o\in O_a$ and $i\in I_a$, define $H_{ao}^{i+1}:F^{i+1}_a\ra P$
by, for $u\in[0,1]$ and $e\in E^i_a$,
\e
H_{ao}^{i+1}(u,e)=\begin{cases} G^i_a(e), & u=0, \\
\mu\bigl((u),G^i_a(e)\bigr), & u\in (0,1), \\
\Pi\ci\uG^i_{ao}(e), & u=1.
\end{cases}
\label{khCeq50}
\e
Write $\bs H_{ao}=(\bs J_a,\bs\ze_a,H_{ao}^j:j\in J_a)$. As in the
proofs of Theorem \ref{kh4thm3} and Proposition \ref{kh4prop}, one
can show $\bs H_{ao}$ is {\it gauge-fixing data} for~$\bigl([0,1]\t
X_a,\bs g_a\bigr)$.

We have $\pd\bigl([0,1]\t X_a\bigr)=\bigl(\{1\}\t X_a\bigr)\amalg
-\bigl(\{0\}\t X_a\bigr)\amalg-\bigl([0,1]\t\pd X_a\bigr)$ in
oriented Kuranishi spaces. The construction of $(\bs J_a,\bs\ze_a)$
and the cases $u=0$, $u=1$ in \eq{khCeq50} ensure that $\bs
H_{ao}\vert_{\{0\}\t X_a}=\bs G_a$ and $\bs H_{ao}\vert_{\{1\}\t
X_a}=\Pi(\ubG_{ao})$, identifying $\{0\}\t X_a\cong X_a\cong\{1\}\t
X_a$. Thus in $KC_k(Y;R)$ we have
\e
\begin{split}
\pd\bigl[[0,1]\t X_a,\bs g_a,\bs H_{ao}\bigr]=\,&[X_a,\bs
f_a,\Pi(\ubG_{ao})]-[X_a,\bs f_a,\bs G_a]\\
&-\bigl[[0,1]\t\pd X_a,\bs g_a\vert_{[0,1]\t\pd X_a},\bs
H_{ao}\vert_{[0,1]\t\pd X_a}\bigr].
\end{split}
\label{khCeq51}
\e

Following the proof of \eq{khCeq19} in equations
\eq{khCeq46}--\eq{khCeq49} above, but replacing $X_a,\bs
f_a,\ubG_{ao}$ by $[0,1]\t X_a,\bs g_a,\bs H_{ao}$, and $X_{\bar
a\bar b}^1$ by $[0,1]\t X_{\bar a\bar b}^1$, and so on, yields
\e
\ts\sum_{a\in A}\sum_{o\in O_a}\rho_a\md{O_a}^{-1} \bigl[[0,1]\t\pd
X_a,\bs g_a\vert_{[0,1]\t\pd X_a},\bs H_{ao}\vert_{[0,1]\t\pd
X_a}\bigr]=0.
\label{khCeq52}
\e
Thus, multiplying \eq{khCeq51} by $\rho_a\md{O_a}^{-1}$, summing
over all $a\in A$ and $o\in O_a$, and using \eq{khCeq52} to cancel
out the contribution of the final term in \eq{khCeq51} gives
\begin{align*}
&\ts\pd\bigl(\sum_{a\in A}\sum_{o\in O_a}\rho_a\ab\md{O_a}^{-1}
\bigl[[0,1]\t X_a,\bs g_a,\bs H_{ao}\bigr]\bigr)=\\
&\ts\sum_{a\in A}\sum_{o\in O_a}\rho_a\ab\md{O_a}^{-1}[X_a,\bs
f_a,\Pi(\ubG_{ao})]-\sum_{a\in A}\rho_a[X_a,\bs f_a,\bs G_a].
\end{align*}
Therefore $\sum_{a\in A}\sum_{o\in O_a}\rho_a\md{O_a}^{-1}[X_a,\bs
f_a,\Pi(\ubG_{ao})]$ is homologous to $\sum_{a\in
A}\rho_a\ab[X_a,\ab\bs f_a,\bs G_a]$, and $\Pi_\ef^\Kh(\be)=\al$.
Hence for each $\al\in KH_k(Y;R)$ we can construct $\be\in
KH_k^\ef(Y;R)$ with $\Pi_\ef^\Kh(\be)=\al$, and
$\Pi_\ef^\Kh:KH_*^\ef(Y;R)\ra KH_*(Y;R)$ is {\it surjective}. This
completes Step~4.

\subsection{Step 5: $\Pi_\ef^\Kh:KH_*^\ef(Y;R)\ra KH_*(Y;R)$ is
injective}
\label{khC5}

As in Step 5, suppose $\be\in KH_k^\ef(Y;R)$ with
$\Pi_\ef^\Kh(\be)=0$, represent $\be$ by the image of a cycle
$\sum_{s\in S}\ze_s\,\tau_s$ in singular homology, and choose
$\sum_{a\in A}\rho_a[X_a,\bs f_a,\bs G_a]$ in $KC_{k+1}(Y;R)$
satisfying \eq{khCeq20}. We now explain how to apply Steps 1--4 to
this $\sum_{a\in A}\rho_a[X_a,\bs f_a,\bs G_a]$, replacing
$k$-cycles by $(k\!+\!1)$-chains, to construct $\sum_{a\in
A}\sum_{o\in O_a}\rho_a\md{O_a}^{-1}[X_a,\bs f_a,\ubG_{ao}]$ in
$KC_{k+1}^\ef(Y;R)$ satisfying \eq{khCeq22}. This will imply that
$\be=0$, so $\Pi_\ef^\Kh:KH_k^\ef(Y;R)\ra KH_k(Y;R)$ is {\it
injective}.

First we discuss Step 1. Introducing connected $(\ti X_d,\bs{\ti
f}_d,\bs{\ti G}_d)$ and $\Ga_d,\eta_d$ for $d\in D$ as in Step 1,
generalizing \eq{khCeq2}, the lift of \eq{khCeq20} to
$\widetilde{KC}_k(Y;R)$ is
\ea
&\ts\sum_{a\in A}\rho_a[\pd X_a,\bs f_a\vert_{\pd X_a},\bs
G_a\vert_{\pd X_a}]=\sum_{s\in S}\ze_s\bigl[\De_k,\tau_s,\bs
G_{\De_k}\bigr]+
\label{khCeq53}\\
&\ts\sum_{d\in D}\eta_d\bigl([\ti X_d/\Ga_d,\bs\pi_*(\bs{\ti
f}_d),\bs\pi_*(\bs{\ti G}_d)]\!-\!\frac{1}{\md{\Ga_d}}[\ti
X_d,\bs{\ti f}_d,\bs{\ti G}_d]\bigr). \nonumber
\ea
Write $\pd(X_a,\bs f_a,\bs G_a)=\coprod_{b=1}^{n_a^1}(X_{ab}^1,\bs
f_{ab}^1,\bs G_{ab}^1)$ as usual. Then \eq{khCeq53} implies that
some of the components $(X_{ab}^1,\bs f_{ab}^1,\bs G_{ab}^1)$ or
$(\ti X_d,\bs{\ti f}_d,\bs{\ti G}_d)$ must be isomorphic to
$(\De_k,\tau_s,\bs G_{\De_k})$ for~$s\in S$.

In Step 1 we choose tent functions $\bs T_a$ for $(X_a,\bs f_a,\bs
G_a)$ for $a\in A$ and $\bs{\ti T}_d$ for $(\ti X_d,\bs{\ti
f}_d,\bs{\ti G}_d)$ for $d\in D$ to cut $X_a,\ti X_d$ into
`arbitrarily small pieces' $X_{ac},\ti X_{df}$ for $c\in C_a$, $f\in
F_d$. The new issue we must deal with here is that when
$(X_{ab}^1,\bs f_{ab}^1,\bs G_{ab}^1)$ or $(\ti X_d,\bs{\ti
f}_d,\bs{\ti G}_d)$ is isomorphic to $(\De_k,\tau_s,\bs G_{\De_k})$,
the tent function $\bs T_a$ or $\bs{\ti T}_d$ also cuts
$(\De_k,\tau_s,\bs G_{\De_k})$ into pieces, so $\sum_{s\in
S}\ze_s\bigl[\De_k,\tau_s,\bs G_{\De_k}\bigr]$ is modified. We must
keep track of these modifications, and in particular ensure that the
modified version of $\sum_{s\in S}\ze_s\bigl[\De_k,\tau_s,\bs
G_{\De_k}\bigr]$ is still the image of a cycle in $KC^\ef_k(Y;R)$
representing~$\be$.

We do this using an idea from Step 4 of \S\ref{khB} in \S\ref{khB4}.
We fix some $N\gg 0$, and choose the tent functions $\bs T_a,\bs{\ti
T}_d$ such that if $(X_{ab}^1,\bs f_{ab}^1,\bs G_{ab}^1)$ or $(\ti
X_d,\bs{\ti f}_d,\bs{\ti G}_d)$ is isomorphic to $(\De_k,\tau_s,\bs
G_{\De_k})$, then $\bs T_a\vert_{X_{ab}^1}$ or $\bs{\ti T}_d$ cuts
$\De_k$ into the $N^{\it th}$ {\it barycentric subdivision} of
$\De_k$. Since we construct the $\bs T_a,\bs{\ti T}_d$ in
\S\ref{khB1} by first choosing tent functions $\bs T_{st}^m$ for
$(X_{st}^m,\bs f_{st}^m,\bs G_{st}^m)$ by induction on decreasing
$m=M,M-1,\ldots,0$, we do this by inserting an extra condition: if
for $0<m\le M$ we have $(X_{st}^m,\bs f_{st}^m,\bs
G_{st}^m)\cong(\De_{k+1-m},\tau,\bs G_{\De_{k+1-m}})$ for some
smooth $\tau:\De_{k+1-m}\ra Y$, then we choose $\bs T_{st}^m$ to be
a standard tent function on $\De_{k+1-m}$ which cuts $\De_{k+1-m}$
into its $N^{\rm th}$ barycentric subdivision.

We must check this extra condition is compatible with conditions
(i)--(iv) near the beginning of \S\ref{khC1}. Part (i) is immediate
since our condition is expressed using isomorphism classes. For
(ii), as $\De_{k+1-m}$ is a simply-connected manifold we cannot
write $(\De_{k+1-m},\tau,\bs G_{\De_{k+1-m}})$ as a quotient
$\bigl(X_{de}^m/\Stab_{\Ga_d}(e),\ab\bs\pi^m_*(\bs f_{de}^m)$ for
$(X_{de}^m,\bs f_{de}^m,\bs G_{de}^m)$ connected and
$\Stab_{\Ga_d}(e)\ne\{1\}$, and if $(X_{de}^m,\bs f_{de}^m,\bs
G_{de}^m)\cong(\De_{k+1-m},\tau,\bs G_{\De_{k+1-m}})$ then
$\Stab_{\Ga_d}(e)$ must act trivially on $(X_{de}^m,\bs f_{de}^m,\bs
G_{de}^m)$ as $\Aut(\De_{k+1-m},\tau,\bs G_{\De_{k+1-m}})=\{1\}$.
Thus $X_{de}^m\cong\De_{k+1-m}\t(\{0\}/\Stab_{\Ga_d}(e))$, so we
deal with (ii) by extending the extra condition to say that if
$(X_{st}^m,\bs f_{st}^m,\bs G_{st}^m)\cong(\De_{k+1-m}\t(\{0\}/\Ga),
\tau,\bs G_{\De_{k+1-m}})$ for $\Ga$ a finite group, then we again
choose $\bs T_{st}^m$ to be a standard tent function on
$\De_{k+1-m}\t(\{0\}/\Ga)$ which cuts $\De_{k+1-m}$ into its $N^{\rm
th}$ barycentric subdivision.

Part (iii) holds automatically since $\pd X_{st}^m$ is the disjoint
union of $k+2-m$ copies of $\De_{k-m}$, and on each of these we have
chosen $\bs T_{st'}^{m+1}$ to be a standard tent function inducing
the $N^{\rm th}$ barycentric subdivision of $\De_{k-m}$ in the
previous inductive step. Part (iv) is satisfied provided $N\gg 0$ is
large enough, since then the $N^{\rm th}$ barycentric subdivision
cuts $X_{st}^m\cong\De_{k+1-m}$ into arbitrarily small pieces. Thus
we can carry out Step 1 in such a way that the $\bs T_a,\bs{\ti
T}_d$ cut each $(\De_k,\tau_s,\bs G_{\De_k})$ into the $N^{\rm th}$
barycentric subdivision of $\De_k$ for some fixed~$N\gg 0$.

We also make one additional change to the definition of the
gauge-fixing data $\bs G_{ac}$ for $(X_{ac},\bs f_{ac})$ for $c\in
C_a$. Here $(X_{ac},\bs f_{ac},\bs G_{ac})$ was defined in
Proposition \ref{khAprop9} using the triple $(Z_{X_a,\bs T_a},\bs
f_a\ci\bs\pi,\bs H_{X_a,\bs T_a})$ constructed from $(X_a,\bs
f_a,\bs G_a)$ and $\bs T_a$ in Definition \ref{khAdef12}. The issue
is that when we use a standard tent function on $\De_k$ to cut
$(\De_k,\tau_s,\bs G_{\De_k})$ into its $N^{\rm th}$ barycentric
subdivision, we obtain $(k!)^N$ small pieces
$(\De_{kl},\tau_s\vert_{\De_{kl}},\bs G_{\De_{kl}})$ for $l$ in an
indexing set $L_{k,N}$, where each $\De_{kl}$ is a small $k$-simplex
with a natural affine diffeomorphism $\De_{kl}\cong\De_k$. However,
because of the way the function $H^{i+1}$ is defined in Definition
\ref{khAdef12}, the diffeomorphism $\De_{kl}\cong\De_k$ does not
identify $\bs G_{\De_{kl}}$ with $\bs G_{\De_k}$, as the functions
$G^k_{\De_{kl}}:\De_{kl}\ra P$ and $G^k_{\De_k}:\De_k\ra P$ in $\bs
G_{\De_{kl}},\bs G_{\De_k}$ are not identified.

We deal this by modifying the definition of the functions
$H^{k+2}_a:W^{k+2}_a\ra P$ in $\bs H_{X_a,\bs T_a}$, by hand, which
then modifies the functions $G^{k+1}_{ac}:V^{k+1}_{ac}\ra P$ in $\bs
G_{ac}$, to identify each $\bs G_{\De_{kl}}$ with $\bs G_{\De_k}$,
and similarly for $\ti X_d,\bs{\ti T}_d$. That is, whenever
$(X_{ab}^m,\bs f_{ab}^m,\bs G_{ab}^m)$ is a component of
$\pd^m(X_a,\bs f_a,\bs G_a)$ isomorphic to some
$(\De_{k+1-m},\tau,\bs G_{\De_{k+1-m}})$, so that $\bs T_{ab}^m$ is
a standard tent function on $X_{ab}^m\cong\De_{k+1-m}$ inducing the
$N^{\rm th}$ subdivision, then on the image in $Z_{X_a,\bs T_a}$
(and hence each $X_{ac}$) of each piece $X_{abl}^m\cong
\De_{(k+1-m)l}$ for $l\in L_{k+1-m,N}$, we define the function
$H^{k+2}_a$ to be identified with $G^{k+1-m}_{\De_{k+1-m}}:
\De_{k+1-m}\ra P$ under the identification~$X_{abl}^m\cong
\De_{(k+1-m)l}\cong\De_{k+1-m}$.

This prescription is compatible on overlaps between the images of
different simplices $X_{abl}^m, X_{ab'l'}^{m'}$, as two such
simplices must intersect in a face of each, and the functions
$G^k_{\De_k}$ restrict to $G^l_{\De_l}$ on any face $\De\cong\De_l$
of $\De_k$. With these changes, following the proof in Step 1 that
$\sum_{a\in A}\sum_{c\in C_a}\rho_a[X_{ac},\bs f_{ac},\bs G_{ac}]$
is homologous to $\sum_{a\in A}\rho_a[X_a,\bs f_a,\bs G_a]$ proves
equation \eq{khCeq21}, as we want. This completes our discussion of
the changes to Step~1.

Step 2 requires essentially no changes. For Step 3, the analogue of
\eq{khCeq14} is
\e
\begin{split}
&\ts\pd\bigl(\sum_{a\in A}\sum_{c\in C_a}\rho_a[X_{ac},\bs
f_{ac},\bs
G_{ac}]\bigr)-\Pi_\rsi^\Kh\ci\Up^N(\sum_{s\in S}\ze_s\,\tau_s)=\\
&\sum_{d\in D}\eta_d\Bigl(\begin{aligned}[t]&\ts\sum_{f\Ga_d\in
F_d/\Ga_d}\bigl[\ti X_{df}/\Stab_{\Ga_d}(f),\bs\pi_*(\bs{\ti
f}_{df}),\bs\pi_*(\bs{\ti
G}_{df})\bigr]\\
&-\md{\Ga_d}^{-1}\ts\sum_{f\in F_d}[\ti X_{df},\bs{\ti
f}_{df},\bs{\ti G}_{df}]\Bigr).
\end{aligned}
\end{split}
\label{khCeq54}
\e
which is a lift of \eq{khCeq21} to $\widetilde{KC}_k(Y;R)$. As for
\eq{khCeq14}, all terms on the l.h.s.\ of \eq{khCeq54} have {\it
trivial stabilizers}, since the $X_a$ which we started with at the
beginning of Step 3 (after changing notation) do, so the $X_{ac}$
do, and the Kuranishi spaces in $\Pi_\rsi^\Kh\ci\Up^N (\sum_{s\in
S}\ze_s\,\tau_s)$ are $k$-simplices $\De_k$, and so are manifolds.

In Step 3, we proved that both sides of \eq{khCeq14} are zero in
$\widetilde{KC}_{k-1}(Y;R)$ provided the $X_{ac},\ti X_{df}$ are
taken sufficiently small, because all terms on the l.h.s.\ of
\eq{khCeq14} have trivial stabilizers. Since this also holds for the
l.h.s.\ of \eq{khCeq54}, the same argument shows that both sides of
\eq{khCeq54} are zero in $\widetilde{KC}_k(Y;R)$ for $X_{ac},\ti
X_{df}$ sufficiently small. Therefore \eq{khCeq21} holds in
$\widetilde{KC}_k(Y;R)$, as well as $KC_k(Y;R)$. This completes the
changes to Step~3.

At the beginning of Step 4 we have a chain $\sum_{a\in
A}\rho_a[X_a,\bs f_a,\bs G_a]$ satisfying
\e
\ts\pd\bigl(\sum_{a\in A}\rho_a[X_a,\bs f_a,\bs G_a]\bigr)=
\Pi_\rsi^\Kh\ci\Up^N(\sum_{s\in S}\ze_s\,\tau_s)
\label{khCeq55}
\e
in $\widetilde{KC}_k(Y;R)$ for $N\gg 0$, where $\sum_{s\in
S}\ze_s\,\tau_s$ is the same as it has been all the way through, but
$\sum_{a\in A}\rho_a[X_a,\bs f_a,\bs G_a]$ has been changed several
times, and now the $X_a$ have trivial stabilizers, the $V^i_a$ are
manifolds, the $(X_a,\bs f_a,\bs G_a)$ are connected, and each
component of $\pd^mX_a$ occurs as the intersection of $m$ distinct
components of $\pd X_a$. Our goal is to choose finite sets $O_a$ and
maps $\ubG_a^i:V_a^i\ra\uP$ for $a\in A$ and $i\in I_a$ such that
\eq{khCeq55} lifts to equation \eq{khCeq22} in~$KH_k^\ef(Y;R)$.

The issue is this: we have to choose these maps $\ubG_a^i$ such that
whenever some component $(X_{ab}^1,\bs f_{ab}^1,\bs G_{ab}^1)$ of
$\pd(X_a,\bs f_a,\bs G_a)$ is isomorphic to $\Pi_\rsi^\Kh(\up)$ for
some $\up:\De_k\ra Y$ occurring in the $N^{\rm th}$ barycentric
subdivision of some $\tau_s:\De_k\ra Y$, then the corresponding
component of $\pd(X_a,\bs f_a,\ubG_{ao})$ must be isomorphic to
$\Pi_\rsi^\ef(\up)$ for all $o\in O_a$. That is, $E_{ab}^{1,k}$ in
$\bs G_{ab}^1$ is a component of $\pd E_a^{k+1}$ in $\bs G_a$, and
the maps $\ubG_{ao}^{k+1}\vert{E_{ab}^{1,k}}:E_{ab}^{1,k}\ra\uP$ and
$\ubG_{\De_k}^k:\De_k\ra\uP$ should agree under the given
diffeomorphism $E_{ab}^{1,k}\cong\De_k$. If we do this, then all
terms on the l.h.s.\ of \eq{khCeq55} isomorphic to terms on the
r.h.s.\ of \eq{khCeq55} lift to terms on the l.h.s.\ of \eq{khCeq22}
isomorphic to terms on the r.h.s.\ of \eq{khCeq22}, and so
\eq{khCeq22} holds.

We ensure this by adding an extra condition in the inductive proof
of Hypothesis \ref{khChyp} in \S\ref{khC4}, when for part (iv) we
chose maps $\ubG_{abo}^{m,i}:E_{ab}^{m,i}\ra\uP$ for $i\in I_{ab}^m$
such that $\ubG_{abo}^m=\bigl(\bs I_{ab}^m,\bs\eta_{ab}^m,
\ubG_{abo}^{m,i}:i\in I_{ab}^m\bigr)$ is effective gauge-fixing data
for $(X_{ab}^m,\bs f_{ab}^m)$. We require that if any $(X_{ab}^m,\bs
f_{ab}^m,\bs G_{ab}^m)$ is isomorphic to $(\De_{k+1-m},\ab\tau,\ab
\bs G_{\De_{k+1-m}})$ for some smooth $\tau:\De_{k+1-m}\ra Y$, then
the corresponding maps $\ubG_{abo}^{m,k+1-m}:E_{ab}^{m,k+1-m}\ra\uP$
for $o\in O_{ab}^m$ are identified with $\ubG_{\De_{k+1-m}}^{k+1-m}:
\De_{k+1-m}\ra\uP$ under the given diffeomorphism $E_{ab}^{m,k+1-m}
\cong\De_{k+1-m}$. Then $(X_{ab}^m,\bs f_{ab}^m,\ubG_{abo}^m)$ is
isomorphic to~$(\De_{k+1-m},\tau,\ubG_{\De_{k+1-m}})$.

We must check that prescribing these particular choices are
compatible with conditions (a)--(e), so that the inductive proof
that we can choose data satisfying Hypothesis \ref{khChyp} is still
valid with this extra condition. For part (a), we claim that if
$(X_{ab}^m,\bs f_{ab}^m,\bs G_{ab}^m)$ is isomorphic to some
$(\De_{k+1-m},\tau,\bs G_{\De_{k+1-m}})$ then $O_{ab}^m$ is a single
point. We prove this by reverse induction on $m=k+1,k,\ldots,0$. For
the inductive step, note that in \S\ref{khC4} we defined $O_{ab}^m=
\Aut(X_{ab}^m,\bs f_{ab}^m,\bs G_{ab}^m)\t D_{ab}^{m1}$. As
$X_{ab}^m\cong\De_{k+1-m}$, it has $k+2-m$ boundary faces
$X_{ab'}^{m+1}$, each of which is a copy of $\De_{k-m}$, and
$D_{ab}^{m1}$ consists of a choice of point in
$O_{\phi^{m+1}(a,b')}^{m+1}$ satisfying some consistency conditions.

By induction, each such $O_{\phi^{m+1}(a,b')}^{m+1}$ is a point, as
$(X_{ab'}^{m+1},\bs f_{ab'}^{m+1},\bs G_{ab'}^{m+1})
\cong(\De_{k-m},\ab\up,\ab\bs G_{\De_{k-m}})$, where $\up=\tau\ci
F_j^{k+1-m}$ for some $j=0,\ldots,k+1-m$, in the notation of
\S\ref{kh41}. Therefore $D_{ab}^{m1}$ is one point, as there is only
one possible choice in $O_{\phi^{m+1}(a,b')}^{m+1}$ for each $b'$.
Also $\Aut(X_{ab}^m,\bs f_{ab}^m,\bs G_{ab}^m)=\{1\}$ as
$\Aut(\De_{k+1-m},\tau,\bs G_{\De_{k+1-m}})=\{1\}$. Hence $O_{ab}^m=
\Aut(X_{ab}^m,\bs f_{ab}^m,\bs G_{ab}^m)\t D_{ab}^{m1}$ is one
point, proving the inductive step.

Since $O_{ab}^m$ is a single point, (a) is trivial as there do not
exist distinct $o,o'$ in $O_{ab}^m$. Part (b) is also trivial, as
$\Aut(X_{ab}^m,\bs f_{ab}^m,\bs G_{ab}^m)=\{1\}$. Part (c) holds
since as above each boundary component $(X_{ab'}^{m+1},\bs
f_{ab'}^{m+1},\bs G_{ab'}^{m+1})$ of $(X_{ab}^m,\bs f_{ab}^m,\bs
G_{ab}^m)$ is isomorphic to $(\De_{k-m},\up,\bs G_{\De_{k-m}})$ for
some $\up$, so in the previous inductive step we chose
$\ubG_{ab'o'}^{m+1}\cong\ubG_{\De_{k-m}}$ for $o'\in O_{ab'}^{m+1}$,
and in \S\ref{kh43} we defined $\ubG_{\De_{k+1-m}}$ so that it
restricted to $\ubG_{\De_{k-m}}$ on each $\De_{k-m}$ face of
$\De_{k+1-m}$. And (d),(e) are not affected as they do not concern
choice of the~$\ubG_{abo}^{m,i}$.

Therefore the inductive proof goes through, and we can choose data
satisfying Hypothesis \ref{khChyp} with this extra condition when
$(X_{ab}^m,\bs f_{ab}^m,\bs G_{ab}^m)\cong(\De_{k+1-m},\ab\tau,\ab
\bs G_{\De_{k+1-m}})$. The proof of \eq{khCeq19} using equations
\eq{khCeq45}--\eq{khCeq49} then generalizes to give \eq{khCeq22},
and Step 5 follows. This completes the proof of
Theorem~\ref{kh4thm2}.

\section{Identification of products $\cup,\cap,\bu$}
\label{khD}

We now prove Theorem \ref{kh4thm6}, that is, we shall show that
$\Pi_\cs^\Kch,\Pi_\cs^\ec,\Pi_\rsi^\Kh, \Pi_\rsi^\ef$ identify the
cup, cap and intersection products $\cup,\cap,\bu$ on
$H^*_\cs,H^\rsi_*(Y;R)$ with those on $KH^*,KH_*(Y;R)$ and
$KH_\ec^*,KH^\ef_*(Y;R)$. Proving this is tricky since we do not
have a good description of $\Pi_\cs^\Kch,\Pi_\cs^\ec$ at the cochain
level. Our approach is to first suppose that $Y$ is oriented, and
prove that $\Pi_\rsi^\Kh:H_*^\rsi(Y;R)\ra KH_*(Y;R)$,
$\Pi_\rsi^\ef:H_*^\rsi(Y;R)\ra KH_*^\ef(Y;R)$ identify the
intersection products $\bu$ on $H_*^\rsi(Y;R)$ and
$KH_*,KH_*^\ef(Y;R)$, using chain level descriptions of $\bu$ on
$C_*^\rsi(Y;R)$ and $KC_*,KC_*^\ef(Y;R)$. Then since Poincar\'e
duality identifies $\bu$ with $\cup,\cap$, we deduce that
$\Pi_\cs^\Kch,\Pi_\cs^\ec,\Pi_\rsi^\Kh,\Pi_\rsi^\ef$ also identify
cup and cap products.

Chain-level definitions of $\bu$ on $KC_*,KC_*^\ef(Y;R)$ are given
in \eq{kh4eq53} and \eq{kh4eq56}. Describing $\bu$ on singular
chains $C_*^\rsi(Y;R)$ is more complex. Two good references are
Bredon \cite[VI.11]{Bred}, and Lefschetz \cite[\S IV]{Lefs}, which
is the original source for much of the theory of the intersection
product on singular homology of manifolds. As these do not discuss
orbifolds, we assume $Y$ is a {\it manifold}.

Let $Y$ be an oriented $n$-manifold without boundary and $\al\in
H_k^\rsi(Y;R)$, $\be\in H_l^\rsi(Y;R)$, so that $\al\bu\be\in
H_{k+l-n}^\rsi(Y;R)$. For Theorem \ref{kh4thm6}(c) we must prove
that $\Pi_\rsi^*(\al\bu\be)=\Pi_\rsi^*(\al)\bu\Pi_\rsi^*(\be)$ in
$KH_{k+l-n}(Y;R)$ or $KH_{k+l-n}^\ef(Y;R)$. First consider the case
in which $\al=[K]$ and $\be=[L]$, where $K,L$ are compact, oriented,
immersed submanifolds of $Y$ without boundary intersecting
transversely in $Y$, with immersions $\io_K:K\ra Y$ and~$\io_L:L\ra
Y$.

Choose effective gauge-fixing data $\ubG_K,\ubG_L$ for
$(K,\io_K),(L,\io_L)$, which is possible by Theorem
\ref{kh3thm5}(a), and set $\bs G_K=\Pi(\ubG_K)$, $\bs
G_L=\Pi(\ubG_L)$. Then $[K,\io_K,\bs G_K]\in KC_k(Y;R)$ with
$\pd[K,\io_K,\bs G_K]=0$ as $\pd K=\es$, so $\bigl[[K,\io_K,\bs
G_K]\bigr]\!\in\!KH_k(Y;R)$. Similarly $\bigl[[K,\io_K,\ubG_K]
\bigr]\!\in\!KH_k^\ef(Y;R)$, $\bigl[[L,\io_L,\bs G_L]\bigr]\!\in\!
KH_l(Y;R)$, and $\bigl[[L,\io_L,\ubG_L]\bigr]\!\in\! KH_l^\ef(Y;R)$.
Considering Step 3 in Appendix \ref{khB}, when we construct an
inverse $(\Pi_\rsi^\ef)^{-1}$ for $\Pi_\rsi^\ef:H_k^\rsi(Y;R)\ra
KH_k^\ef(Y;R)$, we see that
$(\Pi_\rsi^\ef)^{-1}\bigl(\bigl[[K,\io_K,\ubG_K] \bigr]\bigr)$ is
represented, as a singular chain, by a triangulation of $K$ by
$k$-simplices. But this is how we define the class $[K]$ of $K$ in
$H_k^\rsi(Y;R)$. Since $\al=[K]$, this gives
$\bigl[[K,\io_K,\ubG_K]\bigr]=\Pi_\rsi^\ef(\al)$. As
$\bigl[[K,\io_K,\bs G_K]\bigr]=\ab\Pi_\ef^\Kh
\bigl(\bigl[[K,\ab\io_K,\ab\ubG_K]\bigr]\bigr)$ and
$\Pi_\rsi^\Kh\!=\!\Pi_\ef^\Kh\ci\Pi_\rsi^\ef$, we have
$\bigl[[K,\io_K,\bs G_K]\bigr]=\Pi_\rsi^\Kh(\al)$. Similarly
$\bigl[[L,\io_L,\ubG_L]\bigr]=\Pi_\rsi^\ef(\be)$
and~$\bigl[[L,\io_L,\bs G_L]\bigr]=\Pi_\rsi^\Kh(\be)$.

Therefore $\Pi_\rsi^\Kh(\al),\Pi_\rsi^\Kh(\be)$ are represented in
$KC_*(Y;R)$ by $[K,\io_K,\ab\bs G_K],\ab[L,\io_L,\bs G_L]$, so by
\eq{kh4eq53}, $\Pi_\rsi^\Kh(\al)\bu\Pi_\rsi^\Kh(\be)$ is represented
in $KC_{k+l-n}(Y;R)$ by
\e
\begin{split}
\Pi_\Kch^\Kh&\bigl(\Pi_\Kh^\Kch([K,\io_K,\bs G_K])\cup\Pi_\Kh^\Kch
[L,\io_L,\bs G_L]\bigr)\\
&=\Pi_\Kch^\Kh\bigl(\bigl[K^Y,\io_K^Y,\bs C_{\bs G_K}^Y\bigr]\cup
\bigl[L^Y,\io_L^Y,\bs C_{\bs G_L}^Y\bigr]\bigr)\\
&=\Pi_\Kch^\Kh\bigl(\bigl[K^Y\t_YL^Y,\bs\pi_Y,\bs C_{\bs G_K}^Y\t_Y
\bs C_{\bs G_L}^Y\bigr]\bigr)\\
&=\bigl[K^Y\t_YL^Y,\bs\pi_Y,\bs G_{\bs C_{\bs G_K}^Y\t_Y \bs C_{\bs
G_L}^Y}\bigr],
\end{split}
\label{khDeq1}
\e
using the definitions of $\Pi_\Kh^\Kch,\Pi_\Kch^\Kh,\cup$ in
Definitions \ref{kh4def7} and~\ref{kh4def15}.

Since $K,L$ intersect transversely in $Y$, their intersection $K\cap
L$ is a compact immersed $(k\!+\!l\!-\!n)$-submanifold of $Y$, which
can also be considered as a fibre product $K\t_{\io_K,Y,\io_L}L$.
Following Bredon \cite[p.~375]{Bred} we define an {\it
orientation\/} on $K\cap L$, using the orientations on $K,L,Y$. Then
\cite[Th.~VI.11.9]{Bred} shows that $[K]\bu[L]=[K\cap L]$ in
$H_{k+l-n}^\rsi(Y;R)$. Thus, by the argument above,
$\Pi_\rsi^\Kh(\al\bu\be)$ is represented in $KC_{k+l-n}(Y;R)$ by
$[K\cap L,\io_{K\cap L},\bs G_{K\cap L}]$, for some choice of
gauge-fixing data $\bs G_{K\cap L}$ for~$(K\cap L,\io_{K\cap L})$.

Hence, to show that $\Pi_\rsi^\Kh(\al\bu\be)=\Pi_\rsi^\Kh(\al)
\bu\Pi_\rsi^\Kh(\be)$, by \eq{khDeq1} it is enough to prove that
$[K\cap L,\io_{K\cap L},\bs G_{K\cap L}]$ is homologous to
$\bigl[K^Y\t_YL^Y,\bs\pi_Y,\bs G_{\bs C_{\bs G_K}^Y\t_Y \bs C_{\bs
G_L}^Y}\bigr]$ in $KC_{k+l-n}(Y;R)$. Now the topological space
underlying $K^Y\t_YL^Y$ is just $K\cap L$, so $K^Y\t_YL^Y$ is just
$K\cap L$ with a different Kuranishi structure. Examining the
definitions shows that $K^Y\t_YL^Y$ is isomorphic as a Kuranishi
space to $((K\cap L)^Y)^Y$, in the notation of Definition
\ref{kh4def7}, and that $\bs\pi_Y$ agrees with $(\io_{K\cap
L}^Y)^Y$. Comparing our orientation conventions in \S\ref{kh27} with
Bredon's shows that the orientation conventions of $K^Y\t_YL^Y$ also
agree. Now
\e
\begin{split}
\Pi_\Kch^\Kh\ci&\Pi_\Kh^\Kch\ci\Pi_\Kch^\Kh\ci\Pi_\Kh^\Kch\bigl(
[K\cap L,\io_{K\cap L},\bs G_{K\cap L}]\bigr)\\
&=\bigl[((K\cap L)^Y)^Y, (\io_{K\cap L}^Y)^Y,\bs G_{\bs C^Y_{\bs
G_{\bs C^Y_{\bs G_{K\cap L}}}}}\bigr].
\end{split}
\label{khDeq2}
\e

The proof of Theorem \ref{kh4thm3} shows that \eq{khDeq2} is
homologous to $[K\cap L,\io_{K\cap L},\ab\bs G_{K\cap L}]$. But the
right hand side is $\bigl[K^Y\t_YL^Y,\bs\pi_Y,\bs G_{\bs C_{\bs
G_K}^Y\t_Y \bs C_{\bs G_L}^Y}\bigr]$ with an alternative choice of
gauge-fixing data. As $\pd(K^Y\t_YL^Y)=\es$, the homology class does
not depend on the gauge-fixing data. Therefore $[K\cap L,\io_{K\cap
L},\bs G_{K\cap L}]$ is homologous to $\bigl[K^Y\t_YL^Y,\bs\pi_Y,\bs
G_{\bs C_{\bs G_K}^Y\t_Y\bs C_{\bs G_L}^Y}\bigr]$, so
$\Pi_\rsi^\Kh(\al\bu\be)=\Pi_\rsi^\Kh(\al)\bu\Pi_\rsi^\Kh(\be)$.
Similarly $\Pi_\rsi^\ef(\al\bu\be)=\Pi_\rsi^\ef(\al)\bu
\Pi_\rsi^\ef(\be)$.

This proves Theorem \ref{kh4thm6}(c) in the case when $\al,\be$ can
be represented by compact oriented submanifolds $K,L$ intersecting
transversely. For the general case we can use ideas in Lefschetz
\cite{Lefs}, particularly \cite[\S IV]{Lefs}, which is the original
source for much of the theory of the intersection product on
singular homology of manifolds. Lefschetz shows that we can choose a
triangulation $\mathcal T$ of $Y$ by convex polyhedra, with a dual
triangulation ${\mathcal T}^*$, such that if $P,P^*$ are polyhedra
in ${\mathcal T},{\mathcal T}^*$ of dimensions $k,l$ then $P,P^*$
intersect transversely in a convex polyhedron of dimension $k+l-n$,
or~$P\cap P^*=\es$.

We may represent $\al$ by an $R$-linear combination of
$k$-dimensional polyhedra $P$ in $\mathcal T$; to make a singular
chain, we replace each $P$ by its barycentric subdivision into
$k$-simplices. Similarly, we represent $\be$ by an $R$-linear
combination of $l$-dimensional polyhedra $P^*$ in $\mathcal T^*$.
Then $\al\bu\be$ is represented by the corresponding $R$-bilinear
sum of transverse intersections $P\cap P^*$. We then use the
argument above with $K,L$ replaced by $R$-linear combinations of
polygons $P,P^*$, ensuring compatibility at boundaries $\pd P,\pd
P^*$ using the methods of Appendix~\ref{khB}.

This proves Theorem \ref{kh4thm6}(c). Parts (a),(b) when $Y$ is
oriented follow using Poincar\'e duality and Theorems \ref{kh4thm3}
and \ref{kh4thm5}, since Poincar\'e duality identifies
$\cup,\cap,\bu$ on both $H^\rsi_*(Y;R),H^*_\cs(Y;R)$ and
$KH_*(Y;R),KH^*(Y;R)$ and $KH_*^\ef(Y;R),KH^*_\ec(Y;R)$. When $Y$ is
not orientable we deduce (a),(b) using {\it orientation bundles}, as
in Definition \ref{kh4def8}. Part (a) then implies (d), since the
identities in $H^0_\cs(Y;R),KH^0(Y;R),KH^0_\ec(Y;R)$ are
characterized uniquely by the cup products $\cup$. This completes
the proof of Theorem~\ref{kh4thm6}.

\medskip

\noindent{\small\sc The Mathematical Institute, 24-29 St. Giles,
Oxford, OX1 3LB, U.K.}

\noindent{\small\sc E-mail: \tt joyce@maths.ox.ac.uk}

\end{document}